\documentclass[a4paper,10pt]{report}
\usepackage{amsmath}
\usepackage{amssymb}
\usepackage{amsthm}
\usepackage{mathrsfs}
\usepackage{mathtools}
\usepackage{bm}
\usepackage{fullpage}
\usepackage{enumerate}
\usepackage{hyperref}
\usepackage{slashed}
\usepackage{marvosym}
\usepackage{caption}
\usepackage{graphicx}

\newcommand{\upd}{\mathrm{d}}
\newcommand{\D}{\mathscr{D}}
\newcommand{\Lbar}{\underline{L}}
\newcommand{\dVol}{\mathrm{d}\textit{vol}}
\DeclareMathOperator{\tr}{tr}
\DeclareMathOperator{\Div}{div}
\DeclareMathOperator{\Curl}{curl}
\newcommand{\chibar}{\underline{\chi}}
\allowdisplaybreaks[1]

\usepackage[utf8]{inputenc}
\usepackage[backend=biber, style=alphabetic, sorting=anyt]{biblatex}
\renewbibmacro*{in:}{%
	\ifentrytype{article}{}{\printtext{\bibstring{in}\intitlepunct}}}
\newtheorem{theorem}{Theorem}[section]
\newtheorem{lemma}[theorem]{Lemma}
\newtheorem{proposition}[theorem]{Proposition}
\newtheorem{corollary}[theorem]{Corollary}

\newtheorem{conjecture}[theorem]{Conjecture}
\newtheorem{property}[theorem]{Property}

\theoremstyle{definition}
\newtheorem{definition}[theorem]{Definition}

\theoremstyle{remark}
\newtheorem{remark}[theorem]{Remark}

\begin{document}
\pagestyle{empty}

\title{The weak null condition and global existence using the \texorpdfstring{p}{p}-weighted energy method}
\author{Joe Keir \\ \\
{\small DAMTP, Centre for Mathematical Sciences, University of Cambridge}, \\
{\small \sl Wilberforce Road, Cambridge CB3 0WA, UK} \\ \\
{\small Trinity College, Cambridge} \\ \\
\small{j.keir@damtp.cam.ac.uk}}
\maketitle

\setcounter{page}{2}
\pagestyle{plain}

\begin{abstract}
	We prove global existence for solutions arising from small initial data for a large class of quasilinear wave equations satisfying the ``weak null condition'' of Lindblad and Rodnianski, significantly enlarging upon the class of equations for which global existence is known. In addition to the usual weak null condition, we require a certain hierarchical structure in the semilinear terms. Included in this class are the Einstein equations in harmonic coordinates, meaning that a special case of our results is a new proof of the stability of Minkowski space. Our proof also applies to the coupled Einstein-Maxwell system in harmonic coordinates and Lorenz gauge, as well as to various model scalar wave equations which do not satisfy the null condition. Our proof also applies to the Einstein(-Maxwell) equations in wave coordinates (and Lorenz gauge) if, after writing the equations as a set of nonlinear wave equations, we then ``forget'' about the gauge conditions, choosing initial data for the reduced equations which does not satisfy the gauge condition. The methods we use allow us to treat initial data which only has a small ``degenerate energy'', involving a weight that degenerates at null infinity, so the usual (unweighted) energy might be unbounded. We also demonstrate a connection between the weak null condition and \emph{geometric shock formation}, showing that equations satisfying the weak null condition can exhibit ``shock formation at infinity'', of which we provide an explicit example. The methods that we use are very robust and adaptable, including a generalisation of the $p$-weighted energy method of Dafermos and Rodnianski \cite{Dafermos2010b}, adapted to the dynamic geometry using constructions similar to those pioneered by Christodoulou and Klainerman \cite{Christodoulou1993}. This means that our proof applies in a wide range of situations, including those in which the metric remains close to, but never approaches the flat metric in some spatially bounded domain, and those in which the ``geometric'' null infinity and the ``background'' null infinity differ dramatically, for example, when the solution exhibits shock formation at null infinity.
\end{abstract}

\tableofcontents

\chapter{Introduction and overview}

\section{Introduction}

After local existence and uniqueness, one of the most fundamental and foundational issues regarding a PDE is that of \emph{global existence} of its solutions. When the PDE in question is nonlinear, it is often extremely difficult to address this issue in full, but a first step towards a complete understanding can sometimes be made by restricting to suitably \emph{small} initial data. This problem is also important in itself, as it can be interpreted as a kind of ``stability'' of the trivial solution. It addresses the following question: do all suitably small perturbations of the trivial initial data lead to solutions which, at least, exist for all time? In answering this question, we usually have to obtain a good understanding of the solutions and their asymptotics, meaning that stronger kinds of stability (for example, ``asymptotic stability'') are often also proved.

In this paper we examine this problem as it pertains to (systems of) \emph{nonlinear wave equations} in $3+1$ dimensions. These systems are ubiquitous in physics, playing a role in many fields, from fluid dynamics to general relativity. As such, understanding the behaviour of their solutions is of the utmost importance. The (local and global) behaviour of \emph{linear} waves is, of course, well understood, and the solution to the local-in-time problem is also classical, at least in the high-regularity setting (see, for example, \cite{Sogge2008}). In fact, standard methods show that all nonlinear wave equations of the sort we will study in this paper have local-in-time solutions, given sufficiently smooth data. 

In contrast, the \emph{global-in-time} behaviour of solutions to nonlinear wave equations is far more subtle and intricate. For example, the scalar wave equation\footnote{Here, as elsewhere in this work, $\Box$ is the standard flat space wave operator in four dimensions, i.e.\ $\Box = -\partial_t^2 + \Delta$.}
\begin{equation}
\label{equation classical null condition example}
\Box \phi = (\partial_t \phi)^2 - \sum_{i = 1}^3 (\partial_i \phi)^2
\end{equation}
has global solutions for all small initial data, whereas the equation
\begin{equation}
\label{equation John's example}
\Box \phi = (\partial_t \phi)^2
\end{equation}
possesses solutions arising from \emph{arbitrarily small} initial data, but which blow up\footnote{In fact, the situation is even worse: \emph{all} non-trivial, smooth, compactly supported initial data leads to solutions that blow up in a finite time!} in a finite time (see \cite{John1981}). So, while it is true that certain nonlinear wave equations possess global solutions for small initial data, other equations exhibit finite-time blowup, and these equations can look superficially very similar. An important question then becomes the following: what are the conditions on the nonlinearity which guarantee global existence for small initial data?

A small set of particular nonlinear equations can be transformed by some ``trick'' into simpler (e.g. linear) equations (see the example attributed to Nirenberg in \cite{Klainerman1980}) allowing for their solutions to be analysed. However, the first major insight which allowed a whole \emph{class} of nonlinearities to be treated was made by Klainerman (\cite{Klainerman1980}). Here, Klainerman identified a condition on the nonlinearities which was later\footnote{The original paper of Klainerman \cite{Klainerman1980} was limited to spaces with spacial dimension at least six, whereas we are concerned here with equations in $\mathbb{R}\times\mathbb{R}^3$.} proven (independently, in both \cite{Klainerman1986} and \cite{Christodoulou1986}) to \emph{guarantee} global existence for all sufficiently small data. This condition is usually known as the \emph{null condition}, but in this work, in order to avoid confusion with our subject, we will refer to it as the \emph{classical null condition}.

The identification of the classical null condition was important in a number of respects. First, it made it possible to examine a system of equations and to determine that it possesses global solutions for small initial data, simply by examining the structure of the equations, i.e.\ without constructing the solutions or carrying out any detailed estimates. This is of particular importance when the equations in question arise from physics, since, as we have already noted, it frequently has an interpretation as a \emph{global stability} result. Second, it identifies the relevant difference between equations \eqref{equation classical null condition example} and \eqref{equation John's example}: the former satisfies the classical null condition while the latter does not. Third, the null condition can be interpreted physically, and this provides insight into the failure of global existence in equations like equation \eqref{equation John's example}. Specifically, the null condition rules out nonlinear interactions between wave packets that are travelling along the same outgoing null geodesic. This makes sense in view of the fact that wave packets that are localised around the same outgoing null geodesic can interact for an arbitrarily long time, in contrast to those which are travelling in different directions.

Despite its many successes, the classical null condition does not quite give a complete answer to the question ``\emph{what are the structural conditions on a nonlinear wave equation which guarantee global existence for small data?}'' To be precise, although the classical null condition is \emph{sufficient} to guarantee global existence for small data, it is not \emph{necessary}. Indeed, several systems\footnote{Chief among these systems are the \emph{Einstein equations in wave coordinates}, which we will discuss in depth later.} have been discovered which, although they do \emph{not} satisfy the classical null condition, nevertheless possess global solutions for all sufficiently small initial data. All of these systems obey a weaker condition, called the ``\emph{weak null condition}'' (first described in \cite{Lindblad2003}). We postpone the definition of this condition to section \ref{section intro weak null condition}, as it is slightly technical and requires the introduction of some notation.

Thus far, no general proof has been presented proving that \emph{all} nonlinear wave equations satisfying the weak null condition in fact possess global solutions for sufficiently small data. Such a proof would do for the weak null condition what Klainerman's work \cite{Klainerman1980} achieved for the classical null condition. The work presented here takes a significant step in this direction. Specifically, we identify a large class of systems that satisfy, in addition to the weak null condition, an extra \emph{hierarchical} condition on the nonlinearities. This hierarchical condition generalises the algebraic condition on the semilinear terms identified in \cite{Alinhac2006} and extended in \cite{Hidano2017}, in that our semilinear hierarchy is allowed to consist of more than two ``layers'', and (perhaps more significantly) in that we allow for \emph{quasilinear} wave equations.

After introducing the hierarchical weak null condition, we then show that all systems satisfying this condition possess global solutions for all sufficiently small initial data. As far as we are aware, all systems that have been studied so far, and which possess global solutions for all sufficiently small initial data, are actually of this hierarchical form. In particular, the Einstein equations in wave coordinates are a special case of the equations we study. However, we cannot show that \emph{all} equations satisfying the weak null condition are of this hierarchical form. The small-data global existence of solutions obeying the weak null condition but \emph{not} obeying our hierarchical condition (if any such equations exist!) thus remains open.

\vspace{2mm}

One example of a system obeying this weak null hierarchical condition is the Einstein equations in harmonic (or ``wave'') coordinates. This has already been shown to have global solutions for small initial data in the pioneering work of Lindblad and Rodnianski \cite{Lindblad2004}, however, their methods differ significantly from the methods used here. For one thing, we use the $r^p$-weighted energy method, which they did not; this is discussed in more detail below. Also, \cite{Lindblad2004} made significant use of the \emph{wave coordinate condition} itself\footnote{For the ``reduced Einstein equations'', which are a set of nonlinear wave equations, it turns out that, if this condition is imposed on the initial data, then it remains true throughout the solution.}. If this condition holds, then there are significant simplifications in many of the calculations, however, it cannot be assumed to hold for the more general equations we consider in this work.

It is worth pointing out that the simplifications offered by the wave coordinate condition can be utilised in at least two ways. First, it would be possible to obtain a much shorter and simpler proof of global stability, which nevertheless uses the $r^p$-weighted energy method. For example, we can use the extra structure in the Einstein equations to avoid the need to introduce a ``geometric foliation'', and this in turn would simplify our analysis drastically. There are many other issues, for example, involving error terms with ``critical'' decay rates, which can also be handled in a quicker and simpler way using the wave coordinate condition. Overall, this would lead to a much shorter and simpler proof, which would still retain the advantages of the $r^p$-weighted method.

Alternatively, we could retain the geometrical approach and instead use the extra structure to gain additional control over the solutions. For example, when carrying out the $p$-weighted energy estimates, it appears that it would be possible to take larger values of $p$ in the Einstein equations than in the more general systems that we consider, and this would lead to several important consequences, such as the existence of the radiation field and faster decay towards timelike infinity.

We do not carry out either of these programmes in full in this work. Instead, we have chosen to keep the analysis as general as possible, and so we have not assumed that the wave coordinate condition holds. This means that our approach can be applied to a much larger range of problems, and is not limited to the Einstein equations in wave coordinates. For example, we can treat the ``reduced Einstein equations'', which are the Einstein equations expressed in wave coordinates, but we \emph{do not} have to use initial data which satisfies the wave coordinate condition.

For a more complete discussion of the extra structures present in the Einstein equations, and the ways this can be used to help the analysis, see subsection \ref{subsection Einstein equations in wave coords}.

\vspace{2mm}

Before moving on to discuss some of the issues raised above in more detail, we briefly mention some of the methods used in our proof. We have endeavoured, as far as possible, to use methods which are as robust as possible, and which will generalise most easily to other settings. Of particular note is our use of the $r^p$-weighted energy method, introduced in \cite{Dafermos2010b}. This method avoids the use of a large number of approximate symmetries to obtain decay, and instead relies on showing that, through a suitable foliation, the \emph{energy} itself decays. It can therefore be applied in a range of situations in which the approximate symmetries are absent, or are otherwise difficult to make use of due to the presence of fields which do not interact well with these approximate symmetries\footnote{An example of the latter situation is the presence of a \emph{Klein-Gordon} field (i.e.\ a massive scalar field) coupled to the wave equations. For example, one can consider the stability of Minkowski space with a massive scalar field. In this case, the \emph{scaling} vector field does not have good commutation properties with the Klein-Gordon field, however, in this case this difficulty can be overcome (\cite{LeFloch2015, Wang2016})}. 

Furthermore, these methods provide a ``black box'' that can be used in other problems. Specifically, consider wave equations on a manifold made up of two regions: an ``asymptotic region'' which is similar to the space we study, joined to an ``interior region'' which might be significantly more complicated\footnote{Perhaps the most important of such situations are the ``asymptotically flat'' spacetimes of general relativity.}. In such a case, the methods presented in this paper can be easily adapted to handle the asymptotic region, provided that suitable control over the solution is established in the interior region, reducing the problem to one of studying the interior region alone. This kind of approach has already proved extremely fruitful in the context of \emph{linear waves} on black hole backgrounds - see, for example \cite{Dafermos2010a, Dafermos2014, Dafermos2011, Moschidis2015, Moschidis2015a}.

\vspace{10mm}

\begin{flushleft}\textbf{{\LARGE Overview}} \end{flushleft}

\section{The weak null condition}
\label{section intro weak null condition}

\begin{figure}[htb]
	\centering
	\includegraphics[width = 0.9\linewidth, keepaspectratio]{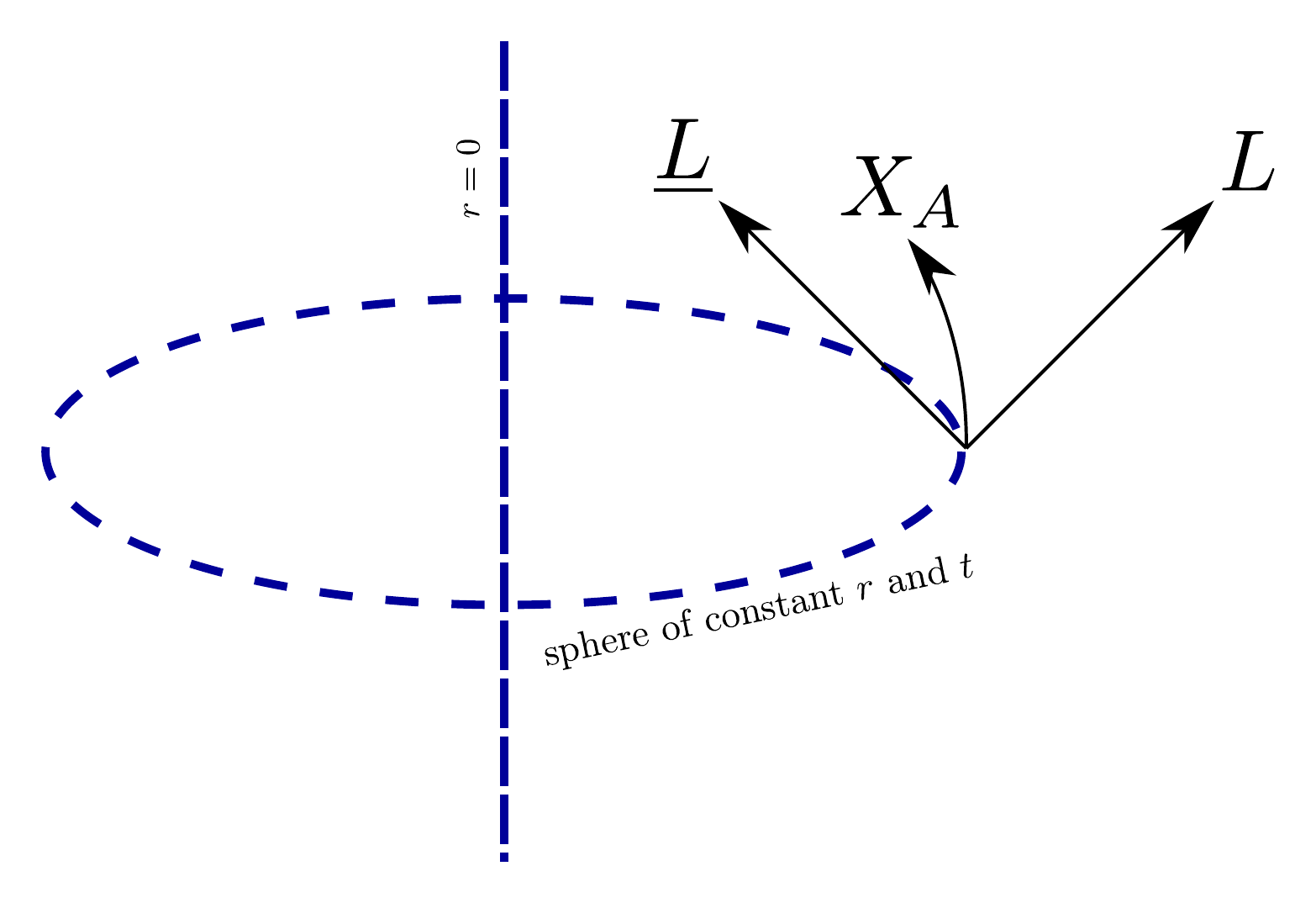}
	\caption{A diagram of the null frame vector fields. The vector fields $L$ and $\Lbar$ are both null, future directed, and orthogonal to the spheres of constant $r$ and $t$. The vector field $L$ points away from the origin $r = 0$ while the vector field $\Lbar$ points toward the origin. The vector fields $X_{A}$ are tangent to the spheres of constant $r$ and $t$. Note that we have suppressed a dimension, so these spheres appear as circles. Note also that in the null frame that we actually use, the time $t$ must be replaced by the ``geometric retarded time'' $\tau$.
	}
	\label{figure null frame}
\end{figure}

Let $\{L, \Lbar, X_A\}$ be a ``null frame'', consisting of the outgoing null vector field\footnote{Here, $t, x^1, x^2, x^3$ are standard rectangular coordinates on $\mathbb{R}^4$, and $r = \sqrt{(x^1)^2 + (x^2)^2 + (x^3)^2}$.} $L = \partial_t + \partial_r$, the ingoing null vector field $\Lbar = \partial_t - \partial_r$ and orthonormal vector fields $X_1$, $X_2$ tangent to the spheres of constant $t$ and $r$ (see figure \ref{figure null frame}). Then the classical null condition rules out any quadratic nonlinear terms proportional to $(\Lbar \phi)^2$. One way to understand why this ensures global existence is to realise that we expect a different decay rate for $(\Lbar \phi)$ compared to the other derivatives. For example, if $\phi$ satisfies the linear wave equation and arises from smooth, compactly supported data, then it is fairly easy to obtain the bounds
\begin{equation*}
\begin{split}
\begin{rcases}
|L\phi| \\ |X_1 \phi| \\ |X_2 \phi|
\end{rcases} \lesssim (1+t)^{-2}
\\
|\Lbar \phi| \lesssim (1+t)^{-1}
\end{split}
\end{equation*}
As such, we refer to the $\Lbar$ derivatives as ``bad derivatives'', and the other null frame derivatives as ``good derivatives''. The classical null condition ensures that, in every pair of derivatives appearing in the quadratic nonlinear terms, at least one derivative is a ``good derivative''. Note that, throughout this work, we are restricting to \emph{derivative nonlinearities}, at least in the quadratic nonlinear terms.

Motivated by these kinds of considerations, H\"ormander introduced \emph{asymptotic systems} (see \cite{Hoermander1987, Hoermander1997}). To construct the asymptotic system corresponding to a given set of nonlinear wave equations, we first ``throw away'' all the terms that involve \emph{only} good derivatives, together with quadratic terms that involve \emph{at least one} good derivative, and also all cubic and higher order terms. In this way, we can obtain transport equations for the variables $(r \Lbar \phi)$ along the outgoing null geodesics. For example, if $\phi$ satisfies a linear wave equation \emph{or} a wave equation obeying the classical null condition, then the corresponding equation for $\phi$ in the asymptotic system is
\begin{equation*}
L(r\Lbar \phi) = 0
\end{equation*}
The solution to this transport equation is trivially obtained by integrating from the surface on which we place initial data.

On the other hand, if $\phi$ satisfies \eqref{equation John's example}, then the corresponding asymptotic system can be written as
\begin{equation*}
L(r\Lbar \phi) = \frac{1}{4}r^{-1} (r\Lbar \phi)^2
\end{equation*}
If we define $u = t-r$ and work in coordinates $(u, r, \theta, \phi)$, then $L = \frac{\partial}{\partial r}$. So, with respect to this coordinate system, we can write the asymptotic system as
\begin{equation*}
\frac{\partial}{\partial \log r} (r\Lbar \phi) = \frac{1}{4}(r\Lbar \phi)^2
\end{equation*}
This has the form of a Riccati equation. Suppose that $(r\Lbar\phi)$ takes the value $(r\Lbar\phi)_0$ at some point on the initial data surface with $r$ coordinate $r_0$, and consider the solution along an integral curve of $L$ originating from this point, parameterised by $r$. Then the solution along this integral curve is given by
\begin{equation*}
(r\Lbar \phi)(r) = \frac{1}{ ((r\Lbar \phi)_0)^{-1} - \frac{1}{4}(\log r - \log r_0) }
\end{equation*}
In particular, if $((r\Lbar \phi)_0)^{-1}$ is ever positive, then the solution will become singular at the point
\begin{equation*}
r = r_0 e^{4((r\Lbar \phi)_0)^{-1}}
\end{equation*}
In particular, $(r\Lbar \phi)_0$ can be arbitrarily small, and yet the solution will still become singular in finite time.

In these two cases, at least, the asymptotic system correctly predicts the behaviour of solutions to the associated wave equations. This motivates the following:
\begin{definition}[The weak null condition]
	Let $\phi_{(a)}$ be a system of nonlinear wave equations such that the corresponding asymptotic system possesses global solutions for all sufficiently small initial data. Suppose also that, for all sufficiently small initial data, the solutions to the asymptotic system obey the bound
	\begin{equation*}
	|(r\Lbar\phi_{(a)})| \lesssim (1+r)^{C\epsilon}
	\end{equation*}
	where $\epsilon = \sup_{(\text{initial data})} |r\Lbar\phi_{(a)}|$ and $C$ is any fixed numerical constant.
	
	Then we say that the system satisfies the \emph{weak null condition}.
\end{definition}
The condition on the growth rate of the solution must be made for technical reasons - see chapter \ref{chapter weak null structure}.

\subsection{The semilinear hierarchy}

A simple example of a system which does \emph{not} obey the classical null condition but which \emph{does} obey the weak null condition is the following:
\begin{equation}
\label{equation weak null example 1}
\begin{split}
\Box \phi_1 &= 0
\\
\Box \phi_2 &= (\partial_t \phi_1)^2
\end{split}
\end{equation}
The presence of the term $(\partial_t \phi_1)^2$ in the equation for $\phi_2$ means that this system fails to satisfy the classical null condition. However, the corresponding asymptotic system is
\begin{equation*}
\begin{split}
\frac{\partial}{\partial \log r} (r\Lbar \phi_1) &= 0
\\
\frac{\partial}{\partial \log r} (r\Lbar \phi_2) &= \frac{1}{4} (r\Lbar \phi_1)^2
\end{split}
\end{equation*}
where, again, we are working in $(u,r,\theta,\phi)$ coordinates, where $u = t-r$. This evidently obeys the weak null condition: $r\Lbar \phi_1$ is just a constant, and $(r\Lbar \phi_2)$ can grow logarithmically. Similarly, it is easy to show that the original set of wave equations \eqref{equation weak null example 1} possesses global solutions for small initial data: these equations can first be solved for $\phi_1$, and then this can be treated as a source term in the equation for $\phi_2$. 

For a slightly more complicated example, we can consider the following system:
\begin{equation}
\label{equation weak null example 2}
\begin{split}
\Box \phi_1 &= 0
\\
\Box \phi_2 &= (\partial_t \phi_1)(\partial_t \phi_2)
\end{split}
\end{equation}
Again, the presence of the term $(\partial_t \phi_1)(\partial_t \phi_2)$ in the equation for $\phi_2$ means that this system fails to satisfy the classical null condition. This time, the corresponding asymptotic system is
\begin{equation*}
\begin{split}
\frac{\partial}{\partial \log r} (r\Lbar \phi_1) &= 0
\\
\frac{\partial}{\partial \log r} (r\Lbar \phi_2) &= \frac{1}{4} (r\Lbar \phi_1)(r\Lbar \phi_2)
\end{split}
\end{equation*}
This asymptotic system predicts that $(r\Lbar \phi_2)$ can grow like $r^{C\epsilon}$, so the system obeys the weak null condition. Again, it is possible to show that the asymptotic system accurately describes the global dynamics of the system\footnote{In other words, the solution to the asymptotic system and the solution to the original system of wave equations have the same asymptotics.}, but the proof of this is significantly more involved.

\begin{figure}[ht]
	\centering
	\includegraphics[width = \linewidth, keepaspectratio]{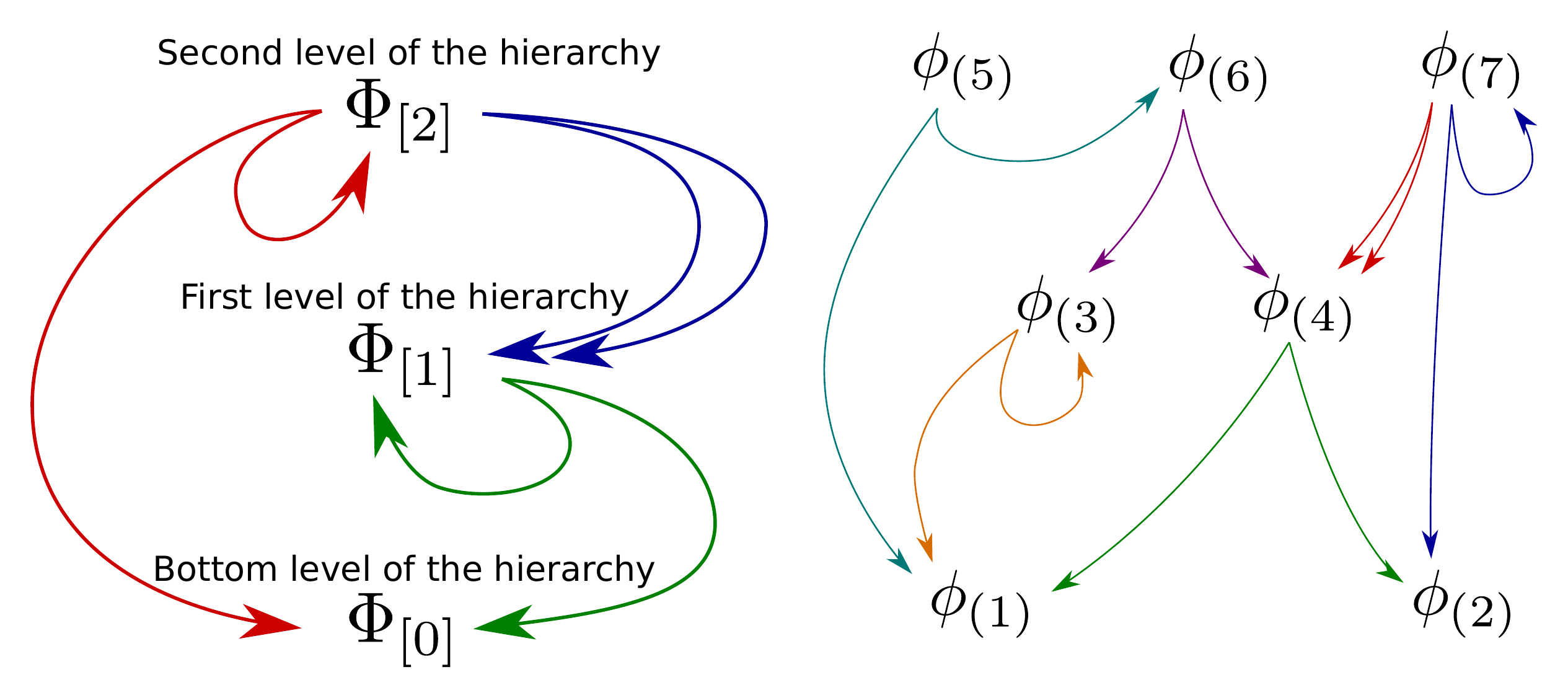}
	\caption{A figure illustrating the structure of the semilinear hierarchy in a case where the hierarchy consists of three levels. On the left, we illustrate the three levels, with the pairs of arrows between the levels representing \emph{possible} pairs of ``bad'' semilinear terms. On the right, we show a diagram of a particular system of wave equations satisfying the hierarchical null condition, also with three levels, and where each pair of ``legs'' represents a bad semilinear term. Note that each pair of legs either has the property that \emph{both} legs point to fields at a lower level of the hierarchy, \emph{or} it has the property that one leg points to a field at the bottom level and one leg points to a field at the same level. The system illustrated is of the form
	\\
	\protect\begin{minipage}{\linewidth}\protect\begin{align*}
	\tilde{\Box}_g \phi_{(1)} &= 0
	\\
	\tilde{\Box}_g \phi_{(2)} &= 0
	\\
	\tilde{\Box}_g \phi_{(3)} &= (\partial \phi_{(1)})(\partial \phi_{(3)})
	\\
	\tilde{\Box}_g \phi_{(4)} &= (\partial \phi_{(1)})(\partial \phi_{(2)})	
	\\
	\tilde{\Box}_g \phi_{(5)} &= (\partial \phi_{(1)})(\partial \phi_{(6)})
	\\
	\tilde{\Box}_g \phi_{(6)} &= (\partial \phi_{(3)})(\partial \phi_{(4)})
	\\
	\tilde{\Box}_g \phi_{(7)} &= (\partial \phi_{(4)})^2 + (\partial \phi_{(2)})(\partial \phi_{(7)})
	\protect\end{align*}\protect\end{minipage}
	}
	\label{figure semilinear hierarchy}
\end{figure}

Common to both of these examples is a certain \emph{hierarchical} structure in the semilinear terms. In fact, it turns out that with a suitable handling of quasilinear terms, \emph{all} examples which have been studied so far exhibit a hierarchical in their semilinear terms. This semilinear hierarchy is of the following form:

\begin{definition}[The semilinear hierarchy]
\hspace{1mm}
\begin{itemize}
	\item Each field $\phi_{(a)}$ can be placed into a certain set of fields $\Phi_{[n]}$, with $[n] = 0, 1, \ldots , N_1$
	\item The quadratic semilinear terms in the wave equations for fields in the set $\Phi_{[0]}$ obey the \emph{classical} null condition
	\item The quadratic semilinear terms in the wave equation for a field in the set $\Phi_{[n]}$ (with $n \geq 1$) that involve \emph{two} $\Lbar$ derivatives involve either a pair of fields from lower in the hierarchy, \emph{or} a pair consisting of one field in the set $\Phi_{[n]}$ together with a field in the set $\Phi_{[0]}$
\end{itemize}
see figure \ref{figure semilinear hierarchy} for a diagrammatic representation of this hierarchy.
\end{definition}

We will refer to such systems as systems satisfying the \emph{hierarchical weak null condition}. It is fairly easy to see that all such systems also satisfy the weak null condition, as defined above. It is equations of this form that we will handle in this paper. Indeed, the main result of this paper is that systems obeying the hierarchical weak null condition admit global solutions for sufficiently small initial data.

\subsection{Quasilinear equations}
In the discussion above, we have focussed exclusively on semilinear equations. However, there are also examples of \emph{quasilinear} wave equations which obey the weak null condition. For example, in \cite{Alinhac2003} Alinhac studied equations of the form
\begin{equation}
\label{equation Alinhac's quasilinear example}
\Box \phi + c(\phi)^2 \slashed{\Delta} \phi = 0
\end{equation}
where $\slashed{\Delta}$ is the Laplacian operator on the sphere of radius $r$ and $c(0) = 1$ is smooth. Also, in \cite{Lindblad2008} Lindblad studied the equation
\begin{equation}
\label{equation Lindblad's quasilinear example}
(g^{-1})^{ab}(\phi) \partial_a \partial_b \phi = 0
\end{equation}
where this equation is given in the standard rectangular coordinates, and where $(g^{-1})^{ab}(\phi) = m^{ab} + h^{ab}(\phi)$, with $m^{ab}$ the (inverse of the) Minkowski metric, and $h^{ab}(\phi)$ some tensor with rectangular components that are \emph{linear}\footnote{In fact, quadratic and higher order terms can be included in the metric, but these lead to \emph{cubic} (and higher order) nonlinearities in the wave equation, which are much easier to handle.} in $\phi$. Finally, Lindblad-Rodnianski studied the Einstein equations in wave coordinates in \cite{Lindblad2004}. In this form, the Einstein equations take the form of a system of quasilinear wave equations for the wave-coordinate components of the metric.

In all of these cases, the authors were able to show global existence for sufficiently small initial data, despite the quasilinear nature of the equations. In fact, these equations were also shown to obey the weak null condition (see \cite{Lindblad1992} for the first two equations, and \cite{Lindblad2003} for the Einstein equations). A key point is that, in the case of a quasilinear system, the wave operator needs to be  \emph{modified} by the addition of certain semilinear terms. This is the reason why, for example, Lindblad studied the equation $(g^{-1})^{ab} \partial_a \partial_b \phi = 0$ rather than $\Box_g \phi = 0$: these two equations differ by certain semilinear terms, related to the Christoffel symbols in rectangular coordinates. This modified wave operator is sometimes called the \emph{reduced} wave operator.

In this paper, we generalise this idea to more complicated systems of quasilinear wave equations, and show how to construct the appropriate ``reduced'' wave operator for these systems. The semilinear hierarchy is then defined with respect to this reduced wave operator, i.e.\ we require that the hierarchical structure is present in the semilinear terms \emph{after} making the suitable modification to the wave operator. Furthermore, the component of the metric perturbation\footnote{Specifically, the metric is given in terms of rectangular indices as $g_{ab} = m_{ab} + h_{ab}$, where $m$ is the Minkowski metric, and the $h_{ab}$ depend on the fields $\phi_{(a)}$. Then the component in question is $h_{LL}$.} which multiplies two ``bad'' derivatives is required to be in the set of fields $\Phi_{[0]}$, that is, it is required to be at the bottom level of the hierarchy. Note that the cases covered by Alinhac and Lindblad are \emph{scalar} wave equations, so there is no question of a hierarchy, whereas the Einstein equations in wave coordinates are rather special, since, along with the wave equations, we \emph{also} have the harmonic coordinate condition itself, which is not a wave equation but a relation between certain first derivatives of the metric components. See chapter \ref{chapter weak null structure} for further details of the general quasilinear structures that we consider, and see subsection \ref{subsection Einstein equations in wave coords} for further discussion of the Einstein equations in wave coordinates.

\subsection{Changing the basis sections}
\label{subsection intro changing basis sections}

A natural question to ask, given the discussion above, is whether \emph{all} systems of equations which satisfy the weak null condition also satisfy the hierarchical weak null condition. Unfortunately, this is not the case, and in fact there exist rather trivial counterexamples. Consider the system
\begin{equation*}
\begin{split}
\Box \phi_1 &= (\partial_t \phi_1 - \partial_t \phi_2)^2
\\
\Box \phi_2 &= (\partial_t \phi_1 - \partial_t \phi_2)^2
\end{split}
\end{equation*}
This system does not exhibit a hierarchy in its semilinear terms. Nevertheless, it does satisfy the weak null condition, and indeed possesses global solutions for sufficiently small data. This is easily seen by defining $\phi_3 := \phi_1 - \phi_2$; then the system for the variables $\phi_3$, $\phi_2$ is identical in form to the system \eqref{equation weak null example 1}.

Of course, our proof also holds in for systems of this kind. In other words, if a system can be put into a form which \emph{does} exhibit a semilinear hierarchy by a redefinition of the fields, then our proof clearly shows that this system possesses global solutions for sufficiently small initial data.

Stated more geometrically, the example above illustrates a ``\emph{change of basis sections}'': see chapter \ref{chapter weak null structure} for a more detailed discussion. Informally, this might be called a ``field redefinition''. Suppose that we have a set of $N_1$ fields satisfying a system of nonlinear wave equations on a manifold $\mathcal{M}$. Then we can consider changing variables to a different set of $N_1$ fields, satisfying equations that can be deduced from the original set of equations. It may be the case that the hierarchy is only evident after such a ``field redefinition''.

\begin{figure}[htpb]
	\centering
	\includegraphics[width = \linewidth, keepaspectratio]{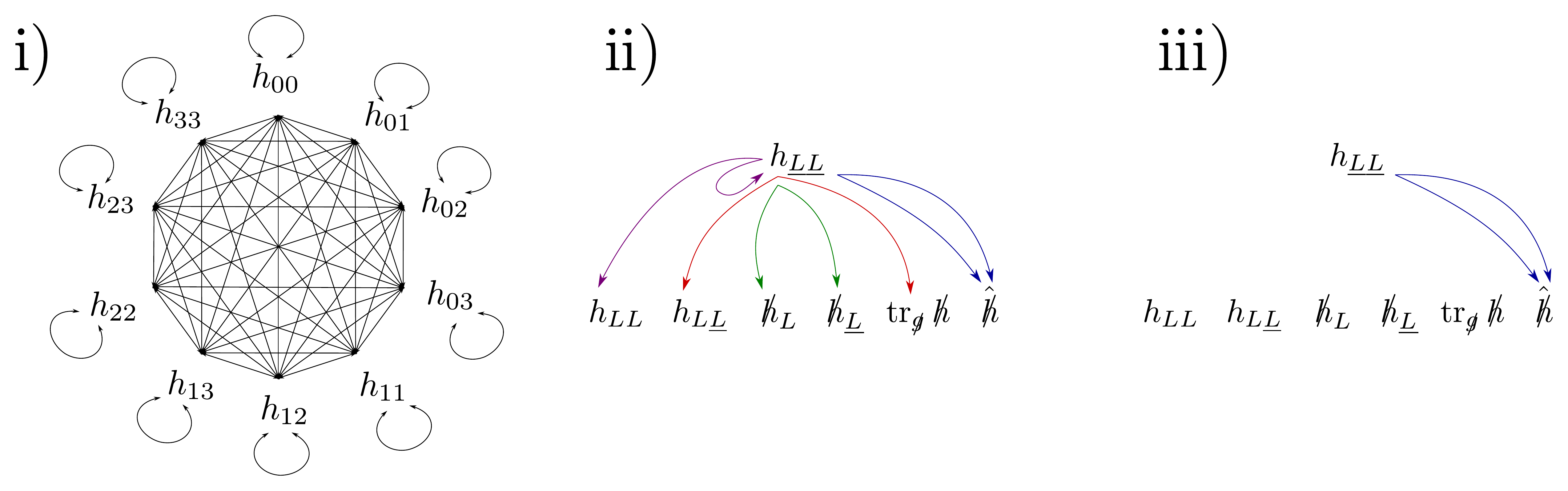}
	\caption{
	A diagram of the semilinear structures present in the Einstein equations in wave coordinates. As before, a pair of arrows originating at the field $\phi$ and pointing to the fields $\psi_1$ and $\psi_2$ means that, in the wave equation for $\phi$, the ``bad'' semilinear term $(\Lbar \psi_1)(\Lbar \psi_2)$ appears.
	\\
	In subfigure i), we see the semilinear structure present in the Einstein equations in wave coordinates \emph{when expressed in terms of the wave-coordinate components} of the metric perturbation $h$. Every possible pair of arrows is present, since every possible combination of ``bad'' terms appears in every equation. In other words, there is no structure present in the semilinear terms: the weak null condition is not evident.
	\\
	In subfigure ii), we show the semilinear structure in the Einstein equations in wave coordinates \emph{expressed in terms of the null frame components of} $h$. Note that it is natural to regard the field $(\slashed{h}_L)_\mu = \slashed{\Pi}_{\mu}^{\protect\phantom{\mu}a} L^b h_{ab}$ as a single field, which, however, takes values in the space of covector fields that are tangent to the spheres. In other words, $\slashed{h}_L$ defines a section of the vector bundle $\mathcal{B}$, rather than a section of $\mathcal{M}\times\mathbb{R}$. The field $\slashed{h}_{\Lbar}$ can be treated identically, while the field $\hat{\slashed{h}}$ takes values in the space of symmetric, trace-free, rank $(0,2)$ tensor fields on the spheres.
	\\
	Note that subfigure ii) shows the presence of a semilinear hierarchy in the Einstein equations. In fact, there are only two levels of the hierarchy, and only the field ($h_{\Lbar\Lbar}$) appears above the bottom level. This is the structure that is present in the Einstein equations in wave coordinates if we simply use the wave coordinate condition to write the Einstein equations as a system of quasilinear wave equations (sometimes called the ``reduced Einstein equations'').
	\\
	Subfigure iii) shows what happens if we also use the wave coordinate condition to rewrite some of the ``bad'' derivatives in terms of ``good derivatives'' of other fields. In other words, we are no longer just dealing with the reduced Einstein equations: we have coupled these equations with the constraint equations and the wave coordinate condition. In this case, there are still two levels in the semilinear hierarchy, but the semilinear structure is much simpler, and there is only one ``bad'' semilinear term.
	}
	\label{figure Einstein semilinear}
\end{figure}

The Einstein equations in wave coordinates offer another example of this kind of phenomenon. In this case, the wave-coordinate components of the metric satisfy a set of nonlinear wave equations, however, in this form the nonlinear terms in these equations do not have a hierarchical structure. A hierarchical structure is, however, present (see \cite{Lindblad2003}); it is made manifest by considering the \emph{null-frame} components of the metric, rather than the wave-coordinate components (see figure \ref{figure Einstein semilinear}). In this case, the transformation from wave-coordinate components to null-frame coordinates is \emph{point dependent}, that is, a different transformation\footnote{This is the case even if the metric $g$ is the flat Minkowski metric $m$ since, for example, we have $L^i = \frac{x^i}{r}$.} needs to be applied at each point in the manifold $\mathcal{M}$. Hence, when expanding the fields with respect to the new basis of sections, the new fields will satisfy wave equations with additional inhomogeneous terms, due to the action of the wave operator on the change-of-basis operator. Handling this needs extra care, particularly when the change-of-basis also depends on the solution itself, as is the case for the Einstein equations - in this case, the null frame is constructed to be null with respect to the metric, but the metric is the solution of the equations! However, we are also able to handle this kind of case with our methods.

In summary, we can extend our theorem to apply to systems of equations which only exhibit a semilinear hierarchy \emph{after a change of basis sections}, and moreover, we can allow (with some technical restrictions) this change-of-basis to depend on the point on the manifold at which we change the basis, and even to depend on the solution to the wave equation.

\vspace{10mm}
We are now ready to state the first version of the theorem proved in this work:

\begin{theorem}[Main theorem, first version]
	\label{theorem main theorem rough statement}
	All systems of nonlinear wave equations with derivative nonlinearities satisfying the hierarchical weak null condition possess global solutions for sufficiently small initial data.
\end{theorem}

Note that \cite{Alinhac2006} considers a special case of this hierarchical structure, where all but one field is at the ``bottom level'', and consequently there are only two levels in the hierarchy. Furthermore, \cite{Alinhac2006} only considered \emph{semilinear} equations. Also, despite the title of this article, it did not establish that ``blowup at infinity'' \emph{actually} occurs - this is done for the first time, as far as we are aware, in appendix \ref{appendix explicit shock formation} of this work.

\subsection{Example systems}
\label{subsection intro example systems}

We now provide some examples of systems to which our analysis applies. We have already mentioned the semilinear systems given in equations \eqref{equation weak null example 1} and \eqref{equation weak null example 1}. More generally, we could consider the semilinear system
\begin{equation*}
\Box \phi_{(a)} = \left(F_{(a)}^{(b)(c)} \right)^{\mu\nu} \partial_\mu \phi_{(b)} \partial_\nu \phi_{(c)}
\end{equation*}
where the tensor fields $F$ are independent of the fields $\phi_{(a)}$ and satisfy
\begin{equation*}
\begin{split}
\left(F_{(a)}^{(b)(c)} \right)^{\mu\nu} &= \left(F_{(a)}^{(c)(b)} \right)^{\mu\nu}
\\
\left(F_{(a)}^{(b)(c)} \right)_{LL} &:= \left(F_{(a)}^{(b)(c)} \right)^{\mu\nu} L_\mu L_\nu
\\
\left(F_{(a)}^{(b)(c)} \right)_{LL} &= 0 \quad \text{if } \begin{cases} a > b \text{ , or}\\ a = b \text{ and } c \geq 1 \end{cases}
\end{split}
\end{equation*}
This generalises systems \eqref{equation weak null example 1} and \eqref{equation weak null example 1} to hierarchies of arbitrarily high order, and to a wider class of nonlinearities.

\vspace{4mm}

We can also handle quasilinear equations, including the scalar quasilinear equations studied by Alinhac \eqref{equation Alinhac's quasilinear example} and Lindblad \eqref{equation Lindblad's quasilinear example}. We can combine these examples with those above, to form quasilinear \emph{systems} of nonlinear wave equations with the semilinear structure: that is, systems of the form
\begin{equation*}
\tilde{\Box}_g \phi_{(a)} = \left(F_{(a)}^{(b)(c)} \right)^{\mu\nu} \partial_\mu \phi_{(b)} \partial_\nu \phi_{(c)}
\end{equation*}
where $\tilde{\Box}_g$ is the ``reduced'' wave operator associated with the metric $g$ (see chapter \ref{chapter weak null structure}) and the tensor fields $\left(F_{(a)}^{(b)(c)} \right)^{\mu\nu}$ are required to satisfy the same conditions as above. Here, $g = m + h$, where $m$ is the Minkowski metric and $h = \mathcal{O}(\phi)$ is such that the field $h_{LL}$ can be expressed in terms of fields $\phi_{(a)}$ that are at the \emph{bottom} level of the hierarchy.

\vspace{4mm}

Additionally, we can handle situations in which a point-dependent change of basis is required to make the semilinear hierarchy manifest. For example, consider the nonlinear covector wave equation
\begin{equation*}
\Box U_\mu = (\Div U) K^\nu (\D_\mu U_\nu)
\end{equation*}
where $K^\nu$ is some fixed (bounded) covector field and $\D$ is the covariant derivative operator. Here, $U_\mu$ is to be treated as a \emph{spacetime} covector. This equation could be ``scalarised'' by expressing it with respect to the usual rectangular vector fields, in which case it takes the form
\begin{equation*}
\Box U_a = (\Div U) K^c (\partial_a U_c)
\end{equation*}
where now the rectangular components of $U$, $U_a$, are a set of scalar fields. The semilinear hierarchy in this system can be made manifest by expressing the system in a null frame, in which it reads
\begin{equation*}
\begin{split}
	(\Box U)_L 
	&= 
	\left( -\frac{1}{2} (LU)_{\Lbar} - \frac{1}{2} (\Lbar U)_{L} + (\slashed{g}^{-1})^{AB} (X_{A} U)_{B} \right) (L U)_K
	\\
	(\Box U)_{\Lbar}
	&=
	\left( -\frac{1}{2} (LU)_{\Lbar} - \frac{1}{2} (\Lbar U)_{L} + (\slashed{g}^{-1})^{AB} (X_{A} U)_{B} \right) (\Lbar U)_K
	\\
	(\Box U)_{A}
	&=
	\left( -\frac{1}{2} (LU)_{\Lbar} - \frac{1}{2} (\Lbar U)_{L} + (\slashed{g}^{-1})^{AB} (X_{A} U)_{B} \right) (X_A U)_K
\end{split}
\end{equation*}
see chapter \ref{chapter preliminaries} for an overview of the notation used here. One important point is that derivatives act on the \emph{rectangular components}, as indicated by the fact that the frame indices are placed outside of brackets. Hence, for example,
\begin{equation*}
(LU)_L := (L U_a) L^a
\end{equation*}

Importantly, the only term involving two ``bad'' derivatives is the term $(\Lbar U)_L (\Lbar U)_K$, which appears in the equation for the field $U_{\Lbar}$. However, this ``bad'' semilinear term involves the field $U_L$, while the equation for the field $U_L$ has the classical null structure. In other words, $U_L$ is at the bottom level of the hierarchy (as are the fields $U_A$), while $U_{\Lbar}$ is at the first level of the hierarchy. A semilinear version of this equation could also be developed, in which the metric is allowed to depend on the fields $U$, but care must be taken in this case to ensure that the principle symbol does not change.

\vspace{4mm}

Next, we mention the Einstein equations in wave coordinates. These were shown to admit global solutions in the pioneering work of \cite{Lindblad2004}. The equations take the form
\begin{equation*}
\tilde{\Box}_g h_{ab} = F_{ab} (\partial h, \partial h) + \mathcal{O}(h (\partial h)^2)
\end{equation*}
where the components of the metric $g$ in wave coordinates are given by
\begin{equation*}
g_{ab} = m_{ab} + h_{ab}
\end{equation*}
where $m_{ab}$ are the components of the Minkowski metric in the standard rectangular coordinates\footnote{i.e.\ $m_{00} = -1$, $m_{ij} = \delta_{ij}$ and $m_{i0} = 0$}. The fields inhomogeneous terms $F_{ab}$ are given by
\begin{equation*}
F_{ab} (\partial h, \partial h)
=
(g^{-1})^{cd}(g^{-1})^{ef} \left( \frac{1}{4} (\partial_a h_{cd}) (\partial_b h_{ef}) - \frac{1}{2} (\partial_a h_{ce}) (\partial_b h_{df}) \right) + Q_{ab} (\partial h, \partial h)
\end{equation*}
where the semilinear terms $Q_{ab}$ satisfy the classical null condition. Again, expressed relative to these wave coordinates, no semilinear hierarchy is present. However, if we express these fields relative to the null frame, then we find that the only term which does not obey the classical null condition is in fact
\begin{equation*}
F_{\Lbar\Lbar}(\partial h, \partial h)
=
-\frac{1}{4} (\Lbar h)_{LL} (\Lbar h)_{\Lbar\Lbar}
-\frac{1}{2}(\Lbar h)_{L\Lbar} (\slashed{g}^{-1})^{\slashed{\mu}\slashed{\nu}} (\Lbar h)_{\slashed{\mu}\slashed{\nu}}
+\frac{1}{2}(\slashed{g}^{-1})^{\slashed{\mu}\slashed{\nu}} (\Lbar h)_{L\slashed{\mu}} (\Lbar h)_{\Lbar\slashed{\nu}}
- \widehat{|\Lbar \slashed{h}|}^2
\end{equation*}
where $\widehat{(\Lbar \slashed{h})}_{\slashed{\mu}\slashed{\nu}}$ is the trace-free part of the tensor $(\Lbar \slashed{h})_{\slashed{\mu}\slashed{\nu}}$ (again, see chapter \ref{chapter preliminaries} for the notation used here). The important point is that this appears as the semilinear term in the equation for $h_{\Lbar\Lbar}$, but there is no term of the form $(\Lbar h)_{\Lbar\Lbar}^2$ in $F_{\Lbar\Lbar}$. In other words, every null frame metric component is at the bottom level of the hierarchy, except for $h_{\Lbar\Lbar}$ which is at the first level of the hierarchy.

In fact, there is a lot more structure in the Einstein equations that can be exploited, which can give improved behaviour relative to most of the other systems discussed above. This additional structure was used extensively throughout the proof of Lindblad-Rodnianski \cite{Lindblad2004}. See subsection \ref{subsection Einstein equations in wave coords} for a more detailed discussion of this extra structure and the ways in which it can be exploited. One important point to note, however, is that these other structures require both the \emph{gauge condition} and the \emph{constraint equations} to be satisfied. We can, instead, consider the Einstein equations in wave coordinates, and then \emph{forget} the gauge condition and the constraint equations. This might actually be important for certain applications - for example, in numerical relativity, it is often impossible to impose a gauge condition and the constraint equations \emph{exactly}. Alternatively, on a more technical note, certain proofs rely on an iterative procedure, where solutions only satisfy the constraints or the gauge conditions in the limit. In these cases, the work done here shows (for the first time) that the resulting equations admit global solutions for sufficiently small data\footnote{Note that, in general, solutions to these equations will have different asymptotics to solutions that also satisfy the constraints and the gauge conditions.}.

\vspace{4mm}

Finally, we mention the coupled Einstein-Maxwell system, in harmonic coordinates and Lorenz gauge. This was first treated in \cite{Zipser2000}, using a framework similar to that used by Christodoulou and Klainerman in \cite{Christodoulou1993}. A treatment similar to that outlined below was later given in \cite{Loizelet2008} (see also \cite{Speck2010}). Here, we work in harmonic coordinates $x^a$ as above, so the functions $x^a$ satisfy
\begin{equation*}
\Box_g x^a = 0
\end{equation*}
We also work in Lorenz gauge, so the field strength $F$ is given by the exterior derivative of a one-form $A$ which in turn satisfies the gauge condition
\begin{equation*}
\Div A = \D^\mu A_\mu = 0
\end{equation*}
Using the Lorenz gauge condition, the Maxwell equations $\D^\mu F_{\mu\nu} = 0$ can be seen to imply the following equation for the one-form $A$:
\begin{equation*}
\Box_g A_\mu - R_\mu^{\phantom{\mu}\nu} A_\nu = 0
\end{equation*}
where $R_{\mu\nu}$ is the Ricci curvature of the manifold. Expanding this in terms of the basis of one forms associated with harmonic coordinates $\upd x^a$, and using the harmonic coordinate condition, we find that the harmonic coordinate components of $A$ satisfy the wave equations
\begin{equation*}
\Box_g A_a - (g^{-1})^{cd} (g^{-1})^{be} (\partial_c h_{ae} + \partial_a h_{ce} - \partial_e h_{ac}) (\partial_d A_b) = 0
\end{equation*}
Similarly, we can deduce the equations satisfied by the harmonic coordinate components of the metric component perturbations $h_{ab}$.

We will call the system obtained in this way the \emph{reduced Einstein-Maxwell equations}. If we express this system relative to the null frame, then we can see that is obeys a semilinear hierarchy. Specifically, dropping terms involving ``good'' derivatives and ignoring numerical constants, we obtain a system of the following form:

\begin{equation*}
\begin{split}
	(\tilde{\Box}_g h)_{LL} &\sim 0 \\
	(\tilde{\Box}_g h)_{L\Lbar} &\sim 0 \\
	(\tilde{\Box}_g \slashed{h})_{L} &\sim 0 \\
	\\
	(\tilde{\Box}_g A)_{L} &\sim (\Lbar A)_{\Lbar} (\Lbar h)_{LL} \\
	(\tilde{\Box}_g A)_{\Lbar} &\sim 
		(\Lbar A)_L (\Lbar h)_{L\Lbar} 
		+ (\Lbar A)_{\Lbar} (\Lbar h)_{LL} 
		+ (\Lbar \slashed{A})\cdot(\Lbar \slashed{h})_{L} \\
	(\tilde{\Box}_g \slashed{A}) &\sim 
		(\Lbar A)_L (\Lbar \slashed{h})_{L}
		+ (\Lbar \slashed{A}) (\Lbar h)_{LL} \\
	\\
	(\tilde{\Box}_g \slashed{h})_{\Lbar} &\sim (\Lbar A)_L (\Lbar \slashed{A}) \\
	(\tilde{\Box}_g \slashed{h}) &\sim \left((\Lbar A)_L\right)^2 \\
	\\
	(\tilde{\Box}_g \slashed{h})_{\Lbar\Lbar } &\sim 
		(\Lbar h)_{LL}(\Lbar h)_{\Lbar\Lbar}
		+ (\Lbar h)_{L\Lbar} (\Lbar \slashed{h})
		+ (\Lbar \slashed{h})_L (\Lbar \slashed{h})_{\Lbar}
		+ (\Lbar \slashed{h})^2
		+ (\Lbar A)_{L} (\Lbar A)_{\Lbar} 
		+ (\Lbar \slashed{A})^2
\end{split}
\end{equation*}
Showing that this system obeys a semilinear hierarchy with \emph{four} levels. Hence our proof shows that this system admits global solutions for all sufficiently small initial data. Note, however, that the gauge conditions can be used again to replace certain ``bad'' derivatives with ``good'' ones. If we perform these substitutions, then the system simplifies considerably and we find
\begin{equation*}
\begin{split}
(\tilde{\Box}_g h)_{LL} &\sim 0 \\
(\tilde{\Box}_g h)_{L\Lbar} &\sim 0 \\
(\tilde{\Box}_g \slashed{h})_{L} &\sim 0 \\
(\tilde{\Box}_g A)_{L} &\sim 0 \\
(\tilde{\Box}_g A)_{\Lbar} &\sim 0\\
(\tilde{\Box}_g \slashed{h})_{\Lbar} &\sim 0 \\
(\tilde{\Box}_g \slashed{h}) &\sim 0 \\
(\tilde{\Box}_g \slashed{h})_{\Lbar\Lbar } &\sim 
	(\widehat{\Lbar \slashed{h}})^2
	+ (\Lbar \slashed{A})^2
\end{split}
\end{equation*}
which exibits a hierarchy with only \emph{two} levels. Note that these two systems will exhibit different asymptotics. The first system can be treated entirely as a system of nonlinear wave equations, while the second system requires us to posit initial data which satisfies both the harmonic coordinate condition and the Lorenz gauge condition. Nevertheless, both systems admit global solutions for sufficiently small initial data.

\section{The \texorpdfstring{$r^p$}{rp}-weighted energy method}
\label{section intro rp weighted method}

Now that we have detailed the kinds of equations to which our proof applies, and placed these in the context of previous results, it is time to turn our attention to the \emph{methods} used to prove global existence. We employ various different techniques to achieve our results, but, conceptually, the central method is the method of $r^p$-\emph{weighted energy estimates}. This was introduced in \cite{Dafermos2010b} in the context of \emph{linear} wave equations, in which context it has proved extremely versatile and useful - see, for example, \cite{Aretakis2011, Holzegel2010, Dafermos2010a, Dafermos2014, Dafermos2014a, Dafermos2016, Keir2016, Johnson2018} (see also \cite{Lindblad2006} for a related, though conceptually different approach using fractional Morawetz estimates). On the other hand, the method has so far seen fairly limited use in the context of nonlinear wave equations, although it has been used for semilinear equations with the classical null condition \cite{Yang2014} as well as certain quasilinear equations in which the quasilinear terms are particularly well behaved (\cite{Yang2013a, Yang2013}).

At the heart of this approach is the idea that, through a suitable foliation, the \emph{energy} of a solution to a wave equation should decay. In this way, the method differs from other vector-field based methods (see \cite{Klainerman1986} for the ``commuting vector fields'' approach, and e.g.\ \cite{Morawetz1968} for the ``vector field multiplier'' approach), which typically seek to show that the energy is conserved (or, possibly, grows slowly). When dealing with nonlinear equations, establishing boundedness is typically insufficient for closing the argument, and we generally find that it is necessary to prove some form of \emph{pointwise decay} for certain quantities. Both the $r^p$-weighted energy method and these older vector field methods rely on \emph{commuting} the equation with a certain set of vector fields, before applying a suitable Sobolev inequality, in order to achieve this decay. They differ in that the older methods rely \emph{entirely} on the use of a set of vector fields which \emph{grow} in a suitable manner to obtain this decay. In contrast, the $r^p$-weighted energy method obtains some of the required decay (specifically, decay towards timelike infinity) simply by establishing that the energy itself decays.

The $r^p$-weighted method has many advantages  over alternative methods. In common with other vector field based methods, it is extremely robust, and can be applied to a variety of different situations with minimal modifications. In particular, it does not rely on any \emph{exact} special properties of the manifold on which we are seeking a solution, which means that it is suitable for use with quasilinear problems, where the metric typically varies with the solution, and the resulting solutions do not posses any exact symmetries. In addition, the $r^p$-weighted method has an advantage over other vector field methods in that it does not require commuting with as many vector fields, since (as mentioned above) we do not need to obtain \emph{all} of the decay through the use of suitably growing vector fields. This, in turn, brings at least three advantages. First, the $r^p$-weighted method is suitable for use on manifolds which do not possess as many approximate\footnote{Here, ``approximate symmetries'' should be understood as meaning that the error terms one encounters if the symmetry is not exact can all be dealt with.} symmetries - for example, black hole spacetimes, which do not admit translation, boost or scaling symmetries. Second, the method does not require an \emph{a priori} assumption that the metric eventually settles down to some given metric (see \cite{Yang2013a}). Third, the $r^p$-weighted method is suitable for handling wave equations that are coupled to other types of equations, in the case where those other equations do not have good commutation properties with all of the vector fields used in the older methods.

Another big advantage of the $r^p$-weighted method, in comparison to other vector field based methods, is that it achieves a certain \emph{decoupling} of problems relating to wave propagation (this point was already appreciated and emphasised in \cite{Dafermos2010b}). It is common to encounter problems where the setting can be divided into two regions: an ``interior'' region exhibiting complicated dynamics, and an ``exterior'' or ``asymptotic'' region which is similar to flat space\footnote{Again, black hole spacetimes provide an illustrative example of this situation.} In these cases, the $r^p$-weighted method can be used to control the solution in the asymptotic region, provided that certain other estimates can be established in the interior region. There are many different ``interior'' regions, but all ``exterior'' regions exhibiting suitable asymptotics can be handled \emph{at once} using the $r^p$-weighted energy estimates. This leaves the only problem as that of establishing suitable estimates in the interior. For this reason, the $r^p$-weighted energy estimates have been referred to as a ``black box'' for handling the asymptotic region.

\subsection{The need for a geometric foliation}
\label{subsection need for geometric foliation}

As stated above, at the core of the $r^p$-weighted method is the idea that, through a suitable foliation, the energy of waves decays. It turns out that a ``suitable'' foliation is one whose leaves are \emph{achronal}\footnote{Strictly speaking, achronality is not required - \cite{Yang2013} used a foliation by leaves which are asymptotically null but might be timelike in places. Nevertheless, there are ``enough'' leaves which are asymptotically achronal in a suitable fashion for this foliation to work.}, \emph{asymptotically null}, and where points are related to each other by flowing along a vector field that is \emph{uniformly timelike}\footnote{This last condition is required in order to rule out foliations like the ``hyperboloidal foliation'' used, for example, in \cite{Klainerman1985, LeFloch2014, LeFloch2015}, in which the leaves of the foliation asymptotically approach each other, while also being asymptotically null. Because the leaves approach one another asymptotically, we cannot expect the energy to decay through such a foliation!}.

In this paper, we are dealing with \emph{quasilinear} wave equations, i.e.\ equations in which the metric depends on the solution to the equations. At the same time, we are also working with equations on $\mathbb{R}\times\mathbb{R}^3$. Hence there are two possible interpretations of a ``null hypersurface'', which we will refer to as the ``background'' and the ``geometric'' interpretations respectively. Since our foliations are only required to be \emph{asymptotically} null, they can be made to agree with one another in some spatially compact region.

On the one hand, we can use the ``background'' structure of the manifold as $\mathbb{R}\times\mathbb{R}^3$ to define ``null'' hypersurfaces. That is, if $(t, x^1, x^2, x^3)$ are the standard coordinates on $\mathbb{R}\times\mathbb{R}^3$ (i.e.\ $t \in \mathbb{R}$ and $(x^1, x^2, x^3) \in \mathbb{R}^3$) then we can define the ``background'' radial coordinate $r := \sqrt{ (x^1)^2 + (x^2)^2 + (x^3)^2 }$. Then, we can define the variable $u_{(\text{background})}$ as $u = t - r$. A possible foliation of our spacetime is then given by joining the level sets of $u_{(\text{background})}$ in the region $r \geq r_0$ to the level sets of $t$ in the region $r \leq r_0$. We refer to this as the ``background null foliation''. In the region $r \geq r_0$ it is null with respect to the Minkowski metric $m$, which, in the $(t, x^1, x^2, x^3)$ coordinate system, has components $m = \text{diag}(-1,1,1,1)$.

\begin{figure}[htbp]
	\centering
	\includegraphics[width = \linewidth, keepaspectratio]{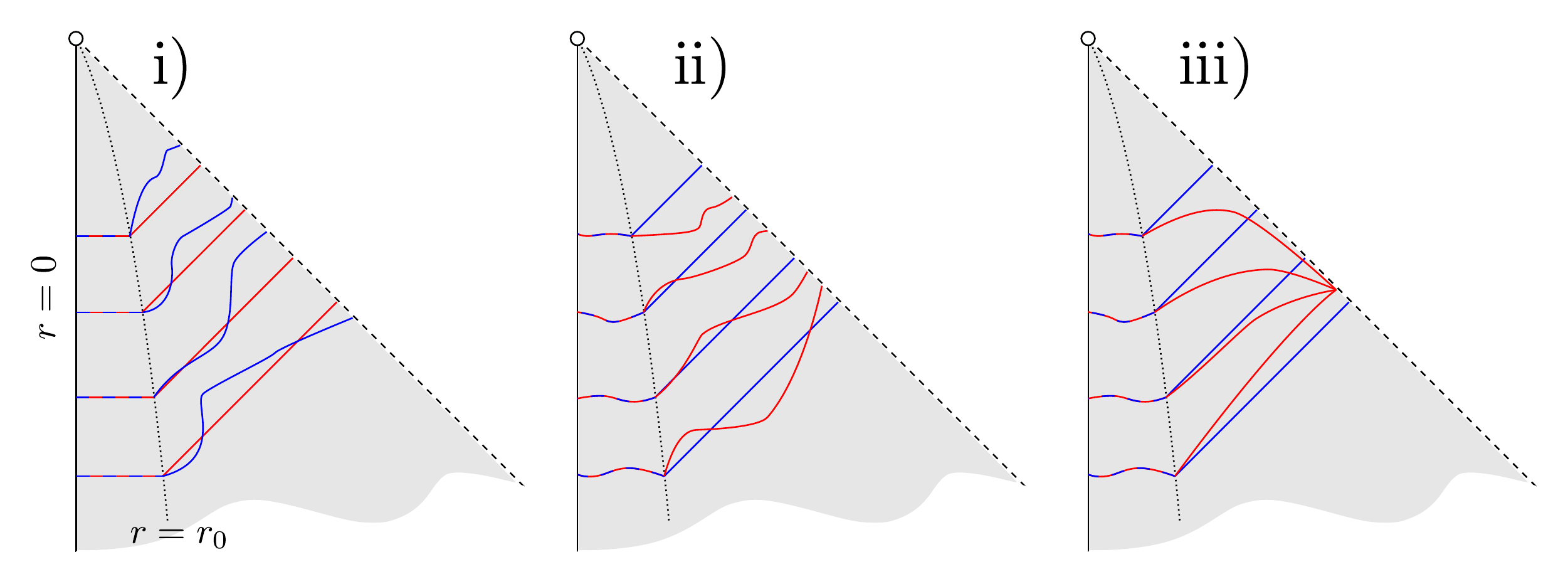}
	\caption{
	``Penrose diagrams'' illustrating the \emph{geometric} and \emph{non-geometric} foliations with respect to different compactifications. In all three figures, the  red lines depict the leaves of the ``background'' foliation, while the blue lines represent the leaves of the ``geometric'' foliation. The solid black line on the left of each figure is the line $r = 0$, and the dotted black line in the interior of each figure represents the hypersurface $r = r_0$. Note that all of the foliations agree in the region $r \leq r_0$.
	\\
	In subfigure i), we have compactified with respect to the ``background'' geometry, so that the lines of constant $(t-r)$ are straight lines at forty five degrees. In other words, this is a conformal compactification with respect to the ``background'' metric $m$, and not the metric $g$. Note that, if we reduce the space to a two dimensional space by quotienting out by the usual action of $SO(3)$ on the spheres of constant $(t-r)$ and $r$, then the ``geometric'' foliation is not invariant under this action. Hence the red lines cannot really represent the geometric foliation - they can be taken either as an indication of the type of behaviour of this foliation, or alternatively, as representing a \emph{single angle} on the sphere.
	\\	
	Subfigure ii) shows the same foliations, but this time we have performed a conformal compactification with respect to the dynamic metric $g$, so that this time the hypersurfaces of constant $u$ are drawn as straight lines at forty five degrees. In other words, this is a true Penrose diagram for the manifold $(\mathcal{M}, g)$. Note that the hypersurfaces of constant $(t-r)$ (i.e.\ the red lines) can be \emph{timelike} in some regions (though not in the region $r \leq r_0$), which is one reason why this foliation is unsuitable. The same comments apply as above: this time, the blue lines can only really be taken as representing a \emph{single angle} on the leaves of the background foliation.
	\\
	Subfigure iii) shows another possible scenario. We have again compactified with respect to the dynamic metric $g$. However, this time, the leaves of constant $t-r$ all approach the \emph{same} value of $u$. We expect to see energy decay through the \emph{geometric} foliation, but, in this case, it is clear that we cannot expect the energy to decay through the background foliation.
	}
	\label{figure penrose diagrams 1}
\end{figure}

On the other hand, the dynamic metric $g$ can be used to define \emph{null geodesics}. In other words, we can consider null geodesics with respect to the metric $g$ that appears in the system of wave equations, and which can depend on the solution to that system. To be more precise: we can define the function $u$ by the requirement that $u = t - r_0$ on the surface $r = r_0$, and that $u$ is constant along outgoing null geodesics\footnote{It is fairly easy to see that such a $u$ satisfies the \emph{eikonal} equation $g^{-1}(\upd u , \upd u) = 0$.} originating from the surface $r = r_0$. Again, we can join the level sets of $u$ in the exterior region $r \geq r_0$ to the level sets of $t$ in the interior region $r \leq r_0$ to define a foliation. We refer to this\footnote{Perhaps we should say that this foliation is only truly ``geometric'' (i.e.\ defined relative to the metric $g$) in the region $r \geq r_0$. However, as long as $g$ and $m$ are suitably close, the leaves of this foliation will be spacelike (with respect to $g$) in the region $r \leq r_0$, which is all that we require.} as the ``geometric null foliation''. See figure \ref{figure penrose diagrams 1} for diagrams of these different foliations. Note that the foliations agree in the region $r \leq r_0$, but can differ in the region $r \geq r_0$.

It is with respect to the \emph{geometric} foliation that we can expect energy decay. This should be obvious from the fact that the wave equations are not defined with reference to the ``background'' structure (i.e.\ the Minkowski metric $m$ or the coordinates $(t,x^1,x^2,x^3,x^4)$) but, instead, they make use of the ``geometric'' structure provided by the metric $g$. To put this another way: we can expect the majority of the energy of the waves to flow along the characteristic curves for the wave equation, which are just the null geodesics of the metric $g$. Therefore, it appears that we need to use the geometric null foliation and not the background null foliation.

We shouldn't give up on the background foliation so easily, however. After all, things are a lot easier if we can work with a foliation that is defined independently of the solution we are seeking. Indeed, when proving the nonlinear stability of Minkowski space, Lindblad and Rodnianski \cite{Lindblad2004} used a ``background'' foliation, whereas Christodoulou and Klainerman \cite{Christodoulou1993} define their foliation in a fully ``geometric'' manner. This is one of main reasons that the proof given in \cite{Lindblad2004} is vastly simpler (and shorter!) than that given in \cite{Christodoulou1993}. Moreover, \cite{Yang2013a} managed to use a ``background'' null foliation, together with $r^p$-weighted energy estimates, to handle certain quasilinear wave equations.

However, the authors of both \cite{Lindblad2004} and \cite{Christodoulou1993} worked with foliations with \emph{uniformly spacelike} leaves, and, unlike the case of null leaves, we can expect a foliation that is spacelike with respect to $m$ to be spacelike with respect to $g$ as well. Furthermore, the equations considered in \cite{Yang2013} involve quasilinear terms that are particularly well behaved. Indeed, from the point of view of energy estimates, these quasilinear equations behave like wave equations with the \emph{classical} null condition, with a few minor additional complications. Indeed, \cite{Yang2013} could have used a foliation by hyperboloidal leaves, which are spacelike and only \emph{asymptotically} null (as used in \cite{Moschidis2015}), which would have avoided the issues associated with timelike portions of the leaves.

Despite the potential advantages, it is in fact \emph{not} possible to use the background null foliation in the case we consider. Nor is it possible to adapt this foliation in a simple way in order to produce a suitable foliation\footnote{In the special case of the Einstein equations in wave coordinates it is actually possible to use the ADM energy to adapt the background foliation into a suitable foliation. This can be done in a way that depends on the initial data, but does not require a full solution to be constructed in the future of the initial data surface. However, the fact that this can be done relies heavily on the additional structures present in the Einstein equations, which are \emph{not present} in a general system of the type we study.}. The reason for this can be easily seen from figure \ref{figure penrose diagrams 1}. The issue is that the leaves of the background null foliation and the leaves of the geometric null foliation can diverge from one another - it turns out that they can diverge at a rate $\sim r^\epsilon$, for initial data of size $\epsilon$. Hence, as shown in figure \ref{figure penrose diagrams 1} iii), it might be the case that the energy of a wave decays with respect to the geometric foliation but does \emph{not} decay with respect to the background foliation.

It is worth comparing this with the case studied in \cite{Yang2013}. There, the quasilinear terms behave better, and it is possible to show that the leaves of a geometric foliation can only move away from the background null leaves by some finite amount. Ultimately, this is responsible for the fact that the background foliation is ``good enough'' to apply the $r^p$-weighted energy method: we can expect energy decay with respect to the background foliation, since this closely follows the geometric foliation. Alternatively, we can compare our case with the special case of the Einstein equations in harmonic coordinates. Here, the leaves of the background foliation and the leaves of the geometric foliation diverge, but they only diverge logarithmically, and (using the special structures present in the Einstein equations) it is possible to correct for this by modifying the background leaves in a manner that depends only on the initial ADM mass of the spacetime. No such modification is available in the general case that we consider.

We therefore use a geometric null foliation in this work. We note that, although this increases the complexity of our treatment, we also gain something from this approach. Namely, this method gives some insight into the geometric origin of various observed phenomena, such as the slow decay towards null infinity, which can arise due to ``shock formation at infinity'' (see section \ref{subsection intro shock formation}). Also, since we are making full use of the geometric structure in the equations, we can expect our estimates to be slightly sharper than they would otherwise be.

We also note here that, in common with other approaches, an important part of our analysis involves foliating the leaves themselves (i.e.\ the level sets of $u$) by ``spheres''. In \cite{Christodoulou1993} these spheres are the intersection of the outgoing null foliation with surfaces of constant ``time'', where these constant time surfaces are themselves chosen to be maximal. In \cite{Speck2016a}, where a ``background'' structure is present, the spheres are the intersection of the surfaces of constant $u$ with the level sets of the ``background'' time coorinate $t$. Finally, in numerous applications (for example \cite{Luk2018}), an ingoing null foliation is constructed, and the ``spheres'' are the intersections of the ingoing and outgoing null leaves. In our case, since we are interested in using the $r^p$-weighted energy method, it is natural to define the spheres as the intersection of the surfaces of constant $u$ and the surfaces of constant $r$. In this sense, our coordinates are similar to the classic ``Bondi coordinates'' near null infinity \cite{Bondi1960}.

\subsection{The need for geometric commutators}
\label{subsection need for geometric commutators}

The necessity of using a geometric foliation immediately leads to the requirement that we use \emph{geometric commutators}. To see why, we can focus on the \emph{angular momentum} vector fields, i.e.\ the vector fields that generate rotations. With respect to the ``background'' geometry, these vector fields can be written as
\begin{equation*}
\Omega_{(ij)} := x^i \partial_j - x^j \partial_i
\end{equation*}
where $i$, $j = 1,2,3$. These vector fields are tangent to the spheres $r = \text{constant}$, $t = \text{constant}$. Vector fields of this kind play a central role in obtaining pointwise bounds, since they can be used in conjunction with the Sobolev inequality on the sphere to provide pointwise bounds with decay in $r$.

\begin{figure}[htpb]
	\centering
	\includegraphics[width = \linewidth, keepaspectratio]{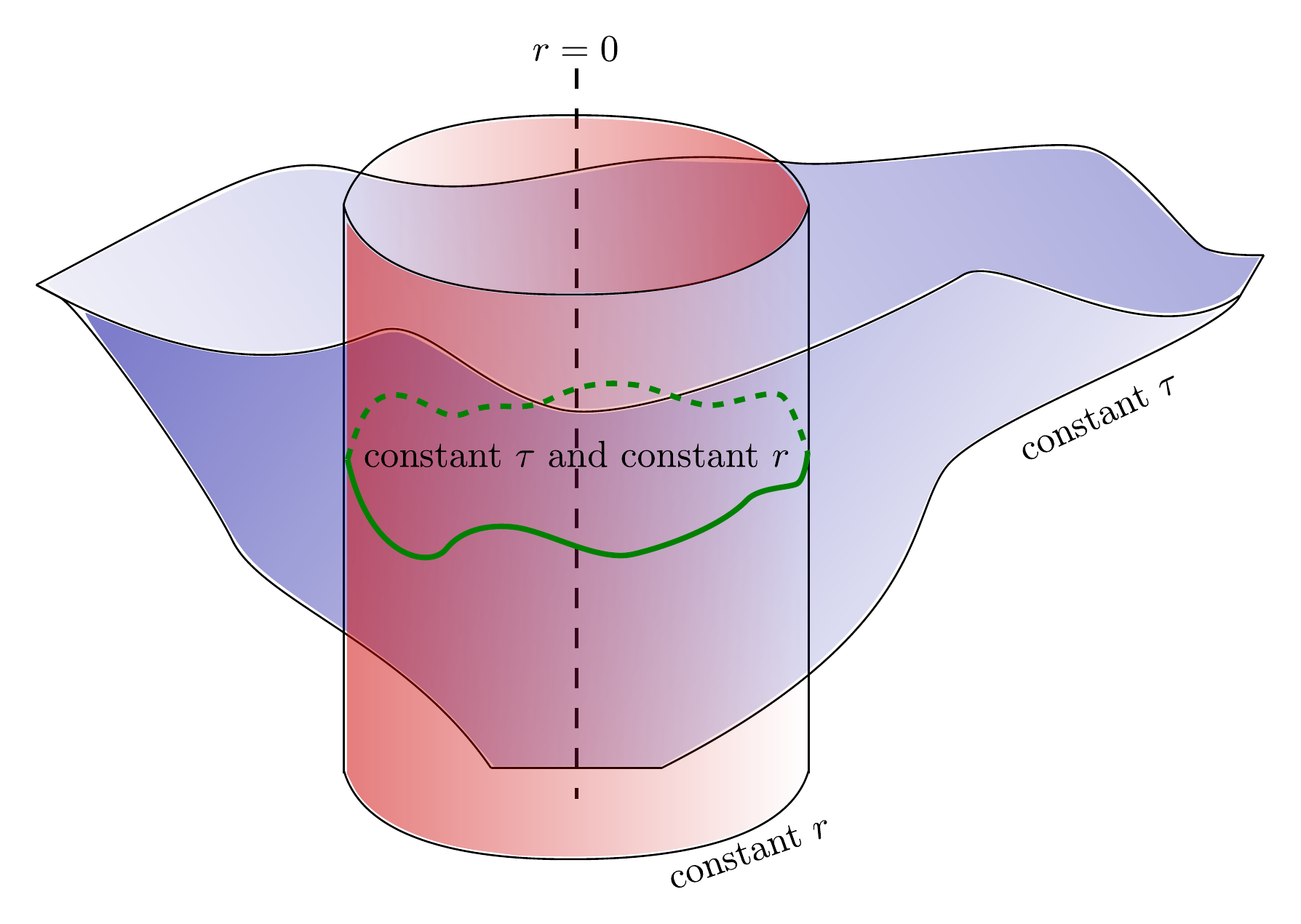}
	\caption{
	A more detailed diagram of our foliation, suppressing only one dimension. We have drawn this from the point of view of the ``background geometry'', i.e.\ the coordinate functions $x^a$ would define a square grid on this diagram.
	\\
	The line $r = 0$ is drawn as a dashed black line in the centre of the diagram. The red cylinder represents a surface of constant $r$, while the blue ``cone'' represents a surface of constant $\tau$. Note that, since the curves of constant $\tau$ are defined ``geometrically'' (in the region $r \geq r_0$), this cone is deformed relative to the standard cone. The blue ``cone'' and the red cylinder intersect along the green curve, which represents a ``sphere'' of constant $\tau$ and $r$. If we had not suppressed one dimension, then the green deformed circle would be replaced by a deformed sphere.
	\\
	The ``background'' angular momentum operators act by rigidly rotating the (red) cylinder of constant $r$ around the axis $r = 0$. It should be clear, from this diagram, that the vector fields which generate these rotations are \emph{not} tangent to the (green) ``sphere''. In other words, this ``sphere'' is not invariant under the action of the background angular momentum operators.
	}
	\label{figure detailed foliation}
\end{figure}

The leaves of our geometric foliation, however, are not themselves foliated by spheres of constant $r$ and $t$. In other words, the vector fields $\Omega_{(ij)}$ are not tangent to the ``spheres'' of constant $u$ and $r$, if $u$ is defined geometrically as above (see figure \ref{figure detailed foliation}). Said yet another way, the leaves of the geometric foliation are not ``spherically symmetric'' from the point of view of the background geometry. If we want to use the Sobolev inequality on these ``spheres'', then we must construct appropriate vector fields which are tangent to them, we cannot simply use the ``background'' angular momentum operators.

The approach taken in \cite{Christodoulou2007e} and subsequently \cite{Speck2016a} is to project the angular momentum operators onto the spheres, using a projection operator which is constructed from the metric $g$. Although we will not take precisely the same approach (see subsection \ref{subsection intro commuting with rnabla slashed} for the details), our methods are similar in spirit. The key point is that the commutation operators are constructed in such a way that they depend on the metric $g$, and hence they depend on the solution to the wave equation. This, in turn, is responsible for a lot of the technical difficulties we must overcome: for example, we are encounter error terms in the energy estimates that depend on the connection coefficients of the null frame, but unlike the Christoffel symbols associated with a ``background'' coordinate system, these connection coefficients cannot be computed directly from the metric perturbations. Instead, we can derive various transport equations for these connection coefficients along outgoing null geodesics, which we then use to estimate the connection coefficients. A much more serious problem is encountered when we try to perform the energy estimates at ``top order'' (that is, after commuting the maximal number of times with the commutation operators), since the estimates we can obtain using the transport equations appear to ``lose derivatives''. In other words, when performing the energy estimates for $\mathscr{Y}^n \phi$, where $\mathscr{Y}$ is any commutation operator, we seem to encounter terms that depend on $\mathscr{Y}^{n+1}\phi$. To overcome this difficult we rely on a number of special features of the null frame connection coefficients, which we detail in subsection \ref{subsection intro recovering sharp decay} below.

\section{Slow decay towards null infinity due to the weak null condition}
\label{section intro slow decay}

We now come to one of the most important differences between systems that only obey the weak null condition, and systems with the classical null condition: the solutions to the former can decay \emph{at a slower rate towards null infinity} compared to solutions of the latter. In fact, if $\phi$ solves a wave equation with the classical null condition, then it is possible to obtain the bound
\begin{equation*}
|\partial \phi| \lesssim (1+r)^{-1}
\end{equation*}
where this bound is also \emph{uniform}\footnote{In fact, it is possible to obtain decay in $u$ as well as in $r$, but to keep this discussion as simple as possible we will ignore the $u$ behaviour for now.} in the parameter $u$. Note that this rate is the same as the sharp rate for linear wave equations. On the other hand, we find that, if $\phi$ solves a wave equation with the hierarchical weak null condition, then we can only obtain a bound of them form
\begin{equation*}
|\partial \phi| \lesssim (1+r)^{-1 + \epsilon}
\end{equation*}
This factor of $r^\epsilon$ makes a fundamental difference in the analysis, as we shall see below. For now, we will concentrate on the origin of this slow decay, and its relationship to the weak null condition.

There are actually two possible sources of this behaviour. First, if we examine the asymptotic system corresponding to the semilinear system
\begin{equation*}
\begin{split}
\Box \phi_1 &= 0 \\
\Box \phi_2 &= (\partial_t \phi_1)(\partial_t \phi_2) \\
\end{split}
\end{equation*}
then we find that it has solutions that behave like
\begin{equation*}
r\Lbar \phi_2 \sim r^\epsilon
\end{equation*}
So, in some situations, we can view the \emph{semilinear} nonlinear terms as giving rise to slow decay towards null infinity. In other words, if we treat the semilinear terms as a forcing term, and substitute for the expected behaviour, then the slow decay of the forcing term towards null infinity naturally leads to slower decay for the solutions towards null infinity. This contrasts with the case of semilinear terms obeying the classical null condition. For example, if $\phi$ obeys the linear wave equation, then we have the behaviour
\begin{equation*}
\begin{split}
|(\partial_t \phi)^2| &\lesssim r^{-2} 
\\
|(\partial_t \phi)^2 - (\partial_r \phi)^2| &\sim r^{-3}
\end{split}
\end{equation*}
So the latter term, which obeys the classical null condition, decays much more rapidly towards null infinity than the former term, which does not obey the classical null condition.

More interestingly, we can also observe this kind of behaviour in \emph{quasilinear} systems, even in cases where the semilinear terms behave well. For example, we can consider the system
\begin{equation*}
\begin{split}
\tilde{\Box}_g \phi &= 0 \\
g_{ab} &= m_{ab} + h_{ab}(\phi)
\end{split}
\end{equation*}
which was previously treated in \cite{Lindblad2008}. This equation effectively has no semilinear terms - all the difficulty is caused by the quasilinear term. However, even in this case, we find that we cannot obtain the rate $r^{-1}$ for \emph{higher} derivatives, i.e.\ we cannot improve on the estimate $\partial\mathscr{Y}^n\phi \sim r^{-1 + \epsilon}$, where $\mathscr{Y}$ stands for any commutation operator and\footnote{Actually, a small improvement can be made, giving only a \emph{logarithmic} loss of decay, but we still do not obtain the rate $r^{-1}$.} $n \geq 1$. The reason for this is that such systems can exhibit ``shock formation at infinity'', which causes their solutions to behave differently from those of the linear wave equation. We explore this issue in more depth in the next section.

Another interesting consequence of this slow decay towards null infinity is that we need to develop estimates which work in the case of ``error terms'' that decay slowly towards null infinity. These error terms can either be generated by nonlinearities (as we have mostly considered above), or alternatively they can be put in ``by hand'', and treated as external ``forcing terms''. For example, we can consider the wave equations
\begin{equation*}
\tilde{\Box}_g \phi_{(A)} = F_{(A)}
\end{equation*}
where the $F_{(A)}$ are some known functions on the manifold. Since our estimates must work with slowly decaying error terms, we are able to allow the $F_{(A)}$ to decay slowly in $r$ too - for example, we can have $F \sim r^{-2-\epsilon}$ if $\phi$ is at the bottom level of the hierarchy, and $F \sim r^{-2+\epsilon}$ if $\phi$ is higher in the hierarchy. Perhaps more interestingly, we can allow the rectangular components of the metric to be given by
\begin{equation*}
g_{ab} = m_{ab} + h^{(0)}_{ab} + h^{(1)}_{ab}(\phi)
\end{equation*}
where $h^{(1)}_{ab}(\phi)$ are some terms depending on the solution to the wave equations (as considered above), but where $h^{(0)}_{ab}$ are some known functions. Then the $h^{(0)}_{ab}$ must be chosen such that the associated null connection components decay suitably in $r$, but this decay can be very slow: for example, we can allow $\tr_{\slashed{g}} \chibar \sim r^{-1+\epsilon}$. Moreover, at least in a region of bounded $r$, we do not require the functions $h^{(0)}_{ab}$ to decay in $t$ at all! This is potentially very important in applications where the ``end state'' is not expected to be trivial, flat space.

\subsection{``Shock formation'' at null infinity and the weak null condition}
\label{subsection intro shock formation}

\begin{figure}[htbp]
	\centering
	\includegraphics[width = \linewidth, keepaspectratio]{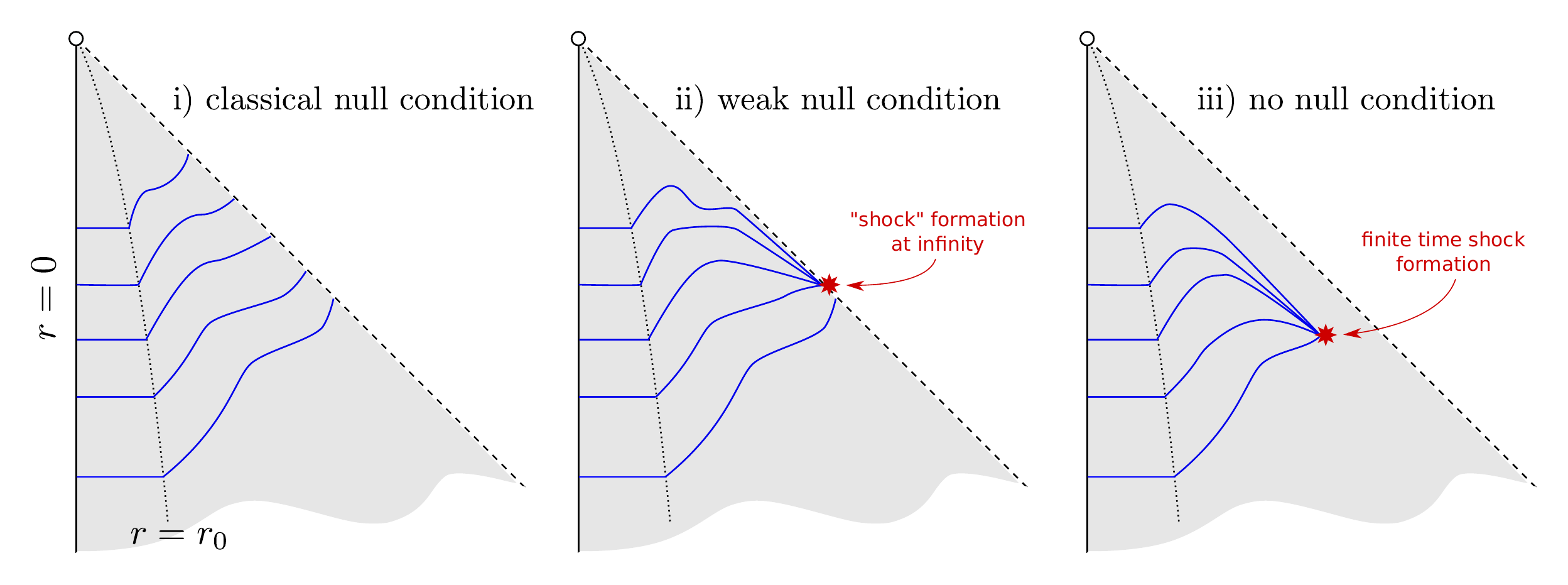}
	\caption{
	This figure shows several conformal compactifications with respect to the \emph{background} geometry, together with (in blue) a geometric foliation at a \emph{single} fixed angle on the sphere (see the comments in the caption of figure \ref{figure penrose diagrams 1}).
	\\
	Subfigure i) shows the kind of foliation which arises if the wave equation actually satisfies the \emph{classical} null condition but is nevertheless quasilinear - see \cite{Yang2013} for examples of this kind of equation. Note that no shocks can form.
	\\
	Subfigure ii) shows a possible behaviour for systems of the kind we are studying, i.e.\ equations satisfying the \emph{weak} null condition. Note that, after some time, the characteristics can converge, signalling that a shock is forming. However, this can only happen \emph{at null infinity}. In other words, the shock cannot form at any finite value of $r$. This is a crucial fact which allows us to prove the global existence result.
	\\
	Finally, subfigure iii) shows the possible behaviour of a system which does not satisfy the weak null condition, and for which shocks form at a \emph{finite} value of $r$. Note that it is not possible to extend the solution classically past the point where the shock forms, since certain derivatives of the solution (actually, the ``bad'' derivatives) blow up there. See \cite{Christodoulou2007e} or \cite{Speck2016a} for more details.
	}
	\label{figure penrose diagrams 2}
\end{figure}

We have already mentioned several times that the leaves of the geometric foliation can diverge from their ``background'' counterparts polynomially. To say this more precisely: we can introduce the ``inverse foliation density'' $\mu$, which measures the separation of the leaves of the foliation relative to the background geometry. When this quantity becomes zero, then the leaves of the foliation come together and a shock forms. It turns out that, in the case of a wave equation satisfying the weak null condition, we can obtain the bounds $(1+r)^{-C\epsilon} \lesssim \mu \lesssim (1+r)^{C\epsilon}$, which means that $\mu$ remains bounded away from zero (and infinity), so shocks cannot form. However, it is possible for the inverse foliation density to \emph{approach} zero at null infinity (see figure \ref{figure penrose diagrams 2}). This is the phenomena that we refer to as ``shock formation at infinity'', and it is responsible for the difference in the decay rates of certain quantities (such as $\Lbar \phi$) in the quasilinear case.

In fact, $\mu$ satisfies a transport equation along the outgoing null geodesics which has the form
\begin{equation*}
\partial_r \log \mu \sim (\Lbar h)_{LL}
\end{equation*}
Here, $(\Lbar h)_{LL}$ is one of the metric components, and so it can be expressed algebraically in terms of the derivatives of the fields $\phi$. It is therefore essential that we can obtain the bound
\begin{equation*}
|(\Lbar h)_{LL}| \lesssim \epsilon (1+r)^{-1}
\end{equation*}
and our definition of the weak null hierarchy ensures this, since the field $h_{LL}$ is required to be at the bottom level of the hierarchy. Substituting this bound into the evolution equation for $\mu$, we recover the bounds claimed above for the inverse foliation density.

It is interesting to compare this to the case where the equations satisfy the classical null condition, and to the case where the weak null condition is not satisfied. In the former case, $\mu$ satisfies a transport equation of the form
\begin{equation*}
\partial_r \log \mu \sim \bar{\partial} h
\end{equation*}
i.e.\ only ``good derivatives'' appear on the right hand side. Since these quantities decay at a rate $\sim \epsilon r^{-2}$, these terms are integrable in $r$ and we find that $\mu$ and its inverse are bounded. In other words, the leaves of the geometric foliation cannot separate or come together in any significant manner (see figure \ref{figure penrose diagrams 2}, subfigure i)).

On the other hand, suppose that we are considering an equation that does \emph{not} satisfy the weak null condition. Then, if we define the rescaled vector field $\check{\Lbar} := \mu \Lbar$, we find that $\mu$ satisfies a transport equation of the form
\begin{equation*}
\partial_r \mu \sim (\check{\Lbar} h)_{LL}
\end{equation*}
Even when the equation does not satisfy the weak null condition, it can still be possible (see \cite{Christodoulou2007e, Speck2016a}) to obtain the bound
\begin{equation*}
|(\check{\Lbar} h)_{LL}| \lesssim (1+r)^{-1}
\end{equation*}
Note two important difference: the transport equation is for the quantity $\mu$ rather than $\log \mu$, and the vector field $\check{\Lbar}$ has taken the place of the vector field $\Lbar$. Note that this is not sufficient to conclude global existence and positivity for $\mu$: instead, we will typically have
\begin{equation*}
|\mu - \mu_0| \sim C\epsilon \log \left( \frac{1 + r}{1+r_0} \right)
\end{equation*}
where $\mu = \mu_0$ when $r = r_0$. We can then expect $\mu$ to vanish, and a shock to form, when $r$ is approximately
\begin{equation*}
r \sim \exp\left(-\frac{1}{C\epsilon}\right)
\end{equation*}
(see figure \ref{figure penrose diagrams 2}, subfigure iii)) and indeed there are many examples of equations (and initial data sets) where this does indeed happen (\cite{Christodoulou2007e, Speck2016a, Holzegel2016}), the most famous of which being the Euler equations in fluid dynamics.

\begin{figure}[htbp]
	\centering
	\includegraphics[width = \linewidth, keepaspectratio]{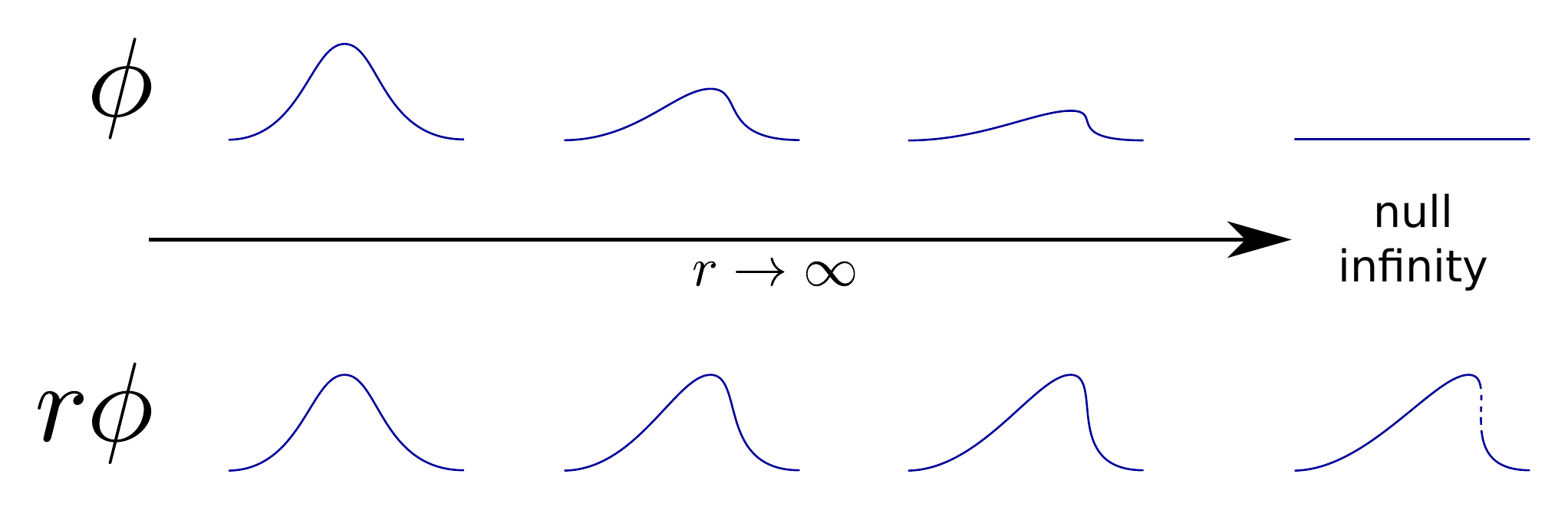}
	\caption{
	A cartoon showing ``shock formation at infinity''. In the upper part of the diagram, we follow a wave packet along an outgoing null geodesic as $r \rightarrow \infty$. Dispersion causes the wave (and its derivatives) to decay, so that in the limit the function $\phi$ is trivial.
	\\
	In the lower part of the diagram, we rescale the wave packet by a factor of $r$, so that we are now following the ``radiation field'' $r\phi$ rather than the field $\phi$. This rescaling counterbalances the dispersion, so the radiation field does not decay as we head towards null infinity. Instead, as $r \rightarrow \infty$ a \emph{shock} can develop in the radiation field, meaning that certain derivatives of the radiation field do not have a finite limit as $r \rightarrow \infty$. In fact, it is the ``bad'' derivative $\Lbar(r\phi)$ which may not have a finite limit.
	}
	\label{figure shock formation at infinity}
\end{figure}

From this point of view, we can view equations with the weak null condition as lying exactly on the borderline between \emph{equations with the classical null condition} and \emph{shock-forming} equations. Shocks cannot form, but they can \emph{almost} form, and, as we head towards null infinity, they get closer and closer to forming. At the same time, dispersion causes the amplitude of the waves to tend to zero as we move out towards null infinity. The practical upshot of this is that certain derivatives of the waves do not behave in the same was as their linear counterparts. In particular, if we rescale the solution by multiplying by a factor of $r$ to compensate for the usual dispersion, then we find that the $\Lbar$ derivative of the solution is not continuous in the limit $r \rightarrow \infty$ (see figure \ref{figure shock formation at infinity}). Another way to say this is that the solutions $\phi$ do not asymptotically approach solutions to the linear wave equation.

\subsection{The need for a sharp Morawetz-type estimate}
\label{subsection intro sharp Morawetz}

As mentioned above, semilinear terms which do not obey the classical null condition decay more slowly towards null infinity. This causes some serious difficulties when attempting to perform the energy estimates. It turns out that the quasilinear case also leads to similar error terms, provided that the ``reduced'' wave operator takes the place of the standard wave operator. To illustrate these difficulties, we will first outline an approach\footnote{This approach is, essentially, the one followed in \cite{Yang2013}.} which works for wave equations satisfying the \emph{classical} null condition Then, we will explain why this approach fails for wave equations that only satisfy the \emph{weak} null condition: essentially, it is because we do not have a \emph{sharp} Morawetz estimate. Then we will explain how this difficulty can be overcome using a more standard vector field method (for example, as it is done in \cite{Lindblad2004}). We are led to the necessity of developing similar techniques for use in conjunction with the $r^p$-weighted energy method, which we will outline in the next subsection.

Suppose that we want to control solutions to a semilinear equation of the form
\begin{equation*}
\Box \phi = (\partial \phi)(\bar{\partial} \phi)
\end{equation*}
where $\bar{\partial}$ stands for one of the ``good'' derivatives (that is, $L$ or one of the angular derivatives) and $\partial$ can be any derivative. This equation obeys the classical null condition. By using the standard energy estimate (for example, multiplying by $\partial_t \phi$ and integrating by parts) we can obtain the estimate
\begin{equation*}
\begin{split}
\int_{\{t = t_1\} \cap \{r \leq r_0\}} |\partial \phi|^2
+ \int_{\{u = t_1 - r_0\} \cap \{r \geq r_0\}} |\bar{\partial} \phi|^2
&\leq
\int_{\{t = t_0\} \cap \{r \leq r_0\}} |\partial \phi|^2
+ \int_{\{u = t_0 - r_0\} \cap \{r \geq r_0\}} |\bar{\partial} \phi|^2
\\
&\phantom{\leq}
+ \int_{\tau = t_0}^{t_1} \left( \int_{\{t = \tau\} \cap \{r \leq r_0\}} |\bar{\partial} \phi||\partial \phi|^2 \right) \upd \tau
\\
&\phantom{\leq}
+ \int_{\tau = t_0}^{t_1} \left( \int_{u = \{\tau - r_0\} \cap \{r \geq r_0\}} |\bar{\partial} \phi||\partial \phi|^2 \right) \upd \tau
\end{split}
\end{equation*}
where the standard volume forms are left implicit. The last two terms are the ``error terms'' - they would be absent if we were looking at the linear wave equation.

Now, we assume that we can obtain the same kinds of pointwise bounds for these quantities as in the linear case. In particular, we can assume that we have the bound
\begin{equation*}
|\bar{\partial}\phi| \lesssim \epsilon (1+r)^{-2}
\end{equation*}
Then, we must control the \emph{spacetime} integral of $\epsilon (1+r)^{-2}|\partial \phi|^2$. Fortunately, if we make the assumption that the energy is bounded by some constant $\mathcal{E}$, then we can also prove the \emph{Morawetz} (or \emph{integrated local energy decay}) estimate, which has the form
\begin{equation}
\label{equation introduction Morawetz}
\int_{\tau = t_0}^{t_1} \left( \int_{\{t = \tau\} \cap \{r \leq r_0\}} (1+r)^{-1-\delta}|\partial \phi|^2 \right) \upd \tau
+ \int_{\tau = t_0}^{t_1} \left( \int_{\{u = \tau - r_0\} \cap \{r \geq r_0\}} (1+r)^{-1-\delta}|\partial \phi|^2 \right) \upd \tau
\lesssim \delta^{-1} \mathcal{E}
\end{equation}
where $\delta > 0$ can be chosen to be any positive constant.

Combining these two estimates, we can verify that the energy is indeed bounded: we have
\begin{equation*}
\int_{\{t = t_1\} \cap \{r \leq r_0\}} |\partial \phi|^2
+ \int_{\{u = t_1 - r_0\} \cap \{r \geq r_0\}} |\bar{\partial} \phi|^2
\leq
\int_{\{t = t_0\} \cap \{r \leq r_0\}} |\partial \phi|^2
+ \int_{\{u = t_0 - r_0\} \cap \{r \geq r_0\}} |\bar{\partial} \phi|^2
+ \epsilon \delta^{-1} \mathcal{E}
\end{equation*}
The first two terms on the right hand side are the initial energy, and the last term can be made very small by taking $\epsilon \ll \delta$. This allows us to \emph{improve} the energy bound, (i.e.\ to obtain a stronger bound than $\mathcal{E}$ on the energy) which allows us to close the argument.

Let us see what happens if we try to repeat this process in the case where the nonlinearity only satisfies the weak null condition. We will again consider the model system
\begin{equation*}
\begin{split}
\Box \phi_1 &= 0 \\
\Box \phi_2 &= (\partial \phi_1)(\partial \phi_2) \\
\end{split}
\end{equation*}
Of course, everything proceeds as usual for the field $\phi_1$. However, when we try to perform the energy estimate for $\phi_2$ we find that we must bound the terms
\begin{equation*}
\int_{\tau = t_0}^{t_1} \left( \int_{\{t = \tau\} \cap \{r \leq r_0\}} |\partial \phi_1| |\partial \phi_2|^2 \right) \upd \tau
+ \int_{\tau = t_0}^{t_1} \left( \int_{\{u = \tau - r_0\} \cap \{r \geq r_0\}} |\partial \phi_1| |\partial \phi_2|^2 \right) \upd \tau
\end{equation*}
If we substitute the bound $|\partial \phi_1| \lesssim \epsilon (1+r)^{-1}$ then we find that we must control terms of the form
\begin{equation*}
\int_{\tau = t_0}^{t_1} \left( \int_{\{t = \tau\} \cap \{r \leq r_0\}} \epsilon (1+r)^{-1} |\partial \phi_2|^2 \right) \upd \tau
+ \int_{\tau = t_0}^{t_1} \left( \int_{\{u = \tau - r_0\} \cap \{r \geq r_0\}} \epsilon (1+r)^{-1} |\partial \phi_2|^2 \right) \upd \tau
\end{equation*}
but in this case \emph{we cannot use the Morawetz estimate}! In particular, we cannot take $\delta \rightarrow 0$ in equation \eqref{equation introduction Morawetz}, since the right hand side includes the factor $\delta^{-1}$. In fact, it is fairly easy to construct counterexamples to the ``$\delta = 0$ Morawetz estimate''\footnote{By this we mean a statement of the form \eqref{equation introduction Morawetz}, with $\delta$ set to $0$ on the left hand side and set to $1$ on the right hand side.}, for example, by using geometric optics or Gaussian beams \cite{Sbierski2013a}.

It is instructive to see how this issue is resolved in the formalism of the ``classical'' vector field method. Here, we use a foliation by surfaces of constant $t$. Also, the decay estimates involving decay in $r$ are replaced by decay estimates with decay in $t$. We find ourselves with the inequality
\begin{equation*}
\int_{t = t_1} |\partial \phi_2|^2
\leq \int_{t = t_0} |\partial \phi_2|^2
+ \int_{\tau = t_0}^{t_1} \left( \int_{\{t = \tau\}}|\partial \phi_1| |\partial \phi_2|^2 \right) \upd \tau
\end{equation*}
Now, after substituting the bound $|\partial \phi_1| \lesssim \epsilon (1+t)^{-1}$ we can use the Gronwall inequality to deduce a bound of the form
\begin{equation*}
\int_{\{t = t_1\}} |\partial \phi_2|^2
\lesssim (1+t_1)^{C\epsilon} \int_{\{t = t_0\}} |\partial \phi_2|^2
\end{equation*}
In other words, the energy of the field $\phi_2$ grows slowly in time\footnote{It is this slow growth of the energy that will eventually be responsible for the worse pointwise behaviour of $\phi_2$. But note that we have not had to assume that $\phi_2 \sim t^{-1}$!}. To say this another way, as $t \rightarrow \infty$ we gradually \emph{lose control} over the energy of the solution. It turns out, however, that we still retain sufficient control to close all of the required estimates.

Evidently we need a similar approach, but adapted to the $r^p$-weighted method. We develop this in the next section.

\subsection{The degenerate energy at null infinity, and a ``sharp'' Morawetz-type estimate}
\label{subsection intro degenerate energy}

We have emphasised above that the solutions to equations which lack the classical null condition can behave differently from their linear counterparts \emph{near null infinity}. There is a conceptual difference, then, between the approach we take - where the leaves of the foliation reach null infinity - and an approach using surfaces of constant $t$, which only approach null infinity in the limit $t \rightarrow \infty$. As we have seen, we can only hope to obtain energy estimates which \emph{degenerate} towards null infinity. Since the leaves of our foliation reach null infinity, this means that we need to work with an energy that degenerates towards null infinity.

Inspired by the above calculations, we could try including a factor of $(1+t)^{-C\epsilon}$ in the energy estimates. However, sine we are using the $r^p$ method, we will mostly working with $u, r$ coordinates, and so it is more natural\footnote{We could, of course, try a weight of the form $(1+r+u)^{-C\epsilon}$, but it turns out that, when we construct $u$ geometrically, certain derivatives of $u$ behave differently than might be expected and consequently this kind of weight does not work well.} to include a weight of the form $(1+r)^{-C\epsilon}$. If we include this decaying weight, then it generates additional terms when we integrate by parts in the energy estimates, and our energy estimates will ``degenerate'' as $r \rightarrow \infty$. It turns out that the some of the additional terms have a ``good sign'', meaning that they contribute spacetime integrals that can be placed on the left hand side of the inequality. Meanwhile, other terms have a ``bad sign'' and must be included on the right hand side of the inequality. It turns out, however, that the terms with the ``good sign'' involve the \emph{bad} derivatives, while the terms with the ``bad sign'' involve only the \emph{good} derivatives. In other words, if we continue to study the model system \eqref{equation weak null example 2} and we try this approach, then we can obtain an estimate of the form
\begin{equation*}
\begin{split}
&\int_{\{t = t_1\} \cap \{r \leq r_0\}} (1+r)^{-C\epsilon} |\partial \phi_2|^2
+ \int_{\{u = t_1 - r_0\} \cap \{r \geq r_0\}} (1+r)^{-C\epsilon}|\bar{\partial} \phi_2|^2
\\
&
+ \int_{\tau = t_0}^{t_1} \left( \int_{\{t = \tau\} \cap \{r \leq r_0\}} C\epsilon (1+r)^{-1-C\epsilon}|\partial \phi_2|^2 \right) \upd \tau
\\
&
+ \int_{\tau = t_0}^{t_1} \left( \int_{\{u = \tau - r_0\} \cap \{r \geq r_0\}} C\epsilon (1+r)^{-1-C\epsilon}|\partial \phi_2|^2 \right) \upd \tau
\\
&\lesssim
\int_{\{t = t_0\} \cap \{r \leq r_0\}} (1+r)^{-C\epsilon}|\partial \phi|^2
+ \int_{\{u = t_0 - r_0\} \cap \{r \geq r_0\}} (1+r)^{-C\epsilon}|\bar{\partial} \phi|^2
\\
&\phantom{\lesssim}
+ \int_{\tau = t_0}^{t_1} \left( \int_{\{t = \tau\} \cap \{r \leq r_0\}} \left( C\epsilon (1+r)^{-1-C\epsilon}|\bar{\partial} \phi_2|^2 + (1+r)^{-C\epsilon}|\partial \phi_1| |\partial \phi_2|^2\right) \right) \upd \tau
\\
&\phantom{\lesssim}
+ \int_{\tau = t_0}^{t_1} \left( \int_{\{u = \tau - r_0\} \cap \{r \geq r_0\}} \left( C\epsilon (1+r)^{-1-C\epsilon}|\bar{\partial} \phi_2|^2 + (1+r)^{-C\epsilon}|\partial \phi_1| |\partial \phi_2|^2 \right) \right) \upd \tau
\end{split}
\end{equation*}
Note that the term which we now have control over (i.e.\ the spacetime integral of $C\epsilon(1+r)^{-1-C\epsilon}|\partial \phi_2|^2$) could also be controlled by choosing $\delta = C\epsilon$ in the Morawetz estimate, so it may appear that we have not gained anything. However - and this is the crucial point - the error terms arising from the nonlinearity are \emph{also} multiplied by the factor $(1+r)^{-C\epsilon}$. This means that, if we substitute $|\partial \phi_1| \lesssim \epsilon (1+r)^{-1}$, and if $C$ is sufficiently large, then the error terms arising from the nonlinearity can now be \emph{absorbed} by the left hand side. This is in marked contrast to the previous approach we described, where the standard energy estimate is combined with the Morawetz estimate.

We still need to control the new error terms on the right hand side, which now only involve the \emph{good} derivatives of $\phi_2$. Note that we can use the standard Morawetz estimate to control these terms in the region $r \leq r_0$. There are several options for dealing with these terms, but since we are using the $r^p$-weighted energy estimates, we can use these to control this term. In fact, if we use the $r^p$ weighted estimates with the choice $p = \delta$, then we find that we can bound a term of the form
\begin{equation*}
\int_{\tau = t_0}^{t_1} \left( \int_{\{u = \tau - r_0\} \cap \{r \geq r_0\}} \left( \delta (1+r)^{-1+\delta}|\bar{\partial} \phi_2|^2 \right) \right) \upd \tau
\end{equation*}
so, for any choice $\delta > 0$, this estimate provides sufficient control. Of course, when performing this estimate there are other error terms which we need to control, but we will postpone our discussion of these terms until the main body of the paper.

In a sense, we have achieved our goal of establishing a ``sharp'' Morawetz estimate. Although we only control the spacetime integral of $(1+r)^{-1-C\epsilon}|\partial \phi|^2$, which could already be done using the Morawetz estimate, all the other terms now appear with a weight of $(1+r)^{-C\epsilon}$, so relative to this weight, we have achieved the ``sharp'' weight $r^{-1}$. In summary, if we allow the implicit constant to depend on $r_0$, and if we set
\begin{equation*}
w := (1+r)^{-C\epsilon}
\end{equation*} 
then we can obtain the bound
\begin{equation}
\label{equation introduction energy estimate}
\begin{split}
&\int_{\{t = t_1\} \cap \{r \leq r_0\}} |\partial \phi|^2
+ \int_{\{u = t_1 - r_0\} \cap \{r \geq r_0\}} \left( w|\bar{\partial} \phi|^2 + r^{-2+\delta}|L(r\phi)|^2 \right)
+ \int_{\tau = t_0}^{t_1} \left( \int_{\{t = \tau\} \cap \{r \leq r_0\}} C\epsilon |\partial \phi|^2 \right) \upd \tau
\\
&
+ \int_{\tau = t_0}^{t_1} \left( \int_{\{u = \tau - r_0\} \cap \{r \geq r_0\}} \left( C\epsilon w (1+r)^{-1}|\partial \phi|^2 + \delta(1+r)^{-1+\delta} |\bar{\partial} \phi|^2 \right) \right) \upd \tau
\\
&\lesssim
\int_{\{t = t_0\} \cap \{r \leq r_0\}} |\partial \phi|^2
+ \int_{\{u = t_0 - r_0\} \cap \{r \geq r_0\}} \left( w|\bar{\partial} \phi|^2 + r^{-2+\delta}|L(r\phi)|^2 \right)
\end{split}
\end{equation}

There are a couple of points to note in the above inequality. The first is that, although we have obtained the ``sharp'' Morawetz estimate, we have had to add the $p$-weighted flux term to the initial energy. This means that, in order to obtain this kind of inequality, we need some additional information about the initial data: in particular, the $L$ derivatives of the initial data need to decay appropriately\footnote{Another consequence of the fact that we need to use the $p$-weighted estimate even to show \emph{boundedness} is that we obtain slightly worse decay towards timelike infinity, i.e.\ slightly worse decay in $u$.}. The second point to be made about equation \eqref{equation introduction energy estimate} is that it concerns the \emph{weighted} energy rather than the standard energy. In other words, we must include the weight $w$ in the energy. This has some important consequences for the later analysis: if we try to obtain pointwise bounds by the usual combination of energy estimates and Sobolev embeddings, then the weight $w$ will lead to slightly worse pointwise decay in $r$. This means that, following this strategy, there is no hope of obtaining the sharp pointwise bounds $(\partial \phi) \sim r^{-1}$, even though we \emph{need} to establish this bound for fields at the bottom level of the hierarchy! The solution to this problem will be outlined in subsection \ref{subsection intro recovering sharp decay}, but note that a similar problem was also encountered when using the ``classical'' vector field method, since here the energy estimate include a factor that grows as $t^{C\epsilon}$, but pointwise bounds of the form $(\partial \phi) \sim t^{-1}$ are required.

\subsection{Slow decay towards timelike infinity due to the weak null condition, upper bounds on the value of \texorpdfstring{p}{p} in the \texorpdfstring{$r^p$}{rp}-weighted energy and the non-existence of a radiation field}

The idea of the $r^p$-weighted energy method is to show that the \emph{energy} decays through a suitable foliation. The rate of decay that can be obtained corresponds to the \emph{maximum} value of $p$ that can be taken in the eponymous energy estimates. For example, in their original paper \cite{Dafermos2010b}, Dafermos and Rodnianski took values of $p$ up to (and including) $p = 2$. This leads to decay rates for the energy of $u^{-2}$, which, in turn, leads to the pointwise bounds
\begin{equation*}
|\phi| \lesssim \begin{cases}
	(1+r)^{-1}(1+u)^{-\frac{1}{2}}
	\\
	(1+r)^{-\frac{1}{2}}(1+u)^{-1}
\end{cases}
\end{equation*}
and related bound for the derivatives of $\phi$. 

In their papers \cite{Angelopoulos2018a} (and \cite{Angelopoulos2018a}), Angelopoulos et al.\  were able to obtain \emph{improved} decay (i.e.\ faster rates of decay in $u$) for higher spherical harmonics in wave equations on black hole spacetimes, using a modified version of the $r^p$-weighted method. This used the fact that the value of $p$ can effectively be taken larger than $2$ for higher spherical harmonics. On the other hand, \cite{Yang2013} was able to show the global existence for solutions to wave equations with the classical null condition using only $p = 1+\delta$ for $\delta > 0$. This, however, leads to the pointwise bounds
\begin{equation*}
|\phi| 
\lesssim  \begin{cases}
	(1+r)^{-1}(1+u)^{-\frac{1}{2}\delta}
	\\
	(1+r)^{-\frac{1}{2}}(1+u)^{-\frac{1}{2} - \frac{1}{2}\delta}
\end{cases}
\end{equation*}
which has significantly worse decay in $u$. Nevertheless, this was found to be sufficient to close all of required bounds.

For most nonlinear equations, the limiting factor for the maximum value of $p$ is the decay in $r$ of the nonlinear terms. Certainly, when the equations satisfy the weak null condition but not the classical null condition, then the nonlinear terms decay slowly towards null infinity and this limits the maximum value of $p$. In fact, we find that the maximum value of $p$ that we can take is only $p = 1 - C\epsilon$. Since we also work with the \emph{degenerate} energy rather than the usual energy, this leads to the bounds
\begin{equation*}
|\phi| 
\lesssim  \begin{cases}
	(1+r)^{-1 + C\epsilon}(1+u)^{\frac{1}{2}C\epsilon}
	\\
	(1+r)^{-\frac{1}{2} + C\epsilon}(1+u)^{-\frac{1}{2} + \frac{1}{2}C\epsilon}
\end{cases}
\end{equation*}
Note that the first inequality actually \emph{grows} in $u$. Importantly, the sum of the exponents is greater than $-1$ in this case, whereas in the other cases it is less than $-1$. The fact that this sum is smaller than $-1$ played a role in the argument of \cite{Yang2013}, so we find that we must modify the argument to take this into account. Nevertheless, we find that we can adapt the argument to work even in the case of this very weak decay in $u$. Indeed, it is a general feature that, while we must be very careful with decay rates in $r$, we have some ``room'' with regards to decay in $u$, and so (in the latter case) we can make do with sub-optimal decay.

A few words should be said at this point about the ``classical'' vector field method of Klainerman. As a direct consequence of commuting with a large set of vector fields, this method is able to achieve the rate
\begin{equation*}
|\phi| \lesssim (1+t)^{-1 + C\epsilon}(1+u)^{-\frac{1}{2}}
\end{equation*}
\emph{even in the case that} $\phi$ \emph{solves an equation without the classical null condition}. Note the improvement in the $u$ decay over that stated previously. The fact that this method achieves improved decay in $u$ reflects a certain advantage in this method over the $r^p$-weighted method, but this comes with a corresponding cost. First, various quantities (for example, the metric components) must be assumed \emph{a priori} to decay in $t$. Decay in $t$, of course, implies decay in $u$. On the other hand, the $r^p$-weighted method can be made to work even if such quantities \emph{never}\footnote{Throughout this work, we have endeavoured to assume the least possible decay in $u$ for various quantities (especially the metric components and the connection components) with a view to future applications. Given the discussion in this section, it may be surprising that we need to assume \emph{any} decay at all in $u$ for these quantities, but for technical reasons certain combinations of connection coefficients are required to decay in $u$. Note, however, that this decay can be very slow - $u^{-\delta}$ is enough - and also that, in the special case of the Einstein equations, we can use the extra structure in the equations to avoid this assumption (see subsection \ref{subsection Einstein equations in wave coords}).} decay in $t$ but are only bounded (see \cite{Yang2013a}), while decaying suitably in $r$. Since the $r^p$-weighted energy method allows for worse behaviour in $u$ for the equations, it is natural that it also obtains weaker estimates in this sense\footnote{Note, however, that we think of the decay towards null infinity (that is, decay in $r$ or $t$) as the \emph{critical} decay, which can cause problems for closing the estimates. The decay in $u$ is typically of secondary importance.}. Second, the fact that the classical vector field method requires commuting with a large set of vector fields can cause problems in a number of situations, a point that we have made several times already.

\section{Other technical issues}
\label{section intro other issues}

There are numerous additional technical issues which we have had to overcome in this work. Most of these issues have already been encountered at various points in the literature, and our strategy for overcoming them follows, to a significant extent, the approaches taken by previous authors. However, certain concessions have to be made to our specific problem, and the approaches must be adapted slightly in order to be made to work. We describe some of these technical issues in the following subsections.

\subsection{Recovering sharp decay rates near null infinity}
\label{subsection intro recovering sharp decay}

Several times above, we have stressed the need to recover the sharp bounds
\begin{equation*}
|\partial \phi| \lesssim \epsilon (1+r)^{-1}
\end{equation*}
at least for the fields at the bottom level of the hierarchy (and especially for the field $h_{LL}$). But it is clear from our discussion above - specifically, our discussion of the degenerate energy - that we will not be able to recover this sharp rate from the usual combination of energy estimates and Sobolev inequalities. In fact, for fields at the $n$-th level in the hierarchy, we begin with the assumption that
\begin{equation*}
|\partial \phi| \lesssim \epsilon (1+r)^{-1+C_{(n)}\epsilon}
\end{equation*}
but, for the same reason, we cannot recover this bound from the combination of energy estimates and Sobolev inequalities, because the degenerate energy always causes worse decay in $r$ than the rate assumed initially. Hence, we actually have the same problem at all levels of the hierarchy.

We solve this problem by using the same kind of argument that was used in, for example, in \cite{Lindblad2008} and \cite{Lindblad2004}\footnote{Note that, in \cite{Lindblad2004}, a different argument was used to recover many of these bounds. In particular, we can recover the improved bound $|\Lbar h|_{LL} \lesssim (1+r)^{-2}$ using additional structure in the Einstein equations - see subsection \ref{subsection Einstein equations in wave coords}.}. This argument consists of integrating the equations of motion along outgoing null geodesics. In other words, we return to the asymptotic system, but this time we reintroduce the terms which were ``dropped'' to form the asymptotic system. However, we treat these ``reintroduced'' terms (all of which, by definition, involve \emph{good} derivatives) as error terms, using the bounds that we have already obtained on the good derivatives.

For example, consider the model equation
\begin{equation*}
\begin{split}
\Box \phi_1 &= 0 \\
\Box \phi_2 &= (\partial_t \phi_1)^2
\end{split}
\end{equation*}
Then we can write the second equation in the null frame as
\begin{equation*}
-L(r\Lbar \phi_2) = -r\slashed{\Delta}\phi_2 - L\phi_2 + r\frac{1}{4}(L\phi_1 + \Lbar \phi_1)^2
\end{equation*}
Assuming that we can show the bound\footnote{In fact, this bound follows easily from using the energy estimates and the Sobolev embedding for $\phi_2$ after commuting sufficiently many times with angular derivatives.} $\slashed{\Delta}\phi_2 \sim \epsilon r^{-3 + C\epsilon}$, and substituting the other pointwise bounds that we have already been able to show using the Sobolev embedding and energy estimates, we find
\begin{equation*}
L(r\Lbar \phi_2) \lesssim \epsilon (1+r)^{-2 + C\epsilon} + \epsilon^2 (1+r)^{-1}
\end{equation*}
So, integrating along the integral curves of $L$, we find that we can obtain the bound
\begin{equation*}
|\Lbar \phi_2| \lesssim \epsilon (1+r)^{-1} \log(2+r)
\end{equation*}
Note that this improves on the rate of decay in $r$ compared to what can be obtained using Sobolev embedding and energy estimates alone.

Similar ideas allow us to improve the rates of decay of other fields. In particular, if $\phi$ is at the bottom level of the hierarchy, then the terms on the right hand side of the transport equation for $r\Lbar \phi$ turn out to be \emph{integrable in} $r$, which means that we can recover the sharp bounds $\partial\phi \sim (1+r)^{-1}$ for these fields. This approach also works for quasilinear equations. For example, in the null frame, the quasilinear equation
\begin{equation*}
\tilde{\Box}_{g(\phi)}\phi = 0
\end{equation*}
can be written as
\begin{equation*}
-L(r\Lbar \phi) = -r\slashed{\Delta} \phi - L\phi + \Gamma(\phi) \cdot \partial \phi
\end{equation*}
where the $\Gamma(\phi)$ are related to the null-frame connection coefficients. In particular, using the ``reduced'' wave operator $\tilde{\Box}_g$ rather than the standard geometric wave operator $\Box_g$ results in these error terms being integrable in $r$. This can be viewed as the motivation for using the reduced wave operator in the definition of the weak null hierarchy.

\subsection{Avoiding a loss of regularity at top-order}

In subsection \ref{subsection need for geometric commutators} we saw that, in our setting, we will need to use commutator vector fields which \emph{depend on the solution}. A consequence of this approach is that some of our estimates appear to \emph{lose derivatives}: when we are estimating the energy of $\mathscr{Y}^n \phi$, where the $\mathscr{Y}$'s are commutation operators, we appear to need knowledge of the energy of $\mathscr{Y}^{n+1} \phi$. This, of course, leads to problems when trying to perform the ``top-order'' energy estimates.

This kind of phenomena was encountered (and successfully dealt with) in the monumental work of Christodoulou and Klainerman \cite{Christodoulou1993}. Before we provide details of the solution to this problem, there are several relevant differences between the approach of \cite{Christodoulou1993} and our present approach, which we will need to sketch first.

The foliation used in \cite{Christodoulou1993} was, in a sense, \emph{entirely geometric}. That is, they constructed an outgoing null foliation, similar to our foliation by outgoing null leaves of constant $u$. At the same time, \cite{Christodoulou1993} also constructed a \emph{geometric} foliation by spacelike surfaces (surfaces of ``constant time''). Note that \cite{Christodoulou1993} made no use at all of any kind of ``background'' structure, so it was necessary to construct the surfaces of constant time solely using the geometry associated with the dynamic metric $g$. The ``spheres of constant time and $u$'', which are the intersections of these two foliations, form an important part of the definition of the null frame, and their geometry is therefore intimately connected to the properties of the null frame connection coefficients.

In contrast, we are making use of the ``background'' geometry, so we can use this to help define our ``spheres''. We have already discussed (see subsection \ref{subsection need for geometric foliation}) the fact that we cannot use a background foliation to define the outgoing ``null'' hypersurfaces, however, we see no need to define the \emph{other} foliation (and thereby the ``spheres'') in a geometric manner. In other words, we could use the ``background'' time coordinate $t$ to define the ``surfaces of constant time'', and then define the ``spheres'' as the intersections of the leaves of this foliation with the leaves of the foliation by outgoing null surfaces (see \cite{Speck2016a} for an example of this approach). However, since we are using the $r^p$-weighted energy method, it is more natural to use the background \emph{radial} coordinate $r = \sqrt{(x^1)^2 + (x^2)^2 + (x^3)^2}$, and then to define the spheres as the intersection of the (geometric) null surfaces with the (background) constant $r$ surfaces.

This choice has some consequences regarding the aforementioned loss of derivatives at top order. If we had chosen to use an entirely ``background geometry'' based approach, then the connection coefficients could be written directly in terms of the derivatives of the metric components, and so there would not be any issues regarding a loss of derivatives. On the other hand, in the ``fully geometric'' approach of \cite{Christodoulou1993}, \emph{all} of the connection coefficients must be controlled by integrating various PDEs satisfied by these quantities. Therefore, a great deal of care must be taken to avoid losing derivatives in this process, and indeed this occupies a large amount of the proof of \cite{Christodoulou1993}.

Our approach (in common with that of \cite{Speck2016a}) occupies a middle ground between these two strategies. Since we use a ``semi-geometric'' approach - using the geometric outgoing null leaves together with the background $r = \text{constant}$ hypersurfaces - it turns out that \emph{some} of the null-frame connection coefficients can be written directly in terms of the derivatives of the (background components of the) metric, while for others we must rely on integrating the associated PDEs. Therefore, we only need to confront the loss-of-regularity issue in a few special cases.

We will not discuss all of the connection coefficients here, but rather, we will take a consider a couple of the most difficult quantities. First, consider the quantity $\tr_{\slashed{g}} \chi_{(\text{small})}$, where $\chi$ is the extrinsic curvature of the spheres $r = \text{constant}$ considered as submanifolds of the hypersurfaces of constant $u$, $\chi_{(\text{small})}$ is $\chi$ minus the Minkowski value, and we take the trace with respect to the metric induced on the sphere by the ambient metric $g$. This quantity satisfies a transport equation along the outgoing null geodesics of the form
\begin{equation*}
 L \left( r^2\tr_{\slashed{g}} \chi_{(\text{small})} \right) = r^2(\slashed{\Delta} h)_{LL} + \ldots
\end{equation*}
where the ellipses stand for terms that are either lower order or that can be written as \emph{exact} $L$ derivatives. If we integrate this equation along the integral curves of $L$, then we find that, in order to control $\tr_{\slashed{g}}\chi$ we need to control \emph{two} derivatives of $h$. This is worrying, since we expect the connection coefficients to behave like \emph{first} derivatives of the metric, and indeed this will eventually lead to a loss of derivatives. Note that integrating this equation leads to the pointwise bound
\begin{equation*}
\tr_{\slashed{g}} \chi_{(\text{small})} \lesssim \epsilon (1+r)^{-2 + C\epsilon}
\end{equation*}

The way around the problem of a loss of derivatives was found in \cite{Christodoulou1993}, and it consists of using the wave equation to rewrite $(\slashed{\Delta} h)_{LL}$ as a perfect $L$ derivative, plus some lower order terms. Specifically, we obtain an equation of the form
\begin{equation*}
L \left( r^2\tr_{\slashed{g}} \chi_{(\text{small})} \right) = r^2(\Box_g h)_{LL} + L\left(r^2(\Lbar h)_{LL} \right) - r(\Lbar h)_{LL} + \ldots
\end{equation*}
Now, $(\Box_g h)_{LL}$ is also lower order by assumption, and the second term on the right hand side can be moved onto the left hand side. In this way we can estimate $\tr_{\slashed{g}} \chi_{(\text{small})}$ in terms of the \emph{first} derivatives of the metric. Note, however, that this approach loses some decay: if we used this equation to obtain a pointwise bound, then we would only be able to show
\begin{equation*}
\tr_{\slashed{g}} \chi_{(\text{small})} \lesssim \epsilon (1+r)^{-1}
\end{equation*}

Another example of a quantity which appears in our estimates and apparently leads to a loss of derivatives is $\slashed{\nabla}^2 \log \mu$, that is, the second angular derivatives of the logarithm of the inverse foliation density. This quantity appears in our energy estimates after we apply one commutation operator, so we should hope that it can be controlled in terms of the \emph{second} derivatives of the metric. However, the transport equation for $\log \mu$ takes the form
\begin{equation*}
L\log \mu \sim (\partial h)
\end{equation*}
so, if we commute this equation twice with angular derivatives, in order to control $\slashed{\Delta} \log \mu$, then we find that we need to control \emph{three} derivatives of the metric. Again, a na\"ive approach here leads to a loss of derivatives. For this reason, \cite{Christodoulou1993} introduced the ``smoothing foliation''. Instead of following this scheme, we take the approach of \cite{Speck2016a}, which is to derive an equation linking the \emph{spherical Laplacian} of $\log \mu$ to the \emph{time derivatives} of $\tr_{\slashed{g}} \chi_{(\text{small})}$. Then, using elliptic estimates, we can control all of the second angular derivatives of $\log \mu$ given control over the spherical Laplacian of $\log \mu$. As above, we find that although this avoids a loss of derivatives, it also leads to worse pointwise bounds, which we discuss in the subsequent subsection.

\subsection{Avoiding decay loss at top-order}

As well as a loss of derivatives, we also encounter an apparent \emph{loss of decay} at the top order in our energy estimates, which, at first sight, appears to ruin our changes of closing the estimates. The origin of this phenomenon is perhaps easiest to understand if we imagine that we are using energy estimates (and commutation operators) based on the \emph{background} geometry, rather than the geometric foliation and geometrically constructed commutation operators. However, the same problem arises in both cases.

First, imagine that we are using an entirely \emph{background}-based approach, Suppose that we commute a maximum number of $n$ times. It is fairly easy to see that, when performing the (degenerate) energy-boundedness estimate for the field $\mathscr{Y}^{m} \phi$, we will encounter an error term of the form
\begin{equation*}
\begin{split}
\int_{\tau = t_0}^{t_1} \left( \int_{\{t = \tau\} \cap \{r \leq r_0\}} \left( w(\mathscr{Y}H^{ab}) \D_a \D_b \mathscr{Y}^{n-1} \phi \right) (\D \mathscr{Y}^n \phi) \right) \upd \tau
\\
+ \int_{\tau = t_0 - r_0}^{t_1 - r_0} \left( \int_{\{u=\tau\} \cap \{r > r_0\}} \left( w (\mathscr{Y}H^{ab}) \D_a \D_b \mathscr{Y}^{n-1} \phi \right) (\D \mathscr{Y}^n \phi) \right) \upd \tau
\end{split}
\end{equation*}
where $(g^{-1})^{ab} = (m^{-1})^{ab} + H^{ab}$. Since $g_{ab} = m_{ab} + h_{ab}(\phi)$, it is possible to express the fields $H^{ab}$ in terms of the fields $\phi$.

Now, it is possible to express all the second derivatives $\D_a \D_b$ in terms of the wave operator, derivatives of time derivatives, derivatives of angular derivatives or lower order terms. For the purposes of this argument, we ignore the terms proportional to the wave operator, as well as the lower order terms (these are much easier to control). This means that, effectively, we can replace $\D_a \D_b \mathscr{Y}^{n-1} \phi$ with $\D \mathscr{Y}^{n} \phi$.

The problem comes when we consider the term $(\mathscr{Y}H)^{ab}$. It should be clear that we at least need to show
\begin{equation*}
\left| (\mathscr{Y}H)^{ab} \right| \lesssim \epsilon (1+r)^{-1}
\end{equation*}
Although we can just about obtain this inequality for the field $\phi$, we \emph{cannot} obtain this kind of pointwise bound after commuting (i.e.\ we cannot obtain $|\mathscr{Y}\phi| \lesssim \epsilon(1+r)^{-1}$), but we find that we necessarily \emph{lose} some decay in $r$ after commuting.

To deal with this kind of problem, \cite{Lindblad2004} used the wave coordinate condition for the Einstein equations (see subsection \ref{subsection Einstein equations in wave coords}). However, in the general case we are considering, this condition cannot be expected to hold. On the other hand, Lindblad (in \cite{Lindblad2008}) gave an alternative approach to dealing with this problem. The idea is \emph{not} to replace $\D_a \D_b \mathscr{Y}^{n-1} \phi$ with $\D \mathscr{Y}^{n} \phi$, but instead to commute \emph{first} with the operators $\D_a$ and only \emph{afterwards} to commute with the operators $\mathscr{Y}$. Suppose that we have already obtained the required bounds for $\mathscr{Y}^{n-1} \phi$. Commuting one more time with $\D_a$, error term appearing in the energy estimate for the field $\D \mathscr{Y}^{n-1} \phi$ is
\begin{equation*}
\begin{split}
\int_{\tau = t_0}^{t_1} \left( \int_{t=\tau \cap \{r \leq r_0\}} \left( w (\partial H^{ab}) \D_a \D_b \mathscr{Y}^{n-1} \phi \right) (\D \mathscr{Y}^n \phi) \right) \upd \tau
\\
+ \int_{\tau = t_0 - r_0}^{t_1 - r_0} \left( \int_{u=\tau \cap \{r > r_0\}} \left( w (\partial \mathscr{Y}H^{ab}) \D_a \D_b \mathscr{Y}^{n-1} \phi \right) (\D \mathscr{Y}^n \phi) \right) \upd \tau
\end{split}
\end{equation*}
In other words, the problematic term $\mathscr{Y} H$ is replaced by the term $\partial H$. Now, we can also write
\begin{equation*}
|(\partial H^{ab}) \D_a \D_b \mathscr{Y}^{n-1} \phi|
\sim
|(\partial H)_{LL} (\D^2 \mathscr{Y}^{n-1} \phi)|
+ \text{good derivatives}
+ \text{lower order terms}
\end{equation*}
The terms involving good derivatives and the lower order terms are comparatively easy to control. It turns out that it \emph{is} possible to obtain the sharp $r^{-1}$ bound for the quantity $|(\partial H)_{LL}|$. Hence we \emph{can} control the energy of the field $\D \mathscr{Y}^{n-1} \phi$. In other words, we can bound the spacetime integral
\begin{equation*}
\int_{\tau = t_0 - r_0}^{t_1 - r_0} \left( \int_{u=\tau \cap \{r > r_0\}} w \epsilon (1+r)^{-1} |\D^2 \mathscr{Y}^{n-1} \phi|^2 \right) \upd \tau
\end{equation*}

Now, we can return to the energy estimate for field $\mathscr{Y}^{n} \phi$. Let us also choose a different weight factor $\tilde{w}$ for this estimate. Now, we can write the error term as
\begin{equation*}
\begin{split}
&\int_{\tau = t_0 - r_0}^{t_1 - r_0} \left( \int_{u=\tau \cap \{r > r_0\}} \left( \tilde{w} (\mathscr{Y}H^{ab}) \D_a \D_b \mathscr{Y}^{n-1} \phi \right) (\D \mathscr{Y}^n \phi) \right) \upd \tau
\\
&\lesssim
\int_{\tau = t_0 - r_0}^{t_1 - r_0} \left( \int_{u=\tau \cap \{r > r_0\}} \left( \tilde{w} \epsilon (1+r)^{-1+C\epsilon} |\D^2 \mathscr{Y}^{n-1} \phi| |\D \mathscr{Y}^n \phi| \right) \right) \upd \tau
\\
&\lesssim
\int_{\tau = t_0 - r_0}^{t_1 - r_0} \left( \int_{u=\tau \cap \{r > r_0\}} \left( \tilde{w} \epsilon (1+r)^{-1+2C\epsilon} |\D^2 \mathscr{Y}^{n-1} \phi|^2 + \tilde{w} \epsilon (1+r)^{-1} |\D \mathscr{Y}^n \phi|^2 \right) \right) \upd \tau
\end{split}
\end{equation*}
The second term is precisely the kind of thing we wanted. On the other hand, if we choose
\begin{equation*}
\tilde{w} = (1+r)^{-2C\epsilon}w
\end{equation*}
then the \emph{first} term has already been controlled!

The discussion above mirrors that in \cite{Lindblad2008}. However, this presupposes a ``background geometry'' based approach. On the other hand, we will make use of both a \emph{geometric} foliation and \emph{geometric} commutator vector fields, so we do not encounter precisely the same error term, nor do we use precisely the same language. Nevertheless, we do encounter an analogue of this difficulty. Specifically, we encounter the error term
\begin{equation*}
\int_{\tau = t_0 - r_0}^{t_1 - r_0} \left( \int_{u=\tau \cap \{r > r_0\}} \left( w |r\slashed{\nabla}\phi| | \slashed{\nabla}^2 \mathscr{Y}^{n-1} \log \mu| |\D \mathscr{Y}^n \phi| \right) \right) \upd \tau
\end{equation*}
where now $u$ is to be understood as the \emph{geometrically} defined retarded time variable. As already mentioned, at \emph{top-order} we can use elliptic estimates to relate $(\slashed{\nabla}^2 \mathscr{Y}^{n-1} \log \mu)$ to $(\D_T \mathscr{Y}^{n-1} \tr_{\slashed{g}}\chi_{(\text{small})})$, where $T$ is the geometric analogue of the vector field $\partial_t$. Furthermore, we can relate $(\D_T \mathscr{Y}^{n-1} \tr_{\slashed{g}}\chi_{(\text{small})})$ to the quantity $(\D \D_T \mathscr{Y}^{n-1} h)_{LL}$, and hence to the quantity\footnote{In this discussion we are ignoring the fact that there might be multiple $\phi$'s, i.e.\ we might be dealing with a system of equations rather than just a single scalar wave equation. This is done entirely to simplify the discussion and the notation: the reader is free to re-insert labels for the different fields $\phi_{(a)}$.} $(\D \D_T \mathscr{Y}^{n-1} \phi)$. By following this sequence of estimates, we find that we must bound the term
\begin{equation*}
\int_{\tau = t_0 - r_0}^{t_1 - r_0} \left( \int_{u=\tau \cap \{r > r_0\}} \left( w |\mathscr{Y}\phi| | \D \D_T \mathscr{Y}^{n-1} \phi| |\D \mathscr{Y}^n \phi| \right) \right) \upd \tau
\end{equation*}
Note that we have written $r\slashed{\nabla} \phi \sim \mathscr{Y} \phi$ since we are commuting with the angular derivatives \emph{weighted by} $r$ (see subsection \ref{subsection intro commuting with rnabla slashed}).

Clearly we will run into exactly the same problem: the quantity $(\mathscr{Y} \phi)$ cannot be shown to decay at the sharp rate $\sim r^{-1}$ but only at a rate $\sim r^{-1+C\epsilon}$. If we were to write $\D_T = \mathscr{Y}$ and then try to estimate this error term, we would find it impossible to bound. The way around this problem is to use the insight of \cite{Lindblad2008}: we should attempt to \emph{first} bound the energy of $\D_T \mathscr{Y}^{n-1} \phi$, and only \emph{afterwards} bound the energy of $\mathscr{Y}^{n} \phi$. It turns out that the error term involving the second angular derivatives of $\log \mu$ does not appear at top-order in the energy estimate for $\D_T \mathscr{Y}^{n-1} \phi$. In fact, if we only commute once with the commutation operators, then the quantity $\slashed{\nabla}^2 \log \mu$ \emph{only} appears when we commute with angular derivatives. Hence, we can actually perform the energy estimate for the field $\D_T \mathscr{Y}^{n-1} \phi$ without running into this problem. Afterwards, as before, we can return to the energy estimate for $\mathscr{Y}^{n} \phi$ and use a slightly more degenerate weight, together with the estimate we have already shown for the energy of $\D_T \mathscr{Y}^{n-1} \phi$ to bound the error term.

One important point to note in the discussion above is that we \emph{only} need to commute with $\D_T$ first, and not all of the translation operators $\D_a$ as in \cite{Lindblad2008}. However, this approach could have been used in \cite{Lindblad2008} too - all second derivatives $\D_a \D_b \phi$ can be written in terms of \emph{first derivatives} of $\D_T \phi$, \emph{first derivatives} of \emph{angular derivatives of} $\phi$, $\Box_g \phi$, and lower order terms\footnote{Note, however, that \cite{Lindblad2008} uses the ``old'' vector field method to obtain decay, so a large number of vector fields is used for commutation. In particular, this already includes the translation vector fields $\partial_a$, so little would be gained by this modification.}. Commuting only with time translation, and avoiding commuting with the spatial translations, matches the spirit of the $r^p$-weighted energy method.

\subsection{The problem of multiple good derivatives}
\label{subsection intro LL problem}

There is another situation in which we seem to lack sufficient decay to close the estimates. However, unlike the cases above, this problem arises \emph{only} as a consequence of our use of the $r^p$-weighted energy estimate and has not yet been encountered in similar problems.

To explain this problem, we will need notation for several different sets of vector fields. Let us write $\mathscr{Z}$ to stand for either the vector field $\D_T$ or the (weighted) covariant angular derivative operators $r\slashed{\nabla}$. If the reader prefers, these latter can be replaced by the angular momentum operators $\Omega_{ij}$. We also write $Z$ to stand for any of the vector fields used in the ``classical'' vector field method: that is, $Z$ can be any of the translation operators $\D_a$, or the angular momentum operators $\Omega_{ij} = x^i \partial_j - x^j \partial_i$, or the boosts $\Omega_{0i} = t\partial_i + x^i \partial_t$, or the scaling operator $S = t\partial_t + x^i \partial_i$.

If we use the ``classical'' vector field method, then we have the bound
\begin{equation}
|\bar{\partial} \phi| \lesssim (1+t)^{-1} |Z \phi|
\end{equation}
for any scalar field $\phi$. In other words, we can exchange a \emph{good} derivative for a commutation operator, and if we do so then we \emph{gain} a decaying factor $\sim t^{-1}$. Similarly, if we encounter a term involving second derivatives, but all of those derivatives are \emph{good}, then we have the bound
\begin{equation}
|\bar{\partial}^2 \phi| \lesssim (1+t)^{-1} |\bar{\partial}Z \phi| + (1+t)^{-1} |\bar{\partial} \phi|
\end{equation}
so again, we can essentially swap one of the ``good'' derivative operators for a vector field $Z$ and a decaying factor $\sim t^{-1}$.

Now, consider what happens if we try to commute \emph{only} with the operators $\mathscr{Z} = \{\D_T , r\slashed{\nabla}\}$ in the spirit of the $r^p$-weighted energy method. Then, if the ``good'' derivative is an angular derivative, we can clearly write
\begin{equation*}
|\slashed{\nabla} \phi| \lesssim r^{-1} |\mathscr{Z} \phi|
\end{equation*}
which is analogous to the first estimate above. On the other hand, we cannot do this if the ``good'' derivative is an $L$ derivative. Nevertheless, it turns out that we can still obtain ``improved'' estimates for the $L$ derivatives (both for $L^2$ bounds and pointwise bounds), primarily through the use of the flux terms in the $r^p$-weighted energy estimates.

However, things are not so straightforward in the case of \emph{second} derivatives. If we encounter a term of the form $(\bar{\partial}^2 \phi)$, then \emph{if at least one of the derivatives is an angular derivative}, we can improve the estimate as above. For example, we have
\begin{equation}
|\bar{\partial} \slashed{\nabla} \phi| \lesssim r^{-1} |\bar{\partial} Z \phi| + r^{-1} |\bar{\partial} \phi|
\end{equation}

This leaves us with the case of \emph{two} $L$ \emph{derivatives}. Here, we cannot play the same game: the best we can do is to write
\begin{equation*}
|LL\phi|
\lesssim
|L T \phi|
+ |\Box_g \phi|
+ r^{-1} |\overline{\partial} \mathscr{Z} \phi|
+ r^{-1} |\partial \phi|
\end{equation*}
where we have used the fact that $T \sim L + \Lbar$. The only dangerous term is the first one. If we ignore the other terms, we have essentially replaced the second $L$ derivative with one of the commutation operators, but we have \emph{not} gained any additional decay in $r$ at all!

This causes problems in handling exactly one of the error terms that arises after commuting. Specifically, when commuting with the angular derivatives and then performing, say, the (weighted) energy boundedness estimate we encounter the error term
\begin{equation*}
\int_{\tau = t_0 - r_0}^{t_1 - r_0} \left( \int_{u=\tau \cap \{r > r_0\}} \left( w |r\slashed{\nabla}\log \mu| |LL \phi| |\D \mathscr{Z} \phi| \right) \right) \upd \tau
\end{equation*}
Similarly, suppose that we consider the error terms after having commuted $n$ times with the commutation operators. Then in the same estimate we encounter the two error terms
\begin{equation}
\label{equation introduction LL problems}
\int_{\tau = t_0 - r_0}^{t_1 - r_0} \left( \int_{u=\tau \cap \{r > r_0\}} w\left( 
	|\mathscr{Z}^n\log \mu| |LL \phi| |\D \mathscr{Z}^n \phi|
	+ |r\slashed{\nabla}\log \mu| |\D_L \D_L \mathscr{Z}^{n-1}\phi| |\D \mathscr{Z}^n \phi|
	\right) \right) \upd \tau
\end{equation}
We will discuss each of these error terms in turn. Each of them causes problems: the first due to a lack of \emph{pointwise} decay, and the second due to a lack of decay in $L^2$.

Consider the first error term in equation \eqref{equation introduction LL problems}. If we follow the logic above, and replace $|LL\phi|$ with $|LT\phi|$, then we must note that \emph{even in the best possible case} we can only obtain the decay rate
\begin{equation*}
|LT\phi| \lesssim \epsilon (1+r)^{-\frac{3}{2}}
\end{equation*}
This is because we are limited by the poor decay towards null infinity (see section \ref{section intro slow decay}) to the range $p \leq 1$ in the $p$-weighted energy estimates. Consequently, even if we imagine that we could actually achieve $p = 1$ (which we can't), we would only obtain the decay rate above. On the other hand, we can expect the behaviour $\mathscr{Z}^n \log \mu \sim (1+r)\slashed{\D} \mathscr{Z}^n \phi$, since we have to integrate in $r$ to obtain bounds for $\mathscr{Z}^n \log \mu$. Hence, we can expect this first error term to behave like
\begin{equation*}
\int_{\tau = t_0 - r_0}^{t_1 - r_0} \left( \int_{u=\tau \cap \{r > r_0\}} 
	w \epsilon (1+r)^{-\frac{1}{2}} |\D \mathscr{Z}^n \phi|^2
\right) \upd \tau
\end{equation*}
but the decay rate of $r^{-\frac{1}{2}}$ is well outside the range we can control. Indeed, even with the help of the ``sharp'' Morawetz estimate (see subsections \ref{subsection intro sharp Morawetz} and \ref{subsection intro degenerate energy}), we need the coefficient to decay at least as fast as $(1+r)^{-1}$. As we have stressed, this is the \emph{sharp} decay rate that we can obtain by this method.

Next, consider the second error term in equation \eqref{equation introduction LL problems}. The best pointwise decay rate that we can obtain for the angular derivatives of the foliation density is $|r\slashed{\nabla} \log \mu| \sim \epsilon r^{C\epsilon}$. Again, if we follow the method outlined above, we end up having to bound the term
\begin{equation*}
\int_{\tau = t_0 - r_0}^{t_1 - r_0} \left( \int_{u=\tau \cap \{r > r_0\}} w\left( 
\epsilon (1+r)^{C\epsilon} |\D_L \D_T \mathscr{Z}^{n-1}\phi| |\D \mathscr{Z}^n \phi|
\right) \right) \upd \tau
\end{equation*}
Here, we have to bound a spacetime integral with a coefficient that \emph{grows} (modulo the weight factor $w$) in $r$. If we try to do this using the Morawetz estimates then we are again out of luck: we need a coefficient that decays at least like $r^{-1}$, not one that grows like $r^{C\epsilon}$.

There are several potential approaches to dealing with these problems, all of which, however, are antithetical to the spirit of the $r^p$-weighted method to a greater or lesser extent. With regard to the first problem (i.e.\ the problem of the pointwise decay of $(LL\phi)$) there appears to be no option but to improve the pointwise decay of $(LL\phi)$. This can be done by commuting a single time, before commuting with anything else, with the vector field $rL$. Of course, in the spirit of the $r^p$-weighted energy method, we would prefer to commute with as few vector fields as possible, but commuting with $rL$ at least once appears inevitable\footnote{Note that this vector field does not generate a symmetry of Minkowski space. Nevertheless, the error terms that are produced by commuting with $rL$ can all be controlled, and some are even beneficial for the $p$-weighted estimates.}.

With regard to the second problem - that of the $L^2$ bound involving a factor that grows in $r$ - there are at least two approaches which can be taken. The first approach is to try to bound the integral in question using the \emph{flux} terms in the $p$-weighted energy estimate. The maximum value of $p$ that we can take in the $p$-weighted energy estimate is $p = 1 - C\epsilon$. If we apply this to the field $\slashed{\D}_T \mathscr{Z}^{n-1} \phi$ then it can be used to produce a bound of the form
\begin{equation*}
\int_{{u = u_1}\cap\{r \geq r_0\}} (1+r)^{-C\epsilon} |\slashed{\D}_L \slashed{\D}_T \mathscr{Z}^{n-1} \phi|^2 \lesssim \mathcal{E}_0
\end{equation*}
for \emph{any} $u_1$ such that the surface $\{u = u_1\} \cap \{r \geq r_0\}$ lies in the future of the initial data surface, and where $\mathcal{E}_0$ is some constant depending on the initial data. Hence, if we choose $w$ appropriately\footnote{Again, we must choose different weight factors for the case where the final commutation operator is a weighted angular derivatives, and for the case where the final commutation operator is $\D_T$, and then we must rely on the fact that these ``bad error terms'' do not appear when commuting with $\D_T$. In general, we can use a ``less degenerate'' energy when we apply the commutation operator $\D_T$ compared with the case when we apply a weighted angular derivative.}, then we can obtain the bound
\begin{equation*}
\begin{split}
&\int_{\tau = t_0 - r_0}^{t_1 - r_0} \left( \int_{\{u=\tau\} \cap \{r > r_0\}} w\left( 
	\epsilon (1+r)^{C\epsilon} |\D_L \D_T \mathscr{Z}^{n-1}\phi| |\D \mathscr{Z}^n \phi|
\right) \right) \upd \tau
\\
&\lesssim
\int_{\tau = t_0 - r_0}^{t_1 - r_0} \left( \int_{\{u=\tau\} \cap \{r > r_0\}} w\left( 
	\epsilon (1+r)^{1 + 2C\epsilon} |\D_L \D_T \mathscr{Z}^{n-1}\phi|^2
	+ \epsilon (1+r)^{-1} |\D \mathscr{Z}^n \phi|^2
\right) \right) \upd \tau
\\
&\lesssim
\epsilon \mathcal{E}_0 (t_1 - t_0)
+ \int_{\tau = t_0 - r_0}^{t_1 - r_0} \left( \int_{\{u=\tau\} \cap \{r > r_0\}} w\left( 
	\epsilon (1+r)^{-1} |\D \mathscr{Z}^n \phi|^2
\right) \right) \upd \tau
\end{split}
\end{equation*}
Although the second term is of the right form, the first term \emph{grows} in the parameter $t_1$, and at an unacceptable rate. This can be corrected \emph{if we presuppose some decay in} $u$ \emph{for the quantity} $(r\slashed{\nabla} \log \mu)$. Again, this is slightly opposed to the spirit of our work: since the $r^p$-weighted method can handle geometries which do not settle down to flat space, we would prefer not to assume pointwise bounds with decay in $\tau$.

An alternative method for controlling these error terms is to commute again with the vector field $rL$, after having commuted with the other vector fields. This approach is much simpler than the approach outlined above. If we write $\mathscr{Y}$ to stand for any of the commutation operators $\slashed{\D}_T$, $r\slashed{\nabla}$, or $r\slashed{\D}_L$, then we can immediately write our error term as\footnote{Here, $\overline{\D}$ stands for a covariant derivative in the direction of one of the \emph{good} derivatives. Actually, we want to commute with covariant derivatives which are \emph{projected} onto the spheres (see subsection \ref{subsection intro commuting with rnabla slashed}), but we will ignore this point for now.}
\begin{equation*}
\int_{\tau = t_0 - r_0}^{\tau_1 - r_0} \left( \int_{\{u = \tau\} \cap \{r > r_0\}} w|\slashed{\nabla}\log\mu| |\overline{\D} \mathscr{Y}^n \phi| |\D \mathscr{Y}^n \phi| \right) \upd \tau
\end{equation*}
if we drop some lower order terms which are easy to control. Now, we do not need the quantity $\slashed{\nabla} \log \mu$ to decay in $u$, and, by making use of the fact that one of the derivatives is a ``good'' derivative, then we can bound this integral. Indeed, we are effectively in the same situation as we would be if we had commuted with the ``full'' set of vector fields $Z$ that were discussed at the beginning of this subsection.

Note that, when we commute with the operator $r\D_L$, we must also control all of the error terms that this produces. It turns out that this requires commuting not just once with the operator $r\D_L$, but $n$ times, where $n$ is the maximum number of times that we commute with the other operators. In other words, we must treat the operator $r\D_L$ and the other commutation operators $\mathscr{Z}$ equally.

Both of the approaches outlined above have advantages and drawbacks. We choose the second approach, since it is slightly simper and since we avoid making the \emph{a priori} assumption of decay\footnote{It turns out that, for very technical reasons, we still need to assumed a very small amount of decay in $u$, at the rate $u^{-\delta}$, for particular combinations of the null frame connection components.} in $u$. Throughout this proof, we have attempted to make our methods as robust and adaptable as possible, with an eye to future applications. In many possible future applications, we cannot hope for the connection coefficients to decay in $u$, since, for example, the spacetime is expected to approach some unknown member of a family of solutions. Since we do not know which member of the family it will approach, we must allow for various quantities to lack decay in $u$. On the other hand, even in spacetimes that are far from ``symmetric''\footnote{By ``symmetric'' in this case, we mean that the error terms produced by commuting with $r\D_L$ are similar to those produced by commuting with $r\D_L$ in Minkowski space. In this case, these ``error terms'' do not vanish, even in the Minkowski space case.}, we can hope to commute with $r\D_L$ at least in the asymptotic region.

It is worth noting that the presence of this error term is closely linked to the behaviour of \emph{angular} derivatives of the foliation density. If we were to change\footnote{Strictly speaking, the ``spheres'' themselves would remain invariant, but the projection operator $\slashed{\Pi}$ changes under the kind of transformation we have in mind.} our ``spheres'' of constant $u$ and $r$, and thereby change the definition of the ``angular derivatives'' (since these are defined as tangent to the spheres), then a judicious choice might relieve us entirely of this problematic error term. In the special case of the Einstein equations, the extra structure present in the equations \emph{does} seem to allow for such a choice, which in turn relieves us of the burden of deciding between \emph{a priori} decay for connection coefficients, and commuting with extra vector fields, both of which are undesirable. However, in the general case we are studying here we have not been able to use such an approach successfully.

\subsection{Commuting with the covariant derivative operator on the spheres and the vector bundle \texorpdfstring{$\mathcal{B}$}{B}}
\label{subsection intro commuting with rnabla slashed}

We now move on to another aspect of our approach which differs slightly from many previous approaches to similar problems. We have already mentioned the importance of commuting with angular derivatives: in conjunction with the Sobolev embedding on the sphere, this is our primary way of obtaining pointwise bounds. Also in subsection \ref{subsection need for geometric commutators} we saw that we necessarily have to use \emph{geometrically} defined angular derivatives.

One option for these commutation operators is to use the background ``angular momentum operators'' $\Omega_{ij} = x^i \partial_j - x^j \partial_i$, and then to project these onto the (geometrically defined) spheres using a projection operator constructed from the metric $g$. This is the approach taken in, for example, \cite{Christodoulou2007e}. However, this strategy can run into some problems if we try to use it in situations which are not \emph{approximately} spherically symmetric.

We have emphasised our desire to use methods which are as robust and adaptable as possible, so that they can easily be applied to a variety of different situations. In keeping with this general approach, we will not commute with the (projected) angular momentum operators, but instead we will commute with the \emph{covariant derivatives on the spheres}. This approach has found success even when the space in question is not even approximately spherically symmetric, for example, in rotating black holes (see \cite{Dafermos2013a}).

In many of the previous uses of this method, the authors treated the Einstein equations (or a related system) in a way which does not explicitly reference the wave equation. Instead of using harmonic coordinates to treat the Einstein equations as a system of quasilinear wave equations, these authors typically choose coordinates to fix various components of the metric, so that the metric takes a particularly simple form. The dynamics of the system is then encoded in a set of \emph{first order} PDEs for the connection coefficients, which link the connection coefficients to the curvature coefficients, together with the Bianchi identities - another set of first order PDEs, but this time relating derivatives of the curvature coefficients to combinations of the curvature coefficients and the connection coefficients. This powerful approach was taken in \cite{Christodoulou1993}, and it has since been refined and used extensively throughout the literature (for example, see \cite{Bieri2010, Dafermos2013a, Dafermos2014b, An2017a, Luk2018, Christodoulou2011}, and references therein).

When pursuing this approach, the quantities involved - typically, the null frame connection coefficients and the null frame components of the curvature tensor - are treated either as scalar fields or as \emph{tensor fields ``on the spheres''}. In other words, we can use the null frame to construct a projection operator which projects onto the spheres, and then these tensor fields are spacetime tensor fields that are invariant under this projection. The operator $\slashed{\nabla}$ can then be viewed as the projection of the covariant derivative operator (with respect to the ambient metric $g$) onto the spheres. This operator raises the order of the tensor field it is applied to: for example, if $\phi$ is a scalar field, then $\slashed{\nabla} \phi$ is a covector field.

An important feature of this kind of method, as practised in the references given above, is that it treats \emph{first order} PDEs. Therefore, commuting with the operator $r\slashed{\nabla}$ can be done to some extent ``by hand'' - it is fairly easy to compute the commutators of first order operators with the covariant derivative, and it is also not too difficult to apply first order differential operators to the projection operators even though they raise the degree of the quantity in question.

On the other hand, we are dealing directly with wave equations, which are \emph{second order} operators. Consequently, commuting with the operator $r\slashed{\nabla}$ is significantly more complicated. Additionally, if $\phi$ is not a scalar field but a higher order tensor field on the spheres, then it is not at all clear what is meant by the wave operator applied to $\phi$. To confront these issues, we develop a more systematic - and more geometric - treatment, which we outline here. Note that an equivalent construction was developed in \cite{Holzegel2015, Holzegel2017}, however, these papers deal only with \emph{linear} equations. By contrast, we develop this formalism for the fully nonlinear problem, in which the metric depends on the solution to the wave equation.

Consider the vector bundle $\mathcal{B}$ whose fibres, at any point $p \in \mathcal{M}$, are given by the cotangent space, at the point $p$, of the (unique) sphere through the point $p$. Clearly, this is a two-dimensional vector bundle over $\mathcal{M}$. Given a section $\phi$ of this vector bundle, we can identify a corresponding section $\varphi$ of the cotangent bundle of $\mathcal{M}$ (i.e.\ a covector field $\varphi$), by requiring that
\begin{enumerate}[(i)]
	\item the restriction of $\varphi$ to a sphere is given by $\phi$
	\item $\varphi$ is invariant under the projection operator
\end{enumerate}
Similarly, given a covector field on $\mathcal{M}$ we can restrict it, at each point, to the sphere though that point, and in this way we form a section of the vector bundle $\mathcal{B}$.

We equip $\mathcal{B}$ with a connection $\slashed{\D}$ as follows: for any vector field $X$ on $\mathcal{M}$, we define $\slashed{\D}_X \phi$ to be the section of $\mathcal{B}$ corresponding (using the above correspondence) to the covector field $\D_X \varphi$, where $\D$ is the standard Levi-Civita connection on the cotangent bundle. This means that to calculate $\slashed{\D}_X \phi$ we can first compute $\D_X \varphi$, and then we restrict this covector field at each point $p$ to the sphere through the point $p$. In a similar way we can deal with higher rank ``tensors'' which are ``tangent to the spheres''.

We stress that this connection $\slashed{\D}$ is \emph{not} the same as the connection on the spheres defined by viewing them as submanifolds of $\mathcal{M}$, and using the Levi-Civita connection associated with the metric induced on the spheres by the ambient metric $g$. This latter connection - which we write as $\slashed{\nabla}$ - is a connection on the (co)tangent bundles of each individual sphere. On the other hand, $\slashed{\D}$ is a connection on the vector bundle $\mathcal{B}$, which is a vector bundle over the whole of the manifold $\mathcal{M}$, and not a vector bundle over a single sphere like the (co)tangent bundle of a sphere. Hence, we can make sense of objects like $\slashed{\D}_L \phi$, though we cannot make sense of objects like $\slashed{\nabla}_L \phi$, since $L$ is not a vector in the tangent space of any specific sphere\footnote{Here our notation differs slightly from that used in, for example \cite{Dafermos2013a}, where the authors \emph{do} write $\slashed{\nabla}_L$ (or rather, $\slashed{\nabla}_4$ since $e_4$ in their work corresponds to our vector $L$). In their work, if the vector field $X$ is not tangent to the spheres, then we are to understand their derivative operator $\slashed{\nabla}_X$ in the same way as our operator $\slashed{\D}_X$. We prefer to keep things as clear as possible, and also to make explicit the fact that $\slashed{\D}$ can be viewed as a connection on a two-dimensional vector bundle over the entire manifold $\mathcal{M}$, a fact which is only implicit in \cite{Dafermos2013a}.}. We also stress that $\slashed{\D}$ is not the same as the covariant derivative $\D$: for example, if $\phi$ is a covector field on $\mathcal{M}$ such that the projection of this covector field onto the spheres vanishes at every point, then $\slashed{\D} \phi = 0$ whereas $\D \phi$ need not vanish. Also note that, given a scalar field $\phi$ on $\mathcal{M}$, $\slashed{\D} \phi = \D \phi = \upd \phi$ is a section of the cotangent bundle, whose fibres are four-dimensional vector spaces, whereas $\slashed{\nabla}\phi$ is actually a section of $\mathcal{B}$, whose fibres are two-dimensional vector spaces.

With this in mind we can take a systematic approach to commuting with the operator $r\slashed{\nabla}$. This operator raises the rank of ``tensor fields that are tangent to the spheres'', which are really sections of $\mathcal{B}$ or its higher-rank analogues. The natural extension of the wave operator to sections of $\mathcal{B}$ is given by $\slashed{\Box}_g \phi := (g^{-1})^{\mu\nu} \slashed{\D}_\mu \slashed{\D}_\nu \phi$, where $\phi$ is a section of $\mathcal{B}$ (or its higher-rank analogues). The error terms that we obtain by commuting with the operator $r\slashed{\nabla}$ typically involve the \emph{curvature} of the connection $\slashed{\D}$, and derivatives of this curvature. As we have tried to emphasise above, the connection $\slashed{\D}$ is not the same as the Levi-Civita connection $\D$, and nor is its curvature. For example, some of the null-frame components of the curvature tensor of $\slashed{\D}$ are nonzero even in the case of flat Minkowski space. On the other hand, the curvature tensor involves two ``spacetime indices'', so it is not the same as the curvature of the connection induced on the spheres.

For a complete discussion of this topic, together with our calculations of the null-frame components of the curvature tensor, we refer the reader to chapter \ref{chapter geometry of vector bundle}.

\subsection{Elliptic estimates in the region near \texorpdfstring{$r = 0$}{r=0}}

Our main tool for deriving pointwise bounds will be the Sobolev inequalities on the spheres, together with commutation with the operator $r\slashed{\nabla}$. However, this approach fails to yield good pointwise estimates in the region close to $r = 0$, since the commutation operator $r\slashed{\nabla}$ is trivial there.

Instead, in the region $r \leq r_0$, we will rely on commutation with the time-translation vector field $T$, together with \emph{elliptic estimates}. This approach has been used previously, for example in \cite{Yang2013}. The main point is that, if we subtract off the terms involving (up to two) time derivatives, then the geometric wave operator $\Box_g$ (restricted to an individual leaf of our foliation) is uniformly elliptic in the region $r \leq r_0$. Then, we can use standard Schauder estimates to control various H\"older norms of the solution.

Importantly, this approach matches our general philosophy: it is extremely robust, and essentially relies on nothing more than that we can commute with a timelike vector field. In fact, unlike some of the other pointwise estimates, the elliptic estimates do not rely on any kind of (weak) null structure in the equations. The only drawback is that this method fails to provide us with any kind of decay in $r$. However, since we only apply these estimates in a region of bounded $r$, this does not cause us any problems.

\subsection{Improved energy decay for time derivatives}

Another aspect of the $r^p$-weighted method is that, suitably modified, it can allow us to recover \emph{improved decay} in $u$ for time derivatives of the solution, as well as for the solution in a spatially bounded region. This modified $r^p$-weighted method was first used in \cite{Schlue2010}, and has recently been refined and employed to great success in \cite{Angelopoulos2018a, Angelopoulos2018b}. In this paper, we show that these \emph{improved decay estimates} can also be used in the \emph{nonlinear} setting.

These improved decay estimates can be obtained in several ways, but all of them rely on commuting with a (possibly weighted) version of the vector field $L$. In the original method \cite{Schlue2010}, we commute with $L$, and then find that the maximum value of $p$ can be increased - for linear equations, it can be increased from $2$ to $4$. Other, more recent methods (\cite{Angelopoulos2018a, Angelopoulos2018b}) rely on commuting with $rL$ or even $r^2L$. In our case, we will struggle to commute with $r^2L$ due to the slow decay of the error terms towards null infinity, but we find that the other approaches are open to us. Indeed, we are already commuting with $rL$ (see the discussion in subsection \ref{subsection intro LL problem}) so this is the approach we take. Note that this has some technical advantages over commuting with $L$ and then taking a higher value of $p$ - for example, we can treat $\phi$ and $rL\phi$ in a similar way in the $p$-weighted energy estimates\footnote{There are additional potential advantages, but these are related to obtaining \emph{sharp} improved estimates, which are not available to us. From this point of view, there are even greater advantages to be had by commuting with $r^2L$}.

One outcome of these estimates is improved decay in $u$ of the \emph{energy} of the field $T\phi$ relative to the field $\phi$. This, of course, can be translated into improved pointwise decay (in $u$) for the field $T\phi$ and its derivatives. Another outcome of these estimates is an improved rate of decay (again, in $u$) in a region of bounded $r$. In fact, using the Morawetz estimate we can already pick out a subset of times (or values of $u$) such that the energy in a region of bounded $r$ decays more rapidly. The improved decay for the energy of $T\phi$ then allows us to interpolate between these times.

The astute reader might think of using these improved estimates to deal with the problems outlined in subsection \ref{subsection intro LL problem}. After all, the problematic term discussed there the lack of sufficient decay in $u$ for a $T$ derivative, which is precisely the kind of thing that we can improve upon by using the improved $p$-weighted estimates. However, in order to employ these methods we have to commute with $rL$, which already fixes the problem discussed in subsection \ref{subsection intro LL problem}! Moreover, in order to actually obtain any improvement in the rates of decay in $u$, we must \emph{also} assume additional \emph{a priori} decay in $u$ for various other quantities.

For a complete discussion of these improved estimates see appendix \ref{appendix improved energy decay}.

\subsection{Semi-global existence and uniqueness}

Although we have not yet been explicit about it, our central argument takes the form of a \emph{bootstrap} or \emph{continuity} argument. In other words, we will begin by making a set of assumptions regarding the solution, which can be imposed initially by choosing suitable initial data, and which should also hold for some (possibly short) amount of ``time''\footnote{Of course, in our case, we will use the ``retarded time'' $u$ in the region $r \geq r_0$. The level sets of this function are null, not spacelike.} beyond the time where the initial data is posed, say, for all times less than $t_1$. These assumptions will dictate the pointwise decay rates for various quantities, as well as some $L^2$-based bounds on various ``energies''. Under these assumptions, we will show that the bounds that we have assumed actually hold with \emph{smaller} constants up to the time $t_1$, a process which is referred to as ``improving'' or ``recovering'' the estimates. Then, by continuity, the original bounds must hold beyond the time $t_1$. It is then straightforward to establish that these bounds hold for all time, and therefore that the solution exists globally.

It is true that the main difficulty of such an argument is usually the establishing of the ``improved'' bounds, and indeed this occupies nearly all of our time in this paper. However, another necessary ingredient in using such an argument to prove global existence is a suitable ``\emph{local existence}'' (and uniqueness) result. By ``suitable'' here, we mean that it is an argument establishing local existence given initial data that satisfies the assumed pointwise and $L^2$ bounds (or perhaps under weaker assumptions). Moreover, the proof should also establish that, while the solution exists, the quantities satisfying the aforementioned bounds are at least \emph{continuous} in time.

For the kinds of nonlinear wave equations we are studying, this local existence result is classical \emph{when the initial data is posed on a spacelike hypersurface}. See, for example, \cite{Sogge2008}. In fact, one does not need to assume any kind of (weak) null structure in order to deduce the required local existence result: even equations like $\Box \phi = (\partial_t \phi)^2$ admit local solutions given initial data satisfying the required bounds\footnote{Solutions to this equation nevertheless blow up in finite time, because the ``required bounds'' do not hold \emph{uniformly} in time.}.

However, we are not placing initial data on an initial spacelike hypersurface - ``Cauchy data'' - but rather (at least in the region $r \geq r_0$) we are placing data on an initial \emph{null} hypersurface - ``characteristic data''. In other words, when we ask for local existence, we need local existence \emph{in} $u$ and not local existence in $t$. This should already be regarded as a ``semi-global'' result: in our coordinate system, it is global in $r$ but local in $u$. This is reflected by the fact that, already, such a result does not hold for the equation $\Box \phi = (\partial_t \phi)^2$, since transverse derivatives of the solution can blow up instantaneously. Note that a detailed result of this form has been proven in \cite{Luk2011}.

We also have to prove this result in the case where the initial data does not have finite energy, but only finite \emph{degenerate} energy together with a finite $p$-weighted flux. Hence the result we need is rather atypical. Nevertheless we find that, by using the weak null hierarchical structure present in the equation, local existence in $u$ can be established. See appendix \ref{appendix semi-global existence} for the full details.

\subsection{Additional structure in the Einstein equations in wave coordinates}
\label{subsection Einstein equations in wave coords}

Several times throughout the introduction we have mentioned that there is extra structure in the Einstein equations (in harmonic coordinates) that can be used to show that solutions to these equations behave significantly better than solutions to the more general equations we are studying. In this section we will detail this extra structure and outline the ways in which it can help to improve various estimates. Note that we do not pursue these estimates fully in this work: taking full advantage of the extra structure means, in some cases, changing our approach in some significant ways. Instead, in the main body of the text we simply point out on occasion how the extra structure can make improve certain estimates.

One major source of additional structure in the Einstein equations in harmonic coordinates is the harmonic coordinate condition itself. In fact, this was used to great effect in \cite{Lindblad2004} to prove the nonlinear stability of Minkowski space in harmonic coordinates. It also provided an essential tool in related work, for example \cite{LeFloch2015, Huneau2016a, Huneau2018a, Wyatt2018}. Since we do not use this structure, the proof that we give provides the first proof of the global existence of solutions to the nonlinear wave equations that constitute the Einstein equations in harmonic coordinates treated \emph{in isolation}; these equations are sometimes called the ``reduced Einstein equations''. On the other hand, if we prescribe initial data that satisfies the harmonic coordinate condition, and which (consequently) also satisfies the constraint equations, then we can couple the system of quasilinear wave equations (the ``reduced Einstein equations'') to the equation defining the harmonic coordinate condition. If we had used that structure, then many parts of the proof would have been significantly easier, and we could also have obtained improved decay estimates for many quantities.

The harmonic coordinate condition means that the map $\mathcal{M} \rightarrow \mathbb{R}^4$ given by the ``background'' or ``rectangular'' coordinates is a \emph{wave map}. This implies that certain first derivatives of the metric fields components are related to other first derivatives of other metric components. Specifically, in terms of rectangular Christoffel symbols, this gives the relation
\begin{equation}
\label{equation introduction wave coordinate condition}
(g^{-1})^{ab} \Gamma^a_{bc} = 0
\end{equation}
It turns out that, if the initial data satisfies this condition and also satisfies the constrain equations, then this ``gauge condition'' propagates.

If we expand equation equation \eqref{equation introduction wave coordinate condition} by first writing $g_{ab} = m_{ab} + h_{ab}$ where $m_{ab}$ are the rectangular components of the Minkowski metric, and then expand the resulting equation in the null frame, we find that it relates certain \emph{bad derivatives} of the metric to other \emph{good derivatives}. Specifically, we have
\begin{equation*}
\begin{rcases}
|L^a L^b \partial h_{ab}| \\
|(\slashed{g}^{-1})^{ab} \partial h_{ab}| \\
|\slashed{\Pi}_a^{\phantom{a}b} L^c \partial h_{bc}| 
\end{rcases}
\sim |\bar{\partial} h|
\end{equation*}
Here, $\slashed{g}$ is the metric induced on the spheres and $\slashed{\Pi}$ is the projection onto the spheres. The important point is that all of these quantities, which we would normally treat as \emph{bad derivatives}, can actually be written in terms of \emph{good derivatives}. Consequently they can be expected to decay much more rapidly in $r$ than they would otherwise be expected to.

This has far-reaching consequences for our proof. For example, the transport equation for the foliation density now reads
\begin{equation*}
L \log \mu \sim (\bar{\partial} h)
\end{equation*}
so that we can expect the inverse foliation density $\mu$ to be \emph{uniformly bounded} in $r$. In other words, ``shocks at infinity'' cannot form, and the equations behave more like quasilinear equations with the \emph{classical} null condition than the general quasilinear equations we consider. 

Additionally, if we substitute for some of the ``bad'' derivatives in the wave equations $\Box_g h_{ab} = F_{ab}$ (see subsection \ref{subsection intro example systems}) and then express these equations relative to the null frame, we again see that the only ``bad'' terms are contained in $F_{\Lbar\Lbar}$. However, now they take the form
\begin{equation*}
F_{\Lbar\Lbar} \sim - \widehat{|\Lbar \slashed{h}|}^2
\end{equation*}
In other words, the Einstein equations are more like the model system
\begin{equation*}
\begin{split}
\Box \phi_1 &= 0 \\
\Box \phi_2 &= (\partial \phi_1)^2
\end{split}
\end{equation*}
than the model system
\begin{equation*}
\begin{split}
\Box \phi_1 &= 0 \\
\Box \phi_2 &= (\partial \phi_1)(\partial \phi_2)
\end{split}
\end{equation*}
which previously appeared to match the structure of the semilinear terms more closely. Note that these two systems have different asymptotics, with $(\partial \phi_2) \sim r^{-1+\epsilon}$ in the first case, and $(\partial \phi_2) \sim r^{-1} \log r$ in the second case.

It turns out that the quantity $L \log \mu$ is also the kind of quantity which often appears as an error term in the energy estimates with a \emph{critical} decay rate. In fact, the only other source of error terms with a \emph{critical} decay rate in $r$ arise from the semilinear terms. With this in mind, we see that, for all the null components of $h$ \emph{apart from} $h_{\Lbar\Lbar}$, we will not have to employ the degenerate energy - we can work with the usual energy instead. 

In a similar vein, it should be possible to take larger values of $p$ (specifically, $p > 1$) in the $p$-weighted energy estimates for all of the metric components apart from $h_{\Lbar\Lbar}$. This not only means that related quantities will decay faster in $u$, but it also means that the radiation fields $\lim_{r \rightarrow \infty} r h$ can be shown to exist for these other metric components. This in turn has further consequences for the ``scattering problem''.

Additionally, we note that the error terms discussed in section \ref{subsection intro LL problem} can potentially be avoided if we make a slightly different definition of the radial coordinate $r$, so that the ``spheres'' are changed. One simple way this might be done is through a conformal transformation - i.e.\ writing the equations with reference to a \emph{conformally rescaled} metric $\Omega^2 g$ instead of the metric $g$. The idea is then to choose $\Omega$ so that this specific error term does not occur. We must be careful, however, to ensure that this conformal transformation leaves the weak null hierarchy intact, and does not introduce additional error terms that we cannot deal with. When we try to carry out this process in the general case, we discover that the additional error terms decay very slowly in $r$, and also that the asymptotic system of the resulting system is altered. Using the wave coordinate condition, however, it appears that neither of these problems arises when dealing with the Einstein equations\footnote{Although the conformal factor required does introduce some extra error terms, and it is not completely clear that these error terms can all be adequately handled, even in the case of the Einstein equations.}. It might therefore be possible to avoid commuting with $rL$ altogether. This would have the consequence of significantly enlarging the initial data set which can be considered to be ``sufficiently small'' to ensure global existence.

Finally, note that there is additional special structure in the semilinear structure of the Einstein equations. Specifically, the only ``bad'' semilinear terms have the form
\begin{equation*}
\widehat{|\partial \slashed{h}|}^2
\end{equation*}
but we can establish the \emph{sharp} pointwise bound $\widehat{|\partial \slashed{h}|} \lesssim \epsilon (1+r)^{-1}$ for these fields. In other words, the only \emph{bad} semilinear terms involve a pair of \emph{good} fields. Because of this, we can expect slightly better pointwise decay\footnote{Specifically, we can expect $|\partial h|_{\Lbar\Lbar} \sim r^{-1}\log r$ instead of $r^{-1+C\epsilon}$}. In addition, and perhaps more importantly, this semilinear structure allows us to perform the energy estimates \emph{directly} on the fields $h_{ab}$, that is, on the rectangular (or wave coordinate) components of the metric, rather than on the null-frame components.

To explain this point, suppose that we want to perform the energy estimates on all of the fields $\mathscr{Y}^n h_{ab}$. Each of these fields is a combination of all of the null frame components of $\mathscr{Y}^n h$. As such, we can expect the semilinear terms to include terms of the form
\begin{equation*}
\widehat{(\Lbar \slashed{h})} \widehat{(\slashed{\D}_{\Lbar} \mathscr{Y}^n \slashed{h})}
\end{equation*}
and in fact these will be the hardest semilinear terms to control. This second term can then be re-expressed in terms of the rectangular components. Substituting the pointwise bound for the first term, we effectively have to control terms of the form
\begin{equation*}
\epsilon (1+r)^{-1} (\slashed{\D} \mathscr{Y}^n h_{ab})
\end{equation*}
and this has just enough decay for us to handle. Hence we see that we only need to distinguish between the various null-frame components of the fields for the purpose of making the \emph{pointwise} bounds - when bounding the energy, we are able to treat the rectangular components directly. This, indeed, is the approach taken in \cite{Lindblad2004}.

On the other hand, we can imagine the case where the semilinear term $F_{\Lbar\Lbar}$ contains a term of the form
\begin{equation*}
(\Lbar h)_{LL} (\Lbar h)_{\Lbar\Lbar}
\end{equation*}
In fact, this term is present if we do not substitute for other terms using the wave coordinate condition. Then, if we try to follow the approach above and estimate the energy of $\mathscr{Y}^n h$, there are two important error terms we will encounter: those of the form
\begin{equation*}
(\Lbar h)_{LL} (\slashed{\D}_{\Lbar} \mathscr{Y}^n h)_{\Lbar\Lbar}
\end{equation*}
which can be treated exactly as above, and also those of the form
\begin{equation*}
(\Lbar h)_{\Lbar\Lbar} (\slashed{\D}_{\Lbar} \mathscr{Y}^n h)_{LL}
\end{equation*}
which cannot. To see why, note that we only expect the pointwise bound $(\Lbar h)_{\Lbar\Lbar} \lesssim (1+r)^{-1+C\epsilon}$. Hence, if we try to estimate the energy of the field $(\mathscr{Y}^n h_{ab})$ directly, then we find a semilinear term that behaves as
\begin{equation*}
\epsilon (1+r)^{-1+C\epsilon} (\slashed{\D} \mathscr{Y}^n h_{ab})
\end{equation*}
the coefficient of which does not have sufficient decay in $r$. However, if we can estimate the energy of the \emph{null-frame components} of $\mathscr{Y}^n h$ instead of the \emph{rectangular} components of $\mathscr{Y}^n h$ then we can handle such terms by using the following strategy: when estimating the energy of the field $(\mathscr{Y}^n h)_{LL}$ there are no such error terms, so there is no problem in estimating the energy of this field. Then, when estimating the energy of the field $(\mathscr{Y}^n h)_{\Lbar\Lbar}$ we can use a ``more degenerate'' weight $w$, in comparison with the weight used in bounding the energy of the field $(\mathscr{Y}^n h)_{LL}$. 

All this means that, in order to work with the more general case, we need to develop a strategy for dealing with the energy of the \emph{null frame} components. In even more generality, we can allow for some \emph{point dependent change of basis sections}, as long as this change-of-basis satisfies appropriate conditions. See subsection \ref{subsection intro changing basis sections} for an overview. Because of the special \emph{semilinear} structure of the Einstein equations in wave coordinates, all of this can be avoided.

\section{Statement of the theorem and structure of the proof}

We can now give a more complete version of our main theorem. For a complete version of the theorem, with all technical details included, see \ref{theorem main theorem}. \ref{theorem main theorem rough statement}

\begin{theorem}[Main theorem, second version]
\label{theorem second version}
Let $\phi_{(a)}$ be a set of scalar fields satisfying the wave equations
\begin{equation*}
\tilde{\Box}_{g(\phi)} \phi_{(a)} = F_{(a)}
\end{equation*}
on $\mathbb{R}^4$, where the inhomogeneous terms $F_{(a)}$ satisfy the \emph{weak null hierarchical condition}, possibly after a change of basis sections for the fields. We also allow for this change-of-basis to depend on the point $p$ on the manifold $\mathcal{M}$, and to depend on the solutions $\phi$, provided that the map which performs this change-of-basis satisfies suitable bounds. The metric $g(\phi)$ is to be written, relative to the standard coordinates on $\mathbb{R}^4$, as
\begin{equation*}
g_{ab} = m_{ab} + h_{ab}(\phi)
\end{equation*}
where $m_{ab} = \text{diag}(-1,1,1,1)$. The ``metric perturbations'' $h_{ab}$ can in turn be written as
\begin{equation*}
h_{ab} = h^{(0)}_{ab} + h^{(1)}_{ab}(\phi)
\end{equation*}
where $h^{(0)}_{ab}$ are some prescribed functions, satisfying some suitable ``smallness'' condition, which, however, does not imply decay towards timelike infinity for most metric quantities. The scalar fields $h^{(1)}_{ab}(\phi)$ are to be linear in the fields $\phi$, plus higher order terms, and the field $h^{(1)}_{LL}$ is to reside at the bottom level of the weak null hierarchy.

Suppose that the initial \emph{degenerate energy} and the initial $p$-weighted energy of the fields $\phi_{(a)}$ is sufficiently small. Here, we choose the weight $w = (1+r)^{-C_{(a)}\epsilon}$ and $p = 1-C_{(a)}\epsilon$, with the constant $C_{(a)}$ depending on the level in the hierarchy of the field $\phi_{(a)}$. Suppose that the same conditions also hold for the fields $\mathscr{Y}^n \phi_{(a)}$, where $\mathscr{Y}$ stands for any of the operators in the set $\{\slashed{\D}_T, r\slashed{\nabla}, r\slashed{\D}_L \}$ and $n \leq N$ for some sufficiently large $N$.

Then there exists a global solution $\{\phi_{(a)}\}$ to the system of wave equations. In addition, this solution obeys various bounds in both $L^2$ and $L^\infty$. The $L^2$ bounds include that the degenerate energy decays along a foliation by outgoing null hypersurfaces. The pointwise bounds are, in general, slightly worse than those that can be obtained for the linear wave equation: in particular, we have the decay rates $|\D \mathscr{Y}^n \phi_{(a)}| \lesssim \epsilon (1+r)^{-1 + C_{(n,a)}\epsilon}$. Importantly, the constant $C_{(0, a)}$ vanishes if $\phi_{(a)}$ is at the bottom level of the weak null hierarchy.

\end{theorem}

\vspace{3mm}

We now outline the structure of the proof of this theorem, as presented in the rest of the paper.
\vspace{2mm}

\textbf{Chapter \ref{chapter preliminaries}} introduces much of the notation we use, and also presents the geometric setting of our proof. For example, we define the foliation that we use, the null frame, we define the ``spheres'' and the geometric coordinates. We also give a more thorough treatment of the wave coordinate condition, although we do not assume that this condition is satisfied for our proof.

In \textbf{chapter \ref{chapter transport eikonal}} we derive certain transport equations along outgoing null geodesics. These are the transport equations which will allow us to control the ``eikonal'' function $u$. As such, they are intimately connected with the behaviour of the foliation by outgoing null geodesics, and their structure is essential for ruling out shock formation (except, possibly, for ``shock formation at infinity'').

\textbf{Chapter \ref{chapter null frame connection coefficients}} defines the rest of the null frame connection coefficients. Algebraic relations between some of these quantities and the first derivatives of $h$ are then deduced - in particular, in the region $r \leq r_0$, \emph{all} of the connection coefficients are related (algebraically) to $h$ and its first derivatives, while in the region $r \geq r_0$ this is true only for a subset of the connection coefficients. We also derive some equations for related quantities: for example, the derivatives of the rectangular components of the null frame, and derivatives of the projection operator. Additionally, we express both the scalar wave operator and the projected wave operator (which acts on sections of $\mathcal{B}$) in terms of the null frame vector fields and connection coefficients.

In \textbf{chapter \ref{chapter weak null structure}} we discuss, in full detail, the \emph{weak null hierarchy}, which is our main structural restriction on the equations. This involved digressions into the definition of the ``reduced wave operator'' $\tilde{\Box}$, as well as changing the basis sections when dealing with a \emph{system} of wave equations, including ``point dependent'' changes of basis sections. We also include a discussion of the action of a conformal rescaling on the structure of the equations, which is essential in the light of the ``radial normalisation'' condition that we impose in the rest of the paper.

\textbf{Chapter \ref{chapter geometry of the null cones}} presents a thorough investigation into the geometry of the outgoing null cones, i.e.\ the surfaces of constant $u$. In particular, we derive numerous equations for those null frame connection coefficients which \emph{cannot} be written algebraically in terms of $h$ and its derivatives in the region $r \geq r_0$. Most of these equations take the form of transport equations along the integral curves of $L$, including transport equations for a specially modified versions of the connection coefficient $\tr_{\slashed{g}} \chi$, which allow us to avoid ``losing a derivative''. There are also elliptic equations for various quantities on the spheres.

\textbf{Chapter \ref{chapter geometry of vector bundle}} is devoted to the geometry of the the vector bundle $\mathcal{B}$. Most importantly, in this chapter we derive explicit expressions for the null frame components of the \emph{curvature} of this vector bundle, defined with respect to the connection $\slashed{\D}$. The null frame components of the curvature are expressed in terms of $h$ and its derivatives, and also the null frame connection coefficients and their derivatives.

\textbf{Chapter \ref{chapter deformation tensors}} presents various calculations relating to \emph{deformation tensors}. These are of general use throughout the rest of the proof, but in particular, in chapter \ref{chapter deformation tensors} we use them to compute quantities associated with using various vector fields as \emph{multipliers} in the energy estimates. A multiplier is the analogue of the vector field $\partial_t$ in the standard energy estimate, which can be obtained by multiplying by $(\partial_t \phi)$ and integrating by parts. In this chapter, we include a prescription for \emph{modifying} the multipliers by lower order terms (necessary for both the Morawetz estimate and the $p$-weighted estimates). We also compute various error terms associated with the multipliers, in particular, the  ``bulk'' error terms, which are spacetime integrals over terms constructed from the deformation tensors. Since the multipliers are defined geometrically, we find that these error terms involve the null frame connection coefficients.

\textbf{Chapter \ref{chapter commuting}} presents the relevant calculations for \emph{commuting} with the commutation operators $\mathscr{Y}$. We begin by computing the commutators of the operators $\mathscr{Y}$ with various other first order operators. Then, we present a systematic treatment of the terms produced by commuting the projected (reduced) wave operator with either a ``vector field''-type operator (such as $\slashed{\D}_T$) or a ``higher rank'' operator (such as $r\slashed{\nabla}$. We then apply this to our set of commutation operators $\mathscr{Y}$. Finally, we present (schematically) the equations satisfied by various geometric quantities (such as the null frame connection coefficients) after having commuted some arbitrary number of times with the commutation operators.

\textbf{Chapter \ref{chapter elliptic estimates and sobolev embedding}} presents the technical tools that we use to obtain pointwise bounds from $L^2$ bounds, namely, elliptic estimates in the region $r \leq r_0$ and Sobolev embedding on the spheres in the region $r \geq r_0$.

\textbf{Chapter \ref{chapter framework for energy estimates}} gives a number of valuable results which we use as part of our energy estimates, including various Hardy inequalities, coarea formulae and estimates for the spherical mean of a function in terms of its degenerate energy. We also compute the various ``boundary'' error terms that are present in the energy estimates - that is, error terms which are integrals over a single leaf of the foliation, rather than spacetime integrals. Again, since both the foliation and the multiplier vector fields are constructed geometrically, these error terms involve the null frame connection coefficients.

In \textbf{chapter \ref{chapter bootstrap}} we finally present the various \emph{bootstrap bounds} which we make. Before this chapter, all of our computations are abstract and the results are completely general. After this chapter, we assume that various quantities are ``small''. We include both pointwise and $L^2$ bootstrap bounds in this chapter.

In \textbf{chapter \ref{chapter energy estimates}} we finally begin to carry out the energy estimates. There are three basic types of energy estimates that we use: the weighted $T$ energy estimate, the Morawetz estimate, and the $p$-weighted estimate. In this chapter, we present each of these energy estimates, along with expressions for the associated error terms, assuming the pointwise bootstrap bounds from the previous chapter.

\textbf{Chapter \ref{chapter boundedness and energy decay}} brings together the energy estimates of the preceding chapter in order to establish first \emph{boundedness} of the degenerate energy, and then degenerate energy \emph{decay}. In particular, we use all three of the basic energy estimates from the previous chapter, and combine them in order to control all of their associated error terms simultaneously. We also obtain boundedness and decay for various spacetime integrals (or ``bulk terms'') in this chapter. Finally, we include a discussion of the case where a \emph{point-dependent change of basis section} is made before carrying out the energy estimates.

\textbf{Chapter \ref{chapter pointwise bounds}} is focussed on establishing \emph{pointwise bounds}. In other words, assuming the pointwise bootstrap bounds hold, we derive other pointwise bounds on various quantities which relate them, in the end, to the \emph{energy} of the fields and their derivatives. We include pointwise bounds on the fundamental quantities - the fields $\phi_{(a)}$ and their derivatives - as well as pointwise bounds on all the relevant geometrical quantities, such as the foliation density or the null frame connection coefficients. We also give pointwise bounds in the region $r \leq r_0$ using elliptic estimates.

While all of the previous analysis has focussed on the abstract equation $\tilde{\slashed{\Box}}_g \phi = F$, in \textbf{chapter \ref{chapter bounds for the inhomogeneous term}} we finally compute the (schematic form of the) inhomogeneous terms $F$ for the specific equations we are dealing with. This involves collecting together all of the error terms coming from both the ``original'' inhomogeneous terms (that is, the terms $F_{(a)}$ in the equations $\tilde{\Box}_g \phi_{(a)} = F_{(a)}$) as well as all of the error terms caused by commuting the commutator operators $\mathscr{Y}$ through the reduced wave operator $\tilde{\Box}_g$. We then use the bootstrap bounds (both the $L^2$ bounds and the $L^\infty$ bounds) to estimate the size of these inhomogeneous terms. In particular, this involves deriving $L^2$ bounds for various \emph{geometric} error terms (for example, error terms involving the null frame Christoffel symbols), which relate these quantities to various $L^2$-based quantities involving the metric and its derivatives. At this point we have to use the equations we have previously derived for these geometric quantities, including the modified transport equations and the elliptic equations in order to avoid a loss of derivatives. Finally, we include in this chapter another discussion of the issue of a point-dependent change of basis, and how this affects the estimates if it is necessary to perform such a change.

Finally, \textbf{chapter \ref{chapter proving the theorem}} brings together all of the estimates of the preceding chapters, to finish the proof of global existence for small initial data.

We also include three appendices. \textbf{Appendix \ref{appendix improved energy decay}} presents the \emph{improved} energy estimates, which, among other things, establish additional decay in $u$ for the energy associated with the $T$ derivative of a field, under the assumption of additional pointwise decay in $u$ for various geometric quantities. 

\textbf{Appendix \ref{appendix semi-global existence}} outlines a proof of the semi-global existence and uniqueness of solutions to the kinds of nonlinear wave equations we are considering, where ``semi-global'' means \emph{local in $u$} but global in $r$. Note that this kind of result already fails for equations like $\Box \phi = (\partial_t \phi)^2$.

Finally, note that our proof is \emph{consistent} with ``shock formation at infinity'', as outlined above, but this does not imply that this kind of behaviour actually occurs, or even that it is possible. To address this issue, in \textbf{appendix \ref{appendix explicit shock formation}}, we presents an explicit example of shock formation at infinity for a particular wave equation obeying the weak null condition, showing that shocks really can form at infinity. In fact, we present a family of initial data for this equation consisting of smooth, compactly supported functions, which can be made arbitrarily small. Nevertheless, \emph{every} member of this family of initial data exhibits shock formation at infinity. Furthermore, the shocks form ``immediately'', that is, they form as $r \rightarrow \infty$ \emph{on the initial data surface}, despite the initial data being trivial\footnote{Perhaps this is not as surprising as it first seems. Recall that shock formation is intimately connected with the behaviour of the \emph{transverse} (``bad'') derivatives. However, when posing initial data on a null hypersurface, we do not prescribe the transverse derivatives; rather, these must be found by solving an ODE. It turns out that the solutions to this ODE can have sufficiently ``bad'' decay towards null infinity, even when the initial data is compactly supported in $r$.} in the region $r \geq r_0$.

In addition to showing that shocks can form at infinity, appendix \ref{appendix explicit shock formation} also gives explicit examples of ``blowup at infinity'' - i.e.\ examples of systems where the solutions \emph{actually} do not have the same asymptotics as the linear equations. This issue has been discussed extensively in the literature, particularly by Alinhac (\cite{Alinhac1995, Alinhac2003, Alinhac2006, Alinhac2012}) but, using our methods, we are able to give the first explicit example of such ``blowup'' for these kinds of equations, albeit arising from (arbitrarily small) initial data that is compactly supported on an initial (asymptotically) \emph{characteristic} slice, rather than from an initial spacelike slice.

\vspace{4mm}

\textbf{Acknowledgements}
\vspace{2mm}

I am very grateful to a large number of people for many helpful discussions and comments on this work. In particular, I would like to thank Mihalis Dafermos, Hans Lindblad and Igor Rodnianski for very valuable discussions and for some comments on the manuscript. I am also indebted to Thomas Johnson for his detailed comments on the introduction. There are many others with whom I have had useful discussions regarding this work, including in particular Georgios Moschidis, Shiwu Yang and Harvey Reall. Finally, I would like to thank my wife Daphne, whose help and support has been essential.

\chapter{Preliminaries}
\label{chapter preliminaries}

\section{The basic geometric set-up}

\subsection{Rectangular coordinates, the eikonal function and foliations}
We begin with a 4-dimensional manifold $\mathcal{M}$ endowed with both a dynamical Lorentzian metric $g$ and a ``background'' Minkowski metric $m$.
\begin{definition}[Rectangular coordinates and the radial function]
We define the standard rectangular coordinates on $\mathbb{R}^4$, $(x^0, x^1, x^2, x^3) := (t, x^1, x^2, x^3)$, with respect to which the Minkowski metric is given by $m = \mathop{\mathrm{diag}}(-1,1,1,1)$. In addition, we define the standard radial function
\begin{equation}
r:= \sqrt{(x^1)^2 + (x^2)^2 + (x^3)^2} 
\end{equation}
\end{definition}

\begin{definition}[The eikonal function]
We define an outgoing eikonal function $u$ with initial data on the hypersurface $r = r_0$, satisfying
\begin{equation}
\begin{split}
 g^{-1}(\upd u,\ \upd u) &= 0 \quad \text{for}\quad r > r_0\\
 u\big|_{r = r_0} &= t - r
\end{split}
\label{eikonal}
\end{equation}
for some fixed radius $r_0$. We extend $u$ to a function on the whole of $\mathcal{M}$ by setting
\begin{equation}
 u := t - r \quad \text{for}\quad r \leq r_0
\end{equation}
\end{definition}

\begin{definition}[The hyperboloidal time]
We define the geometric hyperboloidal time variable $\tau$ by
\begin{equation}
 \tau := \begin{cases}
          t \quad &\text{ if} \quad r \leq r_0 \\
          u + r_0 \quad &\text{ if} \quad r \geq r_0 
         \end{cases}
\end{equation}
\end{definition}

\begin{remark}[Continuity of $\tau$]
Note that the function $\tau$ is continuous, since at $r = r_0$ we have $u = t - r_0$.
\end{remark}

% \begin{definition}[The geometric time coordinate $\tilde{t}$]
%  Recall that the manifold $\mathcal{M}$ is already equipped with the ``background'' or ``Minkowski'' time coordinate $t$. In addition, we define a geometric time coordinate $\tilde{t}$ by
% \begin{equation}
%   \begin{split}
%     \tilde{t} &:= u + r \\
%       &\phantom{:}= \tau - f_{(\alpha)}(r) + r
%   \end{split}
% \end{equation}
% 
% \end{definition}

\begin{definition}[Hypersurfaces and the foliations]
\label{definition hypersurfaces}
Using these coordinates, we foliate the future of some initial hypersurface $\Sigma_{\tau_0} := \left\{ x \in \mathcal{M} | \tau(x) = \tau_0 \right\}$ by the hypersurfaces
\begin{equation}
 \Sigma_\tau := \left\{ x \in \mathcal{M} | \tau(x) = \tau \right\}
\end{equation}
We also define the ``cut-off'' versions of these hypersurfaces:
\begin{equation}
 ^t\Sigma_\tau := \left\{ x \in \mathcal{M} | \tau(x) = \tau, \, t(x) \leq t \right\}
\end{equation}
as well as the ``cut-off'' version of the surface of constant $t$:
\begin{equation}
 _{\tau_0}^{\tau_1}\bar{\Sigma}_t := \left\{ x \in \mathcal{M} | t(x) = t, \tau_0 < \tau(x) \leq \tau_1 \right\}
\end{equation}
% and the ``cylinders'' of constant $r$:
% \begin{equation}
%  C_r := \{x \in \mathcal{M} | r(x) = r \}
% \end{equation}
We define the ``spheres'' of constant $\tau$ and $r$:
\begin{equation}
 S_{\tau, r} := \left\{ x \in \mathcal{M} | \tau(x) = \tau,\ r(x) = r \right\}
\end{equation}
and similarly the ``spheres'' of constant $\tau$ and $t$:
 \begin{equation}
  \bar{S}_{\tau, t} := \left\{ x \in \mathcal{M} | \tau(x) = \tau,\ t(x) = t \right\}
 \end{equation}
Finally, we define the spacetime region
\begin{equation}
 \mathcal{M}_{\tau_1}^{\tau_2} := \{x \in \mathcal{M} | \tau(x) \in [\tau_1, \tau_2] \}
\end{equation}
as well as the ``cut-off'' versions of these spacetime regions:
\begin{equation}
 ^{T}\mathcal{M}_{\tau_1}^{\tau_2} := \{x \in \mathcal{M} | \tau(x) \in [\tau_1, \tau_2] , \, t(x) \leq T \}
\end{equation}

\end{definition}

\begin{figure}[htb]
	\centering
	\includegraphics[width = 0.9\linewidth, keepaspectratio]{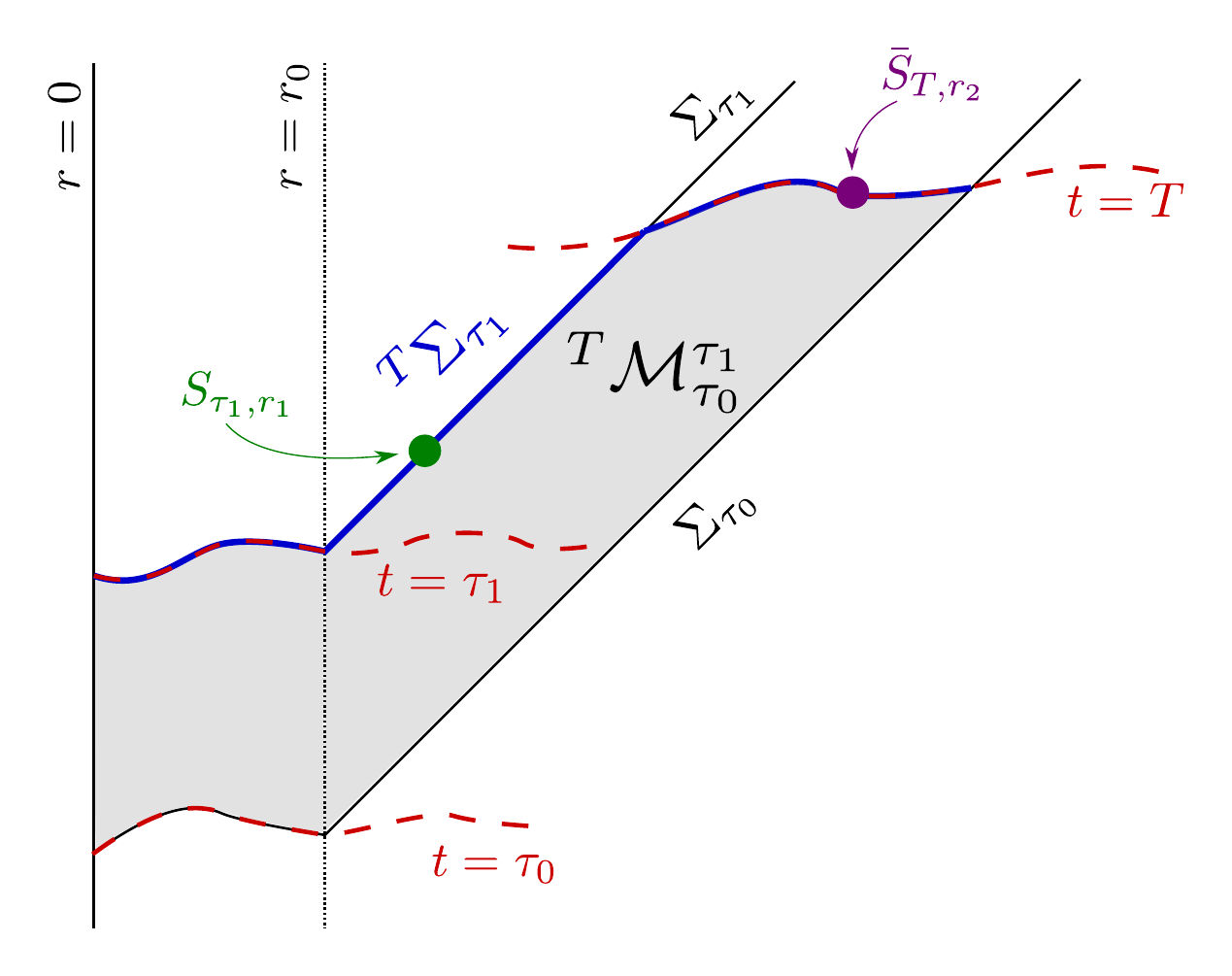}
	\caption{
	A causal diagram of the various spacetime regions. Note that, if we imagine that every point on the diagram represents a sphere $S_{\tau,r}$ (so that this is a sketch of the $(\tau, r)$ plane), then we cannot really draw the curves of constant $t$ in the region $r \geq r_0$, because the sphere $\bar{S_{t, r}}$ does not correspond to any spheres $S_{\tau, r}$. In this sense, we can view the curves of constant $t$ sketched on the diagram, in the region $r \geq r_0$, as actually being curves at some particular \emph{fixed} angle on the sphere.
	}
	\label{figure spacetime regions}
\end{figure}

\begin{definition}[Induced metrics]
Recall that we have a background metric $g$ on $\mathcal{M}$. This metric induces (by restriction) a metric $\underline{g}$ on the hypersurface $\Sigma_\tau$, as well as a metric $\slashed{g}$ on the sphere $S_{\tau,r}$.
\end{definition}

\subsection{Indices}
We make use of the following conventions relating to indices.

\begin{definition}[Abstract indices]
Greek indices $\mu, \nu, \rho, \ldots$ will be used to label \emph{abstract} indices, with the usual convention that raised indices refer to the tangent space and lowered indices refer to the cotangent space. So, for example, in the expression $A^\mu$, the index $\mu$ serves only to indicate that the quantity $A$ is a vector. Consequently, any expression involving only Greek indices will transform covariantly, and repeated Greek indices indicate contractions in the usual way. \end{definition}

\begin{definition}[Rectangular components]
Lower case Latin indices from the start of the alphabet $a, b, c, \ldots$ refer to the rectangular coordinate system $(x^0, x^1, x^2, x^3)$ and the corresponding bases induced on the tangent and cotangent spaces. So, for example, $g_{a, b} := g(\partial_a, \partial_b)$, and for a vector $V$, $V^a := V(x^a) = \upd x^a \cdot V$.

Lower case Latin indices from the middle of the alphabet $i, j, k, \ldots$ refer to the \emph{spatial} rectangular coordinates $(x^1, x^2, x^3)$. For instance, the radial function $r$ satisfies $r^2 = \sum_{i = 1}^3 (x^i)^2$. Note that, unlike other types of indices, we will sum over repeated spatial rectangular indices \emph{regardless of whether those indices are raised or lowered}. For example, we understand $x^i \upd x^i$ to mean $\sum_{i = 1}^3 x^i \upd x^i$.
\end{definition}

\begin{definition}[Angular components]
Upper case Latin indices $A, B, C, \ldots$ refer to the geometric angular coordinates $(\vartheta^1, \vartheta^2)$ defined below. These are coordinates on the spheres $S_{\tau, r}$. Although we shall fix this coordinate system on the spheres (see section \ref{section null frame}) many of the expressions involving these indices transform covariantly and so can be reinterpreted as referring to $S_{\tau, r}$-tensors.
\end{definition}

\begin{remark}[Index conventions]
Unless explicitly stated, we adopt the usual summation convention, i.e.\ repeated indices which refer to a particular coordinate system are summed over.

We use the dynamic metric $g$ and its inverse $g^{-1}$ to lower and raise either abstract indices or rectangular indices. For example, given a vector field $V^\mu$ with rectangular components $V^a$, we define $V_{\mu} := V^\nu g_{\mu\nu}$ and $V_a := V^b g_{ab}$. Likewise, we use the restriction of the metric to the spheres $\slashed{g}$ and its inverse $(\slashed{g})^{-1}$ to lower and raise indices of $S_{\tau, r}$-tensors.
\end{remark}

\begin{definition}[Frame components]
Finally, we will also make use of frame indices: given a vector field $V$, we define \emph{lower} frame indices by contractions, using the metric $g$ if necessary. For example, for a 1-form $\omega$ we define $\omega_V := \omega(V) = \omega_\mu V^\mu$, while for a vector $Z$ we define $Z_V := g(Z, V) = Z^\mu V^\nu g_{\mu\nu}$. Note that we do not actually need to have defined a frame in order to define these lower ``frame'' indices!

On the other hand, given a frame $(V_0, V_1, V_2, V_3)$ of vector fields spanning the tangent space $T_p(\mathcal{M})$ at a point $p$, we can also define \emph{raised} vector indices by expanding vectors in this frame, or by expanding covectors in the basis $(V_0^\flat, V_1^\flat, V_2^\flat, V_3^\flat)$, where we define
\begin{equation}
 (V^\flat)_\mu := V^\nu g_{\mu\nu}
\end{equation}
To be precise, we define raised frame indices for a vector field $Z$ by expanding
\begin{equation*}
 Z := Z^{V_0} V_0 + Z^{V_1} V_1 + Z^{V_2} V_2 + Z^{V_3}V_3
\end{equation*}
while for a covector $\omega$ we define the raised frame indices by
\begin{equation*}
 \omega := \omega^{V_0} V_0^\flat + \omega^{V_1} V_1^\flat + \omega^{V_2} V_2^\flat + \omega^{V_3}V_3^\flat
\end{equation*}
Since the frame consists of linearly independent vector fields spanning the tangent space, if $g$ is non-degenerate then these relations uniquely define the frame components $Z^{V_a}$ and $\omega^{V_a}$ for $a = 0, 1, 2, 3$.

Note that, when frame indices are used for the angular vector fields $X_A$, in order to avoid cluttering the notation we will make the identifications $V^A := V^{X_A}$ and $V_A := V_{X_A}$.

%To denote an arbitrary frame component we will use upper case bold Latin letters, and, as usual, we shall sum over repeated indices. For example, with respect to the frame $(V_0, V_1, V_2, V_3)$ discussed above, we can write the vector $Z$ as
%\begin{equation*}
% Z = Z^{\bm{A}} \bm{A}
%\end{equation*}
%where the frame indices $\bm{A}$ takes the values $V_0, V_1, V_2$ and $V_3$.

\end{definition}

\begin{definition}[Musical notation]
In general, we define musical duals of a vector field $V$ and a covector field $\omega$ by using the metric and its inverse as follows: for \emph{all} vector fields $Z$, we define
\begin{equation}
 \begin{split}
  V^\flat \cdot Z &:= g(V, Z) \\
  g(\omega^\sharp, Z) &:= \omega \cdot Z 
 \end{split}
\end{equation}
and we extend these definitions to general tensors in the obvious way.
\end{definition}

\begin{definition}[Indices outside derivative operators]
 If frame indices are placed outside of a delimiter, then any operators inside the delimiter are to be applied \emph{before} contracting with the frame fields, and they are applied to the \emph{rectangular} indices. For example, given a tensor $h_{\mu\nu}$ and some operator, say the vector field $Z$, we have
\begin{equation*}
 \begin{split}
  (Zh)_{LL} &:= L^a L^b (Z h_{ab}) \\
  |Zh|_{LL} &:= |L^a L^b (Z h_{ab})|
 \end{split}
\end{equation*}
Note that this means that the operator does not operate on the frame fields.
\end{definition}

\begin{definition}[Symmetric and antisymmetric tensors]
 We use round brackets to donate the totally symmetric part of a tensor, and square brackets to donate the totally antisymmetric part of a tensor. For example, given a tensor $T = T_{\mu\nu}$, we define
\begin{equation}
 \begin{split}
  T_{(\mu\nu)} &= \frac{1}{2}\left( T_{\mu\nu} + T_{\nu\mu} \right) \\
  T_{[\mu\nu]} &= \frac{1}{2}\left( T_{\mu\nu} - T_{\nu\mu} \right) \\
 \end{split}
\end{equation}
 We also use vertical lines to indicate that the enclosed indices are \emph{not} included in any symmetrisation or antisymmetrisation, for example
\begin{equation}
 T_{(\mu |\nu| \rho)} = \frac{1}{2}\left( T_{\mu\nu\rho} + T_{\rho\nu\mu} \right)
\end{equation}

\end{definition}

\subsection{Derivative Operators}
We will make use of the standard exterior derivatives and Lie derivatives. In addition, we define several other derivative operators:

\begin{definition}[Covariant derivatives]
We define the covariant derivative operator $\D$ with respect to the metric $g$. Likewise, we define the flat covariant derivative $\nabla$ with respect to the metric $m$, and finally, we define the covariant derivative with respect to the metric $\slashed{g}$, which acts on functions or tensors defined on the spheres $S_{\tau,r}$ and their associated tangent and cotangent bundles.
\end{definition}

\begin{definition}[Second order operators]
We also define the second order operators
\begin{equation}
\begin{split}
 \Box_{g} := (g^{-1})^{\mu\nu} \D^2_{\mu\nu} &:= (g^{-1})^{\mu\nu} \D_\mu \D_\nu \\
 \Box := (m^{-1})^{\mu\nu} \nabla^2_{\mu\nu} &:= (m^{-1})^{\mu\nu} \nabla_\mu \nabla_\nu \\
 \slashed{\Delta} := (\slashed{g}^{-1})^{AB}\slashed{\nabla}^2_{AB} &:= (\slashed{g}^{-1})^{AB}\slashed{\nabla}_A \slashed{\nabla}_B \\
\end{split}
\end{equation}
Note that, although we have not yet defined the coordinates $\{\vartheta^A \}$ which are involved in the definition of $\slashed{\Delta}$ appearing above, this expression is evidently covariant and so it is actually independent of the choice of coordinates $\{\theta^A \}$.

\end{definition}

\begin{definition}[Rectangular derivatives]
 The rectangular derivative vector fields are defined by
\begin{equation}
 \begin{split}
  \partial_t &:= \frac{\partial}{\partial t}\Big|_{x^1, x^2, x^3} \\
  \partial_{x^1} &:= \frac{\partial}{\partial x^1}\Big|_{t, x^2, x^3} \\
  \partial_{x^2} &:= \frac{\partial}{\partial x^2}\Big|_{t, x^1, x^3} \\
  \partial_{x^3} &:= \frac{\partial}{\partial x^3}\Big|_{t, x^1, x^2}
 \end{split}
\end{equation}
In other words, when taking derivatives relative to the rectangular coordinates, it is understood that the other rectangular coordinates are to be held fixed.
\end{definition}

\section{The Null Frame and Geometric Coordinates}
\label{section null frame}

\subsection{Definition of the null frame and geometric coordinates}
\begin{definition}[The outgoing null geodesic generator]
We define the vector field
\begin{equation}
 (L_{(\text{Geo})})^\mu := -(g^{-1})^{\mu\nu} \partial_\nu u \quad \text{for} \quad r \geq r_0 
 \end{equation}
Since $u$ satisfies the eikonal equation \ref{eikonal} in the region $r \geq r_0$, the vector field $L_{(\text{Geo})}$ is null in this region. We can also see that, for $r \geq r_0$, $L_{(\text{Geo})}$ is geodesic:
\begin{equation*}
\begin{split}
 (\D_{L_{(\text{Geo})}} L_{(\text{Geo})})^\mu &= (\D^\nu u) (\D_\nu \D^\mu u) \\
 &= (\D^\nu u) (\D^\mu \D_\nu u) \\
 &= \frac{1}{2}\D^\mu \left( g(L_{(\text{Geo})}, L_{(\text{Geo})}) \right) = 0
\end{split}
\end{equation*}
where we have used the torsion-free property of the Levi-Civita connection.

We extend $L_{(\text{Geo})}$ to the region $r \leq r_0$ by choosing $L$ to be the unique vector field such that, in this region,
\begin{itemize}
	\item $L_{(\text{Geo})}$ is null and future-directed
	\item $L_{(\text{Geo})}$ is orthogonal to the spheres $S_{t,r}$ (using the metric $g$)
	\item $L_{(\text{Geo})}(r) = \mu^{-1}$
\end{itemize}
\end{definition}

\begin{definition}[The inverse foliation density]
We define a very important quantity, the \emph{inverse foliation density} $\mu$ by
\begin{equation}
 \mu^{-1} := -(g^{-1})(\upd r, \upd u)
\end{equation}
Note that this is actually only the inverse foliation density in the region $r \geq r_0$. Nevertheless, we will also write $\mu$ for this quantity in the region $r < r_0$.
\end{definition}

\begin{definition}[The null vector fields $L$ and $\Lbar$]
\label{definition L and Lbar}
Using the foliation density we can define the vector $L$ by
\begin{equation}
 L := \mu L_{(\text{Geo})}
\end{equation}
Note that this vector field is null everywhere and satisfies $L(r) = 1$. Now, we define $\Lbar$ by setting $\Lbar$ to be the unique null vector field satisfying
\begin{itemize}
	\item $\Lbar$ is null and future-directed
	\item $\Lbar$ is orthogonal to the spheres $S_{t,r}$ (using the metric $g$)
	\item $g(L, \Lbar) = -2$
\end{itemize}
\end{definition}

\begin{definition}[Coordinates on the spheres]
We now need to define coordinates on the spheres $S_{\tau, r}$. We fix an atlas $\{ (\mathbb{D}_\alpha, \vartheta^1_\alpha, \vartheta^2_\alpha)\}_{\alpha = 1,2}$ on the sphere $S_{0, r_0}$, where the $\mathbb{D}_\alpha$ are open subsets of $S_{0, r_0}$ which satisfy $S_{0, r_0} = \mathbb{D}_1 \cup \mathbb{D}_2$. We transport these coordinates onto the spheres $S_{0, r}$ for $r < r$ by using the Euclidean vector field $- \frac{x^i}{r}\partial_i$, i.e.\ setting
\begin{equation*}
 - \frac{x^i}{r} \partial_i \vartheta^A_\alpha = 0
\end{equation*}
Now, we transport these coordinates onto the whole of the spacetime tube $r \leq r_0$ by using the Euclidean time translation vector field $\partial_t := \partial_0$, i.e.\ we set
\begin{equation*}
 \frac{\partial \vartheta^A_\alpha}{\partial t}\Big|_{x^1, x^2, x^3} = 0 \quad \text{for}\quad  r \leq r_0
\end{equation*}
note that the vector field $\partial_t$ is tangent to the hypersurface $r = r_0$, so this definition is well-defined. Finally, we extend the coordinates to the region $r > r_0$ by transporting the coordinates using the vector field $L$, i.e.\ we set
\begin{equation*}
 L \vartheta^A_\alpha = 0 \quad \text{for}\quad r > r_0
\end{equation*}
\end{definition}

\begin{definition}[The angular vector fields $X_A$]
In the region now covered by coordinates arising from one of the open sets $\mathbb{D}_\alpha$, we define the vector fields
\begin{equation}
\begin{split}
 X_1 &:= \frac{\partial}{\partial \vartheta^1} \Big|_{\tau, r,  \vartheta^2} \\
 X_2 &:= \frac{\partial}{\partial \vartheta^2} \Big|_{\tau, r,  \vartheta^1}
\end{split}
\end{equation}
\end{definition}

\begin{remark}[Suppression of angular indices]
 There will be occasions in which we will deal with $S_{\tau,r}$ tensors of undetermined rank. In this case, we will suppress the associated angular indices, and use a dot to donate contraction using the metric on $S_{\tau,r}$, and take norms using the metric $\slashed{g}$ and its inverse. For example, if $\phi_\mu := \phi_A \upd \vartheta^A$ is an $S_{\tau,r}$ one-form, we define
 \begin{equation*}
  |\phi|^2 = \phi \cdot \phi = (\slashed{g}^{-1})^{AB} \phi_A \phi_B
  \end{equation*}
\end{remark} 

\begin{remark}[Suppressing the coordinate chart label]
As we have already mentioned, many of our expressions transform covariantly under a change of basis for the tangent space of the $S_{\tau, r}$. In particular, this means that, without ambiguity, we can suppress the coordinate chart label $\alpha$.
\end{remark}

\begin{definition}[The null frame]
 We define the null frame as the set of vector fields $\{L, \Lbar, X_1, X_2\}$.
\end{definition}

\begin{remark}
 From now on we will assume that the null frame spans the tangent space $T_p \mathcal{M}$ at every point $p \in \mathcal{M}$. This will eventually be justified when we close the argument.
\end{remark}

\begin{definition}[Geometric coordinates]
 We refer to the coordinates $(\tau, r, \vartheta^1, \vartheta^2)$ as the \emph{geometric coordinates}
\end{definition}

\begin{remark}[Standard angular coordinates in the region $r < r_0$]
 We may assume that the coordinates $\vartheta^A$ are the standard coordinates $\{\varphi, \theta\}$ in the region $r < r_0$. In particular, this yields the expressions
\begin{equation}
 \begin{split}
  [L, X_A] &= 0 \\
  [\Lbar, X_A] &= 0 \\
  m(X_A, X_B) &= r^2 \mathring{\gamma}_{AB}
 \end{split}
\end{equation}
in the region $r < r_0$, where $\mathring{\gamma}$ is the standard round metric on the unit sphere.
\end{remark}

\begin{definition}[The reference metric on the spheres]
We shall also make use of a reference metric $\mathring{\gamma}$ on the spheres $S_{\tau,r}$. We construct $\mathring{\gamma}$ as follows: for all $r > 0$, there is a natural identification between the sphere $S_{\tau, r}$ and the sphere $\mathbb{S}^2$ (considered as a subset of $\mathbb{R}^3$), defined by the map
\begin{equation}
 \begin{split}
  \pi_{\mathbb{S}^2} &\, : \, S_{\tau, r} \rightarrow \mathbb{S}^2 \subset \mathbb{R}^3 \\
  &(x^0, x^1, x^2, x^3) \mapsto (r^{-1}x^1, r^{-1}x^2, r^{-1}x^3)
 \end{split}
\end{equation}
Note that, if $(x^0, x^1, x^2, x^3) \in S_{\tau,r}$ then we have 
\begin{equation*}
 (x^1)^2 + (x^2)^2 + (x^3)^2 = r^2
\end{equation*}
Note also that, in the region $r \leq r_0$, $x^0$ is constant on the sphere $S_{\tau, r}$, but this is not necessarily the case in the region $r > r_0$.

Now, $\mathbb{S}^2 \subset \mathbb{R}^3$ is equipped with the standard round metric $\gamma$. We define the metric $\mathring{\gamma}$ as the pullback of this metric $\gamma$ to the sphere $S_{\tau,r}$ under the maps given above. Specifically, we define
\begin{equation}
 \mathring{\gamma} := (\pi_{\mathbb{S}^2})^* (\gamma) 
\end{equation}

Similarly, we define $\dVol_{\mathbb{S}^2}$ as the pullback to the sphere $S_{\tau,r}$ of the standard volume form on the sphere $\mathbb{S}^2$ under the map $\pi_{\mathbb{S}^2}$.

\end{definition}

\subsection{The radial vector field and a normalisation condition on the metric}

\begin{definition}[The vector field $R$]
 We define the radial vector field
\begin{equation}
 R^\mu := (g^{-1})^{\mu\nu} \partial_\nu r
\end{equation}
\end{definition}

\begin{property}[The radial component of the metric]
\label{property radial normalisation}
 We shall assume the following condition on the radial component of the metric $g$:
\begin{equation}
 g(R, R) = g^{-1}(\upd r, \upd r) = 1 \quad \text{for } r \geq r_0
\label{equation radial normalisation}
\end{equation}
Note that this implies the following condition on the rectangular components of $g^{-1}$:
\begin{equation}
 \sum_{i,j = 1}^3 (g^{-1})^{ij} x^i x^j = r^2
\end{equation}
\end{property}

\begin{remark}
 Property \ref{property radial normalisation} can always be imposed by means of a conformal transformation. In fact, we will use such a transformation in order to impose this condition on the kinds of equations we shall treat later.
\end{remark}

\subsection{Basic identities involving the null frame and the geometric coordinates}
\begin{proposition}[Metric contractions of the null frame vector fields]
\label{proposition contractions null frame}
In the region $r \geq r_0$, we have the following expressions for the contractions of the null frame vector fields:
\begin{equation}
 \begin{split}
  g(L,L) &= 0 \\
  g(L, \Lbar) &= -2 \\
  g(L, X_A) &= 0 \\
  g(\Lbar, \Lbar) &= 0 \\
  g(\Lbar, X_A) &= 0 \\
  g(X_A, X_B) &= \slashed{g}_{AB}
 \end{split}
\end{equation}
\end{proposition}
\begin{proof}
 These expressions follow from the definitions together with the following calculations:
\begin{equation*}
 \begin{split}
  g(L, X_A) &= \mu g(L_{(\text{Geo})}, X_A) = \mu X_A u = 0 \\
  g(X_A, X_B) &= \slashed{g}(X_A, X_B) = \slashed{g}_{AB}
 \end{split}
\end{equation*}
\end{proof}

\begin{proposition}[Raising and lowering frame indices]
 Given a vector $V$, its upper and lower frame indices are related by
\begin{equation}
 \begin{split}
  V_L &= -2 V^{\Lbar} \\
  V_{\Lbar} &= -2 V^L \\
  V_{A} &= \slashed{g}_{AB} V^B
 \end{split}
\end{equation}
Similar relations hold for tensors of any rank.
\end{proposition}
\begin{proof}
 Expanding the vector $V$ in the null frame, we have
\begin{equation}
 V = V^L L + V^{\Lbar} \Lbar + V^A X_A
\end{equation}
Taking inner products with the null frame proves the proposition.
\end{proof}

\begin{proposition}[Expressing the radial vector field $R$ in the null frame]
\label{proposition R in null frame}
 \begin{equation}
\label{equation R in null frame}
  R = \frac{1}{2}(L - \Lbar)
 \end{equation}
\end{proposition}
\begin{proof}
 $R$ is $g$-orthogonal to the spheres $S_{\tau, r}$, since $g(R, X_A) = X_A(r) = 0$. Hence we can set
\begin{equation*}
 R = R^L L + R^{\Lbar}\Lbar
\end{equation*}
Now, we have
\begin{equation*}
 g(R, R) = 1 = -4R^L R^{\Lbar}
\end{equation*}
and also
\begin{equation*}
\begin{split}
 g(R, L) &= \mu g(R, L_{(\text{Geo})}) \\
 &= -\mu g^{-1}(\upd r, \upd u) \\
 &= 1 = -2R^{\Lbar}
\end{split}
\end{equation*}
Combining the previous two identities yields \eqref{equation R in null frame}
\end{proof}

\begin{proposition}[The action of the null frame on the geometric coordinates]
\label{proposition null frame geo coords}
 Define the $S_{\tau, r}$-tangent vector field $b$ by
\begin{equation}
 b := \left(\Lbar \vartheta^A\right) X_A
\end{equation}
Then
\begin{equation}
 \Lbar u = 2\mu^{-1}
\end{equation}
 and the action of the null frame on the geometric coordinates is given by
\begin{equation}
\label{equation null frame geo coords}
 \begin{split}
  L\tau &= 0 \\
  Lr &= 1 \\
  L\vartheta^A &= 0 \\
  \Lbar\tau &= 2\mu^{-1}  \\
  \Lbar r &= -1 \\
  \Lbar \vartheta^A &= b^A \\
  X_A \tau &= 0 \\
  X_A r &= 0 \\
  X_A \vartheta^B &= \delta_A^B
 \end{split}
\end{equation}
\end{proposition}

\begin{proof}
 The identities in \eqref{equation null frame geo coords} follow immediately from the definitions and the following two computations:
\begin{equation*}
\begin{split}
 Lr &= \mu L_{(\text{Geo})} r \\
 &= -\mu g^{-1}(\upd u, \upd r) = 1
\end{split}
\end{equation*}
Additionally, making use of proposition \ref{proposition R in null frame} we have
\begin{equation*}
\begin{split}
 \Lbar r &= Lr - 2R(r) \\
 &= 1 - 2g(R, R) = -1
\end{split}
\end{equation*}
\end{proof}

Proposition \ref{proposition null frame geo coords} immediately yields the following corollary:

\begin{corollary}[The null frame in terms of coordinate vector fields]
\label{corollary null frame coords}
 Define the coordinate induced vector fields
\begin{equation}
 \begin{split}
  \partial_{\tau} &:= \frac{\partial}{\partial \tau}\Big|_{r, \vartheta^1, \vartheta^2} \\
  \partial_r &:= \frac{\partial}{\partial r}\Big|_{\tau, \vartheta^1, \vartheta^2}
 \end{split}
\end{equation}
Then we have
\begin{equation}
 \begin{split}
  L &= \partial_r \\
  \Lbar &= 2\mu^{-1}\partial_\tau - \partial_r + b^A X_A
 \end{split}
\end{equation}
Conversely, the null frame can be expressed in terms of the coordinate induced vector fields:
\begin{equation}
 \begin{split}
  \partial_\tau &= \frac{1}{2}\mu \left( L + \Lbar - b^A X_A \right) \\
  \partial_r &= L
 \end{split}
\end{equation}

\end{corollary}

\begin{remark}[Geometric interpretation of $b$]
 Note that the vector field $b$ can be viewed as the projection of $\Lbar$ onto the spheres $S_{\tau,r}$ by the natural projection associated with the \emph{coordinate induced vector fields} (not the null frame!). That is, we can define a projection operator $p$ which acts on vector fields $V$ as
 \begin{equation*}
  p(V) := V - V(r) \partial_r - V(\tau) \partial_\tau
 \end{equation*}
 and then $p(\Lbar) = b$.
\end{remark}

\subsection{Commutators of the null frame vector fields}
Corollary \ref{corollary null frame coords} leads to the following expressions for the commutators of the null frame vector fields:
\begin{proposition}[Commutators of the null frame]
\label{proposition commutators null frame}
 The commutators of the null frame vector fields satisfy:
\begin{equation}
 \begin{split}
  [L, X_A] &= 0 \\
  [L, \Lbar] &= -\mu^{-1}(L\mu)(L + \Lbar) + \left(\mu^{-1}(L\mu)b^A + (Lb^A)\right)X_A \\
  [\Lbar, X_A] &= \mu^{-1} (X_A \mu) (L + \Lbar) - \left(\mu^{-1}(X_A\mu)b^B + (X_A b^B) \right)X_B \\
  [X_A, X_B] &= 0
 \end{split}
\end{equation}
\end{proposition}

\begin{remark}[$S_{\tau, r}$-covariant expressions]
The expression for $[\Lbar, X_A]$ can be made explicitly $S_{\tau, r}$-covariant by noting that, in our frame,
\begin{equation*}
 \slashed{\mathcal{L}}_A b = \mathcal{L}_A b = [X_A, b^B X_B] = (X_A b^B)X_B
\end{equation*}
where $\slashed{\mathcal{L}}$ denotes the Lie derivative restricted to the spheres. We have
\begin{equation*}
 [\Lbar, X_A] = \mu^{-1} (X_A \mu) (L + \Lbar) - \left(\mu^{-1}(X_A\mu)b^B + (\slashed{\mathcal{L}}_A b^B) \right)X_B 
\end{equation*}

\end{remark}

\subsection{The projection operators}
\label{subsection projection}
\begin{definition}[The projection operator onto the spheres $S_{\tau,r}$]
\label{definition projection}
We define the projection operator onto the spheres $S_{\tau, r}$ as
\begin{equation}
 \slashed{\Pi}_\nu^{\phantom{\nu}\mu} := \delta_\nu^\mu + \frac{1}{2}L_\nu \Lbar^\mu + \frac{1}{2}\Lbar_{\nu}L^\mu
\end{equation}
\end{definition}
Then we easily see that
\begin{equation}
\label{equation properties of projection}
 \begin{split}
  \slashed{\Pi}_\nu^{\phantom{\nu}\mu} L^{\nu} &= 0 \\
  \slashed{\Pi}_\nu^{\phantom{\nu}\mu} \Lbar^{\nu} &= 0 \\
  \slashed{\Pi}_\nu^{\phantom{\nu}\mu} (X_A)^{\nu} &= (X_A)^\mu \\
  \slashed{\Pi}_\mu^{\phantom{\mu}\rho}\slashed{\Pi}_\nu^{\phantom{\nu}\sigma} g_{\rho \sigma} &= \slashed{g}_{\mu\nu}
\end{split}
\end{equation}

\begin{remark}[Schematic notation for projection operators]
 Given a tensor $\phi$ of indeterminate rank, we shall occasionally write $\slashed{\Pi}(\phi)$ to denote the tensor obtained by projecting all the indices of $\phi$ using the projection operator $\slashed{\Pi}$. For example, if $\phi$ is a one-form, then $(\slashed{\Pi}(\phi))_\mu = \Pi_\mu^{\phantom{\mu}\nu}\phi_\nu$.
\end{remark}

\begin{definition}[Products of projection operators]
 We define the notation
 \begin{equation}
  \slashed{\Pi}_{\mu_1 \ldots \mu_n}^{\nu_1 \ldots \nu_n} := \slashed{\Pi}_{\mu_1}^{\phantom{\mu_1}\nu_1}\slashed{\Pi}_{\mu_2}^{\phantom{\mu_2}\nu_2} \ldots \slashed{\Pi}_{\mu_n}^{\phantom{\mu_n}\nu_n}
 \end{equation}
\end{definition}

We can also define projected derivative operators:

\begin{definition}[The projected exterior derivative]
 Given a function $f$, we define
 \begin{equation}
  \slashed{\upd}_\mu f := \left(\slashed{\upd} f \right)_\mu := \slashed{\Pi}_\mu^{\phantom{\mu}\nu} \left(\upd f \right)_\nu
 \end{equation}
\end{definition}

\begin{definition}[Projected covariant derivatives]
 We define the projected covariant derivative operator $\slashed{\D}_Z$, for any one-form $Z$, by first taking the covariant derivative in the $Z$ direction, $\D_Z$, and then projecting the \emph{lower indices} of the resulting tensor using $\slashed{\Pi}$. For example, given a one-form $\phi_\mu$, we have
 \begin{equation*}
  \slashed{\D}_Z \phi_\mu := \Pi_\mu^{\phantom{\mu}\nu} \D_Z \phi_\nu
 \end{equation*}
 Note that, if $\phi = \phi_{\sigma_1 \ldots \sigma_n}$ is a rank $(0, n)$ $S_{\tau,r}$-tangent tensor, then $\slashed {\D}_A \phi = \slashed{\nabla}_A\phi$. On the other hand, $\slashed{\nabla}_L \phi \equiv 0$ but $\slashed{\D}_L \phi$ does not necessarily vanish.

 The operator $\slashed{\D}$ is also compatible with the metric in the sense that $\slashed{\D}_\mu g_{\nu\rho} = \Pi_{\nu}^{\phantom{\nu}\lambda} \Pi_{\rho}^{\phantom{\rho}\sigma} \D_\mu g_{\lambda \sigma} = 0$. Finally, note that for a scalar field $\phi$, $\slashed{\D}\phi = \D\phi = \upd\phi$.
 
 The operator $\slashed{\D}$ can be viewed as giving rise to a connection on the bundle of $S_{\tau,r}$-tangent tensor fields (see chapter \ref{chapter geometry of vector bundle})
\end{definition}

\begin{definition}[The projected wave operator]
 We also define the projected wave operator $\slashed{\Box}_g$, defined by the formula
 \begin{equation}
 \slashed{\Box}_g \phi := (g^{-1})^{\mu\nu} \slashed{\D}_\mu \slashed{\D}_\nu \phi
 \end{equation}
 i.e.\ we project each covariant derivative and then contract them using $g$. Note that this is not the same as \emph{first} applying the geometric wave operator $\Box_g$ and \emph{then} projecting using $\slashed{\Pi}$: this operator differs from the one we have defined by some first order terms. For example, for a one-form $\phi_\mu$, we have
\begin{equation}
 \slashed{\Box}_g \phi_\alpha = \slashed{\Pi}_\alpha^{\phantom{\alpha}\beta} (g^{-1})^{\mu\nu} \D_\mu \left( \slashed{\Pi}_{\beta}^{\phantom{\beta}\gamma} \D_\nu \phi_\gamma \right)
\end{equation}
which differs from $\slashed{\Pi}_{\alpha}^{\phantom{\alpha}\beta} \Box_g \phi_\beta$ by the first order term $\slashed{\Pi}_\alpha^{\phantom{\alpha}\beta} (g^{-1})^{\mu\nu} \left( \D_\mu  \slashed{\Pi}_{\beta}^{\phantom{\beta}\gamma} \right) \D_\nu \phi_\gamma$.

Note also that the operator $\slashed{\Box}_g$ is \emph{not} the same as applying the operator $(\slashed{g}^{-1})^{\mu\nu}\slashed{\D}_\mu \slashed{\D}_\nu$. Nor is this latter operator the same as the operator $\slashed{\Delta}$: in fact,
\begin{equation*}
\begin{split}
	(\slashed{g}^{-1})^{\mu\nu} \slashed{\D}_\mu \slashed{\D}_\nu \phi
	&= 
	(\slashed{g}^{-1})^{\mu\nu} \slashed{\D}_\mu \left( -\frac{1}{2} \Lbar_\nu \slashed{\D}_L  \phi - \frac{1}{2} L_\nu \slashed{\D}_{\Lbar} \phi + \slashed{\nabla}_\nu \phi \right)
	\\
	&= \slashed{\Delta} \phi - \frac{1}{2} (\slashed{g}^{-1})^{\mu\nu} (\D_\mu \Lbar_\nu) \slashed{\D}_L \phi - \frac{1}{2} (\slashed{g}^{-1})^{\mu\nu} (\D_\mu L_\nu) \slashed{\D}_{\Lbar} \phi 
\end{split}
\end{equation*}

\end{definition}

\begin{proposition}[Frame components of the projected exterior differential of a function]
 For any function $f$, we have
\begin{equation}
 \begin{split}
  \slashed{\upd}_A f & = \left(\slashed{\upd} f \right)_A = X_A (f) \\
  \slashed{\upd}_L f & = \left(\slashed{\upd} f \right)_L = 0 \\
  \slashed{\upd}_{\Lbar} f & = \left(\slashed{\upd} f \right)_{\Lbar} = 0 \\
 \end{split}
\end{equation}
\end{proposition}
\begin{proof}
 This follows from the first three lines of \eqref{equation properties of projection}
\end{proof}

\begin{remark}[$S_{\tau, r}$--tensors and spacetime tensors]
 We slightly abuse notation in order to avoid distinguishing between $S_{\tau, r}$--tensors and tensors embedded in spacetime. For example, we use the notation $\slashed{g}$ to denote both the metric induced on the spheres $S_{\tau, r}$ (which is an $S_{\tau, r}$--tensor) and also to denote the restriction of the metric $g$ to the spheres, i.e.\ $\slashed{g}_{\mu\nu} = \slashed{\Pi}_\mu^{\phantom{\mu}\rho} \slashed{\Pi}_\nu^{\phantom{\nu}\sigma} g_{\rho\sigma}$, which is a spacetime tensor. The context will make clear which is meant.
\end{remark}
% 
% \begin{proposition}[The projection operator onto the hypersurface $\Sigma_{\tau}$]
%  We define the projection operator onto the hypersurface $\Sigma_{\tau}$ as follows:
%  \begin{equation}
%   \underline{\Pi}_\mu^{\phantom{\mu}\nu} := g_\mu^{\phantom{\mu}\nu} + n_\mu n^\nu
%  \end{equation}
%  where $n^\mu$ is the unit normal to the hypersurface $\Sigma_\tau$, i.e.\
%  \begin{equation}
%   n := \frac{1}{\sqrt{2}}\mu^{-\frac{1}{2}}(f'_{(\alpha)})^{-\frac{1}{2}} L - \frac{1}{2\sqrt{2}}\mu^{\frac{1}{2}}(f'_{(\alpha)})^{\frac{1}{2}}(L - \Lbar)
%  \end{equation}
% \end{proposition}

\subsection{Relations between the null frame and rectangular coordinates}

Since the rectangular derivatives span the tangent space to $\mathcal{M}$, we have the standard expressions for the null frame vectors in terms of rectangular derivatives:
\begin{equation}
 \begin{split}
  L &= L^a \partial_a \\
  \Lbar &= \Lbar^a \partial_a \\
  X_A &= (\slashed{\upd}_A x^a) \partial_a
 \end{split}
\end{equation}

In addition, we can express the rectangular derivatives in terms of the null frame vectors:
\begin{equation}
 \partial_a = -\left(\frac{1}{2}\Lbar_a \right) L - \left(\frac{1}{2}L_a \right) \Lbar + (X_A)_a (\slashed{g}^{-1})^{AB} X_B
\end{equation}

Sometimes we will consider sets of fields \emph{labelled} by rectangular indices. The obvious example of such a set is the set of scalar fields $h_{ab}$. In such cases it is often desirable to consider a related set of scalar fields, which are instead labelled by the null frame. An obvious way to do this is simply to project onto the null frame by using the rectangular components of the null frame vector fields. Unfortunately, the rectangular components of the angular vector fields $(X_A)^a$ are expected to grow as $r$. Hence, for example, the (scalar) field $h_{AB}$ will have different asymptotic behaviour to the (scalar) field $h_{ab}$. On the other hand, the fields $h_{AB}$ can be regarded as the components of an $S_{\tau,r}$ tangent tensor field, and this tensor field has the same asymptotics as the scalar field $h$ (when measured using the Riemannian metric on the spheres). To be more precise, we will introduce the \emph{rectangular frame angular one-forms}, which are a set of one-forms, labelled with rectangular indices, which can be used to form $S_{\tau,r}$-tangent tensor fields from sets of fields labelled by rectangular indices.

\begin{definition}[The rectangular frame angular one-forms]
 We define the rectangular frame angular one-forms as the set of one-forms, labelled by rectangular coordinates, given by the fields
 \begin{equation}
  \slashed{\Pi}_\mu^{\phantom{\mu}a} := (\slashed{g}^{-1})^{AB}(X_A)_\mu (X_B)^a
 \end{equation}
\end{definition}

\section{Schematic notation and norms}
We will use the following schematic notation in order to simplify many of the expressions in the following sections.

\begin{definition}[Schematic notation for angular derivatives, and ``good'' and ``bad'' derivatives]
We will use the following notation for angular derivatives: 
\begin{equation}
 \slashed{\nabla}_\mu := \slashed{\Pi}_\mu^{\phantom{\mu}\nu} \slashed{\D}_\nu
\end{equation}
in other words, when we write $\slashed{\nabla}$ we mean the covariant derivative with respect to the metric $\slashed{g}$. For example, for a scalar field $\phi$, $\slashed{\nabla}\phi$ is the one-form $\slashed{\Pi}_\mu^{\phantom{\mu}\nu}\D_\nu \phi = \slashed{\upd}_\mu \phi$. To define norms of the angular derivatives, we take contractions using the inverse metric $\slashed{g}^{-1}$. So, for example, for a scalar field $\phi$ we define
\begin{equation}
 |\slashed{\nabla}\phi| := \sqrt{ (\slashed{g}^{-1})^{\mu\nu} (\slashed{\nabla}_\mu \phi)(\slashed{\nabla}_\nu \phi) }
\end{equation}

We will write ``good'' derivatives schematically as
\begin{equation}
 \bar{\partial} \in \{ L, X_1, X_2 \}
\end{equation}
so, for example, by $\bar{\partial} \phi$ we mean any one of the following: $L\phi$, $X_1 \phi$ or $X_2 \phi$. Similarly, we write
\begin{equation}
 |\bar{\partial} \phi| := \sqrt{ |L\phi|^2 + |\slashed{\nabla}\phi|^2 }
\end{equation}
Likewise, if $\phi$ is an $S_{\tau,r}$ tensor, then corresponding notation will be used, with vector fields replaced by projected covariant derivatives. For example, we have
\begin{equation}
 |\overline{\slashed{\D}} \phi| := \sqrt{ |\slashed{\D}_L \phi|^2 + (\slashed{g}^{-1})^{AB} \slashed{\D}_A \phi \cdot \slashed{\D}_B \phi } 
\end{equation}
Recall that, for $S_{\tau,r}$ tensors, contractions are taken with respect to the metric $\slashed{g}$.

Finally, general (and possibly ``bad'') frame derivatives of a field $\phi$ will be written as $\partial \phi$, i.e.\ $\partial \phi$ may be any of the following: $\Lbar\phi$, $L\phi$, $X_1\phi$ or $X_2\phi$. For scalar fields $\phi$, we write
\begin{equation}
 |\partial \phi| := \sqrt{ |\Lbar\phi|^2 + |L\phi|^2 + |\slashed{\nabla}\phi|^2 }
\end{equation}
and, if $\phi$ is an $S_{\tau,r}$-tensor
\begin{equation}
 |\slashed{\D} \phi| := \sqrt{ |\slashed{\D}_{\Lbar}\phi|^2 + |\slashed{\D}_L\phi|^2 + (\slashed{g}^{-1})^{AB} \slashed{\D}_A \phi \cdot \slashed{\D}_B \phi  }
\end{equation}
Likewise, we can use our frame fields to define a norm of higher rank tensors in the obvious way. For example, for a symmetric $(0,2)$-tensor $h_{\mu\nu}$, we write
\begin{equation}
 |h| := |h_{LL}| + 2|h_{L\Lbar}| + 2|\slashed{h}_{L}| + |h_{\Lbar\Lbar}| + 2|\slashed{h}_{\Lbar}| + |\slashed{h}|
\end{equation}
where we have defined the $S_{\tau,r}$-tangent tensor fields
\begin{equation*}
 \begin{split}
    (\slashed{h}_L)_\mu &:= \slashed{\Pi}_\mu^{\phantom{\mu}\nu} L^\rho h_{\rho \nu} \\
    (\slashed{h}_{\Lbar})_\mu &:= \slashed{\Pi}_\mu^{\phantom{\mu}\nu} \Lbar^\rho h_{\rho \nu} \\
    (\slashed{h})_{\mu\nu} &:= \slashed{\Pi}_\mu^{\phantom{\mu}\rho} \slashed{\Pi}_\nu^{\phantom{\nu}\sigma} h_{\rho \sigma} \\
 \end{split}
\end{equation*}

\end{definition}

\begin{remark}
 One must be careful to distinguish between the schematic notation for general derivatives $\partial \phi$ and the notation we are using for rectangular derivatives $\partial_a \phi$. The schematic notation does not have indices, whereas the notation for the rectangular derivatives includes lower case Latin indices. Note that schematically written quantities such as $|\partial \phi|$ depend on the null frame, whereas the rectangular derivatives do not. It will turn out, however, that $|\partial \phi| \sim \sum_a |\partial_a \phi|$.
\end{remark}

\begin{definition}[Norms of $S_{\tau, r}$--tensors]
 For an $S_{\tau, r}$--tensor, say $\zeta = \zeta^A X_A$, norms are taken by contracting using the metric $\slashed{g}$ and its inverse. For example, we have
\begin{equation*}
 |\zeta| := \sqrt{ \zeta^A \zeta^B \slashed{g}_{AB} }
\end{equation*}

\end{definition}

\begin{definition}[The frame components and frame norm of fields labelled by rectangular indices]
 Let $h_{ab}$ be a set of scalar fields labelled by a symmetric pair of rectangular indices. Then we define the set of frame components of $h$:
 \begin{equation}
  (h)_{(\text{frame})} := \left\{ h_{LL} \, , \, h_{L\Lbar} \, , \, L^a \slashed{\Pi}_\mu^{\phantom{\mu}b} h_{ab} \, , \, h_{\Lbar \Lbar} \, , \, \Lbar^a \slashed{\Pi}_\mu^{\phantom{\mu}b} h_{ab} \, , \, \slashed{\Pi}_\mu^{\phantom{\mu}a} \slashed{\Pi}_\nu^{\phantom{\nu}b} h_{ab} \right\}
 \end{equation}
 Note that $(h)_{(\text{frame})}$ contains scalar fields (e.g.\ $h_{LL}$), $S_{\tau,r}$-tangent one forms (e.g. $L^a \slashed{\Pi}_\mu^{\phantom{\mu}b} h_{ab}$) and a symmetric $S_{\tau,r}$-tangent tensor field (the field $\slashed{\Pi}_\mu^{\phantom{\mu}a} \slashed{\Pi}_\nu^{\phantom{\nu}b} h_{ab}$).
 This notation can be generalised to deal with sets of $S_{\tau,r}$-tangent tensor fields; for example, consider the fields $\slashed{\nabla} h_{ab}$. Then we define
 \begin{equation*}
 (\slashed{\nabla}h)_{(\text{frame})} := \left\{ L^a L^b \nabla_\mu h_{ab} \, , \, L^a \Lbar^b \slashed{\nabla}_\mu h_{ab} \, , \, L^a \slashed{\Pi}_\nu^{\phantom{\nu}b} \slashed{\nabla}_\mu h_{ab} \, , \, \Lbar^a \Lbar^b \slashed{\nabla}_\mu h_{ab} \, , \, \Lbar^a \slashed{\Pi}_\nu^{\phantom{\nu}b}\slashed{\nabla}_\mu h_{ab} \, , \, \slashed{\Pi}_\nu^{\phantom{\nu}a} \slashed{\Pi}_\rho^{\phantom{\rho}b} \slashed{\nabla}_\mu h_{ab} \right\}
 \end{equation*}
 
 We also define the ``frame norm'' of $h$:
 \begin{equation}
  |h|_{(\text{frame})} := \sqrt{ \sum_{X \in (h)_{(\text{frame})}} |X|^2 }
 \end{equation}
 
 Note in particular that, when combining this notation with the other notation defined above for the norms of tensor fields, the transformation from the rectangular frame to the null frame is to be carried out \emph{after} derivatives of the fields labelled by rectangular indices are taken. For example, we have
 \begin{equation*}
  |\partial h|_{(\text{frame})} := \sqrt{ |L h|_{(frame)}^2 + |\Lbar h|_{(frame)}^2 + |\slashed{\nabla} h|_{(frame)}^2 }
 \end{equation*} 
 
\end{definition}

\begin{definition}[The rectangular components and rectangular norm]
 Similarly to the definitions given above, if we have a set of scalar fields labelled by rectangular indices, say $h_{ab}$, then we can write schematically
 \begin{equation}
 h_{(\text{rect})} := \left\{ h_{ab} \ \big| \ a,b \in \{0,1,2,3\} \right\}
 \end{equation}
 and we can define the rectangular norm:
 \begin{equation}
  |h|_{(\text{rect})} := \sqrt{ \sum_{a,b} |h_{ab}|^2 }
 \end{equation}
\end{definition}

\begin{definition}[The notation $\slashed{h}$]
	Sometimes, we will use the notation $\slashed{h}$ or $\slashed{h}_X$ to denote the angular components of the field $h$, which is normally labelled by rectangular indices. Specifically, we define
	\begin{equation*}
	\slashed{h}_{\mu\nu} := \slashed{\Pi}_\mu^{\phantom{\mu}a} \slashed{\Pi}_\nu^{\phantom{\nu}b} h_{ab}
	\end{equation*}
	and, for any vector field $X$ with rectangular components $X^a$, we define
	\begin{equation*}
	\slashed{h}_{X\mu} := X^a \slashed{\Pi}_\mu^{\phantom{\mu}b} h_{ab}
	\end{equation*}
	We define similar notation for derivatives. For example, if $X$ is a vector field, then we define
	\begin{equation*}
	(X\slashed{h})_{\mu\nu} := \slashed{\Pi}_\mu^{\phantom{\mu}a} \slashed{\Pi}_\nu^{\phantom{\nu}b}  (X h_{ab})
	\end{equation*}
	
\end{definition}

%\begin{definition}[Frame components of tensors]
% Given a tensor field, say $h_{\mu}$, we schematically write $(h)_{(\text{frame})}$ to denote any possible contraction of $h$ with the frame fields, i.e.\ $(h)_{(\text{frame})}$ can mean any of the following: $h_\mu L^\mu$, $h_\mu \Lbar^\mu$, $h_\mu (X_1)^\mu$ or $h_\mu (X_2)^\mu$. When taking norms, any angular indices are contracted using the metric $\slashed{g}$ and its inverse. For example, we write
%\begin{equation}
% |h|_{(\text{frame})} := \sqrt{ |h_L|^2 + |h_{\Lbar}|^2 + (\slashed{g}^{-1})^{AB}h_A h_B }
%\end{equation}
%\end{definition}
%
%\begin{remark}
% We can combine the above schematic notation for derivatives and frame components, for example, for a tensor $h_\mu$ we can write
%\begin{equation*}
% \begin{split}
%  \left( |\bar{\partial} h|_{(\text{frame})} \right)^2 &= (|Lh|_L)^2 + (|Lh|_{\Lbar})^2 + (\slashed{g}^{-1})^{AB} (Lh)_A (Lh)_B + (\slashed{g}^{-1})^{AB}(X_A h)_L (X_B h)_L  \\
%  &\phantom{=} + (\slashed{g}^{-1})^{AB} (X_A h)_{\Lbar}(X_B h)_{\Lbar} + (\slashed{g}^{-1})^{AB} (\slashed{g}^{-1})^{CD} (X_A h)_C (X_B h)_D 
% \end{split}
%\end{equation*}
%
%\end{remark}

\begin{definition}[The notation $a \lesssim b$, $a \gtrsim b$ and $a \sim b$]
 Given two quantities $a$ and $b$, the notation $a \lesssim b$ means that there is a \emph{numerical} constant $C$ (independent of all variables which we consider) such that $a \leq C b$. Similarly, the notation $a \gtrsim b$ means that there is a numerical constant $c$ such that $a \geq c b$. Finally, the notation $a \sim b$ means that both $a \lesssim b$ and $a \gtrsim b$.

 We use these same notations with subscripts to donate the fact that the implicit constants can depend on the variables in the subscripts. For example, the notation $a \lesssim_p b$ means that there is a constant $C(p)$, depending on $p$ but independent of all other variables, such that $a \leq C(p) b$.
\end{definition}

\begin{definition}[Cut-off functions]
 Let us define, once and for all, a smooth cut-off function $\chi_0(r)$, which satisfies
\begin{equation}
 \label{equation cut off chi0}
 \chi_0(r) = \begin{cases} 1 \quad \text{ if} \quad r \geq 1 \\ 0 \quad \text{ if} \quad r \leq \frac{1}{2} \end{cases}
\end{equation}
Then, we define the scaled cut-off functions: for $R, R_1, R_2 > 0$ and $R_2 \geq R_1$,
\begin{equation}
\label{equation cut off functions}
 \begin{split}
  \chi_{(R)}(r) &:= \chi_0 (R^{-1}r) \\
  \chi_{(R_1, R_2)}(r) &:= \chi_{(R_1)}(r) \left( 1 - \chi_{(R_2)} \right)
 \end{split}
\end{equation}
\end{definition}

\subsection{Geometric differential operators}
We also define some geometric derivative operators for tensors on the spheres $S_{\tau, r}$:
\begin{definition}[The divergence of $S_{\tau,r}$ vectors and tensors]
 Let $Z = Z^A X_A$ be an $S_{\tau,r}$ vector, and let $Y = Y^{AB} X_A \otimes X_B$ be a symmetric tensor on the sphere $S_{\tau,r}$. Then we define the divergence of $Z$ and $Y$ as follows:
\begin{equation}
 \begin{split}
  \slashed{\Div} X := \slashed{\nabla}_A Z^A = X_A Z^A + \slashed{\Gamma}^A_{AB} Z^B \\
  (\slashed{\Div} Y)^A := \slashed{\nabla}_B Y^{AB} = X_B Y^{AB} + \slashed{\Gamma}^A_{BC} Y^{CB} + \slashed{\Gamma}^B_{BC} Y^{AC}
 \end{split}
\end{equation}
Note that these quantities are covariantly defined with respect to a change of basis for the tangent space of the sphere $S_{\tau,r}$, i.e.\ under such a change of basis $\slashed{\Div} X$ transforms as a scalar and $\slashed{\Div} Y$ transforms as a vector.
\end{definition}

\begin{definition}[The curl of $S_{\tau,r}$ vectors and tensors]
 Again, let $Z = Z^A X_A$ be an $S_{\tau,r}$ vector, and let $Y = Y^{AB} X_A \otimes X_B$ be a symmetric tensor on the sphere $S_{\tau,r}$. Let $\slashed{\varepsilon} = \slashed{\varepsilon}_{AB} \upd\vartheta^A \otimes \upd\vartheta^B$ be the volume form of sphere $S_{\tau,r}$ induced by the metric $\slashed{g}$, normalised by $\slashed{\varepsilon}_{12} = \sqrt{\det \slashed{g}_{AB}}$ (recall that the ordered set of vector fields $(X_1, X_2)$ are chosen to induce the same orientation as the standard Minkowski vector fields $(\partial_{\theta}, \partial_{\phi})$).

 Note that $\slashed{\varepsilon}^{AB}\slashed{\varepsilon}_{AB} = 2$, and also that $\slashed{\nabla}_A \slashed{\varepsilon}_{BC} = 0$. Moreover, note that any antisymmetric tensor on the sphere $S_{\tau, r}$ of rank 2 must be proportional to $\slashed{\varepsilon}$. We define the curl of the vector field $Z$ and the tensor field $Y$ as follows:
\begin{equation}
 \begin{split}
  \slashed{\nabla}_{[A} Z_{B]} &:= \frac{1}{2}(\slashed{\Curl}\, Z) \slashed{\varepsilon}_{AB} \\
  \Rightarrow \slashed{\Curl}\, Z &\phantom{:}= \slashed{\varepsilon}^{AB} \slashed{\nabla}_A Z_B  = \slashed{\varepsilon}^{AB} X_A Z_B \\
  {\slashed{\nabla}_{[A} Y_{B]}}^{C} &:= \frac{1}{2}(\slashed{\Curl}\, Y)^C \slashed{\varepsilon}_{AB} \\
  \Rightarrow (\slashed{\Curl}\, Y)^A &\phantom{:}= \slashed{\varepsilon}^{BC} \slashed{\nabla}_B Y_C^{\phantom{C}A} = \slashed{\varepsilon}^{BC} \left( X_B Y_C^{\phantom{C}A} + \slashed{\Gamma}^A_{BD} Y_C^{\phantom{C}D} \right)
 \end{split}
\end{equation}
 Note that these quantities transform covariantly (as a scalar and a vector, respectively) under \emph{orientation-preserving} transformations of the frame $(X_1, X_2)$.
\end{definition}

\section{Expressions for the metric}
In this section we collect together various expressions for the spacetime metric $g$.
\subsection{Basic expressions for the metric}
\label{subsection basic metric}

\begin{definition}[The tensors $h$ and $H$]
 We define the tensors $h$ and $H$ by the equations
\begin{equation}
 \begin{split}
  g_{\mu\nu} &:= m_{\mu\nu} + h_{\mu\nu} \\
  (g^{-1})^{\mu\nu} &:= (m^{-1})^{\mu\nu} + H^{\mu\nu}
 \end{split}
\end{equation}
\end{definition}
The tensors $h$ and $H$ are to be thought of as ``small'', in a way which shall be made precise later. Using $h$ and $H$ can also establish the following proposition:

\begin{proposition}[The metric and its inverse in rectangular coordinates]
 In rectangular coordinates the metric and its inverse are given by
\begin{equation}
 \begin{split}
  g &= -\upd t^2 + (\upd x^1)^2 + (\upd x^2)^2 + (\upd x^3)^2 + h_{ab} \upd x^a \upd x^b \\
  g^{-1} &= - (\partial_t)^2 + (\partial_{x^1})^2 + (\partial_{x^2 })^2 + (\partial_{x^3})^2 + H^{ab}(\partial_{x^a})(\partial_{x^b})
 \end{split}
\end{equation}
\end{proposition}

\subsection{The metric in geometric coordinates}
In the region $r < r_0$, we only need the basic expressions for the metric given in subsection \ref{subsection basic metric}. However, in the region $r \geq r_0$ we have the following expression for the spacetime metric:
\begin{proposition}
 In the region $r \geq r_0$, the metric is given by
\begin{equation}
\label{equation metric geo coords}
\begin{split}
 g =& -\mu^2 \upd \tau^2 - 2\mu \upd \tau \upd r + \slashed{g}_{AB}\left( \upd \vartheta^A - \frac{1}{2}\mu b^A \upd \tau\right)\left( \upd \vartheta^B - \frac{1}{2}\mu b^B \upd \tau\right)
\end{split}
\end{equation}
\end{proposition}
\begin{proof}
 Using proposition \ref{proposition contractions null frame} and corollary \ref{corollary null frame coords} we can calculate the following metric components:
\begin{equation}
 \begin{split}
  g(\partial_\tau, \partial_\tau) &= -\mu^2 + \frac{1}{4}\mu^2 b^A b^B \slashed{g}_{AB} \\
  g(\partial_\tau, \partial_r) &= -\mu \\
  g(\partial_\tau, X_A) &= -\frac{1}{2}\mu b^B \slashed{g}_{BA} \\
  g(\partial_r, \partial_r) &= 0 \\
  g(\partial_r, X_A) &= 0 \\
  g(X_A, X_B) &= \slashed{g}_{AB}
 \end{split}
\end{equation}
Combining these expressions leads to equation \eqref{equation metric geo coords}.
\end{proof}

\begin{remark}[Relation to Bondi coordinates]
 The geometric coordinates are, in some ways, similar to ``Bondi coordinates'' near infinity: as $r \rightarrow \infty$, $\tau \sim u$. Note, however, that the radial coordinate $r$ is not related in the standard way to the volume of the spheres $S_{\tau, r}$, but is instead simply the ``background'' radial coordinate associated with the rectangular coordinates. It turns out that we could exchange $r$ for the area radius of the spheres without changinve very much - in particular, the decay rates of various quantities would not be qualitatively different. A more serious problem is the factors of $\mu$ in the above expressions, which mean that the usual decay rates for asymptotically flat spacetimes near null infinity are \emph{not} respected, since we expect the behaviour $\mu \sim r^\epsilon$.

% [[Other remarks]]: if \mu is bounded then can transform to something more like real Bondi coordinates (e.g. for Einstein equations). If we ALSO have a better bound on \rho (which also follows from wave coordinate condition!) then these are actually Bondi coordinates after renormalising r and \tau
\end{remark}

\begin{proposition}[The volume form expressed in geometric coordinates]
\label{proposition volume form g}
In the region $r \geq r_0$, the volume form associated with the metric $g$ is given in terms of the geometric coordinates as
\begin{equation}
 \dVol_g = -\mu \sqrt{\det \slashed{g}} \upd \tau \wedge \upd r \wedge \upd \vartheta^1 \wedge \upd \vartheta^2
\end{equation}
where the determinant of $\slashed{g}$ is taken relative to the coordinates $\vartheta^1$, $\vartheta^2$. Alternatively, we regard the tensor $\slashed{g}$ as defining a metric on the spheres $S_{\tau,r}$; then the spacetime volume form is 
\begin{equation}
 \dVol_g = -\mu \upd \tau \wedge \upd r \wedge \dVol_{(\mathbb{S}^2, \slashed{g})}
\end{equation}
where $\dVol_{(\mathbb{S}^2, \slashed{g})}$ is the volume form on $\mathbb{S}^2$ equipped with the induced metric $\slashed{g}$.

\end{proposition}

\begin{proof}
 This follows from a direct computation, beginning with the form of the metric in geometric coordinates given equation \eqref{equation metric geo coords}. The choice of orientation is in accordance with the standard choices. Note that, if $f'_{(\alpha)} = 0$ then we have $G = \slashed{g}$, and the expressions above simplifies.
\end{proof}

\begin{definition}[The scalar $\Omega$]
 We now have two metrics on the spheres $S_{\tau,r}$, given by $\slashed{g}$ and $\mathring{\gamma}$ respectively. We define the scalar $\Omega$ as the square root of the ratio of their respective volume forms, i.e.\
\begin{equation*}
 \dVol_g = -\mu \upd \tau \wedge \upd r \wedge \dVol_{(\mathbb{S}^2, G)} = -\mu \Omega^2 \upd \tau \wedge \upd r \wedge \dVol_{(\mathbb{S}^2, \mathring{\gamma})} 
\end{equation*}
In other words, we define $\Omega$ by
\begin{equation*}
 \sqrt{\det \slashed{g}} = \Omega^2 \sqrt{\det \mathring{\gamma}}
\end{equation*}
Note that this defines a scalar density on the spheres.

\end{definition}

\subsection{The metric in the null frame}
\begin{proposition}[The metric in the null frame]
 In the region $r \geq r_0$, we have the following expressions for the metric and its inverse:
\begin{equation}
\label{equation metric frame components}
 \begin{split}
  g_{\mu\nu} &= -\frac{1}{2}L_\mu \Lbar_\nu - \frac{1}{2} \Lbar_\mu L_\nu + \slashed{g}_{\mu\nu} \\
  g^{\mu\nu} &= -\frac{1}{2}L^\mu \Lbar^\nu - \frac{1}{2} \Lbar^\mu L^\nu + (\slashed{g}^{-1})^{\mu\nu}
 \end{split}
\end{equation}
\end{proposition}
\begin{proof}
 To establish the expression for the metric, we contract the expression given above with each pair of vectors in the null frame $(L, \Lbar, X_1, X_2)$. For example, we have that
\begin{equation*}
 g(L, \Lbar) = -2
\end{equation*}
while on the other hand,
\begin{equation*}
\begin{split}
 \left(-\frac{1}{2}L_\mu \Lbar_\nu - \frac{1}{2} \Lbar_\mu L_\nu + \slashed{g}_{\mu\nu} \right) L^\mu \Lbar^\nu = -\frac{1}{2}\Lbar_\mu L^\mu L_\nu \Lbar^{\nu} = -\frac{1}{2} \left(g(L, \Lbar)\right)^2 = -2
\end{split}
\end{equation*}
In fact, the first line of \eqref{equation metric frame components} holds for all pairs of frame vectors. Hence, as long as the null frame spans the tangent space, the first line of equation \eqref{equation metric frame components} is true. The second line follows by raising indices using the inverse metric $g^{-1}$.
\end{proof}

\subsection{Christoffel symbols}

We also define the rectangular Christoffel symbols:
\begin{definition}[Rectangular Christoffel symbols]
\label{definition Christoffel}
 The rectangular Christoffel symbols are defined by
\begin{equation}
 \D_a \partial_b := \Gamma^c_{ab} \partial_c
\end{equation}
\end{definition}

\begin{proposition}[Rectangular Christoffel symbols in terms of $h$]
\label{proposition Christoffel symbols in terms of h}
 The rectangular Christoffel symbols can be defined in terms of rectangular derivatives of $h$ as follows:
\begin{equation}
 \Gamma^c_{ab} = \frac{1}{2} (g^{-1})^{cd}\left( \partial_a h_{bd} + \partial_b h_{ad} - \partial_d h_{ab} \right)
\end{equation}
\end{proposition}
\begin{proof}
 This proposition follows from the standard expression for the Christoffel symbols in a coordinate induced basis, together with the fact that the coefficients of the Minkowski metric are constants in rectangular coordinates, i.e.\ $\partial_a m_{bc} = 0$
\end{proof}

\section{Second fundamental forms}
In this section we define the second fundamental forms and provide alternative expressions for them.

\begin{definition}[The second fundamental forms $\chi$ and $\chibar$]
 We define the second fundamental forms $\chi$ and $\chibar$ as follows:
\begin{equation}
 \begin{split}
  \chi := \frac{1}{2}\slashed{\mathcal{L}}_L g \\
  \chibar := \frac{1}{2}\slashed{\mathcal{L}}_{\Lbar} g
 \end{split}
\end{equation}
where $\slashed{\mathcal{L}}$ denotes the projection onto the spheres of the Lie derivative. To be clear, in abstract index notation we have
\begin{equation*}
 \begin{split}
  \chi_{\mu\nu} = \frac{1}{2} \slashed{\Pi}_\mu^{\phantom{\mu}\rho} \slashed{\Pi}_\nu^{\phantom{\nu}\sigma} (\mathcal{L}_L g)_{\rho\sigma} \\
  \chibar_{\mu\nu} = \frac{1}{2} \slashed{\Pi}_\mu^{\phantom{\mu}\rho} \slashed{\Pi}_\nu^{\phantom{\nu}\sigma} (\mathcal{L}_{\Lbar} g)_{\rho\sigma}
 \end{split}
\end{equation*}
\end{definition}

We also have alternative expressions for the second fundamental forms, which are often more useful:
\begin{proposition}[Alternative expressions for $\chi$ and $\chibar$]
 The only nonzero frame components of the second fundamental forms $\chi$ and $\chibar$ are given by
\begin{equation}
 \begin{split}
  \chi_{AB} &= g(\D_A L, X_B) \\
  \chibar_{AB} &= g(\D_A \Lbar, X_B)
 \end{split}
\end{equation}
\end{proposition}

\begin{proof}
 We prove the proposition for $\chi$; the $\chibar$ case is similar. First, note that, by equation \eqref{equation properties of projection}, the only frame components of $\chi$ that can be nonzero are $\chi_{AB}$. Now, we have
\begin{equation*}
 \begin{split}
  2\chi_{AB} &= \mathcal{L}_L g(X_A, X_B) \\
  &= L\left( g(X_A, X_B) \right) - g([L, X_A], X_B) - g(X_A, [L,X_B]) \\
  &= g(\D_L X_A, X_B) + g(X_A, \D_L X_B) - g([L, X_A], X_B) - g(X_A, [L,X_B]) \\
  &= g(\D_A L, X_B) + g(X_A, \D_B L) \\
  &= g(\D_A L, X_B) - g(\D_B X_A, L) \\
  &= g(\D_A L, X_B) - g(\D_A X_B, L) \\
  &= 2g(\D_A L, X_B)
 \end{split}
\end{equation*}
Where we have made use of the torsion-free property of the metric to set $[L, X_B] = \D_L X_B - \D_B L$, and we have used of the fact that $[X_A, X_B] = 0$.
\end{proof}

\begin{definition}[Trace and trace-free parts of second fundamental forms]
 We will also decompose the tensors $\chi$ and $\chibar$ into their trace and trace-free parts:
\begin{equation}
 \begin{split}
  \tr_{\slashed{g}}\chi := (\slashed{g}^{-1})^{AB} \chi_{AB} \\
  \hat{\chi}_{AB} := \chi_{AB} - \frac{1}{2} \slashed{g}_{AB} \tr_{\slashed{g}}\chi \\
  \tr_{\slashed{g}}\chibar := (\slashed{g}^{-1})^{AB} \chibar_{AB} \\
  \hat{\chibar}_{AB} := \chibar_{AB} - \frac{1}{2} \slashed{g}_{AB} \tr_{\slashed{g}}\chibar \\
 \end{split}
\end{equation}

\end{definition}

\section{The wave coordinate condition}
\label{section wave coordinate}

One application for the methods developed in this paper are the Einstein equations in ``wave coordinates'' or ``harmonic coordinates''. In this case, the rectangular coordinates themselves obey (scalar) wave equations, that is,
\begin{equation}
 \Box_g x^a = 0
\end{equation}
which is equivalent to the condition on the rectangular Christoffel symbols
\begin{equation}
 (g^{-1})^{bc}\Gamma^a_{bc} = 0
\end{equation}
This can be written in the null frame as
\begin{equation}
 -L^b (\Lbar h_{ab}) - \Lbar^b (L h_{ab}) - 2\slashed{g}^{bc}(\slashed{\nabla}_c h_{ab} ) + ( L^b \Lbar^c - \slashed{g}^{bc}) (\partial_a h_{bc}) = 0
\end{equation}
We can contract this equation with the rectangular components of the null frame $L^a$, $\Lbar^a$ and $\slashed{\Pi}_\mu^{\phantom{\mu}a}$. This leads to the expressions
\begin{equation}
\label{equation wave coordinate condition}
 \begin{split}
  0 &= (\Lbar h)_{LL} + 2\slashed{g}^{\mu\nu} (\slashed{\nabla}_\mu h)_{L\nu} + (\slashed{g}^{-1})^{\mu\nu}(Lh)_{\mu\nu} \\
  0 &= (\slashed{g}^{-1})^{\mu\nu}(\Lbar h)_{\mu\nu} + (Lh)_{\Lbar\Lbar} + 2\slashed{g}^{\mu\nu}(\slashed{\nabla}_\mu h)_{\Lbar\nu} \\
  0 &= \slashed{\Pi}_\mu^{\phantom{\mu}a} L^b (\Lbar h_{ab}) + \slashed{\Pi}_\mu^{\phantom{\mu}a} \Lbar^b (Lh_{ab}) + 2(\slashed{g}^{-1})^{\nu\rho}\slashed{\Pi}_\mu^{\phantom{\mu}a} \slashed{\Pi}_\nu^{\phantom{\nu}b}(\slashed{\nabla}_\rho h_{ab}) + (\slashed{g}^{-1})^{\nu\rho}(\slashed{\nabla}_\mu h)_{\nu\rho} - (\slashed{\nabla}_\mu h)_{L\Lbar} 
 \end{split}
\end{equation}

The Einstein equations, when supplemented with the wave coordinate condition, take the form of a set of nonlinear wave equations for the rectangular components of the metric perturbation $h_{ab}$. Moreover, these equations satisfy the ``weak null condition'' (see chapter \ref{chapter weak null structure}), and so the analysis of this paper can be applied to these equations.

What's more, for many of the calculations which follow, the Einstein equations in wave coordinates actually behave much \emph{better} than a more general wave equation with the weak null condition. The reason for this is that the wave coordinate condition itself implies better behaviour for certain derivatives of the metric, as we can see from the equations in \eqref{equation wave coordinate condition}. Specifically, the quantities
\begin{equation*}
 (\Lbar h)_{LL} \quad , \quad (\slashed{g}^{-1})^{\mu\nu}(\Lbar h)_{\mu\nu} \quad , \quad \slashed{\Pi}_\mu^{\phantom{\mu}a} L^b(\Lbar h_{ab})
\end{equation*}
actually behave like \emph{good} derivatives, despite containing the ``bad derivative'' $\Lbar$. This changes the behaviour of certain quantities dramatically, leading to different asymptotics for the metric components compared with the more general case we have in mind, and making certain computations much easier. This fact was used extensively in the work of Lindblad and Rodnianski \cite{Lindblad2004} and was crucial for their proof of the stability of Minkowksi space. 

One advantage of our approach is that the stability of Minkowksi space can be inferred \emph{without} making additional use of the wave coordinate condition, beyond using it to write the Einstein equations as a set of nonlinear wave equations with the weak null condition. Our approach would therefore allow the use of more general gauge conditions, where one does not expect to gain such good behaviour from the gauge conditions alone.

It is very important to note that we \emph{do not assume that the wave coordinate condition holds} in what follows. However, we will point out the occasions in which using this condition leads to improvements in the behaviour of certain quantities.

\chapter{Transport equations for eikonal quantities}
\label{chapter transport eikonal}
In this chapter we will derive transport equations in the outgoing null direction $L$ for quantities associated with the eikonal function $u$. Note that we shall eventually study a coupled wave equation/eikonal equation system, and so establishing appropriate \emph{a priori} bounds on the eikonal function is crucial for our global existence result. Rather than directly estimating $u$, we will instead estimate its derivatives: $\mu$ as well as the rectangular components of $L$. Recall that $\mu$ is related to the $\Lbar$ derivative of $u$, by $\Lbar u = 2\mu^{-1}$, while the rectangular components of $L$ are related to the rectangular derivatives of $u$ by $L^a = \mu (g^{-1})^{ab} \partial_b u$. 

We first establish the transport equation satisfied by the inverse foliation density $\mu$ in the region $r \geq r_0$ (recall that $\mu = 1$ in the region $r < r_0$).
\begin{proposition}[The transport equation for $\mu$]
\label{proposition transport mu}
 In the region $r \geq r_0$, the inverse foliation density $\mu$ satisfies the following transport equation along integral curves of $L$:
\begin{equation}
\label{equation transport mu}
 \mu^{-1}(L \mu) = \frac{1 - L^i L^i}{r} - \frac{1}{2}(L h)_{L\Lbar} + \frac{1}{4}(Lh)_{LL} + \frac{1}{4}(\Lbar h)_{LL} 
\end{equation}

\end{proposition}
\begin{proof}
 Recall that the vector field $L_{(\text{Geo})}$ is geodesic. Consequently, we find that
\begin{equation}
\label{equation DLL}
 \begin{split}
  \D_{L_{(\text{Geo})}} L_{(\text{Geo})} &= 0 \\
  \Rightarrow \D_L (\mu^{-1}L) &= 0 \\
  \Rightarrow \D_L L &= \mu^{-1} (L\mu) L
 \end{split}
\end{equation}
expanding the equation \eqref{equation DLL} in terms of the rectangular vector fields $\partial_a$, we find
\begin{equation}
\label{equation LL^a first}
 \begin{split}
  \mu^{-1} (L\mu) L &= L^a \D_a (L^b \partial_b) \\
  &= (L L^a)\partial_a + L^a L^b \D_a \partial_b \\
  &= (L L^a) \partial_a + L^a L^b \Gamma^c_{ab} \partial_c
 \end{split}
\end{equation}
Now, we wish to contract this (vector) equation with the 1-form $\upd r = r^{-1} x^i \upd x^i$. We first note a few identities which are obtained by acting on $r$ with the null frame vector fields:
\begin{equation*}
 \begin{split}
  Lr = 1 &\Rightarrow L^i \frac{x^i}{r} = 1 \\
  \Lbar r = -1 &\Rightarrow \Lbar^i \frac{x^i}{r} = -1 \\
  L(Lr) = 0 &\Rightarrow  (L L^i) \frac{x^i}{r} + \frac{L^i L^i}{r} - \frac{L^i x^i}{r^2}  = 0 \\
  &\Rightarrow (L L^i) \frac{x^i}{r} + \frac{L^i L^i}{r} - \frac{1}{r} = 0 \\
  (\slashed{g}^{-1})^{ai} \frac{x^i}{r} &= (\slashed{g}^{-1})^{ab} \partial_b r \\
  &= (g^{-1})^{cd} \slashed{\Pi}_c^{\phantom{c}a} \slashed{\Pi}_d^{\phantom{d}b} \partial_b r \\
  &= (g^{-1})^{cb} \slashed{\Pi}_c^{\phantom{c}a} (\slashed{\upd} r)_b \\
  &= (g^{-1})^{cb} \slashed{\Pi}_c^{\phantom{c}a} \left( -\frac{1}{2}(\slashed{\upd}_{\Lbar} r) L_b -\frac{1}{2}(\slashed{\upd}_L r) \Lbar_b + (\slashed{g}^{-1})^{AB} (\slashed{\upd}_A r) (X_B)_b \right) \\
  &= 0
 \end{split}
\end{equation*}
Returning to equation \eqref{equation LL^a first} and contracting with $\upd r$, we find
\begin{equation*}
 \begin{split}
  \mu^{-1} (L\mu) &=  (L L^i) \frac{x^i}{r} + L^a L^b \Gamma^i_{ab} \frac{x^i}{r}  \\
  &= \frac{1 - L^i L^i}{r} + \frac{1}{2} L^a L^b (g^{-1})^{id}\left( \partial_a h_{bd} + \partial_b h_{ad} - \partial_d h_{ab} \right)\frac{x^i}{r} \\
  &= \frac{1 - L^i L^i}{r} + \frac{1}{2} L^a L^b \left(-\frac{1}{2} L^i \Lbar^d - \frac{1}{2} \Lbar^i L^d + (\slashed{g}^{-1})^{id} \right)  \left( \partial_a h_{bd} + \partial_b h_{ad} - \partial_d h_{ab} \right) \frac{x^i}{r} \\
  &= \frac{1 - L^i L^i}{r} - \frac{1}{4}\left( L^a L^b \Lbar^d - L^a L^b L^d \right)  \left( \partial_a h_{bd} + \partial_b h_{ad} - \partial_d h_{ab} \right) \\
  &= \frac{1 - L^i L^i}{r} - \frac{1}{2}L^a \Lbar^b (L h_{ab}) + \frac{1}{4}L^a L^b (Lh_{ab}) + \frac{1}{4}L^a L^b (\Lbar h_{ab}) 
 \end{split}
\end{equation*}
proving the proposition.

\end{proof}

We can also derive a transport equation for the rectangular components of $L$:

\begin{proposition}[Transport equation for the rectangular coefficients of $L$]
\label{proposition transport La}
 In the region $r \geq r_0$, the rectangular coefficients of $L$ satisfy the following transport equation:
\begin{equation}
 \begin{split}
  LL^a &= \left( \frac{1 - L^i L^i}{r} + \frac{1}{4}(L h)_{LL} \right) L^a + \frac{1}{4} (L h)_{LL} \Lbar^a \\
  & \phantom{=} + (\slashed{g}^{-1})^{AB}\left(\frac{1}{2} (\slashed{\upd}_B h)_{LL} - (L h)_{BL} \right) (X_A)^a 
 \end{split}
\end{equation}
\end{proposition}

\begin{proof}
 We return to equation \eqref{equation LL^a first}, but this time we contract with $\upd x^a$ and substitute for $L(\log\mu)$ to find
\begin{equation}
 \begin{split}
  LL^a &= L(\log\mu) L^a - \Gamma^a_{bc}L^b L^c \\
  &= \left(\frac{1 - L^i L^i}{r} - \frac{1}{2}L^b \Lbar^c (L h_{bc}) + \frac{1}{4}L^b L^c (Lh_{bc}) + \frac{1}{4}L^b L^c (\Lbar h_{bc})\right) L^a \\
 & \phantom{=} - \frac{1}{2} \left(-\frac{1}{2}L^a\Lbar^d - \frac{1}{2}\Lbar^a L^d + (\slashed{g}^{-1})^{ad} \right)L^b L^c \left( \partial_b h_{cd} + \partial_c h_{bd} - \partial_d h_{bc} \right) \\
 & = \left( \frac{1 - L^i L^i}{r} + \frac{1}{4}L^b L^c (L h_{bc}) \right) L^a + \frac{1}{4} L^b L^c (L h_{bc}) \Lbar^a \\
 & \phantom{=} + (\slashed{g}^{-1})^{AB}\left(\frac{1}{2} L^b L^c (\slashed{\upd}_B h_{bc}) - (X_B)^b L^c (L h_{bc}) \right) (X_A)^a \\
 \end{split}
\end{equation}

\end{proof}

\begin{remark}[Redundant equations]
 In the previous proposition we established transport equations for all four rectangular components of $L$. However, these components are not all independent: they are constrained by the relation $L r = 1$, which we could use to express one of the rectangular components of $L$ in terms of the other three.
\end{remark}

\chapter{Null frame connection coefficients}
\label{chapter null frame connection coefficients}

In this chapter we define the null frame connection components and compute the covariant derivatives the vectors in our null frame. We then find relations between some of the null connection components and derivatives of the rectangular components of $h$, and compute derivatives of the rectangular components of the null frame, to complement the transport equations for the rectangular components of $L$ above.

\section{Null frame decomposition of the connection coefficients}

\begin{definition}[The scalar $\omega$]
\label{definition omega}
 We define the scalar quantity $\omega$ as
\begin{equation}
 \omega :=  \begin{cases} \mu^{-1} L\mu = L(\log\mu) \quad &\text{ if} \quad r\geq r_0 \\
 	0 \quad &\text{ if} \quad r < r_0
 	\end{cases}
\end{equation}
\end{definition}

\begin{definition}[The $S_{\tau, r}$-1-form $\zeta$]
 We define the $S_{\tau, r}$ 1-form $\zeta$ by its action on the frame vectors $X_A$:
\begin{equation}
 \zeta(X_A) := g(\D_A L, \Lbar)
\end{equation}
 We can extend $\zeta$ to a 1-form on $\mathcal{M}$ by setting $\zeta(L) = \zeta(\Lbar) = 0$.
\end{definition}

\begin{definition}[The angular connection coefficients]
 We define the angular connection coefficients $\slashed{\Gamma}$ via the following equation:
\begin{equation}
 \slashed{\nabla}_A X_B := \slashed{\Gamma}_{AB}^C X_C
\end{equation}

\end{definition}

\begin{definition}[The $S_{\tau,r}$-tangent vector fields $\sigma_{(A)}$]
 For $(A) \in \{1,2\}$, we define the $S_{\tau, r}$-tangent vector field $\sigma_{(A)}$ as follows:
 \begin{equation}
 \begin{split}
   \sigma_{(A)} &:= \slashed{\mathcal{L}}_{X^A} (\mu b) \\
   &= \slashed{\Pi}\left( \mathcal{L}_{X^A}(\mu b) \right)
  \end{split}
 \end{equation}
 Note that the index labelled as $(A)$ above should not be treated as a tensor index on the spheres $S_{\tau, r}$, i.e.\ it does not transform covariantly under a change of basis for the tangent space of the spheres.
\end{definition}

\begin{proposition}[Null frame connection coefficients]
\label{proposition null connection}
 In the region $r \geq r_0$, the covariant derivative may be decomposed in the null frame as follows:
\begin{equation}
 \begin{split}
  \D_L L &= \omega L \\
  \D_L \Lbar &= -\omega \Lbar - \zeta^A X_A \\
  \D_L X_A &= -\frac{1}{2}\zeta_A L + \chi_A^{\phantom{A}B}X_B \\
  \D_{\Lbar}L &= \omega L + \left(\zeta^A + 2(\slashed{\nabla}^A \log\mu) \right)X_A \\
  \D_{\Lbar}\Lbar &= -\omega \Lbar + 2 (\slashed{\nabla}^A \log\mu) X_A \\
  \D_{\Lbar} X_A &= \mu^{-1}(\slashed{\upd}_A \mu)L + \frac{1}{2}\left(\zeta_A + 2(\slashed{\nabla}_A \log\mu)\right) \Lbar + \left(\chibar_A^{\phantom{A}B} - \mu^{-1}\sigma_{(A)}^{\phantom{(A)}B} \right)X_B \\
  \D_A L &= -\frac{1}{2}\zeta_A L + \chi_A^{\phantom{A}B} X_B \\
  \D_A \Lbar &= \frac{1}{2}\zeta_A \Lbar + \chibar_A^{\phantom{A}B} X_B \\
  \D_A X_B &= \frac{1}{2}\chibar_{AB} L + \frac{1}{2}\chi_{AB} \Lbar + \slashed{\nabla}_A X_B
 \end{split}
\end{equation}
\end{proposition}

\begin{proof}
 Contract each connection coefficient with the null frame components, and make use of the definitions of $\omega$, $\zeta$, $\chi$ and $\chibar$, as well as the commutation identities \ref{proposition commutators null frame}. For example, to compute $\D_L \Lbar$ we start from
 \begin{equation*}
  \D_L \Lbar = -\frac{1}{2}g(\D_L \Lbar, \Lbar) L - \frac{1}{2}g(\D_L \Lbar, L) \Lbar + (\slashed{g}^{-1})^{AB} g(\D_L \Lbar, L) X_B
 \end{equation*}
 and we now find
 \begin{equation*}
  \begin{split}
   g(\D_L \Lbar, \Lbar) &= \frac{1}{2} L(g(\Lbar, \Lbar)) = 0 \\
   g(\D_L \Lbar, L) &= - g(\Lbar, \D_L L) = -\omega \\
   g(\D_L \Lbar, X_A) &= - g(\Lbar, \D_L X_A) = -g(\Lbar, \D_A L) = -\zeta_A
  \end{split}
 \end{equation*}
where we have made use of the identity $[L, X_A] = 0$ as well as the equation $\D_L L = \omega L$, which follows from the fact that $\D_{L_{(\text{Geo})}} L_{(\text{Geo})} = 0$.
 
\end{proof}

\begin{proposition}[Relation between $b^A$ and $\slashed{\upd}_A \mu$]
\label{proposition relation between b and mu}
 Using the connection coefficients, we can relate $L$ derivatives of $b^A$ to angular derivatives of $\mu$. We have
 \begin{equation}
  L b^A + (L\log \mu) b^A = -2 \zeta^A - 2\mu^{-1} (\slashed{g}^{-1})^{AB} (\slashed{\upd}_B \mu)
 \end{equation}
\end{proposition}
\begin{proof}
 We have
 \begin{equation*}
  \begin{split}
   X_A \mu &= -\frac{1}{2} \mu^2 [X_A, \Lbar] u \\
   &= -\mu [X_A, \Lbar]^{\Lbar} \\
   &= \frac{1}{2}\mu g([X_A, \Lbar], L) \\
   &= \frac{1}{2}\mu \left( g(\D_A \Lbar, L) - g(\D_{\Lbar} X_A, L) \right) \\
   &= \frac{1}{2}\mu \left( g(\D_{\Lbar} L, X_A) - g(\D_A L, \Lbar) \right) \\
   &= \frac{1}{2}\mu \left( -\zeta_A + g(\D_L \Lbar, X_A) + g([\Lbar, L] X_A) \right) \\
   &= \frac{1}{2}\mu \left( -2 \zeta_A - g([L, \Lbar], X_A) \right)
  \end{split}
 \end{equation*}
 On the other hand, we have that
 \begin{equation*}
  \begin{split}
   L b^A &= [L, \Lbar]\vartheta^A \\
   &= [L, \Lbar]^A \\
   &= (\slashed{g}^{-1})^{AB} g([L, \Lbar], X_B) - \frac{1}{2} g([L, \Lbar], L)b^A \\
   &= (\slashed{g}^{-1})^{AB} g([L, \Lbar], X_B) - \frac{1}{2} g(\D_L \Lbar, L)b^A \\
   &= (\slashed{g}^{-1})^{AB} g([L, \Lbar], X_B) - \mu^{-1} (L \mu) b^A
  \end{split}
 \end{equation*}
 Combining the previous two equations proves the proposition.
\end{proof}

\section{Recentred variables}
We define several ``recentred variables'' by subtracting off their background Minkowski values.
\begin{definition}[Recentred variables]
 We define the following variables:
\begin{equation}
 \begin{split}
  \big(\chi_{(\text{small})}\big)_{AB} &:= \chi_{AB} - \frac{1}{r}\slashed{g}_{AB} \\
  \tr_{\slashed{g}}\chi_{(\text{small})} &:= \tr_{\slashed{g}}\chi - \frac{2}{r} \\
  \big(\chibar_{(\text{small})}\big)_{AB} &:= \chibar_{AB} + \frac{1}{r}\slashed{g}_{AB} \\
  \tr_{\slashed{g}}\chibar_{(\text{small})} &:= \tr_{\slashed{g}}\chibar + \frac{2}{r} \\
  \mu_{(\text{small})} &:= \mu - 1 \\
  L^i_{(\text{small})} &:= L^i - \frac{x^i}{r} \\
  L^0_{(\text{small})} &:= L^0 - 1\\
  \Lbar^i_{(\text{small})} &:= \Lbar^i + \frac{x^i}{r} \\
  \Lbar^0_{(\text{small})} &:= \Lbar^0 - 1
 \end{split}
\end{equation}

\end{definition}

%\begin{definition}[Schematic notation for connection coefficients]
% Define the following set of variables:
%\begin{equation}
% \bm{\Gamma} := \left\{ \omega, r^{-1}\mu_{(\text{small})}, r^{-1} T(\mu), \tr_{\slashed{g}}\chi_{(\text{small})}, \tr_{\slashed{g}}\chibar_{(\text{small})}, \zeta, \slashed{\upd}\mu, \hat{\chi}, \hat{\chibar}, \slashed{\mathcal{L}}_B (\mu b^A) \right\}
%\end{equation}
%\end{definition}
%This definition will be useful when we do not require detailed knowledge of the structure of error terms. Note that the set $\bm{\Gamma}$ contains scalars, $S_{\tau,r}$-tangent covectors and $S_{(\tau,r)}$-tangent tensors. As such, when we write expressions such as $\bm{\Gamma} \cdot \bm{\Gamma}$, any $S_{\tau,r}$-tangent tensors are contracted in some unspecified way using the metric $\slashed{g}$.

\section{Schematic notation for error terms}

On many occasions the precise form of certain error terms is not important, and we can instead express the error terms schematically. Grouping together expressions which we expect to have similar behaviour, we define our notation below. Note that the sets of fields we define below are nested in such a way that we always include terms with strictly better behaviour alongside terms with the ``typical'' behaviour expected of a given grouping.

\begin{definition}[Error terms]
We define the error terms
\begin{equation}
\bm{\Gamma} := \left\{
\begin{array}{c}
(\partial h_{ab}) \, , \,
(\partial h)_{(\text{frame})} \, , \,
\zeta \, , \,
\slashed{\nabla}\log\mu \, , \,
\tr_{\slashed{g}}\chi_{(\text{small})} \, , \,
\hat{\chi} \, , \,
\omega \, , \,
\tr_{\slashed{g}}\chibar_{(\text{small})} \, , \,
\hat{\chibar} \, , \,
\\
(1+r)^{-1} L^a_{(\text{small})} \, , \,
(1+r)^{-1}\Lbar^a_{(\text{small})} 
\end{array} 
\right\}
\end{equation}
Note that we omit the quantities $\sigma_{(A)}$. This is because we will express all quantities as tensor fields on the spheres $S_{\tau,r}$, and in this case the term $\sigma_{(A)}$ cannot arise, since it does not transform as a tensor with respect to the index $(A)$.
\end{definition}

\begin{definition}[Good error terms]
We define the ``good'' error terms
\begin{equation}
\bm{\Gamma}_{(\text{good})} := \left\{ (\bar{\partial} h_{ab}) \, , \, (\bar{\partial} h)_{(\text{frame})} \, , \, \tr_{\slashed{g}}\chi_{(\text{small})} \, , \, \hat{\chi} \, , \, (1+r)^{-1} L^i_{(\text{small})} \right\}
\end{equation}
These are expected to decay faster toward null infinity compared to the general error terms $\bm{\Gamma}$.
\end{definition}

\begin{remark}[The regularity of the error terms]
 Certain error terms - in particular, $\tr_{\slashed{g}} \chi_{(\text{small})}$, $\hat{\chi}$, $\tr_{\slashed{g}} \chibar_{(\text{small})}$, $\hat{\chibar}$ and $\slashed{\nabla} \log \mu$ - have slightly worse regularity properties than might be desired. Specifically, although we would expect to be able to estimate these quantities in terms of the first derivatives of the metric, in fact we can only estimate these quantities by integrating equations involving the \emph{second} derivatives of the metric. This will cause issues when trying to estimate the highest order error terms.
 
 It turns out that the equations for $\tr_{\slashed{g}} \chi_{(\text{small})}$ and $\tr_{\slashed{g}} \chibar_{(\text{small})}$ can be modified so that these quantities are estimated in terms of the first derivatives of the metric\footnote{Note, however, that in the case of $\tr_{\slashed{g}} \chi_{(\text{small})}$ this results in a loss of \emph{decay} towards infinity.}. At the same time, we find that all of the \emph{good} derivatives of these quantities can be estimated without losing a derivative (i.e.\ the good derivatives can be estimated in terms of the curvature). It turns out that this is sufficient to close our estimates.

\end{remark}

We also define the following schematic notation:
\begin{definition}[Schematic notation for the frame fields]
	We will sometimes make use of the following notation: we use $X_{(\text{frame})}$ to stand for any one of the fields $L^a$, $\Lbar^a$ or $\slashed{\Pi}^a$. Note that there are a total of eight possible scalar fields, and four $S_{\tau,r}$-tangent one-forms, which $X_{(\text{frame})}$ can stand for.
	
	We will also use $\bar{X}_{(\text{frame})}$ to stand for any one of the three frame fields $L^i_{(\text{small})}$.
	
	Finally, we use the notation $X_{(\text{frame, small})}$ to stand for any of the frame fields 
	\begin{equation*}
	\begin{split}
	L^0_{(\text{small})} &:= L^0 - 1 \\
	L^i_{(\text{small})} &:= L^i - r^{-1} x^i \\
	\Lbar^0_{(\text{small})} &:= \Lbar^0 - 1 \\
	\Lbar^i_{(\text{small})} &:= \Lbar^i + r^{-1} x^i \\
	\end{split}
	\end{equation*}
	
\end{definition}

\section{Relations between the connection coefficients and derivatives of \texorpdfstring{$h$}{h}}
It is possible to deduce transport equations for the connection coefficients directly, but this approach will generally lead to regularity problems, i.e.\ we will \emph{lose a derivative} if we try to couple energy estimates for $h$ directly to the transport equations for the connection coefficients. In fact, we will have to overcome this obstacle when dealing with the second fundamental form $\chi$. However, the other connection coefficients can be directly related to derivatives of the rectangular components of $h$ and to the second fundamental form $\chi$, meaning that we only have to overcome this problem once!

\begin{proposition}[$\zeta$ in terms of derivatives of $h$]
\label{proposition zeta}
In the region $r \geq r_0$, we have the following two expressions for the connection coefficient $\zeta_\mu$:
\begin{equation}
 \begin{split}
  \zeta_\mu 
  &=
  2L^i_{(\text{small})} r^{-1} \slashed{\Pi}_\mu^{\phantom{\mu}i}
  - \frac{1}{2} (\Lbar \slashed{h})_{L\mu}
  + \frac{1}{2} (L \slashed{h})_{\Lbar\mu}
  - \frac{1}{2} (\slashed{\nabla}_\mu h)_{LL}
  + \frac{1}{2} (\slashed{\nabla}_\mu h)_{L\Lbar}
  \\
  &=
  2\Lbar^i_{(\text{small})} r^{-1} \slashed{\Pi}_\mu^{\phantom{\mu}i}
  - \frac{1}{2} (\Lbar \slashed{h})_{L\mu}
  + \frac{1}{2} (L \slashed{h})_{\Lbar\mu}
  - \frac{1}{2} (\slashed{\nabla}_\mu h)_{L\Lbar}
  + \frac{1}{2} (\slashed{\nabla}_\mu h)_{\Lbar\Lbar}
 \end{split}
\end{equation}
Schematically, the $S_{\tau,r}$-tangent one-form $\zeta_\mu$ satisfies
\begin{equation}
 |\zeta| \lesssim |\partial h|_{(\text{frame})} + \left| \frac{L^i \slashed{\upd} x^i}{r} \right|
\end{equation}

\end{proposition}

\begin{proof}
Expanding in rectangular coordinates:
\begin{equation}
 \label{equation DAL v1}
 \begin{split}
  -\frac{1}{2}\zeta_A L + \chi_A^{\phantom{A}B} X_B &= \D_A L \\
  &= (X_A)^a \D_a( L^b \partial_b) \\
  &= (\slashed{\upd}_A L^a) \partial_a + (X_A)^a L^b \Gamma^c_{ab} \partial_c
 \end{split}
\end{equation}
Now, we wish to contract this equation with the covector $\upd r$. Note that $X_A (r) = 0$, and also that $(\upd r)_c = \frac{1}{2}(L_c - \Lbar_c)$. Finally, we note that
\begin{equation*}
 \begin{split}
 L^i \frac{x^i}{r} &= 1 \\
 \Rightarrow (\slashed{\upd}_A L^i) \frac{x^i}{r} + L^i \frac{\slashed{\upd}_A x^i}{r} &= 0
 \end{split}
\end{equation*}
Returning to equation \eqref{equation DAL v1} and contracting with $\upd r$ we find
\begin{equation*}
 \begin{split}
  \zeta_A &= 2 L^i \frac{\slashed{\upd}_A x^i}{r} - (X_A)^a L^b (L_c - \Lbar_c) \Gamma^c_{ab} \\
  &= 2 L^i \frac{\slashed{\upd}_A x^i}{r} - \frac{1}{2} (X_A)^a L^b (L^c - \Lbar^c) \left( \partial_a h_{bc} + \partial_b h_{ac} - \partial_c h_{ab} \right) \\
  &= 2 L^i \frac{\slashed{\upd}_A x^i}{r} - \frac{1}{2}(\slashed{\upd}_A x^a) L^b (\Lbar h_{ab}) + \frac{1}{2} (\slashed{\upd}_A x^a) \Lbar^b (L h_{ab}) - \frac{1}{2}L^a (L^b - \Lbar^b) (\slashed{\upd}_A h_{ab})
 \end{split}
\end{equation*}

Alternatively, we can expand the equation $\D_A \Lbar = -\frac{1}{2}\zeta_A \Lbar + \chibar_A^{\phantom{A}B} X_B$ and contract with $\upd r$ to find
\begin{equation*}
 \zeta_A = 2 \Lbar^i \frac{\slashed{\upd}_A x^i}{r} - \frac{1}{2}(\slashed{\upd}_A x^a) L^b (\Lbar h_{ab}) + \frac{1}{2} (\slashed{\upd}_A x^a) \Lbar^b (L h_{ab}) - \frac{1}{2}\Lbar^a (L^b - \Lbar^b) (\slashed{\upd}_A h_{ab})
\end{equation*}
which proves the proposition.

\end{proof}

\begin{remark}
 The proposition above immediately leads to the following identity:
 \begin{equation*}
  2 L^i \frac{\slashed{\upd}_A x^i}{r} - \frac{1}{2}L^a (L^b - \Lbar^b) (\slashed{\upd}_A h_{ab}) = 2 \Lbar^i \frac{\slashed{\upd}_A x^i}{r} - \frac{1}{2}\Lbar^a (L^b - \Lbar^b) (\slashed{\upd}_A h_{ab})
 \end{equation*}
 In fact, this can be deduced from a special case of the following proposition, which uses the condition $g^{-1}(\upd r, \upd r) = 1$ to relate certain derivatives of the metric.
\end{remark}

\begin{remark}[Improved behaviour for $\zeta$ using the wave coordinate condition]
 Recall that the wave coordinate condition leads to improved behaviour of certain derivatives of the rectangular components of the metric, see section \ref{section wave coordinate}. In particular, the quantity $(\Lbar h)_{AL}$ behaves like a \emph{good} derivative, despite the presence of the $\Lbar$ derivative. Thus, if the wave coordinate condition holds, the connection component $\zeta$ behaves like a good derivative of $h$, but in the more general case that we are considering it behaves like an $\Lbar$ derivative.
\end{remark}

%  In the case of the Einstein equations in wave coordinates, we can make use of the wave coordinate condition to improve the behaviour of $\zeta$. Specifically, the wave coordinate condition written in the null frame is
% \begin{equation}
%  -L^b (\Lbar h_{ab}) - \Lbar^b (L h_{ab}) - 2\slashed{g}^{bc}(\slashed{\nabla}_c h_{ab} ) + ( L^b \Lbar^c - \slashed{g}^{bc}) (\partial_a h_{bc}) = 0
% \end{equation}
% We can contract this equation with the null frame, which leads to the equations
% \begin{equation}
%  \begin{split}
%   0 &= (\Lbar h)_{LL} + 2\slashed{g}^{AB} (X_A h)_{BL} + (\slashed{g}^{-1})^{AB}(Lh)_{AB} \\
%   0 &= (\slashed{g}^{-1})^{AB}(\Lbar h)_{AB} + (Lh)_{\Lbar\Lbar} + 2\slashed{g}^{AB}(X_A h)_{B\Lbar} + (Lh)_{\Lbar\Lbar} \\
%   0 &= (\Lbar h)_{LA} + (Lh)_{\Lbar A} + 2(\slashed{g}^{-1})^{BC}(X_B h)_{CA} + (\slashed{g}^{-1})^{BC}(X_A h)_{BC} - (X_A h)_{L\Lbar}
%  \end{split}
% \end{equation}
% 
% 
% 
% One consequence of these equations is that $(\Lbar h)_{LA}$

\begin{proposition}[Relations between first derivatives of the metric]
For any vector field $V$, in the region $r \geq r_0$ we have
\begin{equation}
 -\frac{1}{4}(\mathcal{L}_V g)(L - \Lbar \ ,\ L - \Lbar) + (L - \Lbar)(V(r)) = 0
\end{equation}
\end{proposition}
\begin{proof}
 Since $g^{-1}(\upd r , \upd r) = 1$, we take Lie derivatives in the $V$ direction to find
 \begin{equation*}
  \begin{split}
   0 &= (\mathcal{L}_V g^{-1})(\upd r, \upd r) + 2 g^{-1}(\mathcal{L}_V \upd r , \upd r) \\
   &= - (\mathcal{L}_V g)(R, R) + 2 g^{-1}(\upd \imath_V \upd r , \upd r) \\
   &= - (\mathcal{L}_V g)(R, R) + 2 \imath_R \upd \imath_V \upd r
  \end{split}
 \end{equation*}
where we have made use of Cartan's formula. Note that the identity in the previous remark now follows from setting $V = X_A$ and expanding in rectangular coordinates.
\end{proof}

\begin{proposition}[Relationship between $\chi$, $\chibar$ and derivatives of $h$]
\label{proposition chibar in terms of chi}
In the region $r \geq r_0$, the second fundamental form $\chibar$ can be expressed in terms of the second fundamental form $\chi$ and the first derivatives of $h$ as follows:
\begin{equation}
 \begin{split}
  \chibar_{AB} =& \chi_{AB} - \frac{2}{r}(\slashed{\upd}_A x^i)(\slashed{\upd}_B x^i) + \frac{1}{2}(\Lbar h)_{AB} - \frac{1}{2}(L h)_{AB} \\
  & + \frac{1}{2} (\slashed{\upd}_B h)_{LA} - \frac{1}{2}(\slashed{\upd}_B h)_{\Lbar A}   + \frac{1}{2} (\slashed{\upd}_A h)_{LB} - \frac{1}{2}(\slashed{\upd}_A h)_{\Lbar B}  
 \end{split}
\end{equation}
Schematically, we have
\begin{equation}
 \left|\chibar_{\mu\nu} + \frac{1}{r}\slashed{g}_{\mu\nu} \right| \lesssim \left| \chi_{\mu\nu} - \frac{1}{r}\slashed{g}_{\mu\nu} \right| + \frac{1}{r} \left| \slashed{g}_{\mu\nu} - (\slashed{\upd}_\mu x^i)(\slashed{\upd}_\nu x^i) \right| + |\partial h|_{(\text{frame})}
\end{equation}
 
\end{proposition}
\begin{proof}
 We recall that, using the radial vector field $R$, we have that
\begin{equation*}
 \begin{split}
  \Lbar &= L - 2R \\
  \Rightarrow \chibar_{AB} &= \chi_{AB} - 2g(\D_A R, X_B) \\
 &= \chi_{AB} - 2 (\D_A \upd r)\cdot X_B \\
 &= \chi_{AB} - 2 \D_A \left( \frac{x^i}{r} \upd x^i \right) \cdot X_B \\
 &= \chi_{AB} - \frac{2}{r} (\slashed{\upd}_A x^i)(\slashed{\upd}_B x^i) + \frac{2 x^i}{r} \Gamma^i_{ab} (X_A)^a (X_B)^b \\
 &= \chi_{AB} - \frac{2}{r} (\slashed{\upd}_A x^i)(\slashed{\upd}_B x^i) \\
 & \phantom{=} - \frac{x^i}{r}\left( -\frac{1}{2}L^i \Lbar^c - \frac{1}{2} \Lbar^i L^c + (\slashed{g}^{-1})^{ic} \right) (X_A)^a (X_B)^b (\partial_a h_{bc} + \partial_b h_{ac} - \partial_c h_{ab} ) \\
 &= \chi_{AB} - \frac{2}{r} (\slashed{\upd}_A x^i)(\slashed{\upd}_B x^i) + \frac{1}{2}(L^c - \Lbar^c)(X_A)^a (X_B)^b (\partial_a h_{bc} + \partial_b h_{ac} - \partial_c h_{ab} )
 \end{split}
\end{equation*}
 where we have used the properties of the covariant derivative to write $g(\D_A R, X_B) = (\D_A R^\flat)\cdot X_B = (\D_A \upd r)\cdot X_B$, and we have also used
\begin{equation*}
 0 = (\slashed{g}^{-1})^{Aa}(X_A r) = (\slashed{g}^{-1})^{Aa} (X_A)^i \frac{x^i}{r} = (\slashed{g}^{-1})^{ia}\frac{x^i}{r}
\end{equation*}

\end{proof}

\begin{remark}[Improved behaviour of $\tr_{\slashed{g}}\chibar$ using the wave coordinate condition]
Similarly to $\zeta$, if the wave coordinate condition holds then $\tr_{\slashed{g}}\chibar$ behaves like a \emph{good} derivative of $h$, whereas in the general case it behaves like an $\Lbar$ derivative. To see this, note that the wave coordinate condition leads to an improved estimate on $(\slashed{g}^{-1})^{AB} (\Lbar h)_{AB}$ (see section \ref{section wave coordinate}), which is precisely the ``bad derivative'' term appearing in $\tr_{\slashed{g}}\chibar$.
\end{remark}

\begin{proposition}[Transport equation for the metric on the spheres \texorpdfstring{$\slashed{g}$}{gslashed}]
 The induced metric on the spheres $\slashed{g}$ and its inverse satisfy the following evolution equation along the integral curves of $L$:
\begin{equation}
  \begin{split}
    L\left( r^{-2}\slashed{g}_{AB}\right) &= 2r^{-2}(\chi_{(\text{small})})_{AB} \\
    L\left( r^2 (\slashed{g}^{-1})^{AB} \right) &= -2 r^2 (\chi_{(\text{small})})^{AB}
  \end{split}
\end{equation}

\end{proposition}

\begin{proposition}[Transport equations for the vector fields $\sigma_{(A)}$]
 The vector fields $\sigma_{(A)}$ satisfy the following transport equations along the integral curves of $L$:
 \begin{equation}
  \slashed{\mathcal{L}}_L \sigma_{(A)} = -2 (\slashed{\nabla}_A \mu)\zeta - 2\mu \slashed{\nabla}_A \zeta - 2 \slashed{\nabla}_A\slashed{\nabla}^\sharp  \mu
 \end{equation}
\end{proposition}
\begin{proof}
 We first note that
 \begin{equation}
  \sigma_{(A)} = \slashed{\Pi} [X_A, \mu b] = [X_A, \mu b]
 \end{equation}
 where we have used the fact that $b$ is $S_{[\tau,r]}$-tangent. Now, we have
 \begin{equation}
 \begin{split}
  \slashed{\mathcal{L}}_L \sigma_{(A)} &= \slashed{\Pi}\left( \left[L, [X_A, \mu b] \right] \right) \\
  &= \slashed{\Pi} \left( \left[ X_A, [L, \mu b] \right] + \left[ \mu b, [X_A, L] \right] \right) \\
  &= \slashed{\mathcal{L}}_A \left( \mathcal{L}_L (\mu b) \right)
 \end{split}
 \end{equation}
where in the second line we have used the Jacobi identity, and in the third line we have used the fact that $[X_A, L] = 0$. Now, we make use of proposition \ref{proposition relation between b and mu}:
\begin{equation}
 \begin{split}
   \left(\slashed{\mathcal{L}}_L \sigma_{(A)}\right)^\alpha &= \slashed{\Pi}_\beta^{\phantom{\beta}\alpha}\mathcal{L}_A \left( -2\mu \zeta^\beta -2 \slashed{\nabla}^\beta \mu \right) \\
   &= -2 (\slashed{\nabla}_A \mu) \zeta^\alpha -2\mu (\slashed{\nabla}_A \zeta)^\alpha - 2\slashed{\nabla}_A \slashed{\nabla}^\alpha \mu
  \end{split}
\end{equation}
 
\end{proof}

\section{Derivatives of the rectangular components of the null frame}
\label{section derivatives of rectangular components of null frame}
Recall that in chapter \ref{chapter transport eikonal} we derived transport equations for the rectangular components of $L$. Now that we have defined the null frame connection coefficients, we can also find expressions for some of the other derivatives of the rectangular components of the null frame vector fields. In addition, we can derive expressions for derivatives of the rectangular components of some of the recentred variables $L^i_{(\text{small})}$ and $\Lbar^i_{(\text{small})}$.

\begin{proposition}[Transport equations for the rectangular components of the null frame]
\label{proposition transport rectangular}
The rectangular components null frame vector fields $L^a$ and $\Lbar^b$, and the $S_{\tau,r}$-tangent one forms $\slashed{\Pi}_\mu^{\phantom{\mu}a}$ satisfy the following system of transport equations along the integral curves of $L$:
\begin{equation}
 \begin{split}
  LL^a &= \left( \frac{1 - L^i L^i}{r} + \frac{1}{4}(Lh)_{LL} \right)L^a
  			+ \frac{1}{4}\left((Lh)_{LL}\right) \Lbar^a \\
  		 	&\phantom{=} + (\slashed{g}^{-1})^{\mu\nu} \left( \frac{1}{2}(\slashed{\nabla}_\mu h)_{LL} - \slashed{\Pi}_\mu^{\phantom{\mu}b} L^c (Lh_{bc}) \right)\slashed{\Pi}_\nu^{\phantom{\nu}a}
  \\ \\
  L\Lbar^a &= \frac{1}{4}\left( (L h)_{\Lbar\Lbar} \right) L^a
  			+ \left( \frac{ L^i L^i - 1}{r} + \frac{1}{2}(L h)_{L\Lbar} - \frac{1}{4}(L h)_{LL} \right) \Lbar^a\\
  			&\phantom{=} - \frac{1}{2}(\slashed{g}^{-1})^{\mu\nu}\left(4L^i \frac{\slashed{\nabla}_\mu x^i}{r} + 2\slashed{\Pi}_\mu^{\phantom{\mu}b} L^c (Lh_{bc}) - (\slashed{\nabla}_\mu h)_{LL} \right) \slashed{\Pi}_\nu^{\phantom{\nu}a}
   \\ \\
  \slashed{\D}_L \left( \slashed{\Pi}_\mu^{\phantom{\mu}a} \right)
   			&= \left( -L^i \frac{\slashed{\nabla}_\mu x^i}{r} + \frac{1}{4}(\slashed{\nabla}_\mu h)_{LL} \right) L^a
   			+ \frac{1}{4}(\slashed{\nabla}_\mu h)_{LL} \Lbar^a \\
   			&\phantom{=} + \frac{1}{2}(\slashed{g}^{-1})^{\nu\rho} \left( L^b \slashed{\Pi}_\mu^{\phantom{\mu}c} (\slashed{\nabla}_\nu h_{bc})
   			- L^b \slashed{\Pi}_\nu^{\phantom{\nu}c} (\slashed{\nabla}_\mu h_{bc})
   			- \slashed{\Pi}_\mu^{\phantom{\mu}b}\slashed{\Pi}_\nu^{\phantom{\nu}c} (Lh_{bc}) \right) \slashed{\Pi}_\rho^{\phantom{\rho}a}
 \end{split}
\end{equation}

In addition, the fields $rL^i_{(\text{small})}$ satisfy the transport equations
\begin{equation}
 \begin{split}
    L(r L^i_{(\text{small})} ) &= \frac{1}{4}\left((Lh)_{LL}\right) (rL^i_{(\text{small})})  -\left( L^j_{(\text{small})} L^j_{(\text{small})} \right) L^i + \frac{1}{4}\left(r(L h)_{LL} \right) \Lbar^i_{(\text{small})} \\
     &\phantom{=} + (\slashed{g}^{-1})^{\mu\nu} \left( \frac{1}{2} r(\slashed{\nabla}_\mu h)_{LL} - L^a \slashed{\Pi}_\mu^{\phantom{\mu}b} r(Lh_{ab}) \right) \slashed{\Pi}_\nu^{\phantom{\nu}i} 
  \end{split}
\end{equation}

Schematically, these equations can be written as
\begin{equation}
 \begin{split}
 &\slashed{\D}_L \begin{pmatrix}
                 L^a \\
		 \Lbar^a \\
		 \slashed{\Pi}^a
                \end{pmatrix}
  = \left( (\bar{\partial} h)_{(\text{frame})} + \frac{1}{r^3} (rL^i_{(\text{small})})(rL^i_{(\text{small})}) + \frac{1}{r^2}\slashed{\Pi}^i (rL^i_{(\text{small})}) \right) \cdot 
		\begin{pmatrix}
                 L^a \\
		 \Lbar^a \\
		 \slashed{\Pi}^a
                \end{pmatrix} \\
  & L(rL^i_{(\text{small})}) = (\bar{\partial} h)_{(\text{frame})}(rL^i_{(\text{small})}) + \left( r (\bar{\partial} h)_{(\text{frame})} + \frac{1}{r^2}(rL^j_{(\text{small})})(rL^j_{(\text{small})}) \right)\cdot
		\begin{pmatrix}
                 L^i \\
		 \Lbar^i \\
		 \slashed{\Pi}^i
                \end{pmatrix}
  \end{split}
\end{equation}

\end{proposition}
\begin{proof}
 Each of the expressions above is deduced from expanding the corresponding expressions for the covariant derivatives in rectangular coordinates. The most technically involved calculation is for $\slashed{\D}_L \slashed{\Pi}_\mu^{\phantom{\mu}a}$, for which we can first computed $\D_L \slashed{\Pi}_\mu^{\phantom{\mu}a} = \D_L \left( (\slashed{g}^{-1})^{AB}(X_A)_\mu (X_B)^a\right)$. The only term we now need to calculate is $L (X_A)^a$, which we do as follows:
\begin{equation*}
 \begin{split}
  \D_L X_A &= \D_L \left( (X_A)^a \partial_a \right) = \left( L(X_A)^a + (X_A)^b L^c \Gamma^a_{bc} \right) \partial_a \\
  \Rightarrow L(X_A)^a &= -\frac{1}{2}\zeta_A L^a + \chi_A^{\phantom{A}B}(X_B)^a - (X_A)^b L^c \Gamma^a_{bc} 
 \end{split}
\end{equation*}
and we can now expand the Christoffel symbols $\Gamma^a_{bc}$ using the definition \ref{definition Christoffel}, substitute for $\zeta$ using proposition \ref{proposition zeta} and substitute for $\omega$ using proposition \ref{proposition transport mu}. We also note that
\begin{equation*}
 L^j \frac{x^j}{r} = 1 \Rightarrow L^j_{(\text{small})} \frac{x^j}{r} = 0
\end{equation*}
so
\begin{equation*}
 L^i L^i - 1 = L^i_{(\text{small})}L^i_{(\text{small})}
\end{equation*}
Additionally, we have
\begin{equation*}
 \begin{split}
  L^i \slashed{\nabla}_\mu x^i &= \frac{x^i}{r}\slashed{\nabla}_\mu x^i + L^i_{(\text{small})} \slashed{\nabla}_\mu x^i \\
  &= \frac{1}{2r}\slashed{\nabla}_\mu (x^i x^i) + L^i_{(\text{small})} \slashed{\nabla}_\mu x^i \\
  &= \frac{1}{2r}\slashed{\nabla}_\mu (r^2) + L^i_{(\text{small})} \slashed{\nabla}_\mu x^i \\
  &= L^i_{(\text{small})} \slashed{\nabla}_\mu x^i \\
 \end{split}
\end{equation*}
where we have used the facts that $x^ix^i = r^2$ as well as the fact that $\slashed{\nabla}r = 0$.
\end{proof}

\begin{proposition}[Angular derivatives of rectangular components of the null frame]
\label{proposition angular rectangular}
 The rectangular components of the null frame satisfy the following equations in the region $r \geq r_0$
\begin{equation}
 \begin{split}
  \slashed{\upd}_\mu L^a &= \frac{1}{4}\left(2 L^i \frac{\slashed{\nabla}_\mu x^i}{r} + (\slashed{\nabla}_\mu h)_{LL} \right) L^a
   + \frac{1}{4}\left( (\slashed{\nabla}_\mu h)_{LL} \right) \Lbar^a \\
  &\phantom{=} - \frac{1}{2}(\slashed{g}^{-1})^{\nu\rho} \left( L^b \slashed{\Pi}_\nu^{\phantom{\nu}c}(\slashed{\nabla}_\mu h_{bc}) + \slashed{\Pi}_\mu^{\phantom{\mu}b}\slashed{\Pi}_\nu^{\phantom{\nu}c}(L h_{bc}) - L^b \slashed{\Pi}_\mu^{\phantom{\mu}c}(\slashed{\nabla}_\nu h_{bc}) - 2\chi_{\mu\nu} \right) \slashed{\Pi}_\rho^{\phantom{\rho}a} 
  \\ \\
  \slashed{\upd}_\mu L^i_{(\text{small})} &= \frac{1}{4}\left(2 L^j \frac{\slashed{\nabla}_\mu x^j}{r} + (\slashed{\nabla}_\mu h)_{LL} \right) L^i
   + \frac{1}{4}\left( (\slashed{\nabla}_\mu h)_{LL} \right) \Lbar^i \\
  &\phantom{=} - \frac{1}{2}(\slashed{g}^{-1})^{\nu\rho} \left( L^b \slashed{\Pi}_\nu^{\phantom{\nu}c}(\slashed{\nabla}_\mu h_{bc}) + \slashed{\Pi}_\mu^{\phantom{\mu}b}\slashed{\Pi}_\nu^{\phantom{\nu}c}(L h_{bc}) - L^b \slashed{\Pi}_\mu^{\phantom{\mu}c}(\slashed{\nabla}_\nu h_{bc}) - 2(\chi_{(\text{small})})_{\mu\nu} \right) \slashed{\Pi}_\rho^{\phantom{\rho}i} 
  \\ \\
  \slashed{\upd}_\mu \Lbar^a &= \frac{1}{4}\left( (\slashed{\nabla}_\mu h)_{LL} \right) L^a 
  + \left( \Lbar^i \frac{\slashed{\nabla}_\mu x^i}{r} + \frac{1}{4}(\slashed{\nabla}_\mu h)_{\Lbar\Lbar} \right) \Lbar^a \\
  & \phantom{=} -\frac{1}{2}(\slashed{g}^{-1})^{\nu\rho} \left( \slashed{\Pi}_\mu^{\phantom{\mu}b} \slashed{\Pi}_\nu^{\phantom{\nu}c}(\Lbar h_{bc}) + \Lbar^b \slashed{\Pi}_\nu^{\phantom{\nu}c}(\slashed{\nabla}_\mu h_{bc}) - \Lbar^b\slashed{\Pi}_\mu^{\phantom{\mu}c}(\slashed{\nabla}_\nu h_{bc}) - 2\chibar_{\mu\nu} \right) \slashed{\Pi}_\rho^{\phantom{\rho}a}
  \\
  &= \frac{1}{4}\left( (\slashed{\nabla}_\mu h)_{LL} \right) L^a 
  + \left( \Lbar^i \frac{\slashed{\nabla}_\mu x^i}{r} + \frac{1}{4}(\slashed{\nabla}_\mu h)_{\Lbar\Lbar} \right) \Lbar^a \\
  & \phantom{=} + \left( \chi_{\mu\nu} - \frac{2}{r}\slashed{\Pi}_\mu^{\phantom{\mu}i}\slashed{\Pi}_\nu^{\phantom{\nu}i} - \frac{1}{2} \slashed{\Pi}_\mu^{\phantom{\mu}c} \slashed{\Pi}_\nu^{\phantom{\nu}d}(Lh_{cd}) - \Lbar^c \slashed{\Pi}_\nu^{\phantom{\nu}d} (\slashed{\nabla}_\mu h_{cd}) + \frac{1}{2}L^c \slashed{\Pi}_\mu^{\phantom{\mu}d} (\slashed{\nabla}_\nu h_{cd}) \right) \slashed{\Pi}_\rho^{\phantom{\rho}a}
  \\ \\
  \slashed{\upd}_\mu \Lbar^i_{(\text{small})} &= \frac{1}{4}\left( (\slashed{\nabla}_\mu h)_{LL} \right) L^i 
  + \left( \Lbar^j \frac{\slashed{\nabla}_\mu x^j}{r} + \frac{1}{4}(\slashed{\nabla}_\mu h)_{\Lbar\Lbar} \right) \Lbar^i \\
  & \phantom{=} -\frac{1}{2}(\slashed{g}^{-1})^{\nu\rho} \left( \slashed{\Pi}_\mu^{\phantom{\mu}b} \slashed{\Pi}_\nu^{\phantom{\nu}c}(\Lbar h_{bc}) + \Lbar^b \slashed{\Pi}_\nu^{\phantom{\nu}c}(\slashed{\nabla}_\mu h_{bc}) - \Lbar^b\slashed{\Pi}_\mu^{\phantom{\mu}c}(\slashed{\nabla}_\nu h_{bc}) - 2(\chibar_{(\text{small})})_{\mu\nu} \right) \slashed{\Pi}_\rho^{\phantom{\rho}a}
  \\ \\
  \slashed{\nabla}_\mu \slashed{\Pi}_\nu^{\phantom{\nu}a} &= \frac{1}{4}\left( \Lbar^b \slashed{\Pi}_\nu^{\phantom{\nu}c}(\slashed{\nabla}_\mu h_{bc}) + \Lbar^b \slashed{\Pi}_\mu^{\phantom{\mu}c}(\slashed{\nabla}_\nu h_{bc}) - \slashed{\Pi}_\mu^{\phantom{\mu}b}\slashed{\Pi}_\nu^{\phantom{\nu}c}(\Lbar h_{bc}) + 2\chibar_{\mu\nu} \right)L^a \\
  & \phantom{=} + \frac{1}{4}\left( L^b \slashed{\Pi}_\nu^{\phantom{\nu}c}(\slashed{\nabla}_\mu h_{bc}) + L^b \slashed{\Pi}_\mu^{\phantom{\mu}c}(\slashed{\nabla}_\nu h_{bc}) - \slashed{\Pi}_\mu^{\phantom{\mu}b}\slashed{\Pi}_\nu^{\phantom{\nu}c}(Lh_{bc}) + 2\chi_{\mu\nu} \right)\Lbar^a \\
  & \phantom{=} + \frac{1}{2} (\slashed{g}^{-1})^{\rho\sigma}\left( -\slashed{\Pi}_\nu^{\phantom{\nu}b}\slashed{\Pi}_\rho^{\phantom{\rho}c}(\slashed{\nabla}_\mu h_{bc}) - \slashed{\Pi}_\mu^{\phantom{\mu}b}\slashed{\Pi}_\rho^{\phantom{\rho}c}(\slashed{\nabla}_\nu h_{bc}) + \slashed{\Pi}_\mu^{\phantom{\mu}b}\slashed{\Pi}_\nu^{\phantom{\nu}c}(\slashed{\nabla}_\rho h_{bc}) \right)\slashed{\Pi}_\sigma^{\phantom{\sigma}a}
  \\
  &=
  \frac{1}{4}\left( 2\chi_{\mu\nu} - \frac{4}{r}\slashed{\Pi}_\mu^{\phantom{\mu}i}\slashed{\Pi}_\nu^{\phantom{\nu}i} + L^b \slashed{\Pi}_\nu^{\phantom{\nu}c}(\slashed{\nabla}_\mu h_{bc}) + L^b \slashed{\Pi}_\mu^{\phantom{\nu}c}(\slashed{\nabla}_\nu h_{bc}) \right)L^a \\
  & \phantom{=} + \frac{1}{4}\left( L^b \slashed{\Pi}_\nu^{\phantom{\nu}c}(\slashed{\nabla}_\mu h_{bc}) + L^b \slashed{\Pi}_\mu^{\phantom{\mu}c}(\slashed{\nabla}_\nu h_{bc}) - \slashed{\Pi}_\mu^{\phantom{\mu}b}\slashed{\Pi}_\nu^{\phantom{\nu}c}(Lh_{bc}) + 2\chi_{\mu\nu} \right)\Lbar^a \\
  & \phantom{=} + \frac{1}{2} (\slashed{g}^{-1})^{\rho\sigma}\left( -\slashed{\Pi}_\nu^{\phantom{\nu}b}\slashed{\Pi}_\rho^{\phantom{\rho}c}(\slashed{\nabla}_\mu h_{bc}) - \slashed{\Pi}_\mu^{\phantom{\mu}b}\slashed{\Pi}_\rho^{\phantom{\rho}c}(\slashed{\nabla}_\nu h_{bc}) + \slashed{\Pi}_\mu^{\phantom{\mu}b}\slashed{\Pi}_\nu^{\phantom{\nu}c}(\slashed{\nabla}_\rho h_{bc}) \right)\slashed{\Pi}_\sigma^{\phantom{\sigma}a}
  \\
  \end{split}
\end{equation}

Schematically, we have
\begin{equation}
 \begin{split}
  \slashed{\nabla} \begin{pmatrix}
                    L^a \\
		    \Lbar^a \\
		    \slashed{\Pi}^a
                   \end{pmatrix}
  = \left( (\partial h)_{(\text{frame})} + \frac{1}{r} + \bm{\Gamma} + L^i_{(\text{small})}\frac{\slashed{\nabla}x^i}{r} + \Lbar^i_{(\text{small})}\frac{\slashed{\nabla}x^i}{r} \right)
		   \begin{pmatrix}
                    L^a \\
		    \Lbar^a \\
		    \slashed{\Pi}^a
                   \end{pmatrix} \\ \\
  \slashed{\nabla} \begin{pmatrix}
                    L^i_{(\text{small})} \\
		    \Lbar^i_{(\text{small})} \\
                   \end{pmatrix}
  = \left( (\partial h)_{(\text{frame})} + \bm{\Gamma} + L^i_{(\text{small})}\frac{\slashed{\nabla}x^i}{r} + \Lbar^i_{(\text{small})}\frac{\slashed{\nabla}x^i}{r} \right)
		   \begin{pmatrix}
                    L^a \\
		    \Lbar^a \\
		    \slashed{\Pi}^a
                   \end{pmatrix} \\
 \end{split}
\end{equation}

\end{proposition}

\begin{proof}
 We again expand the corresponding expressions from \ref{proposition null connection}. For the final equation, note that
\begin{equation*}
 \begin{split}
  (X_A)^\mu\slashed{\nabla}_\mu \slashed{\Pi}_\nu^{\phantom{\nu}a} &= \slashed{\D}_A \slashed{\Pi}_\mu^{\phantom{\mu}a} \\
  &= \slashed{\D}_A \left( (\slashed{g}^{-1})^{BC}(X_B)_\mu (X_C)^a \right) \\
  &= \left( X_A (\slashed{g}^{-1})^{BC}\right)(X_B)_\mu (X_C)^a + (\slashed{g}^{-1})^{BC} (\slashed{\D}_A X_B)_\mu (X_C)^a + (\slashed{g}^{-1})^{BC}(X_B)_\mu (X_A (X_C)^a )
 \end{split}
\end{equation*}

For the final line, we rewrite $\chibar$ in terms of $\chi$ using proposition \ref{proposition chibar in terms of chi}.

\end{proof}

Similarly, we can prove the following:
\begin{proposition}[Transport equations for rectangular components in the $\Lbar$ direction]
\label{proposition transport lbar rectangular}
 The rectangular components of the vector fields and the projection operators also  satisfy the following equations in the region $r \geq r_0$:
 \begin{equation}
  \begin{split}
    \Lbar L^a &= \left( -\frac{L^i_{(\text{small})} L^i_{(\text{small})}}{r} + \frac{1}{4}		    (L h)_{\Lbar\Lbar} - \frac{1}{2}(Lh)_{L\Lbar} + \frac{1}{4}(Lh)_{LL} + 		    \frac{1}{4}(\Lbar h)_{LL} \right) L^a
		+ \frac{1}{4}\left( (\Lbar h)_{LL} \right) \Lbar^a\\
		&\phantom{=} + (\slashed{g}^{-1})^{\mu\nu}\bigg(2L^i_{(\text{small})} \frac{\slashed{\nabla}_\mu x^i}{r} - L^b \slashed{\Pi}_\mu^{\phantom{\mu}c}(Lh_{bc}) - \frac{1}{2}(\slashed{\nabla}_\mu h)_{LL} + (\slashed{\nabla}_\mu h)_{L\Lbar} - L^b \slashed{\Pi}_\mu^{\phantom{\mu}c}(\Lbar h_{bc}) \\
		&\phantom{= (\slashed{g}^{-1})^{\mu\nu}\bigg(} + \Lbar^b \slashed{\Pi}_\mu^{\phantom{\mu}c}(Lh_{bc}) + 2\slashed{\nabla}_\mu \log \mu \bigg) \slashed{\Pi}_\nu^{\phantom{\nu}a}
    \\ \\
    \Lbar \Lbar^a &= \frac{1}{4}\left((\Lbar h)_{\Lbar\Lbar} \right)L^a \\
		&\phantom{=} + \left(\frac{L^i_{(\text{small})}L^i_{(\text{small})}}{r} + \frac{1}{2}(Lh)_{L\Lbar} - \frac{1}{4}(Lh)_{LL} - \frac{1}{4}(\Lbar h)_{LL} + \frac{1}{2}(\Lbar h)_{L\Lbar} - \frac{1}{4}(Lh)_{\Lbar\Lbar} \right) \Lbar^a \\
		&\phantom{=} + (\slashed{g}^{-1})^{\mu\nu}\left(\frac{1}{2}(\slashed{\nabla}_\mu h)_{\Lbar\Lbar} - \Lbar^b \slashed{\Pi}_\mu^{\phantom{\mu}c}(\Lbar h_{bc}) + 2(\slashed{\nabla}_\mu \log\mu) \right) \slashed{\Pi}_\nu^{\phantom{\nu}a}
    \\ \\
    \slashed{\D}_{\Lbar} \slashed{\Pi}_\mu^{\phantom{\mu}a} &= 
		\left( \slashed{\nabla}_\mu \log \mu + \frac{1}{4}(\slashed{\nabla}_\mu h)_{\Lbar\Lbar} \right) L^a
		+ \left( L^i_{(\text{small})} \frac{\slashed{\nabla}_\mu x^i}{r} - \frac{1}{4}(\slashed{\nabla}_\mu h)_{LL} + \frac{1}{2}(\slashed{\nabla}_\mu h)_{L\Lbar} + \slashed{\nabla}_\mu \log \mu \right)\Lbar^a \\
		&\phantom{=}- \frac{1}{2}(\slashed{g}^{-1})^{\nu\rho} \left( \slashed{\Pi}_\mu^{\phantom{\mu} b}\slashed{\Pi}_\nu^{\phantom{\nu}c}(\Lbar h_{bc}) + \Lbar^b \slashed{\Pi}_\nu^{\phantom{\nu}c}(\slashed{\nabla}_\mu h_{bc}) - \Lbar^b \slashed{\Pi}_\mu^{\phantom{\nu}c} (\slashed{\nabla}_\nu h_{bc}) \right)\slashed{\Pi}_\rho^{\phantom{\rho}a} 
  \end{split}
 \end{equation}
 Schematically, we have
\begin{equation}
 \slashed{\D}_{\Lbar} \begin{pmatrix} L^a \\ \Lbar^a \\ \slashed{\Pi}^a \end{pmatrix}
  = \left( (\partial h)_{(\text{frame})} + \frac{L^i_{(\text{small})}L^i_{(\text{small})}}{r} + L^i_{(\text{small})}\frac{\slashed{\nabla}x^i}{r} + \slashed{\nabla}\log \mu \right)\begin{pmatrix} L^a \\ \Lbar^a \\ \slashed{\Pi}^a \end{pmatrix}
\end{equation}

\end{proposition}

\section{Derivatives of the projection operator \texorpdfstring{$\slashed{\Pi}$}{Pi}}

On several occasions we will need to take derivatives of the projection operator $\slashed{\Pi}$. These are given in the following set of propositions:

\begin{proposition}[Derivatives of the projection operator]
\label{proposition derivatives of projection}
 The covariant derivatives of the projection operator $\slashed{\Pi}$ are given by
 \begin{equation}
 \begin{split}
  \D_L \slashed{\Pi}_\mu^{\phantom{\mu}\nu} &= -\frac{1}{2}\left( L_\mu \zeta^\nu + \zeta_\mu L^\nu \right) \\
  \D_{\Lbar} \slashed{\Pi}_\mu^{\phantom{\mu}\nu} &= \frac{1}{2}\left(\zeta^\rho + 2(\slashed{\nabla}^\rho \log\mu) \right) \left( \slashed{g}_{\mu\rho} \Lbar^\nu + \Lbar_\mu \slashed{\Pi}_\rho^{\phantom{\rho}\nu} \right) + (\slashed{\nabla}^\rho \log\mu) \left(L_\mu \slashed{\Pi}_\rho^{\phantom{\rho}\nu} + \slashed{g}_{\mu\rho} L^\nu\right) \\
  \slashed{\Pi}_\rho^{\phantom{\rho}\sigma}\D_\sigma \slashed{\Pi}_\mu^{\phantom{\mu}\nu} &= \frac{1}{2}\chi_\rho^{\phantom{\rho}\sigma} \left( \slashed{g}_{\mu\sigma} \Lbar^\nu + \Lbar_\mu \slashed{\Pi}_\sigma^{\phantom{\sigma}\nu} \right) + \frac{1}{2}\chibar_\rho^{\phantom{\rho}\sigma} \left( L_\mu \slashed{\Pi}_\sigma^{\phantom{\sigma}\nu} + \slashed{g}_{\mu\sigma} L^\nu \right) 
 \end{split}
\end{equation}

% \begin{equation}
%  \begin{split}
%   \D_L \slashed{\Pi}_\mu^{\phantom{\mu}\nu} &= -\frac{1}{2}\zeta^A\left( L_\mu (X_A)^\nu + (X_A)_\mu L^\nu \right) \\
%   \D_{\Lbar} \slashed{\Pi}_\mu^{\phantom{\mu}\nu} &= \frac{1}{2}\left(\zeta^A + 2\mu^{-1}(\slashed{\upd}^A \mu) \right) \left( (X_A)_\mu \Lbar^\nu + \Lbar_\mu (X_A)^\nu \right) + \mu^{-1}(\slashed{\upd}^A \mu) \left(L_\mu (X_A)^\nu + (X_A)_\mu L^\nu\right) \\
%   \D_A \slashed{\Pi}_\mu^{\phantom{\mu}\nu} &= \frac{1}{2}\chi_A^{\phantom{A}B} \left( (X_B)_\mu \Lbar^\nu + \Lbar_\mu (X_B)^\nu \right) + \frac{1}{2}\chibar_A^{\phantom{A}B} \left( L_\mu (X_B)^\nu + (X_B)_\mu L^\nu \right) 
%  \end{split}
% \end{equation}

\end{proposition}

Another frequently encountered term is of the form $\slashed{\Pi}\cdot(\D\slashed{\Pi})\D \phi$, where $\phi$ is an $S_{\tau,r}$-tangent tensor field. Despite first appearences, these terms are actually zero-th order in the field $\phi$, i.e.\ they do not depend on any of the derivatives of $\phi$. The following proposition makes this clear.

\begin{proposition}
\label{proposition Pi d Pi}
 Let $\phi$ be a rank one $S_{\tau,r}$-tangent vector field, i.e.\ $\phi_\mu = \phi^A (X_A)_\mu$. Then we have
\begin{equation}
 \begin{split}
  \slashed{\Pi}_\mu^{\phantom{\mu}\nu} (\D_L \slashed{\Pi}_\nu^{\phantom{\nu}\rho}) (\D_L \phi)_\rho
  &= 0 \\
  \slashed{\Pi}_\mu^{\phantom{\mu}\nu} (\D_L \slashed{\Pi}_\nu^{\phantom{\nu}\rho}) (\D_{\Lbar} \phi)_\rho
  &= -\frac{1}{2} \zeta_\mu \left( \zeta^\nu + 2\slashed{\nabla}^\nu \log\mu \right)\phi_\nu \\
  \slashed{\Pi}_\mu^{\phantom{\mu}\nu} (\D_L \slashed{\Pi}_\nu^{\phantom{\nu}\rho}) (\slashed{\Pi}_\lambda^{\phantom{\lambda}\sigma} \D_\sigma \phi_\rho)
  &= - \frac{1}{2} \zeta_\mu \chi_\lambda^{\phantom{\lambda}\sigma} \phi_\sigma \\
  \slashed{\Pi}_\mu^{\phantom{\mu}\nu} (\D_{\Lbar} \slashed{\Pi}_\nu^{\phantom{\nu}\rho}) (\D_L \phi)_\rho
  &= -\frac{1}{2} \left( \zeta_\mu + 2\slashed{\nabla}_\mu \log\mu \right)\zeta^\nu\phi_\nu \\
  \slashed{\Pi}_\mu^{\phantom{\mu}\nu} (\D_{\Lbar} \slashed{\Pi}_\nu^{\phantom{\nu}\rho}) (\D_{\Lbar} \phi)_\rho
  &= \frac{1}{2}\left( \zeta_\mu + 4\slashed{\nabla}_\mu \log\mu \right)\left( \zeta^\nu + 2\slashed{\nabla}^\nu \log\mu \right)\phi_\nu \\
  \slashed{\Pi}_\mu^{\phantom{\mu}\nu} (\D_{\Lbar} \slashed{\Pi}_\nu^{\phantom{\nu}\rho}) (\slashed{\Pi}_\lambda^{\phantom{\lambda}\sigma} \D_\sigma \phi_\rho)
  &= \frac{1}{2} \left( (\zeta_\mu + 2\slashed{\nabla}_\mu \log\mu)\chibar_\lambda^{\phantom{\lambda}\nu} + 2(\slashed{\nabla}_\mu \log\mu) \chi_\lambda^{\phantom{\lambda}\nu} \right)\phi_\nu \\
  \slashed{\Pi}_\mu^{\phantom{\mu}\nu} (\slashed{\Pi}_\sigma^{\phantom{\sigma}\lambda}\D_{\lambda} \slashed{\Pi}_\nu^{\phantom{\nu}\rho}) (\D_L \phi)_\rho
  &= -\frac{1}{2} \chi_{\sigma\mu} \zeta^\nu \phi_\nu \\
  \slashed{\Pi}_\mu^{\phantom{\mu}\nu} (\slashed{\Pi}_\sigma^{\phantom{\sigma}\lambda}\D_{\lambda} \slashed{\Pi}_\nu^{\phantom{\nu}\rho}) (\D_{\Lbar} \phi)_\rho
  &= \frac{1}{2} \left( 2\chi_{\mu\sigma}(\slashed{\nabla}^\nu \log\mu) + \chibar_{\mu\sigma}\left( \zeta^\nu + 2\slashed{\nabla}^\nu \log\mu \right) \right) \phi_\nu \\
  \slashed{\Pi}_\mu^{\phantom{\mu}\nu} (\slashed{\Pi}_\sigma^{\phantom{\sigma}\xi}\D_{\xi} \slashed{\Pi}_\nu^{\phantom{\nu}\rho}) (\slashed{\Pi}_\lambda^{\phantom{\lambda}\kappa}\D_{\kappa} \phi_\rho)
  &= \frac{1}{2}\left( \chi_{\sigma\mu} \chibar_\lambda^{\phantom{\lambda}\nu} + \chibar_{\sigma\mu} \chi_\lambda^{\phantom{\lambda}\nu} \right)\phi_\nu
 \end{split}
\end{equation}

Schematically, we can write
\begin{equation}
 \begin{split}
  \slashed{\Pi}\cdot \begin{pmatrix} \D_L \slashed{\Pi} \\ \D_{\Lbar} \slashed{\Pi} \end{pmatrix} \cdot \begin{pmatrix} \D_L \phi \\ \D_{\Lbar} \phi \end{pmatrix} &= \bm{\Gamma}\cdot \bm{\Gamma}\cdot \phi
  \\
  \slashed{\Pi}\cdot \begin{pmatrix} \D_L \slashed{\Pi} \\ \D_{\Lbar} \slashed{\Pi} \end{pmatrix} \cdot \left((\slashed{\Pi}\cdot \D)\phi \right) &= \bm{\Gamma}\cdot\left( \frac{1}{r} + \bm{\Gamma} \right)\cdot \phi
  \\
  \slashed{\Pi}\cdot \left((\slashed{\Pi}\cdot \D)\Pi \right) \cdot \begin{pmatrix} \D_L \phi \\ \D_{\Lbar} \phi \end{pmatrix} &= \bm{\Gamma}\cdot\left( \frac{1}{r} + \bm{\Gamma} \right)\cdot \phi  
  \\
  \slashed{\Pi}\cdot \left((\slashed{\Pi}\cdot \D)\Pi \right) \cdot \left((\slashed{\Pi}\cdot \D)\phi \right) &= \left( \frac{1}{r} + \bm{\Gamma}_{(\text{good})}\right)\cdot\left( \frac{1}{r} + \bm{\Gamma} \right)\cdot \phi  
  \\
 \end{split}
\end{equation}

\end{proposition}

\begin{proof}
 We derive one of the expressions above, the others follow from very similar calculations. We have
\begin{equation*}
 \begin{split}
  \slashed{\Pi}_\mu^{\phantom{\mu}\nu} (\D_L \slashed{\Pi}_\nu^{\phantom{\nu}\rho}) (\D_{\Lbar} \phi)_\rho
  &= -\frac{1}{2} \zeta^A \slashed{\Pi}_\mu^{\phantom{\mu}\nu}\left( L_\nu (X_A)^\rho + (X_A)_\nu L^\rho \right) (\D_{\Lbar}\phi )_\rho \\
  &= -\frac{1}{2} \zeta^A (X_A)_\mu L^\rho (\D_{\Lbar}\phi )_\rho \\
  &= -\frac{1}{2} \zeta^A (X_A)_\mu (\D_{\Lbar} L)^\rho \phi_\rho \\
  &= -\frac{1}{2}(X_A)_\mu \zeta^A  \left( \zeta^B + 2\mu^{-1}(\slashed{\upd}^B \mu \right) \phi_B
 \end{split}
\end{equation*}
 where we have made use of the expressions in propositions \ref{proposition null connection} and \ref{proposition derivatives of projection}, as well as the fact that $\phi$ is $S_{\tau,r}$-tangent.
\end{proof}

\begin{proposition}[Derivatives of the rectangular components of the projection operators]
\label{proposition derivatives of the rectangular components of Pi}

On occasion we will also need to compute the derivatives of the rectangular components of the projection operators. These are computed in a similar way to the derivatives of the rectangular components of the null frame, and are given by the following equations:
\begin{equation}
 \begin{split}
  L\slashed{\Pi}_a^{\phantom{a}b} &=
  	\left( -L^i_{(\text{small})} \frac{\slashed{\nabla}_\mu x^i}{r} + \frac{1}{4} (\slashed{\nabla}_\mu h)_{LL} \right) \slashed{\Pi}_a^{\phantom{a}\mu} L^b
  	\\
  	&\phantom{=}
  	+ \frac{1}{4}(\slashed{\nabla}_\mu h)_{LL} \slashed{\Pi}_a^{\phantom{a}\mu} \Lbar^b
  	\\
  	&\phantom{=}
  	+ (\slashed{g}^{-1})^{\mu\nu} \left( - L^i\frac{\slashed{\nabla}_\mu x^i}{r}
  	+ \frac{1}{4} (\slashed{\nabla}_\nu h)_{LL}
  	- \frac{1}{2} \Lbar^c \slashed{\Pi}_\nu^{\phantom{\nu}d} (Lh_{cd}) \right) L_a \slashed{\Pi}_\mu^{\phantom{\mu}b}
  	\\
  	&\phantom{=}
  	+ \frac{1}{4} (\slashed{g}^{-1})^{\mu\nu} \left( (\slashed{\nabla}_\nu h)_{LL} - 2L^c \slashed{\Pi}_\nu^{\phantom{\nu}d} (Lh_{cd}) \right) \Lbar_a \slashed{\Pi}_\mu^{\phantom{\mu}b}
  	\\ \\
  	\Lbar \slashed{\Pi}_a^{\phantom{a}b} &=
  	\left( \slashed{\nabla}_\mu \log \mu + \frac{1}{4}(\slashed{\nabla}_\mu h)_{\Lbar\Lbar} \right) \slashed{\Pi}_a^{\phantom{a}\mu} L^b
  	\\
  	&\phantom{=}
  	+ \left( L^i_{(\text{small})}\frac{\slashed{\nabla}_\mu x^i}{r} - \frac{1}{4}(\slashed{\nabla}_\mu h)_{LL}
  	+ \frac{1}{2}(\slashed{\nabla}_\mu h)_{L\Lbar} + \slashed{\nabla}_\mu \log \mu \right) \slashed{\Pi}_a^{\phantom{a}\mu} \Lbar^b
  	\\
  	&\phantom{=}
  	+ (\slashed{g}^{-1})^{\mu\nu} \left(\slashed{\nabla}_\nu \log \mu - \frac{1}{2} \Lbar^c \slashed{\Pi}_\nu^{\phantom{\nu}d} (\Lbar h_{cd}) + \frac{1}{4}(\slashed{\nabla}_\nu h)_{\Lbar\Lbar} \right) L_a \slashed{\Pi}_\mu^{\phantom{\mu}b}
  	\\
  	&\phantom{=}
  	+ (\slashed{g}^{-1})^{\mu\nu} \left( L^i_{(\text{small})}\frac{\slashed{\nabla}_\mu x^i}{r} - \frac{1}{4}(\slashed{\nabla}_\mu h)_{LL} + \frac{1}{2}(\slashed{\nabla}_\mu h)_{L\Lbar} - \frac{1}{2}(\Lbar\slashed{h})_{L\mu}+ \slashed{\nabla}_\mu \log\mu \right) \Lbar_a \slashed{\Pi}_\mu^{\phantom{\mu}b}
    \\ 
    \\
   \slashed{\nabla}_\mu \slashed{\Pi}_a^{\phantom{a}b} &=
    \left( \frac{1}{4} \Lbar^c \slashed{\Pi}_\nu^{\phantom{\nu}d}(\slashed{\nabla}_\mu h_{cd})
    	+ \frac{1}{4} \Lbar^c \slashed{\Pi}_\mu^{\phantom{\mu}d}(\slashed{\nabla}_\nu h_{cd})
    	- \frac{1}{4} \slashed{\Pi}_\mu^{\phantom{\mu}c} \slashed{\Pi}_\nu^{\phantom{\nu}d}(\Lbar h_{cd})
    	+ \frac{1}{2} \chibar_{\mu\nu}
    \right) \slashed{\Pi}_a^{\phantom{a}\nu} L^b
    \\
    &\phantom{=}
    + \left( \frac{1}{4} L^c \slashed{\Pi}_\nu^{\phantom{\nu}d}(\slashed{\nabla}_\mu h_{cd})
    + \frac{1}{4} L^c \slashed{\Pi}_\mu^{\phantom{\mu}d}(\slashed{\nabla}_\nu h_{cd})
    - \frac{1}{4} \slashed{\Pi}_\mu^{\phantom{\mu}c} \slashed{\Pi}_\nu^{\phantom{\nu}d}(L h_{cd})
    + \frac{1}{2} \chi_{\mu\nu}
    \right) \slashed{\Pi}_a^{\phantom{a}\nu} \Lbar^b
    \\
    &\phantom{=}
    + (\slashed{g}^{-1})^{\nu\rho} \left( -\frac{1}{4} \Lbar^c \slashed{\Pi}_\rho^{\phantom{\rho}d}(\slashed{\nabla}_\mu h_{cd})
    + \frac{1}{4} \Lbar^c \slashed{\Pi}_\mu^{\phantom{\mu}d}(\slashed{\nabla}_\rho h_{cd})
    - \frac{1}{4} \slashed{\Pi}_\mu^{\phantom{\mu}c} \slashed{\Pi}_\rho^{\phantom{\rho}d}(\Lbar h_{cd})
    + \frac{1}{2} \chibar_{\mu\rho}
    \right) L_a \slashed{\Pi}_\nu^{\phantom{\nu}b}
    \\
    &\phantom{=}
    + (\slashed{g}^{-1})^{\nu\rho} \left( -\frac{1}{4} L^c \slashed{\Pi}_\rho^{\phantom{\rho}d}(\slashed{\nabla}_\mu h_{cd})
    + \frac{1}{4} L^c \slashed{\Pi}_\mu^{\phantom{\mu}d}(\slashed{\nabla}_\rho h_{cd})
    - \frac{1}{4} \slashed{\Pi}_\mu^{\phantom{\mu}c} \slashed{\Pi}_\rho^{\phantom{\rho}d}(L h_{cd})
    + \frac{1}{2} \chi_{\mu\rho}
    \right) \Lbar_a \slashed{\Pi}_\nu^{\phantom{\nu}b}
 \end{split}
\end{equation}

Schematically, we have
\begin{equation}
 \begin{split}
  L\slashed{\Pi}_a^{\phantom{a}b} = \left( L^i_{(\text{small})}\frac{\slashed{\nabla}x^i}{r} + (\bar{\partial}h)_{(\text{frame})} \right)\cdot \begin{pmatrix} L^a \\ \Lbar^a \end{pmatrix} \cdot (\slashed{\Pi}^a) \\
  \Lbar\slashed{\Pi}_a^{\phantom{a}b} = \left( L^i_{(\text{small})}\frac{\slashed{\nabla}x^i}{r} + (\partial h)_{(\text{frame})} + \slashed{\nabla}\log \mu \right)\cdot \begin{pmatrix} L^a \\ \Lbar^a \end{pmatrix} \cdot (\slashed{\Pi}^a) \\
  \slashed{\nabla}\slashed{\Pi}_a^{\phantom{a}b} = \left( \frac{1}{r} + \bm{\Gamma} \right)\cdot \begin{pmatrix} L^a \\ \Lbar^a \end{pmatrix} \cdot (\slashed{\Pi}^a) \\
 \end{split}
\end{equation}

\end{proposition}

%Occasionally we will want to project a set of scalar fields, laballed by rectangular indices, onto the null frame using the natural projection arising from the rectangular components of the null frame vector fields. For example, for a set of scalar fields $\phi_a$, we can form the scalar fields $\phi_L := L^a \phi_a$, $\phi_{\Lbar} := \Lbar^a \phi_a$ and $\phi_A := (X_A)^a \phi_a$.
%
%A problem arises, however, when we try to use this notation to measure the size of the collection of fields $\phi_a$. Let us suppose, for example, that each field $\phi_a$ is bounded. Then the fields $\phi_L$ and $\phi_{\Lbar}$ will also be bounded, but the fields $\phi_A$ will grow as $|\phi_A| \sim r$, and so the collection of fields $\{ \phi_L, \phi_{\Lbar}, \phi_A\}$ will not have the same behaviour as the collection of fields $\phi_a$.
%
%A better way to project the scalar fields is to make use of ``rectangular frame angular one forms''. These use a version of the projection operators $\slashed{\Pi}$ to form a collection of $S_{\tau,r}$-tangent one-forms from a collection of sets labelled by rectangular indices. Together with the projections associated with the rectangular components of $L$ and $\Lbar$, these one-forms will give us a satisfactory projection of the scalar fields.

\section{Null frame decompositions of the wave operator}
In this section we will establish several different expressions for the scalar wave operator in terms of the null frame and the connection coefficients.

\subsection{Decomposition of the scalar wave operator}

\begin{proposition}[Null decompositions of the scalar wave operator]
\label{proposition scalar wave operator}
 Let $\phi$ be a scalar function. Then in the region $r \geq r_0$ the wave operator can be decomposed as
\begin{equation}
\label{equation wave operator null}
 \begin{split}
  \Box_g \phi &= -\frac{1}{2}L\Lbar \phi - \frac{1}{2}\Lbar L \phi + \slashed{\Delta}\phi - \frac{1}{2}\left(\tr_{\slashed{g}}\chi + \omega\right) \Lbar\phi - \frac{1}{2}\left(\tr_{\slashed{g}}\chibar - \omega \right) L\phi + (\slashed{\nabla}^\alpha \log\mu) \slashed{\upd}_\alpha \phi \\
  &= -L\Lbar\phi + \slashed{\Delta}\phi - \frac{1}{2}(\tr_{\slashed{g}}\chi + 2\omega) \Lbar\phi - \frac{1}{2}\tr_{\slashed{g}}\chibar \, L\phi - \zeta^\alpha \slashed{\upd}_\alpha \phi
 \end{split}
\end{equation}
Additionally, in the same region the $r$-weighted wave operator can be decomposed as
\begin{equation}
 r \Box_g \phi = -L\left( r\Lbar\phi\right) + r\slashed{\Delta}\phi - \frac{1}{2}(\tr_{\slashed{g}}\chi_{(\text{small})} + 2\omega) r \Lbar\phi - \frac{1}{2}(\tr_{\slashed{g}}\chibar )r L\phi - \zeta^\alpha r \slashed{\upd}_\alpha \phi
\end{equation}
%and the $\mu$-weighted wave operator can be decomposed as
%\begin{equation}
% \mu\Box_g \phi = -L(\mu\Lbar\phi) + \mu\slashed{\Delta}\phi - \frac{1}{2}\mu\tr_{\slashed{g}}\chi \, \Lbar\phi - \frac{1}{2}\mu\tr_{\slashed{g}}\chibar \, L\phi - \mu\zeta^A \slashed{\upd}_A \phi
%\end{equation}

\end{proposition}
\begin{proof}
 The wave operator is given by
\begin{equation*}
 \begin{split}
  \Box_g \phi &= (g^{-1})^{\mu\nu} \D_\mu \D_\nu \phi \\
  &= \left(-\frac{1}{2}L^\mu \Lbar^\nu - \frac{1}{2}\Lbar^\mu L^\nu + (\slashed{g}^{-1})^{\mu\nu} \right) \D_\mu \D_\nu \phi \\
  &= -\frac{1}{2}L\Lbar\phi - \frac{1}{2}\Lbar L\phi + (\slashed{g}^{-1})^{AB}X_A X_B \phi + \frac{1}{2}(\D_L \Lbar) \phi + \frac{1}{2}(\D_{\Lbar} L)\phi - (\D_A X_B)\phi
 \end{split}
\end{equation*}
Now, using the expressions in proposition \ref{proposition null connection} the identity
\begin{equation*}
 \slashed{\Delta}\phi = (\slashed{g}^{-1})^{AB} X_A X_B \phi - (\slashed{\nabla}_A X_B) \phi
\end{equation*}
proves the first line of the proposition. To prove the second line, we commute $L$ and $\Lbar$ using the expressions in proposition \ref{proposition commutators null frame}. Finally, to prove the $r$-weighted wave operator decomposition, we recall that $\tr_{\slashed{g}}\chi = 2r^{-1} + \tr_{\slashed{g}}\chi_{(\text{small})}$.
\end{proof}

% \begin{remark}[The $\mu$-weighted wave operator]
%  It is also possible to decompose the $\mu$-weighted wave equation $\mu \Box_g \phi$ in the null frame. In this case, the ``worst'' error term $\omega \Lbar\phi$ can be eliminated. This is the stratagy employed in [[]]
% \end{remark}

\subsection{Decomposition of the projected wave operator}
\begin{proposition}[Null decompositions of the projected wave operator]
\label{proposition tensor wave operator}
 Let $\phi$ be an $S_{\tau,r}$ tensor. Then in the region $r \geq r_0$ the projected wave operator  can be decomposed as

\begin{equation}
  \begin{split}
    \slashed{\Box}_g \phi
      &=
      - \frac{1}{2}\slashed{\D}_L \slashed{\D}_{\Lbar} \phi
      - \frac{1}{2}\slashed{\D}_{\Lbar} \slashed{\D}_L \phi
      + \slashed{\Delta}\phi
      - \frac{1}{2}\left(\tr_{\slashed{g}}\chi + \omega\right) \slashed{\D}_{\Lbar}\phi
      - \frac{1}{2}\left(\tr_{\slashed{g}}\chibar - \omega \right) \slashed{\D}_L\phi
      + (\slashed{\nabla}^\alpha \log\mu) \slashed{\nabla}_\alpha \phi
      \\
      \\
      &=
      - \slashed{\D}_L \slashed{\D}_{\Lbar} \phi
      + \slashed{\Delta}\phi
      - \frac{1}{2}\left(\tr_{\slashed{g}}\chi + 2\omega\right) \slashed{\D}_{\Lbar}\phi
      - \frac{1}{2}\tr_{\slashed{g}}\chibar \slashed{\D}_L\phi
      + \zeta^\alpha \slashed{\nabla}_\alpha \phi
      \\
      &\phantom{=}
      + \frac{1}{2}\left(\slashed{\D}_L \slashed{\D}_{\Lbar} \phi - \slashed{\D}_{\Lbar} \slashed{\D}_L \phi - \slashed{\D}_{[L, \Lbar]}\phi\right)
  \end{split}
\end{equation}

\end{proposition}
\begin{proof}
 This follows in the same way as the proof of the first line of equation \eqref{equation wave operator null}. Note that there is an extra term in the second expression above relative to the second expression in equation \eqref{equation wave operator null}, due to the fact that the covariant derivatives generate curvature factors when commuted past each other. These curvature terms can be expressed in terms of the connection coefficients and the Riemann curvature of $M$, but an alternative description can be given in terms of the curvature of $\slashed{\D}$, interpreted as a connection on the vector bundle of $S_{\tau,r}$-tangent tensor fields - see chapter \ref{chapter geometry of vector bundle}.
\end{proof}

\section{Geometric quantities in the region \texorpdfstring{$r \leq r_0$}{r < r0}}
\label{section geometric quantities in r leq r0}

We will also need to express the various geometric quantities (such as the rectangular components of the null frame, or the null frame connection coefficients) in the region $r \leq r_0$. Note that, in this region, we simply have $\tau = t = x^0$, so that all quantities are ``non-geometric'', and are really defined in terms of the rectangular components of the metric. Hence, we will find that we can express all of these quantities algebraically in terms of the quantities $h_{ab}$.

\begin{proposition}[The rectangular components of the null frame components in the region $r \leq r_0$]
	\label{proposition rectangular components in r < r0}
	In the region $r \leq r_0$, the rectangular components of the frame fields $L$, $\Lbar$, and the rectangular components of the tensor field $\slashed{\Pi}$ can be expressed \emph{algebraically} in terms of the rectangular components of the tensor field $h$. Furthermore, the ``small'' components of these fields can be expressed as a linear combination of the rectangular components of $h$, plus some error term which is at least quadratic in the rectangular components of $h$.
\end{proposition}

\begin{proof}
	Recall that $L$ is defined as an outgoing null vector field, satisfying $L(r) = 1$, and normal to the spheres of constant $r$ and $\tau$. Hence, we can write (in terms of the rectanglar coordinate vector fields)
	\begin{equation*}
	L = A \partial_t + \frac{x^i}{r} \partial_i
	\end{equation*}
	Then, the condition that $L$ is null implies that $A$ satisfies a certain quadratic equation. The condition that $L$ is outgoing picks out one of the solutions to this equation, and so we obtain
	\begin{equation*}
	A = \frac{\frac{x^i}{r} h_{0i} + \sqrt{ \left(\frac{x^i}{r} h_{0i} \right) + (1-h_{00}) \left( 1 + \frac{x^i x^j}{r^2} h_{ij} \right)}}{1-h_{00}}
	\end{equation*}
	
	From this, we also find that $L_{(\text{small})}$ can be expressed as
	\begin{equation*}
	L_{(\text{small})} = A_{(\text{small})} \partial_t
	\end{equation*}
	where
	\begin{equation*}
	\begin{split}
	A_{(\text{small})} 
	&=
	\left( \frac{\frac{x^i}{r} h_{0i} + \sqrt{ \left(\frac{x^i}{r} h_{0i} \right) + (1-h_{00}) \left( 1 + \frac{x^i x^j}{r^2} h_{ij} \right)}}{1-h_{00}} \right) - 1
	\\
	&= \frac{x^i}{r} h_{0i} + \frac{1}{2} \frac{x^i x^j}{r^2} h_{ij} + \frac{1}{2}h_{00} + \mathcal{O}(|h_{(\text{rect})}|^2)
	\end{split}
	\end{equation*}
	
	Next, we recall that $\Lbar$ is defined to be another null vector field normal to the spheres and satisfying $g(L,\Lbar) = -2$. By the condition that $\Lbar$ is null, normal to the spheres and ingoing we find that
	\begin{equation*}
	\Lbar = B\left( A \partial_t - \frac{x^i}{r} \partial_i \right)
	\end{equation*}
	where $A$ is as above, and $B > 0$ is some other function that we need to determine\footnote{Note that, unlike in the region $r \geq r_0$, for reasons of regularity at the origin we avoid imposing $g^{-1}(\upd r, \upd r) = 1$ in the region $r \leq r_0$. Hence, we cannot immediately conclude that $B = 1$.}. By the condition $g(L,\Lbar) = -2$ we find that $B$ is given by
	\begin{equation*}
	\begin{split}
	B 
	&=
	\left( 1 + \frac{1}{2} \frac{x^i x^j}{r^2} h_{ij} + A_{(\text{small})} + \frac{1}{2}(A_{(\text{small})})^2 - A_{(\text{small})} h_{00} - \frac{1}{2}(A_{(\text{small})})^2 h_{00} \right)^{-1}
	\\
	&=
	1 - \frac{x^i}{r} h_{0i} - \frac{x^i x^j}{r^2} h_{ij} - \frac{1}{2}h_{00} + \mathcal{O}(|h_{(\text{rect})}|^2)
	\end{split}
	\end{equation*}
	
	This then allows us to write
	\begin{equation*}
	\Lbar_{(\text{small})}
	=
	(BA - 1) \partial_t - (B-1)\frac{x^i}{r} \partial_i
	\end{equation*}

	Finally, we note that the projection operator $\slashed{\Pi}_\mu^{\phantom{\mu}\nu}$ can be defined in terms of these other fields as
	\begin{equation*}
	\slashed{\Pi}_\mu^{\phantom{\mu}\nu} = \delta_\mu^\nu + \frac{1}{2}L_\mu \Lbar^\nu + \frac{1}{2}\Lbar_\mu L^\nu
	\end{equation*}
	This differs from the ``standard'' (flat or background) projection operator by the addition of another ``small'' tensor field. Specifically, if we define
	\begin{equation*}
	(\slashed{\Pi}_{(\text{background})})_\mu^{\phantom{\mu}\nu}
	:=
	\delta_\mu^\nu + \frac{1}{2}(L - L_{(\text{small})})_\mu (\Lbar - \Lbar_{(\text{small})})^\nu + \frac{1}{2}(\Lbar - \Lbar_{(\text{small})})_\mu (L - L_{(\text{small})})^\nu
	\end{equation*}
	Then we define the tensor field
	\begin{equation*}
	\slashed{\Pi}_{(\text{small})}
	:=
	\slashed{\Pi} - \slashed{\Pi}_{(\text{background})}
	\end{equation*}
	Then it should be clear from the considerations above that the rectangular components of $\slashed{\Pi}_{(\text{small})}$ can be expressed in terms of the rectangular components $h_{ab}$.
	
\end{proof}

\begin{proposition}[The ``inverse foliation density'' in the region \texorpdfstring{$r \leq r_0$}{r < r0}]
	\label{proposition foliation density in r < r0}
	In the region $r \leq r_0$, the funciton $\mu$ can be expressed algebreically in terms of the rectangular components of $h$. In particular, $\mu = 1 + \mu_{(\text{lin})} + \mathcal{O}(|h_{(\text{rect})}|^2)$, where $\mu_{(\text{lin})}$ is linear in the rectangular components of $h$.
\end{proposition}

\begin{proof}
	Recall that the quantity\footnote{We remind the reader that this is not actually the foliation density in the region $r \leq r^0$!} $\mu$ is defined as
	\begin{equation*}
	\mu^{-1} := -(g^{-1})(\upd r , \upd u)
	\end{equation*}
	where, in the region $r \leq r_0$ , we have \emph{defined} $u$ by 
	\begin{equation*}
	u := t - r
	\end{equation*}
	Hence we have
	\begin{equation*}
	\begin{split}
	\mu^{-1} &= (g^{-1})(\upd r , \upd r - \upd t)
	\\
	&= 1 + \frac{x^i x^j}{r^2} H^{ij} - \frac{x^i}{r} H^{0i}
	\end{split}
	\end{equation*}
	where we recall that $g^{-1} = m^{-1} + H$.
\end{proof}

Many geometric quantities are not continuous at $r = r_0$, or else not differentiable here. This will not cause a serious problem: for example, the rectangular components of the frame fields are continuous at $r = r_0$ but not differentiable there. When calculating the derivatives of these quantities in the region $r \geq r_0$, we use the formulae which are derived in section \ref{section derivatives of rectangular components of null frame}. For example, although the the quantity $LL^a$ is not continuous at $r = r_0$, we can use the expression given in proposition \ref{proposition transport La} to calculate $LL^a$ in the region $r \geq r_0$, while in the region $r \leq r_0$ we can simply take the $L$ derivative of given in proposition \ref{proposition rectangular components in r < r0} to express this in terms of derivatives of the metric.

However, the inverse foliation density $\mu$ is not continuous at $r = r_0$, and in the region $r \geq r_0$ we will calculate $\mu$ by integrating the transport equation in proposition \ref{proposition transport mu}. Hence, we must find the correct ``initial value'' for this quantity.

\begin{proposition}[The value of $\mu$ as $r \searrow r_0$]
\label{proposition initial data for mu}
	As $r \rightarrow r_0$ from above, we have
	\begin{equation*}
		\lim_{r \searrow r_0} \mu 
		=
		\frac{1}{ H^{ti}\frac{x^i}{r} + \alpha \left(1 + H^{ij} \frac{x^i x^j}{r^2} \right) }
	\end{equation*}
	where $\alpha$ is defined as
	\begin{equation*}
		\alpha =
			\frac{
				-H^{ti} \frac{x^i}{r}
				-\sqrt{
					1 + \left(H^{ij} + H^{ti} H^{tj}\right) \frac{x^i x^j}{r^2}
				}
			}{1 + H^{ij} \frac{x^i x^j}{r^2}}
	\end{equation*}	
	
\end{proposition}

\begin{proof}
	Recall that in the region $r \geq r_0$, $u$ is defined by
	\begin{equation*}
	\begin{split}
		u \big|_{r = r_0} &= t-r \\
		g^{-1}(\upd u, \upd u) &= 0
	\end{split}
	\end{equation*}
	Hence, the \emph{tangential} derivatives of $u$ on the surface $r = r_0$ are prescribed as part of the ``initial data'', while the \emph{transverse} derivatives of $u$ are determined by the eikonal equation.
	
	We have
	\begin{equation*}
		\upd u =
			\frac{\partial u}{\partial t} \Big|_{r, \vartheta^A} \upd t
			+ \frac{\partial u}{\partial r} \Big|_{t, \vartheta^A} \upd t
			+ \frac{\partial u}{\partial \vartheta^A} \Big|_{r, t} \upd \vartheta^A
	\end{equation*}
	and also
	\begin{equation*}
	\begin{split}
		\left( \frac{\partial u}{\partial t} \Big|_{r, \vartheta^A} \right)\bigg|_{r = r_0} &= 1 \\
		\left( \frac{\partial u}{\partial \vartheta^A} \big|_{r, t} \right)\bigg|_{r = r_0} &= 0
	\end{split}
	\end{equation*}
	by our prescription of the initial data. On the other hand, the $r$ derivative of $u$ can be calculated from the eikonal equation: if we set
	\begin{equation*}
		\left( \frac{\partial u}{\partial r} \Big|_{t, \vartheta^A} \right)\bigg|_{r = r_0} = \alpha
	\end{equation*}
	then the eikonal equation gives us
	\begin{equation*}
		(g^{-1})^{tt} + 2 (g^{-1})^{tr} \alpha + (g^{-1})^{rr} \alpha^2 = 0
	\end{equation*}
	Expanding this in terms of the fields $H^{ab}$, we have
	\begin{equation*}
		-1 + H^{tt} + 2 H^{ti}\frac{x^i}{r} \alpha +  \left( 1 + H^{ij} \frac{x^i x^j}{r^2} \right) \alpha^2 = 0
	\end{equation*}
	Since $u$ is supposed to be a \emph{retarded} solution to the wave equation, we pick the negative root of this quadratic equation for $\alpha$. This gives us
	\begin{equation*}
		\alpha =
			\frac{
				-H^{ti} \frac{x^i}{r}
				-\sqrt{
					1 + \left(H^{ij} + H^{ti} H^{tj}\right) \frac{x^i x^j}{r^2}
				}
			}{1 + H^{ij} \frac{x^i x^j}{r^2}}
	\end{equation*}
	
	Next, we recall that $\mu$ is defined by
	\begin{equation*}
		\mu := \frac{-1}{g^{-1}(\upd u , \upd r)}
	\end{equation*}
	Although $u$ is continuous at $r = r_0$, it is not continuously differentiable. However, its derivatives as $ r \searrow r_0$ are found using the above calculations: we have
	\begin{equation*}
		\lim_{r \searrow r_0} \mu 
		=
		\frac{1}{ H^{ti}\frac{x^i}{r} + \alpha \left(1 + H^{ij} \frac{x^i x^j}{r^2} \right) }
	\end{equation*}
	
	Note that $\alpha = -1 + \mathcal{O}(H)$, so $\lim_{r \searrow r_0} \mu = 1 + \mathcal{O}(H)$. For our global existence result, this is actually all that is needed.
\end{proof}

\begin{proposition}[The connection component \texorpdfstring{$\omega$}{omega} in the region \texorpdfstring{$r \leq r_0$}{r < r0}]
	
	In the region $r \leq r_0$, the connection component $\omega$ can be expressed as
	\begin{equation*}
	\omega = \omega_{(\text{lin})} + \mathcal{O}\left( (h_{(\text{rect})}, \partial h_{(\text{rect})})^2 \right)
	\end{equation*}
	where $\omega_{(\text{lin})}$ is linear in the rectangular components of $h$ and their first derivatives.
	
\end{proposition}

\begin{proof}
	Recall that $\omega$ is defined by the relation
	\begin{equation*}
	\omega = -\frac{1}{2} g(\D_L L, \Lbar)
	\end{equation*}
	Expanding this in rectangular coordinates, we find that
	\begin{equation*}
	\omega = -\frac{1}{2} \Lbar^a g_{ab} (L L^b) - \frac{1}{2} (Lh)_{L\Lbar} + \frac{1}{4}(\Lbar h)_{LL}
	\end{equation*}
	Now, we note that
	\begin{equation*}
	L L^0 = L L^0_{(\text{small})}
	\end{equation*}
	and
	\begin{equation*}
	\begin{split}
	L L^i &= \frac{L^i}{r} - \frac{x^i}{r^2} + L L^i_{(\text{small})}
	\\
	&= \frac{L^i_{(\text{small})}}{r} + L L^i_{(\text{small})}
	\end{split}
	\end{equation*}
	so we can write
	\begin{equation*}
	L L^a = \frac{L^a_{(\text{small})}}{r} - \delta_0^a\frac{L^0_{(\text{small})}}{r} + L L^a_{(\text{small})}
	\end{equation*}
	and in the end we obtain
	\begin{equation*}
	\omega = 
	-\frac{1}{2} \Lbar_a (L L_{(\text{small})}^a)
	- \frac{1}{2} \frac{\Lbar_a L^a_{(\text{small})}}{r}
	+ \frac{1}{2} \frac{\Lbar_0 L^0_{(\text{small})}}{r}
	- \frac{1}{2} (Lh)_{L\Lbar} 
	+ \frac{1}{4}(\Lbar h)_{LL}
	\end{equation*}
	Recalling that $L_{(\text{small})}^a$ can be expressed in terms of the rectangular components of $h$, plus some higher order terms, yields the proposition.
	
\end{proof}

\begin{proposition}[The connection component \texorpdfstring{$\chi$}{chi} in the region \texorpdfstring{$r \leq r_0$}{r < r0}]
	\label{proposition chi in r leq r0}
	The rectangular components of the connection component $\chi$ in the region $r \leq r_0$ can be expressed as
	\begin{equation*}
	\chi_{ab} 
	=
	\frac{2}{r} \slashed{g}_{ab}
	+ (\chi_{(\text{lin})})_{ab}
	+ \mathcal{O}\left(\frac{1}{r}|h_{(\text{rect})}|^2\right)
	+ \mathcal{O}\left(|\partial h_{(\text{rect})}| |h_{(\text{rect})}|\right)
	\end{equation*}
	where $(\chi_{(\text{lin})})_{ab}$ are linear functions of the \emph{derivatives} of the rectangular components of $h$.
\end{proposition}

\begin{proof}
	The connection component $\chi$ is defined as the extrinsic curvature of the spheres $S_{\tau,r}$ with respect to the vector field $L$. Specifically, we have
	\begin{equation*}
	\chi_{\mu\nu} = \frac{1}{2} \slashed{\Pi}_\mu^{\phantom{\mu}\rho} \slashed{\Pi}_\nu^{\phantom{\nu}\sigma} (\mathcal{L}_L g)_{\rho\sigma}
	\end{equation*}
	so, in rectangular coordinates, we have
	\begin{equation*}
	\begin{split}
	\chi_{ab} 
	&=
	\frac{1}{2} \slashed{\Pi}_a^{\phantom{a}c} \slashed{\Pi}_b^{\phantom{b}d} (\mathcal{L}_L g)_{cd}
	\\
	&=
	\frac{1}{2} \slashed{\Pi}_a^{\phantom{a}c} \slashed{\Pi}_b^{\phantom{b}d} \left( Lg_{cd} + (\partial_c L^e) g_{ed} + (\partial_d L^e) g_{ce} \right)
	\\
	&=
	\frac{1}{2} \slashed{\Pi}_a^{\phantom{a}c} \slashed{\Pi}_b^{\phantom{b}d} \left( Lh_{cd} + (\partial_c L^e) ( m_{ed} + h_{ed}) + (\partial_d L^e) (m_{ce} + h_{ce}) \right)
	\\
	&=
	\frac{1}{2} \slashed{\Pi}_a^{\phantom{a}c} \slashed{\Pi}_b^{\phantom{b}d} (Lh_{cd})
	+ (\slashed{\nabla}_a L^e) \slashed{\Pi}_b^{\phantom{b}d}( m_{ed} + h_{ed})
	+ (\slashed{\nabla}_b L^e) \slashed{\Pi}_a^{\phantom{b}d}( m_{ed} + h_{ed})
	\end{split}
	\end{equation*}
	
	Now, we have
	\begin{equation*}
	\slashed{\nabla}_a L^0
	=
	\slashed{\nabla}_a L^0_{(\text{small})}
	\end{equation*}
	and
	\begin{equation*}
	\begin{split}
	\slashed{\nabla}_a L^i
	&=
	\frac{\slashed{\Pi}_a^{\phantom{a}i}}{r}
	+ \slashed{\nabla}_a L^i_{(\text{small})}
	\end{split}
	\end{equation*}
	So, putting these calculations together, we have
	\begin{equation*}
	\slashed{\nabla}_a L^b
	=
	\frac{1}{r} \slashed{\Pi}_a^{\phantom{a}b}
	+ \slashed{\nabla}_a L^b_{(\text{small})}
	- \frac{1}{r} \delta_0^b \slashed{\Pi}_a^{\phantom{a}0}
	\end{equation*}
	and so	
	\begin{equation*}
	\begin{split}
	\chi_{ab} 
	&=
	\frac{2}{r} \slashed{g}_{ab}
	+ \frac{1}{2} \slashed{\Pi}_a^{\phantom{a}c} \slashed{\Pi}_b^{\phantom{b}d} (Lh_{cd})
	+ (\slashed{\nabla}_a L^e_{(\text{small})}) \slashed{\Pi}_b^{\phantom{b}d} (m_{ed} + h_{ed})
	+ (\slashed{\nabla}_b L^e_{(\text{small})}) \slashed{\Pi}_a^{\phantom{a}d} (m_{ed} + h_{ed})
	\\
	&\phantom{=}
	- \frac{1}{r} \slashed{\Pi}_a^{\phantom{a}0} \slashed{\Pi}_b^{\phantom{b}d} (m_{0d} + h_{0d})
	- \frac{1}{r} \slashed{\Pi}_b^{\phantom{b}0} \slashed{\Pi}_a^{\phantom{a}d} (m_{0d} + h_{0d})
	\end{split}
	\end{equation*}
	Note that we actually have $(\slashed{\Pi}_{(\text{background})})_b^{\phantom{b}d} m_{0d} = 0$, and also $(\slashed{\Pi}_{(\text{background})})_b^{\phantom{b}0} = 0$, so these last two terms are actually quadratic in $h$. So, we have
	\begin{equation*}
	\begin{split}
	(\chi_{(\text{small})})_{ab} 
	&=
	\frac{1}{2} \slashed{\Pi}_a^{\phantom{a}c} \slashed{\Pi}_b^{\phantom{b}d} (Lh_{cd})
	+ (\slashed{\nabla}_a L^e_{(\text{small})}) \slashed{\Pi}_b^{\phantom{b}d} (m_{ed} + h_{ed})
	+ (\slashed{\nabla}_b L^e_{(\text{small})}) \slashed{\Pi}_a^{\phantom{a}d} (m_{ed} + h_{ed})
	\\
	&\phantom{=}
	- \frac{1}{r} \slashed{\Pi}_a^{\phantom{a}0} \slashed{\Pi}_b^{\phantom{b}d} (m_{0d} + h_{0d})
	- \frac{1}{r} \slashed{\Pi}_b^{\phantom{b}0} \slashed{\Pi}_a^{\phantom{a}d} (m_{0d} + h_{0d})
	\\
	\\
	&=
	\frac{1}{2} \slashed{\Pi}_a^{\phantom{a}c} \slashed{\Pi}_b^{\phantom{b}d} (Lh_{cd})
	+ (\slashed{\nabla}_a L^e_{(\text{small})}) (\slashed{\Pi}_{(\text{background})})_b^{\phantom{b}d} m_{ed}
	+ (\slashed{\nabla}_b L^e_{(\text{small})}) (\slashed{\Pi}_{(\text{background})})_a^{\phantom{a}d} m_{ed}
	\\
	&\phantom{=}
	+ \mathcal{O}\left(\frac{1}{r}|h_{(\text{rect})}|^2\right)
	+ \mathcal{O}\left(|\partial h_{(\text{rect})}| |h_{(\text{rect})}|\right)
	\end{split}
	\end{equation*}
	Recalling that $L^e_{(\text{small})}$ can be expressed algebraically in terms of $h_{ab}$ in the region $r \leq r_0$ proves the proposition.
	
\end{proof}

\begin{proposition}[The connection component \texorpdfstring{$\chibar$}{chibar} in the region \texorpdfstring{$r \leq r_0$}{r < r0}]
	\label{proposition chibar in r leq r0}
	The rectangular components of the connection component $\chibar$ in the region $r \leq r_0$ can be expressed as
	\begin{equation*}
	\chibar_{ab}
	=
	- \frac{2}{r} \slashed{g}_{ab}
	+ (\chibar_{(\text{lin})})_{ab}
	+ \mathcal{O}\left(\frac{1}{r}|h_{(\text{rect})}|^2\right)
	+ \mathcal{O}\left(|\partial h_{(\text{rect})}| |h_{(\text{rect})}|\right)
	\end{equation*}
	where $(\chibar_{(\text{lin})})_{ab}$ are linear functions of the \emph{derivatives} of the rectangular components of $h$.
\end{proposition}

\begin{proof}
	The connection component $\chibar$ is defined as the extrinsic curvature of the spheres $S_{\tau,r}$ with respect to the vector field $\Lbar$. Specifically, we have
	\begin{equation*}
	\chibar_{\mu\nu} = \frac{1}{2} \slashed{\Pi}_\mu^{\phantom{\mu}\rho} \slashed{\Pi}_\nu^{\phantom{\nu}\sigma} (\mathcal{L}_{\Lbar} g)_{\rho\sigma}
	\end{equation*}
	Hence, we can repeat the steps of proposition \ref{proposition chi in r leq r0}, replacing $L$ with $\Lbar$. This leads to
	\begin{equation*}
	\begin{split}
	\chibar_{ab} 
	&=
	\frac{1}{2} \slashed{\Pi}_a^{\phantom{a}c} \slashed{\Pi}_b^{\phantom{b}d} (\Lbar h_{cd})
	+ (\slashed{\nabla}_a \Lbar^e) \slashed{\Pi}_b^{\phantom{b}d}( m_{ed} + h_{ed})
	+ (\slashed{\nabla}_b \Lbar^e) \slashed{\Pi}_a^{\phantom{b}d}( m_{ed} + h_{ed})
	\end{split}
	\end{equation*}
	
	Now, this time we have
	\begin{equation*}
	\slashed{\nabla}_a \Lbar^0 = \slashed{\nabla}_a \Lbar^0_{(\text{small})}
	\end{equation*}
	and
	\begin{equation*}
	\slashed{\nabla}_a \Lbar^i = -\frac{\slashed{\nabla}_a^{\phantom{a}i}}{r} + \slashed{\nabla}_a \Lbar^i_{(\text{small})}
	\end{equation*}
	so
	\begin{equation*}
	\slashed{\nabla}_a \Lbar^b = -\frac{1}{r}\slashed{\Pi}_a^{\phantom{a}b} + \slashed{\nabla}_a \Lbar^b_{(\text{small})} +\frac{1}{r}\delta_0^b \slashed{\nabla}_a^{\phantom{a}0}
	\end{equation*}
	and we finally obtain
	\begin{equation*}
	\begin{split}
	\chibar_{ab} 
	&=
	-\frac{2}{r} \slashed{g}_{ab}
	+ \frac{1}{2} \slashed{\Pi}_a^{\phantom{a}c} \slashed{\Pi}_b^{\phantom{b}d} (\Lbar h_{cd})
	+ (\slashed{\nabla}_a \Lbar_{(\text{small})}^e) \slashed{\Pi}_b^{\phantom{b}d}( m_{ed} + h_{ed})
	+ (\slashed{\nabla}_b \Lbar_{(\text{small})}^e) \slashed{\Pi}_a^{\phantom{b}d}( m_{ed} + h_{ed})
	\\
	&\phantom{=}
	+ \frac{1}{r} \slashed{\Pi}_a^{\phantom{a}0} \slashed{\Pi}_b^{\phantom{b}d} (m_{0d} + h_{0d})
	+ \frac{1}{r} \slashed{\Pi}_b^{\phantom{b}0} \slashed{\Pi}_a^{\phantom{a}d} (m_{0d} + h_{0d})
	\\
	\\
	&=
	-\frac{2}{r} \slashed{g}_{ab}
	+ \frac{1}{2} \slashed{\Pi}_a^{\phantom{a}c} \slashed{\Pi}_b^{\phantom{b}d} (\Lbar h_{cd})
	+ (\slashed{\nabla}_a \Lbar_{(\text{small})}^e) (\slashed{\Pi}_{(\text{background})})_b^{\phantom{b}d} m_{ed}
	+ (\slashed{\nabla}_b \Lbar_{(\text{small})}^e) (\slashed{\Pi}_{(\text{background})})_a^{\phantom{a}d} m_{ed}
	\\
	&\phantom{=}
	+ \mathcal{O}\left(\frac{1}{r}|h_{(\text{rect})}|^2\right)
	+ \mathcal{O}\left(|\partial h_{(\text{rect})}| |h_{(\text{rect})}|\right)
	\end{split}
	\end{equation*}

\end{proof}

\begin{proposition}[The connection component \texorpdfstring{$\zeta$}{zeta} in the region \texorpdfstring{$r \leq r_0$}{r < r0}]
	\label{proposition zeta in r leq r0}
	
	In the region $r \leq r_0$, the rectangular components of $\zeta$ can be expressed as
	\begin{equation*}
	\zeta_a
	=
	(\zeta_{(\text{lin})})_a
	+ \mathcal{O}\left( |h_{(\text{rect})}| |\partial h_{(\text{rect})}| \right)
	\end{equation*}
	where $(\zeta_{(\text{lin})})_a$ are linear in the \emph{derivatives} of the rectangular components of $h$. 
	
\end{proposition}

\begin{proof}
	Recall that $\zeta$ is defined relative to the null frame as
	\begin{equation*}
	\zeta_A = g(\D_A L , \Lbar)
	\end{equation*}
	Expanding in rectangular coordinates we obtain
	\begin{equation*}
	\begin{split}
	\zeta_A
	&=
	g(\D_A (L^a \partial_a) , \Lbar^b \partial_b)
	\\
	&=
	g\left( (X_A L^a)\partial_a + (X_A)^c L^a \Gamma^d_{ac} \partial_d , \Lbar^b \partial_b \right)
	\\
	&=
	(X_A L^a) \Lbar^b g_{ab}
	+ (X_A)^c L^a \Lbar^b \Gamma^d_{ac} g_{db}
	\\
	&= (X_A L^a) L^b g_{ab} + \frac{1}{2} \left( (Lh)_{A\Lbar} + (X_A h)_{L\Lbar} - (\Lbar h)_{AL} \right)
	\end{split}
	\end{equation*}
	
	We can write this in a form that transforms covariantly under diffeomorphisms of the spheres:
	\begin{equation*}
	\zeta_{\slashed{\alpha}} = (\slashed{\nabla}_{\slashed{\alpha}} L^a) L^b g_{ab} + \frac{1}{2} \left( (Lh)_{\Lbar\slashed{\alpha}} + (\slashed{\nabla}_{\slashed{\alpha}} h)_{L\Lbar} - (\Lbar h)_{L\slashed{\alpha}} \right)
	\end{equation*}
	
	Recalling that the rectangular components of $L$ can be expressed in terms of $h$ in the region $r \leq r_0$ proves the proposition.
	
\end{proof}

\chapter{The reduced wave operator, the weak null structure and non-commutation with null frame}
\label{chapter weak null structure}

Now that we have developed some notation and established some basic properties of the wave operator and the null frame connection coefficients, we are finally in a position to discuss the weak null structure.

The ``weak null structure'' of the equations we are considering in this paper really consists of two conditions: one on the derivatives of the metric component $h_{LL}$, which must be shown to decay at a suitable rate, and another on the semilinear terms, which must obey a suitable hierarchy. In this chapter, we will elaborate on this hierarchy, which will impact the way in which error terms are considered in the deformation tensor calculations of chapter \ref{chapter deformation tensors}, and eventually in the energy estimate of chapters \ref{chapter framework for energy estimates} and \ref{chapter energy estimates}. There are also particular issues for the Einstein equations, which occur due to the fact that the null frame does not commute with the wave operator, which we will address here. In fact, we will develop a more general framework, which can be used to handle situations where the structure of the equations is only evident after a point-dependent change of variables.

Before turning to these issues, we first note that the errors arising from the quasilinear structure are related to the behaviour of the inverse foliation density $\mu$ in our framework: it is evident from proposition \ref{proposition transport mu} that, if $(\Lbar h)_{LL}$ decays at a slower rate than $r^{-1}$ then we cannot hope to control $\mu$. Even with the rate $(\Lbar h)_{LL} \sim \epsilon r^{-1}$ we can expect $\mu$ to grow at a rate $r^\epsilon$, and this is a cause of many of the difficulties encountered in this paper. Note, however, that use of the wave coordinate condition prohibits this growth, since in this case $(\Lbar h)_{LL}$ decays at a rate which is integrable in $r$. Hence $\mu$ is in fact uniformly bounded in $r$ if the wave coordinate condition holds.

\section{The reduced wave operator}
\label{section reduced wave operator}
 In order to elaborate on the hierarchy we require for the semilinear terms, we first introduce the \emph{reduced wave operator} $\tilde{\Box}_g$. This is related to the reduced wave operator $g^{ab}\partial_a \partial_b$ used in \cite{Lindblad2008}. The reason for introducing this operator is that, for wave equations of the form
\begin{equation*}
 \tilde{\Box}_g \phi = 0
\end{equation*}
we can expect the \emph{sharp} decay rate $\partial \phi \sim r^{-1}$. This is not true, in general, for equations of the form
\begin{equation*}
 \Box_g \phi = 0
\end{equation*}
where $g = g(\phi)$, for which we can only obtain decay at a rate $\partial \phi \sim r^{-1 + \epsilon}$. These sharp decay rates will play an important role in defining our semilinear structure.

\begin{definition}[The reduced wave operator]
 We define the \emph{reduced wave operator} as follows: for a scalar field $\phi$,
\begin{equation}
 \tilde{\Box}_g \phi := \Box_g \phi + \chi_{r_0} \omega \Lbar \phi
\end{equation}
Similarly, for an $S_{\tau,r}$-tangent tensor field $\phi$ we define
\begin{equation}
 \tilde{\slashed{\Box}}_g \phi := \slashed{\Box}_g \phi + \chi_{r_0} \omega \slashed{\D}_{\Lbar} \phi
\end{equation}

\end{definition}

Note that we could, instead, have defined the reduced wave operator as the operator
\begin{equation*}
 \Box_g + \chi_{r_0} \frac{1}{4}(\Lbar h)_{LL} \Lbar
\end{equation*}
since, ignoring terms with better behaviour, $\omega = \frac{1}{4}(\Lbar h)_{LL}$. However, we choose to use $\omega$ in our definition in order to simplify various formulae in the following sections, and to keep our discussion as geometric as possible.

\section{The semilinear hierarchy}
\label{section semilinear hierarchy}

Returning to the issue of the semilinear terms, we first observe that systems of wave equations of the kind we are considering can be written as
\begin{equation}
\label{equation set of wave equations}
 \tilde{\Box}_{g(\bm{\phi})} \phi_{(A)} = F_{(A)}(\bm{\phi}, \partial\bm{\phi})
\end{equation}
where $(A)$ labels the different fields, and bold text is used to refer to the vector of fields $(\phi_{(1)}, \phi_{(2)}, \ldots)$. For example, in the case of the the Einstein equations in harmonic coordinates, the fields in question are the rectangular components of the metric perturbation components $h_{ab}$. These are naturally labelled by pairs of rectangular indices, so we can write
\begin{equation}
 \tilde{\Box}_g h_{ab} = F_{ab}
\end{equation}
Recall that lower case latin indices refer to components in the rectangular frame, and so the field $h_{ab}$ is in fact a scalar field, despite appearances.

More formally, let $V$ be a finite dimensional vector space, and let $\bm{\phi}$ be a section of the trivial $V$-bundle over the manifold $\mathcal{M}$. Given a choice of frame for this vector bundle, i.e.\ a choice of vector fields $(v^{(1)}, v^{(2)},\ldots)$ which are linearly independent at every point, we can decompose $\bm{\phi}$ relative to this frame. That is, we can write
\begin{equation*}
 \bm{\phi} = \phi_{(1)} v^{(1)} + \phi_{2} v^{(2)} + \ldots
\end{equation*}
For example, for the Einstein equations in wave coordinates, we can choose the vector space $V$ to be the space of symmetric $4\times 4$ matrices. This space comes equipped with a canonical frame, namely, the matrices $M^{(ij)}$ with components
\begin{equation}
\label{equation frame symmetric matrices}
 M_{ab}^{(ij)} := \delta_a^i \delta_b^j
\end{equation}
and where we can restrict to $i \geq j$.

%In general, the weak null condition requires us to show improved behaviour for the component $h_{LL}$; specifically, we must establish an estimate of the form
%\begin{equation}
%\label{equation good decay of hLL exposition}
% |\Lbar h|_{LL} \lesssim \epsilon(1+r)^{-1}
%\end{equation}
%One way to do this is to appeal directly to the wave coordinate condition: see section \ref{section wave coordinate}, which allows us to express $(\Lbar h)_{LL}$ in terms of \emph{good} derivatives. This was the approach used in [[Lindblad Rodnianski]]. Note that an alternative approach, based on integrating along null geodesics, was used in [[Lindblad]], and this is the approach we will take in the more general case we consider. There are actually two reasons why it is necessary to establish an estimate of the form \eqref{equation good decay of hLL exposition}: first, to handle error terms arising from the quasilinear terms, and second, to handle certain semilinear terms.

The biggest problem arising from the semilinear terms $F_{(A)}$ is that they do not have the classical null structure. In other words, they contain terms of the form $(\partial \phi_{(B)})(\partial \phi_{(C)})$, in addition to the easily controlled ``null forms'' $(\partial \phi_{(B)})(\bar{\partial}\phi_{(C)})$, which contain at least one good derivative. In order to control these terms, we shall require that the fields $\phi_{(A)}$ obey a suitable hierarchy. To be precise, we require that the vector space $V$ splits into a (finite) direct sum of linear subspaces $V_{[i]}$, that is,
\begin{equation}
 V = \bigoplus_{[i] \geq 0} V_{[i]}
\end{equation}
We also define the related linear subspaces
\begin{equation}
W_{[i]} := \bigoplus_{[j] \leq [i]} V_{[j]}
\end{equation}
This set of nested linear subspaces $W_{[i]}$ will provide the hierarchy we need. We require that the subspaces $W_{[i]}$ are respected by the frame in the following sense. Let $\pi_V$ be the canonical projection onto the fibres. Then we require that, for all $(A)$, and all points $x \in \mathcal{M}$,
\begin{equation}
 \pi_V (v^{(A)}, x) \in W_{[i]}
\end{equation}
for some $[i]$, \emph{independent} of $x$. In other words, each frame vector field lies unambiguously in one of the subspaces $W_{[i]}$.

We can now define schematic notation for the fields as follows. We write $\Phi_{[i]}$ for the frame components of a field $\bm{\phi}$ such that
\begin{equation}
 \pi_V (\bm{\phi}, x) \in W_{[i]}
\end{equation}
i.e.\ $\Phi_{[i]}$ represents a set of fields $\phi_{(A)}$ such that $\phi_{(A)}v^{(A)}$ projects into the subspace $W_{[i]}$ at all points. We shall informally write $\phi_{(A)} \in \Phi_{[i]}$ to mean that $\phi_{(A)}v^{(A)} \in W_{[i]}$.

As mentioned above, there are really two components of the ``weak null structure'' which we require, one of which concerns the quasilinear structure (see section \ref{section reduced wave operator}), and the other of which concerns the semilinearities. The structure that we require from the semilinear terms is the following: the ``bad'' semilinear terms $(\partial \phi_{(B)})(\partial \phi_{(C)})$ must respect respect the hierarchy defined by the subspaces $W_{[i]}$. By this, we mean that the wave equation satisfied by a field $\phi_{(A)}$, where $\phi_{(A)}v^{(A)}$ (with no sum over $(A)$) projects into the subspace $W_{[i]}$, can be written schematically as
\begin{equation}
 \tilde{\Box}_{g} \phi_{(A)} = (\partial \Phi_{[i]})(\partial \Phi_{[0]}) + (\partial \Phi_{[i-1]})(\partial \Phi_{[i-1]}) + \text{ null forms } + \text{ cubic terms }
\end{equation}
In other words, the ``bad'' semilinear terms either involve pairs of fields at a lower level in the hierarchy, or else a field in the same level in the hierarchy paired with a field at the bottom level.

We also require a special structure for fields at the bottom level, i.e. the fields $\Phi_{[0]}$. These must satisfy an equation of the form
\begin{equation}
 \tilde{\Box}_{g} \phi_{(A)} = \text{ null forms } + \text{ cubic terms }
\end{equation}
In other words, there are no dangerous semilinear terms on the right hand side.

Finally, having stated the weak null structure as it arises in the structure of the semilinear terms, we can now state the relevant structure which will allow us to deal with the quasilinear terms. In fact, with the notation developed above, this is very easy: for the field $h_{LL}(\bm{\phi})$ we require the condition
\begin{equation}
 h_{LL}(\bm{\phi}) \in \Phi_{[0]}
\end{equation}
In other words, the field $h_{LL}$ can be expressed in terms of fields appearing at the bottom level of the hierarchy.

% There is one exception to this hierarchy, namely, the field\footnote{Assuming the ``non-degeneracy'' condition, that this term is, in fact, linear in the fields $\bm{\phi}$. If, instead, this term is quadratic or lower order in the fields then it does not require any special treatment. Moreover, in that case there are no difficulties with the quasilinear terms, and the entire content of the ``weak null condition'' concerns the semilinear structure.} $h_{LL}[\bm{\phi}]$. For this quantity we shall separately derive an ``independent'' pointwise decay estimate, which requires a slightly different structure for the semilinear terms. For consistency, we can set
% \begin{equation}
%  \Phi_{[-1]} := \{ h_{LL}[\bm{\phi}] \}
% \end{equation}
% Moreover, the semilinear structure is different for this term: we require
% \begin{equation}
% \label{equation schematic wave equation hLL}
%  \Box_{g} \left( h_{LL} \right) = - \frac{1}{4}\left(\Lbar h_{LL} \right)^2 + \text{ null forms } + \text{ cubic terms }
% \end{equation}
% Note that, ignoring better behaved terms $\frac{1}{4}\left(\Lbar h_{LL}\right)^2 = \omega (\Lbar h_{LL})$.

With the help of this hierarchy, we aim to construct a series of energy estimates, which become progressively ``worse'' as we ascend the hierarchy. The fact that these energy estimates become increasingly degenerate as we ascend the hierarchy will eventually lead to increasingly weak pointwise bounds for the fields which are higher in the hierarchy. Thus, when deriving the energy estimates, we must be able to assume worse pointwise bounds for the fields which are higher in the hierarchy. That it is still possible to derive the necessary energy estimates is due precisely to the fact that the semilinear terms respect the hierarchy, as outlined above.

Of particular note in the definition of this hierarchy is the addition of the term $\omega (\Lbar \phi)$ to the wave operator in order to form the reduced wave operator. In particular, this plays an important role when considering the equation for the field $h_{LL}$. At first sight this term appears extremely dangerous, since we have $\omega \sim (\Lbar h)_{LL}$, and wave equations of the form
\begin{equation}
\label{equation John blowup}
 \Box \phi = (\Lbar \phi)^2
\end{equation}
were shown to exhibit finite-time blowup for all nontrivial initial data \cite{John1981}. Note, however, that the wave operator appearing in \eqref{equation John blowup} is the flat wave operator and not the geometric wave operator $\Box_g$. Moreover, the geometric wave operator, expressed in the null frame, already contains terms of this type (see \eqref{equation wave operator null}), and this term is ``cancelled'' by the modification of the wave operator to the reduced wave operator.

\section{Relation to the weak null condition}
The ``weak null condition'' was not originally stated in terms of a hierarchy of the kind described above. Rather, it was has been stated as follows: \emph{the asymptotic system possesses global solutions for small initial data}. To form the asymptotic system, we simply neglect the terms which are expected to have better behaviour, i.e.\ those terms involving good derivatives. This leads to the following conjecture:
\begin{conjecture}[Weak null conjecture, na\"ive version]
	\label{conjecture weak null naive}
	If the asymptotic system corresponding to a system of wave equations admits global solutions for sufficiently small initial data, then so does the system of wave equations.
\end{conjecture}
Actually, this conjecture is almost certainly false, since we can construct systems where the solutions to the asymptotic system grow far too rapidly (see section \ref{section other asymptotic systems}). Instead, the weak null conjecture is often stated as follows:
\begin{conjecture}[Weak null conjecture]
	\label{conjecture weak null}
	If the asymptotic system corresponding to a system of wave equations admits global solutions for sufficiently small initial data, and if those solutions grow no faster than $(r\Lbar \phi) \sim r^{C\epsilon}$, then the system of wave equations also admits global solutions for sufficiently small initial data.
\end{conjecture}

Now, a natural question to ask is the following: \emph{do all systems of nonlinear wave equations which satisfy the weak null condition obey a hierarchy of the form given above?} We can immediately answer this in the negative by providing the following (fairly trivial) counterexample: take a system of wave equations which does obey such a hierarchy, and make a linear change of frame. The resulting system will not obey the hierarchy, despite possessing global solutions.

To be more explicit, consider the following example of a semilinear system obeying the weak null condition:
\begin{equation*}
 \begin{split}
  \Box \phi_{(0)} &= 0 \\
  \Box \phi_{(1)} &= (\Lbar \phi_{(1)})(\Lbar \phi_{(0)})
 \end{split}
\end{equation*}
This clearly obeys a hierarchy of the kind considered above, where we have
\begin{equation*}
 \begin{split}
  V_{[0]} &= \left\{\phi_{(0)} \right\} \\
  V_{[1]} &= \left\{\phi_{(1)} \right\} 
 \end{split}
\end{equation*}
However, if we define the field $\tilde{\phi}_{(0)} := \phi_{(0)} + \phi_{(1)}$ then the system above is equivalent to the system
\begin{equation*}
 \begin{split}
  \Box \tilde{\phi}_{(0)} &= (\Lbar \tilde{\phi}_0)(\Lbar \phi_1) - (\Lbar \phi_1)^2 \\
  \Box \phi_{(1)} &= (\Lbar \tilde{\phi}_0)(\Lbar \phi_1) - (\Lbar \phi_1)^2 \\
 \end{split}
\end{equation*}
This system evidently \emph{does not} obey a hierarchy of the kind discussed above; the bad semilinear terms appearing on the right hand side involve every field. Interestingly, for this system, we have the asymptotic behaviour $(\partial \tilde{\phi}_0) \sim r^{-1 + \epsilon}$ and $(\partial \phi_1) \sim r^{-1 + \epsilon}$. However, one cannot obtain this result through the usual bootstrap type argument; we instead need to ``undo'' the linear transformation and then make bootstrap assumptions on the original fields $\phi_{(0)}$ and $\phi_{(1)}$.

We could overcome this objection by broadening our class of systems to include all systems which can be brought into a suitable heirarchical structure by means of a linear transformation on the frame fields. Of course, our proof easily applies in this case. Note that we can restrict to \emph{linear} transformations of this kind without loss of generality, since we are neglecting cubic and higher order terms.

Even if we broaden the class of allowed systems in this way, we can still find systems of nonlinear wave equations which do not obey a hierarchy of the kind we require: see section \ref{section other asymptotic systems}. Unusually, some of the asymptotic systems associated with the systems we construct in section \ref{section other asymptotic systems} have very rapidly growing solutions, meaning that, if the asymptotic system accurately predicts the behaviour of the fields, then we would expect these fields to grow exponentially or even super-exponentially. Indeed, it is very doubtful in these cases whether the asymptotic system really does give an accurate prediction for the behaviour of the fields, since in these cases cubic and higher order terms would be expected to have large effects. By considering such cases it might be possible to construct a couterexample to the original conjecture \ref{conjecture weak null naive}.

It would be interesting to know whether \emph{all} systems which obey the weak null condition and, in addition, only admit solutions with the asymptotic behaviour $(\partial \phi) \sim r^{-1 + \epsilon}$, can be transformed into a system admitting a hierarchy of the kind described in section \ref{section semilinear hierarchy} by means of a linear transformation.

We first establish that systems with the hierarchy defined in the section above obey the weak null condition, that is, the associated asymptotic systems admit global solutions for all sufficiently small initial data. Moreover, the solutions to these asymptotic systems all predict the behaviour $(\partial \phi_{(A)}) \sim r^{-1 + \epsilon}$ for general $\phi_{(A)}$, and $(\partial \phi_{(A)}) \sim r^{-1}$ for $\phi_{(A)} \in W_{[0]}$.

\begin{proposition}[Systems obeying the hierarchy of section \ref{section semilinear hierarchy} obey the weak null condition]
 Let $\{ \phi_{(A)} \}$ be a collection of fields obeying nonlinear wave equations which admit a hierarchy of the kind described in section \ref{section semilinear hierarchy}. Then the associated asymptotic system admits global solutions for all sufficiently small data, where by ``small'' we mean
\begin{equation}
 \sup_{(A)}\sup_{t = 0} \left( |\Lbar \phi_{(A)}| \, , \, |r\Lbar \phi_{(A)}| \right) \leq \epsilon
\end{equation}
for some sufficiently small $\epsilon$.

 Moreover, the solutions to the asymptotic system obey the bounds
\begin{equation}
 \begin{split}
  |\Lbar \phi_{(A)}| &\lesssim \epsilon(1+r)^{-1} \quad \text{ for} \quad \phi_{(A)} \in V_{[0]} \\
  |\Lbar \phi_{(A)}| &\lesssim \epsilon(1+r)^{-1 + C\epsilon} \quad \text{ for} \quad \phi_{(A)} \in V_{[n]} \text{ , } n \geq 1 \\
 \end{split}
\end{equation}
where $C$ is some suitably large constant, and the initial data satisfies the bound
$|r\Lbar \phi_{(A)}| \big|_{t=0} \leq \epsilon$.
\end{proposition}

\begin{proof}
 Recall that the $r$-weighted wave operator was expressed in proposition \ref{proposition scalar wave operator} as
\begin{equation*}
 r\Box_g \phi = -L(r\Lbar \phi) + r\slashed{\Delta}\phi - \frac{1}{2}(\tr_{\slashed{g}}\chi_{(\text{small})} + 2\omega) r\Lbar \phi - (\tr_{\slashed{g}} \chibar)rL\phi - \zeta^A r\slashed{\nabla}_A \phi 
\end{equation*}
 Hence, if $\tilde{\Box}_g \phi_{(A)} = F_{(A)}$ then we have
\begin{equation*}
 -L(r\Lbar \phi_{(A)}) + r\slashed{\Delta}\phi_{(A)} - \frac{1}{2}(\tr_{\slashed{g}}\chi_{(\text{small})}) r\Lbar \phi_{(A)} - (\tr_{\slashed{g}} \chibar)rL\phi_{(A)} - \zeta^A r\slashed{\nabla}_A \phi_{(A)} = rF_{(A)} 
\end{equation*}
Note that the use of the reduced wave operator has cancelled the term $-\omega r\Lbar \phi$ on the left hand side.

Let us write the ``bad'' semilinear terms appearing in the $F_{(A)}$ as
\begin{equation*}
 F_{(A)} = F_{(A)}^{\phantom{(A)}(BC)}(\Lbar \phi_{(B)})(\Lbar \phi_{(C)}) + \text{null forms} + \text{cubic terms}
\end{equation*}
for some constants $F_{(A)}^{\phantom{(A)}(BC)}$. Let us also define the fields
\begin{equation*}
 r\Lbar \phi_{(A)} := \xi_{(A)}
\end{equation*}

Now, the asymptotic system corresponding to the equation above is a system of transport equations for the fields $\Phi_{(A)}$ along the integral curves of $L$, obtained by dropping second derivatives which contain two good derivatives, as well as all terms involving null forms (or terms with equivalent behaviour, such as $\zeta(\bar{\partial}\phi)$ or $(\tr_{\slashed{g}}\chi)(\partial \phi)$), and additionally we drop all cubic terms. Finally, we allow ourselves to cut-off the semilinear terms for all $r \leq 1$. This leads to an easier treatment of the axis $r = 0$, and can be justified by the availability of elliptic estimates to control the behaviour of fields in this region in the full problem.

Hence, the asymptotic system corresponding to the system above is
\begin{equation*}
 -L(\xi_{(A)}) = r^{-1} \chi_0(r)F_{(A)}^{\phantom{(A)}(BC)}\xi_{(B)}\xi_{(C)}
\end{equation*}
where $\chi_0(r)$ is the smooth cut-off function introduced in equation \eqref{equation cut off chi0}.

If we use the coordinates $(u, r, \vartheta^A)$ then $L = \frac{\partial}{\partial r}$. Defining the new variable
\begin{equation*}
 y := \log r
\end{equation*}
the asymptotic system becomes
\begin{equation}
 \label{equation asymptotic system 1}
 \frac{\partial}{\partial y} \xi_{(A)} = - \chi(e^y) F_{(A)}^{\phantom{(A)}(BC)}\xi_{(B)}\xi_{(C)}
\end{equation}

The hierarchy discussed in section \ref{section semilinear hierarchy} states that, if $\phi_{(A)} \in \Phi_{[i]}$ then $F_{(A)}^{\phantom{(A)}(BC)} = 0$ unless one of the following conditions is satisfied:
\begin{enumerate}
  \item Both $\phi_{(B)}$ and $\phi_{(C)} \in \Phi_{[i-1]}$ \quad \emph{or}
  \item $\phi_{(B)} \in \Phi_{[i]}$ and $\phi_{(C)} \in \Phi_{[0]}$, or $\phi_{(B)} \in \Phi_{[0]}$ and $\phi_{(C)} \in \Phi_{[i]}$
\end{enumerate}
Moreover, if $\phi_{(A)} \in \Phi_{[0]}$ then $F_{(A)}^{\phantom{(A)}(BC)} = 0$ for all $B, C$.

Now, with the help of this hierarchy, we can solve system \eqref{equation asymptotic system 1} with the use of a boostrap. We first set
\begin{equation}
 \tilde{F} := \sup_{A, B, C} \left| F_{(A)}^{(BC)} \right|
\end{equation}
and we consider a particular integral curve of $L$. If this curve intersects the axis $r = 0$ (which corresponds to $y = -\infty$ then regularity of $\Lbar \phi_{(A)}$ requires that $\xi_{(A)} = 0$ at $y = -\infty$. Moreover, this integral curve will enter the region $y \leq 0$, after which the right hand side of equation \eqref{equation asymptotic system 1} vanishes. Hence the unique regular solution along this integral curve is $\xi_{(A)} = 0$ for all $(A)$.

Henceforth we consider an integral curve of $L$ which intersects the initial data surface $t = 0$ at $y = y_0$. The bootstrap assumptions we make along this curve are the following:
\begin{equation*}
 \begin{split}
  |\xi_{(A)}| &\leq 2\inf\{\epsilon e^{y_0} \, , \, \epsilon\} \quad \text{ for} \quad \phi_{(A)} \in \Phi_{[0]} \\
  |\xi_{(A)}| &\leq 2\inf\{\epsilon e^{y_0} \, , \, \epsilon\} + \epsilon \exp\left(\tilde{C}\tilde{F}2^{n}\epsilon (y - y_0)\right) \quad \text{ for} \quad \phi_{(A)} \in \Phi_{[n]} \\
 \end{split}
\end{equation*}
for some large constant $\tilde{C}$ which we determine below. We can check that these bootstrap bounds hold initially. Since the initial data satisfies
\begin{equation*}
 \begin{split}
  |\Lbar \phi_{(A)}|\big|_{t=0} &\leq \epsilon \\
  |r\Lbar \phi_{(A)}|\big|_{t=0} &\leq \epsilon
 \end{split}
\end{equation*}
this implies that
\begin{equation*}
  |\xi_{(A)}|\big|_{t=0} \leq \inf \{ e^y \epsilon \, , \, \epsilon \}
\end{equation*}
so it is easy to see that the bootstrap assumptions hold initially and, by continuity, at least for all sufficiently small values of $y$. 

Now, we can use the equations \eqref{equation asymptotic system 1} to improve the bootstrap. If $\phi_{(A)} \in \Phi_{[0]}$ then the asymptotic system for $\phi_{(A)}$ is simply
\begin{equation*}
 \frac{\partial}{\partial y} \xi_{(A)} = 0
\end{equation*}
and so $\xi_{(A)}$ is just constant and equal to its initial value. In particular this improves the bootstrap assumption by a factor of $1/2$.

Now consider a field $\phi_{(A)} \in \Phi_{[n]}$. For simplicity let us assume that $y_0 \leq 0$; the case $y_0 > 0$ follows similarly. Making use of the bootstrap assumptions, together with the bound $\xi_{(A)} \leq \epsilon$ for $\phi_{(A)} \in \Phi_{[0]}$ obtained above, we find that the derivative of $\xi_{(A)}$ satisfies
\begin{equation*}
\label{equation asymptotic system derivative bound}
 \begin{split}
  \left| \frac{\partial}{\partial y} \xi_{(A)} \right| &\leq
  \tilde{F} \left( 2\epsilon^2 e^{y_0} + \epsilon^2 \exp\left( \tilde{C}\tilde{F}w^n \epsilon(y-y_0) \right)  + \left( 2 \epsilon e^{y_0} + \epsilon \exp\left( \tilde{C}\tilde{F} 2^{n-1} \epsilon(y-y_0) \right) \right)^2 \right) \\
  &\leq \tilde{F} \epsilon^2 \left( (2 e^{y_0} + 4e^{2y_0}) + 2\exp\left( y_0 + \tilde{C} \tilde{F} 2^{n-1} \epsilon(y-y_0) \right) + 2\exp\left( \tilde{C}\tilde{F} 2^n \epsilon(y-y_0) \right) \right)
 \end{split}
\end{equation*}
from which it follows that
\begin{equation*}
 \begin{split}
  |\xi_{(A)}| &\leq \epsilon e^{y_0} + 6\tilde{F}\epsilon^2 (y-y_0) + 4(\tilde{C}^{-1})2^{-n}\epsilon \exp\left( y_0 + \tilde{C} \tilde{F}2^{n-1}\epsilon (y-y_0) \right) \\
  &\phantom{=} + (\tilde{C}^{-1})2^{1-n}\epsilon \exp\left(\tilde{C}\tilde{F}2^n\epsilon(y-y_0)\right)
 \end{split}
\end{equation*}
It is now easy to check that, by taking $\epsilon$ sufficiently small and $\tilde{C}$ sufficiently large, the bootstrap assumptions can be improved.

A continuity argument now shows that the bootstrap assumptions in fact hold for all values of $y$, that is, all along the integral curve of $L$. 

\end{proof}

\section{Other asymptotic systems}
\label{section other asymptotic systems}

As mentioned at the start of this chapter, there are asymptotic systems which do not obey a semilinear hierarchy of the kind described in section \ref{section semilinear hierarchy}, but which nevertheless admit global solutions for suitably small initial data. Some of these can be transformed into a form admitting of a suitable hierarchy by a linear transformation on the fields. However, there are other systems which cannot be transformed in this way, and which lead to very different asymptotic behaviour. In this section we will construct such a system, and show that it leads to the behaviour $\Lbar \phi \sim \epsilon e^{r^{\epsilon}}$. We cannot treat the nonlinear wave equations associated with such systems by our methods, and indeed it is highly doubtful whether they in fact possess global solutions for small data, especially if they include cubic nonlinearities which would be expected to dominate.

\begin{proposition}[A nonlinear wave equation with an asymptotic system admitting exponentially growing solutions]
Consider the following set of (semilinear) nonlinear wave equations:
\begin{equation}
 \begin{split}
  \Box \phi_{(1)} &= 0 \\
  \Box \phi_{(2)} &= -4(T \phi_{(1)})(T \phi_{(2)}) \\
  \Box \phi_{(3)} &= -4(T \phi_{(2)})(T \phi_{(3)})
 \end{split}
\end{equation}
Then the associated asymptotic system admits global solutions for sufficiently small initial data. Using the coordinate $u = t-r$ on the initial surface $t = 0$, suppose that the initial data for the asymptotic system are given in terms of $\xi_{(A)} = r\Lbar \phi_{(A)}$ as
\begin{equation*}
 \xi_{(A)}\big|_{t = 0}(u, \vartheta^A) = \epsilon \mathring{\xi}_{(i)}(u, \vartheta^A)
\end{equation*}
Then the solution of the asymptotic system is the following: for $u > 0$, we have $\xi_{(A)} = 0$, while for $u < 0$ we have
\begin{equation*}
 \begin{split}
  \xi_{(1)}(u,r,\vartheta^A) &= \epsilon\mathring{\xi}_{(1)}(u,\vartheta^A) \\
  \xi_{(2)}(u,r,\vartheta^A) &= \epsilon\mathring{\xi}_{(2)}(u,\vartheta^A)\left(\frac{r}{-u}\right)^{\epsilon \mathring{\xi}_{(1)}(u,\vartheta^A)} \\
  \xi_{(3)}(u,r,\vartheta^A) &= \epsilon\mathring{\xi}_{(3)}(u,\vartheta^A) \exp\left( (-u)^{\epsilon \mathring{\xi}_{(1)}(u, \vartheta^A)}\mathring{\xi}_{(2)}(u, \vartheta^A) \left(\mathring{\xi}_{(1)}(u,\vartheta^A)\right)^{-1} r^{\epsilon \mathring{\xi}_{(1)}(u, \vartheta^A)} \right)
 \end{split}
\end{equation*}
In particular, if there is some point on the initial data surface such that both $\mathring{\xi}_{(1)}$ and $\mathring{\xi}_{(2)}$ are strictly positive, then $|\xi_{(3)}| \sim \epsilon e^{r^\epsilon}$

\end{proposition}

\begin{proof}
 The asymptotic system associated with the system in the proposition is
 \begin{equation*}
  \begin{split}
    L\xi_{(1)} &= 0 \\
    L\xi_{(2)} &= r^{-1} \xi_{(1)}\xi_{(2)} \\
    L\xi_{(3)} &= r^{-1} \xi_{(2)}\xi_{(3)}
  \end{split}
 \end{equation*}
 It is now easy to verify that the solution provided in the proposition solves these equations with initial data which vanishes for integral curves of $L = \partial/\partial r$ which intersect the axis, and with the initial data specified in the proposition at $t = 0$ along integral curves of $L$ which intersect this hypersurface.
\end{proof}

\begin{remark}
It does not appear likely that the nonlinear wave equations with asymptotic systems of this kind would possess global solutions in general. The original motivation for conjecture \ref{conjecture weak null naive} was that both ``good'' derivatives and cubic terms are expected to behave better than the ``bad'' derivatives $\Lbar \phi$. Assuming that the asymptotic system provides a good approximation for the behaviour of the bad derivatives, it is then hoped that energy estimates are sufficient to prove that the good derivatives indeed behave well. One can then justify ``dropping'' these terms from the asymptotic system.

There are several problems for this approach when considering systems of the type considered in this section. First, cubic and higher order terms cannot be expected to behave ``better'', when the field $\phi$ and its derivatives are expected to grow exponentially, rather than to decay. Indeed, by including a term of the form $(\Lbar\phi_{(3)})^3$ in the asymptotic system it is fairly easy to see that solutions will blow up in finite time. Additionally, even if the good derivatives do behave ``better'' by the expected factor of $r^{-1}$, terms which are quadratic in the good derivatives might still large enough to have a dramatic effect on the solution, and possibly to lead to finite-time blowup. For this reason, conjecture \ref{conjecture weak null naive} must be modified to conjecture \ref{conjecture weak null}.
\end{remark}

We have already seen examples of systems which do obey the conditions of conjecture \ref{conjecture weak null}, but which do not obey the hierarchical null condition. However, the examples we have given up to now can all be written in terms of other variables which \emph{do} obey the hierarchical null condition, by changing the basis sections. Below we give an example of a system where this is not possible, and yet the weak null condition is still obeyed.

\begin{proposition}[A system with bounded solutions to the asymptotic system but lacking the weak null hierarchy]
Consider the system
\begin{equation*}
\begin{split}
	\Box \phi_{(1)} &= 4\frac{I_3 - I_2}{I_1} (T\phi_{(2)})(T\phi_{(3)})
	\\
	\Box \phi_{(2)} &= 4\frac{I_1 - I_3}{I_2} (T\phi_{(3)})(T\phi_{(1)})
	\\
	\Box \phi_{(3)} &= 4\frac{I_2 - I_1}{I_3} (T\phi_{(1)})(T\phi_{(2)})
	\\
\end{split}
\end{equation*}
for positive constants $I_1$, $I_2$, $I_3$. Then the corresponding asymptotic system admits global solutions for all small intial data. Moreover, for a generic choice of the constants $I_1$, $I_2$, $I_3$, the system cannot be reduced to one with the semilinear hierarchy by a change of basis sections.
\end{proposition}

\begin{proof}
	If we define $r\Lbar \phi_{(i)} := \xi_{(i)}$, then the asymptotic system is
	\begin{equation*}
	\begin{split}
	I_1 L\xi_{(1)} = (I_3 - I_2) r^{-1} \xi_{(2)}\xi_{(3)}
	\\
	I_2 L\xi_{(2)} = (I_1 - I_3) r^{-1} \xi_{(3)}\xi_{(1)}
	\\
	I_3 L\xi_{(3)} = (I_2 - I_1) r^{-1} \xi_{(1)}\xi_{(2)}
	\end{split}
	\end{equation*}
	Writine $rL = \partial/\partial(\log r)$ we obtain
	\begin{equation*}
	\begin{split}
	I_1 \frac{\partial}{\partial \log r}\xi_{(1)} = (I_3 - I_2) \xi_{(2)}\xi_{(3)}
	\\
	I_2 \frac{\partial}{\partial \log r}\xi_{(2)} = (I_1 - I_3) \xi_{(3)}\xi_{(1)}
	\\
	I_3 \frac{\partial}{\partial \log r}\xi_{(3)} = (I_2 - I_1) \xi_{(1)}\xi_{(2)}
	\end{split}
	\end{equation*}
	These are the Euler equations for the motion of a rigid body with moments of inertia $I_1$, $I_2$, $I_3$ about the principle axes. The solutions are therefore bounded, i.e.\ the asymptotic system predicts the behaviour
	\begin{equation*}
	\Lbar \phi_{(i)} \sim r^{-1}
	\end{equation*}
	On the other hand, for generic $I_1$, $I_2$, $I_3$ there are no closed-form solutions, indicating that the system cannot be transformed into a simple system obeying the weak null condition.
\end{proof}

For this particular system, it may be possible to prove global existence for the corresponding set of wave equations, using the fact that we expect each field to obey the sharp $r^{-1}$ decay rate\footnote{We would also have to prove a kind of stability result for the Euler equations, showing that, if a torque is applied that satisfies particular bounds, then the resulting solutions still remain bounded.}. However, it is clear that the dynamics of systems of ODEs with quadratic nonlinearities can be quite complicated, and we expect that even more intricate systems can be designed, which satisfy the weak null condition but for which our methods cannot be easily adapted. However, we have not been able to construct such a system.

Finally, we remark on the ``quasilinear condition'' $h_{LL} \in \Phi_{[0]}$. This appears necessary in order to avoid the formation of shocks in finite time. Specifically, the inverse foliation density satisfies the transport equation
\begin{equation*}
L \log \mu = (\Lbar h)_{LL} + \text{good terms}
\end{equation*}
The good terms are easy to control. If we only have the bound $(\Lbar h)_{LL} \sim r^{-1+C\epsilon}$, then we can expect the foliation density $\mu$ to grow exponentially in $r$, but then the error terms involving $\mu$ and its derivatives would be completely uncontrollable. We can also compare this with the work of \cite{Speck2016a}, which showed \emph{in the context of a scalar wave equation} that, generically, shocks will form in the case where we cannot obtain a bound like $(\Lbar h)_{LL} \sim r^{-1}$.

\section{Non-commutation with the null frame}

As mentioned in the introduction, one of the applications we have in mind for the methods developed in this paper is to the Einstein equations in wave coordinates. In this case, the fields $\phi_{(A)}$ which obey nonlinear wave equations are the rectangular components of the metric perturbations $h_{ab}$. The associated semilinear terms $F_{ab}$ do obey a hierarchy of the kind described above, but this is only evident when the equations are contracted with the rectangular components of the null frame $L^a$, $\Lbar^a$ and $\slashed{\Pi}_\mu^{\phantom{\mu}a}$.

In the language developed above, this means that we need to change the frame $\{v^{(A)}\}$ from the canonical frame for symmetric matrices \eqref{equation frame symmetric matrices} to one based on the null frame vector fields $L$, $\Lbar$ and $X_A$. We have already mentioned several such changes of basis, but previously in this chapter these changes were assumed to be ``global'' in the sense that the new frame was related to the old frame by some fixed linear map. On the other hand, changing basis to one which depends on the frame vector fields is a ``point dependent'' change of frame, represented by a linear map which depends on the base manifold $\mathcal{M}$. Hence the change-of-frame operator does not commute with the null frame.

To make this more precise, let $w^{(A')} = M_{(A)}^{(A')}(x) v^{(A)}$ be an alternative set of frame vector fields, so that the (point-dependent) change of frame matrix $M_{(A)}^{(A')}$ is non-degenerate at every point $x \in \mathcal{M}$. Then we can write
\begin{equation*}
 \phi_{(A)} = M_{(A)}^{(A')}(x) \phi_{(A')}
\end{equation*}
Since the change of frame matrix $M_{(A)}^{(A')}$ depends on the point $x \in \mathcal{M}$, this change of frame will not commute with the wave operator. Specifically, if $\phi_{(A)}$ satisfies the wave equation \eqref{equation set of wave equations} then $\phi_{(A')}$ will satisfy
\begin{equation*}
 \Box_g \phi_{(A')} = (M^{-1})_{(A')}^{(A)} F_{(A)} - (M^{-1})_{(A')}^{(A)} \left( \Box_g M_{(A)}^{(B')} \right) \phi_{(B')} - 2(M^{-1})_{(A')}^{(A)} \left( \D^{\alpha} M_{(A)}^{(B')} \right) \left( \D_\alpha \phi_{(B')} \right)
\end{equation*}
The second and third terms on the right hand side reflect the fact that the change of frame matrix $M$ does not commute with the wave operator.

We will not ``fully'' change frame in the sense described above, where we use an alternative set of fields $\phi_{(A')}$ to describe the state of the system, but only ``partially'' change frame. The reason for this is that the change of frame we require (namely, a change to the null frame) depends on the fields $\bm{\phi}$ themselves. In other words, the matrix $M_{(A)}^{(B')}$ actually depends on the fields $\bm{\phi}$, in addition to depending directly on the spacetime point $x$. This means that the second order term $\Box_g M_{(A)}^{(B')}$ would introduce additional, leading order error terms which could not be controlled. Instead, we will ``partially'' change frame, which means that we will first take take derivatives of the fields in the original frame $v^{(A)}$, and then change frame afterwards. In other words, we will estimate quantities like
\begin{equation*}
 \left| \left(\partial \phi_{(A)} \right)(M^{-1})_{(A)}^{(A')} \right|^2
\end{equation*}
In order to estimate quantities of this kind, we will have to define a kind of modified energy momentum tensor, see section \ref{section energy estimates involving a point-dependent change of basis}.

We have already seen hints of this kind of structure in our definition of quantities such as $(\partial h)_{LL}$. In this case, we are changing frame from the rectangular frame to the null frame associated with symmetric matrices, but the change of frame occurs outside the derivative operator.

% To aid in making the calculations explicit, from now on we will only consider a change of basis from a frame adapted to symmetric tensors, labelled by the rectangular indices, to a null frame based on the vector fields $L$, $\Lbar$ and the projection operator $\slashed{\Pi}$. However, it should be noted that the methods developed in this paper could in principle be applied to other systems of wave equations, in which the hierarchical structure described at the start of this chapter applies only after a suitable (partial) change of frame, provided that the associated error terms can be controlled.

\section{Compatibility with the radial normalisation condition}
Recall that our geometric set-up required the property \ref{property radial normalisation}, i.e.\ that the radial vector field $R^\mu = (g^{-1})^{\mu\nu} (\upd r)_\nu$ has unit norm with respect to the dynamical metric $g$. However, the systems we have introduced in this section need not satisfy this property, so it would appear that our geometric set-up does not apply.

In this section we will show that any system of quasilinear wave equations of the form given in equation \eqref{equation set of wave equations}, with a metric $g(\bm{\phi})$ not necessarily obeying the radial normalisation property \ref{property radial normalisation}, may be recast (by means of a conformal transformation) into a system of quasilinear wave equations of the same form but with a metric $\mathring{g}(\bm{\phi})$ which \emph{does} obey the radial normalisation property. Moreover, if the original system satisfies the weak null condition, then so does the new system.

Define the scalar field $\Omega_{(R)}$ by
\begin{equation}
  \begin{split}
    (\Omega_{(R)})^{-2} &:= \chi_{(r_0)}(r)(g^{-1})(\upd r, \upd r) + (1 - \chi_{(r_0)}) \\
      &\phantom{:}= 1 -\chi_{(r_0)}\left(\frac{x^i x^j}{r^2} h_{ij} + \mathcal{O}(h^2) \right)
  \end{split}
\end{equation}
where $\chi_{(r_0)}(r)$ is the smooth cut-off function defined in equation \eqref{equation cut off functions}. $\Omega_{(R)}$ is evidently determined by the fields $h_{ij}$, which in turn depend (say, linearly - higher order terms are easier to handle) on the fields $\phi_{(A)}$.

We now define the conformally rescaled metric $\mathring{g}$ by
\begin{equation}
 \mathring{g} := (\Omega_{(R)})^2g
\end{equation}
Note that, by construction, $\mathring{g}(R,R) = 1$ in the region $r \geq r_0$.

To understand the action of this conformal transformation on the \emph{reduced} wave operator $\tilde{\Box}_g$, we must first understand its action on the vector fields $L$, $\Lbar$, the foliation density $\mu$ and the connection coefficient $\omega$. Recall that $L$ is defined as an outgoing null vector field, orthogonal to the spheres, and normalised by $L(r) = 1$. Vector fields which are null with respect to $g$ are also null with respect to $\mathring{g}$, and orthogonality is also preserved, so we find that $L$ is invariant under this transformation.

On the other hand, recall that $\Lbar$ was defined\footnote{Unlike the vector field $L$, it is a \emph{consequence} of the radial normalisation condition that $\Lbar(r) = -1$, so we cannot assume that this is true yet.} by as the future directed, null vector field, orthogonal to the spheres and satisfying $g(L, \Lbar) = -2$. Hence, under a conformal rescaling, we need to use a new vector field $\mathring{\Lbar}$ in place of $\Lbar$, defined by
\begin{equation*}
 \mathring{\Lbar} := (\Omega_{(R)})^{-2} \Lbar
\end{equation*}

Next, recall that $\mu$ is defined by the equation $\mu = -1/(g^{-1}\left(\upd u, \upd r)\right)$. The coordinate function $r$ is invariant under the conformal map since it is defined with respect to the ``background'', independently of the metric $g$. Additionally, the coordinate $u$ is invariant under the conformal map, since $L$ is invariant, and this vector field is all that is required to construct $u$. This means that the new foliation density $\mathring{\mu}$ is related to the old foliation density $\mu$ by
\begin{equation*}
\mathring{\mu} := (\Omega_{(R)})^{2} \mu
\end{equation*}

Finally, we can compute the new connection coefficient $\mathring \omega$. Recall that $\omega$ satisfies $L\log \mu = \omega$. Using the expression for $\mathring{\mu}$, together with the invariance of the vector field $L$, we find
\begin{equation*}
\mathring{\omega} = 2L\log (\Omega_{(R)}) + \omega
\end{equation*}
Note that $L\log (\Omega_{(R)}) \sim L \phi$, so this behaves like a good derivative.

Next, we need to write the wave operator with respect to the rescaled metric $\mathring{g}$ instead of the metric $g$. For this we can use a standard formula: we obtain
\begin{equation}
 \Box_{\mathring{g}} \phi_{(A)} = (\Omega_{(R)})^{-2} \Box_g \phi + 2 \Omega^{-4} (\mathring{g}^{-1})^{\mu\nu}(\D_\mu \log\Omega_{(R)})(\D_\nu \phi_{(A)})
\end{equation}

Putting this all together, we find that if the fields $\phi_{(A)}$ satisfy the quasilinear wave equations $\tilde{\Box}_g \phi_{(A)} = F_{(A)}$, then, written with respect to the conformally rescaled metric, we have
\begin{equation*}
\tilde{\Box}_{\mathring{g}} \phi_{(A)}
=
\Omega^{-2} F_{(A)}
- 2(L\log \Omega_{(R)}) (\mathring{\Lbar} \phi_{(A)})
+ 2(\Omega_{(R)})^{-4} (\mathring{g}^{-1})^{\mu\nu} (\partial_\mu \log \Omega_{(R)}) (\partial_\nu \phi_{(A)})
\end{equation*}

Note that $\Omega_{(R)} = 1 + \mathcal{O}(\bm{\phi})$. This means that all of the additional terms on the right hand side are either cubic terms, or else they involve at most \emph{one} $\mathring{\Lbar}$ derivative. Therefore, if the original system satisfied the weak null condition, then so will the new system. Furthermore, the hierarchical structure of the semilinear terms is also unchanged.

In summary, if we are given a system of wave equations satisfying the hierarchical null condition, then, by writing these equations with respect to a conformally rescaled metric we can impose the radial normalisation condition, and this procedure does not change the null structure of the equations.

% 
% Note that the term $(\Lbar h)_{LL}$ \emph{must} have decay in $r$ which is at least as fast as $r^{-1}$ in order to obtain global existence, since otherwise the inverse foliation density $\mu$ cannot be controlled, and shocks might form in finite time. If we assume that $(\Lbar h)_{LL} \sim \Lbar \phi_{(A)}$, then we cannot hope to obtain such a decay rate from energy estimates alone, even if we have, schematically,
% \begin{equation*}
%  (\Lbar h)_{LL} \sim \Lbar \Phi_{[0]}
% \end{equation*}
% In other words, even if the semilinear terms corresponding to $h_{LL}$ have the classical null structure, we will still generally ``lose'' some decay in the energy estimates. In this case, the loss arises due to the quasilinear structure, rather than the semilinear structure. Nevertheless, if a suitable ``weak null'' condition holds then we can still recover the required decay rate by integrating along null geodesics; see section [[]].

\chapter{Equations governing the geometry of the null cones}
\label{chapter geometry of the null cones}

The nonlinear wave equations we are studying are coupled to a set of equations determining the geometry of the null cones (and the spheres $S_{\tau,r}$) via the metric $g$. In order to close our estimates we must establish equations governing the evolution of certain quantities describing this geometry. Specifically, we require estimates on the null frame connection coefficients described in chapter \ref{chapter null frame connection coefficients}.

Some of the null frame connection coefficients can be written in terms of the rectangular components of the metric perturbation using our normalisation condition, as we have seen in chapters \ref{chapter transport eikonal} and \ref{chapter null frame connection coefficients}. Specifically, we can estimate the scalar field $\omega$ and the vector field $\zeta$ directly in this way; see propositions \ref{proposition transport mu} and \ref{proposition zeta} respectively. Furthermore, the foliation density $\mu$ satisfies the transport equation given in proposition \ref{proposition transport mu}, and so derivatives of $\mu$ can be estimated by commuting this transport equation and then integrating along the integral curves of the vector field $L$.

The remaining null frame connection coefficients are $\tr_{\slashed{g}} \chi$, $\hat{\chi}$, $\tr_{\slashed{g}} \chibar$ and $\hat{\chibar}$. Of these, the latter two can be written in terms of the former two, as shown in proposition \ref{proposition chibar in terms of chi}. Thus, we are left with the problem of establishing equations governing the evolution of $\tr_{\slashed{g}}\chi$ and $\hat{\chi}$.

The approach we take is to first establish an evolution equation for $\tr_{\slashed{g}} \chi$ along the integral curves of $L$. Crucially, this evolution equation \emph{does not lose derivatives}. That is, $\tr_{\slashed{g}}\chi$ can be estimated in terms of the first derivatives of the metric. Furthermore, we can use this equation to establish pointwise estimates on $\tr_{\slashed{g}}\chi_{(\text{small})}$, which, importantly, show that this quantity decays faster than $r^{-1}$.

We will also establish a transport equation for $\hat{\chi}$, which, however, does lose a derivative, in the sense that $L\hat{\chi}$ will behave like second derivatives of the metric. Then, we will make use of the Codazzi equations to show that $\hat{\chi}$ satisfies an elliptic system, which will allow us to estimate derivatives of $\hat{\chi}$ tangent to the spheres, without losing a derivative, and which can also be used to obtain pointwise bounds on $\hat{\chi}$. Combining these two equations, we can estimate the derivatives of $\hat{\chi}$ tangent to the null cones without losing a derivative. Note that, in order to estimate derivatives of $\hat{\chi}$ transverse to the outgoing null cones, we still have to have \emph{a priori} control over certain third derivatives of the metric.

Aditionally, note that the procedure outlined above to estimate derivatives of the foliation density $\mu$ also involves a loss of derivatives, except in the $L$ direction. Specifically, we should expect $\mu$ to behave like a metric component, but in order to estimate it we must integrate equation \eqref{equation transport mu}, which involves first derivatives of the metric. To solve this problem, we will show that $\mu$ also solves an elliptic equation, which couples the behaviour of angular derivatives of $\mu$ to $\tr_{\slashed{g}}\chi$. Hence, we are also able to estimate derivatives of $\mu$ tangent to the outgoing null cones without a loss of derivatives.

Finally, the elliptic systems we establish are only useful\footnote{In fact, this problem can be avoided by making use of the uniformization theorem together with the fact that the elliptic systems we write down are conformally covariant in two dimensions, see \cite{Christodoulou1993}.} given pointwise bounds on the Gauss curvature of the spheres $S_{\tau, r}$. Hence we will also establish an evolution equation for the Gauss curvature of these spheres.

\section{The Riemann curvature tensor}
The equations governing the behaviour of the connection coefficients involve the Riemann curvature tensor. We first state our conventions regarding this tensor, and then examine in more depth certain null frame components of this tensor.

\begin{definition}
We define the Riemann curvature tensor as follows: for any vector fields $V$, $W$, $X$, $Y$ we have
\begin{equation}
 R(V,W,X,Y) := g(\D_V \D_W X - \D_W \D_V X - \D_{[V, W]}X, Y)
\end{equation}
\end{definition}

\begin{proposition}[The Riemann curvature tensor in terms of the rectangular components of the metric]
 We can relate the Riemann curvature tensor to the derivatives of the rectangular components of the metric via the equation
 \begin{equation}
  R_{abc}^{\phantom{abc}d} = \D_a \left( \Gamma^d_{bc}\partial_d \right) - \D_b \left( \Gamma^d_{ac} \partial_d \right)
 \end{equation}
 where we recall that the rectangular Christoffel symbols $\Gamma^a_{bc}$ are defined in terms of the derivatives of $h$ in proposition \ref{proposition Christoffel symbols in terms of h}. Recall that the Christoffel symbols are a set of scalar quantities, labelled by rectangular indices, while the operators $\partial_a$ are a set of vector fields, also labelled by rectangular indices.
\end{proposition}

\begin{proof}
 We have
 \begin{equation*}
  \begin{split}
   R_{abc}^{\phantom{abc}d}\partial_d &= (\D_a \D_b - \D_b \D_a)\partial_c \\
   &= \D_a \left( \Gamma^d_{bc}\partial_d \right) - \D_b \left( \Gamma^d_{ac} \partial_d \right)
  \end{split}
 \end{equation*}
\end{proof}

\begin{proposition}[Expressions for the rectangular Christoffel symbols]
\label{proposition Christoffel symbols null frame}
 We can express the null frame components of the vector fields $\Gamma_{ab}^c \partial_c$ as follows:
\begin{equation}
 \begin{split}
  \Gamma^a_{LL}\partial_a &= \frac{1}{4} \left( (\Lbar h)_{LL} - 2(Lh)_{L\Lbar} 			 \right) L
		   - \frac{1}{4}(Lh)_{LL} \Lbar
		   + \frac{1}{2}(\slashed{g}^{-1})^{\mu\nu} \left( 2L^a \slashed{\Pi}_\mu^{\phantom{\mu}b}(Lh_{ab}) - (\slashed{\nabla}_\mu h)_{LL} \right) \slashed{\nabla}_\nu \\ \\
  \Gamma^a_{L\Lbar}\partial_a &= -\frac{1}{4}(Lh)_{\Lbar\Lbar} L
		   - \frac{1}{4}(\Lbar h)_{LL} \Lbar
		   + \frac{1}{2}(\slashed{g}^{-1})^{\mu\nu}\left( L^a \slashed{\Pi}_\mu^{\phantom{\mu}b}(\Lbar h_{ab}) + \Lbar^a \slashed{\Pi}_\mu^{\phantom{\mu}b}(Lh_{ab}) - (\slashed{\nabla}_\mu h)_{L\Lbar} \right) \slashed{\nabla}_\nu \\ \\ 
  \slashed{\Pi}_\mu^{\phantom{\mu}b}\Gamma^a_{L b}\partial_a &=
		   -\frac{1}{4}\left( \Lbar^a \slashed{\Pi}_\mu^{\phantom{\mu}b}(Lh_{ab}) + (\slashed{\nabla}_\mu h)_{\Lbar L} - L^a \slashed{\Pi}_\mu^{\phantom{\mu}b}(\Lbar h_{ab}) \right) L
		   - \frac{1}{4}(\slashed{\nabla}_\mu h)_{LL} \Lbar \\ 
		   &\phantom{=} + \frac{1}{2}(\slashed{g}^{-1})^{\nu\rho}\left( \slashed{\Pi}_\mu^{\phantom{\mu}a}\slashed{\Pi}_\nu^{\phantom{\nu}b}(Lh_{ab}) + L^a \slashed{\Pi}_\nu^{\phantom{\nu}b}(\slashed{\nabla}_\mu h_{ab}) - L^a \slashed{\Pi}_\mu^{\phantom{\mu}b}(\slashed{\nabla}_\nu h_{ab}) \right)\slashed{\nabla}_\rho \\ \\
  \Gamma^a_{\Lbar\Lbar}\partial_a &= -\frac{1}{4}(\Lbar h)_{\Lbar\Lbar} L
		   + \frac{1}{4}\left( (Lh)_{\Lbar\Lbar} - 2(\Lbar h)_{L\Lbar} \right) \Lbar
		   + \frac{1}{2}(\slashed{g}^{-1})^{\mu\nu} \left( \Lbar^a\slashed{\Pi}_\mu^{\phantom{\mu}b}(\Lbar h_{ab}) - (\slashed{\nabla}_\mu h)_{\Lbar\Lbar} \right) \slashed{\nabla}_\nu \\ \\
  \slashed{\Pi}_\mu^{\phantom{\mu}b}\Gamma^a_{\Lbar b}\partial_a &=
		   -\frac{1}{4}(\slashed{\nabla}_\mu h)_{\Lbar\Lbar} L
		   - \frac{1}{4}\left( L^a \slashed{\Pi}_\mu^{\phantom{\mu}b}(\Lbar h_{ab}) + (\slashed{\nabla}_\mu h)_{L\Lbar} - \Lbar^a \slashed{\Pi}_\mu^{\phantom{\mu}b}(Lh_{ab}) \right) \Lbar \\
		   &\phantom{=} + \frac{1}{2}(\slashed{g}^{-1})^{\nu\rho}\left( \slashed{\Pi}_\mu^{\phantom{\mu}a}\slashed{\Pi}_\nu^{\phantom{\nu}b}(\Lbar h_{ab}) + \Lbar^a \slashed{\Pi}_\nu^{\phantom{\nu}b}(\slashed{\nabla}_\mu h_{ab}) - \Lbar^a \slashed{\Pi}_\mu^{\phantom{\mu}b}(\slashed{\nabla}_\nu h_{ab}) \right) \slashed{\nabla}_\rho \\ \\
  \slashed{\Pi}_\mu^{\phantom{\mu}b}\slashed{\Pi}_\nu^{\phantom{\nu}c}\Gamma^a_{bc}\partial_a &= 
		   -\frac{1}{4}\left( \Lbar^a \slashed{\Pi}_\nu^{\phantom{\nu}b}(\slashed{\nabla}_\mu h_{ab}) + \Lbar^a \slashed{\Pi}_\mu^{\phantom{\mu}b}(\slashed{\nabla}_\nu h_{ab}) - 	\slashed{\Pi}_\mu^{\phantom{\mu}a}\slashed{\Pi}_\nu^{\phantom{\nu}b}(\Lbar h_{ab}) \right)L \\
		   &\phantom{=}- \frac{1}{4}\left( L^a \slashed{\Pi}_\nu^{\phantom{\nu}b}(\slashed{\nabla}_\mu h_{ab}) + L^a \slashed{\Pi}_\mu^{\phantom{\mu}b}(\slashed{\nabla}_\nu h_{ab}) - \slashed{\Pi}_\mu^{\phantom{\mu}a}\slashed{\Pi}_\nu^{\phantom{\nu}b}(Lh_{ab}) \right)\Lbar \\
		   &\phantom{=} + \frac{1}{2}(\slashed{g}^{-1})^{\rho\sigma} \left( \slashed{\Pi}_\nu^{\phantom{\nu}a}\slashed{\Pi}_\rho^{\phantom{\rho}b}(\slashed{\nabla}_\mu h_{ab}) + \slashed{\Pi}_\mu^{\phantom{\mu}a}\slashed{\Pi}_\rho^{\phantom{\rho}b}(\slashed{\nabla}_\nu h_{ab}) - \slashed{\Pi}_\mu^{\phantom{\mu}a}\slashed{\Pi}_\nu^{\phantom{\nu}b}(\slashed{\nabla}_\rho h_{ab}) \right) \slashed{\nabla}_\sigma
 \end{split}
\end{equation}
 Note also that $\Gamma^a_{bc} = \Gamma^a_{cb}$.

\end{proposition}

One particular component of the Riemann tensor has a special structure, meaning that it can be expressed as \emph{perfect} $L$ derivatives, plus some lower order terms. This structure will play a vital role in allowing us to avoid losing derivatives in the estimates for $\tr_{\slashed{g}}\chi$.

\begin{definition}[The tensor field $\alpha$]
 We define the symmetric, $S_{\tau,r}$-tangent tensor field $\alpha$ by its (linear) action on a pair of $S_{\tau, r}$-tangent vector fields $X_A$ and $X_B$:
 \begin{equation}
  \alpha(X_A, X_B) := R_{LALB}
 \end{equation}
 Note that $\alpha$ is symmetric. We can also view $\alpha$ as a spacetime tensor, where we define its action on the frame fields $L$ and $\Lbar$ by $\alpha(L, \cdot) = \alpha(\Lbar, \cdot) = 0$.
\end{definition}

\begin{proposition}[The structure of the tensor $\alpha$]
\label{proposition structure of alpha}
 The tensor $\alpha$ can be expressed as
 \begin{equation}
  \alpha_{AB} = \frac{1}{2}L\left( (Lh)_{AB} - (X_A h)_{BL} - (X_B h)_{AL} \right) + \frac{1}{2}\left( \slashed{\nabla}_A \slashed{\nabla}_B h \right)_{LL} - \frac{1}{4}(\Lbar h)_{LL}\chi_{AB} + \left(\text{Err}_{(\alpha)}\right)_{AB}
 \end{equation}
%  where the error term satisfies
%  \begin{equation}
%   |\text{Err}_{(\alpha)}| \lesssim |\bar{\partial} h|_{(\text{frame})}|\partial h|_{(\text{frame})} + \left( |\chi| + |\zeta| + |\chibar| + \left| \frac{1 - L^i L^i}{r} \right| + \left| L^i \frac{\slashed{\upd} x^i}{r} \right| \right)|\bar{\partial} h|_{(\text{frame})} 
%  \end{equation}
  where the error term is given schematically by
 \begin{equation}
  \text{Err}_{(\alpha)} = (\bar{\partial}h)_{(\text{frame})}\cdot (\partial h)_{(\text{frame})} + \left( \chi + \chibar + \zeta + \frac{L^i_{(\text{small})}L^i_{(\text{small})}}{r} + L^i_{(\text{small})}\frac{\slashed{\nabla}x^i}{r} \right)\cdot (\bar{\partial}h)_{(\text{frame})}
 \end{equation}

\end{proposition}

\begin{proof}
 We begin by noting that
 \begin{equation*}
 \begin{split}
   R_{LALB} &= g\left( L^a (X_A)^b L^c \left( \D_a (\Gamma^d_{bc}\partial_d) - \D_b (\Gamma^d_{ac}\partial_d)\right) \, , \, X_B \right) \\
   &= g\bigg( \D_L \left( \Gamma^a_{AL}\partial_a \right) - \D_A \left( \Gamma^a_{LL}\partial_a \right) - (LL^a)(X_A)^b \Gamma^c_{ab}\partial_c + (X_A L^a)L^b \Gamma^c_{ab}\partial_c \, , \, X_B \bigg)
  \end{split}
 \end{equation*}
where we have used the fat that $[X_A , L] = 0$, and so $X_A L^a = L(X_A)^a$. We now substitute in the expressions from propositions \ref{proposition Christoffel symbols null frame} and propositions \ref{proposition transport rectangular} and \ref{proposition angular rectangular} for the Christoffel symbols and the derivatives of the rectangular components of the null frame respectively. Finally, we also make use of proposition \ref{proposition null connection} for the null connection coefficients.
 
\end{proof}

\begin{proposition}[The structure of $\tr_{\slashed{g}} \alpha$]
\label{proposition structure trace alpha}
 The scalar field $\tr_{\slashed{g}}\alpha$ can be written as a \emph{perfect} $L$ derivative, plus some error terms which will be shown to be lower order. Specifically, we have
 \begin{equation}
 \label{equation trace alpha}
  \begin{split}
   \tr_{\slashed{g}}\alpha &= L \left( \frac{1}{2}(\slashed{g}^{-1})^{\mu\nu}\slashed{\Pi}_\mu^{\phantom{\mu}a}\slashed{\Pi}_\nu^{\phantom{\nu}b} (Lh_{ab}) - (\slashed{g}^{-1})^{\mu\nu}L^a \slashed{\Pi}_\nu^{\phantom{\nu}b}(\slashed{\nabla}_\mu h_{ab}) \right) + \frac{1}{2}(\slashed{\Delta}h)_{LL} - \omega \tr_{\slashed{g}}\chi + \text{Err}_{(\tr_{\slashed{g}}\alpha, \text{low})} \\
   &= L \left( \frac{1}{2}(\slashed{g}^{-1})^{\mu\nu} \slashed{\Pi}_\mu^{\phantom{\mu}a}\slashed{\Pi}_\nu^{\phantom{\nu}b}(Lh_{ab}) - (\slashed{g}^{-1})^{\mu\nu}L^a \slashed{\Pi}_\nu^{\phantom{\nu}b}(\slashed{\nabla}_\mu h_{ab}) + \frac{1}{2}(\Lbar h)_{LL} \right) + \frac{1}{2}(\tilde{\Box}_g h)_{LL} + \text{Err}_{(\tr_{\slashed{g}}\alpha, \text{high})} \\
  \end{split}
 \end{equation}
 where the error terms are given schematically by
 \begin{equation}
  \begin{rcases} \text{Err}_{(\tr_{\slashed{g}}\alpha , \text{low})} \\
   \text{Err}_{(\tr_{\slashed{g}}\alpha, \text{high})}
  \end{rcases}
  = (\bar{\partial}h)_{(\text{frame})}\cdot (\partial h)_{(\text{frame})} + \left( \chi + \chibar + \zeta + \frac{L^i_{(\text{small})}L^i_{(\text{small})}}{r} + L^i_{(\text{small})}\frac{\slashed{\nabla}x^i}{r} \right)\cdot (\bar{\partial}h)_{(\text{frame})}
 \end{equation}
\end{proposition}

\begin{proof}
 The first line of \eqref{equation trace alpha} follows from proposition \ref{proposition structure of alpha}. For the second line, we make use of proposition \ref{proposition scalar wave operator} to express the second angular derivatives of $h$ in terms of $L\Lbar h$ and lower order terms\footnote{Note that the Einstein equations directly imply that $\tr_{\slashed{g}} \alpha = 0$, since $R_{LL} = -\frac{1}{2} R_{LLL\Lbar} - \frac{1}{2} R_{L\Lbar LL} + (\slashed{g}^{-1})^{AB} R_{LALB}$, and the first two terms vanish identically.}.
\end{proof}

\section{The transport equations for \texorpdfstring{$\tr_{\slashed{g}}\chi$ and $\hat{\chi}$}{tr chi and chihat}}

\begin{proposition}[The transport equation for the frame components of $\chi$]
\label{proposition L chi}
 The null frame components of the connection coefficient $\chi$ satisfies the transport equation
\begin{equation}
\label{equation L chi}
 L \chi_{AB} = \alpha_{AB} + \omega \chi_{AB} + \chi_A^{\phantom{A}C}\chi_{BC}
\end{equation}
\end{proposition}

\begin{proof}
 We have
 \begin{equation*}
  \begin{split}
   L\chi_{AB} &= Lg(\D_A L, X_B) \\
   &= g(\D_L \D_A L, X_B) + g(\D_A L, \D_L X_B) \\
   &= R_{LALB} + g(\D_A \D_L L, X_B) + \chi_A^{\phantom{A}C}\chi_{BC} \\
   &= R_{LALB} - g(\D_L L, \D_A X_B) + \chi_A^{\phantom{A}C}\chi_{BC} \\
   &= R_{LALB} + \omega \chi_{AB} + \chi_A^{\phantom{A}C}\chi_{BC}
  \end{split}
 \end{equation*}
\end{proof}

\begin{proposition}[The transport equation for the tensor fields $\chi$ and $\chi_{(\text{small})}$]
\label{proposition transport chi chismall}
 The connection coefficient $\chi_{\mu\nu}$ satisfies the transport equation
\begin{equation}
 \label{equation L chi tensor}
 \slashed{\D}_L \chi_{\mu\nu} = \slashed{\Pi}_\mu^{\phantom{\mu}a}\slashed{\Pi}_\nu^{\phantom{\nu}b} R_{LaLb}+ \omega \chi_{\mu\nu} - \chi_\mu^{\phantom{\mu}\rho}\chi_{\rho\nu}
\end{equation}
and the connection coefficient $(\chi_{(\text{small})})_{\mu\nu}$ satisfies the transport equation
\begin{equation}
 \label{equation L chismall tensor}
 \slashed{\D}_L \left(\chi_{(\text{small})} \right)_{\mu\nu} = \slashed{\Pi}_\mu^{\phantom{\mu}a}\slashed{\Pi}_\nu^{\phantom{\nu}b} R_{LaLb}+ \omega \left(\chi_{(\text{small})}\right)_{\mu\nu} + \omega r^{-1}\slashed{g}_{\mu\nu} - \left(\chi_{(\text{small})}\right)_\mu^{\phantom{\mu}\rho}\left(\chi_{(\text{small})}\right)_{\rho\nu} - 2r^{-1}\left( \chi_{(\text{small})}\right)_{\mu\nu}
\end{equation}

\end{proposition}

\begin{proof}
 We have already established the transport equation for the null frame components of $\chi$, which is given by equation \eqref{equation L chi}. Now, we have
\begin{equation}
 \begin{split}
  (X_A)^\mu (X_B)^\nu \slashed{\D}_L \chi_{\mu \nu} &= (X_A)^\mu (X_B)^\nu \D_L \chi_{\mu \nu} \\
  &= L\chi_{AB} - (\D_L X_A)^\mu (X_B)^\nu \chi_{\mu\nu} - (\D_L X_B)^\mu (X_A)^\nu \chi_{\mu\nu} \\
  &= L\chi_{AB} - 2\chi_A^{\phantom{A}C}\chi_{CB}
 \end{split}
\end{equation}
from which equation \eqref{equation L chi tensor} follows. To prove \eqref{equation L chismall tensor} we simply substitute the definition $\chi_{(\text{small})} = \chi - r^{-1} \slashed{g}$.

\end{proof}

\begin{proposition}[The transport equation for $\tr_{\slashed{g}}\chi$]
\label{proposition transport trace chi}
 $\tr_{\slashed{g}}$ satisfies the following transport equation along the integral curves of the vector field $L$:
 \begin{equation}
 \label{equation L tr chi}
  \begin{split}
   L \tr_{\slashed{g}}\chi &= \tr_{\slashed{g}}\alpha + \omega \tr_{\slashed{g}}\chi - \frac{1}{2}\left( \tr_{\slashed{g}}\chi \right)^2 - \hat{\chi}\cdot \hat{\chi}
  \end{split}
 \end{equation}
\end{proposition}
\begin{proof}
 Recall that $L \slashed{g}_{AB} = 2\chi_{AB}$. From this it follows that
 \begin{equation*}
  L \left( (\slashed{g}^{-1})^{AB} \right) = -2 \chi^{AB}
 \end{equation*}
 Combining this with equation \eqref{equation L chi} proves the proposition.
\end{proof}

\begin{proposition}[Transport equations for the modified versions of $\tr_{\slashed{g}}\chi$]
 By making use of the structure of $\tr_{\slashed{g}}\alpha$ we can find transport equations which either lead to improved decay or improved regularity for $\tr_{\slashed{g}}\chi$. Specifically, we can define the modified versions of $\tr_{\slashed{g}}\chi$:
 \begin{equation}
  \begin{split}
   \mathcal{X}_{(\text{low})} &:= \tr_{\slashed{g}}\chi_{(\text{small})} - \frac{1}{2}(\slashed{g}^{-1})^{\mu\nu}\slashed{\Pi}_\mu^{\phantom{\mu}a}\slashed{\Pi}_\nu^{\phantom{\nu}b}(Lh_{ab}) + (\slashed{g}^{-1})^{\mu\nu} L^a \slashed{\Pi}_\nu^{\phantom{\nu}b}(\slashed{\nabla}_\mu h_{ab}) \\
   \mathcal{X}_{(\text{high})} &:= \tr_{\slashed{g}}\chi_{(\text{small})} - \frac{1}{2}(\slashed{g}^{-1})^{\mu\nu}\slashed{\Pi}_\mu^{\phantom{\mu}a}\slashed{\Pi}_\nu^{\phantom{\nu}b}(Lh_{ab}) + (\slashed{g}^{-1})^{\mu\nu} L^a \slashed{\Pi}_\nu^{\phantom{\nu}b}(\slashed{\nabla}_\mu h_{ab}) - \frac{1}{2}(\Lbar h)_{LL}
  \end{split}
 \end{equation}
 Then these quantities satisfy the transport equations
 \begin{equation}
 \label{equation transport modified chi}
  \begin{split}
   L\left( r^2 \mathcal{X}_{(\text{low})} \right) &= \frac{1}{2}r^2(\slashed{\Delta}h)_{LL} - \frac{1}{2}r^{-2}\left(r^2\mathcal{X}_{(\text{low})} \right)^2 - r^2\hat{\chi}\cdot\hat{\chi} + r^2\text{Err}_{( \mathcal{X}_{(\text{low})} )} \\
   L\left( r^2\mathcal{X}_{(\text{high})} \right) &= \frac{1}{2}r^2(\tilde{\Box}_g h)_{LL} + r^2 \omega \mathcal{X}_{(\text{high})} + 2r\omega - \frac{1}{2}r^2(\mathcal{X}_{(\text{high})})^2 - r^2\hat{\chi}\cdot \hat{\chi} + r^2\text{Err}_{( \mathcal{X}_{(\text{high})} )} \\
  \end{split}
 \end{equation}
 where the error terms are given schematically by
 \begin{equation}
  \begin{rcases} \text{Err}_{( \mathcal{X}_{(\text{low})} )} \\
  \text{Err}_{( \mathcal{X}_{(\text{high})} )}
  \end{rcases}
  = (\bar{\partial}h)_{(\text{frame})}\cdot (\partial h)_{(\text{frame})} + \left( \chi + \chibar + \zeta + \frac{L^i_{(\text{small})}L^i_{(\text{small})}}{r} + L^i_{(\text{small})}\frac{\slashed{\nabla}x^i}{r} \right)\cdot (\bar{\partial}h)_{(\text{frame})}
 \end{equation}

\end{proposition}

\begin{proof}
 We make use of the structure of $\tr_{\slashed{g}}\alpha$, as established in proposition \ref{proposition structure trace alpha}, together with the transport equation for $tr_{\slashed{g}}\chi$ established in proposition \ref{proposition transport trace chi}, and we also substitute for $\omega$ from the definition \ref{definition omega} and proposition $\ref{proposition transport mu}$. Note that the potentially dangerous term $(\partial h)_{LL} \tr_{\slashed{g}} \chi_{(\text{small})}$ cancels.
\end{proof}

\begin{remark}
 Note the additional $r$ weights in the transport equation for $\mathcal{X}_{(\text{low})}$ (the first line of equation \eqref{equation transport modified chi}). These will allow us to show improved decay in $r$ for the quantity $\mathcal{X}_{(\text{low})}$. On the other hand, note that the term $(\slashed{\Delta} h)_{LL}$ appears on the right hand side of this equation, meaning that, in order to estimate $\mathcal{X}_{(\text{low})}$, we already require estimates on the second derivatives of $h$. This is the ``loss of derivatives'' referred to previously. On the other hand, to estimate $\mathcal{X}_{(\text{high})}$ we only need information about the first derivatives of $h$.
\end{remark}

\begin{proposition}[The transport equation for the frame components of $\hat{\chi}$]
 The null frame components of $\hat{\chi}$ satisfy the transport equation
\begin{equation}
 L\hat{\chi}_{AB} = \hat{\alpha}_{AB} + \omega \hat{\chi}_{AB} + \frac{1}{2}\slashed{g}_{AB}\hat{\chi}\cdot \hat{\chi} + \hat{\chi}_A^{\phantom{A}C} \hat{\chi}_{CB}
\end{equation}
\end{proposition}

\begin{proof}
 Recall that we define
\begin{equation*}
 \hat{\chi}_{AB} := \chi_{AB} - \frac{1}{2}\slashed{g}_{AB} \tr_{\slashed{g}}\chi
\end{equation*}
 Also recall that $L \slashed{g}_{AB} = 2\chi_{AB}$. Hence, substituting for $L\chi$ from equation \eqref{equation L chi} and for $L\tr_{\slashed{g}}\chi$ from equation \eqref{equation L tr chi} proves the proposition.
\end{proof}

\begin{proposition}[The transport equation for $\hat{\chi}$]
	\label{proposition transport chihat}
	The $S_{\tau,r}$-tangent tensor field $\hat{\chi}$ satisfies the following transport equation along the integral curves of $L$:
	\begin{equation}
	\slashed{\D}_L \left(r^2 \hat{\chi}_{\mu\nu} \right)
	=
	r^2 \hat{\alpha}_{\mu\nu}
	+ r^2 \omega \hat{\chi}_{\mu\nu}
	+ \frac{1}{2}r^2 \slashed{g}_{\mu\nu}|\hat{\chi}|^2
	- r^2 \hat{\chi}_\mu^{\phantom{\mu}\rho} \hat{\chi}_{\rho \nu}
	- 2 \hat{\chi}_{\mu\nu} \tr_{\slashed{g}}\chi_{(\text{small})}
	\end{equation}
\end{proposition}

\begin{proof}
	We begin by noting that
	\begin{equation*}
	\begin{split}
	L\hat{\chi}_{AB} 
	&= X_A^\mu X_B^\nu \slashed{\D}_L \hat{\chi}_{\mu\nu} + \hat{\chi}_{\mu A} \slashed{\D}_L X_B^\mu + \hat{\chi}_{\mu B} \slashed{\D}_L X_A^\mu \\
	&= X_A^\mu X_B^\nu \left( \slashed{\D}_L \hat{\chi}_{\mu\nu} + \hat{\chi}_{\mu \rho}\chi_\nu^{\phantom{\nu}\rho} + \hat{\chi}_{\nu \rho}\chi_\mu^{\phantom{\mu}\rho} \right) \\
	&= X_A^\mu X_B^\nu \left( \slashed{\D}_L \hat{\chi}_{\mu\nu} + \frac{2}{r}\hat{\chi}_{\mu\nu} + 2\hat{\chi}_{\mu\nu} \tr_{\slashed{g}}\chi_{(\text{small})} + 2\hat{\chi}_{\mu}^{\phantom{\mu}\rho}\hat{\chi}_{\nu\rho} \right) \\
	\end{split}
	\end{equation*}
	
	Inserting the expression from the previous proposition on the left hand side, and then rearranging this expression gives
	\begin{equation*}
	\slashed{\D}_L \hat{\chi}_{\mu\nu} 
	=
	- \frac{2}{r}\hat{\chi}_{\mu\nu} 
	- 2\hat{\chi}_{\mu\nu} \tr_{\slashed{g}}\chi_{(\text{small})} 
	- \hat{\chi}_{\mu}^{\phantom{\mu}\rho}\hat{\chi}_{\nu\rho}
	+ \hat{\alpha}_{\mu\nu}
	+ \omega \hat{\chi}_{\mu\nu}
	+ \frac{1}{2}\slashed{g}_{\mu\nu}|\hat{\chi}|^2
	\end{equation*}
	Now multiplying by $r^2$ proves the proposition.
\end{proof}

%\begin{proposition}[The transport equation for $|\hat{\chi}|^2$]
%	\label{proposition transport chihat squared}
%	We also establish a transport equation for the scalar field $|\hat{\chi}|^2 = \hat{\chi}^{\mu\nu}\hat{\chi}_{\mu\nu}$, which will allow us to accurately estimate the size of $\hat{\chi}$. We have
%	\begin{equation}
%	L\left( r^4 |\hat{\chi}|^2 \right)
%	= 2r^4 \hat{\chi} \cdot \hat{\alpha} 
%	+ 2 r^4 \left(\omega - \tr_{\slashed{g}}\chi_{(\text{small})} \right) |\hat{\chi}|^2
%	-2 r^4 \hat{\chi}^{\mu\nu} \hat{\chi}_\mu^{\phantom{\mu}\rho} \hat{\chi}_{\rho \nu}
%	\end{equation}	
%\end{proposition}
%
%\begin{proof}
%	Since $|\hat{\chi}|^2$ is a scalar field we can perform this computation in the frame $X_A$. Again making use of $L (\slashed{g}^{-1})^{AB} = -2 \chi^{AB}$, we have
%	\begin{equation*}
%	\begin{split}
%	L\left(|\hat{\chi}|^2\right)
%	&=
%	2\hat{\chi}^{AB} L\hat{\chi}_{AB} 
%	-4 \chi^{AB} \hat{\chi}_{AC}\hat{\chi}^{C}_{\phantom{C}B}
%	\\
%	&= 2\hat{\chi} \cdot \hat{\alpha}
%	+ 2 \omega |\hat{\chi}|^2
%	+ 2 \hat{\chi}^{AB} \hat{\chi}_A^{\phantom{A}C} \hat{\chi}_{CB}
%	- 4 \hat{\chi}^{AB} \hat{\chi}_A^{\phantom{A}C} \hat{\chi}_{CB}
%	-4 r^{-1} |\hat{\chi}|^2
%	-2 \tr_{\slashed{g}} \chi_{(\text{small})} |\hat{\chi}|^2
%	\end{split}
%	\end{equation*}
%	and so
%	\begin{equation*}
%	L\left(r^4 |\hat{\chi}|^2 \right)
%	=
%	2r^4 \hat{\chi} \cdot \hat{\alpha}
%	+ 2r^4 \left(\omega - 2 \tr_{\slashed{g}}\chi_{(\text{small})} \right)|\hat{\chi}|^2
%	-2r^4 \hat{\chi}^{\mu\nu} \hat{\chi}_\mu^{\phantom{\mu}\rho} \hat{\chi}_{\rho \nu}
%	\end{equation*}
%\end{proof}

\section{Additional transport equations}

We will also occasionally require equations for both $\Lbar \chi$ and $L \chibar$. These are best given as tensorial equations, rather than equations for the frame components.

\begin{proposition}[An equation for $\slashed{\D}_{\Lbar} \chi_{\mu\nu}$]
	\label{proposition Lbar chi}
	$\chi_{\mu\nu}$ satisfies the transport equation
	\begin{equation}
	\begin{split}
	\slashed{\D}_{\Lbar} \chi_{\mu\nu}
	&=
	\slashed{\Pi}_\mu^{\phantom{\mu}\rho} \slashed{\Pi}_\nu^{\phantom{\nu}\sigma} R_{\Lbar \rho L \sigma}
	+ \frac{1}{2} \slashed{\nabla}_\mu \zeta_\nu
	+ \frac{1}{2} \slashed{\nabla}_\nu \zeta_\mu
	+ 2\slashed{\nabla}^2_{\mu\nu} \log \mu
	+ \frac{1}{2} \zeta_\mu \zeta_\nu
	+ \zeta_\mu \slashed{\nabla}_\nu \log \mu 
	+ \zeta_\nu \slashed{\nabla}_\mu \log \mu
	\\
	&\phantom{=}
	+ 2 (\slashed{\nabla}_\mu \log \mu)(\slashed{\nabla}_\nu \log \mu)
	+ \omega \chi_{\mu\nu}
	- \frac{1}{2} \chi_\mu^{\phantom{\mu}\sigma} \chibar_{\nu \sigma}
	- \frac{1}{2} \chibar_\mu^{\phantom{\mu}\sigma} \chi_{\nu \sigma}
	\end{split}
	\end{equation}
	
\end{proposition}

\begin{proof}
	We find
	\begin{equation*}
	\begin{split}
	\Lbar \chi_{AB} &= R_{\Lbar ALB} + \slashed{\nabla}_A \left(\zeta_B + 2 \slashed{\nabla}_B(\log \mu) \right) + \frac{1}{2}\zeta_A \zeta_B + \zeta_A \slashed{\nabla}_B(\log \mu) + \zeta_B \slashed{\nabla}_A(\log \mu) \\
	&\phantom{=} + 2(\slashed{\nabla}_A \log\mu)(\slashed{\nabla}_B \log\mu) + \omega\chi_{AB} + \chi_A^{\phantom{A}C}\chibar_{BC} + \chi_A^{\phantom{A}C}g([\Lbar, X_B] , X_C) + \chi_B^{\phantom{B}C}g([\Lbar, X_A] , X_C)
	\end{split}
	\end{equation*}
	Symmetrising over $A$ and $B$, and bringing out the factor $(X_A)^\mu (X_B)^\nu$ proves the proposition.
\end{proof}

\begin{proposition}[An equation for $\slashed{\D}_L \chibar_{\mu\nu}$]
	\label{proposition transport chibar}
	The tensor field $\chibar$ satisfies the transport equation
	\begin{equation}
	\slashed{\D}_L \chibar_{\mu\nu}
	=
	\slashed{\Pi}_\mu^{\phantom{\mu}\rho} \slashed{\Pi}_\nu^{\phantom{\nu}\sigma} R_{L \rho \Lbar \sigma}
	- \frac{1}{2} \slashed{\nabla}_\mu \zeta_\nu
	- \frac{1}{2} \slashed{\nabla}_\nu \zeta_\mu
	+ \frac{1}{2} \zeta_\mu \zeta_\nu
	- \omega \chibar_{\mu\nu}
	- \frac{1}{2} \chi_\mu^{\phantom{\mu}\sigma} \chibar_{\sigma\nu}
	- \frac{1}{2} \chibar_\mu^{\phantom{\mu}\sigma} \chi_{\sigma\nu}
	\end{equation}
	
\end{proposition}

\begin{proof}
	Similarly to above, we can calculate
	\begin{equation*}
	L\chibar_{AB}
	=
	R_{LA\Lbar B}
	- \slashed{\nabla}_A \zeta_B
	+ \frac{1}{2} \zeta_A \zeta_B
	- \omega \chibar_{AB}
	+ \chibar_A^{\phantom{A}C}\chi_{BC}
	\end{equation*}
	Now, symmetrising over $A$ and $B$ and bringing out a factor of $(X_A)^\mu (X_B)^\nu$ proves the proposition.
\end{proof}

\section{The elliptic system for \texorpdfstring{$\hat{\chi}$}{chihat}}
In this section we will express the divergence of $\hat{\chi}$ in terms of angular derivatives of $\tr_{\slashed{g}}\chi$, second derivatives of $h$ and some lower order terms. It turns out that this system can be used to establish pointwise estimates of $\hat{\chi}$ as well as higher order estimates on angular derivatives of $\hat{\chi}$. Indeed, combined with proposition \ref{proposition transport chihat} and the elliptic estimates of chapter \ref{chapter elliptic estimates and sobolev embedding} this will allow us to estimate all the \emph{good} derivatives of $\hat{\chi}$ \emph{without losing derivatives}. This plays a crucial role in our eventual bootstrap argument - see chapter \ref{chapter pointwise bounds} and \ref{chapter proving the theorem}.

\begin{proposition}[The divergence of $\hat{\chi}$]
 \label{proposition div chihat}
 $\hat{\chi}$ satisfies the equation
\begin{equation}
 \left(\slashed{\Div}\, \hat{\chi} \right)_\mu = \frac{1}{2}\slashed{\nabla}_\mu \tr_{\slashed{g}}\chi - (\slashed{g}^{-1})^{\nu\rho} R_{\mu \nu L \rho} + \frac{1}{4}\zeta_\mu \tr_{\slashed{g}}\chi - \frac{1}{2}\zeta^\nu \hat{\chi}_{\mu\nu} 
\end{equation}
\end{proposition}

\begin{proof}
 We calculate
\begin{equation*}
 \begin{split}
  X_A \chi_{BC} &= X_A g(\D_B L, X_C) \\
  &= g(\D_A \D_B L, X_C) + g(\D_B L, \D_A X_C) \\
  &= R_{ABLC} + g(\D_B \D_A L, X_C) + g(\D_B L, \D_A X_C) \\
  &= R_{ABLC} + X_B g(\D_A L, X_C) + g(\D_B L, \D_A X_C) - g(\D_A L, \D_B X_C) \\
  &= R_{ABLC} + X_B \chi_{AC} + \frac{1}{2}\zeta_B \chi_{AC} - \frac{1}{2}\zeta_A \chi_{BC} + \slashed{\Gamma}_{AC}^D \chi_{BD} - \slashed{\Gamma}_{BC}^D \chi_{AD}
 \end{split}
\end{equation*}
from which the ``Codazzi equations'' follow:
\begin{equation*}
 \slashed{\nabla}_A \chi_{BC} - \slashed{\nabla}_B \chi_{AC} = R_{ABLC} + \frac{1}{2}\zeta_B \chi_{AC} - \frac{1}{2}\zeta_A \chi_{BC}
\end{equation*}
Contracting with $(\slashed{g}^{-1})^{AC}$ and recalling that $\slashed{\nabla}$ commutes with $\slashed{g}$ and its inverse, and finally decomposing $\chi$ into its trace and trace-free parts, we prove the proposition.

\end{proof}

\section{The spherical laplacian of the foliation density}

We must be careful to avoid a loss of derivatives in terms involving angular derivatives of the foliation density $\mu$. To this end, we will derive an equation\footnote{An alternative approach would be to derive a propagation equation for $\slashed{\Delta} \log \mu$ in the $L$ direction, which gains a derivative relative to the na\"ive estimate. To do this, we would commute the equation given in proposition \ref{proposition transport mu} with $\slashed{\Delta}$, and then use the fact that $\slashed{\Delta} h_{(\text{rect})}$ can be written in terms of a perfect $L$ derivative, $\tilde{\Box}_g h_{(\text{rect})}$ and lower order terms.} for the spherical laplacian of $\mu$, linking its behaviour to that of $T \tr_{\slashed{g}}\chi_{(\text{small})}$. Combined with elliptic estimates and estimates on the Gauss curvature of the spheres of constant $\tau$ and $r$, this will allow us to estimate angular derivatives of $\mu$ at the level of first derivatives of $h$.

\begin{proposition}[The spherical laplacian of $\mu$]
 \label{proposition spherical laplacian mu}
 $\mu$ satisfies the equation
 \begin{equation}
 \begin{split}
  \slashed{\Delta} \left( \log \mu \right) &= 2 T(\tr_{\slashed{g}}\chi) - (\slashed{g}^{-1})^{\mu\nu}R_{\Lbar \mu L \nu} - \tr_{\slashed{g}}\alpha - \slashed{\Div}\zeta - 2|\slashed{\nabla}\log \mu|^2 - 2 \zeta^\mu \slashed{\nabla}_\mu (\log \mu) + \frac{1}{2}(\tr_{\slashed{g}}\chi)^2 \\
  &\phantom{=} + \frac{1}{2}(\tr_{\slashed{g}}\chi)(\tr_{\slashed{g}}\chibar) - \frac{1}{2}|\zeta|^2 + \hat{\chi}\cdot\hat{\chi} + \hat{\chi}\cdot\hat{\chibar}
  \end{split}
 \end{equation}
\end{proposition}
\begin{proof}
 We begin by computing $\Lbar \chi_{AB}$. We find
 \begin{equation*}
  \begin{split}
   \Lbar \chi_{AB} &= R_{\Lbar ALB} + \slashed{\nabla}_A \left(\zeta_B + 2 \slashed{\nabla}_B(\log \mu) \right) + \frac{1}{2}\zeta_A \zeta_B + \zeta_A \slashed{\nabla}_B(\log \mu) + \zeta_B \slashed{\nabla}_A(\log \mu) \\
   &\phantom{=} + 2(\slashed{\nabla}_A \log\mu)(\slashed{\nabla}_B \log\mu) + \omega\chi_{AB} + \chi_A^{\phantom{A}C}\chibar_{BC} + \chi_A^{\phantom{A}C}g([\Lbar, X_B] , X_C) + \chi_B^{\phantom{B}C}g([\Lbar, X_A] , X_C)
  \end{split}
 \end{equation*}
 Additionally, we can compute
 \begin{equation*}
  \begin{split}
   \Lbar (\slashed{g}^{-1})^{AB} &= -(\slashed{g}^{-1})^{AC}(\slashed{g}^{-1})^{BD} \Lbar \slashed{g}_{CD} \\
   &= -2\chibar^{AB} - (\slashed{g}^{-1})^{AC}(\slashed{g}^{-1})^{BD} \left( g([\Lbar, X_C], X_D) + g([\Lbar, X_D], X_C) \right)
  \end{split}
 \end{equation*}
 Combining the above two equations we find that
 \begin{equation*}
  \Lbar \tr_{\slashed{g}}\chi = (\slashed{g}^{-1})^{AB}R_{\Lbar ALB} + \slashed{\Div}\zeta + 2\slashed{\Delta}\log \mu + 2|\slashed{\nabla}\log\mu|^2 + 2\zeta^A \slashed{\nabla}_A (\log \mu) + \frac{1}{2}|\zeta|^2 - \chi\cdot\chibar + \omega\tr_{\slashed{g}}\chi
 \end{equation*}
 Recalling that $2T = L + \Lbar$, we can add equation \eqref{equation L tr chi} to the equation above to proves the proposition.
\end{proof}
% 
% \begin{proposition}[An expression for $\Box_g \tr_{\slashed{g}}\chi$]
%  We have
% \begin{equation}
%  \begin{split}
%   \Box_g \tr_{\slashed{g}}\chi &=
%   - \Lbar \tr_{\slashed{g}}\alpha
%   - \tr_{\slashed{g}}\chi (\Lbar \omega)
%   + \slashed{\Delta}\tr_{\slashed{g}}\chi_{(\text{small})} \\
%   &\phantom{=}- \frac{1}{2}\left( \tr_{\slashed{g}}\chi_{(\text{small})} + \tr_{\slashed{g}}\chibar_{(\text{small})} + 2\omega \right) \left( \tr_{\slashed{g}}\alpha + \omega \tr_{\slashed{g}}\chi - \frac{1}{2}(\tr_{\slashed{g}}\chi)^2 - \hat{\chi}\cdot \hat{\chi} \right) \\
%   &\phantom{=}- 2 \hat{\chi} \cdot (\hat{\alpha} + \omega \hat{\chi}
%   - \tr_{\slashed{g}}\chi) (T \tr_{\slashed{g}}\chi_{(\text{small})})
%   + (\zeta^\alpha + 2\slashed{\nabla}^\alpha \log \mu )\slashed{\nabla}_\alpha \tr_{\slashed{g}}\chi_{(\text{small})}
%  \end{split}
% \end{equation}
% schematically, this is
% \begin{equation*}
%  \Box_g \tr_{\slashed{g}}\chi = \Lbar \tr_{\slashed{g}}\alpha
%   - \tr_{\slashed{g}}\chi (\Lbar \omega)
%   + \slashed{\Delta}\tr_{\slashed{g}}\chi_{(\text{small})}
%   + \begin{pmatrix}
%      \bm{\Gamma}_{(\text{bad})}
%     \end{pmatrix}
% \end{equation*}
% 
% 
% \end{proposition}

\section{The Gauss curvature and its evolution equation}
The Gauss curvature of the spheres $S_{\tau,r}$ plays an important role in various elliptic estimates. In this section we give our conventions for defining the Gauss curvature, and derive an equation for its evolution along the integral curves of $L$.

\begin{definition}[The Gauss curvature]
We define the Gauss curvature $K$ as
\begin{equation}
 K := \frac{1}{2}(\slashed{g}^{-1})^{\mu\rho}(\slashed{g}^{-1})^{\nu\sigma} \slashed{R}_{\mu\nu\rho\sigma}
\end{equation}
where $\slashed{R}$ is the Riemann curvature associated with the metric $\slashed{g}$ on the sphere $S_{\tau,r}$.
\end{definition}

\begin{proposition}[Basic properties of the Gauss curvature]
 The Riemann curvature tensor associated with the metric $\slashed{g}$ can be expressed in terms of the Gauss curvature as
\begin{equation}
 \slashed{R}_{\mu\nu\rho\sigma} = K\left( \slashed{g}_{\mu\rho}\slashed{g}_{\nu\sigma} - \slashed{g}_{\mu\sigma}\slashed{g}_{\nu\rho} \right)
\end{equation}
Similarly, the Ricci curvature tensor can be expressed as
\begin{equation}
 \slashed{R}_{\mu\nu} = K\slashed{g}_{\mu\nu}
\end{equation}
Moreover, the Gauss curvature can be expressed in terms of the Riemann curvature tensor of $(\mathcal{M}, g)$ along with the second fundamental forms $\chi$ and $\chibar$:
\begin{equation}
 K = \frac{1}{2}(\slashed{g}^{-1})^{\mu\rho}(\slashed{g}^{-1})^{\nu\sigma} R_{\mu\nu\rho\sigma} - \frac{1}{2}\hat{\chibar}\cdot \hat{\chi} + \frac{1}{4}\tr_{\slashed{g}}\chi \tr_{\slashed{g}}\chibar
\end{equation}
\end{proposition}

\begin{proof}
 The first two equations in the proposition above follow from standard arguments regarding the symmetries of the Riemann tensor on a two dimensional manifold. For the third equation, note that
 \begin{equation*}
  \begin{split}
   R_{ABCD} &= g( \D_A \D_B X_C - \D_B \D_A X_C \, , \, X_D) \\
   &= g\left( \D_A \left( \frac{1}{2}\chibar_{BC} L + \frac{1}{2}\chi_{BC} \Lbar + \slashed{\nabla}_B X_C \right) - \D_B \left( \frac{1}{2}\chibar_{AC} L + \frac{1}{2}\chi_{AC} \Lbar + \slashed{\nabla}_A X_C \right) \, , \, X_D \right) \\
   &= \frac{1}{4}\chibar_{BC}\chi_{AD} + \frac{1}{4}\chi_{BC}\chibar_{AD} - \frac{1}{4}\chibar_{AC}\chi_{BD} - \frac{1}{4}\chi_{AC}\chibar_{BD} + g\left(\slashed{\nabla}_A \slashed{\nabla}_B X_C - \slashed{\nabla}_B \slashed{\nabla}_A X_C \, , \, X_D \right) \\
   &= \frac{1}{4}\chibar_{BC}\chi_{AD} + \frac{1}{4}\chi_{BC}\chibar_{AD} - \frac{1}{4}\chibar_{AC}\chi_{BD} - \frac{1}{4}\chi_{AC}\chibar_{BD} + \slashed{R}_{ABCD}
  \end{split}
 \end{equation*}
 where we have used the fact that $[X_A, X_B] = 0$. The third part of the proposition above now follows from contracting the equation above with $(\slashed{g}^{-1})^{AC}(\slashed{g}^{-1})^{BD}$.
\end{proof}

\begin{proposition}[The transport equation for the Gauss curvature]
\label{proposition transport Gauss curvature}
 In the region $r \geq r_0$, the Gauss curvature of the spheres $S_{\tau,r}$ satisfies the evolution equation
 \begin{equation}
  LK = \frac{1}{2}\slashed{\Delta} \tr_{\slashed{g}}\chi - \slashed{\nabla}^\mu \slashed{\nabla}^\nu \hat{\chi}_{\mu\nu} - K\tr_{\slashed{g}}\chi
 \end{equation}
 Alternatively, this can be written as
 \begin{equation}
  L(r^2 K) = \frac{1}{2}r^2\slashed{\Delta} \tr_{\slashed{g}}\chi_{(\text{small})} - r^2\slashed{\nabla}^\mu \slashed{\nabla}^\nu \hat{\chi}_{\mu\nu} - r^2K\tr_{\slashed{g}}\chi_{(\text{small})}
 \end{equation}
\end{proposition}

\begin{proof}
 We consider a one parameter family of metrics (with parameter $\rho$) on $\mathbb{S}^2$, defined as follows. Fix some constant $u$, and define
 \begin{equation*}
  \slashed{g}(\rho) = \slashed{g}\big|_{S_{\tau , \rho}}
 \end{equation*}
 i.e.\ we consider the induced metrics on the spheres $S_{\tau, r(\rho)}$ where $\tau$ is a constant and $r(\rho) = \rho$. For the rest of this proof we will use a dot above quantities to denote derivatives with respect to the parameter $\rho$, so, for example, 
 \begin{equation*}
  \dot{\slashed{g}} := \frac{\upd}{\upd \rho} \left(\slashed{g}(\rho) \right)
 \end{equation*}
 Note that, from the expressions \ref{corollary null frame coords}
 
 We first note that
 \begin{equation*}
  \frac{\partial}{\partial \rho}\bigg|_{u, \vartheta^A} = L
 \end{equation*}
 
 Hence, we find that
 \begin{equation*}
  \dot{\slashed{g}}_{AB} = \mathcal{L}_L \slashed{g}_{AB} = 2\chi_{AB}
 \end{equation*}
 A standard computation (using, for example, normal coordinates) reveals that
 \begin{equation*}
  \dot{\slashed{\Gamma}}_{BC}^A = \frac{1}{2}(\slashed{g}^{-1})^{AD}\left( \slashed{\nabla}_B \dot{\slashed{g}}_{CD} + \slashed{\nabla}_C \dot{\slashed{g}}_{BD} - \slashed{\nabla}_D \dot{\slashed{g}}_{BC} \right)
 \end{equation*}
 and so we have
  \begin{equation*}
  \dot{\slashed{\Gamma}}_{BC}^A = \slashed{\nabla}_B \chi_{C}^{\phantom{C}A} + \slashed{\nabla}_C \chi_B^{\phantom{B}A} - \slashed{\nabla}^A \chi_{BC} 
 \end{equation*}
Another standard calculation (which can again be done in normal coordinates) leads to the expression
\begin{equation*}
 \dot{\slashed{R}}_{ABC}^{\phantom{ABC}D} = \slashed{\nabla}_A \dot{\slashed{\Gamma}}_{BC}^D - \slashed{\nabla}_B \dot{\slashed{\Gamma}}_{AC}^D
\end{equation*} 
and hence
\begin{equation*}
 \begin{split}
   \dot{\slashed{R}}_{ACB}^{\phantom{ABC}C} &= \slashed{\nabla}_A \dot{\slashed{\Gamma}}_{CB}^C - \slashed{\nabla}_C \dot{\slashed{\Gamma}}_{AB}^C \\
   &= \slashed{\nabla}_A \slashed{\nabla}_B \tr_{\slashed{g}}\chi - \slashed{\nabla}_C \slashed{\nabla}_A \chi_B^{\phantom{B}C} - \slashed{\nabla}_C \slashed{\nabla}_B \chi_A^{\phantom{A}C} + \slashed{\Delta}\chi_{AB}
  \end{split}
\end{equation*}
Next, we note that
\begin{equation*}
 \frac{\upd}{\upd \rho} \left( \slashed{R}_{ACB}^{\phantom{ACB}C} (\slashed{g}^{-1})^{AB} \right) = (\slashed{g}^{-1})^{AB} \dot{\slashed{R}}_{ACB}^{\phantom{ACB}C} + \slashed{R}_{ACB}^{\phantom{ACB}C} \dot{(\slashed{g}^{-1})}^{AB}
\end{equation*}
and now,
\begin{equation*}
 \dot{(\slashed{g}^{-1})}^{AB} = -2\chi^{AB}
\end{equation*}
Combining the above equations leads to the proposition.
 \end{proof}

\section{The transport equations for the quantities associated with the spheres \texorpdfstring{$S_{\tau,r}$}{S}}

We also need a evolution equations for the metric $\slashed{g}$ on the spheres $S_{\tau,r}$ and the scalar density $\Omega$. Additionally, we require evolution equations for the null frame components of the Christoffel symbols $\slashed{\Gamma}$. These are provided by the following propositions:

\begin{proposition}[The evolution equation for $\slashed{g}$]
\label{proposition transport equation for slashed g}
The metric $\slashed{g}$ satisfies the following transport equation in the outgoing $L$ direction:
\begin{equation}
 \slashed{\mathcal{L}}_L \left( r^{-2} \slashed{g}_{\mu\nu} - \mathring{\gamma}_{\mu\nu} \right) = 2r^{-2} \left( \chi_{(\text{small})} \right)_{\mu\nu}
\end{equation}
\end{proposition}

\begin{proof}
 The induced metric on the spheres $\slashed{g}$ satisfies
\begin{equation*}
 \slashed{\mathcal{L}}_L \slashed{g}_{\mu\nu} = 2\chi_{\mu\nu}
\end{equation*}
Hence we find
\begin{equation*}
 \begin{split}
  \slashed{\mathcal{L}}_L \left( r^{-2} \slashed{g}_{\mu\nu} \right) &= 2 r^{-2} \left(\chi_{\mu\nu} - r^{-1}\slashed{g}_{\mu\nu} \right) \\
  &= r^{-2} \left( \chi_{(\text{small})} \right)_{\mu\nu}
 \end{split}
\end{equation*}

On the other hand, by the definition of the metric $\mathring{\gamma}$, it is Lie-transported along the integral curves of $L$, i.e.\
\begin{equation*}
 \slashed{\mathcal{L}}_L \mathring{\gamma}_{\mu\nu} = 0
\end{equation*}

\end{proof}

\begin{proposition}[Evolution equation for $\Omega$]
\label{proposition transport Omega}
The scalar quantity $\Omega$ satisfies the evolution equation along the integral curves of $L$
\begin{equation}
 \begin{split}
  L \log \Omega &= \frac{1}{2} \tr_{\slashed{g}}\chi \\
  &= \frac{1}{r} + \frac{1}{2}\tr_{\slashed{g}}\chi_{(\text{small})}
 \end{split}
\end{equation}
\end{proposition}

\begin{proof}
Recall that we defined $\Omega$ by
\begin{equation*}
 \det{\slashed{g}} = \Omega^4 \det \mathring{\gamma}
\end{equation*}
Now, taking derivatives along the $L$ direction, and recalling that $\mathring{\gamma}$ is Lie transported along this direction, we find
\begin{equation*}
 \begin{split}
  \left((\slashed{g}^{-1})^{\mu\nu} \mathcal{L}_L \slashed{g}_{\mu\nu} \right) \det{\slashed{g}} &= 4 \Omega^3 (L\Omega) \det \mathring{\gamma} \\
  \Rightarrow (\slashed{g}^{-1})^{\mu\nu} \slashed{\mathcal{L}}_L \slashed{g}_{\mu\nu} &= 4 L \log \Omega
  \Rightarrow 2\tr_{\slashed{g}}\chi = 4 L \log \Omega
 \end{split}
\end{equation*}

\end{proof}

%
%Related to the previous proposition, we will also need to control the Christoffel symbols associated to the metric $\slashed{g}$ on the spheres.
%
%\begin{proposition}[The transport equation for the Christoffel symbols]
% The null frame components of the Christoffel symbols $\slashed{\Gamma}_{AB}^C$ satisfy the transport equations
% \begin{equation}
%  L\slashed{\Gamma}_{AB}^C = R_{LAB}^{\phantom{LAB}C} + \slashed{\nabla}_A \chi_B^{\phantom{B}C} + \frac{1}{2}\zeta^C \chi_{AB} - \frac{1}{2}\zeta_B \chi_A^{\phantom{A}C}
% \end{equation} 
% In particular, we have
% \begin{equation}
%  L \slashed{\Gamma}^B_{AB} = \slashed{\nabla}_A \tr_{\slashed{g}}\chi
% \end{equation}
%\end{proposition}
%
%\begin{proof}
% We compute
% \begin{equation*}
%  \begin{split}
%   L\left( \slashed{\Gamma}_{AB}^C \slashed{g}_{CD} \right)
%   &= L g(\D_A X_B, X_D) \\
%   &= g(\D_L \D_A X_B, X_D) + g(\D_A X_B, \D_L X_D) \\
%   &= R_{LABD} + X_A \chi_{BD} + \frac{1}{2}\zeta_D \chi_{AB} + \slashed{\Gamma}_{AB}^C \chi_{CD} - \frac{1}{2}\zeta_B \chi_{AD} - \slashed{\Gamma}^C_{AD} \chi_{CD}
%  \end{split}
% \end{equation*}
% but on the other hand we have
% \begin{equation*}
%  L\left( \slashed{\Gamma}_{AB}^C \slashed{g}_{CD} \right) = L\left( \slashed{\Gamma}_{AB}^C \right) \slashed{g}_{CD} +  2\slashed{\Gamma}_{AB}^C \chi_{CD}
% \end{equation*}
% which leads to the first part of the proposition. To prove the second part, we contract the indices $B$ and $C$, and use the symmetries of the Riemann tensor.
%\end{proof}

\chapter{The geometry of the vector bundle of \texorpdfstring{$S_{\tau,r}$}{S tau r} tangent tensor fields}
\label{chapter geometry of vector bundle}

We can regard an $S_{\tau,r}$ tangent tensor fields as a section of a particular vector bundle over $\mathcal{M}$, where the fibres are diffeomorphic to the cotangent space of $\mathbb{S}^2$. To put it another way, not only do these tensors live in the cotangent bundle of $\mathcal{M}$, but they live in a subset of the cotangent bundle which can itself be regarded as a vector bundle over $\mathcal{M}$. If we restrict to $S_{\tau,r}$-tangent one-forms, then we will call the associated bundle $\mathcal{B}$.

We can regard the differential operator $\slashed{\D}$ as providing a connection on the vector bundle $\mathcal{B}$. Indeed, it is easy to verify that $\slashed{\D}$ verifies the properties of a connection on this vector bundle: in particular, if $\phi$ is a section of $\mathcal{B}$ then $\slashed{\D}\phi$ is also a section of $\mathcal{B}$ (note that this is not the case for the $\D \phi$, so $\D$ does not provide us with a connection on $\mathcal{B}$, while $\slashed{\D}$ does). Indeed, we have the following elementary proposition:

\begin{proposition}
	$\slashed{\D}$ defines a metric connection on $\mathcal{B}$ with fibre metric $\slashed{g}^{-1}$.	
\end{proposition}

\begin{proof}
	Let $\phi_1$, $\phi_2$ be sections of $\mathcal{B}$. We can use $\slashed{g}^{-1}$ to define the inner product between $\phi_1$ and $\phi_2$ as follows: we first note that
	\begin{equation*}
	\mathcal{B} \subset T^*(\mathcal{M})
	\end{equation*}
	and so we can regard $\phi_1$ and $\phi_2$ as defining the one-forms $(\phi_1)_\mu$ and $(\phi_2)_{\mu}$ in the cotangent spaces of $\mathcal{M}$. Then we define the inner product
	\begin{equation*}
	\langle \phi_1 , \phi_2 \rangle := (\slashed{g}^{-1})^{\mu\nu} (\phi_1)_\mu (\phi_2)_\nu
	\end{equation*}
	Now we note that, since $\phi_1$ and $\phi_2 \in \mathcal{B}$ we can also write
	\begin{equation*}
	\langle \phi_1 , \phi_2 \rangle := (g^{-1})^{\mu\nu} (\phi_1)_\mu (\phi_2)_\nu
	\end{equation*}
	and so
	\begin{equation*}
	\begin{split}
	\upd_\mu \langle \phi_1 , \phi_2 \rangle 
	&=
	\upd_\mu \left(
	(\slashed{g}^{-1})^{\nu\rho} (\phi_1)_\nu (\phi_2)_\rho
	\right)
	\\
	&=(\slashed{g}^{-1})^{\nu\rho} \left( \D_\mu(\phi_1)_\nu \right) (\phi_2)_\rho
	+ (\slashed{g}^{-1})^{\nu\rho} (\phi_1)_\nu \left( \D_\mu(\phi_2)_\rho	\right)
	\\
	&=(\slashed{g}^{-1})^{\nu\rho} \left( \slashed{\D}_\mu(\phi_1)_\nu \right) (\phi_2)_\rho
	+ (\slashed{g}^{-1})^{\nu\rho} (\phi_1)_\nu \left( \slashed{\D}_\mu(\phi_2)_\rho \right)
	\\
	&= \langle \slashed{\D} \phi_1 , \phi_2 \rangle + \langle \phi_1 , \slashed{\D} \phi_2 \rangle
	\end{split}
	\end{equation*}
	Hence $\slashed{\D}$ is a metric connection on $\mathcal{B}$ with the metric $\slashed{g}^{-1}$.
\end{proof}

Note that, despite being a metric connection, $\slashed{\D}$ is \emph{not} simply the Levi-Civita connection associated with $\slashed{g}$, since this is defined on the Riemannian manifolds with metric $\slashed{g}$, namely the spheres $S_{\tau,r}$. Indeed, this connection is precisely $\slashed{\nabla}$. These two connections agree when evaluated in the directions of $S_{\tau, r}$-tangent vector fields, i.e.\
\begin{equation*}
\slashed{\nabla}_X \phi = \slashed{\D}_X \phi \quad \text{if } X \in T(S_{\tau,r})
\end{equation*}
However, if $X$ is not tangent to the spheres then the two connections differ. For example, $\slashed{\nabla}_L \phi = 0$, while in general $\slashed{\D}_L \phi \neq 0$.

When commuting the wave equation, we will encounter terms of the form
\begin{equation*}
[\slashed{\D}_\mu , \slashed{\D}_\nu] \phi
\end{equation*}
where $\phi \in \mathcal{B}$. These terms are naturally expressed in terms of the curvature of the connection $\slashed{\D}$ (which is \emph{not} the Riemann curvature of $\mathcal{M}$, which is the curvature of the connection $\D$!). Moreover, since the basic quantities we will be able to estimate are the rectangular components of the metric perturbations $h_{ab}$, it is essential that we express this curvature in terms of the derivatives of the fields $h_{ab}$.

\begin{proposition}[The connection coefficients in the basis $\slashed{\nabla} x^a$]
	Relative to the basis $\slashed{\nabla}x^a$, the connection coefficients associated with $\slashed{\D}$ can be expressed as
	\begin{equation}
	\omega_b^{\phantom{b}a}
	=
	\slashed{\Pi}_c^{\phantom{c}a} \upd \slashed{\Pi}_b^{\phantom{b}c}
	- \slashed{\Pi}_d^{\phantom{d}a} \Gamma^d_{cb} \upd x^c
	\end{equation}
	so that the connection itself can be written as
	\begin{equation}
	\slashed{\D} \phi = \left( \upd \phi_a + \omega_a^{\phantom{a}b} \phi_b \right) \slashed{\nabla} x^a
	\end{equation}
	for any section $\phi$ of $\mathcal{B}$.
	
\end{proposition}

\begin{proof}
	Let $\phi$ be a section of $\mathcal{B}$. Then we can expand $\phi$ in terms of the sections $\slashed{\nabla} x^a$, where $x^a$ are the rectangular coordinate functions, as follows:
	\begin{equation*}
	\phi = \phi_a \slashed{\nabla} x^a
	\end{equation*}
	
	Note that there are four sections $\slashed{\nabla}x^a$, for $a = 0,1,\ldots 4$ which clearly span the space of sections of $\mathcal{B}$. However, at each point in $\mathcal{M}$ the space of sections is two dimensional, so it would be sufficient to consider a basis consisting of a pair of sections. Nevertheless, the components of $\phi$, $\phi_a$ are uniquely defined by
	\begin{equation*}
	\phi_a := \phi \cdot \partial_a
	\end{equation*}
	
	Note that we can also write
	\begin{equation*}
	\phi = \phi_a \upd x^a 
	\end{equation*}
	and, since $\phi$ is actually $S_{\tau, r}$-tangent, these two expressions are equivalent.
	
	Additionally, note that, since $\slashed{\Pi}_\mu^{\phantom{\mu}\nu} \phi_\nu = \phi_\mu$ it follows that
	\begin{equation*}
	\Pi_a^{\phantom{a}b} \phi_b = \phi_a
	\end{equation*}
	
	Now, we have
	\begin{equation}
	\label{equation defining connection coefficients}
	\begin{split}
	\slashed{\D} \phi
	&=
	\slashed{\D} \left( \phi_a \slashed{\nabla} x^a \right)
	\\
	&= \left( \upd \phi_a + \omega_{a}^{\phantom{a}b} \phi_b \right)\slashed{\nabla} x^a
	\end{split}
	\end{equation}	
	where $\omega_a^{\phantom{a}b}$ are the connection coefficients of $\slashed{\D}$ in the basis of sections $\slashed{\nabla} x^a$; a set of one-forms in $T^*(\mathcal{M})$. Indeed, this expression defines the connection coefficients relative to this basis. The connection coefficients can be found by using the relationship
	\begin{equation*}
	\begin{split}
	(\omega_b^{\phantom{b}a})_\mu \slashed{\nabla}_\nu x^b 
	&=
	\slashed{\D}_\mu \slashed{\nabla}_\nu x^a
	\\
	&=
	\slashed{\Pi}_\nu^{\phantom{\nu} \rho} \D_\mu \left( \slashed{\Pi}_\rho^{\phantom{\rho}\sigma} \partial_\sigma x^a \right)
	\\
	&= 
	\slashed{\Pi}_\nu^{\phantom{\nu} \rho} \D_\mu \left( \slashed{\Pi}_\rho^{\phantom{\rho}a} \right)
	\\
	&=
	\slashed{\Pi}_\nu^{\phantom{\nu} b} (\partial_b)^\rho \D_\mu \left( \slashed{\Pi}_\rho^{\phantom{\rho}a} \right)
	\\
	&=
	\slashed{\Pi}_\nu^{\phantom{\nu}b} \partial_\mu \slashed{\Pi}_b^{\phantom{b}a}
	- \slashed{\Pi}_\nu^{\phantom{\nu}b} \slashed{\Pi}_\rho^{\phantom{\rho}a} \D_\mu (\partial_b)^\rho
	\\
	&=
	\slashed{\Pi}_\nu^{\phantom{\nu}b} \left( \partial_\mu \slashed{\Pi}_b^{\phantom{b}a}
	- \slashed{\Pi}_\rho^{\phantom{\rho}a} (\partial_\mu x^c) \Gamma_{cb}^d (\partial_d)^\rho \right)
	\\
	&=
	\slashed{\Pi}_\nu^{\phantom{\nu}b} \left( \partial_\mu \slashed{\Pi}_b^{\phantom{b}a}
	- \slashed{\Pi}_d^{\phantom{d}a} \Gamma_{cb}^d (\partial_\mu x^c) \right)
	\end{split}
	\end{equation*}
	Note that $\slashed{\Pi}_\nu^{\phantom{\nu}b} = \slashed{\nabla}_\nu x^b$. Hence we can write
	\begin{equation*}
	\omega_b^{\phantom{b}a}
	=
	\upd \slashed{\Pi}_b^{\phantom{b}a}
	- \slashed{\Pi}_d^{\phantom{d}a} \Gamma^d_{cb} \upd x^c
	\end{equation*}
	
	Finally, we note that we are free to project the ``upper'' rectangular index of the connection coefficient using the projection operator $\slashed{\Pi}$ expressed in the rectangular frame; this follows from the defining equation for the connection coefficients (equation \eqref{equation defining connection coefficients}) together with the fact that $\slashed{\Pi}_a^{\phantom{a}b} \phi_b = \phi_a$. Hence, we have
	\begin{equation*}
	\omega_b^{\phantom{b}a}
	=
	\slashed{\Pi}_c^{\phantom{c}a} \upd \slashed{\Pi}_b^{\phantom{b}c}
	- \slashed{\Pi}_d^{\phantom{d}a} \Gamma^d_{cb} \upd x^c
	\end{equation*}
	
\end{proof}

Note that the connection coefficients of $\slashed{\D}$, acting on sections of the cotangent bundle of $\mathcal{M}$, relative to the basis $\upd x^a$ is given by $\Gamma^a_{bc} \upd x^c$. This is evidently \emph{not} the same as the connection coefficients $\omega_b^{\phantom{b}a}$, which is yet another reminder that the connection $\slashed{\D}$ differs from the connection $\D$.

\begin{proposition}[An alternative expression for the connection $\slashed{\D}$]
	We can also express the connection $\slashed{\D}$ as
	\begin{equation}
	\slashed{\D} \phi = \left(\slashed{\Pi}_a^{\phantom{a}c} \upd \phi_c + \tilde{\omega}_a^{\phantom{c}b} \phi_b \right) \slashed{\nabla} x^a
	\end{equation}
	where
	\begin{equation*}
	\tilde{\omega}_a^{\phantom{a}b} := \slashed{\Pi}_a^{\phantom{a}c} \omega_c^{\phantom{c}b}
	\end{equation*}
	are a collection of one-forms on $\mathcal{M}$ which satisfy
	\begin{equation*}
	\tilde{\omega}_a^{\phantom{a}b} = \slashed{\Pi}_a^{\phantom{a}c} \slashed{\Pi}_d^{\phantom{d}b} \tilde{\omega}_c^{\phantom{c}d}
	\end{equation*}
\end{proposition}

\begin{proof}
	We saw above that the connection can be expressed as
	\begin{equation*}
	\slashed{\D} \phi = \left(\upd \phi_a + \omega_a^{\phantom{a}b} \phi_b \right) \slashed{\nabla} x^a	
	\end{equation*}
	where $\phi$ is a section section of $\mathcal{B}$ with components $\phi_a$ given by
	\begin{equation*}
	\phi = \phi_a \slashed{\nabla} x^a
	\end{equation*}
	
	Now, we have
	\begin{equation*}
	\begin{split}
	\left(\upd \phi_a + \omega_a^{\phantom{a}b} \phi_b \right) \slashed{\nabla}_\mu x^a
	&=
	\left(\upd \phi_a + \omega_a^{\phantom{a}b} \phi_b \right) \slashed{\Pi}_\mu^{\phantom{\mu}a}
	\\
	&=
	\left(\upd \phi_a + \omega_a^{\phantom{a}b} \phi_b \right) \slashed{\Pi}_\mu^{\phantom{\mu}c} \slashed{\Pi}_c^{\phantom{c}a}
	\\
	&=
	\left(\slashed{\Pi}_a^{\phantom{a}c} \upd \phi_c + \slashed{\Pi}_a^{\phantom{a}c}\omega_c^{\phantom{c}b} \phi_b \right) \slashed{\Pi}_\mu^{\phantom{\mu}a}
	\\
	&=
	\left(\slashed{\Pi}_a^{\phantom{a}c} \upd \phi_c + \tilde{\omega}_a^{\phantom{c}b} \phi_b \right) \slashed{\nabla}_\mu x^a
	\end{split}
	\end{equation*}
	
\end{proof}

\begin{proposition}[The curvature of the connection $\slashed{\D}$]
	Define the curvature of the connection $\slashed{\D}$ as the collection of two forms on $\mathcal{M}$ with components
	\begin{equation}
	\left(\Omega_a^{\phantom{a}b} \right)_{\mu\nu} \phi_b \slashed{\nabla} x^a := [\slashed{\D}_\mu , \slashed{\D}_\nu]\phi
	\end{equation}
	where $\phi$ is any sufficiently smooth section of $\mathcal{B}$. Then we have
	\begin{equation}
	\Omega_a^{\phantom{a}b} = \slashed{\Pi}_a^{\phantom{a}c} \slashed{\Pi}_d^{\phantom{d}b} \left(\upd \tilde{\omega}_c^{\phantom{c}d} + \upd \slashed{\Pi}_c^{\phantom{c}e} \wedge \upd \slashed{\Pi}_e^{\phantom{e}d} \right) + \tilde{\omega}_a^{\phantom{a}c}\wedge \tilde{\omega}_c^{\phantom{c}b}
	\end{equation}
\end{proposition}

\begin{proof}
	We can extend the operator $\slashed{\D}$ to act on sections of the cotangent bundle of $\mathcal{M}$ by simply using the operator $\D$ in this case\footnote{In fact, the curvature of $\slashed{\D}$, is independent of the choice of connection on the cotangent bundle.}. Then we find
	\begin{equation*}
	\begin{split}
	\slashed{\D}_\mu \slashed{\D}_\nu \phi
	&=
	\slashed{\D}_\mu \left( \left(\slashed{\Pi}_a^{\phantom{a}b} \D_\nu \phi_b + \left(\tilde{\omega}_a^{\phantom{a}b} \right)_\nu \phi_b \right) \slashed{\nabla} x^a \right)
	\\
	&=
	\left(\slashed{\Pi}_a^{\phantom{a}b} \D_\mu \left(\slashed{\Pi}_b^{\phantom{b}c} \D_\nu \phi_c + \left(\tilde{\omega}_b^{\phantom{b}c} \right)_\nu \phi_c \right)  
	+ \left(\tilde{\omega}_a^{\phantom{a}c}\right)_\mu \left(\slashed{\Pi}_c^{\phantom{c}b} \D_\nu \phi_b + \left(\tilde{\omega}_c^{\phantom{c}b} \right)_\nu \phi_b \right) \right) \slashed{\nabla} x^a 
	\\
	&=
	\bigg( \slashed{\Pi}_a^{\phantom{a}b} \D_\mu \left(\slashed{\Pi}_b^{\phantom{b}c} \right) \D_\nu \phi_c
	+ \slashed{\Pi}_a^{\phantom{a}c} \D^2_{\mu\nu}\phi_c
	+ \slashed{\Pi}_a^{\phantom{a}b} \left( \D_\mu \left( \tilde{\omega}_b^{\phantom{b}c}\right)_\nu \right) \phi_c
	+ \left( \tilde{\omega}_a^{\phantom{a}c} \right)_\nu \D_\mu \phi_c
	+ \left( \tilde{\omega}_a^{\phantom{a}c} \right)_\mu \D_\nu \phi_c
	\\
	&\phantom{=}
	+ \left( \tilde{\omega}_a^{\phantom{a}c} \right)_\mu \left( \tilde{\omega}_c^{\phantom{c}b} \right)_\nu \phi_b \bigg) \slashed{\nabla} x^a
	\end{split}
	\end{equation*}
	and so it is fairly easy to see that
	\begin{equation*}
	[\slashed{\D}_\mu, \slashed{\D}_\nu]\phi
	= \left(
	\slashed{\Pi}_a^{\phantom{a}c} \upd \slashed{\Pi}_c^{\phantom{c}b} \wedge \upd \phi_b
	+ \slashed{\Pi}_a^{\phantom{a}c} \upd \tilde{\omega}_c^{\phantom{c}b} \phi_b
	+ \tilde{\omega}_a^{\phantom{a}c} \wedge \tilde{\omega}_c^{\phantom{c}b} \phi_b \right)_{\mu\nu} \slashed{\nabla} x^a
	\end{equation*}
	We can write the first term on the right hand side in a more convenient form as follows: we have
	\begin{equation*}
	\begin{split}
	\phi_a &= \slashed{\Pi}_a^{\phantom{a}b} \phi_b \\
	\Rightarrow \upd \phi_a &= \left(\upd \slashed{\Pi}_a^{\phantom{a}b}\right) \phi_b + \slashed{\Pi}_a^{\phantom{a}b} \upd\phi_b \\
	\Rightarrow \slashed{\Pi}_a^{\phantom{a}b}\upd \phi_b &= \slashed{\Pi}_a^{\phantom{a}c} \left(\upd \slashed{\Pi}_c^{\phantom{c}b}\right) \phi_b + \slashed{\Pi}_a^{\phantom{a}b} \upd \phi_b \\
	\Rightarrow \slashed{\Pi}_a^{\phantom{a}c} \left(\upd \slashed{\Pi}_c^{\phantom{c}b}\right) \phi_b &= 0 \\
	\Rightarrow \left( \upd \slashed{\Pi}_a^{\phantom{a}c} \wedge \upd \slashed{\Pi}_c^{\phantom{c}b} \right) \phi_b + \slashed{\Pi}_a^{\phantom{a}c} \upd \phi_b \wedge \upd \slashed{\Pi}_c^{\phantom{c}b} &= 0
	\end{split}
	\end{equation*}
	so we have
	\begin{equation*}
	[\slashed{\D}_\mu , \slashed{\D}_\nu]\phi 
	=
	\left( \upd \slashed{\Pi}_a^{\phantom{a}c} \wedge \upd \slashed{\Pi}_c^{\phantom{c}b} 
	+ \slashed{\Pi}_a^{\phantom{a}c} \upd \tilde{\omega}_c^{\phantom{c}b}
	+ \tilde{\omega}_a^{\phantom{a}c} \wedge \tilde{\omega}_c^{\phantom{c}b} \right) \phi_b \slashed{\nabla} x^a
	\end{equation*}
	Finally, we note that we are free to use the rectangular components of the  projection operators to ``project'' the index $b$ (since $\phi_b = \slashed{\Pi}_b^{\phantom{b}c} \phi_c$) and the index $a$ (since $\slashed{\nabla} x^a = \slashed{\Pi}_c^{\phantom{c}a} \slashed{\nabla} x^c$), proving the proposition.
	
\end{proof}

%
%
%\begin{proposition}[The curvature of the connection $\slashed{\D}$]
%	
%Expressed in terms of the sections $\slashed{\nabla} x^a$, the curvature of the connection $\slashed{\D}$ is given by
%\begin{equation*}
%	\Omega_b^{\phantom{b}a}
%	=
%	-\upd(\slashed{\Pi}_d^{\phantom{d}a} \Gamma^d_{cb}) \wedge \upd x^c
%	+ \left( \upd \slashed{\Pi}_b^{\phantom{b}c}
%	- \slashed{\Pi}_d^{\phantom{d}c} \Gamma^d_{be} \upd x^e \right)
%	\wedge
%	\left( \upd \slashed{\Pi}_c^{\phantom{c}a}
%	- \slashed{\Pi}_f^{\phantom{f}a} \Gamma^f_{cg} \upd x^g \right)
%\end{equation*}
%	
%\end{proposition}
%
%\begin{proof}
%	In terms of the connection coefficients, we can write the curvature of $\slashed{\D}$ as the collection of two-forms on $\mathcal{M}$
%	\begin{equation*}
%	\Omega_b^{\phantom{b}a} := \upd \omega_b^{\phantom{b}a} + \omega_b^{\phantom{b}c} \wedge \omega_c^{\phantom{c}a}
%	\end{equation*}
%	and, in this case, for $\phi$ a section of $\mathcal{B}$ we have
%	\begin{equation*}
%	[\slashed{\D}_\mu, \slashed{\D}_\nu] \phi = \left(\Omega_b^{\phantom{b}a} \right)_{\mu\nu} \phi_a \slashed{\nabla} x^b
%	\end{equation*}
%	
%	Using the expression above for the connection coefficients, we find that
%	\begin{equation*}
%	\Omega_b^{\phantom{b}a}
%	=
%	-\upd(\slashed{\Pi}_d^{\phantom{d}a} \Gamma^d_{cb}) \wedge \upd x^c
%	+ \left( \upd \slashed{\Pi}_b^{\phantom{b}c}
%		- \slashed{\Pi}_d^{\phantom{d}c} \Gamma^d_{be} \upd x^e \right)
%	\wedge
%	\left( \upd \slashed{\Pi}_c^{\phantom{c}a}
%	- \slashed{\Pi}_f^{\phantom{f}a} \Gamma^f_{cg} \upd x^g \right)
%	\end{equation*}
%	
%\end{proof}

\begin{proposition}[The connection coefficients $\tilde{\omega}$ expressed in terms of derivatives of the rectangular components of the metric]
	\label{proposition connection coefficients omega in terms of the metric}
	Using the notation
	\begin{equation*}
	\begin{split}
	(X\slashed{h})_{Y\mu} &:= (Xh_{ab})Y^a \slashed{\Pi}_\mu^{\phantom{\mu}b}
	\\
	(X\slashed{h})_{\mu \nu} &:= (Xh_{ab}) \slashed{\Pi}_\mu^{\phantom{\mu}a} \slashed{\Pi}_\nu^{\phantom{\mu}b}
	\\
	(\slashed{\nabla}_\mu \slashed{h})_{X\nu} &:= (\slashed{\nabla}_\mu h_{ab})X^a \slashed{\Pi}_\nu^{\phantom{\nu}b}
	\\
	(\slashed{\nabla}_\mu \slashed{h})_{\nu \rho} &:= (\slashed{\nabla}_\mu h_{ab}) \slashed{\Pi}_\nu^{\phantom{\nu}a} \slashed{\Pi}_\rho^{\phantom{\rho}b}
	\end{split}
	\end{equation*}
	for any vector fields $X$ and $Y \in T(\mathcal{M})$, we find
	\begin{equation}
	\begin{split}
	\left(\tilde{\omega}_b^{\phantom{b}a}\right)_\mu
	&=
	L_\mu \slashed{\Pi}_b^{\phantom{b}\nu} \slashed{\Pi}_\rho^{\phantom{\rho}a} (\slashed{g}^{-1})^{\rho \sigma} \left(
		\frac{1}{4}(\slashed{\nabla}_\nu \slashed{h})_{\Lbar \sigma}
		+ \frac{1}{4}(\Lbar \slashed{h})_{\nu\sigma}
		- \frac{1}{4}(\slashed{\nabla}_\sigma \slashed{h})_{\Lbar \nu}
	\right)
	\\
	&\phantom{=}
	+ \Lbar_\mu \slashed{\Pi}_b^{\phantom{b}\nu} \slashed{\Pi}_{\rho}^{\phantom{\rho}a} (\slashed{g}^{-1})^{\rho\sigma} \left(
	\frac{1}{4}(\slashed{\nabla}_\nu \slashed{h})_{L\sigma}
	+ \frac{1}{4}(L \slashed{h})_{\nu\sigma} 
	- \frac{1}{4}(\slashed{\nabla}_\sigma \slashed{h})_{L\nu}
	\right)
	\\
	&\phantom{=}
	+ \slashed{\Pi}_\mu^{\phantom{\mu}c} \slashed{\Pi}_b^{\phantom{b}\nu} \slashed{\Pi}_\rho^{\phantom{\rho} a} (\slashed{g}^{-1})^{\rho\sigma} \slashed{\Pi}_c^{\phantom{c}\lambda} \left(
		- \frac{1}{2}(\slashed{\nabla}_\nu \slashed{h})_{\lambda\sigma}
		- \frac{1}{2}(\slashed{\nabla}_\lambda \slashed{h})_{\nu\sigma}
		+ \frac{1}{2}(\slashed{\nabla}_\sigma \slashed{h})_{\lambda\nu}
		\right)
	\end{split}
	\end{equation}

\end{proposition}

\begin{proof}
	We begin with the expression
	\begin{equation*}
	\left(\tilde{\omega}_b^{\phantom{b}a}\right)_\mu = \slashed{\Pi}_b^{\phantom{b}c} \slashed{\Pi}_d^{\phantom{d}a} \left( \partial_\mu \slashed{\Pi}_c^{\phantom{c}d} - \Gamma^d_{ce} \partial_\mu x^e \right)
	\end{equation*}
	Note that, from proposition \ref{proposition derivatives of the rectangular components of Pi}, the first term vanishes, since $\upd \slashed{\Pi}_a^{\phantom{a}b}$ vanishes when we project \emph{both} its indices using $\slashed{\Pi}$.
	
	Expanding in terms of the null frame, we have
	\begin{equation*}
	\begin{split}
	\left(\omega_b^{\phantom{b}a}\right)_\mu 
	&=
	\frac{1}{2}L_\mu \Lbar^e \slashed{\Pi}_b^{\phantom{b}c}\slashed{\Pi}_d^{\phantom{d}a} \Gamma^d_{ce}
	+ \frac{1}{2}\Lbar_\mu L^e \slashed{\Pi}_b^{\phantom{b}c}\slashed{\Pi}_d^{\phantom{d}a} \Gamma^d_{ce}
	- \slashed{\Pi}_\mu^{\phantom{\mu}c} \slashed{\Pi}_b^{\phantom{b}d}\slashed{\Pi}_e^{\phantom{e}a} \Gamma^{e}_{cd}
	\end{split}
	\end{equation*}
	Expanding the rectangular Christoffel symbols in terms of derivatives of $h$ proves the proposition.

\end{proof}

\begin{proposition}[An expression for the curvature of the connection $\slashed{\D}$]
	\label{proposition expression for Omega}
The curvature of the connection $\slashed{\D}$, expressed relative to the basis of sections of $\mathcal{B}$ given by $\slashed{\nabla} x^a$ and in terms of derivatives of the rectangular components of $h$ when possible, has the following form:
\begin{equation}
	\begin{split}
	\left( \Omega_b^{\phantom{b}a} \right)_{\mu\nu}
	&=
	\bigg( 
		\frac{1}{8} r^{-1} (\slashed{g}^{-1})^{\lambda\sigma}\slashed{\Pi}_b^{\phantom{b}\rho} \slashed{\Pi}_\sigma^{\phantom{\sigma}a} 
		\Big( 
			\left( \slashed{\D}_L (r\slashed{\nabla}_\rho \slashed{h}) \right)_{\Lbar \lambda} 
			- \left( \slashed{\D}_L (r\slashed{\nabla}_\lambda \slashed{h}) \right)_{\Lbar \rho}
			- \left(\slashed{\D}_{\Lbar} (r\slashed{\nabla}_\rho \slashed{h}) \right)_{L\lambda}
			\\
			&\phantom{=	\bigg( \frac{1}{8} r^{-1} (\slashed{g}^{-1})^{\lambda\sigma}\slashed{\Pi}_b^{\phantom{b}\rho} \slashed{\Pi}_\sigma^{\phantom{\sigma}a} \Big( }
			+  \left(\slashed{\D}_{\Lbar} (r\slashed{\nabla}_\lambda \slashed{h}) \right)_{L\rho}
		\Big)
		+ \textit{l.o.t.}_{(\Omega , L\wedge\Lbar)}
		\bigg)L_\mu \wedge \Lbar_\nu
	\\
	&\phantom{=}
	+ \bigg(
		\frac{1}{4}r^{-1} \slashed{\Pi}_c^{\phantom{c}\rho} \slashed{\Pi}_b^{\phantom{b}\sigma} \slashed{\Pi}_d^{\phantom{d}\delta} (\slashed{g}^{-1})^{ad} \Big(
			\left(\slashed{\nabla}_\rho (r\slashed{\nabla}_\delta \slashed{h})\right)_{\Lbar\sigma}
			- \left(\slashed{\nabla}_\rho (r\slashed{\nabla}_\sigma \slashed{h}) \right)_{\Lbar\delta}
			- \left(\slashed{\D}_{\Lbar} (r\slashed{\nabla}_\delta \slashed{h})\right)_{\rho\sigma}
			\\
			&\phantom{= + \bigg(
				\frac{1}{4}r^{-1} \slashed{\Pi}_c^{\phantom{c}\rho} \slashed{\Pi}_b^{\phantom{b}\sigma} \slashed{\Pi}_d^{\phantom{d}\delta} (\slashed{g}^{-1})^{ad} \Big(}
			+ \left(\slashed{\D}_{\Lbar} (r\slashed{\nabla}_\sigma \slashed{h})\right)_{\rho\delta}
		\Big)
		+ \textit{l.o.t.}_{(\Omega , L\wedge\slashed{\Pi})}
		\bigg)L_\mu \wedge \slashed{\Pi}_\nu^{\phantom{\nu}c}
	\\
	&\phantom{=}
	+ \bigg(
		\frac{1}{4}r^{-1} \slashed{\Pi}_c^{\phantom{c}\rho} \slashed{\Pi}_b^{\phantom{b}\sigma} \slashed{\Pi}_d^{\phantom{d}\delta} (\slashed{g}^{-1})^{ad} \Big(
			\left(\slashed{\nabla}_\rho (r\slashed{\nabla}_\delta \slashed{h})\right)_{L\sigma}
			- \left(\slashed{\nabla}_\rho (r\slashed{\nabla}_\sigma \slashed{h})\right)_{L\delta}
			- \left(\slashed{\D}_{L} (r\slashed{\nabla}_\delta \slashed{h})\right)_{\rho\sigma}
			\\
			&\phantom{=	+ \bigg( \frac{1}{4}r^{-1} \slashed{\Pi}_c^{\phantom{c}\rho} \slashed{\Pi}_b^{\phantom{b}\sigma} \slashed{\Pi}_d^{\phantom{d}\delta} (\slashed{g}^{-1})^{ad} \Big(}
			+ \left(\slashed{\D}_{L} (r\slashed{\nabla}_\sigma \slashed{h})\right)_{\rho\delta}
		\Big)
		+ \textit{l.o.t.}_{(\Omega , \Lbar\wedge\slashed{\Pi})}
		\bigg)\Lbar_\mu \wedge \slashed{\Pi}_\nu^{\phantom{\nu}c}
	\\
	&\phantom{=}
	+ \bigg(
		\frac{1}{2}r^{-1} \slashed{\Pi}_c^{\phantom{c}\rho} \slashed{\Pi}_b^{\phantom{b}\sigma} \slashed{\Pi}_e^{\phantom{e}\delta} \slashed{\Pi}_d^{\phantom{d}\lambda} (\slashed{g}^{-1})^{ae} \left(
		\left(\slashed{\nabla}_\lambda (r\slashed{\nabla}_\sigma \slashed{h})\right)_{\rho\delta}
		- \left(\slashed{\nabla}_\lambda (r\slashed{\nabla}_\delta \slashed{h})\right)_{\rho\sigma}
		\right)
		+ \textit{l.o.t.}_{(\Omega , \slashed{\Pi}\wedge\slashed{\Pi})}
		\bigg)\slashed{\Pi}_\mu^{\phantom{\mu}c} \wedge \slashed{\Pi}_\nu^{\phantom{\nu}d}
	\end{split}
\end{equation}
where the lower order terms are given schematically as
\begin{equation*}
	\begin{split}
	\textit{l.o.t.}_{(\Omega , L\wedge\Lbar)}
	&=
	r^{-1} (\partial h)_{(\text{frame})}
	+ \begin{pmatrix}
		(\partial h)_{(\text{frame})} \\
		\slashed{\nabla} \log \mu \\
		r^{-1} L^i_{(\text{small})} \slashed{\Pi}^i
		\end{pmatrix}
	\cdot \begin{pmatrix}
		(\partial h)_{(\text{frame})} \\
		\slashed{\nabla} \log \mu \\
		r^{-1} L^i_{(\text{small})} \slashed{\Pi}^i
		\end{pmatrix}
	\\
	\\
	\textit{l.o.t.}_{(\Omega , L\wedge\slashed{\Pi})}
	&=
	r^{-1} (\slashed{\nabla} h)_{(\text{frame})}
	+ \begin{pmatrix}
		(\slashed{\nabla} h)_{(\text{frame})} \\
		(\Lbar h)_{(\text{frame})} \\
		(\slashed{\nabla} \log \mu) \\
		\chibar \\
		r^{-1}\slashed{\Pi}^i \otimes \slashed{\Pi}^i \\
		r^{-1} L^i_{(\text{small})} \slashed{\Pi}^i
		\end{pmatrix}
	\cdot \begin{pmatrix}
		(\slashed{\nabla} h)_{(\text{frame})} \\
		(\Lbar h)_{(\text{frame})} \\
		\end{pmatrix}
	+ \chibar \begin{pmatrix}
		\slashed{\nabla} \log \mu \\
		r^{-1} L^i_{(\text{small})} \slashed{\Pi}^i
		\end{pmatrix}
	\\
	&\phantom{=}
	+ (\slashed{\nabla} \log \mu) (r^{-1}\slashed{\Pi}^i \otimes \slashed{\Pi}^i)
	\\
	\\
	\textit{l.o.t.}_{(\Omega , \Lbar\wedge\slashed{\Pi})}
	&=
	r^{-1}(\slashed{\nabla} h)_{(\text{frame})}
	+ \begin{pmatrix}
		(\bar{\partial} h)_{(\text{frame})} \\		
		(r^{-1} L^i_{(\text{small})} \slashed{\Pi}^i)
	\end{pmatrix} \cdot
	\begin{pmatrix}
		(\bar{\partial} h)_{(\text{frame})} \\
		\chi \\
		r^{-1} \slashed{\Pi}^i \otimes \slashed{\Pi}^i
	\end{pmatrix}
	\\
	\\
	\textit{l.o.t.}_{(\Omega , \slashed{\Pi}\wedge\slashed{\Pi})}
	&=
	\begin{pmatrix}
		(\bar{\partial} h)_{(\text{frame})} \\
		\chi \\
		r^{-1} \slashed{\Pi}^i \otimes \slashed{\Pi}^i
	\end{pmatrix} \cdot
	\begin{pmatrix}
		(\bar{\partial} h)_{(\text{frame})} \\
		\chi
	\end{pmatrix}
	\end{split}
\end{equation*}

\end{proposition}

\begin{proof}
	Recall the expression
	\begin{equation*}
	\left(\Omega_a^{\phantom{a}b} \right)_{\mu\nu}
	= \slashed{\Pi}_a^{\phantom{a}c} \slashed{\Pi}_d^{\phantom{d}b} \left( \upd \tilde{\omega}_c^{\phantom{c}d} + \upd \slashed{\Pi}_c^{\phantom{c}e} \wedge \upd \slashed{\Pi}_e^{\phantom{e}d} + \tilde{\omega}_c^{\phantom{c}e} \wedge \tilde{\omega}_e^{\phantom{e}d} \right)
	\end{equation*}
	We will focus on each of the three terms on the right in turn. First, we calculate $\upd \omega$, which includes the highest order terms. We have
	\begin{equation*}
	\begin{split}
	\left( \upd \tilde{\omega}_b^{\phantom{b}a} \right)_{\mu\nu}
	&=
	(\upd L^{\flat})_{\mu\nu} \slashed{\Pi}_b^{\phantom{b}\rho} \slashed{\Pi}_\sigma^{\phantom{\sigma}a} (\slashed{g}^{-1})^{\sigma\lambda} \left(
	\frac{1}{4}(\slashed{\nabla}_\rho \slashed{h})_{\Lbar \lambda}
	+ \frac{1}{4}(\Lbar\slashed{h})_{\rho\lambda}
	- \frac{1}{4}(\slashed{\nabla}_\lambda \slashed{h})_{\Lbar \rho}
	\right)
	\\
	&\phantom{=}
	+ (\upd \Lbar^{\flat})_{\mu\nu} \slashed{\Pi}_b^{\phantom{b}\rho} \slashed{\Pi}_{\sigma}^{\phantom{\sigma}a} (\slashed{g}^{-1})^{\sigma\lambda} \left(
	\frac{1}{4}(\slashed{\nabla}_\rho \slashed{h})_{L\lambda}
	+ \frac{1}{4}(L\slashed{h})_{\rho\lambda} 
	- \frac{1}{4}(\slashed{\nabla}_\lambda \slashed{h})_{L\rho}
	\right)
	\\
	&\phantom{=}
	+ (\upd \slashed{\Pi}^{c})_{\mu\nu} \slashed{\Pi}_b^{\phantom{b}\rho} \slashed{\Pi}_\sigma^{\phantom{\sigma} a} (\slashed{g}^{-1})^{\sigma\lambda} \slashed{\Pi}_c^{\phantom{c}\delta} \left(
	\frac{1}{2}(\slashed{\nabla}_\rho \slashed{h})_{\delta\lambda}
	+ \frac{1}{2}(\slashed{\nabla}_\delta \slashed{h})_{\rho\lambda}
	- \frac{1}{2}(\slashed{\nabla}_\lambda \slashed{h})_{\delta\rho}
	\right)
	\\
	&\phantom{=}
	+ \frac{1}{2} L \left( \slashed{\Pi}_b^{\phantom{b}\rho} \slashed{\Pi}_\sigma^{\phantom{\sigma}a} (\slashed{g}^{-1})^{\sigma \lambda} \left(
	\frac{1}{4}(\slashed{\nabla}_\rho \slashed{h})_{\Lbar \lambda}
	+ \frac{1}{4}(\Lbar \slashed{h})_{\rho\lambda}
	- \frac{1}{4}(\slashed{\nabla}_\lambda \slashed{h})_{\Lbar \rho}
	\right) \right) L_\mu \wedge \Lbar_\nu
	\\
	&\phantom{=}
	- \slashed{\Pi}_c^{\phantom{c}\rho} \slashed{\nabla}_\rho \left( \slashed{\Pi}_b^{\phantom{b}\sigma} \slashed{\Pi}_\lambda^{\phantom{\lambda}a} (\slashed{g}^{-1})^{\lambda \delta} \left(
	\frac{1}{4}(\slashed{\nabla}_\sigma \slashed{h})_{\Lbar \delta}
	+ \frac{1}{4}(\Lbar\slashed{h})_{\sigma\delta}
	- \frac{1}{4}(\slashed{\nabla}_\delta \slashed{h})_{\Lbar \sigma}
	\right) \right) L_\mu \wedge \slashed{\Pi}_\nu^{\phantom{\nu}c}
	\\
	&\phantom{=}
	-\frac{1}{2} \Lbar \left( \slashed{\Pi}_b^{\phantom{b}\rho} \slashed{\Pi}_{\sigma}^{\phantom{\rho}a} (\slashed{g}^{-1})^{\sigma\lambda} \left(
	\frac{1}{4}(\slashed{\nabla}_\rho \slashed{h})_{L\lambda}
	+ \frac{1}{4}(L\slashed{h})_{\rho\lambda}
	- \frac{1}{4}(\slashed{\nabla}_\lambda \slashed{h})_{L\rho}
	\right) \right) L_\mu \wedge \Lbar_\nu
	\\
	&\phantom{=}
	- \slashed{\Pi}_c^{\phantom{c}\rho} \slashed{\nabla}_\rho \left( \slashed{\Pi}_b^{\phantom{b}\sigma} \slashed{\Pi}_{\lambda}^{\phantom{\rho}a} (\slashed{g}^{-1})^{\lambda\delta} \left(
	\frac{1}{4}(\slashed{\nabla}_\sigma \slashed{h})_{L\delta} 
	+ \frac{1}{4}(L\slashed{h})_{\sigma\delta}
	- \frac{1}{4}(\slashed{\nabla}_\delta \slashed{h})_{L\sigma}
	\right) \right)	\Lbar_\mu \wedge \slashed{\Pi}_\nu^{\phantom{\nu}c}
	\\
	&\phantom{=}
	+ \frac{1}{2} \Lbar \left( \slashed{\Pi}_b^{\phantom{b}\rho} \slashed{\Pi}_\sigma^{\phantom{\sigma} a} (\slashed{g}^{-1})^{\sigma\lambda} \slashed{\Pi}_c^{\phantom{c}\delta} \left(
	\frac{1}{2}(\slashed{\nabla}_\rho \slashed{h})_{\delta\lambda}
	+ \frac{1}{2}(\slashed{\nabla}_\delta \slashed{h})_{\rho\lambda}
	- \frac{1}{2}(\slashed{\nabla}_\lambda \slashed{h})_{\delta\rho}
	\right) \right) L_\mu \wedge \slashed{\Pi}_\nu^{\phantom{\nu}c}
	\\
	&\phantom{=}
	+ \frac{1}{2} L \left( \slashed{\Pi}_b^{\phantom{b}\rho} \slashed{\Pi}_\sigma^{\phantom{\sigma} a} (\slashed{g}^{-1})^{\sigma\lambda} \slashed{\Pi}_c^{\phantom{c}\delta} \left(
	\frac{1}{2}(\slashed{\nabla}_\rho \slashed{h})_{\delta\lambda}
	+ \frac{1}{2}(\slashed{\nabla}_\delta \slashed{h})_{\rho\lambda}
	- \frac{1}{2}(\slashed{\nabla}_\lambda \slashed{h})_{\delta\rho}
	\right) \right) \Lbar_\mu \wedge \slashed{\Pi}_\nu^{\phantom{\nu}c}
	\\
	&\phantom{=}
	-  \slashed{\Pi}_d^{\phantom{d}\kappa}\slashed{\nabla}_\kappa \left( \slashed{\Pi}_b^{\phantom{b}\rho} \slashed{\Pi}_\sigma^{\phantom{\sigma} a} (\slashed{g}^{-1})^{\sigma\lambda} \slashed{\Pi}_c^{\phantom{c}\delta} \left(
	\frac{1}{2}(\slashed{\nabla}_\rho \slashed{h})_{\delta\lambda}
	+ \frac{1}{2}(\slashed{\nabla}_\delta \slashed{h})_{\rho\lambda}
	- \frac{1}{2}(\slashed{\nabla}_\lambda \slashed{h})_{\delta\rho}
	\right) \right) \slashed{\Pi}_\mu^{\phantom{\mu}d} \wedge \slashed{\Pi}_\nu^{\phantom{\nu}c}
	\end{split}
	\end{equation*}
	
	We do not have to expand all of the terms above, since we will eventually project both of the rectangular indices ($a$ and $b$ in the formula above) using $\slashed{\Pi}$. Hence, for example, any term proportional to $L^a$ or $\Lbar_b$ can safely be ignored.
	
	Now, using propositions \ref{proposition null connection}, \ref{proposition transport rectangular}, \ref{proposition angular rectangular}, \ref{proposition transport lbar rectangular} and \ref{proposition derivatives of projection} we find that
	\begin{equation*}
	\begin{split}
	\left(\upd L^\flat \right)_{\mu\nu}
	&=
	\frac{1}{2}\omega L_\mu \wedge \Lbar_\nu 
	- \slashed{\Pi}_c^{\phantom{c}\rho} (\slashed{\nabla}_\rho \log \mu) L_\mu \wedge \slashed{\Pi}_\nu^{\phantom{\nu}c}
	\\
	&=
	\left( -\frac{1}{2} \frac{L^i_{(\text{small})}L^i_{(\text{small})}}{r}
		- \frac{1}{4}(L h)_{L\Lbar}
		+ \frac{1}{8}(Lh)_{LL}
		+ \frac{1}{8} (\Lbar h)_{LL}
		\right) L_\mu \wedge \Lbar_\nu
		\\
	&\phantom{=}
	- \slashed{\Pi}_c^{\phantom{c}\rho} (\slashed{\nabla}_\rho \log \mu) L_\mu \wedge \slashed{\Pi}_\nu^{\phantom{\nu}c}
	\\
	\\
	\left(\upd \Lbar^\flat \right)_{\mu\nu}
	&=
	\frac{1}{2}\omega L_\mu \wedge \Lbar_\nu
	- \slashed{\Pi}_c^{\phantom{c}\rho} (\slashed{\nabla}_\rho \log \mu) L_\mu \wedge \slashed{\Pi}_\nu^{\phantom{\nu}c}
	\\
	&=
	\left( -\frac{1}{2} \frac{L^i_{(\text{small})}L^i_{(\text{small})}}{r}
		- \frac{1}{4}(L h)_{L\Lbar}
		+ \frac{1}{8}(Lh)_{LL}
		+ \frac{1}{8} (\Lbar h)_{LL}
		\right) L_\mu \wedge \Lbar_\nu
		\\
	&\phantom{=}
	- \slashed{\Pi}_c^{\phantom{c}\rho} (\slashed{\nabla}_\rho \log \mu) L_\mu \wedge \slashed{\Pi}_\nu^{\phantom{\nu}c}
	\\
	\\
	\\
	\left(\upd \slashed{\Pi}^c \right)_{\mu\nu} \slashed{\Pi}_c^{\phantom{c}\delta}
	&=
	-\frac{1}{2}(\zeta^\delta + \slashed{\nabla}^\delta \log \mu) L_\mu \wedge \Lbar_\nu
	\\
	&\phantom{=}
	- \left(\frac{1}{2} \chibar_c^{\phantom{c}\delta} 
		+ \frac{1}{4} (\slashed{g}^{-1})^{\sigma\delta} \slashed{\Pi}_c^{\phantom{c}\lambda} \left(
			(\slashed{\nabla}_\sigma \slashed{h})_{\Lbar\lambda}
			- (\slashed{\nabla}_\lambda \slashed{h})_{\Lbar\sigma}
			- (\Lbar\slashed{h})_{\sigma\lambda}
			\right)
		\right)
	L_\mu \wedge \slashed{\Pi}_\nu^{\phantom{\nu}c}
	\\
	&\phantom{=}
	- \left( \frac{1}{2} \chi_c^{\phantom{c}\delta}
		+ \frac{1}{4} (\slashed{g}^{-1})^{\sigma\delta} \slashed{\Pi}_c^{\phantom{c}\lambda} \left(
			(\slashed{\nabla}_\sigma \slashed{h})_{L\lambda}
			- (\slashed{\nabla}_\lambda \slashed{h})_{L\sigma}
			- (L\slashed{h})_{\sigma\lambda}
			\right)
		\right)
	\Lbar_\mu \wedge \slashed{\Pi}_\nu^{\phantom{\nu}c}
	\\
	&\phantom{=}
	+ \frac{1}{2} (\slashed{g}^{-1})^{\rho\delta} \slashed{\Pi}_d^{\phantom{d}\kappa}\slashed{\Pi}_c^{\phantom{c}\lambda}
	\left(
		-(\slashed{\nabla}_\lambda \slashed{h})_{\rho\kappa}
		-(\slashed{\nabla}_\kappa \slashed{h})_{\rho\lambda}
		+(\slashed{\nabla}_\rho \slashed{h})_{\lambda\kappa}
	\right)
	\slashed{\Pi}_\mu^{\phantom{\mu}c}\wedge\slashed{\Pi}_\nu^{\phantom{\nu}d}
	\\
	\end{split}
	\end{equation*}
	where, for the final lines, we have also made use of the calculation
	\begin{equation*}
	\begin{split}
	\D_\mu \slashed{\Pi}_\nu^{\phantom{\nu}c}
	&=
	\D_\mu \left( \slashed{\Pi}_\nu^{\phantom{\nu}\rho} \slashed{\Pi}_\rho^{\phantom{\rho}c} \right)
	\\
	&=
	\slashed{\Pi}_\rho^{\phantom{\rho}c} \D_\mu \slashed{\Pi}_\nu^{\phantom{\nu}\rho} + \slashed{\D}_\mu \slashed{\Pi}_\nu^{\phantom{\nu}c}
	\end{split}
	\end{equation*}
	
	Next, we note the following calculations, which will help us to expand the expressions above:
	\begin{equation*}
	\begin{split}
	\slashed{\Pi}_c^{\phantom{c}a} \slashed{\D}_L \left(\slashed{\Pi}_\mu^{\phantom{\mu}c} \right)
	&=
	\frac{1}{2} (\slashed{g}^{-1})^{\rho\nu} \slashed{\Pi}_\rho^{\phantom{\rho}a} \left( \slashed{\nabla}_\nu \slashed{h})_{L\mu} - (\slashed{\nabla}_\mu \slashed{h})_{L\nu} - (L\slashed{h})_{\mu\nu} \right) 
	\\
	\\
	\slashed{\Pi}_b^{\phantom{b}c} \slashed{\D}_L \left(\slashed{\Pi}_c^{\phantom{c}\mu} \right)
	&= \slashed{\D}_L \left( g_{bc} (\slashed{g}^{-1})^{\mu\nu} \slashed{\Pi}_\nu^{\phantom{\nu}c} \right)
	\\
	&= (\slashed{g}^{-1})^{\mu\nu} \slashed{g}_{bc} (\slashed{\D}_L \slashed{\Pi}_\nu^{\phantom{\nu}c}) + (\slashed{g}^{-1})^{\mu\nu} \slashed{\Pi}_b^{\phantom{b}\rho} (L\slashed{h})_{\nu\rho}
	\\
	&= \frac{1}{2}(\slashed{g}^{-1})^{\mu\nu} \slashed{\Pi}_b^{\phantom{b}\rho} \left( (L\slashed{h})_{\nu\rho} + (\slashed{\nabla}_\rho \slashed{h})_{L\nu} - (\slashed{\nabla}_\nu \slashed{h})_{L\rho} \right)
	\\ 
	\\
	\slashed{\Pi}_c^{\phantom{c}a} \slashed{\D}_{\Lbar} \left(\slashed{\Pi}_\mu^{\phantom{\mu}c} \right)
	&=
	\frac{1}{2} (\slashed{g}^{-1})^{\nu\rho} \slashed{\Pi}_\rho^{\phantom{\rho}a} \left( (\slashed{\nabla}_\nu \slashed{h})_{\Lbar\mu} - (\Lbar \slashed{h})_{\mu\nu} - (\slashed{\nabla}_\mu \slashed{h})_{\Lbar\nu} \right)
	\\
	\\
	\slashed{\Pi}_b^{\phantom{b}c} \slashed{\D}_{\Lbar} \left(\slashed{\Pi}_c^{\phantom{c}\mu} \right)
	&=
	(\slashed{g}^{-1})^{\mu\nu} \slashed{g}_{bc} (\slashed{\D}_{\Lbar} \slashed{\Pi}_\nu^{\phantom{\nu}c}) + (\slashed{g}^{-1})^{\mu\nu} \slashed{\Pi}_b^{\phantom{b}\rho} (\Lbar\slashed{h})_{\nu\rho}
	\\
	&=
	\frac{1}{2}(\slashed{g}^{-1})^{\mu\nu} \slashed{\Pi}_b^{\phantom{b}\rho} \left( (\Lbar \slashed{h})_{\nu\rho} + (\slashed{\nabla}_\rho \slashed{h})_{\Lbar\nu} - (\slashed{\nabla}_\nu \slashed{h})_{\Lbar\rho} \right)
	\\
	\\
	\slashed{\Pi}_c^{\phantom{c}a} \slashed{\nabla}_{\rho} \left(\slashed{\Pi}_\mu^{\phantom{\mu}c} \right)
	&=
	\frac{1}{2} (\slashed{g}^{-1})^{\sigma\lambda} \slashed{\Pi}_\lambda^{\phantom{\lambda}a}
	\left( (\slashed{\nabla}_\sigma \slashed{h})_{\rho\mu} - (\slashed{\nabla}_\rho \slashed{h})_{\mu\sigma} - (\slashed{\nabla}_\mu \slashed{h})_{\rho\sigma} \right)
	\\
	\\
	\slashed{\Pi}_b^{\phantom{b}c} \slashed{\nabla}_{\rho} \left(\slashed{\Pi}_c^{\phantom{c}\mu} \right)
	&=
	(\slashed{g}^{-1})^{\mu\nu} \slashed{g}_{bc} (\slashed{\nabla}_{\rho} \slashed{\Pi}_\nu^{\phantom{\nu}c}) + (\slashed{g}^{-1})^{\mu\nu} \slashed{\Pi}_b^{\phantom{b}\sigma} (\slashed{\nabla}_\rho \slashed{h})_{\nu\sigma}
	\\
	&=
	\frac{1}{2} (\slashed{g}^{-1})^{\mu\nu} \slashed{\Pi}_b^{\phantom{b}\sigma}
	\left( (\slashed{\nabla}_\rho \slashed{h})_{\nu\sigma} + (\slashed{\nabla}_\sigma \slashed{h})_{\nu\rho} - (\slashed{\nabla}_\nu \slashed{h})_{\rho\sigma} \right)
	\end{split}
	\end{equation*}
	
	Finally, we will need to deal with terms involving second derivatives. These also need to be rewritten in a convenient way, to involve derivatives of $\mathscr{Z}h$ for commutation operators $\mathscr{Z}$, or to make use of the fact that $\tilde{\Box}_g h_{ab} = F_{ab}$. We have
	\begin{equation*}
	\begin{split}
	\slashed{\D}_L \left( (\slashed{\nabla}_\rho \slashed{h})_{\Lbar\lambda} \right)
	&=
	r^{-1} \left(\slashed{\D}_L \left( r\slashed{\nabla}_\rho \slashed{h} \right) \right)_{\Lbar\lambda}
	- r^{-1} (\slashed{\nabla}_\rho \slashed{h})_{\Lbar\lambda}
	+ \left(\frac{L^i L^i - 1}{r}\right) (\slashed{\nabla}_\rho \slashed{h})_{\Lbar\lambda}
	\\
	&\phantom{=}
	- 2(\slashed{g}^{-1})^{\sigma\kappa} \left(L^i_{(\text{small})} \frac{\slashed{\nabla}_\kappa x^i}{r} \right) (\slashed{\nabla}_\rho \slashed{h})_{\sigma\lambda}
	- \left(L^i_{(\text{small})} \frac{\slashed{\nabla}_\lambda x^i}{r} \right) (\slashed{\nabla}_\rho h)_{L\Lbar}
	\\
	&\phantom{=}
	+ \frac{1}{4}(\slashed{\nabla}_\rho \slashed{h})_{L\lambda} (Lh)_{\Lbar\Lbar}
	+ \frac{1}{2}(\slashed{\nabla}_\rho \slashed{h})_{\Lbar\lambda} (Lh)_{L\Lbar}
	- \frac{1}{4}(\slashed{\nabla}_\rho \slashed{h})_{\Lbar\lambda} (Lh)_{LL}
	\\
	&\phantom{=}
	+ 2(\slashed{g}^{-1})^{\sigma\kappa} (\slashed{\nabla}_\rho \slashed{h})_{\sigma\lambda} (L\slashed{h})_{L\kappa}
	- (\slashed{g}^{-1})^{\sigma\kappa} (\slashed{\nabla}_\rho \slashed{h})_{\sigma\lambda} (\slashed{\nabla}_\kappa h)_{LL}
	+ \frac{1}{4} (\slashed{\nabla}_\rho h)_{L\Lbar} (\slashed{\nabla}_\lambda h)_{LL}
	\\
	&\phantom{=}
	+ \frac{1}{4} (\slashed{\nabla}_\rho h)_{\Lbar\Lbar} (\slashed{\nabla}_\lambda h)_{LL}
	+ \frac{1}{2} (\slashed{g}^{-1})^{\sigma\kappa} (\slashed{\nabla}_\rho \slashed{h})_{\Lbar\sigma} (\slashed{\nabla}_\kappa \slashed{h})_{L\lambda}
	- \frac{1}{2} (\slashed{g}^{-1})^{\sigma\kappa} (\slashed{\nabla}_\rho \slashed{h})_{\Lbar\sigma} (\slashed{\nabla}_\lambda \slashed{h})_{L\kappa}
	\\
	&\phantom{=}
	- \frac{1}{2} (\slashed{g}^{-1})^{\sigma\kappa} (\slashed{\nabla}_\rho \slashed{h})_{\Lbar\sigma} (L \slashed{h})_{\kappa\lambda}
	\end{split}
	\end{equation*}
	
	\begin{equation*}
	\begin{split}
	\slashed{\D}_L \left( \left( \slashed{\nabla}_\rho \slashed{h} \right)_{\delta\lambda} \right)
	&=
	r^{-1} \left(\slashed{\D}_L \left(r\slashed{\nabla}_\rho \slashed{h} \right)\right)_{\delta\lambda}
	- r^{-1} \left(\slashed{\nabla}_\rho \slashed{h} \right)_{\delta\lambda}
	- \left(L^i_{(\text{small})} \frac{ \slashed{\nabla}_\delta x^i}{r} \right) (\slashed{\nabla}_\rho \slashed{h})_{L\lambda}
	\\
	&\phantom{=}
	- \left(L^i_{(\text{small})} \frac{ \slashed{\nabla}_\lambda x^i}{r} \right) (\slashed{\nabla}_\rho \slashed{h})_{L\delta}
	+ \frac{1}{4}(\slashed{\nabla}_\rho \slashed{h})_{L\delta} (\slashed{\nabla}_\lambda h)_{LL}
	+ \frac{1}{4}(\slashed{\nabla}_\rho \slashed{h})_{\Lbar\delta} (\slashed{\nabla}_\lambda h)_{LL}
	\\
	&\phantom{=}
	+ (\slashed{g}^{-1})^{\sigma\kappa} (\slashed{\nabla}_\rho \slashed{h})_{\delta\sigma} (\slashed{\nabla}_\kappa \slashed{h})_{L\lambda}
	+ (\slashed{g}^{-1})^{\sigma\kappa} (\slashed{\nabla}_\rho \slashed{h})_{\delta\sigma} (\slashed{\nabla}_\lambda \slashed{h})_{L\kappa}
	- (\slashed{g}^{-1})^{\sigma\kappa} (\slashed{\nabla}_\rho \slashed{h})_{\delta\sigma} (L \slashed{h})_{\lambda\kappa}
	\\
	&\phantom{=}
	+ \frac{1}{4}(\slashed{\nabla}_\rho \slashed{h})_{L\lambda} (\slashed{\nabla}_\delta h)_{LL}
	+ \frac{1}{4}(\slashed{\nabla}_\rho \slashed{h})_{\Lbar\lambda} (\slashed{\nabla}_\delta h)_{LL}
	+ (\slashed{g}^{-1})^{\sigma\kappa} (\slashed{\nabla}_\rho \slashed{h})_{\lambda\sigma} (\slashed{\nabla}_\kappa \slashed{h})_{L\delta}
	\\
	&\phantom{=}
	+ (\slashed{g}^{-1})^{\sigma\kappa} (\slashed{\nabla}_\rho \slashed{h})_{\lambda\sigma} (\slashed{\nabla}_\delta \slashed{h})_{L\kappa}
	- (\slashed{g}^{-1})^{\sigma\kappa} (\slashed{\nabla}_\rho \slashed{h})_{\lambda\sigma} (L \slashed{h})_{\delta\kappa}
	\end{split}
	\end{equation*}
	
	\begin{equation*}
	\begin{split}
	\slashed{\D}_{\Lbar} \left( \left( \slashed{\nabla}_\rho \slashed{h} \right)_{L\lambda} \right)
	&=
	r^{-1} \left( \slashed{\D}_{\Lbar} \left( r\slashed{\nabla}_\rho \slashed{h} \right) \right)_{L\lambda}
	+ r^{-1} \left(\slashed{\nabla}_\rho \slashed{h} \right)_{L\lambda}
	- \left( \frac{L^i_{(\text{small})}L^i_{(\text{small})}}{r} \right) (\slashed{\nabla}_\rho \slashed{h})_{\Lbar\lambda}
	\\
	&\phantom{=}
	+ 2 (\slashed{g}^{-1})^{\sigma\kappa} \left( L^i_{(\text{small})} \frac{\slashed{\nabla}_\kappa x^i}{r} \right) (\slashed{\nabla}_\rho \slashed{h})_{\sigma\lambda}
	+ 2 (\slashed{g}^{-1})^{\sigma\kappa} (\slashed{\nabla}_\kappa \log \mu) (\slashed{\nabla}_\rho \slashed{h})_{\sigma\lambda}
	\\
	&\phantom{=}
	+ (\slashed{\nabla}_\lambda \log \mu)(\slashed{\nabla}_\rho h)_{LL}
	+ \left( L^i_{(\text{small})} \frac{\slashed{\nabla}_\lambda x^i}{r} \right) (\slashed{\nabla}_\rho h)_{L\Lbar}
	+ (\slashed{\nabla}_\lambda \log \mu)(\slashed{\nabla}_\rho h)_{L\Lbar}
	\\
	&\phantom{=}
	+ \frac{1}{4}(\slashed{\nabla}_\rho \slashed{h})_{L\lambda} (Lh)_{\Lbar\Lbar}
	- \frac{1}{2}(\slashed{\nabla}_\rho \slashed{h})_{L\lambda} (Lh)_{L\Lbar}
	+ \frac{1}{4}(\slashed{\nabla}_\rho \slashed{h})_{L\lambda} (Lh)_{LL}
	+ \frac{1}{4}(\slashed{\nabla}_\rho \slashed{h})_{L\lambda} (\Lbar h)_{LL}
	\\
	&\phantom{=}
	+ \frac{1}{4}(\slashed{\nabla}_\rho \slashed{h})_{\Lbar\lambda} (\Lbar h)_{LL}
	- (\slashed{g}^{-1})^{\sigma\kappa} (\slashed{\nabla}_\rho \slashed{h})_{\sigma\lambda} (L\slashed{h})_{L\kappa}
	- \frac{1}{2}(\slashed{g}^{-1})^{\sigma\kappa} (\slashed{\nabla}_\rho \slashed{h})_{\sigma\lambda} (\slashed{\nabla}_\kappa h)_{LL}
	\\
	&\phantom{=}
	+ (\slashed{g}^{-1})^{\sigma\kappa} (\slashed{\nabla}_\rho \slashed{h})_{\sigma\lambda} (\slashed{\nabla}_\kappa h)_{L\Lbar}
	- (\slashed{g}^{-1})^{\sigma\kappa} (\slashed{\nabla}_\rho \slashed{h})_{\sigma\lambda} (\Lbar\slashed{h})_{L\kappa}
	+ (\slashed{g}^{-1})^{\sigma\kappa} (\slashed{\nabla}_\rho \slashed{h})_{\sigma\lambda} (L\slashed{h})_{\Lbar\kappa}
	\\
	&\phantom{=}
	+ \frac{1}{4}(\slashed{\nabla}_\rho h)_{LL} (\slashed{\nabla}_\lambda h)_{\Lbar\Lbar}
	- \frac{1}{4}(\slashed{\nabla}_\rho h)_{L\Lbar} (\slashed{\nabla}_\lambda h)_{LL}
	+ \frac{1}{2}(\slashed{\nabla}_\rho h)_{L\Lbar} (\slashed{\nabla}_\lambda h)_{L\Lbar}
	\\
	&\phantom{=}
	-\frac{1}{2} (\slashed{g}^{-1})^{\sigma\kappa} (\slashed{\nabla}_\rho \slashed{h})_{L\sigma} (\Lbar \slashed{h})_{\lambda\kappa}
	-\frac{1}{2} (\slashed{g}^{-1})^{\sigma\kappa} (\slashed{\nabla}_\rho \slashed{h})_{L\sigma} (\slashed{\nabla}_\lambda \slashed{h})_{\Lbar\kappa}
	+\frac{1}{2} (\slashed{g}^{-1})^{\sigma\kappa} (\slashed{\nabla}_\rho \slashed{h})_{L\sigma} (\slashed{\nabla}_\kappa \slashed{h})_{\Lbar\lambda}
	\end{split}
	\end{equation*}

	\begin{equation*}
	\begin{split}
	\slashed{\D}_{\Lbar} \left( (\slashed{\nabla}_\rho \slashed{h})_{\delta\lambda} \right)
	&=
	r^{-1} \left( \slashed{\D}_{\Lbar} \left( r\slashed{\nabla}_\rho \slashed{h} \right) \right)_{\delta\lambda}
	+ r^{-1} (\slashed{\nabla}_\rho \slashed{h})_{\delta\lambda}
	+ (\slashed{\nabla}_\lambda \log \mu)(\slashed{\nabla}_\rho \slashed{h})_{L\delta}
	+ (\slashed{\nabla}_\delta \log \mu)(\slashed{\nabla}_\rho \slashed{h})_{L\lambda}
	\\
	&\phantom{=}
	+ \left( L^i_{(\text{small})} \frac{ \slashed{\nabla}_\lambda x^i}{r} \right) (\slashed{\nabla}_\rho \slashed{h})_{\Lbar\delta}
	+ \left( L^i_{(\text{small})} \frac{ \slashed{\nabla}_\delta x^i}{r} \right) (\slashed{\nabla}_\rho \slashed{h})_{\Lbar\lambda}
	+ (\slashed{\nabla}_\lambda \log \mu)(\slashed{\nabla}_\rho \slashed{h})_{\Lbar\delta}
	\\
	&\phantom{=}
	+ (\slashed{\nabla}_\delta \log \mu)(\slashed{\nabla}_\rho \slashed{h})_{\Lbar\lambda}
	+ \frac{1}{4}(\slashed{\nabla}_\rho \slashed{h})_{L\delta} (\slashed{\nabla}_\lambda h)_{\Lbar\Lbar}
	- \frac{1}{4}(\slashed{\nabla}_\rho \slashed{h})_{\Lbar\delta} (\slashed{\nabla}_\lambda h)_{LL}
	+ \frac{1}{2}(\slashed{\nabla}_\rho \slashed{h})_{\Lbar\delta} (\slashed{\nabla}_\lambda h)_{L\Lbar}
	\\
	&\phantom{=}
	+ \frac{1}{2}(\slashed{g}^{-1})^{\sigma\kappa} (\slashed{\nabla}_\rho \slashed{h})_{\delta\sigma} (\Lbar\slashed{h})_{\lambda\kappa}
	+ \frac{1}{2}(\slashed{g}^{-1})^{\sigma\kappa} (\slashed{\nabla}_\rho \slashed{h})_{\delta\sigma} (\slashed{\nabla}_\lambda\slashed{h})_{\Lbar\kappa}
	- \frac{1}{2}(\slashed{g}^{-1})^{\sigma\kappa} (\slashed{\nabla}_\rho \slashed{h})_{\delta\sigma} (\slashed{\nabla}_\kappa\slashed{h})_{\Lbar\lambda}
	\\
	&\phantom{=}
	+ \frac{1}{4}(\slashed{\nabla}_\rho \slashed{h})_{L\lambda} (\slashed{\nabla}_\delta h)_{\Lbar\Lbar}
	- \frac{1}{4}(\slashed{\nabla}_\rho \slashed{h})_{\Lbar\lambda} (\slashed{\nabla}_\delta h)_{LL}
	+ \frac{1}{2}(\slashed{\nabla}_\rho \slashed{h})_{\Lbar\lambda} (\slashed{\nabla}_\delta h)_{L\Lbar}
	\\
	&\phantom{=}
	+ \frac{1}{2}(\slashed{g}^{-1})^{\sigma\kappa} (\slashed{\nabla}_\rho \slashed{h})_{\lambda\sigma} (\Lbar\slashed{h})_{\delta\kappa}
	+ \frac{1}{2}(\slashed{g}^{-1})^{\sigma\kappa} (\slashed{\nabla}_\rho \slashed{h})_{\lambda\sigma} (\slashed{\nabla}_\delta\slashed{h})_{\Lbar\kappa}
	- \frac{1}{2}(\slashed{g}^{-1})^{\sigma\kappa} (\slashed{\nabla}_\rho \slashed{h})_{\lambda\sigma} (\slashed{\nabla}_\kappa\slashed{h})_{\Lbar\delta}
	\end{split}
	\end{equation*}

	\begin{equation*}
	\begin{split}
	\slashed{\nabla}_\rho \left( (\slashed{\nabla}_\sigma \slashed{h})_{L\delta} \right)
	&=
	r^{-1} \left( \slashed{\nabla}_\rho ( r\slashed{\nabla}_\sigma \slashed{h}) \right)_{L\delta}
	+ \frac{1}{2} \left(L^i_{(\text{small})} \frac{\slashed{\nabla}_\rho x^i}{r} \right) (\slashed{\nabla}_\sigma \slashed{h})_{L\slashed{\delta}}
	+ \chi_\rho^{\phantom{\rho}\kappa} (\slashed{\nabla}_\sigma \slashed{h})_{\delta\kappa}
	+ \frac{1}{2} \chi_{\rho\delta} (\slashed{\nabla}_\sigma h)_{LL}
	\\
	&\phantom{=}
	- \left( \frac{1}{r} \slashed{\Pi}_\rho^{\phantom{\rho}i} \slashed{\Pi}_\sigma^{\phantom{\sigma}i} \right) (\slashed{\nabla}_\sigma h)_{LL} 
	+ \frac{1}{2} \chi_{\rho\delta} (\slashed{\nabla}_\sigma h)_{L\Lbar}
	+ \frac{1}{4}(\slashed{\nabla}_\sigma \slashed{h})_{L\delta} (\slashed{\nabla}_\rho h)_{LL}
	+ \frac{1}{4}(\slashed{\nabla}_\sigma \slashed{h})_{\Lbar\delta} (\slashed{\nabla}_\rho h)_{LL}
	\\
	&\phantom{=}
	- \frac{1}{2}(\slashed{g}^{-1})^{\kappa\lambda} (\slashed{\nabla}_\sigma \slashed{h})_{\delta\kappa} (\slashed{\nabla}_\rho \slashed{h})_{L\lambda}
	- \frac{1}{2}(\slashed{g}^{-1})^{\kappa\lambda} (\slashed{\nabla}_\sigma \slashed{h})_{\delta\kappa} (L\slashed{h})_{\rho\lambda}
	+ \frac{1}{2}(\slashed{g}^{-1})^{\kappa\lambda} (\slashed{\nabla}_\sigma \slashed{h})_{\delta\kappa} (\slashed{\nabla}_\lambda \slashed{h})_{L\rho}
	\\
	&\phantom{=}
	+ \frac{1}{4}(\slashed{\nabla}_\sigma h)_{LL} (\slashed{\nabla}_\mu \slashed{h})_{L\nu}
	+ \frac{1}{4}(\slashed{\nabla}_\sigma h)_{LL} (\slashed{\nabla}_\nu \slashed{h})_{L\mu}
	- \frac{1}{4}(\slashed{\nabla}_\sigma h)_{L\Lbar} (\slashed{\nabla}_\rho \slashed{h})_{L\delta}
	- \frac{1}{4}(\slashed{\nabla}_\sigma h)_{L\Lbar} (\slashed{\nabla}_\delta \slashed{h})_{L\rho}
	\\
	&\phantom{=}
	+ \frac{1}{4}(\slashed{\nabla}_\sigma h)_{L\Lbar} (L \slashed{h})_{\rho\delta}
	- \frac{1}{2}(\slashed{g}^{-1})^{\kappa\lambda} (\slashed{\nabla}_\sigma \slashed{h})_{L\kappa} (\slashed{\nabla}_\rho \slashed{h})_{\delta\lambda}
	- \frac{1}{2}(\slashed{g}^{-1})^{\kappa\lambda} (\slashed{\nabla}_\sigma \slashed{h})_{L\kappa} (\slashed{\nabla}_\delta \slashed{h})_{\rho\lambda}
	\\
	&\phantom{=}
	+ \frac{1}{2}(\slashed{g}^{-1})^{\kappa\lambda} (\slashed{\nabla}_\sigma \slashed{h})_{L\kappa} (\slashed{\nabla}_\lambda \slashed{h})_{\rho\delta}
	\end{split}
	\end{equation*}

	\begin{equation*}
	\begin{split}
	\slashed{\nabla}_\rho \left( (\slashed{\nabla}_\sigma \slashed{h})_{\Lbar\delta} \right)
	&=
	r^{-1} \left( \slashed{\nabla}_\rho ( r\slashed{\nabla}_\sigma \slashed{h}) \right)_{\Lbar\delta}
	+ \left( \Lbar^i_{(\text{small})} \frac{\slashed{\nabla}_\rho x^i}{r}\right)(\slashed{\nabla}_\sigma \slashed{h})_{\Lbar\delta}
	+ \left(\frac{1}{r} \slashed{\Pi}_\rho^{\phantom{\rho}i} \slashed{\Pi}_\delta^{\phantom{\delta}i} \right) (\slashed{\nabla}_\sigma h)_{\Lbar\Lbar}
	+ \chibar_\rho^{\phantom{\rho}\kappa} (\slashed{\nabla}_\sigma \slashed{h})_{\delta\kappa}
	\\
	&\phantom{=}
	+ \frac{1}{2} \chibar_{\rho\delta} (\slashed{\nabla}_\sigma h)_{\Lbar\Lbar}
	+ \frac{1}{4}(\slashed{\nabla}_\sigma \slashed{h})_{L\delta} (\slashed{\nabla}_\rho \slashed{h})_{LL}
	+ \frac{1}{4}(\slashed{\nabla}_\sigma \slashed{h})_{\Lbar\delta}(\slashed{\nabla}_\rho \slashed{h})_{\Lbar\Lbar}
	- \frac{1}{2}(\slashed{g}^{-1})^{\kappa\lambda}(\slashed{\nabla}_\sigma \slashed{h})_{\delta\kappa}(\Lbar \slashed{h})_{\rho\lambda}
	\\
	&\phantom{=}
	- \frac{1}{2}(\slashed{g}^{-1})^{\kappa\lambda}(\slashed{\nabla}_\sigma \slashed{h})_{\delta\kappa}(\slashed{\nabla}_\rho \slashed{h})_{\Lbar\lambda}
	+ \frac{1}{2}(\slashed{g}^{-1})^{\kappa\lambda}(\slashed{\nabla}_\sigma \slashed{h})_{\delta\kappa}(\slashed{\nabla}_\lambda \slashed{h})_{\Lbar\rho}
	+ \frac{1}{4}(\slashed{\nabla}_\sigma h)_{L\Lbar} (\slashed{\nabla}_\rho \slashed{h})_{\Lbar\delta}
	\\
	&\phantom{=}
	+ \frac{1}{4}(\slashed{\nabla}_\sigma h)_{L\Lbar} (\slashed{\nabla}_\delta \slashed{h})_{\Lbar\rho}
	- \frac{1}{4}(\slashed{\nabla}_\sigma h)_{L\Lbar} (\Lbar \slashed{h})_{\rho\delta}
	- \frac{1}{4}(\slashed{\nabla}_\sigma h)_{\Lbar\Lbar} (\Lbar \slashed{h})_{\rho\delta}
	+ \frac{1}{4}(\slashed{\nabla}_\sigma h)_{\Lbar\Lbar} (\slashed{\nabla}_\rho \slashed{h})_{\Lbar\delta}
	\\
	&\phantom{=}
	+ \frac{1}{4}(\slashed{\nabla}_\sigma h)_{\Lbar\Lbar} (\slashed{\nabla}_\delta \slashed{h})_{\Lbar\rho}
	- \frac{1}{2}(\slashed{g}^{-1})^{\kappa\lambda} (\slashed{\nabla}_\sigma \slashed{h})_{\Lbar\kappa} (\slashed{\nabla}_\rho \slashed{h})_{\lambda\delta}
	- \frac{1}{2}(\slashed{g}^{-1})^{\kappa\lambda} (\slashed{\nabla}_\sigma \slashed{h})_{\Lbar\kappa} (\slashed{\nabla}_\delta \slashed{h})_{\lambda\rho}
	\\
	&\phantom{=}
	+ \frac{1}{2}(\slashed{g}^{-1})^{\kappa\lambda} (\slashed{\nabla}_\sigma \slashed{h})_{\Lbar\kappa} (\slashed{\nabla}_\lambda \slashed{h})_{\rho\delta}
	\end{split}
	\end{equation*}
	
	\begin{equation*}
	\begin{split}
	\slashed{\nabla}_\kappa \left( (\slashed{\nabla}_\rho \slashed{h})_{\delta\lambda} \right)
	&=
	r^{-1} \left( \slashed{\nabla}_\rho ( r\slashed{\nabla}_\sigma \slashed{h}) \right)_{\Lbar\delta}
	+ \frac{1}{2} \chi_{\lambda\kappa} (\slashed{\nabla}_\rho \slashed{h})_{L\delta}
	+ \frac{1}{2} \chi_{\delta\kappa} (\slashed{\nabla}_\rho \slashed{h})_{L\lambda}
	+ \frac{1}{2} \chi_{\lambda\kappa} (\slashed{\nabla}_\rho \slashed{h})_{\Lbar\delta}
	+ \frac{1}{2} \chi_{\delta\kappa} (\slashed{\nabla}_\rho \slashed{h})_{\Lbar\lambda}
	\\
	&\phantom{=}
	- \left( \frac{1}{r}\slashed{\Pi}_\kappa^{\phantom{\kappa}i} \slashed{\Pi}_\lambda^{\phantom{\lambda}i} \right) (\slashed{\nabla}_\rho \slashed{h})_{L\delta}
	- \left( \frac{1}{r}\slashed{\Pi}_\kappa^{\phantom{\kappa}i} \slashed{\Pi}_\delta^{\phantom{\delta}i} \right) (\slashed{\nabla}_\rho \slashed{h})_{L\lambda}
	+ \frac{1}{4}(\slashed{\nabla}_\rho \slashed{h})_{L\delta} (\slashed{\nabla}_\kappa \slashed{h})_{L\lambda}
	\\
	&\phantom{=}
	+ \frac{1}{4}(\slashed{\nabla}_\rho \slashed{h})_{L\delta} (\slashed{\nabla}_\lambda \slashed{h})_{L\kappa}
	- \frac{1}{4}(\slashed{\nabla}_\rho \slashed{h})_{\Lbar\delta} (\slashed{\nabla}_\kappa \slashed{h})_{L\lambda}
	- \frac{1}{4}(\slashed{\nabla}_\rho \slashed{h})_{\Lbar\delta} (\slashed{\nabla}_\lambda \slashed{h})_{L\kappa}
	+ \frac{1}{4}(\slashed{\nabla}_\rho \slashed{h})_{\Lbar\delta} (L \slashed{h})_{\lambda\kappa}
	\\
	&\phantom{=}
	- \frac{1}{2}(\slashed{g}^{-1})^{\zeta\xi} (\slashed{\nabla}_\rho \slashed{h})_{\delta\zeta} (\slashed{\nabla}_\kappa \slashed{h})_{\lambda\xi}
	- \frac{1}{2}(\slashed{g}^{-1})^{\zeta\xi} (\slashed{\nabla}_\rho \slashed{h})_{\delta\zeta} (\slashed{\nabla}_\lambda \slashed{h})_{\kappa\xi}
	+ \frac{1}{2}(\slashed{g}^{-1})^{\zeta\xi} (\slashed{\nabla}_\rho \slashed{h})_{\delta\zeta} (\slashed{\nabla}_\xi \slashed{h})_{\lambda\kappa}
	\\
	&\phantom{=}
	+ \frac{1}{4}(\slashed{\nabla}_\rho \slashed{h})_{L\lambda} (\slashed{\nabla}_\kappa \slashed{h})_{L\delta}
	+ \frac{1}{4}(\slashed{\nabla}_\rho \slashed{h})_{L\lambda} (\slashed{\nabla}_\lambda \slashed{h})_{L\kappa}
	- \frac{1}{4}(\slashed{\nabla}_\rho \slashed{h})_{\Lbar\lambda} (\slashed{\nabla}_\kappa \slashed{h})_{L\delta}
	- \frac{1}{4}(\slashed{\nabla}_\rho \slashed{h})_{\Lbar\lambda} (\slashed{\nabla}_\delta \slashed{h})_{L\kappa}
	\\
	&\phantom{=}
	+ \frac{1}{4}(\slashed{\nabla}_\rho \slashed{h})_{\Lbar\lambda} (L \slashed{h})_{\delta\kappa}
	- \frac{1}{2}(\slashed{g}^{-1})^{\zeta\xi} (\slashed{\nabla}_\rho \slashed{h})_{\lambda\zeta} (\slashed{\nabla}_\kappa \slashed{h})_{\delta\xi}
	- \frac{1}{2}(\slashed{g}^{-1})^{\zeta\xi} (\slashed{\nabla}_\rho \slashed{h})_{\lambda\zeta} (\slashed{\nabla}_\delta \slashed{h})_{\kappa\xi}
	\\
	&\phantom{=}
	+ \frac{1}{2}(\slashed{g}^{-1})^{\zeta\xi} (\slashed{\nabla}_\rho \slashed{h})_{\lambda\zeta} (\slashed{\nabla}_\xi \slashed{h})_{\delta\kappa}
	\end{split}
	\end{equation*}

	\begin{equation*}
	\begin{split}
	\slashed{\D}_L \left((\Lbar \slashed{h})_{\rho\lambda} \right)
	- \slashed{\D}_{\Lbar} \left((L \slashed{h})_{\rho\lambda} \right)
	&=
	-\left( L^i_{(\text{small})} \frac{\slashed{\nabla}_\rho x^i}{r} \right) (\Lbar\slashed{h})_{L\lambda}
	-\left( L^i_{(\text{small})} \frac{\slashed{\nabla}_\lambda x^i}{r} \right) (\Lbar\slashed{h})_{L\rho}
	\\
	&\phantom{=}
	-\left( L^i_{(\text{small})} \frac{\slashed{\nabla}_\rho x^i}{r} \right) (L\slashed{h})_{\Lbar\lambda}
	-\left( L^i_{(\text{small})} \frac{\slashed{\nabla}_\lambda x^i}{r} \right) (L\slashed{h})_{\Lbar\rho}
	- (\slashed{\nabla}_\rho \log \mu) (L\slashed{h})_{L\lambda}
	\\
	&\phantom{=}
	- (\slashed{\nabla}_\lambda \log \mu) (L\slashed{h})_{L\rho}
	- (\slashed{\nabla}_\rho \log \mu) (L\slashed{h})_{\Lbar\lambda}
	- (\slashed{\nabla}_\lambda \log \mu) (L\slashed{h})_{\Lbar\rho}
	\\
	&\phantom{=}
	+ \frac{1}{4}(\slashed{\nabla}_\rho \slashed{h})_{LL} (\Lbar \slashed{h})_{L\lambda}
	+ \frac{1}{4}(\Lbar\slashed{h})_{\Lbar\lambda} (\slashed{\nabla}_\rho h)_{LL}
	+ \frac{1}{2}(\slashed{g}^{-1})^{\sigma\kappa} (\Lbar\slashed{h})_{\lambda\sigma} (\slashed{\nabla}_\kappa \slashed{h})_{L\rho}
	\\
	&\phantom{=}
	- \frac{1}{2}(\slashed{g}^{-1})^{\sigma\kappa} (\Lbar\slashed{h})_{\lambda\sigma} (\slashed{\nabla}_\rho \slashed{h})_{L\kappa}
	+ \frac{1}{4}(\slashed{\nabla}_\lambda \slashed{h})_{LL} (\Lbar \slashed{h})_{L\rho}
	+ \frac{1}{4}(\Lbar\slashed{h})_{\Lbar\rho} (\slashed{\nabla}_\lambda h)_{LL}
	\\
	&\phantom{=}
	+ \frac{1}{2}(\slashed{g}^{-1})^{\sigma\kappa} (\Lbar\slashed{h})_{\rho\sigma} (\slashed{\nabla}_\kappa \slashed{h})_{L\lambda}
	- \frac{1}{2}(\slashed{g}^{-1})^{\sigma\kappa} (\Lbar\slashed{h})_{\rho\sigma} (\slashed{\nabla}_\lambda \slashed{h})_{L\kappa}
	\\
	&\phantom{=}
	- \frac{1}{4}(L\slashed{h})_{L\lambda} (\slashed{\nabla}_\rho h)_{\Lbar\Lbar}
	+ \frac{1}{4}(L\slashed{h})_{\Lbar\lambda} (\slashed{\nabla}_\rho h)_{LL}
	- \frac{1}{2}(L\slashed{h})_{\Lbar\lambda} (\slashed{\nabla}_\rho h)_{L\Lbar}
	\\
	&\phantom{=}
	- \frac{1}{2}(\slashed{g}^{-1})^{\sigma\kappa} (L\slashed{h})_{\lambda\sigma} (\slashed{\nabla}_\kappa \slashed{h})_{\Lbar\rho}
	+ \frac{1}{2}(\slashed{g}^{-1})^{\sigma\kappa} (L\slashed{h})_{\lambda\sigma} (\slashed{\nabla}_\rho \slashed{h})_{\Lbar\kappa}
	\\
	&\phantom{=}
	- \frac{1}{4}(L\slashed{h})_{L\rho} (\slashed{\nabla}_\lambda h)_{\Lbar\Lbar}
	+ \frac{1}{4}(L\slashed{h})_{\Lbar\rho} (\slashed{\nabla}_\lambda h)_{LL}
	- \frac{1}{2}(L\slashed{h})_{\Lbar\rho} (\slashed{\nabla}_\lambda h)_{L\Lbar}
	\\
	&\phantom{=}
	- \frac{1}{2}(\slashed{g}^{-1})^{\sigma\kappa} (L\slashed{h})_{\rho\sigma} (\slashed{\nabla}_\kappa \slashed{h})_{\Lbar\lambda}
	+ \frac{1}{2}(\slashed{g}^{-1})^{\sigma\kappa} (L\slashed{h})_{\rho\sigma} (\slashed{\nabla}_\lambda \slashed{h})_{\Lbar\kappa}
	\end{split}
	\end{equation*}
	
	\begin{equation*}
	\begin{split}
	\slashed{\D}_L \left( (\slashed{\nabla}_\rho \slashed{h})_{\sigma\delta} \right)
	- \slashed{\nabla}_\rho \left( (L\slashed{h})_{\sigma\delta} \right)
	&=
	-\left( L^i \frac{\slashed{\nabla}_\delta x^i}{r} \right) (\slashed{\nabla}_\rho \slashed{h})_{L\sigma}
	-\left( L^i \frac{\slashed{\nabla}_\sigma x^i}{r} \right) (\slashed{\nabla}_\rho \slashed{h})_{L\delta}
	-\frac{1}{2}\chi_{\delta\rho} (L\slashed{h})_{L\sigma}
	\\
	&\phantom{=}
	-\frac{1}{2}\chi_{\sigma\rho} (L\slashed{h})_{L\delta}
	-\frac{1}{2}\chi_{\delta\rho} (L\slashed{h})_{\Lbar\sigma}
	-\frac{1}{2}\chi_{\sigma\rho} (L\slashed{h})_{\Lbar\delta}
	- \chi_{\rho}^{\phantom{\rho}\lambda} (\slashed{\nabla}_\lambda \slashed{h})_{\sigma\delta}
	\\
	&\phantom{=}
	+ \left( \frac{1}{r} \slashed{\Pi}_\rho^{\phantom{\rho}i} \slashed{\Pi}_\delta^{\phantom{\rho}i} \right) (L\slashed{h})_{L\sigma}
	+ \left( \frac{1}{r} \slashed{\Pi}_\rho^{\phantom{\rho}i} \slashed{\Pi}_\sigma^{\phantom{\rho}i} \right) (L\slashed{h})_{L\delta}
	+ \frac{1}{4}(\slashed{\nabla}_\rho \slashed{h})_{L\sigma} (\slashed{\nabla}_\delta h)_{LL}
	\\
	&\phantom{=}
	+ \frac{1}{4}(\slashed{\nabla}_\rho \slashed{h})_{\Lbar\sigma} (\slashed{\nabla}_\delta h)_{LL}
	+ \frac{1}{2}(\slashed{g}^{-1})^{\lambda\kappa} (\slashed{\nabla}_\rho \slashed{h})_{\sigma\lambda} (\slashed{\nabla}_\kappa \slashed{h})_{L\delta}
	\\
	&\phantom{=}
	- \frac{1}{2}(\slashed{g}^{-1})^{\lambda\kappa} (\slashed{\nabla}_\rho \slashed{h})_{\sigma\lambda} (\slashed{\nabla}_\delta \slashed{h})_{L\kappa}
	- \frac{1}{4}(L\slashed{h})_{L\sigma} (\slashed{\nabla}_\rho \slashed{h})_{L\delta}
	- \frac{1}{4}(L\slashed{h})_{L\sigma} (\slashed{\nabla}_\delta \slashed{h})_{L\rho}
	\\
	&\phantom{=}
	- \frac{1}{4}(L\slashed{h})_{\Lbar\sigma} (\slashed{\nabla}_\rho \slashed{h})_{L\delta}
	- \frac{1}{4}(L\slashed{h})_{\Lbar\sigma} (\slashed{\nabla}_\delta \slashed{h})_{L\rho}
	+ \frac{1}{4}(L\slashed{h})_{\Lbar\sigma} (L\slashed{h})_{\delta\rho}
	\\
	&\phantom{=}
	+ \frac{1}{2}(\slashed{g}^{-1})^{\lambda\kappa} (L\slashed{h})_{\sigma\lambda} (\slashed{\nabla}_\delta \slashed{h})_{\kappa\rho}
	- \frac{1}{2}(\slashed{g}^{-1})^{\lambda\kappa} (L\slashed{h})_{\sigma\lambda} (\slashed{\nabla}_\kappa \slashed{h})_{\delta\rho}
	\\
	&\phantom{=}
	+ \frac{1}{4}(\slashed{\nabla}_\rho \slashed{h})_{L\delta} (\slashed{\nabla}_\sigma h)_{LL}
	+ \frac{1}{4}(\slashed{\nabla}_\rho \slashed{h})_{\Lbar\delta} (\slashed{\nabla}_\sigma h)_{LL}
	+ \frac{1}{2}(\slashed{g}^{-1})^{\lambda\kappa} (\slashed{\nabla}_\rho \slashed{h})_{\delta\lambda} (\slashed{\nabla}_\kappa \slashed{h})_{L\sigma}
	\\
	&\phantom{=}
	- \frac{1}{2}(\slashed{g}^{-1})^{\lambda\kappa} (\slashed{\nabla}_\rho \slashed{h})_{\delta\lambda} (\slashed{\nabla}_\sigma \slashed{h})_{L\kappa}
	- \frac{1}{4}(L\slashed{h})_{L\delta} (\slashed{\nabla}_\rho \slashed{h})_{L\sigma}
	- \frac{1}{4}(L\slashed{h})_{L\delta} (\slashed{\nabla}_\sigma \slashed{h})_{L\rho}
	\\
	&\phantom{=}
	- \frac{1}{4}(L\slashed{h})_{\Lbar\delta} (\slashed{\nabla}_\rho \slashed{h})_{L\sigma}
	- \frac{1}{4}(L\slashed{h})_{\Lbar\delta} (\slashed{\nabla}_\sigma \slashed{h})_{L\rho}
	+ \frac{1}{4}(L\slashed{h})_{\Lbar\delta} (L\slashed{h})_{\sigma\rho}
	\\
	&\phantom{=}
	+ \frac{1}{2}(\slashed{g}^{-1})^{\lambda\kappa} (L\slashed{h})_{\delta\lambda} (\slashed{\nabla}_\sigma \slashed{h})_{\kappa\rho}
	- \frac{1}{2}(\slashed{g}^{-1})^{\lambda\kappa} (L\slashed{h})_{\delta\lambda} (\slashed{\nabla}_\kappa \slashed{h})_{\sigma\rho}
	\end{split}
	\end{equation*}

	\begin{equation*}
	\begin{split}
	\slashed{\D}_{\Lbar} \left( (\slashed{\nabla}_\rho \slashed{h})_{\sigma\delta} \right)
	- \slashed{\nabla}_\rho \left( (\Lbar \slashed{h})_{\sigma\delta} \right)
	&=
	(\slashed{\nabla}_\rho \log \mu) (L\slashed{h})_{\sigma\delta}
	+ (\slashed{\nabla}_\rho \log \mu) (\Lbar\slashed{h})_{\sigma\delta}
	+ (\slashed{\nabla}_\delta \log \mu) (\slashed{\nabla}_\rho \slashed{h})_{L\sigma}
	\\
	&\phantom{=}
	+ (\slashed{\nabla}_\sigma \log \mu) (\slashed{\nabla}_\rho \slashed{h})_{L\delta}
	+ (\slashed{\nabla}_\delta \log \mu) (\slashed{\nabla}_\rho \slashed{h})_{\Lbar\sigma}
	+ (\slashed{\nabla}_\sigma \log \mu) (\slashed{\nabla}_\rho \slashed{h})_{\Lbar\delta}
	\\
	&\phantom{=}
	+ \left( L^i_{(\text{small})}\frac{\slashed{\nabla}_\delta x^i}{r} \right) (\slashed{\nabla}_\rho \slashed{h})_{\Lbar \sigma}
	+ \left( L^i_{(\text{small})}\frac{\slashed{\nabla}_\sigma x^i}{r} \right) (\slashed{\nabla}_\rho \slashed{h})_{\Lbar \delta}
	\\
	&\phantom{=}
	- \left( \frac{1}{r} \slashed{\Pi}_\rho^{\phantom{\rho}i} \slashed{\Pi}_\delta^{\phantom{\delta}i} \right) (\Lbar\slashed{h})_{\Lbar\sigma}
	- \left( \frac{1}{r} \slashed{\Pi}_\rho^{\phantom{\rho}i} \slashed{\Pi}_\delta^{\phantom{\sigma}i} \right) (\Lbar\slashed{h})_{\Lbar\delta}
	- \chibar_\rho^{\phantom{\rho}\kappa} (\slashed{\nabla}_\kappa \slashed{h})_{\sigma\delta}
	\\
	&\phantom{=}
	- \frac{1}{2} \chibar_{\rho\delta} (\Lbar \slashed{h})_{\Lbar\sigma}
	- \frac{1}{2} \chibar_{\rho\sigma} (\Lbar \slashed{h})_{\Lbar\delta}
	+ \frac{1}{4}(\slashed{\nabla}_\rho \slashed{h})_{L\sigma} (\slashed{\nabla}_\delta h)_{\Lbar\Lbar}
	\\
	&\phantom{=}
	- \frac{1}{4}(\slashed{\nabla}_\rho \slashed{h})_{\Lbar\sigma} (\slashed{\nabla}_\delta h)_{LL}
	- \frac{1}{4}(\slashed{\nabla}_\rho \slashed{h})_{\Lbar\sigma} (\slashed{\nabla}_\delta h)_{L\Lbar}
	- \frac{1}{2}(\slashed{g}^{-1})^{\kappa\lambda} (\slashed{\nabla}_\rho \slashed{h})_{\sigma\kappa} (\Lbar \slashed{h})_{\delta\lambda}
	\\
	&\phantom{=}
	- \frac{1}{2}(\slashed{g}^{-1})^{\kappa\lambda} (\slashed{\nabla}_\rho \slashed{h})_{\sigma\kappa} (\slashed{\nabla}_\delta \slashed{h})_{\Lbar\lambda}
	+ \frac{1}{2}(\slashed{g}^{-1})^{\kappa\lambda} (\slashed{\nabla}_\rho \slashed{h})_{\sigma\kappa} (\slashed{\nabla}_\lambda \slashed{h})_{\Lbar\delta}
	\\
	&\phantom{=}
	- \frac{1}{4}(\Lbar\slashed{h})_{L\sigma} (\slashed{\nabla}_\rho \slashed{h})_{\Lbar\delta}
	- \frac{1}{4}(\Lbar\slashed{h})_{L\sigma} (\slashed{\nabla}_\delta \slashed{h})_{\Lbar\rho}
	+ \frac{1}{4}(\Lbar\slashed{h})_{L\sigma} (\Lbar \slashed{h})_{\rho\delta}
	\\
	&\phantom{=}
	+ \frac{1}{4}(\Lbar\slashed{h})_{\Lbar\sigma} (\Lbar \slashed{h})_{\rho\delta}
	- \frac{1}{4}(\Lbar\slashed{h})_{\Lbar\sigma} (\slashed{\nabla}_\rho \slashed{h})_{\Lbar\delta}
	- \frac{1}{4}(\Lbar\slashed{h})_{\Lbar\sigma} (\slashed{\nabla}_\delta \slashed{h})_{\Lbar\rho}
	\\
	&\phantom{=}
	+ \frac{1}{2}(\slashed{g}^{-1})^{\kappa\lambda} (\Lbar \slashed{h})_{\sigma\kappa} (\slashed{\nabla}_\rho \slashed{h})_{\delta\lambda}
	+ \frac{1}{2}(\slashed{g}^{-1})^{\kappa\lambda} (\Lbar \slashed{h})_{\sigma\kappa} (\slashed{\nabla}_\delta \slashed{h})_{\rho\lambda}
	\\
	&\phantom{=}
	- \frac{1}{2}(\slashed{g}^{-1})^{\kappa\lambda} (\Lbar \slashed{h})_{\sigma\kappa} (\slashed{\nabla}_\lambda \slashed{h})_{\rho\delta}
	+ \frac{1}{4}(\slashed{\nabla}_\rho \slashed{h})_{L\delta} (\slashed{\nabla}_\sigma h)_{\Lbar\Lbar}
	- \frac{1}{4}(\slashed{\nabla}_\rho \slashed{h})_{\Lbar\delta} (\slashed{\nabla}_\sigma h)_{LL}
	\\
	&\phantom{=}
	- \frac{1}{4}(\slashed{\nabla}_\rho \slashed{h})_{\Lbar\delta} (\slashed{\nabla}_\sigma h)_{L\Lbar}
	- \frac{1}{2}(\slashed{g}^{-1})^{\kappa\lambda} (\slashed{\nabla}_\rho \slashed{h})_{\delta\kappa} (\Lbar \slashed{h})_{\sigma\lambda}
	\\
	&\phantom{=}
	- \frac{1}{2}(\slashed{g}^{-1})^{\kappa\lambda} (\slashed{\nabla}_\rho \slashed{h})_{\delta\kappa} (\slashed{\nabla}_\sigma \slashed{h})_{\Lbar\lambda}
	+ \frac{1}{2}(\slashed{g}^{-1})^{\kappa\lambda} (\slashed{\nabla}_\rho \slashed{h})_{\delta\kappa} (\slashed{\nabla}_\lambda \slashed{h})_{\Lbar\sigma}
	\\
	&\phantom{=}
	- \frac{1}{4}(\Lbar\slashed{h})_{L\delta} (\slashed{\nabla}_\rho \slashed{h})_{\Lbar\sigma}
	- \frac{1}{4}(\Lbar\slashed{h})_{L\delta} (\slashed{\nabla}_\sigma \slashed{h})_{\Lbar\rho}
	+ \frac{1}{4}(\Lbar\slashed{h})_{L\delta} (\Lbar \slashed{h})_{\rho\sigma}
	\\
	&\phantom{=}
	+ \frac{1}{4}(\Lbar\slashed{h})_{\Lbar\delta} (\Lbar \slashed{h})_{\rho\sigma}
	- \frac{1}{4}(\Lbar\slashed{h})_{\Lbar\delta} (\slashed{\nabla}_\rho \slashed{h})_{\Lbar\sigma}
	- \frac{1}{4}(\Lbar\slashed{h})_{\Lbar\delta} (\slashed{\nabla}_\sigma \slashed{h})_{\Lbar\rho}
	\\
	&\phantom{=}
	+ \frac{1}{2}(\slashed{g}^{-1})^{\kappa\lambda} (\Lbar \slashed{h})_{\delta\kappa} (\slashed{\nabla}_\rho \slashed{h})_{\sigma\lambda}
	+ \frac{1}{2}(\slashed{g}^{-1})^{\kappa\lambda} (\Lbar \slashed{h})_{\delta\kappa} (\slashed{\nabla}_\sigma \slashed{h})_{\rho\lambda}
	\\
	&\phantom{=}
	- \frac{1}{2}(\slashed{g}^{-1})^{\kappa\lambda} (\Lbar \slashed{h})_{\delta\kappa} (\slashed{\nabla}_\lambda \slashed{h})_{\rho\sigma}
	\end{split}
	\end{equation*}
	
	Note that, in the formulae above, we have used equation \ref{proposition chibar in terms of chi} to replace $\chibar - \frac{1}{2} (\Lbar \slashed{h})$ by $\chi$ and \emph{good} derivatives, or, in some cases, to replace $\chi - \frac{1}{2} (L \slashed{h})$ by $\chibar$ and a combination of $\Lbar$ and angular derivatives. Although it is easy to see why we have made the former substitution ($\chi$ and the good derivatives of the metric generally behave better than $\chibar$ and the bad derivatives of the metric), the latter requires some explanation. The reason for this unusual substitution is that the most difficult error term we will encounter takes the form $\Div \Omega_a^{\phantom{a}b}$ and, when calculating this quantity, we will need to take $L$ derivatives of some quantities, and $\Lbar$ derivatives of others. For those quantities where we will later take an $\Lbar$ derivative, it is preferable to write the quantity in terms of $L$ derivatives or $\chi$, while if we will later take an $L$ derivative, then it is preferable to write the quantity in terms of $\Lbar$ derivatives and $\chibar$.

	Putting this all together, we have the following expression:
	\begin{equation*}
	\begin{split}
	&\slashed{\Pi}_c^{\phantom{c}a} \slashed{\Pi}_b^{\phantom{b}d} (\upd \tilde{\omega}_d^{\phantom{d}c})_{\mu\nu}
	\\
	&=
	\bigg( 
		\frac{1}{8} r^{-1} (\slashed{g}^{-1})^{\lambda\sigma}\slashed{\Pi}_b^{\phantom{b}\rho} \slashed{\Pi}_\sigma^{\phantom{\sigma}a}  \left( 
			\left( \slashed{\D}_L (r\slashed{\nabla}_\rho \slashed{h}) \right)_{\Lbar \lambda} 
			- \left( \slashed{\D}_L (r\slashed{\nabla}_\lambda \slashed{h}) \right)_{\Lbar \rho}
			- \left(\slashed{\D}_{\Lbar} (r\slashed{\nabla}_\rho \slashed{h}) \right)_{L\lambda}
			+  \left(\slashed{\D}_{\Lbar} (r\slashed{\nabla}_\lambda \slashed{h}) \right)_{L\rho}
		\right)
	\\
	&\phantom{= \bigg(}
		+ \textit{l.o.t.}_{(\upd\omega , L\wedge\Lbar)}
	\bigg)L_\mu \wedge \Lbar_\nu
	\\
	&\phantom{=}
	+ \bigg(
		\frac{1}{4}r^{-1} \slashed{\Pi}_c^{\phantom{c}\rho} \slashed{\Pi}_b^{\phantom{b}\sigma} \slashed{\Pi}_d^{\phantom{d}\delta} (\slashed{g}^{-1})^{ad} \left(
			\left(\slashed{\nabla}_\rho (r\slashed{\nabla}_\delta \slashed{h})\right)_{\Lbar\sigma}
			- \left(\slashed{\nabla}_\rho (r\slashed{\nabla}_\sigma \slashed{h})\right)_{\Lbar\delta}
			- \left(\slashed{\D}_{\Lbar} (r\slashed{\nabla}_\delta \slashed{h})\right)_{\rho\sigma}
			+ \left(\slashed{\D}_{\Lbar} (r\slashed{\nabla}_\sigma \slashed{h})\right)_{\rho\delta}
		\right)
	\\
	&\phantom{= \bigg(}
	+ \textit{l.o.t.}_{(\upd\omega , L\wedge\slashed{\Pi})}
	\bigg)L_\mu \wedge \slashed{\Pi}_\nu^{\phantom{\nu}c}
	\\
	&\phantom{=}
	+ \bigg(
		\frac{1}{4}r^{-1} \slashed{\Pi}_c^{\phantom{c}\rho} \slashed{\Pi}_b^{\phantom{b}\sigma} \slashed{\Pi}_d^{\phantom{d}\delta} (\slashed{g}^{-1})^{ad} \left(
			\left(\slashed{\nabla}_\rho (r\slashed{\nabla}_\delta \slashed{h})\right)_{L\sigma}
			- \left(\slashed{\nabla}_\rho (r\slashed{\nabla}_\sigma \slashed{h})\right)_{L\delta}
			- \left(\slashed{\D}_{L} (r\slashed{\nabla}_\delta \slashed{h})\right)_{\rho\sigma}
			+ \left(\slashed{\D}_{L} (r\slashed{\nabla}_\sigma \slashed{h})\right)_{\rho\delta}
		\right)
	\\
	&\phantom{= \bigg(}
	+ \textit{l.o.t.}_{(\upd\omega , \Lbar\wedge\slashed{\Pi})}
	\bigg)\Lbar_\mu \wedge \slashed{\Pi}_\nu^{\phantom{\nu}c}
	\\
	&\phantom{=}
	+ \bigg(
	\frac{1}{2}r^{-1} \slashed{\Pi}_c^{\phantom{c}\rho} \slashed{\Pi}_b^{\phantom{b}\sigma} \slashed{\Pi}_e^{\phantom{e}\delta} \slashed{\Pi}_d^{\phantom{d}\lambda} (\slashed{g}^{-1})^{ae} \left(
	\left(\slashed{\nabla}_\lambda (r\slashed{\nabla}_\sigma \slashed{h})\right)_{\rho\delta}
	- \left(\slashed{\nabla}_\lambda (r\slashed{\nabla}_\delta \slashed{h})\right)_{\rho\sigma}
	\right)
	+ \textit{l.o.t.}_{(\upd\omega , \slashed{\Pi}\wedge\slashed{\Pi})}
	\bigg)\slashed{\Pi}_\mu^{\phantom{\mu}c} \wedge \slashed{\Pi}_\nu^{\phantom{\nu}d}
	\end{split}
	\end{equation*}
	where the lower order terms are given schematically by
	\begin{equation*}
	\begin{split}
	 \textit{l.o.t.}_{(\upd\omega , L\wedge\Lbar)}
	 &=
	 (\partial h)_{(\text{frame})} \left( r^{-1} + (\partial h)_{(\text{frame})} + \slashed{\nabla} \log \mu + r^{-1} L^i_{(\text{small})} \slashed{\Pi}^i \right)
	 \\
	 \\
	 \textit{l.o.t.}_{(\upd\omega , L\wedge\slashed{\Pi})}
	 &=
	 r^{-1}(\slashed{\nabla} h)_{(\text{frame})}
	 + (\slashed{\nabla} \log \mu) \begin{pmatrix}
	 (\slashed{\nabla} h)_{(\text{frame})} \\
	 (\Lbar h)_{(\text{frame})}
	 \end{pmatrix}
	 + \begin{pmatrix}
	 r^{-1} \Lbar^i_{(\text{small})} \slashed{\Pi}^i  \\
	 r^{-1} L^i_{(\text{small})} \slashed{\Pi}^i
	 \end{pmatrix} (\slashed{\nabla} h)_{(\text{frame})}
	 \\
	 &\phantom{=}
	 + (r^{-1} \slashed{\Pi}^i \otimes \slashed{\Pi}^i) \begin{pmatrix}
	 (\slashed{\nabla} h)_{(\text{frame})} \\
	 (\Lbar h)_{(\text{frame})}
	 \end{pmatrix}
	 + \chibar \begin{pmatrix}
	 (\slashed{\nabla} h)_{(\text{frame})} \\
	 (\Lbar h)_{(\text{frame})}
	 \end{pmatrix}
	 + (\slashed{\nabla} h)_{(\text{frame})}(\slashed{\nabla} h)_{(\text{frame})}
	 \\
	 &\phantom{=}
	 + (\slashed{\nabla} h)_{(\text{frame})}(\Lbar h)_{(\text{frame})}
	 + (\Lbar h)_{(\text{frame})}(\Lbar h)_{(\text{frame})}
	 \\
	 \\
	 \textit{l.o.t.}_{(\upd\omega , \Lbar\wedge\slashed{\Pi})}
	 &=
	 r^{-1}(\slashed{\nabla} h)_{(\text{frame})}
	 + (r^{-1} L^i_{(\text{small})} \slashed{\Pi}^i) (\slashed{\nabla} h)_{(\text{frame})}
	 + (r^{-1} \slashed{\Pi}^i \otimes \slashed{\Pi}^i) (L h)_{(\text{frame})}
	 \\
	 &\phantom{=}
	 + \chi (\bar{\partial} h)_{(\text{frame})}
	 + (\bar{\partial} h)_{(\text{frame})}(\bar{\partial} h)_{(\text{frame})}
	 \\
	 \\
	 \textit{l.o.t.}_{(\upd\omega , \slashed{\Pi}\wedge\slashed{\Pi})}
	 &=
	 (r^{-1} \slashed{\Pi}^i \otimes \slashed{\Pi}^i)(\slashed{\nabla} h)_{(\text{frame})}
	 + \chi (\slashed{\nabla} h)_{(\text{frame})}
	 + (\slashed{\nabla} h)_{(\text{frame})} (\bar{\partial} h)_{(\text{frame})}
	\end{split}
	\end{equation*} 
	Note the important cancellation in the terms $\textit{l.o.t.}_{(\upd\omega , L\wedge\slashed{\Pi})}$ of terms proportional to $(\slashed{\nabla} \log \mu) (Lh)_{(\text{frame})}$, which would otherwise have caused serious problems.

	Next, we shall compute the terms
	\begin{equation*}
	\slashed{\Pi}_a^{\phantom{a}c} \slashed{\Pi}_d^{\phantom{d}b} \upd \slashed{\Pi}_c^{\phantom{c}e} \wedge \upd \slashed{\Pi}_e^{\phantom{e}d}
	\end{equation*}
	From proposition \ref{proposition derivatives of the rectangular components of Pi} we find that
	\begin{equation*}
	\begin{split}
	\slashed{\Pi}_a^{\phantom{a}c} \partial_\mu \slashed{\Pi}_c^{\phantom{c}e}
	&=
	L_\mu \bigg( 
		\left( 
			-\frac{1}{2}(\slashed{\nabla}_\nu \log \mu)
			- \frac{1}{8} (\slashed{\nabla}_\nu h)_{\Lbar\Lbar}
			\right)
		\slashed{\Pi}_a^{\phantom{a}\nu} L^e
		\\
		&\phantom{= L_\mu \bigg(}
		+ \left(
			-\frac{1}{2}r^{-1} L^i_{(\text{small})} \slashed{\Pi}_\nu^{\phantom{\nu}i}
			- \frac{1}{2}\slashed{\nabla}_\nu \log \mu
			+ \frac{1}{8}(\slashed{\nabla}_\nu h)_{LL}
			- \frac{1}{8}(\slashed{\nabla}_\nu h)_{L\Lbar}
			\right)
		\slashed{\Pi}_a^{\phantom{a}\nu} \Lbar^e
	\bigg)
	\\
	&\phantom{=}
	+ \Lbar_\mu \bigg(
		\left( 
			\frac{1}{2}r^{-1} L^i_{(\text{small})} \slashed{\Pi}_\nu^{\phantom{\nu}i}
			- \frac{1}{8}(\slashed{\nabla}_\nu h)_{LL}
			\right)
		\slashed{\Pi}_a^{\phantom{a}\nu} L^e
		+ \left(
			-\frac{1}{8} (\slashed{\nabla}_\nu h)_{LL}
			\right)
		\slashed{\Pi}_a^{\phantom{a}\nu} \Lbar^e
	\bigg)
	\\
	&\phantom{=}
	+ \slashed{\Pi}_\mu^{\phantom{\mu}c} \bigg(
		\left(
			\frac{1}{2} \chibar_{ac}
			+ \frac{1}{4}(\slashed{\nabla}_c \slashed{h})_{\Lbar a}
			+ \frac{1}{4}(\slashed{\nabla}_a \slashed{h})_{\Lbar c}
			- \frac{1}{4}(\Lbar \slashed{h})_{ac}
			\right)
		L^e
		\\
		&\phantom{= + \slashed{\Pi}_\mu^{\phantom{\mu}c} \bigg(}
		+ 
		\left(
		\frac{1}{2} \chi_{ac}
		+ \frac{1}{4}(\slashed{\nabla}_c \slashed{h})_{L a}
		+ \frac{1}{4}(\slashed{\nabla}_a \slashed{h})_{L c}
		- \frac{1}{4}(L \slashed{h})_{ac}
		\right)
		\Lbar^e
	\bigg)
	\end{split}
	\end{equation*}
	and also that
	\begin{equation*}
	\begin{split}
	\slashed{\Pi}_d^{\phantom{d}b} \partial_\mu \slashed{\Pi}_e^{\phantom{e}d}
	&=
	L_\mu \bigg( 
		\left(
			-\frac{1}{2} (\slashed{\nabla}_\rho \log \mu)
			+ \frac{1}{4} (\Lbar \slashed{h})_{\Lbar \rho}
			- \frac{1}{8} (\slashed{\nabla}_\rho h)_{\Lbar\Lbar}
		\right)
		L_e (\slashed{g}^{-1})^{\rho\sigma} \slashed{\Pi}_\sigma^{\phantom{\sigma}b}
	\\
	&\phantom{=}
		+ \left(
			-\frac{1}{2}r^{-1} L^i_{(\text{small})} \slashed{\Pi}_\rho^{\phantom{\rho}i}
			- \frac{1}{2}\slashed{\nabla}_\rho \log \mu
			+ \frac{1}{8}(\slashed{\nabla}_\rho h)_{LL}
			- \frac{1}{4}(\slashed{\nabla}_\rho h)_{L\Lbar}
			+ \frac{1}{4}(\Lbar \slashed{h})_{L\rho}
		\right)
		\Lbar_e (\slashed{g}^{-1})^{\rho\sigma}\slashed{\Pi}_\sigma^{\phantom{\sigma}b}
	\bigg)
	\\
	&\phantom{=}
	+ \Lbar_\mu \bigg(
		\left( 
			\frac{1}{2}r^{-1} L^i_{(\text{small})} \slashed{\Pi}_\rho^{\phantom{\nu}i}
			- \frac{1}{8}(\slashed{\nabla}_\rho h)_{LL}
			+ \frac{1}{4} (L\slashed{h})_{\Lbar\rho}
			\right)
		L_e (\slashed{g}^{-1})^{\rho\sigma} \slashed{\Pi}_\sigma^{\phantom{\sigma}b}
		\\
		&\phantom{=	+ \Lbar_\mu \bigg(}
		+ \left(
			-\frac{1}{8} (\slashed{\nabla}_\rho h)_{LL}
			+ \frac{1}{4} (L\slashed{h})_{L\rho}
			\right)
		\Lbar_e (\slashed{g}^{-1})^{\rho\sigma}\slashed{\Pi}_\sigma^{\phantom{\sigma}b}
	\bigg)
	\\
	&\phantom{=}
	+ \bigg(
		\left(
			\frac{1}{2}\chibar_{\mu\rho}
			- \frac{1}{4}(\slashed{\nabla}_\mu \slashed{h})_{\Lbar\rho}
			+ \frac{1}{4}(\slashed{\nabla}_\rho \slashed{h})_{\Lbar\mu}
			- \frac{1}{4}(\Lbar \slashed{h})_{\mu\rho}
			\right)
		L_e (\slashed{g}^{-1})^{\nu\rho} \slashed{\Pi}_\nu^{\phantom{\nu}d}
		\\
		&\phantom{= + \slashed{\Pi}_\mu^{\phantom{\mu}c} \bigg(}
		+ 
		\left(
			\frac{1}{2}\chi_{\mu\rho}
			- \frac{1}{4}(\slashed{\nabla}_\mu \slashed{h})_{L\rho}
			+ \frac{1}{4}(\slashed{\nabla}_\rho \slashed{h})_{L\mu}
			- \frac{1}{4}(L \slashed{h})_{\mu\rho}
			\right)
		\Lbar_e (\slashed{g}^{-1})^{\nu\rho} \slashed{\Pi}_\nu^{\phantom{\nu}d}
	\bigg)	
	\end{split}
	\end{equation*}
	and so we have, schematically
	\begin{equation*}
	\begin{split}
	\slashed{\Pi}_a^{\phantom{a}c} \slashed{\Pi}_d^{\phantom{d}b} \upd \slashed{\Pi}_c^{\phantom{c}e} \wedge \upd \slashed{\Pi}_e^{\phantom{e}d}
	&=
	\bigg(
		(\slashed{\nabla}\log \mu)(r^{-1} L^i_{(\text{small})}\slashed{\Pi}^i)
		+ (r^{-1} L^i_{(\text{small})}\slashed{\Pi}^i)(r^{-1} L^i_{(\text{small})}\slashed{\Pi}^i)
		+ (\slashed{\nabla}\log\mu)(\bar{\partial} h)_{(\text{frame})}
		\\
		&\phantom{= \bigg(}
		+ (r^{-1} L^i_{(\text{small})}\slashed{\Pi}^i)(\bar{\partial} h)_{(\text{frame})}
		+ (\slashed{\nabla} h)_{(\text{frame})}(\partial h)_{(\text{frame})}
		\bigg)
	L_\mu \wedge \Lbar_\nu
	\\
	&\phantom{=}
	+ \bigg(
		\chibar (\slashed{\nabla} \log \mu)
		+ \chibar (r^{-1} L^i_{(\text{small})} \slashed{\Pi}^i)
		+ (\slashed{\nabla} \log \mu)(r^{-1} \slashed{\Pi}^i \otimes \slashed{\Pi}^i)
		\\
		&\phantom{= + \bigg(}
		+ \chibar \left( (\slashed{\nabla} h)_{(\text{frame})} + (\Lbar h)_{(\text{frame})} \right)
		+ (\slashed{\nabla} \log \mu) \left( (\slashed{\nabla} h)_{(\text{frame})} + (\Lbar h)_{(\text{frame})} \right)
		\\
		&\phantom{= + \bigg(}
		+ (r^{-1} L^i_{(\text{small})} \slashed{\Pi}^i) \left( (\slashed{\nabla} h)_{(\text{frame})} + (\Lbar h)_{(\text{frame})} \right)
		\\
		&\phantom{= + \bigg(}
		+ (r^{-1} \slashed{\Pi}^i \otimes \slashed{\Pi}^i) \left( (\slashed{\nabla} h)_{(\text{frame})} + (\Lbar h)_{(\text{frame})} \right)
		\\
		&\phantom{= + \bigg(}
		+ \left( (\slashed{\nabla} h)_{(\text{frame})} + (\Lbar h)_{(\text{frame})} \right) \left( (\slashed{\nabla} h)_{(\text{frame})} + (\Lbar h)_{(\text{frame})} \right)	
	\bigg) L_\mu \wedge \slashed{\Pi}_\nu^{\phantom{\nu}c}
	\\
	&\phantom{=}
	+ \bigg(
		\chi (r^{-1} L^i_{(\text{small})} \slashed{\Pi}^i)
		+ (r^{-1} L^i_{(\text{small})} \slashed{\Pi}^i) (r^{-1} \slashed{\Pi}^i \otimes \slashed{\Pi}^i)
		+ \chi (\bar{\partial} h)_{(\text{frame})}
		\\
		&\phantom{= + \bigg(}
		+ (r^{-1} L^i_{(\text{small})} \slashed{\Pi}^i) (\bar{\partial} h)_{(\text{frame})}
		+ (r^{-1} \slashed{\Pi}^i \otimes \slashed{\Pi}^i) (\bar{\partial} h)_{(\text{frame})}
		\\
		&\phantom{= + \bigg(}
		+ (\bar{\partial} h)_{(\text{frame})}(\bar{\partial} h)_{(\text{frame})}
	\bigg) \Lbar_\mu \wedge \slashed{\Pi}_\nu^{\phantom{\nu}c}
	\\
	&\phantom{=}
	+ \bigg(
		\chi \cdot \chi
		+ \chi (r^{-1} \slashed{\Pi}^i \otimes \slashed{\Pi}^i)
		+ \chi (\bar{\partial} h)_{(\text{frame})}
		+ (r^{-1} \slashed{\Pi}^i \otimes \slashed{\Pi}^i)(\bar{\partial} h)_{(\text{frame})}
		\\
		&\phantom{= + \bigg(}
		+ (\bar{\partial} h)_{(\text{frame})}(\bar{\partial} h)_{(\text{frame})}
	\bigg) \slashed{\Pi}_\mu^{\phantom{\mu}c} \wedge \slashed{\Pi}_\nu^{\phantom{\nu}d}
	\end{split}
	\end{equation*}

	Finally, we must compute the term
	\begin{equation*}
	\slashed{\Pi}_a^{\phantom{a}c} \slashed{\Pi}_d^{\phantom{d}b} \tilde{\omega}_c^{\phantom{c}e} \wedge \tilde{\omega}_e^{\phantom{e}d}
	\end{equation*}
	From proposition \ref{proposition connection coefficients omega in terms of the metric} we see that, schematically, we have
	\begin{equation*}
	\begin{split}
	\slashed{\Pi}_a^{\phantom{a}c} \slashed{\Pi}_d^{\phantom{d}b} \tilde{\omega}_c^{\phantom{c}e} \wedge \tilde{\omega}_e^{\phantom{e}d}
	&=
	(\partial h)_{(\text{frame})} (\bar{\partial} h)_{(\text{frame})} L_\mu \wedge \Lbar_\nu
	+ (\slashed{\nabla} h)_{(\text{frame})} \left( (\slashed{\nabla} h)_{(\text{frame})} + (\Lbar h)_{(\text{frame})} \right) L_\mu \wedge \slashed{\Pi}_\nu^{\phantom{\nu}c}
	\\
	&\phantom{=}
	+ (\slashed{\nabla} h)_{(\text{frame})} (\bar{\partial} h)_{(\text{frame})} \Lbar_\mu \wedge \slashed{\Pi}_\nu^{\phantom{\nu}c}
	+ (\slashed{\nabla} h)_{(\text{frame})}(\slashed{\nabla} h)_{(\text{frame})} \slashed{\Pi}_\mu^{\phantom{\mu}c} \wedge \slashed{\Pi}_\nu^{\phantom{\nu}d}
	\end{split}
	\end{equation*}

\end{proof}

\begin{definition}[Abstract indices for sections of $\mathcal{B}$]
	We will now use slashed greek indices to label abstract fibre indices associated with the vector bundle $\mathcal{B}$. So, for example, if $\phi$ is a section of $\mathcal{B}$ then we write
	\begin{equation*}
	\phi_{\slashed{\mu}} := \phi_a \slashed{\nabla}_{\slashed{\mu}} x^a
	\end{equation*}
	
	Since $\slashed{\D}$ is a metric connection with metric $\slashed{g}$, we can make use of $\slashed{g}$ and its inverse to make sense of \emph{raised} indices referring to the vector bundle $\mathcal{B}$, for example, we have
	\begin{equation*}
	\phi^{\slashed{\mu}} := (\slashed{g}^{-1})^{\slashed{\mu}\slashed{\nu}} \phi_a \slashed{\nabla}_{\slashed{\nu}} x^a
	\end{equation*}
	
	Similarly, we can use this notation to make sense sections of vector bundles over $\mathcal{M}$ whose fibres are \emph{products} of the fibres\footnote{The fibres of $\mathcal{B}$ being, of course, the vector spaces consisting of the set of $S_{\tau,r}$-tangent one-forms at a point.} of $\mathcal{B}$ or the appropriate dual vector spaces. Then, the covariant derivative $\slashed{\D}$ can be extended to act on these vector bundles in the obvious way. So, for example, the curvature $\Omega$ can be expressed as the quantity $\Omega_{\slashed{\alpha} \phantom{\slashed{\beta}} \mu\nu}^{\phantom{\slashed{\alpha}} \slashed{\beta}}$. Finally, note that we can extend the covariant derivative $\slashed{\D}$ to act on this kind of object by using the standard covariant derivative $\D$ when applied to indices referring to the tangent or cotangent bundle of $\mathcal{M}$.
	
	Finally, note that we can pass freely between spacetime indices $\mu$, $\nu$ and the fibre indices $\slashed{\mu}$, $\slashed{\nu}$ using the projection operators $\slashed{\Pi}$. For example, given a spacetime covector we can define a section of $\mathcal{B}$ by projecting the covector using $\slashed{\Pi}$. Likewise, given a section $\phi$ of $\mathcal{B}$, we can extend this to a spacetime covector by the prescription $\phi(L) = \phi(\Lbar) = 0$.
\end{definition}

\begin{proposition}[The quantities $T^\nu \Omega_{\slashed{\alpha}\slashed{\beta} \mu\nu}$, $r \Omega_{\slashed{\alpha}\slashed{\beta} \mu \slashed{\gamma}}$ and $rL^\nu \Omega_{\slashed{\alpha}\slashed{\beta} \mu\nu}$]
	\label{proposition TOmega rOmega and rLOmega}
	
	The quantity $T^\nu \Omega_{\slashed{\alpha}\slashed{\beta} \mu\nu}$ is given by
	\begin{equation*}
	\begin{split}
	T^\nu \Omega_{\slashed{\alpha}\slashed{\beta} \mu\nu}
	&=
	-\frac{1}{4}L_\mu \Omega_{\slashed{\alpha}\slashed{\beta} \Lbar L}
	+ \frac{1}{4}\Lbar_\mu \Omega_{\slashed{\alpha}\slashed{\beta} \Lbar L}
	- \frac{1}{2} \slashed{\Pi}_\mu^{\phantom{\mu}\nu} \Omega_{\slashed{\alpha}\slashed{\beta} L\nu}
	- \frac{1}{2} \slashed{\Pi}_\mu^{\phantom{\mu}\nu} \Omega_{\slashed{\alpha}\slashed{\beta} \Lbar\nu}
	\\
	&=
	\frac{1}{8}r^{-1} \left(
		-\left(\slashed{\D}_L (r\slashed{\nabla}_{\slashed{\alpha}} \slashed{h})\right)_{\Lbar\slashed{\beta}}
		+ \left(\slashed{\D}_L (r\slashed{\nabla}_{\slashed{\beta}}\slashed{h})\right)_{\Lbar\slashed{\alpha}}
		+ \left(\slashed{\D}_{\Lbar} (r\slashed{\nabla}_{\slashed{\alpha}}\slashed{h})\right)_{L\slashed{\beta}}
		- \left(\slashed{\D}_{\Lbar} (r\slashed{\nabla}_{\slashed{\beta}}\slashed{h})\right)_{L\slashed{\alpha}}
	\right) L_\mu
	\\
	&\phantom{=}
	+ \frac{1}{8}r^{-1} \left(
		\left(\slashed{\D}_L (r\slashed{\nabla}_{\slashed{\alpha}}\slashed{h})\right)_{\Lbar\slashed{\beta}}
		- \left(\slashed{\D}_L (r\slashed{\nabla}_{\slashed{\beta}}\slashed{h})\right)_{\Lbar\slashed{\alpha}}
		- \left(\slashed{\D}_{\Lbar} (r\slashed{\nabla}_{\slashed{\alpha}}\slashed{h})\right)_{L\slashed{\beta}}
		+ \left(\slashed{\D}_{\Lbar} (r\slashed{\nabla}_{\slashed{\beta}}\slashed{h})\right)_{L\slashed{\alpha}}
	\right) \Lbar_\mu
	\\
	&\phantom{=}
	+ \frac{1}{4}r^{-1} \left(
		\left( \slashed{\nabla}_{\slashed{\gamma}} (r \slashed{\nabla}_{\slashed{\beta}} \slashed{h})\right)_{L\slashed{\alpha}}
		- \left( \slashed{\nabla}_{\slashed{\gamma}} (r \slashed{\nabla}_{\slashed{\alpha}} \slashed{h}) \right)_{L\slashed{\beta}}
		- \left( \slashed{\D}_{L} (r \slashed{\nabla}_{\slashed{\beta}} \slashed{h}) \right)_{\slashed{\alpha}\slashed{\gamma}}
		+ \left( \slashed{\D}_{L} (r \slashed{\nabla}_{\slashed{\alpha}} \slashed{h}) \right)_{\slashed{\beta}\slashed{\gamma}}
	\right) \slashed{\Pi}_\mu^{\phantom{\mu}\slashed{\gamma}}
	\\
	&\phantom{=}
	+ \frac{1}{4}r^{-1} \left(
		\left( \slashed{\nabla}_{\slashed{\gamma}} (r \slashed{\nabla}_{\slashed{\beta}} \slashed{h})\right)_{\Lbar\slashed{\alpha}}
		- \left( \slashed{\nabla}_{\slashed{\gamma}} (r \slashed{\nabla}_{\slashed{\alpha}} \slashed{h}) \right)_{\Lbar\slashed{\beta}}
		- \left( \slashed{\D}_{\Lbar} (r \slashed{\nabla}_{\slashed{\beta}} \slashed{h}) \right)_{\slashed{\alpha}\slashed{\gamma}}
		+ \left( \slashed{\D}_{\Lbar} (r \slashed{\nabla}_{\slashed{\alpha}} \slashed{h}) \right)_{\slashed{\beta}\slashed{\gamma}}
	\right) \slashed{\Pi}_\mu^{\phantom{\mu}\slashed{\gamma}}
	\\
	&\phantom{=}
	+ \textit{l.o.t.}_{(\Omega, L\wedge\Lbar)} (L_\mu + \Lbar_\mu)
	+ \left( \textit{l.o.t.}_{(\Omega, L\wedge\slashed{\Pi})} + \textit{l.o.t.}_{(\Omega, \Lbar\wedge\slashed{\Pi})}\right) \slashed{\Pi}_\mu^{\phantom{\mu}\slashed{\gamma}}
	\end{split}
	\end{equation*}
	and the quantity
	\begin{equation*}
	r \Omega_{\slashed{\alpha}\slashed{\beta} \mu \slashed{\gamma}}
	=
	r \slashed{\Pi}_{\slashed{\gamma}}^{\phantom{\slashed{\gamma}}\nu} \Omega_{\slashed{\alpha}\slashed{\beta} \mu \nu}
	\end{equation*}
	is given by
	\begin{equation*}
	\begin{split}
	r \Omega_{\slashed{\alpha}\slashed{\beta} \mu \slashed{\gamma}}
	&=
	\frac{1}{4} \left(
		\left( \slashed{\nabla}_{\slashed{\gamma}} (r \slashed{\nabla}_{\slashed{\beta}} \slashed{h})\right)_{\Lbar\slashed{\alpha}}
		- \left( \slashed{\nabla}_{\slashed{\gamma}} (r \slashed{\nabla}_{\slashed{\alpha}} \slashed{h}) \right)_{\Lbar\slashed{\beta}}
		- \left( \slashed{\D}_{\Lbar} (r \slashed{\nabla}_{\slashed{\beta}} \slashed{h}) \right)_{\slashed{\alpha}\slashed{\gamma}}
		+ \left( \slashed{\D}_{\Lbar} (r \slashed{\nabla}_{\slashed{\alpha}} \slashed{h}) \right)_{\slashed{\beta}\slashed{\gamma}}
	\right) L_\mu
	\\
	&\phantom{=}
	+ \frac{1}{4} \left(
		\left( \slashed{\nabla}_{\slashed{\gamma}} (r \slashed{\nabla}_{\slashed{\beta}} \slashed{h})\right)_{L\slashed{\alpha}}
		- \left( \slashed{\nabla}_{\slashed{\gamma}} (r \slashed{\nabla}_{\slashed{\alpha}} \slashed{h}) \right)_{L\slashed{\beta}}
		- \left( \slashed{\D}_{L} (r \slashed{\nabla}_{\slashed{\beta}} \slashed{h}) \right)_{\slashed{\alpha}\slashed{\gamma}}
		+ \left( \slashed{\D}_{L} (r \slashed{\nabla}_{\slashed{\alpha}} \slashed{h}) \right)_{\slashed{\beta}\slashed{\gamma}}
	\right) \Lbar_\mu
	\\
	&\phantom{=}
	+ \frac{1}{2} \left(
		\left( \slashed{\nabla}_{\slashed{\gamma}} (r\slashed{\nabla}_{\slashed{\alpha}} \slashed{h}) \right)_{\slashed{\beta}\mu}
		- \left( \slashed{\nabla}_{\slashed{\gamma}} (r\slashed{\nabla}_{\slashed{\beta}} \slashed{h}) \right)_{\slashed{\alpha}\mu}
		- \left( \slashed{\nabla}_{\mu} (r\slashed{\nabla}_{\slashed{\alpha}} \slashed{h}) \right)_{\slashed{\beta}\slashed{\gamma}}
		+ \left( \slashed{\nabla}_{\mu} (r\slashed{\nabla}_{\slashed{\beta}} \slashed{h}) \right)_{\slashed{\alpha}\slashed{\gamma}}
	\right)
	\\
	&\phantom{=}
	+ r \cdot \left(\textit{l.o.t.}_{(\Omega, L\wedge\slashed{\Pi})} \right)L_\mu
	+ r \cdot \left(\textit{l.o.t.}_{(\Omega, \Lbar\wedge\slashed{\Pi})} \right)\Lbar_\mu
	+ r \cdot \left(\textit{l.o.t.}_{(\Omega, \slashed{\Pi}\wedge\slashed{\Pi})} \right) \slashed{\Pi}_\mu^{\phantom{\mu}\slashed{\gamma}}
	\end{split}
	\end{equation*}
	
	Finally, the quantity $rL^\nu \Omega_{\slashed{\alpha}\slashed{\beta}\mu\nu}$ is given by
	\begin{equation*}
	\begin{split}
	rL^\nu \Omega_{\slashed{\alpha}\slashed{\beta}\mu\nu}
	&=
	\bigg( \frac{1}{4} (\slashed{g}^{-1})^{\lambda\sigma}\slashed{\Pi}_b^{\phantom{b}\rho} \slashed{\Pi}_\sigma^{\phantom{\sigma}a}  \left( 
	\left( \slashed{\D}_L (r\slashed{\nabla}_\rho \slashed{h}) \right)_{\Lbar \lambda} 
	- \left( \slashed{\D}_L (r\slashed{\nabla}_\lambda \slashed{h}) \right)_{\Lbar \rho}
	- \left(\slashed{\D}_{\Lbar} (r\slashed{\nabla}_\rho \slashed{h}) \right)_{L\lambda}
	+  \left(\slashed{\D}_{\Lbar} (r\slashed{\nabla}_\lambda \slashed{h}) \right)_{L\rho}
	\right)
	\\
	&\phantom{= \bigg(}
	+ r\textit{l.o.t.}_{(\Omega , L\wedge\Lbar)}
	\bigg)L_\mu
	\\
	&\phantom{=} + \bigg(
	\frac{1}{2}r^{-1} \slashed{\Pi}_c^{\phantom{c}\rho} \slashed{\Pi}_b^{\phantom{b}\sigma} \slashed{\Pi}_d^{\phantom{d}\delta} (\slashed{g}^{-1})^{ad} \Big(
	- \left(\slashed{\nabla}_\rho (r\slashed{\nabla}_\delta \slashed{h})\right)_{L\sigma}
	+ \left(\slashed{\nabla}_\rho (r\slashed{\nabla}_\sigma \slashed{h})\right)_{L\delta}
	+ \left(\slashed{\D}_{L} (r\slashed{\nabla}_\delta \slashed{h})\right)_{\rho\sigma}
	\\
	&\phantom{= + \bigg(
	\frac{1}{2} \slashed{\Pi}_c^{\phantom{c}\rho} \slashed{\Pi}_b^{\phantom{b}\sigma}
	\slashed{\Pi}_d^{\phantom{d}\delta} (\slashed{g}^{-1})^{ad} \Big(}
	- \left(\slashed{\D}_{L} (r\slashed{\nabla}_\sigma \slashed{h})\right)_{\rho\delta} \Big)
	+ r\textit{l.o.t.}_{(\Omega , \Lbar\wedge\slashed{\Pi})}
	\bigg)\slashed{\Pi}_\nu^{\phantom{\nu}c}
	\end{split}
	\end{equation*}
	
\end{proposition}

\begin{proof}
	This follows immediately from proposition \ref{proposition expression for Omega} and the definition $T^\mu = \frac{1}{2}(L^\mu + \Lbar^\mu)$.
\end{proof}

\begin{proposition}[The divergences of $T^\nu \Omega_{\slashed{\alpha}\slashed{\beta} \mu\nu}$, $r \Omega_{\slashed{\alpha}\slashed{\beta} \mu \slashed{\gamma}}$ and $rL^\nu \Omega_{\slashed{\alpha}\slashed{\beta} \mu\nu}$]
\label{proposition div Z Omega}
The divergences of the $S_{\tau,r}$-tangent tensor fields $T^\nu \Omega_{\slashed{\alpha}\slashed{\beta} \mu\nu}$, $r \Omega_{\slashed{\alpha}\slashed{\beta} \mu \slashed{\gamma}}$ and $rL^\nu \Omega_{\slashed{\alpha}\slashed{\beta} \mu\nu}$ satisfy the schematic equations
\begin{equation}
\begin{split}
\slashed{\D}^\mu \left( T^\nu \Omega_{\slashed{\alpha}\slashed{\beta} \mu\nu} \right)
&=
r^{-1} \left(\slashed{\D} \slashed{\D}_T \left( r\slashed{\nabla} h \right) \right)_{(\text{frame})}
+ r^{-2} \left(\slashed{\D} \left( r\slashed{\nabla} \left( r\slashed{\nabla} h \right) \right) \right)_{(\text{frame})}
+ r^{-1} \tilde{\slashed{\Box}}_g (r\slashed{\nabla}h)_{(\text{frame})}
\\
&\phantom{=}
+ \bm{\Gamma}_{(-1)} \cdot \begin{pmatrix}
	(\partial T h)_{(\text{frame})} \\
	\left(\D (r\slashed{\nabla} h)\right)_{(\text{frame})} \\
	(\tilde{\Box}_g h)_{(\text{frame})} \\
	\D_T \chi_{(\text{small})} \\
	r\slashed{\nabla} \chi_{(\text{small})} \\
	r^{-1} \slashed{\D}_T (r\slashed{\nabla} \log \mu) \\
	\slashed{\nabla} (r\slashed{\nabla} \log \mu) \\
	\bm{\Gamma}_{(-1)} \cdot \bm{\Gamma}_{(-1, \text{small})}
	\end{pmatrix}
\end{split}
\end{equation}
and
\begin{equation}
\slashed{\D}^\mu \left( r \Omega_{\slashed{\alpha}\slashed{\beta} \mu \slashed{\gamma}}\right)
=
r^{-1} \left(\slashed{\D} \left( r\slashed{\nabla} (r\slashed{\nabla} h) \right) \right)_{(\text{frame})}
+ \left( \tilde{\slashed{\Box}}_g (r\slashed{\nabla} h) \right)_{(\text{frame})}
+ \bm{\Gamma}_{(-1)} \cdot \begin{pmatrix}
	\left(\slashed{\D}(r\slashed{\nabla} h)\right)_{(\text{frame})} \\
	r(\tilde{\Box}_g h)_{(\text{frame})} \\
	r^{-1} (r\slashed{\nabla})^2 \log \mu \\
	r\slashed{\nabla} \chi_{(\text{small})} \\
	r\bm{\Gamma}_{(-1)} \cdot \bm{\Gamma}_{(-1, \text{small})}
\end{pmatrix}
\end{equation}
and also
\begin{equation}
\begin{split}
\slashed{\D}^\mu \left( T^\nu \Omega_{\slashed{\alpha}\slashed{\beta} \mu\nu} \right)
&=
\left(\slashed{\D}_L \left( \slashed{\D}_T (r\slashed{\nabla} h) \right) \right)_{(\text{frame})}
+ r^{-1} \left(\slashed{\nabla} \left( r\slashed{\nabla} (r\slashed{\nabla} h) \right) \right)_{(\text{frame})}
+ \left( \tilde{\slashed{\Box}}_g (r\slashed{\nabla} h) \right)_{(\text{frame})}
\\
&\phantom{=}
+ r\bm{\Gamma}_{(-1)} \cdot \begin{pmatrix}
\left(L (T h)\right)_{(\text{frame})} \\
(\tilde{\Box}_g h)_{(\text{frame})} \\
r^{-1} \left( \slashed{\D} (r\slashed{\nabla} h) \right)_{(\text{frame})} \\
\slashed{\nabla} \chi_{(\text{small})} \\
\bm{\Gamma}_{(-1)} \cdot \bm{\Gamma}_{(-1, \text{small})}
\end{pmatrix}
\end{split}
\end{equation}

where we have defined
\begin{equation*}
\begin{split}
\bm{\Gamma}_{(-1)} &:= \left\{
r^{-1} \ , \ (\partial h)_{(\text{frame})} \ , \ \chi \ , \ \slashed{\nabla}\log\mu \ , \ r^{-1}L^i \ , \ r^{-1} \Lbar^i \ , \ r^{-1} \slashed{\Pi}^i \otimes \slashed{\Pi}^i
\right\}
\\
\bm{\Gamma}_{(-1, \text{small})} &:= \left\{
(\partial h)_{(\text{frame})} \ , \ \chi_{(\text{small})} \ , \ \slashed{\nabla}\log\mu \ , \ r^{-1}L^i_{(\text{small})} \ , \ r^{-1} \Lbar^i_{(\text{small})}
\right\}
\\
\end{split}
\end{equation*}
and where it is possible that additional factors of $r^{-1}$ might be present\footnote{Note that we will only be considering this kind of error term away from the origin, and so terms with additional factors of $r^{-1}$ are \emph{better} behaved.}
	
\end{proposition}

\begin{proof}
	We begin by computing
	\begin{equation*}
	\slashed{\D}^\mu \left( f^L L_\mu + f^{\Lbar} \Lbar_\mu + f_{\slashed{\gamma}} \slashed{\Pi}_\mu^{\phantom{\mu}\slashed{\gamma}} \right)
	\end{equation*}
	where $f^L$ and $f^{\Lbar}$ are $S_{\tau,r}$-tangent tensor fields, as is $f_{\slashed{\gamma}} \slashed{\Pi}_\mu^{\phantom{\mu}\slashed{\gamma}} (X_A)^\mu $. Then we find
	\begin{equation*}
	\slashed{\D}^\mu \left( f^L L_\mu + f^{\Lbar} \Lbar_\mu + f_{\slashed{\gamma}} \slashed{\Pi}_\mu^{\phantom{\mu}\slashed{\gamma}} \right)
	=
	\slashed{\D}_L f^L + \slashed{\D}_{\Lbar} f^{\Lbar} + \slashed{\nabla}^{\slashed{\gamma}} f_{\slashed{\gamma}}
	+ (\tr_{\slashed{g}}\chi + \omega) f^L
	+ (\tr_{\slashed{g}}\chibar - \omega) f^{\Lbar}
	- \left( \slashed{\nabla}^{\slashed{\gamma}} \log \mu \right) f_{\slashed{\gamma}}
	\end{equation*}

	Now we turn to the quantity
	\begin{equation*}
	\slashed{\D}^\mu \left( T^\nu \Omega_{\slashed{\alpha}\slashed{\beta} \mu\nu} \right)
	\end{equation*}
	Using the expression given in \ref{proposition TOmega rOmega and rLOmega} we find that
	\begin{equation*}
	\begin{split}
	\slashed{\D}^\mu \left( T^\nu \Omega_{\slashed{\alpha}\slashed{\beta} \mu\nu} \right)
	&=
	- \frac{1}{8}r^{-1} \left(\slashed{\D}_L\slashed{\D}_L (r\slashed{\nabla}_{\slashed{\alpha}}\slashed{h})\right)_{\Lbar\slashed{\beta}}
	+ \frac{1}{8}r^{-1} \left(\slashed{\D}_L\slashed{\D}_L (r\slashed{\nabla}_{\slashed{\beta}}\slashed{h})\right)_{\Lbar\slashed{\alpha}}
	+ \frac{1}{8}r^{-1} \left(\slashed{\D}_L\slashed{\D}_{\Lbar} (r\slashed{\nabla}_{\slashed{\alpha}}\slashed{h})\right)_{L\slashed{\beta}}
	\\
	&\phantom{=}
	- \frac{1}{8}r^{-1} \left(\slashed{\D}_L\slashed{\D}_{\Lbar} (r\slashed{\nabla}_{\slashed{\beta}}\slashed{h})\right)_{L\slashed{\alpha}}
	+ \frac{1}{8}r^{-1}	\left(\slashed{\D}_{\Lbar}\slashed{\D}_L (r\slashed{\nabla}_{\slashed{\alpha}}\slashed{h})\right)_{\Lbar\slashed{\beta}}
	- \frac{1}{8}r^{-1} \left(\slashed{\D}_{\Lbar}\slashed{\D}_L (r\slashed{\nabla}_{\slashed{\beta}}\slashed{h})\right)_{\Lbar\slashed{\alpha}}
	\\
	&\phantom{=}
	- \frac{1}{8}r^{-1} \left(\slashed{\D}_{\Lbar}\slashed{\D}_{\Lbar} (r\slashed{\nabla}_{\slashed{\alpha}}\slashed{h})\right)_{L\slashed{\beta}}
	+ \frac{1}{8}r^{-1}\left(\slashed{\D}_{\Lbar}\slashed{\D}_{\Lbar} (r\slashed{\nabla}_{\slashed{\beta}}\slashed{h})\right)_{L\slashed{\alpha}}
	+ \frac{1}{4}r^{-1} \left( \slashed{\Delta} (r \slashed{\nabla}_{\slashed{\beta}} \slashed{h})\right)_{L\slashed{\alpha}}
	\\
	&\phantom{=}
	- \frac{1}{4}r^{-1} \left( \slashed{\Delta} (r \slashed{\nabla}_{\slashed{\alpha}} \slashed{h}) \right)_{L\slashed{\beta}}
	- \frac{1}{4}r^{-1} \left( \slashed{\nabla}^{\slashed{\gamma}} \slashed{\D}_{L} (r \slashed{\nabla}_{\slashed{\beta}} \slashed{h}) \right)_{\slashed{\alpha}\slashed{\gamma}}
	+ \frac{1}{4}r^{-1} \left( \slashed{\nabla}^{\slashed{\gamma}} \slashed{\D}_{L} (r \slashed{\nabla}_{\slashed{\alpha}} \slashed{h}) \right)_{\slashed{\beta}\slashed{\gamma}}
	\\
	&\phantom{=}
	+ \frac{1}{4}r^{-1}	\left( \slashed{\Delta} (r \slashed{\nabla}_{\slashed{\beta}} \slashed{h})\right)_{\Lbar\slashed{\alpha}}
	- \frac{1}{4}r^{-1} \left( \slashed{\Delta} (r \slashed{\nabla}_{\slashed{\alpha}} \slashed{h}) \right)_{\Lbar\slashed{\beta}}
	- \frac{1}{4}r^{-1} \left( \slashed{\nabla}^{\slashed{\gamma}}\slashed{\D}_{\Lbar} (r \slashed{\nabla}_{\slashed{\beta}} \slashed{h}) \right)_{\slashed{\alpha}\slashed{\gamma}}
	\\
	&\phantom{=}
	+ \frac{1}{4}r^{-1} \left( \slashed{\nabla}^{\slashed{\gamma}}\slashed{\D}_{\Lbar} (r \slashed{\nabla}_{\slashed{\alpha}} \slashed{h}) \right)_{\slashed{\beta}\slashed{\gamma}}
	+ \ldots
	\\
	\\
	&=
	\frac{1}{4}r^{-1} \left(\slashed{\D}_L \slashed{\D}_T (r\slashed{\nabla}_{\slashed{\beta}} \slashed{h}) \right)_{\Lbar\slashed{\alpha}}
	+ \frac{1}{4}r^{-1} \left( \tilde{\slashed{\Box}}_g (r\slashed{\nabla}_{\slashed{\beta}} \slashed{h}) \right)_{\Lbar\slashed{\alpha}}
	- \frac{1}{4}r^{-1} \left(\slashed{\D}_L \slashed{\D}_T (r\slashed{\nabla}_{\slashed{\alpha}} \slashed{h}) \right)_{\Lbar\slashed{\beta}}
	\\
	&\phantom{=}
	- \frac{1}{4}r^{-1} \left( \tilde{\slashed{\Box}}_g (r\slashed{\nabla}_{\slashed{\alpha}} \slashed{h}) \right)_{\Lbar\slashed{\beta}}
	+ \frac{1}{4}r^{-1} \left(\slashed{\D}_{\Lbar} \slashed{\D}_T (r\slashed{\nabla}_{\slashed{\beta}} \slashed{h}) \right)_{L\slashed{\alpha}}
	+ \frac{1}{4}r^{-1} \left( \tilde{\slashed{\Box}}_g (r\slashed{\nabla}_{\slashed{\beta}} \slashed{h})\right)_{L\slashed{\alpha}}
	\\
	&\phantom{=}
	- \frac{1}{4}r^{-1} \left(\slashed{\D}_{\Lbar} \slashed{\D}_T (r\slashed{\nabla}_{\slashed{\alpha}} \slashed{h}) \right)_{L\slashed{\beta}}
	- \frac{1}{4}r^{-1} \left( \tilde{\slashed{\Box}}_g (r\slashed{\nabla}_{\slashed{\alpha}} \slashed{h})\right)_{L\slashed{\beta}}
	- \frac{1}{4}r^{-2} \left(\slashed{\D}_L (r^2 \slashed{\nabla}^{\slashed{\gamma}} \slashed{\nabla}_{\slashed{\beta}} \slashed{h}) \right)_{\slashed{\alpha}\slashed{\gamma}}
	\\
	&\phantom{=}
	- \frac{1}{4}r^{-2} \left(\slashed{\D}_{\Lbar} (r^2 \slashed{\nabla}^{\slashed{\gamma}} \slashed{\nabla}_{\slashed{\beta}} \slashed{h}) \right)_{\slashed{\alpha}\slashed{\gamma}}
	+ \frac{1}{4}r^{-2} \left(\slashed{\D}_L (r^2 \slashed{\nabla}^{\slashed{\gamma}} \slashed{\nabla}_{\slashed{\alpha}} \slashed{h}) \right)_{\slashed{\beta}\slashed{\gamma}}
	+ \frac{1}{4}r^{-2} \left(\slashed{\D}_{\Lbar} (r^2 \slashed{\nabla}^{\slashed{\gamma}} \slashed{\nabla}_{\slashed{\alpha}} \slashed{h}) \right)_{\slashed{\beta}\slashed{\gamma}}
	\\
	&\phantom{=}
	+ \textit{l.o.t.}_{\left(\Div(T\Omega)\right)}
	\end{split}
	\end{equation*}
	where the ellipsis in the first equality stand for some lower order terms which, however, are different from the lower order terms denoted by $\textit{l.o.t.}_{\left(\Div(T\Omega)\right)}$ (though schematically they are the same). These lower order terms are given schematically by
	\begin{equation*}
	\textit{l.o.t.}_{\left(\Div(T\Omega)\right)}
	= (\slashed{\D} \bm{\Gamma}_{(-1, \text{small})}) \cdot \bm{\Gamma}_{(-1)}
	+ (\slashed{\D} \bm{\Gamma}_{(-1)}) \cdot \bm{\Gamma}_{(-1, \text{small})}
	\end{equation*}
	with possible additional factors of $r^{-1}$.
	
	Note that we have, schematically,
	\begin{equation*}
	\slashed{\D} \bm{\Gamma}_{(-1)}
	=
	\begin{pmatrix}
	\bm{\Gamma}_{(-1)} \cdot \bm{\Gamma}_{(-1)} \\
	(\partial T h)_{(\text{frame})} \\
	\left(\slashed{\D} (r\slashed{\nabla} h) \right)_{(\text{frame})} \\
	(\tilde{\Box}_g h)_{(\text{frame})} \\
	\mathscr{Z} \chi_{(\text{small})} \\
	r^{-1}\slashed{\D}_T (r\slashed{\nabla} \log\mu) \\
	\slashed{\nabla} (r\slashed{\nabla} \log\mu) \\
	\end{pmatrix}
	\end{equation*}
	and
	\begin{equation*}
	\slashed{\D} \bm{\Gamma}_{(-1, \text{small})}
	=
	\begin{pmatrix}
	\bm{\Gamma}_{(-1, \text{small})} \cdot \bm{\Gamma}_{(-1)} \\
	(\partial T h)_{(\text{frame})} \\
	\left(\slashed{\D} (r\slashed{\nabla} h) \right)_{(\text{frame})} \\
	(\tilde{\Box}_g h)_{(\text{frame})} \\
	\mathscr{Z} \chi_{(\text{small})} \\
	r^{-1}\slashed{\D}_T (r\slashed{\nabla} \log\mu) \\
	\slashed{\nabla} (r\slashed{\nabla} \log\mu) \\
	\end{pmatrix}
	\end{equation*}
	where again, additional factors of $r^{-1}$ might be present, in addition to factors of the form $L^a$, $\Lbar^a$ or $\slashed{\Pi}^a$ which are expected to behave like constants.

	Next, we turn to the quantity
	\begin{equation*}
	\slashed{\D}^\mu \left( r\Omega_{\slashed{\alpha}\slashed{\beta} \mu \slashed{\gamma}} \right)
	=
	\slashed{\D}^\mu \left( r\slashed{\Pi}_{\slashed{\gamma}}^{\phantom{\slashed{\gamma}}\nu}\Omega_{\slashed{\alpha}\slashed{\beta} \mu \nu} \right)
	\end{equation*}
	This requires a more delicate treatment, since the factor of $r$ can potentially cause a slower decay in the radial direction. This could cause serious problems: when commuting we will encounter error terms involving this quantity, and to close our estimates we require it to decay no slower than $r^{-2+\epsilon}$. In particular, even terms of the form $r(\partial h)_{(\text{frame})} (\bar{\partial} h)_{(\text{frame})}$ in $r\Omega_{\slashed{\alpha}\slashed{\beta} \mu \slashed{\gamma}}$ can be problematic: if we take a derivative of this quantity then we can expect terms of the form $r(\partial \mathscr{Z} h)_{(\text{frame})} (\bar{\partial} h)_{(\text{frame})}$, which does not have the required decay.
	
	Using proposition \ref{proposition TOmega rOmega and rLOmega} we find that
	\begin{equation*}
	\begin{split}
	\slashed{\D}^\mu \left( r\Omega_{\slashed{\alpha}\slashed{\beta} \mu \slashed{\gamma}} \right)
	&=
	\frac{1}{4}r^{-1} \left( \slashed{\D}_L \left( r\slashed{\nabla}_{\slashed{\gamma}} (r\slashed{\nabla}_{\slashed{\beta}} \slashed{h}) \right) \right)_{\Lbar\slashed{\alpha}}
	- \frac{1}{4}r^{-1} \left( \slashed{\D}_L \left( r\slashed{\nabla}_{\slashed{\gamma}} (r\slashed{\nabla}_{\slashed{\alpha}} \slashed{h}) \right) \right)_{\Lbar\slashed{\beta}}
	\\
	&\phantom{=}
	+ \frac{1}{4}r^{-1} \left( \slashed{\D}_{\Lbar} \left( r\slashed{\nabla}_{\slashed{\gamma}} (r\slashed{\nabla}_{\slashed{\beta}} \slashed{h}) \right) \right)_{L\slashed{\alpha}}
	- \frac{1}{4}r^{-1} \left( \slashed{\D}_{\Lbar} \left( r\slashed{\nabla}_{\slashed{\gamma}} (r\slashed{\nabla}_{\slashed{\alpha}} \slashed{h}) \right) \right)_{L\slashed{\beta}}
	\\
	&\phantom{=}
	+ \frac{1}{2} r^{-1} \left( \slashed{\nabla}^{\slashed{\mu}} \left( r\slashed{\nabla}_{\slashed{\gamma}} (r\slashed{\nabla}_{\slashed{\alpha}} \slashed{h}) \right) \right)_{\slashed{\beta}\slashed{\mu}}
	- \frac{1}{2} r^{-1} \left( \slashed{\nabla}^{\slashed{\mu}} \left( r\slashed{\nabla}_{\slashed{\gamma}} (r\slashed{\nabla}_{\slashed{\beta}} \slashed{h}) \right) \right)_{\slashed{\alpha}\slashed{\mu}}
	\\
	&\phantom{=}
	+ \frac{1}{2}\left( \tilde{\slashed{\Box}}_g (r\slashed{\nabla}_{\slashed{\beta}} \slashed{h}) \right)_{\slashed{\alpha} \slashed{\gamma}}
	- \frac{1}{2}\left( \tilde{\slashed{\Box}}_g (r\slashed{\nabla}_{\slashed{\beta}} \slashed{h}) \right)_{\slashed{\alpha} \slashed{\gamma}}
	+ \textit{l.o.t.}_{(\Div(\slashed{\Omega}))}
	\end{split}
	\end{equation*}
	where the lower order terms are given schematically by
	\begin{equation*}
	\begin{split}
	\textit{l.o.t.}_{(\Div(\slashed{\Omega}))}
	&=
	\slashed{\D}_L \left( r\left( \textit{l.o.t.}_{(\Omega, L\wedge\slashed{\Pi})} \right) \right)
	+ \slashed{\D}_{\Lbar} \left( r \left( \textit{l.o.t.}_{(\Omega, \Lbar\wedge\slashed{\Pi})} \right) \right)
	+ \slashed{\nabla} \left( r\left( \textit{l.o.t.}_{(\Omega, \slashed{\Pi}\wedge\slashed{\Pi})} \right) \right)
	\\
	&\phantom{=}
	+ \bm{\Gamma}_{(-1)} \begin{pmatrix}
		\left(\slashed{\D}(\mathscr{Z} h) \right)_{(\text{frame})} \\
		r\left(\textit{l.o.t.}_{(\Omega, L\wedge\slashed{\Pi})} \right)\\
		r\left(\textit{l.o.t.}_{(\Omega, \Lbar\wedge\slashed{\Pi})} \right)\\
		r\left(\textit{l.o.t.}_{(\Omega, \slashed{\Pi}\wedge\slashed{\Pi})} \right)
	\end{pmatrix}
	\end{split}
	\end{equation*}
	Using the expressions for the lower order terms given in proposition \ref{proposition expression for Omega} we have (again schematically)
	\begin{equation*}
	\begin{split}
	\slashed{\D}_L \left( r\left( \textit{l.o.t.}_{(\Omega, L\wedge\slashed{\Pi})} \right) \right)
	&=
	\bm{\Gamma}_{(-1)} \cdot \begin{pmatrix}
		\left(\mathscr{\D}_L(r\slashed{\nabla} h)\right)_{(\text{frame})} \\
		\left(L(r\Lbar h)\right)_{(\text{frame})} \\
		\slashed{\D}_L \left(r\slashed{\nabla} \log \mu \right) \\
		\slashed{\D}_L \left(r\chibar_{(\text{small})} \right) \\
		r\bm{\Gamma}_{(-1, \text{small})} \cdot \bm{\Gamma}_{(-1)}
		\end{pmatrix}
	\end{split}
	\end{equation*}
	We now note that we have
	\begin{equation*}
	\begin{split}
	\left(L(r\Lbar h)\right)_{(\text{frame})}
	&=
	- r(\tilde{\Box}_g h)_{(\text{frame})}
	+ r(\slashed{\Delta}h)_{(\text{frame})}
	- \frac{1}{2}r \tr_{\slashed{g}}\chi_{(\text{small})} (\Lbar h)_{(\text{frame})}
	- \frac{1}{2} r \tr_{\slashed{g}} \chibar (L h)_{(\text{frame})}
	\\
	&\phantom{=}
	- r\zeta^\alpha (\slashed{\nabla}_\alpha h)_{(\text{frame})}
	\\
	\\
	&=
	r(\tilde{\Box}_g h)_{(\text{frame})}
	+ \left(\slashed{\nabla} (r\slashed{\nabla} h)\right)_{(\text{frame})}
	+ r\bm{\Gamma}_{(-1)}(\partial h)_{(\text{frame})}
	\end{split}
	\end{equation*}
	where in the second line we have given a schematic expression.
	
	We also note that
	\begin{equation*}
	\begin{split}
	\slashed{\D}_L (r\slashed{\nabla}_{\alpha} \log \mu)
	&=
	r\slashed{\nabla}_{\alpha} (L\log\mu)
	- r \left(\chi_{(\text{small})}\right)_\alpha^{\phantom{\alpha}\beta} \slashed{\nabla}_\beta \log \mu
	\\ \\
	&=
	r\slashed{\nabla}_{\alpha} \left( -r^{-1} L^i_{(\text{small})}L^i_{(\text{small})}
		- \frac{1}{2}(Lh)_{L\Lbar}
		+ \frac{1}{4}(Lh)_{LL}
		+ \frac{1}{4}(\Lbar h)_{LL}
		\right)
	\\
	&\phantom{=}
	- r \left(\chi_{(\text{small})}\right)_\alpha^{\phantom{\alpha}\beta} \slashed{\nabla}_\beta \log \mu
	\\ \\
	&=
	\left(\slashed{\D}(r\slashed{\nabla}h)\right)_{(\text{frame})}
	+ r\bm{\Gamma}_{(-1)}\bm{\Gamma}_{(-1, \text{small})}
	\end{split}
	\end{equation*}
	and
	\begin{equation*}
	\begin{split}
	\slashed{\D}_L \left( r(\chibar_{(\text{small})})_{\slashed{\mu}\slashed{\nu}} \right)
	&=
	r \slashed{\Pi}_{\slashed{\mu}}^{\phantom{\slashed{\mu}}\slashed{\rho}}\slashed{\Pi}_{\slashed{\nu}}^{\phantom{\slashed{\nu}}\slashed{\sigma}} R_{L\slashed{\rho}\Lbar\slashed{\sigma}}
	- \frac{1}{2}r \slashed{\nabla}_{\slashed{\mu}} \zeta_{\slashed{\nu}}
	- \frac{1}{2}r \slashed{\nabla}_{\slashed{\nu}} \zeta_{\slashed{\mu}}
	+ \frac{1}{2}r \zeta_{\slashed{\mu}} \zeta_{\slashed{\nu}}
	- r\omega \chibar_{\slashed{\mu}\slashed{\nu}}
	+ \chi_{\slashed{\mu}\slashed{\nu}}
	\\ \\
	&=
	\left(L(r\Lbar h)\right)_{(\text{frame})}
	+ \left(\slashed{\nabla}(r\slashed{\nabla}h)\right)_{(\text{frame})}
	+ \left(\slashed{\D}_L(r\slashed{\nabla}h)\right)_{(\text{frame})}
	+ \left(\slashed{\D}_{\Lbar}(r\slashed{\nabla}h)\right)_{(\text{frame})}
	\\
	&\phantom{=}
	+ r\bm{\Gamma}_{(-1)} \cdot \bm{\Gamma}_{(-1, \text{small})}
	\end{split}	
	\end{equation*}
	where in the second line we have expanded the Riemann tensor in terms of the rectangular derivatives of $h$, and also substituted for $\zeta$, $\omega$ and $\chibar$ using propositions \ref{proposition zeta}, \ref{proposition transport mu} and \ref{proposition chibar in terms of chi}.
	
	Note that the term $r R_{L\slashed{\mu}\Lbar\slashed{\nu}}$ is actually \emph{better} behaved than a term of the form $r R_{L\slashed{\mu}L\slashed{\nu}}$, since we have
	\begin{equation*}
	r(LLh)_{(\text{frame})} = \frac{1}{2}r(L T h)_{(\text{frame})} - r(L\Lbar h)_{(\text{frame})}
	\end{equation*}
	The second term can be estimated in a similar way to previously, but the first term behaves only like $r^{-\delta}$, rather than $r^{-1+\epsilon}$ as required.
	
	In summary, we have
	\begin{equation*}
	\slashed{\D}_L \left( r\left( \textit{l.o.t.}_{(\Omega, L\wedge\slashed{\Pi})} \right)\right)
	=
	\bm{\Gamma}_{(-1)} \cdot \begin{pmatrix}
		\left(\slashed{\D} (r\slashed{\nabla}h) \right)_{(\text{frame})} \\
		r(\tilde{\Box}_g h)_{(\text{frame})} \\
		r\bm{\Gamma}_{(-1)} \cdot \bm{\Gamma}_{(-1, \text{small})}
		\end{pmatrix}
	\end{equation*}
	where, possibly, additional factors of $r^{-1}$ might be present.

	Next, we examine the term
	\begin{equation*}
	\slashed{\D}_{\Lbar} \left( r\left( \textit{l.o.t.}_{(\Omega, \Lbar\wedge\slashed{\Pi})} \right)\right)
	\end{equation*}
	From proposition \ref{proposition expression for Omega} we have, schematically,
	\begin{equation*}
	\begin{split}
	\slashed{\D}_{\Lbar} \left( r\left( \textit{l.o.t.}_{(\Omega, \Lbar\wedge\slashed{\Pi})} \right)\right)
	&=
	\bm{\Gamma}_{(-1)} \left(\slashed{\D}_{\Lbar} (r\slashed{\nabla}h)\right)_{(\text{frame})}
	+ \bm{\Gamma}_{(-1)} \left(\Lbar (rLh)\right)_{(\text{frame})}
	+ \bm{\Gamma}_{(-1)} \left(\slashed{\D}_{\Lbar} \left(L^i_{(\text{small})} \slashed{\Pi}^i \right) \right)
	\\
	&\phantom{=}
	+ \bm{\Gamma}_{(-1, \text{small})} \left(\slashed{\D}_{\Lbar} (r\chi) \right)
	+ \bm{\Gamma}_{(-1, \text{small})} \left(\slashed{\D}_{\Lbar} \left(\slashed{\Pi}^i \otimes \slashed{\Pi}^i \right) \right)
	\end{split}
	\end{equation*}
	
	Now, we have, schematically,
	\begin{equation*}
	\left( \Lbar( rLh) \right)_{(\text{frame})}
	=
	(\tilde{\Box}_g h)_{(\text{frame})}
	+ r\bm{\Gamma}_{(-1)}(\partial h)_{(\text{frame})}
	\end{equation*}
	and
	\begin{equation*}
	\begin{split}
	\slashed{\D}_{\Lbar} \left( L^i_{(\text{small})} \slashed{\Pi}^i \right)
	&=
	r\bm{\Gamma}_{(-1)} \cdot \bm{\Gamma}_{(-1, \text{small})}
	\\
	\slashed{\D}_{\Lbar} \left(\slashed{\Pi}^i \otimes \slashed{\Pi}^i \right)
	&=
	r\bm{\Gamma}_{(-1)} \cdot \bm{\Gamma}_{(-1, \text{small})}
	\end{split}
	\end{equation*}
	
	Also, we note that
	\begin{equation*}
	\begin{split}
	\slashed{\D}_{\Lbar} (r\chi_{\slashed{\mu}\slashed{\nu}})
	&=
	\slashed{\D}_{\Lbar} \left(r(\chi_{(\text{small})})_{\slashed{\mu}\slashed{\nu}} + \slashed{g}_{\slashed{\mu}\slashed{\nu}} \right)
	\\
	&=
	r\slashed{\Pi}_{\slashed{\mu}}^{\phantom{\slashed{\mu}}\slashed{\rho}} \slashed{\Pi}_{\slashed{\nu}}^{\phantom{\slashed{\nu}}\slashed{\sigma}}
	R_{\Lbar\slashed{\rho}L\slashed{\sigma}}
	+ \frac{1}{2}r\slashed{\nabla}_{\slashed{\mu}} \zeta_{\slashed{\nu}}
	+ \frac{1}{2}r\slashed{\nabla}_{\slashed{\nu}} \zeta_{\slashed{\mu}}
	+ 2r\slashed{\nabla}^2_{\slashed{\mu}\slashed{\nu}} \log \mu
	+ \frac{1}{2}r \zeta_\mu \zeta_\nu
	\\
	&\phantom{=}
	+ \zeta_\mu r\slashed{\nabla}_\nu \log \mu
	+ \zeta_\nu r\slashed{\nabla}_\mu \log \mu
	+ 2r (\slashed{\nabla}_\mu \log \mu) (\slashed{\nabla}_\nu \log \mu)
	+ r\omega \chi_{\mu\nu}
	+ (\chi_{(\text{small})}))_{\slashed{\mu}\slashed{\nu}}
	\\
	&\phantom{=}
	- \frac{1}{2}r (\chi_{(\text{small})}))_{\slashed{\mu}}^{\phantom{\slashed{\mu}}\slashed{\sigma}} (\chibar_{(\text{small})}))_{\slashed{\sigma}\slashed{\nu}}
	- \frac{1}{2}r (\chibar_{(\text{small})}))_{\slashed{\mu}}^{\phantom{\slashed{\mu}}\slashed{\sigma}} (\chi_{(\text{small})}))_{\slashed{\sigma}\slashed{\nu}}
	\end{split}
	\end{equation*}
	where we have used proposition \ref{proposition Lbar chi}. Again, note the presence of the component of the Riemann tensor $R_{\Lbar \slashed{\rho} L \slashed{\sigma}}$, in which the highest order terms either involve at least one angular derivative, or a term of the form $(\Lbar L h)_{(\text{frame})}$, meaning that these terms have additional decay in $r$ as required. Indeed, we find that, schematically,
	\begin{equation*}
	\begin{split}
	\slashed{\D}_{\Lbar} (r\chi_{\slashed{\mu}\slashed{\nu}})
	=
	\left(\slashed{\D}(r\slashed{\nabla} h)\right)_{(\text{frame})}
	+ \left(\Lbar(r L h)\right)_{(\text{frame})}
	+ r^{-1} (r\slashed{\nabla})^2 \log \mu
	+ r\bm{\Gamma}_{(-1)} \cdot \bm{\Gamma}_{(0,\text{small})}
	\end{split}
	\end{equation*}
	Note that these error terms first show up when we commute for a \emph{second} time, i.e.\ when we are in the process of estimating \emph{third} derivatives of the metric. Thus we can esimate a term of the form $r^{-1} (r\slashed{\nabla})^2 \log \mu$ by using its transport equation in the $L$ direction without encountering any regularity issues.
	
	In summary, we have, schematically,
	\begin{equation*}
	\slashed{\D}_{\Lbar} \left(r\left( \textit{l.o.t.}_{(\Omega, \Lbar\wedge\slashed{\Pi})}\right)\right)
	=
	\bm{\Gamma}_{(-1)} \begin{pmatrix}
		\left(\slashed{\D}(r\slashed{\nabla}h)\right)_{(\text{frame})} \\
		r\left(\tilde{\Box}_g h\right)_{(\text{frame})} \\
		r^{-1} (r\slashed{\nabla})^2 \log \mu \\
		r\bm{\Gamma}_{(-1)}\cdot\bm{\Gamma}_{(0,\text{small})}
	\end{pmatrix}
	\end{equation*}

	Finally, we turn to the angular derivative term
	\begin{equation*}
	\slashed{\nabla} \left( r \cdot \left(\textit{l.o.t.}_{(\Omega,\slashed{\Pi}\wedge\slashed{\Pi})} \right)\right)
	\end{equation*}
	Substituting again for this term using proposition \ref{proposition expression for Omega} we find
	\begin{equation*}
	\slashed{\nabla} \left( r \cdot \left(\textit{l.o.t.}_{(\Omega,\slashed{\Pi}\wedge\slashed{\Pi})} \right)\right)
	=
	\bm{\Gamma}_{(-1)} \cdot \begin{pmatrix}
		\left(\overline{\slashed{\D}} (r\slashed{\nabla} h)\right)_{(\text{frame})} \\
		r\slashed{\nabla}\chi_{(\text{small})} \\
		r\bm{\Gamma}_{(-1)} \cdot \bm{\Gamma}_{(-1, \text{small})}	
	\end{pmatrix}
	\end{equation*}
	Again we note that $r\slashed{\nabla} \chi = r\slashed{\nabla}\chi_{(\text{small})}$ can be estimated by using its transport equation in the $L$ direction, which requires us to have control over third derivatives of the metric.

	Finally, note that the computation of the divergence of $rL^\nu \Omega_{\slashed{\alpha}\slashed{\beta} \mu\nu}$ can be calculated using very similar computations to the divergence of $T^\nu \Omega_{\slashed{\alpha}\slashed{\beta} \mu\nu}$.

\end{proof}

\chapter{Deformation tensor calculations}
\label{chapter deformation tensors}

In this chapter we will calculate the deformation tensors associated with various vector fields. These expressions will later be used both for performing energy estimates and estimating the error terms produced when commuting.

\section{Basic deformation tensor calculations and the multiplier vector fields}

\begin{definition}[Deformation tensors]
 We define the deformation tensor $^{(Z)}\pi$ associated with a vector field $Z$ as follows:
\begin{equation}
 ^{(Z)}\pi_{\mu\nu} := \D_\mu Z_\nu + \D_\nu Z_\mu
\end{equation}
\end{definition}

\begin{proposition}[The deformations tensor $^{(L)}\pi$ and $^{(\Lbar)}\pi$]
\label{proposition deformation L Lbar}
 In the region $r \geq r_0$ the deformation tensor associated with the null vector $L$ is given by
\begin{equation}
 \begin{split}
  ^{(L)}\pi_{\mu\nu} &= -\omega L_\mu L_\nu - \frac{1}{2}\omega L_\mu \Lbar_\nu - \frac{1}{2}\omega \Lbar_\mu L_\nu - L_\mu (\zeta_\nu + \mu^{-1}\slashed{\upd}_\nu \mu )
  \\
  & \phantom{=} - (\zeta_\mu + \mu^{-1}\slashed{\upd}_\mu \mu ) L_\nu + 2\chi_{\mu\nu}
 \end{split}
\end{equation}
while in the region $r < r_0$ the rectangular components of the deformation tensor are
\begin{equation}
  ^{(L)}\pi_{ab} = 2r \mathring{\gamma}_{ab} + L h_{ab} - 2L^c \partial_{(a}h_{b)c}
\end{equation}

% \begin{equation}
%  \begin{split}
%   ^{(L)}\pi_{\mu\nu} &= \frac{1}{2} \left(L^b \Lbar^c \Lbar_a \Gamma^a_{bc} \right) L_\mu L_\nu + \frac{1}{4}\left( L^b L^c \Lbar_a + \Lbar^b L^c L_a \right) \Gamma^a_{bc}  (L_\mu \Lbar_\nu + \Lbar_\mu L_\nu) \\
%   & \phantom{=} - \frac{1}{2} (\slashed{g}^{-1})^{AB} \left( L^b (X_A)^c \Lbar_a + \Lbar^b L^c (X_A)_a \right) \Gamma^a_{bc} (L_\mu (X_B)_\nu + (X_B)_\mu L_\nu ) \\
%   & \phantom{=} + \frac{1}{4} (L^b L^c L_a \Gamma^a_{bc}) \Lbar_\mu \Lbar_\nu - \frac{1}{2}(\slashed{g}^{-1})^{AB}( L^b (X_A)^c L_a + L^b L^c (X_A)_a )\Gamma^a_{bc} \left(L_\mu (X_B)_\nu + (X_B)_\mu L_\nu \right) \\
%   & \phantom{=} + (\slashed{g}^{-1})^{AC} (\slashed{g}^{-1})^{BD} \left( g((\slashed{\upd}_A L^a)\partial_a, X_B) + (X_A)^b L^c (X_B)_a \Gamma^a_{bc} \right) (X_C)_\mu (X_D)_\nu
%  \end{split}
% \end{equation}

For $r \geq r_0$ the deformation tensor associated with the null vector $\Lbar$ is
\begin{equation}
 \begin{split}
  ^{(\Lbar)}\pi_{\mu\nu} &= \frac{1}{2}\omega L_\mu \Lbar_\nu + \frac{1}{2}\omega \Lbar_\mu L_\nu -  \mu^{-1} L_\mu (\slashed{\upd}_\nu \mu) - \mu^{-1} (\slashed{\upd}_\mu \mu) L_\nu + \omega \Lbar_\mu \Lbar_\nu \\
  & \phantom{=} + \Lbar_\mu \zeta_\nu + \zeta_\mu \Lbar_\nu  + 2 \chibar_{\mu\nu}
 \end{split}
\end{equation}
\end{proposition}
while in the region $r < r_0$ the deformation tensor associated with $\Lbar$ has the rectangular components
\begin{equation}
  ^{(\Lbar)}\pi_{ab} = -2r \mathring{\gamma}_{ab} + \Lbar h_{ab} - 2\Lbar^c \partial_{(a}h_{b)c} 
\end{equation}

\begin{proof}
 We decompose the deformation tensors in the null frame and then contract with each vector field in the null frame. That is, we set
\begin{equation*}
 \begin{split}
  ^{(Z)}\pi &= {^{(Z)}\pi}^{LL}LL + {^{(Z)}\pi}^{L\Lbar}L\Lbar + {^{(Z)}\pi}^{\Lbar L} \Lbar L + {^{(Z)}\pi}^{LA} L X_A + {^{(Z)}\pi}^{AL} X_A L + {^{(Z)}\pi}^{\Lbar\Lbar}\Lbar\Lbar + {^{(Z)}\pi}^{\Lbar A} \Lbar X_A \\
  & \phantom{=} + {^{(Z)}\pi}^{A \Lbar} X_A \Lbar + {^{(Z)}\pi}^{AB} X_A X_B
 \end{split}
\end{equation*}
where, for example $^{(Z)}\pi^{LA} = -\frac{1}{2}(\slashed{g}^{-1})^{AB}\left(g(\D_B Z, \Lbar) + g(\D_{\Lbar} Z, X_B) \right)$. We then use proposition \ref{proposition null connection} to compute each one of these inner products.

In the region $r < r_0$, we make use of the identity
\begin{equation}
 ^{(Z)}\pi_{ab} = \partial_a Z_b + \partial_b Z_a + 2\Gamma^c_{ab}Z_c
\end{equation}
as well as the commutation identities of the null frame in the region $r < r_0$.

\end{proof}

We now introduce the \emph{multiplier vector fields} and compute their associated deformation tensors. Before doing this, we first define the vector field $T$, which will be used in the definition of the multiplier vector fields.

\begin{definition}[The vector field $T$]
 We define the vector field $T$ as
\begin{equation}
 T:= \frac{1}{2}(L + \Lbar)
\end{equation}
\end{definition}

\begin{definition}[Multiplier vector fields]
 We define the set of multiplier vector fields, $\mathcal{Z}$ as follows:
\begin{equation}
 \mathcal{Z} := \left\{ wT\ ,\quad w f_R(r)R\ ,\quad f_L(r) r^p L\ \right\}
\end{equation}
for $w$, $f_R$, and $f_L$ functions of $r$ to be determined below, and for $p$ a positive constant. We will sometimes refer to the function $w$ as a ``weight function''.
\end{definition}

\begin{remark}[$\mu$ Weighted multipliers]
 It would be possible to use modified, ``weighted'' multiplier vector fields which incorperate powers of $\mu$ in order to eliminate the terms of the form $\omega (\Lbar \phi)$ from the energy estimates. Specifically, if we define
 \begin{equation*}
  \begin{split}
   \check{T} &:= \frac{1}{2}\left( L + \mu \Lbar \right) \\
   \check{R} &:= \frac{1}{2}\left(L - \mu \Lbar \right)
  \end{split}
 \end{equation*}
 then the worst terms $\frac{1}{2}\omega \Lbar_\mu \Lbar_\nu$ in the deformation tensors are eliminated, and the corresponding error terms in the energy estimates will also be removed.

 Note, however, that there are additional error terms in the energy estimates which arise from the semilinear structure, and these can be comparable to the errors removed in this way. For example, let $\phi$ be a field satisfying the ``reduced wave equation''
 \begin{equation*}
  \Box_g \phi + \omega \Lbar \phi = 0
 \end{equation*}
 Then the additional error terms are precisely of the same form as those ``removed'' by the use of the weighted multipliers. In order to remove these error terms, we should instead use the weighted multipliers of the form
 \begin{equation*}
  \check{T}^{(-1)} := \frac{1}{2}\left( L + \mu^{-1} \Lbar \right)
 \end{equation*}
 
 More generally, defining
 \begin{equation*}
  \check{T}^{(s)} := \frac{1}{2}\left( L + \mu^{s} \Lbar \right)
 \end{equation*}
 and choosing $s$ appropriately, we can remove error terms of the form $\omega (\Lbar\phi)^2$ in the energy estimates. Note, however, that this stratagy can only be used to remove error terms of the form $(\Lbar h)_{LL} (\Lbar \phi)^2$ in the energy estimates; error terms involving other fields, such has $(\Lbar \phi_{(1)})(\Lbar \phi)^2$ where $\phi_{(1)} \neq h_{LL}$ cannot be removed in the same way. Note also that modifying the multiplier vector fields in this way produces additional error terms in the energy estimates, of which the most dangerous has the form
 \begin{equation*}
  \mu^{-1}(T\mu)|\slashed{\nabla}\phi|^2
 \end{equation*}
 If we were able to show that $\mu$ is bounded in $r$ (as in the case of the Einstein equations), and if we also know that $T\mu$ decays at a rate which is integrable in $\tau$ (which can possibly be shown using the improved decay estimates for $T$ derivatives, in appendix \ref{appendix improved energy decay}) then we can handle such terms. However, as mentioned in the introduction we wish to close our estimates without making use of the improved decay estimates, and we also cannot generally show that $T\mu$ is bounded in $r$; it might grow like $r^\epsilon$. For these reasons we will not use the modified vector fields.
\end{remark}

\begin{proposition}[Deformation tensors of the multiplier vector fields]
\label{proposition deformation multipliers}
 The deformation tensors associated with the multiplier vector fields $Z \in \mathcal{Z}$ are as follows:

In the region $r \geq r_0$, the deformation tensor of the weighted vector field $wT$ is
\begin{equation}
  \begin{split}
    ^{(wT)}\pi_{\mu\nu} &= \frac{1}{2}(w\omega - w') \Lbar_\mu \Lbar_\nu -\frac{1}{2}(w\omega - w') L_\mu L_\nu 
    - \frac{1}{2}w(\slashed{g}^{-1})^{AB}\left(\zeta_A + 2\mu^{-1}\slashed{\upd}_A \mu \right)\left( L_\mu (X_B)_\nu + (X_B)_\mu L_\nu \right) \\
    &\phantom{=}+ \frac{1}{2}w(\slashed{g}^{-1})^{AB}\zeta_A \left(\Lbar_\mu (X_B)_\nu + (X_B)_\mu \Lbar_\nu \right) + w\left( (\chi_{(\text{small})})_{\mu\nu} + (\chibar_{(\text{small})})_{\mu\nu} \right)
  \end{split}
\end{equation}

while in the region $r \leq r_0$ the deformation tensor of the vector field $wT$, expressed in rectangular indices, is
\begin{equation}
 ^{(wT)}\pi_{ab} = wT h_{ab} - 2wT^c \partial_{(a}h_{b)c} + \frac{1}{2}w'(\partial_a r)(L_b + \Lbar_b) + \frac{1}{2}w'(\partial_b r)(L_a + \Lbar_a)
\end{equation}

In the region $r \geq r_0$ the deformation tensor of the vector field $w f_R(r)R$ (the ``weighted Morawetz vector field'') is
\begin{equation}
 \begin{split}
  ^{(w f_R R)}\pi_{\mu\nu} &= \frac{1}{2}(wf'_R + w' f_R - \omega w f_R)L_\mu L_\nu
  -\frac{1}{2}\left( wf'_R + w' f_R + \omega w f_R \right)\left( L_\mu \Lbar_\nu + \Lbar_\mu L_\nu \right) \\
  &\phantom{=}
  + \frac{1}{2}\left( wf'_R + w' f_R - \omega wf_R \right) \Lbar_\mu \Lbar_\nu 
  - \frac{1}{2}\zeta^A wf_R \left( L_\mu (X_A)_\nu + (X_A)_\mu L_\nu \right) \\
  &\phantom{=}
  - \frac{1}{2}\zeta^A wf_R \left( \Lbar_\mu (X_A)_\nu + (X_A)_\mu \Lbar_\nu \right)
  + wf_R\left(\chi_{\mu\nu} - \chibar_{\mu\nu}\right)
 \end{split}
\end{equation}

while in the region $r \leq r_0$ the associated deformation tensor is
\begin{equation}
 \begin{split}
  ^{(wf_R R)}\pi_{ab} &= \frac{1}{2}(w f_R' + w' f_R) \left( L_a - \Lbar_a\right) \left( L_b - \Lbar_b \right) + wf_R r\mathring{\gamma}_{ab} - \frac{1}{2} wf_R \left( Lh_{ab} - \Lbar h_{ab}\right) - wf_R R^c \partial_{(a}h_{b)c}
 \end{split}
\end{equation}

In the region $r \geq r_0$ the deformation tensor of the vector field $ f_L(r) r^p L$ is
\begin{equation}
 \begin{split}
  ^{(f_L r^p L)}\pi_{\mu\nu} &= \left(p r^{p-1} f_L + r^p f'_L -\omega r^p f_L \right) L_\mu L_\nu \\
  & \phantom{=} - \frac{1}{2}\left(pr^{p-1} f_L + r^p f'_L + \omega r^p f_L \right)\left( L_\mu \Lbar_\nu + \Lbar_\mu L_\nu \right) \\
  & \phantom{=} - r^p f_L (\slashed{g}^{-1})^{AB}(\zeta_A + 2\mu^{-1}\slashed{\upd}_A \mu )\left( L_\mu (X_B)_\nu + (X_B)_\mu L_\nu \right) \\
  & \phantom{=} + 2r^p f_L\chi_{\mu\nu}
 \end{split}
\end{equation}

\end{proposition}

\begin{proof}
 These computations are similar to those in proposition \ref{proposition deformation L Lbar}. We also make use of the fact that 
 \begin{equation*}
  ^{(f(r)Z)}\pi_{\mu\nu} = f(r) ^{(Z)}\pi_{\mu\nu} + f'(r)\frac{1}{2}(L_\mu - \Lbar_\mu) Z_\nu + f'(r)\frac{1}{2} Z_\mu (L_\nu - \Lbar_\nu) 
 \end{equation*}
 for any function $f(r)$ and vector field $Z$.
\end{proof}

\begin{remark}
 We will later pick the function $f_L(r)$ to be supported only in the range $r \geq r_0$, so we will not need any expressions relating to the vector field $f_L r^p L$ in the region $r < r_0$.
\end{remark}

\section{Energy momentum tensors and compatible currents}

\begin{definition}[The energy momentum tensor]
 For any $S_{\tau,r}$-tangent tensor field $\phi$, we define the associated energy momentum tensor $Q$ as
 \begin{equation}
\label{equation definition energy momentum}
  Q_{\mu\nu}[\phi] := (\slashed{\D}_\mu \phi)\cdot (\slashed{\D}_\nu \phi) - \frac{1}{2}g_{\mu\nu} (g^{-1})^{\alpha \beta}(\slashed{\D}_\alpha \phi)\cdot(\slashed{\D}_\beta \phi)
 \end{equation}
\end{definition}

\begin{definition}[The energy momentum tensors after a point-dependent change of basis for the fields]
 As mentioned above, we will also need to consider a point-dependent change of basis for the fields $\phi$, and the associated energy momentum tensors require some modification.

 Let $\phi_{(a)}$ be a collection of $N$ $S_{\tau, r}$-tangent tensor fields labelled by the index $(a)$. Let $M_{(A)}^{\phantom{(A)}(a)}(x)$ be a point-dependent change of basis matrix, that is, a collection of scalar fields such that, at every point $x \in \mathcal{M}$, the matrix $M_{(A)}^{(a)}$ has rank $N$. Define the fields
\begin{equation}
 \phi_{(A)} := M_{(A)}^{(a)}\phi_{(a)}
\end{equation}
where, as usual, repeated pairs of indices are summed over if one indix appears in the raised position and the other appears in the lowered position. Associated with the fields $\phi_{(A)}$ we define the energy momentum tensors
\begin{equation}
 (Q_{(A)})_{\mu\nu}[\phi] := M_{(A)}^{\phantom{(A)}(a)} M_{(A)}^{\phantom{(A)}(b)} \left( (\slashed{\D}_\mu \phi_{(a)})\cdot (\slashed{\D}_\nu \phi_{(b)}) - \frac{1}{2}g_{\mu\nu} (g^{-1})^{\alpha \beta}(\slashed{\D}_\alpha \phi_{(a)})\cdot(\slashed{\D}_\beta \phi_{(b)}) \right) 
\end{equation}
\end{definition}

\begin{remark}[The paradigm case for a point-dependent change of basis]
 When analysing the Einstein equations in harmonic coordinates, in order to make the semilinear structure apparent we need to change from rectangular coordinates to null frame coordinates. That is, we begin by considering the fields $h_{ab}$, labelled by the symmetric pair of rectangular indices $(ab)$ since these are the fields which satisfy a set of nonlinear wave equations. However, in order to make use of the weak null structure we must change to fields labelled by null frame indices. In other words, we take the matrices $M$ to be
\begin{equation*}
 \begin{split}
  M_{(LL)}^{\phantom{LL}(ab)} &:= L^a L^b \\
  M_{(L\Lbar)}^{\phantom{L\Lbar}(ab)} &:= L^a \Lbar^b \\
  M_{(La)}^{\phantom{La}(bc)} &:= L^b \slashed{\Pi}_a^{\phantom{a}c} \\
  M_{(\Lbar\Lbar)}^{\phantom{\Lbar\Lbar}(ab)} &:= \Lbar^a \Lbar^b \\
  M_{(\Lbar a)}^{\phantom{\Lbar a}(bc)} &:= \Lbar^b \slashed{\Pi}_a^{\phantom{a}c} \\
  M_{(ab)}^{\phantom{ab}(cd)} &:= \slashed{\Pi}_a^{\phantom{a}c}\slashed{\Pi}_b^{\phantom{b}d} \\
 \end{split}
\end{equation*}
 Since the rectangular components of the null frame fields $L$, $\Lbar$ and $\slashed{\Pi}$ are functions of the metric perturbations $h_{ab}$, which are themselves scalar fields on $\mathcal{M}$, this is a point-dependent change of basis. It is important to remember that, from the point of view of regularity, the scalar fields $M$ are the same order as the metric perturbations $h$.
\end{remark}

\begin{proposition}[The divergence of the energy momentum tensor]
 \label{proposition divergence energy momentum}
  Let $\phi_{\sigma_1 \ldots \sigma_n}$ be a rank $(0,n)$ tensor. Then the associated energy momentum tensor satisfies the following equation:
 \begin{equation}
    \D_\nu Q^{\nu}_{\phantom{\nu}\mu}[\phi] = \left( \slashed{\Box}_g \phi \right) \cdot \slashed{\D}_\mu \phi + (\slashed{\D}^\nu \phi)\cdot\left( [\slashed{\D}_\mu \, , \, \slashed{\D}_\nu]\phi\right)
 \end{equation}
\end{proposition}

\begin{proof}
Taking the divergence of equation \eqref{equation definition energy momentum} proves the proposition. Note that the commutator term vanishes in the case that $\phi$ is a scalar field. Note also that this term can be expressed in terms of the curvature $\Omega_{\mu\nu}$.
\end{proof}

\begin{proposition}[The divergence of the energy momentum tensors after a point-dependent change of basis]
\label{proposition divergence of point dependent energy momentum}
 We borrow notation from the notation already developped for the null components of derivatives of the rectangular components of the metric perturbations ($(\partial h)_{LL}$ for example), in which indices outside of parentheses are contracted \emph{after} differential operators are first applied to the rectangular components. Thus, for example, for a set of scalar fields $\phi_{(a)}$ we write
\begin{equation*}
 (L\phi)_{(A)} := M_{(A)}^{(a)}(L\phi_{(a)})
\end{equation*}
% We also use the notation
% \begin{equation}
%  |\phi|_{(\text{frame})} := \sqrt{ \phi_{LL}^2 + \phi_{L\Lbar}^2 + \phi_{\Lbar\Lbar}^2 + |\slashed{\Pi}_\mu^{\phantom{\mu}a} L^b \phi_{ab} |^2 + |\slashed{\Pi}_\mu^{\phantom{\mu}a} \Lbar^b \phi_{ab} |^2 + |\slashed{\Pi}_\mu^{\phantom{\mu}a} \slashed{\Pi}_\nu^b \phi_{ab} |^2  }
% \end{equation}
% and similar schematic notation $(\phi)_{(\text{frame})}$ to refer to the collection of null frame decomposed components of $\phi$.

Then, for a set of $S_{\tau, r}$-tangent tensor field $\phi_{(a)}$, we have

\begin{equation}
 \begin{split}
  Z^\mu \D_\nu & (Q_{(A)})^\nu_{\phantom{\nu}\mu}[\phi] =
  \left( \slashed{\Box}_g \phi \right)_{(A)} \cdot \left(\slashed{\D}_Z \phi\right)_{(A)}
  + (\slashed{\D}^\nu \phi)_{(A)}\cdot\left( [\slashed{\D}_\mu \, , \, \slashed{\D}_\nu]\phi\right)_{(A)}\\
  &
  + \frac{1}{2}g(Z,L) \left(L M_{(A)}^{\phantom{(A)}(a)} \right) (\slashed{\D}_{\Lbar}\phi_{(a)}) \cdot (\slashed{\D}_{\Lbar}\phi)_{(A)}
  + \frac{1}{2}g(Z,\Lbar) \left(L M_{(A)}^{\phantom{(A)}(a)} \right) (\slashed{\nabla}\phi_{(a)}) \cdot (\slashed{\nabla}\phi)_{(A)}\\
  &
  - \frac{1}{2} \left(L M_{(A)}^{\phantom{(A)}(a)} \right) (\slashed{\D}_{\Lbar}\phi_{(a)}) \cdot (\slashed{\nabla}_{\slashed{\Pi}(Z)}\phi)_{(A)}
  - \frac{1}{2} \left(L M_{(A)}^{\phantom{(A)}(a)} \right) (\slashed{\nabla}_{\slashed{\Pi}(Z)}\phi_{(a)}) \cdot (\slashed{\D}_{\Lbar} \phi)_{(A)}\\
  &
  + \frac{1}{2}g(Z,\Lbar) \left(\Lbar M_{(A)}^{\phantom{(A)}(a)} \right) (\slashed{\D}_{L}\phi_{(a)}) \cdot (\slashed{\D}_{L}\phi)_{(A)}
  + \frac{1}{2}g(Z,L) \left(\Lbar M_{(A)}^{\phantom{(A)}(a)} \right) (\slashed{\nabla}\phi_{(a)}) \cdot (\slashed{\nabla}\phi)_{(A)}\\
  &
  - \frac{1}{2} \left(\Lbar M_{(A)}^{\phantom{(A)}(a)} \right) (\slashed{\D}_L \phi_{(a)}) \cdot (\slashed{\nabla}_{\slashed{\Pi}(Z)}\phi)_{(A)}
  - \frac{1}{2} \left(\Lbar M_{(A)}^{\phantom{(A)}(a)} \right) (\slashed{\nabla}_{\slashed{\Pi}(Z)}\phi_{(a)}) \cdot (\slashed{\D}_L \phi)_{(A)}\\
  & 
  - \frac{1}{2} g(Z,\Lbar) \left(\slashed{\nabla}^\alpha M_{(A)}^{\phantom{(A)}(a)} \right) (\slashed{\nabla}_\alpha \phi_{(a)}) \cdot (\slashed{\D}_L \phi)_{(A)}
  - \frac{1}{2} g(Z,\Lbar) \left(\slashed{\nabla}^\alpha M_{(A)}^{\phantom{(A)}(a)} \right) (\slashed{\D}_L \phi_{(a)}) \cdot (\slashed{\nabla}_\alpha \phi)_{(A)}\\
  &
  - \frac{1}{2} g(Z,L) \left(\slashed{\nabla}^\alpha M_{(A)}^{\phantom{(A)}(a)} \right) (\slashed{\nabla}_\alpha \phi_{(a)}) \cdot (\slashed{\D}_{\Lbar} \phi)_{(A)}
  - \frac{1}{2} g(Z,L) \left(\slashed{\nabla}^\alpha M_{(A)}^{\phantom{(A)}(a)} \right) (\slashed{\D}_{\Lbar} \phi_{(a)}) \cdot (\slashed{\nabla}_\alpha \phi)_{(A)}\\
  &
  + \left(\slashed{\nabla}^\alpha M_{(A)}^{\phantom{(A)}(a)} \right) (\slashed{\nabla}_\alpha \phi_{(a)}) \cdot (\slashed{\nabla}_{\slashed{\Pi}(Z)}\phi)_{(A)}
  + \left(\slashed{\nabla}^\alpha M_{(A)}^{\phantom{(A)}(a)} \right) (\slashed{\nabla}_{\slashed{\Pi}(Z)}\phi_{(a)}) \cdot (\slashed{\nabla}_\alpha \phi)_{(A)} \\
  &
  + \frac{1}{2} \left( \slashed{\nabla}_{\slashed{\Pi}(Z)} M_{(A)}^{\phantom{(A)}(a)} \right) (\slashed{\D}_L \phi_{(a)})\cdot (\slashed{\D}_{\Lbar} \phi)_{(A)}
  + \frac{1}{2} \left( \slashed{\nabla}_{\slashed{\Pi}(Z)} M_{(A)}^{\phantom{(A)}(a)} \right) (\slashed{\D}_{\Lbar} \phi_{(a)})\cdot (\slashed{\D}_L \phi)_{(A)} \\
  &
  - \left( \slashed{\nabla}_{\slashed{\Pi}(Z)} M_{(A)}^{\phantom{(A)}(a)} \right) (\slashed{\nabla} \phi_{(a)})\cdot (\slashed{\nabla} \phi)_{(A)}
 \end{split}
\end{equation}

\end{proposition}

\begin{remark}[Further remarks on proposition \ref{proposition divergence of point dependent energy momentum}]
There are a large number of additional error terms appearing in the formulae given in proposition \ref{proposition divergence of point dependent energy momentum}, each of which involves a derivative of the change of basis matrix $\partial M$, multiplied by a derivative of the field in the original basis $\slashed{\D}\phi_{(a)}$, and then multiplied by a derivative of the field in the new basis $(\slashed{\D}\phi)_{(A)}$.

The important point about these error terms is that not all of them are necessarily ``small'', i.e.\ they do not all need to come with a factor of $\epsilon$. This is because we can seperately estimate the terms involving the field in the original basis, and indeed we can perform the energy estimates in such a way that the energy of these terms is already suitably small.

This is best illustrated with an example, for which we shall turn to the paradigmatic example of fields $\phi_{(ab)}$ labelled by rectangular indices $(ab)$. let $X$, $Y \in \{ L, \Lbar, \slashed{\Pi}_a \}$. Then, when performing energy estimates associated with the fields $\phi_{XY}$, and using the multiplier vector field $Z$, we will encounter additional error terms of the form
\begin{equation*}
 \int_{\mathcal{M}} \left( \frac{1}{r} (\slashed{\nabla} \phi)_{XY} \cdot (\slashed{\D}_Z \phi)_{(\text{frame})} + \frac{1}{r} (\slashed{\nabla} \phi)_{(\text{frame})} \cdot (\slashed{\D}_Z \phi)_{XY} \right) \dVol_g 
\end{equation*}

The idea we use to estimate these terms is to exploit the fact that at least one of the derivatives involved is a \emph{good} derivative. Moreover, since the rectangular components of the null frame are of order one, we have
\begin{equation*}
 (\partial \phi)_{(\text{frame})} \sim (\partial \phi)_{(\text{rect})}
\end{equation*}
The rectangular components of $\phi$ are themselves functions of the fields $\phi_{(A)}$, and so can be assumed to behave similarly to the worst of these fields, i.e.\ the fields at the highest level of the semilinear hierarchy. Hence, if $Z$ is a multiplier of order unity (for example, if $Z = T$) then we estimate
\begin{equation*}
 \begin{split}
  &\int_{\mathcal{M}} \left( \frac{1}{r} (\slashed{\nabla} \phi)_{XY} \cdot (\slashed{\D} \phi)_{(\text{frame})} + \frac{1}{r} (\slashed{\nabla} \phi)_{(\text{frame})} \cdot (\slashed{\D} \phi)_{XY} \right) \dVol_g \\
  &\lesssim
  \int_{\mathcal{M}} \left( \epsilon r^{-1 + \delta}|\slashed{\nabla} \phi|^2_{XY} + \epsilon^{-1} r^{-1-\delta}|\slashed{\D} \phi|_{(\text{frame})} + \epsilon^{-1} r^{-1 + \delta} |\slashed{\nabla} \phi|^2_{(\text{frame})} + \epsilon r^{-1-\delta} |\slashed{\D} \phi|^2_{XY} \right) \dVol_g
 \end{split}
\end{equation*}
Assuming that we have \emph{already} improved the energy estimates at the highest level of the hierarchy, we can estimate the terms involving $\phi_{XY}$.

On the other hand, if we are doing the $p$-weighted energy estimates, where $Z = r^p L$, we instead estimate
\begin{equation*}
 \begin{split}
  &\int_{\mathcal{M}} \left( r^{p-1} (\slashed{\nabla} \phi)_{XY} \cdot (\slashed{\D}_L \phi)_{(\text{frame})} + r^{p-1}(\slashed{\nabla} \phi)_{(\text{frame})} \cdot (\slashed{\D}_L \phi)_{XY} \right) \dVol_g \\
  &\lesssim
  \int_{\mathcal{M}} \left( \epsilon r^{p-1}|\slashed{\nabla} \phi|^2_{XY} + \epsilon^{-1} r^{p-1}|\slashed{\D}_L \phi|_{(\text{frame})} + \epsilon^{-1} r^{p-1 - \delta} \tau^{2\delta} |\slashed{\nabla} \phi|^2_{(\text{frame})} + \epsilon r^{p-1+\delta}\tau^{-2\delta} |\slashed{\D}_L \phi|^2_{XY} \right) \dVol_g
 \end{split}
\end{equation*}
Importantly, not only can these terms now be estimated, but these error terms are consistent with a taking a larger value for $p$ for some of the null-decomposed components of $\phi$. In particular, it may be possible, in certain cases, to take $p > 1$ for certain null decomposed components of the equations. Whether is in fact possible will depend on the structure of the inhomogeneous term $F_{XY}$.

\end{remark}

\begin{definition}[Compatible currents]
 Associated to a scalar field $\phi$ and a vector field $Z$ we have the following two ``currents'':
 \begin{equation}
  \begin{split}
   ^{(Z)}J^\mu[\phi] &:= Z_\nu Q^{\mu\nu}[\phi] \\
   ^{(Z)}K[\phi] &:= \frac{1}{2}{^{(Z)}\pi}^{\mu\nu} Q_{\mu\nu}[\phi]
  \end{split}
 \end{equation}
\end{definition}
 
\begin{proposition}[Compatible current identity]
\label{proposition compatible current identity}
 The divergence of the current $^{(Z)}J$ satisfies the following identity:
 \begin{equation}
  \begin{split}
    \D_\mu {^{(Z)}J^\mu}[\phi] &= {^{(Z)}K}[\phi] - \omega (\slashed{\D}_{\Lbar}\phi)\cdot(\slashed{\D}_Z \phi) + \tilde{\slashed{\Box}}_g \phi \cdot (\slashed{\D}_Z \phi) + \left( [\slashed{\D}_\mu \, , \, \slashed{\D}_\nu]\phi \right) \cdot (\slashed{\D}^\mu \phi) Z^\nu 
  \end{split}
 \end{equation}
\end{proposition}
\begin{proof}
 This is a straightforward application of the previous proposition.
\end{proof}

\begin{definition}[Modified compatible currents]
\label{definition modified currents}
 Given an $S_{\tau,r}$-tangent tensor field $\phi$, a vector field $Z$ and a function $f_Z$ we associate the following modified compatible currents:
 \begin{equation}
  \begin{split}
   ^{(Z,\, f_Z)}\tilde{J}^\mu[\phi] &:= {^{(Z)}J}^\mu[\phi] + \frac{1}{2}f_Z (g^{-1})^{\mu\nu} \phi \cdot (\slashed{\D}_\nu \phi) - \frac{1}{4}(g^{-1})^{\mu\nu}(\D_\nu f_Z)|\phi|^2 \\
   ^{(Z,\, f_Z)}\tilde{K}[\phi] &:= {^{(Z)}K}[\phi] + \frac{1}{2}f_Z (g^{-1})^{\mu\nu} \slashed{\D}_\mu \phi \cdot \slashed{\D}_\nu \phi - \frac{1}{4}(\Box_g f_Z) |\phi|^2
  \end{split}
 \end{equation}
\end{definition}

\begin{proposition}[Modified compatible current identity]
\label{proposition modified compatible current identity}
 An easy calculation yields the following identity:
 \begin{equation}
  \begin{split}
    \D_\mu \left( {^{(Z,\, f_Z)}\tilde{J}}^\mu[\phi] \right) &= {^{(Z,\, f_Z)}\tilde{K}}[\phi] - \omega (\slashed{\D}_{\Lbar}\phi)\cdot\left( (\slashed{\D}_Z \phi) + \frac{1}{2}f_Z \phi \right)
    + \tilde{\slashed{\Box}}_g \phi \cdot \left((\slashed{\D}_Z\phi) + \frac{1}{2}f_Z \phi \right) \\
    &\phantom{=} + \left( [\slashed{\D}_\mu \, , \, \slashed{\D}_\nu]\phi \right) \cdot (\slashed{\D}^\mu \phi) Z^\nu
  \end{split}
 \end{equation}
\end{proposition}

\begin{proposition}[Modified compatible current identity after a point-dependent change of basis]
 We can also calculate
 \begin{equation}
  \begin{split}
    \D_\mu \left( {^{(Z,\, f_Z)}\tilde{J}}^\mu[\phi]_{(A)} \right) &= 
    {^{(Z,\, f_Z)}\tilde{K}}[\phi]_{(A)} - \omega (\slashed{\D}_{\Lbar}\phi)_{(A)}\cdot\left( (\slashed{\D}_Z \phi)_{(A)} + \frac{1}{2}f_Z \phi_{(A)} \right)\\
    &\phantom{=} + (\tilde{\slashed{\Box}}_g \phi)_{(A)} \cdot \left((\slashed{\D}_Z\phi)_{(A)} + \frac{1}{2}f_Z \phi_{(A)} \right)
    + \left( [\slashed{\D}_\mu \, , \, \slashed{\D}_\nu]\phi \right)_{(A)} \cdot (\slashed{\D}^\mu \phi)_{(A)} Z^\nu \\
    &\phantom{=}
    - \frac{1}{4} (L M_{(A)}^{\phantom{(A)}(a)}) \left( f_z\phi_{(a)} \cdot (\slashed{\D}_{\Lbar} \phi)_{(A)} + f_z\phi_{(A)} \cdot (\slashed{\D}_{\Lbar} \phi_{(a)}) - (\Lbar f_Z) \phi_{(A)}\cdot \phi_{(a)} \right) \\
    &\phantom{=}
    - \frac{1}{4} (\Lbar M_{(A)}^{\phantom{(A)}(a)}) \left( f_z\phi_{(a)} \cdot (\slashed{\D}_{L} \phi)_{(A)} + f_z\phi_{(A)} \cdot (\slashed{\D}_{L} \phi_{(a)}) - (L f_Z) \phi_{(A)}\cdot \phi_{(a)} \right) \\
    &\phantom{=}
    + \frac{1}{2} (\slashed{\nabla}^\alpha M_{(A)}^{\phantom{(A)}(a)}) \left( f_z\phi_{(a)} \cdot (\slashed{\nabla}_{\alpha} \phi)_{(A)} + f_z\phi_{(A)} \cdot (\slashed{\nabla}_{\alpha} \phi_{(a)}) - (\slashed{\nabla}_{\alpha} f_Z) \phi_{(A)}\cdot \phi_{(a)} \right)
  \end{split}
 \end{equation}
\end{proposition}

\section{The (modified) compatible multiplier currents}
To each vector field $Z \in \mathcal{Z}$ we shall associate (modified) compatible currents.

\begin{definition}[The weighted $T$-energy currents]
 To the weighted vector field $wT$ we associate the following (unmodified) compatible currents:
 \begin{equation}
  \begin{split}
   ^{(wT)}J^\mu[\phi] &= wT_\nu Q^{\mu\nu}[\phi] \\
   ^{(wT)}K[\phi] &= \frac{1}{2}{^{(wT)}\pi}^{\mu\nu} Q_{\mu\nu}[\phi] \\
  \end{split}
 \end{equation}
\end{definition}

\begin{proposition}[An expression for the bulk weighted $T$-energy current term $^{(wT)}K$]
\label{proposition bulk weighted T current}
The bulk weighted $T$-energy current ${^{(wT)}K}$ can be decomposed as follows:
 \begin{equation}
  \begin{split}
   {^{(wT)}K}[\phi] - w\omega (\slashed{\D}_{\Lbar} \phi)\cdot (\slashed{\D}_T \phi) := -\frac{1}{2}w' |\slashed{\D}_{\Lbar}\phi|^2 + w\textit{Err}_{(wT,\, \text{bulk})}[\phi]
  \end{split}
 \end{equation}
 where, in the region $r \geq r_0$,
\begin{equation}
 \begin{split}
  \textit{Err}_{(wT,\, \text{bulk})}[\phi] &=
  \frac{1}{2}\partial_r(\log w)|\slashed{\D}_L \phi|^2
  -\frac{1}{4}\omega |\slashed{\D}_{\Lbar} \phi|^2
  - \frac{1}{4}\omega |\slashed{\D}_L\phi|^2 \\
  &\phantom{=}
  + \frac{1}{4}\left( \tr_{\slashed{g}}\chi_{(\text{small})} + \tr_{\slashed{g}}\chibar_{(\text{small})} -2\omega \right) (\slashed{\D}_L\phi)\cdot(\slashed{\D}_{\Lbar}\phi) \\
  &\phantom{=} 
  - \frac{1}{2}(\zeta^\alpha + 2\slashed{\nabla}^\alpha \log\mu )(\slashed{\D}_L \phi)\cdot(\slashed{\nabla}_\alpha \phi)
  + \frac{1}{2} \zeta^\alpha (\slashed{\D}_{\Lbar} \phi) \cdot (\slashed{\nabla}_\alpha \phi) \\
  &\phantom{=} + \frac{1}{2}\left( (\hat{\chi})^{\alpha\beta} + (\hat{\chibar})^{\alpha\beta} \right) (\slashed{\nabla}_\alpha \phi )\cdot(\slashed{\nabla}_\beta \phi)
 \end{split}
\end{equation}
which can be expressed schematically as
\begin{equation}
 \begin{split}
  \textit{Err}_{(wT,\, \text{bulk})}[\phi] &= 
  \partial_r (\log w) |\overline{\slashed{\D}}\phi|^2
  +\omega (\slashed{\D}\phi)^2
  + \bm{\Gamma} (\slashed{\D}\phi)\cdot(\overline{\slashed{\D}}\phi) 
 \end{split}
\end{equation}

while in the region $r < r_0$, using the fact that $T^0 = 1$ and $T^i = 0$, we have
\begin{equation}
  \left|\textit{Err}_{(T,\, \text{bulk})}[\phi]\right| \lesssim 
  |\partial_r(\log w)| |\slashed{\D} \phi|^2
  + |\partial h|_{(\text{rect})} |g|_{(\text{rect})}^2 |\slashed{\D} \phi|^2 
\end{equation}

\end{proposition}

\begin{proof}
This follows easily from the expressions given in \ref{proposition deformation multipliers}. In order to obtain the expressions for the error terms, we also need to use the expressions for $\omega$, $\zeta$ and $\chibar$ given in propositions \ref{proposition transport mu}, \ref{proposition zeta} and \ref{proposition chibar in terms of chi}, as well as the equations
\begin{equation*}
 \begin{split}
  1 - L^i L^i &= -2 L^i L^i_{(\text{small})} - L^i_{(\text{small})}L^i_{(\text{small})} \\
  L^i \slashed{\upd}_A x^i &= L^i_{(\text{small})} \slashed{\upd}_A x^i
 \end{split}
\end{equation*}
where the second equality follows from the fact that 
\begin{equation*}
 x^i \slashed{\upd}_A x^i = \frac{1}{2}\slashed{\upd}_A (x^i x^i) = r \slashed{\upd}_A r = 0 
\end{equation*}

\end{proof}

\begin{definition}[The weighted Morawetz currents]
 We define the weighted Morawetz currents as the modified compatible currents ${^{(Z,\, f_Z)}\tilde{J}}$ and ${^{(Z,\, f_Z)}\tilde{K}}$ with the choices $Z = w f_R R$ and $f_Z = 2w r^{-1} f_R$. We abbreviate these as follows:
 \begin{equation}
  \begin{split}
   {^{(wR)}\tilde{J}}[\phi] &:= {^{(wf_R R,\, 2r^{-1} wf_R)}\tilde{J}}[\phi] \\
   {^{(wR)}\tilde{K}}[\phi] &:= {^{(wf_R R,\, 2r^{-1} wf_R)}\tilde{K}}[\phi]
  \end{split}
 \end{equation}

\end{definition}

\begin{proposition}[An expression for the bulk weighted Morawetz current ${^{(wR)}\tilde{K}}$]
\label{proposition bulk weighted Morawetz current}
 The bulk weighted Morawetz current ${^{(wR)}\tilde{K}}$ can be decomposed as follows:
 \begin{equation}
  \begin{split}
   &{^{(wR)}\tilde{K}}[\phi] - \omega w f_R (\slashed{\D}_{\Lbar}\phi)\cdot \left( \slashed{\D}_R \phi + r^{-1}\phi \right) \\
   &= \frac{1}{4}w f'_R |\slashed{\D}_L\phi|^2 + \frac{1}{4} w f'_R |\slashed{\D}_{\Lbar}\phi|^2 + \left(r^{-1} wf_R - \frac{1}{2} wf'_R\right)|\slashed{\nabla}\phi|^2 - \frac{1}{2}r^{-1} wf_R'' |\phi|^2 + w\textit{Err}_{(R,\, \text{bulk})}
  \end{split}
 \end{equation}
 where, in the region $r \geq r_0$, the error term $\textit{Err}_{(wR,\, \text{bulk})}$ is given by
 \begin{equation}
  \begin{split}
   \textit{Err}_{(wR,\, \text{bulk})}[\phi] &:= 
    \frac{1}{4}\omega  f_R |\slashed{\D}_{\Lbar}\phi|^2
   - \frac{1}{4}\omega  f_R |\slashed{\D}_L\phi|^2
   - \frac{1}{2} f_R \zeta^\alpha (\slashed{\D}_L\phi)\cdot(\slashed{\nabla}_\alpha \phi)
   - \frac{1}{2} f_R \zeta^\alpha (\slashed{\D}_{\Lbar}\phi)\cdot (\slashed{\nabla}_\alpha \phi) \\
   &\phantom{:=}
   - \frac{1}{2}\omega f_R |\slashed{\nabla}\phi|^2 
   + \frac{1}{2} f_R \left( \hat{\chi}^{\alpha\beta} - \hat{\chibar}^{\alpha\beta} \right) (\slashed{\nabla}_\alpha \phi) \cdot(\slashed{\nabla}_\beta \phi) \\
   &\phantom{:=}
   + \frac{1}{4} f_R \left( \tr_{\slashed{g}} \chi_{(\text{small})} - \tr_{\slashed{g}} \chibar_{(\text{small})} -2\omega \right) (\slashed{\D}_L\phi)\cdot(\slashed{\D}_{\Lbar}\phi) \\
   &\phantom{:=}
   - \frac{1}{4} (r^{-1} f'_R - r^{-2}w f_R) \left(2\omega - \tr_{\slashed{g}}\chibar_{(\text{small})} + \tr_{\slashed{g}}\chi_{(\text{small})} \right)|\phi|^2 \\
   &\phantom{:=}
   -\omega r^{-1} f_R (\slashed{\D}_{\Lbar}\phi)\cdot \phi
   + \frac{1}{4}\partial_r (\log w) f_R |\slashed{\D}_L \phi|^2
   + \frac{1}{4}\partial_r (\log w) f_R |\slashed{\D}_{\Lbar}\phi|^2\\
   &\phantom{:=}
   - \frac{1}{2}\partial_r (\log w) f_R |\slashed{\nabla}\phi|^2 
   - \frac{1}{2}r^{-1} \left( \frac{w''}{w} f_R - \partial_r (\log w) f_R' \right)|\phi|^2\\
   &\phantom{:=}
   - \frac{1}{4}r^{-1} \partial_r (\log w) f_R \left( 2\omega - \tr_{\slashed{g}}\chibar_{(\text{small})} + \tr_{\slashed{g}}\chi_{(\text{small})} \right)|\phi|^2
  \end{split}
 \end{equation}
 which can be expressed schematically as
\begin{equation}
 \begin{split}
  \textit{Err}_{(R,\, \text{bulk})}[\phi] & = \omega f_R (\slashed{\D} \phi)^2  + f_R \bm{\Gamma}(\slashed{\D}\phi)\cdot(\overline{\slashed{\D}}\phi) + (r^{-1}f'_R - r^{-2}f_R)\bm{\Gamma}|\phi|^2 \\
  &\phantom{=}+ f_R\omega r^{-1} (\slashed{\D}_{\Lbar}\phi)\cdot \phi
  + \partial_r(\log w) f_R |\slashed{\D}\phi|^2
  + r^{-1}\left( \frac{w''}{w} f_R + \partial_r(\log w) f_R'\right)|\phi|^2\\
  &\phantom{=}
  + r^{-1} \partial_r(\log w) f_R \bm{\Gamma}|\phi|^2
 \end{split}
\end{equation}

and in the region $r < r_0$ the error term satisfies
\begin{equation}
 \textit{Err}_{(R,\, \text{bulk})} \lesssim \left(|w f_R| + |r^{-1}w f_R| + |r^{-1} w' f_R'| + |r^{-1} w'' f_R| \right) \left( |h| + |\partial h|_{(\text{rect})}\right)\left(  |\slashed{\D}\phi|^2 + |\phi|^2 \right)
\end{equation}

\end{proposition}

\begin{proof}
 This calculation makes use of propositions \ref{proposition deformation multipliers} and the definition of the modified current \ref{definition modified currents}. Note also that, for a function $f_R = f_R(r)$, we have
 \begin{equation*}
  \Box_g f_R = f''_R + \omega f'_R - \frac{1}{2}\left( \tr_{\slashed{g}}\chibar - \tr_{\slashed{g}}\chi \right) f'_R
 \end{equation*}
\end{proof}

\begin{definition}[The currents $^{(L,p)}\tilde{J}$ and $^{(L,p)}\tilde{K}$]
 We define the $p$-weighted energy currents as the modified compatible currents ${^{(Z,\, f_Z)}\tilde{J}}$ and ${^{(Z,\, f_Z)}\tilde{K}}$ with the choices $Z = f_L r^p L$ and $f_Z = 2r^{p-1} f_L$. We abbreviate these as follows:
 \begin{equation}
  \begin{split}
   {^{(L,p)}\tilde{J}}[\phi] &:= {^{(f_L r^p L,\,  2r^{p-1} f_L)}\tilde{J}}[\phi] \\
   {^{(L,p)}\tilde{K}}[\phi] &:= {^{(f_L r^p L,\,  2r^{p-1} f_L)}\tilde{K}}[\phi]
  \end{split}
 \end{equation}
\end{definition}

\begin{remark}[Alternative forms for the $p$-weighted estimates]
 We could have used a modified version of the $p$ weighted estimates, incorperating the weight $\Omega$ instead of weights $r$. Specifically, we could have chosen $Z = f_L r^p L$ and $f_Z = 2 f_L r^p L(\log \Omega)$. Then, instead of defining $\psi := r \phi$ we would define $\psi := \Omega \phi$. This leads to some cancellations among the error terms, which are important if we wish to take $p \geq 1$. However, due to other error terms (which are not removable in this way) we are restricted to the range $p < 1$ in any case, so we choose to use the simpler versions of the $p$ weighted estimates.
\end{remark}

A fairly straightforward calculation yields the following two proposition:

\begin{proposition}[An expression for the bulk $p$-weighted energy current for $p \leq 1$]
\label{proposition bulk p current p<1}
 In the region $r \geq R$, the bulk $p$-weighted energy current $^{(L,p)}\tilde{K}$ can be decomposed as follows:
\begin{equation}
 \begin{split}
  &^{(L, p)}\tilde{K}[\phi] - \omega f_L r^p(\slashed{\D}_{\Lbar}\phi)\cdot \left(\slashed{\D}_L \phi + r^{-1}\phi \right) \\
  &= \frac{1}{2} p r^{p-1} f_L |\slashed{\D}_L\phi|^2 + \frac{1}{2}(2-p)r^{p-1} f_L |\slashed{\nabla}\phi|^2 + \frac{1}{2} p(1-p)r^{p-3}f_L |\phi|^2 + \textit{Err}_{(L, p, \text{bulk})}
 \end{split}
\end{equation}
 where $\textit{Err}_{(L, p, \text{bulk})}$ is given by the following expression:
\begin{equation}
 \begin{split}
  \textit{Err}_{(L, p, \text{bulk})} &:= \frac{1}{2}\left( r^p f_L' - r^p f_L \omega \right) |\slashed{\D}_L\phi|^2 \\
  &\phantom{:=}
  + \frac{1}{2}r^p f_L \left(\tr_{\slashed{g}}\chi_{(\text{small})} - 2\omega\right) (\slashed{\D}_L\phi)\cdot(\slashed{\D}_{\Lbar}\phi) \\
  & \phantom{:=} - r^p f_L (\zeta^\alpha + 2 \slashed{\nabla}^\alpha \log\mu)(\slashed{\D}_L\phi)\cdot(\slashed{\nabla}_\alpha \phi)
  + r^p f_L \hat{\chi}^{\alpha\beta}(\slashed{\nabla}_\alpha \phi)\cdot(\slashed{\nabla}_\beta \phi) \\
  & \phantom{:=} - \frac{1}{2}(r^p f_L' - 2 f_L r^p \omega)|\slashed{\nabla}\phi|^2 \\
  & \phantom{:=} - r^{p-2}\left(\frac{1}{2}r f_L'' + pf_L' + \frac{1}{4}\left( (p-1)f_L + rf_L'\right) \left( 2\omega - \tr_{\slashed{g}}\chibar_{(\text{small})} + \tr_{\slashed{g}}\chi_{(\text{small})} \right) \right)  |\phi|^2 \\
  &\phantom{:=} - r^{p-1}f_L \omega (\slashed{\D}_{\Lbar}\phi)\cdot\phi
 \end{split}
\end{equation}

Schematically, the error term is given by
\begin{equation}
 \begin{split}
  \textit{Err}_{(L, p, \text{bulk})} - \frac{1}{2}(1-p)f_L r^{p-2}\omega |\phi|^2 &= 
  r^p(f_L' + f_L\omega) |\slashed{\D}_L \phi|^2
  + r^p f_L \left( \tr_{\slashed{g}}\chi_{(\text{small})} - 2\omega \right)(\slashed{\D}_L \phi) \cdot (\slashed{\D}_{\Lbar}\phi) \\
  &\phantom{=} + r^p f_L \begin{pmatrix} \zeta \\ \slashed{\nabla}\log\mu \end{pmatrix} (\slashed{\D}_L \phi)\cdot (\slashed{\nabla}\phi )
  + \begin{pmatrix} r^p f_L\hat{\chi} \\ r^p f'_L \\ r^p f_L\omega \end{pmatrix} (\slashed{\nabla}\phi)^2 \\
  &\phantom{=} + \begin{pmatrix} r^{p-1}f_L'' \\ r^{p-2}f_L' \\ \end{pmatrix} (\phi)^2 
  + r^{p-2}\begin{pmatrix} (p-1)f_L \\ rf'_L \end{pmatrix} \cdot \begin{pmatrix} \tr_{\slashed{g}}\chibar_{(\text{small})} \\ \tr_{\slashed{g}}\chibar_{(\text{small})} \end{pmatrix} (\phi)^2
 \end{split}
\end{equation}

\end{proposition}

\begin{proposition}[An expression for the bulk $p$-weighted energy current in terms of $\psi$]
\label{proposition bulk p current in terms of psi}
 Define the $S_{\tau,r}$-tangent tensor field
\begin{equation}
 \psi := r\phi
\end{equation}
 Then, in the region $r \geq R$, the bulk $p$-weighted energy current $^{(L,p)}\tilde{K}$ can be decomposed as follows:
\begin{equation}
 \begin{split}
  &^{(L, p)}\tilde{K}[\phi] - \omega f_L r^p(\slashed{\D}_{\Lbar}\phi)\cdot \left(\slashed{\D}_L \phi + r^{-1}\phi \right) \\
  &= \frac{1}{2}f_L r^{p-3}\left( p |\slashed{\D}_L\psi|^2 + (2-p) |\slashed{\nabla}\psi|^2\right) - \frac{1}{2}\mu^{-1}r^{-2} L\left( \mu p f_L r^p |\phi|^2 \right) + \frac{1}{2}p r^{p-2} f_L' |\phi|^2 \\
  &\phantom{=} + \textit{Err}_{(L, p, \text{bulk})} + \frac{1}{2}p \omega  r^{p-2} f_L |\phi|^2
 \end{split}
\end{equation}

\end{proposition}

\begin{remark}
 Note that in proposition \ref{proposition bulk p current p<1} we are limited to the range $p \leq 1$ in order for the coefficient of $|\phi|^2$ to be positive. In fact, the presence of semilinear terms which do not satisfy the classical null condition also limits us to the range $p \leq 1$.

 The most dangerous error term in the above expressions appears to be $r^p f_L (\tr_{\slashed{g}}\chi_{\text{small}} - 2\omega) (\slashed{\D}_L\phi)\cdot(\slashed{\D}_{\Lbar}\phi)$. This involves a ``bad derivative'' $\slashed{\D}_{\Lbar}\phi$, which can only be controlled with the help of the Morawetz vector field, which can control such a term but only with a decaying weight in $r$. Hence, this term limits the maximum value of $p$ which can be picked. Note that, interestingly, the quantity $\left(\tr_{\slashed{g}}\chi_{(\text{small})} - 2\omega \right)$ has improved regularity compared to $\tr_{\slashed{g}}\chi_{(\text{small})}$, but worse decay.

 Note that we are also limited to the range $p \leq 2$ in order to ensure that the coefficient of $|\slashed{\nabla}\psi|^2$ has a good sign; this restriction holds even in the linear case (see \cite{Dafermos2010b, Moschidis2015}). Note that, in general, the semilinear terms lead to more difficulties, and it is these which will actually limit the maximum value for $p$.

 Another dangerous error term in proposition \ref{proposition bulk p current in terms of psi} is the term $\frac{1}{2}\omega r^{p-2}|\phi|^2$. If we imagine, for a moment, that we could take $p>1$, then it is only possible to control this term if we assume that the part of $\omega$ which behaves like $r^{-1}$ either decays in $\tau$ or is uniformly positive. Note again that, if the wave coordinate condition holds, then $\omega$ decays faster than $r^{-1}$ anyway.

 One might think that the term proportional to $\tr_{\slashed{g}}\chibar r^{p-2}|\phi|^2$ is even more dangerous, since $\tr_{\slashed{g}}\chibar$ has worse decay in $r$ than $\omega$. However, if we perform the $p$ weighted energy estimates with a weight $\Omega^p$ instead of $r^p$, then this term vanishes. On the other hand, this introduces other error terms, which require us to have improved control over the foliation density $\mu$.
 
 Note, however, that both of these last terms discussed above can be handled in the case $p < 1$ without too much difficulty, and this is the approach that we will take.
\end{remark}

\chapter{Commuting}
\label{chapter commuting}

In order to obtain $L^\infty$ bounds we shall need to commute the wave equation with various operators. A standard approach is to commute with some set of vector fields, but following \cite{Dafermos2013a} we shall commute with both vector fields and covariant derivatives. To be precise, we will commute with the operator $\slashed{\D}_{T}$, where $T = \frac{1}{2}(L + \Lbar)$, as well as the differential operators $r \slashed{\nabla}$. Note that the former may be applied to a tensor field without changing the rank of that field, whereas the latter is an operator which maps tensors of rank $(0,n)$ to tensors of rank $(0, n+1)$. As discussed in detail in chapter \ref{chapter geometry of vector bundle}, these tensor fields can either be interpreted as sections of the cotangent bundle, or as sections of the bundle $\mathcal{B}$.

In this section, we will also compute the terms created by commuting with the vector field $r\slashed{\D}_L$. These calculations can be used (see appendix \ref{appendix improved energy decay}) to give \emph{improved} decay in $\tau$ for the field $\slashed{\D}_T \phi$, and for the original field $\phi$ in the region $r \leq r_0$. Although these improved decay estimates do not form part of our scheme for proving global existence, they are interesting consequences of commuting with the operator $r\slashed{\D}_L$. More importantly, we will also be able to prove $p$-weighted energy estimates for the field $r\slashed{\D}_L \phi$, that is, we will be able to control terms of the form
\begin{equation*}
\int_{\Sigma_\tau} r^{p} \left(\slashed{\D}_L (r^2\slashed{\D}_L \phi)\right)^2 \dVol_{\Sigma_\tau}
+ \int_{\mathcal{M}_\tau^{\tau_1}} r^{p-1} \left(\overline{\slashed{\D}} (r\slashed{\D}_L \phi)\right)^2 \dVol_g
\end{equation*}
Estimates of these quantities (in particular, the second term appearing above) turn out to be crucial in closing our scheme. In particular, we will find that commuting with the operator $r\slashed{\nabla}$ generates an error term which can only be controlled using an estimate of this form.

We first define the \emph{commutation operators}:
\begin{definition}[The commutation operators]
	We define the set of commutation operators
	\begin{equation*}
	\mathscr{Z} := \left\{ \slashed{\D}_T \ , \ r\slashed{\nabla} \right\}
	\end{equation*}
	We use the schematic notation $\mathscr{Z}^n f$ to mean a quantity involving $n$ operators from the set $\mathscr{Z}$ applied to the field $f$.
	
\end{definition}

We also commute with the operator $r\slashed{\D}_L$, however, we will obtain slightly different results for fields after commuting with this operator. In particular, we will not show boundedness or decay for the (weighted) $T$-energy of fields after commuting with $r\slashed{\D}_L$, however, we will still obtain the $p$-weighted energy estimates. With this in mind, we also define the \emph{extended} set of commutation operators:
\begin{definition}[The extended set of commutation operators]
	We define the set
	\begin{equation*}
	\mathscr{Y} := \left\{ r\slashed{\D}_L \ , \ \slashed{\D}_T \ , \ r\slashed{\nabla} \right\} = \mathscr{Z} \cup \{ r\slashed{\D}_L \}
	\end{equation*}
	We use the schematic notation $\mathscr{Y}^n f$ to mean a quantity involving $n$ operators from the set $\mathscr{Y}$ applied to the field $f$.
	
\end{definition}

\section{Some useful identities for computing error terms from commutators}

In this section we collect together several useful identities which can be used to simplify expressions which arise when computing commutators.

\begin{definition}[Rectangular components of second derivatives]
 Let $\phi$ be a scalar field. Then we define
\begin{equation}
 \D^2_{ab}\phi := \partial_a \partial_b \phi - \Gamma_{ab}^c \partial_c \phi
\end{equation}
Note that $\D^2_{ab}\phi$ are the rectangular components of the tensor field $\D_\mu \D_\nu \phi$.
\end{definition}

\begin{proposition}[A more explicit expression for the Riemann tensor]
\label{proposition explicit expression for Riemann}
 The rectangular components of the Riemann tensor can be expressed as
\begin{equation}
 \begin{split}
  R_{abcd} &= \frac{1}{2}\left( \D^2_{ac} h_{bd} + \D^2_{bd}h_{ac} - \D^2_{ad} h_{bc} - \D^2_{bc} h_{ad} \right) \\
  &\phantom{=} + \frac{1}{4}(g^{-1})^{ef} \bigg(
  (\partial_a h_{fd})(\partial_e h_{bc}) + (\partial_a h_{ce})(\partial_b h_{fd}) + (\partial_a h_{cf})(\partial_d h_{be})
  - (\partial_a h_{fd})(\partial_b h_{ce}) \\
  &\phantom{= + \frac{1}{4}(g^{-1})^{ef} \bigg(}
  - (\partial_a h_{fd})(\partial_c h_{be}) - (\partial_a h_{ed})(\partial_f h_{bc})
  +(\partial_b h_{fd})(\partial_c h_{ae}) + (\partial_b h_{cf})(\partial_e h_{ad})\\
  &\phantom{= + \frac{1}{4}(g^{-1})^{ef} \bigg(}
  + (\partial_b h_{ed})(\partial_f h_{ac})
  - (\partial_b h_{ce})(\partial_f h_{ad}) - (\partial_b h_{fd})(\partial_e h_{ac}) - (\partial_b h_{cf})(\partial_d h_{ae}) \\
  &\phantom{= + \frac{1}{4}(g^{-1})^{ef} \bigg(}
  + (\partial_c h_{ae})(\partial_f h_{bd}) + (\partial_c h_{bf})(\partial_e h_{ad}) + (\partial_c h_{af})(\partial_d h_{be})
  - (\partial_c h_{be})(\partial_f h_{ad}) \\
  &\phantom{= + \frac{1}{4}(g^{-1})^{ef} \bigg(}
  - (\partial_c h_{bf})(\partial_d h_{ae}) - (\partial_c h_{af})(\partial_e h_{bd})
  + (\partial_d h_{be})(\partial_f h_{ac}) + (\partial_d h_{ae})(\partial_f h_{bc}) \\
  &\phantom{= + \frac{1}{4}(g^{-1})^{ef} \bigg(}
  + (\partial_d h_{be})(\partial_f h_{ac}) - (\partial_d h_{ae})(\partial_f h_{bc})
  + (\partial_e h_{bc})(\partial_f h_{ad}) + (\partial_e h_{ad})(\partial_f h_{bc}) \\
  &\phantom{= + \frac{1}{4}(g^{-1})^{ef} \bigg(}
  + (\partial_e h_{bd})(\partial_f h_{ac})
  - (\partial_e h_{ac})(\partial_f h_{bd}) - (\partial_e h_{bd})(\partial_f h_{ac}) - (\partial_e h_{ad})(\partial_f h_{bc}) \bigg) 
 \end{split}
\end{equation}
Note that the rectangular components of the metric perturbation $h_{ab}$ are \emph{scalar fields}, while $\D^2_{ab}$ represents the rectangular components of the second covariant derivative.
\end{proposition}

The proposition above will allow us to express the components of the Riemann tensor in terms of derivatives of the rectangular components of $h$. We also need the following proposition, which can be used to express the second covariant derivatives of the rectangular components of $h$ in terms of derivatives of ``commuted'' field, i.e.\ the fields $\mathcal{Z}h_{ab}$ where $\mathcal{Z}$ is either the operator $\D_T$ or the operator $r\slashed{\nabla}$.

\begin{proposition}[Expressing second covariant derivatives in terms of the commutation operators]
 \label{proposition second derivatives as commutators}
 We can express all second derivatives of a scalar field in terms of first derivatives of fields after commutation operators have been applied, plus some lower order terms.

 Let $\phi$ be a scalar field. Then its second derivatives can be expressed as
\begin{equation}
 \begin{split}
  L^a L^b \D^2_{ab} \phi &= 2LT\phi - \slashed{\Delta}\phi + \tilde{\Box}_g \phi + \frac{1}{2} \left( \tr_{\slashed{g}}\chi \right) \Lbar \phi + \frac{1}{2}\left( \tr_{\slashed{g}}\chibar - 2\omega \right)L\phi + \zeta^\alpha \slashed{\nabla}_\alpha \phi \\ \\
  L^a \Lbar^b \D^2_{ab} \phi &= \slashed{\Delta}\phi - \tilde{\Box}_g \phi - \frac{1}{2}(\tr_{\slashed{g}}\chi - 2\omega) \Lbar\phi - \frac{1}{2}(\tr_{\slashed{g}}\chibar) L\phi \\ \\
  L^a \slashed{\Pi}_\mu^{\phantom{\mu}b} \D^2_{ab} \phi &= r^{-1} \slashed{\D}_L(r\slashed{\nabla}_\mu \phi) + \frac{1}{2}\zeta_\mu L\phi - r^{-1}\slashed{\nabla}_\mu \phi \\ \\
  \Lbar^a \Lbar^b \D^2_{ab}\phi &= 2\Lbar T\phi - \slashed{\Delta}\phi + \tilde{\Box}_g \phi + \frac{1}{2}(\tr_{\slashed{g}}\chibar - 2\omega)L\phi + \frac{1}{2}(\tr_{\slashed{g}}\chi )\Lbar \phi - \zeta^\alpha \slashed{\nabla}_\alpha \phi \\ \\
  \Lbar^a \slashed{\Pi}_\mu^{\phantom{\mu}b} \D^2_{ab}\phi &= r^{-1} \slashed{\D}_{\Lbar} (r\slashed{\nabla}_\mu \phi) - (\slashed{\nabla}_\mu \log\mu) L\phi - \frac{1}{2}\left( \zeta_\mu + 2(\slashed{\nabla}_\mu \log\mu) \right) \Lbar\phi + r^{-1}\slashed{\nabla}_\mu \phi \\ \\
  \slashed{\Pi}_\mu^{\phantom{\mu}a}\slashed{\Pi}_\nu^{\phantom{\nu}b} \D^2_{ab} \phi &= \slashed{\nabla}^2_{\mu\nu}\phi - \frac{1}{2}\chibar_{\mu\nu} L\phi - \frac{1}{2}\chi_{\mu\nu} \Lbar \phi
 \end{split}
\end{equation}

\end{proposition}

\begin{definition}[The $S_{\tau,r}$-tangent tensor fields $\slashed{R}$, $\slashed{R}_Z$ and $\slashed{R}_{ZW}$]
We define the $S_{\tau,r}$-tangent tensor field
\begin{equation}
 \left(\slashed{R} \right)_{\mu\nu\rho\sigma} := \slashed{\Pi}_{\mu\nu\rho\sigma}^{\alpha\beta\gamma\delta} R_{\alpha\beta\gamma\delta}
\end{equation}
Given a vector field $Z$ (which is not necessarily $S_{\tau,r}$-tangent itself) we define the $S_{\tau,r}$-tangent tensor field
\begin{equation}
 \left(\slashed{R}_Z \right)_{\mu\nu\rho} := \slashed{\Pi}_{\mu\nu\rho}^{\alpha\beta\gamma} Z^\sigma R_{\sigma \alpha\beta\gamma}
\end{equation}
Similarly, given a pair of vector fields $Z$, $W$ we define the $S_{\tau,r}$-tangent tensor field
\begin{equation}
 \left(\slashed{R}_{ZW} \right)_{\mu\nu} := \slashed{\Pi}_{\mu\nu}^{\alpha\beta} Z^\sigma W^\rho R_{\sigma\rho\alpha\beta}
\end{equation}
\end{definition}

\begin{proposition}[An expression for $R_{(\text{frame})}$]
 The frame components of the Riemann tensor are given, schematically, by
\begin{equation}
 R_{(\text{frame})} = (\partial T h)_{(\text{frame})} + \left(\partial (r \slashed{\nabla} h)\right)_{(\text{frame})} + (\tilde{\Box}_g h)_{(\text{frame})} + \bm{\Gamma} \cdot (\partial h)_{(\text{frame})} + (\partial h)_{(\text{frame})} \cdot (\partial h)_{(\text{frame})} 
\end{equation}

\end{proposition}

\begin{proposition}[The curvature components $\slashed{R}$, $\slashed{R}_L$, $\slashed{R}_{\Lbar}$ and $\slashed{R}_{L\Lbar}$]
\label{proposition R RL}
 The curvature components $\slashed{R}$, $\slashed{R}_L$, $\slashed{R}_{\Lbar}$ and $\slashed{R}_{L\Lbar}$ frequently occur in error terms. These terms are given schematically by:
\begin{equation}
 \begin{split}
  \slashed{R} &= \bigg( r^{-1}\slashed{\nabla}(r\slashed{\nabla} h)+ r^{-1}|\partial h| + |\bm{\Gamma}_{(\text{good})}| |\partial h| + |\bm{\Gamma}| |\bar{\partial} h| \bigg)_{(\text{frame})} \\
  \\
  \slashed{R}_L &= \bigg( r^{-1}\slashed{\D}_L (r\slashed{\nabla} h) + r^{-1}\slashed{\nabla}(r\slashed{\nabla} h)+ r^{-1}|\partial h| + |\bm{\Gamma}_{(\text{good})}| |\partial h| + |\bm{\Gamma}| |\bar{\partial} h| \bigg)_{(\text{frame})} \\
  \\
  \slashed{R}_{\Lbar} &= \bigg( r^{-1}\slashed{\D}_{\Lbar} (r\slashed{\nabla} h) + r^{-1}\slashed{\nabla}(r\slashed{\nabla} h )+ r^{-1}|\partial h| + |\bm{\Gamma}||\partial h| \bigg)_{(\text{frame})} \\
  \\
 \slashed{R}_{L\Lbar} &= \bigg( r^{-1}\slashed{\D}_{\Lbar} (r\slashed{\nabla} h) + r^{-1}\slashed{\D}_L(r\slashed{\nabla} h) + r^{-1}|\partial h|  + |\bm{\Gamma}||\partial h| \bigg)_{(\text{frame})} \\
  \\
 \end{split}
\end{equation}
\end{proposition}

\begin{proof}
 Since $\slashed{\D}^2_{ab}\phi$ are the rectangular components of the tensor field $\D_\mu \D_\nu \phi$ we can calculate, for example,
\begin{equation*}
 \begin{split}
  L^a \slashed{\Pi}_\mu^{\phantom{\mu}b} \D^2_{ab} \phi &= L^\nu \slashed{\Pi}_\mu^{\phantom{\mu}\rho} \D^2_{\nu\rho} \phi \\
  &= \slashed{\Pi}_\mu^{\phantom{\mu}\rho}\D_L \left( \slashed{\Pi}_\rho^{\phantom{\rho}\sigma} \D_\sigma \phi \right) - \slashed{\Pi}_\mu^{\phantom{\mu}\rho} \left( \D_L \slashed{\Pi}_\rho^{\phantom{\rho}\sigma} \right) \D_\sigma \phi \\
  &= \slashed{\D}_L \slashed{\nabla}_\mu \phi + \frac{1}{2}\zeta_\mu \D_L \phi \\
  &= r^{-1}\slashed{\D}_L \left( r\slashed{\nabla}_\mu \phi \right) + \frac{1}{2}\zeta_\mu \D_L \phi - r^{-1} \slashed{\nabla}_\mu \phi 
 \end{split}
\end{equation*}

\end{proof}

\begin{proposition}[First derivatives of the projection operator]
 The first derivatives of the projection operator $\slashed{\Pi}$ are given by
\begin{equation}
  \begin{split}
      \slashed{\Pi}_{\sigma}^{\phantom{\sigma}\lambda} \D_\mu \left( \slashed{\Pi}_\lambda^{\phantom{\lambda}\rho} \right)
	&= \frac{1}{4}\zeta_\sigma \Lbar_\mu L^\rho 
	- \frac{1}{4}\left( \zeta_\sigma + 2(\slashed{\nabla}_\sigma \mu) \right) L_\mu \Lbar^\rho
	- \frac{1}{2}(\slashed{\nabla}_\sigma \mu)L_\mu L^\rho \\
	& \phantom{=} + \frac{1}{2}\chi_{\mu\sigma} \Lbar^\rho
	+ \frac{1}{2}\chibar_{\mu\sigma} L^\rho
  \end{split}
\end{equation}

\end{proposition}

\begin{proof}
 This is a short computation making use of proposition \ref{proposition derivatives of projection}. Note that the tensor computed above is ``transverse'' in the sense that
\begin{equation}
 (X_A)_\rho \left( \slashed{\Pi}_{\sigma}^{\phantom{\sigma}\lambda} \D_\mu \left( \slashed{\Pi}_\lambda^{\phantom{\lambda}\rho} \right) \right) = 0  
\end{equation}

\end{proof}

By a similar (though more lengthy) computation we can prove the following related proposition:

\begin{proposition}[Applying the wave operator to the projection operator]
 The projection operator $\slashed{\Pi}$ satisfies the following identity:
\begin{equation}
 \begin{split}
  \slashed{\Pi}_\lambda^{\phantom{\lambda}\nu}\Box_g \slashed{\Pi}_\nu^{\phantom{\nu}\rho}
  &= \left( -\frac{1}{2}\slashed{\D}_L \slashed{\nabla}_\lambda \log\mu + \frac{1}{4}\slashed{\D}_{\Lbar}\zeta_\lambda + \frac{1}{2}\slashed{\Div}\chibar_\lambda - \omega (\slashed{\nabla}_\lambda \log \mu) - \frac{1}{2} \omega \zeta_\lambda + \frac{1}{4}\zeta_\lambda \tr_{\slashed{g}}\chibar - \frac{1}{4}\zeta^\alpha \chibar_{\alpha \lambda} \right) L^\rho \\
  &\phantom{=} + \left( -\frac{1}{4}\slashed{\D}_L\left( \zeta_\lambda + 2\slashed{\nabla}_\lambda \log \mu \right)  + \frac{1}{2}\slashed{\Div}\chi_{\lambda} - \frac{1}{4}\zeta_\lambda \tr_{\slashed{g}}\chi - \frac{1}{2} (\slashed{\nabla}_\lambda \log \mu)\tr_{\slashed{g}}\chi + \frac{1}{4}\zeta^\alpha \chibar_{\alpha \lambda} \right) \Lbar^\rho \\
  &\phantom{=} + \left( \zeta_\lambda \zeta^\alpha + (\slashed{\nabla}_\lambda \log \mu)\zeta^\alpha + \zeta_\lambda (\slashed{\nabla}^\alpha \log \mu) + \chi^{\alpha\beta}\chibar_{\beta\lambda} + \chibar^{\alpha\lambda}\chi_{\beta\lambda} \right) \slashed{\Pi}_\alpha^{\phantom{\alpha}\rho}
 \end{split}
\end{equation} 
which is given schematically by
\begin{equation}
  \slashed{\Pi}_\lambda^{\phantom{\lambda}\nu} \Box_g \slashed{\Pi}_\nu^{\phantom{\nu}\rho} =
  \begin{pmatrix} \slashed{\nabla}\omega \\ \slashed{\D}_{\Lbar}\zeta \\ \slashed{\Div}\chibar_{(\text{small})} \\ \bm{\Gamma}\cdot \left(r^{-1} + \bm{\Gamma} \right)\end{pmatrix} L^\rho
  + \begin{pmatrix} \slashed{\D}_L \zeta \\ \slashed{\nabla}\omega \\ \slashed{\Div}\chi_{(\text{small})} \\ \bm{\Gamma}\cdot \left(r^{-1} + \bm{\Gamma}\right) \end{pmatrix} \Lbar^\rho
  + \begin{pmatrix} r^{-2} \\ r^{-1}\bm{\Gamma} \\ \bm{\Gamma}\cdot\bm{\Gamma} \end{pmatrix}\slashed{\Pi}^\rho
\end{equation}
where we have suppressed the index $\lambda$ on the right hand side.
\end{proposition}

\begin{proposition}[Commuting $L$ and $\Lbar$ through the reduced wave operator]
 Let $\phi_\alpha$ be an $S_{\tau,r}$-tangent one-form. Then we have, schematically,
\begin{equation}
 \begin{split}
  L^\alpha \tilde{\Box}_g \phi_\alpha &= \begin{pmatrix} r^{-1} \\ \bm{\Gamma} \end{pmatrix} \overline{\slashed{\D}}\phi + \begin{pmatrix} \slashed{\D}_L \zeta \\ \slashed{\nabla}\omega \\ \slashed{\Div}\chi_{(\text{small})} \end{pmatrix}\cdot \phi + \bm{\Gamma} \cdot \begin{pmatrix} r^{-1} \\ \bm{\Gamma} \end{pmatrix} \cdot \phi
  \\ \\
  \Lbar^\alpha \tilde{\Box}_g \phi_\alpha &= \begin{pmatrix} r^{-1} \\ \bm{\Gamma} \end{pmatrix} \slashed{\D}\phi + \begin{pmatrix} \slashed{\D}_{\Lbar} \zeta \\ \slashed{\nabla}\omega \\ \slashed{\Div}\chibar_{(\text{small})} \end{pmatrix}\cdot \phi + \bm{\Gamma} \cdot \begin{pmatrix} r^{-1} \\ \bm{\Gamma} \end{pmatrix} \cdot \phi
 \end{split}
\end{equation}
\end{proposition}

% 
% 
%\begin{equation}
%  \begin{split}
%    \slashed{\Pi}_\sigma^{\phantom{\sigma}\lambda} \Box_g \left( \slashed{\Pi}_\lambda^{\phantom{\lambda}\rho} \right)
%    &= \Bigg( - \frac{1}{4}L\left( \zeta^A + 2\mu^{-1}(\slashed{\upd}^A \mu) \right) 
%      - \frac{1}{4}\left(\zeta^A + 2\mu^{-1}(\slashed{\upd}^A \mu) \right)(\chi_B^{\phantom{B}A} - \omega\delta_B^A )
%      + \frac{1}{2}\slashed{\nabla}^B \chi_B^{\phantom{B}A} \\
%      & \phantom{= + \Bigg( } + \frac{1}{4}\left( \chi^{AB}\zeta_B - (\tr_{\slashed{g}}\chi)\left( \zeta^A + 2\mu^{-1}(\slashed{\upd}^A \mu) \right)\right)
%      \Bigg) (X_A)_\sigma \Lbar^\rho \\
%    &\phantom{=} + \Bigg( - \frac{1}{2}L\left( \mu^{-1}(\slashed{\upd}^A \mu) \right)
%      + \frac{1}{4}\Lbar(\zeta^A)
%      + \frac{1}{2}\slashed{\nabla}^B \chi_A^{\phantom{A}B}
%      - \frac{1}{2}\mu^{-1}(\slashed{\upd}^B \mu) \left( \chi_B^{\phantom{B}A} + \omega\delta_B^A \right) \\
%      & \phantom{= + \Bigg(} - \frac{1}{4}\zeta^B \left( \mu^{-1}(\slashed{\upd}_B \mu)\rho^A - \slashed{\mathcal{L}}_B \rho^A \right)
%      - \frac{1}{4}(\tr_{\slashed{g}}\chi)\left( \zeta^A - 2\mu^{-1}(\slashed{\upd}^A \mu) \right) \Bigg)(X_A)_\sigma L^\rho \\
%    &\phantom{=} + \left( -\zeta^A \left( \zeta^B + 2\mu^{-1}(\slashed{\upd}^B \mu) \right) + \chibar^{AC}\chi_C^{\phantom{C}B} \right)\left( (X_A)_\sigma (X_B)^\rho + (X_B)_\sigma (X_A)^\rho \right)
%  \end{split}
%\end{equation}

\section{Commuting with first order operators}

In many situations we simply need to commute the commutation operators with the first order operators $\slashed{\D}_L$, $\slashed{\D}_{\Lbar}$ and $\slashed{\nabla}$. We have

\begin{proposition}[Commuting $\slashed{\D}_T$ with first order operators]
	\label{proposition commuting DT with first order operators}
	We have the following schematic expressions:
	\begin{equation}
	\begin{split}
	[\slashed{\D}_L \, , \slashed{\D}_T] \phi &=
	- \frac{1}{2}\omega \slashed{\D}_L \phi
	- \frac{1}{2}\omega \slashed{\D}_{\Lbar} \phi
	- (\zeta^\alpha + \slashed{\nabla}^\alpha \log \mu) \slashed{\nabla}_\alpha \phi
	+ \frac{1}{2}L^\mu \Lbar^\nu [\slashed{\D}_\mu \, , \slashed{\D}_\nu]\phi
	\\ \\
	[\slashed{\D}_{\Lbar} \, , \slashed{\D}_T] \phi &=
	\frac{1}{2}\omega \slashed{\D}_L \phi
	+ \frac{1}{2}\omega \slashed{\D}_{\Lbar} \phi
	+ (\zeta^\alpha + \slashed{\nabla}^\alpha \log \mu) \slashed{\nabla}_\alpha \phi
	- \frac{1}{2}L^\mu \Lbar^\nu [\slashed{\D}_\mu \, , \slashed{\D}_\nu]\phi
	\\ \\
	[\slashed{\nabla}_{\slashed{\alpha}} \, , \slashed{\D}_T] \phi &=
	- \frac{1}{2}(\slashed{\nabla}_{\slashed{\alpha}}\log\mu) \slashed{\D}_L \phi
	- \frac{1}{2}(\slashed{\nabla}_{\slashed{\alpha}}\log\mu) \slashed{\D}_{\Lbar} \phi
	+ \frac{1}{2}\left( (\chi_{(\text{small})})_{\slashed{\alpha}}^{\phantom{\slashed{\alpha}}\slashed{\beta}} + (\chibar_{(\text{small})})_{\slashed{\alpha}}^{\phantom{\slashed{\alpha}}\slashed{\beta}} \right) \slashed{\nabla}_{\slashed{\beta}}\phi
	\\
	&\phantom{=}
	+ \slashed{\Pi}_{\slashed{\alpha}}^{\phantom{\slashed{\alpha}}\mu} T^\nu [\slashed{\D}_\mu, \slashed{\D}_\nu]\phi
	\end{split}
	\end{equation}
	where the terms involving the curvature of $\mathcal{B}$ (i.e.\ terms involving $[\slashed{\D}_\mu, \slashed{\D}_\nu]\phi$) are absent if $\phi$ is a scalar field.
\end{proposition}

\begin{proposition}[Commuting $\slashed{\D}_{rL}$ with first order operators]
	\label{proposition commuting rL with first order operators}
	We have
	\begin{equation}
	\begin{split}
	[\slashed{\D}_L \, , \slashed{\D}_{rL}]\phi
	&=
	\slashed{\D}_L \phi
	\\ \\
	[\slashed{\D}_{\Lbar} \, , \slashed{\D}_{rL}]\phi
	&=
	(-1 + r\omega)\slashed{\D}_L \phi
	+ r\omega \slashed{\D}_{\Lbar} \phi
	+ 2r(\zeta^\alpha + \slashed{\nabla}^\alpha \log \mu) \slashed{\nabla}_\alpha \phi
	-r L^\mu \Lbar^\nu [\slashed{\D}_\mu, \slashed{\D}_\nu]\phi
	\\ \\
	[\slashed{\nabla}_{\slashed{\alpha}} \, , \slashed{\D}_{rL}]\phi
	&=
	r\chi_{\slashed{\alpha}}^{\phantom{\slashed{\alpha}}\slashed{\beta}} \slashed{\nabla}_{\slashed{\beta}}
	+ r\slashed{\Pi}_{\slashed{\alpha}}^{\phantom{\slashed{\alpha}} \mu} L^\nu [\slashed{\D}_\mu \, , \slashed{\D}_\nu]\phi
	\end{split}
	\end{equation}
	where, again, the terms involving $[\slashed{\D}_\mu, \slashed{\D}_\nu]\phi$ are absent if $\phi$ is a scalar field.
\end{proposition}

\begin{proposition}[Commuting $r\slashed{\nabla}$ with first order operators]
	\label{proposition commuting rnabla with first order operators}
	We have
	\begin{equation}
	\begin{split}
	[\slashed{\D}_L \, , r\slashed{\nabla}_{\slashed{\alpha}}] \phi
	&=
	-r(\chi_{(\text{small})})_{\slashed{\alpha}}^{\phantom{\slashed{\alpha}}\slashed{\beta}} \slashed{\nabla}_{\slashed{\beta}} \phi
	+ rL^\mu \slashed{\Pi}_{\slashed{\alpha}}^{\phantom{\slashed{\alpha}}\nu} [\slashed{\D}_\mu \slashed{\D}_\nu]\phi
	\\ \\
	[\slashed{\D}_{\Lbar} \, , r\slashed{\nabla}_{\slashed{\alpha}}] \phi
	&=
	r(\slashed{\nabla}_{\slashed{\alpha}} \log \mu) \slashed{\D}_L \phi
	+ r(\slashed{\nabla}_{\slashed{\alpha}} \log \mu) \slashed{\D}_{\Lbar} \phi
	- r(\chibar_{(\text{small})})_{\slashed{\alpha}}^{\phantom{\slashed{\alpha}}\slashed{\beta}} \slashed{\nabla}_{\slashed{\beta}}\phi
	+ r\Lbar^\mu \slashed{\Pi}_{\slashed{\alpha}}^{\phantom{\slashed{\alpha}}\nu} [\slashed{\D}_\mu \slashed{\D}_\nu]\phi
	\\ \\
	[\slashed{\nabla}_{\slashed{\alpha}} \, , r\slashed{\nabla}_{\slashed{\beta}}] \phi
	&=
	r\slashed{\Pi}_{\slashed{\alpha}}^{\phantom{\slashed{\alpha}}\mu} \slashed{\Pi}_{\slashed{\beta}}^{\phantom{\slashed{\beta}}\nu} [\slashed{\D}_\mu \slashed{\D}_\nu]\phi
	\end{split}
	\end{equation}
	where again, the terms which are linear in $\phi$ are absent if $\phi$ is a scalar field. Note that in many cases we will not need to make use of the final equality, since we have
	\begin{equation*}
	r\slashed{\nabla}_{\alpha} \slashed{\nabla}_{\beta} \phi
	= \slashed{\nabla}_{\alpha} \left( r\slashed{\nabla}_{\beta} \phi \right)
	\end{equation*}
	and so, if we do not care about the order in which the derivatives are applied, we can avoid commuting them. Note also that the final term can be expressed either in terms of the Gauss curvature $K$, or the curvature of $\mathcal{B}$, $\Omega$.
	
\end{proposition}

\section{Preliminary commutation calculations}

Throughout this chapter we will require the deformation tensors associated with the commutation operators $\slashed{\D}_T$ and $r\slashed{\nabla}$. Here, $T = \frac{1}{2}(L + \Lbar)$ is the same vector field which has already appeared as a multiplier. The deformation tensor of $T$ was already calculated in proposition \ref{proposition deformation multipliers}, so in the following proposition we will calculate the deformation tensor associated with the operator $r \slashed{\nabla}$.

Note that, as in the case of the multiplier vector fields, it would be possible to modify the commutation operators (in particular, we could use a $\mu$-weighted version of the vector field $T$) in order to remove some of the worst terms. However, as before we expect terms which are at least this bad to appear due to the semilinear terms, and so we prefer to treat all of these terms together, making use of the semilinear hierarchy to close our estimates. Moreover, modifying the vector fields in this way introduces additional error terms which are not possible to control without improved estimates on $T\mu$, which we cannot, in general prove\footnote{Once again, if the wave coordinate condition holds, then we can in fact prove improved estimates on $\mu$, by commuting the equation for $(\Lbar h)_{LL}$ in \eqref{equation wave coordinate condition}}.

% \begin{proposition}[The deformation tensor of the vector field $\check{T}$]
%  Recall that we have defined
% \begin{equation}
%  \check{T} := \frac{1}{2}( L + \mu\Lbar )
% \end{equation}
% Then $\check{T}$ has the associated deformation tensor in the region $r \geq r_0$:
% \begin{equation}
%   \begin{split}
%     ^{(T)}\pi_{\mu\nu} &= -\frac{1}{2}\omega L_\mu L_\nu 
%     + \frac{1}{4}\left( \omega\mu_{(\text{small})} - (\Lbar \mu) \right)(L_\mu \Lbar_\nu + \Lbar_\mu L_\nu ) \\
%     & \phantom{=} - \frac{1}{2}(\slashed{g}^{-1})^{AB}\left(\zeta_A + \mu^{-1}\slashed{\upd}_A \mu + \slashed{\upd}_A \mu \right)\left( L_\mu (X_B)_\nu + (X_B)_\mu L_\nu \right) \\
%     & \phantom{=} + \frac{1}{2}(\slashed{g}^{-1})^{AB}\left(\mu\zeta_A + \slashed{\upd}_A \mu \right)\left(\Lbar_\mu (X_B)_\nu + (X_B)_\mu \Lbar_\nu \right) \\
%     & \phantom{=} + \left( (\chi_{(\text{small})})^{AB} + \mu(\chibar_{(\text{small})})^{AB} - r^{-1}\mu_{(\text{small})} (\slashed{g}^{-1})^{AB} \right)(X_A)_\mu (X_B)_\nu
%   \end{split}
% \end{equation}
% \end{proposition}

\begin{definition}[Higher order deformation tensors]
 We can also define the higher order deformation tensor $^{(\Xi)}\pi$ associated with higher order operators. In particular, let $\Xi$ be a field which defines, at every point, a vector in the vector space formed by taking the Cartesian product of the tangent space at that point, with the space of $S_{\tau,r}$-tangent one-forms at that point. In terms of abstract indices, we can write this tensor as $\Xi_{\slashed{\alpha}}^{\phantom{\slashed{\alpha}}\mu}$. Then, we associate to $\Xi$ the deformation tensor:
\begin{equation}
\begin{split}
	^{(\Xi)}\pi_{\slashed{\alpha} \mu \nu} &:= \slashed{\Pi}_{\slashed{\alpha}}^{\phantom{\slashed{\alpha}}\slashed{\beta}} \left( \D_\mu \Xi_{\slashed{\beta} \nu} + \D_\nu \Xi_{\slashed{\beta} \mu} \right)\\
	&=
	\slashed{\D}_\mu \Xi_{\slashed{\alpha} \nu} + \slashed{\D}_\nu \Xi_{\slashed{\alpha} \mu}
\end{split}
\end{equation}
where in the second line we must remember that the covariant derivative $\slashed{\D}$ ``projects'' all the \emph{slashed} indices using $\slashed{\Pi}$.

\end{definition}

\begin{proposition}[The deformation tensor associated with the operator $r\slashed{\nabla}$]
 The deformation tensor of the operator $r\slashed{\nabla}$ is given by
\begin{equation}
  \begin{split}
    ^{(r\slashed{\Pi})}\pi_{\slashed{\alpha} \mu \nu} &= 
      -\frac{1}{2}r (\slashed{\nabla}_\slashed{\alpha} \log\mu)(L_\mu \Lbar_\nu + \Lbar_\mu L_\nu )
      - r (\slashed{\nabla}_\slashed{\alpha} \log\mu)L_\mu L_\nu \\
      &\phantom{=} + \frac{1}{2}r \left( (\chi_{(\text{small})})_{\slashed{\alpha}\mu} \Lbar_\nu + (\chi_{(\text{small})})_{\slashed{\alpha}\nu}\Lbar_\mu \right)
      + \frac{1}{2}r \left( (\chibar_{(\text{small})})_{\slashed{\alpha}\mu}L_\nu + (\chibar_{(\text{small})})_{\slashed{\alpha}\nu} L_\mu \right)
    \bigg)
  \end{split}
\end{equation}
\end{proposition}

\begin{remark}[The trace $^{(r\slashed{\Pi})}\pi_{\slashed{\alpha} \rho}^{\phantom{\slashed{\alpha} \rho}\rho}$]
 Note that the trace of the deformation tensor satisfies
\begin{equation}
  ^{(r\slashed{\Pi})}\pi_{\slashed{\alpha} \rho}^{\phantom{\slashed{\alpha} \rho}\rho} = 2r\mu^{-1} \slashed{\nabla}_\slashed{\alpha} \mu
\end{equation}

\end{remark}

\begin{definition}[Commutation currents]
 Given a vector field $Z$, we define the associated \emph{commutation current} $^{(Z)}\mathscr{J}[\phi]$ as follows:
\begin{equation}
 {^{(Z)}\mathscr{J}^\mu}[\phi] := {^{(Z)}\pi}^{\mu\nu}\cdot\slashed{\D}_\nu \phi - \frac{1}{2}\left( \tr_g {^{(Z)}\pi} \right) \cdot\slashed{\D}^\mu \phi
\end{equation}
 Likewise, given a tensor field $\Xi_{\slashed{\alpha}}^{\phantom{\slashed{\alpha}}\mu}$ we can define the associated commutation current:
\begin{equation}
 ^{(\Xi)}{\mathscr{J}_{\slashed{\alpha}}^{\phantom{\slashed{\alpha}}\mu}}[\phi] := {^{(\Xi)}\pi}_{\slashed{\alpha}}^{\phantom{\slashed{\alpha}}\mu\nu}\cdot\slashed{\D}_\nu \phi - \frac{1}{2}\left( (g^{-1})^{\rho\sigma}\ {^{(\Xi)}\pi_{\slashed{\alpha} \rho\sigma}} \right)\cdot \slashed{\D}^\mu \phi
\end{equation}

\end{definition}

\begin{proposition}[Basic commutation identity with vector fields]
\label{proposition commute vector field}
 The projected reduced wave operator satisfies the following commutation identity with a vector field operator $\slashed{\D}_Z$:
\begin{equation}
\label{equation commute vector field}
 \begin{split}
    \tilde{\slashed{\Box}}_g (\slashed{\D}_Z \phi) &=
      \slashed{\D}_Z( \tilde{\slashed{\Box}}_g \phi)
      + \slashed{\D}_\mu {^{(Z)}\mathscr{J}[\phi]^\mu}
      + \frac{1}{2}\left(\tr_g {^{(Z)}\pi}\right) \tilde{\slashed{\Box}}_g \phi
      - (Z\omega)\slashed{\D}_{\Lbar}\phi
      + \omega \slashed{\D}_{[\Lbar, Z]}\phi
      - \frac{1}{2}\omega \left(\tr_{\slashed{g}} {^{(Z)}\pi}\right) \slashed{\D}_{\Lbar} \phi \\
      &\phantom{=}
      + [\slashed{\D}_\mu \, , \slashed{\D}_\nu]\left( Z^\nu \slashed{\D}^\mu \phi \right) 
      + \slashed{\D}^\mu\left(Z^\nu [\slashed{\D}_\mu \, , \slashed{\D}_\nu]\phi \right)
      + \omega \Lbar^\mu Z^\nu [\slashed{\D}_\mu \, , \slashed{\D}_\nu] \phi
 \end{split}
\end{equation}

\end{proposition}

\begin{proof}
 This follows by direct, but fairly long computation, making use of the definition of the commutation current ${^{(Z)}\mathscr{J}}[\phi]$. We have
 \begin{equation*}
 \begin{split}
 	\slashed{\Box}_g \slashed{\D}_Z \phi
 	&=
 	\slashed{\D}_\mu \slashed{\D}^\mu \left( Z^\nu \slashed{\D}_\nu \phi \right)
 	\\
 	&=
 	\slashed{\D}_\mu \left( (\D^\mu Z^\nu)\slashed{\D}_\nu \phi + Z^\nu \slashed{\D}^\mu \slashed{\D}_\nu \phi \right)
 	\\
 	&=
 	(\Box_g Z^\nu) \slashed{\D}_\nu \phi
 	+ 2(\D^\mu Z^\nu) \slashed{\D}_\mu \slashed{\D}_\nu \phi
 	+ Z^\nu \slashed{\D}_\mu \slashed{\D}^\mu \slashed{\D}_\nu \phi
 	\\
 	&=
 	\slashed{\D}_Z \left(\slashed{\Box}_g \phi\right)
 	+ \left( \D_\mu {^{(Z)}\pi}^{\mu\nu} \right) \slashed{\D}_\nu \phi
 	- (\D_\mu \D_\nu Z^\mu) \slashed{\D}^\nu \phi
 	+ {^{(Z)}\pi}^{\mu\nu} \slashed{\D}_\mu \slashed{\D}_\nu \phi
 	+ (\D^\mu Z^\nu) [\slashed{\D}_\mu, \slashed{\D}_\nu]\phi
 	\\
 	&\phantom{=}
 	+ Z^\nu [\slashed{\D}_\mu, \slashed{\D}_\nu](\slashed{\D}_\mu \phi)
 	+ Z^\nu \slashed{\D}^\mu [\slashed{\D}_\mu, \slashed{\D}_\nu]\phi
 	\\
 	&=
 	\slashed{\D}_Z \left(\slashed{\Box}_g \phi\right)
 	+ \left( \D_\mu {^{(Z)}\pi}^{\mu\nu} \right) \slashed{\D}_\nu \phi
 	- \frac{1}{2}\left( \D_\mu \tr_g{^{(Z)}\pi} \right) \slashed{\D}^\mu \phi
 	+ {^{(Z)}\pi}^{\mu\nu} \slashed{\D}_\mu \slashed{\D}_\nu \phi
 	+ [\slashed{\D}_\mu, \slashed{\D}_\nu] \left( Z^\nu \slashed{\D}^\mu \phi \right)
 	\\
 	&\phantom{=}
 	+ \slashed{\D}^\mu \left( Z^\nu [\slashed{\D}_\mu, \slashed{\D}_\nu]\phi \right)
 	\\
 	&=
 	\slashed{\D}_Z \left(\slashed{\Box}_g \phi\right)
 	+ \slashed{\D}_\mu \left( {^{(Z)}\pi}^{\mu\nu} \slashed{\D}_\nu \phi
 		- \frac{1}{2} \left( \tr_g{^{(Z)}\pi} \right) \slashed{\D}^\mu \phi \right)
 	+ \frac{1}{2} \left( \tr_g{^{(Z)}\pi} \right) \slashed{\Box}_g \phi
 	\\
 	&\phantom{=}
 	+ [\slashed{\D}_\mu, \slashed{\D}_\nu] \left( Z^\nu \slashed{\D}^\mu \phi \right)
 	+ \slashed{\D}^\mu \left( Z^\nu [\slashed{\D}_\mu, \slashed{\D}_\nu]\phi \right)
 \end{split}
 \end{equation*}
 and we also have
 \begin{equation*}
 \omega \slashed{\D}_{\Lbar} \slashed{\D}_Z \phi
 =
  \slashed{\D}_Z \left( \omega \slashed{\D}_{\Lbar} \phi \right)
 + \omega \left( \slashed{\D}_{\Lbar} \slashed{\D}_Z \phi - \slashed{\D}_Z \slashed{\D}_{\Lbar} \phi - \slashed{\D}_{[\Lbar, Z]}\phi \right)
 - (Z\omega) \slashed{\D}_{\Lbar} \phi
 + \omega \slashed{\D}_{[\Lbar, Z]} \phi
 \end{equation*}
\end{proof}

\begin{remark}
 Note that, in the case that $\phi$ is a scalar field, all of the terms involving commutators of the covariant derivatives vanish\footnote{Of course, the term $\omega \slashed{\D}_{[\Lbar, Z]}\phi$ does not necessarily vanish!}. This is evidently the fact when the commutator is applied directly to the field $\phi$, and we also have
 \begin{equation*}
 \begin{split}
 [\D_\mu, \D_\nu] \left( Z^\nu \D^\mu \phi \right)
 &=
 R_{\mu\nu \phantom{\nu} \rho}^{\phantom{\mu\nu} \nu} Z^\rho \D^\mu \phi
 + R_{\mu\nu \phantom{\mu} \rho}^{\phantom{\mu\nu} \nu} Z^\nu \D^\rho \phi
 \\
 &=
 \text{Ric}[g]_{\mu\nu} Z^\mu \D^\nu \phi
 - \text{Ric}[g]_{\nu\mu} Z^\mu \D^\nu \phi
 \\
 &= 0
 \end{split}
 \end{equation*}
\end{remark}

\begin{proposition}[Basic commutation identity with tensorial operators]
\label{proposition commute tensor field}
 The operator $\Xi_{\slashed{\alpha}}^{\phantom{\slashed{\alpha}}\mu} \slashed{\D}_\mu$ satisfies the following commutation identity:
\begin{equation}
\label{equation commute tensor field}
 \begin{split}
 \tilde{\slashed{\Box}}_g (\Xi_{\slashed{\alpha}}^{\phantom{\slashed{\alpha}}\mu}\slashed{\D}_\mu \phi) &=
 \Xi_{\slashed{\alpha}}^{\phantom{\slashed{\alpha}}\mu}\slashed{\D}_\mu( \tilde{\slashed{\Box}}_g \phi)
 + \slashed{\D}_\mu {^{(\Xi)}\mathscr{J}[\phi]_{\slashed{\alpha}}^{\phantom{\slashed{\alpha}}\mu}}
 + \frac{1}{2}\left( {^{(\Xi)}\pi_{\slashed{\alpha}\mu}^{\phantom{\slashed{\alpha}\mu}\mu}}\right) \tilde{\slashed{\Box}}_g \phi
 - (\Xi_{\slashed{\alpha}}^{\phantom{\slashed{\alpha}}\mu} \D_\mu \omega)\slashed{\D}_{\Lbar}\phi
 \\
 &\phantom{=}
 + \omega \left( (\slashed{\D}_{\Lbar} \Xi_{\slashed{\alpha}}^{\phantom{\slashed{\alpha}}\mu})
 	- \Xi_{\slashed{\alpha}}^{\phantom{\slashed{\alpha}}\nu}(\D_\nu \Lbar^\mu) \right)
 	\slashed{\D}_\mu \phi
 - \frac{1}{2} \omega \left({^{(\Xi)}\pi_{\slashed{\alpha}\mu}^{\phantom{\slashed{\alpha}\mu}\mu}}\right) \slashed{\D}_{\Lbar} \phi
 + [\slashed{\D}_\mu \, , \slashed{\D}_\nu]\left( \Xi_{\slashed{\alpha}}^{\phantom{\slashed{\alpha}}\nu} \slashed{\D}^\mu \phi \right) 
 \\
 &\phantom{=}
 + \slashed{\D}^\mu\left(\Xi_{\slashed{\alpha}}^{\phantom{\slashed{\alpha}}\nu} [\slashed{\D}_\mu \, , \slashed{\D}_\nu]\phi \right)
 + \omega \Lbar^\mu \Xi_{\slashed{\alpha}}^{\phantom{\slashed{\alpha}}\nu} [\slashed{\D}_\mu \, , \slashed{\D}_\nu] \phi
 \end{split}
\end{equation}

\end{proposition}

\begin{proof}
 This follows from a computation which is very similar to the proof of proposition \ref{proposition commute vector field}. One must be careful to distinguish between different types of indices when applying the operators $\slashed{\D}$.
\end{proof}

\begin{remark}
 Unlike the expression \eqref{equation commute vector field}, the term $[\slashed{\D}_\mu \, , \slashed{\D}_\nu]\left( \Xi_{\slashed{\alpha}}^{\phantom{\slashed{\alpha}}\nu} \slashed{\D}^\mu \phi \right)$ does \emph{not} in general vanish, even if $\phi$ is a scalar field. In fact, in this case we have
 \begin{equation*}
 [\slashed{\D}_\mu \, , \slashed{\D}_\nu]\left( \Xi_{\slashed{\alpha}}^{\phantom{\slashed{\alpha}}\nu} \D^\mu \phi \right)
 =
 \Omega_{\mu\nu\slashed{\alpha}}^{\phantom{\mu\nu\slashed{\alpha}}\slashed{\beta}} \Xi_{\slashed{\beta}}^{\phantom{\slashed{\beta}}\nu} \slashed{\D}^\mu \phi
 \end{equation*}

%[[rewrite this]]
%At first sight the term in the final line appears to be related to the deformation tensor $^{(\Xi)}\pi$, however, the index labeled $\gamma$ in the above expression is \emph{not} projected using the projection operator $\slashed{\Pi}$, as it would need to be in order to express this term in terms of $^{(\Xi)}\pi$. Indeed, for a tensor $\Xi$ satisfying
%\begin{equation}
% L^\alpha \Xi_\alpha^{\phantom{\alpha}\beta} = \Lbar^{\alpha} \Xi_\alpha^{\phantom{\alpha}\beta} = 0
%\end{equation}
%we have, schematically,
%\begin{equation}
%  \begin{split}
%    \D^\mu \Xi_\gamma^{\phantom{\gamma}\nu} + \D^\nu \Xi_\gamma^{\phantom{\gamma}\mu}  - \D_\rho \Xi_\gamma^{\phantom{\gamma}\rho} (g^{-1})^{\mu\nu}
%    &= {^{(\Xi)}\pi_\gamma^{\phantom{\gamma}\mu\nu}} - \frac{1}{2}{^{(\Xi)}\pi_{\gamma \rho}^{\phantom{\gamma \rho}\rho}} (g^{-1})^{\mu\nu}
%    + \left( \left(r^{-1} + \bm{\Gamma} \right)\cdot \Xi \right)_\gamma^{\phantom{\gamma}\mu\nu}
%  \end{split}
%\end{equation}
%Making use of proposition \ref{proposition derivatives of projection} as well as the expression above, we find that the final term in equation \eqref{equation commute tensor field} is schematically of the form
%\begin{equation}
% \mu (r^-1 + \bm{\Gamma})\cdot \left( {^{(\Xi)}\pi} + \bm{\Gamma}\cdot \Xi \right)\cdot \D\phi + \mu r^{-2} \Xi \cdot \slashed{\nabla}\phi
%\end{equation}

\end{remark}

\section{Null frame decomposition of the commutation currents}

The commutation currents $^{(T)}\mathscr{J}[\phi]$, $^{(rL)}\mathscr{J}[\phi]$ and $^{(r\slashed{\Pi})}\mathscr{J}[\phi]$ have the following null frame decompositions:

\begin{proposition}[Null frame decomposition of $^{(T)}\mathscr{J}{[\phi]}$]
 The null frame components of the commutation current $^{(T)}\mathscr{J}[\phi]$ are given by
\begin{equation}
 \begin{split}
    ^{(T)}\mathscr{J}_L[\phi] &= - \frac{1}{2}\left( \tr_{\slashed{g}}\chi_{(\text{small})} + \tr_{\slashed{g}}\chibar_{(\text{small})} \right) \slashed{\D}_L \phi
      -\omega \slashed{\D}_{\Lbar} \phi
      - \zeta^\alpha \slashed{\nabla}_\alpha \phi \\ \\
    ^{(T)}\mathscr{J}_{\Lbar}[\phi] &= \omega \slashed{\D}_L \phi
      -\frac{1}{2}\left( \tr_{\slashed{g}}\chi_{(\text{small})} + \tr_{\slashed{g}}\chibar_{(\text{small})} \right)\slashed{\D}_{\Lbar} \phi
      + \left( \zeta^\alpha + 2\slashed{\nabla}^\alpha \log\mu \right)\slashed{\nabla}_\alpha \phi \\ \\
    \slashed{\Pi}_\mu^{\phantom{\mu}\nu} \left({^{(T)}\mathscr{J}_{\mu}}[\phi]\right) &= -\frac{1}{2}\left( \zeta_\mu + 2\slashed{\nabla}_\mu \log\mu \right) \slashed{\D}_L \phi
      + \frac{1}{2}\zeta_\mu \slashed{\D}_{\Lbar}\phi
      + \left( \hat{\chi}_\mu^{\phantom{\mu}\nu} + \hat{\chibar}_\mu^{\phantom{\mu}\nu} \right) \slashed{\nabla}_\nu \phi 
 \end{split}
\end{equation}

\end{proposition}

\begin{proposition}[Null frame decomposition of $^{(rL)}\mathscr{J}{[\phi]}$]
 The null frame components of the commutation current $^{(T)}\mathscr{J}[\phi]$ are given by
\begin{equation}
 \begin{split}
    ^{(rL)}\mathscr{J}_L[\phi] &= -\left(2 + r (\tr_{\slashed{g}}\chi_{(\text{small})} ) \right)\slashed{\D}_L \phi \\ \\
    ^{(rL)}\mathscr{J}_{\Lbar}[\phi] &= -2(1-r\omega) \slashed{\D}_L \phi
      -\left( 2 + r(\tr_{\slashed{g}}\chi_{(\text{small})})\right) \slashed{\D}_{\Lbar} \phi
      + 2r\left( \zeta^\mu + \slashed{\nabla}^\mu \log\mu \right)\slashed{\nabla}_\mu \phi \\ \\
    \slashed{\Pi}_\mu^{\phantom{\mu}\nu} \left({^{(rL)}\mathscr{J}}_{\nu}[\phi] \right) &= -r\left( \zeta_\mu + \slashed{\nabla}_\mu \log\mu \right) \slashed{\D}_L \phi
      - \left( 1 + r\omega \right) \slashed{\nabla}_\mu \phi 
      + 2r \hat{\chi}_\mu^{\phantom{\mu}\nu}\slashed{\nabla}_\nu \phi 
 \end{split}
\end{equation}

\end{proposition}

\begin{proposition}[Null frame decomposition of $^{(r\slashed{\Pi})}\mathscr{J}{[\phi]}$]
 The null frame components of the higher order commutation current $^{(r\slashed{\Pi})}\mathscr{J}[\phi]$ are given by
\begin{equation}
 \begin{split}
    ^{(r\slashed{\Pi})}\mathscr{J}_{\slashed{\alpha} L}[\phi] &= -r(\chi_{(\text{small})})_{\slashed{\alpha}}^{\phantom{\slashed{\alpha}}\slashed{\beta}} \slashed{\nabla}_{\slashed{\beta}} \phi \\
    \\
    ^{(r\slashed{\Pi})}\mathscr{J}_{\slashed{\alpha} \Lbar}[\phi] &= 2r (\slashed{\nabla}_{\slashed{\alpha}} \log\mu) \slashed{\D}_L \phi
      -r(\chibar_{(\text{small})})_{\slashed{\alpha}}^{\phantom{\slashed{\alpha}}\slashed{\beta}} \slashed{\nabla}_{\slashed{\beta}} \phi \\
    \\
    \slashed{\Pi}_\mu^{\phantom{\mu}\nu} \left({^{(r\slashed{\Pi})}\mathscr{J}}_{\slashed{\alpha} \nu}[\phi]\right) &= \frac{1}{2}r(\chibar_{(\text{small})})_{\slashed{\alpha}\mu} \slashed{\D}_{L}\phi
      +\frac{1}{2}r(\chi_{(\text{small})})_{\slashed{\alpha}\mu} \slashed{\D}_{\Lbar}\phi - r(\slashed{\nabla}_\slashed{\alpha} \log \mu)\slashed{\nabla}_\mu \phi
 \end{split}
\end{equation}

\end{proposition}

\begin{proposition}[Decomposition of the commutation current of the vector field $T$]
\label{proposition commute T}
 Define the commutation current
\begin{equation}
 ^{(T)}\mathscr{K}[\phi] := \slashed{\D}_\mu \left( {^{(T)}\mathscr{J}[\phi]} \right)^\mu + \frac{1}{2}\left( \tr_{\slashed{g}} {^{(T)}\pi} \right) \tilde{\slashed{\Box}}_g \phi - (T\omega)\slashed{\D}_{\Lbar}\phi - \frac{1}{2}\omega \left( \tr_{\slashed{g}} {^{(T)}\pi} \right) \slashed{\D}_{\Lbar} \phi + \omega \slashed{\D}_{[\Lbar, T]} \phi
\end{equation}
 Then $^{(T)}\mathscr{K}[\phi]$ can be decomposed as
\begin{equation}
 \begin{split}
    ^{(T)}\mathscr{K}[\phi] &= {^{(T)}\mathscr{K}}_{(\pi, \Lbar)}[\phi] + {^{(T)}\mathscr{K}}_{(\pi, L)}[\phi] + {^{(T)}\mathscr{K}}_{(\pi, \slashed{\Pi})}[\phi] + {^{(T)}\mathscr{K}}_{(\pi, \text{elliptic})}[\phi] + {^{(T)}\mathscr{K}}_{(\pi, \text{good})}[\phi]  \\
    & \phantom{=} + {^{(T)}\mathscr{K}}_{(\phi)}[\phi] +  {^{(T)}\mathscr{K}}_{(\text{low})}[\phi]
 \end{split}
\end{equation}
where the various terms are defined by
\begin{equation}
 \begin{split}
  {^{(T)}\mathscr{K}}_{(\pi, \Lbar)}[\phi] &:= \frac{1}{2}(\slashed{\Div}\zeta) \slashed{\D}_{\Lbar}\phi \\ \\
  {^{(T)}\mathscr{K}}_{(\pi, L)}[\phi] &:= \left( \frac{1}{2}(T \tr_{\slashed{g}}\chi_{(\text{small})}) + \frac{1}{2}(T \tr_{\slashed{g}}\chibar_{(\text{small})}) - \frac{1}{2}(\slashed{\Div}\zeta) \right)\slashed{\D}_L \phi \\ \\
  {^{(T)}\mathscr{K}}_{(\pi, \slashed{\Pi})}[\phi] &:= \left( \frac{1}{2}(\slashed{\D}_T \zeta^{\alpha}) - (\slashed{\nabla}^\alpha \omega) \right)\slashed{\nabla}_\alpha \phi \\ \\
  {^{(T)}\mathscr{K}}_{(\pi, \text{elliptic})}[\phi] &:= - (\slashed{\Delta}\log \mu)\slashed{\D}_L \phi + \left( (\slashed{\Div}\, \hat{\chi})^\alpha + (\slashed{\Div}\, \hat{\chibar})^\alpha \right)\slashed{\nabla}_\alpha \phi \\ \\
  {^{(T)}\mathscr{K}}_{(\pi, \text{good})}[\phi] &:= \left( \frac{1}{4}(L\tr_{\slashed{g}}\chi_{(\text{small})}) + \frac{1}{4}(L\tr_{\slashed{g}}\chibar_{(\text{small})}) - \frac{1}{2}(L\omega) \right)\slashed{\D}_{\Lbar}\phi \\
    &\phantom{:=} + \left( -\frac{1}{4}(L\tr_{\slashed{g}}\chi_{(\text{small})}) - \frac{1}{4}(L\tr_{\slashed{g}}\chibar_{(\text{small})}) - \frac{1}{2}(L\omega) \right) \slashed{\D}_L \phi
    - (\slashed{\D}_L \zeta^\alpha)\slashed{\nabla}_\alpha \phi \\ \\
  {^{(T)}\mathscr{K}}_{(\phi)}[\phi] &:= \omega \slashed{\D}_{\Lbar} \slashed{\D}_T \phi
    + \zeta^\alpha r^{-1} \slashed{\D}_{\Lbar} (r\slashed{\nabla}_\alpha \phi)
    - \omega \slashed{\D}_L \slashed{\D}_T \phi
    - (\zeta^\alpha + 2 \slashed{\nabla}^\alpha \log \mu)r^{-1} \slashed{\D}_L (r\slashed{\nabla}_\alpha \phi) \\
    &\phantom{:=} + \left( \hat{\chi}^{\alpha\beta} + \hat{\chibar}^{\alpha\beta} + \frac{1}{2}\left( \tr_{\slashed{g}}\chi_{(\text{small})} + \tr_{\slashed{g}}\chibar_{(\text{small})} \right)(\slashed{g}^{-1})^{\alpha\beta} \right)r^{-1} \slashed{\nabla}_\alpha (r\slashed{\nabla}_\beta \phi) \\ \\
 \end{split}
\end{equation}

and ${^{(T)}\mathscr{K}}_{(\text{low})}$ is given schematically by
\begin{equation}
 \label{equation TKlow}
 \begin{split}
    {^{(T)}\mathscr{K}}_{(\text{low})}[\phi]_{\sigma_1 \ldots \sigma_n} &=
     \bm{\Gamma}\cdot \begin{pmatrix}r^{-1} \\ \bm{\Gamma} \end{pmatrix}\cdot (\slashed{\D}\phi)
     + \begin{pmatrix} \omega L^\mu \Lbar^\nu \\ (\zeta^\mu + 2\slashed{\nabla}^\mu \log \mu)L^\nu \\ \zeta^\mu \Lbar^\nu \end{pmatrix} [\slashed{\D}_\mu \, , \slashed{\D}_\nu]\phi
 \end{split}
\end{equation}

Note that, importantly, the terms involving the curvature $\Omega$ in the above expression are \emph{not} present if $\phi$ is a scalar field.

\end{proposition}

\begin{proposition}[Decomposition of the commutation current of the vector field $rL$]
\label{proposition commute rL}
 The commutation current
\begin{equation}
 ^{(rL)}\mathscr{K}[\phi] := \slashed{\D}_\mu \left( {^{(rL)}\mathscr{J}[\phi]} \right)^\mu - (rL\omega)\slashed{\D}_{\Lbar}\phi  + \frac{1}{2}\left( \tr_{\slashed{g}} {^{(rL)}\pi} \right) \tilde{\slashed{\Box}}_g \phi + \omega\slashed{\D}_{[\Lbar, rL]}\phi - \frac{1}{2}\omega ( \tr {^{(rL)}\pi})\slashed{\D}_{\Lbar}\phi
\end{equation}
 can be decomposed as
\begin{equation}
 \begin{split}
    ^{(rL)}\mathscr{K}[\phi] &= {^{(rL)}\mathscr{K}}_{(\text{large})}[\phi] + {^{(rL)}\mathscr{K}}_{(\pi, \Lbar)}[\phi] + {^{(rL)}\mathscr{K}}_{(\pi, L)}[\phi] + {^{(rL)}\mathscr{K}}_{(\pi, \slashed{\Pi})}[\phi] + {^{(rL)}\mathscr{K}}_{(\phi)}[\phi] \\
      &\phantom{=} + {^{(rL)}\mathscr{K}}_{(\text{low})}[\phi]
 \end{split}
\end{equation}
where the various terms are defined by
\begin{equation}
 \begin{split}
  {^{(rL)}\mathscr{K}}_{(\text{large})}[\phi] &:= \slashed{\Delta}\phi 
  + r^{-1} \slashed{\D}_L \left( r\slashed{\D}_L \phi \right)
  + r^{-1} \slashed{\D}_L \phi
  \\ \\
  {^{(rL)}\mathscr{K}}_{(\pi, \Lbar)}[\phi] &:= \left( \frac{1}{2}rL(\tr_{\slashed{g}}\chi_{(\text{small})}) - r(L\omega) \right)\slashed{\D}_{\Lbar}\phi
  \\ \\
  {^{(rL)}\mathscr{K}}_{(\pi, L)}[\phi] &:= \left( \frac{1}{2}r(\Lbar \tr_{\slashed{g}}\chi_{(\text{small})}) - r(\slashed{\Div}\zeta) - r(\slashed{\Delta}\log\mu) - r(L\omega) \right)\slashed{\D}_L \phi
  \\ \\
  {^{(rL)}\mathscr{K}}_{(\pi, \slashed{\Pi})}[\phi] &:= \left( -r(\slashed{\D}_L \zeta^\alpha) - 2r(\slashed{\nabla}^\alpha \omega) + 2r(\slashed{\Div}\hat{\chi})^\alpha  \right)\slashed{\nabla}_\alpha \phi
  \\ \\
  {^{(rL)}\mathscr{K}}_{(\phi)}[\phi] &:= 
  (\tr_{\slashed{g}}\chi_{(\text{small})} - \omega)r\slashed{\Delta}\phi 
  + 2r\hat{\chi}^{\mu\nu}\slashed{\nabla}_\mu \slashed{\nabla}_\nu \phi
  - 2(\zeta^\mu + \slashed{\nabla}^\mu \log\mu )\slashed{\D}_L (r \slashed{\nabla}_\mu \phi)
  - \omega \slashed{\D}_L \left(r\slashed{\D}_L \phi\right)
 \end{split}
\end{equation}

and where ${^{(rL)}\mathscr{K}}_{(\text{low})}$ is given schematically by
\begin{equation}
\begin{split}
  {^{(rL)}\mathscr{K}}_{(\text{low})}[\phi] 
  &= 
  r \omega \tilde{\slashed{\Box}}_g \phi
  + r\bm{\Gamma}\cdot \bm{\Gamma} \cdot (\overline{\slashed{\D}} \phi)
  + \begin{pmatrix} 1 \\ r\omega \\ r\tr_{\slashed{g}}\chi_{(\text{small})} \end{pmatrix} \begin{pmatrix} \omega \\ \tr_{\slashed{g}}\chi_{(\text{small})} \end{pmatrix} \slashed{\D}_{\Lbar} \phi
  \\
  &\phantom{=}
  + r(\zeta^\mu + \slashed{\nabla}^\mu \log \mu)L^\nu [\slashed{\D}_\mu \, , \slashed{\D}_\nu]\phi
\end{split}
\end{equation}

% 
% 
% \begin{equation}
%   \begin{split}
%     {^{(rL)}\mathscr{K}}_{(\text{large})}[\phi] &:= \slashed{\Delta}\phi + r^{-1}\slashed{\D}_L \slashed{\D}_L (r\phi) - 2\slashed{\Box}_g \phi \\ \\
%     {^{(rL)}\mathscr{K}}_{(\pi, \Lbar)}[\phi] &:= \frac{1}{2}L\left( r\tr_{\slashed{g}}\chi_{(\text{small})} \right) \slashed{\D}_{\Lbar}\phi \\ \\
%     {^{(rL)}\mathscr{K}}_{(\pi, L)}[\phi] &:= \left( r\mu^{-1}\check{T}(\tr_{\slashed{g}}\chi) - \frac{1}{2}\mu^{-1}L(r\tr_{\slashed{g}}\chi) - rL(\omega) - r\slashed{\Div}\zeta - r\mu^{-1}\slashed{\Delta}\mu \right)\slashed{\D}_L \phi \\ \\
%     {^{(rL)}\mathscr{K}}_{(\pi, A)}[\phi] &:= \left( -L\left( r\zeta^A + r\mu^{-1}(\slashed{\upd}^A\mu) \right) + 2r(\slashed{\Div}\hat{\chi})^A - r (\slashed{\nabla}^A \omega) \right)\slashed{\nabla}_A \phi \\ \\
%     {^{(rL)}\mathscr{K}}_{(\phi)}[\phi] &:= -\left( r\tr_{\slashed{g}}\chi_{(\text{small})} + r \omega \mu_{(\text{small})} \right)\slashed{\Delta}\phi
%     	+ 2r\hat{\chi}^{AB} \slashed{\nabla}^2_{AB} \phi
%     	-\omega\slashed{\D}_L\left(r \slashed{\D}_L \phi\right) \\
%     	&\phantom{:=} -2(\zeta^A + \mu^{-1}\slashed{\upd}^A\mu)(X_A)^\mu \slashed{\D}_L (r\slashed{\nabla}_\mu \phi)
% 	- r(\tr_{\slashed{g}}\chi_{(\text{small})})\slashed{\Box}_g \phi    
%   \end{split}
% \end{equation}
% and ${^{(rL)}\mathscr{K}}_{(\text{low})}[\phi]$ satisfies the bound
% \begin{equation}
%   \begin{split}
%     \left|{^{(rL)}\mathscr{K}}_{(\text{low})}[\phi]\right| &\lesssim \left|(1+\mu^{-1})\bm{\Gamma}\cdot(1+r\bm{\Gamma})\cdot \slashed{\D}\phi \right|
%   \end{split}
% \end{equation}
% 

\end{proposition}

\begin{proposition}[Decomposition of the commutation current of the differential operator $r\slashed{\nabla}$]
\label{proposition commute rnabla}
 Define the commutation current
\begin{equation}
 \begin{split}
  ^{(r\slashed{\Pi})}\mathscr{K}[\phi]_{\slashed{\alpha}} &:= \slashed{\D}_\mu \left( {^{(r\slashed{\Pi})}\mathscr{J}[\phi]} \right)_{\slashed{\alpha}}^{\phantom{\slashed{\alpha}}\mu}
  - r(\slashed{\nabla}_{\slashed{\alpha}} \omega) \slashed{\D}_{\Lbar}\phi
  + \omega \left(\slashed{\D}_{\Lbar} (r\slashed{\Pi}_{\slashed{\alpha}}^{\phantom{\slashed{\alpha}}\mu})\right)\slashed{\D}_\mu \phi 
  - \omega \left( r \slashed{\nabla}_{\slashed{\alpha}} \Lbar^\mu \right) \slashed{\D}_\mu \phi
  \\
  &\phantom{:=} - \frac{1}{2}\left({^{(r\slashed{\Pi})}\pi}_{\slashed{\alpha}\beta}^{\phantom{\slashed{\alpha}\beta}\beta} \right)\omega \slashed{\D}_{\Lbar} \phi
  + \frac{1}{2}\left( {^{(r\slashed{\Pi})}\pi}_{\slashed{\alpha}\beta}^{\phantom{\slashed{\alpha}\beta}\beta} \right)\tilde{\slashed{\Box}}_g \phi
 \end{split}
\end{equation}
 Then $^{(r\slashed{\Pi})}\mathscr{K}[\phi]_\slashed{\alpha}$ can be decomposed as
\begin{equation}
 \begin{split}
    ^{(r\slashed{\Pi})}\mathscr{K}[\phi]_{\slashed{\alpha}} &= {^{(r\slashed{\Pi})}\mathscr{K}}_{(\pi, \Lbar)}[\phi]_{\slashed{\alpha}} + {^{(r\slashed{\Pi})}\mathscr{K}}_{(\pi, L)}[\phi]_{\slashed{\alpha}} + {^{(r\slashed{\Pi})}\mathscr{K}}_{(\pi, \slashed{\Pi})}[\phi]_{\slashed{\alpha}} + {^{(r\slashed{\Pi})}\mathscr{K}}_{(\pi, \text{elliptic})}[\phi]_{\slashed{\alpha}} \\
    & \phantom{=} + {^{(r\slashed{\Pi})}\mathscr{K}}_{(\pi, \text{good})}[\phi]_{\slashed{\alpha}} + {^{(r\slashed{\Pi})}\mathscr{K}}_{(\phi)}[\phi]_{\slashed{\alpha}} + {^{(r\slashed{\Pi})}\mathscr{K}}_{(r\slashed{\D}_L)}[\phi]_{\slashed{\alpha}} + {^{(r\slashed{\Pi})}\mathscr{K}}_{(\text{low})}[\phi]_{\slashed{\alpha}} 
 \end{split}
\end{equation}
where the various terms are defined by
\begin{equation}
 \begin{split}
    {^{(r\slashed{\Pi})}\mathscr{K}}_{(\pi, \Lbar)}[\phi]_{\slashed{\alpha}} &:= \left( \frac{1}{4}r \left(\slashed{\nabla}_{\slashed{\alpha}} \tr_{\slashed{g}}\chi_{(\text{small})} \right) - r\slashed{\nabla}_{\slashed{\alpha}} \omega \right)\slashed{\D}_{\Lbar}\phi \\
    \\
    {^{(r\slashed{\Pi})}\mathscr{K}}_{(\pi, L)}[\phi]_{\slashed{\alpha}} &:= \left( \frac{1}{4}r \left( \slashed{\nabla}_{\slashed{\alpha}} \tr_{\slashed{g}}\chibar_{(\text{small})} \right) - r\slashed{\nabla}_{\slashed{\alpha}} \omega \right) \slashed{\D}_L \phi \\
    \\
    {^{(r\slashed{\Pi})}\mathscr{K}}_{(\pi, \slashed{\Pi})}[\phi]_{\slashed{\alpha}} &:= \left( \frac{1}{2}r \left(\slashed{\D}_{\Lbar} \chi_{(\text{small})}\right)_{\slashed{\alpha}}^{\phantom{\slashed{\alpha}}\beta}
    + \frac{1}{2}r \left(\slashed{\D}_L \chibar_{(\text{small})}\right)_{\slashed{\alpha}}^{\phantom{\slashed{\alpha}}\beta}
    \right)\slashed{\nabla}_{\slashed{\beta}} \phi \\ 
    \\
    {^{(r\slashed{\Pi})}\mathscr{K}}_{(\pi, \text{elliptic})}[\phi]_{\slashed{\alpha}} &:= \left(\frac{1}{2}r (\slashed{\Div}\,\hat{\chi})_{\slashed{\alpha}} \right) \slashed{\D}_{\Lbar}\phi
      + \left(\frac{1}{2}r (\slashed{\Div}\,\hat{\chibar})_{\slashed{\alpha}} \right)\slashed{\D}_L \phi 
      - \left(r \slashed{\nabla}^\beta \slashed{\nabla}_{\slashed{\alpha}} \log\mu\right) \slashed{\nabla}_\beta \phi \\
    \\
    {^{(r\slashed{\Pi})}\mathscr{K}}_{(\phi)}[\phi]_{\slashed{\alpha}} &:=
    (\chi_{(\text{small})})_{\slashed{\alpha}}^{\phantom{\slashed{\alpha}}\beta} \slashed{\D}_{\Lbar} \left( r\slashed{\nabla}_\beta \phi \right)
    + (\chibar_{(\text{small})})_{\slashed{\alpha}}^{\phantom{\slashed{\alpha}}\beta} \slashed{\D}_L \left( r\slashed{\nabla}_\beta \phi \right) 
     - 2(\slashed{\nabla}_{\slashed{\alpha}} \log \mu) \slashed{\nabla}^\beta (r\slashed{\nabla}_\beta \phi) \\
    \\
    {^{(r\slashed{\Pi})}\mathscr{K}}_{(r\slashed{\D}_L)}[\phi]_{\slashed{\alpha}} &:=
    - (\slashed{\nabla}_{\slashed{\alpha}} \log \mu)\slashed{\D}_L \left( r\slashed{\D}_L \phi \right)
 \end{split}
\end{equation}

and the current ${^{(r\slashed{\Pi})}\mathscr{K}}_{(\text{low})}[\phi]$ is given schematically by
\begin{equation}
 \begin{split}
  {^{(r\slashed{\Pi})}\mathscr{K}}_{(\text{low})}[\phi] &=
  r\bm{\Gamma}\cdot \begin{pmatrix} r^{-1} \\ \bm{\Gamma} \end{pmatrix}\cdot \slashed{\D}\phi
  + \begin{pmatrix}
  	r(\chibar_{\text{small}})_{\slashed{\alpha}}^{\phantom{\slashed{\alpha}}\mu} L^\nu \\
  	r(\chi_{\text{\small} })_{\slashed{\alpha}}^{\phantom{\slashed{\alpha}}\mu} \Lbar^\nu \\  	
  \end{pmatrix} \cdot [\slashed{\D}_\mu \, , \slashed{\D}_\nu]\phi
 \end{split}
\end{equation}

Note that, as before, the terms involving the curvature components do \emph{not} appear in the case that $\phi$ is a scalar field.

\end{proposition}

\begin{remark}[Explanation of the various error terms]
	From the point of view of regularity, the most dangerous error terms are those which involve top order derivatives of the connection coefficients, i.e.\ the terms given in $\mathscr{K}_{(\pi, \Lbar)}$, $\mathscr{K}_{(\pi, L)}$, and $\mathscr{K}_{(\pi, A)}$. In particular, terms involving $T(\tr_{\slashed{g}}\chi_{(\text{small})})$ and $r\slashed{\nabla}_\alpha (\tr_{\slashed{g}} \chi_{(\text{small})})$ are the most difficult to control. Of these, the terms appearing in $\mathscr{K}_{(\pi, \Lbar)}$ require the most delicate handling, since these are critical from the point of view of regularity and \emph{also} require good decay properties.
	
	The terms given in $\mathscr{K}_{(\pi, \text{elliptic})}$ can be estimated in terms of other quantities by using elliptic estimates, while the terms appearing in $\mathscr{K}_{(\pi, \text{good})}$ can be estimated straightforwardly using the transport equations satisfied by the connection coefficients. The terms appearing in $^{(T)}\mathscr{K}_{(\text{low})}$ are lower order and are easily controlled.
	
	The terms appearing in $\mathscr{K}_{(\phi)}$ are relatively easy to control, with the important exception of the term $\omega \slashed{\D}_{\Lbar} \slashed{\D}_T \phi$ in ${^{(T)}\mathscr{K}_{(\phi)}}$. This is ``critical'' from the point of view of decay, and also contains a bad derivative of the field in question, namely $\slashed{\D}_T \phi$. It is this term which is responsible for the energy growth for higher order energies in quasilinear equations, even if there are no semilinear terms.
	
	The terms in ${^{(r\slashed{\Pi})}\mathscr{K}}_{(\pi, \slashed{\Pi})}[\phi]_{\slashed{\alpha}}$ also require some special handling: to control these we must use the transport equations satisfied by $\chi$ and $\chibar$ \emph{and} the structure of the curvature component $R_{L\slashed{\alpha}\Lbar\slashed{\beta}}$.
	
	The term in ${^{(r\slashed{\Pi})}\mathscr{K}_{(r\slashed{\D}_L)}}$ needs some special treatment: see the following remark.
	
	Finally, we reiterate that all terms involving the curvature tensor $\Omega$ do not arise if $\phi$ is a scalar field. This is important, since the components of this tensor are second order in the metric (and hence second order in the fields), but since they do not arise for scalar fields they can never be above leading order.
\end{remark}

\begin{remark}[The error term ${^{(r\slashed{\Pi})}\mathscr{K}_{(r\slashed{\D}_L)}}$]
	\label{remark error term needing commutation with rL}
	
	The term in ${^{(r\slashed{\Pi})}\mathscr{K}_{(r\slashed{\D}_L)}}$ is the \emph{only} term arising from commuting with $r\slashed{\nabla}$ and $\slashed{\D}_T$ which also requires us to commute with $r\slashed{\D}_L$. In other words, if this term did not arise here, then we could get away with commuting with the operators $r\slashed{\nabla}$ and $\slashed{\D}_T$ alone; there would be no need to commute with $r\slashed{\D}_L$.
	
	One could ask the following question: is it possible to change our set-up so that this term does not arise? For example, we could try changing the way the spheres $S_{(\tau,r)}$ are defined, so that the operator $r\slashed{\nabla}$ has slightly different commutation properties with the wave operator. In fact, if we choose our normalisation so that
	\begin{equation*}
	g^{-1}(\upd r , \upd r) = \mu^{-1}
	\end{equation*}
	instead of $g^{-1}(\upd r , \upd r) = \mu^{-1}$, then this error term can be avoided. However, the conformal transformation needed to enforce this condition will introduce new inhomogeneous terms of the form
	\begin{equation*}
	\slashed{\D}^\mu (\log \mu) \slashed{\D}_\mu \phi
	\end{equation*}
	These terms involve a term of the form
	\begin{equation*}
	(\Lbar h)_{LL} (\slashed{\D}_{\Lbar} \phi)
	\end{equation*}
	which does not have the null structure. Even worse, there is a term of the form
	\begin{equation*}
	(\Lbar \log \mu) (\slashed{\D}_L \phi)
	\end{equation*}
	which behaves like $r^{\epsilon} (\slashed{\D}_L \phi)$. An error term of this form cannot be controlled. Note, however, that the additional structure arising from the wave coordinate condition might allow us to control such a term, provided that we can arrange, for example, that $T \log \mu \sim \tau^{-1-\delta}$. In other words, we need sufficient decay in $\tau$ (and no growth in $r$) for this quantity.
	
\end{remark}

\section{Additional error terms in the commutators}

In the expressions for the commutators with a vector field (in proposition \ref{proposition commute vector field}) or with a tensor field (in proposition \ref{proposition commute tensor field}) there are additional error terms which have not been considered in the section above. These error terms involve the curvature of the vector bundle $\mathcal{B}$, and they vanish if $\phi$ is a scalar field with the exception of a single term which appears when commuting with $r\slashed{\nabla}$.

\begin{proposition}[Additional error terms from commuting with $T$]
	\label{proposition additional error term commuting with T}
	Let $\phi$ be an $S_{\tau,r}$-tangent tensor field. Then the additional error terms when commuting with $T$ are given schematically by
	\begin{equation}
	\begin{split}
	&[\slashed{\D}_\mu \, , \slashed{\D}_\nu] \left(T^\nu \slashed{\D}^\mu \phi \right)
	+ \omega \Lbar^\mu T^\nu [\slashed{\D}_\mu \, , \slashed{\D}_\nu] \phi
	+ \slashed{\D}^\mu \left( T^\nu [\slashed{\D}_\mu \, , \slashed{\D}_\nu]\phi \right)
	\\
	&=
	\begin{pmatrix}
	r^{-1} \left(\slashed{\D}(\mathscr{Z} h)\right)_{(\text{frame})} \\
	\bm{\Gamma}_{(-1)}\cdot\bm{\Gamma}_{(-1, \text{small})} \\
	\end{pmatrix} \cdot \slashed{\D} \phi
	+ \begin{pmatrix}
	r^{-1} \left( \tilde{\slashed{\Box}}_g (\mathscr{Z} h) \right)_{(\text{frame})} \\
	r^{-1} \left(\slashed{\D}(\mathscr{Z}^2 h)\right)_{(\text{frame})} \\
	\bm{\Gamma}_{(-1)} \cdot (\slashed{\D} \mathscr{Z} h)_{(\text{frame})} \\
	\bm{\Gamma}_{(-1)} \cdot (\tilde{\Box}_g h)_{(\text{frame})} \\
	\bm{\Gamma}_{(-1)} \cdot \mathscr{Z} \chi_{(\text{small})} \\
	r^{-1}\bm{\Gamma}_{(-1)} \mathscr{Z}^2 \log \mu \\
	\bm{\Gamma}_{(-1)}\cdot\bm{\Gamma}_{(-1)}\cdot\bm{\Gamma}_{(-1, \text{small})}
	\end{pmatrix} \cdot \phi
	\end{split}
	\end{equation}
	where additional factors of $r^{-1}$ may be present. Additionally, if $\phi$ is in fact a scalar field, then these terms vanish identically.

\end{proposition}

\begin{proof}
	We have
	\begin{equation*}
	[\slashed{\D}_\mu \, , \slashed{\D}_\nu] \left(T^\nu \slashed{\D}^\mu \phi_{\slashed{\alpha}_1 \ldots \slashed{\alpha}_n} \right)
	=
	T^\nu \Omega_{\slashed{\alpha}_1 \phantom{\slashed{\beta}}\mu\nu}^{\phantom{\slashed{\alpha_1}}\slashed{\beta}} \slashed{\D}^\mu \phi_{\slashed{\beta} \slashed{\alpha}_2 \ldots \slashed{\alpha}_n}
	+ \ldots
	+ T^\nu \Omega_{\slashed{\alpha}_n \phantom{\slashed{\beta}}\mu\nu}^{\phantom{\slashed{\alpha_n}}\slashed{\beta}} \slashed{\D}^\mu \phi_{\slashed{\alpha}_1 \ldots \slashed{\alpha}_{n-1} \slashed{\beta}}
	\end{equation*}
	So, using proposition \ref{proposition TOmega rOmega and rLOmega} we have, schematically,
	\begin{equation*}
	[\slashed{\D}_\mu \, , \slashed{\D}_\nu] \left(T^\nu \slashed{\D}^\mu \phi\right)
	=
	\left( r^{-1} \left(\slashed{\D}(\mathscr{Z} h)\right)_{(\text{frame})} 
	+ \bm{\Gamma}_{(-1)}\cdot\bm{\Gamma}_{(-1, \text{small})} \right) \cdot \slashed{\D} \phi
	\end{equation*}
	Similarly, we find that, schematically,
	\begin{equation*}
	\omega \Lbar^\mu T^\nu [\slashed{\D}_\mu \, , \slashed{\D}_\nu] \phi
	=
	\left( r^{-1} \bm{\Gamma}_{(-1, \text{small})} \left(\slashed{\D}(\mathscr{Z} h)\right)_{(\text{frame})} 
	+ \bm{\Gamma}_{(-1)}\cdot\bm{\Gamma}_{(-1)}\cdot\bm{\Gamma}_{(-1, \text{small})} \right) \cdot \phi
	\end{equation*}
	
	Finally, using proposition \ref{proposition div Z Omega} we have
	\begin{equation*}
	\slashed{\D}^\mu \left( T^\nu [\slashed{\D}_\mu \, , \slashed{\D}_\nu]\phi \right)
	=
	\begin{pmatrix}
	r^{-1} \left(\slashed{\D}(\mathscr{Z} h)\right)_{(\text{frame})} \\
	\bm{\Gamma}_{(-1)}\cdot\bm{\Gamma}_{(-1, \text{small})} \\
	\end{pmatrix} \cdot \slashed{\D} \phi
	+ \begin{pmatrix}
	r^{-1} \left( \tilde{\slashed{\Box}}_g (\mathscr{Z} h) \right)_{(\text{frame})} \\
	r^{-1} \left(\slashed{\D}(\mathscr{Z}^2 h)\right)_{(\text{frame})} \\
	\bm{\Gamma}_{(-1)} \cdot (\slashed{\D} \mathscr{Z} h)_{(\text{frame})} \\
	\bm{\Gamma}_{(-1)} \cdot (\tilde{\Box}_g h)_{(\text{frame})} \\
	\bm{\Gamma}_{(-1)} \cdot \mathscr{Z} \chi_{(\text{small})} \\
	r^{-1} \mathscr{Z}^2 \log \mu \\
	\bm{\Gamma}_{(-1)}\cdot\bm{\Gamma}_{(-1)}\cdot\bm{\Gamma}_{(-1, \text{small})}
	\end{pmatrix} \cdot \phi
	\end{equation*}
	where there may be additional factors of $r^{-1}$.

\end{proof}

\begin{proposition}[Additional error terms from commuting with $r\slashed{\nabla}$]
	\label{proposition additional error term commuting with rnabla}
	Let $\phi$ be an $S_{\tau,r}$-tangent tensor field. Then the additional error terms when commuting with $r\slashed{\nabla}$ are given schematically by
	\begin{equation}
	\begin{split}
	&[\slashed{\D}_\mu \, , \slashed{\D}_\nu] \left(r\slashed{\Pi}_{\slashed{\alpha}}^{\phantom{\slashed{\alpha}}\nu} \slashed{\D}^\mu \phi \right)
	+ r\omega \Lbar^\mu \slashed{\Pi}_{\slashed{\alpha}}^{\phantom{\slashed{\alpha}}\nu} [\slashed{\D}_\mu \, , \slashed{\D}_\nu] \phi
	+ \slashed{\D}^\mu \left( r\slashed{\Pi}_{\slashed{\alpha}}^{\phantom{\slashed{\alpha}}\nu} [\slashed{\D}_\mu \, , \slashed{\D}_\nu]\phi \right)
	\\
	&=
	\begin{pmatrix}
	(\slashed{\nabla} \mathscr{Z} h)_{(\text{frame})} \\
	r\bm{\Gamma}_{(-1)} \cdot \bm{\Gamma}_{(-1)}
	\end{pmatrix} \cdot (\slashed{\D} \phi)
	+ \begin{pmatrix}
	r^{-1} \slashed{\D}(\mathscr{Z}^2 h)_{(\text{frame})} \\
	\tilde{\slashed{\Box}}_g(\mathscr{Z}h)_{(\text{frame})} \\
	\bm{\Gamma}_{(-1)}\cdot\left(\slashed{\D}(\mathscr{Z}h) \right)_{(\text{frame})} \\
	r\bm{\Gamma}_{(-1)}\cdot (\tilde{\Box}_g h)_{(\text{frame})} \\
	r^{-1} \bm{\Gamma}_{(-1)}\cdot \mathscr{Z}^2 \log \mu \\
	\bm{\Gamma}_{(-1)}\cdot\mathscr{Z} \chi_{(\text{small})} \\
	r\bm{\Gamma}_{(-1)}\cdot\bm{\Gamma}_{(-1)} \cdot \bm{\Gamma}_{(-1)}
	\end{pmatrix} \cdot \phi
	\end{split}
	\end{equation}
	where additional factors of $r^{-1}$ may be present.
	
	If $\phi$ is a scalar field, then the only nonzero terms are given schematically by
	\begin{equation*}
	[\slashed{\D}_\mu \, , \slashed{\D}_\nu](\slashed{\Pi}_{\slashed{\alpha}}^{\phantom{\slashed{\alpha}}\nu} \D^\mu \phi)
	=
	(\slashed{\nabla} \mathscr{Z} h)_{(\text{frame})} \cdot (\D \phi)
	+ r\bm{\Gamma}_{(-1)} \cdot \bm{\Gamma}_{(-1)} (\D\phi)
	\end{equation*}
	
\end{proposition}

\begin{proof}
	We have
	\begin{equation*}
	\begin{split}
	[\slashed{\D}_\mu \, , \slashed{\D}_\nu] \left(r\slashed{\Pi}_{\slashed{\alpha}}^{\phantom{\slashed{\alpha}}\nu} \slashed{\D}^\mu \phi_{\slashed{\alpha}_1 \ldots \slashed{\alpha}_n} \right)
	&=
	r\slashed{\Pi}_{\slashed{\beta}}^{\phantom{\slashed{\beta}}\nu} \Omega_{\slashed{\alpha}\phantom{\slashed{\beta}}\mu\nu}^{\phantom{\slashed{\alpha}}\slashed{\beta}} \slashed{\D}^\mu \phi_{\slashed{\alpha}_1 \ldots \slashed{\alpha}_n}
	+ r\slashed{\Pi}_{\slashed{\alpha}}^{\phantom{\slashed{\alpha}}\nu} \Omega_{\slashed{\alpha_1}\phantom{\slashed{\beta}}\mu\nu}^{\phantom{\slashed{\alpha_1}}\slashed{\beta}} \slashed{\D}^\mu \phi_{\slashed{\beta} \slashed{\alpha}_{2} \ldots \slashed{\alpha}_{n}}
	+ \ldots
	\\ 
	&\phantom{=}
	\ldots
	+ r\slashed{\Pi}_{\slashed{\alpha}}^{\phantom{\slashed{\alpha}}\nu} \Omega_{\slashed{\alpha}_n\phantom{\slashed{\beta}}\mu\nu}^{\phantom{\slashed{\alpha}_n}\slashed{\beta}} \slashed{\D}^\mu \phi_{\slashed{\alpha}_1 \ldots \slashed{\alpha}_{n-1} \slashed{\beta}}
	\\ 
	&=
	r\Omega_{\slashed{\alpha}\phantom{\slashed{\beta}}\mu\slashed{\beta}}^{\phantom{\slashed{\alpha}}\slashed{\beta}} \slashed{\D}^\mu \phi_{\slashed{\alpha}_1 \ldots \slashed{\alpha}_n}
	+ r \Omega_{\slashed{\alpha}_1 \phantom{\slashed{\beta}}\mu\slashed{\alpha}}^{\phantom{\slashed{\alpha_1}}\slashed{\beta}} \slashed{\D}^\mu \phi_{\slashed{\beta} \slashed{\alpha}_{2} \ldots \slashed{\alpha}_{n}}
	+ \ldots
	+ r \Omega_{\slashed{\alpha}_n \phantom{\slashed{\beta}}\mu\slashed{\alpha}}^{\phantom{\slashed{\alpha_n}}\slashed{\beta}} \slashed{\D}^\mu \phi_{\slashed{\alpha}_1 \ldots \slashed{\alpha}_{n-1} \slashed{\beta}} 
	\end{split}
	\end{equation*}
	So, using proposition \ref{proposition TOmega rOmega and rLOmega} we have, schematically,
	\begin{equation*}
	\begin{split}
	[\slashed{\D}_\mu \, , \slashed{\D}_\nu] \left(r\slashed{\Pi}_{\slashed{\alpha}}^{\phantom{\slashed{\alpha}}\nu} \slashed{\D}^\mu \phi_{\slashed{\alpha}_1 \ldots \slashed{\alpha}_n} \right)
	&=
	(\slashed{\nabla} \mathscr{Z} h)_{(\text{frame})} \cdot (\slashed{\D} \phi)
	+ r\bm{\Gamma}_{(-1)} \cdot \bm{\Gamma}_{(-1)} (\slashed{\D}\phi)
	\end{split}
	\end{equation*}
	Similarly, we have
	\begin{equation*}
	 r\omega \Lbar^\mu \slashed{\Pi}_{\slashed{\alpha}}^{\phantom{\slashed{\alpha}}\nu} [\slashed{\D}_\mu \, , \slashed{\D}_\nu] \phi
	 =
	 \bm{\Gamma}_{(-1)} \cdot (\slashed{\nabla} \mathscr{Z} h)_{(\text{frame})} \cdot \phi
	 + r\bm{\Gamma}_{(-1)} \cdot \bm{\Gamma}_{(-1)} \cdot \bm{\Gamma}_{(-1)} \cdot \phi	 
	\end{equation*}
	Finally, using proposition \ref{proposition div Z Omega} we have
	\begin{equation*}
	\slashed{\D}^\mu \left( r\slashed{\Pi}_{\slashed{\alpha}}^{\phantom{\slashed{\alpha}}\nu} [\slashed{\D}_\mu \, , \slashed{\D}_\nu]\phi \right)
	=
	\begin{pmatrix}
	(\slashed{\nabla} \mathscr{Z} h)_{(\text{frame})} \\
	r\bm{\Gamma}_{(-1)} \cdot \bm{\Gamma}_{(-1)}
	\end{pmatrix} \cdot (\slashed{\D} \phi)
	+ \begin{pmatrix}
	r^{-1} \slashed{\D}(\mathscr{Z}^2 h)_{(\text{frame})} \\
	\tilde{\slashed{\Box}}_g(\mathscr{Z}h)_{(\text{frame})} \\
	\bm{\Gamma}_{(-1)}\cdot\left(\slashed{\D}(\mathscr{Z}h) \right)_{(\text{frame})} \\
	r\bm{\Gamma}_{(-1)}\cdot (\tilde{\Box}_g h)_{(\text{frame})} \\
	r^{-1} \bm{\Gamma}_{(-1)}\cdot \mathscr{Z}^2 \log \mu \\
	\bm{\Gamma}_{(-1)}\cdot\mathscr{Z} \chi_{(\text{small})} \\
	r\bm{\Gamma}_{(-1)}\cdot\bm{\Gamma}_{(-1)} \cdot \bm{\Gamma}_{(-1)}
	\end{pmatrix} \cdot \phi
	\end{equation*}
	where additional factors of $r^{-1}$ might be present.
	
	If $\phi$ is a scalar field, then the only nonzero terms in all of the above calculations are
	\begin{equation*}
	r\Omega_{\slashed{\alpha}\phantom{\slashed{\beta}}\mu\slashed{\beta}}^{\phantom{\slashed{\alpha}}\slashed{\beta}}\D^\mu \phi
	\end{equation*}
\end{proof}

\begin{proposition}[Additional error terms from commuting with $rL$]
	\label{proposition additional error term commuting with rL}
	Let $\phi$ be an $S_{\tau,r}$-tangent tensor field. Then the additional error terms when commuting with $rL$ are given schematically by
	\begin{equation}
	\begin{split}
	&[\slashed{\D}_\mu \, , \slashed{\D}_\nu] \left(rL^\nu \slashed{\D}^\mu \phi \right)
	+ \omega \Lbar^\mu rL^\nu [\slashed{\D}_\mu \, , \slashed{\D}_\nu] \phi
	+ \slashed{\D}^\mu \left( rL^\nu [\slashed{\D}_\mu \, , \slashed{\D}_\nu]\phi \right)
	\\
	&=
	\begin{pmatrix}
	\left(\slashed{\D}(\mathscr{Z} h)\right)_{(\text{frame})} \\
	r\bm{\Gamma}_{(-1)}\cdot\bm{\Gamma}_{(-1, \text{small})} \\
	\end{pmatrix} \cdot \overline{\slashed{\D}} \phi
	+ \begin{pmatrix}
	r^{-1} \left( \slashed{\D}_L \left( r \slashed{\D}_L \mathscr{Z} h \right) \right)_{(\text{frame})} \\
	\bm{\Gamma}_{(-1)} \left( L \left( r L h \right) \right)_{(\text{frame})} \\
	(\tilde{\slashed{\Box}}_g \mathscr{Z} h)_{(\text{frame})} \\
	r\bm{\Gamma}_{(-1)} (\tilde{\Box}_g h)_{(\text{frame})} \\
	r^{-1} \left( \slashed{\D} \mathscr{Z}^2 h \right)_{(\text{frame})} \\
	\bm{\Gamma}_{(-1)} \left( \slashed{\D} \mathscr{Z} h \right)_{(\text{frame})} \\
	\bm{\Gamma}_{(-1)} \mathscr{Z} \chi_{(\text{small})} \\
	r^{-1} \bm{\Gamma}_{(-1)} \mathscr{Z}^2 \log \mu \\
	r \bm{\Gamma}_{(-1)} \bm{\Gamma}_{(-1)} \bm{\Gamma}_{(-1, \text{small})}
	\end{pmatrix} \cdot \phi
	\end{split}
	\end{equation}
\end{proposition}

\begin{proof}
	We have
	\begin{equation*}
	\begin{split}
	[\slashed{\D}_\mu \, , \slashed{\D}_\nu] \left(rL^\nu \slashed{\D}^\mu \phi_{\slashed{\alpha}_1 \ldots \slashed{\alpha}_n} \right)
	&=
	rL^\nu \Omega_{\slashed{\alpha}_1 \phantom{\slashed{\beta}}\mu\nu}^{\phantom{\slashed{\alpha_1}}\slashed{\beta}} \slashed{\D}^\mu \phi_{\slashed{\beta} \slashed{\alpha}_2 \ldots \slashed{\alpha}_n}
	+ \ldots
	+ rL^\nu \Omega_{\slashed{\alpha}_n \phantom{\slashed{\beta}}\mu\nu}^{\phantom{\slashed{\alpha_n}}\slashed{\beta}} \slashed{\D}^\mu \phi_{\slashed{\alpha}_1 \ldots \slashed{\alpha}_{n-1} \slashed{\beta}}
	\\
	& =
	-\frac{1}{2} r \Omega_{\slashed{\alpha}_1 \phantom{\slashed{\beta}}L\Lbar}^{\phantom{\slashed{\alpha_1}}\slashed{\beta}} \slashed{\D}_L \phi_{\slashed{\beta} \slashed{\alpha}_2 \ldots \slashed{\alpha}_n}
	+ r\Omega_{\slashed{\alpha}_1 \phantom{\slashed{\beta}}L}^{\phantom{\slashed{\alpha_1}}\slashed{\beta} \phantom{L} \slashed{\gamma}} \slashed{\nabla}_{\slashed{\gamma}} \phi_{\slashed{\beta} \slashed{\alpha}_2 \ldots \slashed{\alpha}_n}
	+ \ldots
	\\
	&\phantom{=} \ldots
	-\frac{1}{2} r \Omega_{\slashed{\alpha_n} \phantom{\slashed{\beta}}L\Lbar}^{\phantom{\slashed{\alpha_n}}\slashed{\beta}} \slashed{\D}_L \phi_{\slashed{\alpha}_1 \ldots \slashed{\alpha}_{n-1} \slashed{\beta}}
	+ r\Omega_{\slashed{\alpha_n} \phantom{\slashed{\beta}}L}^{\phantom{\slashed{\alpha_n}}\slashed{\beta} \phantom{L} \slashed{\gamma}} \slashed{\nabla}_{\slashed{\gamma}} \phi_{\slashed{\alpha}_1 \ldots \slashed{\alpha}_{n-1} \slashed{\beta}}
	\end{split}
	\end{equation*}
	So, using proposition \ref{proposition TOmega rOmega and rLOmega} we have, schematically,
	\begin{equation*}
	[\slashed{\D}_\mu \, , \slashed{\D}_\nu] \left(rL^\nu \slashed{\D}^\mu \phi\right)
	=
	\left( \left(\slashed{\D}(\mathscr{Z} h)\right)_{(\text{frame})} 
	+ r\bm{\Gamma}_{(-1)}\cdot\bm{\Gamma}_{(-1, \text{small})} \right) \cdot \overline{\slashed{\D}} \phi
	\end{equation*}

	Similarly, we find that, schematically
	\begin{equation*}
	r\omega \Lbar^\mu L^\nu [\slashed{\D}_\mu \, , \slashed{\D}_\nu] \phi
	=
	\omega \left( \left(\slashed{\D}(\mathscr{Z} h)\right)_{(\text{frame})} 
	+ r\bm{\Gamma}_{(-1)}\cdot\bm{\Gamma}_{(-1, \text{small})} \right) \cdot \phi
	\end{equation*}
	
	Next, we consider
	\begin{equation*}
	\slashed{\D}^\mu \left(rL^\nu [\slashed{\D}_\mu \, , \slashed{\D}_\nu] \phi \right)
	\end{equation*}
	We have
	\begin{equation*}
	\slashed{\D}^\mu \left( r L^\nu \Omega_{\slashed{\alpha} \phantom{\slashed{\beta}} \mu \nu}^{\phantom{\slashed{\alpha}}\slashed{\beta}} \right)
	=
	\frac{1}{2}\slashed{\D}_L \left( r\Omega_{\slashed{\alpha} \phantom{\slashed{\beta}} L\Lbar}^{\phantom{\slashed{\alpha}}\slashed{\beta}} \right)
	+ r \slashed{\nabla}^{\slashed{\gamma}} \Omega_{\slashed{\alpha} \phantom{\slashed{\beta}} \slashed{\gamma} L}^{\phantom{\slashed{\alpha}}\slashed{\beta}}
	\end{equation*}
	Schematically, this is
	\begin{equation*}
	\slashed{\D}_L \left( \left(\slashed{\D}(\mathscr{Z} h)\right)_{(\text{frame})} 
	+ r\bm{\Gamma}_{(-1)}\cdot\bm{\Gamma}_{(-1, \text{small})} \right)
	+ r\slashed{\nabla} \left( r^{-1}\left(\slashed{\D}(\mathscr{Z} h)\right)_{(\text{frame})} 
	+ \bm{\Gamma}_{(-1)}\cdot\bm{\Gamma}_{(-1, \text{small})} \right)
	\end{equation*}
	Recalling the definitions of $\bm{\Gamma}_{(-1)}$ and $\bm{\Gamma}_{(-1,\text{small})}$ we have
	\begin{equation*}
	r\slashed{\nabla} \bm{\Gamma}_{(-1)}
	=
	\begin{pmatrix}
	(\slashed{\D} \mathscr{Z} h)_{(\text{frame})} \\
	\mathscr{Z} \chi_{(\text{small})} \\
	r^{-1} \mathscr{Z}^2 \log \mu \\
	r\bm{\Gamma}_{(-1)} \cdot \bm{\Gamma}_{(-1, \text{small})}
	\end{pmatrix}
	\end{equation*}
	and
	\begin{equation*}
	r\slashed{\nabla} \left(\slashed{\D}(\mathscr{Z}h)\right)_{(\text{frame})}
	=
	\left(\slashed{\D} (\mathscr{Z}^2 h) \right)_{(\text{frame})}
	+ r\bm{\Gamma}_{(-1)} \left(\slashed{\D} (\mathscr{Z} h) \right)_{(\text{frame})}
	\end{equation*}
	
	Next, we note that
	\begin{equation*}
	\begin{split}
	\slashed{\D}_L \left( \slashed{\D} (\mathscr{Z} h) \right)_{(\text{frame})}
	&=
	\begin{pmatrix}
	\left( \slashed{\D}_L \slashed{\D}_{\Lbar} \mathscr{Z} h \right)_{(\text{frame})} \\
	\left( \slashed{\D}_L \slashed{\nabla} \mathscr{Z} h \right)_{(\text{frame})} \\
	\left( \slashed{\D}_L \slashed{\D}_L \mathscr{Z} h \right)_{(\text{frame})} \\
	\bm{\Gamma}_{(-1)} \left( \slashed{\D} (\mathscr{Z} h) \right)_{(\text{frame})}
	\end{pmatrix}
	\\
	&=
	\begin{pmatrix}
	\left( \tilde{\slashed{\Box}}_g \mathscr{Z} h \right)_{(\text{frame})} \\
	r^{-1} \left( \overline{\slashed{\D}} \mathscr{Z}^2 h \right)_{(\text{frame})} \\
	r^{-1} \left( \slashed{\D}_L \left( r\slashed{\D}_L \mathscr{Z} h \right) \right)_{(\text{frame})} \\
	\bm{\Gamma}_{(-1)} \left( \slashed{\D} (\mathscr{Z} h) \right)_{(\text{frame})}
	\end{pmatrix}
	\end{split}
	\end{equation*}
	and also
	\begin{equation*}
	\slashed{\D}_L \bm{\Gamma}_{(-1)}
	=
	\begin{pmatrix}
	r^{-2} \\
	(L \partial h)_{(\text{frame})} \\
	\slashed{\D}_L \chi \\
	r^{-1} \bm{\Gamma}_{(-1)}
	\end{pmatrix}
	\end{equation*}
	Now, using proposition \ref{proposition transport chi chismall} we have
	\begin{equation*}
	\slashed{\D}_L \chi
	=
	\begin{pmatrix}
	\alpha \\
	r^{-1} \omega \\
	r^{-1} \chi_{(\text{small})} \\
	\omega \chi_{(\text{small})} \\
	\chi_{(\text{small})} \chi_{(\text{small})}
	\end{pmatrix}	
	\end{equation*}
	and we can estimate
	\begin{equation*}
	\alpha = r^{-1} \left(L (r L h) \right)_{(\text{frame})}
	+ r^{-1} \left( \overline{\slashed{\D}} \mathscr{Z} h \right)_{(\text{frame})}
	+ \bm{\Gamma}_{(-1)}\bm{\Gamma}_{(-1, \text{small})}
	\end{equation*}
	
	In summary, we have, schematically,
	\begin{equation*}
	\slashed{\D}^\mu \left( r L^\nu \Omega_{\slashed{\alpha} \phantom{\slashed{\beta}} \mu \nu}^{\phantom{\slashed{\alpha}}\slashed{\beta}} \right)
	=
	\begin{pmatrix}
	r^{-1} \left( \slashed{\D}_L \left( r \slashed{\D}_L \mathscr{Z} h \right) \right)_{(\text{frame})} \\
	\bm{\Gamma}_{(-1)} \left( L \left( r L h \right) \right)_{(\text{frame})} \\
	(\tilde{\slashed{\Box}}_g \mathscr{Z} h)_{(\text{frame})} \\
	r\bm{\Gamma}_{(-1)} (\tilde{\Box}_g h)_{(\text{frame})} \\
	r^{-1} \left( \slashed{\D} \mathscr{Z}^2 h \right)_{(\text{frame})} \\
	\bm{\Gamma}_{(-1)} \left( \slashed{\D} \mathscr{Z} h \right)_{(\text{frame})} \\
	\bm{\Gamma}_{(-1)} \mathscr{Z} \chi_{(\text{small})} \\
	r^{-1} \bm{\Gamma}_{(-1)} \mathscr{Z}^2 \log \mu \\
	r \bm{\Gamma}_{(-1)} \bm{\Gamma}_{(-1)} \bm{\Gamma}_{(-1, \text{small})}
	\end{pmatrix}
	\end{equation*}

\end{proof}

\section{Notation for commuted fields}

Given a field $\phi$ we use the following notation for repeatedly commuted fields: we write $\mathscr{Z}^N \phi$ to mean
\begin{equation}
 \mathscr{Z}^N \phi := \mathscr{Z}_{(1)}\ldots\mathscr{Z}_{(N)} \phi
\end{equation}
where, for each $n$, $\mathscr{Z}_{(n)}$ is either $\slashed{\D}_{T}$ or $r \slashed{\nabla}$. Similarly, we write
\begin{equation}
\mathscr{Y}^N \phi := \mathscr{Y}_{(1)}\ldots\mathscr{Y}_{(N)} \phi
\end{equation}
where, for each $n$, $\mathscr{Y}_{(n)}$ is either $\slashed{\D}_{T}$, $r \slashed{\nabla}$, \emph{or} $r\slashed{\D}_L$.

Note that the rank of these quantities depends on the number of times $r\slashed{\nabla}$ appears in the commuting operators, but if $\phi$ is a scalar field the maximum rank of $\mathscr{Z}^N \phi$ or $\mathscr{Y}^N \phi$ is $N$. We define
\begin{equation}
 |\mathscr{Z}^N \phi| := \sum |\mathscr{Z}_{(1)} \ldots \mathscr{Z}_{(N)} \phi|
\end{equation}
where the sum is taken over all ordered sets of commutation operators $(\mathscr{Z}_{(1)}, \ldots, \mathscr{Z}_{(N)})$, where for all $n \leq N$, we have $\mathscr{Z}_{(n)} \in \{ \slashed{\D}_{T}, r\slashed{\nabla} \}$. Similarly, we can define $\slashed{\D}\mathscr{Z}^N \phi$ and $|\slashed{\D}\mathscr{Z}^N \phi|$ for the derivatives of the commuted field $\phi$, and finally we can define $\overline{\slashed{\D}}\mathscr{Z}^N \phi$ and $|\overline{\slashed{\D}}\mathscr{Z}^N \phi|$ for the \emph{good} derivatives of the commuted field. The same notation will be used for the set $\mathscr{Y}$ instead of the set $\mathscr{Z}$, with the obvious modification that we allow operators of the form $r\slashed{\D}_L$ to appear.

\section{Commuted equations for geometric quantities}

We will need to establish suitable equations which will allow us to estimate various geometric quantities (such as the foliation density $\mu$, the connection coefficients $\bm{\Gamma}$ etc.) \emph{after} having commuted some number of times with the operators $\mathscr{Y} \in \{ \slashed{\D}_T , r\slashed{\nabla} , r\slashed{\D}_L \}$. These are presented in this section.

In the previous sections, we have only encountered the connection coefficients $\bm{\Gamma}$ or a low number of derivatives of these. Since these quantities have only appeared with a small number of derivatives in the previous expressions, they will mostly be estimated in $L^{\infty}$, in which case our previous division of the connection coefficients into the sets $\bm{\Gamma}$ and $\bm{\Gamma}_{(\text{good})}$ was sufficient for our purposes. However, from this point onwards we need a more detailed breakdown of the error terms.

\begin{definition}[Detailed error terms]
We define the following sets, which will be used schematically from now on:
\begin{equation}
 \begin{split}
  \tilde{\bm{\Gamma}}_{(-1, \text{large})}
  &:= 
  \left\{
  r^{-1}  \ , \
  r^{-1} X_{(\text{frame})} \ , \
  r^{-1} x^a
  \right\}
  \\ \\
  \tilde{\bm{\Gamma}}_{(0, \text{large})}
  &:= 
  \left\{
  1 \ , \
  X_{(\text{frame})}
  \right\} \cup \tilde{\bm{\Gamma}}_{(-1, \text{large})}
  \\ \\
  \tilde{\bm{\Gamma}}_{(C_{(n)}\epsilon, \text{large})}^{(n)}
  &:=
  \left\{
  \mathscr{Y}^{\leq n} X_{(\text{frame})}
  \right\} \cup \tilde{\bm{\Gamma}}_{(0, \text{large})}
  \\ \\
  \tilde{\bm{\Gamma}}_{(-1-\frac{3}{2}\delta)}^{(n)} 
  &:= 
  \begin{Bmatrix} (\overline{\mathscr{\D}} \mathscr{Y}^{\leq n} h)_{(\text{frame})} \ , \
  \overline{\mathscr{\D}} \mathscr{Y}^{\leq n} h_{(rect)} \ , \
  \mathscr{Y}^{\leq n}\tr_{\slashed{g}} \chi_{(\text{small})} \vspace{4mm} \\
  \mathscr{Y}^{\leq n}\hat{\chi} \ , \
  r^{-1}\bar{X}_{(\text{frame})}
  \end{Bmatrix}
  \\ \\ 
  \tilde{\bm{\Gamma}}_{(-1)}^{(n)}
  &:= 
  \left\{
  (\partial h)_{LL} \ , \ 
  \omega \ , \
  r^{-1}X_{(\text{frame, small})}
  \right\} \cup \tilde{\bm{\Gamma}}_{(-1-2\delta)}^{(n)}
  \\ \\
  \tilde{\bm{\Gamma}}_{(-1 + C_{(n)} \epsilon)}^{(n)}
  &:=
  \begin{Bmatrix}
  (\slashed{\D} \mathscr{Y}^{\leq n} h)_{(\text{frame})} \ , \
  \slashed{\D} \mathscr{Y}^{\leq n} h_{(\text{rect})} \ , \
  \mathscr{Y}^{\leq n} \tr_{\slashed{g}}\chibar_{(\text{small})} \ , \
  \mathscr{Y}^{\leq n} \hat{\chibar} \vspace{4mm} \\
  \mathscr{Y}^{\leq n-1} \slashed{\nabla} \log \mu \ , \
  r^{-1} \mathscr{Y}^{\leq n} X_{(\text{frame, small})}
  \end{Bmatrix} \cup \tilde{\bm{\Gamma}}_{(-1)}^{(n)}
  \\ \\
  \tilde{\bm{\Gamma}}_{\left(-\frac{1}{2} + \delta \right)}^{(n)}
  &:=
  \left\{
  (\mathscr{Y}^{\leq n} h)_{(\text{frame})} \ , \
  \mathscr{Y}^{\leq n} h_{(\text{rect})}
  \right\} \cup \tilde{\bm{\Gamma}}_{(-1 + C_{(n)} \epsilon)}^{(n)}
  \\ \\
  \tilde{\bm{\Gamma}}_{(C_{(n)}\epsilon)}^{(n)}
  &:=
  \left\{
  \mathscr{Y}^{\leq n} \log \mu
  \right\} \cup 
  \tilde{\bm{\Gamma}}_{\left(-\frac{1}{2} + \delta \right)}^{(n)} \cup
  \left(r \tilde{\bm{\Gamma}}_{(-1 + C_{(n)}\epsilon)}^{(n)} \right)
 \end{split}
\end{equation}

Now, from these basic ingredients we build up the sets $\bm{\Gamma}^{(n)}_{(a)}$ by applying the following set of rules:
\begin{equation*}
 \begin{split}
  \tilde{\bm{\Gamma}}^{(n)}_{(a)} &\subset \bm{\Gamma}^{(m)}_{(b)} \text{\quad \quad if } n \leq m \text{ and } a \leq b \\
  \tilde{\bm{\Gamma}}^{(n)}_{(a, \text{large})} &\subset \bm{\Gamma}^{(m)}_{(b, \text{large})} \text{\quad \quad if } n \leq m \text{ and } a \leq b \\
  \tilde{\bm{\Gamma}}^{(n_1)}_{(a)} \cdot \tilde{\bm{\Gamma}}^{(n_2)}_{(b)} &\subset \bm{\Gamma}^{(\max\{n_1, n_2\})}_{(a+b)} \\
  \tilde{\bm{\Gamma}}^{(n_1)}_{(a)} \cdot \tilde{\bm{\Gamma}}^{(n_2)}_{(b, \text{large})} &\subset \bm{\Gamma}^{(\max\{n_1, n_2\})}_{(a+b)} \\
  \tilde{\bm{\Gamma}}^{(n_1)}_{(a, \text{large})} \cdot \tilde{\bm{\Gamma}}^{(n_2)}_{(b, \text{large})} &\subset \bm{\Gamma}^{(\max\{n_1, n_2\})}_{(a+b, \text{large})}
 \end{split}
\end{equation*}

The idea is that the upper index indicates the number commutation operators appearing, while the lower index indicates the expected behaviour as a function of $r$. So, for example, we would expect to be able to estimate the quantities in $\bm{\Gamma}^{(2)}_{-1-\delta}$ in terms of quantities involving two $\mathscr{Y}$ operators, and the expected behaviour is $\bm{\Gamma}^{(2)}_{-1-\delta} \sim \epsilon r^{-1-\delta}$. Those quantities which are  labelled ``large'' behave the same way, except that these quantities are not expected to come with the small factor of $\epsilon$.

Note that, since $C_{(n)} \gg C_{(n-1)}$ we have, for example,
\begin{equation*}
 \bm{\Gamma}^{(n)}_{(C_{(n)}\epsilon)} \cdot \bm{\Gamma}^{(n)}_{(C_{(n)}\epsilon)} \subset \bm{\Gamma}^{(n)}_{(C_{(n+1)}\epsilon)}
\end{equation*}

Note also that the reason for including the factor $\frac{3}{2}$ in the definition of the term $\bm{\Gamma}^{(n)}_{(-1-\frac{3}{2}\delta)}$ so that we have, for example,
\begin{equation*}
\bm{\Gamma}^{(n)}_{(-1-\frac{3}{2}\delta)} \bm{\Gamma}^{(n)}_{(C_{(n)}\epsilon, \text{large})} = \bm{\Gamma}^{(n)}_{(-1-\delta)}
\end{equation*}
In other words, we will not have to count factors of $\epsilon$ when dealing with error terms of this form. However, this will require us to eventually improve the rate of decay in the bootstrap bounds so that, for example, $(\bar{\partial} h)_{(\text{rect})} \sim r^{-1-\frac{3}{2}\delta}$.

\end{definition}

We now need some preparatory calculations, which allow us to commute with the operators $\mathscr{Z}$ and $\mathscr{Y}$ an arbitrary number of times.

\begin{proposition}
 \label{proposition commute Zn D}
Let $\phi$ be an $S_{\tau,r}$-tangent tensor field, and suppose that $C_{(n)} \gg C_{(n-1)}$. Then we have, schematically,
\begin{equation*}
  [\mathscr{Z}^n, \slashed{\D}] \phi
 =
 \sum_{j+k \leq n-1} \bm{\Gamma}_{(C_{(j+1)} \epsilon)}^{(j+1)} \overline{\slashed{\D}} \mathscr{Z}^k \phi 
 + \sum_{j+k \leq n-1} \bm{\Gamma}_{(C_{(j)} \epsilon)}^{(j)} \slashed{\D} \mathscr{Z}^k \phi 
 + \sum_{j+k \leq n-1} \bm{\Gamma}_{(-1 + C_{(j+1)} \epsilon)}^{(j+1)} \mathscr{Z}^k \phi 
\end{equation*}
and
\begin{equation*}
 \begin{split}
 [\mathscr{Z}^n , \overline{\slashed{\D}}] \phi
 &=
 \sum_{j+k \leq n-1} \bm{\Gamma}^{(j)}_{(-1+C_{(j)}\epsilon)} \slashed{\D} \mathscr{Z}^k \phi 
 + \sum_{j+k \leq n-1} \bm{\Gamma}^{(j+1)}_{(-\delta)} \overline{\slashed{\D}} \mathscr{Z}^k \phi 
 + \sum_{j+k \leq n-1} \bm{\Gamma}_{(-1 + C_{(j+1)} \epsilon)}^{(j+1)} \mathscr{Z}^k \phi 
 \\
 &=
 \sum_{j+k \leq n-1} \bm{\Gamma}_{(-\delta)}^{(j)} \slashed{\D} \mathscr{Z}^k \phi 
 + \sum_{j+k \leq n-1} \bm{\Gamma}_{(-1 + C_{(j+1)} \epsilon)}^{(j+1)} \mathscr{Z}^k \phi 
 \end{split}
\end{equation*}

Additionally, if $\phi$ is a scalar field, then we have
\begin{equation*}
[\mathscr{Z}^n, \slashed{\D}] \phi
=
\sum_{j+k \leq n-1} \bm{\Gamma}_{(C_{(j+1)} \epsilon)}^{(j+1)} \overline{\slashed{\D}} \mathscr{Z}^k \phi 
+ \sum_{j+k \leq n-1} \bm{\Gamma}_{(C_{(j)} \epsilon)}^{(j)} \slashed{\D} \mathscr{Z}^k \phi 
+ \sum_{j+k \leq n-1} \bm{\Gamma}_{(-1 + C_{(j)} \epsilon)}^{(j)} \mathscr{Z}^k \phi 
\end{equation*}
and
\begin{equation*}
[\mathscr{Z}^n , \overline{\slashed{\D}}] \phi
=
\sum_{j+k \leq n-1} \bm{\Gamma}_{(-\delta)}^{(j+1)} \slashed{\D} \mathscr{Z}^k \phi 
+ \sum_{j+k \leq n-1} \bm{\Gamma}_{(-1 + C_{(j)} \epsilon)}^{(j)} \mathscr{Z}^k \phi 
\end{equation*}

\end{proposition}

\begin{proof}
 Using propositions \ref{proposition commuting DT with first order operators} and \ref{proposition commuting rnabla with first order operators} as well as proposition \ref{proposition expression for Omega} we have, schematically
\begin{equation*}
\begin{split}
  \mathscr{Z} \slashed{\D} \phi 
  &=
  \slashed{\D} \mathscr{Z} \phi
  + \bm{\Gamma}_{(C_{(0)}\epsilon)}^{(0)} \slashed{\D}\phi
  + \bm{\Gamma}_{(C_{(1)}\epsilon)}^{(1)} \overline{\slashed{\D}}\phi
  + \begin{pmatrix}
     \Omega_{L\Lbar} \\
     r\Omega_{L\slashed{\alpha}} \\
     r\Omega_{\Lbar\slashed{\alpha}} \\
  \end{pmatrix}
  \cdot \phi
  \\
  &=
  \slashed{\D} \mathscr{Z} \phi
  + \bm{\Gamma}_{(C_{(0)}\epsilon)}^{(0)} \slashed{\D}\phi
  + \bm{\Gamma}_{(C_{(1)}\epsilon)}^{(1)} \overline{\slashed{\D}}\phi
  + \begin{pmatrix}
  r^{-1} (\slashed{\D} \mathscr{Z} h)_{(\text{frame})} \\
  (\slashed{\nabla} \mathscr{Z} h)_{(\text{frame})} \\
  r\bm{\Gamma}_{(-1)}\cdot\bm{\Gamma}_{(-1,\text{small})} \\
  \end{pmatrix}
  \cdot \phi
\end{split}
\end{equation*}
where, if $\phi$ is a scalar field, then the terms which are linear in $\phi$ (rather than its derivatives) are absent. Note that, in the case that $\slashed{\D} = \slashed{\nabla}$ we \emph{do not} commute the derivatives, but instead view $r\slashed{\nabla}$ as a commutation operator $\mathscr{Z}$. Hence we avoid picking up a factor of $\bm{\Gamma}^{(0)}_{(-1, \text{large})}$. Making use of proposition \ref{proposition R RL} and the definitions above, we find that we can in fact write
\begin{equation*}
	\mathscr{Z} \slashed{\D} \phi =
	\slashed{\D} \mathscr{Z} \phi
	+ \bm{\Gamma}^{(0)}_{(C_{(0)}\epsilon)} \slashed{\D} \phi
	+ \bm{\Gamma}^{(1)}_{(C_{(1)}\epsilon)} \overline{\slashed{\D}} \phi
	+ \bm{\Gamma}_{(-1 + C_{(1)}\epsilon)}^{(1)} \phi
\end{equation*}

On the other hand, commuting with the ``good'' derivatives, we find
\begin{equation*}
\begin{split}
 \mathscr{Z} \overline{\slashed{\D}} \phi 
 &=
  \overline{\slashed{\D}} \mathscr{Z} \phi
  + \bm{\Gamma}^{(0)}_{(-1 + C_{(0)}\epsilon)} \slashed{\D}\phi
  + \bm{\Gamma}^{(1)}_{(-\delta)} \overline{\slashed{\D}}\phi
  + \begin{pmatrix}
     \Omega_{L\Lbar} \\
     r\Omega_{L\slashed{\alpha}}
  \end{pmatrix}
  \cdot \phi
  \\
 &=
  \overline{\slashed{\D}} \mathscr{Z} \phi
  + \bm{\Gamma}^{(0)}_{(-1 + C_{(0)}\epsilon)} \slashed{\D}\phi
  + \bm{\Gamma}^{(1)}_{(-\delta)} \overline{\slashed{\D}}\phi
  + \begin{pmatrix}
  r^{-1}(\slashed{\D} \mathscr{Z}h)_{(\text{frame})} \\
  (\slashed{\nabla} \mathscr{Z}h)_{(\text{frame})} \\
  \bm{\Gamma}_{(-1)} \cdot \bm{\Gamma}_{(-1)} \\
  r\bm{\Gamma}_{(-1)} \cdot (\bar{\partial}h)_{(\text{frame})}
  \end{pmatrix}
  \cdot \phi
  \\
  &=
  \overline{\slashed{\D}} \mathscr{Z} \phi
  + \bm{\Gamma}^{(0)}_{(-1 + C_{(0)}\epsilon)} \slashed{\D}\phi
  + \bm{\Gamma}^{(1)}_{(-\delta)} \overline{\slashed{\D}}\phi
  + \bm{\Gamma}^{(1)}_{(-1-\delta)} \cdot \phi
\end{split}
\end{equation*}

Now we suppose that, for all $n \leq N$, we have
\begin{equation*}
	\mathscr{Z}^n \slashed{\D} \phi
	=
	\slashed{\D} \mathscr{Z}^n \phi
	+ \sum_{j+k \leq n-1} \bm{\Gamma}_{(C_{(j+1)} \epsilon)}^{(j+1)} \overline{\slashed{\D}} \mathscr{Z}^k \phi
	+ \sum_{j+k \leq n-1} \bm{\Gamma}_{(C_{(j)} \epsilon)}^{(j)} \slashed{\D} \mathscr{Z}^k \phi
	+ \sum_{j+k \leq n-1} \bm{\Gamma}_{(-1 + C_{(j+1)} \epsilon)}^{(j+1)} \mathscr{Z}^k \phi 
\end{equation*}
and
\begin{equation*}
 \mathscr{Z}^n \overline{\slashed{\D}} \phi
 =
 \overline{\slashed{\D}} \mathscr{Z}^n \phi
 + \sum_{j+k \leq n-1} \bm{\Gamma}^{(j)}_{(-1+C_{(j)}\epsilon)} \slashed{\D} \mathscr{Z}^k \phi 
 + \sum_{j+k \leq n-1} \bm{\Gamma}^{(j+1)}_{(-\delta)} \overline{\slashed{\D}} \mathscr{Z}^k \phi
 + \sum_{j+k \leq n-1} \bm{\Gamma}_{(-1 - \delta)}^{(j+1)} \mathscr{Z}^k \phi 
\end{equation*}
These equations evidently hold for $N = 1$, as shown above.

Note that this means that acting with an operator $\mathscr{Z}$ on a term of the form $\bm{\Gamma}^{(n)}_{(-\delta)}$ or $\bm{\Gamma}^{(n)}_{(C_{(n)}\epsilon)}$ serves to add one to the index $n$, at least for $n \leq N-1$. This is either immediate from the definitions of these terms (invlolving $n$ commutation operators) or follows from commuting the commutation operators through $\slashed{\D}$ or $\overline{\slashed{\D}}$. For example,
\begin{equation*}
	\begin{split}
	\mathscr{Z}^n \slashed{\D} h_{(\text{rect})}
	&=
	\slashed{\D} \mathscr{Z}^n h_{(\text{rect})}
	+ \sum_{j+k \leq n-1} \bm{\Gamma}_{(C_{(j+1)} \epsilon)}^{(j+1)} \overline{\slashed{\D}} \mathscr{Z}^k h_{(\text{rect})}
	+ \sum_{j+k \leq n-1} \bm{\Gamma}_{(C_{(j)} \epsilon)}^{(j)} \slashed{\D} \mathscr{Z}^k h_{(\text{rect})}
	\\
	&\phantom{=}
	+ \sum_{j+k \leq n-1} \bm{\Gamma}_{(-1 + C_{(j+1)} \epsilon)}^{(j+1)} \mathscr{Z}^k h_{(\text{rect})}
	\\
	\\
	&=
	\bm{\Gamma}^{(n)}_{(-1 + C_{(n)}\epsilon)}
	+ \sum_{j+k \leq n-1} \bm{\Gamma}_{(C_{(j+1)} \epsilon)}^{(j+1)} \bm{\Gamma}^{(k)}_{(-1-\frac{3}{2}\delta)}
	+ \sum_{j+k \leq n-1} \bm{\Gamma}_{(C_{(j)} \epsilon)}^{(j)} \bm{\Gamma}_{(-1 + C_{(k)} \epsilon)}^{(k)}
	\\
	&\phantom{=}
	+ \sum_{j+k \leq n-1} \bm{\Gamma}_{(-1 + C_{(j+1)} \epsilon)}^{(j+1)} \bm{\Gamma}^{(k)}_{(-\frac{1}{2} + \delta)}
	\\
	&=
	\bm{\Gamma}^{(n)}_{(-1 + C_{(n)}\epsilon)}
	\end{split}
\end{equation*}
where for the final line it is important to note that, for all $j + k$ such that $j+k \leq n-1$, we have
\begin{equation*}
C_{(n)} \gg C_{(j)} + C_{(k)}
\end{equation*}

Now, we can compute
\begin{equation*}
 \begin{split}
  \mathscr{Z}^{N+1} \slashed{\D} \phi
  &= \mathscr{Z} \left( \mathscr{Z}^N \slashed{\phi} \right) \\
  \\
  &= \mathscr{Z} \bigg( \slashed{\D} \mathscr{Z}^N \phi
  + \sum_{j+k \leq N-1} \bm{\Gamma}_{(C_{(j+1)} \epsilon)}^{(j+1)} \overline{\slashed{\D}} \mathscr{Z}^k \phi
  + \sum_{j+k \leq N-1} \bm{\Gamma}_{(C_{(j)} \epsilon)}^{(j)} \slashed{\D} \mathscr{Z}^k \phi
  \\
  &\phantom{= \mathscr{Z} \bigg(}
  + \sum_{j+k \leq N-1} \bm{\Gamma}_{(-1 + C_{(j+1)} \epsilon)}^{(j+1)} \mathscr{Z}^k \phi 
  \bigg)
  \\
  \\
  &= \slashed{\D} \mathscr{Z}^{N+1} \phi
  + \bm{\Gamma}^{(0)}_{(C_{(0)}\epsilon)} \slashed{\D} \mathscr{Z}^N \phi
  + \bm{\Gamma}^{(1)}_{(-1+C_{(1)}\epsilon)} \mathscr{Z}^N \phi
  + \sum_{j+k \leq N-1} \bm{\Gamma}^{(j+2)}_{(C_{(j+2)}\epsilon)} \overline{\slashed{\D}} \mathscr{Z}^k \phi
  \\
  &\phantom{=}
  + \sum_{j+k \leq N-1} \bm{\Gamma}^{(j+1)}_{(C_{(j+1)}\epsilon)} \overline{\slashed{\D}} \mathscr{Z}^{k+1} \phi
  + \sum_{j+k \leq N-1} \bm{\Gamma}^{(j+1)}_{(C_{(j+1)}\epsilon)} \slashed{\D}\mathscr{Z}^k \phi \\
  &\phantom{=}
  + \sum_{j+k \leq N-1} \bm{\Gamma}^{(j)}_{(C_{(j)}\epsilon)} \left( \slashed{\D} \mathscr{Z}^{k+1} \phi + \bm{\Gamma}^{(0)}_{(C_{(0)}\epsilon)}\slashed{\D} \mathscr{Z}^k \phi + \bm{\Gamma}^{(1)}_{(C_{(1)}\epsilon)} \mathscr{Z}^k \phi \right) \\
  &\phantom{=}
  + \sum_{j+k \leq N-1} \left( \bm{\Gamma}_{(-1 + C_{(j+2)} \epsilon)}^{(j+2)} \mathscr{Z}^k \phi + \bm{\Gamma}_{(-1 + C_{(j+1)} \epsilon)}^{(j+1)} \mathscr{Z}^{k+1} \phi \right) \\
  \\
  &= \slashed{\D} \mathscr{Z}^{N+1} \phi + \sum_{j+k\leq N} \bm{\Gamma}^{(j)}_{(C_{(j)}\epsilon)} \slashed{\D} \mathscr{Z}^k \phi
  + \sum_{j+k \leq N} \bm{\Gamma}^{(j+1)}_{(-1 + C_{(j+1)}\epsilon)} \mathscr{Z}^k \phi
 \end{split}
\end{equation*}
where in the last line we have used $\bm{\Gamma}^{(j)}_{(C_{(j)}\epsilon)} \bm{\Gamma^{(0)}}_{(C_{(0)}\epsilon)} = \bm{\Gamma}^{(j+1)}_{(C_{(j)}\epsilon)}$, and other similar equations, which all follow from the fact that $C_{(j)}$ is sufficiently large compared with $C_{(j-1)}$.

Hence, the claim holds also for all $n \leq N+1$, and so, by induction, it holds for all $n$.

The proof for the case of good derivatives follows from an almost identical calculation.
\end{proof}

\begin{proposition}
	\label{proposition commute D Yn}
	Let $\phi$ be an $S_{\tau,r}$-tangent tensor field, and suppose that $C_{(n)} \gg C_{(n-1)}$. Then we have, schematically,
	\begin{equation}
	[\mathscr{Y}^n, \slashed{\D}] \phi
	=
	\sum_{j+k \leq n-1} \left(1 + \bm{\Gamma}_{(C_{(j+1)} \epsilon)}^{(j+1)} \right) \overline{\slashed{\D}} \mathscr{Y}^k \phi 
	+ \sum_{j+k \leq n-1} \bm{\Gamma}_{(C_{(j)} \epsilon)}^{(j)} \slashed{\D} \mathscr{Y}^k \phi 
	+ \sum_{j+k \leq n-1} \bm{\Gamma}_{(-1 + C_{(j+1)} \epsilon)}^{(j+1)} \mathscr{Y}^k \phi 
	\end{equation}
	On the other hand, if we commute with a ``good'' derivative then we have
	\begin{equation*}
	\begin{split}
	[\mathscr{Y}^n , \overline{\slashed{\D}}] \phi
	&=
	\sum_{j+k \leq n-1} \bm{\Gamma}^{(j)}_{(-1+C_{(j)}\epsilon)} \slashed{\D} \mathscr{Y}^k \phi 
	+ \sum_{j+k \leq n-1} \left(1 + \bm{\Gamma}^{(j+1)}_{(-\delta)}\right) \overline{\slashed{\D}} \mathscr{Y}^k \phi 
	+ \sum_{j+k \leq n-1} \bm{\Gamma}_{(-1 + C_{(j+1)} \epsilon)}^{(j+1)} \mathscr{Y}^k \phi 
	\\
	&=
	\sum_{j+k \leq n-1} \bm{\Gamma}_{(-\delta)}^{(j)} \slashed{\D} \mathscr{Y}^k \phi 
	+ \sum_{j+k \leq n-1} \bm{\Gamma}_{(-1 + C_{(j+1)} \epsilon)}^{(j+1)} \mathscr{Y}^k \phi 
	\end{split}
	\end{equation*}

\end{proposition}
\begin{proof}
	Using proposition \ref{proposition commuting rL with first order operators} we have, schematically,
	\begin{equation*}
	[r\slashed{\D}_L , \slashed{\D}] \phi
	=
	\slashed{\D}_L \phi
	+ \bm{\Gamma}^{(1)}_{(C_{(1)}\epsilon)} \overline{\slashed{\D}} \phi
	+ \bm{\Gamma}^{(0)}_{(0)} \slashed{\D} \phi
	+ \begin{pmatrix} r\Omega_{L\Lbar} \\ r\slashed{\Omega}_L \end{pmatrix} \phi
	\end{equation*}
	and so, making use of the proposition above, we have
	\begin{equation*}
	[\mathscr{Y} , \slashed{\D}] \phi
	=
	\slashed{\D}_L \phi
	+ \bm{\Gamma}^{(1)}_{(C_{(1)}\epsilon)} \overline{\slashed{\D}} \phi
	+ \bm{\Gamma}^{(C_{(0)}\epsilon)}_{(0)} \slashed{\D} \phi
	+ \begin{pmatrix} r\Omega_{L\Lbar} \\ r\slashed{\Omega}_L \end{pmatrix} \phi
	\end{equation*}
	Relative to the proposition above, the only difference is the appearence of the term $\slashed{\D}_L \phi$ (with the ``large'' coefficient)	and so this proposition follows straightforwardly. The same is true when we commute with good derivatives.

\end{proof}

\begin{proposition}
\label{proposition commute Zn L}
Let $\phi$ be any $S_{\tau,r}$-tangent tensor field. Then we have
\begin{equation*}
  \mathscr{Z}^{n} \slashed{\D}_L \phi =
 \slashed{\D}_L \mathscr{Z}^{n} \phi
 + \bm{\Gamma}^{(0)}_{(-1)} \mathscr{Z}^{n} \phi
 + \sum_{j+k \leq n-1} \bm{\Gamma}^{(j+1)}_{(-1+C_{(j+1)}\epsilon)} \overline{\slashed{\D}} \mathscr{Z}^{k} \phi
 + \sum_{j+k \leq n-1} \bm{\Gamma}^{(j+1)}_{(-1-\delta)} \mathscr{Z}^k \phi
\end{equation*}
\end{proposition}

\begin{proof}
Using propositions \ref{proposition commuting DT with first order operators} and \ref{proposition commuting rnabla with first order operators} we find that, schematically, for any $S_{\tau,r}$-tangent tensor field $\phi$,
\begin{equation*}
 \begin{split}
  [\slashed{\D}_L \, , \slashed{\D}_T] \phi 
  &= 
  \omega \slashed{\D}_T \phi
  + \bm{\Gamma}^{(1)}_{(-1+C_{(1)}\epsilon)} \overline{\slashed{\D}}\phi
  + \Omega_{L\Lbar} \cdot \phi 
  \\
  &=
  \bm{\Gamma}^{(0)}_{(-1)} \cdot \slashed{\D}_T \phi
  + \bm{\Gamma}^{(1)}_{(-1+C_{(1)}\epsilon)} \overline{\slashed{\D}}\phi
  + \bm{\Gamma}^{(1)}_{(-2+C_{(1)}\epsilon)} \cdot \phi
  \\
  \\
  [\slashed{\D}_L \, , r\slashed{\nabla}] \phi 
  &=
  \chi_{(\text{small})} \cdot r\slashed{\nabla} \phi 
  + r\Omega_{L\slashed{\alpha}} \cdot \phi
  \\
  &=
  \bm{\Gamma}^{(0)}_{(-1-\delta)} \cdot (r\slashed{\nabla}\phi)
  + \bm{\Gamma}^{(1)}_{(-1-\delta)} \cdot \phi
 \end{split}
\end{equation*}
where we have also made use of the structure of the structure of the tensor $\Omega$ given in proposition \ref{proposition expression for Omega}. So, putting these together, we can write
\begin{equation*}
  [\slashed{\D}_L \, , \mathscr{Z}] \phi 
  =
  \bm{\Gamma}^{(0)}_{(-1)} \mathscr{Z}\phi
  + \bm{\Gamma}^{(1)}_{(-1+C_{(1)}\epsilon)} \overline{\slashed{\D}}\phi + \bm{\Gamma}_{(-1 - \delta)}^{(1)}\phi
\end{equation*}

Now, suppose that, for all $n \leq N$ we have
\begin{equation*}
 \mathscr{Z}^n \slashed{\D}_L \phi =
 \slashed{\D}_L \mathscr{Z}^n \phi
 + \bm{\Gamma}^{(0)}_{(-1)} \mathscr{Z}^n \phi
 + \sum_{j+k \leq n-1} \bm{\Gamma}^{(j+1)}_{(-1+C_{(j+1)}\epsilon)} \overline{\slashed{\D}} \mathscr{Z}^{k} \phi
 + \sum_{j+k \leq n-1} \bm{\Gamma}^{(j+1)}_{(-1-\delta)} \mathscr{Z}^k \phi
\end{equation*}
Clearly this holds for $N=1$, by the previous calculation.

Now, we can compute
\begin{equation*}
 \begin{split}
  \mathscr{Z}^{N+1} \slashed{\D}_L &\phi
  =
  \mathscr{Z} \slashed{\D}_L \mathscr{Z}^N \phi
 + \bm{\Gamma}^{(1)}_{(-1+ C_{(1)}\epsilon)} \mathscr{Z}^N \phi
 + \bm{\Gamma}^{(0)}_{(-1)} \mathscr{Z}^{N+1} \phi
 + \sum_{j+k \leq N-1} \bm{\Gamma}^{(j+2)}_{(-1+C_{(j+2)}\epsilon)} \overline{\slashed{\D}} \mathscr{Z}^{k} \phi \\
 &\phantom{=}
 + \sum_{j+k \leq N-1} \bm{\Gamma}^{(j+1)}_{(-1+C_{(j+1)}\epsilon)} \mathscr{Z}\overline{\slashed{\D}} \mathscr{Z}^{k} \phi
 + \sum_{j+k \leq N} \bm{\Gamma}^{(j+1)}_{(-1-\delta)} \mathscr{Z}^k \phi
 \\ \\
 &=
 \slashed{\D}_L \mathscr{Z}^{N+1} \phi
 + \bm{\Gamma}^{(0)}_{(-1)} \mathscr{Z}^{N+1} \phi
 + \bm{\Gamma}^{(1)}_{(-1+C_{(1)}\epsilon)}\overline{\slashed{\D}} \mathscr{Z}^N \phi
 + \bm{\Gamma}^{(1)}_{(-1-\delta)} \mathscr{Z}^N \phi \\
 &\phantom{=}
 + \sum_{j+k \leq N} \bm{\Gamma}^{(j+1)}_{(-1+C_{(j+1)}\epsilon)} \overline{\slashed{\D}} \mathscr{Z}^{k} \phi 
 + \sum_{j+k \leq N} \bm{\Gamma}^{(j+1)}_{(-1-\delta)} \mathscr{Z}^k \phi \\
 &\phantom{}
 + \sum_{j+k \leq N-1} \bm{\Gamma}^{(j+1)}_{(-1+C_{(j+1)}\epsilon)} \left( \overline{\slashed{\D}} \mathscr{Z}^{k+1} \phi
   + \sum_{l+m \leq k-1} \bm{\Gamma}^{(l)}_{(-\delta)}\slashed{\D}\mathscr{Z}^m \phi
   + \sum_{l+m \leq k-1} \bm{\Gamma}^{(l+1)}_{(-1+C_{(l+1)}\epsilon)}\mathscr{Z}^m \phi \right)
 \end{split}
\end{equation*}
Now, we can write
\begin{equation*}
 \overline{\slashed{\D}} \mathscr{Z}^m \phi = \slashed{\D}_L \mathscr{Z}^m \phi + r^{-1}\mathscr{Z}^{m+1}\phi
\end{equation*}
and so, in the end, we find
\begin{equation*}
 \mathscr{Z}^{N+1} \slashed{\D}_L \phi =
 \slashed{\D}_L \mathscr{Z}^{N+1} \phi
 + \bm{\Gamma}^{(0)}_{(-1)} \mathscr{Z}^{N+1} \phi
 + \sum_{j+k \leq N} \bm{\Gamma}^{(j)}_{(-1+C_{(j)}\epsilon)} \overline{\slashed{\D}} \mathscr{Z}^{k} \phi
 + \sum_{j+k \leq N} \bm{\Gamma}^{(j+1)}_{(-1-\delta)} \mathscr{Z}^k \phi
\end{equation*}
so the proposition holds also for $n = N+1$, and so it holds for all $n$.
\end{proof}

\begin{proposition}
	\label{proposition commute Yn L}
	Let $\phi$ be any $S_{\tau,r}$-tangent tensor field. Then we have
	\begin{equation*}
	\begin{split}
	\mathscr{Y}^{n} \slashed{\D}_L \phi 
	&=
	\sum_{j \leq n} \slashed{\D}_L \mathscr{Y}^j \phi
	+ \sum_{j \leq n} \bm{\Gamma}^{(0)}_{(-1)} \mathscr{Y}^{j} \phi
	+ \sum_{j+k \leq n-1} \bm{\Gamma}^{(j+1)}_{(-1+C_{(j+1)}\epsilon)} \overline{\slashed{\D}} \mathscr{Y}^{k} \phi
	+ \sum_{j+k \leq n-1} \bm{\Gamma}^{(j+1)}_{(-1-\delta)} \mathscr{Y}^k \phi
	\end{split}	
	\end{equation*}
\end{proposition}

\begin{proof}
	Using proposition \ref{proposition commuting rL with first order operators} we find that
	\begin{equation*}
	r\slashed{\D}_L \slashed{\D}_L \phi
	=
	\slashed{\D}_L (r\slashed{\D}_L \phi)
	- \slashed{\D}_L \phi
	\end{equation*}
	Combining this with the calculations in the previous proposition, we find that
	\begin{equation*}
	[\slashed{\D}_L \, , \mathscr{Y}] \phi
	=
	\slashed{\D}_L \phi + \bm{\Gamma}^{(0)}_{(-1)} \mathscr{Y}\phi
	+ \bm{\Gamma}^{(1)}_{(-1+C_{(1)}\epsilon)} \overline{\slashed{\D}}\phi
	+ \bm{\Gamma}^{(1)}_{(-1-\delta)} \phi
	\end{equation*}
	
	Now, we suppose that, for all $n \leq N$, schematically we have
	\begin{equation*}
	\mathscr{Y}^{n} \slashed{\D}_L \phi 
	=
	\sum_{j \leq n} \slashed{\D}_L \mathscr{Y}^j \phi
	+ \bm{\Gamma}^{(0)}_{(-1)} \mathscr{Y}^{n} \phi
	+ \sum_{j+k \leq n-1} \bm{\Gamma}^{(j+1)}_{(-1+C_{(j+1)}\epsilon)} \overline{\slashed{\D}} \mathscr{Y}^{k} \phi
	+ \sum_{j+k \leq n-1} \bm{\Gamma}^{(j+1)}_{(-1-\delta)} \mathscr{Y}^k \phi
	\end{equation*}
	Applying one more operator $\mathscr{Y}$ we have
	\begin{equation*}
	\begin{split}
	\mathscr{Y}^{N+1} \slashed{\D}_L \phi 
	&=
	\sum_{j \leq N} \mathscr{Y}\slashed{\D}_L \mathscr{Y}^j \phi
	+ \bm{\Gamma}^{(1)}_{(-1 + C_{(1)}\epsilon)} \mathscr{Y}^{N} \phi
	+ \bm{\Gamma}^{(0)}_{(-1)} \mathscr{Y}^{N+1} \phi
	+ \sum_{j+k \leq N-1} \bm{\Gamma}^{(j+2)}_{(-1+C_{(j+2)}\epsilon)} \overline{\slashed{\D}} \mathscr{Y}^{k} \phi
	\\
	&\phantom{=}
	+ \sum_{j+k \leq N-1} \bm{\Gamma}^{(j+1)}_{(-1+C_{(j+1)}\epsilon)} \mathscr{Y} \overline{\slashed{\D}} \mathscr{Y}^{k} \phi
	+ \sum_{j+k \leq N-1} \bm{\Gamma}^{(j+2)}_{(-1-\delta)} \mathscr{Y}^k \phi
	+ \sum_{j+k \leq n-1} \bm{\Gamma}^{(j+1)}_{(-1-\delta)} \mathscr{Y}^{k+1} \phi
	\\
	&=
	\sum_{j \leq N} \left( \slashed{\D}_L \mathscr{Y}^{j+1} \phi
		+ \bm{\Gamma}^{(0)}_{(-1)} \mathscr{Y}^{j+1} \phi
		+ \bm{\Gamma}^{(1)}_{(-1 + C_{(1)})} \overline{\slashed{\D}} \mathscr{Y}^{j} \phi
		+ \bm{\Gamma}^{(1)}_{(-1 - \delta)} \mathscr{Y}^{j} \phi
		\right)
	+ \bm{\Gamma}^{(1)}_{(-1 + C_{(1)}\epsilon)} \mathscr{Y}^n \phi
	\\
	&\phantom{=}
	+ \bm{\Gamma}^{(1)}_{(-1)} \mathscr{Y}^{N+1} \phi
	+ \sum_{j+k \leq N-1} \bm{\Gamma}^{(j+1)}_{(-1+C_{(j+1)}\epsilon)} \overline{\slashed{\D}} \mathscr{Y}^k \phi
	+ \sum_{j+k \leq N-1} \bm{\Gamma}^{(j)}_{(-1+C_{(j)}\epsilon)} \overline{\slashed{\D}} \mathscr{Y}^{k+1} \phi
	\\
	&\phantom{=}
	+ \sum_{j+k \leq N-1} \bm{\Gamma}^{(j+2)}_{(-1+C_{(j)}\epsilon)} \overline{\slashed{\D}} \mathscr{Y}^{k+1} \phi
	\end{split}
	\end{equation*}
	where we have made use of the previous propositions to commute $\mathscr{Y}$ and $\slashed{\D}$, and to commute $\mathscr{Y}$ and $\overline{\slashed{\D}}$. Collecting terms shows that, schematically
	\begin{equation*}
	\mathscr{Y}^{N+1} \slashed{\D}_L \phi 
	=
	\sum_{j \leq N+1} \slashed{\D}_L \mathscr{Y}^j \phi
	+ \bm{\Gamma}^{(0)}_{(-1)} \mathscr{Y}^{N+1} \phi
	+ \sum_{j+k \leq N} \bm{\Gamma}^{(j+1)}_{(-1+C_{(j+1)}\epsilon)} \overline{\slashed{\D}} \mathscr{Y}^{k} \phi
	+ \sum_{j+k \leq N} \bm{\Gamma}^{(j+1)}_{(-1-\delta)} \mathscr{Y}^k \phi
	\end{equation*}
	
\end{proof}

\begin{proposition}[An alternative expression for $\lbrack\mathscr{Y}^n , \slashed{\D}_L \rbrack \phi$]
	\label{proposition alternative expression for commutator L Yn}

	Let $\phi$ be any $S_{\tau,r}$-tangent tensor field. Suppose that the operator $\slashed{\D}_T$ appears $k$ times in the expansion of $\mathscr{Y}^n$. Then we have
	\begin{equation*}
	[\slashed{\D}_L, \mathscr{Y}^n] \phi
	- k\omega \mathscr{Y}^n \phi
	=
	\sum_{j + k \leq n-1} \bm{\Gamma}^{(j+1)}_{(-1-\delta)} \mathscr{Y}^{k+1} \phi
	+ \sum_{\substack{j+k \leq n-1 \\ k \leq n-2}} \bm{\Gamma}^{(j)}_{(-1+C_{(j)}\epsilon)} \mathscr{Y}^{k+1}\phi 
	\end{equation*}
	
\end{proposition}

\begin{proof}
	This follows from the previous proposition, together with the fact that $\overline{\slashed{\D}} \sim r^{-1} \mathscr{Y}$.
\end{proof}

We will also occasionally need to commute the operators $\mathscr{Y}$ and $\mathscr{Z}$. Clearly, the commutator $[\mathscr{Y}, \mathscr{Z}]$ can be expressed in terms of the operators $\mathscr{Z}$ \emph{unless} $\mathscr{Y} = r\slashed{\D}_L$. Thus, we only need to consider $[r\slashed{\D}_L , \mathscr{Z}]$.

\begin{proposition}[Commuting $\mathscr{Z}$ with $(r\slashed{\D}_L)$].
	\label{proposition commuting Z with rL}
	For any $S_{\tau,r}$-tangent tensor field $\phi$, we have, schematically
	\begin{equation}
	[\mathscr{Z} , r\slashed{\D}_L]\phi
	=
	\bm{\Gamma}^{(0)}_{(0)} \slashed{\D}_T \phi
	+ \left( \bm{\Gamma}^{(1)}_{(-1 + C_{(1)}\epsilon)} + r(\chi_{(\text{small})}) \right) \mathscr{Z} \phi
	+ \left(r\Omega_{L\Lbar} + r^2\slashed{\Omega}_L \right) \phi
	\end{equation}
	
\end{proposition}

\begin{proof}
	We consider the two cases $\mathscr{Z} = \slashed{\D}_T$ and $\mathscr{Z} = r\slashed{\nabla}$. First, we have
	\begin{equation*}
	\begin{split}
	[\slashed{\D}_T, r\slashed{\D}_L] \phi
	&=
	\frac{1}{2}r\omega \slashed{\D}_L \phi
	+ \frac{1}{2}r\omega \slashed{\D}_{\Lbar} \phi
	+ r(\zeta^\alpha + \slashed{\nabla}^\alpha \log \mu) \slashed{\nabla}_\alpha \phi
	- \frac{1}{2} r \Omega_{L\Lbar} \cdot \phi
	\\
	&=
	\bm{\Gamma}^{(0)}_{(0)} \slashed{\D}_T \phi
	+ \bm{\Gamma}^{(1)}_{(-1 + C_{(1)}\epsilon)} \mathscr{Z} \phi
	+ r \Omega_{L\Lbar} \cdot \phi
	\end{split}
	\end{equation*}
	and
	\begin{equation*}
	\begin{split}
	[r\slashed{\nabla}, r\slashed{\D}_L] \phi
	&=
	r^2(\chi_{(\text{small})}) \cdot \slashed{\nabla} \phi
	+ r^2 \slashed{\Omega}_L \cdot \phi
	\\
	&=
	r(\chi_{(\text{small})}) \cdot \mathscr{Z} \phi
	+ r^2 \slashed{\Omega}_L \cdot \phi
	\end{split}
	\end{equation*}
	so, putting these together we have
	\begin{equation*}
	[\mathscr{Z}, r\slashed{\D}_L] \phi
	=
	\bm{\Gamma}^{(1)}_{(0)}\mathscr{Z} \phi
	+ \bm{\Gamma}^{(1)}_{(-\delta)} \phi	
	\end{equation*}
	Now, suppose that for all $n \leq N$ we have
	\begin{equation*}
	[\mathscr{Z}, (r\slashed{\D}_L)^n] \phi
	=
	\sum_{j+k \leq n-1} \bm{\Gamma}^{(j+1)}_{(C_{(j+1)}\epsilon)} \mathscr{Y}^{k+1} \phi
	\end{equation*}
	This clearly holds in the case $N = 1$, by the calculation above. Then we have
	\begin{equation*}
	\begin{split}
	[\mathscr{Z}, (r\slashed{\D}_L)^{N+1}]\phi
	&=
	[\mathscr{Z}, (r\slashed{\D}_L)^N] (r\slashed{\D}_L \phi)
	+ (r\slashed{\D}_L)^N [\mathscr{Z}, r\slashed{\D}_L] \phi
	\\
	&=
	\sum_{j+k \leq N-1} \bm{\Gamma}^{(j+1)}_{(C_{(j+1)}\epsilon)} \mathscr{Y}^{k+1} (r\slashed{\D}_L \phi)
	+ (r\slashed{\D}_L)^N \left( 
		\bm{\Gamma}^{(1)}_{(0)} \mathscr{Z} \phi
		+ \bm{\Gamma}^{(1)}_{(-\delta)} \phi \right)
	\\
	&=
	\sum_{j+k \leq N-1} \bm{\Gamma}^{(j+1)}_{(C_{(j+1)}\epsilon)} \mathscr{Y}^{k+2} \phi
	+ \sum_{j+k \leq N} \bm{\Gamma}^{(j+1)}_{(C_{(j+1)}\epsilon)} \mathscr{Y}^{k+1} \phi
	+ \sum_{j+k \leq N} \bm{\Gamma}^{(j+1)}_{(-\delta)} \mathscr{Y}^{k} \phi
	\\
	&=
	\sum_{j+k \leq N} \bm{\Gamma}^{(j+1)}_{(C_{(j+1)}\epsilon)} \mathscr{Y}^{k+1} \phi
	\end{split}
	\end{equation*}
	proving the inductive step.
\end{proof}

\begin{proposition}[Transport equations for commuted rectangular components of the frame fields]
	\label{proposition transport Yn rectangular}
	After applying the operators $\mathscr{Y}$ $n$ times, the rectangular components of the frame fields satisfy transport equations of the form
	\begin{equation*}
	\begin{split}
	\slashed{\D}_L \mathscr{Y}^n X_{(\text{frame})}
	&=
	\bm{\Gamma}^{(0)}_{(-1)} \mathscr{Y}^n X_{(\text{frame})}
	+ \bm{\Gamma}^{(0)}_{(0, \text{large})} (\overline{\slashed{\D}} \mathscr{Y}^n h)_{(\text{frame})}
	+ \bm{\Gamma}^{(0)}_{(-2, \text{large})} (r\mathscr{Y}^n \bar{X}_{(\text{frame})})
	\\
	&\phantom{=}
	+ \bm{\Gamma}^{(n-1)}_{(-1 + C_{(n-1)}\epsilon)}
	\\
	\\
	\slashed{\D}_L \left( r\mathscr{Y}^n \bar{X}_{(\text{frame})} \right)
	&=
	\bm{\Gamma}^{(0)}_{(-1)} \mathscr{Y}^n (r\bar{X}_{(\text{frame})})
	+ \bm{\Gamma}^{(n)}_{(-1+C_{(n)}\epsilon)} \bm{\Gamma}^{(0)}_{(-\frac{3}{2}\delta)}
	+ r \bm{\Gamma}^{(0)}_{(0, \text{large})} (\overline{\slashed{\D}} \mathscr{Y}^n h)_{(\text{frame})}
	+ r\bm{\Gamma}^{(n-1)}_{(-1-\delta)}
	\end{split}
	\end{equation*}
	
\end{proposition}

\begin{proof}
	From proposition \ref{proposition transport rectangular} we have
	\begin{equation*}
	\begin{split}
	\slashed{\D}_L X_{(\text{frame})} &= (\bar{\partial}h_{(\text{rect})})X_{(\text{frame})}^3 + r^{-2} X_{(\text{frame})}^2 (r\bar{X}_{(\text{frame})}) + r^{-3} X_{(\text{frame})} (r\bar{X}_{(\text{frame})})^2 \\
	L (r\bar{X}_{(\text{frame})}) &= (\bar{\partial}h_{(\text{rect})}) X_{(\text{frame})}^2 (r\bar{X}_{(\text{frame})}) + r(\bar{\partial} h_{(\text{rect})}) X_{(\text{frame})}^3 + r^{-2}X_{(\text{frame})} (r\bar{X}_{(\text{frame})})^2
	\end{split}
	\end{equation*}
	which we can write schematically as
	\begin{equation*}
	\begin{split}
	\slashed{\D}_L X_{(\text{frame})} &=
	\bm{\Gamma}^{(0)}_{(-1-\delta)} X_{(\text{frame})} \\
	L (r\bar{X}_{(\text{frame})}) &=
	\bm{\Gamma}^{(0)}_{(-1-\delta)} (r\bar{X}_{(\text{frame})})
	+ r(\bar{\partial} h_{(\text{rect})}) X_{(\text{frame})}^3
	\end{split}
	\end{equation*}
	Applying the operator $\mathscr{Y}^n$ to the first of these equations and using proposition \ref{proposition commute Yn L}, we find
	\begin{equation*}
	\begin{split}
	\slashed{\D}_L \mathscr{Y}^n X_{(\text{frame})}
	&=
	[\slashed{\D}_L , \mathscr{Y}^n] X_{(\text{frame})}
	+ \sum_{j+k \leq n} \bm{\Gamma}^{(j)}_{(-1-\delta)} \mathscr{Y}^k X_{(\text{frame})} \\ \\
	&=
	\sum_{j \leq n} \bm{\Gamma}^{(0)}_{(-1)}\mathscr{Y}^j X_{(\text{frame})}
	+ \sum_{j+k \leq n-1} \bm{\Gamma}^{(j+1)}_{(-1 + C_{(j+1)}\epsilon)} \overline{\slashed{\D}} \mathscr{Y}^k X_{(\text{frame})} \\
	&\phantom{=}
	+ \sum_{j+k \leq n} \bm{\Gamma}^{(j)}_{(-1-\delta)} \mathscr{Y}^k X_{(\text{frame})} \\
	\end{split}
	\end{equation*}
	We can further decompose
	\begin{equation*}
	\overline{\slashed{\D}} \mathscr{Y}^k X_{(\text{frame})}
	=
	r^{-1} \mathscr{Y}^{k+1} X_{(\text{frame})}
	\end{equation*}
	where we have used the fact that $\overline{\slashed{\D}} \phi = r^{-1} \mathscr{Y}\phi$. Hence we have
	\begin{equation*}
	\begin{split}
	\slashed{\D}_L \mathscr{Y}^n X_{(\text{frame})}
	&=
	\sum_{j \leq n} \bm{\Gamma}^{(0)}_{(-1)}\mathscr{Y}^n X_{(\text{frame})}
	+\sum_{j+k \leq n-1} \bm{\Gamma}^{(j+1)}_{(-2+C_{(j+1)}\epsilon)} \mathscr{Y}^{k+1} X_{(\text{frame})}
	\\
	&\phantom{=}
	+ \sum_{j+k \leq n} \bm{\Gamma}^{(j)}_{(-1-\delta)} \mathscr{Y}^k X_{(\text{frame})}
	\end{split}
	\end{equation*}
	so we conclude that
	\begin{equation*}
	\slashed{\D}_L \mathscr{Y}^n X_{(\text{frame})}
	=
	\bm{\Gamma}^{(0)}_{(-1)} \mathscr{Y}^n X_{(\text{frame})}
	+ \bm{\Gamma}^{(n)}_{(-1-\delta)}
	+ \bm{\Gamma}^{(n-1)}_{(-1 + C_{(n-1)}\epsilon)}
	\end{equation*}
	If we pay attention to the highest order terms, then we actually find that
	\begin{equation*}
	\slashed{\D}_L \mathscr{Y}^n X_{(\text{frame})}
	=
	\bm{\Gamma}^{(0)}_{(-1)} \mathscr{Y}^n X_{(\text{frame})}
	+ \bm{\Gamma}^{(0)}_{(0, \text{large})} (\overline{\slashed{\D}} \mathscr{Y}^n h)_{(\text{frame})}
	+ \bm{\Gamma}^{(0)}_{(-2, \text{large})} (r\mathscr{Y}^n \bar{X}_{(\text{frame})})
	+ \bm{\Gamma}^{(n-1)}_{(-1 + C_{(n-1)}\epsilon)}
	\end{equation*}

	Similarly, applying the operator $\mathscr{Y}^n$ to the equation for $L(r\bar{X}_{(\text{frame})})$, we have
	\begin{equation*}
	\begin{split}
	\slashed{\D}_L \left( r\mathscr{Y}^n \bar{X}_{(\text{frame})} \right)
	&=
	[\slashed{\D}_L , \mathscr{Y}^n] \left( r \bar{X}_{(\text{frame})} \right)
	+ \sum_{j+k \leq n} \bm{\Gamma}^{(j)}_{(-1-\delta)} \left( r\mathscr{Y}^k \bar{X}_{(\text{frame})} \right)
	+ r(\overline{\slashed{\D}} \mathscr{Y}^n h)_{(\text{frame})} X_{(\text{frame})} \\
	&\phantom{=}
	+ r \bm{\Gamma}^{(0)}_{(-1-\delta)} (\mathscr{Y}^n X_{(\text{frame})})
	+ r \bm{\Gamma}^{(n-1)}_{(-1-\delta)}
	\\ \\
	&= \sum_{j \leq n} \bm{\Gamma}^{(0)}_{(-1)} \mathscr{Y}^j (r\bar{X}_{(\text{frame})})
	+ \sum_{j+k \leq n-1} \bm{\Gamma}^{(j+1)}_{(-1+C_{(j+1)}\epsilon)}\left( \overline{\slashed{\D}} \left( r\mathscr{Y}^k \bar{X}_{(\text{frame})} \right) \right) \\
	&\phantom{=}
	+ \sum_{j+k \leq n-1} \bm{\Gamma}^{(j+1)}_{(-1-\delta)}(r\mathscr{Y}^k\bar{X}_{(\text{frame})})
	+ \sum_{j+k \leq n} \bm{\Gamma}^{(j)}_{(-1-\delta)} \left( r\mathscr{Y}^k \bar{X}_{(\text{frame})} \right) \\
	&\phantom{=}
	+ r(\overline{\slashed{\D}} \mathscr{Y}^n h)_{(\text{frame})} X_{(\text{frame})}
	+ r \bm{\Gamma}^{(n-1)}_{(-1-\delta)}
\end{split}
\end{equation*}
Again, we can substitute the good derivatives for $r^{-1}\mathscr{Y}$ to find
\begin{equation*}
	\begin{split}
	\slashed{\D}_L \left( r\mathscr{Y}^n \bar{X}_{(\text{frame})} \right)
	&=
	\bm{\Gamma}^{(0)}_{(-1)} \mathscr{Y}^n (r\bar{X}_{(\text{frame})})
	+ \bm{\Gamma}^{(n)}_{(-1+C_{(n)}\epsilon)} \bm{\Gamma}^{(0)}_{(-\frac{3}{2}\delta)}
	+ r \bm{\Gamma}^{(0)}_{(0, \text{large})} (\overline{\slashed{\D}} \mathscr{Y}^n h)_{(\text{frame})}
	+ r\bm{\Gamma}^{(n-1)}_{(-1-\delta)}	
	\end{split}
\end{equation*}

\end{proof}

\begin{proposition}[Transport equations for $\mathscr{Y}^n X_{(\text{frame, small})}$]
	\label{proposition transport Yn rectangular small}
	After commuting with the operators $\mathscr{Z}$ $n$ times, the ``small'' null frame components $X_{(\text{frame, small})}$ satisfy equations of the form
	\begin{equation}
	\slashed{\D}_L \mathscr{Z}^n X_{(\text{frame, small})}
	=
	\bm{\Gamma}^{0}_{(-1)} \mathscr{Z}^n X_{(\text{frame})} + \bm{\Gamma}^{(n)}_{(-1-\frac{1}{2}\delta)} + \bm{\Gamma}^{(n-1)}_{(-1+C_{(n-1)}\epsilon)}
	\end{equation}
\end{proposition}

\begin{proof}
	We have
	\begin{equation*}
	\begin{split}
	L_{(\text{small})}^0 = L^0 - 1 \\
	L_{(\text{small})}^i = L^i - \frac{x^i}{r} \\
	\Lbar_{(\text{small})}^0 = \Lbar^0 - 1 \\
	\Lbar_{(\text{small})}^i = \Lbar^i + \frac{x^i}{r} \\	
	\end{split}
	\end{equation*}
	So the evolution equations for $L^0_{(\text{small})}$ and $\Lbar^0_{(\text{small})}$ obey identical transport equations to those for $L^0$ and $\Lbar^0$. On the other hand, we have
	\begin{equation*}
	\mathscr{Y} \left(\frac{x^i}{r}\right) = \begin{pmatrix} r^{-1} \\ 1 \end{pmatrix} X^i_{(\text{frame})}
	\end{equation*}
	so that, schematically,
	\begin{equation*}
	\slashed{\D}_L (\mathscr{Z}^n X_{(\text{frame, small})}) = \slashed{\D}_L \left( \mathscr{Z}^n X_{(\text{frame})} + \mathscr{Z}^{n-1} X_{(\text{frame})}\right) 
	\end{equation*}
	Hence these quantities, too, obey transport equations which are schematically of the same form as those for $\mathscr{Z}^n X_{(\text{frame})}$.
	
\end{proof}

\begin{proposition}[Transport equation for $\mathscr{Y}^n \mu$]
 \label{proposition transport Yn mu}
 The commuted versions of the foliation density $\mu$ satisfy transport equations of the form
\begin{equation}
 \begin{split}
  \slashed{\D}_L \left( \mathscr{Y}^n \log \mu \right)
  &=\bm{\Gamma}^{(0)}_{(-1)} \mathscr{Y}^n \log \mu
  + \sum_{j+k\leq n-1} \bm{\Gamma}^{(j+1)}_{(-1-\delta)} (\mathscr{Y}^k \log \mu)
  + \bm{\Gamma}^{(0)}_{(-1, \text{large})} \left( \mathscr{Y}^n \bar{X}_{(\text{frame})} \right) \\
  &\phantom{=}
  + (\slashed{\D} \mathscr{Y}^n h)_{(\text{frame})}
  + \bm{\Gamma}^{(0)}_{(-1 -\delta)} (\mathscr{Y}^n X_{(\text{frame})})
  + \bm{\Gamma}^{(n-1)}_{(-1+3C_{(n-1)}\epsilon)}
 \end{split}
\end{equation}

\end{proposition}

\begin{proof}

Recall proposition \ref{proposition transport mu}, which we can write schematically as
\begin{equation*}
 \begin{split}
  L\log \mu &= r^{-2} x^i L^i_{(\text{small})} + r^{-1}L^i_{(\text{small})}L^i_{(\text{small})} + (\partial h)_{LL} + (\bar{\partial} h)_{(\text{frame})} \\
  &= r^{-2} x^i \bar{X}_{(\text{frame})} + r^{-1} (\bar{X}_{(\text{frame})})^2 + (\partial h)_{LL} + (\bar{\partial} h)_{(\text{frame})}
  \end{split}
\end{equation*}
Commuting with $\mathscr{Y}^n$, we find
\begin{equation*}
 \begin{split}
  \slashed{\D}_L \left( \mathscr{Y}^n \log \mu \right) &= 
  [\slashed{\D}_L , \mathscr{Y}^n]\log \mu
  + \bm{\Gamma}^{(0)}_{(-1)} \left( r^{-1} \mathscr{Y}^n \bar{X}_{(\text{frame})} \right)
  + \sum_{j+k \leq n-1} \bm{\Gamma}^{(j)}_{(-1 + C_{(j)}\epsilon , \text{large})}
  (\mathscr{Y}^{k+1} \bar{X}_{(\text{frame})}) \\
  &\phantom{=}
  + r^{-3}\sum_{j+k \leq n-1} \mathscr{Y}^{j} (r\bar{X}_{(\text{frame})}) \mathscr{Y}^{k} (r\bar{X}_{(\text{frame})})
  + (\slashed{\D} \mathscr{Y}^n h)_{(\text{frame})} \\
  &\phantom{=}
  + \bm{\Gamma}^{(0)}_{(-1 -\delta)} (\mathscr{Y}^n X_{(\text{frame})})
  + \bm{\Gamma}^{(n-1)}_{(-1+C_{(n-1)}\epsilon)}
  + \sum_{j+k \leq n} \left([\slashed{\D}, \mathscr{Y}^j] h_{(\text{rect})} \right) \left( \mathscr{Y}^{k} X_{(\text{frame})} \right)^2 
  \\
  \\
  \\
  &=
  \bm{\Gamma}^{(0)}_{(-1)} \mathscr{Y}^n \log \mu
  + \sum_{j+k \leq n-1} \bm{\Gamma}^{(j+1)}_{(-1+C_{(j+1)}\epsilon)} \mathscr{Y}^{k+1} \log \mu
  + \bm{\Gamma}^{(0)}_{(-1, \text{large})} \left( \mathscr{Z}^n \bar{X}_{(\text{frame})} \right)
  \\
  &\phantom{=}
  + \sum_{j+k \leq n-1} \bm{\Gamma}^{(j)}_{(C_{(j)}\epsilon , \text{large})}
  (r^{-1}\mathscr{Y}^{k+1} \bar{X}_{(\text{frame})})
  + r^{-3}\sum_{j+k \leq n-1} \mathscr{Y}^{j}(r\bar{X}_{(\text{frame})}) \mathscr{Y}^{k}(r\bar{X}_{(\text{frame})})
  \\
  &\phantom{=}
  + (\slashed{\D} \mathscr{Y}^n h)_{(\text{frame})}
  + \bm{\Gamma}^{(0)}_{(-1 -\delta)} (\mathscr{Y}^n X_{(\text{frame})})
  + \bm{\Gamma}^{(n-1)}_{(-1+C_{(n-1)}\epsilon)} \\
  &\phantom{=}
  + \sum_{j+k \leq n} \left( \bm{\Gamma}^{(j)}_{(C_{(j)}\epsilon)}\bm{\Gamma}^{(0)}_{(-1-\delta)}
  	+ \bm{\Gamma}^{(j-1)}_{(-1 + 2C_{(j-1)}\epsilon)} \right) \bm{\Gamma}^{(k)}_{(C_{(k)}\epsilon)} 
  \\
  \\
  \\
  &=\bm{\Gamma}^{(0)}_{(-1)} \mathscr{Y}^n \log \mu
  + \sum_{j+k\leq n-1} \bm{\Gamma}^{(j+1)}_{(-1-\delta)} (\mathscr{Y}^k \log \mu)
  + \bm{\Gamma}^{(0)}_{(-1, \text{large})} \left( \mathscr{Y}^n \bar{X}_{(\text{frame})} \right) \\
  &\phantom{=}
  + (\slashed{\D} \mathscr{Y}^n h)_{(\text{frame})}
  + \bm{\Gamma}^{(0)}_{(-1 -\delta)} (\mathscr{Y}^n X_{(\text{frame})})
  + \bm{\Gamma}^{(n-1)}_{(-1+3C_{(n-1)}\epsilon)}
 \end{split}
\end{equation*}

\end{proof}

\begin{remark}
The expression given in the proposition above ``loses'' a derivative in a certain sense: the foliation density $\mu$ is expected to behave like a metric components $h$, since its derivative $\omega$ is a connection coefficient. However, by integrating the expression for $\slashed{\D}_L \left(\mathscr{Z}^n \log \mu \right)$ above, we see that we can estimate $\mathscr{Z}^n \log \mu$ on the level of $(\slashed{\D} \mathscr{Z}^n h)$. In other words, $\log\mu$ can be estimated in terms of the \emph{first derivatives} of the metric, rather than the metric itself, and for this reason we say that this expression loses derivatives.
\end{remark}

\begin{proposition}[An expression for $\mathscr{Y}^n \omega$]
\label{proposition Yn omega}
 We have, schematically,
\begin{equation}
 \begin{split}
   \mathscr{Y}^n \omega
   &=
   \left( r^{-1}\bm{\Gamma}^{(0)}_{(0, \text{large})} + \bm{\Gamma}^{(0)}_{(-1+C_{(0)}\epsilon)} \right)\left( \mathscr{Z}^n \bar{X}_{(\text{frame})} \right)
   + \bm{\Gamma}^{(0)}_{(-1-\delta)} \left( \mathscr{Z}^n X_{(\text{frame})} \right)
   + (\slashed{\D} \mathscr{Z}^n h)_{LL}
   \\
   &\phantom{=}
   + (\overline{\slashed{\D}} \mathscr{Z}^n h)_{(\text{frame})}
   + \bm{\Gamma}^{(n)}_{(C_{(n)}\epsilon)}\bm{\Gamma}^{(0)}_{(-1-\delta)}
   + \bm{\Gamma}^{(n-1)}_{(-1+ 3C_{(n-1)}\epsilon)}
 \end{split}
\end{equation}

\end{proposition}

\begin{proof}
 Recall definition \ref{definition omega}, which gives $\omega = L\log\mu$. Combining this with proposition \ref{proposition transport mu} we have, schematically,
\begin{equation*}
 \omega = r^{-2} x^i \bar{X}_{(\text{frame})} + r^{-1}\bar{X}_{(\text{frame})}^2 + (\partial h)_{LL} + (\bar{\partial} h)_{(\text{frame})}
\end{equation*}

Now, applying the operator $\mathscr{Y}^n$ we have
\begin{equation*}
 \begin{split}
  \mathscr{Y}^n \omega
  &=
  \bm{\Gamma}^{(0)}_{(0, \text{large})} \left( r^{-1}\mathscr{Y}^n \bar{X}_{(\text{frame})} \right)
  +\sum_{j+k \leq n-1} \bm{\Gamma}^{(j)}_{(C_{(j)}\epsilon, \text{large})} \bm{\Gamma}^{(k)}_{(-1-\delta)} 
  \\
  &\phantom{=}
  + r^{-3}\sum_{j+k \leq n} \mathscr{Y}^j \left( r\bar{X}_{(\text{frame})} \right) \mathscr{Y}^k \left(r \bar{X}_{(\text{frame})} \right) 
  \\
  &\phantom{=}
  + \bm{\Gamma}^{(0)}_{(-1+C_{(0)}\epsilon)} (\mathscr{Y}^k L^a)
  + \bm{\Gamma}^{(0)}_{(-1-\delta)} (\mathscr{Y}^k X_{(\text{frame})})
  \\
  &\phantom{=}
  + \sum_{\substack{j+k \leq n \\ k \leq n-1}} \left( (\slashed{\D} \mathscr{Y}^j h_{(\text{rect})}) + [\slashed{\D}, \mathscr{Y}^j] h_{(\text{rect})}) \right) (\mathscr{Y}^k X_{(\text{frame})})^2
 \end{split}
\end{equation*}
Now we note that, schematically,
\begin{equation*}
[\mathscr{Y}^n , \slashed{\D}] h_{(\text{rect})}
=
\bm{\Gamma}^{(n)}_{(C_{(n)}\epsilon)}\bm{\Gamma}^{(0)}_{(-1-\delta)} + \bm{\Gamma}^{(n-1)}_{(-1+ 2C_{(n-1)}\epsilon)}
\end{equation*}
Additionally, we note that
\begin{equation*}
\begin{split}
	L^0 &= 1 + L^0_{(\text{small})} \\
	L^i &= \frac{x^i}{r} + L^i_{(\text{small})} 
\end{split}
\end{equation*}
and so, schematically,
\begin{equation*}
\mathscr{Y}^n L^a = \mathscr{Y}^n \bar{X}_{(\text{frame})} + \mathscr{Y}^{n-1} X_{(\text{frame})}
\end{equation*}
and so finally
\begin{equation*}
\begin{split}
	\mathscr{Y}^n \omega
	&=
	\left( r^{-1}\bm{\Gamma}^{(0)}_{(0, \text{large})} + \bm{\Gamma}^{(0)}_{(-1+C_{(0)}\epsilon)} \right)\left( \mathscr{Z}^n \bar{X}_{(\text{frame})} \right)
	+ \bm{\Gamma}^{(0)}_{(-1-\delta)} \left( \mathscr{Z}^n X_{(\text{frame})} \right)
	+ (\slashed{\D} \mathscr{Z}^n h)_{LL}
	\\
	&\phantom{=}
	+ (\overline{\slashed{\D}} \mathscr{Z}^n h)_{(\text{frame})}
	+ \bm{\Gamma}^{(n)}_{(C_{(n)}\epsilon)}\bm{\Gamma}^{(0)}_{(-1-\delta)}
	+ \bm{\Gamma}^{(n-1)}_{(-1+ 3C_{(n-1)}\epsilon)}
\end{split}
\end{equation*}

\end{proof}

\begin{proposition}[An expression for $\mathscr{Y}^n \zeta$]
 \label{proposition Yn zeta}
 The derivatives of the connection coefficient $\zeta$ satisfy an equation of the following schematic form:
\begin{equation}
\begin{split}
 \mathscr{Y}^n \zeta
 &= (\slashed{\D} \mathscr{Y}^n h)_{(\text{frame})}
 + \bm{\Gamma}^{(0)}_{(-1 + C_{(0)}\epsilon)} (\mathscr{Y}^n X_{(\text{frame})})
 + r^{-1} \bm{\Gamma}^{(0)}_{(0, \text{large})} (\mathscr{Y}^n \bar{X}_{(\text{frame})})
 + \bm{\Gamma}^{(n)}_{(C_{(n)}\epsilon)} \bm{\Gamma}^{(0)}_{(-1-\delta)}
 \\
 &\phantom{=}
 + \bm{\Gamma}^{(n-1)}_{(-1+2C_{(n-1)}\epsilon)}
\end{split}
\end{equation}
\end{proposition}

\begin{proof}
Recall proposition \ref{proposition zeta}, which we can write schematically as
\begin{equation*}
 \zeta = (\partial h)_{(\text{frame})} + r^{-1} X_{(\text{frame})} \bar{X}_{(\text{frame})}
\end{equation*}
where we have made use of the fact that
\begin{equation*}
 L^i \frac{ \slashed{\nabla}_\mu x^i}{r} = r^{-1} L^i_{(\text{small})} \slashed{\Pi}_\mu^{\phantom{\mu}i}
\end{equation*}
which in turn follows from the fact that $\slashed{\nabla} r^2 = \slashed{\nabla}(x^i x^i) = 0$.

Applying an operator $\mathscr{Y}$ $n$ times we find
\begin{equation*}
 \begin{split}
  \mathscr{Z}^n \zeta
  &=
  \sum_{j+k+l \leq n} \left( (\slashed{\D} \mathscr{Y}^j h_{(\text{rect})}) + [\slashed{\D} ,  \mathscr{Y}^j] h_{(\text{rect})}  \right)(\mathscr{Y}^k X_{(\text{frame})})(\mathscr{Z}^l X_{(\text{frame})}) \\
  &\phantom{=}
  + r^{-1} \sum_{j+k \leq n} (\mathscr{Y}^{j} X_{(\text{frame})})(\mathscr{Y}^{k} \bar{X}_{(\text{frame})})
  \\
  \\
  &= (\slashed{\D} \mathscr{Y}^n h)_{(\text{frame})}
  + \bm{\Gamma}^{(0)}_{(-1 + C_{(0)}\epsilon)} (\mathscr{Y}^n X_{(\text{frame})})
  + r^{-1} \bm{\Gamma}^{(0)}_{(0, \text{large})} (\mathscr{Z}^n \bar{X}_{(\text{frame})})
  + \bm{\Gamma}^{(n)}_{(C_{(n)}\epsilon)} \bm{\Gamma}^{(0)}_{(-1-\delta)}
  \\
  &\phantom{=}
  + \bm{\Gamma}^{(n-1)}_{(-1+2C_{(n-1)}\epsilon)}
 \end{split}
\end{equation*}

\end{proof}

\begin{proposition}[A transport equation for $\mathscr{Y}^n \tr_{\slashed{g}} \chi$ for a low number of derivatives]
	\label{proposition Yn tr chi low}
The commuted trace of the renormalised second fundamental form $\mathscr{Z}^n \tr_{\slashed{g}} \chi_{(\text{small})}$ satisfies a transport equation along the integral curves of $L$ which can be given schematically by the following system:
\begin{equation}
 \begin{split}
  \mathscr{Y}^n \mathcal{X}_{(\text{low})}
  &=
  \mathscr{Y}^n \tr_{\slashed{g}}\chi_{(\text{small})}
  + (\overline{\slashed{\D}} \mathscr{Y}^n h)_{(\text{frame})}
  + \bm{\Gamma}^{(0)}_{(-1-\delta)} (\mathscr{Z}^n X_{(\text{frame})})
  + \bm{\Gamma}^{(n)}_{(-\delta)} \bm{\Gamma}^{(0)}_{(-1 + C_{(0)}\epsilon)}
  \\
  &\phantom{=}
  + \bm{\Gamma}^{(0)}_{(-\delta)} \bm{\Gamma}^{(n)}_{(-1 + C_{(n)}\epsilon)} 
  + \bm{\Gamma}^{(n-1)}_{(-1-\delta)}
  \\
  \\
  \slashed{\D}_L \left( r^2 \mathscr{Y}^n \mathcal{X}_{(\text{low})} \right)
  &= \bm{\Gamma}^{(0)}_{(-1)} \left( r^2 \mathscr{Z}^n \mathcal{X}_{(\text{low})} \right)
  + r(\overline{\slashed{\D}} \mathscr{Y}^{n+1} h)_{(\text{frame})}
  + r(\overline{\slashed{\D}} \mathscr{Y}^{n} h_{(\text{rect})}) \bm{\Gamma}^{(1)}_{(C_{(1)}\epsilon, \text{large})} 
  \\
  &\phantom{=}
  + \bm{\Gamma}^{(0)}_{(1-\delta)} (\slashed{\D} \mathscr{Y}^n h)_{(\text{frame})}
  + \bm{\Gamma}^{(1)}_{(-\delta)} (\mathscr{Y}^n X_{(\text{frame})})
  + \bm{\Gamma}^{(0)}_{(-\frac{1}{2}+\delta)} \bm{\Gamma}^{(n)}_{(C_{(n)}\epsilon)}
  \\
  &\phantom{=}
  + r^2 \bm{\Gamma}^{(0)}_{(-1-\delta)} \bm{\Gamma}^{(n)}_{(-1-\delta)}
  + \bm{\Gamma}^{(n-1)}_{(-\delta)}
 \end{split}
\end{equation}

\end{proposition}

\begin{proof}
Recall proposition \ref{proposition transport trace chi}, which we can write schematically as the following system:
\begin{equation*}
 \begin{split}
  \mathcal{X}_{(\text{low})} &= \tr_{\slashed{g}}\chi_{(\text{small})} + (\bar{\partial} h)_{(\text{frame})} 
  \\ 
  \\
  \\
  L \left( r^2 \mathcal{X}_{(\text{low})} \right)
  &= r^{-2} (r^2 \mathcal{X}_{(\text{low})})^2
  + r(\overline{\slashed{\D}}\mathscr{Z} h)_{(\text{frame})}
  + r^2 (\hat{\chi})^2
  + r^2 (\bar{\partial} h)_{(\text{frame})} (\partial h)_{(\text{frame})} \\
  &\phantom{=}
  + r^2 \bm{\Gamma} (\bar{\partial} h)_{(\text{frame})}
  + r X_{(\text{frame})} \bar{X}_{(\text{frame})} (\bar{\partial} h)_{(\text{frame})}
  \\
  \\
  &= r^{-2} (r^2 \mathcal{X}_{(\text{low})})^2
  + r(\overline{\slashed{\D}}\mathscr{Z} h)_{(\text{frame})}
  + r^2 (\partial h)_{(\text{frame)}} (\bar{\partial} h)_{(\text{frame)}}
  + r^2 \bm{\Gamma}^{(0)}_{(-1-\delta)}\bm{\Gamma}^{(0)}_{(-1-\delta)}
 \end{split}
\end{equation*}
where are making use of the expression in terms of $\mathcal{X}_{(\text{low})}$ rather than $\mathcal{X}_{(\text{high})}$.

If we commute $n$ times with a commutation operator $\mathscr{Y}$, and make use of the fact that
\begin{equation*}
[\mathscr{Y}^n \, , \overline{\slashed{\D}}] h_{(\text{rect})}
=
\bm{\Gamma}^{(n)}_{(-1 + C_{(n)}\epsilon)} \bm{\Gamma}^{(0)}_{(-\frac{1}{2} + \delta)}
+ \bm{\Gamma}^{(n-1)}_{(-1-\delta)}
\end{equation*}
(which follows from proposition \ref{proposition commute D Yn}), then we obtain the following system:
\begin{equation*}
 \begin{split}
  \mathscr{Y}^n \mathcal{X}_{(\text{low})}
  &=
  \mathscr{Y}^n \tr_{\slashed{g}}\chi_{(\text{small})}
  + (\overline{\slashed{\D}} \mathscr{Y}^n h)_{(\text{frame})}
  + \bm{\Gamma}^{(0)}_{(-1-\delta)} (\mathscr{Z}^n X_{(\text{frame})})
  \\
  &\phantom{=}
  + \bm{\Gamma}^{(0)}_{(-\frac{1}{2} + \delta)} \bm{\Gamma}^{(n)}_{(-1 + C_{(n)}\epsilon)} 
  + \bm{\Gamma}^{(n-1)}_{(-1-\delta)}
  \\
  \\
  \\
  \slashed{\D}_L \left( r^2 \mathscr{Y}^n \mathcal{X}_{(\text{low})} \right)
  &=
  [\slashed{\D}_L , \mathscr{Y}^n] r^2 \mathcal{X}_{(\text{low})}
  + r^{-2} \sum_{j+k\leq n} \left( \mathscr{Y}^j r^2 \mathcal{X}_{(\text{low})} \right)\left( \mathscr{Y}^k r^2 \mathcal{X}_{(\text{low})} \right) \\
  &\phantom{=}
  + r \sum_{j+ \leq n} \left( (\overline{\slashed{\D}} \mathscr{Y}^{j+1} h_{(\text{rect})}) + [\overline{\slashed{\D}} , \mathscr{Y}^{j}] \mathscr{Z} h_{(\text{rect})}) \right) \bm{\Gamma}^{(k)}_{(C_{(k)}\epsilon, \text{large})} \\
  &\phantom{=}
  + r \sum_{j+k+\ell \leq n} 
  	\left( (\overline{\slashed{\D}} \mathscr{Y}^{j+1} h_{(\text{rect})}) + [\overline{\slashed{\D}} , \mathscr{Y}^{j}] \mathscr{Z} h_{(\text{rect})}) \right)
  	\\
  	&\phantom{= + r \sum_{j+k+\ell \leq n}}
  	\times \left( (\slashed{\D} \mathscr{Y}^{k+1} h_{(\text{rect})}) + [\slashed{\D} , \mathscr{Y}^{k}] \mathscr{Z} h_{(\text{rect})}) \right)
  	\bm{\Gamma}^{(\ell)}_{(2C_{(\ell)}\epsilon, \text{large})} 
  \\
  &\phantom{=}
  + r^2 \sum_{j+k \leq n} \bm{\Gamma}^{(j)}_{(-1-\delta)}\bm{\Gamma}^{(k)}_{(-1-\delta)}
  \\
  \\
  &= \sum_{j \leq n} \bm{\Gamma}^{(0)}_{(-1)} \left( r^2 \mathscr{Y}^j \mathcal{X}_{(\text{low})} \right)
  + \sum_{j+k \leq n-1} \bm{\Gamma}^{(j+1)}_{(-1+C_{(j+1)}\epsilon)} \left(\overline{\slashed{\D}} (r^2 \mathscr{Y}^k \mathcal{X}_{(\text{low})}) \right) \\
  &\phantom{=}
  + r^2 \sum_{j+k \leq n-1} \bm{\Gamma}^{(j+1)}_{(-1-\delta)} \mathscr{Y}^k \mathcal{X}_{(\text{low})}
  + \bm{\Gamma}^{(0)}_{(-1-\delta)} \left( r^2 \mathscr{Y}^n \mathcal{X}_{(\text{low})} \right)
  + r(\overline{\slashed{\D}} \mathscr{Y}^{n+1} h)_{(\text{frame})}
  \\
  &\phantom{=}
  + r(\overline{\slashed{\D}} \mathscr{Y}^{n} h_{(\text{rect})}) \bm{\Gamma}^{(1)}_{(C_{(1)}\epsilon, \text{large})}
  + \bm{\Gamma}^{(0)}_{(1-\delta)} (\slashed{\D} \mathscr{Y}^n h)_{(\text{frame})}
  + \bm{\Gamma}^{(1)}_{(-\delta)} (\mathscr{Y}^n X_{(\text{frame})})
  \\
  &\phantom{=}
  + \bm{\Gamma}^{(0)}_{(-\frac{1}{2}+\delta)} \bm{\Gamma}^{(n)}_{(C_{(n)}\epsilon)}
  + r^2\bm{\Gamma}^{(n)}_{(-1-\delta)} \bm{\Gamma}^{(0)}_{(-1-\delta)}
  + \bm{\Gamma}^{(n-1)}_{(-\delta)}
  \\
  \\
  &= \bm{\Gamma}^{(0)}_{(-1)} \left( r^2 \mathscr{Z}^n \mathcal{X}_{(\text{low})} \right)
  + r(\overline{\slashed{\D}} \mathscr{Y}^{n+1} h)_{(\text{frame})}
  + r(\overline{\slashed{\D}} \mathscr{Y}^{n} h_{(\text{rect})}) \bm{\Gamma}^{(1)}_{(C_{(1)}\epsilon, \text{large})} 
  \\
  &\phantom{=}
  + \bm{\Gamma}^{(0)}_{(1-\delta)} (\slashed{\D} \mathscr{Y}^n h)_{(\text{frame})}
  + \bm{\Gamma}^{(1)}_{(-\delta)} (\mathscr{Y}^n X_{(\text{frame})})
  + \bm{\Gamma}^{(0)}_{(-\frac{1}{2}+\delta)} \bm{\Gamma}^{(n)}_{(C_{(n)}\epsilon)}
  \\
  &\phantom{=}
  + r^2\bm{\Gamma}^{(n)}_{(-1-\delta)}\bm{\Gamma}^{(0)}_{(-1-\delta)}
  + \bm{\Gamma}^{(n-1)}_{(-\delta)}
 \end{split}
\end{equation*}

\end{proof}

The proposition above also has the undesirable feature that is \emph{loses a derivative}. That is, we can use it to estimate $\mathscr{Z}^n \tr_{\slashed{g}}\chi_{(\text{small})}$ in terms of $(\overline{\slashed{\D}} \mathscr{Z}^{n+1} h)_{(\text{frame})}$, together with other, lower order quantities. This is not a problem when we are engaged in proving pointwise bounds, where some loss of derivatives is unavoidable in any case, occurring, for example, in the Sobolev inequalities.

However, the loss of derivatives will lead to problems when we are required to estimate $\mathscr{Z}^n \tr_{\slashed{g}} \chi_{(\text{small})}$ in $L^2$. Since $\tr_{\slashed{g}} \chi_{(\text{small})}$ is one of the connection coefficients, we can hope to estimate it on the same level as the derivatives of the metric; i.e.\ we could hope to be able to estimate $\mathscr{Z}^n \tr_{\slashed{g}}\chi_{(\text{small})}$ in terms of $(\slashed{\D} \mathscr{Z}^{n} h)_{(\text{frame})}$. It turns out that we require such a bound in order to close our estimates, however, it cannot be obtained by using the previous expression. Instead, we must use the following proposition.

\begin{proposition}[A transport equation for $\mathscr{Y}^n \tr_{\slashed{g}} \chi$ for a high number of derivatives]
	\label{proposition Yn tr chi high}
Suppose that the rectangular components of the metric satisfy
\begin{equation*}
 \tilde{\Box}_g h_{ab} = F_{ab}
\end{equation*}
Then the quantity $\mathscr{Y}^n \tr_{\slashed{g}} \chi_{(\text{small})}$ obeys a transport equation along the integral curves of $L$ which can be written in the following form:
\begin{equation*}
 \begin{split}
  \mathscr{Y}^n \mathcal{X}_{(\text{high})}
  &=
  \mathscr{Y}^n \tr_{\slashed{g}}\chi_{(\text{small})}
  + (\slashed{\D}\mathscr{Y}^n h)_{LL}
  + (\overline{\slashed{\D}}\mathscr{Y}^n h)_{(frame)}
  + \bm{\Gamma}^{(0)}_{(-1+C_{(0)}\epsilon)} (\mathscr{Y}^n \bar{X}_{(\text{frame})})
  \\
  &\phantom{=}
  + \bm{\Gamma}^{(0)}_{(-1-\delta)} (\mathscr{Y}^n X_{(\text{frame})})
  + \bm{\Gamma}^{(n-1)}_{(-1+2C_{(n-1)}\epsilon)}
  \\
  \\
  \\
  \slashed{\D}_L \left( r^2 \mathscr{Y}^n\mathcal{X}_{(\text{high})} \right)
  &= 
  \bm{\Gamma}^{(0)}_{(-1)} \left( r^2 \mathscr{Y}^n \mathcal{X}_{(\text{high})} \right)
  + r^2 (\mathscr{Y}^n F)_{LL}
  + r^2 (F)_{(\text{frame})} (\mathscr{Y}^n \bar{X}_{(\text{frame})})
  \\
  &\phantom{=}
  + r^2 \sum_{\substack{j+k \leq n \\ j,k \leq n-1}} \bm{\Gamma}^{(j)}_{(C_{(j)}\epsilon, \text{large})} (\mathscr{Y}^{k} F)_{(\text{frame})}
  + r(\mathscr{Y}^n \omega)
  + r^2 \bm{\Gamma}^{(0)}_{(-1-\delta)} (\slashed{\D} \mathscr{Y}^n h)_{(\text{frame})}
  \\
  &\phantom{=}
  + r^2 \bm{\Gamma}^{(0)}_{(-1+C_{(0)}\epsilon)} (\overline{\slashed{\D}} \mathscr{Y}^n h)_{(\text{frame})}
  + \bm{\Gamma}^{(0)}_{(-\frac{1}{2} + \delta)} \bm{\Gamma}^{(n)}_{(C_{(n)}\epsilon)}
  + r^2 \bm{\Gamma}^{(n)}_{(-1-\delta)} \bm{\Gamma}^{(0)}_{(-1-\delta)}
  \\
  &\phantom{=}
  + r^2 \sum_{\substack{j+k \leq n \\ j,k \leq n-1}} \bm{\Gamma}^{(j)}_{(-1+C_{(j)}\epsilon)} \bm{\Gamma}^{(k)}_{(-1+C_{(k)}\epsilon)}
 \end{split}
\end{equation*}

\end{proposition}

\begin{proof}
 Using the second expression given in proposition \ref{proposition transport trace chi} we have the following system for $\tr_{\slashed{g}}\chi_{(\text{small})}$, given schematically:
\begin{equation*}
 \begin{split}
  \mathcal{X}_{(\text{high})} &= \tr_{\slashed{g}} \chi_{(\text{small})} + (\partial h)_{LL} + (\bar{\partial} h)_{(\text{frame})}
  \\
  \\
  \\
  L\left(r^2 \mathcal{X}_{(\text{high})} \right)
  &=
  r^2 (\tilde{\Box}_g h)_{LL}
  + \omega \left(r^2 \mathcal{X}_{(\text{high})} \right)
  + r\omega
  + r^{-2} \left(r^2 \mathcal{X}_{(\text{high})} \right)^2
  + r^2 (\hat{\chi})^2
  + r^2 (\bar{\partial} h)_{(\text{frame})} (\partial h)_{(\text{frame})} 
  \\
  &\phantom{=}
  + r^2 \left(\bm{\Gamma} + r^{-1} \bar{X}_{(\text{frame})}^2 + r^{-1} X_{(\text{frame})} \bar{X}_{(\text{frame})} \right) (\bar{\partial} h)_{(\text{frame})}
  \\
  \\
  &=
  r^2 F_{LL}
  + \bm{\Gamma}^{(0)}_{(-1)}\left( r^2 \mathcal{X}_{(\text{high})} \right)
  + r\omega
  + r^{-2} \left(r^2 \mathcal{X}_{(\text{high})}\right)^2
  + r^2 (\partial h)_{(\text{frame})}\cdot (\bar{\partial} h)_{(\text{frame})}
  \\
  &\phantom{=}
  + r^2 \bm{\Gamma}^{(0)}_{(-1-\delta)} \cdot \bm{\Gamma}^{(0)}_{(-1-\delta)}
 \end{split}
\end{equation*}

Now commuting $n$ times with an operator $\mathscr{Y}$, we obtain a system of the form
\begin{equation*}
 \begin{split}
  \mathscr{Y}^n \mathcal{X}_{(\text{high})}
  &=
  \mathscr{Y}^n \tr_{\slashed{g}}\chi_{(\text{small})}
  + (\slashed{\D}\mathscr{Y}^n h)_{LL}
  + (\overline{\slashed{\D}}\mathscr{Y}^n h)_{(frame)}
  + \bm{\Gamma}^{(0)}_{(-1+C_{(0)}\epsilon)} (\mathscr{Y}^n \bar{X}_{(\text{frame})})
  \\
  &\phantom{=}
  + \bm{\Gamma}^{(0)}_{(-1-\delta)} (\mathscr{Y}^n X_{(\text{frame})})
  + \bm{\Gamma}^{(n-1)}_{(-1+2C_{(n-1)}\epsilon)}
  \\
  \\
  \\
  \slashed{\D}_L \left( r^2 \mathscr{Y}^n\mathcal{X}_{(\text{high})} \right)
  &=
  [\slashed{\D}_L , \mathscr{Y}^n] r^2 \mathcal{X}_{(\text{high})}
  + r^2(\mathscr{Y}^n F)_{LL}
  + r^2(F)_{(\text{frame})} (\mathscr{Y}^n \bar{X}_{(\text{frame})})
  \\
  &\phantom{=}
  + r^2 \sum_{ \substack{j+k+l \leq n \\ j,k,l \leq n-1}} (\mathscr{Y}^j F_{(\text{rect})}) (\mathscr{Y}^k X_{(\text{frame})})(\mathscr{Y}^l X_{(\text{frame})})
  + r \sum_{j\leq n} (\mathscr{Y}^j \omega)
  \\
  &\phantom{=}
  + \bm{\Gamma}^{(0)}_{(-1)} r^2 (\mathscr{Y}^n \mathcal{X}_{(\text{high})})
  + r^2 \bm{\Gamma}^{(0)}_{(-1-\delta)} (\slashed{\D} \mathscr{Y}^n h)_{(\text{frame})}
  + r^2 \bm{\Gamma}^{(0)}_{(-1+C_{(0)}\epsilon)} (\overline{\slashed{\D}} \mathscr{Y}^n h)_{(\text{frame})}
  \\
  &\phantom{=}
  + r^2 \bm{\Gamma}^{(n)}_{(-1-\delta)} \bm{\Gamma}^{(0)}_{(-1-\delta)}
  + r^2 \sum_{\substack{j+k\leq n \\ j \leq n-1 \\ k \leq n-1}} \bm{\Gamma}^{(j)}_{(-1+C_{(j)}\epsilon)} \bm{\Gamma}^{(k)}_{(-1+C_{(k)}\epsilon)} 
  \\
  \\
  &= 
  \bm{\Gamma}^{(0)}_{(-1)} \left( r^2 \mathscr{Y}^n \mathcal{X}_{(\text{high})} \right)
  + r^2 (\mathscr{Y}^n F)_{LL}
  + r^2 (F)_{(\text{frame})} (\mathscr{Y}^n \bar{X}_{(\text{frame})})
  \\
  &\phantom{=}
  + r^2 \sum_{\substack{j+k \leq n \\ j,k \leq n-1}} \bm{\Gamma}^{(j)}_{(C_{(j)}\epsilon, \text{large})} (\mathscr{Y}^{k} F)_{(\text{frame})}
  + r(\mathscr{Y}^n \omega)
  + r^2 \bm{\Gamma}^{(0)}_{(-1-\delta)} (\slashed{\D} \mathscr{Y}^n h)_{(\text{frame})}
  \\
  &\phantom{=}
  + \bm{\Gamma}^{(0)}_{(-\frac{1}{2} + \delta)} \bm{\Gamma}^{(n)}_{(C_{(n)}\epsilon)}
  + r^2 \bm{\Gamma}^{(0)}_{(-1+C_{(0)}\epsilon)} (\overline{\slashed{\D}} \mathscr{Y}^n h)_{(\text{frame})}
  + r^2 \bm{\Gamma}^{(n)}_{(-1-\delta)} \bm{\Gamma}^{(0)}_{(-1-\delta)}
  \\
  &\phantom{=}
  + r^2 \sum_{\substack{j+k \leq n \\ j,k\leq n-1}} \bm{\Gamma}^{(j)}_{(-1+C_{(j)}\epsilon)} \bm{\Gamma}^{(k)}_{(-1+C_{(k)}\epsilon)}
 \end{split}
\end{equation*}
where we have again made the substitution $\overline{\slashed{\D}} = r^{-1}\mathscr{Y}$.

\end{proof}

\begin{proposition}[A transport equation for the commuted tensor field $\hat{\chi}$]
\label{proposition transport Yn chihat}
After applying $n$ of the operators $\mathscr{Y}$ to the trace-free tensor field $\hat{\chi}$, the resulting fields satisfy transport equations along the integral curves of the vector field $L$ which can be written in the form
\begin{equation}
\begin{split}
\mathscr{Y}^n \hat{\mathcal{X}}
&=
\mathscr{Y}^n \hat{\chi}
+ (\overline{\slashed{\D}} \mathscr{Y}^n h)_{(\text{frame})}
+ (\bar{\partial} h_{(\text{rect})}) \bm{\Gamma}^{(0)}_{(0, \text{large})} (\mathscr{Y}^n X_{(\text{frame})})
+ \bm{\Gamma}^{(n)}_{(-\delta)} \bm{\Gamma}^{(0)}_{(-1+C_{(0)}\epsilon)}
\\
&\phantom{=}
+ \bm{\Gamma}^{(0)}_{(-\delta)} \bm{\Gamma}^{(n)}_{(-1+C_{(n)}\epsilon)}
+ \bm{\Gamma}^{(n-1)}_{(-1-\delta)}
\\
\\
\\
\slashed{\D}_L \left( r^2 \mathscr{Y}^n \hat{\mathcal{X}} \right)
&=
r(\overline{\slashed{\D}} \mathscr{Y}^{n+1} h)_{(\text{frame})}
+ r (\overline{\slashed{\D}}\mathscr{Z}h_{(\text{rect})}) \bm{\Gamma}^{(0)}_{(0, \text{large})} (\mathscr{Y}^n X_{(\text{frame})})
+ r \bm{\Gamma}^{(1)}_{(C_{(1)}\epsilon, \text{large})} (\overline{\slashed{\D}} \mathscr{Y}^n h_{(\text{rect})})
\\
&\phantom{=}
+ \bm{\Gamma}^{(0)}_{(-1)} (r^2 \mathscr{Y}^n \hat{\mathcal{X}})
+ r^2 \bm{\Gamma}^{(0)}_{(-1-\delta)} (\slashed{\D} \mathscr{Y}^n h)_{(\text{frame})}
+ \bm{\Gamma}^{(0)}_{(-\frac{1}{2} + \delta)} \bm{\Gamma}^{(n)}_{(C_{(n)}\epsilon)}
+ r^2 \bm{\Gamma}^{(0)}_{(-1-\delta)} \bm{\Gamma}^{(n)}_{(-1-\delta)}
\\
&\phantom{=}
+ r \bm{\Gamma}^{(n-1)}_{(-1-\delta)}
\end{split}
\end{equation}

\end{proposition}

\begin{proof}
	Recall proposition \ref{proposition transport chihat}, which, when combined with proposition \ref{proposition structure of alpha} can be written schematically in the form
	\begin{equation*}
	\begin{split}
	\hat{\mathcal{X}}
	&=
	\hat{\chi}
	+ (\bar{\partial} h)_{(\text{frame})}
	\\
	\\
	\slashed{\D}_L \left(r^2 \hat{\mathcal{X}} \right)
	&=
	r (\overline{\slashed{\D}}\mathscr{Z} h)_{(\text{frame})}
	+ r^2 \omega \hat{\mathcal{X}}
	+ r^{-2} (r^2 \hat{\mathcal{X}})^2
	+ r^2 (\partial h)_{(\text{frame})} (\bar{\partial} h)_{(\text{frame})}
	+ r^2 \bm{\Gamma}^{(0)}_{(-1-\delta)} \bm{\Gamma}^{(0)}_{(-1-\delta)}
	\end{split}
	\end{equation*}
	
	Now, commuting with $\mathscr{Y}^n$ we obtain a system of the form
	\begin{equation*}
	\begin{split}
	\mathscr{Y}^n \hat{\mathcal{X}}
	&=
	\mathscr{Y}^n \hat{\chi}
	+ (\overline{\slashed{\D}} \mathscr{Y}^n h)_{(\text{frame})}
	+ (\bar{\partial} h_{(\text{rect})}) \bm{\Gamma}^{(0)}_{(0, \text{large})} (\mathscr{Y}^n X_{(\text{frame})})
	+ \bm{\Gamma}^{(0)}_{(-\frac{1}{2}+\delta)} \bm{\Gamma}^{(n)}_{(-1+C_{(n)}\epsilon)}
	\\
	&\phantom{=}
	+ \bm{\Gamma}^{(n-1)}_{(-1-\delta)}
	\\
	\\
	\\
	\slashed{\D}_L \left( r^2 \mathscr{Y}^n \hat{\mathcal{X}} \right)
	&=
	[\slashed{\D}_L , \mathscr{Y}^n] r^2 \hat{\mathcal{X}}
	+ r(\overline{\slashed{\D}} \mathscr{Y}^{n+1} h)_{(\text{frame})}
	+ r([\overline{\slashed{\D}} , \mathscr{Y}^n] \mathscr{Z} h)_{(\text{frame})}
	\\
	&\phantom{=}
	+ r (\overline{\slashed{\D}}\mathscr{Z} h_{(\text{rect})}) \bm{\Gamma}^{(0)}_{(0, \text{large})} (\mathscr{Y}^n X_{(\text{frame})})
	+ r \bm{\Gamma}^{(1)}_{(C_{(1)}\epsilon, \text{large})} (\overline{\slashed{\D}} \mathscr{Y}^n h_{(\text{rect})})
	+ \bm{\Gamma}^{(0)}_{(-1)} (r^2 \mathscr{Y}^n \hat{\mathcal{X}})
	\\
	&\phantom{=}
	+ r^2 \bm{\Gamma}^{(0)}_{(-1-\delta)} \bm{\Gamma}^{(n)}_{(-1-\delta)}
	+ r^2 \bm{\Gamma}^{(0)}_{(-1-\delta)} (\slashed{\D} \mathscr{Y}^n h)_{(\text{frame})}
	+ r^2 \bm{\Gamma}^{(0)}_{(-1+C_{(0)}\epsilon)} (\overline{\slashed{\D}} \mathscr{Y}^n h)_{(\text{frame})}
	\\
	&\phantom{=}
	+ \bm{\Gamma}^{(0)}_{(-\frac{1}{2} + \delta)} \bm{\Gamma}^{(n)}_{(C_{(n)}\epsilon)}
	+ r \bm{\Gamma}^{(n-1)}_{(-1-\delta)}
	\\
	\\
	&=
	r(\overline{\slashed{\D}} \mathscr{Y}^{n+1} h)_{(\text{frame})}
	+ r (\overline{\slashed{\D}}\mathscr{Z}h_{(\text{rect})}) \bm{\Gamma}^{(0)}_{(0, \text{large})} (\mathscr{Y}^n X_{(\text{frame})})
	+ r \bm{\Gamma}^{(1)}_{(C_{(1)}\epsilon, \text{large})} (\overline{\slashed{\D}} \mathscr{Y}^n h_{(\text{rect})})
	\\
	&\phantom{=}
	+ \bm{\Gamma}^{(0)}_{(-1)} (r^2 \mathscr{Y}^n \hat{\mathcal{X}})
	+ r^2 \bm{\Gamma}^{(0)}_{(-1-\delta)} (\slashed{\D} \mathscr{Y}^n h)_{(\text{frame})}
	+ \bm{\Gamma}^{(0)}_{(-\frac{1}{2} + \delta)} \bm{\Gamma}^{(n)}_{(C_{(n)}\epsilon)}
	+ r^2 \bm{\Gamma}^{(0)}_{(-1-\delta)} \bm{\Gamma}^{(n)}_{(-1-\delta)}
	\\
	&\phantom{=}
	+ r \bm{\Gamma}^{(n-1)}_{(-1-\delta)}
	\end{split}
	\end{equation*}

\end{proof}

The proposition above is useful for obtaining information about the $L$ derivative of $\hat{\chi}$. However, this equation also ``loses a derivative'': in order to estimate $\mathscr{Y}^n \hat{\chi}$ using this equation, we must already be in possession of information regarding $(\overline{\slashed{\D}} \mathscr{Y}^{n+1} h)$. On the other hand, we would expect $\hat{\chi}$ to behave similarly to $(\bar{\partial} h)$. The resulting problems can be overcome with the use of elliptic estimates, which involve estimating the divergence of $\hat{\chi}$, that is, the tensor field 
\begin{equation*}
 \slashed{\nabla}^\nu \hat{\chi}_{\nu \mu}
\end{equation*}
In order to estimate sufficiently many derivatives of $\hat{\chi}$, we will in fact have to commute the equations for this quantity with the operators $\mathscr{Y}$. This is achieved in the next few propositions.

%
%
%\begin{proposition}
%Let $\phi$ be an $S_{\tau,r}$-tangent tensor field. Then we have, schematically
%\begin{equation}
%	[\slashed{g}^{-1} \slashed{\nabla} , \mathscr{Z}^n]\phi
%	\sum_{j+k \leq n-1} \bm{\Gamma}^{(j)}_{(-1+C_{(j)}\epsilon)} (\slashed{\D}\mathscr{Z}^k\phi)
%	+ \sum_{j+k \leq n-1} \bm{\Gamma}^{(j+1)}_{(-1+C_{(j+1)}\epsilon)}(\mathscr{Z}^k \phi)
%\end{equation}
%
%\end{proposition}
%
%\begin{proof}
%Using proposition \ref{proposition commuting rnabla with first order operators} we have, schematically,
%\begin{equation*}
%	[\slashed{g}^{-1} \slashed{\nabla} , \mathscr{Z}]\phi
%	= \bm{\Gamma} \cdot \slashed{\D} \phi
%	+ \begin{pmatrix} r^{-1} \\ \bm{\Gamma} \\ r^{-1}\bm{\Gamma} \\ \bm{\Gamma}\cdot\bm{\Gamma} \\ r\bm{\Gamma}\cdot\bm{\Gamma}_{(\text{good})} \\ r\slashed{R} \\ \slashed{R}_L \\ \slashed{R}_{\Lbar} \end{pmatrix} \cdot \phi
%\end{equation*} 
%for any $S_{\tau,r}$ tangent field $\phi$, where we have used $\slashed{\nabla} \slashed{g} = 0$ and $\slashed{\D}_T \slashed{g} = \bm{\Gamma}^{(0)}_{(-1 + C_{(0)}\epsilon)}$. Now, we claim that
%\begin{equation*}
%	\slashed{g}^{-1} \slashed{\nabla} \mathscr{Z}^n \phi
%	= \mathscr{Z}^n \left(\slashed{g}^{-1} \slashed{\nabla} \phi\right) 
%	+ \sum_{j+k \leq n-1} \bm{\Gamma}^{(j)}_{(-1+C_{(j)}\epsilon)} (\slashed{\D} \mathscr{Z}^k\phi) + \sum_{j+k \leq n-1} \bm{\Gamma}^{(j+1)}_{(-1 + C_{(j+1)}\epsilon)} (\mathscr{Z}^k \phi)
%\end{equation*}
%	
%To prove this claim, we follow almost exactly the same steps as in proposition \ref{proposition commute Zn D}.
%	
%\end{proof}

\begin{proposition}[Commuting the equation for the divergence of $\hat{\chi}$]
	\label{proposition commuted div chihat}
	Let $\phi_{\mu\nu}$ be a symmetric, trace-free, $S_{\tau,r}$-tangent tensor field. Define
	\begin{equation*}
	\slashed{\Div} \left(  r^n \slashed{\nabla}^n \phi \right)_{\mu \alpha_1 \ldots \alpha_n}
	:=
	(\slashed{g}^{-1})^{\sigma\rho} \slashed{\nabla}_\rho \left( r^n \slashed{\nabla}_{\alpha_1} \ldots \slashed{\nabla}_{\alpha_n} \phi_{\mu\sigma} \right)
	\end{equation*}
	
	Then the $S_{\tau,r}$-tangent tensor fields $\mathscr{Y}^n \hat{\chi}$ obey equations of the form
	\begin{equation*}
	\begin{split}
	\slashed{\Div} (\mathscr{Y}^n \hat{\chi})
	&=
	r^{-1} (\overline{\slashed{\D}} \mathscr{Y}^{n+1} h)_{(\text{frame})}
	+ r^{-1} \mathscr{Y}^{(n+1)} \tr_{\slashed{g}}\chi_{(\text{small})}
	+ \sum_{j \leq n} r^{-1} \bm{\Gamma}^{(j)}_{(-1+C_{(j)}\epsilon)}
	\\
	&\phantom{=}
	+ \sum_{j+k \leq n} \bm{\Gamma}^{(j)}_{(-1+C_{(j)}\epsilon)} \bm{\Gamma}^{(k)}_{(-1-\delta)}
	\end{split}
	\end{equation*}
	
\end{proposition}

\begin{proof}
	Recall proposition \ref{proposition div chihat}, which we can write schematically as
	\begin{equation*}
	\slashed{\Div} \hat{\chi} 
	= 
	\slashed{\nabla} \tr_{\slashed{g}} \chi_{(\text{small})}
	+ r^{-1} (\overline{\slashed{\D}} \mathscr{Z} h)_{(\text{frame})}
	+ r^{-1} \bm{\Gamma}^{(0)}_{(-1+C_{(0)}\epsilon)}
	+ \bm{\Gamma}^{(0)}_{(-1+C_{(0)}\epsilon)}\cdot \bm{\Gamma}^{(0)}_{(-1-\delta)}
	\end{equation*}
	Commuting with an operator $\mathscr{Y}$ $n$ times, we obtain an equation of the form
	\begin{equation*}
	\begin{split}
	\slashed{\Div} \mathscr{Y}^n \hat{\chi}
	&=
	[\slashed{g}^{-1}\slashed{\nabla}, \mathscr{Y}^n] \hat{\chi}
	+ \sum_{j \leq n+1} r^{-1} \mathscr{Y}^{j} \tr_{\slashed{g}}\chi_{(\text{small})}
	+ r^{-1} (\overline{\slashed{\D}} \mathscr{Y}^{n+1} h)_{(\text{frame})}
	+ r^{-1} \bm{\Gamma}^{(n)}_{(-1+C_{(n)}\epsilon)}
	\\
	&\phantom{=}
	+ \sum_{j+k \leq n} \bm{\Gamma}^{(j)}_{(-1+C_{(j)}\epsilon)} \bm{\Gamma}^{(k)}_{(-1-\delta)}
	\end{split}
	\end{equation*}
	Recall that $\slashed{\D}$ is a metric connection with metric $\slashed{g}$, so $\mathscr{Y}$ commutes with $\slashed{g}$ and $\slashed{g}^{-1}$. Hence, making use of the propositions above we have
	\begin{equation*}
	\begin{split}
	\slashed{\Div} \mathscr{Y}^n \hat{\chi}
	&=
	\sum_{j \leq n+1} r^{-1} \mathscr{Y}^{j} \tr_{\slashed{g}}\chi_{(\text{small})}
	+ r^{-1} (\overline{\slashed{\D}} \mathscr{Y}^{n+1} h)_{(\text{frame})}
	+ r^{-1} \bm{\Gamma}^{(n)}_{(-1+C_{(n)}\epsilon)}
	\\
	&\phantom{=}
	+ \sum_{j+k \leq n} \bm{\Gamma}^{(j)}_{(-1+C_{(j)}\epsilon)} \bm{\Gamma}^{(k)}_{(-1-\delta)}
	\end{split}
	\end{equation*}
	
\end{proof}

\begin{proposition}[Commuting the equation for $\slashed{\D}_T \mathscr{Y}^n \hat{\chi}$]
	\label{proposition DT Yn chi hat}
\end{proposition}

\begin{proof}
	Proposition \ref{proposition Lbar chi} gives
	\begin{equation*}
	\begin{split}
	\slashed{\D}_{\Lbar} \hat{\chi}_{\mu\nu}
	&=
	\slashed{\Pi}_\mu^{\phantom{\mu}\rho} \slashed{\Pi}_\nu^{\phantom{\nu}\sigma} R_{\Lbar \rho L \sigma}
	- \frac{1}{2} \slashed{g}_{\mu\nu} (\slashed{g}^{-1})^{\rho\sigma} R_{\Lbar \rho L \sigma}
	+ \frac{1}{2} \slashed{\nabla}_\mu \zeta_\nu
	+ \frac{1}{2} \slashed{\nabla}_\nu \zeta_\mu
	- \frac{1}{4} \slashed{g}_{\mu\nu} \slashed{\Div} \zeta
	+ 2 \slashed{\nabla}^2_{\mu\nu} \log \mu
	\\
	&\phantom{=}
	- \slashed{g}_{\mu\nu} \slashed{\Delta} \log \mu
	+ \frac{1}{2} \zeta_\mu \zeta_\nu
	- \frac{1}{4} \slashed{g}_{\mu\nu} |\zeta|^2
	+ \zeta_\mu \slashed{\nabla}_\nu \log \mu
	+ \zeta_\nu \slashed{\nabla}_\mu \log \mu
	- \frac{1}{2} \slashed{g}_{\mu\nu} \left(\zeta \cdot \slashed{\nabla} \log \mu\right)
	\\
	&\phantom{=}
	+ 2(\slashed{\nabla}_\mu \log \mu)(\slashed{\nabla}_\nu \log \mu)
	- \slashed{g}_{\mu\nu} |\slashed{\nabla}\log \mu|^2
	+ \omega \hat{\chi}_{\mu\nu}
	- \frac{1}{2} \hat{\chi}_\mu^{\phantom{\mu} \sigma} \hat{\chibar}_{\nu\sigma}
	- \frac{1}{2} \hat{\chi}_\nu^{\phantom{\nu} \sigma} \hat{\chibar}_{\mu\sigma}
	- r^{-1} \hat{\chibar}_{\mu\nu}
	+ r^{-1} \hat{\chi}_{\mu\nu}
	\\
	&\phantom{=}
	- \frac{1}{2} \tr_{\slashed{g}}\chi_{(\text{small})} \hat{\chibar}_{\mu\nu}
	- \frac{1}{2} \tr_{\slashed{g}}\chibar_{(\text{small})} \hat{\chi}_{\mu\nu}
	+ \frac{1}{4} \slashed{g}_{\mu\nu} \hat{\chi}^{\rho\sigma} \hat{\chibar}_{\rho\sigma}
	\end{split}
	\end{equation*}
	so, schematically, we have
	\begin{equation*}
	\begin{split}
	\slashed{\D}_{\Lbar} \hat{\chi}
	=
	r^{-1} (\slashed{\D} \mathscr{Y} h)_{(\text{frame})}
	+ \slashed{\nabla}^2 \log \mu
	+ r^{-1} \bm{\Gamma}^{(0)}_{(-1+C_{(0)}\epsilon)}
	+ \bm{\Gamma}^{(1)}_{(-1+C_{(1)}\epsilon)} \bm{\Gamma}^{(1)}_{(-1+C_{(1)}\epsilon)}
	\end{split}
	\end{equation*}
	Now, commuting with $\mathscr{Y}^{n-1}$ and using proposition \ref{proposition commute D Yn} we find that, schematically,
	\begin{equation*}
	\begin{split}
	\slashed{\D}_{\Lbar} \mathscr{Y}^{n-1} \hat{\chi}
	&=
	[\slashed{\D} , \mathscr{Y}^{n-1}] \hat{\chi}_{\mu\nu}
	+ \sum_{j \leq n} r^{-1} (\slashed{\D} \mathscr{Y}^j h)_{(\text{frame})}
	\\
	&\phantom{=}
	+ \sum_{j+k+\ell \leq n-1} r^{-1} ([\slashed{\D} , \mathscr{Y}^j] \mathscr{Y} h_{(\text{rect})}) (\mathscr{Y}^k X_{(\text{frame})}) (\mathscr{Y}^\ell X_{(\text{frame})})
	+ \slashed{\nabla}^2 \mathscr{Z}^{n-1} \log \mu
	\\
	&\phantom{=}
	+ [\slashed{\nabla}^2 , \mathscr{Y}^{n-1}] \log \mu
	+ r^{-1} \bm{\Gamma}^{(n-1)}_{(-1+C_{(n-1)}\epsilon)}
	+ \sum_{j+k \leq n-1} \bm{\Gamma}^{(j+1)}_{(-1+C_{(j+1)}\epsilon)} \bm{\Gamma}^{(k+1)}_{(-1+C_{(k+1)}\epsilon)}
	\\
	\\
	&=
	\bm{\Gamma}^{(1)}_{(-1+C_{(1)}\epsilon)} \bm{\Gamma}^{(n)}_{(-1+C_{(n)}\epsilon)}
	+ \bm{\Gamma}^{(n-1)}_{(-1-\delta)}
	+ r^{-1} (\slashed{\D} \mathscr{Y}^{n} h)_{(\text{frame})}
	\\
	&\phantom{=}
	+ \bm{\Gamma}^{(1)}_{(-2+C_{(1)}\epsilon)} (\mathscr{Y}^{n-1} \bar{X}_{(\text{frame})})
	+ \slashed{\nabla}^2 \mathscr{Z}^{n-1} \log \mu
	\end{split}
	\end{equation*}

	Combining this with proposition \ref{proposition transport Yn chihat} we find that
	\begin{equation*}
	\begin{split}
	\slashed{\D}_{T} \mathscr{Y}^{n-1} \hat{\chi}
	&=
	[\slashed{\D} , \mathscr{Y}^{n-1}] \hat{\chi}
	+ \sum_{j \leq n} r^{-1} (\slashed{\D} \mathscr{Y}^j h)_{(\text{frame})}
	\\
	&\phantom{=}
	+ \sum_{j+k+\ell \leq n-1} r^{-1} ([\slashed{\D} , \mathscr{Y}^j] \mathscr{Y} h_{(\text{rect})}) (\mathscr{Y}^k X_{(\text{frame})}) (\mathscr{Y}^\ell X_{(\text{frame})})
	+ \slashed{\nabla}^2 \mathscr{Z}^{n-1} \log \mu
	\\
	&\phantom{=}
	+ [\slashed{\nabla}^2 , \mathscr{Y}^{n-1}] \log \mu
	+ r^{-1} \bm{\Gamma}^{(n-1)}_{(-1+C_{(n-1)}\epsilon)}
	+ \sum_{j+k \leq n-1} \bm{\Gamma}^{(j+1)}_{(-1+C_{(j+1)}\epsilon)} \bm{\Gamma}^{(k+1)}_{(-1+C_{(k+1)}\epsilon)}
	\\
	\\
	&=
	\bm{\Gamma}^{(1)}_{(-1+C_{(1)}\epsilon)} \bm{\Gamma}^{(n)}_{(-1+C_{(n)}\epsilon)}
	+ \bm{\Gamma}^{(n-1)}_{(-1-\delta)}
	+ r^{-1} (\slashed{\D} \mathscr{Y}^{n} h)_{(\text{frame})}
	\\
	&\phantom{=}
	+ \bm{\Gamma}^{(1)}_{(-2+C_{(1)}\epsilon)} (\mathscr{Y}^{n-1} \bar{X}_{(\text{frame})})
	+ \slashed{\nabla}^2 \mathscr{Z}^{n-1} \log \mu
	\end{split}
	\end{equation*}
\end{proof}

\begin{proposition}[An expression for $\mathscr{Y}^n \tr_{\slashed{g}} \chibar$ in terms of other quantities]
	\label{proposition Yn tr chibar}
 The null frame connection component $\mathscr{Y}^n \tr_{\slashed{g}}\chibar_{(\text{small})}$ are related to various other quantities via equations with the following schematic form:
\begin{equation*}
\begin{split}
\mathscr{Y}^n \tr_{\slashed{g}} \chibar_{(\text{small})}
&=
\mathscr{Y}^n \tr_{\slashed{g}}\chi_{(\text{small})}
+ (\slashed{\D} \mathscr{Y}^n h)_{(\text{frame})}
+ \bm{\Gamma}^{(0)}_{(-1+C_{(0)}\epsilon)} (\mathscr{Y}^n X_{(\text{frame})})
+ \bm{\Gamma}^{(1)}_{(C_{(1)}\epsilon)} (\overline{\slashed{\D}} \mathscr{Y}^n h)_{(\text{frame})}
\\
&\phantom{=}
+ \bm{\Gamma}^{(0)}_{(-\delta)} \bm{\Gamma}^{(n)}_{(-1+C_{(n)}\epsilon)}
+ \bm{\Gamma}^{(n-1)}_{(-1-\delta)}
\end{split}
\end{equation*}

\end{proposition}

\begin{proof}
 Recall proposition \ref{proposition chibar in terms of chi}, which we can take the trace of and then write schematically as
\begin{equation*}
  \tr_{\slashed{g}} \chibar_{(\text{small})} 
  =
  \tr_{\slashed{g}} \chi_{(\text{small})}
  + r^{-1} \left( 2 - \sum_{i=1}^3 \slashed{\Pi}_\mu^{\phantom{\mu}i} \slashed{\Pi}^{\mu i} \right) + (\partial h)_{(\text{frame})} 
\end{equation*}
and we find that we can write
\begin{equation*}
 \begin{split}
  2 - \sum_{i=1}^3 \slashed{\Pi}_\mu^{\phantom{\mu} i} \slashed{\Pi}^{\mu i}
  &=
  2 - \sum_{i=1}^3 \left( \delta_\mu^i + \frac{1}{2} L_\mu \Lbar^i + \frac{1}{2} \Lbar_\mu L^i \right) \left( (g^{-1})^{\mu i} + \frac{1}{2} L^\mu \Lbar^i + \frac{1}{2} \Lbar^\mu L^i \right) \\
  &=
  2 - \sum_{i=1}^3 \left( (g^{-1})^{ii} + L^i \Lbar^i \right) \\
  &=
  \sum_{i=1}^3 \left( -H^{ii} + \frac{x^i}{r} L^i_{(\text{small})} - \frac{x^i}{r} \Lbar^i_{(\text{small})} - L^i_{(\text{small})}\Lbar^i_{(\text{small})} \right)
 \end{split}
\end{equation*}
where we recall the definition of the tensor field $H$, which we have so far not made much use of:
\begin{equation*}
 (g^{-1})^{\mu\nu} := (m^{-1})^{\mu\nu} + H^{\mu\nu}
\end{equation*}
where $m$ is the Minkowski metric, with values relative to the rectangular indices
\begin{equation*}
 \begin{split}
  m_{ij} &= \delta_{ij} \\
  m_{00} &= -1 \\
  m_{0i} = m_{i0} &= 0
 \end{split}
\end{equation*}
Hence, we have $H^{ij} = -h_{ij} + \mathcal{O}\left( (h_{(\text{rect})})^2 \right)$. Putting this together, we find, schematically
\begin{equation*}
 \begin{split}
  \tr_{\slashed{g}}\chibar_{(\text{small})}
  &=
  \tr_{\slashed{g}} \chi_{(\text{small})}
  + (\partial h)_{(\text{frame})} \\
  &\phantom{=}
  + r^{-1} \left( h_{(\text{rect})} + \frac{x^i}{r} X_{(\text{frame, small})} + X_{(\text{frame, small})} \bar{X}_{(\text{frame})} + \mathcal{O}\left( (h_{(\text{rect})})^2 \right) \right)
  \\
  &=
  \tr_{\slashed{g}}\chi_{(\text{small})}
  + (\partial h)_{(\text{frame})}
  + r^{-2} x^i X_{(\text{frame, small})}
  + \bm{\Gamma}^{(0)}_{(-1, \text{large})} \bm{\Gamma}^{(0)}_{(-\delta)}
  + \bm{\Gamma}^{(0)}_{(-1-\delta)}
 \end{split}
\end{equation*}

Applying the operators $\mathscr{Y}^n$, we find (again, schematically)
\begin{equation*}
	\begin{split}
	\mathscr{Y}^n \tr_{\slashed{g}} \chibar_{(\text{small})}
	&=
	\mathscr{Y}^n \tr_{\slashed{g}}\chi_{(\text{small})}
	+ (\slashed{\D} \mathscr{Y}^n h)_{(\text{frame})}
	+ ([\slashed{\D} , \mathscr{Y}^n] h)_{(\text{frame})} \\
	&\phantom{=}
	+ \bm{\Gamma}^{(0)}_{(-1+C_{(0)}\epsilon)} (\mathscr{Y}^n X_{(\text{frame})})
	+ r^{-2} x^i (\mathscr{Y}^n X_{(\text{frame, small})})
	+ \bm{\Gamma}^{(n)}_{(-1-\delta)}
	\end{split}
\end{equation*}

Now, using proposition \ref{proposition commute D Yn} proves the proposition.

\end{proof}

\begin{proposition}[An expression for $\mathscr{Y}^n\hat{\chibar}$]
	\label{proposition Yn chibar hat}
The tensor fields $\mathscr{Y}^n \hat{\chibar}$ can be expressed, schematically, in terms of other quantities as
\begin{equation*}
\begin{split}
\mathscr{Y}^n \hat{\chibar}
&=
\mathscr{Y}^n \hat{\chi}
+ (\slashed{\D} \mathscr{Y}^n h)_{(\text{frame})}
+ \bm{\Gamma}^{(0)}_{(-1 + C_{(0)}\epsilon)} (\mathscr{Y}^n X_{(\text{frame})})
+ \bm{\Gamma}^{(0)}_{(-1, \text{large})} (\mathscr{Y}^n h)_{(\text{frame})}
+ \bm{\Gamma}^{(0)}_{(-\delta)} \bm{\Gamma}^{(n)}_{(-1+C_{(n)}\epsilon)}
\\
&\phantom{=}
+ \bm{\Gamma}^{(n-1)}_{(-1 + 2C_{(n-1)}\epsilon)}
\end{split}
\end{equation*}

\end{proposition}

\begin{proof}
Recall proposition \ref{proposition chibar in terms of chi}. Taking the trace-free part of the expression given in that proposition, we find that $\hat{\chibar}$ obeys an equation of the form
\begin{equation*}
 \hat{\chibar} = \hat{\chi} + (\partial h)_{(\text{frame})} + r^{-1} \sum_{i=1}^3 \left( \slashed{\Pi}^i \otimes \slashed{\Pi}^i - \frac{1}{2} \slashed{g} (\slashed{g}^{-1})^{ii} \right)
\end{equation*}
To estimate this last term, we first denote by $T_{(\text{rect})}$ the vector field with values relative to the rectangular coordinate system
\begin{equation*}
 \begin{split}
  T_{(\text{rect})}^0 &= 1 \\
  T_{(\text{rect})}^i &= 0
 \end{split}
\end{equation*}
In other words, relative to the rectangular coordinates system $x^a$, we have
\begin{equation*}
 T_{(\text{rect})} = \partial_{x^0} = \partial_t
\end{equation*}

Then, we can write
\begin{equation*}
 \begin{split}
  \sum_{i=1}^{3} \slashed{\Pi}_\mu^{\phantom{\mu}i} \slashed{\Pi}_\nu^{\phantom{\nu}i} - \frac{1}{2} \slashed{g}_{\mu\nu} (\slashed{g}^{-1})^{ii}
  &=
  \left( \slashed{g}_{\mu \alpha} \slashed{g}_{\nu \beta} - \frac{1}{2}\slashed{g}_{\mu\nu} \slashed{g}_{\alpha\beta} \right) \left( (m^{-1})^{\alpha\beta} + T_{(\text{rect})}^\alpha T_{(\text{rect})}^\beta \right) \\
  &=
  \left( \slashed{g}_{\mu \alpha} \slashed{g}_{\nu \beta} - \frac{1}{2}\slashed{g}_{\mu\nu} \slashed{g}_{\alpha\beta} \right) \left( (g^{-1})^{\alpha\beta} - H^{\alpha\beta} + T_{(\text{rect})}^\alpha T_{(\text{rect})}^\beta \right) \\
  &=
  \left( \slashed{g}_{\mu \alpha} \slashed{g}_{\nu \beta} - \frac{1}{2}\slashed{g}_{\mu\nu} \slashed{g}_{\alpha\beta} \right) \left( (\slashed{g}^{-1})^{\alpha\beta} - H^{\alpha\beta} + T_{(\text{rect})}^\alpha T_{(\text{rect})}^\beta \right) \\
  &=
  -\hat{H}_{\mu\nu}
  + \left( \slashed{g}_{\mu \alpha} \slashed{g}_{\nu \beta} - \frac{1}{2}\slashed{g}_{\mu\nu} \slashed{g}_{\alpha\beta} \right) T^\alpha _{(\text{rect})} T^\beta_{(\text{rect})}
 \end{split}
\end{equation*}
Now, we can decompose the tensor fields $\slashed{g}$ relative to the rectangular one-forms $\upd x^a$. That is, we write 
\begin{equation*}
 \slashed{g}_{\mu\nu} = \slashed{g}_{ab} (\upd x^a)_\mu (\upd x^b)_\nu 
\end{equation*}
and we also recall that
\begin{equation*}
 \begin{split}
  \slashed{g}_{\mu\nu} &= g_{\mu\nu} + \frac{1}{2}L_\mu \Lbar_\nu + \frac{1}{2}\Lbar_\mu L_\nu \\
  &= m_{\mu\nu} + h_{\mu\nu} + \frac{1}{2}L_\mu \Lbar_\nu + \frac{1}{2}\Lbar_\mu L_\nu \\
 \end{split}
\end{equation*}
Hence, we have
\begin{equation*}
 \begin{split}
  \slashed{g}_{\mu \alpha} T^\alpha_{(\text{rect})}
  &=
  \left( g_{a0} + \frac{1}{2} g_{a0}(\Lbar_0 L^0 + L_0 \Lbar^0) + \frac{1}{2}g_{ai} (\Lbar_0 L^i + L_0 \Lbar^i) \right) (\upd x^a)_{\mu} 
  \\
  \\
  &= \bigg( g_{a0}\left(1 - \Lbar^0 L^0 + \frac{1}{2} h_{0b}(L^0 \Lbar^b + \Lbar^0 L^b) \right) - \frac{1}{2} g_{ai} (\Lbar^0 L^i + L^0 \Lbar^i) \\
  &\phantom{=\bigg(}
  + \frac{1}{2} g_{ai} h_{0a}(\Lbar^a L^i + L^a \Lbar^i) \bigg)(\upd x^a)_{\mu} 
  \\
  \\
  &= \bigg( g_{a0}\left(- \Lbar^0_{(\text{small})} - L^0_{(\text{small})} - \Lbar^0_{(\text{small})}L^0_{(\text{small})} + \frac{1}{2} h_{0b}(L^0 \Lbar^b + \Lbar^0 L^b) \right) \\
  &\phantom{=\bigg(}
  - \frac{1}{2} g_{ai} \Big( r^{-1} x^i \Lbar^0_{(\text{small})} - r^{-1} x^i L^0_{(\text{small})} + L^i_{(\text{small})} + \Lbar^i_{(\text{small})} + \Lbar^0_{(\text{small})} L^i_{(\text{small})} 
  \\
  &\phantom{=\bigg( - \frac{1}{2} g_{ai} \Big(}
  + L^0_{(\text{small})}\Lbar^i_{(\text{small})} \Big)
  + \frac{1}{2}g_{ai} h_{0a} (\Lbar^a L^i + L^a \Lbar^i)
  \bigg)(\upd x^a)_{\mu}
 \end{split}
\end{equation*}
and so, schematically,
\begin{equation*}
 \begin{split}
  \slashed{g}_{\mu\alpha} T^\alpha_{(\text{rect})}
  &=
  g_{(\text{rect})}\bigg( X_{(\text{frame, small})} + (X_{(\text{frame, small})})^2 + h_{(\text{rect})} (X_{(\text{frame})})^2 \\
  &\phantom{= g_{(\text{rect})}\bigg(} 
  + r^{-1} x^i X_{(\text{frame, small})}\bigg) 
 \end{split}
\end{equation*}

Putting this all together, we have, schematically,
\begin{equation*}
 \begin{split}
  \hat{\chibar}
  &=
  \hat{\chi}
  + (\partial h)_{(\text{frame})}
  + r^{-1} h_{(\text{frame})}
  + r^{-1} (g_{(\text{rect})})^2 (X_{(\text{frame, small})})^2
  + r^{-1} (g_{(\text{rect})})^2 (X_{(\text{frame, small})})^4 
  \\
  &\phantom{=}
  + r^{-1} (g_{(\text{rect})})^2 (h_{(\text{rect})})^2 (X_{(\text{frame})})^4 
  + r^{-1} (g_{(\text{rect})})r^{-2} (x^i X_{(\text{frame, small})})^2
  + r^{-1} \mathscr{O}\left((h_{(\text{frame})})^2 \right)
  \\
  \\
  &= \hat{\chi}
  + (\partial h)_{(\text{frame})}
  + r^{-1} h_{(\text{frame})}
  + r^{-1} \bm{\Gamma}^{(0)}_{(0, \text{large})} (X_{(\text{frame, small})})^2
  + r^{-1} \bm{\Gamma}^{(0)}_{(0, \text{large})} h_{(\text{rect})}
 \end{split}
\end{equation*}

Now, applying the operators $\mathscr{Y}^n$ we find that, schematically,
\begin{equation*}
 \begin{split}
  \mathscr{Y}^n \hat{\chibar}
  &=
  \mathscr{Y}^n \hat{\chi}
  + (\slashed{\D} \mathscr{Y}^n h)_{(\text{frame})}
  + ([\slashed{\D} , \mathscr{Y}^n] h)_{(\text{frame})}
  + \bm{\Gamma}^{(0)}_{(-1 + C_{(0)}\epsilon)} (\mathscr{Y}^n X_{(\text{frame})})
  + \bm{\Gamma}^{(0)}_{(-1, \text{large})} (\mathscr{Z}^n h)_{(\text{frame})}
  \\
  &\phantom{=}
  + \bm{\Gamma}^{(n-1)}_{(-1 + 2C_{(n-1)}\epsilon)}
 \end{split}
\end{equation*}
Substituting for the third term using proposition \ref{proposition commute Zn D} proves the proposition.

\end{proof}

\begin{proposition}[Commuting the expression for $\slashed{\Delta} \mu$]
\label{proposition commute laplacian mu}
 Suppose that the rectangular components of $h$ satisfy equations of the form
\begin{equation*}
 \tilde{\Box}_g h_{ab} = F_{ab}
\end{equation*}
for some scalar fields $F_{ab}$. Then we have the following expression for the commuted version of $\slashed{\Delta} \log \mu$:
\begin{equation}
\begin{split}
\slashed{\Delta} \mathscr{Z}^{n-1} \log \mu
&=
\bm{\Gamma}^{(0)}_{(-1-\delta)} \mathscr{Y}^n \log \mu
+ (1+\bm{\Gamma}^{(0)}_{(C_{(0)}\epsilon)}) \mathscr{Y}^{n} \tr_{\slashed{g}} \chi_{(\text{small})}
+ r^{-1} \mathscr{Y}^n \zeta
+ r^{-1} (\slashed{\D} \mathscr{Y}^n h)_{(\text{frame})}
\\
&\phantom{=}
+ \sum_{j+k \leq n-1} (\mathscr{Y}^{j} F)_{(\text{rect})} \bm{\Gamma}^{(k)}_{(C_{(k)}\epsilon, \text{large})}
+ \bm{\Gamma}^{(n)}_{(-1+C_{(n)}\epsilon)} \bm{\Gamma}^{(0)}_{(-1+C_{(0)}\epsilon)}
\\
&\phantom{=}
+ \sum_{\substack{j+k \leq n \\ j,k \leq n-1}} \bm{\Gamma}^{(j)}_{(-1+C_{(j)}\epsilon)}\bm{\Gamma}^{(k)}_{(-1+C_{(k)}\epsilon)}
\end{split}
\end{equation}

\end{proposition}

\begin{proof}
Recall proposition \ref{proposition spherical laplacian mu} for the spherical laplacian of $\mu$, which we can write schematically as
\begin{equation*}
  \slashed{\Delta} \log \mu
  =
  T \tr_{\slashed{g}} \chi
  + (\slashed{g}^{-1})^{\mu\nu} R_{L \mu \Lbar \nu}
  + \tr_{\slashed{g}} \alpha
  + r^{-1} \mathscr{Z} \zeta
  + \bm{\Gamma}^{(0)}_{(-1+C_{(0)}\epsilon)} \cdot \bm{\Gamma}^{(0)}_{(-1+C_{(0)}\epsilon)}
\end{equation*}
and, using propositions \ref{proposition explicit expression for Riemann} and \ref{proposition second derivatives as commutators} we have
\begin{equation*}
 \begin{split}
  (\slashed{g}^{-1})^{\mu\nu} R_{L \mu \Lbar \nu}
  + \tr_{\slashed{g}} \alpha
  &=
  r^{-1} (\slashed{\D} \mathscr{Y} h)_{(\text{frame})} + (\tilde{\Box}_g h)_{(\text{frame})} + \begin{pmatrix} r^{-1} \\ \bm{\Gamma} \\ (\partial h)_{(\text{frame})} \end{pmatrix} (\partial h)_{(\text{frame})}
  \\
  \\
  &=
  r^{-1}(\slashed{\D} \mathscr{Y} h)_{(\text{frame})}
  + F_{(\text{frame})}
  + \bm{\Gamma}^{(1)}_{(-1+C_{(1)}\epsilon)} \cdot \bm{\Gamma}^{(0)}_{(-1+C_{(0)}\epsilon)}
 \end{split}
\end{equation*}
so, schematically we have
\begin{equation*}
 \begin{split}
  \slashed{\Delta} \log \mu
  &=
  T \tr_{\slashed{g}} \chi_{(\text{small})}
  + r^{-1} \mathscr{Z} \zeta
  + r^{-1} (\slashed{\D}\mathscr{Y} h)_{(\text{frame})}
  + F_{(\text{frame})}
  + \bm{\Gamma}^{(0)}_{(-1+C_{(0)}\epsilon)} \cdot \bm{\Gamma}^{(1)}_{(-1+C_{(1)}\epsilon)}
 \end{split}
\end{equation*}

This time, we only need to commute $(n-1)$ times with the operators $\mathscr{Z}$. We find
\begin{equation*}
 \begin{split}
  \slashed{\Delta} \mathscr{Z}^{n-1} \log \mu
  &=
  [\slashed{\Delta} , \mathscr{Z}^{n-1}] \log \mu
  + \slashed{\D}_T \mathscr{Y}^{n-1} \tr_{\slashed{g}} \chi_{(\text{small})}
  + [\slashed{\D}, \mathscr{Y}^{n-1}] \tr_{\slashed{g}} \chi_{(\text{small})}
  + r^{-1} \mathscr{Y}^n \zeta
  \\
  &\phantom{=}
  + r^{-1} (\slashed{\D} \mathscr{Y}^n h)_{(\text{frame})}
  + \sum_{j+k \leq n-1} (\mathscr{Y}^{j} F)_{(\text{rect})} \bm{\Gamma}^{(k)}_{(C_{(k)}\epsilon, \text{large})}
  + \bm{\Gamma}^{(n)}_{(-1+C_{(n)}\epsilon)} \bm{\Gamma}^{(0)}_{(-1+C_{(0)}\epsilon)}
  \\
  &\phantom{=}
  + \sum_{\substack{j + k \leq n \\ j,k \leq n-1}} \bm{\Gamma}^{(j)}_{(-1+C_{(j)}\epsilon)}\bm{\Gamma}^{(k)}_{(-1+C_{(k)}\epsilon)}
 \end{split}
\end{equation*}

Before we can make further progress, we need to compute $[\slashed{\Delta}, \mathscr{Z}^n] \log \mu$. We again remark that $\slashed{\D}$ is a metric connection on $\mathcal{B}$, with fibre metric $\slashed{\D}$. Thus we find
\begin{equation*}
\begin{split}
	[\slashed{\Delta} \, , \mathscr{Z}^{n-1}]
	&=
	(\slashed{g}^{-1})^{\mu\nu} [\slashed{\nabla}_\mu \slashed{\nabla}_\nu \, , \mathscr{Z}^{n-1}]
	\\
	&= (\slashed{g}^{-1})^{\mu\nu} \left( \slashed{\nabla}_\mu [\slashed{\nabla}_\nu \, , \mathscr{Z}^{n-1}] + [\slashed{\nabla}_\mu \, , \mathscr{Z}^{n-1}] \slashed{\nabla}_\nu \right)
\end{split}
\end{equation*}
so, schematically,
\begin{equation*}
	[\slashed{\Delta} \, , \mathscr{Z}^{n-1}] \log \mu
	= r^{-1} \mathscr{Z} [\slashed{\nabla}, \mathscr{Z}^{n-1}] \log \mu
	+ [\slashed{\nabla}, \mathscr{Z}^{n-1}] \left( r^{-1} \mathscr{Z} \log \mu \right)
\end{equation*}
and so, using propositions \ref{proposition commuting DT with first order operators} and \ref{proposition commuting rnabla with first order operators} we obtain
\begin{equation*}
[\slashed{\nabla} \, , \mathscr{Y}^{n-1}] \log \mu
=
\bm{\Gamma}^{(1)}_{(-1+C_{(1)}\epsilon)} \mathscr{Y}^{n-1} \log \mu
+ \bm{\Gamma}^{(n-2)}_{(-1+2C_{(n-2)}\epsilon)}
\end{equation*}
and hence
\begin{equation*}
[\slashed{\Delta} \, , \mathscr{Y}^{n-1}] \log \mu
=
r^{-1}\bm{\Gamma}^{(1)}_{(-1+C_{(1)}\epsilon)} \mathscr{Y}^{n} \log \mu
+ \bm{\Gamma}^{(n-1)}_{(-2+2C_{(n-2)}\epsilon)}
\end{equation*}

Putting this together with the expression for $[\slashed{\D}_T, \mathscr{Y}^{n-1}]$ given in proposition \ref{proposition commute D Yn}, we find

\begin{equation}
\begin{split}
\slashed{\Delta} \mathscr{Z}^{n-1} \log \mu
&=
r^{-1}\bm{\Gamma}^{(1)}_{(-1+C_{(1)}\epsilon)} \mathscr{Y}^n \log \mu
+ (1+\bm{\Gamma}^{(0)}_{(C_{(0)}\epsilon)}) \mathscr{Y}^{n} \tr_{\slashed{g}} \chi_{(\text{small})}
+ r^{-1} \mathscr{Y}^n \zeta
+ r^{-1} (\slashed{\D} \mathscr{Y}^n h)_{(\text{frame})}
\\
&\phantom{=}
+ \sum_{j+k \leq n-1} (\mathscr{Y}^{j} F)_{(\text{rect})} \bm{\Gamma}^{(k)}_{(C_{(k)}\epsilon, \text{large})}
+ \bm{\Gamma}^{(n)}_{(-1+C_{(n)}\epsilon)} \bm{\Gamma}^{(0)}_{(-1+C_{(0)}\epsilon)}
\\
&\phantom{=}
+ \sum_{\substack{j+k \leq n \\ j,k \leq n-1}} \bm{\Gamma}^{(j)}_{(-1+C_{(j)}\epsilon)}\bm{\Gamma}^{(k)}_{(-1+C_{(k)}\epsilon)}
\end{split}
\end{equation}

\end{proof}

\chapter{Elliptic estimates and Sobolev embedding}
\label{chapter elliptic estimates and sobolev embedding}

In this chapter we provide estimates for several different systems of the form $\mathscr{O}(\Phi) = F$, for some differential operator $\mathscr{O}$, and some inhomogeneous term $F$. The systems are elliptic in the following sense: if the operator $\mathscr{O}$ is an $n$-th order differential operator, then we obtain estimates for $n$ derivatives of $\Phi$ in terms of $F$.

Note that these estimates are only elliptic in the above sense \emph{assuming a lower bound on the Gauss curvature} $K$ of the spheres $S_{\tau,r}$. Thus in order to make use of these bounds, we will have to couple these equations to the equation for the evolution of the Gauss curvature given in proposition \ref{proposition transport Gauss curvature}.

\begin{proposition}[Elliptic estimates for solutions to Poisson's equation]
\label{proposition Poisson estimate}
Let $\Phi$ be a scalar field on the sphere $S_{\tau,r}$. Then $\Phi$ satisfies the following equation:
\begin{equation}
 |\slashed{\nabla}^2 \Phi|^2 + K|\slashed{\nabla}\Phi|^2 = (\slashed{\Delta}\Phi)^2 + \slashed{\nabla}_\mu \left( (\slashed{\nabla}^\mu \slashed{\nabla}^\nu \Phi)(\slashed{\nabla}_\nu \Phi) - (\slashed{\nabla}^\nu \Phi)(\slashed{\Delta}\Phi) \right)
\end{equation}
and so $\Phi$ satisfies the following estimate:
\begin{equation}
 \int_{S_{\tau,r}} \left(|\slashed{\nabla}^2 \Phi|^2 + K|\slashed{\nabla}\Phi|^2 \right) \dVol_{\slashed{g}} = \int_{S_{\tau,r}} (\slashed{\Delta}\Phi)^2 \dVol_{\slashed{g}}
\end{equation}

On the other hand, if $\Phi$ is a higher order $S_{\tau, r}$-tangent tensor field, then we have
\begin{equation}
	\int_{S_{\tau,r}} \left(|\slashed{\nabla}^2 \Phi|^2 \right) \dVol_{\slashed{g}} \lesssim \int_{S_{\tau,r}} \left(  (\slashed{\Delta}\Phi)^2 + K|\slashed{\nabla}\Phi|^2 + |\slashed{\nabla} K| |\slashed{\nabla}\Phi| |\Phi| \right) \dVol_{\slashed{g}}
\end{equation}

\end{proposition}

\begin{proof}
 We begin from the equation
\begin{equation*}
 |\slashed{\nabla}^2 \Phi|^2 = (\slashed{\nabla}_\mu \slashed{\nabla}_\nu \Phi ) \cdot (\slashed{\nabla}^\mu \slashed{\nabla}^\nu \Phi )
\end{equation*}
We now integrate by parts, and use the fact that, since the sphere $S_{\tau,r}$ is two-dimensional, the Riemann curvature of $S_{\tau,r}$ may be expressed in terms of the Gauss curvature as
\begin{equation}
 \slashed{R}_{\mu\nu\rho\sigma} = K\left( \slashed{g}_{\mu\rho}\slashed{g}_{\nu\sigma} - \slashed{g}_{\mu\sigma} \slashed{g}_{\nu\rho} \right)
\end{equation}
For the case of a scalar field, we need to keep track of the signs of the various terms, but for the general case we simply commute the derivatives and integrate by parts.

\end{proof}

\begin{proposition}[Elliptic estimates for div-curl systems of vector fields]
\label{proposition elliptic vector estimate}
 Let $\Phi$ be a vector field on the sphere $S_{\tau,r}$. Then $\Phi$ satisfies the following equation
\begin{equation}
 |\slashed{\nabla}\Phi|^2 + K|\Phi|^2 = (\slashed{\Curl}\, \Phi)^2 + (\slashed{\Div}\, \Phi)^2 + \slashed{\nabla}_\mu \left( (\slashed{\nabla}^\nu \Phi^\mu)\Phi_\nu - (\slashed{\Div}\, \Phi)\Phi^\mu \right)
\end{equation}
and so we have the following estimate:
\begin{equation}
 \int_{S_{\tau,r}} \left( |\slashed{\nabla}\Phi|^2 + K|\Phi|^2 \right) \dVol_{\slashed{g}} = \int_{S_{\tau,r}} \left( (\slashed{\Curl}\, \Phi)^2 + (\slashed{\Div}\, \Phi)^2 \right) \dVol_{\slashed{g}}
\end{equation}

On the other hand, if $\Phi$ is a higher order tensor field, and we define
\begin{equation*}
\begin{split}
	(\slashed{\Div}\Phi)_{\alpha_1 \ldots \alpha_n} &:= \slashed{\nabla}^\mu \Phi_{\mu \alpha_1 \ldots \alpha_n}
	\\
	\frac{1}{2}\slashed{\epsilon}_{\mu\nu} (\slashed{\Curl} \Phi)_{\alpha_1 \ldots \alpha_n} &:=
	\slashed{\nabla}_{[\mu} \Phi_{\nu]\alpha_1 \ldots \alpha_n}
\end{split}
\end{equation*}
then we have
\begin{equation}
	\int_{S_{\tau,r}} \left( |\slashed{\nabla}\Phi|^2 \right) \dVol_{\slashed{g}} \lesssim \int_{S_{\tau,r}} \left( (\slashed{\Curl}\, \Phi)^2 + (\slashed{\Div}\, \Phi)^2 + K |\Phi|^2 \right) \dVol_{\slashed{g}}
\end{equation}

\end{proposition}

\begin{proof}
 We begin by decomposing the derivatives of $\Phi$ into their symmetric and antisymmetric parts:
\begin{equation*}
 \begin{split}
  |\slashed{\nabla}\Phi|^2 &= \left(\slashed{\nabla}^\mu \Phi^\nu \right)\left( \frac{1}{2}\slashed{\nabla}_\mu \Phi_\nu + \frac{1}{2}\slashed{\nabla}_\nu \Phi_\mu + \slashed{\nabla}_{[\mu} \Phi_{\nu]} \right) \\
  \Rightarrow |\slashed{\nabla}\Phi|^2 &= \left( \slashed{\nabla}^\mu \Phi^\nu \right) \left( \slashed{\nabla}_\nu \Phi_\mu \right) + (\slashed{\Curl}\, \Phi)^2
 \end{split}
\end{equation*}
We now integrate by parts to deal with the first term on the right hand side in a similar way to proposition \ref{proposition Poisson estimate}.

The proof in the case of higher order fields is almost identical, except that some additional terms involving the Gauss curvature are produced when commuting derivatives.
\end{proof}

\begin{remark}
We will actually not make any use of the proposition above regarding div-curl systems, but we include it here for completeness and because, with a different choice of normalisation, we would have had to use this kind of estimate. These elliptic estimates are useful for avoiding derivative loss when estimating vector fields such as $\zeta$, however, due to our choice of normalisation of $\upd r$, we have been able to write $\zeta$ directly in terms of the derivatives of $h$. This would not have been the case if we had chosen, for example (following Christodoulou-Klainerman) $r$ to be an affine parameter along the integral curves of $L$, and in that case we would have to turn to elliptic estimates in order to avoid a loss of derivatives associated with terms of the form $\slashed{\nabla} \mathscr{Z}^{n} \zeta$. 
\end{remark}

\begin{proposition}[Elliptic estimates for Hodge systems of symmetric trace-free tensor fields]
	\label{proposition elliptic estimates Hodge}
 Now let $\Phi$ be a tensor field on the sphere $S_{\tau,r}$ which is \emph{symmetric} and \emph{trace-free}, i.e.\ $\Phi^{\mu\nu} = \Phi^{\nu\mu}$ and $\Phi^{\mu\nu} \slashed{g}_{\mu\nu} = 0$. Then $\Phi$ satisfies the following equation:
\begin{equation}
 |\slashed{\nabla}\Phi|^2 + 2K|\Phi|^2 = 2|\slashed{\Div}\, \Phi|^2 + \slashed{\nabla}_\mu \left( \Phi^{\nu\rho}\slashed{\nabla}_\nu \Phi^\mu_{\phantom{\mu}\rho} - \Phi^\mu_{\phantom{\mu}\rho} (\slashed{\Div}\, \Phi)^\rho \right)
\end{equation}
and so $\Phi$ satisfies the following estimate:
\begin{equation}
 \int_{S_{\tau,r}} \left( |\slashed{\nabla}\Phi|^2 + 2K|\Phi|^2 \right) \dVol_{\slashed{g}} = \int_{S_{\tau,r}} \left( 2|\slashed{\Div}\, \Phi|^2 \right) \dVol_{\slashed{g}}
\end{equation}

On the other hand, if $\Phi$ is a higher order $S_{\tau,r}$-tangent tensor field which is symmteric and trace-free in its first two indices, then, defining the divergence in the obvious way, we have
\begin{equation}
\int_{S_{\tau,r}} \left( |\slashed{\nabla}\Phi|^2 \right) \dVol_{\slashed{g}} \lesssim \int_{S_{\tau,r}} \left( |\slashed{\Div}\, \Phi|^2 + K|\Phi|^2 \right) \dVol_{\slashed{g}}
\end{equation}

\end{proposition}

\begin{proof}
 We proceed in a similar way to the proof of proposition \ref{proposition elliptic vector estimate}. This time, additional terms arise when commuting the covariant derivatives, and so it is important to make use of the fact that $\Phi$ is trace-free and symmetric. In this way it is fairly easy to derive the equation
\begin{equation*}
 |\slashed{\nabla}\Phi|^2 + 2K|\Phi|^2 = |\slashed{\Div}\, \Phi|^2 + |\slashed{\Curl}\, \Phi|^2 + \slashed{\nabla}_\mu \left( \Phi^{\nu\rho}\slashed{\nabla}_\nu \Phi^\mu_{\phantom{\mu}\rho} - \Phi^\mu_{\phantom{\mu}\rho} (\slashed{\Div}\, \Phi)^\rho \right)
\end{equation*}
where we recall that
\begin{equation*}
 (\slashed{\Curl}\, \Phi)_\mu := \slashed{\epsilon}^{\nu\rho}\slashed{\nabla}_\nu \Phi_{\rho\mu}
\end{equation*}
We claim that
\begin{equation*}
 (\slashed{\Curl}\, \Phi)_\mu = {^*(\slashed{\Div}\, \Phi)}_\mu
\end{equation*}
where $*$ denotes the Hodge star; i.e.
\begin{equation*}
 {^*(\slashed{\Div}\, \Phi)}_\mu = \slashed{\epsilon}_\mu^{\phantom{\mu}\nu}(\slashed{\Div}\, \Phi)_\nu
\end{equation*}
To prove this, note that
\begin{equation*}
 \slashed{\epsilon}_\mu^{\phantom{\mu}\nu}\slashed{\nabla}^\rho \Phi_{\nu\rho} = \slashed{\nabla}^\rho \left(\epsilon_\mu^{\phantom{\mu}\nu}\Phi_{\nu\rho}\right)
\end{equation*}
Now, we decompose the tensor field $\slashed{\epsilon}_\mu^{\phantom{\mu}\nu}\phi_{\nu\rho}$ into its irreducible parts:
\begin{equation*}
 \slashed{\epsilon}_\mu^{\phantom{\mu}\nu}\phi_{\nu\rho} = (\Phi_{(1)})_{\mu\rho} + (\Phi_{(2)})\slashed{g}_{\mu\rho} + (\Phi_{(3)})\slashed{\epsilon}_{\mu\rho}
\end{equation*}
where $\Phi_{(1)}$ is symmetric and trace-free, and $\Phi_{(2)}$ and $\Phi_{(3)}$ are given by
\begin{equation*}
 \begin{split}
  \Phi_{(2)} &:= \slashed{\epsilon}^{\mu\nu}\phi_{\mu\nu} \\
  \Phi_{(3)} &:= \frac{1}{2}\slashed{\epsilon}^{\mu\rho} \left( \slashed{\epsilon}_\mu^{\phantom{\mu}\nu}\phi_{\nu\rho} \right) \\
  &= \frac{1}{2} \slashed{g}^{\mu\nu}\phi_{\mu\nu}
 \end{split}
\end{equation*}
These expressions are valid for any tensor field $\Phi$. However, in our case $\Phi$ is symmetric and trace-free, so $\Phi_{(2)} = \Phi_{(3)} = 0$. Hence
\begin{equation*}
 \epsilon_\mu^{\phantom{\mu}\nu}\phi_{\nu\rho} = \epsilon_\rho^{\phantom{\rho}\nu}\phi_{\mu\rho}
\end{equation*}
and so we find that
\begin{equation*}
 \begin{split}
  {^*(\slashed{\Div}\, \Phi)}_\mu &= \slashed{\nabla}^\rho \left(\slashed{\epsilon}_\mu^{\phantom{\mu}\nu}\Phi_{\nu\rho}\right) \\
  &= \slashed{\nabla}^\rho \left(\slashed{\epsilon}_\rho^{\phantom{\rho}\nu}\Phi_{\nu\mu}\right) \\
  &= \slashed{\epsilon}^{\rho\nu} \slashed{\nabla}_\rho \Phi_{\nu\mu} \\
  &= (\slashed{\Curl}\, \Phi)_\mu
 \end{split}
\end{equation*}
Moreover, it is easy to show that
\begin{equation*}
 |{^*(\slashed{\Div}\, \Phi)}| = |(\slashed{\Div}\, \Phi)|
\end{equation*}
and so
\begin{equation*}
 |(\slashed{\Curl}\, \Phi)| = |(\slashed{\Div}\, \Phi)|
\end{equation*}

As usual, the calculations are similar in the case of a higher order field, with some additional terms involving the Gauss curvature produced when commuting derivatives.

\end{proof}

\begin{proposition}[Sobolev embedding]
	\label{proposition Sobolev}
 Let $\phi$ be an $S_{\tau,r}$-tangent tensor field, and suppose that the Gauss curvature of the sphere $S_{\tau,r}$ satisfies the bound
 \begin{equation*}
 |K - r^{-2}| \lesssim \epsilon r^{-2-\delta} 
 \end{equation*}
 Suppose also that the rectangular components of the metric satisfy
 \begin{equation*}
 \begin{split}
	|\slashed{g}_{ab}| &\lesssim 1 \\
	|(\slashed{g}^{-1})^{ab}| &\lesssim 1 \\
	|\slashed{\Gamma}^a_{bc}| &\lesssim r^{-1} \\
	|\slashed{\nabla} \slashed{\Gamma}^a_{bc}| &\lesssim r^{-2} \\
 \end{split}
 \end{equation*}
 where $\slashed{\Gamma}^a_{bc}$ are the rectangular components of the Christoffel symbols associated with the covariant derivative operator $\slashed{\nabla}$. In other words, we have
 \begin{equation*}
 \slashed{\nabla}_a \partial_b = \slashed{\Gamma}_{ab}^c \partial_c
 \end{equation*}
 
 Then
 \begin{equation}
	r^2 \left(||f||_{L^{\infty}(S_{\tau,r})}\right)^2 \lesssim \int_{S_{\tau,r}} \left( |f|^2 + r^2|\slashed{\nabla} f|^2 + r^4|\slashed{\nabla}^2 f|^2 \right) \dVol_{\slashed{g}}
 \end{equation}
\end{proposition}

\begin{proof}
	Define the rescaled metric on the sphere $S_{\tau,r}$
	\begin{equation*}
	r^2\mathring{\slashed{g}} := \slashed{g}
	\end{equation*}
	For a tensor field $\phi_{\mu_1 \ldots \mu_n}$ we write
	\begin{equation*}
	|\phi|_{\mathring{\slashed{g}}} := \phi_{\mu_1 \ldots \mu_n} \phi_{\nu_1 \ldots \nu_n} (\mathring{\slashed{g}}^{-1})^{\mu_1 \nu_1} \ldots (\mathring{\slashed{g}}^{-1})^{\mu_n \nu_n}
	\end{equation*}
	
	First, we shall deal with the case of a scalar field $f$. Then we have
	\begin{equation*}
	\int_{S_{\tau,r}} \left( |f|^2 + r^2|\slashed{\nabla} f|^2 + r^4|\slashed{\nabla}^2 f|^2 \right) \dVol_{\slashed{g}}
	=\int_{S_{\tau,r}} \left( |f|^2 + |\slashed{\nabla} f|_{\mathring{\slashed{g}}}^2 + |\slashed{\nabla}^2 f|_{\mathring{\slashed{g}}}^2 \right) r^2 \dVol_{\mathring{\slashed{g}}}
	\end{equation*}
	where we note that the Levi-Civita connections of $\slashed{g}$ and $\mathring{\slashed{g}}$ agree, since the metrics themselves differ only by a conformal factor $r^2$ which is \emph{constant} on the sphere $S_{\tau,r}$.
	
	Now, we can compute the Gauss cuvature of $\mathring{\slashed{g}}$, denoted by $\mathring{K}$. We find that
	\begin{equation*}
	\mathring{K} = r^2 K
	\end{equation*}
	So, using the hypothesis of the proposition, we have
	\begin{equation*}
	|\mathring{K} - 1| \lesssim \epsilon r^{-\delta}
	\end{equation*}
	
	In particular, $\mathring{K}$ is always positive, and is bounded above and below. Indeed, for sufficiently small $\epsilon$, we have
	\begin{equation*}
	\frac{1}{2} \leq \mathring{K} \leq \frac{3}{2}
	\end{equation*}
	Under these conditions on the Gauss curvature, it can be shown (see e.g.\ \cite{Aubin1998}) that an inequality of the following form holds, with an implicit constant which is uniformly bounded above and below: for any scalar field $f$,
	\begin{equation}
	\left(||f||_{L^{\infty}(S_{\tau,r})}\right)^2 \lesssim \int_{S_{\tau,r}} \left( |f|^2 + |\slashed{\nabla} f|_{\mathring{\slashed{g}}}^2 + |\slashed{\nabla}^2 f|_{\mathring{\slashed{g}}}^2 \right) \dVol_{\mathring{\slashed{g}}}
	\end{equation}
	
	Now, for a higher order tensor field $f_{\mu_1 \ldots \mu_n}$ we have
	\begin{equation*}
	|f|_{\mathring{\slashed{g}}} = r^n |f|
	\end{equation*}
	and also, using the bounds on $\mathring{\slashed{g}}_{ab}$ and $(\mathring{\slashed{g}}^{-1})^{ab}$, we can show
	\begin{equation*}
	\begin{split}
	|f|_{\mathring{\slashed{g}}}
	& \sim r^n \sum_{a_1, \ldots, a_n} |f_{a_1 \ldots a_n}|
	\\
	|\slashed{\nabla} f|_{\mathring{\slashed{g}}} 
	& \sim r^{n+1} \left( \sum_{a_1, \ldots, a_{n}} |\slashed{\nabla} f_{a_1 \ldots a_n}| 
	+ \sum_{a,b,c}|\Gamma^a_{bc}| |f|_{\mathring{\slashed{g}}} \right)
	\\
	|\slashed{\nabla}^2 f|_{\mathring{\slashed{g}}}
	& \sim r^{n+2} \left( \sum_{a_1, \ldots, a_{n}} |\slashed{\nabla}^2 f_{a_1 \ldots a_n}| 
	+ \sum_{a,b,c}|\slashed{\Gamma}^a_{bc}| |\slashed{\nabla}f|_{\mathring{\slashed{g}}}
	+ \sum_{a,b,c}|\slashed{\nabla} \slashed{\Gamma}^a_{bc}|_{\mathring{\slashed{g}}} |f|_{\mathring{\slashed{g}}} \right)
	\end{split}
	\end{equation*}
	where we recall that $\slashed{\Gamma}^a_{bc}$ are the rectangular components of the Christoffel symbols associated with the metric $\slashed{g}$. Hence, using the Sobolev inequality for scalar fields given above, we have
	\begin{equation*}
	\left( || f_{a_1 \ldots a_n} ||_{L^\infty(S_{\tau,r})} \right)^2
	\lesssim
	\int_{S_{\tau,r}} \left( |f_{a_1 \ldots a_n}|^2 + |\slashed{\nabla} f_{a_1 \ldots a_n}|_{\mathring{\slashed{g}}}^2 + |f_{a_1 \ldots a_n}|_{\mathring{\slashed{g}}}^2 \right) \dVol_{\mathring{\slashed{g}}}	
	\end{equation*}
	Summing over the rectangular indices, and using the expressions above and the bounds for the Christoffel symbols and their derivatives, we find
	\begin{equation*}
	\begin{split}
	||f||_{L^\infty(S_{\tau,r})}
	&\lesssim
	\int_{S_{\tau,r}} r^{-n} \left(
	|f|^2_{\mathring{\slashed{g}}}
	+ |\slashed{\nabla}f|^2_{\mathring{\slashed{g}}}
	+ |\slashed{\nabla}^2 f|^2_{\mathring{\slashed{g}}}
	\right)\dVol_{\mathring{\slashed{g}}}
	\\
	&\lesssim
	\int_{S_{\tau,r}} \left(
	|f|^2
	+ r^{2}|\slashed{\nabla}f|^2
	+ r^{4}|\slashed{\nabla}^2 f|^2
	\right)r{-2}\dVol_{\slashed{g}}
	\end{split}	
	\end{equation*}
	Multiplying by $r^{2}$ proves the proposition in the case of a tensor field.
\end{proof}

\begin{remark}[An alternative proof of the Sobolev embedding]
	We note here that we could give an alternative proof of the proposition above, using the fact that we have a map $S_{\tau,r} \rightarrow \mathbb{S}^2$, i.e.\ a map from $S_{\tau,r}$ to the unit sphere, given by flowing along the integral curves of $L$. We also have a second metric $\bm{\gamma}$ on $S_{\tau,r}$, given by the pull-back of the standard round metric on the unit sphere by this mapping. Using this second metric, we can rewrite the derivatives $\slashed{\nabla}$ in terms of the \emph{standard} angular derivatives (i.e.\ using the Levi-Civita connection of $\gamma$), with an error term given by the difference between the two connections. This difference can be expressed in terms of a tensor field on the spheres, which can itself be expressed in terms of quantities of the form $\slashed{\nabla} \gamma$. These terms can then finally be controlled by commuting the equation $\slashed{\mathcal{L}} \gamma = 0$ with the operator $\slashed{\nabla}$, which is the covariant derivative with respect to $\slashed{g}$.
\end{remark}

\begin{proposition}[Derivatives of spherical integrals]
\label{proposition derivatives of spherical integrals}
 Let $\phi$ be any scalar field. Then we have the following two equations for the derivatives of the spherical integrals of $\phi$:
\begin{equation}
 \begin{split}
  \frac{\upd}{\upd r} \left( \int_{S_{\tau,r}} \phi \, \dVol_{\slashed{g}} \right) &= \int_{S_{\tau,r}} \Bigg( \left(1-\frac{1}{2}f'_{(\alpha)} \mu \right) \left(L\phi + (\tr_{\slashed{g}}\chi)\phi \right) - \frac{1}{2}f'_{(\alpha)} \mu \left(\Lbar\phi + (\tr_{\slashed{g}}\chibar)\phi \right) \\
  & \phantom{= \int_{S_{\tau,r}} \Bigg(}
  + \frac{1}{2}f'_{(\alpha)} \mu \left(b^A\slashed{\upd}_A \phi + \mu^{-1}\slashed{\nabla}_A (\mu b^A) \phi \right) \Bigg) \dVol_{\slashed{g}} \\ \\
  \frac{\upd}{\upd \tau} \left( \int_{S_{\tau,r}} \phi \, \dVol_{\slashed{g}} \right) &= \int_{S_{\tau,r}} \frac{1}{2}\mu \bigg( \left(L\phi + (\tr_{\slashed{g}}\chi)\phi \right) + \left(\Lbar\phi + (\tr_{\slashed{g}}\chibar)\phi \right) - \\
  & \phantom{= \int_{S_{\tau,r}} \frac{1}{2}\mu \bigg( }
  \left(b^A\slashed{\upd}_A \phi + \mu^{-1}\slashed{\nabla}_A (\mu b^A) \phi \right) \bigg) \dVol_{\slashed{g}} 
 \end{split}
\end{equation}

\end{proposition}

\begin{proof}
 In this proof, we shall use the notation $\partial_r = \partial_r\big|_{\tau,\vartheta_1,\vartheta_2}$, i.e.\ it is to be understood that the vector field $\partial_r$ is defined with respect to the $\{\tau,r,\vartheta^A\}$ coordinate system.

 Let $\varphi_{(\lambda)}$ be the diffeomorphism defined by flowing along the integral curves of the vector field $\partial_r$ by an affine distance $\lambda$, i.e.\ given a point $p \in \mathbb{R}^4$, we have $\frac{\upd}{\upd \lambda} \varphi_{(\lambda)}(p) = \partial_r \circ \varphi_{(\lambda)}(p)$, and $\varphi_{(0)}$ is the identity map on $\mathbb{R}^4$. Then 

\begin{equation*}
 \int_{S_{\tau, r + \lambda}} \phi \, \dVol_{\slashed{g}} = \int_{S_{\tau, r}} \varphi_{(\lambda)}^* (\phi) \dVol_{\varphi_{(\lambda)}^*(\slashed{g})}
\end{equation*}
 The standard relationship between the derivatives of the determinant of a matrix and its trace then allows us to see that
\begin{equation*}
 \frac{\upd}{\upd \lambda}\Big|_{\lambda = 0} \dVol_{\varphi_{(\lambda)}^*(\slashed{g})} = \frac{1}{2} (\slashed{g}^{-1})^{AB} \left(\slashed{\mathcal{L}}_{\partial_r} \slashed{g}\right)_{AB} = \frac{1}{2} \tr_{\slashed{g}} \left({^{(\partial_r)}\slashed{\pi}}\right)
\end{equation*}
Now, we use corollary \ref{corollary null frame coords} to find that
\begin{equation*}
 \partial_r = \left(1-\frac{1}{2}\mu f'_{(\alpha)} \right) L - \frac{1}{2} \mu f'_{(\alpha)} \Lbar + \frac{1}{2} \mu f' b^A X_A
\end{equation*}
which, together with proposition \ref{proposition deformation L Lbar} implies
\begin{equation*}
 {^{(\partial_r)}\slashed{\pi}}_{AB} = 2\left(1-\frac{1}{2}\mu f'_{(\alpha)} \right) \chi_{AB} - \mu f'_{(\alpha)} \chibar_{AB} + f' \slashed{\nabla}_{(A} \left( \mu b_{B)}\right)
\end{equation*}

In a similar way, we can prove the second equality given in the proposition, where we now make use of the relation
\begin{equation*}
 \frac{\partial}{\partial \tau}\Big|_{r,\vartheta_1,\vartheta_2} = \frac{1}{2}\mu \left( L + \Lbar - b^A X_A \right)
\end{equation*}

\end{proof}

\chapter{The framework for energy estimates}
\label{chapter framework for energy estimates}

In this chapter we present the abstract framework in which we perform energy estimates, which are of fundamental importance to our analysis, and perform some preliminary calculations related to the energy estimates. We first present several useful propositions which allow us to quantify the energies that we consider and to compare them with the standard $L^2$ norms. We will also present propositions which will allow us to deal with lower order terms, as well as propositions which allow us to estimate $L^2$ norms on the spheres $S_{\tau,r}$ in terms of the energy. Combining these with Sobolev embedding on the spheres will lead to pointwise bounds in terms of the energy. Finally, we present the fundamental energy flux identity in an abstract form. This will be used in subsequent chapters, along with specific choices for the associated vector fields and an $L^\infty$ bootstrap in order to prove various $L^2$ bounds.

Note that many of the propositions in this section make use of the various definitions of hypersurfaces given in definition \ref{definition hypersurfaces}.

\section{Preliminary calculations relating to the energy estimates}

Throughout this section we will need to use the expressions
\begin{equation}
\label{equation d dtau and ddr on sigma}
 \begin{split}
  \frac{\partial}{\partial \tau}\bigg|_{r,\vartheta^1,\vartheta^2} &= \frac{1}{2}\mu \left( L + \Lbar - b^A X_A \right) \\
  \frac{\partial}{\partial r}\bigg|_{\tau,\vartheta^1,\vartheta^2} &= L
 \end{split}
\end{equation}

\begin{proposition}[The volume form on $\Sigma_\tau$]
\label{proposition volume form gbar}
 In the region $r \geq r_0$, the hypersurface $\Sigma_\tau$ is null. We choose the volume form on the surfaces $\Sigma_{\tau}$ to be
\begin{equation}
  \dVol_{\underline{g}} =  \Omega^2 \upd r \wedge \dVol_{\mathbb{S}^2}
\end{equation}

% \begin{equation}
%  \dVol_{\underline{g}} =  \frac{\sqrt{2}}{2}\mu^{\frac{1}{2}}(f'_{(\alpha)})^{\frac{1}{2}}\sqrt{1-\frac{1}{2}\mu f'_{(\alpha)} }\sqrt{-\det g} \, \varepsilon_{abcd} (Z_{(\tau)})^a (Z_{(r)})^b (X_1)^c (X_2)^d \upd r \wedge \upd\vartheta^1 \wedge \upd\vartheta^2
% \end{equation}

\end{proposition}
% 
% \begin{proof}
%  We first note that 
% \begin{equation}
%  \begin{split}
%   (\upd \tau)_{\mu} &= - \mu^{-1} L_\mu + \frac{1}{2}f'_{(\alpha)} \left( L_\mu - \Lbar_\mu \right) \\
%   \Rightarrow (g^{-1})^{\mu\nu} (\upd \tau)_\mu (\upd \tau)_\nu &= - 2\mu^{-1} f'_{(\alpha)} + (f'_{(\alpha)})^2
%  \end{split}
% \end{equation}
% Together with proposition \ref{proposition volume form g} this proves the proposition.
% \end{proof}

\begin{proposition}[The induced volume form on the spheres $S_{\tau,r}$]
 \label{proposition volume form slashed g}
 In the region $r \geq r_0$, the induced volume form on the sphere $S_{\tau,r}$ is given by
 \begin{equation}
  \dVol_{\slashed{g}} = \dVol_{(\mathbb{S}^2, \slashed{g})} = \Omega^2 \dVol_{\mathbb{S}^2}
 \end{equation}
\end{proposition}

% \begin{proof}
%  We need to calculate the norm of the covector field $\upd r$ \emph{with respect to the metric} $\underline{g}$. We have
%  \begin{equation*}
%   (\underline{g}^{-1})^{\mu\nu} (\upd r)_\mu (\upd r)_\nu = (g^{-1})^{\mu\nu}\underline{\Pi}_\mu^{\phantom{\mu}\rho} \underline{\Pi}_\nu^{\phantom{\nu}\sigma} (\upd r)_\rho (\upd r)_\sigma
%   \end{equation*}
% The restriction of the one-form $\upd r$ to the surface $\Sigma_\tau$ is given by
% \begin{equation*}
%  \underline{\Pi}_\mu^{\phantom{\mu}\nu} (\upd r)_\nu = \left( -\frac{1}{4} + \frac{1}{2}\mu^{-1}(f'_{(\alpha)})^{-1} \right) L_\mu + \left( -\frac{1}{4} - \frac{1}{8}\mu f'_{(\alpha)} \right) \Lbar_{\mu}
% \end{equation*}
% and so we find
% \begin{equation*}
%  (\underline{g}^{-1})^{\mu\nu} (\upd r)_\mu (\upd r)_\nu = \frac{1}{2}\mu^{-1}(f'_{(\alpha)})^{-1} - \frac{1}{8}\mu f'_{(\alpha)}
% \end{equation*}
%  This, together with proposition \ref{proposition volume form gbar} proves the proposition.
%   
% \end{proof}

\begin{proposition}[The induced volume form on hypersurfaces of constant $t$]
 In the region $r \geq r_0$, the induced volume form on a hypersurface of constant $t$ is given by
\begin{equation}
  \begin{split}
    \dVol_{(\bar{\Sigma}_t)} &:= \left( \frac{2}{L^0 + \Lbar^0 - b^0} \right)^{-1} \sqrt{L^0 \Lbar^0 - (\slashed{g}^{-1})^{AB}(X_A)^0(X_B)^0} \Omega^2 \upd r \wedge \dVol_{\mathbb{S}^2}
  \end{split}
\end{equation}

 On the other hand, in the region $r \leq r_0$ the induced volume form on a hypersurface of constant $t$ is given by
\begin{equation}
  \begin{split}
    \dVol_{(\Sigma_t)} &= (g^{-1})^{00}\sqrt{\det g} \upd x^1 \wedge \upd x^2 \wedge \upd x^3
  \end{split} 
\end{equation}

\end{proposition}

\begin{proof}
 Note that
\begin{equation*}
 \left.\frac{\partial t}{\partial \tau}\right|_{r, \vartheta^1, \vartheta^2} = \frac{1}{2}\mu \left(L^0 + \Lbar^0 - b^A (X_A)^0 \right)
\end{equation*}
 and also that 
\begin{equation*}
 |\upd t| = \sqrt{L^0\Lbar^0 - (\slashed{g}^{-1})^{AB}(X_A)^0 (X_B)^0}
\end{equation*}

\end{proof}

\begin{proposition}[Estimating the spherical mean of a tensor field in terms of its energy]
\label{proposition spherical mean in terms of energy}
Let $\phi$ be an $S_{\tau,r}$-tangent tensor field. Let $t$ be such that the sphere $S_{\tau,r}$ is in the interior of the sphere $\bar{S}_{\tau, t}$, where both spheres are being considered as subsets of the hypersurface $\Sigma_\tau$. 

Let $\alpha < 1$, and let $r \geq r_0$.
Then we have
\begin{equation}
  \begin{split}
    \int_{S_{\tau,r}}|\phi|^2 \dVol_{\mathbb{S}^2} &\lesssim \frac{1}{(1-\alpha)}\frac{1}{r^{1-\alpha}}\int_{^t\Sigma_\tau} (1+r')^{-\alpha}\left( |\slashed{\D}_L \phi|^2 \right) (r')^2 \upd r' \wedge \dVol_{\mathbb{S}^2} \\
    &\phantom{\lesssim} + \int_{\bar{S}_{\tau,t}}|\phi|^2 \dVol_{\mathbb{S}^2}
  \end{split}
\end{equation}
where the implicit constant depends only on $\alpha$ and $r_0$.

Similarly we can estimate
\begin{equation}
  \begin{split}
    &\int_{\bar{S}_{\tau_1,t}}|\phi|^2 \dVol_{\mathbb{S}^2} \lesssim \int_{\bar{S}_{\tau_0,t}}|\phi|^2 \dVol_{\mathbb{S}^2} \\
    & + \frac{1}{(1-\alpha)}\sup_{x \in \bar{S}_{\tau_1,t}} \left\{\frac{1}{r(x)^{1-\alpha}}\right\} \int_{_{\tau_0}^{\tau_1}\bar{\Sigma}_t} (1+r')^{-\alpha}\left( |\slashed{\D}_{\Lbar} \phi|^2 + |\slashed{\D}_L \phi|^2 + \textit{Err}_{(t-\partial_r)}[\phi] \right) (r')^2 \upd r' \wedge \dVol_{\mathbb{S}^2}\\
  \end{split}
\end{equation}

Finally, by combining the above two inequalities we can estimate
\begin{equation}
 \begin{split}
    &\int_{S_{\tau,r}}|\phi|^2 \dVol_{\mathbb{S}^2} \lesssim \frac{1}{(1-\alpha)}\frac{1}{r^{1-\alpha}}\int_{^t\Sigma_\tau} (1+r')^{-\alpha}\left( |\slashed{\D}_L \phi|^2 \right) (r')^2 \upd r' \wedge \dVol_{\mathbb{S}^2} \\
     & + \frac{1}{(1-\alpha)}\sup_{x \in \bar{S}_{\tau_1,t}} \left\{\frac{1}{r(x)^{1-\alpha}}\right\}\int_{_{\tau_0}^{\tau_1}\bar{\Sigma}_t} (1+r')^{-\alpha}\left( |\slashed{\D}_{\Lbar} \phi|^2 + |\slashed{\D}_L \phi|^2 + \textit{Err}_{(t-\partial_r)}[\phi] \right) (r')^2 \upd r' \wedge \dVol_{\mathbb{S}^2} \\
    & + \int_{\bar{S}_{\tau_0,t}}|\phi|^2 \dVol_{\mathbb{S}^2} 
 \end{split}
\end{equation}

The error term in the above inequalities satisfy: 
\begin{equation}
  \begin{split}
    &\left| \textit{Err}_{(t-\partial_r)}[\phi] \right| \lesssim \frac{|L^0_{(\text{small})}| + |\Lbar^0_{(\text{small})}| + |b^0|}{\left|L^0 + \Lbar^0 - b^0\right|} \left( |\slashed{\D}_L \phi|^2 + |\slashed{\D}_{\Lbar} \phi|^2 \right) + \frac{|L^0| |b^0|}{\left| L^0 + \Lbar^0 - b^0 \right|}|\slashed{\nabla}\phi|^2
  \end{split}
\end{equation}

\end{proposition}

\begin{proof}
 To prove the first inequality, let $R_t(\tau,\vartheta^1, \vartheta^2)$ be the value of $r$ such that, along the line $\tau(x) = \tau$, $\vartheta^1(x)=\vartheta^1$, $\vartheta^2(x) = \vartheta^2$, we have $r(x) = R_t(\tau,\vartheta^1, \vartheta^2)$ when $t(x) = t$. We abbreviate this as $R_t$. Then we have
 \begin{equation*}
  \begin{split}
   \int_{S_{\tau,r}} |\phi|^2 \dVol_{\mathbb{S}^2}
    &= \int_{\mathbb{S}^2} \left( \int_{r' = r}^{R_t} \frac{ \partial |\phi|}{\partial r'} (\tau, r', \vartheta^1, \vartheta^2) \upd r \right)^2 \dVol_{\mathbb{S}^2} 
    - \int_{\bar{S}_{\tau, t}} |\phi|^2 \dVol_{\mathbb{S}^2} \\
    &\phantom{=} + \int_{\mathbb{S}^2} 2|\phi(\tau, r, \vartheta^1, \vartheta^2)| |\phi(\tau, R_t, \vartheta^1, \vartheta^2)| \dVol_{\mathbb{S}^2} \\
    &\lesssim \int_{\mathbb{S}^2} \left( \int_{r' = r}^{R_t} \frac{ \partial |\phi|}{\partial r'} (\tau, r', \vartheta^1, \vartheta^2) \upd r \right)^2 \dVol_{\mathbb{S}^2}
    + \int_{\bar{S}_{\tau, t}} |\phi|^2 \dVol_{\mathbb{S}^2} \\
  \end{split}
 \end{equation*}
By Cauchy-Schwarz we have
\begin{equation}
  \begin{split}
    &\int_{\mathbb{S}^2} \left( \int_{r' = r}^{R_t} \frac{ \partial |\phi|}{\partial r'} (\tau, r', \vartheta^1, \vartheta^2) \upd r \right)^2 \dVol_{\mathbb{S}^2} \\
    &\lesssim \left(\int_{r' = r}^{\sup_{(\vartheta_1, \vartheta_2)} R_t} (r')^{-2 + \alpha} \upd r' \right)\left( \int_{^t\Sigma_{\tau}} (1+r)^{-\alpha}\left( \frac{\partial |\phi|}{\partial r} \right)^2 r^2 \upd r \wedge \dVol_{\mathbb{S}^2}  \right) \\
    &\lesssim \frac{1}{(1-\alpha)}\frac{1}{r^{1-\alpha}} \int_{^t\Sigma_{\tau}}(1+r)^{-\alpha} \left( \frac{\partial |\phi|}{\partial r} \right)^2 r^2 \upd r \wedge \dVol_{\mathbb{S}^2} 
  \end{split}
\end{equation}
where we have assumed that $\alpha < 1$ and $r \geq r_0$.

Now, we recall the expression for $\partial_r$ given in \eqref{equation d dtau and ddr on sigma}. Additionally, for any vector field $X$ we have
\begin{equation*}
 \begin{split}
  \left(X|\phi|\right)^2 &= \left(X\sqrt{\phi_{\sigma_1 \ldots \sigma_n} \phi^{\sigma_1 \ldots \sigma_n}} \right)^2 \\
  &= \left( \frac{1}{2}|\phi|^{-1} X\left( \phi_{\sigma_1 \ldots \sigma_n} \phi^{\sigma_1 \ldots \sigma_n} \right) \right)^2 \\
  &= |\phi|^{-2} \left( \phi^{\sigma_1 \ldots \sigma_n} \D_X \phi_{\sigma_1 \ldots \sigma_n} \right)^2 \\
  &= |\phi|^{-2} \left( \phi^{\sigma_1 \ldots \sigma_n} \slashed{\D}_X \phi_{\sigma_1 \ldots \sigma_n} \right)^2 \\
  &\leq |\slashed{\D}_X \phi|^2
 \end{split}
\end{equation*}
where the penultimate line follows from the fact that $\phi$ is $S_{\tau,r}$-tangent. The considerations above lead to the bound
 \begin{equation*}
  \begin{split}
   \int_{S_{\tau,r}} |\phi|^2 \dVol_{\mathbb{S}^2}
    &\lesssim \frac{1}{(1-\alpha)}\frac{1}{r^{1-\alpha}} \int_{\mathbb{S}^2} \int_{r' = r}^{R_t} (1+r)^{-\alpha}\left(  |\slashed{\D}_L \phi|^2 \right)r^2 \upd r \wedge \dVol_{\mathbb{S}^2} 
    + \int_{\bar{S}_{\tau, t}} |\phi|^2 \dVol_{\mathbb{S}^2} \\
  \end{split}
 \end{equation*}

The second and third inequalities are proven similarly, where we also need to make use of the calculation
\begin{equation*}
 \begin{split}
  \left.\frac{\partial}{\partial r}\right|_{t, \vartheta^A} &= \left( \frac{\Lbar^0 - b^0 }{L^0 + \Lbar^0 - b^0}\right) L - \left(\frac{L^0}{L^0 + \Lbar^0 - b^0}\right)\Lbar + \left(\frac{L^0}{L^0 + \Lbar^0 - b^0}\right) b^A X_A \\
  &= \frac{1}{2}(L - \Lbar) + \frac{ \left(\Lbar^0_{(\text{small})} - L^0_{(\text{small})} - b^0 \right)}{2(L^0 + \Lbar^0 - b^0)}\left( L - \Lbar \right) + \left(\frac{L^0}{L^0 + \Lbar^0 - b^0}\right) b^A X_A 
 \end{split}
\end{equation*}
\end{proof}

\begin{proposition}[Estimating the spherical mean of a field defined by a point dependent change of basis]
\label{proposition spherical mean after change of basis}
Let $\phi_{(A)} = M_{(A)}^{\phantom{(A)}(a)} \phi_{(a)}$, where $M$ is a (possibly point dependent) change of basis matrix, and the $\phi_{(a)}$ are some $S_{\tau,r}$-tangent tensor fields which all have the same rank. Write
\begin{equation*}
 |\phi_{(\text{orig})}| := \sup_{(a)} |\phi_{(a)}|
\end{equation*}
and similarly define quantities involving derivatives, e.g.
\begin{equation*}
 |\slashed{\D}\phi_{(\text{orig})}| := \sup_{(a)} |\slashed{\D}\phi_{(a)}|
\end{equation*}
Finally, we define
\begin{equation*}
 |M_{(A)}| := \sup_{(a)} |M_{(A)}^{\phantom{(A)} (a)}|
\end{equation*}
and similarly for derivatives of the change of basis matrix $M$.

Then we have
\begin{equation}
 \begin{split}
    &\int_{S_{\tau,r}}|\phi_{(A)}|^2 \dVol_{\mathbb{S}^2} \lesssim \frac{1}{(1-\alpha)}\frac{1}{r^{1-\alpha}}\int_{^t\Sigma_\tau} (1+r)^{-\alpha}\left( |\slashed{\D}_L \phi|_{(A)}^2 + |L M_{(A)}|^2 |\phi_{(\text{orig})}|^2 \right) (r')^2 \upd r' \wedge \dVol_{\mathbb{S}^2} \\
     & + \frac{1}{(1-\alpha)}\sup_{x \in \bar{S}_{\tau_1,t}} \left\{\frac{1}{r(x)^{1-\alpha}}\right\}\int_{_{\tau_0}^{\tau_1}\bar{\Sigma}_t} (1+r)^{-\alpha}\bigg( |\slashed{\D}_{\Lbar} \phi|_{(A)}^2 + |\slashed{\D}_L \phi|_{(A)}^2 + |\partial M_{(A)}|^2 |\phi_{(\text{orig})}|^2 \\
     &\phantom{+ \frac{1}{(1-\alpha)}\sup_{x \in \bar{S}_{\tau_1,t}} \left\{\frac{1}{r(x)^{1-\alpha}}\right\}\int_{_{\tau_0}^{\tau_1}\bar{\Sigma}_t} (1+r)^{-\alpha}\bigg( }
     + \textit{Err}_{(t-\partial_r)}[\phi]_{(A)} \bigg) (r')^2 \upd r' \wedge \dVol_{\mathbb{S}^2} \\
    & + \int_{\bar{S}_{\tau_0,t}}|\phi_{(A)}|^2 \dVol_{\mathbb{S}^2} 
 \end{split}
\end{equation}
where the error term satisfies
The error term in the above inequalities satisfy: 
\begin{equation}
  \begin{split}
    \left| \textit{Err}_{(t-\partial_r)}[\phi]_{(A)} \right| &\lesssim \frac{|L^0_{(\text{small})}| + |\Lbar^0_{(\text{small})}| + |b^0|}{\left|L^0 + \Lbar^0 - b^0\right|} \left( |\slashed{\D}_L \phi|_{(A)}^2 + |\slashed{\D}_{\Lbar} \phi|_{(A)}^2 + |\partial M_{(A)}|^2 |\phi_{(\text{orig})}|^2 \right) \\
    &+ \frac{|L^0| |b^0|}{\left| L^0 + \Lbar^0 - b^0 \right|}\left( |\slashed{\nabla}\phi|_{(A)}^2 + |\slashed{\nabla}M_{(A)}|^2 |\phi_{(\text{orig})}|^2 \right) 
  \end{split}
\end{equation}

\end{proposition}

\begin{proof}
 By proposition \ref{proposition spherical mean in terms of energy} we have
 \begin{equation*}
  \begin{split}
    & \int_{S_{\tau,r}}|\phi_{(A)}|^2 \dVol_{\mathbb{S}^2} \lesssim \frac{1}{(1-\alpha)}\frac{1}{r^{1-\alpha}}\int_{^t\Sigma_\tau} (1+r)^{-\alpha}\left( |\slashed{\D}_L \phi_{(A)}|^2 \right) (r')^2 \upd r' \wedge \dVol_{\mathbb{S}^2} \\
    & + \frac{1}{(1-\alpha)}\sup_{x \in \bar{S}_{\tau_1,t}} \left\{\frac{1}{r(x)^{1-\alpha}}\right\}\int_{_{\tau_0}^{\tau_1}\bar{\Sigma}_t} (1+r)^{-\alpha}\left( |\slashed{\D}_{\Lbar} \phi_{(A)}|^2 + |\slashed{\D}_L \phi_{(A)}|^2 + \textit{Err}_{(t-\partial_r)}[\phi_{(A)}] \right) (r')^2 \upd r' \wedge \dVol_{\mathbb{S}^2} \\
    & + \int_{\bar{S}_{\tau_0,t}}|\phi_{(A)}|^2 \dVol_{\mathbb{S}^2} 
 \end{split}
 \end{equation*}
 Now, we have
 \begin{equation*}
  \begin{split}
   |\slashed{\D}_L \phi_{(A)}|^2 &= |\slashed{\D}_L \left( M_{(A)}^{\phantom{(A)}(a)} \phi_{(a)} \right)|^2 \\
   &= \left| (\slashed{\D}_L \phi)_{(A)} + \left( L M_{(A)}^{\phantom{(A)}(a)} \right) \phi_{(a)} \right|^2 \\
   &\lesssim | |\slashed{\D}_L \phi|_{(A)}^2  + | L M_{(A)}|^2 | \phi_{(\text{orig})} |^2
  \end{split}
 \end{equation*}
 and the other terms can be treated similarly.
 
\end{proof}

\begin{remark}[The notation $|\phi|_{(A)}$]
 We will sometimes use the notation $|\phi|_{(A)}$ instead of $|\phi_{(A)}|$. That is, when referring to a field defined by a point dependent change of basis, we might place the index $(A)$ outside of the delimiter, even when there is no derivative operator inside the delimiter. In such a case, there is no difference between $|\phi_{(A)}|$ and $|\phi|_{(A)}$, unlike, for example $|\slashed{\D}_L \phi_{(A)}|$ and $|\slashed{\D}_L \phi_{(A)}|$, which denote different objects: in the former, the matrix $M$ is applied to the field $\phi$ and then the derivative operator $\slashed{\D}$ is applied, while in the latter the derivative operator is applied before the matrix $M$. 
\end{remark}

\begin{proposition}[Weighted Hardy inequality on the surface $^t\Sigma_\tau$]
\label{proposition Hardy}
 We can use the following weighted version of Hardy's inequality to estimate integrals of $|\phi|^2$ in terms of its derivatives, as long as the term involving $|\phi|^2$ is multiplied by a factor which decays at a rate at least as fast as $(1+r)^{-2}$.

  Let $0 \leq \alpha < 1$. Then  the following inequality holds: for all $S_{\tau,r}$-tangent tensor fields $\phi$,
\begin{equation}
 \begin{split}
  \int_{^t\Sigma_\tau}(1+r)^{-1-\alpha}r^{-1} |\phi|^2 r^2 \upd r \wedge \dVol_{\mathbb{S}^2} 
  &\leq \frac{1}{(1-\alpha)^2}
  \int_{^t\Sigma_\tau} (1+r)^{-\alpha}|\slashed{\D}_L\phi|^2 \, r^2\upd r \wedge \dVol_{\mathbb{S}^2} \\
  &\phantom{\leq \frac{1}{(1-\alpha)^2} \Bigg(}
  + \frac{1}{(1-\alpha)} \int_{\bar{S}_{\tau, t}}(1+r)^{-\alpha}|\phi|^2 r \dVol_{\mathbb{S}^2}
 \end{split}
\end{equation}

Additionally, if $f(r)$ is a compactly supported function, supported away from the origin, and $\alpha \neq 1$, then we have
\begin{equation*}
\begin{split}
 \int_{^t\Sigma_\tau} r^{-\alpha}f |\phi|^2 \upd r \wedge \dVol_{\mathbb{S}^2} 
 &\leq
 \frac{1}{(1-\alpha)^2}\int_{^t\Sigma_\tau} fr^{-\alpha} |\slashed{\D}_L\phi|^2 \, r^2 \upd r \wedge \dVol_{\mathbb{S}^2}
 \\
 &\phantom{\leq}
 + \frac{1}{|1-\alpha|}\int_{^t\Sigma_\tau \cap \, \mathrm{supp}(f')} r^{1-\alpha} f' |\phi|^2 \upd r \wedge \dVol_{\mathbb{S}^2}
\end{split}
\end{equation*}

Finally, note that if $\alpha > 1$ then we have
\begin{equation*}
\int_{^t\Sigma_\tau \cap \{r \geq r_0\}} r^{-\alpha} |\phi|^2 \upd r \wedge \dVol_{\mathbb{S}^2}
\lesssim
\frac{1}{(\alpha - 1)^2} \int_{\Sigma_\tau \cap \{r \geq r_0\}} r^{2-\alpha} |\slashed{\D}_L\phi|^2 \upd r \wedge \dVol_{\mathbb{S}^2}
+ \frac{1}{\alpha -1} \int_{S_{\tau,r_0}} r^{1-\alpha} |\phi|^2 \dVol_{\mathbb{S}^2}
\end{equation*}

\end{proposition}

\begin{proof}
 First consider the case $0 < \alpha < 1$. We have
\begin{equation*}
 \begin{split}
  &\int_{^t\Sigma_\tau}(1+r)^{-1-\alpha}|\phi|^2 r \upd r \wedge \dVol_{\mathbb{S}^2} \\
  &= \int_{^t\Sigma_\tau} \partial_r\left( \frac{(1+r)^{-\alpha}(1+\alpha r)}{\alpha (1-\alpha)} - \frac{1}{\alpha(1-\alpha)} \right)|\phi|^2 \upd r \wedge \dVol_{\mathbb{S}^2} \\
  & \leq \frac{1}{1-\alpha}\Bigg( \int_{^t\Sigma_\tau} (1+r)^{-1-\alpha}\left( \left.\frac{\partial |\phi|^2}{\partial r}\right|_{\tau,\vartheta^1,\vartheta^2}\right) r^2 \upd r \wedge \dVol_{\mathbb{S}^2} + \int_{\bar{S}_{\tau,r}} (1+r)^{-\alpha}|\phi|^2 r \dVol_{\mathbb{S}^2} \Bigg)
 \end{split}
\end{equation*}
where in the last line we have integrated by parts. Now, making use of the second expression in equation \eqref{equation d dtau and ddr on sigma}, for any $\delta > 0$ we have
\begin{equation*}
 \begin{split}
  &\int_{^t\Sigma_\tau}(1+r)^{-\alpha}|\phi|^2 r \upd r \wedge \dVol_{\mathbb{S}^2} \\
  & \leq \frac{1}{1-\alpha}\Bigg( \int_{^t\Sigma_\tau} \left( \delta (1+r)^{-1-\alpha}r|\phi|^2 + \delta^{-1}(1+r)^{-1-\alpha}r^3|\slashed{\D}_L\phi|^2 \right)  \upd r \wedge \dVol_{\mathbb{S}^2} \\
  &\phantom{\leq \frac{1}{1-\alpha}\Bigg( } + \int_{\bar{S}_{\tau,r}} (1+r)^{-\alpha}|\phi|^2 r \dVol_{\mathbb{S}^2} \Bigg)
 \end{split}
\end{equation*}
Choosing $\delta (1-\alpha)^{-1}$ sufficiently small so that the first term on the right hand side can be absorbed by the left hand side proves the proposition in the case $0 < \alpha < 1$.

Now consider the case $\alpha = 0$. We have
\begin{equation*}
 \begin{split}
  &\int_{^t\Sigma_\tau}(1+r)^{-1}|\phi|^2 r \upd r \wedge \dVol_{\mathbb{S}^2} \\
  &= \int_{^t\Sigma_\tau} \partial_r \left( r - \log(1+r) \right)|\phi|^2 \upd r \wedge \dVol_{\mathbb{S}^2} \\
  &\lesssim \int_{^t\Sigma_\tau} \frac{r^2}{(1+r)} \left( \left.\frac{\partial |\phi|^2}{\partial r}\right|_{\tau,\vartheta^1,\vartheta^2}\right) \upd r \wedge \dVol_{\mathbb{S}^2} + \int_{\bar{S}_{\tau,r}} |\phi|^2 r \, \dVol_{\mathbb{S}^2} \\
  &\lesssim \int_{^t\Sigma_\tau}  \left( \delta\frac{r}{(1+r)}|\phi|^2 + \delta^{-1}\frac{r^3}{(1+r)}|\slashed{\D}_L \phi|^2  \right) \upd r \wedge \dVol_{\mathbb{S}^2} + \int_{\bar{S}_{\tau,r}} |\phi|^2 r \, \dVol_{\mathbb{S}^2} \\
  \end{split}
 \end{equation*}
Once again we can pick $\delta$ sufficiently small that the first term on the right hand side can be absorbed by the left hand side.

The proof of the second part of the proposition follows in almost exactly the same way; the condition on the support of $f$ ensures that the integrals are well-defined. Finally, for the third part of the proposition we note that taking $\alpha > 1$ means that the boundary term on $\bar{S}_{t,r}$ has the ``right'' sign and so can be ignored.
\end{proof}

\begin{proposition}[Weighted Hardy inequality on the surface $^t\Sigma_\tau$ after a point dependent change of basis]
\label{proposition Hardy after change of basis}

 Let $0 \leq \alpha < 1$, and let $\phi_{(a)}$ be a collection of $S_{\tau,r}$-tangent tensor fields. Let $\phi_{(A)} = M_{(A)}^{\phantom{(A)}(a)}$, and use the same notational conventions as above. Then
\begin{equation}
 \begin{split}
  &\int_{^t\Sigma_\tau}(1+r)^{-1-\alpha}r^{-1} |\phi|_{(A)}^2 r^2 \upd r \wedge \dVol_{\mathbb{S}^2} \\
  &\lesssim \frac{1}{(1-\alpha)^2}
  \int_{^t\Sigma_\tau} (1+r)^{-\alpha} \left( |\slashed{\D}_L\phi|_{(A)}^2 + |L M_{(A)}|^2 |\phi_{(\text{orig})}|^2 \right) \, r^2\upd r \wedge \dVol_{\mathbb{S}^2} \\
  &\phantom{\leq \frac{1}{(1-\alpha)^2} \Bigg(}
  + \frac{1}{(1-\alpha)} \int_{\bar{S}_{\tau, t}}(1+r)^{-\alpha}|\phi|_{(A)}^2 r \dVol_{\mathbb{S}^2}
 \end{split}
\end{equation}

Additionally, if $f(r)$ is a compactly supported function, supported away from the origin, and $\alpha \neq 1$, then we have
\begin{equation*}
 \begin{split}
 &\int_{^t\Sigma_\tau} r^{-\alpha}f |\phi|_{(A)}^2 \upd r \wedge \dVol_{\mathbb{S}^2} \\
 &\lesssim \frac{1}{(1-\alpha)^2}\int_{^t\Sigma_\tau} fr^{-\alpha} \left( |\slashed{\D}_L\phi|_{(A)}^2 + |L M_{(A)}|^2 |\phi_{(\text{orig})}|^2\right) \, r^2 \upd r \wedge \dVol_{\mathbb{S}^2} \\
 &\phantom{\lesssim}
 + \frac{1}{|1-\alpha|}\int_{^t\Sigma_\tau \cap \, \mathrm{supp}(f')} r^{1-\alpha} f' |\phi|_{(A)}^2 \upd r 
 \end{split}
\end{equation*}

\end{proposition}

\begin{proof}
To prove this proposition, we begin by applying proposition \ref{proposition Hardy} to the field $\phi_{(A)}$. We then use 
\begin{equation*}
\begin{split}
 \slashed{\D}_L (\phi_{(A)}) &= M_{(A)}^{\phantom{(A)}(a)} \slashed{\D}_L \phi_{(a)} + (LM_{(A)}^{\phantom{(A)}(a)}) \phi_{(a)} \\
 &= (\slashed{\D}_L \phi)_{(A)} + (LM_{(A)}^{\phantom{(A)}(a)}) \phi_{(a)} \\
 \end{split}
\end{equation*}
\end{proof}

\begin{proposition}[Weighted Hardy inequality on the surface $\bar{\Sigma}_t$]
 \label{proposition Hardy on constant t}
 Similarly to proposition \ref{proposition Hardy}, let $0 \leq \alpha \leq 1$. Then the following inequality holds: for all $S_{\tau,r}$-tangent tensor fields $\phi$,
\begin{equation}
 \begin{split}
  &\int_{_{\tau_0}^{\tau_1}\Sigma_t}(1+r)^{-1-\alpha} |\phi|^2 r \upd r \wedge \dVol_{\mathbb{S}^2} \\ 
  &\lesssim \frac{1}{(1-\alpha)^2}
  \int_{_{\tau_0}^{\tau_1}\Sigma_t} (1+r)^{-\alpha} \left(|\slashed{\D}_L\phi|^2 + |\slashed{\D}_{\Lbar}\phi|^2 + \textit{Err}_{(t-\partial_r)}[\phi]\right) r^2\upd r \wedge \dVol_{\mathbb{S}^2} \\
  &\phantom{\leq \frac{1}{1-\alpha} }
  + \frac{1}{(1-\alpha)} \int_{\bar{S}_{\tau_0, t}}(1+r)^{-\alpha}|\phi|^2 r \dVol_{\mathbb{S}^2}
 \end{split}
\end{equation}
\end{proposition}
\begin{proof}
 The proof is almost identical to the proof of proposition \ref{proposition Hardy}. Note that an additional term appears at the ``inner'' boundary $\bar{S}_{\tau_1, t}$ when integrating by parts, but this term has a good sign and so can be dropped.
\end{proof}

\begin{proposition}[Weighted Hardy inequality on the surface $\bar{\Sigma}_t$ after a point dependent change of basis]
 \label{proposition Hardy on constant t after change of basis}
 Similarly to proposition \ref{proposition Hardy after change of basis}, let $0 \leq \alpha \leq 1$. Let $\phi_{(a)}$ be a collection of $S_{\tau,r}$-tangent tensor fields of the same rank, and let $\phi_{(A)} := M_{(A)}^{\phantom{(A)}(a)}\phi_{(a)}$. Then
\begin{equation}
 \begin{split}
  &\int_{_{\tau_0}^{\tau_1}\Sigma_t}(1+r)^{-1-\alpha} |\phi|_{(A)}^2 r \upd r \wedge \dVol_{\mathbb{S}^2} \\ 
  &\lesssim \frac{1}{(1-\alpha)^2}
  \int_{_{\tau_0}^{\tau_1}\Sigma_t} (1+r)^{-\alpha} \left(|\slashed{\D}_L\phi|_{(A)}^2 + |\slashed{\D}_{\Lbar}\phi|_{(A)}^2 + |\partial M_{(A)}|^2 |\phi_{(\text{orig})}|^2 + \textit{Err}_{(t-\partial_r)}[\phi]_{(A)}\right) r^2\upd r \wedge \dVol_{\mathbb{S}^2} \\
  &\phantom{\leq \frac{1}{1-\alpha} }
  + \frac{1}{(1-\alpha)} \int_{\bar{S}_{\tau_0, t}}(1+r)^{-\alpha}|\phi|_{(A)}^2 r \dVol_{\mathbb{S}^2}
 \end{split}
\end{equation}
\end{proposition}

Sometimes we will encounter lower order terms on surfaces of constant $t$, whose coefficients decay in $\tau$ rather than in $r$. Since $\tau$ also increases as we move outwards along these hypersurfaces, it is also possible to prove a Hardy inequality involving these kinds of weights.

\begin{proposition}[Weighted Hardy inequality on the surface $\bar{\Sigma}_t$ with decaying weights in $\tau$]
 \label{proposition Hardy on constant t with tau weight}
 Similarly to proposition \ref{proposition Hardy}, let $0 \leq \alpha \leq 1$. Then the following inequality holds: for all $S_{\tau,r}$-tangent tensor fields $\phi$,
\begin{equation}
 \begin{split}
  &\int_{_{\tau_0}^{\tau_1}\Sigma_t}(1+r)^{-\alpha} (1+\tau)^{-2} |\phi|^2 \upd r \wedge \dVol_{\mathbb{S}^2} \\ 
  &\lesssim 
  \int_{_{\tau_0}^{\tau_1}\Sigma_t} (1+r)^{-\alpha} \left(|\mu|^2 |\slashed{\D}\phi|^2 + \textit{Err}_{(t-\partial_{\tau})}[\phi]\right) r^2\upd r \wedge \dVol_{\mathbb{S}^2} \\
  &\phantom{\leq \frac{1}{1-\alpha} }
  + \int_{\bar{S}_{\tau_0, t}}(1+r)^{-\alpha} (1+\tau)^{-1} |\phi|^2 r^2 \dVol_{\mathbb{S}^2}
 \end{split}
\end{equation}
where
\begin{equation}
 \begin{split}
  |\textit{Err}_{(t-\partial_{\tau})}[\phi]|
  &\lesssim |\mu|^2(1+r)^{-\alpha} (1+\alpha^2)  \left| \frac{L^0 + \Lbar^0 - b^0}{2 L^0} \right|^2 ((1+r)^{-2} + r^{-2} )  |\phi|^2 \\
  &\phantom{\lesssim}
  + |\mu|^2(1+r)^{-\alpha} \left( \left| \frac{\Lbar^0 - L^0 - b^0}{L^0} \right|^2 + |b|^2 \right) |\slashed{\D}\phi|^2 
 \end{split}
\end{equation}

\end{proposition}
\begin{proof}
We can calculate
\begin{equation*}
 \left. \frac{\partial}{\partial \tau}\right|_{t, \vartheta^1, \vartheta^2} = -\frac{1}{2}\mu \left( \frac{ \Lbar^0 - b^0}{L^0} \right) L + \frac{1}{2}\mu \Lbar - \frac{1}{2}\mu b^A X_A
\end{equation*}
and so, in particular, we have
\begin{equation*}
 \left. \frac{\partial r}{\partial \tau}\right|_{t, \vartheta^1, \vartheta^2} = -\mu \left( \frac{L^0 + \Lbar^0 - b^0}{2 L^0} \right)
\end{equation*}
Hence we have
\begin{equation*}
 \begin{split}
  &\int_{_{\tau_0}^{\tau_1}\Sigma_t}(1+r)^{-\alpha} (1+\tau)^{-2} |\phi|^2 r^2 \upd r \wedge \dVol_{\mathbb{S}^2} \\
  &= \int_{_{\tau_0}^{\tau_1}\Sigma_t} -(1+r)^{-\alpha} \left(\frac{\partial}{\partial \tau} (1+\tau)^{-1}\right) |\phi|^2 r^2 \upd r \wedge \dVol_{\mathbb{S}^2} \\
  &= \int_{_{\tau_0}^{\tau_1}\Sigma_t} (1+r)^{-\alpha} \bigg( 
  \left( \mu \left( \frac{L^0 + \Lbar^0 - b^0}{2 L^0} \right) (\alpha(1+r)^{-1} - 2r^{-1} ) (1+\tau)^{-1}\right) |\phi|^2 \\
  &\phantom{= \int_{_{\tau_0}^{\tau_1}\Sigma_t}}
  + 2(1+\tau)^{-1} \mu \left( -\frac{1}{2} \left( \frac{ \Lbar^0 - b^0}{L^0} \right) \slashed{\D}_L \phi + \frac{1}{2} \slashed{\D}_{\Lbar} \phi  - \frac{1}{2} b^A \slashed{\nabla}_A \phi \right) \cdot \phi \bigg)
  r^2 \upd r \wedge \dVol_{\mathbb{S}^2} \\
  &- \int_{\bar{S}_{\tau_1, t}} (1+r)^{-\alpha} (1+\tau)^{-1} |\phi|^2 r^2 \upd r \wedge \dVol_{\mathbb{S}^2}
  + \int_{\bar{S}_{\tau_0, t}} (1+r)^{-\alpha} (1+\tau)^{-1} |\phi|^2 r^2 \upd r \wedge \dVol_{\mathbb{S}^2} \\
  &\lesssim \int_{_{\tau_0}^{\tau_1}\Sigma_t}
  \bigg( (1+\alpha^2)|\mu|^2 \left| \frac{L^0 + \Lbar^0 - b^0}{2 L^0} \right|^2 ((1+r)^{-2} + r^{-2} ) (1+r)^{-\alpha} |\phi|^2 \\
  &\phantom{= \int_{_{\tau_0}^{\tau_1}\Sigma_t}}
  + (1+r)^{-\alpha} |\mu|^2 \left( 1 + \left| \frac{\Lbar^0 - b^0}{L^0} \right|^2 + |b|^2 \right) |\slashed{\D}\phi|^2
  \bigg) r^2 \upd r \wedge \dVol_{\mathbb{S}^2} \\
  &
  + \int_{\bar{S}_{\tau_0, t}} (1+r)^{-\alpha} (1+\tau)^{-1} |\phi|^2 r^2 \upd r \wedge \dVol_{\mathbb{S}^2} \\
  \end{split}
\end{equation*}

\end{proof}

By using proposition \ref{proposition spherical mean in terms of energy} we can estimate the spherical mean of a field $\phi$ in terms of an ``energy'' type quantity. The estimate in proposition \ref{proposition spherical mean in terms of energy} is also useful when estimating the spherical mean of a field which is ``weighted'' by some factor of $r^{k}$, but only when $k < 1$. However, we sometimes need to estimate the spherical mean of a field weighted by a higher power of $r$, for which we use the following pair of propositions.

\begin{proposition}[An estimate for higher weighted spherical integrals of $|\phi|^2$]
	\label{proposition higher weighted spherical integral}
	Let $\psi = r\phi$, and let $r \geq r_0$. Then the spherical mean of $\phi$ satisfies
	\begin{equation*}
	\int_{S_{\tau,r}} |\phi|^2 \dVol_{\mathbb{S}^2}
	\lesssim
	\left(\frac{r_0}{r} \right)^2 \int_{S_{\tau,r_0}} |\phi|^2 \dVol_{\mathbb{S}^2}
	+ \frac{1}{1-p} r^{-1-p}\int_{\Sigma_\tau \cap \{r_0 \leq r' \leq r\}} (r')^p |\slashed{\D}_L\psi|^2   \upd r' \wedge \dVol_{\mathbb{S}^2}
	\end{equation*}
	
\end{proposition}

\begin{proof}
	Following calculations which are almost identical to those in proposition \ref{proposition spherical mean in terms of energy} we have, for $p < 1$,
	\begin{equation*}
	\int_{S_{\tau,r}} |\psi|^2 \dVol_{\mathbb{S}^2}
	\lesssim
	\int_{S_{\tau,r_0}} |\psi|^2 \dVol_{\mathbb{S}^2}
	+ \frac{1}{1-p} r^{1-p}\int_{\Sigma_\tau \cap \{r_0 \leq r' \leq r\}} (r')^p |\slashed{\D}_L\psi|^2   \upd r' \wedge \dVol_{\mathbb{S}^2}
	\end{equation*}
	and so
	\begin{equation*}
	\int_{S_{\tau,r}} |\phi|^2 \dVol_{\mathbb{S}^2}
	\lesssim
	\left(\frac{r_0}{r} \right)^2 \int_{S_{\tau,r_0}} |\phi|^2 \dVol_{\mathbb{S}^2}
	+ \frac{1}{1-p} r^{-1-p}\int_{\Sigma_\tau \cap \{r_0 \leq r' \leq r\}} (r')^p |\slashed{\D}_L\psi|^2   \upd r' \wedge \dVol_{\mathbb{S}^2}
	\end{equation*}
\end{proof}

\begin{proposition}[An estimate for higher weighted integrals of $|\phi|^2$]
	\label{proposition higher weighted spatial integral}
	
	Define
	\begin{equation*}
	\begin{split}
	\mathcal{E}^{(T,\alpha)}
	&=
	\lim_{t \rightarrow \infty}
	\bigg(
	\int_{^t\Sigma_\tau} (1+r)^{-\alpha}\left( |\slashed{\D}_L \phi|^2 \right) (r)^2 \upd r \wedge \dVol_{\mathbb{S}^2} \\
	&\phantom{= \lim_{t \rightarrow \infty} \bigg(}
	+ \int_{_{\tau_0}^{\tau_1}\bar{\Sigma}_t} (1+r)^{-\alpha}\left( |\slashed{\D}_{\Lbar} \phi|^2 + |\slashed{\D}_L \phi|^2 + \textit{Err}_{(t-\partial_r)}[\phi] \right) (r)^2 \upd r \wedge \dVol_{\mathbb{S}^2} \\
	&\phantom{=\lim_{t \rightarrow \infty} \bigg(}
	+ \int_{\bar{S}_{\tau_0,t}}|\phi|^2 \dVol_{\mathbb{S}^2}
	\bigg)
	\\
	\\
	\mathcal{E}^{(L,p)} &= \int_{\Sigma_\tau \cap \{r \geq r_0\}} r^p |L\psi|^2 \upd r \wedge \dVol_{\mathbb{S}^2}
	\end{split}
	\end{equation*}
	where, as usual, $\psi = r\phi$. 
	
	Then we have
	\begin{equation*}
	\int_{\Sigma_\tau \cap \{r \geq r_0\}} r^{p - \delta} |\phi|^2 \upd r \wedge \dVol_{\mathbb{S}^2}
	\lesssim
	(r_0)^{1 + \alpha} \frac{1}{(1-\alpha)(1-p+\delta)} \mathcal{E}^{(T, \alpha)} + \delta^{-1} (r_0)^{-\delta} \mathcal{E}^{(L,p)}
	\end{equation*}
	
\end{proposition}

\begin{proof}
	Again, following calculations which are almost identical to those in proposition \ref{proposition spherical mean in terms of energy} we have, for $p < 1$,
	\begin{equation*}
	\int_{S_{\tau,r}} |\psi|^2 \dVol_{\mathbb{S}^2}
	\lesssim
	\int_{S_{\tau,r_0}} |\psi|^2 \dVol_{\mathbb{S}^2}
	+ \frac{1}{1-p} r^{1-p}\int_{\Sigma_\tau \cap \{r_0 \leq r' \leq r\}} (r')^p |L\psi|^2   \upd r' \wedge \dVol_{\mathbb{S}^2}
	\end{equation*}
	The first term is can be bounded by bounded by $\mathcal{E}^{(T, \alpha)}$ (using proposition \ref{proposition spherical mean in terms of energy})  and the integrand in the second term is bounded by $\mathcal{E}^{(L,p)}$.
	
	Now, multiplying by $r^{-2 + p - \delta}$ and integrating from $r = r_0$ to infinity, we obtain
	\begin{equation*}
	\int_{\Sigma_\tau \cap \{r \geq r_0 \}} r^{-2+p-\delta}|\psi|^2 \upd r \wedge \dVol_{\mathbb{S}^2}
	\lesssim
	(r_0)^{1 + \alpha} \frac{1}{(1-\alpha)(1-p+\delta)} \mathcal{E}^{(T, \alpha)} + \delta^{-1} (r_0)^{-\delta} \mathcal{E}^{(L,p)}
	\end{equation*}

\end{proof}

It is often useful to have an expression relating spacetime integrals to spatial integrals, which are then integrated over time. In particular, such an expression is necessary when controlling certain error terms by the corresponding flux terms on the hypersurfaces $\Sigma_\tau$, with the help of Gronwall's inequality. As such, we have the following expression, which is an easy application of the coarea formula.

\begin{proposition}[The coarea formula for spacetime integrals]
\label{proposition coarea}
 Let $f$ be any smooth function on $\mathcal{M}$. Then
 \begin{equation}
 \label{equation coarea spacetime}
  \int_{^t\mathcal{M}_{\tau_0}^{\tau_1} } f\, \dVol_g = \int_{\tau_0}^{\tau_1} \left(
  \int_{^t\Sigma_\tau \cap r< r_0} f (g^{-1})^{00}\sqrt{-\det g} \, \upd x^1 \wedge \upd x^2 \wedge \upd x^3
   + \int_{^t\Sigma_\tau \cap r\geq r_0}  f \Omega^2 \upd r \wedge \dVol_{\mathbb{S}^2}  \right) \upd \tau
 \end{equation}
 where the determinant $\det g$ is calculated relative to the rectangular coordinates $x^a$.

\end{proposition}

We also have a coarea formula relating spacetime integrals to integrals over the level sets of the function $t$:

\begin{proposition}[The coarea formula for spacetime integrals]
\label{proposition coarea t}
 Let $f$ be any smooth function on $\mathcal{M}$. Then
 \begin{equation}
 \label{equation coarea t}
  \int_{^t\mathcal{M}_{\tau_0}^{\tau_1} } f\, \dVol_g = \int_{t' = \tau_0}^{t} \left(
  \int_{^{\tau_1}_{\tau_0}\bar{\Sigma}_{t'} \cap r< r_0} (f g^{-1})^{00} \sqrt{-\det g} \, \upd x^1 \wedge \upd x^2 \wedge \upd x^3
   + \int_{^t\Sigma_\tau \cap r\geq r_0}  f \dVol_{(\bar{\Sigma}_{t'})}  \right) \upd t'
 \end{equation}
 where the determinant $\det g$ is calculated relative to the rectangular coordinates $x^a$.

\end{proposition}

\begin{proposition}[Gronwall inequality]
 \label{proposition Gronwall}
  We shall prove the following form of Gronwall's inequality: Let $f(t)>0$, $h(t)>0$ be continuous functions, let $g(t) \geq 0 \in L^1_{\text{loc}}$ and let $G(t, T) > 0$ be a $C^1$, nondecreasing function of its first argument, satisfying $G(T,T)\geq 0$ for any $T$. Suppose $f$ satisfies the integral inequality
\begin{equation}
\label{equation Gronwall inequality}
 h(T) + f(T) \leq \int_{T_0}^T \bigg( \epsilon g(t)f(t) \bigg)\upd t + f(T_0) + h(T_0) + G(T, T_0)
\end{equation}
for all $T \geq 0$.

Then $f$ satisfies
\begin{equation}
\label{equation Gronwall conclusion}
 h(T) + f(T) \leq \exp\left(\int_{T_0}^T \epsilon g(t)\upd t \right) \bigg( f(T_0) + h(T_0) + G(t, T_0) \bigg)
\end{equation}
for all $T \geq 0$.
\end{proposition}

\begin{proof}
 First, we assume that, for all $T_0 \leq T \leq T_{\text{max}}$ and for some $\delta_0 > 0$, we have
\begin{equation}
 \label{equation Gronwall internal 1}
 f(T) + h(T) \leq (1+\delta_0) \exp\left(\int_{T_0}^T (1+\delta_0) \epsilon g(t)\upd t \right) \bigg( f(T_0) + h(T_0) + G(T, T_0) \bigg)
\end{equation}
Since $\delta_0 > 0$ this clearly holds (by continuity) for some $T_{\text{max}} > T_0$. Now, from the inequality \eqref{equation Gronwall inequality} satisfied by $f$, we find that, for all $T \leq T_{\text{max}}$,
\begin{equation*}
  \begin{split}
    f(T) + h(T) &\leq \int_{T_0}^T (1+\delta_0)\epsilon g(t) \exp\left(\int_{T_0}^t (1+\delta_0) \epsilon g(t')\upd t' \right) \bigg( f(T_0) + h(T_0) + G(t, T_0) \bigg) \upd t \\
    & \phantom{\leq} + f(T_0) + h(T_0) + G(T, T_0)
  \end{split}
\end{equation*}
Now, we calculate
\begin{equation*}
  \begin{split}
   &\int_{T_0}^T (1+\delta_0)\epsilon g(t) \exp\left(\int_{T_0}^t (1+\delta_0)\epsilon g(t')\upd t' \right) \bigg( f(T_0) + h(T_0) + G(t, T_0) \bigg) \upd t 
   + f(T_0) + h(T_0) + G(T, T_0) \\ \\
   &= \int_{T_0}^T \bigg( \partial_t \left( \exp\left(\int_{T_0}^t (1+\delta_0) \epsilon g(t')\upd t' \right) \bigg( f(T_0) + h(T_0) + G(t, T_0) \bigg) \right) 
   \\
   &\phantom{= \int_{T_0}^T \bigg(}
   - \exp\left(\int_{T_0}^t (1+\delta_0)\epsilon g(t')\upd t' \right) \partial_t G(t, T_0) \bigg) \upd t
   + f(T_0) + h(T_0) + G(T, T_0) \\ \\
   &= \exp\left(\int_{T_0}^T (1+\delta_0)\epsilon g(t)\upd t \right) \bigg( f(T_0) + h(T_0) + G(T, T_0) \bigg)
   - \int_{T_0}^T \exp\left( \int_{T_0}^t (1+\delta_0)\epsilon g(t')\upd t' \right) \partial_t G(t, T_0) \upd t \\ \\
   &\leq \exp\left(\int_{T_0}^T (1+\delta_0)\epsilon g(t)\upd t \right) \bigg( f(T_0) + h(T_0) + G(T, T_0) \bigg)
  \end{split}
\end{equation*}
where in the last line we have used the facts that $\partial_t G(t, T_0) \geq 0$.

Hence, the original bound on $f$ (equation \eqref{equation Gronwall internal 1}) actually holds \emph{without} the first factor of $(1+\delta_0)$ up to the time $T_{\text{max}}$. Hence, the bound \eqref{equation Gronwall internal 1} actually holds for $T$ in the range $T \in [T_0, T_{\text{max}} + \delta_1]$ for some $\delta_1 > 0$. Hence we can take $T_{\text{max}} = \infty$.

Now, we have shown that the bound in equation \eqref{equation Gronwall internal 1} holds for all $T$ and for all $\delta > 0$. Suppose, for the sake of contradiction, that it \emph{does not} hold for $\delta = 0$. Then there is some time $T_1 > T_0$ such that
\begin{equation*}
	f(T_1) + h(T_1) > \exp\left(\int_{T_0}^{T_1} \epsilon g(t)\upd t \right) \bigg( f(T_0) + h(T_0) + G(t, T_0) \bigg)
\end{equation*}
Since this is a strict inequality, for all sufficiently small $\delta$ it must be that
\begin{equation*}
	f(T_1) + h(T_1) > (1+\delta_0)\exp\left(\int_{T_0}^{T_1} (1+\delta_0)\epsilon g(t)\upd t \right) \bigg( f(T_0) + h(T_0) + G(t, T_0) \bigg)
\end{equation*}
since the right hand side is continuous in $\delta_0$. But this contradicts the bounds which we have already established. Hence we can set $\delta_0 = 0$ in equation \eqref{equation Gronwall internal 1}.

\end{proof}

\begin{proposition}[An integral inequality]
	\label{proposition integral inequality}
	
	Let $f(t)$ be an integrable function satisfying the bound
	\begin{equation*}
	\int_{\tau}^{\tau_1} f(t) \upd t \leq C(1+\tau)^{\alpha}
	\end{equation*}
	for all $\tau \geq \tau_0$ and for all $\tau_1 \geq \tau$, and where $C > 0$ and $\alpha$ are constants.
	
	Then, if $\beta \geq 0$ is non-negative constant such that $\alpha + \beta < 1$, we have
	\begin{equation*}
	\int_{\tau}^{\tau_1} (1+t)^\beta f(t) \upd t
	\leq C\left( 1 + \frac{\beta}{\alpha + \beta} \right) (1+\tau)^{\alpha + \beta}
	\end{equation*}
	
\end{proposition}

\begin{proof}
	If $\alpha + \beta < 1$ then we have
	\begin{equation*}
	\begin{split}
	\int_{\tau}^{\tau_1} (1+t)^\beta f(t) \upd t
	&=
	\int_{\tau}^{\tau_1} (1+t)^\beta \frac{\partial}{\partial t} \left( -\int_{t}^{\tau_1} f(t') \upd t' \right) \upd t
	\\
	&=
	\int_{\tau}^{\tau_1} \beta (1+t)^{-1 + \beta} \left( \int_{t}^{\tau_1} f(t') \upd t' \right) \upd t
	+ (1+\tau)^\beta \int_{\tau}^{\tau_1} f(t) \upd t
	\\
	&\leq C \beta \int_{\tau}^{\tau_1} (1+t)^{-1 + \alpha + \beta} \upd t
	+ C (1+\tau)^{\alpha + \beta}
	\\
	&\leq C\left( 1 + \frac{\beta}{\alpha + \beta} \right) (1+\tau)^{\alpha + \beta}
	\end{split}
	\end{equation*}

\end{proof}

% \begin{equation}
% \begin{split}
%  \int_{r = r_1}^{\infty} \left( \int_{S_{\tau,r}} |\phi|^2 r^{p-2} \dVol_{\slashed{g}(\tau,r)} \right) \upd r
%  &\lesssim \frac{1}{p-1} (r_1)^{p-1} \int_{S_{\tau,r_1}} |\phi|^2 \dVol_{\slashed{g}(\tau,r)} \\
%  &+ \frac{1}{(p-1)^2} \int_{r = r_1}^{\infty} \left(\int_{S_{\tau,r}} |\slashed{\D}_L\psi|^2 r^{p-2} \dVol_{\slashed{g}(\tau,r)} \right) \upd r \\
%  &+ \textit{Err}_{(p-\text{Hardy})}
% \end{split}
% \end{equation}
% where $\psi := r\phi$ and the error term satisfies
% \begin{equation}
% \begin{split}
% 	|Err_{(p-\text{Hardy})}| &\lesssim \int_{r=r_1}^{\infty} \Bigg( \int_{S_{\tau,r}} \bigg(
% 	|f'_{(\alpha)}|^2 \mu^2 \left( r^{p}|\slashed{\D}_{\Lbar}\phi|^2 + r^{p - 2}|\rho^A \slashed{\nabla}_A \psi|^2 \right) \\
% 	&\phantom{\lesssim}
% 	+\bigg( \left(1-\frac{1}{2}f'_{(\alpha)}\mu \right)\tr_{\slashed{g}}\chi_{(\text{small})}
% 	+ f'_{(\alpha)} \mu r^{p-2} + f'_{(\alpha)}\mu r^{p-1}\tr_{\slashed{g}}\chibar_{(\text{small})} \\
% 	&\phantom{\lesssim}
% 	+ f'_{(\alpha)} \slashed{\nabla}_A(\mu \rho^A) r^{p-1} \bigg)|\phi|^2 \bigg) \dVol_{\slashed{g}(\tau,r)}  \Bigg) \upd r 
% \end{split}
% \end{equation}
 
% \begin{proof}
% We begin by writing
% \begin{equation*}
%  \int_{r = r_1}^{\infty} \left(\int_{S_{\tau,r}} |\phi|^2 r^{p-2} \dVol_{\slashed{g}(\tau,r)} \right) \upd r
%  = \int_{r=r_1}^{\infty} \frac{1}{p-1} \partial_r (r^{p-1}) \left( \int_{S_{\tau,r}} |\psi|^2 r^{-2}\dVol_{\slashed{g}(\tau,r)} \right) \upd r
% \end{equation*}
% Integrating by parts an proceeding as in proposition \ref{proposition Hardy Sigma tau} yields the proof.
% \end{proof}

\section{Boundary terms in the energy estimates}
Given a vector field $Z \in \mathcal{Z}$, there will be a corresponding boundary term in the energy estimate of the form
\begin{equation*}
 \int_{\Sigma_\tau} \imath_{({^{(Z)}J})} \dVol_{g}
\end{equation*}
Similarly, given a modified energy current ${^{(Z,f_Z)}\tilde{J}}$ we will find that the energy estimate involves a boundary term of the form
\begin{equation*}
  \int_{\Sigma_\tau} \imath_{({^{(Z,f_Z)}\tilde{J}})} \dVol_{g}
\end{equation*}
In this section, we will provide expressions for the boundary terms corresponding to the (modified) energy currents ${^{(Z)}J}$ (or ${^{(Z,f_Z)}\tilde{J}}$) corresponding to each of the multiplier vector fields $Z \in \mathcal{Z}$.

%We first define a quantity, related to the determinant of the metric in geometric coordinates, which occurs frequently in the expressions in this section.
% 
%\begin{definition}[The quantity $G$]
% We define the scalar density function $\sqrt{G}$ as function relating the volume forms $\dVol_{(\mathbb{S}^2, G)}$ and $\dVol_{\mathbb{S}^2}$, which are both volume forms on $\mathbb{S}^2$ and hence are proportional. Specifically, we set
%\begin{equation}
% \dVol_{(\mathbb{S}^2, G)} = \sqrt{G} \dVol_{\mathbb{S}^2}
%\end{equation}
%Note that this is technically a scalar density on the spheres $S_{\tau,r}$, rather than a scalar field.

%\end{definition}

\begin{proposition}[The weighted $T$-energy boundary term]
\label{proposition boundary weighted T}
 The boundary term on $^t\Sigma_\tau$ in the weighted $T$-energy is given by the following expression:

 \begin{equation}
  \begin{split}
  &\int_{^t\Sigma_\tau} \imath_{{^{(wT)}J}} \dVol_g = \int_{^t\Sigma_{\tau}} \frac{1}{2}w\left( \left(1-\frac{1}{2}\mu f'_{(\alpha)} \right)(\slashed{\D}_L\phi)^2 + |\slashed{\nabla}\phi|^2 \right) \Omega^2 \upd r \wedge \dVol_{\mathbb{S}^2} \\
  \end{split}
 \end{equation}

In addition, the boundary term on a surface of constant $t$ in the $T$-energy is given by

\begin{equation}
 \begin{split}
  \int_{\{t = t_0\}} \imath_{{^{(T)}J}} \dVol_g &= \int_{\{t = t_0\}} \frac{1}{2} \frac{1}{(L^0 + \Lbar^0 - b^0)} w\left( |\slashed{\D}_L \phi|^2 + |\slashed{\D}_{\Lbar} \phi|^2 + 2|\slashed{\nabla}\phi|^2 + \textit{Err}_{(T,t_0\text{-bdy})} \right) \Omega^2\upd r \wedge \dVol_{\mathbb{S}^2}
 \end{split}
\end{equation}

 where $\textit{Err}_{(T,t_0\text{-bdy})}$ satisfies
 \begin{equation}
 	|\textit{Err}_{(T,t_0\text{-bdy})}| \lesssim \left( |(L_{(\text{small})})^0| + |(\Lbar_{(\text{small})})^0| + |\slashed{\nabla} t| \right) \left( |\slashed{\D}_{\Lbar} \phi|^2 + |\slashed{\D}_L \phi|^2 + |\slashed{\nabla}\phi|^2 \right)
 \end{equation}

\end{proposition}

\begin{proof}
 We begin by noting that
 \begin{equation*}
  {^{(wT)}J^\mu}[\phi] = \frac{1}{2}w(\slashed{\D}_L\phi + \slashed{\D}_{\Lbar}\phi)\cdot \slashed{\D}^\mu \phi - \frac{1}{4}w(L^\mu + \Lbar^\mu) \left( -(\slashed{\D}_L\phi)\cdot(\slashed{\D}_{\Lbar}\phi) + |\slashed{\nabla}\phi|^2 \right)
 \end{equation*}
%  A similar formula, but without the $\mu$ weight, holds for the unweighted $\check{T}$ energy. We also need to compute
 \begin{equation*}
  (g^{-1})^{\mu\nu}(\upd\tau)_\nu = -\mu^{-1}L^\mu
 \end{equation*}
 and we also make use of the fact that $L\tau + \Lbar\tau = 2\mu^{-1}$, together with the expressions in propositions \ref{proposition volume form g} and \ref{proposition volume form gbar}. Additionally, we make use of the fact that the volume form can be written as
\begin{equation*}
 \dVol_g := -\frac{2}{(L^0 + \Lbar^0 - b^0)} \upd t \wedge \upd r \wedge \dVol_{(\mathbb{S}^2, G)}
\end{equation*}
and that
\begin{equation*}
 (\upd t)_\mu := -\frac{1}{2}\Lbar^0 L_\mu -\frac{1}{2}L^0 \Lbar_\mu + \slashed{\nabla}_\mu t
\end{equation*}

\end{proof}
% 
% \begin{proposition}[The weighted $T$ energy boundary terms]
% The boundary terms associated with the weighted $T$ energy (i.e.\ the energy associated with the vector field $wT$) are exactly the same as those appearing in the $T$ energy, but multiplied by the weight $w$.

% \end{proposition}

\begin{proposition}[The weighted Morawetz boundary terms]
\label{proposition boundary weighted Morawetz}
 The boundary term in the weighted Morawetz energy estimate corresponding to the modified energy current ${^{(wR)}\tilde{J}}[\phi]$ is given by

 \begin{equation}
  \begin{split}
   \int_{^t\Sigma_\tau} \imath_{{^{(wR)}\tilde{J}}} \dVol_g &=
   \int_{^t\Sigma_\tau} \frac{1}{2} w\Bigg( f_R|\slashed{\D}_L\phi|^2 
   - f_R  |\slashed{\nabla}\phi|^2
   + 4r^{-1} f_R \phi\cdot(\slashed{\D}_L\phi) \\
   & \phantom{=\int_{^t\Sigma_\tau} \frac{1}{2} \Bigg(}
   + 2(r^{-2}f_R - r^{-1}f'_R - r^{-1}f_R \partial_r(\log w) )|\phi|^2 \Bigg) \Omega^2 \upd r \wedge \dVol_{\mathbb{S}^2}
  \end{split}
 \end{equation}

 In addition, the boundary term on a surface of constant $t$ in the Morawetz estimate is given by
 \begin{equation}
    \begin{split}
      \int_{\{t = t_0\}} \imath_{{^{(R)}\tilde{J}}} \dVol_g = \int_{\{t = t_0\}}\frac{2}{(L^0 + \Lbar^0 - b^0)} w\textit{Err}_{(wR,t_0 \text{-bdy})} \Omega^2 \upd r \wedge \dVol_{\mathbb{S}^2}
    \end{split}
 \end{equation}
 where the error term satisfies
 \begin{equation}
  \begin{split}
    |\textit{Err}_{(wR,t_0 \text{-bdy})}| &\lesssim |f_R| \left(1 + |(L_{(\text{small})})^0| + |(\Lbar_{(\text{small})})^0| + |\slashed{\nabla} t| \right) \Big( |\slashed{\D}_L\phi|^2 +|\slashed{\D}_{\Lbar}\phi|^2 + |\slashed{\nabla}\phi|^2 \\
    & \phantom{\lesssim} + \left(r^{-1} |f'_R||f_R|^{-1} + r^{-2} + r^{-1}f_R \partial_r(\log w)\right)|\phi| \Big)
  \end{split}
 \end{equation}

\end{proposition}

\begin{proof}
 Recall that the modified energy current associated to the Morawetz vector field is
 \begin{equation*}
  {^{(wR)}\tilde{J}}^\mu[\phi] = {^{(wf_R R)}J^\mu}[\phi] + 2r^{-1}wf_R\phi\cdot \slashed{\D}^\mu \phi - \D^\mu \left( wr^{-1} f_R \right) |\phi|^2
 \end{equation*}
 We contract this with the one-form $(\upd \tau)_\mu$, and make use of proposition \ref{proposition volume form g}. A slightly involved but straightforward calculation leads to the expressions in the proposition.
\end{proof}

\begin{proposition}[The $p$-weighted boundary term]
\label{proposition boundary p}
 The boundary term in the $p$-weighted estimate, corresponding to the modified energy current ${^{(L,p)}\tilde{J}}[\phi]$ is given by
 \begin{equation}
 \label{equation p weighted boundary}
  \begin{split}
   &\int_{^t\Sigma_\tau} \imath_{{^{(L,p)}\tilde{J}}} \dVol_g =
   \int_{^t\Sigma_\tau} \bigg( f_L r^{p-2} (\slashed{\D}_L\psi)^2 
   - \frac{1}{2} r^{-2} \frac{\partial}{\partial r} \Big|_{\tau,\vartheta^1,\vartheta^2} \left( f_L r^{p+1} |\phi|^2 \right) \bigg) \Omega^2\upd r \wedge \dVol_{\mathbb{S}^2} \\ \\
   &= \int_{^t\Sigma_\tau} \left(
   f_L r^{p-2}(\slashed{\D}_L\psi)^2 
    + \textit{Err}_{(\text{L,p-bdy})} \right) \Omega^2 \upd r \wedge \dVol_{\mathbb{S}^2}  -  \int_{\bar{S}_{\tau,t}}\frac{1}{2} f_L r^{p-1} |\phi|^2 \Omega^2 \, \dVol_{\mathbb{S}^2} \\
   &\phantom{=} \\
  \end{split}
 \end{equation}
 Where the error term satisfies
 \begin{equation}
 \begin{split}
  |\textit{Err}_{(\text{L,p-bdy})} | &\lesssim
   \left|  \left.\frac{\partial}{\partial r}\right|_{\tau, \vartheta^1, \vartheta^2}\log \Omega - r^{-1} \right|r^{p-1} |\phi|^2 
  \end{split}
 \end{equation}
\end{proposition}
% 
% Finally, the boundary term on a surface of constant $r$ is
% \begin{equation}
%  \begin{split}
%     &\int_{C_r} \imath_{{^{(L,p)}\tilde{J}}} \dVol_g \\
%     &= \int_{C_r} \frac{1}{4}\mu \bigg( f_L r^{p} |\slashed{\D}_L \phi|^2 - f_L r^p |\slashed{\nabla}\phi|^2 + f_L r^{p-1} \phi \cdot \left( \slashed{\D}_L \phi - \slashed{\D}_{\Lbar}\phi \right) \\
%     &\phantom{= \int_{C_r} }
%     - \frac{1}{2}\left( (p-1)f_L r^{p-2} + r^{p-1}f_L \right) |\phi|^2 \bigg) \upd \tau \wedge \upd \vartheta^1 \wedge \upd \vartheta^2
%  \end{split}
% \end{equation}

\begin{proof}
 Recall that the modified energy current ${^{(L,p)}\tilde{J}}[\phi]$ is given by
 \begin{equation*}
 ({^{(L,p)}\tilde{J}}[\phi])^\mu = ({^{(f_L r^p L)}}J)^\mu + r^{p-1} f_L \phi \D^\mu \phi - \frac{1}{2} \partial^\mu\left( r^{p-1} f_L \right) \phi^2
 \end{equation*}
A long computation, beginning by taking the inner product of this energy current with the one-form $\upd \tau$, leads to the first expression in \eqref{equation p weighted boundary}. To obtain the second expression, we integrate by parts, picking up a term from the derivative of the quantity $\Omega$.
% 
% To compute the boundary term on a surface of constant $t$ we must make use of the fact that
% \begin{equation*}
%  \left.\frac{\partial}{\partial r}\right|_{t, \vartheta^1, \vartheta^2} = \frac{1}{(L^0 + \Lbar^0 - b^0)} \left( (\Lbar^0 - b^0)L - L^0 \Lbar + b^A X_A \right)
% \end{equation*}

\end{proof}

\section{The (modified) energy identity}

The fundamental identity which we will use to obtain control over $L^2$-type quantities is the \emph{energy identity} (suitably modified by lower order terms). 

\begin{proposition}[The (modified) energy identity]
\label{proposition energy identity}
 Let $Z$, $f_Z$ be a smooth vector field and a smooth function on the manifold $\mathcal{M}$. Let $\mathcal{N}\subset\mathcal{M}$ be any precompact subset of $\mathcal{M}$, and assume that $\phi$ is a smooth function on $\mathcal{M}$. Then the following energy identity is a consequence of the divergence theorem and proposition \ref{proposition modified compatible current identity}
 \begin{equation}
 \begin{split}
    &\int_{\mathcal{N}} \left(\Div {^{(Z,f_Z)}\tilde{J}}[\phi] \right) \dVol_{g} \\
    &= \int_{ \mathcal{N}} \bigg( 
      {^{(Z,\, f_Z)}\tilde{K}}[\phi] - \omega (\slashed{\D}_{\Lbar} \phi)\cdot \left( (\slashed{\D}_Z \phi) + \frac{1}{2} f_Z \phi \right) + \left( \tilde{\slashed{\Box}}_g \phi \right) \cdot \left((\slashed{\D}_Z\phi) + \frac{1}{2}f_Z \phi \right) \\
      &\phantom{=\int_{ \mathcal{N}} \bigg(} + \left( [\slashed{\D}_\mu \, , \, \slashed{\D}_\nu]\phi \right) \cdot (\slashed{\D}^\mu \phi) Z^\nu 
	  \bigg) \dVol_g \\
 	&= \int_{\partial \mathcal{N}} \imath_{{^{(Z,f_Z)}\tilde{J}[\phi]}} \dVol_g
 \end{split}
 \end{equation}
\end{proposition}

\begin{remark}
 The smoothness conditions on $Z$ and $f_Z$ can obviously be relaxed, by taking smooth approximations to some non-smooth vector field and function. In addition, the precompactness condition can also be dropped if we already know that $\phi$ is in a suitable Sobolev space, and if we have suitable bounds on the vector field $Z$ and the function $f_Z$ so that the integrals in the above expression can be bounded by the Sobolev norm in question.
\end{remark}

\begin{proposition}[The (modified) energy identity after a point-dependent change of basis]
\label{proposition energy identity after change of basis}
 Let $Z$, $f_Z$ be a smooth vector field and a smooth function on the manifold $\mathcal{M}$. Let $\mathcal{N}\subset\mathcal{M}$ be any precompact subset of $\mathcal{M}$, and assume that $\phi$ is a smooth function on $\mathcal{M}$. 

 Let $\phi_{(a)}$ be a collection of scalar fields labelled by the index $(a)$, and let $M_{(A)}^{\phantom{(A)}(a)}$ be a (possibly point-dependent) change-of-basis matrix, i.e.\ a collection of scalar fields such that, at every point, the rank of the matrix $M_{(A)}^{\phantom{(A)}(a)}$ is equal to the number of scalar fields $\phi_{(a)}$.

 Then we have the following energy identity:
 \begin{equation}
 \begin{split}
      \int_{\mathcal{N}} \left(\Div {^{(Z,f_Z)}\tilde{J}}[\phi]_{(A)} \right) \dVol_{g} &= \int_{ \mathcal{N}} \bigg( 
      {^{(Z,\, f_Z)}\tilde{K}}[\phi]_{(A)} - \omega (\slashed{\D}_{\Lbar} \phi)_{(A)}\cdot \left( (\slashed{\D}_Z \phi)_{(A)} + \frac{1}{2} f_Z \phi_{(A)} \right)\\
      & \phantom{ = \int_{ \mathcal{N}} \bigg( }
      + \left( \tilde{\slashed{\Box}}_g \phi \right)_{(A)} \cdot \left((\slashed{\D}_Z\phi)_{(A)} + \frac{1}{2}f_Z \phi_{(A)} \right) \\
      & \phantom{ = \int_{ \mathcal{N}} \bigg( }
      + \left( [\slashed{\D}_\mu \, , \, \slashed{\D}_\nu]\phi \right)_{(A)} \cdot (\slashed{\D}^\mu \phi)_{(A)} Z^\nu 
      + \textit{Err}_{(\partial M)}[\phi]_{(A)}
      \bigg) \dVol_g \\
      &= \int_{\partial \mathcal{N}} \imath_{{^{(Z,f_Z)}\tilde{J}[\phi]_{(A)}}} \dVol_g
 \end{split}
 \end{equation}
\end{proposition}
where, as usual, indices outside parentheses or other delimiters are understood as being contracted \emph{after} the relevant differential operators have been applied. For example,
\begin{equation}
 \begin{split}
  {^{(Z,f_Z)}\tilde{J}}[\phi]_{(A)} &= 
  Z_\nu Q^{\mu\nu}[\phi]_{(A)} + \frac{1}{2}f_Z (g^{-1})^{\mu\nu} \phi_{(A)} \cdot (\slashed{\D}_\nu \phi)_{(A)} - \frac{1}{4}(g^{-1})^{\mu\nu}(\D_\nu f_Z)|\phi_{(A)}|^2 \\
  &= (\slashed{\D}_\mu \phi)_{(A)} \cdot (\slashed{\D}_Z \phi)_{(A)}
  - \frac{1}{2}Z_\mu \left( -(\slashed{\D}_L \phi)_{(A)} \cdot (\slashed{\D}_{\Lbar}\phi)_{(A)} + (\slashed{\nabla}\phi)_{(A)}\cdot (\slashed{\nabla}\phi)_{(A)} \right) \\
  &\phantom{=} + \frac{1}{2}f_Z (g^{-1})^{\mu\nu} \phi_{(A)} \cdot (\slashed{\D}_\nu \phi)_{(A)} - \frac{1}{4}(g^{-1})^{\mu\nu}(\D_\nu f_Z)|\phi_{(A)}|^2 \\
  &= M_{(A)}^{\phantom{(A)}(a)}M_{(A)}^{\phantom{(A)}(b)} \bigg(
  (\slashed{\D}_\mu \phi_{(a)}) \cdot (\slashed{\D}_Z \phi_{(b)})
  - \frac{1}{2}Z_\mu \left( -(\slashed{\D}_L \phi_{(a)})\cdot (\slashed{\D}_{\Lbar}\phi_{(b)}) + (\slashed{\nabla}\phi_{(a)})\cdot (\slashed{\nabla}\phi_{(b)}) \right) \\
  &\phantom{=} + \frac{1}{2}f_Z (g^{-1})^{\mu\nu} \phi_{(a)} \cdot (\slashed{\D}_\nu \phi_{(b)}) - \frac{1}{4}(g^{-1})^{\mu\nu}(\D_\nu f_Z) \phi_{(a)} \cdot \phi_{(b)} \bigg) \\
 \end{split}
\end{equation}

We use the schematic notation $\phi_{(\text{orig})}$ to mean the supremum over the values of the fields in the original basis, that is
\begin{equation}
 |\phi_{(\text{orig})}| := \sup_{(a)} |\phi_{(a)}|
\end{equation}
We extend this notation in the obvious ways for derivatives of the $\phi_{(a)}$, so, for example,
\begin{equation*}
 \begin{split}
  |\slashed{\D}\phi_{(\text{orig})}| &:= \sup_{(a)} |\slashed{\D}\phi_{(a)}| \\
  |M_{(A)}| &:= \sup_{(a)} |M_{(A)}^{\phantom{(A)}(a)}|
 \end{split}
\end{equation*}

In terms of this notation, the error term associated with the derivatives of $M$ satisfies
\begin{equation}
 \begin{split}
  \left| \textit{Err}_{(\partial M)}[\phi]_{(A)} \right|
  &\lesssim  
  |g(Z,L)| \bigg( |L M_{(A)}| |\slashed{\D}\phi_{(\text{orig})}| |\slashed{\D} \phi|_{(A)} 
  + |\Lbar M_{(A)}| |\slashed{\nabla}\phi_{(\text{orig})}| |\slashed{\nabla}\phi|_{(A)} \\
  &\phantom{\lesssim |g(Z,L)| \bigg(}
  + |\slashed{\nabla} M_{(A)}| |\slashed{\nabla}\phi_{(\text{orig})}| |\slashed{\D}\phi|_{(A)}
  + |\slashed{\nabla} M_{(A)}| |\slashed{\D}\phi_{(\text{orig})}| |\slashed{\nabla}\phi|_{(A)} \bigg) \\
  &\phantom{\lesssim}
  + |g(Z,\Lbar)| \bigg( |L M_{(A)}| |\slashed{\nabla}\phi_{(\text{orig})}| |\slashed{\nabla} \phi|_{(A)} 
  + |\Lbar M_{(A)}| |\slashed{\D}_L\phi_{(\text{orig})}| |\slashed{\D}_L\phi|_{(A)} \\
  &\phantom{\lesssim |g(Z,\Lbar)| \bigg(}
  + |\slashed{\nabla} M_{(A)}| |\slashed{\nabla}\phi_{(\text{orig})}| |\slashed{\D}_L\phi|_{(A)}
  + |\slashed{\nabla} M_{(A)}| |\slashed{\D}_L\phi_{(\text{orig})}| |\slashed{\nabla}\phi|_{(A)} \bigg) \\
  &\phantom{\lesssim}
  + |\slashed{\Pi}(Z)| \bigg( |LM_{(A)}| |\slashed{\D}\phi_{(\text{orig})}| |\slashed{\nabla}\phi|_{(A)} + |LM_{(A)}| |\slashed{\nabla}\phi_{(\text{orig})}| |\slashed{\D}\phi|_{(A)} \\
  &\phantom{\lesssim + |\slashed{\Pi}(Z)| \bigg(}
  + |\Lbar M_{(A)}| |\slashed{\D}_L\phi_{(\text{orig})}| |\slashed{\nabla}\phi|_{(A)} + |\Lbar M_{(A)}| |\slashed{\nabla}\phi_{(\text{orig})}| |\slashed{\D}_L\phi|_{(A)} \\
  &\phantom{\lesssim + |\slashed{\Pi}(Z)| \bigg(}
  + |\slashed{\nabla} M_{(A)}| |\slashed{\nabla}\phi_{(\text{orig})}| |\slashed{\nabla}\phi|_{(A)}
  + |\slashed{\nabla} M_{(A)}| |\slashed{\D}\phi_{(\text{orig})}| |\slashed{\D}_L \phi|_{(A)} \\
  &\phantom{\lesssim + |\slashed{\Pi}(Z)| \bigg(}
  + |\slashed{\nabla} M_{(A)}| |\slashed{\D}_L \phi_{(\text{orig})}| |\slashed{\D}\phi|_{(A)} \bigg) \\
  &\phantom{\lesssim} 
  +  |f_Z| \bigg( |LM_{(A)}| |\slashed{\D}\phi_{(\text{orig})}| |\phi|_{(A)}| + |LM_{(A)}| |\phi|_{(\text{orig})}| |\slashed{\D}\phi_{(\text{orig})}|_{(A)} \\
  &\phantom{\lesssim + |\slashed{\Pi}(Z)| \bigg(}
  +  |\Lbar M_{(A)}| |\slashed{\D}_L\phi_{(\text{orig})}| |\phi|_{(A)}| + |\Lbar M_{(A)}| |\phi|_{(\text{orig})}| |\slashed{\D}_L\phi_{(\text{orig})}|_{(A)} \\
  &\phantom{\lesssim + |\slashed{\Pi}(Z)| \bigg(}
  + |\slashed{\nabla}M_{(A)}| |\slashed{\nabla}\phi_{(\text{orig})}| |\phi|_{(A)}| + |\slashed{\nabla}M_{(A)}| |\phi|_{(\text{orig})}| |\slashed{\nabla}\phi_{(\text{orig})}|_{(A)} \bigg) \\
  &\phantom{\lesssim}
  + |\Lbar f_Z| |LM_{(A)}| |\phi_{(\text{orig})}| |\phi|_{(A)}
  + |L f_Z| |\Lbar M_{(A)}| |\phi_{(\text{orig})}| |\phi|_{(A)} \\
  &\phantom{\lesssim}
  + |\slashed{\nabla} f_Z| |\slashed{\nabla}M_{(A)}| |\phi_{(\text{orig})}| |\phi|_{(A)}
 \end{split}
\end{equation}

\chapter{The bootstrap}
\label{chapter bootstrap}

In this chapter we initiate the bootstrap, i.e.\ we state a series of bounds which we shall later improve. These bounds will be pointwise, or $L^{\infty}$ bounds for the lower derivatives, as well as $L^2$-based bounds for higher derivatives. For the sake of clarity, we will state the pointwise bounds for all relevant quantities, including those quantities which can be expressed in terms of other quantities - in other words, our pointwise bounds are not all independent of each other.

\section{Constants}

In order to state the bootstrap assumptions, we first need to define a whole set of large and small constants, and also establish the relationships between them.

We begin by defining a collections of constants:
\begin{equation}
  0 \ll C_{(0)} \ll C_{(1)} \ll C_{(2)} \ll \ldots
\end{equation}

We further ``subdivide'' the intervals between the constants by introducing another set of constants, with two numbered indices, satisfying (for all $n$)
\begin{equation}
  C_{(n-1)} \ll C_{(n, 1)} \ll C_{(n, 2)} \ll \ldots \ll C_{(n)}
\end{equation}

Finally, we introduce the very large constant $\mathring{C}$ which obeys
\begin{equation}
\mathring{C} \gg C_{(n)}
\end{equation}
for all $n$.

While the constants $C_{(n)}$ and $C_{(n,m)}$ are involved in the \emph{pointwise} bounds, we also need some constants in order to specify the $L^2$ bootstrap bounds. These will eventually be related to the constants $C_{(n)}$ and $C_{(n,m)}$, though they are not precisely the same. We introduce some set of constants
\begin{equation}
\begin{split}
0 \ll C_{[0]} \ll C_{[1]} \ll C_{[2]} \ll \ldots \\
C_{[n-1]} \ll C_{[n, 1]} \ll C_{[n, 2]} \ll \ldots \ll C_{[n]}
\end{split}
\end{equation}

In contrast, $\delta$ is some small constant, which will later be chosen suitably small, although the very small constant $\epsilon$ satisfies
\begin{equation*}
 \epsilon \ll \delta
\end{equation*}

Note that the notation $A \ll B$ means that, for all large constants $C$ which appear when comparing $A$ and $B$, we have $CA < B$. In particular, we have $C_{(n,m)}\epsilon < \delta$.

Finally, $\beta$ is a fixed constant in the range
\begin{equation*}
 0 < \beta < \frac{1}{2}
\end{equation*}
We should think of the constants $\beta$, $\delta$ and $\epsilon$ as satisfying the relationship
\begin{equation*}
 \epsilon \ll \delta \ll \beta
\end{equation*}

\section{Pointwise bootstrap bounds}
\label{section pointwise bootstrap}

We shall assume that the following pointwise bounds hold at all times $\tau$ satisfying $\tau_0 \leq \tau \leq \tau_{\text{max}}$.

Pointwise bounds on the fields: for all $n \leq M_1$, for all $\phi_{(A)}$, in the region $r \geq  \frac{1}{2}r_0$ we have
\begin{equation}
 \label{equation bootstrap fields}
 |\mathscr{Y}^n \phi_{(A)}| \leq \epsilon (1+r)^{-\frac{1}{2} + \delta}
\end{equation}

We also have the following pointwise bounds on the fields and their first derivatives, yielding improved decay in $r$ at the expense of worse behaviour at large $\tau$:
\begin{equation}
\label{equation bootstrap improved bound fields}
|\phi_{(A)}| + |\mathscr{Y} \phi_{(A)}| \leq \epsilon (1+r)^{-1 + C_{[N_1]}\epsilon}(1+\tau)^{C{(N_1)}\delta}
\end{equation}

Additionally, we have the following pointwise bounds on the fields and their first derivatives, giving decay in $\tau$ and a better constant for the lowest order quantities: for all $\phi_{(A)}$ and for all values of $r$ we have
\begin{equation}
\begin{split}
\label{equation bootstrap fields no commute}
|\phi_{(A)}| &\leq \epsilon^3 (1+r)^{-\frac{1}{2} + \delta}(1+\tau)^{-\beta} \\
|\mathscr{Y}\phi_{(A)}| &\leq \epsilon^5 (1+r)^{-\frac{1}{2} + \delta}(1+\tau)^{-\beta} \\
|\partial \phi_{(A)}| &\leq \epsilon^5 (1+r)^{-1 + C_{(0,m)}}(1+\tau)^{-\beta} \\
|\partial \mathscr{Y} \phi_{(A)}| &\leq \epsilon^5 (1+r)^{-1 + C_{(1,m)}} \\
\end{split}
\end{equation}
Note that, at the lowest order, we have even smaller quantities (of order $\epsilon^3$ or $\epsilon^5$ rather than $\epsilon$).

Pointwise bounds on the derivatives of the fields: for all $n \leq M_1$ and for any field $\phi_{(A)}$ in the region $r \geq \frac{1}{2}r_0$ we have
\begin{equation}
\label{equation bootstrap bad derivs of fields}
|\slashed{\D} \mathscr{Y}^n \phi_{(A)}| \leq \epsilon\left( (1+r)^{-1} + (1+r)^{-1 + \delta}(1+\tau)^{-\beta}\right) 
\end{equation}

We also have the following bounds on the derivatives of the fields, giving more detailed information about decay in $r$: for all $n \leq M_1$ and for $\phi_{(A)} \in \Phi_{[m]}$, in the region $r \geq \frac{1}{2}r_0$ we have
\begin{equation}
\label{equation bootstrap bad derivs of fields strong r}
|\slashed{\D} \mathscr{Y}^n \phi_{(A)}| \leq \epsilon (1+r)^{-1 + C_{(n,m)}\epsilon}
\end{equation}

Additionally, we have the following pointwise bounds on the derivatives of the fields, giving decay also in $\tau$ for the lowest order quantities: for $\phi_{(A)} \in \Phi_{[m]}$ and for all values of $r$ we have
\begin{equation}
\label{equation bootstrap derivs of fields no commute}
|\partial \phi_{(A)}| \leq \epsilon (1+r)^{-1 + C_{(0,m)}\epsilon}(1+\tau)^{-\beta}
\end{equation}
In particular, the bad derivatives of the ``good fields'' obey the following bounds: for $\phi_{(A)} \in \Phi_{[0]}$, in the region $r \geq \frac{1}{2}r_0$
\begin{equation}
\label{equation bootstrap bad derivs of good fields}
|\partial \phi_{(A)}| \leq \epsilon\left( (1+r)^{-1-\delta} + (1+r)^{-1}(1+\tau)^{-\beta} \right)
\end{equation}

We also have the following pointwise bounds on the ``good'' derivatives of the fields: for all $n \leq M_1$, in the region $r \geq \frac{1}{2}r_0$
\begin{equation}
\label{equation bootstrap good derivs of fields}
|\overline{\slashed{\D}} \mathscr{Y}^n \phi_{(A)}| \leq \epsilon  (1+r)^{-1 - \delta}
\end{equation}

Again, for the lowest order derivatives we have an improved estimate which gives decay also in $\tau$: for all fields $\phi_{(A)}$, for all values of $r$ we have
\begin{equation}
\label{equation bootstrap good derivs of fields no commute}
|\overline{\slashed{\D}} \phi_{(A)}| \leq \epsilon  (1+r)^{-1 - \delta}(1+\tau)^{-\beta}
\end{equation}

Pointwise bounds on the derivatives of the fields in the region $r \leq \frac{1}{2}r_0$: for all $n \leq M_1$
\begin{equation}
 \label{equation bootstrap derivs of fields r<r0}
 |\partial^{n + 1} \phi_{(A)}| \leq \epsilon
\end{equation}

%Pointwise bounds on the derivatives of the ``good fields'' in the region $\frac{1}{2}r_0 \leq r \leq r_0$: for all $n \leq M_1$ and for $\phi_{(A)} \in \Phi_{[0]}$
%\begin{equation}
%	\label{equation bootstrap derivs of good fields r<r0}
%	|\Lbar \phi_{(A)}| \leq \epsilon(1+\tau)^{-\beta}
%\end{equation}

Pointwise bounds on the metric components: for all $n \leq M_1$, in the region $r \geq \frac{1}{2}r_0$
\begin{equation}
 \begin{split}
  \label{equation bootstrap metric}
  |\mathscr{Y}^n h_{ab}| &\leq \epsilon (1+r)^{-\frac{1}{2} + \delta} \\
  |\mathscr{Y}^n h|_{(\text{frame})} &\leq \epsilon (1+r)^{-\frac{1}{2} + \delta} 
 \end{split}
\end{equation}

Pointwise bounds on the derivatives of the metric components: for all $n \leq M_1$, in the region $r \geq \frac{1}{2}r_0$
\begin{equation}
\begin{split}
\label{equation bootstrap bad derivs of metric}
|\slashed{\D} \mathscr{Y}^n h_{ab}| &\leq \epsilon \left( (1+r)^{-1} + (1+r)^{-1+\delta}(1+\tau)^{-\beta} \right) \\
|\slashed{\D} \mathscr{Y}^n h|_{(\text{frame})} &\leq \epsilon \left( (1+r)^{-1} + (1+r)^{-1 + \delta}(1+\tau)^{-\beta} \right) \\
\end{split}
\end{equation}

We also have the following estimates, giving stronger control over the behaviour of the derivatives in the $r$ direction: in the region $r \geq \frac{1}{2}r_0$, for all $n \leq M_1$
\begin{equation}
\begin{split}
\label{equation bootstrap bad derivs of metric strong r}
|\slashed{\D} \mathscr{Y}^n h_{ab}| &\leq \epsilon (1+r)^{-1 + C_{(n)}\epsilon} \\
|\slashed{\D} \mathscr{Y}^n h|_{(\text{frame})} &\leq \epsilon (1+r)^{-1 + C_{(n)}\epsilon}
\end{split}
\end{equation}

Additionally, we have the following estimates on the derivatives of the ``good'' components of the metric:
\begin{equation}
\label{equation bootstrap bad derivs of good components of metric}
|\partial h|_{LL} \leq \epsilon \left( (1+r)^{-1-\delta} + \epsilon(1+r)^{-1}(1+\tau)^{-\beta} \right)
\end{equation}

%\begin{equation}
% \begin{split}
%  \label{equation bootstrap bad derivs of metric high order}
%  |\slashed{\D} \mathscr{Z}^n h_{ab}| &\leq \epsilon \left( (1+r)^{-1-\delta} + (1+r)^{-1 + \delta}(1+\tau)^{-\beta} \right) \\
%  |\slashed{\D} \mathscr{Z}^n h|_{(\text{frame})} &\leq \epsilon \left( (1+r)^{-1-\delta} + (1+r)^{-1 + \delta}(1+\tau)^{-\beta} \right) \\
% \end{split}
%\end{equation}
%
%Detailed pointwise bounds on the derivatives of the metric components: for all $m \leq M_1 - 3$

%we also have, for all $m \leq M_1 -4$,
%\begin{equation}
% \begin{split}
%  \label{equation bootstrap DT bad derivs of metric low order}
%  |\slashed{\D} \slashed{\D}_T \mathscr{Z}^n h_{ab}| &\leq \epsilon \left( (1+r)^{-1-\delta} + (1+r)^{-1 + C_{(n,T)}\epsilon}(1+\tau)^{-\beta} \right) \\
%  |\slashed{\D} \slashed{\D}_T \mathscr{Z}^n h|_{(\text{frame})} &\leq \epsilon \left( (1+r)^{-1-\delta} + (1+r)^{-1 + C_{(n,T)}\epsilon}(1+\tau)^{-\beta} \right) \\
% \end{split}
%\end{equation}

Pointwise bounds on the derivatives of the metric components in the region $r \leq \frac{1}{2}r_0$: for all $m \leq M_1 - 1$
\begin{equation}
 \begin{split}
  \label{equation bootstrap derivs metric r<r0}
  |\partial^{m + 1} h_{ab}| \leq \epsilon \\
  |\partial^{m + 1} h|_{(\text{frame})} \leq \epsilon 
 \end{split}
\end{equation}
%
%Additionally, in the region the region $\frac{1}{2}r_0 \leq r \leq r_0$ we have
%\begin{equation}
%	\label{equation bootstrap derivs good component of metric r<r_0}
%	|\Lbar h|_{LL} \leq \epsilon (1+\tau)^{-\beta}
%\end{equation}

Pointwise bounds on the ``good'' derivatives of the metric components: for all $m \leq M_1$, in the region $r \geq \frac{1}{2}r_0$
\begin{equation}
 \begin{split}
  \label{equation bootstrap good derivs of metric}
  |\overline{\slashed{\D}} \mathscr{Y}^m h_{ab}| &\leq \epsilon  (1+r)^{-1 - \delta} \\
  |\overline{\slashed{\D}} \mathscr{Y}^m h|_{(\text{frame})} &\leq \epsilon  (1+r)^{-1 - \delta} \\
 \end{split}
\end{equation}
%
%Pointwise bounds on bad derivatives of the ``good metric component'': in the region $r \geq \frac{1}{2}r_0$
%\begin{equation}
%\label{equation bootstrap bad derivs of good metric component}
% |\partial h|_{LL} \leq \epsilon\left( (1+r)^{-1-\delta} + (1+r)^{-1}(1+\tau)^{-\beta}\right)
%\end{equation}

Pointwise bounds on the foliation density:
\begin{equation*}
\label{equation bootstrap foliation density}
	\begin{split}
	|\mu| &\leq (1+\epsilon)(1+r)^{C_{(0)}\epsilon} \\
  	|\mu^{-1}| &\leq (1+\epsilon)(1+r)^{C_{(0)}\epsilon}
\end{split}
\end{equation*}
Note that these bounds are to hold in both the regions $r \leq r_0$ and $r \geq r_0$.

We also have the following bounds on the derivatives of the foliation density: for all $0 \leq n \leq N_1$,
\begin{equation}
\begin{split}
\label{equation bootstrap derivs foliation density}
|\mathscr{Y}^n \log\mu| &\leq \epsilon \left( 1 + (1+r)^{\delta}(1+\tau)^{-\beta} \right) \\
|r\slashed{\nabla} \mathscr{Z}^{n-1} \log\mu| &\leq \epsilon \left( r^{-1} + (1+r)^{-1+\delta}(1+\tau)^{-\beta} \right) \\
\end{split}
\end{equation}

Again, we can give the following bounds on derivatives of the foliation density, which give stronger control in $r$: for all $0 \leq n \leq M_1$
\begin{equation}
\label{equation bootstrap derivs foliation density strong r}
|\mathscr{Y}^n \log\mu| \leq \epsilon (1+r)^{C_{(n)}\epsilon}
\end{equation}

%\begin{equation}
% \label{equation bootstrap foliation density}
% \begin{split}
%  |\mu| &\leq (1+\epsilon)(1+r)^{C_{(0)}\epsilon} \\
%  |\mu^{-1}| &\leq (1+\epsilon)(1+r)^{C_{(0)}\epsilon}
% \end{split}
%\end{equation}
%
%%we also need the bound
%%\begin{equation}
%%\label{equation bootstrap first derivs foliation density}
%% |r\slashed{\nabla} \log \mu| \leq \epsilon \left( 1 + (1+r)^{C_{(1)}\epsilon}(1+\tau)^{-\beta} \right)
%%\end{equation}
%
%additionally, 

Pointwise bounds on the connection coefficients: $\omega$ satisfies the ``zero-th order'' bounds:
\begin{equation}
 \label{equation bootstrap omega}
  |\omega| \leq \epsilon \left( (1+r)^{-1-\delta} + (1+r)^{-1}(1+\tau)^{-\beta} \right)
\end{equation}
additionally, for all $0 \leq n \leq M_1$, the ``bad'' connection coefficients satisfy the bounds
\begin{equation}
 \label{equation bootstrap bad connection components}
 \begin{split}
  |\mathscr{Y}^n \omega| &\leq \epsilon \left( (1+r)^{-1} + (1+r)^{-1 + \delta}(1+\tau)^{-\beta} \right) \\
  |\mathscr{Y}^n \zeta| &\leq \epsilon \left( (1+r)^{-1} + (1+r)^{-1 + \delta}(1+\tau)^{-\beta} \right) \\
  |\mathscr{Y}^n \tr_{\slashed{g}} \chibar_{(\text{small})}| &\leq \epsilon \left( (1+r)^{-1} + (1+r)^{-1 + \delta}(1+\tau)^{-\beta} \right) \\
  |\mathscr{Y}^n \hat{\chibar} | &\leq \epsilon \left( (1+r)^{-1} + (1+r)^{-1 + \delta}(1+\tau)^{-\beta} \right)
 \end{split}
\end{equation}
as well as the following bounds, which give stronger bounds in $r$:
\begin{equation}
\label{equation bootstrap bad connection components strong r}
\begin{split}
|\mathscr{Y}^n \omega| &\leq \epsilon (1+r)^{-1+ C_{(n)}\epsilon}\\
|\mathscr{Y}^n \zeta| &\leq \epsilon (1+r)^{-1+ C_{(n)}\epsilon} \\
|\mathscr{Y}^n \tr_{\slashed{g}} \chibar_{(\text{small})}| &\leq \epsilon (1+r)^{-1+ C_{(n)}\epsilon} \\
|\mathscr{Y}^n \hat{\chibar} | &\leq \epsilon (1+r)^{-1+ C_{(n)}\epsilon}
\end{split}
\end{equation}

%\begin{equation}
% \label{equation bootstrap bad connection components highest regulariy}
% \begin{split}
%  |\mathscr{Z}^n \omega| &\leq \epsilon \left( (1+r)^{-1} + (1+r)^{-1 + \delta}(1+\tau)^{-\beta} \right) \\
%  |\mathscr{Z}^n \zeta| &\leq \epsilon \left( (1+r)^{-1} + (1+r)^{-1 + \delta}(1+\tau)^{-\beta} \right)
% \end{split}
%\end{equation}

%and for all $0 \leq n \leq M_1 - 3$
%\begin{equation}
% \label{equation bootstrap bad connection components low regulariy}
% \begin{split}
%  |\mathscr{Z}^n \tr_{\slashed{g}} \chibar_{(\text{small})}| & \leq \epsilon \left( (1+r)^{-1} + (1+r)^{-1 + C_{(n+1)}\epsilon}(1+\tau)^{-\beta} \right) \\
%  |\mathscr{Z}^n \hat{\chibar} | &\leq \epsilon (1+r)^{-1 + C_{(n+1)}\epsilon}
% \end{split}
%\end{equation}

Additionally, the ``good'' connection coefficients obey the following estimates: for all $n \leq M_1$,
\begin{equation}
 \label{equation bootstrap good connection components}
 \begin{split}
  |\mathscr{Y}^n \tr_{\slashed{g}} \chi_{(\text{small})}| & \leq \epsilon (1+r)^{-1 - \delta} \\
  |\mathscr{Y}^n \hat{\chi} | &\leq \epsilon (1+r)^{-1 - \delta}
 \end{split}
\end{equation}

Finally, certain combinations of connection coefficients satisfy the following bounds, giving uniform decay in $\tau$:
\begin{equation*}
\begin{split}
	|\slashed{\nabla} \log \mu| &\lesssim (1+r)^{-1+C_{(1)}\epsilon} (1+\tau)^{-C^* \delta} \\
	|\zeta| &\lesssim (1+r)^{-1+C_{(0)}\epsilon} (1+\tau)^{-C^* \delta} \\
	|\chi_{(\text{small})} + \chibar_{(\text{small})}| &\lesssim (1+r)^{-1+C_{(0)}\epsilon} (1+\tau)^{-C^* \delta}
\end{split}
\end{equation*}
Note that these are the only quantities related to the metric which are required to decay \emph{uniformly} in $\tau$. Note, however, that the uniform decay in $\tau$ is at a very slow rate: we only require decay like $\tau^{-C^* \delta}$.

Note that all of the null frame connection coefficients are considered only in the region $r \geq r_0$; in the region $r \leq r_0$ we deal directly with derivatives of the rectangular components of $h$.

Pointwise bounds on the rectangular components of the frame fields: in the region $r \geq \frac{1}{2}r_0$ we have
\begin{equation}
\label{equation bootstrap rectangular components zeroth order}
\begin{split}
	|L^a_{(\text{small})}| &\leq \epsilon \\
	|\Lbar^a_{(\text{small})}| &\leq \epsilon \\
	|\slashed{\Pi}_\mu^{\phantom{\mu} a}| &\leq 1 + \epsilon
\end{split}
\end{equation}
and for all $n \leq M_1$
%\begin{equation}
%\label{equation bootstrap rectangular components}
%  \begin{split}
%    |\mathscr{Z}^n L^a| &\leq (1+r)^{-\delta} + (1+r)^{\delta}(1+\tau)^{-\beta} \\
%    |\mathscr{Z}^n\Lbar^a| &\leq (1+r)^{-\delta} + (1+r)^{\delta}(1+\tau)^{-\beta} \\
%    |\mathscr{Z}^n \slashed{\Pi}_\mu^{\phantom{\mu} a}| &\leq (1+r)^{-\delta} + (1+r)^{\delta}(1+\tau)^{-\beta} \\
%    |\mathscr{Z}^n L^a_{(\text{small})}| &\leq \epsilon(1+r)^{-\delta} + \epsilon(1+r)^{\delta}(1+\tau)^{-\beta} \\
%	|\mathscr{Z}^n\Lbar^a_{(\text{small})}| &\leq \epsilon(1+r)^{-\delta} + \epsilon(1+r)^{\delta}(1+\tau)^{-\beta} \\
%  \end{split}
%\end{equation}
%We also have the following bounds, giving improved bounds in $r$ on certain components:
\begin{equation}
\label{equation bootstrap rectangular components strong r}
	\begin{split}
	|\mathscr{Y}^n L^a| &\leq 2(1+r)^{C_{(n)}\epsilon} \\
	|\mathscr{Y}^n\Lbar^a| &\leq 2(1+r)^{C_{(n)}\epsilon} \\
	|\mathscr{Y}^n \slashed{\Pi}_\mu^{\phantom{\mu} a}| &\leq 2(1+r)^{C_{(n)}\epsilon} \\
	|\mathscr{Y}^n L^a_{(\text{small})}| &\leq \epsilon(1+r)^{C_{(n)}\epsilon} \\
	|\mathscr{Y}^n\Lbar^a_{(\text{small})}| &\leq \epsilon(1+r)^{C_{(n)}\epsilon} \\
	\end{split}
\end{equation}
and the following bounds on the ``good'' rectangular components of the frame fields: for $n \leq M_1$
\begin{equation}
\label{equation bootstrap Li}
	|\mathscr{Y}^n L^i_{(\text{small})}| \leq \epsilon(1+r)^{-\delta}
\end{equation}
Recall that the fields $L$ and $\Lbar$ can be expressed directly in terms of the metric components in the region $r \leq r_0$.

Pointwise bounds on other geometric quantities: for $n \leq M_1$, in the region $r \geq r_0$
\begin{equation}
 \label{equation bootstrap b}
 |\mathscr{Y}^n b| \leq \epsilon + \epsilon(1+r)^{C_{(n)}\epsilon}(1+\tau)^{-\beta} \\
\end{equation}

%For $n \leq M_1$
%\begin{equation}
% \label{equation bootstrap K}
% \left|\mathscr{Z}^n \left( r^2 K - 1 \right) \right| \leq \epsilon
%\end{equation}
We also have
\begin{equation}
 \begin{split}
  \left| \sqrt{-\det g} - 1 \right| \leq \epsilon
  \left|r^2 K - 1\right| \leq \epsilon
 \end{split}
\end{equation}

For $n \leq M_1$
\begin{equation}
 \left|\mathscr{Y}^n \left(r^{-2}\slashed{g}_{\mu\nu} - \mathring{\gamma}_{\mu\nu} \right)\right| \lesssim \epsilon
\end{equation}

The scalar $\Omega$ obeys
\begin{equation}
 \begin{split}
  1 - \epsilon \leq \frac{\Omega}{r} &\leq 1 + \epsilon \\
 \end{split}
\end{equation}

The rectangular coordinate $t$ can be expressed in terms of the geometric coordinates, in terms of which it satisfies
\begin{equation}
 |t - \tau - r| \lesssim \epsilon r
\end{equation}

\begin{remark}[Long range, non-decaying perturbations]
	We remark here that we have allowed many quantities, including derivatives of the metric components, to include a part which decays in $r$ but not in $\tau$. In fact, our bootstrap assumptions are such that, at \emph{fixed} $r$, no geometric quantity is assumed \emph{a priori} to decay in $\tau$, except for a few specific combinations of connection coefficients.
	
%	, with the exception of the quantity $(\Lbar h)_{LL}$. This is required to decay in time in a neighbourhood of $r = r_0$, but still, it might not decay uniformly in $\tau$ away from this region. Note that this is required specifically to ensure that the quantity $\lim_{r\rightarrow\infty} r(\Lbar h)_{\Lbar\Lbar}$ tends to zero as $\tau \rightarrow \infty$.
	
	Note also that many of the geometric perturbations are allowed to include ``long ranged'' and \emph{non decaying} (in $\tau$) parts. That is, the non-decaying in $\tau$ part is allowed to decay at the rate $r^{-1}$. On the other hand, the decaying (in $\tau$) parts of many geometric quantities quantities (for example, $\zeta$ or $\tr_{\slashed{g}}\chibar_{(\text{small})}$) is allowed to be ``super long ranged'', that is, to decay at a rate even slower than $r^{-1}$.
	
	In fact, we will recover bounds which include decay in $\tau$. Thus, for example, we will eventually recover a bound of the form
	\begin{equation*}
		|\partial \phi_{(A)}| \leq \frac{1}{2}\epsilon (1+r)^{-1+ C_{(0)}\epsilon} (1+\tau)^{-\beta}
	\end{equation*}
	from which it follows trivially that
	\begin{equation*}
	|\partial \phi_{(A)}| \leq \frac{1}{2}\epsilon \left( (1+r)^{-1} + (1+r)^{-1+C_{(0)}\epsilon} (1+\tau)^{-\beta} \right)
	\end{equation*}
	improving the bootstrap bound in question. However, we retain the non-decaying terms in our bootstrap argument, partly because our methods can accommodate them without serious difficulty, and partly because there are important problems for which a non-decaying source term might be present (in particular, a long-ranged source: see section \ref{section intro slow decay} of the introduction).
		
	Finally, note that all of the assumptions involving the operator $r\slashed{\D}_L$ applied to some field are required \emph{only} so that we can control the term $\slashed{\nabla}\log\mu (\slashed{\D}_L (r\slashed{\D}_L \phi)$ appearing in the commutator with $r\slashed{\nabla}$, and the additional error terms appearing when commuting with $r\slashed{\D}_L$. If we could find another way to control this single term, then it would be unnecessary to commute with $r\slashed{\D}_L$, and consequently we could replace the bootstrap bounds with the corresponding bounds after having replaced $\mathscr{Y}$ with $\mathscr{Z}$. Note that this is in fact the case for the Einstein equations in harmonic coordinates.
	
	%If the equation for $\log \mu$ is such that, even after commuting with the operator $r\slashed{\nabla}$, we can prove some better decay (in $r$, say) for this quantity, then we may not need to commute with $r\slashed{\D}_L$, and consequently we would not need to make the bootstrap assumptions involving this operator. This is the case, for example, for the Einstein equations in harmonic coordinates, where the harmonic coordinate condition implies that $L\log\mu \sim (\bar{\partial}h)$, instead of $L\log\mu \sim (\bar{\partial}h) + (\partial h)_{(LL)}$ as in the general case.
	
\end{remark}

\section{\texorpdfstring{$L^2$}{L2} based bootstrap bounds}
\label{section L2 bootstrap bounds}

$L^2$ bootstrap assumptions for the fields: for $\phi_{(a)} \in \Phi_{[m]}$, and for $n \leq N_2$, and for any $\tau_1 \geq \tau \geq \tau_0$ we assume the bootstrap bounds
\begin{equation}
\begin{split}
	\int_{\Sigma_\tau} 
		(1+r)^{-C_{[n,m]}} |\overline{\slashed{\D}} \mathscr{Y}^n\phi_{(a)}|^2
	\dVol_{\Sigma_\tau}
	+ \int_{\Sigma_\tau \cap \{r \leq r_0\}} 
	|\slashed{\D} \mathscr{Y}^n\phi_{(a)}|^2
	\dVol_{\Sigma_\tau}
	&\leq \epsilon^{2(N_2 + 2 - n)} (1+\tau)^{-1 + C_{(n)}\delta}
	\\
	\int_{\Sigma_\tau \cap \{r \geq r_0\}} 
		\delta (1+r)^{\frac{1}{2}\delta} |\slashed{\D}_L \mathscr{Y}^n\phi_{(a)}|^2
	\dVol_{\Sigma_\tau}
	&\leq \epsilon^{2(N_2 + 2 - n)} (1+\tau)^{-1 + C_{(n)}\delta}
	\\
	\int_{\Sigma_\tau \cap \{r \geq r_0\}} 
		\delta (1+r)^{1-C_{[n,m]}} |\slashed{\D}_L (r\mathscr{Y}^n\phi_{(a)})|^2
	r^{-2}\dVol_{\Sigma_\tau}
	&\leq \epsilon^{2(N_2 + 2 - n)} (1+\tau)^{C_{(n)}\delta}
	\\
	\int_{\mathcal{M}_\tau^{\tau_1}} 
		C_{[n,m]}\epsilon (1+r)^{-1-C_{[n,m]}} |\slashed{\D} \mathscr{Y}^n\phi_{(a)}|^2
	\dVol_{\Sigma_\tau}
	&\leq \epsilon^{2(N_2 + 2 - n)} (1+\tau)^{-1 + C_{(n)}\delta}
	\\
	\int_{\mathcal{M}_\tau^{\tau_1}} 
		\left( \delta^2 (1+r)^{-1-\delta} |\slashed{\D} \mathscr{Y}^n\phi_{(a)}|^2
		+ \delta^2 (1+r)^{-3-\delta} |\mathscr{Y}^n\phi_{(a)}|^2
	\right)\dVol_{g}
	&\leq \epsilon^{2(N_2 + 2 - n)} (1+\tau)^{-1 + C_{(n)}\delta}
	\\
	\int_{\mathcal{M}_\tau^{\tau_1}} 
		\delta (1+r)^{-1+ \frac{1}{2}\delta} |\overline{\slashed{\D}} \mathscr{Y}^n\phi_{(a)}|^2
	\dVol_{g}
	&\leq \epsilon^{2(N_2 + 2 - n)} (1+\tau)^{-1 + C_{(n)}\delta}
\end{split}
\end{equation}

Similarly, if $\phi_{(A)} = M_{(A)}^{(a)} \phi_{(a)}$, where $M_{(A)}^{(a)}$ is a change-of-basis matrix, then we assume that, if $\phi_{(A)} \in \bm{\Phi}_{[m]}$, then we have the bootstrap bounds
\begin{equation}
\begin{split}
	\int_{\Sigma_\tau} 
		(1+r)^{-C_{[n,m]}} |\overline{\slashed{\D}} \mathscr{Y}^n\phi|_{(A)}^2
	\dVol_{g}
	&\leq \epsilon^{2(N_2 + 2 - n)} (1+\tau)^{-1 + C_{(n)}\delta}
	\\
	\int_{\Sigma_\tau} 
		\delta (1+r)^{\frac{1}{2}\delta} |\slashed{\D}_L \mathscr{Y}^n\phi|_{(A)}^2
	\dVol_{g}
	&\leq \epsilon^{2(N_2 + 2 - n)} (1+\tau)^{-1 + C_{(n)}\delta}
	\\
	\int_{\Sigma_\tau} 
		\delta (1+r)^{1-C_{[n,m]}} |\slashed{\D}_L (r\mathscr{Y}^n\phi)|_{(A)}^2
	r^{-2} \dVol_{g}
	&\leq \epsilon^{2(N_2 + 2 - n)} (1+\tau)^{C_{(n)}\delta}
	\\
	\int_{\mathcal{M}_\tau^{\tau_1}} 
		C_{[n,m]}\epsilon (1+r)^{-1-C_{[n,m]}} |\slashed{\D} \mathscr{Y}^n\phi|_{(A)}^2
	\dVol_{g}
	&\leq \delta^{-1} \epsilon^{2(N_2 + 2 - n)} (1+\tau)^{-1 + C_{(n)}\delta}
	\\
	\int_{\mathcal{M}_\tau^{\tau_1}} 
		\delta (1+r)^{-1-\delta} |\slashed{\D} \mathscr{Y}^n\phi|_{(A)}^2
	\dVol_{g}
	&\leq \delta^{-1} \epsilon^{2(N_2 + 2 - n)} (1+\tau)^{-1 + C_{(n)}\delta}
	\\
	\int_{\mathcal{M}_\tau^{\tau_1}} 
		\delta (1+r)^{-1+ \frac{1}{2}\delta} |\overline{\slashed{\D}} \mathscr{Y}^n\phi|_{(A)}^2
	\dVol_{g}
	&\leq \delta^{-1} \epsilon^{2(N_2 + 2 - n)} (1+\tau)^{-1 + C_{(n)}\delta}
\end{split}
\end{equation}

We also assume the following $L^2$ bootstrap assumptions for the metric components: 
\begin{equation}
\begin{split}
	\int_{\Sigma_\tau} 
		(1+r)^{-C_{[n]}} |\overline{\slashed{\D}} \mathscr{Y}^n h_{(\text{rect})}|^2
	\dVol_{g}
	&\leq \epsilon^{2(N_2 + 2 - n)} (1+\tau)^{-1 + C_{(n)}\delta}
	\\
	\int_{\Sigma_\tau} 
		\delta (1+r)^{\frac{1}{2}\delta} |\slashed{\D}_L \mathscr{Y}^n h_{(\text{rect})}|^2
	\dVol_{g}
	&\leq \epsilon^{2(N_2 + 2 - n)} (1+\tau)^{-1 + C_{(n)}\delta}
	\\
	\int_{\Sigma_\tau} 
		\delta (1+r)^{1-C_{[n]}} |\slashed{\D}_L (r\mathscr{Y}^n h_{(\text{rect})})|^2
	r^{-2}\dVol_{g}
	&\leq \epsilon^{2(N_2 + 2 - n)} (1+\tau)^{C_{(n)}\delta}
	\\
	\int_{\mathcal{M}_\tau^{\tau_1}} 
		C_{[n]}\epsilon (1+r)^{-1-C_{[n]}} |\slashed{\D} \mathscr{Y}^n h_{(\text{rect})}|^2
	\dVol_{g}
	&\leq \epsilon^{2(N_2 + 2 - n)} (1+\tau)^{-1 + C_{(n)}\delta}
	\\
	\int_{\mathcal{M}_\tau^{\tau_1}} 
		\delta (1+r)^{-1-\delta} |\slashed{\D} \mathscr{Y}^n h_{(\text{rect})}|^2
	\dVol_{g}
	&\leq \epsilon^{2(N_2 + 2 - n)} (1+\tau)^{-1 + C_{(n)}\delta}
	\\
	\int_{\mathcal{M}_\tau^{\tau_1}} 
		\delta (1+r)^{-1+ \frac{1}{2}\delta} |\overline{\slashed{\D}} \mathscr{Y}^n h_{(\text{rect})}|^2
	\dVol_{g}
	&\leq \epsilon^{2(N_2 + 2 - n)} (1+\tau)^{-1 + C_{(n)}\delta}
\end{split}
\end{equation}
and
\begin{equation}
\begin{split}
	\int_{\Sigma_\tau} 
		(1+r)^{-C_{[n]}} |\overline{\slashed{\D}} \mathscr{Y}^n h|_{(\text{frame})}^2
	\dVol_{g}
	&\leq \epsilon^{2(N_2 + 2 - n)} (1+\tau)^{-1 + C_{(n)}\delta}
	\\
	\int_{\Sigma_\tau} 
		\delta (1+r)^{\frac{1}{2}\delta} |\slashed{\D}_L \mathscr{Y}^n h|_{(\text{frame})}^2
	\dVol_{g}
	&\leq \epsilon^{2(N_2 + 2 - n)} (1+\tau)^{-1 + C_{(n)}\delta}
	\\
	\int_{\Sigma_\tau} 
		\delta (1+r)^{1-C_{[n,m]}} |\slashed{\D}_L (r\mathscr{Y}^n h)|_{(\text{frame})}^2
	r^{-2}\dVol_{g}
	&\leq \epsilon^{2(N_2 + 2 - n)} (1+\tau)^{C_{(n)}\delta}
	\\
	\int_{\mathcal{M}_\tau^{\tau_1}} 
		C_{[n]}\epsilon (1+r)^{-1-C_{[n]}} |\slashed{\D} \mathscr{Y}^n h|_{(\text{frame})}^2
	\dVol_{g}
	&\leq \epsilon^{2(N_2 + 2 - n)} (1+\tau)^{-1 + C_{(n)}\delta}
	\\
	\int_{\mathcal{M}_\tau^{\tau_1}} 
		\delta (1+r)^{-1-\delta} |\slashed{\D} \mathscr{Y}^n h|_{(\text{frame})}^2
	\dVol_{g}
	&\leq \epsilon^{2(N_2 + 2 - n)} (1+\tau)^{-1 + C_{(n)}\delta}
	\\
	\int_{\mathcal{M}_\tau^{\tau_1}} 
		\delta (1+r)^{-1+ \frac{1}{2}\delta} |\overline{\slashed{\D}} \mathscr{Y}^n h|_{(\text{frame})}^2
	\dVol_{g}
	&\leq \epsilon^{2(N_2 + 2 - n)} (1+\tau)^{-1 + C_{(n)}\delta}
\end{split}
\end{equation}
and finally
\begin{equation}
\begin{split}
	\int_{\Sigma_\tau} 
		(1+r)^{-C_{[n,0]}} |\overline{\slashed{\D}} \mathscr{Y}^n h|_{LL}^2
	\dVol_{g}
	&\leq \epsilon^{2(N_2 + 2 - n)} (1+\tau)^{-1 + C_{(n)}\delta}
	\\
	\int_{\Sigma_\tau} 
		\delta (1+r)^{\frac{1}{2}\delta} |\slashed{\D}_L \mathscr{Y}^n h|_{LL}^2
	\dVol_{g}
	&\leq \epsilon^{2(N_2 + 2 - n)} (1+\tau)^{-1 + C_{(n)}\delta}
	\\
	\int_{\Sigma_\tau} 
		\delta (1+r)^{1-C_{[n,0]}} |\slashed{\D}_L (r\mathscr{Y}^n h)|_{LL}^2
	r^{-2}\dVol_{g}
	&\leq \epsilon^{2(N_2 + 2 - n)} (1+\tau)^{C_{(n)}\delta}
	\\
	\int_{\mathcal{M}_\tau^{\tau_1}} 
		C_{[n,0]}\epsilon (1+r)^{-1-C_{[n,0]}} |\slashed{\D} \mathscr{Y}^n h|_{LL}^2
	\dVol_{g}
	&\leq \epsilon^{2(N_2 + 2 - n)} (1+\tau)^{-1 + C_{(n)}\delta}
	\\
	\int_{\mathcal{M}_\tau^{\tau_1}} 
		\delta (1+r)^{-1-\delta} |\slashed{\D} \mathscr{Y}^n h|_{LL}^2
	\dVol_{g}
	&\leq \epsilon^{2(N_2 + 2 - n)} (1+\tau)^{-1 + C_{(n)}\delta}
	\\
	\int_{\mathcal{M}_\tau^{\tau_1}} 
		\delta (1+r)^{-1+ \frac{1}{2}\delta} |\overline{\slashed{\D}} \mathscr{Y}^n h|_{LL}^2
	\dVol_{g}
	&\leq \epsilon^{2(N_2 + 2 - n)} (1+\tau)^{-1 + C_{(n)}\delta}
\end{split}
\end{equation}

\chapter{Energy estimates}
\label{chapter energy estimates}
In this chapter we will use the bootstrap estimates from the previous chapter to prove various energy estimates. Our strategy will be to first write out the energy estimates with the error terms which follow from the bootstrap. These will then be combined in order to prove a boundedness and Morewetz (or ``integrated local energy decay'') estimate, as well as several different $p$-weighted energy estimates.

In this chapter we present the three basic energy estimates, arising from the use of the currents associated with the vector fields $T$ (the ``$T$-energy'' estimate), $R$ (the ``Morawetz'' energy estimate) and $r^p L$ (the ``$p$-weighted'' energy estimate). Each of these generates error terms, which we will later control by adding together suitable combinations of the basic energy inequalities, and using the Gronwall inequality.

Note that, from now on, the value of $r_0$ is considered fixed, and so all implicit constants are allowed to depend on $r_0$.

Note also that these estimates can be proved under the bootstrap assumptions of chapter \ref{chapter bootstrap}, but with the constants $C_{(n)}$ and $C_{(n,m)}$ replaced the larger constant $\mathring{C}$. These assumptions are therefore strictly weaker than those of chapter \ref{chapter bootstrap}; we only need the more detailed bootstrap assumptions when estimating the inhomogenous terms $F$ and when improving some of the pointwise bounds.

\begin{remark}[The constant $C_{(\phi)}$]
 The constant $C_{(\phi)}$ must be chosen suitably for each field $\phi_{(A)}$ appearing in the system of wave equations we are studying. Additionally, different constants must be chosen for the commuted fields; that is, we can consider the set of fields $\mathscr{Y}^n \phi_{(A)}$ for various choices of $n$, $(A)$ and different collections of differential operators $\mathscr{Y}$, and this constant must be allowed to take different values depending on which field we are estimating. However, $C_{(\phi)}$ is not a function of the \emph{value} of $\phi$: it is a constant, taking the same value at all points in $\mathcal{M}$. Rather, it depends on where the field $\phi$ fits into the structure of the set of wave equations under consideration.
\end{remark}

\section{The basic weighted \texorpdfstring{$T$}{T} energy estimates}

\begin{proposition}[The basic weighted $T$-energy estimate]
\label{proposition basic weighted T energy}
 Let $\phi$ be an $S_{\tau, r}$-tangent tensor field satisfying the equation
\begin{equation*}
 \tilde{\slashed{\Box}}_g \phi = F
\end{equation*}
for some $S_{\tau, r}$-tangent tensor field $F$. Define the weighted $T$ energy of $\phi$
\begin{equation}
 \begin{split}
  \mathcal{E}^{(wT)}[\phi](\tau, t, \tau_0) &:= 
   \int_{^t\Sigma_\tau \cap \{ r \leq r_0\} } w\left( |\slashed{\D}_L \phi|^2 + |\slashed{\D}_{\Lbar} \phi|^2 + |\slashed{\nabla}\phi|^2 \right) r^2 \upd r \wedge \dVol_{\mathbb{S}^2} \\
  &\phantom{:=} + \int_{^t\Sigma_\tau \cap \{r \geq r_0\} } w\left( |\slashed{\D}_L \phi|^2 + |\slashed{\nabla}\phi|^2 \right) r^2\upd r \wedge \dVol_{\mathbb{S}^2} \\
  &\phantom{:=}+ \int_{^\tau_{\tau_0}\bar{\Sigma}_t} w\left( |\slashed{\D}_L \phi|^2 + |\slashed{\D}_{\Lbar} \phi|^2 + |\slashed{\nabla}\phi|^2 \right) r^2 \upd r \wedge \dVol_{\mathbb{S}^2}
 \end{split}
\end{equation}
for $\tau > \tau_0$. Assume that all the bootstrap assumptions from chapter \ref{chapter bootstrap} hold. Choose the weight function $w$ to be
\begin{equation}
 w = (1+r)^{-C_{(\phi)}\epsilon}
\end{equation}

Then, for all sufficiently large constants $C_{(\phi)}$ and for any $\alpha > 0$ we have the following bound for the energy $\mathcal{E}^{(wT)}[\phi]$:
\begin{equation}
\label{equation basic weighted T energy estimate}
 \begin{split}
  &\mathcal{E}^{(wT)}[\phi](\tau, t, \tau_0)
  + \int_{^t\mathcal{M}_{\tau_0}^\tau} C_{(\phi)}\epsilon (1+r)^{-1} w|\slashed{\D}\phi|^2 \dVol_g  
   \lesssim
  \mathcal{E}^{(wT)}[\phi](\tau_0, t, \tau_0) \\
  &\phantom{\lesssim}
  + \int_{^t\mathcal{M}_{\tau_0}^{\tau}} w\bigg(
  \epsilon(1+\tau)^{-1-\delta}|\overline{\slashed{\D}}\phi|^2
  + C_{(\phi)}\epsilon(1+r)^{-1}|\overline{\slashed{\D}}\phi|^2
  + \epsilon (1+r)^{-3} |\phi|^2 \\
  &\phantom{\lesssim + \int_{^t\mathcal{M}_{\tau_0}^{\tau_1}} \bigg(}
  + \epsilon (1+r)^{-2 -\delta}(1+\tau)^{-1-\delta} |\phi|^2
  + \epsilon^{-1} (1+r) |F|^2
  + \epsilon (1+r)^{-1-\delta}|\slashed{\D}\phi|^2
   \bigg) \dVol_g 
 \end{split}
\end{equation}
where, if $\phi$ is in fact a scalar field, then the terms involving only $\phi$ (and not its derivatives) on the right hand side of the inequality above are not present.

\end{proposition}

\begin{proof}
 Applying the energy estimate \ref{proposition energy identity} to the field $\phi$, using the multiplier $wT$, in the spacetime region $^t\mathcal{M}_{\tau_0}^{\tau}$, the compatible current identity \ref{proposition compatible current identity} gives
\begin{equation*}
 \begin{split}
  &\int_{^t\Sigma_\tau} \imath_{^{(wT)}J[\phi]} \dVol_g
  + \int_{^{\tau_0}_{\tau}\bar{\Sigma}_t} \imath_{^{(wT)}J[\phi]} \dVol_g
  - \int_{^t\Sigma_{\tau_0}} \imath_{^{(wT)}J[\phi]} \dVol_g \\
  &= \int_{^t\mathcal{M}_{\tau_0}^{\tau}} \bigg( {^{(wT)}K}[\phi] - w\omega (\slashed{\D}_{\Lbar}\phi)\cdot(\slashed{\D}_T \phi) + w F\cdot (\slashed{\D}_T \phi) + wT^\mu ([\slashed{\D}_\mu , \slashed{\D}_\nu]\phi)\cdot (\slashed{\D}^\nu \phi)
  \bigg)\dVol_g
 \end{split}
\end{equation*}
where, if $\phi$ is a scalar field, then the final term on the right hand side clearly vanishes.

We will first show that the terms on the right hand side (the ``bulk terms'') are bounded by the spacetime integral given in the proposition, and then we will deal with the boundary terms. We will concentrate on the region $r \geq r_0$ below; in the region $r \leq r_0$ it is easy to see that suitable bounds hold. In particular, in this region $r$ is bounded, so no decay in $r$ is necessary.

Now, from proposition \ref{proposition bulk weighted T current} we have, schematically,
\begin{equation*}
 \begin{split}
  &{^{(wT)}K}[\phi] - w\omega (\slashed{\D}_{\Lbar}\phi)\cdot(\slashed{\D}_T \phi) - \frac{1}{2}C_{(\phi)}\epsilon (1+r)^{-1} w|\slashed{\D}\phi|^2 \\
  &= C_{(\phi)}\epsilon (1+r)^{-1} w |\overline{\slashed{\D}}\phi|^2 + w\omega |\slashed{\D}\phi|^2 + w\bm{\Gamma}(\slashed{\D}\phi)\cdot(\overline{\slashed{\D}}\phi)
 \end{split}
\end{equation*}
Note that the bounds in chapter \ref{chapter bootstrap} imply, in particular,
\begin{equation*}
\begin{split}
 \left| \bm{\Gamma} \right| &\leq \epsilon \left( (1+r)^{-1} + (1+r)^{-1+\delta}(1+\tau)^{-\beta} \right)
 \\
 |\omega| &\leq \epsilon (1+r)^{-1}
\end{split}
\end{equation*}
and so, making use of these bounds we have
\begin{equation*}
 \begin{split}
  &\left| {^{(wT)}K}[\phi] - w\omega (\slashed{\D}_{\Lbar}\phi)\cdot(\slashed{\D}_T \phi) - \frac{1}{2}C_{(\phi)}\epsilon (1+r)^{-1} w|\slashed{\D}\phi|^2 \right| \\
  &\lesssim \epsilon(1+r)^{-1}w|\slashed{\D}\phi|^2
  + \epsilon(1+r)^{-1 + \delta}(1+\tau)^{-\beta}w|\slashed{\D}\phi| |\overline{\slashed{\D}}\phi| 
  + C_{(\phi)}\epsilon (1+r)^{-1}w|\overline{\slashed{\D}}\phi|^2
  \\ \\
  &\lesssim \epsilon(1+r)^{-1}w|\slashed{\D}\phi|^2 
  + \epsilon(1+\tau)^{-1-2\beta+2\delta} w|\overline{\slashed{\D}}\phi|^2
  + C_{(\phi)}\epsilon (1+r)^{-1}w|\overline{\slashed{\D}}\phi|^2
 \end{split}
\end{equation*}

The next term on the right hand side is estimated as
\begin{equation*}
 \begin{split}
 \left| wF \cdot (\slashed{\D}_T \phi) \right| &\lesssim \epsilon(1+r)^{-1}w|\slashed{\D}\phi|^2 + \epsilon^{-1}(1+r)w|F|^2
  \end{split}
\end{equation*}

These are the only terms which arise if $\phi$ is a scalar field; the terms discussed below are only present if $\phi$ is a higher rank tensor field.

The final term is given schematically as
\begin{equation*}
 wT^\mu \left( [\slashed{\D}_\mu, \slashed{\D}_\nu]\phi \right) \cdot \left( \slashed{\D}^\nu \phi \right)
  = w\begin{pmatrix} \Omega_{L\Lbar} \cdot \slashed{\D}\phi \\ \slashed{\Omega}_{L} \cdot \slashed{\nabla} \phi \\ \slashed{\Omega}_{\Lbar} \cdot \slashed{\nabla} \phi \end{pmatrix} \cdot \phi
\end{equation*}
where we remind the reader that $\slashed{\Omega}_L$ refers to the $S_{\tau,r}$-tangent tensor field $\slashed{\Pi}_\mu^{\nu} \Omega_{L \nu \slashed{\alpha} \slashed{\beta}}$, and similarly for $\slashed{\Omega}_{\Lbar}$.

Let us first consider the region $r \geq r_0$. Using the expressions for the components of the curvature $\Omega$ given in chapter \ref{chapter geometry of vector bundle} together with the bootstrap we easily find that
\begin{equation*}
\begin{split}
wT^\mu \left( [\slashed{\D}_\mu, \slashed{\D}_\nu]\phi \right) \cdot \left( \slashed{\D}^\nu \phi \right)
&\lesssim
\epsilon w \left( (1+r)^{-2} + (1+r)^{-2+\delta}(1+\tau)^{-\beta} \right)|\phi| |\slashed{\D}\phi|
\\
&\lesssim
\epsilon (1+r)^{-1}w|\slashed{\D}\phi|^2
+ \epsilon (1+r)^{-3}w|\phi|^2
+ \epsilon (1+r)^{-3+2\delta}(1+\tau)^{-2\beta}w|\phi|^2
\end{split}
\end{equation*}

We can further estimate
\begin{equation*}
	\epsilon w (1+r)^{-3+2\delta} (1+\tau)^{-2\beta}|\phi|^2
	\lesssim
	\epsilon w (1+r)^{-3 - 2\delta} |\phi|^2
	+ \epsilon w (1+r)^{-2 - 2\delta} (1+\tau)^{-1-\delta}|\phi|^2
\end{equation*}
where the first quantitiy on the right hand side is sufficient to bound the left hand side in the region $r \leq \tau$, and the second quantity bounds the right hand side in the region $r \geq \tau$, and we have used the fact that $\beta \gg \delta$.

In the region $r \leq r_0$ we can bound the same terms by simply using the schematic expression
\begin{equation*}
 \begin{split}
  T^\mu \left( [\D_\mu, \D_\nu]\phi \right) \cdot \left( \slashed{\D}^\nu \phi \right) &\sim \Omega_{abcd} T^a \phi \cdot (\slashed{\D}\phi) \\
  &\sim \left(\partial^2 h_{(\text{rect})} + (\partial h_\text{rect})\cdot(\partial h_\text{rect}) \right)  \phi \cdot (\slashed{\D}\phi) \\
  &\lesssim \epsilon |\phi|^2 + \epsilon |\slashed{\D}\phi|^2
 \end{split}
\end{equation*}
where we have used the bootstrap bounds in the region $r \leq r_0$. Note that all of the components of $\Omega$ are ``small'' (meaning that they vanish if $h = 0$) except for the \emph{angular} components, which behave as
\begin{equation*}
	\Omega_{\slashed{\alpha}\slashed{\beta}} \sim r^{-2}
\end{equation*}

We now put together all of the above calculations, and choose $C_{(\phi)}$ sufficiently large to absorb all of the error terms of the form $\epsilon(1+r)^{-1}|\slashed{\D}\phi|^2$ by the positive bulk term on the left hand side.

Now we need to estimate the boundary terms. From proposition \ref{proposition boundary weighted T} we have
\begin{equation*}
 \begin{split}
  \int_{^t\Sigma_\tau} \imath_{{^{(wT)}J}} \dVol_g &= \int_{^t\Sigma_{\tau}} \frac{1}{2}w\left( |\slashed{\D}_L\phi|^2 + |\slashed{\nabla}\phi|^2 \right) \Omega^2 \upd r \wedge \dVol_{\mathbb{S}^2} 
 \end{split}
\end{equation*}
Now, the bootstrap bound on $\Omega$ implies
\begin{equation*}
\Omega^2 \sim r^2
\end{equation*}
and so
\begin{equation*}
 \int_{^t\Sigma_\tau} \imath_{{^{(wT)}J}} \dVol_g \sim \int_{^t\Sigma_\tau} w\left( |\slashed{\D}_L \phi|^2 + |\slashed{\nabla}\phi|^2 \right) r^2 \upd r \wedge \dVol_{\mathbb{S}^2} 
\end{equation*}

Finally the boundary term on a surface of constant $t$ in the weighted $T$-energy is given by
\begin{equation*}
 \begin{split}
  &\int_{\{t = t_0\}} \imath_{{^{(wT)}J}} \dVol_g \\
  &= \int_{\{t = t_0\}} \frac{1}{2} \frac{1}{(L^0 + \Lbar^0 - b^0)} w\left( |\slashed{\D}_L \phi|^2 + |\slashed{\D}_{\Lbar} \phi|^2 + 2|\slashed{\nabla}\phi|^2 + \textit{Err}_{(wT,t_0\text{-bdy})} \right) \upd r \wedge \dVol_{(\mathbb{S}^2, G)} \\
  &= \int_{\{t = t_0\}} \frac{1}{2} \frac{1}{(L^0 + \Lbar^0 - b^0)} w\left( |\slashed{\D}_L \phi|^2 + |\slashed{\D}_{\Lbar} \phi|^2 + 2|\slashed{\nabla}\phi|^2 + \textit{Err}_{(wT,t_0\text{-bdy})} \right) \Omega^2 \upd r \wedge \dVol_{\mathbb{S}^2}
 \end{split}
\end{equation*}
where we have the bound
\begin{equation*}
 |\textit{Err}_{(wT, t_0-\text{bdy})}| \lesssim \left( |(L_{(\text{small})})^0| + |(\Lbar_{(\text{small})})^0| + |\slashed{\nabla}t| \right)|\slashed{\D}\phi|^2
\end{equation*}
Note that $\slashed{\nabla}_\mu t = \slashed{\Pi}_\mu^{\phantom{\mu}0}$. Hence the bootstrap assumptions imply that
\begin{equation*}
 |\textit{Err}_{(T, t_0-\text{bdy})}| \lesssim \epsilon |\slashed{\D}\phi|^2
\end{equation*}
Moreover, we have
\begin{equation*}
 |L^0 + \Lbar^0 - b^0 - 2| \lesssim \epsilon
\end{equation*}
hence we can conclude that 
\begin{equation*}
 \int_{\{t = t_0\}} \imath_{{^{(wT)}J}} \dVol_g \sim \int_{\{t = t_0\}} w \left( |\slashed{\D}\phi|^2 \right) r^2 \dVol_{\mathbb{S}^2}
\end{equation*}

Putting together all of the calculations above proves the proposition.

\end{proof}

\section{The basic weighted Morawetz energy estimate}

\begin{proposition}[The basic Morawetz energy estimate]
\label{proposition basic Morawetz estimate}
 Let $\phi$ be an $S_{\tau,r}$-tangent tensor field satisfying the equation
\begin{equation*}
 \tilde{\slashed{\Box}}_g \phi = F
\end{equation*}
for some $S_{\tau,r}$-tangent tensor field $F$. Choose the weight function
\begin{equation*}
 w := (1+r)^{-C_{(\phi)}\epsilon}
\end{equation*}
 Assume that all the bootstrap assumptions from chapter \ref{chapter bootstrap} hold. Then for all sufficiently small $\epsilon$ we have
\begin{equation}
\label{equation basic weighted Morawetz estimate}
 \begin{split}
  &\int_{^t\mathcal{M}_{\tau_0}^{\tau}} \left( \delta(1+r)^{-1-\frac{1}{2}\delta}|\slashed{\D}\phi|^2 + \chi_{(r_0)}(r)(1+r)^{-1 - C_{(\phi)}\epsilon}|\slashed{\nabla}\phi|^2 + \delta\left(1+ \frac{1}{2}\delta\right) r^{-1}(1+r)^{-2-\frac{1}{2}\delta}|\phi|^2 \right) \dVol_g \\
  &\lesssim \mathcal{E}^{(T, C_{(\phi)}\epsilon)}[\phi](\tau,t,\tau_0) + \mathcal{E}^{(T, C_{(\phi)}\epsilon)}[\phi](\tau_0, t, \tau_0) 
  \\
  &\phantom{\lesssim}
  + \int_{^t\mathcal{M}_{\tau_0}^{\tau}} w\bigg( 
  (1 +C_{(\phi)})\epsilon(1+r)^{-1}|\slashed{\D}\phi|^2 
  + \epsilon(1+\tau)^{-1-\delta}|\overline{\slashed{\D}}\phi|^2 
  + \epsilon r^{-1}(1+r)^{-1 - \delta}(1+\tau)^{-1-\delta}|\phi|^2 \\
  &\phantom{\lesssim + \int_{^t\mathcal{M}_{\tau_0}^{\tau_1}} \bigg(}
  + (1+C_{(\phi)})\epsilon r^{-1}(1+r)^{-2} |\phi|^2
  + \epsilon^{-1} (1+r) |F|^2
   \bigg) \dVol_g \\
  &\phantom{=}
  + \int_{\bar{S}_{t,\tau_0}}(1+r)^{-C_{(\phi)}\epsilon} r|\phi|^2 \dVol_{\mathbb{S}^2}
  + \int_{\bar{S}_{t,\tau}}(1+r)^{-C_{(\phi)}\epsilon} r|\phi|^2 \dVol_{\mathbb{S}^2}
 \end{split}
\end{equation}
where $\chi_{(r_0)}(r)$ is the smooth monotone cut-off function defined in equation \eqref{equation cut off functions}.

Moreover, if the field $\phi$ is in fact a scalar field, then the terms
\begin{equation*}
 \int_{^t\mathcal{M}_{\tau_0}^{\tau}} \epsilon^{-1}(1+r)^{-3} |\phi|^2 \dVol_g
\end{equation*}
and
\begin{equation*}
 \int_{^t\mathcal{M}_{\tau_0}^{\tau} \cap \{r < r_0\}} \epsilon^{-1} r^{-2}|\phi|^2 \dVol_g
\end{equation*}
are absent.
\end{proposition}

\begin{proof}
 We apply the energy estimate to the field $\phi$ in the spacetime region $^t\mathcal{M}^\tau_{\tau_0}$, with the modified energy current $^{(wR)}\tilde{J}[\phi]$, and with the choices
\begin{equation*}
 \begin{split}
  f_R(r) &:= 1 - (1+r)^{-\frac{1}{2}\delta + C_{(\phi)}\epsilon} \\
  w(r) &:= (1+r)^{-C_{(\phi)}\epsilon}
 \end{split}
\end{equation*}
Using the modified compatible current identity \ref{proposition modified compatible current identity} We obtain
\begin{equation}
\label{equation weighted Morawetz estimate 1}
 \begin{split}
  &\int_{^t\Sigma_\tau} \imath_{^{(wR)}\tilde{J}[\phi]} \dVol_g
  + \int_{^{\tau_0}_{\tau}\bar{\Sigma}_t} \imath_{^{(wR)}\tilde{J}[\phi]} \dVol_g
  - \int_{^t\Sigma_{\tau_0}} \imath_{^{(wR)}\tilde{J}[\phi]} \dVol_g \\
  &=\int_{^t\mathcal{M}^\tau_{\tau_0}} \bigg(
  {^{(wR)}\tilde{K}}[\phi]
  - \omega w f_R (\slashed{\D}_{\Lbar}\phi)\cdot \left( \slashed{\D}_R \phi 
  + r^{-1}\phi \right) 
  + w f_R F \cdot \left( \slashed{\D}_R \phi + r^{-1}\phi \right) \\
  &\phantom{=\int_{^t\mathcal{M}^\tau_{\tau_0}} \bigg(}
  + w f_R R^\nu \left([\slashed{\D}_\mu, \slashed{\D}_\nu]\phi \right)\cdot \left(\slashed{\D}^\mu \phi \right) \bigg) \dVol_g
 \end{split}
\end{equation}
As before we shall first estimate the bulk terms, and then turn to the boundary terms.

Making use of proposition \ref{proposition bulk weighted Morawetz current}, together with the facts that
\begin{equation*}
 \begin{split}
  f'_R &= \left(\frac{1}{2}\delta - C_{(\phi)}\epsilon \right)(1+r)^{-1-\frac{1}{2}\delta + C_{(\phi)}\epsilon} \\
  f''_R &= -\left((1+\frac{1}{2}\delta - C_{(\phi)}\epsilon\right)\left(\frac{1}{2}\delta - C_{(\phi)}\epsilon\right)(1+r)^{-2-\frac{1}{2}\delta + C_{(\phi)}\epsilon} \\
  w' &= -C_{(\phi)}\epsilon (1+r)^{-1-C_{(\phi)}\epsilon} \\
  w'' &= C_{(\phi)}\epsilon \left( 1 + C_{(\phi)}\epsilon \right)(1+r)^{-2-C_{(\phi)}\epsilon}
 \end{split}
\end{equation*}
we have
\begin{equation}
\label{equation weighted Morawetz estimate 2}
 \begin{split}
  &{^{(wR)}\tilde{K}}[\phi]
  - \omega w f_R (\slashed{\D}_{\Lbar}\phi)\cdot \left( \slashed{\D}_R \phi 
  + r^{-1}\phi \right) \\
  &= \frac{1}{8}(\delta - 2C_{(\phi)}\epsilon)(1+r)^{-1-\frac{1}{2}\delta} \left( |\slashed{\D}_L \phi|^2 + |\slashed{\D}_{\Lbar} \phi|^2 \right) \\
  &\phantom{=}
  + \left(\frac{(1+r)^{1+\frac{1}{2}\delta-C_{(\phi)}\epsilon} - (1+r) - \frac{1}{4}(\delta - 2C_{(\phi)}\epsilon) r}{r(1+r)^{1+\frac{1}{2}\delta}} \right) |\slashed{\nabla}\phi|^2
  + \frac{1}{4}\delta \left(1+\frac{1}{2}\delta\right) r^{-1}(1+r)^{-2-\frac{1}{2}\delta}|\phi|^2 \\
  &\phantom{=}
  + w\textit{Err}_{(wR, \text{bulk})}
 \end{split}
\end{equation}
The error term will be dealt with below; we first turn to the other terms. We first note that, if $\epsilon$ is chosen sufficiently small relative to $\delta$ and $C_{(\phi)}$, then we have $(\delta - 2C_{(\phi)}\epsilon) \geq \frac{1}{2}\delta$. The only term which needs additional consideration is the coefficient of $|\slashed{\nabla}\phi|^2$; in the region $r \geq r_0$ this satisfies
\begin{equation*}
 \frac{(1+r)^{1+\frac{1}{2}\delta-C_{(\phi)}\epsilon} - (1+r) - \frac{1}{4}(\delta - 2C_{(\phi)}\epsilon) r}{r(1+r)^{1+\frac{1}{2}\delta}} \gtrsim w(1+r)^{-1}
\end{equation*}
where the implicit constant depends only on $r_0$. In the region $r < r_0$ we instead have the bound
\begin{equation*}
 \frac{(1+r)^{1+\frac{1}{2}\delta-C_{(\phi)}\epsilon} - (1+r) - \frac{1}{4}(\delta - 2C_{(\phi)}\epsilon) r}{r(1+r)^{1+\frac{1}{2}\delta}} \gtrsim (\delta - 2C_{(\phi)}\epsilon)(1+r)^{-1 - \frac{1}{2}\delta}
\end{equation*}
Choosing $\epsilon$ sufficiently small compared to $\delta$, we have $\delta - 2C_{(\phi)}\epsilon \geq \frac{1}{2}\delta$. Hence we have the inequality
\begin{equation*}
 \frac{(1+r)^{1+\frac{1}{2}\delta-C_{(\phi)}\epsilon} - (1+r) - \frac{1}{4}(\delta - 2C_{(\phi)}\epsilon) r}{r(1+r)^{1+\frac{1}{2}\delta}} \gtrsim \chi_{(r_0)}(r) w(1+r)^{-1} + \delta(1+r)^{-1-\frac{1}{2}\delta}
\end{equation*}

The error term in \eqref{equation weighted Morawetz estimate 2} can be estimated, using the bootstrap bounds of chapter \ref{chapter bootstrap} together with the expression given for the error term in proposition \ref{proposition bulk weighted Morawetz current}, as
\begin{equation*}
 \begin{split}
  \left| \textit{Err}_{(wR, \text{bulk})} \right| &\lesssim
  (1 + C_{(\phi)})\epsilon (1+r)^{-1}|\slashed{\D}\phi|^2
  + \epsilon(1+r)^{-1 + \delta}(1+\tau)^{-\beta}|\slashed{\D}\phi||\overline{\slashed{\D}}\phi| \\
  &\phantom{\lesssim}
  + (1 + C_{(\phi)}\epsilon)\epsilon r^{-1}(1+r)^{-2 + \delta}(1+\tau)^{-\beta}|\phi|^2
  + (1 + C_{(\phi)}\epsilon)\epsilon r^{-1}(1+r)^{-2}|\phi|^2 \\
  &\phantom{\lesssim}
  + \epsilon (1+r)^{-2}|\slashed{\D}\phi||\phi|
  \\
  &\lesssim
  C_{(\phi)}\epsilon (1+r)^{-1}|\slashed{\D}\phi|^2
  + \epsilon(1+\tau)^{-1 - \delta}|\overline{\slashed{\D}}\phi|^2
  + \epsilon r^{-1}(1+r)^{-1-\delta}(1+\tau)^{-1-\delta}|\phi|^2 \\
  &\phantom{\lesssim}
  + C_{(\phi)}\epsilon r^{-1}(1+r)^{-2}|\phi|^2  
 \end{split}
\end{equation*}

The next term can be bounded as follows:
\begin{equation*}
 \left| w\omega f_R (\slashed{\D}_{\Lbar} \phi) \cdot \left( \slashed{\D}_R\phi + r^{-1}\phi \right) \right|
 \lesssim \epsilon(1+r)^{-1} w|\slashed{\D}\phi|^2 + \epsilon(1+r)^{-3}w|\phi|^2 
\end{equation*}
where we have used the fact that, as $r \rightarrow 0$, $f_R \sim \delta r$.

In a similar way, we can bound
\begin{equation*}
 \begin{split}
  \left| w f_R F \left(\slashed{\D}_R \phi + r^{-1}\phi \right) \right|
  &\lesssim \epsilon(1+r)^{-1}w|\slashed{\D}\phi|^2
  + \epsilon(1+r)^{-3}w|\phi|^2 
  + \epsilon^{-1}(1+r)w|F|^2
 \end{split}
\end{equation*}

Next, we turn to those terms which are absent if $\phi$ is a scalar field. First, in the region $r \geq r_0$ we can bound
\begin{equation*}
 \begin{split}
  \left| w f_R R^\nu \left( [\slashed{\D}_\mu, \slashed{\D}_\nu]\phi \right) \cdot (\slashed{\D}^\mu \phi) \right|
  &\lesssim \epsilon(1+r)^{-1}w|\slashed{\D}\phi|^2 + \epsilon(1+r)^{-3}w|\phi|^2
  + \epsilon (1+r)^{-2-2\delta} (1+\tau)^{-1-\delta}w|\phi|^2
 \end{split}
\end{equation*}
where the computations are almost identical to those appearing in the $T$-energy estimate. In the region $r \leq r_0$ we can similarly bound
\begin{equation*}
  \left| f_R R^\nu \left( [\slashed{\D}_\mu, \slashed{\D}_\nu]\phi \right) \cdot (\slashed{\D}^\mu \phi) \right| \lesssim \left( |\partial^2 h| + |\partial h|^2 \right)|\phi||\partial \phi|
\end{equation*}
Since, in this region, $r$ is bounded, and the bootstrap bounds $|\partial h|\leq \epsilon$ and $|\partial^2 h|\leq \epsilon$ are sufficient to give the same estimate in this region.

Next, we need to estimate the boundary terms. Unlike in the $T$-energy estimate, these do not come with any particular sign, and we only aim to show that they can be bounded by the terms appearing in the $T$-energy.

From proposition \ref{proposition boundary weighted Morawetz}, using the fact that $|f_R| < 1$, we have
\begin{equation*}
 \int_{^t\Sigma_\tau} \imath_{^{(wR)}\tilde{J}}\dVol_g 
  \lesssim \int_{^t\Sigma_\tau} w\left( |\overline{\slashed{\D}}\phi|^2 + r^{-1}(1+r)^{-1}|\phi|^2 \right)\Omega^2 \upd r \wedge \dVol_{\mathbb{S}^2}
\end{equation*}
Using the bootstrap bounds on $\Omega$ we find
\begin{equation*}
 \begin{split}
  \int_{^t\Sigma_\tau} \imath_{^{(wR)}\tilde{J}}\dVol_g 
  &\lesssim \int_{^t\Sigma_\tau} w\left( |\overline{\slashed{\D}}\phi|^2 + r^{-1}(1+r)^{-1}|\phi|^2 \right)r^2 \upd r \wedge \dVol_{\mathbb{S}^2} \\
  &\lesssim \mathcal{E}^{(T, C_{(\phi)}\epsilon)}[\phi](\tau,t,\tau) + \int_{^t\Sigma_\tau} (1+r)^{-1 - C_{(\phi)}\epsilon}|\phi|^2 r \upd r \wedge \dVol_{\mathbb{S}^2}
 \end{split}
\end{equation*}
Now, we can use the Hardy inequality of proposition \ref{proposition Hardy} together with the fact that $C_{(\phi)}\epsilon \ll 1$ to show
\begin{equation*}
 \int_{^t\Sigma_\tau} (1+r)^{-1} w|\phi|^2 r \upd r \wedge \dVol_{\mathbb{S}^2}
 \lesssim \int_{^t\Sigma_\tau} w|\slashed{\D}_L \phi|^2 r^2 \upd r \wedge \dVol_{\mathbb{S}^2}
  + \int_{\bar{S}_{r,t}} w|\phi|^2 r \dVol_{\mathbb{S}^2} 
\end{equation*}
Hence we find that
\begin{equation*}
 \int_{^t\Sigma_\tau} (1+r)^{-1}w|\phi|^2 r \upd r \wedge \dVol_{\mathbb{S}^2}
 \lesssim \mathcal{E}^{(T, C_{(\phi)}\epsilon)}[\phi](\tau,t,\tau)
  + \int_{\bar{S}_{r,t}} w|\phi|^2 r \dVol_{\mathbb{S}^2}
\end{equation*}

Next, we turn to the boundary term on $^\tau_{\tau_0}\bar{\Sigma}_t$, i.e.\ the part of the boundary on which $t$ is constant. Again, using proposition \ref{proposition boundary weighted Morawetz} together with the bootstrap bounds of chapter \ref{chapter bootstrap} we find
\begin{equation*}
 \int_{^\tau_{\tau_0}\bar{\Sigma}_t} \imath_{^{(wR)}\tilde{J}}\dVol_g
  \lesssim
  \int_{^\tau_{\tau_0}\bar{\Sigma}_t} w\left( |\slashed{\D}\phi|^2 + r^{-1}(1+r)^{-1} \right)|\phi|^2 r^2 \upd r \wedge \dVol_{\mathbb{S}^2}
\end{equation*}

We can now use the other version of the Hardy inequality, given in proposition \ref{proposition Hardy on constant t}, to estimate the lower order term. Making use of the bootstrap bounds, it is easy to show
\begin{equation*}
  \int_{^\tau_{\tau_0}\bar{\Sigma}_t} w(1+r)^{-1}|\phi|^2 r \upd r \wedge \dVol_{\mathbb{S}^2}
  \lesssim \int_{^\tau_{\tau_0}\bar{\Sigma}_t}w |\slashed{\D}\phi|^2 r^2 \upd r \wedge \dVol_{\mathbb{S}^2} + \int_{\bar{S}_{\tau_0,t}} w |\phi|^2 r \dVol_{\mathbb{S}^2}
\end{equation*}

This finishes the calculation of the boundary terms. Putting all of the calculations above together, and absorbing terms when appropriate (assuming, in particular, that $\epsilon \ll \delta$) proves the proposition.

\end{proof}

\section{The basic \texorpdfstring{$p$}{p}-weighted energy estimates}

\begin{proposition}
\label{proposition basic p weighted estimate}

Let $\phi$ be an $S_{\tau, r}$-tangent tensor field satisfying the equation
\begin{equation*}
 \tilde{\slashed{\Box}}_g \phi = F_1 + F_2 + F_3
\end{equation*}
for some $S_{\tau, r}$-tangent tensor fields $F_1$, $F_2$, $F_3$. Let $\psi := r\phi$. Define the $p$-weighted energy of $\phi$
\begin{equation}
 \begin{split}
  \mathcal{E}^{(L,p)}[\phi](\tau, R) &:= 
  \int_{^t\Sigma_\tau} \chi_{(2r_0, R)}r^{p} |\slashed{\D}_L \psi|^2 \, \upd r \wedge \dVol_{\mathbb{S}^2}
 \end{split}
\end{equation}
for some constant $R$, and where the cut-off function $\chi_{(2r_0, R)}$ is defined in equation \eqref{equation cut off functions}.

Assume that all the bootstrap assumptions from chapter \ref{chapter bootstrap} hold. Let $0 < p \leq \delta$, and choose $t$ sufficiently large relative to $R$ and $\tau$ so that 
\begin{equation*}
 \text{supp}(\chi_{r_0, R}) \cap {_{\tau_0}^{\tau}\bar{\Sigma}}_t = \emptyset
\end{equation*}

Then, for all sufficiently small $\delta$, and for all sufficiently small $\epsilon$ (depending on $\delta$) we have 
\begin{equation}
\label{equation basic p weighted estimate low p}
 \begin{split}
  &\mathcal{E}^{(L,p)}[\phi](\tau, R) 
  + \int_{^t\mathcal{M}_{\tau_0}^{\tau}} \chi_{(2r_0, R)}\left( pr^{p-1}|\slashed{\D}_L \phi|^2 + (2-p)r^{p-1}|\slashed{\nabla}\phi|^2 + p(1-p)r^{p-3}|\phi|^2 \right) \dVol_g \\
  &\lesssim
  \mathcal{E}^{(L,p)}[\phi](\tau_0, R) \\
  &\phantom{\lesssim}
  + \int_{^t\mathcal{M}_{\tau_0}^{\tau}} \bigg(
  \epsilon (1+r)^{-1-\delta}|\slashed{\D}\phi|^2
  + \epsilon \chi_{(2r_0, R)} r^{p-2} (1+\tau)^{-1-\delta}|\slashed{\D}_L \psi|^2
  \\
  &\phantom{\lesssim + \int_{^t\mathcal{M}_{\tau_0}^{\tau}} }
  + \epsilon \chi_{(2r_0,R)} (p-1)r^{p-2 - \delta}(1+\tau)^{-1-\delta} |\phi|^2 
  + \epsilon^{-1}\chi_{(2r_0, R)} r^p (1+\tau)^{1+\delta}|F_1|^2\\
  &\phantom{\lesssim + \int_{^t\mathcal{M}_{\tau_0}^{\tau}} }
  + \epsilon^{-1}\chi_{(2r_0, R)} r^{p+1-2\delta}(1+\tau)^{6\delta}|F_2|^2
  + \epsilon^{-1}\chi_{(2r_0, R)} r^{p+1}|F_3|^2
  + \epsilon\chi_{(2r_0, R)} r^{p-3}|\phi|^2
   \bigg) \dVol_g \\
  &\phantom{\lesssim}
  + \int_{^t\mathcal{M}_{\tau_0}^{\tau} \cap \{r_0 \leq r \leq 2r_0\}} \left( |\slashed{\D}\phi|^2 + |\phi|^2 \right) \dVol_g \\
  &\phantom{\lesssim}
  + \int_{^t\Sigma_{\tau}} \epsilon (1+r)^{p-\delta} |\slashed{\D}_L \phi|^2 r^2 \upd r \wedge \dVol_{\mathbb{S}^2}
  + \int_{\bar{S}_{t, \tau}} \epsilon(1+r)^{p + 1 - \delta } |\phi|^2 \dVol_{\mathbb{S}^2} \\
  &\phantom{\lesssim}
  + \int_{^t\Sigma_{\tau_0}} \epsilon \frac{1}{(1-p+\delta)^2}(1+r)^{p-\delta} |\slashed{\D}_L \phi|^2 r^2 \upd r \wedge \dVol_{\mathbb{S}^2}
  + \int_{\bar{S}_{t, \tau_0}} \epsilon \frac{1}{(1-p+\delta)}(1+r)^{p + 1 - \delta } |\phi|^2 \dVol_{\mathbb{S}^2} \\
  &\phantom{\lesssim} 
  + \int_{\tau_0}^\tau \left( \int_{^t\Sigma_{\tau'}\cap\left\{ \frac{1}{2}R \leq r \leq R \right\}} p^{-2} r^{p-1}|\overline{\slashed{\D}}\phi|^2 r^2 \upd r \wedge \dVol_{\mathbb{S}^2} 
  + \int_{S_{\tau,R}} p^{-1} r^p |\phi|^2 \dVol_{\mathbb{S}^2} \right) \upd \tau
 \end{split}
\end{equation}
where, if $\phi$ is in fact a scalar field, then the term in the spacetime integral involving 
\begin{equation*}
 \epsilon\chi_{(r_0, R)} r^{p-3}|\phi|^2
\end{equation*}
is not present.

On the other hand, for all positive values of $p < 1$, if $t$ is sufficiently large compared to $R$ and $(\tau - \tau_0)$, then we can define the cut-off $p$-weighted energy
\begin{equation}
 \begin{split}
  \mathcal{E}^{(L,p)}[\phi](\tau, R) &:= 
  \int_{^t\Sigma_\tau} \chi_{(2r_0, R)} r^{p-2} |\slashed{\D}_L \psi|^2 \upd r \wedge \dVol_{\mathbb{S}^2}
 \end{split}
\end{equation}
Then, for all sufficiently small $\delta$ and for all sufficiently small $\epsilon$ (depending on $\delta$) we have
\begin{equation}
\label{equation basic p weighted estimate high p}
 \begin{split}
  &\mathcal{E}^{(L,p)}[\phi](\tau, R) 
  + \int_{^t\mathcal{M}_{\tau_0}^{\tau}} \chi_{(2r_0, R)}\left( pr^{p-1}|\slashed{\D}_L \phi|^2 + (2-p)r^{p-1}|\slashed{\nabla}\phi|^2 + p(1-p)r^{p-3}|\phi|^2 \right) \dVol_g \\
  &\lesssim
  \mathcal{E}^{(L,p)}[\phi](\tau_0, R)
  \\
  &\phantom{\lesssim}
  + \int_{^t\mathcal{M}_{\tau_0}^{\tau}}\bigg( 
  \epsilon \chi_{(2r_0, R)}(1+r)^{p-2} (1+\tau)^{1-2\beta} |\slashed{\D}\phi|^2
  + \epsilon \chi_{(2r_0, R)} r^{p-2} (1+\tau)^{-1-C_{(\phi)}\epsilon}|\slashed{\D}_L \psi|^2
  \\
  &\phantom{\lesssim + \int_{^t\mathcal{M}_{\tau_0}^{\tau}} \bigg(}
  + \epsilon \chi_{(2r_0, R)} r^{p-2-\frac{1}{2}\delta}(1+\tau)^{-1-\delta} |\phi|^2  
  + \epsilon^{-1}\chi_{(r_0, R)} r^p (1+\tau)^{1+\delta}|F_1|^2
  \\
  &\phantom{\lesssim + \int_{^t\mathcal{M}_{\tau_0}^{\tau}} \bigg(}
  + \epsilon^{-1}\chi_{(r_0, R)} r^{p+1-4\delta}(1+\tau)^{6\delta}|F_2|^2
  + \epsilon^{-1}\chi_{(2r_0, R)} r^{p+1}|F_3|^2
  + \epsilon\chi_{(r_0, R)} r^{p-3}|\phi|^2
   \bigg) \dVol_g \\
  &\phantom{\lesssim} 
  + \int_{^t\mathcal{M}_{\tau_0}^{\tau} \cap \{r_0 \leq r \leq 2r_0\}} \left( |\slashed{\D}\phi|^2 + |\phi|^2 \right) \dVol_g \\
  &\phantom{\lesssim}
  + \int_{^t\Sigma_{\tau}} \epsilon \frac{1}{(1-p+\delta)^2} r^{p-\delta} |\slashed{\D}_L \psi|^2 \upd r \wedge \dVol_{\mathbb{S}^2}
  + \int_{S_{\tau, R}} \epsilon \frac{1}{(1-p+\delta)} r^{p + 1 - \delta } |\phi|^2 \dVol_{\mathbb{S}^2} \\
  &\phantom{\lesssim}
  + \int_{^t\Sigma_{\tau_0}} \epsilon \frac{1}{(1-p+\delta)^2} r^{p-\delta} |\slashed{\D}_L \psi|^2 \upd r \wedge \dVol_{\mathbb{S}^2}
  + \int_{S_{\tau_0, R}} \epsilon \frac{1}{(1-p+\delta)} r^{p + 1 - \delta } |\phi|^2 \dVol_{\mathbb{S}^2} \\
  &\phantom{\lesssim} 
  + \int_{\tau_0}^\tau \left( \int_{^t\Sigma_{\tau'}\cap\left\{ \frac{1}{2}R \leq r \leq R \right\}} p^{-2} r^{p-1}|\overline{\slashed{\D}}\phi|^2 r^2 \upd r \wedge \dVol_{\mathbb{S}^2} 
  + \int_{S_{\tau,R}} p^{-1} r^p |\phi|^2 \dVol_{\mathbb{S}^2} \right) \upd \tau
 \end{split}
\end{equation}
Again, if $\phi$ is a scalar field then the final term in the spacetime integral, involving the term 
\begin{equation*}
 \epsilon\chi_{(r_0)} r^{p-3}|\phi|^2
\end{equation*}
is not present.

\end{proposition}

\begin{proof}
 We apply the energy estimate to the field $\phi$ in the spacetime region $^t\mathcal{M}^\tau_{\tau_0}$, with the modified energy current $^{(L,p)}\tilde{J}[\phi]$. We also choose the function $f_L(r)$ to be
\begin{equation*}
 f_L(r) = \chi_{(2r_0, R)}(r)
\end{equation*}

 Using the modified compatible current identity \ref{proposition modified compatible current identity} we obtain
\begin{equation}
 \begin{split}
  &\int_{^t\Sigma_\tau} \imath_{^{(L,p)}\tilde{J}[\phi]} \dVol_g
  + \int_{^{\tau_0}_{\tau}\bar{\Sigma}_t} \imath_{^{(L,p)}\tilde{J}[\phi]} \dVol_g
  - \int_{^t\Sigma_{\tau_0}} \imath_{^{(L,p)}\tilde{J}[\phi]} \dVol_g \\
  &=\int_{^t\mathcal{M}^\tau_{\tau_0}} \bigg(
  {^{(L,p)}\tilde{K}}[\phi]
  - \omega f_L (\slashed{\D}_{\Lbar}\phi)\cdot \left( r^p\slashed{\D}_L \phi 
  + r^{p-1}\phi \right) 
  + f_L (F_1 + F_2 + F_3) \cdot \left( r^p\slashed{\D}_L \phi + r^{p-1}\phi \right) \\
  &\phantom{=\int_{^t\mathcal{M}^\tau_{\tau_0}} \bigg(}
  + f_L r^p L^\nu \left([\slashed{\D}_\mu, \slashed{\D}_\nu]\phi \right)\cdot \left(\slashed{\D}^\mu \phi \right) \bigg)\dVol_g
 \end{split}
\end{equation}
As before we shall first estimate the bulk terms, and then turn to the boundary terms.

First, using proposition \ref{proposition bulk p current p<1} we have
 \begin{equation*}
  \begin{split}
   {^{(L,p)}\tilde{K}}[\phi]
   - \omega f_L (\slashed{\D}_{\Lbar}\phi)\cdot \left( r^p\slashed{\D}_L \phi 
   + r^{p-1}\phi \right)
   &= \frac{1}{2}f_L r^{p-1} \left( p|\slashed{\D}_L \phi|^2 + (2-p)|\slashed{\nabla}\phi|^2 + p(1-p)r^{-2}|\phi|^2 \right) \\
   &\phantom{=}+ \textit{Err}_{(L,p,\text{bulk})}
  \end{split}
 \end{equation*}
We estimate the error terms in two different ways, depending on whether $p < \delta$ or $p \geq \delta$. First, for very small values of $p$ (specifically, $p \leq \delta$), we can use the bootstrap bounds from chapter \ref{chapter bootstrap} to bound
\begin{equation}
\label{equation p weighted internal 1}
 \begin{split}
  \left| \textit{Err}_{(L,p,\text{bulk})} \right|
  &\lesssim r^p f_L' |\overline{\slashed{\D}}\phi|^2
  + \left( r^{p-1}f_L'' + r^{p-2}f_L' \right)|\phi|^2
  + \epsilon f_L r^{p-1}(1+\tau)^{-\beta}|\slashed{\D}_L \phi| |\slashed{\D}\phi| \\
  &\phantom{\lesssim}
  + \epsilon f_L r^{p - 1 - \delta}|\slashed{\D}_L \phi| |\slashed{\D}\phi|
  + \epsilon f_L r^{p-1 + \delta}(1+\tau)^{-\beta}|\slashed{\D}_L \phi| |\slashed{\nabla}\phi|
  + \epsilon f_L r^{p-1}|\slashed{\nabla}\phi|^2 \\
  &\phantom{\lesssim}
  + \epsilon f_L (p-1)r^{p-3 + \delta} (1+\tau)^{-\beta} |\phi|^2
  + \epsilon f_L (p-1)r^{p-3} |\phi|^2
 \end{split}
\end{equation}

Now, we estimate the error term as follows:
\begin{equation*}
 \begin{split}  
  \left| \textit{Err}_{(L,p,\text{bulk})} \right|
  &\lesssim r^p f_L'|\overline{\slashed{\D}}\phi|^2 + (r^{p-1}f_L'' + r^{p-2}f_L')|\phi|^2 
  + \epsilon (1+r)^{-1-\delta}|\slashed{\D}\phi|^2 \\
  &\phantom{\lesssim}
  + \epsilon f_L r^{2p - 1 + \delta}(1+\tau)^{-2\beta}|\slashed{\D}_L \phi|^2
  + \epsilon f_L r^{2p - 1 - \delta}|\slashed{\D}_L \phi|^2
  + \epsilon f_L r^{p - 1}|\slashed{\nabla}\phi|^2 \\
  &\phantom{\lesssim}
  + \epsilon f_L r^{p-1+2\delta}(1+\tau)^{-\beta}|\slashed{\D}_L \phi|^2
  + \epsilon f_L (p-1)r^{p-3 + \delta} (1+\tau)^{-\beta} |\phi|^2
  + \epsilon f_L (p-1)r^{p-3} |\phi|^2
 \end{split}
\end{equation*}
We further decompose these error terms as follows. We have
\begin{equation*}
 \epsilon f_L r^{2p - 1 + \delta}(1+\tau)^{-2\beta}|\slashed{\D}_L \phi|^2
 \leq \epsilon f_L r^{2p - 1 - \delta}|\slashed{\D}_L \phi|^2
 + \epsilon f_L r^{2p - \delta}(1+\tau)^{-1-\delta}|\slashed{\D}_L \phi|^2
\end{equation*}
where, in fact, the inequality holds without the second term in the region $r \leq \tau$ (the ``interior''), and without the first term in the region $r \geq \tau$ (the ``exterior''), and where we have used that $\delta \ll \beta$. Moreover, we have $p \leq \delta$, so we actually obtain
\begin{equation*}
 \epsilon f_L r^{2p - 1 + \delta}(1+\tau)^{-2\beta}|\slashed{\D}_L \phi|^2
 \leq 
 \epsilon f_L r^{p - 1}|\slashed{\D}_L \phi|^2
 + \epsilon f_L r^{p}(1+\tau)^{-1-\delta}|\slashed{\D}_L \phi|^2
\end{equation*}

We also have
\begin{equation*}
  \epsilon f_L r^{2p -\delta}(1+\tau)^{-1-\delta}|\slashed{\D}_L \phi|^2
  \lesssim f_L r^{2p -2 -\delta}(1+\tau)^{-1-\delta}|\slashed{\D}_L \psi|^2 + f_L r^{2p -2 -\delta}(1+\tau)^{-1-\delta}|\phi|^2
\end{equation*}

The next error term satisfies
\begin{equation*}
  \epsilon f_L r^{2p - 1 - \delta}|\slashed{\D}_L \phi|^2
  \leq \epsilon f_L r^{p - 1}|\slashed{\D}_L \phi|^2
\end{equation*}
where we have used the fact that $p \leq \delta$.

Finally, we estimate
\begin{equation*}
   \epsilon f_L (p-1)r^{p-3 + \delta} (1+\tau)^{-\beta} |\phi|^2
   \leq \epsilon f_L (p-1)r^{p-3 -\delta}|\phi|^2 + \epsilon f_L (p-1)r^{p-2-\delta} (1+\tau)^{-1-\delta} |\phi|^2 
\end{equation*}
where, again, the first term on the right hand side is all that is needed in the interior, while the second term is all that is needed in the exterior, and we have used $\delta \ll \beta$.

This completes the calculations needed for the term $\textit{Err}_{(L,p,\text{bulk})}$ for very small values of $p$. In summary, we have shown that, for $p \leq \delta$,
\begin{equation*}
 \begin{split}  
  \left| \textit{Err}_{(L,p,\text{bulk})} \right|
  &\lesssim r^p f_L'|\overline{\slashed{\D}}\phi|^2 + (r^{p-1}f_L'' + r^{p-2}f_L')|\phi|^2
  + \epsilon (1+r)^{-1-\delta}|\slashed{\D}\phi|^2 \\
  &\phantom{\lesssim}
  + \epsilon f_L r^{p-1}|\slashed{\D}_L \phi|^2
  + \epsilon f_L r^{p-2} (1+\tau)^{-1-\delta}|\slashed{\D}_L \psi|^2
  + \epsilon f_L (p-1)r^{p-3} |\phi|^2 \\
  &\phantom{\lesssim}
  + \epsilon f_L r^{p-1}|\slashed{\nabla}\phi|^2
  + \epsilon f_L (p-1)r^{p-2 - \delta}(1+\tau)^{-1-\delta} |\phi|^2
 \end{split}
\end{equation*}

For larger values of $p$, we can estimate the error term \eqref{equation p weighted internal 1} as
\begin{equation*}
 \begin{split}
  \left| \textit{Err}_{(L, p, \text{bulk})} \right|
  &\lesssim r^p f_L' |\overline{\slashed{\D}}\phi|^2
  + (r^{p-1}f''_L + r^{p-2}f'_L )|\phi|^2
  + \epsilon f_L r^{p-2}(1+\tau)^{1-\delta}|\slashed{\D}\phi|^2 \\
  &\phantom{\lesssim}
  + \epsilon f_L r^p (1+\tau)^{-1-\beta}|\slashed{\D}_L\phi|^2
  + \epsilon f_L r^{p-2\delta} (1+\tau)^{-1+\delta} |\slashed{\D}_L\phi|^2 \\
  &\phantom{\lesssim}
  + \epsilon f_L r^{p-1+2\delta}(1+\tau)^{-2\beta}|\slashed{\D}_L\phi|^2
  + \epsilon f_L r^{p-1}|\slashed{\nabla}\phi|^2
  + \epsilon f_L (p-1)r^{p-3+\delta}(1+\tau)^{-\beta}|\phi|^2 \\
  &\phantom{\lesssim}
  + \epsilon f_L (p-1)r^{p-3}|\phi|^2
 \end{split}
\end{equation*}

% \begin{equation*}
%  \begin{split}  
%   \left| \textit{Err}_{(L,p,\text{bulk})} \right|
%   &\lesssim r^p f_L'|\overline{\slashed{\D}}\phi|^2 + (r^{p-1}f_L'' + r^{p-2}f_L')|\phi|^2 
%   + \epsilon (1+r)^{-1-\delta} (1+\tau)^{1-2\beta} |\slashed{\D}\phi|^2 \\
%   &\phantom{\lesssim}
%   + \epsilon f_L r^{2p - 1 + \delta}(1+\tau)^{-1}|\slashed{\D}_L \phi|^2
%   + \epsilon (1+r)^{-1-\delta} (1+\tau)^{1 - \delta} |\slashed{\D}\phi|^2 \\
%   &\phantom{\lesssim}
%   + \epsilon f_L r^{2p - 1 - (2\mathring{C} - 3)\delta} (1+\tau)^{-1 + \delta} |\slashed{\D}_L \phi|^2 
%   + \epsilon f_L r^{p-1}|\slashed{\nabla}\phi|^2 \\
%   &\phantom{\lesssim}
%   + \epsilon f_L r^{p-1+2\delta}(1+\tau)^{-2\beta}|\slashed{\D}_L \phi|^2
%   + \epsilon f_L (p-1)r^{p-3 + \delta} (1+\tau)^{-\beta} |\phi|^2
%   + \epsilon f_L (p-1)r^{p-3} |\phi|^2
%  \end{split}
% \end{equation*}

We also need to further decompose these terms. We have
\begin{equation*}
 \begin{split}
  \epsilon r^{p-2\delta}(1+\tau)^{-1+\delta}|\slashed{\D}_L\phi|^2
  &\lesssim \epsilon r^{p-1-\delta}|\slashed{\D}_L\phi|^2 
  + \epsilon r^{p}(1+\tau)^{-1-\delta}|\slashed{\D}_L\phi|^2
 \end{split}
\end{equation*}
where the first term on the right hand side bounds the left hand side in the interior ($r \leq \tau$) and the second term bounds it in the exterior ($r \geq \tau$).

Similarly, we can bound
\begin{equation*}
 \begin{split}
  r^{p-1+2\delta}(1+\tau)^{-2\beta}|\slashed{\D}_L\phi|^2
  &\lesssim \epsilon r^{p-1+2\delta - 2\beta}|\slashed{\D}_L\phi|^2
  + r^{p-2\beta + 3\delta}(1+\tau)^{-1 - \delta}|\slashed{\D}_L\phi|^2 \\
  &\lesssim \epsilon r^{p-1}|\slashed{\D}_L\phi|^2
  + r^{p}(1+\tau)^{-1 - \delta}|\slashed{\D}_L\phi|^2
 \end{split}
\end{equation*}

In each of the calculations above, we have a term of the form $r^p(1+\tau)^{-1-\delta}|\slashed{\D}_L\phi|^2$. We can further estimate this as
\begin{equation*}
 r^p(1+\tau)^{-1-\delta}|\slashed{\D}_L\phi|^2
 \lesssim r^{p-2}(1+\tau)^{-1-\delta}|\slashed{\D}_L\psi|^2 + r^{p-2}(1+\tau)^{-1-\delta}|\phi|^2
\end{equation*}
In fact, the lower order term with the critical decay in $r$, that is, $r^{p-2}(1+\tau)^{-1-\delta}|\phi|^2$, is absent. This term can only arise from the terms involving $\omega$, which decays at the critical rate $r^{-1}$. However, the relevant terms are
\begin{equation*}
 -\omega r^p f_L (\slashed{\D}_L\phi)\cdot(\slashed{\D}_{\Lbar}\phi) - \omega r^{p-1}f_l (\slashed{\D}_{\Lbar}\phi)\cdot \phi
 = -\omega r^{p-1} f_L (\slashed{\D}_L \psi)\cdot (\slashed{\D}_{\Lbar}\phi)
\end{equation*}
so the lower order term is absent. The remaining terms involving $\phi$ and not its derivatives have improved decay in $r$.

Finally, we estimate the term
\begin{equation*}
 \begin{split}
  r^{p-3+\delta}(1+\tau)^{-\beta}|\phi|^2
  & \lesssim r^{p-3+\delta-\beta}|\phi|^2
  + r^{p-2-\beta} (1+\tau)^{-1-\delta}|\phi|^2 \\
  & \lesssim r^{p-3}|\phi|^2
  + r^{p-2}(1+\tau)^{-1-\delta}|\phi|^2
 \end{split}
\end{equation*}

In summary, we have shown that
\begin{equation*}
 \begin{split}
  \left| \textit{Err}_{(L,p,\text{bulk})} \right|
  &\lesssim r^p f_L' |\overline{\slashed{\D}}\phi|^2
  + (r^{p-1}f''_L + r^{p-2}f'_L )|\phi|^2
  + \epsilon f_L r^{p-2}(1+\tau)^{1-\delta}|\slashed{\D}\phi|^2 \\
  &\phantom{\lesssim}
  + \epsilon f_L r^{p-2} (1+\tau)^{-1-\delta}|\slashed{\D}_L\psi|^2 
  + \epsilon f_L r^{p-2 - \frac{1}{2}\delta} (1+\tau)^{-1-\delta}|\phi|^2
  + \epsilon f_L r^{p-1}|\slashed{\D}_L\phi|^2 \\
  &\phantom{\lesssim}
  + \epsilon f_L r^{p-1}|\slashed{\nabla}\phi|^2
  + \epsilon f_L (p-1)r^{p-3}|\phi|^2
 \end{split}
\end{equation*}

We also need to estimate the terms involving the derivatives of $f_L$. We do this in the same way, regardless of the value of $p$. Using the coarea formula \ref{proposition coarea} we have
\begin{equation*}
 \begin{split}
  &\int_{^t\mathcal{M}_{\tau_0}^\tau} \left( r^p f_L'|\overline{\slashed{\D}}\phi|^2 + (r^{p-1}f_L'' + r^{p-2}f_L')|\phi|^2 \right)\dVol_g \\
  &= \int_{^t\mathcal{M}_{\tau_0}^\tau \cap \{r_0 \leq r \leq 2r_0\} } \left( r^p  \chi'_{(2r_0)} |\overline{\slashed{\D}}\phi|^2 + (r^{p-1} \chi''_{(2r_0)} + r^{p-2} \chi'_{(2r_0)})|\phi|^2 \right)\dVol_g \\
  &\phantom{=}
  + \int_{\tau_0}^\tau \left( \int_{^t\Sigma_{\tau'}} \left( -r^p \chi'_{(R)}|\overline{\slashed{\D}}\phi|^2 - (r^{p-1}\chi''_{(R)} + r^{p-2}\chi'_{(R)})|\phi|^2 \right) \Omega^2 \upd r \wedge \dVol_{\mathbb{S}^2} \right) \upd \tau' \\
  &\lesssim \int_{^t\mathcal{M}_{\tau_0}^\tau \cap \{r_0 \leq r \leq 2r_0\}} \left( |\overline{\slashed{\D}}\phi|^2 + |\phi|^2 \right) \dVol_g \\
  &\phantom{\lesssim} + \int_{\tau_0}^\tau \left( \int_{^t\Sigma_{\tau'}\cap \left\{\frac{1}{2}R \leq r \leq R \right\}} \left(r^{p-1}|\overline{\slashed{\D}}\phi|^2 + r^{p-3}|\phi|^2 \right) r^2 \upd r \wedge \dVol_{\mathbb{S}^2} \right) \upd \tau'
 \end{split}
\end{equation*}
where in the last line we have used the fact that $\Omega \sim r$ as well as
\begin{equation*}
 \begin{split}
  |\chi'_{(R)}| &\lesssim R^{-1} \\
  |\chi''_{(R)}| &\lesssim R^{-2}
 \end{split}
\end{equation*}
which follows from the definition of the cut off functions in equation \eqref{equation cut off functions}. Moreover, the derivatives of the cut-off funciton $\chi_{(R)}$ are supported only in the region $\frac{1}{2}R \leq r \leq R$. Hence we have, for example,
\begin{equation*}
 \chi'_{(R)}(r) \lesssim r^{-1}
\end{equation*}

We further estimate the zero-th order term by using the Hardy inequality (proposition \ref{proposition Hardy}) and also using the fact that $p < 1$:
\begin{equation*}
\begin{split}
  &\int_{\tau_0}^\tau \left( \int_{^t\Sigma_{\tau'}\cap \left\{\frac{1}{2}R \leq r \leq R \right\}}  r^{p-1}|\phi|^2  \upd r \wedge \dVol_{\mathbb{S}^2} \right) \upd \tau' \\
  &\lesssim \int_{\tau_0}^\tau \left( \int_{^t\Sigma_{\tau'}\cap \left\{\frac{1}{2}R \leq r \leq R \right\}}  p^{-2} r^{p-1}|\slashed{\D}_L \phi|^2 r^2 \upd r \wedge \dVol_{\mathbb{S}^2} + \int_{S_{\tau', R}} p^{-1} r^p |\phi|^2 \dVol_{\mathbb{S}^2} \right) \upd \tau'
 \end{split}
\end{equation*}
% and we bound this final term using proposition \ref{proposition spherical mean in terms of energy}, choosing $\alpha= \delta$:
% \begin{equation*}
%  \begin{split}
%   \int_{\tau_0}^\tau \left(  \int_{S_{\tau', R}} r^p |\phi|^2 \dVol_{\mathbb{S}^2} \right) \upd \tau'
%   &\lesssim \int_{\tau_0}^\tau \Bigg(  R^{p - 1 + \delta} \int_{^t\Sigma_{\tau'}} (1+r)^{-\delta} |\slashed{\D}_L\phi|^2 r^2 \upd r \wedge \dVol_{\mathbb{S}^2} \\
%   &\phantom{\lesssim \int_{\tau_0}^\tau \Bigg(}
%   + R^{p - 1 + \delta} \int_{{_{\tau_0}^{\tau'}\bar{\Sigma}}_{t}} (1+r)^{-\delta} |\slashed{\D}\phi|^2 r^2 \upd r \wedge \dVol_{\mathbb{S}^2} \\
%   &\phantom{\lesssim \int_{\tau_0}^\tau \Bigg(}
%   + \int_{\bar{S}_{\tau_0, t}} r^p |\phi|^2 \dVol_{\mathbb{S}^2}
%  \Bigg) \upd \tau'
%  \end{split}
% \end{equation*}

Next, we estimate the term involving the inhomogeneity $F_1 + F_2 + F_3$. We have
\begin{equation*}
 \begin{split}
 f_L (F_1 + F_2 + F_3) \cdot \left( r^p\slashed{\D}_L \phi + r^{p-1}\phi \right)
 &\lesssim
 \epsilon f_L r^{p-3 + 4 \delta} (1+\tau)^{-6\delta} \left| \slashed{\D}_L \psi \right|^2
 + \epsilon f_L r^{p-2}(1+\tau)^{-1-\delta} |\slashed{\D}_L \psi|^2
 \\
 &\phantom{\lesssim}
 + \epsilon f_L r^{p-1} |\slashed{\D}_L \phi|^2
 + \epsilon f_L r^{p-3} |\phi|^2
 + \epsilon^{-1} f_L r^p (1+\tau)^{1+\delta} |F_1|^2
 \\
 &\phantom{\lesssim}
 + \epsilon^{-1} f_L r^{p + 1 - 4 \delta} (1+\tau)^{6\delta} |F_2|^2
 + \epsilon^{-1} f_L r^{p+1} |F_3|^2
 \end{split}
\end{equation*}
Similarly to before, we have
\begin{equation*}
 \begin{split}
 \epsilon f_L r^{p-3 + 4\delta} (1+\tau)^{-6\delta} \left| \slashed{\D}_L \psi \right|^2
 &\lesssim \epsilon f_L r^{p - 1}|\slashed{\D}_L \phi|^2
 + \epsilon f_L r^{p-3-\delta}|\phi|^2 \\
 &\phantom{\lesssim}+ \epsilon f_L r^{p-2}(1+\tau)^{-1-\delta}|\slashed{\D}_L \psi|^2
 + \epsilon f_L r^{p-2-\delta}(1+\tau)^{-1-\delta}|\phi|^2
 \end{split}
\end{equation*}

Next, we have to estimate the error term involving curvature terms, which is not present if $\phi$ is a scalar field. If $r \geq r_0$ then we have
\begin{equation*}
 f_L r^p L^\nu \left( [\slashed{\D}_\mu , \slashed{\D}_\nu]\phi \right) \cdot \slashed{\D}^\mu \phi
 \lesssim f_L r^p |\Omega_{L\Lbar}| |\slashed{\D}_L \phi| |\phi|
 + f_L r^p |\Omega_{L\slashed{\alpha}}| |\slashed{\nabla} \phi| |\phi|
\end{equation*}
Note the absence of a ``bad derivative'' term $\slashed{\D}\phi$, which follows from the fact that
\begin{equation*}
 L^\mu L^\nu [\slashed{\D}_\mu , \slashed{\D}_\nu]\phi = 0
\end{equation*}
Hence, using the expressions in chapter \ref{chapter geometry of vector bundle} we obtain
\begin{equation*}
f_L r^p L^\nu \left( [\slashed{\D}_\mu , \slashed{\D}_\nu]\phi \right) \cdot \slashed{\D}^\mu \phi
\lesssim
\epsilon f_L r^{p-1}|\overline{\slashed{\D}}\phi|^2
+ \epsilon f_L r^{p-3}|\phi|^2 
+ \epsilon f_L r^{p-2-\delta}(1+\tau)^{-1-\delta}|\phi|^2
\end{equation*}

where we have used the fact that $\delta \ll \beta$.

We also have to deal with this term in the region $\frac{1}{2}r_0 \leq r \leq r_0$ (recall that $f_L \equiv 0$ in the region $r \leq \frac{1}{2}r_0$). Here, the bootstrap bounds lead easily to the estimate
\begin{equation*}
 f_L r^p L^\nu \left( [\D_\mu , \D_\nu]\phi \right) \cdot \slashed{\D}^\mu \phi
\lesssim \epsilon \left( |\phi|^2 + |\slashed{\D}\phi|^2 \right)
\end{equation*}

This concludes the required estimates for the bulk terms; we must now estimate the boundary terms arising in the $p$-weighted estimates. Using proposition \ref{proposition boundary p} together with the fact that $f_L$ is supported away from $^{\tau}_{\tau_0}\bar{\Sigma}_t$ we have
\begin{equation*}
 \begin{split}
  \int_{^t\Sigma_\tau} \imath_{^{(L,p)}\tilde{J}[\phi]} \dVol_g
  &= \int_{^t\Sigma_\tau} f_L r^{p-2}\left( |\slashed{\D}_L \psi|^2 + \textit{Err}_{(L,p-\text{bdy})} \right) \Omega^2\upd r \wedge \dVol_{\mathbb{S}^2} \\
 \end{split}
\end{equation*}
where 
\begin{equation*}
 \begin{split}
  \left| \textit{Err}_{(L,p-\text{bdy})} \right|
  &\lesssim 
  + \left| \left.\frac{\partial}{\partial r}\right|_{\tau,\vartheta^1, \vartheta^2} \log \Omega - r^{-1} \right| r|\phi|^2 
 \end{split}
\end{equation*}
Now, we have
\begin{equation*}
 \begin{split}
  \left.\frac{\partial}{\partial r}\right|_{\tau,\vartheta^1, \vartheta^2} \log \Omega - r^{-1}
  &= L\log \Omega - r^{-1} \\
  &= \frac{1}{2}\tr_{\slashed{g}}\chi - r^{-1} \\
  &= \frac{1}{2}\tr_{\slashed{g}}\chi_{(\text{small})}
 \end{split}
\end{equation*}
so, using the bootstrap bounds we find that
\begin{equation*}
 + \left| \left.\frac{\partial}{\partial r}\right|_{\tau,\vartheta^1, \vartheta^2} \log \Omega - r^{-1} \right| r|\phi|^2
 \lesssim \epsilon r (1+r)^{-1-\delta} |\phi|^2
\end{equation*}

Now, we need to further estimate the corresponding error term, which (using the support of $f_L$) can be bounded by
\begin{equation*}
 \int_{^t \Sigma_\tau} \epsilon f_L  r^{p -\delta} |\phi|^2 \upd r \wedge \dVol_{\mathbb{S}^2}
\end{equation*}
Again, this term will be estimated in different ways depending on whether $p < \delta$ or $p > \delta$.

If $p < \delta$ then we estimate this term in terms of the $T$ energy. Specifically, using the fact that $0 < f_L < 1$ we can use the first part of proposition \ref{proposition Hardy} to show
\begin{equation*}
 \begin{split}
 \int_{^t \Sigma_\tau} \epsilon f_L  r^{p -\delta} |\phi|^2 \upd r \wedge \dVol_{\mathbb{S}^2}
 &\lesssim 
 \frac{1}{(1 - p + \delta)^2} \int_{^t \Sigma_\tau} \epsilon (1+r)^{p -\delta} |\slashed{\D}_L\phi|^2 r^2 \upd r \wedge \dVol_{\mathbb{S}^2} \\
 &\phantom{\lesssim}
 +  \frac{\epsilon}{(1 - p + \delta)}\int_{\bar{S}_{\tau,t}}(1+r)^{p + 1 -\delta}|\phi|^2 \, \dVol_{\mathbb{S}^2}
 \end{split}
\end{equation*}

% The boundary term on $\bar{S}_{\tau,r}$ can itself be bounded by the $T$ energy. Using the bootstrap bounds and proposition \ref{proposition spherical mean in terms of energy} we have
% \begin{equation*}
%  \int_{\bar{S}_{\tau,t}}(1+r)^{p-\mathring{C}\delta}|\phi|^2 r \, \dVol_{\mathbb{S}^2}
%  \lesssim \int_{^\tau_{\tau_0}\bar{\Sigma}_t} (1+r)^{p-\mathring{C}\delta}|\slashed{\D}\phi|^2 r^2 \upd r \wedge \dVol_{\mathbb{S}^2} + \int_{\bar{S}_{\tau_0, t}}(1+r)^{p-\mathring{C}\delta}|\phi|^2 \dVol_{\mathbb{S}^2}
% \end{equation*}

On the other hand, if $p > \delta$ then we make use of the second part of proposition \ref{proposition Hardy}, applied to the field $\psi = r\phi$. We choose $\alpha = 2 + \delta - p$. Since $p < 1$, we have $\alpha > 1$. We obtain
\begin{equation*}
\begin{split}
 &\int_{^t\Sigma_\tau} \epsilon f_L  r^{p - 2 -\delta} |\psi|^2 \upd r \wedge \dVol_{\mathbb{S}^2} \\
 &\lesssim
 \int_{^t\Sigma_\tau \cap \{r \leq R\}} \epsilon \left(1 - \chi_{(2r_0)}\right) r^{p -2- \delta} | \psi|^2  \upd r \wedge \dVol_{\mathbb{S}^2} \\
 &\lesssim
 \int_{^t\Sigma_\tau \cap \left\{r_0 \leq r \leq 2r_0 \right\} } \epsilon (r_0)^{-1} |\chi'_{(2r_0)}| r^{p +1 - \delta} |\phi|^2  \upd r \wedge \dVol_{\mathbb{S}^2} \\
 &\phantom{\lesssim} 
 + \frac{1}{(1 - p + \delta)^2}\int_{^t\Sigma_\tau \cap \left\{ r \leq R \right\} } \epsilon r^{p - \delta} |\slashed{\D}_L\psi|^2  \upd r \wedge \dVol_{\mathbb{S}^2}
 +  \frac{1}{(1 - p + \delta)}\int_{S_{\tau,R}} \epsilon r^{p+1 - \delta} |\phi|^2 \dVol_{\mathbb{S}^2} 
 \end{split}
\end{equation*}

%Since $r_0$ is some fixed constant, we can use the Hardy inequality to bound
%\begin{equation*}
% \int_{^t\Sigma_\tau \cap \left\{\frac{1}{2}r_0 \leq r \leq r_0 \right\} } \epsilon (r_0)^{-1} |\chi'_{(2r_0)}| r^{p +1 - \mathring{C}\delta} |\phi|^2  \upd r \wedge \dVol_{\mathbb{S}^2}
% \lesssim \epsilon \, \mathcal{E}^{(T)}[\phi] (\tau, t, \tau_0) + \epsilon \int_{\bar{S}_{\tau_0, t}} |\phi|^2 \dVol_{\mathbb{S}^2}
%\end{equation*}
%where, as in many of the  previous inequalities, the implicit constant depends on $r_0$.

\end{proof}

\chapter{Boundedness and energy decay}
\label{chapter boundedness and energy decay}

In this chapter we will combine the basic energy estimates (established in the previous chapter) to prove bounds on solutions to inhomogenous wave equations on manifolds whose metric components and connection coefficients obey the bootstrap bounds of chapter \ref{chapter bootstrap}. Note that these bounds also apply to linear wave equations on such manifolds.

The approach we take will be to prove progressively better and better decay in $\tau$ for a \emph{degenerate} energy, that is, a quantity which is similar to the energy, but includes an additional factor of the form $(1+r)^{-C\epsilon}$, meaning that this energy degenerates as $r \rightarrow \infty$. The necessity for using this kind of degenerate energy comes from the presence of certain error terms in the energy estimates, which decay at the critical rate $(1+r)^{-1}$, and can only be controlled by the use of a Gronwall inequality, rather than being absorbed by some ``bulk terms'' associated with the energy currents.

As mentioned above, we will aim to prove progressively stronger decay in $\tau$ for this degenerate energy. We will begin by using the weighted $T$ energy estimate to prove that the degenerate energy grows at most exponentially fast in $\tau$. Note that, in order to prove this kind of bound for a linear wave equation with respect to a foliation by \emph{uniformly spacelike leaves} (e.g. with respect to the foliation given by surfaces of constant $t$) it is only necessary that the metric components and their derivatives be bounded, and moreover, it is the full energy, rather than the degenerate energy, that can be controlled\footnote{Note, however, that in this case the energy grows like $t^\epsilon$, where $t$ parameterises the leaves of the foliation.}. This is because the error terms arising in the $T$ energy estimate can be controlled with the help of the boundary terms in the $T$ energy, which involve all the derivatives of $\phi$ when evaluated on a spacelike hypersurface, together with a Gronwall inequality. However, we require our estimates to be given with respect to an (asymptotically) null foliation, and the $T$ energy evaluated on these leaves only involves the ``good derivatives''. Hence it cannot be used to bound bulk error terms involving the bad derivatives.

If the equations we were dealing with obeyed a stronger version of the null condition then these error terms could be handled with the use of an integrated local energy decay estimate. For example, if the equations obeyed the classical null condition of \cite{Klainerman1980} or if the coefficient of the ``good metric term'' $h_{LL}$ depended only on the \emph{good derivatives} of $\phi$, then the bulk error terms in the $T$ energy would be of the form $(1+r)^{-1-\delta}|\slashed{\D}\phi|^2$, and these can be controlled by the Morawetz energy current. In our case, however, these terms decay at the critical rate of $(1+r)^{-1}$, and so extra work is required even to conclude that the energy grows at most exponentially. In fact, this is accomplished by the inclusion of the weight $w = (1+r)^{-C\epsilon}$, which generates a positive bulk term in the weighted $T$ energy estimate involving the bad derivatives term $|\slashed{\D}_{\Lbar}\phi|^2$. Note that there are also other bulk error terms which cannot be controlled by the Morawetz estimate, since they decay at a supercritical rate $(1+r)^{-1 + C\delta}$. Fortunately, all of these error terms only involve the \emph{good derivatives}, so these can also be controlled by appealing to the Gronwall inequality. However, this leads to exponential growth in $\tau$ for the degenerate energy! Such an estimate is very far from providing sufficient control to recover the bootstrap bounds, which requires the energy to decay at least at the rate $(1+\tau)^{-2\beta}$.

Next, we need to incorporate the $p$ weighted energy estimates in order to improve the decay in $\tau$ of the degenerate energy. Note that the $p$-weighted estimates contribute a positive bulk term involving the good derivatives, and with a weight which behaves like $(1+r)^{-1 + p}$. Hence, by adding the $p$-weighted estimate with $p \sim C\delta$ to the Morawetz and $T$ energy estimates, we can absorb the error terms which previously led to exponential growth in $\tau$. Moreover, there are certain error terms in the $p$-weighted energy estimates which can be controlled if the degenerate energy is known \emph{a priori} to be finite (this is the reason for first considering only the weighted $T$ estimate). Putting these estimates together will lead to boundedness of the degenerate energy. Additionally, these estimates will allow us to prove a version of ``integrated local energy decay''; we will be able to show that a spacetime integral of a weighted, energy-type quantity (including, importantly, the ``bad derivatives'') is also bounded by uniformly in $\tau$.

These estimates are still not sufficiently powerful to be useful in closing our bootstrap bounds. Indeed, we need to show that the energy decays at a rate $(1+\tau)^{-2\beta}$, and so far we have only shown that it is bounded. In order to upgrade our estimates, we need to include the $p$-weighted estimates with higher values of $p$. Note that there are error terms in the $p$-weighted estimates, with $p = P \geq C\delta$, which tend to zero \emph{if} the $p$-weighted energy with $p = P - C\delta$ is already known to be finite. Hence we need to prove the $p$-weighted estimates for all values of $p$ in the interval $(C\delta, 1 - C\epsilon)$. Once we have proved this, these estimates can be used to show that the degenerate energy decays at a rate  of the form $(1+\tau)^{-1 + C\epsilon}$, which is sufficiently fast to be useful in closing the bootstrap bounds.

\section{Exponential growth of the degenerate energy and integrated local energy}

In this section we will use the weighted $T$ energy estimate to prove that the degenerate energy grows at most exponentially in $\tau$.

\begin{lemma}[Exponential growth of the degenerate energy and integrated local energy]
\label{lemma exponential growth}
 Let $\phi$ be an $S_{\tau,r}$-tangent tensor field satisfying
\begin{equation*}
 \tilde{\slashed{\Box}}_g \phi = F
\end{equation*}
for some $S_{\tau, r}$ tangent tensor field $F$. Suppose that $F$ satisfies the bounds
\begin{equation*}
 \begin{split}
  &\int_{^t\mathcal{M}^\tau_{\tau_0}} w \left(\epsilon^{-1}(1+r) |F|^2\right)\dVol
  \lesssim \tilde{\mathcal{E}} \exp \left( C_{[\phi]}\epsilon (\tau - \tau_0) \right)
 \end{split}
\end{equation*}
where
\begin{equation*}
 w = (1+r)^{-C_{[\phi]}\epsilon}
\end{equation*}
and $\tilde{\mathcal{E}}$ is some constant.

%Moreover, if the rank of $\phi$ is at least one (so $\phi$ is not a scalar field), then assume that
%\begin{equation*}
% \begin{split}
%  &\int_{^t\mathcal{M}^\tau_{\tau_0}} w\left( \epsilon (1+r)^{-2 + \delta} (1+\tau)^{-1-\delta}|\phi|^2 + \epsilon^{-1}(1+r)^{-3} |\phi|^2 \right) \dVol_g
% \lesssim \mathcal{E}_0 \exp \left( \tilde{C}\epsilon (\tau - \tau_0) \right)
% \end{split}
%\end{equation*}
%for some constant $\mathcal{E}_0$.

Suppose that all the bootstrap bounds of chapter \ref{chapter bootstrap} are satisfied. Suppose, moreover, that the initial energy for $\phi$ satisfies
\begin{equation*}
 \mathcal{E}^{(T, C_{[\phi]}\epsilon)}(\tau_0, t, \tau_0) \lesssim \mathcal{E}_0
\end{equation*}

Finally, suppose $\tilde{C}$ is sufficiently large compared to $C_{[\phi]}$.

Then, for all sufficiently small $\delta$, for all sufficiently small $\epsilon$ and for all sufficiently large $\tilde{C}$ we have

\begin{equation}
 \mathcal{E}^{(T, C_{[\phi]}\epsilon)}(\tau, t, \tau_0)
 + \int_{^t\mathcal{M}^\tau_{\tau_0}} C_{[\phi]}\epsilon w(1+r)^{-1}|\slashed{\D}\phi|^2 \dVol_g
 \lesssim \left(\mathcal{E}_0 + \tilde{\mathcal{E}}\right) \exp \left( 2\tilde{C}\epsilon (\tau - \tau_0) \right)
\end{equation}
\end{lemma}

\begin{proof}

Using the basic weighted $T$ energy estimate \ref{proposition basic weighted T energy}, for sufficiently large $C_{[\phi]}$, together with the bounds assumed in the lemma, we have
\begin{equation*}
 \begin{split}
  &\mathcal{E}^{(T, C_{[\phi]}\epsilon)}[\phi](\tau, t, \tau_0)
  + \int_{^t\mathcal{M}^\tau_{\tau_0}} C_{[\phi]}\epsilon w(1+r)^{-1}|\slashed{\D}\phi|^2 \dVol_g \\
  &\lesssim \mathcal{E}^{(T, C_{[\phi]}\epsilon)}[\phi](\tau_0, t, \tau_0)
  + \int_{^t\mathcal{M}^\tau_{\tau_0}} \left( C_{[\phi]} \epsilon w(1+r)^{-1 + \delta}|\overline{\slashed{\D}}\phi|^2 \right) \dVol_g
  + (\mathcal{E}_0 + \tilde{\mathcal{E}})\exp\left( C_{[\phi]} \epsilon(\tau - \tau_0) \right)
 \end{split}
\end{equation*}

Now, we can write
\begin{equation*}
 \begin{split}
  &\int_{^t\mathcal{M}^\tau_{\tau_0}} \left( C_{[\phi]} \epsilon w(1+r)^{-1 + \delta}|\overline{\slashed{\D}}\phi|^2 \right) \dVol_g \\
  &= \int_{\tau_0}^{\tau} \left( \int_{^t\Sigma_{\tau'}} C_{[\phi]}\epsilon (1+r)^{-1 +\delta - C_{[\phi]}\epsilon }|\overline{\slashed{\D}}\phi|^2 \Omega^2 \upd r \wedge \dVol_{\mathbb{S}^2} \right) \upd \tau' \\
  &\lesssim \int_{\tau_0}^{\tau} \left( \int_{^t\Sigma_{\tau'}} C_{[\phi]}\epsilon (1+r)^{-1 +\delta - C_{[\phi]}\epsilon }|\overline{\slashed{\D}}\phi|^2 r^2 \upd r \wedge \dVol_{\mathbb{S}^2} \right) \upd \tau' \\
  &\lesssim \int_{\tau_0}^{\tau} \left( \int_{^t\Sigma_{\tau'}} C_{[\phi]}\epsilon w |\overline{\slashed{\D}}\phi|^2 r^2 \upd r \wedge \dVol_{\mathbb{S}^2} \right) \upd \tau'
 \end{split}
\end{equation*}

Hence we have
\begin{equation*}
 \begin{split}
  &\int_{^t\Sigma_{\tau}}  w |\overline{\slashed{\D}}\phi|^2 r^2 \upd r \wedge \dVol_{\mathbb{S}^2}
  + \int_{^\tau_{\tau_0}\bar{\Sigma}_t}  w |\slashed{\D}\phi|^2 r^2 \upd r \wedge \dVol_{\mathbb{S}^2}
  + \int_{^t\mathcal{M}^\tau_{\tau_0}} C_{[\phi]}\epsilon w(1+r)^{-1}|\slashed{\D}\phi|^2 \dVol_g \\
  &\lesssim
  \int_{^t\Sigma_{\tau}}  w |\overline{\slashed{\D}}\phi|^2 r^2 \upd r \wedge \dVol_{\mathbb{S}^2}
  + \int_{\tau_0}^{\tau} \left( \int_{^t\Sigma_{\tau'}} \epsilon w |\overline{\slashed{\D}}\phi|^2 r^2 \upd r \wedge \dVol_{\mathbb{S}^2} \right) \upd \tau'
  + (\mathcal{E}_0 + \tilde{\mathcal{E}})\exp\left( C_{[\phi]} \epsilon (\tau - \tau_0) \right)
 \end{split}
\end{equation*}
so, using the Gronwall inequality of proposition \ref{proposition Gronwall} with the choices
\begin{equation*}
 \begin{split}
  f(\tau) &= \int_{^t\Sigma_{\tau'}} \epsilon w |\overline{\slashed{\D}}\phi|^2 r^2 \upd r \wedge \dVol_{\mathbb{S}^2} \\
  h(\tau) &= \int_{^\tau_{\tau_0}\bar{\Sigma}_t}  w |\slashed{\D}\phi|^2 r^2 \upd r \wedge \dVol_{\mathbb{S}^2} + \int_{^t\mathcal{M}^\tau_{\tau_0}} C_{[\phi]}\epsilon w(1+r)^{-1}|\slashed{\D}\phi|^2 \dVol_g\\
  g(\tau) &= C \\
	G(\tau, \tau_0) &= (\mathcal{E}_0 + \tilde{\mathcal{E}})\exp\left( \tilde{C} \epsilon (\tau - \tau_0) \right)
 \end{split}
\end{equation*}
where $C$ is the implicit constant in the preceding inequality. We find
\begin{equation*}
 \begin{split}
  &\int_{^t\Sigma_{\tau}}  w |\overline{\slashed{\D}}\phi|^2 r^2 \upd r \wedge \dVol_{\mathbb{S}^2}
  + \int_{^\tau_{\tau_0}\bar{\Sigma}_t}  w |\slashed{\D}\phi|^2 r^2 \upd r \wedge \dVol_{\mathbb{S}^2}
  + \int_{^t\mathcal{M}^\tau_{\tau_0}} C_{[\phi]}\epsilon w(1+r)^{-1}|\slashed{\D}\phi|^2 \dVol_g
  \\
  &\lesssim e^{C\epsilon (\tau - \tau_0)}\left( \mathcal{E}^{(T, C_{[\phi]}\epsilon)}[\phi](\tau_0, t, \tau_0) + \tilde{\mathcal{E}}e^{\tilde{C}\epsilon(\tau - \tau_0)} \right) \\
  &\lesssim \left( \mathcal{E}_0 + \tilde{\mathcal{E}} \right) \exp\left( (C + \tilde{C})\epsilon(\tau - \tau_0) \right)
 \end{split}
\end{equation*}
so, choosing $\tilde{C} \geq 2C$ proves the lemma.

\end{proof}

\section{Boundedness of the degenerate energy, the integrated local energy, and the \texorpdfstring{$p$}{p}-weighted energy estimate for very small \texorpdfstring{$p$}{p}}

\begin{lemma}[Energy boundedness, the integrated local energy decay estimate and the $p$-weighted energy estimate for very small $p$]
\label{lemma boundedness}
 Let $\phi$ be an $S_{\tau,r}$-tangent tensor field satisfying
\begin{equation*}
 \tilde{\slashed{\Box}}_g \phi = F = F_1 + F_2 + F_3
\end{equation*}
for some $S_{\tau, r}$ tangent tensor fields $F_1$, $F_2$ and $F_3$. Let $\tau_2 \geq \tau_1$. Suppose that these satisfy the bound
\begin{equation*}
 \begin{split}
  &\int_{^t\mathcal{M}^{\tau_2}_{\tau_1}} \epsilon^{-1}\bigg( 
  	(1+r)^{1-C_{[\phi]}\epsilon} |F|^2
  	+ (1+r)^{\frac{1}{2}\delta}(1+\tau)^{1+\delta}|F_1|^2 
  	+ (1+r)^{1-3\delta}(1+\tau)^{6\delta}|F_2|^2
  	\\
  	&\phantom{\int_{^t\mathcal{M}^{\tau_2}_{\tau_1}} \epsilon^{-1}\bigg(}
  	+ (1+r)^{1+\frac{1}{2}\delta}|F_3|^2 
  	\bigg)\dVol_g
  \lesssim \tilde{\mathcal{E}}_1 \\
 \end{split}
\end{equation*}

%Moreover, if the rank of $\phi$ is at least one (so $\phi$ is not a scalar field), then assume that
%\begin{equation*}
% \begin{split}
%  \int_{^t\mathcal{M}^{\tau^2}_{\tau_1}} \epsilon^{-1}\left( r^{-1}(1+r)^{-2 + C_{[\phi]}\epsilon}|\phi|^2 + \chi_{(r_0, R)} r^{-3+\frac{1}{2}\delta}|\phi|^2 \right) \dVol_g
%  \lesssim \mathcal{E}_1
% \end{split}
%\end{equation*}
%and
%\begin{equation*}
% \int_{^t\mathcal{M}^{\tau^2}_{\tau_1}\cap\{r \leq r_0\}} r^{-2}\left( \epsilon^{-1}|\phi|^2 + \epsilon |\slashed{\D}_T \phi|^2 \right)\dVol_g
% \lesssim \mathcal{E}_1
%\end{equation*}

Choose the weight function
\begin{equation*}
 w = (1+r)^{-C_{[\phi]}\epsilon}
\end{equation*}

Suppose, moreover, that the initial energy of $\phi$ satisfies
\begin{equation*}
 \mathcal{E}^{(wT)}(\tau_2, t, \tau_1) \lesssim \mathcal{E}_1
\end{equation*}
and
\begin{equation*}
 \delta^3 \mathcal{E}^{(L, \frac{1}{2}\delta)}(\tau_1,R) \lesssim \mathcal{E}_1
\end{equation*}
where $t$ is sufficiently large relative to $R$, $\tau_1$ and $\tau_2$ so that
\begin{equation}
 \label{equation relationship between R and t}
 \{r \leq R\} \cap \left\{^{\tau_2}_{\tau_1}\bar{\Sigma}_t \right\} = \emptyset
\end{equation}

Define the modified weight
\begin{equation*}
 \tilde{w} = (1+r)^{-\frac{1}{2}C_{[\phi]}\epsilon}
\end{equation*}
Suppose additionally that there is some $\tau_0 \leq \tau_1$ such that
\begin{equation*}
 \mathcal{E}^{(\tilde{w}T)}(\tau_0, t, \tau_0) \lesssim \mathcal{E}_0
\end{equation*}
Furthermore, on the initial hypersurface $^t\Sigma_{\tau_0}$ suppose that we have
\begin{equation*}
 \int_{\bar{S}_{t, r}}|\phi|^2 \dVol_{\mathbb{S}^2} \lesssim \mathcal{E}_0 (t - \tau_0)^{-1 + \frac{1}{2}C_{[\phi]}\epsilon}
\end{equation*}

Finally, suppose that all the bootstrap bounds of chapter \ref{chapter bootstrap} are satisfied. 

Then, for all sufficiently small $\delta$, for all sufficiently small $\epsilon$ we have
\begin{equation}
 \begin{split}
  &\delta^3 \mathcal{E}^{(L, \frac{1}{2}\delta)}[\phi](\tau_2, R)
  + \mathcal{E}^{(wT)}[\phi](\tau_2, t, \tau_1) \\
  &
  + \int_{^t\mathcal{M}^{\tau_2}_{\tau_1}} \bigg( \delta^2 (1+r)^{-1-\frac{1}{2}\delta}|\slashed{\D}\phi|^2 + \delta \chi_{(r_0)}(r)(1+r)^{-1 - C_{[\phi]}\epsilon}|\slashed{\nabla}\phi|^2 + \delta^2 r^{-1}(1+r)^{-2-\frac{1}{2}\delta}|\phi|^2 \\
  &\phantom{+ \int_{^t\mathcal{M}^{\tau_2}_{\tau_1}} \bigg(}
  + C_{[\phi]}\epsilon(1+r)^{-1-C_{[\phi]}\epsilon}|\slashed{\D}\phi|^2 \\
  &\phantom{+ \int_{^t\mathcal{M}^{\tau_2}_{\tau_1}} \bigg(}
  + \chi_{(2r_0, R)}\left( \delta^4 r^{-1+\frac{1}{2}\delta}|\slashed{\D}_L \phi|^2 + \delta^3 r^{-1 + \frac{1}{2}\delta|}\slashed{\nabla}\phi|^2 + \delta^3 r^{-3 + \frac{1}{2}\delta}|\phi|^2 \right) \bigg)\dVol_g \\
  &\lesssim
  \delta^{-1} \mathcal{E}_1 
  + \delta^{-1} \tilde{\mathcal{E}}_1
  \\
  &\phantom{\lesssim}
  + \left( \mathcal{E}_0 + \tilde{\mathcal{E}}_0 \right)\left( (t - \tau)^{-\frac{1}{2}C_{[\phi]}\epsilon} + (\tau - \tau_1)(t - \tau_0)^{-1 + \frac{3}{2}\delta} + \epsilon^{-1}R^{-1 + \frac{3}{2}\delta}e^{\tilde{C}\epsilon(\tau - \tau_0)} \right)
 \end{split}
\end{equation}

In particular, suppose that the conditions of the lemma hold in the limit $R \rightarrow \infty$, $t \rightarrow \infty$, where the limit is taken such that \eqref{equation relationship between R and t} is true. For example, we could take $t = 2(\tau_2 + R)$ which, under the bootstrap assumptions, can be seen to imply \eqref{equation relationship between R and t}. We could then take the limit $R \rightarrow \infty$, with $t$ considered a function of $R$, and $\tau_0$, $\tau_1$, $\tau_2$ fixed (finite) constants. Then we define the energies
\begin{equation}
 \begin{split}
  \mathcal{E}^{(wT)}[\phi](\tau) &:= \lim_{t \rightarrow \infty} \mathcal{E}^{(wT)}[\phi](\tau, t, \tau) \\
  \mathcal{E}^{(L, \alpha)}[\phi](\tau) &:= \lim_{R \rightarrow \infty} \mathcal{E}^{(L, \alpha)}[\phi](\tau, R) \\ 
 \end{split}
\end{equation}
Then we have the bound
\begin{equation}
 \begin{split}
  &\delta^3 \mathcal{E}^{(L, \frac{1}{2}\delta)}[\phi](\tau_2)
  + \mathcal{E}^{(wT)}[\phi](\tau_2) \\
  &
  + \int_{\mathcal{M}^{\tau_2}_{\tau_1}} \bigg( \delta^2 (1+r)^{-1-\frac{1}{2}\delta}|\slashed{\D}\phi|^2 + \delta \chi_{(r_0)}(r)(1+r)^{-1 - C_{[\phi]}\epsilon}|\slashed{\nabla}\phi|^2 + \delta^2 r^{-1}(1+r)^{-2-\frac{1}{2}\delta}|\phi|^2 \\
  &\phantom{+ \int_{^t\mathcal{M}^{\tau_2}_{\tau_1}} \bigg(}
  + C_{[\phi]}\epsilon(1+r)^{-1-C_{[\phi]}\epsilon}|\slashed{\D}\phi|^2 \\
  &\phantom{+ \int_{^t\mathcal{M}^{\tau_2}_{\tau_1}} \bigg(}
  + \chi_{(2r_0)}\left(\delta^4 r^{-1+\frac{1}{2}\delta}|\slashed{\D}_L \phi|^2 + \delta^3 r^{-1 + \frac{1}{2}\delta|}\slashed{\nabla}\phi|^2 + \delta^3 r^{-3 + \frac{1}{2}\delta}|\phi|^2 \right) \bigg)\dVol_g \\
  &\lesssim
  \delta^{-1} \mathcal{E}_1 
  + \delta^{-1} \tilde{\mathcal{E}}_1
 \end{split}
\end{equation}
In particular, this implies the following three estimates:
\begin{enumerate}
 \item Degenerate energy boundedness:
 \begin{equation*}
  \mathcal{E}^{(wT)}[\phi](\tau_2)
  \lesssim
  \delta^2 \mathcal{E}^{(L, \frac{1}{2}\delta)}[\phi](\tau_1)
  + \delta^{-1} \mathcal{E}^{(wT)}[\phi](\tau_1)
  + \delta^{-1} \tilde{\mathcal{E}}_1
 \end{equation*}
 \item Integrated local energy decay:
 \begin{equation*}
  \begin{split}
    &\int_{\mathcal{M}^{\tau_2}_{\tau_1}} \bigg( \delta^2 (1+r)^{-1-\frac{1}{2}\delta}|\slashed{\D}\phi|^2 + \delta \chi_{(r_0)}(r)(1+r)^{-1-C_{[\phi]}\epsilon}|\slashed{\nabla}\phi|^2 + \delta^2 r^{-1}(1+r)^{-2-\frac{1}{2}\delta}|\phi|^2 \\
    &\phantom{+ \int_{^t\mathcal{M}^{\tau_2}_{\tau_1}} \bigg(}
    + C_{[\phi]}\epsilon(1+r)^{-1-C_{[\phi]}\epsilon}|\slashed{\D}\phi|^2 \bigg)\dVol_g \\
    &\lesssim
    \delta^2 \mathcal{E}^{(L, \frac{1}{2}\delta)}[\phi](\tau_1)
    + \delta^{-1} \mathcal{E}^{(wT)}[\phi](\tau_1)
    + \delta^{-1} \tilde{\mathcal{E}}_1
  \end{split}
 \end{equation*}
 \item $p$-weighted energy estimate with $p = \frac{1}{2}\delta$
 \begin{equation*}
  \begin{split}
    & \mathcal{E}^{(L, \frac{1}{2}\delta)}[\phi](\tau_2) \\
    &
    + \int_{\mathcal{M}^{\tau_2}_{\tau_1}} \bigg(
    \chi_{(2r_0)}\left(\delta r^{-1+\frac{1}{2}\delta}|\slashed{\D}_L \phi|^2 + r^{-1 + \frac{1}{2}\delta|}\slashed{\nabla}\phi|^2 + r^{-3 + \frac{1}{2}\delta}|\phi|^2 \right) \bigg)\dVol_g \\
    &\lesssim
    \delta^{-1} \mathcal{E}^{(L, \frac{1}{2}\delta)}[\phi](\tau_1)
    + \delta^{-4}\mathcal{E}^{(wT)}[\phi](\tau_1)
    + \delta^{-4}\tilde{\mathcal{E}}_1
 \end{split}
 \end{equation*}

\end{enumerate}

\end{lemma}

\begin{proof}
 We begin by choosing the weight function $w = (1+r)^{-C_{[\phi]}\epsilon}$. We multiply the basic weighted Morawetz energy estimate (equation \eqref{equation basic weighted Morawetz estimate}) by $\delta$. We also multiply the basic $p$-weighted estimate \eqref{equation basic p weighted estimate low p}, by $\delta^3$, with the choice $p = \frac{1}{2}\delta$. We add together these two estimates, and also add the basic weighted $T$ energy estimate (equation \eqref{equation basic weighted T energy estimate}). We also assume some numerical bound on $\delta$, so that $\delta \ll 1$. Performing these energy estimates in the region $^t\mathcal{M}^{\tau_2}_{\tau_1}$, using the facts that $\epsilon \ll \delta$ we obtain the lengthy expression

\begin{equation}
\label{equation long energy estimate 1}
 \begin{split}
  &\mathcal{E}^{(wT)}[\phi](\tau_2, t, \tau_1) 
  + \int_{^t\mathcal{M}_{\tau_1}^{\tau_2}} C_{[\phi]}\epsilon (1+r)^{-1 - C_{[\phi]}\epsilon} |\slashed{\D}\phi|^2 \dVol_g \\
  & + \int_{^t\mathcal{M}_{\tau_1}^{\tau_2}} \left( \delta^2(1+r)^{-1-\frac{1}{2}\delta}|\slashed{\D}\phi|^2 + \delta\chi_{(r_0)}(r)(1+r)^{-1-C_{[\phi]}\epsilon}|\slashed{\nabla}\phi|^2 + \delta^2 r^{-1}(1+r)^{-2-\frac{1}{2}\delta}|\phi|^2 \right) \dVol_g \\
  & + \delta^3 \mathcal{E}^{(L,\frac{1}{2}\delta)}[\phi](\tau_2, R) \\
  & + \int_{^t\mathcal{M}_{\tau_1}^{\tau_2}} \chi_{(2r_0, R)}\left(\delta^4 r^{-1+\frac{1}{2}\delta}|\slashed{\D}_L \phi|^2 + \delta^3 r^{-1+\frac{1}{2}\delta}|\slashed{\nabla}\phi|^2 +  \delta^3 r^{-3 + \frac{1}{2}\delta}|\phi|^2 \right) \dVol_g \\
  &\lesssim
  \mathcal{E}^{(wT)}[\phi](\tau_2, t, \tau_1) 
  + \delta^3 \mathcal{E}^{(L,\frac{1}{2}\delta)}[\phi](\tau_1, R) \\
  &\phantom{\lesssim}
  + \int_{^t\mathcal{M}_{\tau_1}^{\tau_2}} w\bigg(
  \epsilon(1+\tau)^{-1-\delta}|\overline{\slashed{\D}}\phi|^2
  + (1+C_{[\phi]})\epsilon(1+r)^{-1}|\overline{\slashed{\D}}\phi|^2
  + \epsilon (1+C_{[\phi]}) r^{-1}(1+r)^{-2} |\phi|^2 \\
  &\phantom{\lesssim + \int_{^t\mathcal{M}_{\tau_1}^{\tau_2}} \bigg(}
  + \epsilon r^{-1}(1+r)^{-1-\delta}(1+\tau)^{-1-\delta}|\phi|^2
  + \epsilon^{-1} (1+r) |F|^2
   \bigg) \dVol_g \\
  &\phantom{\lesssim}
  + \int_{^t\mathcal{M}_{\tau_1}^{\tau_2}} \bigg(
  \delta^3\epsilon \chi_{(2r_0, R)} r^{\frac{1}{2}\delta-2} (1+\tau)^{-1-\delta}|\slashed{\D}_L \psi|^2
  + \delta^3\epsilon \chi_{(2r_0,R)} r^{\frac{1}{2}\delta-2 - \frac{1}{2}\beta}(1+\tau)^{-1-\delta} |\phi|^2 
  \\
  &\phantom{\lesssim + \int_{^t\mathcal{M}_{\tau_1}^{\tau_2}} }
  + \delta^3\epsilon^{-1}\chi_{(2r_0, R)} r^{1-\delta}(1+\tau)^{2\beta}|F|^2
  \bigg) \dVol_g 
  \\
  &\phantom{=}
  + \int_{\bar{S}_{t,\tau_1}}\delta (1+r)^{-C_{[\phi]}\epsilon} r|\phi|^2 \dVol_{\mathbb{S}^2}
  + \int_{\bar{S}_{t,\tau_2}}\delta (1+r)^{-C_{[\phi]}\epsilon} r|\phi|^2 \dVol_{\mathbb{S}^2} \\
  &\phantom{\lesssim}
  + \int_{^t\Sigma_{\tau_1}} \epsilon (1+r)^{-\frac{1}{2}\delta} |\slashed{\D}_L \phi|^2 r^2 \upd r \wedge \dVol_{\mathbb{S}^2}
  + \int_{\bar{S}_{t, \tau_1}} \epsilon(1+r)^{1 - \frac{1}{2}\delta } |\phi|^2 \dVol_{\mathbb{S}^2} \\
  &\phantom{\lesssim}
  + \int_{^t\Sigma_{\tau_2}} \epsilon (1+r)^{-\frac{1}{2}\delta} |\slashed{\D}_L \phi|^2 r^2 \upd r \wedge \dVol_{\mathbb{S}^2}
  + \int_{\bar{S}_{t, \tau_2}} \epsilon(1+r)^{1 - \frac{1}{2}\delta } |\phi|^2 \dVol_{\mathbb{S}^2} \\
  &\phantom{\lesssim} 
  + \int_{\tau_1}^{\tau_2} \left( \int_{^t\Sigma_{\tau'}\cap\left\{ \frac{1}{2}R \leq r \leq R \right\}} \delta^5 r^{\delta-1}|\overline{\slashed{\D}}\phi|^2 r^2 \upd r \wedge \dVol_{\mathbb{S}^2} 
  + \int_{S_{\tau,R}} \delta^4 |\phi|^2 \dVol_{\mathbb{S}^2} \right) \upd \tau
 \end{split}
\end{equation}
where, if $\phi$ is a scalar field, then the terms
\begin{equation*}
\int_{^t\mathcal{M}_{\tau_0}^{\tau_1}} \left(  w(1+\delta C_{[\phi]}) \epsilon r^{-1}(1+r)^{-2} |\phi|^2 
  + \delta^2 \epsilon \chi_{(r_0, R)} r^{\frac{1}{2}\delta - 3}|\phi|^2
  \right) \dVol_g
\end{equation*}
are absent.

Note that we have already absorbed many of the error terms on the right hand side by the positive terms on the left hand side. For example, most of the ``bulk'' error terms arising in the $T$-energy estimate have been absorbed by the bulk terms appearing in the Morawetz estimate. At the same time, the ``boundary'' terms in the Morawetz estimate have been absorbed by the boundary terms in the $T$-energy estimate. It is possible to accomplish both these tasks at the same time because $\epsilon \ll \delta \ll 1$.

Many of the remaining error terms on the right hand side of the expression above can absorbed by the left hand side in the region $r \leq \frac{1}{2}R$, using the facts that $\epsilon$ is very small compared to all the other parameters and that $\delta$ is small compared to the other parameters except for $\epsilon$. Moreover, the constant $\mathring{C}$ is large, but we still have $\mathring{C}\epsilon \ll \frac{1}{2}\delta$. Hence we have, for example, in the region $r \leq \frac{1}{2}R$
\begin{equation*}
 (1+C_{[\phi]})\epsilon (1+r)^{-1}|\overline{\slashed{\D}}\phi|^2
 \ll \mathring{C} \delta^3 \chi_{(2r_0, R)} r^{-1 + \mathring{C}\delta} |\overline{\slashed{\D}}\phi|^2 + \delta^2 (1+r)^{-1-\delta}|\slashed{\D}\phi|^2
\end{equation*}
where the first term on the right hand side bounds the left hand side in $2r_0 \leq r \leq \frac{1}{2}R$, and the second term bounds it in the region $0 \leq r \leq 2r_0$. Note that we cannot bound these terms in this way in the region $r \geq \frac{1}{2}R$. However, their contribution in this region can be absorbed into the final line of the expression \eqref{equation long energy estimate 1}.

By also making use of the Hardy inequalities, we can deal with the terms
\begin{equation*}
 \int_{^t\mathcal{M}_{\tau_1}^{\tau_2}} \epsilon r^{-1}(1+r)^{-2}|\phi|^2 \dVol_g
\end{equation*}
in a similar way (similarly for terms which are ``better'' than this, in that they have additional decay in $r$). The additional boundary term obtained when applying the Hardy inequality can be absorbed into the final term on the right hand side of \eqref{equation long energy estimate 1}.

Additionally, since $\delta \ll \beta$, then we can use the first part of the Hardy inequality \ref{proposition Hardy} to bound the term
\begin{equation*}
 \begin{split}
  &\int_{^t\mathcal{M}_{\tau_1}^{\tau_2}} \delta^3 \epsilon \chi_{(2r_0, R)} r^{\frac{1}{2}\delta - 2 - \beta}(1+\tau)^{-1-\delta} |\phi|^2 \dVol_g \\
  &\lesssim
  \int_{^t\mathcal{M}_{\tau_1}^{\tau_2}} \beta^{-2} \delta^3 \epsilon (1+r)^{\frac{1}{2}\delta - \beta}(1+\tau)^{-1-\delta} |\slashed{\D}_L\phi|^2 \dVol_g \\
  &\phantom{\lesssim}
  + \int_{\tau_1}^{\tau_2} \left( \int_{\bar{S}_{t,\tau} } \beta^{-1} \delta^3 \epsilon r^{\frac{1}{2}\delta + 1 - \beta} (1+\tau)^{-1-\delta} |\phi|^2 \dVol_{\mathbb{S}^2}\right) \upd \tau
 \end{split}
\end{equation*}

Furthermore, we have the bound
\begin{equation*}
 \int_{^t\mathcal{M}_{\tau_1}^{\tau_2} \cap \{r < r_0\}} \delta^3 \left( |\slashed{\D}\phi|^2 + |\phi|^2 \right)  \dVol_g
 \ll
 \int_{^t\mathcal{M}_{\tau_1}^{\tau_2} \cap \{r < r_0\}} \delta^2\left( (1+r)^{-1-\delta} |\slashed{\D}\phi|^2 + r^{-1}(1+r)^{-2-\delta} |\phi|^2 \right)  \dVol_g
\end{equation*}
using the fact that $\delta$ is sufficiently small compared to $r_0$ (which, as before, is considered a fixed constant now, so that all other variables can depend implicitly on $r_0$).

Next, we observe that, since $C_{[\phi]}\epsilon \ll \delta$ we have
\begin{equation*}
 \int_{^t\Sigma_{\tau_2}} \epsilon (1+r)^{-\frac{1}{2}\delta} |\slashed{\D}_L \phi|^2 r^2 \upd r \wedge \dVol_{\mathbb{S}^2}
 \lesssim \epsilon \mathcal{E}^{(wT)}(\tau_2, t, \tau_2)
\end{equation*}
and similarly for the boundary term on $^t\Sigma_{\tau_1}$.

Now, we shall deal with the error terms on the spheres $\bar{S}_{\tau, R}$, for various values of $\tau$. To bound these, we need to use lemma \ref{lemma exponential growth}. The most difficult to control terms are those with the largest $r$-weights, which in fact arise from the Morawetz current, and are of the form
\begin{equation*}
 \int_{\bar{S}_{t,\tau_1}} (1+r)^{-C_{[\phi]}\epsilon} r|\phi|^2 \dVol_{\mathbb{S}^2}
\end{equation*}
We handle these by making use of proposition \ref{proposition spherical mean in terms of energy}. Indeed, we have
\begin{equation*}
 \begin{split}
  &\int_{\bar{S}_{t,\tau_1}} (1+r)^{-C_{[\phi]}\epsilon} r|\phi|^2 \dVol_{\mathbb{S}^2}
  \leq \sup_{r \in \bar{S}_{t, \tau_1}} \left\{ r^{1-C_{[\phi]}\epsilon} \right\} \int_{\bar{S}_{t,\tau_1}} |\phi|^2 \dVol_{\mathbb{S}^2} \\
  &\lesssim \sup_{r \in \bar{S}_{t, \tau_1}} \left\{ r^{1-C_{[\phi]}\epsilon} \right\} 
  \left( \sup_{r \in \bar{S}_{t, \tau_1}} \left\{ r^{-1+\frac{1}{2}C_{[\phi]}\epsilon} \right\} 
  \int_{ ^{\tau_1}_{\tau_0}\bar{\Sigma}_t } (1+r)^{-\frac{1}{2}C_{[\phi]}\epsilon}|\slashed{\D}\phi|^2 r^2 \upd r \wedge \dVol_{\mathbb{S}^2}
  + \int_{\bar{S}_{t, \tau_0}} |\phi|^2 \dVol_{\mathbb{S}^2}
  \right)
 \end{split}
\end{equation*}
To obtain bounds on the behaviour of $r$ on the spheres $\bar{S}_{t, \tau}$, we observe that for $r \geq r_0$
\begin{equation*}
 \left. \frac{\partial t}{\partial r} \right|_{\tau, \vartheta^1, \vartheta^2} = L(t) = L^0 = 1 + \mathcal{O}(\epsilon (1+r)^{-\delta})
\end{equation*}
while at $r = r_0$ we have $t = \tau + r_0$. Hence we obtain the bound
\begin{equation*}
 |t - r - \tau| \lesssim \epsilon \left (1+r)^{1-\delta} - (1+r_0)^{1-\delta} \right)
\end{equation*}
so the reverse triangle inequality implies, in particular,
\begin{equation*}
 r \sim (t - \tau) \quad \text{ for} \quad r \geq r_0
\end{equation*}
Eventually, we have obtained the bound
\begin{equation*}
 \begin{split}
  \int_{\bar{S}_{t,\tau_1}} (1+r)^{-C_{[\phi]}\epsilon} r|\phi|^2 \dVol_{\mathbb{S}^2}
  &\lesssim (t - \tau_1)^{-\frac{1}{2}C_{[\phi]}\epsilon}
  \int_{ ^{\tau_1}_{\tau_0}\Sigma_t } (1+r)^{-\frac{1}{2}C_{[\phi]}\epsilon}|\slashed{\D}\phi|^2 r^2 \upd r \wedge \dVol_{\mathbb{S}^2} \\
  &\phantom{\lesssim}
  + (t-\tau_1)^{1-C_{[\phi]}\epsilon} \int_{\bar{S}_{t, \tau_0}} |\phi|^2 \dVol_{\mathbb{S}^2}
  \\
  &\lesssim (t - \tau_1)^{-\frac{1}{2}C_{[\phi]}\epsilon}\mathcal{E}^{(\tilde{w}T)}(\tau_0, t, \tau_0) + (t-\tau_1)^{1-C_{[\phi]}\epsilon} \int_{\bar{S}_{t, \tau_0}} |\phi|^2 \dVol_{\mathbb{S}^2}
 \end{split}
\end{equation*}
In exactly the same way, we can bound the error terms on the sphere $\bar{S}_{t, \tau_2}$.

There is also an error term arising from the \emph{integral over} $\tau$ of boundary terms on the spheres $S_{\tau,r}$. We bound this as
\begin{equation*}
 \begin{split}
  \int_{\tau_1}^{\tau_2} \left( \int_{S_{\tau,R}} |\phi|^2 \dVol_{\mathbb{S}^2} \right) \upd \tau
  &\lesssim \int_{\tau_1}^{\tau_2} \left( R^{-1+\frac{1}{2}C_{[\phi]}\epsilon}\mathcal{E}^{\tilde{w}T}(\tau, t, \tau_0) + \int_{\bar{S}_{t, \tau_0}} |\phi|^2 \dVol_{\mathbb{S}^2} \right)\upd \tau \\
  &\lesssim R^{-1 + \frac{1}{2}C_{[\phi]}\epsilon} \left( \int_{\tau_1}^{\tau_2} \mathcal{E}^{\tilde{w}T}(\tau, t, \tau_0)\right) + (\tau_2 - \tau_1) \int_{\bar{S}_{t, \tau_0}} |\phi|^2 \dVol_{\mathbb{S}^2} \\
  &\lesssim R^{-1 + \frac{1}{2}C_{[\phi]}\epsilon} \left( \int_{\tau_1}^{\tau_2} \mathcal{E}^{\tilde{w}T}(\tau, t, \tau_0)\right) + (\tau_2 - \tau_1)\int_{\bar{S}_{t, \tau_0}} |\phi|^2 \dVol_{\mathbb{S}^2}
 \end{split}
\end{equation*}

Next, we will deal with the error terms arising from the spacetime integral over the region $\frac{1}{2}R \leq r \leq R$. Since in this region $r \sim R$, we have
\begin{equation*}
 \begin{split}
  &\int_{\tau_1}^{\tau_2} \left( \int_{^t\Sigma_{\tau'}\cap\left\{ \frac{1}{2}R \leq r \leq R \right\}}  r^{\delta-1}|\overline{\slashed{\D}}\phi|^2 r^2 \upd r \wedge \dVol_{\mathbb{S}^2} \right) \upd \tau \\
  &\lesssim R^{-1 + \delta + \frac{1}{2}C_{[\phi]}} \int_{\tau_1}^{\tau_2} \left( \int_{^t\Sigma_{\tau'}\cap\left\{ \frac{1}{2}R \leq r \leq R \right\}} r^{-\frac{1}{2}C_{[\phi]}}|\overline{\slashed{\D}}\phi|^2 r^2 \upd r \wedge \dVol_{\mathbb{S}^2} \right) \upd \tau \\
  &\lesssim R^{-1 + 2\delta} \int_{\tau_1}^{\tau_2} \mathcal{E}^{\tilde{w}T}(\tau, t, \tau_0) \\
 \end{split}
\end{equation*}

Putting together all of the considerations above, we are able to simplify equation \eqref{equation long energy estimate 1} to obtain
\begin{equation}
\label{equation long energy estimate 2}
 \begin{split}
  &\mathcal{E}^{(wT)}[\phi](\tau_2, t, \tau_1) 
  + \int_{^t\mathcal{M}_{\tau_1}^{\tau_2}} C_{[\phi]}\epsilon (1+r)^{-1 - C_{[\phi]}\epsilon} |\slashed{\D}\phi|^2 \dVol_g \\
  & + \int_{^t\mathcal{M}_{\tau_1}^{\tau_2}} \left( \delta^2(1+r)^{-1-\frac{1}{2}\delta}|\slashed{\D}\phi|^2 + \delta\chi_{(r_0)}(r)(1+r)^{-1-C_{[\phi]}\epsilon}|\slashed{\nabla}\phi|^2 + \delta^2 r^{-1}(1+r)^{-2-\frac{1}{2}\delta}|\phi|^2 \right) \dVol_g \\
  & + \delta^3 \mathcal{E}^{(L,\frac{1}{2}\delta)}[\phi](\tau_2, R) \\
  & + \int_{^t\mathcal{M}_{\tau_1}^{\tau_2}} \chi_{(2r_0, R)}\left(  \delta^4 r^{-1+\frac{1}{2}\delta}|\slashed{\D}_L \phi|^2 + \delta^3 r^{-1+\frac{1}{2}\delta}|\slashed{\nabla}\phi|^2 + \delta^3 r^{-3 + \frac{1}{2}\delta}|\phi|^2 \right) \dVol_g \\
  &\lesssim
  \mathcal{E}^{(wT)}[\phi](\tau_2, t, \tau_1) 
  + \delta^3 \mathcal{E}^{(L,\frac{1}{2})}[\phi](\tau_1, R) 
  + (t - \tau_2)^{-\frac{1}{2}C_{[\phi]}\epsilon}\mathcal{E}^{(\tilde{w}T)}(\tau_0, t, \tau_0)  
  \\
  &\phantom{\lesssim}
  + \int_{^t\mathcal{M}_{\tau_1}^{\tau_2}} w\bigg(
   \epsilon^{-1}r^{-1}(1+r)^{-2} |\phi|^2 
  + \epsilon^{-1} (1+r) |F|^2
   \bigg) \dVol_g \\
  &\phantom{\lesssim}
  + \int_{^t\mathcal{M}_{\tau_1}^{\tau_2}} \bigg(
  \delta^3\epsilon \chi_{(2r_0, R)} r^{\frac{1}{2}\delta-2} (1+\tau)^{-1-\delta}|\slashed{\D}_L \psi|^2
  + w \epsilon (1+\tau)^{-1-\delta} |\overline{\slashed{\D}} \phi|^2 \\
  &\phantom{\lesssim + \int_{^t\mathcal{M}_{\tau_1}^{\tau_2}} }
  + \delta^3\epsilon^{-1}\chi_{(2r_0, R)} r^{1-\delta}(1+\tau)^{2\beta}|F|^2
   \bigg) \dVol_g 
  \\
  &\phantom{=}
  + (t-\tau_0)^{1-C_{[\phi]}\epsilon} \int_{\bar{S}_{t, \tau_0}} |\phi|^2 \dVol_{\mathbb{S}^2}
  + R^{-1 + 2\delta} \left( \int_{\tau_1}^{\tau_2} \mathcal{E}^{\tilde{w}T}(\tau, t, \tau_0)\right) + (\tau_2 - \tau_1) \int_{\bar{S}_{t, \tau_0}} |\phi|^2 \dVol_{\mathbb{S}^2}
 \end{split}
\end{equation}

Using the conditions on the inhomogeneous term $F$ and the conditions on the initial data, we can further simplify this to find
\begin{equation}
\label{equation long energy estimate 3}
 \begin{split}
  &\mathcal{E}^{(wT)}[\phi](\tau_2, t, \tau_1) 
  + \int_{^t\mathcal{M}_{\tau_1}^{\tau_2}} C_{[\phi]}\epsilon (1+r)^{-1 - C_{[\phi]}\epsilon} |\slashed{\D}\phi|^2 \dVol_g \\
  & + \int_{^t\mathcal{M}_{\tau_1}^{\tau_2}} \left( \delta^2(1+r)^{-1-\frac{1}{2}\delta}|\slashed{\D}\phi|^2 + \delta\chi_{(r_0)}(r)(1+r)^{-1-C_{[\phi]}\epsilon}|\slashed{\nabla}\phi|^2 + \delta^2 r^{-1}(1+r)^{-2-\frac{1}{2}\delta}|\phi|^2 \right) \dVol_g \\
  & + \delta^3 \mathcal{E}^{(L,\frac{1}{2}\delta)}[\phi](\tau_2, R) \\
  & + \int_{^t\mathcal{M}_{\tau_1}^{\tau_2}} \chi_{(2r_0, R)}\left( \delta^4 r^{-1+\frac{1}{2}\delta}|\slashed{\D}_L \phi|^2 + \delta^3 r^{-1+\frac{1}{2}\delta}|\slashed{\nabla}\phi|^2 + \delta^3 r^{-3 + \frac{1}{2}\delta}|\phi|^2 \right) \dVol_g \\
  &\lesssim
  \tilde{\mathcal{E}}_1  
  + \mathcal{E}_1
  + (t - \tau_2)^{-\frac{1}{2}C_{[\phi]}\epsilon}\mathcal{E}_0
  + (\tau_2 - \tau_1) (t - \tau_0)^{-1 + \frac{1}{2}C_{[\phi]}\epsilon}\mathcal{E}_0
  \\
  &\phantom{\lesssim}
  + \int_{^t\mathcal{M}_{\tau_1}^{\tau_2}} \bigg(
  \delta^3\epsilon \chi_{(2r_0, R)} r^{\frac{1}{2}\delta-2} (1+\tau)^{-1-\delta}|\slashed{\D}_L \psi|^2
  + \delta^3\epsilon \chi_{(2r_0,R)} r^{\frac{1}{2}\delta-2 - \frac{1}{2}\beta}(1+\tau)^{-1-\delta} |\phi|^2 \\
  &\phantom{\lesssim + \int_{^t\mathcal{M}_{\tau_1}^{\tau_2}} \bigg(}
  + w\epsilon(1+\tau)^{-1-\delta}|\overline{\slashed{\D}}\phi|^2
   \bigg) \dVol_g 
  \\
  &\phantom{\lesssim} 
  + R^{-1 + 2\delta} \left( \int_{\tau_1}^{\tau_2} \mathcal{E}^{\tilde{w}T}(\tau, t, \tau_0)\right) 
 \end{split}
\end{equation}

Now, to handle this final term, we appeal to lemma \ref{lemma exponential growth}. If $C_{[\phi]}$ is sufficiently large, then we find that we can bound
\begin{equation*}
 \mathcal{E}^{\tilde{w}T}(\tau, t, \tau_0) \lesssim \left(\mathcal{E}_0 + \tilde{\mathcal{E}}_0\right) e^{\tilde{C}\epsilon (\tau - \tau_0)}
\end{equation*}
and so we have
\begin{equation*}
R^{-1 + 2\delta} \left( \int_{\tau_1}^{\tau_2} \mathcal{E}^{\tilde{w}T}(\tau, t, \tau_0)\right) 
\lesssim \epsilon^{-1} R^{-1 + 2\delta} \left(\mathcal{E}_0 + \tilde{\mathcal{E}}_0\right) e^{\tilde{C}\epsilon (\tau_2 - \tau_0)}
\end{equation*}

We can also use the Hardy inequality to estimate the term
\begin{equation*}
 \begin{split}
  &\int_{^t\mathcal{M}_{\tau_1}^{\tau_2}} \delta^3\epsilon \chi_{(2r_0,R)} r^{\frac{1}{2}\delta-2 - \frac{1}{2}\beta}(1+\tau)^{-1-\delta} |\phi|^2 \dVol_g \\
  &\lesssim \epsilon \int_{\tau_1}^{\tau_2} (1+\tau)^{-1-\delta} \left(
 \int_{^t\Sigma_\tau} (1+r)^{-\frac{1}{2}\beta + \frac{1}{2}\delta} |\slashed{\D}_L \phi|^2 r^2 \upd r \wedge \dVol_{\mathbb{S}^2} \right) \upd \tau \\
  &\phantom{\lesssim}
  + \epsilon \int_{\tau_1}^{\tau_2} (1+\tau)^{-1-\delta} \left(\int_{\bar{S}_{t, \tau}} r^{1 -\frac{1}{2}\beta + \frac{1}{2}\delta} |\phi|^2 \dVol_{\mathbb{S}^2}
  \right) \upd \tau
 \end{split}
\end{equation*}
and the second term on the right hand side of the above inequality can be estimated as
\begin{equation*}
 \begin{split}
  &\epsilon \int_{\tau_1}^{\tau_2} (1+\tau)^{-1-\delta} \left(\int_{\bar{S}_{t, \tau}} r^{1 -\frac{1}{2}\beta + \frac{1}{2}\delta} |\phi|^2 \dVol_{\mathbb{S}^2}
  \right) \upd \tau \\
  &\lesssim \epsilon \int_{\tau_1}^{\tau_2} (1+\tau)^{-1-\delta} (t - \tau)^{-\frac{1}{2}\beta + \frac{1}{2}\delta + \frac{1}{2}C_{[\phi]}\epsilon} \left(\int_{^{\tau}_{\tau_0}\bar{\Sigma}_t} (1+r)^{-\frac{1}{2}C_{[\phi]}\epsilon} |\phi|^2 \dVol_{\mathbb{S}^2} \right) \upd \tau \\
  &\lesssim (1+\tau_1)^{-1-\delta} (t-\tau_2)^{-\frac{1}{2}\beta + \frac{1}{2}\delta + \frac{1}{2}C_{[\phi]}\epsilon} \left(\mathcal{E}_0 + \tilde{\mathcal{E}}_0\right) e^{\tilde{C}\epsilon(\tau_2 - \tau_0)}
 \end{split}
\end{equation*}
Although this quantity grows rapidly in $\tau_2$, for fixed $\tau_2$ and $\tau_1$ it tends to zero in the limit $t \rightarrow \infty$. This is the limit we will later take.

We have now arrived at the inequality
\begin{equation*}
 \begin{split}
  &\mathcal{E}^{(wT)}[\phi](\tau_2, t, \tau_1) 
  + \int_{^t\mathcal{M}_{\tau_1}^{\tau_2}} C_{[\phi]}\epsilon (1+r)^{-1 - C_{[\phi]}\epsilon} |\slashed{\D}\phi|^2 \dVol_g \\
  & + \int_{^t\mathcal{M}_{\tau_1}^{\tau_2}} \left( \delta^2(1+r)^{-1-\frac{1}{2}\delta}|\slashed{\D}\phi|^2 + \delta\chi_{(r_0)}(r)(1+r)^{-1-C_{[\phi]}\epsilon}|\slashed{\nabla}\phi|^2 + \delta^2 r^{-1}(1+r)^{-2-\frac{1}{2}\delta}|\phi|^2 \right) \dVol_g \\
  & + \delta^3 \mathcal{E}^{(L,\frac{1}{2}\delta)}[\phi](\tau_2, R) \\
  & + \int_{^t\mathcal{M}_{\tau_1}^{\tau_2}} \chi_{(2r_0, R)}\left(  \delta^4 r^{-1+\frac{1}{2}\delta}|\slashed{\D}_L \phi|^2 + \delta^3 r^{-1+\frac{1}{2}\delta}|\slashed{\nabla}\phi|^2 +  \delta^3 r^{-3 + \frac{1}{2}\delta}|\phi|^2 \right) \dVol_g \\
  &\lesssim
  \mathcal{E}_1
  + \tilde{\mathcal{E}}_1
  + (t - \tau_2)^{-\frac{1}{2}C_{[\phi]}\epsilon}\left(\mathcal{E}_0 + \tilde{\mathcal{E}}_0\right)
  + (\tau_2 - \tau_1) (t - \tau_0)^{-1 + 2\delta}\left(\mathcal{E}_0 + \tilde{\mathcal{E}}_0\right)
  \\
  &\phantom{\lesssim}
  + \int_{^t\mathcal{M}_{\tau_1}^{\tau_2}} \bigg(
  \delta^3\epsilon \chi_{(2r_0, R)} r^{\frac{1}{2}\delta-2} (1+\tau)^{-1-\delta}|\slashed{\D}_L \psi|^2
  + \epsilon (1+\tau)^{-1-\delta}|\overline{\slashed{\D}}\phi|^2 
   \bigg) \dVol_g 
  \\
  &\phantom{\lesssim} 
  + \left(\mathcal{E}_0 + \tilde{\mathcal{E}}_0\right) e^{\tilde{C}\epsilon(\tau_2 - \tau_0)} \left( \epsilon^{-1}R^{-1+\frac{3}{2}\delta} + (t - \tau_2)^{-\frac{1}{2}\beta + \frac{1}{2}\delta + \frac{1}{2}C_{[\phi]}\epsilon} \right)
 \end{split}
\end{equation*}

Now, we can appeal to the Gronwall inequality \ref{proposition Gronwall}, this time with the choices
\begin{equation*}
 \begin{split}
  h(\tau) &= \int_{^t\Sigma_\tau} w|\overline{\slashed{\D}}\phi|^2 r^2 \upd r \wedge \dVol_{\mathbb{S}^2}
  + \int_{^{\tau}_{\tau_0}\bar{\Sigma}_t} w|\slashed{\D}\phi|^2 r^2 \upd r \wedge \dVol_{\mathbb{S}^2}
  \\
  &\phantom{=} 
  + \int_{^t\mathcal{M}_{\tau_1}^{\tau}}\bigg(
  C_{[\phi]}\epsilon (1+r)^{-1 - C_{[\phi]}\epsilon} |\slashed{\D}\phi|^2 \delta^2(1+r)^{-1-\frac{1}{2}\delta}|\slashed{\D}\phi|^2 
  + \delta\chi_{(r_0)}(r)(1+r)^{-1-C_{[\phi]}\epsilon}|\slashed{\nabla}\phi|^2 \\
  &\phantom{=\int_{^t\mathcal{M}_{\tau_1}^{\tau}}\bigg(}
  + \delta^2 r^{-1}(1+r)^{-2-\frac{1}{2}\delta}|\phi|^2 \bigg) \dVol_g \\ \\
  f(\tau) &= \int_{^t\Sigma_\tau} \left( \delta^3 \chi_{(2r_0,R)}r^{\frac{1}{2}\delta}|\slashed{\D}_L \psi|^2 + w\epsilon|\overline{\slashed{\D}}\phi|^2 \right) \upd r \wedge \dVol_{\mathbb{S}^2} \\ \\
  G(\tau, \tau_1) &= \mathcal{E}_1 + \tilde{\mathcal{E}}_1 + \left(\mathcal{E}_0 + \tilde{\mathcal{E}}_0 \right)\bigg( (t - \tau)^{-\frac{1}{2}C_{[\phi]}\epsilon} + (\tau - \tau_1)(t - \tau_0)^{-1 + \frac{3}{2}\delta} + \epsilon^{-1}R^{-1 + \frac{3}{2}\delta}e^{\tilde{C}\epsilon(\tau - \tau_0)} \\
  &\phantom{= \tilde{\mathcal{E}}\bigg(}
  + (t - \tau)^{-\frac{1}{2}\beta + \frac{1}{2}\delta + \frac{1}{2}C_{[\phi]}\epsilon}e^{\tilde{C}\epsilon(\tau - \tau_0)} \bigg) \\
  g(\tau) &= \epsilon(1+\tau)^{-1-\delta}
 \end{split}
\end{equation*}
Note that here $\tau_0$, $R$ and $t$ are considered fixed constants. In particular, the ``initial time'' appearing as the second argument of $G$ is $\tau_1$, not $\tau_0$.

Since $g(\tau)$ is integrable, and since $\epsilon \ll \delta$, we obtain, for all $\tau_2 \geq \tau_1 \geq \tau_0$
\begin{equation*}
 f(\tau_2) + h(\tau_2) \lesssim f(\tau_1) + h(\tau_1) + G(\tau_2, \tau_1)
\end{equation*}

Note in particular that, for any fixed $\tau_1$ and $\tau_2$, if we take the limits $t \rightarrow \infty$ and $R \rightarrow \infty$ then $G(\tau_2, \tau_1) \rightarrow \tilde{\mathcal{E}}_1$, implying the final part of the lemma.

\end{proof}

\section{Decay of the degenerate energy and integrated local energy, and the \texorpdfstring{$p$}{p}-weighted energy estimate for large \texorpdfstring{$p$}{p}}

In the previous section we established boundedness of the degenerate energy, as well as some additional estimates (integrated local energy decay and the $p$-weighted energy estimate with very small $p$). Although these estimates already give us some control over the solution $\phi$, they are insufficient to close the bootstrap assumptions of chapter \ref{chapter bootstrap}. For this, we need to show that the degenerate energy \emph{decays}, and that is the subject of this section. 

In order to prove decay, we will again use the $p$-weighted energy estimate, but this time with $p \sim 1$ instead of $p \sim \delta$. Then, this can be combined with the estimates of the previous sections (in particular the degenerate energy boundedness estimate) to prove degenerate energy decay.

\begin{lemma}[The $p$ weighted energy estimates for large $p$]
\label{lemma p weighted}

Let $\phi$ be an $S_{\tau,r}$-tangent tensor field satisfying
\begin{equation*}
 \tilde{\slashed{\Box}}_g \phi = F = F_1 + F_2 + F_3 = F_4 + F_5 + F_6
\end{equation*}
for some $S_{\tau, r}$-tangent tensor fields $F$, $F_1 , \ldots , F_6$. Let these tensor fields satisfy
\begin{equation*}
 \begin{split}
  &\int_{\mathcal{M}^\tau_{\tau_0}} \bigg( \epsilon^{-1} \chi_{(r_0)} \Big(
  r^{1-C_{[\phi]}\epsilon}(1+\tau)^{1+\delta}|F_4|^2
  + r^{2-C_{[\phi]}\epsilon-2\delta}(1+\tau)^{6\delta}|F_5|^2
  + r^{2-C_{[\phi]}\epsilon}|F_6|^2
  \Big)\bigg) \dVol_g 
  \\
  &\lesssim \tilde{\mathcal{E}}_0 (1+\tau)^{C_{[\phi]}\delta}
  \\
  \\
  &\int_{\mathcal{M}^{\tau_1}_{\tau}} \epsilon^{-1}\bigg( \Big(
  (1+r)^{1-C_{[\phi]}\epsilon}|F|^2
  + (1+r)^{\frac{1}{2}\delta}(1+\tau)^{1+\delta}|F_1|^2
  + (1+r)^{1-\frac{7}{2}\delta}(1+\tau)^{6\delta}|F_2|^2
  \\
  &\phantom{\int_{\mathcal{M}^{\tau_1}_{\tau}} \epsilon^{-1}\bigg( \Big(}
  + (1+r)^{1+\frac{1}{2}\delta}|F_3|^2
  \Big) \bigg) \dVol_g
  \lesssim \tilde{\mathcal{E}}_0(1+\tau)^{-1 + C_{[\phi]}\delta}
 \end{split}
\end{equation*}

Define
\begin{equation*}
 \begin{split}
    \tilde{w} &= (1+r)^{-\frac{1}{2}C_{[\phi]}\epsilon}
 \end{split}
\end{equation*}

Let the initial data of $\phi$ satisfy
\begin{equation*}
 \begin{split}
  \mathcal{E}^{(\tilde{w}T)}[\phi](\tau_0) &\lesssim \mathcal{E}_0 \\
  \mathcal{E}^{(L, 1-C_{[\phi]}\epsilon)}[\phi](\tau_0) &\lesssim \mathcal{E}_0 \\
  \int_{\bar{S}_{t,r}}|\phi|^2 \dVol_{\mathbb{S}^2} &\lesssim \mathcal{E}_0(t - \tau_0)^{-1+\frac{1}{2}C_{[\phi]}\epsilon}
 \end{split}
\end{equation*}

Finally, suppose that all of the bootstrap bounds of chapter \ref{chapter bootstrap} hold.
 
Then, for sufficiently small $\epsilon$, for sufficiently small $\delta$, and for sufficiently large constants $C_{[\phi]}$, for all $\tau \geq \tau_0$ we have
\begin{equation*}
 \begin{split}
  &\mathcal{E}^{(L, 1-C_{[\phi]}\epsilon)}[\phi](\tau)
  + \int_{\mathcal{M}^{\tau}_{\tau_0}} \chi_{(2r_0)}r^{-C_{[\phi]}\epsilon} \left( |\slashed{\D}_L \phi|^2 + |\slashed{\nabla}\phi|^2 + C_{[\phi]}\epsilon r^{-2}|\phi|^2 \right) \dVol_g
  \\
  &\lesssim \delta^{-8\delta^{-1}} (\mathcal{E}_0 + \tilde{\mathcal{E}}_0) (1+\tau)^{C_{[\phi]}\delta}
 \end{split}
\end{equation*}

\end{lemma}

\begin{proof}
In order to establish decay of the integrated local energy, and to increase $p$, we will rely on a double continuity argument. Suppose that we have
\begin{equation}
\label{equation large p internal bootstrap 1}
 \int_{^t\mathcal{M}^{\tau}_{\tau_0}} \epsilon (1+\tau)^{1-\delta} (1+r)^{-1-C_{[\phi]}\epsilon} |\slashed{\D}\phi|^2 \dVol_g \lesssim (C_{[\phi]})^{-1} \delta^{-8\delta^{-1}\kappa}(\mathcal{E}_0 + \tilde{\mathcal{E}}_0) (1+\tau)^{1-\kappa}
\end{equation}
for all $\kappa \in [\delta, \kappa_{(\text{max})}]$, for some $\kappa_{\text{max}} \leq 1$. Note that, by lemma \ref{lemma boundedness} we can at least take $\kappa_{(\text{max})} = \delta$. Additionally, suppose that for all $p \in [\delta, p_{(\text{max})}]$, and for all $R$, we have
\begin{equation}
\label{equation large p internal bootstrap 2}
 \begin{split}
  &\mathcal{E}^{(L,p)}[\phi](\tau,R) 
  + \int_{^t\mathcal{M}_{\tau_0}^{\tau}} \chi_{(2r_0, R)}\left( pr^{p-1}|\slashed{\D}_L \phi|^2 + (2-p)r^{p-1}|\slashed{\nabla}\phi|^2 + p(1-p)r^{p-3}|\phi|^2 \right) \dVol_g \\
  &\lesssim
  (C_{[\phi]})^{-1}\delta^{-8\delta^{-1}\kappa}(\mathcal{E}_0 + \tilde{\mathcal{E}}_0) (1+\tau)^{1-\kappa}
 \end{split}
\end{equation}
Again, this certainly holds for $p_{(\text{max})} = \delta$ and for $\kappa = \delta$. We will refer to equations \eqref{equation large p internal bootstrap 1} and \eqref{equation large p internal bootstrap 2} as the ``internal bootstrap assumptions''. These will be closed in the course of this proof.

Now we return to equation \eqref{equation basic p weighted estimate high p}. For $p \leq 1 - C_{[\phi]}\epsilon$ we have 
\begin{equation*}
 \begin{split}
  &\mathcal{E}^{(L,p)}[\phi](\tau, R) 
  + \int_{^t\mathcal{M}_{\tau_0}^{\tau}} \chi_{(2r_0, R)}\left( pr^{p-1}|\slashed{\D}_L \phi|^2 + (2-p)r^{p-1}|\slashed{\nabla}\phi|^2 + p(1-p)r^{p-3}|\phi|^2 \right) \dVol_g \\
  &\lesssim
  \mathcal{E}^{(L,p)}[\phi](\tau_0, R)
  \\
  &\phantom{\lesssim}
  + \int_{^t\mathcal{M}_{\tau_0}^{\tau}}\bigg( 
  \epsilon \chi_{(2r_0, R)}(1+r)^{-1 - C_{[\phi]}\epsilon } (1+\tau)^{1-\delta} |\slashed{\D}\phi|^2
  + \epsilon \chi_{(2r_0, R)} r^{p-2} (1+\tau)^{-1-\delta}|\slashed{\D}_L \psi|^2
  \\
  &\phantom{\lesssim + \int_{^t\mathcal{M}_{\tau_0}^{\tau}} \bigg(}
  + \epsilon \chi_{(2r_0, R)} r^{p-2-\frac{1}{2}\delta}(1+\tau)^{-1-\delta} |\phi|^2
  + \epsilon \chi_{(2r_0, R)} r^{p-3} |\phi|^2
  + \epsilon^{-1}\chi_{(r_0, R)} r^p (1+\tau)^{1+\delta}|F_1|^2
  \\
  &\phantom{\lesssim + \int_{^t\mathcal{M}_{\tau_0}^{\tau}} \bigg(}
  + \epsilon^{-1}\chi_{(r_0, R)} r^{p+1-4\delta}(1+\tau)^{6\delta}|F_2|^2
  + \epsilon^{-1}\chi_{(r_0, R)} r^{p+1}|F_3|^2
  \bigg) \dVol_g \\
  &\phantom{\lesssim}
  + \int_{^t\Sigma_{\tau}} \epsilon \frac{1}{(1-p+\delta)^2} r^{p-\delta} |\slashed{\D}_L \psi|^2 \upd r \wedge \dVol_{\mathbb{S}^2}
  + \int_{S_{\tau, R}} \epsilon \frac{1}{(1-p+\delta)} r^{p + 1 - \delta } |\phi|^2 \dVol_{\mathbb{S}^2} \\
  &\phantom{\lesssim}
  + \int_{^t\Sigma_{\tau_0}} \epsilon \frac{1}{(1-p+\delta)^2} r^{p-\delta} |\slashed{\D}_L \psi|^2 \upd r \wedge \dVol_{\mathbb{S}^2}
  + \int_{S_{\tau_0, R}} \epsilon \frac{1}{(1-p+\delta)} r^{p + 1 - \delta } |\phi|^2 \dVol_{\mathbb{S}^2} \\
  &\phantom{\lesssim} 
  + \int_{\tau_0}^\tau \left( \int_{^t\Sigma_{\tau'}\cap\left\{ \frac{1}{2}R \leq r \leq R \right\}} p^{-2} r^{-\epsilon}|\overline{\slashed{\D}}\phi|^2 r^2 \upd r \wedge \dVol_{\mathbb{S}^2} 
  + \int_{S_{\tau',R}} p^{-1} r^{1-\epsilon} |\phi|^2 \dVol_{\mathbb{S}^2} \right) \upd \tau
 \end{split}
\end{equation*}

Making use of the internal bootstrap assumptions made at the beginning of the proof, together with the assumptions of the inhomogeneity $F$, we find that, for $\frac{1}{2}\delta \leq p \leq p_{(\text{max})} + \delta$ we have
\begin{equation*}
 \begin{split}
  &\mathcal{E}^{(L,p)}[\phi](\tau, R) 
  + \int_{^t\mathcal{M}_{\tau_0}^{\tau}} \chi_{(2r_0, R)}\left( pr^{p-1}|\slashed{\D}_L \phi|^2 + (2-p)r^{p-1}|\slashed{\nabla}\phi|^2 + p(1-p)r^{p-3}|\phi|^2 \right) \dVol_g \\
  &\lesssim
  (C_{[\phi]})^{-1}\delta^{-8\delta^{-1}\kappa}(\mathcal{E}_0 + \tilde{\mathcal{E}}_0) \bigg( 
  	\delta^{-1+8\delta^{-1}\kappa}(1+\tau)^{C_{[\phi]}\delta}
  	+ (1+\tau)^{1-\kappa}
  	+ \frac{\epsilon}{(1-p+\delta)^2}(1+\tau)^{1-\kappa}
  	\\
  	&\phantom{\lesssim (C_{[\phi]})^{-1}\delta^{-8\delta^{-1}\kappa}(\mathcal{E}_0 + \tilde{\mathcal{E}}_0) \bigg( }
  	+ \frac{\epsilon}{(1-p+\delta)^2}(1+\tau_0)^{1-\kappa}
  \bigg)
  \\
  &\phantom{\lesssim}
  + \int_{^t\mathcal{M}_{\tau_0}^{\tau}}\bigg( 
  \epsilon \chi_{(2r_0, R)} r^{p-2} (1+\tau)^{-1-\delta}|\slashed{\D}_L \psi|^2
  + \epsilon \chi_{(2r_0, R)} r^{p-2 - \frac{1}{2}\delta}(1+\tau)^{-1-\delta} |\phi|^2  
   \bigg) \dVol_g \\
  &\phantom{\lesssim} 
  + \int_{S_{\tau, R}} \epsilon \frac{1}{(1-p+\delta)} r^{p + 1 - \delta } |\phi|^2 \dVol_{\mathbb{S}^2} \\
  &\phantom{\lesssim} 
  + \int_{\tau_0}^\tau \left( \int_{^t\Sigma_{\tau'}\cap\left\{ \frac{1}{2}R \leq r \leq R \right\}} p^{-2} r^{-\epsilon}|\overline{\slashed{\D}}\phi|^2 r^2 \upd r \wedge \dVol_{\mathbb{S}^2} 
  + \int_{S_{\tau',R}} p^{-1} r^{1-\epsilon} |\phi|^2 \dVol_{\mathbb{S}^2} \right) \upd \tau
 \end{split}
\end{equation*}

Using lemma \ref{lemma boundedness} we have
\begin{equation*}
\int_{^t\Sigma_{\tau'}\cap\left\{ \frac{1}{2}R \leq r \leq R \right\}} p^{-2} r^{-\delta}|\overline{\slashed{\D}}\phi|^2 r^2 \upd r \wedge \dVol_{\mathbb{S}^2} \lesssim p^{-2} R^{-\delta + C_{[\phi]}\epsilon} \delta^{-1}\mathcal{E}_0
\end{equation*}
and so we easily obtain
\begin{equation*}
 \int_{\tau_0}^{\tau} \left( \int_{^t\Sigma_{\tau'}\cap\left\{ \frac{1}{2}R \leq r \leq R \right\}} p^{-2} r^{-\delta}|\overline{\slashed{\D}}\phi|^2 r^2 \upd r \wedge \dVol_{\mathbb{S}^2} \right) \upd \tau' \lesssim p^{-2} (\tau - \tau_0) R^{-\frac{1}{2}\delta} \delta^{-1}\mathcal{E}_0
\end{equation*}
Importantly, this quantity tends to zero as $R \rightarrow \infty$.

Next, we consider the boundary terms on the spheres. First, we note that, for any $p' < 1$ we have
\begin{equation*}
 \begin{split}
  \int_{S_{\tau, R}} |\psi|^2 \dVol_{\mathbb{S}^2}
  &= \int_{r_0}^{R} \left( \int_{S_{\tau, r}} \frac{\partial |\psi|}{\partial r} \dVol_{\mathbb{S}^2} \right)^2 \upd r
  - \int_{S_{\tau, r_0}} |\psi|^2 \dVol_{\mathbb{S}^2}
  + 2 \int_{\mathbb{S}^2} |\psi (\tau, R, \vartheta)| |\psi (\tau, r_0, \vartheta)| \dVol_{\mathbb{S}^2} \\
  &\lesssim \int_{r_0}^{R} \left( \int_{S_{\tau, r}} \frac{\partial |\psi|}{\partial r} \dVol_{\mathbb{S}^2} \right)^2 \upd r
  + \int_{S_{\tau, r_0}} |\psi|^2 \dVol_{\mathbb{S}^2} \\
  &\lesssim  R^{1-p'} \int_{r_0}^{R} \left( \int_{S_{\tau, r}} r^{p'} |\slashed{\D}_L \psi|^2 \dVol_{\mathbb{S}^2} \right)\upd r
  + (r_0)^2 \int_{S_{\tau, r_0}} |\phi|^2 \dVol_{\mathbb{S}^2} \\
  &\lesssim  R^{1-p'} \int_{\Sigma_\tau \cap \{r_0 \leq r \leq R\}} r^{p'} |\slashed{\D}_L \psi|^2 \upd r \wedge \dVol_{\mathbb{S}^2}
  + (r_0)^2 \mathcal{E}_0 \\
 \end{split}
\end{equation*}
where in the last line we have used the first part of proposition \ref{proposition spherical mean in terms of energy}.

Hence we have
\begin{equation*}
  \begin{split}
    \int_{S_{\tau, R}} r^{p + 1 - \delta } |\phi|^2 \dVol_{\mathbb{S}^2}
    &=
    \int_{S_{\tau, R}} r^{p - 1 - \delta } |\psi|^2 \dVol_{\mathbb{S}^2} \\
    & \lesssim  R^{p-p' - \delta} \int_{\Sigma_\tau \cap \{r_0 \leq r \leq R\}} r^{p'} |\slashed{\D}_L \psi|^2 \upd r \wedge \dVol_{\mathbb{S}^2}
  + R^{p-1-\delta}(r_0)^2 \mathcal{E} \\
  \end{split}
\end{equation*}

Now, we make the further restriction $p \leq p_{(\text{max})} + \frac{1}{2}\delta$. Then we can choose $p' = p_{(\text{max})}$. We find
\begin{equation*}
  \begin{split}
    \int_{S_{\tau, R}} r^{p + 1 - \delta } |\phi|^2 \dVol_{\mathbb{S}^2}
    &=
    \int_{S_{\tau, R}} r^{p - 1 - \delta } |\psi|^2 \dVol_{\mathbb{S}^2} \\
    & \lesssim  R^{p-p_{(\text{max})} - \delta} \int_{r_0}^{R} \int_{S_{\tau, r}} r^{p_{(\text{max})}} |\slashed{\D}_L \psi|^2 \dVol_{\mathbb{S}^2} \upd r
  + R^{p-1-\delta}(r_0)^2 \mathcal{E}_0 \\
    &\lesssim R^{-\frac{1}{2}\delta} \left( \mathcal{E}^{(L, p_{(\text{max})})}[\phi](\tau, 2R) + \mathcal{E}_0\right) \\
    &\lesssim R^{-\frac{1}{2}\delta} \delta^{-1}(\mathcal{E}_0 + \tilde{\mathcal{E}}_0)\left( (1+\tau)^{C_{[\phi]}\delta} + \delta^{1-8\delta^{-1}\kappa}(1+\tau)^{1-\kappa} \right) \\
  \end{split}
\end{equation*}
where we again note that we are allowing implicit constants to depend on $r_0$. Note that this quantity tends to zero as $R \rightarrow \infty$.

We also need to deal with the bulk zero-th order terms. Using the second part of the Hardy inequality \ref{proposition Hardy} we have
\begin{equation*}
 \begin{split}
   &\int_{^t\mathcal{M}_{\tau_0}^{\tau}}\bigg( 
   \chi_{(2r_0, R)} r^{p-2 - \frac{1}{2}\delta}(1+\tau)^{-1-\delta} |\phi|^2  
   \bigg) \dVol_g \\
   &\lesssim 
   \delta^{-2} \int_{^t\mathcal{M}_{\tau_0}^{\tau}}\bigg( 
   \chi_{(2r_0, R)} r^{p-2- \frac{1}{2}\delta}(1+\tau)^{-1-\delta} |\slashed{\D}_L \psi|^2  
   \bigg) \dVol_g 
   + \delta^{-1}\int_{^t\mathcal{M}_{\tau_0}^{\tau} \cap \{r_0 \leq r \leq 2r_0\} } (1+\tau)^{-1-\delta}|\phi|^2 \dVol_g \\
   &\phantom{\lesssim} 
   + \delta^{-1}\int_{^t\mathcal{M}_{\tau_0}^{\tau} \cap \{\frac{1}{2}R \leq r \leq R \} } r^{p-3+\delta}(1+\tau)^{-1-\delta} |\phi|^2 \dVol_g
 \end{split}
\end{equation*}
The first term on the right hand side is bounded by $\delta^{-2-8\delta^{-1}\kappa}(\mathcal{E}_0 + \tilde{\mathcal{E}}_0)(1+\tau)^{1-\kappa-\delta}$ by the (internal) bootstrap bound, since $p - 3 + \delta \leq p_{(\text{max})} - 2$. Similarly, the second term is bounded by $\delta^{-2}\mathcal{E}$ by the use of the integrated local energy decay estimate. Note that these terms also appear with a factor of $\epsilon$, and $\epsilon/(\delta^3) \ll 1$.

 For the third term, using the previous calculation we have
\begin{equation*}
  \int_{S_{\tau, r'}} |\psi|^2 \dVol_{\mathbb{S}^2} \lesssim  (r')^{1-p'} \int_{\Sigma_\tau \cap \{r_0 \leq r \leq r'\}} r^{p'} |\slashed{\D}_L \psi|^2 \upd r \wedge \dVol_{\mathbb{S}^2}
  + (r_0)^2 \mathcal{E}_0
\end{equation*}
 Multiplying by $(r')^{p-3+\delta}$ and integrating from $r' = r_0$ to $R$, we obtain
 \begin{equation*}
  \int_{r' = r_0}^{R} \left( \int_{S_{\tau, r'}} |\phi|^2 r^{p-1+\delta} \dVol_{\mathbb{S}^2} \right) \upd r' \lesssim  R^{p - p' - 1 + \delta} \int_{\Sigma_\tau \cap \{r_0 \leq r \leq r'\}} r^{p'} |\slashed{\D}_L \psi|^2 \upd r \wedge \dVol_{\mathbb{S}^2}
  + (r')^{p-2+\delta}(r_0)^2 \mathcal{E}_0 
 \end{equation*}
 we now choose $p' = p_{(\text{max})}$, to obtain
 \begin{equation*}
  \begin{split}
  &\int_{\Sigma_\tau \cap \{r_0 \leq r \leq R\}} |\phi|^2 r^{p - 3 + \delta} r^2 \upd r \wedge \dVol_{\mathbb{S}^2} \\
  &\lesssim  R^{p - p_{(\text{max})} - 1 + \delta} \int_{\Sigma_\tau \cap \{r_0 \leq r \leq r'\}} r^{p_{(\text{max})}} |\slashed{\D}_L \psi|^2 \upd r \wedge \dVol_{\mathbb{S}^2}
  + R^{p-2+\delta} (r_0)^2 \mathcal{E}_0 \\
  &\lesssim  R^{p - p_{(\text{max})} - 1 + \delta} (C_{[\phi]})^{-1} \delta^{-8\delta^{-1}\kappa} (\mathcal{E}_0 + \tilde{\mathcal{E}}_0) \left( (1+\tau)^{C_{[\phi]}\delta} + (1+\tau)^{1-\kappa} \right)
  + R^{p-2+\delta} \mathcal{E}_0
  \end{split}
 \end{equation*}
 Since $p < 1$ and $p \leq p_{(\text{max})} + \frac{1}{2}\delta$
 both of these terms tend to zero as $R \rightarrow \infty$. Note that both of these terms also appear inside an integral over $\tau$, with the weight $(1+\tau)^{-1-\delta}$.
 
 Finally, using \ref{proposition spherical mean in terms of energy} we can bound
 \begin{equation*}
  \int_{\tau_0}^{\tau} \left( \int_{S_{\tau', R}} p^{-1}r^{1-\delta}|\phi|^2 \dVol_{\mathbb{S}^2} \right) \upd \tau'
  \lesssim (\tau - \tau_0) R^{-\frac{1}{2}\delta}\delta^{-1}\mathcal{E}_0
 \end{equation*}
 
 In summary, we have
 \begin{equation*}
 \begin{split}
  &\mathcal{E}^{(L,p)}[\phi](\tau, R) 
  + \int_{^t\mathcal{M}_{\tau_0}^{\tau}} \chi_{(2r_0, R)}\left( pr^{p-1}|\slashed{\D}_L \phi|^2 + (2-p)r^{p-1}|\slashed{\nabla}\phi|^2 + p(1-p)r^{p-3}|\phi|^2 \right) \dVol_g \\
  &\lesssim
  \delta^{-8\delta^{-1}\kappa}(C_{[\phi]})^{-1}(\mathcal{E}_0 + \tilde{\mathcal{E}}_0) \bigg(
	  \delta^{-1+8\delta^{-1}\kappa}(1+\tau)^{C_{[\phi]}\delta} 
	  + (1+\tau)^{1-\kappa} 
	  + \frac{\epsilon}{(1-p+\delta)^2}(1+\tau)^{1-\kappa} 
	  \\
	  &\phantom{\lesssim
	  	\delta^{-1}(C_{[\phi]})^{-1}}
	  + \frac{\epsilon}{(1-p+\delta)^2}(1+\tau_0)^{1-\kappa} + \tilde{f}(\tau, R) \bigg)
  \int_{^t\mathcal{M}_{\tau_0}^{\tau}}\bigg( 
  \epsilon \chi_{(2r_0, R)} r^{p-2} (1+\tau)^{-1-\delta}|\slashed{\D}_L \psi|^2   \bigg) \dVol_g \\
 \end{split}
\end{equation*}
where
\begin{equation*}
 |\tilde{f}(\tau,R)| \lesssim R^{-1+\frac{3}{2}\delta}(1+\tau)^{1-\kappa-\delta} + R^{-1+\delta}(1+\tau)^{-\delta} + R^{-\frac{1}{2}\delta} (1+\tau)
\end{equation*}
Now, we can appeal to the Gronwall inequality (proposition \ref{proposition Gronwall}). We find that, for $p \leq p_{(\text{max})} + \frac{1}{2}\delta$ and $p < 1 - C_{[\phi]}\epsilon$, we have
\begin{equation*}
 \begin{split}
  &\mathcal{E}^{(L,p)}[\phi](\tau, R) 
  + \int_{^t\mathcal{M}_{\tau_0}^{\tau}} \chi_{(2r_0, R)}\left( pr^{p-1}|\slashed{\D}_L \phi|^2 + (2-p)r^{p-1}|\slashed{\nabla}\phi|^2 + p(1-p)r^{p-3}|\phi|^2 \right) \dVol_g \\
  &\lesssim
  \delta^{-8\delta^{-1}\kappa}(C_{[\phi]})^{-1}(\mathcal{E}_0 + \tilde{\mathcal{E}}_0) \left(  \delta^{-1+8\delta^{-1}\kappa}(1+\tau)^{C_{[\phi]}\delta} + (1+\tau)^{1-\kappa} + \frac{\epsilon}{(1-p+\delta)^2}(1+\tau)^{1-\kappa} + \tilde{f}(\tau,R) \right)
  \\
 \end{split}
\end{equation*}
In fact, the factor $\tilde{f}(\tau, R)$ can be dropped from the right hand side. This is because we can take the limit $R \rightarrow \infty$. The left hand side is monotonically non-decreasing as $R$ increases, while $\tilde{f}(\tau,R) \rightarrow 0$ as $R \rightarrow \infty$. Hence we have, in fact,
\begin{equation*}
 \begin{split}
  &\mathcal{E}^{(L,p)}[\phi](\tau, R) 
  + \int_{^t\mathcal{M}_{\tau_0}^{\tau}} \chi_{(2r_0, R)}\left( pr^{p-1}|\slashed{\D}_L \phi|^2 + (2-p)r^{p-1}|\slashed{\nabla}\phi|^2 + p(1-p)r^{p-3}|\phi|^2 \right) \dVol_g \\
  &\lesssim
  \delta^{-8\delta^{-1}\kappa}(C_{[\phi]})^{-1}(\mathcal{E}_0 + \tilde{\mathcal{E}}_0)\left( \delta^{-1+8\delta^{-1}\kappa}(1+\tau)^{C_{[\phi]}\delta} + (1+\tau)^{1-\kappa}  \right)
 \end{split}
\end{equation*}
meaning that the range of $p$ can in fact be extended to $p \leq \sup\{ 1-C_{[\phi]}\epsilon, p_{(\text{max})} + \frac{1}{2}\delta\}$. By iterating this argument, we can evidently increase $p$ up to the limit $p = (1 - C_{[\phi]}\epsilon)$. In other words, the second internal bootstrap bound \ref{equation large p internal bootstrap 2} holds for all $p \leq 1-C_{[\phi]}\epsilon$.

We also need to show that $\kappa$ can be increased. Suppose that we have
\begin{equation*}
\kappa \leq 1-C_{[\phi]}\delta
\end{equation*}
Let us also assume that
\begin{equation*}
8\delta^{-1}\kappa \geq 1
\end{equation*}
Then, by combining the result of the calculation above with lemma \ref{lemma boundedness} we find that we can bound
\begin{equation*}
 \begin{split}
  &\mathcal{E}^{(L, 1-C_{[\phi]}\epsilon)}[\phi](\tau)
  + \int_{\mathcal{M}^\tau_{\tau_0}} \bigg( 
  C_{[\phi]}\epsilon (1+r)^{-1- C_{[\phi]}\epsilon}|\slashed{\D}\phi|^2 + \delta^2 (1+r)^{-1- \delta}|\slashed{\D}\phi|^2 \\
  & \phantom{\int_{\mathcal{M}^\tau_{\tau_0}} \bigg(}
  + \chi_{(2r_0)}\left( r^{3-C_{[\phi]}\epsilon}|\slashed{\D}_L \phi|^2 + r^{-3C_{[\phi]}\epsilon }|\slashed{\nabla}\phi|^2 + C_{[\phi]}\epsilon r^{-2 - 3C_{[\phi]}\epsilon}|\phi|^2 \right) \bigg) \dVol_g \\
  & \lesssim \delta^{-8\delta^{-1}\kappa}(C_{[\phi]})^{-1}(\mathcal{E}_0 + \tilde{\mathcal{E}}_0) (1+\tau)^{1-\kappa}
 \end{split}
\end{equation*}
for all $\tau$. Now, by the pigeonhole principle, we can pick a diadic sequence of times $\tau_n$ such that
\begin{equation*}
 \begin{split}
  &\int_{\Sigma_{\tau_n} } \bigg( 
  C_{[\phi]}\epsilon (1+r)^{-1- C_{[\phi]}\epsilon}|\slashed{\D}\phi|^2 + \delta^2 (1+r)^{-1- \delta}|\slashed{\D}\phi|^2 \\
  & \phantom{\int_{\mathcal{M}^\tau_{\tau_0}} \bigg(}
  + \chi_{(2r_0)}\left( r^{-3C_{[\phi]}\epsilon}|\slashed{\D}_L \phi|^2 + r^{-3C_{[\phi]}\epsilon}|\slashed{\nabla}\phi|^2 + C_{[\phi]}\epsilon r^{-2 - 3C_{[\phi]}\epsilon}|\phi|^2 \right) \bigg) r^2 \upd r \wedge \dVol_{\mathbb{S}^2} \\
  & \lesssim \delta^{-8\delta^{-1}\kappa}(C_{[\phi]})^{-1}(\mathcal{E}_0 + \tilde{\mathcal{E}}_0) (1+\tau_n)^{-\kappa}
 \end{split}
\end{equation*}
In particular, at these times we have
\begin{equation*}
  \mathcal{E}^{(wT)}[\phi](\tau_n) \lesssim \delta^{-8\delta^{-1}\kappa}(C_{[\phi]})^{-1}(\mathcal{E}_0 + \tilde{\mathcal{E}}_0) (1+\tau_n)^{-\kappa}
\end{equation*}

Moreover, the above calculations also lead to the conclusion that for all $\tau \geq \tau_0$ we have
\begin{equation*}
 \mathcal{E}^{(L, 1-C_{[\phi]}\epsilon)}[\phi](\tau) \lesssim \delta^{-8\delta^{-1}\kappa}(C_{[\phi]})^{-1}(\mathcal{E}_0 + \tilde{\mathcal{E}}_0) (1+\tau_n)^{1-\kappa}
\end{equation*}

We need to interpolate between these two inequalities in order to obtain decay in $\tau$ for the quantity $\mathcal{E}^{(L,\frac{1}{2}\delta)}[\phi](\tau)$. Using H\"older's inequality, at the times $\tau_n$ we have
\begin{equation*}
  \begin{split}
    \mathcal{E}^{(L,\frac{1}{2}\delta)}[\phi](\tau_n) 
    &= \int_{\Sigma_{\tau_n}} \chi_{(2r_0)} r^{\frac{1}{2}\delta} |\slashed{\D}_L \psi|^2 \upd r \wedge \dVol_{\mathbb{S}^2} \\
    &\lesssim \left( \int_{\Sigma_{\tau_n}} \chi_{(2r_0)} r^{1-C_{[\phi]}\epsilon } |\slashed{\D}_L \psi|^2 \upd r \wedge \dVol_{\mathbb{S}^2} \right)^{\frac{1}{p}} \left( \int_{\Sigma_{\tau_n}} \chi_{(2r_0)} r^{-C_{[\phi]}\epsilon} |\slashed{\D}_L \psi|^2 \upd r \wedge \dVol_{\mathbb{S}^2} \right)^\frac{1}{q}
  \end{split}
\end{equation*}
where, as usual, we have set $\psi := r \phi$, and where the exponents are chosen as
\begin{equation*}
  \begin{split}
    p &= \frac{1}{\frac{1}{2}\delta + C_{[\phi]}\epsilon} \\
    q &= \frac{1}{1 - \frac{1}{2}\delta - C_{[\phi]}\epsilon}
  \end{split}
\end{equation*}
Moreover, we have
\begin{equation*}
 \begin{split}
  \int_{\Sigma_{\tau_n}} \chi_{(2r_0)} r^{-C_{[\phi]}\epsilon} |\slashed{\D}_L \psi|^2 \upd r \wedge \dVol_{\mathbb{S}^2}
  &\lesssim \int_{\Sigma_{\tau_n}} w \left( |\slashed{\D} \phi|^2 + r^{-1}(1+r)^{-1}|\phi|^2 \right) r^2 \upd r \wedge \dVol_{\mathbb{S}^2} \\
  &\lesssim \mathcal{E}^{(wT)}[\phi](\tau_n)
 \end{split}
\end{equation*}
where in the last line we have used the Hardy inequality.

Consequently, we find that
\begin{equation*}
 \mathcal{E}^{(L,\frac{1}{2}\delta)}[\phi](\tau_n)
 \lesssim \delta^{-8\delta^{-1}\kappa}(C_{[\phi]})^{-1}(\mathcal{E}_0 + \tilde{\mathcal{E}}_0) (1+\tau_n)^{- \kappa'} 
\end{equation*}
where $\kappa'$ is defined as
\begin{equation*}
  \kappa'  := \kappa - \frac{1}{2}\delta - C_{[\phi]}\epsilon
\end{equation*}
Note that, since $\epsilon \ll \delta \ll 1$, we have
\begin{equation*}
 \kappa' \geq \kappa - \frac{3}{4}\delta
\end{equation*}

Hence, \emph{at the times} $\tau_n$, we have
\begin{equation*}
  \mathcal{E}^{(wT)}[\phi](\tau_n) + \mathcal{E}^{(L,\frac{1}{2}\delta)}[\phi](\tau_n)
 \lesssim \delta^{-8\delta^{-1}\kappa}(\mathcal{E}_0 + \tilde{\mathcal{E}}_0) (1+\tau_n)^{-\kappa + \frac{3}{4}\delta}
\end{equation*}
Now, we can apply lemma \ref{lemma boundedness} at the intermediate times $\tau \in [\tau_n , \tau_{n+1}]$. In particular, we find that
\begin{equation*}
 \int_{\mathcal{M}^\tau_{\tau_n}} \epsilon (1+r)^{-1-C_{[\phi]}\epsilon} |\slashed{\D}\phi|^2 \dVol_g \lesssim (C_{[\phi]})^{-1} \delta^{-1-8\delta^{-1}\kappa}(\mathcal{E}_0 + \tilde{\mathcal{E}}_0) (1+\tau_n)^{-\kappa + \frac{3}{4}\delta}
\end{equation*}
Since the sequence $\{\tau_n\}$ is diadic, and $\tau_n \leq \tau \leq \tau_{n+1}$, in fact we have, for all $\tau_1 \geq \tau \geq \tau_0$
\begin{equation*}
 \int_{\mathcal{M}^\tau_{\tau_1}} \epsilon (1+r)^{-1-C_{[\phi]}\epsilon} |\slashed{\D}\phi|^2 \dVol_g \lesssim (C_{[\phi]})^{-1}\delta^{-1-8\delta^{-1}\kappa}(\mathcal{E}_0 + \tilde{\mathcal{E}}_0) (1+\tau_n)^{-\kappa + \frac{3}{4}\delta}
\end{equation*}

Now, by again decomposing the interval $[\tau_0, \tau_1]$ into a sum of diadic intervals, we have
\begin{equation*}
 \begin{split}
  &\int_{\mathcal{M}^{\tau}_{\tau_0}} \epsilon (1+\tau)^{1-C_{[\phi]}\epsilon} (1+r)^{-1-C_{[\phi]}\epsilon} |\slashed{\D}\phi|^2 \dVol_g \\
  &\leq (C_{[\phi]})^{-1}\sum_{n=0}^{\lceil \log (\tau - \tau_0) \rceil} \int_{\mathcal{M}^{\tau_0 + e^n}_{\tau_0 + e^{(n-1)}}} C_{[\phi]}\epsilon (1+\tau)^{1-C_{[\phi]}\epsilon} (1+r)^{-1-C_{[\phi]}\epsilon} |\slashed{\D}\phi|^2 \dVol_g \\
  &\lesssim (C_{[\phi]})^{-1}\sum_{n=0}^{\lceil \log (\tau - \tau_0) \rceil} \delta^{-1-8\delta^{-1}\kappa}\mathcal{E}_0 (1+\tau_0 + e^{(n-1)})^{1- \kappa - \frac{1}{4}\delta}
 \end{split}
\end{equation*}
Recall that we are assuming that
\begin{equation*}
 \kappa \leq 1 - C_{[\phi]}\delta
\end{equation*}
The sum above is bounded by a geometric series, and we find that
\begin{equation*}
 \begin{split}
  \int_{\mathcal{M}^{\tau}_{\tau_0}} \epsilon (1+\tau)^{1-C_{[\phi]}\epsilon} (1+r)^{-1-C_{[\phi]}\epsilon} |\slashed{\D}\phi|^2 \dVol_g
  \lesssim (C_{[\phi]})^{-2}(\delta)^{-2-8\delta^{-1}\kappa}\mathcal{E}_0 (1+\tau)^{1-\kappa-\frac{1}{4}\delta}
 \end{split}
\end{equation*}
Now, we note that
\begin{equation*}
\delta^{-2 -8\delta^{-1}\kappa} = \delta^{-8\delta^{-1}(\kappa + \frac{1}{4}\delta)}
\end{equation*}

Hence we find that, as long as $\kappa \leq 1 - C_{[\phi]}\delta$, we can improve $\kappa$ to $\kappa + \frac{1}{4}\delta$. In other words, beginning with the assumption that the integrated local energy decays like $\tau^{1-\kappa}$, we have shown that, in fact, the integrated local energy decays like $\tau^{1-\kappa - \frac{1}{4}\delta}$. Recall that lemma \ref{lemma boundedness} ensures that the inequality holds for $\kappa = \frac{1}{2}\delta$.

Iterating the argument above, we can clearly increase $\kappa$ up to a point where
\begin{equation*}
 \kappa = 1 - C_{[\phi]}\delta
\end{equation*}
Provided only that $C_{[\phi]}$ is sufficiently large. In other words, the internal bootstrap bounds \ref{equation large p internal bootstrap 1} and \ref{equation large p internal bootstrap 2} can be improved.

%We can now follow through the steps of the argument one more time, although this time we do not improve the decay by the full $\frac{1}{4}\delta$, since when we come to sum up the integrals we obtain
%\begin{equation*}
% \begin{split}
%  &\int_{\mathcal{M}^{\tau}_{\tau_0}} \epsilon (1+\tau)^{1-\delta} (1+r)^{-1-C_{[\phi]}\epsilon} |\slashed{\D}\phi|^2 \dVol_g \\
%  &\lesssim (C_{[\phi]})^{-1}\sum_{n=0}^{\lceil \log (\tau - \tau_0) \rceil} \delta^{-2}\mathcal{E}_0 (1+\tau_0 + e^{(n-1)})^{-\kappa + \frac{3}{4}\delta}
% \end{split}
%\end{equation*}
%This time, the series above is summable, and so we simply obtain
%\begin{equation*}
% \int_{\mathcal{M}^{\tau}_{\tau_0}} \epsilon (1+\tau)^{1-\delta} (1+r)^{-1-C_{[\phi]}\epsilon} |\slashed{\D}\phi|^2 \dVol_g \\
% \lesssim \delta^{-1}\mathcal{E}_0 
%\end{equation*}
%which is now uniform in $\tau$. In other words, we can increase $\kappa$ to set $\kappa = 1$.

\end{proof}

Inspection of the proof above also yields the following two decay estimates, the first of which will be crucial in obtaining pointwise decay in $\tau$ and eventually closing the bootstrap bounds of chapter \ref{chapter bootstrap}. In fact, the following two corollaries hold with decay at the rate $(1+\tau)^{-1 + \frac{3}{4}\delta}$, but to prevent the notation from becoming cluttered we state them with the decay rate $(1+\tau)^{-1 + \delta}$, which is sufficient for our purposes.

\begin{corollary}[Degenerate energy decay]
\label{corollay energy decay}
 Let the conditions of lemma \ref{lemma p weighted} hold. Let
 \begin{equation*}
  w = (1+r)^{-C_{[\phi]}\epsilon}
 \end{equation*}
 Then the degenerate energy decays in $\tau$ as
 \begin{equation*}
  \mathcal{E}^{(wT)}[\phi](\tau) \lesssim \delta^{-8\delta^{-1}}(\mathcal{E}_0 + \tilde{\mathcal{E}}_0) (1+\tau)^{-1 + C_{[\phi]}\delta}
 \end{equation*}
 for all $\tau \geq \tau_0$.
\end{corollary}

\begin{corollary}[Decay of the integrated local energy]
\label{corollary ILED decay}
Let the conditions of lemma \ref{lemma p weighted} hold.

 Then integrated local energy decays in $\tau$ as
 \begin{equation*}
   \begin{split}
    &\int_{\mathcal{M}^{\tau_1}_{\tau}} \bigg( \delta^2 (1+r)^{-1-\frac{1}{2}\delta}|\slashed{\D}\phi|^2 + \delta \chi_{(r_0)}(r)(1+r)^{-1-C_{[\phi]}\epsilon}|\slashed{\nabla}\phi|^2 + \delta^2 r^{-1}(1+r)^{-2-\frac{1}{2}\delta}|\phi|^2 \\
    &\phantom{+ \int_{^t\mathcal{M}^{\tau_2}_{\tau_1}} \bigg(}
    + C_{[\phi]}\epsilon(1+r)^{-1-C_{[\phi]}\epsilon}|\slashed{\D}\phi|^2 \bigg)\dVol_g
    \lesssim
    \delta^{-8\delta^{-1}}(\mathcal{E}_0 + \tilde{\mathcal{E}}_0) (1+\tau)^{-1 + C_{[\phi]}\delta}
   \end{split}
 \end{equation*}
 for all $\tau \geq \tau_0$.
 
\end{corollary}

Note that this second corollary follows from the ``internal boostrap'' argument employed in the proof of the lemma, with $\kappa = 1-\delta$. Strictly speaking, we have only shown this for the last term in the integrand, but it is easy to modify the proof above to include the other spacetime integral terms given above.

\begin{corollary}[Decay of the $p$ weighted energy for $p = \frac{1}{2}\delta$]
	
	\label{corollary small p decay}
	Let the conditions of lemma \ref{lemma p weighted} hold.
	
	Then the $p$-weighted energy with $p = \frac{1}{2}\delta$ decays in $\tau$ as
	\begin{equation*}
	\begin{split}
	&\mathcal{E}^{(L, \frac{1}{2}\delta)}[\phi] (\tau)
	+ \int_{\mathcal{M}^{\tau_1}_\tau} \chi_{2r_0} r^{-1 + \frac{1}{2}\delta} \left( |\slashed{\D}_L \phi|^2 + |\slashed{\nabla}\phi|^2 + \delta r^{-3 + \frac{1}{2}\delta}|\phi|^2 \right) \dVol_g
	\\
	&\lesssim
	\delta^{-8\delta^{-1}}(\mathcal{E}_0 + \tilde{\mathcal{E}}_0) (1+\tau)^{-1 + C_{[\phi]}\delta}
	\end{split}
	\end{equation*}
	for all $\tau \geq \tau_0$.
	
\end{corollary}

\section{Energy estimates involving a point-dependent change of basis}

\label{section energy estimates involving a point-dependent change of basis}

The energy estimates proved above are suitable for analysing solutions to systems of wave equations of the form
\begin{equation*}
 \Box_g \phi_{(a)} = F_{(a)}
\end{equation*}
However, we also want to consider systems of wave equations which require a point dependent change of basis in order to reveal their structure, e.g. wave equations of the same form but where the weak null structure is only apparent after transforming to a new basis of fields $\phi_{(A)}$, defined by
\begin{equation*}
 \phi_{(A)} := M_{(A)}^{\phantom{(A)}(a)}\phi_(a)
\end{equation*}
and where the change-of-basis matrix $M$ can depend on the point on the manifold $\mathcal{M}$, possible through a dependence on the fields $\phi_{(a)}$ themselves. An example of such a system is the Einstein equations in harmonic coordinates, for which we must change from rectangular coordinates to a system of coordinates based on the null frame. That is, in order to exploit the semilinear structure in the Einstein equations, we must transform from the fields $h_{ab}$, which are the rectangular components of the metric perturbation, to the fields $X^a Y^b h_{ab}$, where $X^a$ and $Y^a$ are the rectangular components of the null frame fields $L$, $\Lbar$ and $\slashed{\Pi}_a$.

In such a case, the approach we take will be to \emph{first} prove energy estimates for the original components of $\phi$ (i.e.\ the components $\phi_{(a)}$, which correspond in the Einstein equations to the rectangular components $h_{ab}$, and for which we can use the energy estimates above) and \emph{then} to prove energy estimates for fields in the new basis (corresponding, in the case of the Einstein equations, to the null frame components of the fields $h$). Hence, for the estimates in this section, we will allow error terms involving the energy of the field in the original basis.

\begin{lemma}[Energy boundedness, the integrated local energy decay estimate and the $p$-weighted energy estimate for very small $p$ after a point-dependent change of basis]
\label{lemma boundedness after change of basis}
 Let $\phi_{(a)}$ be a set of $S_{\tau,r}$-tangent tensor field labelled by $(a)$, satisfying
\begin{equation*}
 \tilde{\slashed{\Box}}_g \phi_{(a)} = F_{(a)} = F_{(a,1)} + F_{(a,2)} = = F_{(a,4)} + F_{(a,5)} + F_{(a,6)}
\end{equation*}
for some collection of $S_{\tau, r}$ tangent tensor fields $F_{(a)}$, $F_{(a,1)} \ldots F_{(a,6)}$. Let $\tau_2 \geq \tau_1$. Suppose that $F_{(a)}$, $F_{(a,1)}$ and $F_{(a,2)}$ satisfy the bounds
\begin{equation*}
 \begin{split}
  &\int_{^t\mathcal{M}^{\tau_2}_{\tau_1}} \epsilon^{-1}\bigg(
  	(1+r)^{1-C_{(a,\text{orig})}\epsilon} |F_{(a)}|^2
  	+ (1+r)^{\frac{1}{2}\delta}(1+\tau)^{1 + \delta} |F_{(a,1)}|^2
  	+ (1+r)^{1-\frac{7}{2}\delta}(1+\tau)^{6\delta}|F_{(a,2)}|^2
  	\\
  	&\phantom{\int_{^t\mathcal{M}^{\tau_2}_{\tau_1}} \epsilon^{-1}\bigg(}
  	+ (1+r)^{1+\frac{1}{2}\delta}|F_{(a,3)}|^2
  \bigg)\dVol_g
  \leq \tilde{\mathcal{E}}_{(\text{orig})} (1+\tau)^{-1+C_{[\phi]}\delta}
 \end{split}
\end{equation*}

Choose the weight functions associated to the $\phi_{(a)}$
\begin{equation*}
 w_{(a,\text{orig})} = (1+r)^{-C_{(a,\text{orig})}\epsilon}
\end{equation*}

Suppose, moreover, that the initial energy of $\phi_{(a)}$ satisfies
\begin{equation*}
 \mathcal{E}^{(w_{(a,\text{orig})}T)}[\phi_{(a)}](\tau_2, t, \tau_1) \leq \mathcal{E}_{(\text{orig})}
\end{equation*}
and
\begin{equation*}
 \mathcal{E}^{(L, \frac{1}{2}\delta)}[\phi_{(a)}](\tau_1,R) \leq \mathcal{E}_{(\text{orig})}
\end{equation*}
where $t$ is sufficiently large relative to $R$, $\tau_1$ and $\tau_2$ so that
\begin{equation}
 \label{equation relationship between R and t 2}
 \{r \leq R\} \cap \left\{^{\tau_2}_{\tau_1}\bar{\Sigma}_t \right\} = \emptyset
\end{equation}

Define the modified weight
\begin{equation*}
 \tilde{w}_{(a,\text{orig})} = (1+r)^{-\frac{1}{2}C_{(a,\text{orig})}\epsilon}
\end{equation*}
Suppose additionally that there is some $\tau_0 \leq \tau_1$ such that
\begin{equation*}
 \mathcal{E}^{(\tilde{w}_{(a,\text{orig})}T)}(\tau_0, t, \tau_0) \leq \tilde{\mathcal{E}}_{(\text{orig})}
\end{equation*}
Furthermore, on the initial hypersurface $^t\Sigma_{\tau_0}$ suppose that we have
\begin{equation*}
 \int_{\bar{S}_{t, r}}|\phi_{(a)}|^2 \dVol_{\mathbb{S}^2} \leq \tilde{\mathcal{E}}_{(\text{orig})} (t - \tau_0)^{-1 + \frac{1}{2}C_{[\phi]}\epsilon}
\end{equation*}

Let $\phi_{(A)} = M_{(A)}^{\phantom{(A)}(a)}\phi_{(a)}$, where the change-of-basis matrix $M_{(A)}^{\phantom{(A)}(a)}$ is of maximal rank. Moreover, let the matrix $M$ satisfy
\begin{equation}
 \begin{split}
   |M_{(A)}^{\phantom{(A)}(a)}| &\lesssim 1 \\
   |\bar{\partial} M_{(A)}^{\phantom{(A)}(a)}| &\lesssim (1+r)^{-1 - \delta} \\
   |\partial M_{(A)}^{\phantom{(A)}(a)}| &\lesssim (1+r)^{-1} + (1+r)^{-1 + \delta}(1+\tau)^{-\beta}
 \end{split}
\end{equation}

Denote
\begin{equation*}
 \begin{split}
  F_{(A)} &:= M_{(A)}^{\phantom{(A)}(a)}F_{(a)} \\
  F_{(A)} &= F_{(A,1)} + F_{(A,2)} + F_{(A,3)} = F_{(A,4)} + F_{(A,5)} + F_{(A,6)}
 \end{split}
\end{equation*}
 and suppose that $F_{(A)}$, $F_{(A,1)}$, $F_{(A,2)}$ and $F_{(A,3)}$ satisfy the bounds
\begin{equation*}
 \begin{split}
  &\int_{^t\mathcal{M}^{\tau_2}_{\tau_1}} \epsilon^{-1}\bigg(
  	(1+r)^{1-C_{(A,\text{new})}\epsilon} |F_{(A)}|^2
  	+ (1+r)^{\frac{1}{2}\delta}(1+\tau)^{1 + \delta} |F_{(A,1)}|^2
  	+ (1+r)^{1-\frac{7}{2}\delta}(1+\tau)^{6\delta}|F_{(A,2)}|^2
  	\\
  	&\phantom{\int_{^t\mathcal{M}^{\tau_2}_{\tau_1}} \epsilon^{-1}\bigg(}
  	+ (1+r)^{1+\frac{1}{2}\delta}|F_{(A,3)}|^2
  \bigg)\dVol_g
  \leq \tilde{\mathcal{E}}_{(\text{new})} (1+\tau)^{-1+C_{[\phi]}\delta}
 \end{split}
\end{equation*}

Choose a new weight function, associated with the new fields $\phi_{(A)}$
\begin{equation*}
 w_{(A,\text{new})} = (1+r)^{-C_{(A,\text{new})}\epsilon}
\end{equation*}

Suppose, moreover, that the initial energy of the fields $\phi_{(A)}$ satisfy
\begin{equation*}
 \mathcal{E}^{(w_{(A,\text{new})}T)}[\phi]_{(A)}(\tau_2, t, \tau_1) \leq \mathcal{E}_{(\text{new})}
\end{equation*}
and
\begin{equation*}
 \mathcal{E}^{(L, \frac{1}{2}\delta)}[\phi]_{(A)}(\tau_1,R) \leq \mathcal{E}_{(\text{new})}
\end{equation*}

Define the a new modified weight
\begin{equation*}
 \tilde{w}_{(A,\text{new})} = (1+r)^{-\frac{1}{2}C_{(A,\text{new})}\epsilon}
\end{equation*}
Suppose additionally that at $\tau = \tau_0$ we have
\begin{equation*}
 \mathcal{E}^{(\tilde{w}_{(A,\text{new})}T)}(\tau_0, t, \tau_0) \leq \tilde{\mathcal{E}}_{(\text{new})}
\end{equation*}
Furthermore, on the initial hypersurface $^t\Sigma_{\tau_0}$ suppose that we have
\begin{equation*}
 \int_{\bar{S}_{t, r}}|\phi_{(A)}|^2 \dVol_{\mathbb{S}^2} \leq \tilde{\mathcal{E}}_{(\text{new})} (t - \tau_0)^{-1 + \frac{1}{2}C_{(A, \text{new})}\epsilon}
\end{equation*}

Suppose that all the bootstrap bounds of chapter \ref{chapter bootstrap} are satisfied.

Then, for all sufficiently small $\delta$, for all sufficiently small $\epsilon$ we have

\begin{equation}
 \begin{split}
  & \delta^3 \mathcal{E}^{(L, \frac{1}{2}\delta)}[\phi]_{(A)}(\tau_2, R)
  + \mathcal{E}^{(w_{(A, \text{new})}T)}[\phi]_{(A)}(\tau_2, t, \tau_1) \\
  &
  + \int_{^t\mathcal{M}^{\tau_2}_{\tau_1}} \bigg( 
  	\delta^2 (1+r)^{-1-\delta}|\slashed{\D}\phi|_{(A)}^2 
  	+ \delta \chi_{(r_0)}(r)(1+r)^{-1 - C_{(A, \text{new})}\epsilon}|\slashed{\nabla}\phi|_{(A)}^2 
  	+ \delta^2 r^{-1}(1+r)^{-2-\delta}|\phi|_{(A)}^2 
  \\
  &\phantom{+ \int_{^t\mathcal{M}^{\tau_2}_{\tau_1}} \bigg(}
  + C_{(A, \text{new})}\epsilon(1+r)^{-1-C_{(A, \text{new})}\epsilon}|\slashed{\D}\phi|_{(A)}^2 \\
  &\phantom{+ \int_{^t\mathcal{M}^{\tau_2}_{\tau_1}} \bigg(}
  + \chi_{(2r_0, R)}\left( \delta^4 r^{-1+\frac{1}{2}\delta}|\slashed{\D}_L \phi|_{(A)}^2 + \delta^3 r^{-1 + \frac{1}{2}\delta|}\slashed{\nabla}\phi|_{(A)}^2 + \delta^3 r^{-3 + \frac{1}{2}\delta}|\phi|_{(A)}^2 \right) \bigg)\dVol_g \\
  &\lesssim
  \delta^2 \mathcal{E}^{(L, \frac{1}{2}\delta)}[\phi]_{(A)}(\tau_1, R)
  + \sum_{(a)} \mathcal{E}^{(L, \frac{1}{2}\delta)}[\phi_{(a)}](\tau_1, R)
  + \delta^{-8\delta^{-1}}\mathcal{E}^{(w_{(A, \text{new})}T)}[\phi]_{(A)}(\tau_1, t, \tau_1)
  \\
  &\phantom{\lesssim}
  + \delta^{-2-8\delta^{-1}} \sum_{(a)}\mathcal{E}^{(w_{(A, \text{orig})}T)}[\phi_{(a)}](\tau_1, t, \tau_1)
  + \delta^{-8\delta^{-1}} \tilde{\mathcal{E}}_{(\text{new})} + \delta^{-2-8\delta^{-1}}\tilde{\mathcal{E}}_{(\text{orig})} \\
  &\phantom{\lesssim}
  + \epsilon^{-1} \left(\mathcal{E}_{(\text{new})} + \delta^{-2}\mathcal{E}_{(\text{orig})} \right) \left( (t - \tau)^{-\frac{1}{2}C_{(a, \text{orig})}\epsilon} + (\tau - \tau_1)(t - \tau_0)^{-1 + \frac{3}{2}\delta} + \epsilon^{-1}R^{-1 + \frac{3}{2}\delta}e^{\tilde{C}\epsilon(\tau - \tau_0)} \right) \\
  &\phantom{\lesssim}
  + \left(\mathcal{E}_{(\text{new})} + \delta^{-2}\mathcal{E}_{(\text{orig})} \right)\left( (t - \tau)^{-\frac{1}{2}C_{(A, \text{new})}\epsilon} + (\tau - \tau_1)(t - \tau_0)^{-1 + \frac{3}{2}\delta} + \epsilon^{-1}R^{-1 + \frac{3}{2}\delta}e^{\tilde{C}\epsilon(\tau - \tau_0)} \right)
 \end{split}
\end{equation}

In particular, suppose that the conditions of the lemma hold in the limit $R \rightarrow \infty$, $t \rightarrow \infty$, where the limit is taken such that \eqref{equation relationship between R and t 2} is true. For example, we could take $t = 2(\tau_2 + R)$ which, under the bootstrap assumptions, can be seen to imply \eqref{equation relationship between R and t}. We could then take the limit $R \rightarrow \infty$, with $t$ considered a function of $R$, and $\tau_0$, $\tau_1$, $\tau_2$ fixed (finite) constants. Then we define the energies
\begin{equation}
 \begin{split}
  \mathcal{E}^{(w_{(A, \text{new})}T)}[\phi]_{(A)}(\tau) &:= \lim_{t \rightarrow \infty} \mathcal{E}^{(w_{(A, \text{new})}T)}[\phi]_{(A)}(\tau, t, \tau) \\
  \mathcal{E}^{(L, \alpha)}[\phi]_{(A)}(\tau) &:= \lim_{R \rightarrow \infty} \mathcal{E}^{(L, \alpha)}[\phi]_{(A)}(\tau, R) \\ 
 \end{split}
\end{equation}
Then we have the bound
\begin{equation}
 \begin{split}
  &\delta^3 \mathcal{E}^{(L, \frac{1}{2}\delta)}[\phi]_{(A)}(\tau_2)
  + \mathcal{E}^{(w_{(A, \text{new})}T)}[\phi]_{(A)}(\tau_2) \\
  &
  + \int_{\mathcal{M}^{\tau_2}_{\tau_1}} \bigg( \delta^2 (1+r)^{-1-\delta}|\slashed{\D}\phi|_{(A)}^2 + \delta \chi_{(r_0)}(r)(1+r)^{-1 - C_{(\phi)}\epsilon}|\slashed{\nabla}\phi|_{(A)}^2 + \delta^2 r^{-1}(1+r)^{-2-\delta}|\phi|_{(A)}^2 \\
  &\phantom{+ \int_{^t\mathcal{M}^{\tau_2}_{\tau_1}} \bigg(}
  + C_{(\phi)}\epsilon(1+r)^{-1-C_{(\phi)}\epsilon}|\slashed{\D}\phi|_{(A)}^2 \\
  &\phantom{+ \int_{^t\mathcal{M}^{\tau_2}_{\tau_1}} \bigg(}
  + \chi_{(2r_0)}\left(\delta^4 r^{-1+\frac{1}{2}\delta}|\slashed{\D}_L \phi|_{(A)}^2 + \delta^3 r^{-1 + \frac{1}{2}\delta|}\slashed{\nabla}\phi|_{(A)}^2 + \delta^3 r^{-3 + \frac{1}{2}\delta}|\phi|_{(A)}^2 \right) \bigg)\dVol_g \\
  &\lesssim
  \delta^2 \mathcal{E}^{(L, \frac{1}{2}\delta)}[\phi]_{(A)}(\tau_1)
  + \sum_{(a)} \mathcal{E}^{(L, \frac{1}{2}\delta)}[\phi_{(a)}](\tau_1)
  + \delta^{-8\delta^{-1}}\mathcal{E}^{(w_{(A, \text{new})}T)}[\phi]_{(A)}(\tau_1)
  \\
  &\phantom{\lesssim}
  + \delta^{-2-8\delta^{-1}} \sum_{(a)}\mathcal{E}^{(w_{(A, \text{orig})}T)}[\phi_{(a)}](\tau_1)
  + \delta^{-8\delta^{-1}} \tilde{\mathcal{E}}_{(\text{new})} 
  + \delta^{-2-8\delta^{-1}}\tilde{\mathcal{E}}_{(\text{orig})}
 \end{split}
\end{equation}
In particular, this implies the following three estimates:
\begin{enumerate}
 \item Degenerate energy boundedness:
 \begin{equation*}
 \begin{split}
  \mathcal{E}^{(w_{(A, \text{new})}T)}[\phi]_{(A)}(\tau_2)
  &\lesssim
  \delta^2 \mathcal{E}^{(L, \frac{1}{2}\delta)}[\phi]_{(A)}(\tau_1)
  + \sum_{(a)} \mathcal{E}^{(L, \frac{1}{2}\delta)}[\phi_{(a)}](\tau_1)
  + \delta^{-1}\mathcal{E}^{(w_{(A, \text{new})}T)}[\phi]_{(A)}(\tau_1)
  \\
  &\phantom{\lesssim}
  + \delta^{-2 -8\delta^{-1}} \sum_{(a)}\mathcal{E}^{(w_{(A, \text{orig})}T)}[\phi_{(a)}](\tau_1)
  + \delta^{-8\delta^{-1}} \tilde{\mathcal{E}}_{(\text{new})} 
  + \delta^{-2 -8\delta^{-1}}\tilde{\mathcal{E}}_{(\text{orig})}
 \end{split}
 \end{equation*}
 \item Integrated local energy decay:
 \begin{equation*}
  \begin{split}
    &\int_{\mathcal{M}^{\tau_2}_{\tau_1}} \bigg( \delta^2 (1+r)^{-1-\delta}|\slashed{\D}\phi|_{(A)}^2 + \delta \chi_{(r_0)}(r)(1+r)^{-1-C_{(\phi)}\epsilon}|\slashed{\nabla}\phi|_{(A)}^2 + \delta^2 r^{-1}(1+r)^{-2-\delta}|\phi|_{(A)}^2 \\
    &\phantom{+ \int_{^t\mathcal{M}^{\tau_2}_{\tau_1}} \bigg(}
    + C_{(\phi)}\epsilon(1+r)^{-1-C_{(\phi)}\epsilon}|\slashed{\D}\phi|_{(A)}^2 \bigg)\dVol_g \\
    &\lesssim
    \delta^2 \mathcal{E}^{(L, \frac{1}{2}\delta)}[\phi]_{(A)}(\tau_1)
    + \sum_{(a)} \mathcal{E}^{(L, \frac{1}{2}\delta)}[\phi_{(a)}](\tau_1)
    + \delta^{-8\delta^{-1}}\mathcal{E}^{(w_{(A, \text{new})}T)}[\phi]_{(A)}(\tau_1)
    \\
    &\phantom{\lesssim}
    + \delta^{-2-8\delta^{-1}} \sum_{(a)}\mathcal{E}^{(w_{(A, \text{orig})}T)}[\phi_{(a)}](\tau_1)
    + \delta^{-8\delta^{-1}} \tilde{\mathcal{E}}_{(\text{new})} 
    + \delta^{-2-8\delta^{-1}}\tilde{\mathcal{E}}_{(\text{orig})}
  \end{split}
 \end{equation*}
 \item $p$-weighted energy estimate with $p = \frac{1}{2}\delta$
 \begin{equation*}
  \begin{split}
    & \mathcal{E}^{(L, \frac{1}{2}\delta)}[\phi]_{(A)}(\tau_2) \\
    &
    + \int_{\mathcal{M}^{\tau_2}_{\tau_1}} \bigg(
    \chi_{(2r_0)}\left(\delta r^{-1+\frac{1}{2}\delta}|\slashed{\D}_L \phi|_{(A)}^2 + r^{-1 + \frac{1}{2}\delta|}\slashed{\nabla}\phi|_{(A)}^2 + r^{-3 + \frac{1}{2}\delta}|\phi|_{(A)}^2 \right) \bigg)\dVol_g \\
    &\lesssim
    \delta^{-8\delta^{-1}} \mathcal{E}^{(L, \frac{1}{2}\delta)}[\phi]_{(A)}(\tau_1)
    + \sum_{(a)} \delta^{-3} \mathcal{E}^{(L, \frac{1}{2}\delta)}[\phi_{(a)}](\tau_1)
    + \delta^{-3-8\delta^{-1}}\mathcal{E}^{(w_{(A, \text{new})}T)}[\phi]_{(A)}(\tau_1)
    \\
    &\phantom{\lesssim}
    + \delta^{-5-8\delta^{-1}} \sum_{(a)}\mathcal{E}^{(w_{(A, \text{orig})}T)}[\phi_{(a)}](\tau_1)
    + \delta^{-3-8\delta^{-1}} \tilde{\mathcal{E}}_{(\text{new})} 
    + \delta^{-5-8\delta^{-1}}\tilde{\mathcal{E}}_{(\text{orig})}
 \end{split} 
 \end{equation*}

\end{enumerate}

Additionally, if the inhomogeneous terms $F_{(A,4)}$, $F_{(A,5)}$ and $F_{(A,6)}$ satisfy
\begin{equation*}
 \begin{split}
  &\int_{\mathcal{M}^\tau_{\tau_0}} \epsilon^{-1} \chi_{(r_0)} \bigg(
  	r^{1-C_{(A, \text{new})}\epsilon}(1+\tau)^{1+\delta}|F_{(A,4)}|^2
  	+ r^{2 - C_{(A)}\epsilon -2\delta} (1+\tau)^{6\delta}|F_{(A,5)}|^2
  	\\
  	&\phantom{\int_{\mathcal{M}^\tau_{\tau_0}} \epsilon^{-1} \chi_{(r_0)} \bigg(}
  	+ r^{2-C_{(A, \text{new})}\epsilon}|F_{(A,6)}|^2
  \bigg) \dVol_g 
  \lesssim \tilde{\mathcal{E}}_0 (1+\tau)^{C_{[\phi]}\delta}
 \end{split}
\end{equation*}
and we suppose also that the initial data of $\phi_{(A)}$ and $\phi_{(a)}$ also satisfy
\begin{equation*}
 \begin{split}
  \mathcal{E}^{(L, 1-C_{(A)}\epsilon)}[\phi]_{(A)}(\tau_0) &\lesssim \mathcal{E}_0 \\
  \mathcal{E}^{(L, 1-C_{(a)}\epsilon)}[\phi_{(a)}](\tau_0) &\lesssim \mathcal{E}_0 \\
  \int_{\bar{S}_{t,r}}|\phi|_{(A)}^2 \dVol_{\mathbb{S}^2} &\lesssim \mathcal{E}_0(t - \tau_0)^{-1+\frac{1}{2}C_{(A)}\epsilon} \\
  \int_{\bar{S}_{t,r}}|\phi_{(a)}|^2 \dVol_{\mathbb{S}^2} &\lesssim \mathcal{E}_0(t - \tau_0)^{-1+\frac{1}{2}C_{(a)}\epsilon}
 \end{split} 
\end{equation*}

Then, for sufficiently small $\epsilon$, for sufficiently small $\delta$, and for sufficiently large constants $C_{(A)}$, $C_{(a)}$, for all $\tau \geq \tau_0$ we have
\begin{equation*}
 \begin{split}
  &\mathcal{E}^{(L, 1-C_{(A)}\epsilon)}[\phi]_{(A)}(\tau)
  + \int_{\mathcal{M}^{\tau}_{\tau_0}} \chi_{(2r_0)}r^{-C_{(A)}\epsilon} \left( |\slashed{\D}_L \phi|_{(A)}^2 + |\slashed{\nabla}\phi|_{(A)}^2 + C_{(A)}\epsilon r^{-2}|\phi|_{(A)}^2 \right) \dVol_g 
  \\
  &\lesssim
  \mathcal{E}^{(L, 1-C_{(A)}\epsilon)}[\phi]_{(A)}(\tau_1)
  + \sum_{(a)} \delta^{-2} \mathcal{E}^{(L, 1-C_{(a)}\epsilon)}[\phi_{(a)}](\tau_1)
  + \delta^{-3-8\delta^{-1}}\mathcal{E}^{(w_{(A, \text{new})}T)}[\phi]_{(A)}(\tau_1)
  \\
  &\phantom{\lesssim}
  + \delta^{-5-8\delta^{-1}} \sum_{(a)}\mathcal{E}^{(w_{(A, \text{orig})}T)}[\phi_{(a)}](\tau_1)
  + \delta^{-3-8\delta^{-1}} \tilde{\mathcal{E}}_{(\text{new})} 
  + \delta^{-5-8\delta^{-1}}\tilde{\mathcal{E}}_{(\text{orig})}
 \end{split}
\end{equation*}

Moreover, the degenerate energy decays in $\tau$ as
	\begin{equation*}
	\begin{split}
	(1+\tau)^{1-C_{[\phi]}\delta}\mathcal{E}^{(w_{(A, \text{new})}T)}[\phi]_{(A)}(\tau)
	&\lesssim
	\delta^2 \mathcal{E}^{(L, 1-C_{(A)}\epsilon)}[\phi]_{(A)}(\tau_0)
	+ \sum_{(a)} \mathcal{E}^{(L,  1-C_{(a)}\epsilon)}[\phi_{(a)}](\tau_0)
	\\
	&\phantom{\lesssim}
	+ \delta^{-8\delta^{-1}}\mathcal{E}^{(w_{(A, \text{new})}T)}[\phi]_{(A)}(\tau_0)
	+ \delta^{-2-8\delta^{-1}} \sum_{(a)}\mathcal{E}^{(w_{(A, \text{orig})}T)}[\phi_{(a)}](\tau_0)
	\\
	&\phantom{\lesssim}
	+ \delta^{-8\delta^{-1}} \tilde{\mathcal{E}}_{(\text{new})} 
	+ \delta^{-2-8\delta^{-1}}\tilde{\mathcal{E}}_{(\text{orig})}
	\end{split}
	\end{equation*}
and the integrated local energy decays in $\tau$ as:
	\begin{equation*}
	\begin{split}
	&(1+\tau)^{1-C_{[\phi]}\delta}\int_{\mathcal{M}^{\tau}_{\tau_1}} \bigg( 
		\delta^2 (1+r)^{-1-\delta}|\slashed{\D}\phi|_{(A)}^2 
		+ \delta \chi_{(r_0)}(r)(1+r)^{-1-C_{(\phi)}\epsilon}|\slashed{\nabla}\phi|_{(A)}^2
		\\
		&\phantom{(1+\tau)^{1-C_{[\phi]}\delta}\int_{\mathcal{M}^{\tau}_{\tau_1}} \bigg(}
		+ \delta^2 r^{-1}(1+r)^{-2-\delta}|\phi|_{(A)}^2 
		+ C_{(\phi)}\epsilon(1+r)^{-1-C_{(\phi)}\epsilon}|\slashed{\D}\phi|_{(A)}^2
	\bigg)\dVol_g \\
	&\lesssim
	\delta^{2} \mathcal{E}^{(L, 1-C_{(A)}\epsilon)}[\phi]_{(A)}(\tau_0)
	+ \sum_{(a)} \mathcal{E}^{(L,  1-C_{(a)}\epsilon)}[\phi_{(a)}](\tau_0)
	+ \delta^{-8\delta^{-1}}\mathcal{E}^{(w_{(A, \text{new})}T)}[\phi]_{(A)}(\tau_0)
	\\
	&\phantom{\lesssim}
	+ \delta^{-2-8\delta^{-1}} \sum_{(a)}\mathcal{E}^{(w_{(A, \text{orig})}T)}[\phi_{(a)}](\tau_0)
	+ \delta^{-8\delta^{-1}} \tilde{\mathcal{E}}_{(\text{new})} 
	+ \delta^{-2-8\delta^{-1}}\tilde{\mathcal{E}}_{(\text{orig})}
	\end{split}
	\end{equation*}
and finally, the $p$-weighted energy estimate with $p = \frac{1}{2}\delta$ decays in time as
	\begin{equation*}
	\begin{split}
	& (1+\tau)^{1-C_{[\phi]}\delta}\mathcal{E}^{(L, \frac{1}{2}\delta)}[\phi]_{(A)}(\tau)
	\\
	&
	+ (1+\tau)^{1-C_{[\phi]}\delta}\int_{\mathcal{M}^{\tau}_{\tau_1}} \bigg(
	\chi_{(2r_0)}\left(\delta r^{-1+\frac{1}{2}\delta}|\slashed{\D}_L \phi|_{(A)}^2 + r^{-1 + \frac{1}{2}\delta|}\slashed{\nabla}\phi|_{(A)}^2 + r^{-3 + \frac{1}{2}\delta}|\phi|_{(A)}^2 \right) \bigg)\dVol_g \\
	&\lesssim
	\delta^{-1} \mathcal{E}^{(L, \frac{1}{2}\delta)}[\phi]_{(A)}(\tau_1)
	+ \sum_{(a)} \delta^{-3} \mathcal{E}^{(L, \frac{1}{2}\delta)}[\phi_{(a)}](\tau_1)
	+ \delta^{-3-8\delta^{-1}}\mathcal{E}^{(w_{(A, \text{new})}T)}[\phi]_{(A)}(\tau_1)
	\\
	&\phantom{\lesssim}
	+ \delta^{-5-8\delta^{-1}} \sum_{(a)}\mathcal{E}^{(w_{(A, \text{orig})}T)}[\phi_{(a)}](\tau_1)
	+ \delta^{-3-8\delta^{-1}} \tilde{\mathcal{E}}_{(\text{new})} 
	+ \delta^{-5-8\delta^{-1}}\tilde{\mathcal{E}}_{(\text{orig})}
	\end{split} 
	\end{equation*}

\end{lemma}

\begin{proof}
 The proof proceeds in almost the same way as lemmas \ref{lemma boundedness} and \ref{lemma p weighted} and corollaries \ref{corollay energy decay} and \ref{corollary ILED decay}, with the only new error terms arising from derivatives of $M$. Hence we will only sketch the proof here, indicating the points at which the proof differs from those mentioned above.
 
 The new error terms in the energy estimates can be read off from proposition \ref{proposition energy identity after change of basis}, specifically the expression for the error term $\textit{Err}_{(\partial M)}[\phi]_{(A)}$, with the choices 
 \begin{equation*}
  (Z, f_Z) = \begin{cases} (w_{(A, \text{new})}T, 0) \\
  (w_{(A, \text{new})}f_R R, 2 r^{-1} w_{(A, \text{new})} f_R  ) \\
  (f_L r^{\frac{1}{2}\delta} L, 2f_L r^{-1 + \frac{1}{2}\delta})
 \end{cases}
 \end{equation*}
 with $f_{R} = 1 - (1+r)^{-\delta + C_{(A, \text{new})}}$ and $f_L = \chi_{(2r_0, R)}(r)$. Note that, for all these vector fields, $\slashed{\Pi}(Z)$ vanishes. 
 
 We can now attempt to repeat the estimates of lemma \ref{lemma boundedness}, but this time incorporating the change of basis. Specifically, we need to repeat the weighted $T$ energy estimate, the weighted Morawetz estimate, and the $p$-weighted estimate with $p = \frac{1}{2}\delta$.
 
 When performing the weighted $T$-energy estimate and the weighted Morawetz estimate, we encounter new error terms of the form
 
\begin{equation*}
 \begin{split}
  \int_{\mathcal{M}} w_{(A, new)}\Bigg(  
  &  
  (1+r)^{-1 - \delta} |\slashed{\D}\phi_{(\text{orig})}| |\slashed{\D}\phi|_{(A)}
  + \left( (1+r)^{-1} + (1+r)^{-1+\delta}(1+\tau)^{-\beta}\right) |\overline{\slashed{\D}}\phi_{(\text{orig})}| |\overline{\slashed{\D}}\phi|_{(A)} \\
  &
  + \delta (1+r)^{-2 - \delta} |\slashed{\D}\phi_{(\text{orig})}| |\phi|_{(A)}
  + \delta (1+r)^{-2 - \delta} |\slashed{\D}\phi|_{(A)} |\phi_{(\text{orig})} | \\
  &
  + \delta(1+r)^{-2+\delta}(1+\tau)^{-\beta} |\phi_{(\text{orig})}| |\slashed{\D}_L\phi|_{(A)}
  + \delta (1+r)^{-2+\delta}(1+\tau)^{-\beta} |\slashed{\D}_L\phi_{(\text{orig})}| |\phi|_{(A)} \\
  &
  + \delta^2 \left( (1+r)^{-3} + (1+r)^{-3+\delta}(1+\tau)^{-\beta}\right) |\phi_{(\text{orig})}| |\phi|_{(A)}
  \bigg) \dVol_{g}
\end{split}
\end{equation*} 

These terms can then be bounded by the following collection of terms, for any small constant $c$:
\begin{equation*}
 \begin{split}
 \int_{\mathcal{M}} \Bigg( 
  &
  c\delta (1+r)^{-1 - \delta} |\slashed{\D}\phi|^2_{(A)}
  + (c\delta)^{-1} (1+r)^{-1-\delta} |\slashed{\D}\phi_{(\text{orig})}|^2
  + c\delta (1+r)^{-1} |\overline{\slashed{\D}}\phi|^2_{(A)} \\
  &
  + (c\delta) (1+r)^{-1} |\overline{\slashed{\D}}\phi_{(\text{orig})}|^2
  + c\delta (1+r)^{-1 + \delta} (1+\tau)^{-\beta} |\overline{\slashed{\D}}\phi|^2_{(A)}\\
  &
  + (c\delta)^{-1} (1+r)^{-1 + \delta} (1+\tau)^{-\beta} |\overline{\slashed{\D}}\phi_{(\text{orig})}|^2
  + c\delta (1+r)^{-3-\delta} |\phi|_{(A)}^2
  + c^{-1}\delta (1+r)^{-3+\delta} |\phi_{(\text{orig})}|^2 \\
  &
  + c\delta (1+r)^{-1 + 3\delta} (1+\tau)^{-2\beta}|\slashed{\D}_L \phi|_{(A)}^2
  + c\delta^3 (1+r)^{-3 + \delta} (1+\tau)^{-\beta} |\phi|^2_{(A)} \\
  &
  + c^{-1}\delta (1+r)^{-3 + \delta} (1+\tau)^{-\beta} |\phi_{(\text{orig})}|^2
  \bigg) \dVol_g
 \end{split}
\end{equation*}

The terms involving $\phi_{(\text{orig})}$ and its derivatives have all appeared before in lemma \ref{lemma boundedness}, the only difference being that they include the large numerical factor $(c\delta)^{-1}$ where before there was either a coefficient of order $1$ or a small coefficient of order $\delta$. Hence, they cannot be absorbed by the left hand side of the estimates. However, the terms involving $\phi_{(\text{orig})}$ and its derivatives can be estimated \emph{without} the change of basis, in other words, for these fields we can simply appeal directly to lemma \ref{lemma boundedness}, and conclude the bounds of the relevant terms (which must then be multiplied by $\delta^{-1}$ or, at worse, $\delta^{-2}$).
  
On the other hand, the terms involving the new field $\phi_{(A)}$ come with the small coefficient $c\delta$ or, for those terms involving the field itself rather than its derivative, $c\delta^2$. These terms are otherwise of precisely the same form as terms which appeared as error terms in lemma \ref{lemma boundedness}. In fact, we can see from propositions \ref{proposition basic weighted T energy}, \ref{proposition basic Morawetz estimate} and \ref{proposition basic p weighted estimate} that these error terms can be absorbed the the bulk terms in the weighted $T$ energy, the Morawetz estimate or the $p$-weighted estimate with $p = \frac{1}{2}\delta$, as long as $c$ is sufficiently small.

Hence, we can try to control these terms in the same way as was done in lemma \ref{lemma boundedness}. An additional complication arises when we try to use the Hardy inequalities, where additional error terms appear; see propositions \ref{proposition Hardy after change of basis} and \ref{proposition Hardy on constant t after change of basis}. However, using the bounds assumed  on the derivatives of $M$, we find that proposition \ref{proposition Hardy after change of basis} yields

\begin{equation*}
 \begin{split}
  &\int_{^t\Sigma_\tau}(1+r)^{-1-\alpha}r^{-1} |\phi|_{(A)}^2 r^2 \upd r \wedge \dVol_{\mathbb{S}^2} \\
  &\lesssim \frac{1}{(1-\alpha)^2}
  \int_{^t\Sigma_\tau} (1+r)^{-\alpha} \left( |\slashed{\D}_L\phi|_{(A)}^2 + (1+r)^{-2-2\delta} |\phi_{(\text{orig})}|^2 \right) \, r^2\upd r \wedge \dVol_{\mathbb{S}^2} \\
  &\phantom{\leq \frac{1}{(1-\alpha)^2} \Bigg(}
  + \frac{1}{(1-\alpha)} \int_{\bar{S}_{\tau, t}}(1+r)^{-\alpha}|\phi|_{(A)}^2 r \dVol_{\mathbb{S}^2} \\
  &\lesssim \frac{1}{(1-\alpha)^2}
  \int_{^t\Sigma_\tau} (1+r)^{-\alpha} \left( |\slashed{\D}_L\phi|_{(A)}^2 + (1+r)^{-2\delta} |\slashed{\D}_L\phi_{(\text{orig})}|^2 \right) \, r^2\upd r \wedge \dVol_{\mathbb{S}^2} \\
  &\phantom{\leq }
  + \frac{1}{(1-\alpha)} \int_{\bar{S}_{\tau, t}}(1+r)^{-\alpha}|\phi|_{(A)}^2 r \dVol_{\mathbb{S}^2} 
  + \frac{1}{(1-\alpha)^2(1-2\delta-\alpha)} \int_{\bar{S}_{\tau, t}}(1+r)^{-2\delta -\alpha}|\phi|_{(\text{orig})}^2 r \dVol_{\mathbb{S}^2} 
  \\
 \end{split}
\end{equation*}
Here, the additional error terms involving $\phi_{(\text{orig})}$ have been dealt with by appealing again to the Hardy inequality of proposition \ref{proposition Hardy}.

Similarly, proposition \ref{proposition Hardy on constant t after change of basis} gives
\begin{equation*}
  \begin{split}
  &\int_{_{\tau_0}^{\tau_1}\Sigma_t}(1+r)^{-1-\alpha} |\phi|_{(A)}^2 r \upd r \wedge \dVol_{\mathbb{S}^2} \\ 
  &\lesssim \frac{1}{(1-\alpha)^2}
  \int_{_{\tau_0}^{\tau_1}\Sigma_t} (1+r)^{-\alpha} \bigg(|\slashed{\D}_L\phi|_{(A)}^2 + |\slashed{\D}_{\Lbar}\phi|_{(A)}^2 + (1+r)^{-2} |\phi_{(\text{orig})}|^2 \\
  &\phantom{ lesssim \frac{1}{(1-\alpha)^2} \int_{_{\tau_0}^{\tau_1}\Sigma_t} (1+r)^{-\alpha} \bigg(}
  + (1+r)^{-2 + 2\delta}(1+\tau)^{-4\beta}|\phi_{(\text{orig})}|^2 + \textit{Err}_{(t-\partial_r)}[\phi]_{(A)}\bigg) r^2\upd r \wedge \dVol_{\mathbb{S}^2} \\
  &\phantom{\leq \frac{1}{1-\alpha} }
  + \frac{1}{(1-\alpha)} \int_{\bar{S}_{\tau_0, t}}(1+r)^{-\alpha}|\phi|_{(A)}^2 r \dVol_{\mathbb{S}^2}
 \end{split}
\end{equation*}
Again, we would like to make use of the Hardy inequality \ref{proposition Hardy on constant t} to deal with the new error terms involving $\phi_{(\text{orig})}$. There is a problem, however, with the integral of the term
\begin{equation*}
 (1+r)^{-2 + 2\delta} (1+\tau)^{-2\beta} |\phi_{(\text{orig})}|^2
\end{equation*}
which does not have sufficient decay in $r$ to appeal to proposition \ref{proposition Hardy on constant t}. On the other hand, we have
\begin{equation*}
 (1+r)^{-2 + 2\delta} (1+\tau)^{-2\beta} |\phi_{(\text{orig})}|^2
 \lesssim (1+r)^{-2 + 2\delta - 2\beta} |\phi_{(\text{orig})}|^2
  + (1+r)^{-2\beta + 2\delta} (1+\tau)^{-2} |\phi_{(\text{orig})}|^2
\end{equation*}
where the first term provides the bound in the interior region $r \leq \tau$, and the second in the exterior region $\tau \leq r$. The first term can be bounded by making use of the Hardy inequality \ref{proposition Hardy on constant t}, and the second term by making use of the alternative Hardy inequality \ref{proposition Hardy on constant t with tau weight}, together with the bootstrap bounds.

In order to finish this proof, we also need to repeat the $p$-weighted estimate for $p = \frac{1}{2}\delta$, but including the error terms arising from derivatives of $M$. These can also be read off from proposition \ref{proposition energy identity after change of basis}, and are of the form
\begin{equation*}
 \begin{split}
  \int_{\mathcal{M}}\chi_{(2r_0, R)} r^p \bigg( &
  (1+r)^{-1-\delta} |\overline{\slashed{\D}}\phi_{(\text{orig})}| |\overline{\slashed{\D}}\phi|_{(A)}
  + (1+r)^{-1 + \delta} (1+\tau)^{-\beta} |\slashed{\D}_L \phi_{(\text{orig})}| |\slashed{\D}_L \phi|_{(A)} \\
  &
  + (1+r)^{-2 - \delta} |\slashed{\D}\phi_{(\text{orig})}| |\phi|_{(A)}
  + (1+r)^{-2 - \delta} |\slashed{\D}\phi|_{(A)} |\phi_{(\text{orig})}| \\
  &
  + (1+r)^{-2 + \delta}(1+\tau)^{-\beta} |\slashed{\D}_L \phi_{(\text{orig})}| |\phi|_{(A)}
  + (1+r)^{-2 + \delta}(1+\tau)^{-\beta} |\slashed{\D}_L \phi|_{(A)} |\phi_{(\text{orig})}| \\
  &
  + (1+r)^{-3 + \delta}(1+\tau)^{-\beta} |\phi_{(\text{orig})}| |\phi|_{(A)}
  \bigg) \dVol_g
 \end{split}
\end{equation*}
and they can be bounded by the terms
\begin{equation*}
\begin{split}
  \int_{\mathcal{M}}\chi_{(2r_0, R)} r^p \bigg( &
  c\delta (1+r)^{-1-\delta} |\overline{\slashed{\D}}\phi|_{(A)}^2
  + (c\delta)^{-1} (1+r)^{-1-\delta} |\overline{\slashed{\D}}\phi_{(\text{orig})}|^2
  + c\delta (1+r)^{-1 + \delta} (1+\tau)^{-\beta}|\slashed{\D}_L \phi|_{(A)}^2
  \\
  &
  + (c\delta)^{-1} (1+r)^{-1 + \delta} (1+\tau)^{-\beta} |\slashed{\D}_L \phi_{(\text{orig})}|
  + c\delta (1+r)^{-3 - \delta} |\phi|_{(A)}^2
  \\
  &
  + (c\delta)^{-1} (1+r)^{-1-\delta} |\slashed{\D} \phi_{(\text{orig})}|^2
  + c\delta (1+r)^{-1-\delta} |\slashed{\D} \phi|_{(A)}  
  + (c\delta)^{-1} (1+r)^{-3 - \delta} |\phi_{(\text{orig})}|^2 
  \\
  &
  + c\delta (1+r)^{-1 + 3\delta}(1+\tau)^{-2\beta} |\slashed{\D}_L \phi|_{(A)}^2
  + (c\delta)^{-1} (1+r)^{-1 + 3\delta}(1+\tau)^{-2\beta} |\slashed{\D}_L \phi_{(\text{orig})}|^2
  \\
  &
  + (c\delta)^{-1} (1+r)^{-3 + 3\delta}(1+\tau)^{-2\beta} |\phi_{(\text{orig})}|^2
  \bigg) \dVol_g
 \end{split}
\end{equation*}
Once again, the terms involving $\phi_{(A)}$ and its derivatives can be absorbed by the bulk terms in the Morawetz estimate and the $p$-weighted estimate, as long as $p \gg c\delta$. In particular, this holds with $p \geq \frac{1}{2}\delta$, as long as $c \ll \frac{1}{2}$. On the other hand, the terms involving $\phi_{(\text{orig})}$ have already been bounded. In particular, the terms involving $\phi_{(A)}$ (and \emph{not} its derivatives) can be bounded by proposition \ref{proposition higher weighted spatial integral} in terms of the $T$-energy and the $p$-weighted energy, for all sufficiently small $c$. Note that these terms come with the large factor $(c\delta)^{-1}$, whereas the quantities we can control come with the a factor $\delta$.

Finally we note that, when performing the $p$-weighted estimates, it is necessary to treat one of the error terms on the boundary as a total $L$ derivative, and this also requires some extra care, because $(\slashed{\D}_L \phi)_{(A)}$ is \emph{not} a perfect $L$ derivative, since the matrix $M$ appears outside of the derivative operator. In fact, the additional error terms can be read off from proposition \ref{proposition boundary p}:
\begin{equation*}
 \begin{split}
  &\int_{^t\Sigma_\tau} \chi_{(2r_0 R)} r^{p-1} |LM_{(A)}| |\phi|_{(A)} |\phi_{(\text{orig})}| \Omega^2 \upd r \wedge \dVol_{\mathbb{S}^2} \\
  &\lesssim \int_{^t\Sigma_\tau} \chi_{(2r_0 R)} r^{p-2-\delta}  \left( c\delta |\phi|^2_{(A)} + (c\delta)^{-1} |\phi_{(\text{orig})}|^2 \right) \Omega^2 \upd r \wedge \dVol_{\mathbb{S}^2} \\
 \end{split}
\end{equation*}
and this term is similar to other error terms already dealt with.

In order to prove the second part of the lemma, we also need to repeat the $p$-weighted energy estimates with larger values of $p$, however, once again we can use the estimates above, splitting the error terms into terms involving $\phi_{(A)}$ which are of exactly the same form as those dealt with in lemma \ref{lemma p weighted}, and error terms involving $\phi_{(\text{orig})}$ which are also of the same form as the error terms in lemma \ref{lemma p weighted} but with the large factor $(c\delta)^{-1}$ in place of the small factor $\delta$.

In summary, picking $c$ sufficiently small so that all of the conditions above hold, we have proved the lemma.

\end{proof}

\chapter{Pointwise bounds}
\label{chapter pointwise bounds}

In this chapter we will prove pointwise bounds on fields $\phi$ and its derivatives, where the associated degenerate energy is bounded. In fact, in order to use the Sobolev inequalities, we will need \emph{higher order} energies, involving not just the field $\phi$ but the fields $\mathcal{Z}^n \phi$, where we recall that $\mathcal{Z}$ either represents differentiating with respect to the vector field $T$ or applying the operator $r\slashed{\nabla}$. In fact, in the region $r \geq r_0$ we will rely on commuting with $r\slashed{\nabla}$, together with Sobolev inequalities on the sphere $S_{\tau_r}$, whereas in the region $r \leq r_0$ we rely on commuting with $\slashed{\D}_T$.

\section{Pointwise bounds in the region \texorpdfstring{$r \geq r_0$}{r geq r0}}

\subsection{Pointwise bounds on the field}

Using the proposition above, we have the following bound on the fields themselves:
\begin{proposition}[Pointwise bounds on the field in terms of its energy]
\label{proposition pointwise bound on field in terms of energy}
 Let $\phi$ be an $S_{\tau,r}$-tangent tensor field satisfying
\begin{equation*}
 \sum_{n=0}^{2} \mathcal{E}^{(w_{n}T)}[ (r\slashed{\nabla})^n\phi](\tau) \lesssim \tilde{\mathcal{E}}
\end{equation*}
where $w_n = (1+r)^{-C_{\left((r\slashed{\nabla}^n\phi \right)}\epsilon}$.

Additionally, suppose that the initial data on $\Sigma_{\tau_0}$ satisfies
\begin{equation*}
 \lim_{t \rightarrow \infty} \sum_{n=0}^2 \int_{\bar{S}_{\tau_0, t}} |(r\slashed{\nabla})^n \phi|^2 = 0
\end{equation*}
where $\tau_0 \leq \tau$, and suppose that all of the bootstrap bounds of chapter \ref{chapter bootstrap} hold.

Then, for all $r \geq r_0$ we have
\begin{equation}
 |\phi| \lesssim \sqrt{\tilde{\mathcal{E}}} (1+r)^{-\frac{1}{2} + \frac{1}{2}C^*\epsilon}
\end{equation}
where 
\begin{equation*}
C^* := \sup_{0 \leq n \leq 2} C_{\left[(r\slashed{\nabla}^n)\phi \right]}
\end{equation*}
\end{proposition}

\begin{proof}
 From proposition \ref{proposition spherical mean in terms of energy} we have 
\begin{equation*}
 \int_{S_{\tau,r}} \sum_{m=0}^{2} r^{2m} |\slashed{\nabla}^m \phi|^2 \dVol_{\mathbb{S}^2}
 \lesssim \sum_{m=0}^{2} \left( r^{-1 + C_{[\phi]}\epsilon} \mathcal{E}^{(w_m T)}(\tau, t, \tau_0)[(r\slashed{\nabla})^m \phi]
  + \int_{\bar{S}_{\tau_0, t}}|(r\slashed{\nabla})^m\phi|^2 \dVol_{\mathbb{S}^2} \right)
\end{equation*}
Taking the limit $t \rightarrow \infty$, we find
\begin{equation*}
 \int_{S_{\tau,r}} \sum_{m=0}^{2} r^{2m} |\slashed{\nabla}^m \phi|^2 \dVol_{\mathbb{S}^2}
 \lesssim \sum_{m=0}^{2} r^{-1 + C_{[\phi]}\epsilon} \mathcal{E}^{(wT)}(\tau)[(r\slashed{\nabla})^m\phi] \lesssim r^{-1 + C_{(\phi)}\epsilon} \tilde{\mathcal{E}}
\end{equation*}
Now, appealing to proposition \ref{proposition Sobolev}, recalling that $C^*$ is at least as large the other constants and that $\dVol_{\slashed{g}} \sim r^2 \dVol_{\mathbb{S}^2}$, we prove the proposition.
\end{proof}

Combining this with corollary \ref{corollay energy decay} we obtain the following result:

\begin{proposition}[Pointwise decay of fields satisfying the wave equation]
\label{proposition pointwise decay in r and tau}
 Let the conditions of lemma \ref{lemma p weighted} hold. In particular, this means that $\phi$ satisfies the wave equation
\begin{equation*}
 \tilde{\slashed{\Box}}_g \phi = F
\end{equation*}
where $F$, $g$ and the initial data for $\phi$ satisfy the conditions given in lemma \ref{lemma p weighted}. Additionally, let the same conditions be satisfied by the fields $r\slashed{\nabla}\phi$ and $(r\slashed{\nabla})^2\phi$, not necessarily with the same inhomogenous terms $F$. In other words, let
\begin{equation*}
 \tilde{\slashed{\Box}}_g \left( (r\slashed{\nabla})^m \phi \right) = F_{(m)}
\end{equation*}
for $m = 0,1,2$, and where the inhomogenous terms $F_{(m)}$ satisfy the conditions satisfied by $F$. Moreover, the constants $C_{[\phi]}$ can depend on the integer $m$, i.e.\ we have the constants  $C_{[\phi]}$ $C_{[r\slashed{\nabla}\phi]}$ and $C_{\left[(r\slashed{\nabla})^2 \phi\right]}$. Additionally, we set
\begin{equation*}
 C^* = \max_{0 \leq n \leq 2} \left\{ C_{\left[ (r\slashed{\nabla})^n \phi \right]} \right\}
\end{equation*}

Then we have
\begin{equation}
 |\phi| \lesssim \sqrt{\tilde{\mathcal{E}}} (1+\tau)^{-\frac{1}{2} + \frac{1}{2}C^*\delta} (1+r)^{-\frac{1}{2} + \frac{1}{2}C^*\epsilon}
\end{equation}
\end{proposition}

\begin{proof}
 We first set
\begin{equation*}
 w_{(m)} = (1+r)^{-C_{\left[ (r\slashed{\nabla})^m \phi \right]}}
\end{equation*}

 Applying corollary \ref{corollay energy decay} we find
\begin{equation*}
 \mathcal{E}^{(w_{(m)}T)}[(r\slashed{\nabla})^m \phi](\tau) \lesssim \tilde{\mathcal{E}}(1+\tau)^{-1+C_{\left[ (r\slashed{\nabla})^n \phi \right]}\delta}
\end{equation*}
for $m = 0, 1, 2$. Now, applying proposition \ref{proposition pointwise bound on field in terms of energy} we find that
\begin{equation*}
 |\phi| \lesssim \sqrt{\tilde{\mathcal{E}}} (1+\tau)^{-\frac{1}{2} + \frac{1}{2}C^*} (1+r)^{-\frac{1}{2} + \frac{1}{2}C^*\epsilon}
\end{equation*}

\end{proof}

We also have an additional way to bound the field $\phi$, which makes use of the boundedness of the $p$-weighted energy to obtain improved decay in $r$ relative to the proposition above, at the cost of slower decay in $\tau$.

\begin{proposition}[Improved pointwise bounds on the field in terms of the $p$-weighed energy]
	\label{proposition improved pointwise bounds in terms of p energy}
	Let $\phi$ be an $S_{\tau,r}$-tangent tensor field, and let $r \geq r_0$. Let $w_n(r) = (1+r)^{-C_{(n)}\epsilon}$ for any $C_{(n)} \geq 0$. Let $p < 1$. Then we have
	\begin{equation}
	|\phi|^2 \lesssim \sum_{n=0}^2 \left( 
	r^{-2} \mathcal{E}^{(w_n T)}[(r\slashed{\nabla})^n \phi](\tau)
	+ \frac{1}{1-p} r^{-1-p} \mathcal{E}^{(L, p)}[(r\slashed{\nabla})^n \phi](\tau)
	\right)
	\end{equation}
\end{proposition}

\begin{proof}
	Using proposition \ref{proposition higher weighted spherical integral} we have
	\begin{equation*}
	\int_{S_{\tau,r}} |\phi|^2 \dVol_{\mathbb{S}^2}
	\lesssim
	\left( \frac{r_0}{r} \right)^2 \int_{S_{\tau,r_0}} |\phi|^2 \dVol_{\mathbb{S}^2}
	+ \frac{1}{1-p} r^{-1-p} \mathcal{E}^{(L, p)}[\phi][\tau]
	\end{equation*}
	We can control the boundary term at $r = r_0$ by using proposition \ref{proposition spherical mean in terms of energy}, which allows us to obtain
	\begin{equation*}
	\int_{S_{\tau,r}} |\phi|^2 \dVol_{\mathbb{S}^2}
	\lesssim
	r^{-2} \mathcal{E}^{(wT)}[\phi](\tau)
	+ \frac{1}{1-p} r^{-1-p} \mathcal{E}^{(L, p)}[\phi](\tau)
	\end{equation*}
	Commuting twice with $r\slashed{\nabla}$, using $\Omega \sim r$ and the Sobolev inequality on the sphere proves the proposition.
\end{proof}

\subsection{Pointwise bounds on derivatives}
\label{subsection pointwise bounds on derivatives}
Proposition \ref{proposition pointwise decay in r and tau} establishes decay in both $\tau$ and $r$ for solutions to the wave equation. In particular, the field $\phi$ itself, which is the field satisfying the wave equation, is shown to decay. Although this is important to close our bootstrap, we have also made use of bootstrap assumptions on the decay of \emph{derivatives} of fields satisfying the wave equation. Indeed, in order to close the bootstrap arguments we must show that the derivatives of the field decay faster (in $r$) than the field itself, and the ``good'' derivatives must be shown to decay faster still.

We first need to prove a few preliminary propositions:

\begin{proposition}[An estimate for spatial integrals of bad derivatives]
\label{proposition estimate for spatial integrals of bad derivs}
Let $\phi$ be an $S_{\tau,r}$-tangent tensor field, and assume that the pointwise bootstrap bounds of chapter \ref{chapter bootstrap} hold. Then
\begin{equation}
 \int_{^t\Sigma_\tau} (1+r)^{-1-\delta} |\slashed{\D}\phi|^2 \upd r \wedge \dVol_{\mathbb{S}^2}
 \lesssim \int_{^t\mathcal{M}_\tau^{\tau + 1}} \left( \sum_{n = 0}^1 (1+r)^{-1-\delta}|\slashed{\D} \mathscr{Z}^n\phi|^2 + \epsilon (1+r)^{-3-\delta}|\phi|^2 \right) \upd r \wedge \dVol_{\mathbb{S}^2}
\end{equation}
\end{proposition}
\begin{proof}
 We have
 \begin{equation*}
   \begin{split}
     &\int_{^t\Sigma_{\tau_1}} (1+r)^{-1-\delta} |\slashed{\D}\phi|^2 r^2 \upd r \wedge \dVol_{\mathbb{S}^2}
     \lesssim
     \int_{\tau_1}^{\tau_1 + 1} -\partial_{\tau} \left( \int_{^t\Sigma_\tau} (1+r)^{-1-\delta} |\slashed{\D}\phi|^2 \Omega^2 \upd r \wedge \dVol_{\mathbb{S}^2} \right) \upd \tau \\
     &\lesssim
     \int_{^t\mathcal{M}^{\tau_1 + 1}_{\tau_1}} (1+r)^{-1-\delta}\left( \left. \frac{\partial}{\partial \tau} \right|_{r, \vartheta^1, \vartheta^2} |\slashed{\D}\phi|^2 \right) \mu^{-1} \dVol_g \\
     &\lesssim
     \int_{^t\mathcal{M}^{\tau_1 + 1}_{\tau_1}} (1+r)^{-1-\delta} \left| ( L + \Lbar - b^A X_A) |\slashed{\D}\phi|^2 \right|  \dVol_g \\
     &\lesssim
     \int_{^t\mathcal{M}^{\tau_1 + 1}_{\tau_1}} (1+r)^{-1-\delta} \left( | \slashed{\D}_T\slashed{\D}\phi| + r^{-1} |b| |r\slashed{\nabla} \slashed{\D} \phi| \right) |\slashed{\D}\phi| \dVol_g \\
     &\lesssim
     \int_{^t\mathcal{M}^{\tau_1 + 1}_{\tau_1}} (1+r)^{-1-\delta} \bigg( 
       |\slashed{\D} \slashed{\D}_T \phi| |\slashed{\D}\phi| 
       + r^{-1}|b| |\slashed{\D} (r\slashed{\nabla})\phi| |\slashed{\D}\phi|
       + (1+|b|)|\bm{\Gamma}| |\slashed{\D}\phi|^2
       + |\Omega_{L\Lbar}| |\phi| |\slashed{\D}\phi|
       \\
       &\phantom{\lesssim \int_{^t\mathcal{M}^{\tau_1 + 1}_{\tau_1}} (1+r)^{-1-\delta} \bigg(}
       + (1+|b|)\left( |\slashed{\Omega}_{L}| + |\slashed{\Omega}_{\Lbar}| \right) |\phi| |\slashed{\D}\phi|
       \bigg) \dVol_g \\
     &\lesssim
     \int_{^t\mathcal{M}^{\tau_1 + 1}_{\tau_1}} (1+r)^{-1-\delta}\bigg( 
       |\slashed{\D} \mathscr{Z}\phi| 
       + \epsilon (1+r)^{-1 + 2\delta}|\slashed{\D}\phi|
       + \epsilon (1+r)^{-2 + 2\delta}|\phi| \bigg) |\slashed{\D}\phi| \dVol_g \\
     &\lesssim
     \int_{^t\mathcal{M}^{\tau_1 + 1}_{\tau_1}} (1+r)^{-1-\delta}\left( |\slashed{\D} \mathscr{Z}\phi|^2 + |\slashed{\D}\phi|^2 + \epsilon(1+r)^{-3-\delta}|\phi|^2 \right)  \dVol_g \\
   \end{split}
 \end{equation*}
 where we have used the fact that $\dVol_g = -\mu \Omega^2 \upd \tau \wedge \upd r \wedge \dVol_{\mathbb{S}^2}$, and that
 \begin{equation*}
   \left. \frac{\partial}{\partial \tau} \right|_{r, \vartheta^1, \vartheta^2} = \frac{1}{2}\mu \left( L + \Lbar - b^A X_A \right)
 \end{equation*}
 and also the calculations in propositions \ref{proposition commute T} and \ref{proposition commuting rnabla with first order operators}.
\end{proof}

\begin{corollary}[Decay for the spherical mean of the bad derivatives along a subsequence]
\label{corollary decay for spherical mean along a subsequence}
Let $\phi$ be an $S_{\tau,r}$-tangent tensor field, and assume that the pointwise bootstrap bounds of chapter \ref{chapter bootstrap} hold. Suppose also that the following integral over the entire spacetime region $\Sigma_{\tau}$ (not just the ``cut-off'' spacetime region $^t\Sigma_\tau$) is finite:
\begin{equation*}
  \int_{\Sigma_\tau} \left( \sum_{n = 0}^1 \delta^2 (1+r)^{-1-\delta} |\slashed{\D} \mathscr{Z}^n\phi|^2 + \delta^2 (1+r)^{-3-\delta}|\phi|^2 \right) \upd r \wedge \dVol_{\mathbb{S}^2}
  \leq \mathcal{E}
\end{equation*}
Then there is a diadic sequence of radii $r_n$, with $r_n \rightarrow \infty$ as $n \rightarrow \infty$ such that
\begin{equation*}
  \int_{S_{\tau,r_n}} \left((1+r)^{1 - \delta} |\slashed{\D}\phi|^2 + (1+r)^{-3-\delta}|\phi|^2 \right) \dVol_{\mathbb{S}^2} \leq (r_n)^{-1} \delta^{-2} \mathcal{E}
\end{equation*}

\end{corollary}

\begin{proposition}[An expression for $\slashed{\D}_L \left(\mathcal{Z}^n \slashed{\D}_{\Lbar} \phi\right)$]
\label{proposition expression for LLbar phi}
Let $\phi$ be a scalar field satisfying
\begin{equation*}
 \tilde{\slashed{\Box}}_g \phi = F
\end{equation*}
for some scalar $F$, and suppose that all of the bootstrap bounds of chapter \ref{chapter bootstrap} hold.

Then for all $r \geq r_0$ we have

\begin{equation*}
\begin{split}
  |\slashed{\D}_L \left(\mathscr{Z}^n \Lbar \phi\right)|
  &\lesssim \epsilon (1+r)^{-1} (\slashed{\D} \mathscr{Z}^n \phi)
  + r^{-1} |\overline{\slashed{\D}} \mathscr{Z}^{n+1}\phi|
  + \sum_{j+k \leq n-1} |\bm{\Gamma}^{(j)}_{(-1+C_{(j)}\epsilon)}| |\slashed{\D}\mathscr{Z}^k \phi|
  \\
  &\phantom{=}
  + \sum_{j+k \leq n} |\bm{\Gamma}^{(j)}_{(-1+C_{(j)}\epsilon)}| |\overline{\slashed{\D}}\mathscr{Z}^k \phi|
  + \sum_{j+k \leq n} |\bm{\Gamma}^{(j)}_{(-1-\frac{1}{2}\delta)}| |\slashed{\D}\mathscr{Z}^k \phi|
  + \sum_{j+k \leq n} |\bm{\Gamma}^{(j)}_{(-2-\delta)}| |\mathscr{Z}^k \phi|
 \end{split}
\end{equation*}

Similarly, if $\phi$ is a higher valence $S_{\tau,r}$-tangent field, then
\begin{equation*}
\begin{split}
|\slashed{\D}_L \left(\mathscr{Z}^n \slashed{\D}_{\Lbar} \phi\right)|
&\lesssim \epsilon (1+r)^{-1} (\slashed{\D} \mathscr{Z}^n \phi)
+ r^{-1} |\overline{\slashed{\D}} \mathscr{Z}^{n+1}\phi|
+ \sum_{j+k \leq n-1} |\bm{\Gamma}^{(j)}_{(-1+C_{(j)}\epsilon)}| |\slashed{\D}\mathscr{Z}^k \phi|
\\
&\phantom{=}
+ \sum_{j+k \leq n} |\bm{\Gamma}^{(j)}_{(-1+C_{(j)}\epsilon)}| |\overline{\slashed{\D}}\mathscr{Z}^k \phi|
+ \sum_{j+k \leq n} |\bm{\Gamma}^{(j)}_{(-1-\frac{1}{2}\delta)}| |\slashed{\D}\mathscr{Z}^k \phi|
\\
&\phantom{=}
+ \sum_{j+k \leq n} |\bm{\Gamma}^{(j+1)}_{(-2 + C_{(j+1)}\epsilon)}| |\mathscr{Z}^k \phi|
\end{split}
\end{equation*}

\end{proposition}

\begin{proof}
 We begin by noting that, if $\phi$ is a scalar field, then
\begin{equation}
\label{equation LLbar phi scalar case}
 L\Lbar \phi = -\tilde{\Box}_g \phi + \slashed{\Delta}\phi - \frac{1}{2}\tr_{\slashed{g}}\chi (\Lbar \phi) - \frac{1}{2}\tr_{\slashed{g}}\chibar (L \phi) - \zeta^\alpha \slashed{\nabla}_\alpha \phi
\end{equation}
which follows from proposition \ref{proposition scalar wave operator} together with the definition $\tilde{\Box}_g \phi = \Box_g \phi + \omega \Lbar \phi$.

On the other hand, if $\phi$ is a higher rank $S_{\tau,r}$-tangent tensor field, then we have, schematically,
\begin{equation}
\label{equation LLbar phi tensor case}
 \slashed{\D}_L \slashed{\D}_{\Lbar} \phi = -\tilde{\slashed{\Box}}_g \phi + \slashed{\Delta}\phi - \frac{1}{2}\tr_{\slashed{g}}\chi (\Lbar \phi) - \frac{1}{2}\tr_{\slashed{g}}\chibar (L \phi) - \zeta^\alpha \slashed{\nabla}_\alpha \phi + \Omega_{L\Lbar} \cdot \phi
\end{equation}
which follows from proposition \ref{proposition tensor wave operator}.

Hence, if $\phi$ is a scalar field then we have
\begin{equation*}
  |L\Lbar \phi| \lesssim
  |F|
  + r^{-1}| \slashed{\nabla} (r\slashed{\nabla}\phi) |
  + \left( r^{-1} + |\tr_{\slashed{g}}\chi_{(\text{small})}| \right)| \partial \phi| 
  + \left( r^{-1} + |\zeta| + |\tr_{\slashed{g}}\chibar_{(\text{small})}| \right)| \bar{\partial}\phi|
\end{equation*}
while if $\phi$ is higher rank, then
\begin{equation*}
 \begin{split}
  \left| \slashed{\D}_L \left( \slashed{\D}_{\Lbar} \phi \right) \right|
  &\lesssim
  |F|
  + r^{-1}| \slashed{\nabla} (r\slashed{\nabla}\phi) |
  + \left( r^{-1} + |\tr_{\slashed{g}}\chi_{(\text{small})}| \right)| \partial \phi| 
  + \left( r^{-1} + |\zeta| + |\tr_{\slashed{g}}\chibar_{(\text{small})}| \right)| \bar{\partial}\phi| \\
  &\phantom{\lesssim}
  + |\Omega_{L\Lbar}| |\phi|
 \end{split} 
\end{equation*}

To prove the proposition we need to commute equations \eqref{equation LLbar phi scalar case} and \eqref{equation LLbar phi tensor case} with the operators $\mathscr{Z}^n$. Using propositions \ref{proposition commute Zn L} and \ref{proposition commute Zn D} we find, for a scalar field $\phi$, schematically
\begin{equation*}
\begin{split}
	\slashed{\D}_L \left(\mathscr{Z}^n \Lbar \phi\right)
	&=
	\bm{\Gamma}^{(0)}_{(-1)}(\mathscr{Z}^n \Lbar \phi)
	+ \sum_{j+k \leq n-1} \bm{\Gamma}^{(j)}_{(-1+C_{(j)}\epsilon)} (\overline{\slashed{\D}}\mathscr{Z}^k \Lbar\phi)
	+ \sum_{j+k \leq n-1} \bm{\Gamma}^{(j+1)}_{(-1-\delta)} (\mathscr{Z}^k \Lbar \phi)
	+ (\mathscr{Z}^n F)
	\\
	&\phantom{=}
	+ r^{-1} (\mathscr{Z}^n \overline{\slashed{\D}} \mathscr{Z} \phi)
	+ \sum_{j+k \leq n} \bm{\Gamma}^{(j)}_{(-1-\delta)} (\mathscr{Z}^k \Lbar \phi)
	+ \sum_{j+k \leq n} \bm{\Gamma}^{(j)}_{(-1+C_{j}\epsilon)} (\mathscr{Z}^k L \phi) \\
	&\phantom{=}
	+ \sum_{j+k \leq n} \bm{\Gamma}^{(j)}_{(-1+C_{j}\epsilon)} (\mathscr{Z}^k \slashed{\nabla} \phi)
	\\
	\\
	&=
	\bm{\Gamma}^{(0)}_{(-1)} (\slashed{\D} \mathscr{Z}^n \phi)
	+ \sum_{j+k \leq n-1} \bm{\Gamma}^{(j)}_{(-1+C_{(j)}\epsilon)} (\slashed{\D} \mathscr{Z}^k \phi)
	+ \sum_{j+k \leq n-1} \bm{\Gamma}^{(j)}_{(-2+C_{(j+1)}\epsilon)} (\mathscr{Z}^k \phi)
	\\
	&\phantom{=}
	+ \sum_{j+k \leq n-1} \bm{\Gamma}^{(j)}_{(-2+2C_{(j)}\epsilon)} (\slashed{\D} \mathscr{Z}^{k+1}\phi)
	+ \sum_{j+k \leq n-1} \bm{\Gamma}^{(j)}_{(-3+2C_{(j)}\epsilon)} (\mathscr{Z}^{k+1}\phi)
	\\
	&\phantom{=}
	+ \sum_{j+k \leq n-1} \bm{\Gamma}^{(j+1)}_{(-1-\frac{1}{2}\delta)} (\slashed{\D} \mathscr{Z}^k \phi)
	+ \sum_{j+k \leq n-1} \bm{\Gamma}^{(j+1)}_{(-2-\frac{1}{2}\delta)} (\mathscr{Z}^k \phi)
	+ (\mathscr{Z}^n F)
	+ r^{-1} (\overline{\slashed{\D}} \mathscr{Z}^{n+1}\phi)
	\\
	&\phantom{=}
	+ \sum_{j+k \leq n} \bm{\Gamma}^{(j)}_{(-1+C_{j}\epsilon)} (\overline{\slashed{\D}}\mathscr{Z}^k \phi)
	+ \sum_{j+k \leq n} \bm{\Gamma}^{(j)}_{(-2-\frac{1}{2}\delta)} (\slashed{\D}\mathscr{Z}^k \phi)
	\\
	\\
	&= \bm{\Gamma}^{(0)}_{(-1)} (\slashed{\D} \mathscr{Z}^n \phi)
	+ r^{-1} (\overline{\slashed{\D}} \mathscr{Z}^{n+1}\phi)
	+ \sum_{j+k \leq n-1} \bm{\Gamma}^{(j)}_{(-1+C_{(j)}\epsilon)} (\slashed{\D}\mathscr{Z}^k \phi)
	+ (\mathscr{Z}^n F)
	\\
	&\phantom{=}
	+ \sum_{j+k \leq n} \bm{\Gamma}^{(j)}_{(-1+C_{(j)}\epsilon)} (\overline{\slashed{\D}}\mathscr{Z}^k \phi)
	+ \sum_{j+k \leq n} \bm{\Gamma}^{(j)}_{(-1-\frac{1}{2}\delta)} (\slashed{\D}\mathscr{Z}^k \phi)
	+ \sum_{j+k \leq n} \bm{\Gamma}^{(j)}_{(-2-\frac{1}{2}\delta)} (\mathscr{Z}^k \phi)
\end{split}
\end{equation*}

On the other hand, if $\phi$ is a higher rank field, then we have additional lower order terms, and we find
\begin{equation*}
\begin{split}
	\slashed{\D}_L \left(\mathscr{Z}^n \Lbar \phi\right)
	&= \bm{\Gamma}^{(0)}_{(-1)} (\slashed{\D} \mathscr{Z}^n \phi)
	+ r^{-1} (\overline{\slashed{\D}} \mathscr{Z}^{n+1}\phi)
	+ \sum_{j+k \leq n-1} \bm{\Gamma}^{(j)}_{(-1+C_{(j)}\epsilon)} (\slashed{\D}\mathscr{Z}^k \phi)
	+ (\mathscr{Z}^n F)
	\\
	&\phantom{=}
	+ \sum_{j+k \leq n} \bm{\Gamma}^{(j)}_{(-1+C_{(j)}\epsilon)} (\overline{\slashed{\D}}\mathscr{Z}^k \phi)
	+ \sum_{j+k \leq n} \bm{\Gamma}^{(j)}_{(-1-\frac{1}{2}\delta)} (\slashed{\D}\mathscr{Z}^k \phi)
	+ \sum_{j+k \leq n} \bm{\Gamma}^{(j+1)}_{(-2 + C_{(j+1)}\epsilon)} (\mathscr{Z}^k \phi)
\end{split}
\end{equation*}

\end{proof}

\begin{proposition}[A pointwise estimate for the bad derivatives]
	\label{proposition pointwise bound lbar}
 Let $\phi$ be an $S_{\tau,r}$-tangent tensor field satisfying
\begin{equation*}
 \tilde{\slashed{\Box}}_g \phi = F
\end{equation*}
for some $S_{\tau,r}$-tangent tensor field $F$ of the same rank as $\phi$, and suppose that all of the bootstrap bounds of chapter \ref{chapter bootstrap} hold. Suppose also that
\begin{equation}
  \begin{split}
    &\int_{\Sigma_{\tau}\cap\{r \geq r'\}}
    \bigg(
    r^{1 - 2\delta} \sum_{0 \leq n \leq 2} |(r\slashed{\nabla})^n F|^2
    + r^{-1 - \delta} \sum_{0 \leq n \leq 3} |\overline{\slashed{\D}} (r\slashed{\nabla})^n \phi|^2 
    \bigg) r^2 \upd r \wedge \dVol_{\mathbb{S}^2}  \\
    &  
    + \int_{\mathcal{M}_\tau^{\infty} \cap \{r' \leq r\}}
    r^{-1-\delta} \sum_{n = 0}^3 |\slashed{\D} \mathscr{Z}^n \phi|^2 \dVol_g 
    \leq \mathcal{E}(1+\tau)^{-2\kappa}
  \end{split}
\end{equation}
for some $\kappa \geq 0$.

Then, for $r \geq r_0$ we have
\begin{equation}
  |\slashed{\D}_{\Lbar} \phi| \lesssim \sqrt{\mathcal{E}} r^{-1 + \delta} (1+\tau)^{-\kappa^*}
\end{equation}
where
\begin{equation*}
\kappa^* := \begin{cases}
\kappa \quad &\text{if } \kappa \leq 1 \\
1 \quad &\text{if } \kappa > 1
\end{cases}
\end{equation*}

\end{proposition}

\begin{proof}
	
	We apply (a very slight modification\footnote{Specifically, we integrate up to some maximum value of $r$ rather than a maximum value of $t$.} of) proposition \ref{proposition spherical mean in terms of energy} applied to the field $\left((r\slashed{\nabla})^m \slashed{\D}_{\Lbar} \phi\right)$, with the choice $\alpha = -1 + 2\delta$, and integrating up to some maximum value $r = R$. Summing the resulting inequalities for $m = 0$, $1$, $2$ we find
	\begin{equation*}
	\begin{split}
	&\sum_{m = 0}^2 \int_{S_{\tau, r'}} |(r\slashed{\nabla})^m \slashed{\D}_{\Lbar}\phi|^2 \dVol_{\mathbb{S}^2}
	\\
	&\lesssim
	(r')^{-2 + 2\delta} \int_{\Sigma_{\tau}\cap\{r' \leq r \leq R\}}
	r^{1 - 2\delta}
	\bigg(
	\epsilon^2 (1+r)^{-2} |\slashed{\D} \mathscr{Z}^{\leq 2} \phi|^2
	+ r^{-2} |\overline{\slashed{\D}} \mathscr{Z}^{\leq 3}\phi|^2
	+ \sum_{j+k \leq 1} |\bm{\Gamma}^{(j)}_{(-1+C_{(j)}\epsilon)}|^2 |\slashed{\D}\mathscr{Z}^k \phi|^2
	\\
	&\phantom{\lesssim
		(r')^{-2 + 2\delta} \int_{\Sigma_{\tau}\cap\{r' \leq r \leq R\}}
		r^{1 - 2\delta}
		\bigg(}
	+ \sum_{j+k \leq 2} |\bm{\Gamma}^{(j)}_{(-1+C_{(j)}\epsilon)}|^2 |\overline{\slashed{\D}}\mathscr{Z}^k \phi|^2
	+ \sum_{j+k \leq 2} |\bm{\Gamma}^{(j)}_{(-1-\frac{1}{2}\delta)}|^2 |\slashed{\D}\mathscr{Z}^k \phi|^2
	\\
	&\phantom{\lesssim
		(r')^{-2 + 2\delta} \int_{\Sigma_{\tau}\cap\{r' \leq r \leq R\}}
		r^{1 - 2\delta}
		\bigg(}
	+ \sum_{j+k \leq 2} |\bm{\Gamma}^{(j+1)}_{(-2 + C_{(j+1)}\epsilon)}|^2 |\mathscr{Z}^k \phi|^2
	+ \sum_{n \leq 2} |\mathscr{Z}^n F|^2
	\bigg) r^2\upd r \wedge \dVol_{\mathbb{S}^2}
	\\
	&\phantom{\lesssim}
	+ \int_{S_{\tau, R}} |(r\slashed{\nabla})^{\leq 2} \slashed{\D}_{\Lbar} \phi |^2 \dVol_{\mathbb{S}^2}
	\\
	\\
	&\lesssim (r')^{-2+2\delta} \int_{\Sigma_\tau \cap \{r' \leq r \leq R\}} \bigg(
	\sum_{n=0}^2 \epsilon^2 (1+r)^{-1-2\delta + 2C_{(2)}\epsilon}|\slashed{\D} \mathscr{Z}^{n} \phi|^2
	+ \sum_{n=0}^3 (1+r)^{-1-2\delta} |\overline{\slashed{\D}} \mathscr{Z}^n \phi|^2
	\\
	&\phantom{\lesssim (r')^{-2+2\delta} \int_{\Sigma_\tau \cap \{r' \leq r \leq R\}} \bigg(}
	+ \sum_{n=0}^2 \epsilon^2 (1+r)^{-3-2\delta + 2C_{(3)}\epsilon} |\mathscr{Z}^n \phi|^2
	+ \sum_{n=0}^2 (1+r)^{1-2\delta}|\mathscr{Z}^n F|^2
	\bigg) r^2 \upd r \wedge \dVol_{\mathbb{S}^2}
	\\
	&\phantom{\lesssim}
	+ \int_{S_{\tau, R}} |(r\slashed{\nabla})^{\leq 2} \slashed{\D}_{\Lbar} \phi |^2 \dVol_{\mathbb{S}^2}
	\\
	\\
	&\lesssim (r')^{-2+2\delta} \int_{\Sigma_\tau \cap \{r' \leq r \leq R\}} \bigg(
	\sum_{n=0}^2 \epsilon^2 (1+r)^{-1-\delta}|\slashed{\D} \mathscr{Z}^{n} \phi|^2
	+ \sum_{n=0}^3 (1+r)^{-1-\delta} |\overline{\slashed{\D}} \mathscr{Z}^n \phi|^2
	\\
	&\phantom{\lesssim (r')^{-2+2\delta} \int_{\Sigma_\tau \cap \{r' \leq r \leq R\}} \bigg(}
	+ \sum_{n=0}^2 \epsilon^2 (1+r)^{-3-\delta} |\mathscr{Z}^n \phi|^2
	+ \sum_{n=0}^2 (1+r)^{1-2\delta}|\mathscr{Z}^n F|^2
	\bigg) r^2 \upd r \wedge \dVol_{\mathbb{S}^2}
	\\
	&\phantom{=}
	+ \int_{S_{\tau, R}} \sum_{n=0}^2 \left( \epsilon^2 (1+r)^{2C_{(1)}\epsilon} |\slashed{\D} \mathscr{Z}^n \phi|^2 + \epsilon^2 (1+r)^{-2 + C_{(2)}\epsilon} |\mathscr{Z}^n \phi |^2 \right) \dVol_{\mathbb{S}^2}
	\end{split}
	\end{equation*}
	
	We need some additional control on the terms
	\begin{equation*}
	\int_{\Sigma_\tau \cap \{r' \leq r \leq R \}} \bigg(
	\sum_{n=0}^2 \epsilon^2 (1+r)^{-1-\delta}|\slashed{\D} \mathscr{Z}^{n} \phi|^2
	+ \sum_{n=0}^2 \epsilon^2 (1+r)^{-3-\delta} |\mathscr{Z}^n \phi|^2 \bigg) r^2\upd r \wedge \dVol_{\mathbb{S}^2}
	\end{equation*}
	as well as the terms given as an integral over the sphere $S_{\tau,R}$.

	First, we note that, for $\tau' > \tau$ we have
	\begin{equation*}
	\begin{split}
	&\int_{\Sigma_\tau \cap \{r \geq r'\}} \bigg(
	\sum_{n=0}^2 \epsilon^2 (1+r)^{-1-\delta}|\slashed{\D} \mathscr{Z}^{n} \phi|^2
	+ \sum_{n=0}^2 \epsilon^2 (1+r)^{-3-\delta} |\mathscr{Z}^n \phi|^2 \bigg) r^2\upd r \wedge \dVol_{\mathbb{S}^2}
	\\
	&=\int_{\Sigma_{\tau'} \cap \{r \geq r'\}} \bigg(
	\sum_{n=0}^2 \epsilon^2 (1+r)^{-1-\delta}|\slashed{\D} \mathscr{Z}^{n} \phi|^2
	+ \sum_{n=0}^2 \epsilon^2 (1+r)^{-3-\delta} |\mathscr{Z}^n \phi|^2 \bigg) r^2\upd r \wedge \dVol_{\mathbb{S}^2}
	\\
	&\phantom{=}
	+ \int_{\tau}^{\tau'} \frac{\partial}{\partial \tilde{\tau}}\Bigg( \int_{\Sigma_{\tilde{\tau}}\cap \{r \geq r'\}}
	\sum_{n=0}^2 \epsilon^2 (1+r)^{-1-\delta}|\slashed{\D} \mathscr{Z}^{n} \phi|^2
	+ \sum_{n=0}^2 \epsilon^2 (1+r)^{-3-\delta} |\mathscr{Z}^n \phi|^2 \bigg) r^2\upd r \wedge \dVol_{\mathbb{S}^2} \Bigg) \upd \tilde{\tau}
	\\
	\\
	&\lesssim \int_{\Sigma_{\tau'}} \bigg(
	\sum_{n=0}^2 \epsilon^2 (1+r)^{-1-\delta}|\slashed{\D} \mathscr{Z}^{n} \phi|^2
	+ \sum_{n=0}^2 \epsilon^2 (1+r)^{-3-\delta} |\mathscr{Z}^n \phi|^2 \bigg) r^2\upd r \wedge \dVol_{\mathbb{S}^2}
	\\
	&\phantom{\lesssim}
	+ \int_{\mathcal{M}_{\tau}^{\tau'}} \bigg(
	\sum_{n=0}^2 \epsilon^2 (1+r)^{-1-\delta}|\mathscr{Z}\slashed{\D} \mathscr{Z}^{n} \phi| |\slashed{\D}\mathscr{Z}^n\phi| 
	+ \sum_{n=0}^2 \epsilon^2 (1+r)^{-3-\delta} |\mathscr{Z}^{n+1} \phi||\mathscr{Z}^n \phi| \bigg) \dVol_g
	\\
	\\
	&\lesssim \int_{\Sigma_{\tau'}} \bigg(
	\sum_{n=0}^2 \epsilon^2 (1+r)^{-1-\delta}|\slashed{\D} \mathscr{Z}^{n} \phi|^2
	+ \sum_{n=0}^2 \epsilon^2 (1+r)^{-3-\delta} |\mathscr{Z}^n \phi|^2 \bigg) r^2\upd r \wedge \dVol_{\mathbb{S}^2}
	\\
	&\phantom{=}
	+ \int_{\mathcal{M}_{\tau}^{\tau'}} \bigg(
	\sum_{n=0}^3 \epsilon^2 (1+r)^{-1-\frac{1}{2}\delta}|\slashed{\D}\mathscr{Z}^n\phi|^2 
	+ \sum_{n=0}^3 \epsilon^2 (1+r)^{-3-\frac{1}{2}\delta} |\mathscr{Z}^{n} \phi|^2 \bigg) \dVol_g
	\end{split}
	\end{equation*}
	where we have used the fact that
	\begin{equation*}
	\left.\frac{\partial}{\partial \tau} \right|_{r, \vartheta^A} = \mu T - \frac{1}{2}\mu b^\mu \slashed{\nabla}_\mu 
	\end{equation*}
	together with the fact that $\dVol_g = -\mu\Omega^2 \upd \tau \wedge \upd r \wedge \dVol_{\mathbb{S}^2} \sim -\mu r^2 \upd \tau \wedge \upd r \wedge \dVol_{\mathbb{S}^2}$.

	Now, the conditions in given in the proposition imply that we can pick a dyadic sequence of times $\tau_n \rightarrow \infty$ such that
	\begin{equation*}
	\int_{\Sigma_{\tau_n}} \bigg(
	\sum_{n=0}^2 \epsilon^2 (1+r)^{-1-\delta}|\slashed{\D} \mathscr{Z}^{n} \phi|^2
	+ \sum_{n=0}^2 \epsilon^2 (1+r)^{-3-\delta} |\mathscr{Z}^n \phi|^2 \bigg) r^2\upd r \wedge \dVol_{\mathbb{S}^2}
	\lesssim
	\epsilon^2 \delta^{-1} \mathcal{E} (1+\tau_n)^{-1}
	\end{equation*}
	Since this sequence is dyadic, we can choose $\tau'$ to be a member of this sequence and, at the same time, to satisfy $\tau' \sim \tau$. Hence, choosing this value for $\tau'$, and using the assumptions of the proposition we find that
	\begin{equation*}
	\int_{\Sigma_\tau \cap} \bigg(
	\sum_{n=0}^2 \epsilon^2 (1+r)^{-1-\delta}|\slashed{\D} \mathscr{Z}^{n} \phi|^2
	+ \sum_{n=0}^2 \epsilon^2 (1+r)^{-3-\delta} |\mathscr{Z}^n \phi|^2 \bigg) r^2\upd r \wedge \dVol_{\mathbb{S}^2}
	\lesssim \epsilon^2 \delta^{-1}\mathcal{E} (1+\tau)^{-1}		
	\end{equation*}
	In particular, the restriction of the integrand to the region $r' \leq r \leq R$ is also bounded by $\epsilon^2 \delta^{-1}\mathcal{E} (1+\tau)^{-1}$.

	Finally, to control the terms which are integrals over the sphere $S_{\tau,R}$, we use corollary \ref{corollary decay for spherical mean along a subsequence} and choose $R = r_n$ to be one of the radii given in that corollary. Hence:
	\begin{equation*}
	\begin{split}
	&\int_{S_{\tau, r_n}} \sum_{n=0}^2 \left( \epsilon^2 (1+r)^{2C_{(1)}\epsilon} |\slashed{\D} \mathscr{Z}^n \phi|^2 + \epsilon^2 (1+r)^{-2 + C_{(2)}\epsilon} |\mathscr{Z}^n \phi |^2 \right) \dVol_{\mathbb{S}^2}
	\\
	&\lesssim
	\epsilon^2 \delta^{-2} \mathcal{E} (1+r_n)^{-1 + \frac{1}{2}\delta + 2C_{(2)}\epsilon}(1+\tau)^{-1}
	\\
	&\lesssim
	\epsilon^2 \delta^{-2} \mathcal{E} (1+r_n)^{-1 + \delta}(1+\tau)^{-1}
	\end{split}
	\end{equation*}
	Recall also that we can choose a sequence of such $r_n$'s such that $r_n \rightarrow \infty$, in which this final term tends to zero. Hence, using the conditions of the proposition to control the integrals of $\mathscr{Z}^n F$ and the the good derivatives $\overline{\slashed{\D}} \mathscr{Z}^n \phi$,we have
	\begin{equation*}
	\begin{split}
	\sum_{m = 0}^2 \int_{S_{\tau, r'}} |(r\slashed{\nabla})^m \slashed{\D}_{\Lbar}\phi|^2 \dVol_{\mathbb{S}^2}
	&\lesssim
	(r')^{-2+2\delta} \mathcal{E} \left( (1+\tau)^{-2\beta} + (1+\tau)^{-1} \right)
	\\
	&\lesssim
	(r')^{-2+2\delta} \mathcal{E} (1+\tau)^{-2\beta}
	\end{split}
	\end{equation*}

\end{proof}

We remark here that, in the proposition above, we have made our first use of the ``improved'' pointwise decay bootstrap assumptions, i.e.\ bootstrap assumptions of the form
\begin{equation*}
\bm{\Gamma}^{(n)}_{(-1+C_{(n)}\epsilon)} \lesssim \epsilon (1+r)^{-1+C_{(n)}\epsilon}
\end{equation*}
Previously we had only used the bootstrap assumptions of the form
\begin{equation*}
\bm{\Gamma}^{(n)}_{(-1+C_{(n)}\epsilon)} \lesssim \epsilon (1+r)^{-1+\delta}(1+\tau)^{-\beta} + (1+r)^{-1-\delta}
\end{equation*}
which does not have such good \emph{uniform} decay in $r$.

We also need to prove that the ``good'' derivatives of the field, $\slashed{\nabla}\phi$ and $\slashed{\D}_L \phi$, obey improved pointwise estimates.

\begin{proposition}[Improved pointwise decay for angular derivatives]
	\label{proposition improved pointwise bound angular}
Let $\phi$ be an $S_{\tau,r}$-tangent tensor field satisfying
\begin{equation*}
 \sum_{n=0}^{3} \mathcal{E}^{(w_{n}T)}[ (r\slashed{\nabla})^n \phi](\tau) \lesssim \tilde{\mathcal{E}}
\end{equation*}
where $w_{n} = (1+r)^{-C_{\left[r\slashed{\nabla}^n \phi \right]}\epsilon}$.

Suppose also that
\begin{equation*}
\sum_{n=0}^{2} \mathcal{E}^{(L, 2\delta)}[ (r\slashed{\nabla})^n (r\slashed{\D}_L \phi)](\tau) \lesssim \tilde{\mathcal{E}}
\end{equation*}

Additionally, suppose that the initial data on $\Sigma_{\tau_0}$ satisfies
\begin{equation*}
 \lim_{t \rightarrow \infty} \sum_{n=0}^3 \int_{\bar{S}_{\tau_0, t}} |(r\slashed{\nabla})^n \phi|^2 = 0
\end{equation*}
where $\tau_0 \leq \tau$, and suppose that all of the bootstrap bounds of chapter \ref{chapter bootstrap} hold.

Then, for all $r \geq r_0$ we have
\begin{equation}
\begin{split}
	|\slashed{\nabla}\phi| 
	&\lesssim
	\sqrt{\tilde{\mathcal{E}}} (1+r)^{-\frac{3}{2} + \frac{1}{2}C^*\epsilon} 
\end{split}
\end{equation}
\end{proposition}
where 
\begin{equation*}
C^* = \sup_{0\leq n \leq 3} C_{\left[ (r\slashed{\nabla})^n \phi \right]}
\end{equation*}

\begin{proof}
This proposition is easily proved by applying proposition \ref{proposition pointwise bound on field in terms of energy} to the field $(r\slashed{\nabla}\phi)$.
\end{proof}

\begin{proposition}[Improved pointwise decay for $L$ derivatives]
	\label{proposition improved pointwise bound L}
	
	Let $\phi$ be an $S_{\tau,r}$-tangent tensor field satisfying
	\begin{equation*}
	\sum_{n=0}^2 \mathcal{E}^{(L, p)} [(r\slashed{\nabla})^n (r\slashed{\D}_L \phi)](\tau) \lesssim \tilde{\mathcal{E}}
	\end{equation*}
	
	Suppose also that $\phi$ satisfies the wave equation
	\begin{equation*}
	\tilde{\slashed{\Box}}_g \phi = F
	\end{equation*}
	and that
	\begin{equation*}
	\begin{split}
	&\int_{\Sigma_{\tau}\cap\{r \geq r'\}}
	\bigg(
	r^{1 - 2\delta} \sum_{0 \leq n \leq 2} |(r\slashed{\nabla})^n F|^2
	+ r^{-1 - \delta} \sum_{0 \leq n \leq 3} |\overline{\slashed{\D}} (r\slashed{\nabla})^n \phi|^2 
	\bigg) r^2 \upd r \wedge \dVol_{\mathbb{S}^2}  \\
	&  
	+ \int_{\mathcal{M}_\tau^{\infty} \cap \{r' \leq r\}}
	r^{-1-\delta} \sum_{n = 0}^3 |\slashed{\D} \mathscr{Z}^n \phi|^2 \dVol_g 
	\leq \mathcal{E}_1
	\end{split}
	\end{equation*}
	and
	\begin{equation*}
	\sum_{n=0}^3 \mathcal{E}^{(w_n T)}[\mathscr{Z}^n \phi] \lesssim \tilde{\mathcal{E}}_1
	\end{equation*}
	where $w_n = (1+r)^{-C_{(n)}\epsilon}$.
	
	Finally, suppose that for some $\tau_0 \leq \tau$ we have
	\begin{equation*}
	\lim_{t\rightarrow\infty} \sum_{n=0}^3 \int_{\bar{S}_{\tau_0, t}} |\mathscr{Z}^n \phi|^2 = 0
	\end{equation*}
	
	Then for $r \geq r_0$ we have
	\begin{equation*}
	|\slashed{\D}_L \phi|
	\lesssim
	r^{-4} \mathcal{E}_1
	+ \frac{1}{1-p} r^{-3-p} \mathcal{E}
	\end{equation*}

\end{proposition}

\begin{proof}
	Using proposition \ref{proposition higher weighted spherical integral} applied to the field $(r\slashed{\nabla})^n r\slashed{\D}_L \phi$ we obtain
	\begin{equation*}
	\sum_{n=0}^2 \int_{S_{\tau,r}} |(r\slashed{\nabla})^n \slashed{\D}_L \phi|^2 \dVol_{\mathbb{S}^2}
	\lesssim
	r^{-4} \sum_{n=0}^2 \int_{S_{\tau,r_0}} |(r\slashed{\nabla})^n \slashed{\D}_L \phi|^2 \dVol_{\mathbb{S}^2}
	+ \frac{1}{1-p} r^{-3-p} \sum_{n=0}^2 \mathcal{E}^{(L,p)}[(r\slashed{\nabla})^n r\slashed{\D}_L \phi](\tau)
	\end{equation*}
	
	Now, we can estimate
	\begin{equation*}
	\begin{split}
	\sum_{n=0}^2 \int_{S_{\tau,r_0}} |(r\slashed{\nabla})^n \slashed{\D}_L \phi|^2 \dVol_{\mathbb{S}^2}
	&\lesssim
	\sum_{n=0}^2 \int_{S_{\tau,r_0}}\left(
		|(r\slashed{\nabla})^n \slashed{\D}_T \phi|^2
		+ |(r\slashed{\nabla})^n \slashed{\D}_{\Lbar} \phi|^2
	\right)	\dVol_{\mathbb{S}^2}
	\end{split}
	\end{equation*}
	The first term is bounded during the proof of \ref{proposition pointwise bound on field in terms of energy}, and the second term is bounded during the proof of proposition \ref{proposition pointwise bound lbar}. We obtain
	\begin{equation*}
	\sum_{n=0}^2 \int_{S_{\tau,r_0}}
	|(r\slashed{\nabla})^n \slashed{\D}_T \phi|^2 
	\dVol_{\mathbb{S}^2}
	\lesssim
	(1+r)^{-\frac{1}{2} + \frac{1}{2} C^* \epsilon} \sum_{n=0}^3 \mathcal{E}^{(w_n T)}[\mathscr{Z}^n \phi](\tau)
	\lesssim
	\tilde{\mathcal{E}}_1(1+r)^{-\frac{1}{2} + \frac{1}{2} C^* \epsilon}
	\end{equation*}
	and
	\begin{equation*}
	\sum_{n=0}^2 \int_{S_{\tau,r_0}}
	|(r\slashed{\nabla})^n \slashed{\D}_{\Lbar} \phi|^2
	\dVol_{\mathbb{S}^2}
	\lesssim
	\sqrt{\tilde{\mathcal{E}}_1} r^{-1+\delta}
	\end{equation*}
	Combining these bounds proves the proposition.

\end{proof}

\begin{remark}[An alternative way to obtain pointwise bounds for $L$ derivatives]
	
	We remark here that an alternative way of bounding the $L$ derivatives is available to us, based on using the one-dimensional Sobolev inequality on the half-line $r \geq r_0$. In order to use this, we need to have $L^2$ bounds for both $r^{1+\delta} L\phi$ and $L\left(r^{1+\delta} L\phi \right)$ which, in fact, are available to us in the form of the $p$-weighted energy estimates. In fact, this approach\footnote{Note that this is the approach taken in \cite{Yang2013}.} can lead to improved decay for the $L$ derivatives \emph{without} commuting with the operator $r\slashed{\D}_L$, but, taking this approach, we would not simultaneously obtain decay in $\tau$. Note that our bootstrap assumptions do not actually require decay at a rate which is simultaneously better than $r^{-1}$ \emph{and} decaying in $\tau$!
	
\end{remark}

\subsection{Improved pointwise bounds on bad derivatives}
If we have \emph{already} produced bounds on the ``good'' derivatives, then we can use these to produce improved pointwise bounds on the ``bad'' derivatives by integration along the integral curves of $L$. These are a kind of $L^\infty$-$L^\infty$ bound, since they make direct use of $L^\infty$ bounds of the kind we have proved above, in order to improve certain other $L^\infty$ bounds.

\begin{proposition}
	\label{proposition improved pointwise bound Lbar}
 Let $\phi$ be an $S_{\tau,r}$-tangent tensor field satisfying
\begin{equation*}
 \tilde{\Box}_g \phi = F
\end{equation*}
for some $S_{\tau,r}$-tangent tensor field $F$. Suppose that the pointwise bootstrap assumptions of chapter \ref{chapter bootstrap} hold, and, moreover, suppose that the ``good'' derivatives of $\phi$ satisfy
\begin{equation*}
 \sum_{n=0}^{1} |\overline{\slashed{\D}}\mathscr{Z}\phi| \lesssim \sqrt{\mathcal{E}} r^{-1-\delta}
\end{equation*}
in the region $r \geq r_0$. If $\phi$ is not a scalar field but a higher rank tensor field, then suppose also that
\begin{equation*}
 |\phi| \lesssim \sqrt{\mathcal{E}} r^{-\frac{1}{2} + \delta}
\end{equation*}
in the region $r \geq r_0$. Finally, in all cases, suppose that
\begin{equation*}
 \sup_{r = r_0} |\slashed{\D}_{\Lbar} \phi| \lesssim \sqrt{\mathcal{E}}
\end{equation*}

Then, if $F$ satisfies
\begin{equation*}
 |F| \lesssim \epsilon r^{-1-\delta}|\slashed{\D}_{\Lbar} \phi| + \sqrt{\mathcal{E}} r^{-2-\delta}
\end{equation*}
in the region $r \geq r_0$, then for all $r \geq r_0$ we have
\begin{equation}
 |\slashed{\D}_{\Lbar} \phi| \lesssim \delta^{-1} \sqrt{\mathcal{E}} r^{-1}
\end{equation}

On the other hand, if $F$ satisfies
\begin{equation*}
 |F| \lesssim \epsilon r^{-1}|\slashed{\D}_{\Lbar} \phi| + \sqrt{\mathcal{E}} r^{-2+2C_{(F)}\epsilon}
\end{equation*}
in the region $r \geq r_0$, then for all $r \geq r_0$ we have
\begin{equation}
 |\slashed{\D}_{\Lbar} \phi| \lesssim \tilde{C}^{-1} \epsilon^{-1} \sqrt{\mathcal{E}} r^{-1 + \tilde{C}\epsilon}
\end{equation}
where the constant $\tilde{C}$ is chosen sufficiently large compared with the constant $C_{(F)}$ and the implicit constants in the other inequalities.

\end{proposition}

\begin{proof}
 First, consider the case when $\phi$ is a scalar field. From proposition \ref{proposition scalar wave operator} we have
\begin{equation*}
 L(r\Lbar \phi) = -r\tilde{\Box}_g \phi  + r\slashed{\Delta}\phi - \frac{1}{2}(\tr_{\slashed{g}}\chi_{(\text{small})}) (r\Lbar \phi) + \left( \frac{1}{r} - \frac{1}{2}\tr_{\slashed{g}}\chibar_{(\text{small})} \right) rL\phi - \zeta^\alpha r\slashed{\nabla}\phi
\end{equation*}
and so, we have
\begin{equation*}
 \begin{split}
  &|(r\Lbar \phi)(\tau, r, \vartheta^1, \vartheta^2)| \leq 
  |(r\Lbar \phi)(\tau, r_0, \vartheta^1, \vartheta^2)| \\
  &\phantom{\lesssim}
  + \int_{r_0}^r \bigg(
  r'|F|
  + (r')^{-1}|\slashed{\nabla}r(\slashed{\nabla})\phi|
  + |\tr_{\slashed{g}}\chi_{(\text{small})}| |r'\Lbar \phi|
  + \left( (r')^{-1} + |\tr_{\slashed{g}}\chibar_{(\text{small})}| \right) r'|L\phi|
  + |\zeta| r'|\slashed{\nabla}\phi|
  \bigg) \upd r' \\ \\
 \end{split}
\end{equation*}
Substituting in the assumed bounds on the good derivatives of $\phi$, together with the pointwise bounds on the connection coefficients assumed in chapter \ref{chapter bootstrap}, we obtain
\begin{equation*}
 \begin{split}
  |(r\Lbar \phi)(\tau, r, \vartheta^1, \vartheta^2)|
  &\lesssim r_0|(\Lbar \phi)(\tau, r_0, \vartheta^1, \vartheta^2)|
  + \int_{r_0}^r \bigg(
  r'|F|
  + \epsilon r^{-1-\delta}|r'\Lbar\phi|
  + \sqrt{\mathcal{E}} r^{-1-\delta+\mathring{C}\epsilon}
  \bigg) \upd r'
 \end{split}
\end{equation*}
Now, if in the region $r \geq r_0$, $F$ satisfies the bound
\begin{equation*}
 |F| \lesssim \epsilon r^{-1-\delta}|\Lbar \phi| + \sqrt{\mathcal{E}} r^{-2 - \delta}
\end{equation*}
then we can estimate
\begin{equation*}
  |(r\Lbar \phi)(\tau, r, \vartheta^1, \vartheta^2)|
  \lesssim
  r_0|(\Lbar \phi)(\tau, r_0, \vartheta^1, \vartheta^2)|
  + \int_{r_0}^r \bigg(
  \epsilon r^{-1-\delta}|r'\Lbar\phi|
  + \sqrt{\mathcal{E}} r^{-1-\delta+\mathring{C}\epsilon}
  \bigg) \upd r'
\end{equation*}
and now, using the Gronwall inequality, we find
\begin{equation*}
 |\Lbar \phi|
 \lesssim
 r^{-1} \left( |(\Lbar \phi)(\tau, r_0, \vartheta^1, \vartheta^2)| + \delta^{-1}\sqrt{\mathcal{E}} \right)
\end{equation*}

If instead $F$ satisfies the bound
\begin{equation*}
 |F| \lesssim \epsilon r^{-1}|\Lbar \phi| + \sqrt{\mathcal{E}} r^{-2 + 2C_{(F)}\epsilon}
\end{equation*}
then we can estimate
\begin{equation*}
  |(r\Lbar \phi)(\tau, r, \vartheta^1, \vartheta^2)|
  \lesssim
  r_0|(\Lbar \phi)(\tau, r_0, \vartheta^1, \vartheta^2)|
  + \int_{r_0}^r \bigg(
  \epsilon r^{-1}|r'\Lbar\phi|
  + C_{(F)}^{-1} \epsilon^{-1} \sqrt{\mathcal{E}} r^{-1-2C_{(F)}\epsilon}
  \bigg) \upd r'
\end{equation*}
for some constant $C_{(F)} > 1$. This time the Gronwall inequality yields
\begin{equation*}
 |\Lbar \phi|
 \lesssim
 r^{-1 + \tilde{C}\epsilon} \left( |(\Lbar \phi)(\tau, r_0, \vartheta^1, \vartheta^2)| + \sqrt{\mathcal{E}} \right)
\end{equation*}
where the constant $\tilde{C}$ must be chosen sufficiently large compared with both $C_{(F)}$ and the implicit constant in the preceeding inequality.

Now, we need to generalise these results to the case when $\phi$ is an $S_{\tau,r}$-tangent tensor field rather than a scalar field. In this case before we are able to integrate, we first need to ``scalarise'' the equations by contracting with the vector fields $\slashed{\Pi}_a^{\phantom{a}\mu}$. 

We first note that, from proposition \ref{proposition tensor wave operator} that, schematically,
\begin{equation}
\label{equation schematic L r Lbar phi}
 \begin{split}
  \slashed{\D}_L \left( r \slashed{\D}_{\Lbar} \phi \right)
  &=
  -r\tilde{\slashed{\Box}}_g \phi
  + r\slashed{\Delta}\phi
  - \frac{1}{2}r\tr_{\slashed{g}}\chi_{(\text{small})} \slashed{\D}_{\Lbar}\phi
  + \left(1 - \frac{1}{2}r\tr_{\slashed{g}}\chibar \right) \slashed{\D}_L \phi
  + r\zeta \slashed{\nabla}\phi \\
  &\phantom{=}
  + r\Omega_{L\Lbar}\phi
 \end{split}
\end{equation}

Now, making use of proposition \ref{proposition transport rectangular} together with the fact that
\begin{equation*}
 \slashed{\Pi}_\mu^{\phantom{\mu}a} \slashed{\Pi}^\mu_{\phantom{\mu}b} = \delta_a^b
\end{equation*}
we find that
\begin{equation*}
  \left|\slashed{\D}_L \left( \slashed{\Pi}^\mu_{\phantom{\mu}a} \right) \right|
  \lesssim \left( |\bar{\partial} h|_{(\text{frame})} + r^{-1} |L^i_{(\text{small})} \slashed{\Pi}^i| \right)\cdot \left( |L_a| + |\Lbar_a| + |\slashed{\Pi}_a| \right)
\end{equation*}
and so, using the pointwise bootstrap assumptions, we have that
\begin{equation*}
 \left|\slashed{\D}_L \left( \slashed{\Pi}^\mu_{\phantom{\mu}a} \right) \right|
  \lesssim \epsilon r^{-1-\delta}
\end{equation*}
in the region $r \geq r_0$.

Hence, contracting equation \eqref{equation schematic L r Lbar phi} with the vector fields $\slashed{\Pi}^{\mu_1}_{\phantom{\mu_1} a_1}$, \ldots $\slashed{\Pi}^{\mu_n}_{\phantom{\mu_n} a_n}$, and taking $\phi$ to be a rank $n$ tensor satisfying $\tilde{\Box}_g \phi = F$, we obtain
\begin{equation*}
 \begin{split}
  \left| L \left( \left( r\slashed{\D}_{\Lbar}\phi \right)_{a_1 \ldots a_n} \right) \right|
  &\lesssim
  r|F_{a_1 \ldots a_n}|
  + \left|\left(\slashed{\nabla} (r\slashed{\nabla})\phi\right)_{a_1 \ldots a_n} \right|
  + \epsilon r^{-1-\delta} \left| \left(r\slashed{\D}_{\Lbar} \phi \right)_{a_1 \ldots a_n} \right| \\
  &\phantom{\lesssim}
  + \left( 1 + \epsilon r^{C_{(1)}\epsilon} \right) \left|\left(\overline{\slashed{\D}} \phi \right)_{a_1 \ldots a_n} \right|
  + \epsilon r^{-1 + 2C_{(1)}\epsilon} |\phi|
 \end{split}
\end{equation*}
where we have also made use of the pointwise bootstrap bounds for the connection coefficients and the curvature. Now, making use of the bounds assumed in the proposition for $|\phi|$ and $|\overline{\slashed{\D}}\phi|$, we obtain
\begin{equation*}
 \left| L \left( \left( r\slashed{\D}_{\Lbar}\phi \right)_{a_1 \ldots a_n} \right) \right|
  \lesssim
  r|F_{a_1 \ldots a_n}|
  + \epsilon^2 r^{-1-\delta}
  + \epsilon r^{-1-\delta} \left| \left( r\slashed{\D}_{\Lbar}\phi \right)_{a_1 \ldots a_n} \right|
\end{equation*}
The rest of the proof follows in exactly the same way as the scalar case. Note that, strictly speaking, this establishes bounds on the \emph{scalar fields} $\left( r\slashed{\D}_{\Lbar}\phi \right)_{a_1 \ldots a_n}$. However, since these bounds are established for \emph{all} sets of rectangular indices, and since also $|\slashed{\Pi}_{\mu}^{\phantom{\mu}a}| \lesssim 1$ and $|(\slashed{g}^{-1})^{ab}| \sim 1$, this actually establishes bounds on the tensor field $r\slashed{\D}_{\Lbar}\phi$.

\end{proof}

\section{Summary of pointwise estimates for solutions of the wave equation}

In the following lemma we bring together the pointwise bounds we have proved above, to give an overall summary of pointwise decay results for solutions of the wave equation. Note that this can be viewed as a result establishing pointwise bounds on solutions to the \emph{linear} inhomogeneous wave equation on manifolds whose geometry obeys the pointwise bounds given in chapter \ref{chapter bootstrap}.

\begin{lemma}[Pointwise bounds for solutions of the wave equation]
	\label{lemma pointwise bounds summary}
	
	Suppose that the pointwise bootstrap bounds of chapter \ref{chapter bootstrap} hold. Suppose additionally that the fields $h_{ab}$ satisfy
	\begin{equation}
	|\tilde{\Box}_g \mathscr{Z}^n h| \lesssim (1+r)^{-1-2\delta}
	\end{equation}
	for $n \leq 2$.
	
	Let $\phi$ be an $S_{\tau,r}$-tangent tensor field satisfying
	\begin{equation*}
	\tilde{\slashed{\Box}}_g \phi = F = F_1 + F_2
	\end{equation*}
	for some $S_{\tau,r}$-tangent tensor field $F$. Suppose also that
	\begin{equation*}
	\tilde{\slashed{\Box}}_g \left( (r\slashed{\nabla})^m \phi \right) = F_{(m)} = F_{(m,1)} + F_{(m,2)}
	\end{equation*}
	for $0 \leq m \leq 4$.
	
	Also, suppose that, for $0 \leq m \leq 2$ we have
	\begin{equation*}
	\tilde{\slashed{\Box}}_g (r\slashed{\D}_L (r\slashed{\nabla})^m \phi)
	=
	\slashed{\Delta} (r\slashed{\nabla})^m \phi
	+ r^{-1} \slashed{\D}_L \left( r\slashed{\D}_L \left( (r\slashed{\nabla})^m \phi \right) \right)
	+ r^{-1} \slashed{\D}_L \left( (r\slashed{\nabla})^m \phi \right)
	+ F_{(m, rL, 1)}
	+ F_{(m, rL, 2)}
	\end{equation*}
	We write $F_{(m, rL)} := F_{(m, rL, 1)} + F_{(m, rL, 2)}$.
	
	Suppose that the inhomogenous terms satisfy
	\begin{equation*}
	\begin{split}
	&\sum_{m=0}^4 \int_{\mathcal{M}_{\tau_0}^\tau} \bigg(
	\epsilon^{-1} \chi_{(r_0)} r^{1-C_{\left[ (r\slashed{\nabla})^m \phi \right]}\epsilon} (1+\tau)^{1+\delta} |F_{(m,1)}|^2
	\\
	&\phantom{\sum_{m=0}^4 \int_{\mathcal{M}_{\tau_0}^\tau} \bigg(}
	+ \epsilon^{-1} \chi_{(r_0)} r^{2-C_{\left[ (r\slashed{\nabla})^m \phi \right]}\epsilon - 2\delta} (1+\tau)^{2\beta}|F_{(m,2)}|^2 \bigg) \dVol_g \lesssim \mathcal{E}_0
	\\
	\\
	&\sum_{m=0}^4 \int_{\mathcal{M}_{\tau}^{\tau_1}} \bigg(
	\epsilon^{-1}(1+r)^{1-C_{\left[ (r\slashed{\nabla})^m \phi \right]}\epsilon} \left(|F_{(m)}|^2 + |\mathscr{Z}^m F|^2 \right)
	+ \epsilon^{-1}(1+r)^{\frac{1}{2}\delta}(1+\tau)^{1+\delta}|F_{(m,1)}|^2
	\\
	&\phantom{\sum_{m=0}^4 \int_{\mathcal{M}_{\tau_0}^\tau} \bigg(}
	+ \epsilon^{-1}(1+r)^{1-\delta}(1+\tau)^{2\beta}|F_{(m,1)}|^2
	\bigg)\dVol_g
	\lesssim
	\mathcal{E}_0 (1+\tau)^{-1}
	\\
	\\
	&\sum_{m=0}^2 \int_{\mathcal{M}^{\tau_1}_{\tau}} \left(
		\epsilon^{-1} \chi_{(2r_0)} r^{1-\frac{1}{2}\delta} (1+\tau)^{2\beta}|F_{(m, rL)}|^2
	\right) \dVol_g 
	\lesssim \mathcal{E}_0 (1+\tau)^{-1+C_{(m)}\delta}
	\\
	\\
	&\sum_{m=0}^2 \int_{\mathcal{M}^\tau_{\tau_0}} \bigg(
		\epsilon^{-1} \chi_{(2r_0)} r^{1-C_{\left[ (r\slashed{\nabla})^m \phi \right]}\epsilon} (1+\tau)^{1+\delta}|F_{(m, rL,1)}|^2
		\\
		&\phantom{\sum_{m=0}^2 \int_{\mathcal{M}^\tau_{\tau_0}} \bigg(}
		+ \epsilon^{-1} \chi_{(2r_0)} r^{2-2\delta+C_{\left[ (r\slashed{\nabla})^m \phi \right]}\epsilon} (1+\tau)^{2\beta}|F_{(m, rL,2)}|^2
	\bigg) \dVol_g 
	\lesssim \mathcal{E}_0
	\end{split}
	\end{equation*}
	for all $\tau \geq \tau_0$ and for all $\tau_1 \geq \tau$.
	
	Define
	\begin{equation*}
	\begin{split}
	w_{(m)} &= (1+r)^{-C_{\left[ (r\slashed{\nabla})^m \phi \right]}}
	\\
	\tilde{w}_{(m)} &= (1+r)^{-\frac{1}{2}C_{\left[ (r\slashed{\nabla})^m \phi \right]}}
	\end{split}
	\end{equation*}
	
	Suppose that the initial data for $\phi$ at $\tau = \tau_0$ satisfies
	\begin{equation*}
	\begin{split}
	\int_{\bar{S}_{\tau_0, t}} \sum_{m = 0}^4 |\mathscr{Z}^m\phi|^2 \dVol_{\mathbb{S}^2}
	&\lesssim
	\tilde{\mathcal{E}}(t-\tau_0)^{-1+\frac{1}{2}C_{[\mathscr{Z}^m\phi]}\epsilon}
	\\
	\sum_{m=0}^{4} \mathcal{E}^{(\tilde{w}_m T)}[\mathscr{Z}^m \phi](\tau_0)
	&\lesssim
	\tilde{\mathcal{E}}
	\\
	\delta^{-1}\sum_{m=0}^{4} \mathcal{E}^{(L, 1-C_{[\mathscr{Z}^m \phi]}\epsilon)}[\mathscr{Z}^m \phi](\tau_0)
	&\lesssim
	\tilde{\mathcal{E}}
	\\
	\sum_{m=0}^2 \mathcal{E}^{(L, 1-C_{[\mathscr{Z}^m \phi]}\epsilon)}[r\slashed{\D}_L \mathscr{Z}^m \phi] (\tau_0)
	&\lesssim
	\tilde{\mathcal{E}}
	\end{split}
	\end{equation*}

	Define
	\begin{equation*}
	C^* = \max_{0 \leq m \leq 4} C_{\left[ (r\slashed{\nabla})^m \phi \right]}
	\end{equation*}
	
	Then, if $\epsilon$ is sufficiently small, in the region $r \geq r_0$ we have
	\begin{equation*}
	\begin{split}
	|\phi| &\lesssim
	\sqrt{\tilde{\mathcal{E}}}(1+r)^{-\frac{1}{2} + \frac{1}{2} C^* \epsilon} (1+\tau)^{-\frac{1}{2} + \frac{1}{2}C^* \delta}
	\\
	|\slashed{\D}_{\Lbar} \phi| &\lesssim
	\sqrt{\tilde{\mathcal{E}}}(1+r)^{-\frac{1}{2} + \delta} (1+\tau)^{-\frac{1}{2} + \frac{1}{2}C^* \delta}
	\\
	|\slashed{\nabla}\phi| &\lesssim
	\sqrt{\tilde{\mathcal{E}}}(1+r)^{-\frac{3}{2} + \frac{1}{2} C^* \epsilon} (1+\tau)^{-\frac{1}{2} + \frac{1}{2}C^* \delta}
	\\
	|\slashed{\D}_L \phi| &\lesssim
	\sqrt{\tilde{\mathcal{E}}}(1+r)^{-\frac{3}{2} + \delta} (1+\tau)^{-\frac{1}{2} + \frac{1}{2}C^* \delta}
	\end{split}
	\end{equation*}
	In particular, for sufficiently small $\delta$ we have
	\begin{equation*}
	\begin{split}
	|\slashed{\D} \phi| &\lesssim
	\sqrt{\tilde{\mathcal{E}}}\left( (1+r)^{-\frac{1}{2} + \frac{1}{2} C^* \epsilon} (1+\tau)^{-\beta} + (1+r)^{-\frac{3}{2} + \delta} \right)
	\\
	|\overline{\slashed{\D}} \phi| &\lesssim
	\sqrt{\tilde{\mathcal{E}}}(1+r)^{-\frac{3}{2} + \delta}
	\end{split}
	\end{equation*}
	
	Additionally, suppose that, in the region $r \geq r_0$, $F$ satisfies the bound
	\begin{equation*}
	|F| \lesssim \epsilon r^{-1-\delta} |\slashed{\D}_{\Lbar} \phi| + \sqrt{\mathcal{E}} r^{-2-\delta} 
	\end{equation*}
		
	Then we have the following bound, giving improved decay in $r$ for the bad derivatives:
	\begin{equation*}
	|\slashed{\D}_{\Lbar} \phi| \lesssim \sqrt{\mathcal{E}} (1+r)^{-1}
	\end{equation*}
	
	On the other hand, suppose instead that in the region $r \geq r_0$, $F$ satisfies the bound
	\begin{equation*}
	|F| \lesssim \epsilon r^{-1} |\slashed{\D}_{\Lbar} \phi| + \sqrt{\mathcal{E}} r^{-2 + C_{(F)}\epsilon} 
	\end{equation*}
	for some constant $C_{(F)}$.
	
	Then we have the following bound, giving improved decay in $r$ for the bad derivatives:
	\begin{equation*}
	|\slashed{\D}_{\Lbar} \phi| \lesssim \sqrt{\mathcal{E}} (1+r)^{-1 + \tilde{C}\epsilon}
	\end{equation*}
	where $\tilde{C}$ is some constant which is sufficiently large compared with $C_{(F)}$.

\end{lemma}

\begin{proof}
	This result follows from combining all the pointwise bounds proven in this chapter with the energy estimates of the previous chapter.
\end{proof}

\begin{remark}
	The precise number of operators $\mathscr{Z}$ in the proposition above is, in a certain sense, wasteful. To be precise: if we wish only to bound $|\phi|$ pointwise, then we need only commute with the operators $\mathscr{Z}$ twice. In order to obtain the first set of bounds on the derivatives $|\slashed{\D} \phi|$, we need only to commute with the operators $\mathscr{Z}$ three times. Commuting four times is only required in order to obtain the improved pointwise bounds on the bad derivatives stated at the end of the lemma. Commuting \emph{after already applying the operator} $r\slashed{\D}_L$ is, of course, only necessary in order to control the term $|\slashed{\D}_L (r\slashed{\D}_L\phi)|$
\end{remark}

\section{Pointwise decay estimates for \texorpdfstring{$\mathscr{Y}^n \phi$}{Yn phi}}

Note that the summary above follows from bounds on the (weighted) $T$-energy, together with some bounds on the $p$-weighted energies. To put this another way, we assume that the field $\phi$ satisfies a wave equation with an inhomogeneous term which is suitably small. As such, we can use these bounds to give pointwise bounds on fields like $\mathscr{Z}^n \phi$ and $\slashed{\D} \mathscr{Z}^n \phi$. Note, however, that we also commute with the operator $r\slashed{\D}_L$, and the resulting fields do not satisfy an equation of the required form. On the other hand, we can still establish $p$-weighted energy estimates after commuting with $r\slashed{\D}_L$, so we need estimates for the various quantities which only makes use of the $p$-weighted energy. These is given in this subsection.

In some of the following estimates, we encounter boundary terms at $r = r_0$. Our strategy for dealing with these will be to exchange operators $\mathscr{Y}$ for a combination of operators of the form $\mathscr{Z}$ and $\tilde{\slashed{\Box}}_g$, making use of the fact that, at $r = r_0$, we do not have to worry about factors of $r$.

\begin{proposition}[Exchanging derivatives of the sphere $r = r_0$]
	\label{proposition exchanging derivatives at r0}
	Let $\phi$ satisfy
	\begin{equation*}
	\tilde{\slashed{\Box}} \mathscr{Z}^\ell \phi = F_{(\ell)}
	\end{equation*}
	and we write $F = F_{(0)}$. Then, if we evaluate all the terms on the sphere $r = r_0$, we have the following schematic expression for $\mathscr{Y}^n \phi$: for all $n \geq 0$,
	\begin{equation}
	\mathscr{Y}^n \phi 
	=
	\mathscr{Z}^{\leq n} \phi
	+ \slashed{\D} \mathscr{Z}^{\leq n-1} \phi
	+ \mathscr{Y}^{(\leq n-2)} F
	+ \sum_{j+k \leq n-1} \bm{\Gamma}^{(j)}_{(0)} \slashed{\D} \mathscr{Z}^{k} \phi
	+ \sum_{j+k+\ell \leq n-1} \bm{\Gamma}^{(j)}_{(0)} \mathscr{Y}^{k-1} F_{(\ell)}
	\end{equation}
\end{proposition}

\begin{proof}
	Suppose that, for all $n \leq N$, on the sphere $S_{\tau,r_0}$ we have the schematic expression
	\begin{equation*}
	\mathscr{Y}^n \phi 
	=
	\mathscr{Z}^{\leq n} \phi
	+ \slashed{\D} \mathscr{Z}^{\leq n-1} \phi
	+ \mathscr{Y}^{(\leq n-2)} F
	+ \sum_{j+k \leq n-1} \bm{\Gamma}^{(j)}_{(0)} \slashed{\D} \mathscr{Z}^{k} \phi
	+ \sum_{j+k+\ell \leq n-1} \bm{\Gamma}^{(j)}_{(0)} \mathscr{Y}^{k-1} F_{(\ell)}
	\end{equation*}
	This is evidently true for $N = 1$, in which case we actually have 
	\begin{equation*}
	\mathscr{Y} \phi = \begin{pmatrix} \mathscr{Z} \phi \\ r_0 \slashed{\D}_L \phi \end{pmatrix}
	\end{equation*}
	Note that we do not need to keep track of the decay of various quantities in $r$, since we are evaluating these quantities at $r = r_0$.
	
	Now, we apply one more operator $\mathscr{Y}$. Using the inductive hypothesis we have
	\begin{equation*}
	\begin{split}
	\mathscr{Y}^{N+1} \phi
	&=
	\mathscr{Y} \mathscr{Z}^{\leq N}\phi
	+ \mathscr{Y} \slashed{\D} \mathscr{Z}^{\leq N-1} \phi
	+ \mathscr{Y}^{(\leq N-1)} F
	+ \sum_{j+k \leq n-1} \bm{\Gamma}^{(j+1)}_{(0)} \slashed{\D} \mathscr{Z}^{k} \phi
	+ \sum_{j+k \leq n-1} \bm{\Gamma}^{(j)}_{(0)} \mathscr{Y} \slashed{\D} \mathscr{Z}^{k} \phi
	\\
	&\phantom{=}
	+ \sum_{j+k+\ell \leq n-1} \bm{\Gamma}^{(j+1)}_{(0)} \mathscr{Y}^{k-1} F_{(\ell)}
	+ \sum_{j+k+\ell \leq n-1} \bm{\Gamma}^{(j)}_{(0)} \mathscr{Y}^{k} F_{(\ell)}
	\end{split}
	\end{equation*}
	
	The third, fourth sixth and seventh terms on the right hand side are already of the required form. For the first term, we note that either $\mathscr{Y} = \mathscr{Z}$ or $\mathscr{Y} = r\slashed{\D}_L$, so when $r = r_0$ we have
	\begin{equation*}
	\mathscr{Y} \mathscr{Z}^{\leq N}\phi
	=
	\begin{pmatrix}
	\mathscr{Z}^{\leq N+1}\phi \\
	\slashed{\D}_L \mathscr{Z}^{\leq N}\phi
	\end{pmatrix}
	\end{equation*}
	
	For the remaining terms we treat separately the cases $\mathscr{Y} = \mathscr{Z}$ and $\mathscr{Y} = r\slashed{\D}_L$. In the former case we have
	\begin{equation*}
	\mathscr{Z} \slashed{\D} \mathscr{Z}^{\leq N-1} \phi
	=
	\slashed{\D} \mathscr{Z}^{\leq N} \phi
	+
	\bm{\Gamma}^{(1)}_{(0)} \slashed{\D} \mathscr{Z}^{\leq N-1} \phi
	\end{equation*}
	while in the latter case we make use of the computation
	\begin{equation*}
	\begin{split}
	\slashed{\D}_L \slashed{\D} \phi
	&=
	\begin{pmatrix}
	\slashed{\D}_L \slashed{\D}_L \phi \\
	\slashed{\D}_L \slashed{\D}_{\Lbar} \phi \\
	\slashed{\D}_L \slashed{\nabla} \phi \\
	\end{pmatrix}
	=
	\begin{pmatrix}
	\slashed{\D}_L \slashed{\D}_T \phi \\
	\tilde{\slashed{\Box}}_g \phi \\
	\slashed{\nabla} (r\slashed{\nabla} \phi) \\
	\slashed{\D} \phi \\
	\bm{\Gamma}^{(1)}_{(0)} \slashed{\D} \phi \\
	\slashed{\D}_L (r\slashed{\nabla} \phi) \\
	\end{pmatrix}
	=
	\begin{pmatrix}
	\slashed{\D} \mathscr{Z} \phi \\
	\slashed{\D} \phi \\
	\bm{\Gamma}^{(1)}_{(0)} \slashed{\D} \phi \\
	F
	\end{pmatrix}
	\end{split}
	\end{equation*}
	
	Combining the calculations above we find that
	\begin{equation*}
	\mathscr{Y}^{N+1} \phi
	=
	\mathscr{Z}^{\leq N+1} \phi
	+ \slashed{\D} \mathscr{Z}^{\leq N} \phi
	+ \mathscr{Y}^{(\leq N)}F
	+ \sum_{j+k \leq n} \bm{\Gamma}^{(j)}_{(0)} \slashed{\D} \mathscr{Z}^k \phi 
	+ \sum_{j+k \leq n-1} \bm{\Gamma}^{(j)}_{(0)} F_{(k)}
	+ \sum_{j+k+\ell \leq n} \bm{\Gamma}^{(j)}_{(0)} \mathscr{Y}^k F_{(\ell)}
	\end{equation*}
	which proves the inductive step.

\end{proof}

\begin{proposition}[Pointwise bound on $\mathscr{Y}^n \phi$]
	\label{proposition pointwise bound Yn phi}
	Set
	\begin{equation*}
	w^* := (1+r)^{-C^* \epsilon}
	\end{equation*}
	then for any $p < 1$ and for $r \geq r_0$ we have
	\begin{equation*}
	\begin{split}
	|\mathscr{Y}^n \phi|^2
	&\lesssim
	\frac{1}{1-p} r^{-1-p} \mathcal{E}^{(L,p)}[\mathscr{Y}^{\leq n+2} \phi](\tau)
	+ r^{-\frac{5}{2} + C^* \epsilon} \mathcal{E}^{(w^*T)}[\mathscr{Z}^{\leq n+2} \phi](\tau) 
	\\
	&\phantom{\lesssim}
	+ r^{-2}\int_{S_{\tau,r_0}} \left(
	|\mathscr{Y}^{\leq n} F|^2
	+ \sum_{j+k+\ell \leq n+1} |\bm{\Gamma}^{(j)}_{(0)}|^2 |\mathscr{Y}^{k-1} F_{(\ell)}|^2
	\right) \dVol_{\mathbb{S}^2}
	\\
	&\phantom{\lesssim}
	+ r^{-2}\left( \sup_{\substack{x \in S_{\tau,r_0} \\ j \leq \lfloor \frac{n}{2} \rfloor}} |\bm{\Gamma}^{(j)}_{(0)}|^2 \right) \int_{S_{\tau, r_0}} |\slashed{\D} \mathscr{Z}^\leq n-1 \phi|^2 \dVol_{\mathbb{S}^2}
	\\
	&\phantom{\lesssim}
	+ r^{-2}\left( \sup_{\substack{x \in S_{\tau,r_0} \\ j \leq \lfloor \frac{n}{2} \rfloor}} |\slashed{\D} \mathscr{Z}^j \phi|^2 \right) \int_{S_{\tau, r_0}} |\bm{\Gamma}^{(j)}_{(0)}|^2 \dVol_{\mathbb{S}^2} 
	\end{split}
	\end{equation*}
	
\end{proposition}

\begin{proof}
	Following the calculations in proposition \ref{proposition improved pointwise bounds in terms of p energy} (applied to the field $\mathscr{Y}^n \phi$) we find that
	\begin{equation*}
	\int_{S_{\tau,r}} |\mathscr{Y}^n \phi|^2 \dVol_{\mathbb{S}^2}
	\lesssim
	r^{-2} \int_{S_{\tau,r_0}} |\mathscr{Y}^n \phi|^2 \dVol_{\mathbb{S}^2}
	+ \frac{1}{1-p} r^{-1-p} \mathcal{E}^{(L,p)}[\mathscr{Y}^n \phi](\tau)
	\end{equation*}
	
	Next, we use proposition \ref{proposition exchanging derivatives at r0} to evaluate the boundary term at $r_0$. We have
	\begin{equation*}
	\begin{split}
	\int_{S_{\tau,r}} |\mathscr{Y}^n \phi|^2 \dVol_{\mathbb{S}^2}
	&\lesssim
	\frac{1}{1-p} r^{-1-p} \mathcal{E}^{(L,p)}[\mathscr{Y}^n \phi](\tau)
	\\
	&\phantom{\lesssim}
	+ r^{-2} \int_{S_{\tau,r_0}} \Bigg( 
		|\mathscr{Z}^{\leq n} \phi|^2 
		+ \sum_{j+k \leq n-1} |\bm{\Gamma}^{(j)}_{(0)}|^2 |\slashed{\D} \mathscr{Z}^k \phi|^2
		+ |\mathscr{Y}^{\leq n-2} F|^2
		\\
		&\phantom{\lesssim + r^{-2} \int_{S_{\tau,r_0}} \Bigg(}
		+ \sum_{j+k+\ell \leq n-1} |\bm{\Gamma}^{(j)}_{(0)}|^2 |\mathscr{Y}^{k-1} F_{(\ell)}|^2
	\Bigg) \dVol_{\mathbb{S}^2}
	\end{split}
	\end{equation*}
	We bound the third term on the right hand side by a combination of $L^2$ and $L^\infty$ estimates, where we take the term with the lowest number of derivatives in $L^\infty$. We find
	\begin{equation*}
	\begin{split}
	\int_{S_{\tau,r}} |\mathscr{Y}^n \phi|^2 \dVol_{\mathbb{S}^2}
	&\lesssim
	\frac{1}{1-p} r^{-1-p} \mathcal{E}^{(L,p)}[\mathscr{Y}^n \phi](\tau)
	+ r^{-\frac{3}{2} + C^* \epsilon} \mathcal{E}^{(w^*T)}[\mathscr{Z}^{\leq n} \phi](\tau) 
	\\
	&\phantom{\lesssim}
	+ r^{-2}\int_{S_{\tau,r_0}} \left(
	|\mathscr{Y}^{\leq n-2} F|^2
	+ \sum_{j+k+\ell \leq n-1} |\bm{\Gamma}^{(j)}_{(0)}|^2 |\mathscr{Y}^{k-1} F_{(\ell)}|^2
	\right) \dVol_{\mathbb{S}^2}
	\\
	&\phantom{\lesssim}
	+ r^{-2}\left( \sup_{\substack{x \in S_{\tau,r_0} \\ j \leq \lfloor \frac{n-2}{2} \rfloor}} |\bm{\Gamma}^{(j)}_{(0)}|^2 \right) \int_{S_{\tau, r_0}} |\slashed{\D} \mathscr{Z}^\leq n-1 \phi|^2 \dVol_{\mathbb{S}^2}
	\\
	&\phantom{\lesssim}
	+ r^{-2}\left( \sup_{\substack{x \in S_{\tau,r_0} \\ j \leq \lfloor \frac{n-2}{2} \rfloor}} |\slashed{\D} \mathscr{Z}^j \phi|^2 \right) \int_{S_{\tau, r_0}} |\bm{\Gamma}^{(j)}_{(0)}|^2 \dVol_{\mathbb{S}^2} 
	\end{split}
	\end{equation*}
	Repeating the same calculation with $\phi$ replaced by $\mathscr{Y}^{\leq 2} \phi$ and remembering that $\mathscr{Y}$ can be $r\slashed{\nabla}$, and finally using the Sobolev inequality on the spheres proves the proposition.
	
\end{proof}

\begin{corollary}[Pointwise bounds on $\overline{\slashed{\D}} \mathscr{Y}^n \phi$]
	\label{corollary Dbar Yn phi}
	Again, let
	\begin{equation*}
	w^* := (1+r)^{-C^* \epsilon}
	\end{equation*}
	then for any $p < 1$ and for $r \geq r_0$ we have
	\begin{equation*}
	\begin{split}
	|\overline{\slashed{\D}}\mathscr{Y}^n \phi|^2
	&\lesssim
	\frac{1}{1-p} r^{-3-p} \mathcal{E}^{(L,p)}[\mathscr{Y}^{\leq n+3} \phi](\tau)
	+ r^{-\frac{9}{2} + C^* \epsilon} \mathcal{E}^{(w^*T)}[\mathscr{Z}^{\leq n+3} \phi](\tau) 
	\\
	&\phantom{\lesssim}
	+ r^{-4}\int_{S_{\tau,r_0}} \left(
	|\mathscr{Y}^{\leq n+1} F|^2
	+ \sum_{j+k+\ell \leq n+2} |\bm{\Gamma}^{(j)}_{(0)}|^2 |\mathscr{Y}^{k-1} F_{(\ell)}|^2
	\right) \dVol_{\mathbb{S}^2}
	\\
	&\phantom{\lesssim}
	+ r^{-4}\left( \sup_{\substack{x \in S_{\tau,r_0} \\ j \leq \lfloor \frac{n+1}{2} \rfloor}} |\bm{\Gamma}^{(j)}_{(0)}|^2 \right) \int_{S_{\tau, r_0}} |\slashed{\D} \mathscr{Z}^\leq n-1 \phi|^2 \dVol_{\mathbb{S}^2}
	\\
	&\phantom{\lesssim}
	+ r^{-4}\left( \sup_{\substack{x \in S_{\tau,r_0} \\ j \leq \lfloor \frac{n+1}{2} \rfloor}} |\slashed{\D} \mathscr{Z}^j \phi|^2 \right) \int_{S_{\tau, r_0}} |\bm{\Gamma}^{(j)}_{(0)}|^2 \dVol_{\mathbb{S}^2} 
	\end{split}
	\end{equation*}
	
\end{corollary}

\begin{proof}
	This follows immediately \ref{proposition pointwise bound Yn phi} together with the fact that $|\overline{\slashed{\D}\phi}| \lesssim r^{-1} |\mathscr{Y} \phi|$
\end{proof}

\begin{proposition}[An estimate for $\slashed{\D}_{\Lbar} (\mathscr{Y}^n \phi)$]
	
	Suppose that $\phi$ is an $S_{\tau,r}$-tangent tensor field satisfying
	\begin{equation*}
	\tilde{\slashed{\Box}}_g \mathscr{Y}^n \phi = F_{(\mathscr{Y}, n)}
	\end{equation*}
	Suppose also that
	\begin{equation*}
	|\tilde{\Box}_g h|_{(\text{frame})}
	\lesssim
	\epsilon (1+r)^{-2 + C^* \epsilon}
	\end{equation*}
	Then we have
	\begin{equation*}
	\begin{split}
	|\slashed{\D}_{\Lbar} \mathscr{Y}^n \phi|
	&\lesssim
	\sum_{j + k \leq n} (1+r)^{jC_{(1)}\epsilon} |\slashed{\D}_{\Lbar} \mathscr{Z}^k \phi|
	+ \sum_{j + k \leq n} (1+r)^{1 + jC_{(1)}\epsilon} |F_{(\mathscr{Y}, k-1)}|
	\\
	&\phantom{\lesssim}
	+ \frac{1}{1-p} r^{-2-p + \delta} \mathcal{E}^{(L,p)}[\mathscr{Y}^{\leq n+2} \phi](\tau)
	+ r^{-1-\frac{5}{2} + \delta} \mathcal{E}^{(w^*T)}[\mathscr{Z}^{\leq n+2} \phi](\tau) 
	\\
	&\phantom{\lesssim}
	+ r^{-3 + \delta}\int_{S_{\tau,r_0}} \left(
	|\mathscr{Y}^{\leq n} F|^2
	+ \sum_{j+k+\ell \leq n+1} |\bm{\Gamma}^{(j)}_{(0)}|^2 |\mathscr{Y}^{k-1} F_{(\ell)}|^2
	\right) \dVol_{\mathbb{S}^2}
	\\
	&\phantom{\lesssim}
	+ r^{-3 + \delta}\left( \sup_{\substack{x \in S_{\tau,r_0} \\ j \leq \lfloor \frac{n}{2} \rfloor}} |\bm{\Gamma}^{(j)}_{(0)}|^2 \right) \int_{S_{\tau, r_0}} |\slashed{\D} \mathscr{Z}^\leq n-1 \phi|^2 \dVol_{\mathbb{S}^2}
	\\
	&\phantom{\lesssim}
	+ r^{-3 + \delta}\left( \sup_{\substack{x \in S_{\tau,r_0} \\ j \leq \lfloor \frac{n}{2} \rfloor}} |\slashed{\D} \mathscr{Z}^j \phi|^2 \right) \int_{S_{\tau, r_0}} |\bm{\Gamma}^{(j)}_{(0)}|^2 \dVol_{\mathbb{S}^2} 
	\end{split}
\end{equation*}
	where the final term is absent if $\phi$ is a scalar field.
	
\end{proposition}

\begin{proof}
	
	We are considering $\Lbar \mathscr{Y}^n \phi$. Recall that the $\mathscr{Y}$'s can either be $\mathscr{Z}$ or $r\slashed{\D}_L$. If all the $\mathscr{Y}$'s are $\mathscr{Z}$'s then we are done. Otherwise, there is at least one instance of the operator $r\slashed{\D}_L$, and we commute the $\mathscr{Y}$'s so that an operator $r\slashed{\D}_L$ appears on the left. Using proposition \ref{proposition commuting Z with rL} we can write
	\begin{equation*}
	\begin{split}
	\slashed{\D}_{\Lbar} \mathscr{Y}^n \phi
	&=
	\slashed{\D}_{\Lbar} \left( \mathscr{Z}^n \phi \right)
	+ \slashed{\D}_{\Lbar} \left( (r\slashed{\D}_L)\mathscr{Y}^{n-1} \phi \right)
	\\
	&\phantom{=}
	+ \slashed{\D}_{\Lbar} \bigg( 
		\bm{\Gamma}^{(0)}_{(0)} \slashed{\D}_T \mathscr{Y}^{n-2} \phi
		+\left( \bm{\Gamma}^{(1)}_{(-1+C_{(1)}\epsilon)} + r\chi_{(\text{small})} \right) \mathscr{Y}^{n-1} \phi
		+ \left( r\Omega_{L\Lbar} + r^2 \slashed{\Omega}_L \right) \mathscr{Y}^{n-2} \phi		
	\bigg)
	\\
	&=
	\slashed{\D}_{\Lbar} \left( \mathscr{Z}^n \phi \right)
	+ \slashed{\D}_{\Lbar} \left( (r\slashed{\D}_L)\mathscr{Y}^{n-1} \phi \right)
	\\
	&\phantom{=}
	+ \bm{\Gamma}^{(1)}_{(C_{(1)}\epsilon)} \slashed{\D} \mathscr{Y}^{n-2} \phi
	+ \bm{\Gamma}^{(0)}_{(0)} \slashed{\D} \mathscr{Y}^{n-1} \phi
	+ \bm{\Gamma}^{(1)}_{(-\delta)} \slashed{\D} \mathscr{Y}^{n-1} \phi
	\\
	&\phantom{=}
	+ \left( \slashed{\D}_{\Lbar} \bm{\Gamma}^{(1)}_{(-1+C_{(1)}\epsilon)} + \slashed{\D}_{\Lbar} (r\chi_{(\text{small})}) \right) \mathscr{Y}^{(n-1)} \phi
	+ \bm{\Gamma}^{(1)}_{(C_{(1)}\epsilon)} \slashed{\D} \mathscr{Y}^{n-2} \phi
	\\
	&\phantom{=}
	+ \slashed{\D}_{\Lbar} \left( r\Omega_{L\Lbar} + r^2 \slashed{\Omega}_L \right) \mathscr{Y}^{n-2} \phi
	\end{split}
	\end{equation*}

	Using the expression for the wave equation \ref{proposition tensor wave operator} we have, schematically,
	\begin{equation*}
	\begin{split}
	\slashed{\D}_{\Lbar} \left(r\slashed{\D}_L (\mathscr{Y}^{n-1}\phi) \right)
	&=
	r\tilde{\slashed{\Box}}_g (\mathscr{Y}^{n-1}\phi)
	+ (\overline{\slashed{\D}} (\mathscr{Y}^{n}\phi))
	+ (1+\bm{\Gamma}^{(0)}_{(-\delta)})(\slashed{\D}_{\Lbar}(\mathscr{Y}^{n-1}\phi))
	\\
	&\phantom{=}
	+ \bm{\Gamma}^{(1)}_{(C_{(1)}\epsilon)}(\overline{\slashed{\D}}(\mathscr{Y}^{n-1}\phi))
	+ r\Omega_{L\Lbar} (\mathscr{Y}^{n-1}\phi)
	\\
	&=
	rF_{(\mathscr{Y}, n-1)}
	+ (\overline{\slashed{\D}} (\mathscr{Y}^{n}\phi))
	+ (1+\bm{\Gamma}^{(0)}_{(-\delta)})(\slashed{\D}_{\Lbar}(\mathscr{Y}^{n-1}\phi))
	\\
	&\phantom{=}
	+ \bm{\Gamma}^{(1)}_{(C_{(1)}\epsilon)}(\overline{\slashed{\D}}(\mathscr{Y}^{n-1}\phi))
	+ \bm{\Gamma}^{(1)}_{(-1+C_{(1)}\epsilon)} (\mathscr{Y}^{n-1}\phi)
	\end{split}
	\end{equation*}
	
	Next, we note that
	\begin{equation*}
	\slashed{\D}_{\Lbar} \bm{\Gamma}^{(1)}_{(-1 + C_{(1)}\epsilon)}
	=
	\bm{\Gamma}^{(2)}_{(-1 + C_{(2)}\epsilon)}
	\end{equation*}
	which follows from writing $\slashed{\D}_{\Lbar} = 2\slashed{\D}_T - r^{-1} r\slashed{\D}_L$.
	
	We can use proposition \ref{proposition Lbar chi} together with the pointwise bootstrap bounds to bound
	\begin{equation*}
	|\slashed{\D}_{\Lbar} \chi_{(\text{small})}|
	\lesssim
	\epsilon \bm{\Gamma}^{(1)}_{(-2 + C_{(1)}\epsilon)}
	+ |\tilde{\Box}_g h|_{(\text{frame})}
	\end{equation*}
	
	Finally, using proposition \ref{proposition expression for Omega} we have
	\begin{equation*}
	\slashed{\D}_{\Lbar} \left( r\Omega_{L\Lbar} + r^2 \slashed{\Omega}_L \right)
	=
	\bm{\Gamma}^{(1)}_{(-1 + C_{(1)}\epsilon)}
	\end{equation*}
	
	Putting this all together, we have
	\begin{equation*}
	\begin{split}
	|\slashed{\D}_{\Lbar} \mathscr{Y}^n \phi|
	&\lesssim
	|\slashed{\D}_{\Lbar} \mathscr{Z}^n \phi|
	+ r|F_{(\mathscr{Y}, n-1)}|
	+ |\overline{\slashed{\D}}\mathscr{Y}^{n} \phi|
	+ |\slashed{\D}_{\Lbar} \mathscr{Y}^{n-1} \phi|
	+ \epsilon (1+r)^{C_{(1)}\epsilon} |\overline{\slashed{\D}}\mathscr{Y}^{n-1} \phi|
	\\
	&\phantom{=}
	+ \epsilon (1+r)^{-1 + C_{(2)}\epsilon} |\mathscr{Y}^{n-1} \phi|
	+ r|\tilde{\Box}_g h|_{(\text{frame})} |\mathscr{Y}^{n-1} \phi|
	+ \epsilon (1+r)^{C_{(1)}\epsilon} |\slashed{\D}_{\Lbar} \mathscr{Y}^{n-1} \phi|
	\\
	&\phantom{=}
	+ \epsilon (1+r)^{-1 + C_{(1)}\epsilon} |\mathscr{Y}^{n-2} \phi|
	\end{split}
	\end{equation*}
	
	Iterating this bound, we find that
	\begin{equation*}
	\begin{split}
	|\slashed{\D}_{\Lbar} \mathscr{Y}^n \phi|
	&\lesssim
	\sum_{j + k \leq n} (1+r)^{jC_{(1)}\epsilon} |\slashed{\D}_{\Lbar} \mathscr{Z}^k \phi|
	+ \sum_{j + k \leq n} (1+r)^{1 + jC_{(1)}\epsilon} |F_{(\mathscr{Y}, k-1)}|
	\\
	&\phantom{\lesssim}
	+ \sum_{j + k \leq n} (1+r)^{(j+1)C_{(1)} \epsilon} |\overline{\slashed{\D}} \mathscr{Y}^{k-1} \phi|
	+ \epsilon \sum_{j + k \leq n} (1+r)^{-1 + C_{(2)}\epsilon + jC_{(1)}\epsilon} |\mathscr{Y}^{k-1} \phi|
	\\
	&\phantom{\lesssim}
	+ \sum_{j + k \leq n} (1+r)^{jC_{(1)}\epsilon} |\tilde{\Box}_g h|_{(\text{frame})}|\mathscr{Y}^{k-1} \phi|
	\end{split}
	\end{equation*}
	Now, substituting the bounds from proposition \ref{proposition pointwise bound Yn phi} and from corollary \ref{corollary Dbar Yn phi}, together with the bound we have assumed on $|\tilde{\Box}_g h|_{(\text{frame})}$, and using the fact that $\delta \gg \epsilon$ proves the proposition.

\end{proof}

\section{Pointwise estimates for other geometric quantities}

\label{section pointwise bounds for geometric quantities}

We will eventually be considering solutions to wave equations in which the metric itself depends on the solution to the wave equation. Hence, the components of the metric perturbation (both in the rectangular and null frames) are assumed to satisfy suitable bounds, which will in fact follow in a simple way from the relationship between these quantities and the solutions of some wave equations. More complex relationships hold between the metric components and other geometric quantities, such as the connection coefficients and the curvature. The goal of this section is to deduce bounds on these quantities, under the assumption that the metric components and their derivatives behave similarly to solutions to the wave equation.

\begin{proposition}[Pointwise bounds on the rectangular components of the null frame]
\label{proposition pointwise bounds rectangular}
Suppose that the metric components satisfy the bounds
\begin{equation*}
\begin{split}
  |\bar{\partial} h|_{(\text{frame})} &\lesssim \sqrt{\mathcal{E}} (1+r)^{-1-2\delta}
  \\
  |h_{ab}| \big|_{r = r_0} &\lesssim \sqrt{\mathcal{E}}
\end{split}
\end{equation*}
then, for all sufficiently small $\mathcal{E}$ (compared to $\delta$), we have
\begin{equation*}
 \begin{split}
  |L^a| &\leq 1 + C_{(0)}\delta^{-1}\sqrt{\mathcal{E}} \\
  |\Lbar^a| &\leq 1 + C_{(0)}\delta^{-1}\sqrt{\mathcal{E}} \\
  |\slashed{\Pi}^a| &\leq 1 + C_{(0)}\delta^{-1}\sqrt{\mathcal{E}} \\
  |L_{(\text{small})}^i| &\leq C_{(0)}\delta^{-1}\sqrt{\mathcal{E}}r^{-2\delta}
 \end{split}
\end{equation*}
in the region $r \geq r_0$. In particular, this implies the bound
\begin{equation*}
  |X_{(\text{frame, small})}| \leq C_{(0)}\sqrt{\mathcal{E}}
\end{equation*}
\end{proposition}

\begin{proof}
We prove this proposition by a continuity argument: suppose that, for all $r_0 \leq r \leq R$, we have the bounds
\begin{equation*}
 \begin{split}
  \max\{|L^a| , |\Lbar^a| , |\slashed{\Pi}^a| \} &\leq 1 + C_{(0)}\delta^{-1}\sqrt{\mathcal{E}} \\
  L^i_{(\text{small})} & \leq C_{(0)}\delta^{-1}\sqrt{\mathcal{E}} r^{-\delta}
 \end{split}
\end{equation*}
for some sufficiently large constant $C_{(0)}$.

Then, for all $r_0 \leq r \leq R$, proposition \ref{proposition transport rectangular} shows that the rectangular components of the null frame satisfy a system of the form
\begin{equation*}
 \begin{split}
  \left| \slashed{\D}_L \begin{pmatrix} L^a \\ \Lbar^a \\ \slashed{\Pi}^a \end{pmatrix} \right|
  &\lesssim 
  \left( \sqrt{\mathcal{E}}r^{-1-2\delta} + C_{(0)}^2\delta^{-2}\mathcal{E} r^{-1-4\delta}\right)
 \end{split}
\end{equation*}
We wish to integrate these equations along the integral curves of $L$, but we first need to explain how to integrate the tensorial equation for $\slashed{\D}_L \slashed{\Pi}^a$. We have
\begin{equation*}
 L(|\slashed{\Pi}^a|^2) = L\left( \slashed{\Pi}_\mu^{\phantom{\mu}a} \slashed{\Pi}^{\mu a} \right) = 2\slashed{\Pi}^{\mu a} \slashed{\D}_L\left( \slashed{\Pi}_\mu^{\phantom{\mu}a} \right) 
\end{equation*}
and so
\begin{equation*}
  |\slashed{\Pi}^a|^2(\tau, r, \vartheta^1, \vartheta^2)
  \leq
  |\slashed{\Pi}^a|^2(\tau, r_0, \vartheta^1, \vartheta^2)
  + \int_{r_0}^r 2\slashed{\Pi}^{\mu a} \slashed{\D}_L\left( \slashed{\Pi}_\mu^{\phantom{\mu}a} \right) (\tau, r', \vartheta^1, \vartheta^2) \upd r'
\end{equation*}
Also, we noting (using proposition \ref{proposition rectangular components in r < r0}) that, at $r = r_0$, we have
\begin{equation*}
 \begin{split}
  |L^a| &\leq 1 + \sqrt{\mathcal{E}}\\
  |\Lbar^a| &\leq 1 + \sqrt{\mathcal{E}}\\
  |\slashed{\Pi}^a| &\leq 1 + \sqrt{\mathcal{E}}\\
 \end{split}
\end{equation*}
Hence, integrating these equations along the integral curves of $L$ we find that\footnote{Note a bound of the form $|\slashed{\Pi}^a|^2 \leq 1 + \sqrt{\mathcal{E}}$, together with the fact that $\sqrt{\mathcal{E}} \ll 1$, implies the bound $|\slashed{\Pi}^a| \leq 1 + \sqrt{\mathcal{E}}$}
\begin{equation*}
 \max\{|L^a| , |\Lbar^a| , |\slashed{\Pi}^a| \} \leq 1 + \sqrt{\mathcal{E}} + C\delta^{-1}\sqrt{\mathcal{E}} + C_{(0)}^2\delta^{-3}\mathcal{E}
\end{equation*}
where $C$ is a constant related to the implicit constants in the transport equations. Hence, since $\sqrt{\mathcal{E}} \ll \delta \ll 1$, if we choose $C_{(0)}$ sufficiently large then we can improved the bootstrap bounds, and so the bounds on $|L^a|$, $|\Lbar^a|$ and $|\slashed{\Pi}^a|$ actually hold for all $r_0 \leq r \leq R + \tilde{\epsilon}$, for some $\tilde{\epsilon} > 0$.

Now, we examine the transport equation for $rL^i_{(\text{small})}$, which is also presented in proposition \ref{proposition transport rectangular}. We find that, for $r \leq R$ we have
\begin{equation*}
 \left| L\left( rL^i_{(\text{small})} \right) \right|
 \lesssim 
 \left( \epsilon \sqrt{\mathcal{E}} r^{-3\delta} + \left(\sqrt{\mathcal{E}} + C_{(0)}^2 \delta^{-2} \mathcal{E} \right) r^{-2\delta} \right)
\end{equation*}
and so, integrating this from $r = r_0$ (where $|L^i_{(\text{small})}| \lesssim \sqrt{\mathcal{E}}$ - see proposition \ref{proposition rectangular components in r < r0}) and using the relationships between the various constants, we find that, in the region $r \geq r_0$,
\begin{equation*}
 |L^i_{(\text{small})}| \lesssim \delta^{-1}\sqrt{\mathcal{E}} r^{-2\delta}
\end{equation*}
Now, this is an improvement over the bound assumed for $L^i_{(\text{small})}$, if $C_{(0)}$ is sufficiently large.

Hence, all of the bounds which were assumed to hold up to $r = R$ actually hold for slightly larger values of $r$, and so, by continuity, they hold for all values of $r$, proving the proposition.

\end{proof}

\begin{proposition}[A pointwise bound on the commuted rectangular components of the null frame fields]
	\label{proposition pointwise bounds Yn Xframe}
	Suppose that the pointwise bounds of chapter \ref{chapter bootstrap} hold. Additionally, we suppose that certain \emph{improved} versions of the pointwise bootstrap bounds of chapter \ref{chapter bootstrap} hold, in the following senses:
	\begin{itemize}
		\item when bounding a term involving $m$ commutation operators $\mathscr{Y}$, where $m \leq n$ we will assume the bounds
		\begin{equation*}
		\begin{split}
			|\bm{\Gamma}^{(m-1)}_{(-1+C_{(m-1)}\epsilon)}| &\lesssim \sqrt{\mathcal{E}} (1+r)^{-1+C_{(m-1)}\epsilon}
			\\
			|\bm{\Gamma}^{(m-1)}_{(-1-\delta)}| &\lesssim \sqrt{\mathcal{E}} (1+r)^{-1-(2-c_{[m]})\delta}
		\end{split}
		\end{equation*}
		where the constants $c_{[m]}$ are such that $C_{[m]} \gg c_{[m]} \gg C_{(m)}$.
		\item when bounding terms that directly involve the metric perturbation $h$ or its derivatives, we assume the bounds
		\begin{equation*}
		\begin{split}
			|\mathscr{Y}^n h_{(\text{rect})}| &\lesssim \sqrt{\mathcal{E}} (1+r)^{-2\delta}
			\\
			|\slashed{\D} \mathscr{Y}^n h_{(\text{rect})}| &\lesssim \sqrt{\mathcal{E}} (1+r)^{-1+\delta}
			\\
			|\slashed{\D} \mathscr{Y}^n h|_{(\text{frame})} &\lesssim \sqrt{\mathcal{E}} (1+r)^{-1+\delta}
			\\
			|\overline{\slashed{\D}} \mathscr{Y}^n h_{(\text{rect})}| &\lesssim \sqrt{\mathcal{E}} (1+r)^{-1-2\delta}
			\\
			|\overline{\slashed{\D}} \mathscr{Y}^n h|_{(\text{frame})} &\lesssim \sqrt{\mathcal{E}} (1+r)^{-1-2\delta}
		\end{split}
		\end{equation*}
	\end{itemize}

	Finally, suppose that
	\begin{equation*}
	\sqrt{\mathcal{E}} \ll \epsilon
	\end{equation*}
	
	Then, for sufficiently small $\mathcal{E}$ and $\epsilon$ we have
	\begin{equation}
	\begin{split}
	|\mathscr{Y}^n X_{(\text{frame})}| 
	&\lesssim (1+r)^{C_{(n-1)} \epsilon} \\
	\\
	|\mathscr{Y}^n \bar{X}_{(\text{frame})}|
	&\lesssim \delta^{-1}\sqrt{\mathcal{E}}(1+r)^{-2\delta}
	\end{split}
	\end{equation}
\end{proposition}

\begin{proof}
	Recall proposition \ref{proposition transport Yn rectangular}. Combining this with the bounds we have assumed above, in the region $r \geq r_0$ we have the following bounds:
	\begin{equation*}
	\begin{split}
	|\slashed{\D}_L \mathscr{Y}^n X_{(\text{frame})}|
	&\lesssim
	\sqrt{\mathcal{E}} r^{-1} |\mathscr{Y}^n X_{(\text{frame})}|
	+ \sqrt{\mathcal{E}} r^{-1+C_{(n-1)}\epsilon}
	+ r^{-2} \left| r\mathscr{Y}^n \bar{X}_{(\text{frame})}\right|
	\\
	\left|\slashed{\D}_L \left( r\mathscr{Y}^n \bar{X}_{(\text{frame})} \right) \right|
	&\lesssim
	\sqrt{\mathcal{E}} r^{-1} \left|r\mathscr{Y}^n \bar{X}_{(\text{frame})} \right|
	+ \sqrt{\mathcal{E}} r^{-2\delta}
	\end{split}
	\end{equation*}
	
	Using the expressions satisfied by the derivatives of the frame fields given in propositions \ref{proposition transport rectangular}, \ref{proposition transport lbar rectangular} and \ref{proposition angular rectangular} we have the bounds
	\begin{equation*}
	\begin{split}
	\sup_{S_{\tau, r_0}} |\mathscr{Y}^n X_{(\text{frame})}| &\lesssim 1 \\
	\sup_{S_{\tau, r_0}} |r\mathscr{Y}^n \bar{X}_{(\text{frame})}| &\lesssim \sqrt{\mathcal{E}}
	\end{split}
	\end{equation*}
	
	So, using the Gronwall inequality, we find that, for $r \geq r_0$
	\begin{equation*}
	|r\mathscr{Y}^n \bar{X}_{(\text{frame})}|
	\lesssim 
	\delta^{-1}\sqrt{\mathcal{E}}(1+r)^{1-2\delta}
	\end{equation*}
	and so
	\begin{equation*}
	|\mathscr{Y}^n \bar{X}_{(\text{frame})}|
	\lesssim 
	\delta^{-1}\sqrt{\mathcal{E}}(1+r)^{-2\delta}
	\end{equation*}
	Now, subsituting this in to the bound for $|\slashed{\D}_L \mathscr{Y}^n X_{(\text{frame})}|$ we obtain
	\begin{equation*}
	|\slashed{\D}_L\mathscr{Y}^n X_{(\text{frame})}|
	\lesssim
	\sqrt{\mathcal{E}} r^{-1} |\mathscr{Y}^n X_{(\text{frame})}|
	+ \sqrt{\mathcal{E}} r^{-1+C_{(n-1)}\epsilon}
	+ \delta^{-1}\sqrt{\mathcal{E}} r^{-1-\delta}
	\end{equation*}
	So, again, the Gronwall inequality gives us
	\begin{equation*}
	|\mathscr{Y}^n X_{(\text{frame})}|
	\lesssim
	\left(1 + \frac{\sqrt{\mathcal{E}}}{C_{(n-1)}\epsilon} + \frac{\sqrt{\mathcal{E}}}{\delta^2} \right) (1+r)^{C_{(n-1)}\epsilon}
	\end{equation*}
	assuming that $\sqrt{\mathcal{E}} \ll C_{(n-1)}\epsilon$. Again, using this assumption together with the fact that $\epsilon \ll \delta$ we actually have
	\begin{equation*}
	|\mathscr{Y}^n X_{(\text{frame})}|
	\lesssim
	(1+r)^{C_{(n-1)}\epsilon}
	\end{equation*}
	
\end{proof}

\begin{proposition}[Pointwise bounds on $\mathscr{Y}^n X_{(\text{frame, small})}$]
	\label{proposition pointwise bound rectangular small}
	Suppose that the same conditions hold as in proposition \ref{proposition pointwise bounds Yn Xframe}. Then we have
	\begin{equation}
	|\mathscr{Y}^n X_{(\text{frame, small})}| 
	\lesssim
	\frac{\sqrt{\mathcal{E}}}{C_{(n-1)}\epsilon} (1+r)^{C_{(n-1)}\epsilon}
	\end{equation}
%	We \emph{also} have
%	\begin{equation}
%	|\mathscr{Y}^n X_{(\text{frame, small})}| 
%	\leq \delta^{-1}\sqrt{\mathcal{E}}\left( (1+r)^{\frac{1}{2}\delta + C_{(n)}\sqrt{\mathcal{E}}}(1+\tau)^{-\beta} 
%	+ (1+r)^{-2\delta + C_{(n)}\sqrt{\mathcal{E}}} \right)
%	\end{equation}
\end{proposition}

\begin{proof}
	The proof of this proposition is very similar to the proof of the previous proposition. From proposition \ref{proposition transport Yn rectangular small}, in the region $r \geq r_0$ we have
	\begin{equation*}
	|\slashed{\D}_L\mathscr{Y}^n X_{(\text{frame, small})}|
	\lesssim
	\sqrt{\mathcal{E}} r^{-1} |\mathscr{Y}^n X_{(\text{frame, small})}|
	+ \sqrt{\mathcal{E}} r^{-1+C_{(n-1)}\epsilon}
	+ \delta^{-1}\sqrt{\mathcal{E}} r^{-1-\delta}
	\end{equation*}
	This time, the initial data for $\mathscr{Y}^n X_{(\text{frame, small})}$ is bounded by $\sqrt{\mathcal{E}}$ (see proposition \ref{proposition rectangular components in r < r0}). Hence, using the Gronwall inequality we can prove that, in fact,
	\begin{equation*}
	|\mathscr{Y}^n X_{(\text{frame, small})}|
	\lesssim
	\left(\frac{\sqrt{\mathcal{E}}}{C_{(n-1)}\epsilon} + \frac{\sqrt{\mathcal{E}}}{\delta^2} \right) (1+r)^{C_{(n-1)}\epsilon}
	\end{equation*}
	which, \emph{a fortiori} proves the proposition.
\end{proof}

\begin{proposition}[A pointwise bound for $\mu$]
	\label{proposition pointwise bound mu}
	Suppose that 
	\begin{equation*}
	\begin{split}
	|h_{(\text{rect})}| &\lesssim \sqrt{\mathcal{E}} \\
	|\bar{\partial}h|_{(\text{frame})} &\lesssim \sqrt{\mathcal{E}}(1+r)^{-1-2\delta} \\
	|\partial h|_{LL} &\lesssim \sqrt{\mathcal{E}}(1+r)^{-1}
	\end{split}
	\end{equation*}
	
	Then, for all sufficiently small $\mathcal{E}$, there is some numerical constant $C$ such that we have
	\begin{equation}
	\begin{split}
	|\mu| &\leq 2(1+r)^{C\sqrt{\mathcal{E}}} \\
	|\mu^{-1}| &\leq 2(1+r)^{C\sqrt{\mathcal{E}}} \\
	\end{split}
	\end{equation}
\end{proposition}

\begin{proof}
	Recall proposition \ref{proposition transport mu}. Combining this with the result of proposition \ref{proposition pointwise bounds rectangular} we have
	\begin{equation*}
	\begin{split}
	|L\log \mu| 
	&\lesssim
	r^{-1} \left| \frac{x^i}{r} \right| |\bar{X}_{(\text{small})}|
	+ r^{-1} |\bar{X}_{(\text{small})}|^2
	+ |\bar{\partial} h|_{(\text{frame})}
	+ |\partial h|_{LL}
	\\
	&\lesssim
	\sqrt{\mathcal{E}} \left( \delta^{-1}(1+r)^{-1-2\delta} + (1+r)^{-1} \right)
	\end{split}
	\end{equation*}
	Integrating from $r=r_0$ (where $\mu = 1 + \mathcal{O}(|h_{(\text{rect})}|) = 1 + \mathcal{O}(\sqrt{\mathcal{E}})$ - see proposition \ref{proposition initial data for mu}), we have
	\begin{equation*}
	|\log \mu| \lesssim \sqrt{\mathcal{E}} \left( \delta^{-2}(1+r)^{-2\delta} + \log (1+r) \right)
	\end{equation*}
	from which it follows that
	\begin{equation*}
	\left.\begin{aligned} |\mu| \\ |\mu^{-1}|
	\end{aligned}\right\rbrace 
	\leq \exp(C\sqrt{\mathcal{E}}\delta^{-2})(1+r)^{C\sqrt{\mathcal{E}}}
	\leq 2(1+r)^{C\sqrt{\mathcal{E}}}
	\end{equation*}
	where $C$ is some numerical constant, and where we have chosen $\mathcal{E}$ sufficiently small relative to $C \delta^{-2}$.
\end{proof}

\begin{proposition}[A pointwise bound for $\mathscr{Y}^n \mu$]
	\label{proposition pointwise bound Yn mu}
	Suppose that the same conditions hold as in proposition \ref{proposition pointwise bounds Yn Xframe}. Then we have the bounds
	\begin{equation*}
	|\mathscr{Y}^n \log \mu|
	\lesssim
	\delta^{-2}\sqrt{\mathcal{E}} (1+r)^{\delta}
	\end{equation*}
	
	On the other hand, if we have
	\begin{equation*}
	|\slashed{\D} \mathscr{Y}^n h|_{LL} \lesssim \sqrt{\mathcal{E}} (1+r)^{-1 + C_{(n)}\epsilon}
	\end{equation*}
	then
	\begin{equation*}
	|\mathscr{Y}^n \log \mu|
	\lesssim
	C_{(n)}^{-1} \epsilon^{-1} \sqrt{\mathcal{E}} (1+r)^{C_{(n)}\epsilon}
	\end{equation*}
	
\end{proposition}

\begin{proof}
	Recall proposition \ref{proposition transport Yn mu}. Substituting the first pointwise bounds we have assumed, along with the bounds for the frame components we have derived above, we obtain
	\begin{equation*}
	\begin{split}
	|\slashed{\D}_L \left(\mathscr{Y}^n \log \mu \right)|
	&\lesssim
	\sqrt{\mathcal{E}} (1+r)^{-1} (\mathscr{Y}^n\log \mu)
	+ \delta^{-1}\sqrt{\mathcal{E}} (1+r)^{-1 - \delta}
	+ \sqrt{\mathcal{E}} (1+r)^{-1 + \delta}
	\end{split}
	\end{equation*}
	where we have used the fact that $\delta \gg \epsilon$. Making use of Gronwall's inequality, together with the fact that $|\mathscr{Y}^n \log \mu| \lesssim \sqrt{\mathcal{E}}$ at $r = r_0$ (see proposition \ref{proposition initial data for mu}), we obtain the bound
	\begin{equation*}
	|\mathscr{Y}^n \log \mu|
	\lesssim
	\delta^{-2}\sqrt{\mathcal{E}}(1+r)^{\delta}
	\end{equation*}

	On the other hand, making use of the alternative bound
	\begin{equation*}
	|\slashed{\D} \mathscr{Y}^n h|_{(\text{frame})}
	\lesssim
	\mathcal{E}(1+r)^{-1+C_{(n)}\epsilon}
	\end{equation*}
	we instead obtain the bound
	\begin{equation*}
	|\mathscr{Y}^n \log \mu|
	\lesssim
	\left(\delta^{-2}\sqrt{\mathcal{E}} + C_{(n)}^{-1} \epsilon^{-1} \sqrt{\mathcal{E}} \right)(1+r)^{C_{(n)}\epsilon}
	\end{equation*}
	Making use of the fact that $\epsilon \ll \delta$ proves the second part of the proposition.
\end{proof}

\begin{proposition}[A pointwise bound on $\omega$]
	\label{proposition pointwise bound omega} 
 Suppose that the metric perturbations satisfy
\begin{equation*}
 \begin{split}
  |h_{(\text{rect})}| &\lesssim \sqrt{\mathcal{E}} \\
  |\bar{\partial} h|_{(\text{frame})} &\lesssim \sqrt{\mathcal{E}} (1+r)^{-1-2\delta} \\
  |\partial h|_{LL} &\lesssim \sqrt{\mathcal{E}}
 \end{split}
\end{equation*}
Then the connection coefficient $\omega$ satisfies
\begin{equation}
 \omega \lesssim \delta^{-1}\sqrt{\mathcal{E}}(1+r)^{-1}
\end{equation}

%Additionally, suppose that the metric perturbations satisfy
%\begin{equation*}
% |\partial \mathscr{Z}^n h| \lesssim \sqrt{\mathcal{E}} (1+r)^{-1 + C_{(n)}\epsilon}
%\end{equation*}
%then $\omega$ satisfies
%\begin{equation*}
% \mathscr{Z}^n\omega \lesssim \sqrt{\mathcal{E}}(1+r)^{-1 + C_{(n)}\epsilon}
%\end{equation*}

\end{proposition}

\begin{proof}
 This follows directly from definition \ref{definition omega}, and the bounds we have already obtained in proposition \ref{proposition pointwise bounds Yn Xframe}.
\end{proof}

\begin{proposition}[A pointwise bound for $\mathscr{Y}^n \omega$]
	\label{proposition pointwise bound Yn omega}
	Suppose that the same conditions hold as in proposition \ref{proposition pointwise bounds Yn Xframe}. Then we have the bound
	\begin{equation*}
	|\mathscr{Y}^n \omega|
	\lesssim
	\sqrt{\mathcal{E}} (1+r)^{-1 + \delta}
	\end{equation*}
	
	On the other hand, if we also have the bound
	\begin{equation*}
	|\slashed{\D}\mathscr{Y}^n h|_{LL}
	\lesssim
	\sqrt{\mathcal{E}} (1+r)^{-1+C_{(n)}\epsilon}
	\end{equation*}
	then  we have the bound
	\begin{equation*}
	|\mathscr{Y}^n \omega|
	\lesssim
	\sqrt{\mathcal{E}} (1+r)^{-1 + C_{(n)}\epsilon}
	\end{equation*}
	
\end{proposition}

\begin{proof} 
	Recall proposition \ref{proposition Yn omega}. This proposition follows immediately from the expression given in proposition \ref{proposition Yn omega} together with the bootstrap bounds we have assumed, the bounds in proposition \ref{proposition pointwise bounds Yn Xframe} and the fact that $\sqrt{\mathcal{E}} \ll \epsilon \ll 1$.
\end{proof}

\begin{proposition}[A pointwise bound on $\zeta$]
	\label{proposition pointwise bound zeta}
 Suppose that the metric perturbations satisfy
\begin{equation*}
 \begin{split}
  |\partial h|_{(\text{frame})} 
  &\lesssim
  \sqrt{\mathcal{E}} (1+r)^{-1+ \delta}
  \\
  |\bar{\partial}h|_{(\text{frame})} 
  &\lesssim
  \sqrt{\mathcal{E}} (1+r)^{-1-2\delta}
 \end{split}
\end{equation*}
and suppose that $\mathcal{E}$ is sufficiently small.

Then the connection coefficient $\zeta$ satisfies
\begin{equation}
 |\zeta| \lesssim
 \sqrt{\mathcal{E}} (1+r)^{-1+\delta}
\end{equation}
in the region $r \geq r_0$.

On the other hand, if we have the bound
\begin{equation*}
	|\partial h|_{(\text{frame})} 
	\lesssim
	\sqrt{\mathcal{E}} (1+r)^{-1+C_{(0)}\epsilon}
\end{equation*}
Then the connection coefficient $\zeta$ satisfies
\begin{equation}
	|\zeta|
	\lesssim
	\delta^{-1}\sqrt{\mathcal{E}} (1+r)^{-1+ \epsilon}
\end{equation}
in the region $r \geq r_0$.

\end{proposition}

\begin{proof}
From proposition \ref{proposition zeta} we have
\begin{equation*}
  |\zeta| \lesssim |\partial h|_{(\text{frame})} + r^{-1}|L^i_{(\text{small})}| |\slashed{\Pi}^i|
\end{equation*}
and from proposition \ref{proposition pointwise bounds rectangular} we have
\begin{equation*}
 \begin{split}
  |L^i_{(\text{small})}| &\lesssim \delta^{-1}\sqrt{\mathcal{E}}r^{-\delta} \\
  |\slashed{\Pi}^i| &\leq 2
 \end{split}
\end{equation*}
assuming $\sqrt{\mathcal{E}}$ is small enough. Hence, using the bounds assumed on the derivatives of $h$ we can prove the proposition.
\end{proof}

\begin{proposition}[A pointwise bound on $\mathscr{Y}^n \zeta$]
	\label{proposition pointwise bound Yn zeta}
	Suppose that the same conditions hold as in proposition \ref{proposition pointwise bounds Yn Xframe}. Then we have the bound
	\begin{equation*}
	|\mathscr{Y}^n \zeta|
	\lesssim
	\sqrt{\mathcal{E}} (1+r)^{-1 + \delta}
	\end{equation*}
	
	On the other hand, if we also have the bound
	\begin{equation*}
	|\slashed{\D} \mathscr{Y}^n h|_{(\text{frame})}
	\lesssim
	\sqrt{\mathcal{E}} (1+r)^{-1 + C_{(n)}\epsilon}
	\end{equation*}
	then we have the bound
	\begin{equation*}
	|\mathscr{Z}^n \zeta|
	\lesssim
	\sqrt{\mathcal{E}} (1+r)^{-1 + C_{(n)}\epsilon}
	\end{equation*}
	
\end{proposition}

\begin{proof}
	Recall proposition \ref{proposition Yn zeta}. Substituting in the pointwise bounds we have assumed, together with the bounds on the rectangular components of the null frame obtained in proposition \ref{proposition pointwise bounds Yn Xframe}, and using the fact that $\epsilon \ll \delta$ we obtain the results of the proposition.
\end{proof}

\begin{proposition}[A pointwise bound for $\tr_{\slashed{g}}\chi_{(\text{small})}$]
	\label{proposition pointwise bound trchi small}
	Suppose that the bootstrap bounds hold. Suppose in addition that the metric perturbations satisfy
	\begin{equation*}
	\begin{split}
	|\partial h|_{(\text{frame})} 
	&\lesssim
	\sqrt{\mathcal{E}} (1+r)^{-1+ \delta}
	\\
	|\bar{\partial}h|_{(\text{frame})} 
	&\lesssim 
	\sqrt{\mathcal{E}} (1+r)^{-1-3\delta} 
	\\
	|\overline{\slashed{\D}} \mathscr{Z} h|_{(\text{frame})} 
	&\lesssim 
	\sqrt{\mathcal{E}} (1+r)^{-1-2\delta} \\
	\end{split}
	\end{equation*}
	
	Then we have
	\begin{equation}
	|\mathcal{X}_{(\text{low})}| \lesssim \delta^{-1}\sqrt{\mathcal{E}} (1+r)^{-1-2\delta}
	\end{equation}
	in the region $r \geq r_0$.
	
	Additionally, we have
	\begin{equation}
	|\tr_{\slashed{g}}\chi_{(\text{small})}| \lesssim \delta^{-1}\sqrt{\mathcal{E}} (1+r)^{-1 - 2\delta}
	\end{equation}
		
\end{proposition}

\begin{proof}
	Substituting the bounds we have assumed above into the transport equation for $r^2 \mathcal{X}_{(\text{low})}$ given in proposition \ref{proposition transport trace chi} we find
	\begin{equation*}
	\left|\D_L \left(r^2 \mathcal{X}_{(\text{low})} \right) \right|
	\lesssim
	r^{-2} \left|r^2 \mathcal{X}_{(\text{low})} \right|
	+ \sqrt{\mathcal{E}} (1+r)^{-2\delta}
	\end{equation*}
	Now, since at $r = r_0$ we can express $\chi$ in terms of the rectangular components of $h$ and their derivatives (see proposition \ref{proposition chi in r leq r0}), the initial data satisfies
	\begin{equation*}
	\left|r^2 \mathcal{X}_{(\text{low})}\right| \big|_{r = r_0} \lesssim \sqrt{\mathcal{E}}
	\end{equation*}
	Hence, using Gronwall's inequality, we find
	\begin{equation*}
	\left|r^2 \mathcal{X}_{(\text{low})}\right|
	\lesssim
	\delta^{-1}\sqrt{\mathcal{E}} (1+r)^{1-2\delta}
	\end{equation*}
	which easily proves the first part of the proposition. To prove the second part, we note that, from proposition \ref{proposition transport trace chi} we have
	\begin{equation*}
	|\mathcal{X}_{(\text{low})} - \tr_{\slashed{g}}\chi_{(\text{small})}|
	\lesssim
	|\bar{\partial}h|_{(\text{frame})}
	\end{equation*}
	Using the bounds we have assumed for the derivatives of $h$ and the triangle inequality we can prove the second part of the proposition.
\end{proof}

\begin{proposition}[A pointwise bound for $\mathscr{Y}^n \tr_{\slashed{g}}\chi_{(\text{small})}$]
	\label{proposition pointwise bound Yn tr chi}
	Suppose that the same conditions hold as in proposition \ref{proposition pointwise bounds Yn Xframe}. Suppose in addition that the ``good'' derivatives of the metric components satisfy the following bounds, which include $n+1$ commutation operators:
	\begin{equation*}
	|\overline{\slashed{\D}} \mathscr{Y}^{n+1} h|_{(\text{frame})} \lesssim \sqrt{\mathcal{E}}(1+r)^{-1-2\delta}
	\end{equation*}

	Then we have the bounds
	\begin{equation}
	\begin{split}
	|\mathscr{Y}^n \mathcal{X}_{(\text{low})}| 
	&\lesssim
	\delta^{-1}\sqrt{\mathcal{E}} (1+r)^{-1-2\delta}
	\\
	|\mathscr{Y}^n \tr_{\slashed{g}}\chi_{(\text{small})}|
	&\lesssim
	\delta^{-1}\sqrt{\mathcal{E}} (1+r)^{-1-2\delta + c_{[n-1]}\epsilon}
	\end{split}
	\end{equation}
	
\end{proposition}

\begin{proof}
	Recall proposition \ref{proposition Yn tr chi low}. Substituting the bounds assumed in the proposition, we find
	\begin{equation*}
	\left|\slashed{\D}_L \left(r^2 \mathscr{Y}^n \mathcal{X}_{(\text{low})} \right) \right|
	\lesssim
	\sqrt{\mathcal{E}} (1+r)^{-1} \left|r^2 \mathscr{Y}^n \mathcal{X}_{(\text{low})} \right|
	+ \delta^{-1}\sqrt{\mathcal{E}}(1+r)^{-2\delta}
	\end{equation*}
	Next we use the Gronwall inequality, together with the bounds on the initial data
	\begin{equation*}
	|\mathscr{Y}^n \mathcal{X}_{(\text{low})}| \big|_{r = r_0} \lesssim \sqrt{\mathcal{E}}
	\end{equation*}
	which follows from the fact that $\mathcal{X}_{(\text{low})}$ can be expressed in terms of the rectangular derivatives of $h$ at $r = r_0$.
	
	Next we note that, by substituting the bounds assumed in the proposition into the expression for $\mathscr{Z}^n \mathcal{X}_{(\text{low})}$ given in proposition \ref{proposition Yn tr chi low} we obtain
	\begin{equation*}
	\left|\mathscr{Y}^n \mathcal{X}_{(\text{low})} - \mathscr{Y}^n \tr_{\slashed{g}}\chi_{(\text{small})}\right|
	\lesssim
	\delta^{-1}\sqrt{\mathcal{E}}(1+r)^{-1-2\delta + C_{(n-1)}\epsilon}
	\end{equation*}
	Now, the triangle inequality gives the required bound on $\mathscr{Y}^n \tr_{\slashed{g}}\chi_{(\text{small})}$.
\end{proof}

\begin{proposition}[A pointwise bound for $\hat{\chi}$]
	\label{proposition pointwise bound chihat}
	Suppose that the pointwise bootstrap bounds of chapter \ref{chapter bootstrap} hold. Additionally, suppose that the derivatives of the metric satisfy
	\begin{equation*}
	\begin{split}
	|\partial h|_{(\text{frame})} &\lesssim \sqrt{\mathcal{E}} (1+r)^{-1+\delta}
	\\
	|\partial h|_{LL} &\lesssim \sqrt{\mathcal{E}}(1+r)^{-1}
	\\
	|\bar{\partial} h|_{(\text{frame})} &\lesssim \sqrt{\mathcal{E}}(1+r)^{-1-3\delta}
	\\
	|\overline{\slashed{\D}} \mathscr{Z} h|_{(\text{frame})} &\lesssim \sqrt{\mathcal{E}}(1+r)^{-1-2\delta}
	\end{split}
	\end{equation*}
	
	Then, for sufficiently small $\mathcal{E}$, in the region $r \geq r_0$, we have
	\begin{equation}
	|\hat{\chi}| \lesssim \sqrt{\mathcal{E}} (1+r)^{-1 - 2\delta}
	\end{equation}
	
\end{proposition}

\begin{proof}
	Recall propositions \ref{proposition transport chihat} and \ref{proposition structure of alpha}. Multiplying by the equation given in proposition \ref{proposition transport chihat} by $r^2$, and substituting the expression for $\alpha$ given in \ref{proposition structure of alpha}, we find
	\begin{equation*}
	|\slashed{\D}_L \left( r^2 \hat{\chi} \right)|
	\lesssim
	|\tr_{\slashed{g}}\chi_{(\text{small})}| |r^2 \hat{\chi}|
	+ r^{-2} |r^2 \hat{\chi}|^2
	+ r|\overline{\slashed{\D}}\mathscr{Z} h|_{(\text{frame})}
	+ |\partial h|_{LL} |r^2\hat{\chi}|
	+ r^2 |\bar{\partial} h|_{(\text{frame})} \bm{\Gamma}^{(0)}_{(-1 + C_{(0)}\epsilon)}
	\end{equation*}
	Substituting in the bounds we have assumed, we find, for $r \geq r_0$,
	\begin{equation*}
	|\slashed{\D}_L \left( r^2 \hat{\chi} \right)|
	\lesssim
	\epsilon(1+r)^{-1-\delta} |r^2 \hat{\chi}|
	+ \sqrt{\mathcal{E}}(1+r)^{-1}|r^2\hat{\chi}|
	+ \sqrt{\mathcal{E}} (1+r)^{-2\delta}
	\end{equation*}
	We also have that
	\begin{equation*}
	|\hat{\chi}| \big|_{r = r_0} \lesssim \sqrt{\mathcal{E}}
	\end{equation*}
	since $\hat{\chi}$ can be expressed in terms of the derivatives of $h$ at $r = r_0$. Hence, applying Gronwall's inequality, we have
	\begin{equation*}
	|r^2 \hat{\chi}|^2 
	\lesssim
	\sqrt{\mathcal{E}}(1+r)^{1-2\delta}
	\end{equation*}
	for sufficiently small $\mathcal{E}$, and for $r \geq r_0$.
	
\end{proof}

\begin{proposition}[A pointwise bound for $\mathscr{Y}^n \hat{\chi}$]
	\label{proposition pointwise bound Yn chihat}
	Suppose that the same conditions hold as in proposition \ref{proposition pointwise bounds Yn Xframe}. Suppose in addition that the ``good'' derivatives of the metric components satisfy the following bounds, which include $n+1$ commutation operators:
	\begin{equation*}
	|\overline{\slashed{\D}} \mathscr{Y}^{n+1} h|_{(\text{frame})} \lesssim \sqrt{\mathcal{E}}(1+r)^{-1-2\delta}
	\end{equation*}

	Then we have the bound
	\begin{equation}
	|\mathscr{Y}^n \hat{\chi}| \lesssim \delta^{-1}\sqrt{\mathcal{E}} (1+r)^{-1-2\delta+c_{[n-1]}\epsilon}
	\end{equation}
	
\end{proposition}

\begin{proof}

Recall proposition \ref{proposition transport Yn chihat}. Substituting the bounds assumed in the proposition, we find
\begin{equation*}
\left|\slashed{\D}_L \left(r^2 \mathscr{Y}^n \hat{\chi} \right) \right|
\lesssim
\sqrt{\mathcal{E}} (1+r)^{-1} \left|r^2 \mathscr{Y}^n \hat{\chi} \right|
+ \sqrt{\mathcal{E}}(1+r)^{-2\delta}
\end{equation*}
Next we use the Gronwall inequality, together with the bounds on the initial data
\begin{equation*}
|\mathscr{Y}^n \mathcal{X}_{(\text{low})}| \big|_{r = r_0} \lesssim \sqrt{\mathcal{E}}
\end{equation*}
which follows from the fact that $\mathcal{X}_{(\text{low})}$ can be expressed in terms of the rectangular derivatives of $h$ at $r = r_0$ (see proposition \ref{proposition chi in r leq r0}).

\end{proof}

\begin{proposition}[A pointwise bound for $\tr_{\slashed{g}}\chibar_{(\text{small})}$]
	\label{proposition pointwise bound tr chibar}	
	Suppose that the bootstrap bounds hold. Suppose in addition that the metric perturbations satisfy
	\begin{equation*}
	\begin{split}
	|h_{(\text{rect})}| &\lesssim \sqrt{\mathcal{E}}(1+r)^{-2\delta} \\
	|\partial h|_{(\text{frame})} &\lesssim
	\sqrt{\mathcal{E}} (1+r)^{-1+\delta} \\
	\\
	|\bar{\partial}h|_{(\text{frame})} &\lesssim \sqrt{\mathcal{E}} (1+r)^{-1-3\delta} \\
	|\overline{\slashed{\D}} \mathscr{Z} h|_{(\text{frame})} &\lesssim \sqrt{\mathcal{E}} (1+r)^{-1-2\delta} \\
	\end{split}
	\end{equation*}
	Then we have
	\begin{equation*}
	|\tr_{\slashed{g}}\chibar_{(\text{small})}|
	\lesssim
	\delta^{-1}\sqrt{\mathcal{E}}(1+r)^{-1 + \delta}
	\end{equation*}
	
	On the other hand, if we assume the bound
	\begin{equation*}
	|\partial h|_{(\text{frame})} \lesssim \sqrt{\mathcal{E}} (1+r)^{-1 + C_{(0)}\epsilon}
	\end{equation*}
	then we have the bound
	\begin{equation*}
	|\tr_{\slashed{g}}\chibar_{(\text{small})}|
	\lesssim
	\delta^{-1}\sqrt{\mathcal{E}}(1+r)^{-1 + C_{(0)}\epsilon}
	\end{equation*}
	
\end{proposition}

\begin{proof}
	Recall proposition \ref{proposition chibar in terms of chi}. Taking the trace of the expression given in this proposition, we have
	\begin{equation*}
	|\tr_{\slashed{g}}\chibar_{(\text{small})}|
	\lesssim
	|\tr_{\slashed{g}}\chi_{(\text{small})}|
	+ |\partial h|_{(\text{frame})}
	+ r^{-1} \left|2 - \sum_{i=1}^3 \slashed{\Pi}_\mu^{\phantom{\mu}i} \slashed{\Pi}^{\mu i} \right|
	\end{equation*}
	We deal with this final term by following the calculations in proposition \ref{proposition Yn tr chibar}. Specifically, we have
	\begin{equation*}
	\left|2 - \sum_{i=1}^3 \slashed{\Pi}_\mu^{\phantom{\mu}i} \slashed{\Pi}^{\mu i} \right|
	\lesssim
	|X_{(\text{frame, small})}|
	+ |h_{(\text{rect})}|
	+ |\mathcal{O}(h_{(\text{rect})})^2|
	\end{equation*}
	Now, using the bounds on $h_{(\text{rect})}$ assumed in the proposition, along with the bounds on $X_{(\text{frame, small})}$ from proposition \ref{proposition pointwise bounds rectangular}, we have
	\begin{equation*}
	\left|2 - \sum_{i=1}^3 \slashed{\Pi}_\mu^{\phantom{\mu}i} \slashed{\Pi}^{\mu i} \right|
	\lesssim
	1
	\end{equation*}
	where we have assumed that $\mathcal{E} \ll \delta$.
	
	Now, substituting the first bounds we have assumed for $\partial h$, together with the bound for $\tr_{\slashed{g}}\chi_{(\text{small})}$ from proposition \ref{proposition pointwise bound trchi small} proves the first part of the proposition. On the other hand, substituting the second bound for $\partial h$ we obtain the second part of the proposition.
\end{proof}

\begin{proposition}[A pointwise bound for $\mathscr{Y}^n\tr_{\slashed{g}}\chibar_{(\text{small})}$]
	\label{proposition pointwise bound Yn tr chibar}	
	Suppose that the same bounds hold as in proposition \ref{proposition pointwise bounds Yn Xframe}. Suppose additionally that
	\begin{equation*}
	|\overline{\slashed{\D}} \mathscr{Y}^{n+1} h|_{(\text{frame})}
	\lesssim
	\sqrt{\mathcal{E}} (1+r)^{-1-2\delta}
	\end{equation*}
	
	Then we have
	\begin{equation*}
	|\mathscr{Y}^n\tr_{\slashed{g}}\chibar_{(\text{small})}|
	\lesssim
	\delta^{-1}\sqrt{\mathcal{E}}(1+r)^{-1 + \delta}
	\end{equation*}
	
	If, in addition, we have the bound
	\begin{equation*}
	|\slashed{\D} \mathscr{Y}^n h|_{(\text{frame})}
	\lesssim
	\sqrt{\mathcal{E}} (1+r)^{-1 + C_{(n)}\epsilon}
	\end{equation*}
	then we can obtain the bound
	\begin{equation*}
	|\mathscr{Y}^n\tr_{\slashed{g}}\chibar_{(\text{small})}|
	\lesssim
	\delta^{-1}\sqrt{\mathcal{E}} (1+r)^{-1 + C_{(n)}\epsilon}
	\end{equation*}
	
\end{proposition}

\begin{proof}
	Recall proposition \ref{proposition Yn tr chibar}. Substituting the bounds we have assumed, along with the bounds on $\mathscr{Y}^n X_{(\text{frame})}$, $\mathscr{Z}^n X_{(\text{frame, small})}$ and $\mathscr{Y}^n \tr_{\slashed{g}}\chi_{(\text{small})}$ which were proved above, we obtain the required bound.
\end{proof}

\begin{proposition}[A pointwise bound for $\hat{\chibar}$]
	\label{proposition pointwise bound chibar hat}	
	Suppose that
	\begin{equation*}
	\begin{split}
	|h|_{(\text{rect})} &\lesssim  \sqrt{\mathcal{E}}(1+r)^{-2\delta} \\
	|\partial h|_{(\text{frame})} &\lesssim	\sqrt{\mathcal{E}} (1+r)^{-1+\delta} \\
	|\bar{\partial}h|_{(\text{frame})} &\lesssim \sqrt{\mathcal{E}} (1+r)^{-1-3\delta} \\
	|\overline{\slashed{\D}} \mathscr{Z} h|_{(\text{frame})} &\lesssim \sqrt{\mathcal{E}} (1+r)^{-1-2\delta} \\
	\end{split}
	\end{equation*}
	
	Then we have
	\begin{equation}
	|\hat{\chibar}|
	\lesssim
	\delta^{-1}\sqrt{\mathcal{E}}(1+r)^{-1 + \delta}
	\end{equation}
	
	If, in addition, we have the bound
	\begin{equation*}
	|\partial h|_{(\text{frame})} \lesssim	\sqrt{\mathcal{E}} (1+r)^{-1+C_{(0)}\epsilon}
	\end{equation*}
	then
	\begin{equation}
	|\hat{\chibar}|
	\lesssim
	\delta^{-1}\sqrt{\mathcal{E}}(1+r)^{-1 + C_{(0)}\epsilon}
	\end{equation}
	
\end{proposition}

\begin{proof}
	Recall proposition \ref{proposition chibar in terms of chi}. Subtracting the trace from the expression given in this proposition, we have
	\begin{equation*}
	|\hat{\chibar}|
	\lesssim
	|\hat{\chi}|
	+ |\partial h|_{(\text{frame})}
	+ r^{-1}\left| \sum_{i=1}^3 \left( \slashed{\Pi}_\mu^{\phantom{\mu}i} \slashed{\Pi}_\nu^{\phantom{\nu}i} - \frac{1}{2} \slashed{g}_{\mu\nu} \slashed{\Pi}_\rho^{\phantom{\rho}i} \slashed{\Pi}^{\rho i} \right) \right|
	\end{equation*}
	To control this second term, we follow the calculations in proposition \ref{proposition Yn chibar hat}, which lead to the estimate
	\begin{equation*}
	\left| \sum_{i=1}^3 \left( \slashed{\Pi}_\mu^{\phantom{\mu}i} \slashed{\Pi}_\nu^{\phantom{\nu}i} - \frac{1}{2} \slashed{g}_{\mu\nu} \slashed{\Pi}_\rho^{\phantom{\rho}i} \slashed{\Pi}^{\rho i} \right) \right|
	\lesssim
	(1+|X_{(\text{frame})}|)
	\left(|h_{(\text{rect})}|
	+ |X_{(\text{frame, \text{small}})}|
	+ \mathcal{O}(|h_{(\text{rect})}|^2)
	\right)
	\end{equation*}
	So, substituting the bounds for $h_{(\text{rect})}$ and $(\partial h)_{(\text{frame})}$ assumed in the proposition, along with the bounds on $X_{(\text{frame})}$ and $X_{(\text{frame, small})}$ from the previous propositions, we have proved the proposition.
\end{proof}

\begin{proposition}[A pointwise bound for $\mathscr{Y}^n \hat{\chibar}$]
	\label{proposition pointwise bound Yn chibar hat}	
	Suppose that the same bounds hold as in proposition \ref{proposition pointwise bounds Yn Xframe}. Suppose, additionally, that
	\begin{equation*}
	|\overline{\slashed{\D}}\mathscr{Y}^{n+1} h|_{(\text{frame})} \lesssim \sqrt{\mathcal{E}} (1+r)^{-1-2\delta}
	\end{equation*}
	Then we have
	\begin{equation*}
	|\mathscr{Y}^n \hat{\chibar}|
	\lesssim
	\delta^{-1}\sqrt{\mathcal{E}}(1+r)^{-1 + \delta}
	\end{equation*}
	
	If, in addition, we have
	\begin{equation*}
	|\slashed{\D} \mathscr{Y}^n h|_{(\text{frame})}
	\lesssim
	\sqrt{\mathcal{E}} (1+r)^{-1+C_{(n)}\epsilon}
	\end{equation*}
	then
	\begin{equation*}
	|\mathscr{Y}^n \hat{\chibar}|
	\lesssim
	\delta^{-1}\sqrt{\mathcal{E}}(1+r)^{-1 + C_{(n)}\epsilon}
	\end{equation*}
	
\end{proposition}

\begin{proof}
Recall proposition \ref{proposition Yn chibar hat}. Substituting the bounds we have assumed, along with the bounds on $\mathscr{Y}^n X_{(\text{frame})}$, $\mathscr{Y}^n X_{(\text{frame, small})}$ and $\mathscr{Y}^n \tr_{\slashed{g}}\chi_{(\text{small})}$ which were proved above, we obtain the required bound.
\end{proof}

\begin{proposition}[A pointwise bound on $\Omega$]
	\label{proposition pointwise bound on Omega}
	Suppose that the same bounds hold as in proposition \ref{proposition pointwise bound Yn tr chi}. Suppose also that
	\begin{equation*}
		|h|_{(\text{frame})} \lesssim \sqrt{\mathcal{E}} (1+r)^{-\frac{1}{2} + \mathring{C}\delta}
	\end{equation*}
	Then
\end{proposition}

\begin{proof}
	Recall proposition \ref{proposition transport Omega}. Substituting the bound that we have already obtained for $\tr_{\slashed{g}} \chi_{(\text{small})}$ we have
	\begin{equation*}
	\left|L\log\Omega - \frac{1}{r}\right| \lesssim \delta^{-1} \sqrt{\mathcal{E}} (1+r)^{-1-2\delta}
	\end{equation*}
	Note also that, at $r = r_0$ we have
	\begin{equation*}
	\frac{\det \slashed{g}}{\det \mathring{\gamma}} \lesssim (r_0)^2 (1+|h|_{(\text{rect})})
	\end{equation*}
	So, integrating from $r = r_0$ we obtain
	\begin{equation*}
	e^{-\delta^{-2} \sqrt{\mathcal{E}}} \lesssim \left| \frac{\Omega}{r} \right| \lesssim e^{\delta^{-2} \sqrt{\mathcal{E}}}
	\end{equation*}
\end{proof}

Note that we only need to provide a bound on $\Omega$ and not its higher derivatives. Similarly, we need such a bound on the Gauss curvature of the spheres $K$:

\begin{proposition}[A pointwise bound on the Gauss curvature $K$]
	\label{proposition pointwise bound on Gauss curvature}
	Suppose that the Gauss curvature $K$ satisfies
	\begin{equation*}
	|r^2 K| \lesssim 1 + \epsilon
	\end{equation*}
	Suppose also that the same conditions hold as in proposition \ref{proposition pointwise bounds Yn Xframe}, and also that
	\begin{equation*}
	|\mathscr{Z}^2 h| \lesssim \sqrt{\mathcal{E}} (1+r)^{-\frac{1}{2} + \mathring{C}\delta}
	\end{equation*}
	
	Then
	\begin{equation*}
	|r^2 K| \lesssim 1 + \delta^{-2} \sqrt{\mathcal{E}}
	\end{equation*}
\end{proposition}

\begin{proof}
	Recall proposition \ref{proposition transport Gauss curvature}. Substituting the bounds that we have already obtained for $\tr_{\slashed{g}} \chi_{(\text{small})}$ and its first two derivatives, and also for $\hat{\chi}$ and its first two derivatives, we have
	\begin{equation*}
	\left| L(r^2 K) \right|
	\lesssim
	\delta^{-1} \sqrt{\mathcal{E}} (1+r)^{-1-2\delta + c_{[1]}\epsilon}
	+ \delta^{-1} \sqrt{\mathcal{E}} (1+r)^{-1-2\delta} |r^2 K|
	\end{equation*}
	Now, if we substitute the bound for $K$ in the proposition, and also note that at $r = r_0$ we have
	\begin{equation*}
	|r^2 K| \big|_{r = r_0} \lesssim (r_0)^2 (1+\sqrt{\mathcal{E}})
	\end{equation*}
	then we obtain the bound
	\begin{equation*}
	|r^2 K| \lesssim \delta^{-2} \sqrt{\mathcal{E}}
	\end{equation*}
\end{proof}

\section{Pointwise bounds in the region \texorpdfstring{$r \leq r_0$}{r leq r0}}
\label{section pointwise bounds in r < r0}

All of the pointwise bounds we have derived above apply exclusively to the region $r \geq r_0$. In order to obtain pointwise bounds in the region $r \leq r_0$ we will instead use elliptic estimates.

\begin{proposition}[A uniformly elliptic operator in the region $r \leq r_0$]
	Suppose that
	\begin{equation*}
	|h_{(\text{rect})}| \lesssim \epsilon
	\end{equation*}
	Then for all sufficiently small $\epsilon$, in the region $r \leq r_0$, the restriction of the operator
	\begin{equation}
	\Delta_g
	:=
	\tilde{\Box}_g
	- (g^{-1})^{00}T^2
	- 2(g^{-1})^{00}\partial_i T
	- (g^{-1})^{0c} (g^{-1})^{ab}\left(\partial_a h_{bc} - \frac{1}{2}\partial_c h_{ab} \right) T
	\end{equation}
	to $\Sigma_\tau$ is a uniformly elliptic operator on functions 
	\begin{equation*}
	\phi \ : \ \Sigma_\tau \cap \{r \leq r_0\} \rightarrow \mathbb{R}
	\end{equation*}
	Moreover, the constant appearing in the uniformity condition can be chosen to be independent of $\tau$.
	
%	Similarly, the restriction of the operator
%	\begin{equation}
%	\Delta_g
%	:=
%	\slashed{\tilde{\Box}}_g
%	- \left((g^{-1})^{00} T^{\mu}T^\nu 
%	- 2(g^{-1})^{00}\slashed{\D}
%	
%	- (g^{-1})^{0c} (g^{-1})^{ab}\left(\partial_a h_{bc} - \frac{1}{2}\partial_c h_{ab} \right) T
%	\end{equation}
	
\end{proposition}

\begin{proof}
	Note that in the region $r \leq r_0$, we have
	\begin{equation*}
	\tilde{\Box}_g = \Box_g
	\end{equation*}
	since $\chi_{r_0} \cdot \omega \slashed{\D}_{\Lbar} \phi = 0$ here. Expanding the wave operator in rectangular coordinates, we find, for any scalar field $\phi$,
	\begin{equation*}
	\begin{split}
	\Box_g \phi
	&=
	(g^{-1})^{ab} \partial_a \partial_b \phi
	+ (g^{-1})^{ab}(g^{-1})^{cd} \left( \partial_a h_{bd} - \frac{1}{2}\partial_d h_{ab} \right) \partial_c \phi \\
	\\
	&=
	(g^{-1})^{00} T^2 \phi
	+ 2(g^{-1})^{0i} \partial_i T \phi
	+ (g^{-1})^{ij} \partial_i \partial_j \phi
	+ (g^{-1})^{0c}(g^{-1})^{ab} \left( \partial_a h_{bc} - \frac{1}{2}\partial_c h_{ab} \right) T \phi 
	\\
	&\phantom{=}
	+ (g^{-1})^{ic}(g^{-1})^{ab} \left( \partial_a h_{bc} - \frac{1}{2}\partial_c h_{ab} \right) \partial_i \phi 
	\end{split}
	\end{equation*}
	Hence we can express the operator $\Delta_g$ acting on a scalar field $\phi$ as
	\begin{equation*}
	\Delta_g \phi
	=
	(g^{-1})^{ij} \partial_i \partial_j \phi
	+ (g^{-1})^{ic}(g^{-1})^{ab} \left( \partial_a h_{bc} - \frac{1}{2}\partial_c h_{ab} \right) \partial_i \phi 
	\end{equation*}
	Note that the restriction of this operator to the leaf $\Sigma_\tau \cap \{r \leq r_0\}$ can be expressed by an identical formula.
	
	The principle symbol of $\Delta_g$ is therefore
	\begin{equation*}
	\sigma(x, \xi) :=
	(g^{-1})^{ij}(x) \xi_i \xi_j
	\end{equation*}
	for any vector field $\xi$ on $\Sigma_\tau \cap \{r \leq r_0\}$. To check that $\Delta_g$ is uniformly elliptic, we need to show that
	\begin{equation*}
	\sigma(x, \xi) \gtrsim |\xi|^2 
	\end{equation*}
	with some implicit constant which is uniform in $x$. But we can calculate, for sufficiently small $\epsilon$
	\begin{equation*}
	\begin{split}
	(g^{-1})^{ij}(x) \xi_i \xi_j
	&=
	\left((m^{-1})^{ij} + H^{ij}(x) \right) \xi_i \xi_j \\
	&= |\xi|^2 + H^{ij}(x)\xi_i \xi_j \\
	&\geq \frac{1}{2} |\xi|^2
	\end{split}
	\end{equation*}
	where the last line follows from
	\begin{equation*}
	H^{ij} = -h^{ij} + \mathcal{O}(h_{(\text{rect})}^2)
	\end{equation*}
	and the fact that the rectangular components of $h$ are $\mathcal{O}(\epsilon)$. Hence, for all sufficiently small $\epsilon$, the operator $\Delta_g$ is uniformly elliptic. Moreover, the bounds we have used can all be chosen to be independent of $\tau$, since we have the (uniform in $\tau$) bound $h_{(\text{rect})} \lesssim \epsilon$
\end{proof}

\begin{proposition}[$H^2$ and $C^{0, \frac{1}{2}}$ bounds]
	\label{proposition H2 and C0 bounds}
		Let $\phi$ be a scalar field satisfying
		\begin{equation*}
		\tilde{\Box}_g \phi = F
		\end{equation*}
		
		Suppose that
		\begin{equation*}
		\int_{\Sigma_\tau \cap \{r \leq r_0\}} \left( (\partial T \phi)^2 + (T \phi)^2 + |\phi|^2 + F^2 \right) \dVol_{\Sigma_g} \lesssim \mathcal{E}
		\end{equation*}
		
		Suppose also that 
		\begin{equation*}
		\begin{split}
		||h_{\text{rect}}||_{L^\infty[\Sigma_\tau \cap \{r \leq r_0\}]} &\lesssim \epsilon
		\\
		||\partial h_{\text{rect}}||_{L^\infty[\Sigma_\tau \cap \{r \leq r_0\}]} &\lesssim \epsilon
		\end{split}
		\end{equation*}
		
		Then we have
		\begin{equation*}
		\begin{split}
		||\phi||_{H^2[\Sigma_\tau \cap \{r \leq \frac{3}{4} r_0\}]} &\lesssim \sqrt{\mathcal{E}}
		\\
		||\phi||_{C^{0,\frac{1}{2}}[\Sigma_\tau \cap \{r \leq \frac{3}{4}r_0\}]} &\lesssim \sqrt{\mathcal{E}}
		\end{split}
		\end{equation*}
		
\end{proposition}

\begin{proof}
	Since $\phi$ satisfies the wave equation, we have
	\begin{equation*}
	\Delta_g \phi
	=
	F
	- (g^{-1})^{00} T^2 \phi
	- 2 (g^{-1})^{0i} \partial_i T \phi
	- (g^{-1})^{0c}(g^{-1})^{ab} \left(\partial_a h_{bc} - \frac{1}{2} \partial_c h_{ab} \right) T\phi
	\end{equation*}
	Moreover, $\Delta_g$ is uniformly elliptic and $g_{ab}$, $\partial g_{ab} \in L^\infty$ (since $g_{ab} = m_{ab} + h_{ab}$ and the $m_{ab}$ are just constants). Hence it is easy to show (for example, by considering the elliptic equation satisfied by the cut-off function $\chi_{(\frac{3}{4}r_0, r_0)}\cdot \phi$) that
	\begin{equation*}
	\int_{\Sigma_\tau \cap \{r \leq \frac{3}{4}r_0\}}
	\sum_{i,j}|\partial_i \partial_j \phi|^2 \dVol_g
	\lesssim
	\int_{\Sigma_\tau \cap \{r \leq r_0\}}
	\left( |\Delta_g \phi|^2 + |\phi|^2 \right) \dVol_g
	\end{equation*}
	and so, since the lower order term is controlled by assumption, we have
	\begin{equation*}
	||\phi||^2_{H^2[\Sigma_\tau \cap \{r \leq \frac{3}{4}r_0\}]} \lesssim \mathcal{E}
	\end{equation*}
	
	Now, using the Sobolev inequality we can immediately show the $C^{0, \frac{1}{2}}$ bound.
\end{proof}

\begin{proposition}[$H^3$ and $C^{1,\frac{1}{2}}$ bounds]
	\label{proposition H3 and C1 bounds}
	Let $\phi$ be a scalar field satisfying
	\begin{equation*}
	\tilde{\Box}_g \phi = F
	\end{equation*}
	where
	\begin{equation*}
	\begin{split}
	||h_{\text{rect}}||_{L^\infty[\Sigma_\tau \cap \{ r \leq r_0\}]} \lesssim \epsilon
	\\
	||\partial h_{\text{rect}}||_{L^\infty[\Sigma_\tau \cap \{ r \leq r_0\}]} \lesssim \epsilon
	\end{split}
	\end{equation*}
	and suppose that
	\begin{equation*}
	\begin{split}
	||\phi||^2_{H^1[\Sigma_{\tau} \cap \{r \leq \frac{3}{4}r_0\}]}
	& \lesssim \mathcal{E} \\
	|| T \phi ||^2_{H^{2}[\Sigma_{\tau} \cap \{r \leq \frac{3}{4}r_0\}]}
	& \lesssim \mathcal{E} \\
	|| T^2 \phi ||^2_{H^{1}[\Sigma_{\tau} \cap \{r \leq \frac{3}{4}r_0\}]}
	& \lesssim \mathcal{E} \\
	||F||^2_{H^{1}[\Sigma_{\tau} \cap \{r \leq \frac{3}{4}r_0\}]} 
	& \lesssim \mathcal{E} \\
	||h||^2_{H^{2}[\Sigma_{\tau} \cap \{r \leq \frac{3}{4}r_0\}]}
	& \lesssim \mathcal{E} \\
	||Th||^2_{H^{1}[\Sigma_{\tau} \cap \{r \leq \frac{3}{4}r_0\}]}
	& \lesssim \mathcal{E} \\
	\end{split}
	\end{equation*}
	
	Then we have
	\begin{equation*}
	\begin{split}
	||\phi||_{H^3[\Sigma_{\tau} \cap \{r \leq \frac{2}{3}r_0\}]}
	\lesssim \sqrt{\mathcal{E}} \\
	||\phi||_{C^{1, \frac{1}{2}}[\Sigma_{\tau} \cap \{r \leq \frac{2}{3}r_0\}]}
	\lesssim \sqrt{\mathcal{E}}
	\end{split}	
	\end{equation*}
	
\end{proposition}

\begin{proof}
	As before, we have
	\begin{equation}
	\label{equation elliptic internal 1}
	\Delta_g \phi
	=
	F
	- (g^{-1})^{00} T^{2} \phi
	- 2 (g^{-1})^{0i} \partial_i T \phi
	- (g^{-1})^{0c}(g^{-1})^{ab} \left(\partial_a h_{bc} - \frac{1}{2} \partial_c h_{ab} \right) T\phi
	\end{equation}

	Next, we commute equation \eqref{equation elliptic internal 1} with the vector fields $\partial_i$. Using the notation $\underline{\partial}$ to stand for any of the derivatives $\partial_i$, $i \in \{1,2,3\}$, it is easy to obtain the inequality
	\begin{equation*}
	\begin{split}
	|\Delta_g  \underline{\partial} \phi|
	&\lesssim
	|\underline{\partial} F|
	+ |\underline{\partial} T^2 \phi |
	+ |\underline{\partial}^2 T \phi|
	+ \left( |\underline{\partial} h| + |T h| \right) |\underline{\partial} T \phi|
	\\
	&\phantom{\lesssim}
	+ |\underline{\partial} h| \left( |T^2 \phi| + |\underline{\partial} T \phi| + |\underline{\partial}^2 \phi| \right)
	+ \left( |\underline{\partial}^2 h| + |\underline{\partial} T h| \right) |T \phi|
	\end{split}
	\end{equation*}

	Now, the elliptic estimates tell us that
	\begin{equation*}
	\begin{split}
	||\phi||^2_{H^{3}[\Sigma_{\tau} \cap \{r \leq \frac{2}{3}r_0\}]}
	&\lesssim
	||F||^2_{H^{1}[\Sigma_{\tau} \cap \{r \leq \frac{3}{4}r_0\}]}
	+ || T^2 \phi ||^2_{H^{1}[\Sigma_{\tau} \cap \{r \leq \frac{3}{4}r_0\}]}
	+ || T \phi ||^2_{H^{2}[\Sigma_{\tau} \cap \{r \leq \frac{3}{4}r_0\}]}
	\\
	&\phantom{\lesssim}
	+ ||h||^2_{H^{2}[\Sigma_{\tau} \cap \{r \leq \frac{3}{4}r_0\}]}
	+ ||Th||^2_{H^{1}[\Sigma_{\tau} \cap \{r \leq \frac{3}{4}r_0\}]}
	+ ||\phi||^2_{H^1[\Sigma_{\tau} \cap \{r \leq \frac{3}{4}r_0\}]}
	\end{split}
	\end{equation*}
	and so, by the conditions of the proposition, together with the Sobolev inequality we have
	\begin{equation*}
	||\phi||^2_{C^{1, \frac{1}{2}}[\Sigma_{\tau} \cap \{r \leq \frac{2}{3}r_0\}]}
	\end{equation*}

\end{proof}

\begin{proposition}[$C^{k, \frac{1}{2}}$ estimates]
	\label{proposition Ck estimates}
	Let $\phi$ be a scalar field satisfying
	\begin{equation*}
	\tilde{\Box}_g \phi = F
	\end{equation*}
	and suppose that
	\begin{equation*}
	\begin{split}
	||F||_{C^{k, \frac{1}{2}}[\Sigma_\tau \cap \{ r \leq \frac{2}{3}r_0\}]}
	&\lesssim \sqrt{\mathcal{E}}
	\\
	||h||_{C^{k+1, \frac{1}{2}}[\Sigma_\tau \cap \{ r \leq \frac{2}{3}r_0\}]}
	&\lesssim \sqrt{\mathcal{E}}
	\\
	||Th||_{C^{k, \frac{1}{2}}[\Sigma_\tau \cap \{ r \leq \frac{2}{3}r_0\}]}
	&\lesssim \sqrt{\mathcal{E}}
	\\
	||T^2 \phi||_{C^{k, \frac{1}{2}}[\Sigma_\tau \cap \{ r \leq \frac{2}{3}r_0\}]}
	&\lesssim \sqrt{\mathcal{E}}
	\\
	||T \phi||_{C^{k+1, \frac{1}{2}}[\Sigma_\tau \cap \{ r \leq \frac{2}{3}r_0\}]}
	&\lesssim \sqrt{\mathcal{E}}
	\\
	||\phi||_{C^{0}[\Sigma_\tau \cap \{ r \leq \frac{2}{3}r_0\}]}
	&\lesssim \sqrt{\mathcal{E}}
	\end{split}
	\end{equation*}
	for some $k \geq 0$.
	
	Then we have
	\begin{equation*}
	\begin{split}
	||\phi||_{C^{k+2, \frac{1}{2}}[\Sigma_{\tau} \cap \{r \leq \frac{1}{2}r_0\}]}
	\lesssim \sqrt{\mathcal{E}}
	\end{split}	
	\end{equation*}
	
\end{proposition}

\begin{proof}
	This follows immediately from the interior Schauder estimates, together with the expression in equation \eqref{equation elliptic internal 1} for the elliptic equation satisfied by $\phi$.
\end{proof}

\begin{remark}[The need for the $C^{0, \frac{1}{2}}$ and $C^{1, \frac{1}{2}}$ bounds]
	If we want to set $k = 0$ in the higher order interior Schauder estimates above, then we need to already have $C^{1, \frac{1}{2}}$ bounds on $h$ and $T\phi$, as well as $C^{0, \frac{1}{2}}$ bounds on $Th$ and $T^2 \phi$, and a $C^0$ bound on $\phi$. Since we will treat equations in which $h$ depends (smoothly) on the fields $\phi$, these requirements reduce to $C^{1, \frac{1}{2}}$ bounds on $\phi$ and $T\phi$ and $C^{0, \frac{1}{2}}$ bounds on $T^2\phi$. Hence, we need to first establish these bounds, for which we use the other propositions.
\end{remark}

When we perform the elliptic estimates, it is necessary to restrict the estimates to a sequence of progressively smaller balls. Hence, we need additional control over various quantities in the region $\frac{1}{2} r_0 \leq r \leq r_0$.

\begin{proposition}[Pointwise bounds in the region $\frac{1}{2} r_0 \leq r \leq r_0$]
	\label{proposition pointwise bounds 1/2 r0 leq r leq r0}
	Suppose that the pointwise bootstrap bounds of chapter \ref{chapter bootstrap} hold. Suppose, in addition, that $\phi$ is an $S_{\tau,r}$-tangent tensor field satisfying
	\begin{equation*}
	\tilde{\slashed{\Box}}_g (r\slashed{\nabla})^n \phi = F_n
	\end{equation*}
	$n = 0,1,2,3$, and for some $S_{\tau,r}$ tangent tensor fields $F_n$ satisfying
	
	\begin{equation*}
	\int_{\Sigma_\tau \cap \{r \leq r_0\}} |F_n|^2 \lesssim \sqrt{\mathcal{E}}
	\end{equation*}
	
	Suppose additionally that $\phi$ satisfies
	\begin{equation*}
	\sum_{m=0}^3 \int_{\Sigma_\tau \cap \{\frac{1}{2} r_0 \leq r \leq r_0\}}
	\left( |\slashed{\D} (r\slashed{\nabla})^m \phi|^2 + |(r\slashed{\nabla})^m \phi|^2 \right) \upd r \wedge \dVol_{\mathbb{S}^2}
	\lesssim
	\mathcal{E}
	\end{equation*}
	
	Then we have
	\begin{equation*}
	\sup_{S_{\tau,r}} \left( |\phi|^2 + |\slashed{\D}\phi|^2 \right)
	\lesssim
	\mathcal{E}
	+ \int_{S_{\tau,r_0}} \sum_{m=0}^3 \left(|\mathscr{Z}^m\phi|^2 + |\slashed{\D}\mathscr{Z}^m \phi|^2 \right) \dVol_{\mathbb{S}^2}
	\end{equation*}
	
\end{proposition}

\begin{proof}
	Following similar steps to proposition \ref{proposition spherical mean in terms of energy}, but this time choosing $r \in [\frac{1}{2}r_0, r_0]$, we find that for a scalar field $\phi$
	\begin{equation*}
	\begin{split}
	\int_{S_{\tau, r}} |\phi|^2 \dVol_{\mathbb{S}^2}
	&=
	\int_{\mathbb{S}^2} \left( \int_{r' = r}^{r_0} \frac{\partial \phi}{\partial r'} (\tau, r, \vartheta^1, \vartheta^2) \upd r' \right)^2 \dVol_{\mathbb{S}^2}
	- \int_{S_{\tau, r_0}} |\phi|^2 \dVol_{\mathbb{S}^2}
	\\
	&\phantom{=}
	+ \int_{\mathbb{S}^2} 2\phi(\tau, r, \vartheta^1, \vartheta^2)\phi(\tau, r_0, \vartheta^1, \vartheta^2) \dVol_{\mathbb{S}^2}
	\\
	\\
	&\lesssim 
	r_0 \int_{\mathbb{S}^2} \int_{r' = r}^{r_0} \left(\frac{\partial \phi}{\partial r'} \right)^2 \upd r' \dVol_{\mathbb{S}^2}
	+ \int_{S_{\tau, r_0}} |\phi|^2 \dVol_{\mathbb{S}^2}
	\end{split}
	\end{equation*}
	The computation for a higher rank $S_{\tau,r}$-tangent tensor field is similar: we can consider the sum of the squares of the rectangular components, and use the bootstrap bounds to pass back to an estimate for $|\phi|^2$.
	
	Under the pointwise bootstrap assumptions the first term is bounded by the degenerate energy. Hence, for $r \in [\frac{1}{2}r_0, r_0]$ we have
	\begin{equation*}
	\int_{S_{\tau, r}} |\phi|^2 \dVol_{\mathbb{S}^2}
	\lesssim
	\mathcal{E} + \int_{S_{\tau, r_0}} |\phi|^2 \dVol_{\mathbb{S}^2}
	\end{equation*}
	
	Applying the same computation to the fields $\mathscr{Z}^m \phi$, and recalling that $\mathscr{Z}$ can be $r\slashed{\nabla}$, we can use the Sobolev embedding on the sphere (proposition \ref{proposition Sobolev}) to show the pointwise bound in the region $\frac{1}{2}r_0 \leq r \leq r_0$
	\begin{equation*}
	\sup_{S_{\tau, r}}  \left( |\phi|^2 |\slashed{\D}_T \phi|^2 + |\slashed{\nabla} \phi|^2 \right)
	\lesssim
	\mathcal{E} + \sum_{m=0}^3 \int_{S_{\tau, r_0}} |\mathscr{Z}^m\phi|^2 \dVol_{\mathbb{S}^2}
	\end{equation*}
	
	Note that we cannot take $r$ to be arbitrarily small in this way, due to the factor of $r$ in the Sobolev inequality on the sphere. This is the reason for using elliptic estimates for small $r$.
	
	Next, we want to estimate the derivatives of $\phi$ in the region $\frac{1}{2}r_0 \leq r \leq r_0$. We note that, in this region
	\begin{equation*}
	\begin{split}
	|\slashed{\D}_T\phi| = |\mathscr{Z}\phi| \\
	|\slashed{\nabla} \phi| \lesssim |\mathscr{Z}\phi|
	\end{split}
	\end{equation*}
	So, if we are prepared to apply one more operator $\mathscr{Z}$ then we immediately obtain pointwise control of the $T$ derivatives and the angular derivatives. The only derivative we do not gain control over is the $\partial_r$ derivative. 
	
	Recall the notation $R$ for the vector field $\partial_r$. We note that, in the region $r \leq r_0$, for an $S_{\tau,r}$-tangent tensor field $\phi_{\alpha_1 \ldots \alpha_n}$ we can estimate
	\begin{equation*}
	\begin{split}
	\tilde{\slashed{\Box}}_g \phi_{\alpha_1 \ldots \alpha_n}
	&=
	(m^{-1})^{ab} \slashed{\D}^2_{ab} \phi_{\alpha_1 \ldots \alpha_n}
	+ H^{ab} \slashed{\D}^2_{ab} \phi_{\alpha_1 \ldots \alpha_n}
	\\
	\\
	&=
	- \slashed{\D}_T \slashed{\D}_T \phi_{\alpha_1 \ldots \alpha_n}
	+ \slashed{\D}_R \slashed{\D}_R \phi_{\alpha_1 \ldots \alpha_n}
	+ \frac{2}{r} \slashed{\D}_R \phi_{\alpha_1 \ldots \alpha_n}
	+ \slashed{\Delta}_{(m)} \phi_{\alpha_1 \ldots \alpha_n}
	+ (\D_T T)^\mu \slashed{\D}_\mu \phi_{\alpha_1 \ldots \alpha_n}
	\\
	&\phantom{=}
	- (\D_R R)^\mu \slashed{\D}_\mu \phi_{\alpha_1 \ldots \alpha_n}
	+ \left(\slashed{\D}_T \slashed{\Pi}_{\alpha_1 \ldots \alpha_n}^{\beta_1 \ldots \beta_n} \right) \D_T \phi_{\beta_1 \ldots \beta_n}
	- \left(\slashed{\D}_R \slashed{\Pi}_{\alpha_1 \ldots \alpha_n}^{\beta_1 \ldots \beta_n} \right) \D_R \phi_{\beta_1 \ldots \beta_n}
	\\
	&\phantom{=}
	- (\slashed{m}^{-1})^{\mu\nu}\left(\slashed{\nabla}_\mu \slashed{\Pi}_{\alpha_1 \ldots \alpha_n}^{\beta_1 \ldots \beta_n} \right) \D_\nu \phi_{\beta_1 \ldots \beta_n}
	+ H^{ab} \slashed{\D}^2_{ab} \phi_{\alpha_1 \ldots \alpha_n}
	\end{split}
	\end{equation*}
	where the operator $\slashed{\Delta}_{(m)}$ is the standard laplacian operator on the spheres, i.e.\
	\begin{equation*}
	\slashed{\Delta}_{(m)} := (\slashed{m}^{-1})^{\mu\nu} \slashed{\nabla}_\mu \slashed{\nabla}_\nu
	\end{equation*}
	and where $\slashed{m}$ is the restriction of the Minkowski metric $m$ to the spheres $S_{\tau,r}$. The error terms are bounded by the bootstrap assumptions. Hence we have
	\begin{equation*}
	|\slashed{\D}_R \slashed{\D}_R \phi|
	\lesssim
	|F|
	+ |\slashed{\D} \mathscr{Z} \phi|
	+ |\slashed{\D} \phi|
	+ |\phi|
	+ |H_{(\text{rect})}| |\slashed{\D}^2 \phi|
	\end{equation*}
	in the region $r \leq r_0$. Using the bootstrap bounds on  $|H_{(\text{rect})}|$, this last term can be absorbed into the others. Then, following similar calculations to those given above, we can show that
	\begin{equation*}
	\sup_{S_{\tau, r}}  |\slashed{\D}_R \phi|^2
	\lesssim
	\mathcal{E} + \sum_{m=0}^3 \int_{S_{\tau, r_0}} |\mathscr{Z}^m\phi|^2 \dVol_{\mathbb{S}^2}
	\end{equation*}
	
	Since $R$, $T$ and the angular derivatives $\slashed{\nabla}$ span the tangent space of $\mathcal{M}$, this completes the required bounds for $\slashed{\D}\phi$.
	
\end{proof}

\chapter{Estimates for the inhomogeneous term \texorpdfstring{$F$}{F}}
\label{chapter bounds for the inhomogeneous term}

In the previous chapters we have proved energy decay subject to certain pointwise bounds, and we have also recovered comparable pointwise bounds from the decay of certain higher order energies. However, we are not yet in a position to close the bootstrap, since the energy estimates (and some of the pointwise bounds) were only proved \emph{subject to certain bounds on the inhomogeneous term} $F$. The object of this chapter is to show how to recover these bounds on $F$.

For the first time in this chapter, instead of considering a general system of the form
\begin{equation*}
\tilde{\slashed{\Box}}_g\phi = F
\end{equation*}
for some $S_{\tau,r}$-tangent tensor field $\phi$ and some unspecified $S_{\tau,r}$-tangent tensor field $F$, we shall instead consider systems of the form
\begin{equation}
\label{equation system of wave equations}
\begin{split}
	\tilde{\Box}_g\phi_{(A)} &= F_{(A,0)}
	\\
	F_{(A,0)} &= F_{(A,0)}^{(0)} + \left(F_{(A,0)}^{(BC)}\right)^{\mu\nu}(\partial_\mu \phi_B)(\partial_\nu \phi_C) + \mathcal{O}\left(\phi (\partial\phi)^2 \right)
\end{split}
\end{equation}
where the $\phi_{(A)}$ are a set of \emph{scalar} fields. Note that we can easily add additional higher order terms, for example cubic terms of the form $(\partial \phi)^3$, but for clarity we will not include these terms.

In the expression given above, the term $F_{(A,0)}^{(0)}$ is some fixed function which obeys all of the relevant bounds. On the other hand, the coefficients
\begin{equation*}
\left(F_{(A,0)}^{(BC)}\right)_{LL}
\end{equation*}
are required to satisfy certain conditions, related to the weak null condition.

Under these conditions, together with the pointwise \emph{and} $L^2$ bootstrap bounds on the $\phi_{(A)}$ and their derivatives, we will show that $F_{(A,0)}$ satisfies all of the bounds required of the inhomogeneous term, as discussed above.

Next, we commute this system with the operators $\mathscr{Z}^n$, creating a system of the form
\begin{equation*}
	\tilde{\slashed{\Box}}_g \mathscr{Z}^n \phi_{(A)} = F_{(A,n)} 
\end{equation*}
We now need to show that $F_{(A,n)}$ \emph{also} obeys all of the bounds required of the inhomogeneous term. In order to do this, we need to first construct $F_{(A,n)}$, and then show that all of the terms in $F_{(A,n)}$ can be suitably controlled. This is the main computation carried out in this chapter.

First, we summarise the required bounds on the inhomogeneous term $F$. We need $L^2$ bounds of the following forms:
\begin{equation*}
\begin{split}
&\int_{\mathcal{M}_{\tau_0}^\tau} \bigg(
\epsilon^{-1} \chi_{(r_0)} r^{1-C_{(\phi)}\epsilon} (1+\tau)^{1+\delta} |F|^2 \bigg) \dVol_g \lesssim \mathcal{E} (1+\tau)^{C_{[\phi]}\delta}
\\
\\
&\int_{\mathcal{M}_{\tau}^{\tau_1}} \bigg(
\epsilon^{-1}(1+r)^{1-C_{(\phi)}\epsilon} |F|^2
+ \epsilon^{-1}(1+r)^{1-\delta}(1+\tau)^{2\beta}|F|^2
\bigg)\dVol_g
\lesssim
\mathcal{E}(1+\tau)^{-1}
\end{split}
\end{equation*}
These $L^2$ bounds are required to hold for all $F_{(A,n)}$, $n \leq M_2$.

Additionally, for some smaller number $M_3 < M_2$, the following pointwise bounds are required to hold in the region $r \geq r_0$: if $\phi_{(A)} \in \Phi_{[0]}$ then we require
\begin{equation*}
|F_{(A,n)}| \lesssim \epsilon r^{-1-\delta}|\slashed{\D}_{\Lbar}\phi_{(A)}| + \sqrt{\mathcal{E}}r^{-2-\delta}
\end{equation*}
on the other hand, if $\phi_{(A)} \in \Phi_{[m]}$ then we require
\begin{equation*}
|F_{(A,n)}| \lesssim \epsilon r^{-1}|\slashed{\D}_{\Lbar}\phi_{(A)}| + \sqrt{\mathcal{E}}r^{-2 + C_{(A,n)}\epsilon}
\end{equation*}

\section{Bounds on the inhomogeneous terms before commuting}
\label{section bounds on the inhomogeneous terms before commuting}

We aim to show that the inhomogeneous terms $F_{(A,n)}$ satisfy all the bounds above. This will be achieved via an induction argument, so we first need to show that the bounds hold for the ``base case'', that is, \emph{before} commuting the equations.

\begin{proposition}[$L^2$ bounds for the inhomogeneous terms before commuting]
	\label{proposition L2 bounds for F before commuting}
	Let $\phi_{(A)}$ be a set of scalar fields satisfying the equations
	\begin{equation*}
	\begin{split}
	\tilde{\Box}_g\phi_{(A)} &= F_{(A,0)}
	\\
	F_{(A,0)} &= F_{(A,0)}^{(0)} + \left(F_{(A,0)}^{(BC)}\right)^{\mu\nu}(\partial_\mu \phi_B)(\partial_\nu \phi_C) + \mathcal{O}\left(\phi (\partial\phi)^2\right)
	\end{split}
	\end{equation*}
	where we further decompose
	\begin{equation*}
	\begin{split}
	F_{(A,0)}^{(0)}
	&=
	F_{(A,0,1)}^{(0)}
	+ F_{(A,0,2)}^{(0)}
	+ F_{(A,0,3)}^{(0)}
	\\
	&=
	F_{(A,0,4)}^{(0)}
	+ F_{(A,0,5)}^{(0)}
	+ F_{(A,0,6)}^{(0)}
	\end{split}	
	\end{equation*}
	
	We require that the tensor fields $F_{(A,0)}^{(BC)}$ have \emph{constant rectangular components}\footnote{This condition can easily be weakened to the condition that the rectangular components are bounded.}. Also, we suppose that they satisfy\footnote{Again, these conditions can be weakened: in a bounded region, the frame coefficients of $F_{(A,0)}^{(BC)}$ need only be bounded, and the coefficients which we set to $0$ can actually be nonzero but decaying sufficiently fast in $r$.} the structural equations
	\begin{equation*}
	\begin{split}
	\left(F_{(A,0)}^{(BC)}\right)^{\mu\nu} 
	&= 	\left(F_{(A,0)}^{(CB)}\right)^{\mu\nu}
	\\
	\left(F_{(A,0)}^{(BC)}\right)_{LL}
	&= 0 \quad \text{if}\quad \phi_{(A)} \in \Phi_{[0]}
	\\
	\left(F_{(A,0)}^{(BC)}\right)_{LL}
	&= 0 \quad \text{if}\quad \phi_{(A)} \in \Phi_{[n]} \text{\, and either } \begin{cases}
	\phi_{(B)} \in \Phi_{[n+1]} \\
	\phi_{(B)} \in \Phi_{[n]} \text{\, and \,} \phi_{(C)} \in \Phi_{[m]} \text{ , } m \geq 1
	\end{cases}
	\end{split}
	\end{equation*}
	
	Suppose moreover that the terms $F^{(0)}_{(A,0)}$ satisfy the following conditions: if $\phi_{(A)} \in \Phi_{[0]}$, then 
	\begin{equation*}
	\begin{split}
	&\int_{\mathcal{M}_{\tau_0}^\tau} \epsilon^{-1} \bigg(
		(1+r)^{1-C_{[0,0]}\epsilon}|F^{(0)}_{(A,0)}|^2
		+ (1+r)^{\frac{1}{2}\delta}(1+\tau)^{1+\delta} |F^{(0)}_{(A,0,1)}|^2
		+ (1+r)^{1-3\delta}(1+\tau)^{2\beta} |F^{(0)}_{(A,0,2)}|^2
		\\
		&\phantom{\int_{\mathcal{M}_{\tau_0}^\tau} \epsilon^{-1} \bigg(}
		+ (1+r)^{1+\frac{1}{2}\delta} |F^{(0)}_{(A,0,3)}|^2
	\bigg)\dVol_g
	\lesssim
	\frac{1}{C_{[0,0]}} \tilde{\mathcal{E}}(1+\tau)^{-1}
	\\
	\\
	&\int_{\mathcal{M}_{\tau_0}^\tau \cap \{r \geq r_0\}} \epsilon^{-1}\bigg(
		r^{1-C_{[0,0]}\epsilon} (1+\tau)^{1+\delta} |F^{(0)}_{(A,0,4)}|^2
		+ r^{2-C_{[0,0]}\epsilon-2\delta} (1+\tau)^{2\beta} |F^{(0)}_{(A,0,5)}|^2
		\\
		&\phantom{\int_{\mathcal{M}_{\tau_0}^\tau \cap \{r \geq r_0\}} \epsilon^{-1}\bigg(}
		+ r^{2-C_{[0.0]}\epsilon} |F^{(0)}_{(A,0,6)}|^2
	\bigg)\dVol_g
	\lesssim
	\frac{1}{C_{[0,0]}} \tilde{\mathcal{E}}
	\end{split}
	\end{equation*}
	
	On the other hand, if $\phi_{(A)} \in \Phi_{[m]}$, then suppose that
	\begin{equation*}
	\begin{split}
		&\int_{\mathcal{M}_{\tau_0}^\tau} \epsilon^{-1} \bigg(
			(1+r)^{1-C_{[0,m]}\epsilon}|F^{(0)}_{(A,0)}|^2
			+ (1+r)^{\frac{1}{2}\delta}(1+\tau)^{1+\delta} |F^{(0)}_{(A,0,1)}|^2
			+ (1+r)^{1-3\delta}(1+\tau)^{2\beta} |F^{(0)}_{(A,0,2)}|^2
			\\
			&\phantom{\int_{\mathcal{M}_{\tau_0}^\tau} \epsilon^{-1} \bigg(}
			+ (1+r)^{1+\frac{1}{2}\delta} |F^{(0)}_{(A,0,3)}|^2
		\bigg)\dVol_g
		\lesssim
		\frac{1}{C_{[0,m]}} \tilde{\mathcal{E}}(1+\tau)^{-1+C_{(0,m)}\delta}
		\\
		\\
		&\int_{\mathcal{M}_{\tau_0}^\tau \cap \{r \geq r_0\}} \epsilon^{-1}\bigg(
			r^{1-C_{[0,m]}\epsilon} (1+\tau)^{1+\delta} |F^{(0)}_{(A,0,4)}|^2
			+ r^{2-C_{[0,m]}\epsilon-2\delta} (1+\tau)^{2\beta} |F^{(0)}_{(A,0,5)}|^2
			\\
			&\phantom{\int_{\mathcal{M}_{\tau_0}^\tau \cap \{r \geq r_0\}} \epsilon^{-1}\bigg(}
			+ r^{2-C_{[0.m]}\epsilon} |F^{(0)}_{(A,0,6)}|^2
		\bigg)\dVol_g
		\lesssim
		\frac{1}{C_{[0,m]}} \tilde{\mathcal{E}}(1+\tau)^{C_{(0,m)}\delta}
	\end{split}
	\end{equation*}

	Furthermore, suppose that both the pointwise bounds \emph{and} the $L^2$ bounds of chapter \ref{chapter bootstrap} hold.
	
	Then, for all sufficiently small $\epsilon$, we can decompose $F_{(A,0)}$ as
	\begin{equation*}
	\begin{split}
	F_{(A,0)}
	&=
	F_{(A,0,1)} + F_{(A,0,2)} + F_{(A,0,3)}
	\\
	&=
	F_{(A,0,4)} + F_{(A,0,5)} + F_{(A,0,6)}
	\end{split}	
	\end{equation*}
	where, if $\phi_{(A)} \in \Phi_{[0]}$, then 
	\begin{equation*}
	\begin{split}
	&\int_{\mathcal{M}_{\tau_0}^\tau} \epsilon^{-1} \bigg(
		(1+r)^{1-C_{[0,0]}\epsilon}|F_{(A,0)}|^2
		+ (1+r)^{\frac{1}{2}\delta}(1+\tau)^{1+\delta} |F_{(A,0,1)}|^2
		+ (1+r)^{1-3\delta}(1+\tau)^{2\beta} |F_{(A,0,2)}|^2
		\\
		&\phantom{\int_{\mathcal{M}_{\tau_0}^\tau} \epsilon^{-1} \bigg(}
		+ (1+r)^{1+\frac{1}{2}\delta} |F_{(A,0,3)}|^2
	\bigg)\dVol_g
	\lesssim
	\frac{1}{C_{[0,0]}} \tilde{\mathcal{E}}(1+\tau)^{-1}
	\\
	\\
	&\int_{\mathcal{M}_{\tau_0}^\tau \cap \{r \geq r_0\}} \epsilon^{-1}\bigg(
		r^{1-C_{[0,0]}\epsilon} (1+\tau)^{1+\delta} |F_{(A,0,4)}|^2
		+ r^{2-C_{[0,0]}\epsilon-2\delta} (1+\tau)^{2\beta} |F_{(A,0,5)}|^2
		\\
		&\phantom{\int_{\mathcal{M}_{\tau_0}^\tau \cap \{r \geq r_0\}} \epsilon^{-1}\bigg(}
		+ r^{2-C_{[0.0]}\epsilon} |F_{(A,0,6)}|^2
	\bigg)\dVol_g
	\lesssim
	\frac{1}{C_{[0,0]}} \tilde{\mathcal{E}}
	\end{split}
	\end{equation*}
	
	On the other hand, if $\phi_{(A)} \in \Phi_{[m]}$, then
	\begin{equation*}
	\begin{split}
	&\int_{\mathcal{M}_{\tau_0}^\tau} \epsilon^{-1} \bigg(
		(1+r)^{1-C_{[0,m]}\epsilon}|F_{(A,0)}|^2
		+ (1+r)^{\frac{1}{2}\delta}(1+\tau)^{1+\delta} |F_{(A,0,1)}|^2
		+ (1+r)^{1-3\delta}(1+\tau)^{2\beta} |F_{(A,0,2)}|^2
		\\
		&\phantom{\int_{\mathcal{M}_{\tau_0}^\tau} \epsilon^{-1} \bigg(}
		+ (1+r)^{1+\frac{1}{2}\delta} |F_{(A,0,3)}|^2
	\bigg)\dVol_g
	\lesssim
	\frac{1}{C_{[0,m]}} \tilde{\mathcal{E}}(1+\tau)^{-1+C_{(0,m)}\delta}
	\\
	\\
	&\int_{\mathcal{M}_{\tau_0}^\tau \cap \{r \geq r_0\}} \epsilon^{-1}\bigg(
		r^{1-C_{[0,m]}\epsilon} (1+\tau)^{1+\delta} |F_{(A,0,4)}|^2
		+ r^{2-C_{[0,m]}\epsilon-2\delta} (1+\tau)^{2\beta} |F_{(A,0,5)}|^2
		\\
		&\phantom{\int_{\mathcal{M}_{\tau_0}^\tau \cap \{r \geq r_0\}} \epsilon^{-1}\bigg(}
		+ r^{2-C_{[0,m]}\epsilon} |F_{(A,0,6)}|^2
	\bigg)\dVol_g
	\lesssim
	\frac{1}{C_{[0,m]}} \tilde{\mathcal{E}}(1+\tau)^{C_{(0,m)}\delta}
	\end{split}
	\end{equation*}

\end{proposition}

\begin{proof}

	We consider each of the terms in $F_{(A,0)}$ in turn, and we show that they obey the required bounds. Note that $F_{(A,0)}^{(0)}$ satisfies the $L^2$ bounds by assumption. In addition to assuming the pointwise bootstrap bounds of chapter \ref{chapter bootstrap}, we will also make use of the $L^2$ bootstrap bounds from sections \ref{section L2 bootstrap bounds} and \ref{section L2 bounds for geometric error terms}.
	
	First, we shall consider the bounds on $F$ which require decay in $\tau$. These bounds are related to the $T$-energy estimate, the Morawetz estimate and the $p$-weighted energy estimate with small $p$. First, we suppose that $\phi_{(A)} \in \Phi_{[0]}$. Then we have
	\begin{equation*}
	\begin{split}
	&\int_{\mathcal{M}_{\tau_0}^\tau} \epsilon^{-1}(1+r)^{1-C_{[0,0]}\epsilon} \left| \left( F^{(BC)}_{(A,0)} \right)^{\mu\nu} (\partial_\mu \phi_B) (\partial_\nu \phi_C) \right|^2 \dVol_g
	\\
	&\lesssim
	\int_{\mathcal{M}_{\tau_0}^\tau} \epsilon^{-1}(1+r)^{1-C_{[0,0]}\epsilon}\left(
		|\partial \phi|^2 |\bar{\partial} \phi|^2 
		\right) \dVol_g
	\end{split}
	\end{equation*}
	where we write $\phi$ to stand for any of the fields $\phi_{(A)}$. Substituting the pointwise bounds for $(\bar{\partial}\phi)$ and then the $L^2$ bootstrap bounds from chapter \ref{chapter bootstrap} we have
	\begin{equation*}
	\begin{split}
	&\int_{\mathcal{M}_{\tau_0}^\tau} \epsilon^{-1}(1+r)^{1-C_{[0,0]}\epsilon} \left| \left( F^{(BC)}_{(A,0)} \right)^{\mu\nu} (\partial_\mu \phi_B) (\partial_\nu \phi_C) \right|^2 \dVol_g
	\\
	&\lesssim
	\int_{\mathcal{M}_{\tau_0}^\tau} \epsilon(1+r)^{-1-2\delta-C_{[0,0]}\epsilon} |\partial \phi|^2 \dVol_g
	\\
	&\lesssim
	\tilde{\mathcal{E}} \epsilon \delta^{-1} (1+\tau)^{-1+C_{(0,0)}\delta}
	\end{split}
	\end{equation*}
	Since $\epsilon \ll \delta$ this gives the required bound in the case $\phi_{(A)} \in \Phi_{[0]}$.
	
	Similarly, we can use the pointwise bounds on $\phi$ and $\partial \phi$ to bound
	\begin{equation*}
	\begin{split}
	&\int_{\mathcal{M}_{\tau_0}^\tau} \epsilon^{-1}(1+r)^{1-C_{[0,0]}\epsilon} |\phi|^2 |\partial\phi|^4 \dVol_g
	\\
	&\lesssim
	\int_{\mathcal{M}_{\tau_0}^\tau} \epsilon^3 (1+r)^{-2+2\delta-C_{[0,0]}\epsilon} |\partial \phi|^2 \dVol_g
	\\
	&\lesssim
	\tilde{\mathcal{E}} \epsilon^3 \delta^{-1} (1+\tau)^{-1+C_{(0,0)}\delta}
	\end{split}
	\end{equation*}

	Next, suppose that $\phi_{(A)} \in \Phi_{[n]}$. Writing $\phi$ to stand for any of the fields $\phi_{(A)}$, and $\phi_{[n]}$ to stand for any of the fields in $\Phi_{[n]}$, we note that
	\begin{equation*}
	\begin{split}
	&\int_{\mathcal{M}_{\tau_0}^\tau} \epsilon^{-1}(1+r)^{1-C_{(\phi_{(A)})}\epsilon} \left| \left( F^{(BC)}_{(A,0)} \right)^{\mu\nu} (\partial_\mu \phi_B) (\partial_\nu \phi_C) \right|^2 \dVol_g
	\\
	&\lesssim
	\int_{\mathcal{M}_{\tau_0}^\tau} \epsilon^{-1}(1+r)^{1-C_{(\phi_{(A)})}\epsilon}\left( |\partial \phi|^2 |\bar{\partial} \phi|^2 + |\partial \phi_{[n-1]}|^4 + |\partial \phi_{[n]}|^2|\partial \phi_{[0]}|^2 \right) \dVol_g
	\end{split}
	\end{equation*}

	Substituting for the pointwise bounds from chapter \ref{chapter bootstrap} we have\footnote{There are several possible choices regarding \emph{which} terms to estimate in $L^\infty$ here.}
	\begin{equation*}
	\begin{split}
	&\int_{\mathcal{M}_{\tau_0}^\tau} \epsilon^{-1}(1+r)^{1-C_{[n]}\epsilon} \left| \left( F^{(BC)}_{(A,0)} \right)^{\mu\nu} (\partial_\mu \phi_B) (\partial_\nu \phi_C) \right|^2 \dVol_g
	\\
	&\lesssim
	\int_{\mathcal{M}_{\tau_0}^\tau} \epsilon\bigg(
	(1+r)^{-1-2\delta - C_{[n]}\epsilon}|\partial \phi|^2
	+ (1+r)^{-1 - C_{[n]}\epsilon + 2C_{(n-1)}\epsilon}|\partial \phi_{[n-1]}|^2
	\\
	&\phantom{\int_{\mathcal{M}_{\tau_0}^\tau} \epsilon\bigg(}
	+ (1+r)^{-1 - C_{[n]}\epsilon} |\partial \phi_{[n]}|^2	
	\bigg) \dVol_g
	\end{split}
	\end{equation*}
	and now making use of the $L^2$ bounds we find that, if
	\begin{equation*}
	C_{[n]} > C_{[n-1]} + 2C_{(n-1)}
	\end{equation*}
	then we have
	\begin{equation*}
	\begin{split}
	&\int_{\mathcal{M}_{\tau_0}^\tau} \epsilon^{-1}(1+r)^{1-C_{[n]}\epsilon} \left| \left( F^{(BC)}_{(A,0)} \right)^{\mu\nu} (\partial_\mu \phi_B) (\partial_\nu \phi_C) \right|^2 \dVol_g
	\\
	&\lesssim
	\left(\epsilon \delta^{-2} + \frac{1}{C_{[n-1]}} + \frac{1}{C_{[n]}} \right) \tilde{\mathcal{E}} (1+\tau)^{-1+C_{(0,0)}\delta}
	\end{split}
	\end{equation*}
	Hence, if $C_{[n-1]}$ is sufficiently large then we have the required bound.

	To bound the other terms, we make the following choices: terms involving at least one ``good'' derivative, together with cubic terms, are included in $F_{(A,0,2)}$ and $F_{(A,0,5)}$ respectively, while we include the ``bad'' terms, involving only bad derivatives, in $F_{(A,0,2)}$ and $F_{(A,0,4)}$. To bound the term corresponding to $F_{(A,0,2)}$ we need to bound
	\begin{equation*}
		\int_{M_{\tau}^{\tau_1}} \epsilon^{-1} (1+r)^{1-3\delta} (1+\tau)^{2\beta} \left| \left( F^{(BC)}_{(A,0)} \right)^{\mu\nu} (\partial_\mu \phi_B) (\partial_\nu \phi_C) \right|^2 \dVol_g
	\end{equation*}
	
	Again, we begin with the case $\phi_{(A)} \in \Phi_{[0]}$. Then we have
	\begin{equation*}
	\begin{split}
	&\int_{\mathcal{M}_{\tau}^{\tau_1}} \epsilon^{-1} (1+r)^{1-3\delta} (1+\tau)^{2\beta} \left| \left( F^{(BC)}_{(A,0)} \right)^{\mu\nu} (\partial_\mu \phi_B)(\partial_\nu \phi_C) \right|^2 \dVol_g
	\\
	&\lesssim
	\int_{\mathcal{M}_{\tau}^{\tau_1}} \epsilon^{-1} (1+r)^{1-3\delta} (1+\tau)^{2\beta}|\partial \phi|^2 |\bar{\partial} \phi|^2 \dVol_g
	\\
	&\lesssim
	\int_{\mathcal{M}_{\tau}^{\tau_1}} \epsilon (1+r)^{-1-5\delta} |\bar{\partial} \phi|^2 \dVol_g
	\end{split}
	\end{equation*}
	Now, using the bootstrap bounds we have
	\begin{equation*}
	\begin{split}
	&\int_{\mathcal{M}_{\tau}^{\tau_1}} \epsilon^{-1} (1+r)^{1-3\delta} (1+\tau)^{2\beta} \left| \left( F^{(BC)}_{(A,0)} \right)^{\mu\nu} (\partial_\mu \phi_B)(\partial_\nu \phi_C) \right|^2 \dVol_g
	\\
	&\lesssim
	\epsilon \tilde{\mathcal{E}} \delta^{-1} (1+\tau)^{-1}
	\end{split}
	\end{equation*}
	
	We need a similar bound for the case $\phi_{(A)} \in \Phi_{[n]}$, with $n \geq 1$. This time we have
	\begin{equation*}
	\begin{split}
	&\int_{\mathcal{M}_{\tau}^{\tau_1}} \epsilon^{-1} (1+r)^{1-3\delta} (1+\tau)^{2\beta} \left| \left( F^{(BC)}_{(A,0)} \right)^{\mu\nu} (\partial_\mu \phi_B)(\partial_\nu \phi_C) \right|^2 \dVol_g
	\\
	&\lesssim
	\int_{\mathcal{M}_{\tau}^{\tau_1}} \epsilon^{-1} (1+r)^{1-3\delta} (1+\tau)^{2\beta} \left(|\partial \phi|^2 |\bar{\partial} \phi|^2
	+ |\partial \phi_{[n-1]}|^4
	+ |\partial \phi_{[n]}|^2 |\partial \phi_{[0]}|^2
	\right)
	\dVol_g
	\\
	&\lesssim
	\int_{\mathcal{M}_{\tau}^{\tau_1}} \epsilon^{-1} (1+r)^{-1-\delta} |\partial \phi|^2
	\dVol_g
	\end{split}
	\end{equation*}
	This term can be dealt with in exactly the same way as before, so we have
	\begin{equation*}
	\begin{split}
	&\int_{\mathcal{M}_{\tau}^{\tau_1}} \epsilon^{-1} (1+r)^{1-3\delta} (1+\tau)^{2\beta} \left| \left( F^{(BC)}_{(A,0)} \right)^{\mu\nu} (\partial_\mu \phi_B)(\partial_\nu \phi_C) \right|^2 \dVol_g
	\lesssim
	\epsilon \tilde{\mathcal{E}} \delta^{-1} (1+\tau)^{-1 + C_{(0,n)}\delta}
	\end{split}
	\end{equation*}
	
	Note also that we can bound the cubic terms this way: we have
	\begin{equation*}
	\begin{split}
	&\int_{\mathcal{M}_{\tau}^{\tau_1}} \epsilon^{-1} (1+r)^{1-3\delta} (1+\tau)^{2\beta} |\phi|^2 |\partial \phi|^4 \dVol_g
	\\
	&\lesssim
	\int_{\mathcal{M}_{\tau}^{\tau_1}} \epsilon^{3} (1+r)^{-1-3\delta} |\bar{\partial} \phi|
	\dVol_g
	\\
	&\lesssim
	\epsilon^3 \tilde{\mathcal{E}} \delta^{-1} (1+\tau)^{-1+C_{(0,n)}\delta}
	\dVol_g
	\end{split}
	\end{equation*}
	where the factor of $(1+\tau)^{C_{(0,n)}\delta}$ is absent if $\phi \in \Phi_{[0]}$.

	Next, we bound the terms corresponding to $F_{(A, 0, 5)}$. Note that, as mentioned above, we only include ``good'' terms in $F_{(A,0,5)}$. First, using the pointwise bootstraps we have
	\begin{equation*}
	\int_{\mathcal{M}_\tau^{\tau_1}} \epsilon^{-1} (1+r)^{2-3\delta} (1+\tau)^{2\beta} |\partial \phi|^2 |\bar{\partial}\phi|^2 \dVol_g
	\lesssim
	\int_{\mathcal{M}_\tau^{\tau_1}} \epsilon (1+r)^{-\delta} |\bar{\partial}\phi|^2 \dVol_g
	\end{equation*}
	and so using the $L^2$ bootstrap bounds we have
	\begin{equation*}
	\int_{\mathcal{M}_\tau^{\tau_1}} \epsilon^{-1} (1+r)^{2-3\delta} (1+\tau)^{2\beta} |\partial \phi|^2 |\bar{\partial}\phi|^2 \dVol_g
	\lesssim
	\epsilon \tilde{\mathcal{E}}
	\end{equation*}
	
	Similarly, we can bound
	\begin{equation*}
	\int_{\mathcal{M}_\tau^{\tau_1}} \epsilon^{-1} (1+r)^{2-3\delta} (1+\tau)^{2\beta} |\phi|^2 |\partial \phi|^4 \dVol_g
	\lesssim
	\int_{\mathcal{M}_\tau^{\tau_1}} \epsilon^3 (1+r)^{-1-\delta} |\partial\phi|^2 \dVol_g
	\end{equation*}
	which can be bounded as before.
	
	Finally, we turn to those terms which can be bounded in $F_{(A,0,4)}$. Using the pointwise bootstrap bounds we have
	\begin{equation*}
	\begin{split}
	&\int_{\mathcal{M}_{\tau}^{\tau_1}} \epsilon^{-1} (1+r)^{1-C_{[0,m]}\epsilon} (1+\tau)^{1+\delta} \left( |\partial \phi_{[0]}|^2|\partial \phi_{[m]}|^2 + |\partial \phi_{[m-1]}|^2|\partial \phi_{[m-1]}|^2 \right) \dVol_g
	\\
	&\lesssim
	\int_{\mathcal{M}_{\tau}^{\tau_1}} \epsilon^{-1} (1+r)^{1-C_{[0,m]}\epsilon} (1+\tau)^{1+\delta} \left(
		\epsilon^2 (1+r)^{-2} |\partial \phi_{[m]}|^2
		+ \epsilon^2 (1+r)^{-2 + 2C_{(0,m-1)}\epsilon} |\partial \phi_{[m-1]}|^2
	\right) \dVol_g
	\\
	&\lesssim
	\int_{\mathcal{M}_{\tau}^{\tau_1}} (1+\tau)^{1+\delta} \epsilon   \bigg(
		(1+r)^{-1-C_{[0,m]}\epsilon} |\partial \phi_{[m]}|^2
		+ (1+r)^{-1-C_{[0,m]}\epsilon + 2C_{(0,m-1)}\epsilon} |\partial \phi_{[m-1]}|^2
	\bigg) \dVol_g
	\end{split}
	\end{equation*}
	Now, we use the fact that $C_{[0,m]} - 2C_{(0,m-1)} \geq C_{[0,m-1]}$ and break the interval $[\tau, \tau_1]$ into diadic pieces to obtain
	\begin{equation*}
	\begin{split}
	&\int_{\mathcal{M}_{\tau}^{\tau_1}} \epsilon^{-1} (1+r)^{1-C_{[0,m]}\epsilon} (1+\tau)^{1+\delta} \left( |\partial \phi_{[0]}|^2|\partial \phi_{[m]}|^2 + |\partial \phi_{[m-1]}|^2|\partial \phi_{[m-1]}|^2 \right) \dVol_g
	\\
	&\lesssim
	\sum_{j=0}^{\lceil \log(\tau_1 - \tau + 1) \rceil - 1} \int_{\mathcal{M}_{\tau + e^j - 1}^{\tau + e^{j+1} - 1}} (1+\tau)^{1+\delta} \epsilon   \bigg(
		(1+r)^{-1-C_{[0,m]}\epsilon} |\partial \phi_{[m]}|^2
		\\
		&\phantom{\lesssim \sum_{j=0}^{\lceil \log(\tau_1 - \tau + 1) \rceil - 1} \int_{\mathcal{M}_{\tau + e^j - 1}^{\tau + e^{j+1} - 1}} (1+\tau)^{1+\delta} \epsilon \bigg(}
		+ (1+r)^{-1-C_{[0,m]}\epsilon + 2C_{(0,m-1)}\epsilon} |\partial \phi_{[m-1]}|^2
	\bigg) \dVol_g
	\\
	&\lesssim
	\sum_{j=0}^{\lceil \log(\tau_1 - \tau + 1) \rceil - 1}
		(\tau + e^{j+1}) \int_{\mathcal{M}_{\tau + e^j - 1}^{\tau + e^{j+1} - 1}} \epsilon \bigg(
			(1+r)^{-1-C_{[0,m]}\epsilon} |\partial \phi_{[m]}|^2
			\\
			&\phantom{\lesssim \sum_{j=0}^{\lceil \log(\tau_1 - \tau + 1) \rceil - 1}	(\tau + e^{j+1}) \int_{\mathcal{M}_{\tau + e^j - 1}^{\tau + e^{j+1} - 1}} \epsilon \bigg(}
			+ (1+r)^{-1-C_{[0,m-1]}\epsilon} |\partial \phi_{[m-1]}|^2
		\bigg) \dVol_g
	\end{split}
	\end{equation*}
	and so, using the $L^2$ bootstrap bounds we have
	\begin{equation*}
	\begin{split}
	&\int_{\mathcal{M}_{\tau}^{\tau_1}} \epsilon^{-1} (1+r)^{1-C_{[0,m]}\epsilon} (1+\tau)^{1+\delta} \left( |\partial \phi_{[0]}|^2|\partial \phi_{[m]}|^2 + |\partial \phi_{[m-1]}|^2|\partial \phi_{[m-1]}|^2 \right) \dVol_g
	\\
	&\lesssim
	\epsilon \tilde{\mathcal{E}} \sum_{j=0}^{\lceil \log(\tau_1 - \tau + 1) \rceil - 1} (\tau + e^{j+1}) (\tau + e^{j}-1)^{-1+C_{(0,m)}\delta}
	\\
	&\lesssim
	\epsilon \tilde{\mathcal{E}} \sum_{j=0}^{\lceil \log(\tau_1 - \tau + 1) \rceil - 1} (\tau + e^{j+1})^{C_{(0,m)}\delta}
	\\
	&\lesssim
	\frac{1}{C_{(0,m)}\delta}\epsilon \tilde{\mathcal{E}} (1 + \tau_1)^{C_{(0,m)}\delta}
	\end{split}
	\end{equation*}
		
\end{proof}

\section{Expressions for the inhomogeneous terms after commuting}
\label{section expressions for the inhomogeneous terms after commuting}

In this section, we will begin with a system of equations satisfying
\begin{equation*}
\tilde{\Box}_g \phi_{(a)} = F_{(a)}
\end{equation*}
We then consider the system formed by commuting these equations with the operators $\mathscr{Z}$ or $\mathscr{Y}$ some number of times, forming a system of equations of the form
\begin{equation*}
\tilde{\Box}_g \mathscr{Z}^n \phi_{(a)} = \tilde{F}_{(a,n)}
\end{equation*}
or
\begin{equation*}
\tilde{\Box}_g \mathscr{Y}^n \phi_{(a)} 
	- \slashed{\Delta} \mathscr{Y}^{n-1} \phi_{(a)}
	- (2^k-1)r^{-1} \slashed{\D}_L (r\slashed{\D}_L \mathscr{Y}^{n-1} \phi_{(a)})
	- (2^k-1)r^{-1} \slashed{\D}_L ( \mathscr{Y}^{n-1} \phi_{(a)})
	= \tilde{F}_{(rL, a, n)}
\end{equation*}
where $k$ labels the number of times the operator $r\slashed{\D}_L$ appears in the operator expansion of $\mathscr{Y}^n$.

We are aiming to show that the $\tilde{F}_{(a,n)}$ and $\tilde{F}_{(rL, a, n)}$ satisfy suitable $L^2$ bounds, so that we can use the energy estimates proved in chapter \ref{chapter boundedness and energy decay} and appendix \ref{appendix improved energy decay}. Before we can do this, we first need to discover the structure of the inhomogeneous terms $\tilde{F}_{(n)}$ and $\tilde{F}_{(rL, n)}$.

\begin{proposition}[The structure of the inhomogeneous terms after commuting once with $\slashed{\D}_T$]
	\label{proposition inhomogeneous terms commute with T}
	Suppose that $\phi$ is an $S_{\tau,r}$-tangent tensor field satisfying the equation
	\begin{equation*}
	\tilde{\Box}_g \phi = F
	\end{equation*}
	
	Then, if $\phi$ is a scalar field, $T\phi$ satisfies the following equation, given schematically:
	\begin{equation}
	\begin{split}
	\tilde{\Box}_g (T \phi) 
	&=
	(T F)
	+ \begin{pmatrix}
		r^{-1} (\slashed{\D} \mathscr{Z} h)_{(\text{frame})}  \\
		(\overline{\slashed{\D}} \mathscr{Z} h)_{(\text{frame})}  \\
		(\tilde{\Box}_g h)_{(\text{frame})} \\
		\bm{\Gamma}^{(1)}_{(-2+2C_{(1)}\epsilon)}
	\end{pmatrix} (\D \phi)
	+ \begin{pmatrix}
		\mathscr{Z} \tr_{\slashed{g}} \chi_{(\text{small})} \\
		(\slashed{\D} \mathscr{Z} h)_{(\text{frame})}
	\end{pmatrix}(\overline{\slashed{\D}} \phi)
	\\
	&\phantom{=}
	+ \bm{\Gamma}^{(0)}_{(-1)} \left(\D (T \phi) \right)
	+ r^{-1} \begin{pmatrix}
		\slashed{\nabla} \log \mu \\
		\zeta \\
		(\chi_{(\text{small})} + \chibar_{(\text{small})})
	\end{pmatrix} (\slashed{\D}\mathscr{Z} \phi)
	\end{split}
	\end{equation}
	
	On the other hand, if $\phi$ is a higher order $S_{\tau,r}$-tangent tensor field, then $\slashed{\D}_T \phi$ satisfies the following equation, also given schematically:
	\begin{equation}
	\begin{split}
	\tilde{\Box}_g (\slashed{\D}_T \phi) 
	&=
	(\slashed{\D}_T F)
	+ \begin{pmatrix}
	r^{-1} (\slashed{\D} \mathscr{Z} h)_{(\text{frame})}  \\
	(\overline{\slashed{\D}} \mathscr{Z} h)_{(\text{frame})}  \\
	(\tilde{\Box}_g h)_{(\text{frame})} \\
	\bm{\Gamma}^{(1)}_{(-2+2C_{(1)}\epsilon)}
	\end{pmatrix} (\slashed{\D} \phi)
	+ \begin{pmatrix}
	\mathscr{Z} \tr_{\slashed{g}} \chi_{(\text{small})} \\
	(\slashed{\D} \mathscr{Z} h)_{(\text{frame})}
	\end{pmatrix}(\overline{\slashed{\D}} \phi)
	+ \bm{\Gamma}^{(0)}_{(-1)} (\slashed{\D}\slashed{\D}_T \phi)
	\\
	&\phantom{=}
	+ r^{-1} \begin{pmatrix}
	\slashed{\nabla} \log \mu \\
	\zeta \\
	(\chi_{(\text{small})} + \chibar_{(\text{small})})
	\end{pmatrix} (\slashed{\D}\mathscr{Z} \phi)
	+ \begin{pmatrix}
	r^{-1} \left( \tilde{\slashed{\Box}}_g (\mathscr{Z} h) \right)_{(\text{frame})} \\
	r^{-1} \left(\slashed{\D}(\mathscr{Z}^2 h)\right)_{(\text{frame})} \\
	\bm{\Gamma}^{(1)}_{(-1+C_{(1)}\epsilon)} \cdot (\slashed{\D} \mathscr{Z} h)_{(\text{frame})} \\
	\bm{\Gamma}^{(1)}_{(-1+C_{(1)}\epsilon)} \cdot (\tilde{\Box}_g h)_{(\text{frame})} \\
	\bm{\Gamma}^{(1)}_{(-1+C_{(1)}\epsilon)} \cdot \mathscr{Z} \chi_{(\text{small})} \\
	\bm{\Gamma}^{(1)}_{(-2+C_{(1)}\epsilon)} \mathscr{Z}^2 \log \mu \\
	\bm{\Gamma}^{(1)}_{(-3 + 3C_{(1)}\epsilon)}
	\end{pmatrix} \phi	
	\end{split}
	\end{equation}

\end{proposition}

\begin{proof}
	First, recall proposition \ref{proposition commute vector field}. This immediately lead to the  expression
	\begin{equation*}
	\begin{split}
	\tilde{\slashed{\Box}}_g \left( \slashed{\D}_T \phi \right)
	&=
	\slashed{\D}_T F
	+ \slashed{\D}_\mu \left( {^{(T)}\mathscr{J}}^\mu \right)
	+ \frac{1}{2}(\tr_{g} {^{(T)}\pi}) \tilde{\slashed{\Box}}_g \phi
	- (T\omega) \slashed{\D}_{\Lbar} \phi
	+ \omega \slashed{\D}_{[\Lbar, T]} \phi
	- \frac{1}{2}\omega (\tr_{g} {^{(T)}\pi}) \slashed{\D}_{\Lbar} \phi
	\\
	&\phantom{=}
	+ [\slashed{\D}_\mu \, , \slashed{\D}_\nu] (T^\nu \slashed{\D}^\mu \phi)
	+ \slashed{\D}^\mu \left( T^\nu [\slashed{\D}_\mu \, , \slashed{\D}_\nu]\phi \right)
	+ \omega \Lbar^\mu T^\nu [\slashed{\D}_\mu \, , \slashed{\D}_\nu] \phi
	\end{split}
	\end{equation*}
	
	Now, using proposition \ref{proposition commute T} we have, schematically,
	\begin{equation*}
	\begin{split}
	&
	\slashed{\D}_\mu \left( {^{(T)}\mathscr{J}}^\mu \right)
	+ \frac{1}{2}(\tr_{g} {^{(T)}\pi}) \tilde{\slashed{\Box}}_g \phi
	- (T\omega) \slashed{\D}_{\Lbar} \phi
	+ \omega \slashed{\D}_{[\Lbar, T]} \phi
	- \frac{1}{2}\omega (\tr_{g} {^{(T)}\pi}) \slashed{\D}_{\Lbar} \phi
	\\
	&=
	\begin{pmatrix}
		r^{-1} (\slashed{\D} \mathscr{Z} h)_{(\text{frame})}  \\
		L \tr_{\slashed{g}} \chi_{(\text{small})} \\
		\bm{\Gamma}^{(1)}_{(-2+2C_{(1)}\epsilon)}
	\end{pmatrix} (\slashed{\D} \phi)
	+ \begin{pmatrix}
		\mathscr{Z} \tr_{\slashed{g}} \chi_{(\text{small})} \\
		\slashed{\Delta} \log \mu \\
		\slashed{\Div} \hat{\chi} \\
		(\slashed{\D} \mathscr{Z} h)_{(\text{frame})}
	\end{pmatrix}(\overline{\slashed{\D}} \phi)
	+ \bm{\Gamma}^{(0)}_{(-1)} (\slashed{\D}\slashed{\D}_T \phi)
	\\
	&\phantom{=}
	+ r^{-1}\begin{pmatrix}
		|\slashed{\nabla} \log \mu| \\
		|\zeta| \\
		|\chi_{(\text{small})} + \chibar_{(\text{small})}|
	\end{pmatrix}(\slashed{\D}\mathscr{Z} \phi)
	+ \begin{pmatrix}
		\bm{\Gamma}^{(1)}_{(-2+C_{(1)}\epsilon)} (\slashed{\D} \mathscr{Z} h)_{(\text{frame})} \\
		\bm{\Gamma}^{(1)}_{(-3+3C_{(1)}\epsilon)}
	\end{pmatrix} \phi
	\end{split}
	\end{equation*}
	Some of these terms can be further simplified. From proposition \ref{proposition transport trace chi} it is fairly easy to see that, schematically,
	\begin{equation*}
	L\tr_{\slashed{g}}\chi_{(\text{small})}
	=
	(\overline{\slashed{\D}} \mathscr{Z} h)_{(\text{frame})}
	+ (\tilde{\Box}_g h)_{(\text{frame})}
	+ \bm{\Gamma}^{(1)}_{(-2 -\delta)}
	\end{equation*}
	Also, from proposition \ref{proposition spherical laplacian mu} we have
	\begin{equation*}
	\slashed{\Delta} \log \mu
	=
	\mathscr{Z} \tr_{\slashed{g}}\chi_{(\text{small})}
	+ (\tilde{\Box}_g h)_{(\text{frame})}
	+ r^{-1} (\slashed{\D} \mathscr{Z} h)_{(\text{frame})}
	+ \bm{\Gamma}^{(1)}_{(-2+2C_{(1)}\epsilon)}
	\end{equation*}
	Finally, from proposition \ref{proposition div chihat} we have
	\begin{equation*}
	\slashed{\Div} \hat{\chi}
	=
	r^{-1} \mathscr{Z} \tr_{\slashed{g}}\chi_{(\text{small})}
	+ r^{-1} (\slashed{\D} \mathscr{Z} h)_{(\text{frame})}
	+ \bm{\Gamma}^{(1)}_{(-2 - \delta)}
	\end{equation*}
	Putting this all together leads to
	\begin{equation*}
	\begin{split}
	&
	\slashed{\D}_\mu \left( {^{(T)}\mathscr{J}}^\mu \right)
	+ \frac{1}{2}(\tr_{g} {^{(T)}\pi}) \tilde{\slashed{\Box}}_g \phi
	- (T\omega) \slashed{\D}_{\Lbar} \phi
	+ \omega \slashed{\D}_{[\Lbar, T]} \phi
	- \frac{1}{2}\omega (\tr_{g} {^{(T)}\pi}) \slashed{\D}_{\Lbar} \phi
	\\
	&=
	\begin{pmatrix}
	r^{-1} (\slashed{\D} \mathscr{Z} h)_{(\text{frame})}  \\
	(\overline{\slashed{\D}} \mathscr{Z} h)_{(\text{frame})}  \\
	(\tilde{\Box}_g h)_{(\text{frame})} \\
	\bm{\Gamma}^{(1)}_{(-2+2C_{(1)}\epsilon)}
	\end{pmatrix} (\slashed{\D} \phi)
	+ \begin{pmatrix}
	\mathscr{Z} \tr_{\slashed{g}} \chi_{(\text{small})} \\
	(\slashed{\D} \mathscr{Z} h)_{(\text{frame})}
	\end{pmatrix}(\overline{\slashed{\D}} \phi)
	+ \bm{\Gamma}^{(0)}_{(-1)} (\slashed{\D}\slashed{\D}_T \phi)
	\\
	&\phantom{=}
	+ r^{-1}\begin{pmatrix}
		|\slashed{\nabla} \log \mu| \\
		|\zeta| \\
		|\chi_{(\text{small})} + \chibar_{(\text{small})}|
	\end{pmatrix}(\slashed{\D}\mathscr{Z} \phi)
	+ \begin{pmatrix}
	\bm{\Gamma}^{(0)}_{(-2+C_{(1)}\epsilon)} (\slashed{\D} \mathscr{Z} h)_{(\text{frame})} \\
	\bm{\Gamma}^{(0)}_{(-3+3C_{(1)}\epsilon)}
	\end{pmatrix} \phi
	\end{split}
	\end{equation*}
	Note that the terms involving $\phi$ (and not its derivatives) are absent if $\phi$ is a scalar field.
	
	The remaining terms also vanish if $\phi$ is a scalar field. Using proposition \ref{proposition additional error term commuting with T} we have
	\begin{equation}
	\begin{split}
	&[\slashed{\D}_\mu \, , \slashed{\D}_\nu] \left(T^\nu \slashed{\D}^\mu \phi \right)
	+ \omega \Lbar^\mu T^\nu [\slashed{\D}_\mu \, , \slashed{\D}_\nu] \phi
	+ \slashed{\D}^\mu \left( T^\nu [\slashed{\D}_\mu \, , \slashed{\D}_\nu]\phi \right)
	\\
	&=
	\begin{pmatrix}
	r^{-1} \left(\slashed{\D}(\mathscr{Z} h)\right)_{(\text{frame})} \\
	\bm{\Gamma}^{(1)}_{(-2 + 2C_{(1)}\epsilon)} \\
	\end{pmatrix} (\slashed{\D} \phi)
	+ \begin{pmatrix}
	r^{-1} \left( \tilde{\slashed{\Box}}_g (\mathscr{Z} h) \right)_{(\text{frame})} \\
	r^{-1} \left(\slashed{\D}(\mathscr{Z}^2 h)\right)_{(\text{frame})} \\
	\bm{\Gamma}^{(1)}_{(-1+C_{(1)}\epsilon)} \cdot (\slashed{\D} \mathscr{Z} h)_{(\text{frame})} \\
	\bm{\Gamma}^{(1)}_{(-1+C_{(1)}\epsilon)} \cdot (\tilde{\Box}_g h)_{(\text{frame})} \\
	\bm{\Gamma}^{(1)}_{(-1+C_{(1)}\epsilon)} \cdot \mathscr{Z} \chi_{(\text{small})} \\
	\bm{\Gamma}^{(1)}_{(-2+C_{(1)}\epsilon)} \mathscr{Z}^2 \log \mu \\
	\bm{\Gamma}^{(1)}_{(-3 + 3C_{(1)}\epsilon)}
	\end{pmatrix} \phi
	\end{split}
	\end{equation}
	
	Putting these calculations together proves the proposition.
	
\end{proof}

\begin{proposition}[The structure of the inhomogeneous terms after commuting once with $r\slashed{\nabla}$]
	\label{proposition inhomogeneous terms commute with rnabla}
	Suppose that $\phi$ is an $S_{\tau,r}$-tangent tensor field satisfying the equation
	\begin{equation*}
	\tilde{\Box}_g \phi = F
	\end{equation*}
	
	Then, if $\phi$ is a scalar field, $r\slashed{\nabla}\phi$ satisfies the following equation, given schematically:
	\begin{equation*}
	\begin{split}
	\tilde{\slashed{\Box}}_g (r\slashed{\nabla}\phi) 
	&= 
	r\slashed{\nabla}F
	+ \begin{pmatrix}
		\mathscr{Z} \tr_{\slashed{g}}\chi_{(\text{small})} \\
		(\slashed{\D} \mathscr{Z} h)_{LL} \\
		(\overline{\slashed{\D}} \mathscr{Z} h)_{(\text{frame})} \\
		\bm{\Gamma}^{(1)}_{(-1+C_{(1)}\epsilon)}
	\end{pmatrix} (\partial \phi)
	+ \begin{pmatrix}
		r(\tilde{\Box}_g h)_{(\text{frame})} \\
		r\slashed{\nabla}^2 \log \mu \\
		\bm{\Gamma}^{(1)}_{(-1+2C_{(1)}\epsilon)}
	\end{pmatrix} (\slashed{\nabla}\phi)
	+ \bm{\Gamma}^{(0)}_{(-1-\delta)} (\slashed{\D} \mathscr{Z}\phi) 
	\\
	&\phantom{=}
	+ \bm{\Gamma}^{(1)}_{(-1+C_{(1)}\epsilon)} (\overline{\slashed{\D}} \mathscr{Z}\phi)
	+ \bm{\Gamma}^{(1)}_{(-1+C_{(1)}\epsilon)} (L(rL\phi))
	\end{split}
	\end{equation*}
	
	On the other hand, if $\phi$ is a higher rank tensor field, then $r\slashed{\nabla}$ satisfies the an equation of the following schematic form:
	\begin{equation*}
	\begin{split}
	\tilde{\slashed{\Box}}_g (r\slashed{\nabla}\phi) 
	&= 
	r\slashed{\nabla}F
	+ \begin{pmatrix}
		\mathscr{Z} \tr_{\slashed{g}}\chi_{(\text{small})} \\
		(\slashed{\D} \mathscr{Z} h)_{LL} \\
		(\overline{\slashed{\D}} \mathscr{Z} h)_{(\text{frame})} \\
		\bm{\Gamma}^{(1)}_{(-1+2C_{(1)}\epsilon)}
	\end{pmatrix} (\slashed{\D} \phi)
	+ \begin{pmatrix}
		r(\tilde{\Box}_g h)_{(\text{frame})} \\
		r\slashed{\nabla}^2 \log \mu \\
		\bm{\Gamma}^{(1)}_{(-1+2C_{(1)}\epsilon)}
	\end{pmatrix} (\slashed{\nabla}\phi)
	+ \bm{\Gamma}^{(1)}_{(-1-\delta)} (\slashed{\D} \mathscr{Z}\phi) 
	\\
	&\phantom{=}
	+ \bm{\Gamma}^{(1)}_{(-1+C_{(1)}\epsilon)} (\overline{\slashed{\D}} \mathscr{Z}\phi)
	+ \bm{\Gamma}^{(1)}_{(-1+C_{(1)}\epsilon)} (\slashed{\D}_L(r\slashed{\D}_L\phi))
	+ \begin{pmatrix}
		r^{-1} \slashed{\D}(\mathscr{Z}^2 h)_{(\text{frame})} \\
		\tilde{\slashed{\Box}}_g (\mathscr{Z}h)_{(\text{frame})} \\
		\bm{\Gamma}^{(1)}_{(-1+C_{(1)}\epsilon)} (\tilde{\Box}_g h)_{(\text{frame})} \\
		\bm{\Gamma}^{(1)}_{(-1+C_{(1)}\epsilon)} \mathscr{Z}^2 \log \mu \\
		\bm{\Gamma}^{(1)}_{(-1+C_{(1)}\epsilon)} \mathscr{Z} \chi_{(\text{small})} \\
		\bm{\Gamma}^{(1)}_{(-2+3C_{(1)}\epsilon)} 
	\end{pmatrix} \cdot (\phi)
	\end{split}
	\end{equation*}
	
\end{proposition}

\begin{proof}
	We recall proposition \ref{proposition commute tensor field}, which gives the expression
	\begin{equation*}
	\begin{split}
	\tilde{\slashed{\Box}}_g \left( r\slashed{\nabla}_{\slashed{\alpha}} \phi \right)
	&=
	r\slashed{\nabla}_{\slashed{\alpha}} F
	+ \slashed{\D}_\mu \left( {^{(r\slashed{\Pi})}\mathscr{J}[\phi]_{\slashed{\alpha}}^{\phantom{\slashed{\alpha}}\mu}} \right)
	+ \frac{1}{2}\left( {^{(r\slashed{\Pi})}\pi_{\slashed{\alpha}\mu}^{\phantom{\slashed{\alpha}\mu}\mu}}\right) \tilde{\slashed{\Box}}_g \phi
	- (r\slashed{\nabla}_{\slashed{\alpha}} \omega)\slashed{\D}_{\Lbar}\phi
	\\
	&\phantom{=}
	+ \omega \left( \left(\slashed{\D}_{\Lbar} (r\slashed{\Pi}_{\slashed{\alpha}}^{\phantom{\slashed{\alpha}}\mu}) \right)
	- r\slashed{\nabla}_{\slashed{\alpha}} \Lbar^\mu \right)
	\slashed{\D}_\mu \phi
	- \frac{1}{2} \omega \left({^{(r\slashed{\Pi})}\pi_{\slashed{\alpha}\mu}^{\phantom{\slashed{\alpha}\mu}\mu}}\right) \slashed{\D}_{\Lbar} \phi
	+ [\slashed{\D}_\mu \, , \slashed{\D}_\nu]\left( r\slashed{\Pi}_{\slashed{\alpha}}^{\phantom{\slashed{\alpha}}\nu} \slashed{\D}^\mu \phi \right) 
	\\
	&\phantom{=}
	+ \slashed{\D}^\mu\left(r\slashed{\Pi}_{\slashed{\alpha}}^{\phantom{\slashed{\alpha}}\nu} [\slashed{\D}_\mu \, , \slashed{\D}_\nu]\phi \right)
	+ r\omega \Lbar^\mu \Pi_{\slashed{\alpha}}^{\phantom{\slashed{\alpha}}\nu} [\slashed{\D}_\mu \, , \slashed{\D}_\nu] \phi 
	\end{split}
	\end{equation*}
	
	Now, we can use proposition \ref{proposition commute rnabla} to express the main error terms arising from commuting. First, suppose that $\phi$ is a scalar field. Then we have
	\begin{equation*}
	\begin{split}
	&
	\slashed{\D}_\mu \left( {^{(r\slashed{\Pi})}\mathscr{J}[\phi]_{\slashed{\alpha}}^{\phantom{\slashed{\alpha}}\mu}} \right)
	+ \frac{1}{2}\left( {^{(r\slashed{\Pi})}\pi_{\slashed{\alpha}\mu}^{\phantom{\slashed{\alpha}\mu}\mu}}\right) \tilde{\Box}_g \phi
	- (r\slashed{\nabla}_{\slashed{\alpha}} \omega) (\Lbar \phi)
	+ \omega \left( \left(\slashed{\D}_{\Lbar} (r\slashed{\Pi}_{\slashed{\alpha}}^{\phantom{\slashed{\alpha}}\mu}) \right)
	- r\slashed{\nabla}_{\slashed{\alpha}} \Lbar^\mu \right) \D_\mu \phi
	\\
	&
	- \frac{1}{2} \omega \left({^{(r\slashed{\Pi})}\pi_{\slashed{\alpha}\mu}^{\phantom{\slashed{\alpha}\mu}\mu}}\right) (\Lbar \phi)
	\\
	&=
	\begin{pmatrix}
	\mathscr{Z} \tr_{\slashed{g}}\chi_{(\text{small})} \\
	r (\slashed{\Div} \hat{\chi}) \\
	(\slashed{\D} \mathscr{Z} h)_{LL} \\
	(\overline{\slashed{\D}} \mathscr{Z} h)_{(\text{frame})} \\
	\bm{\Gamma}^{(1)}_{(-1+2C_{(1)}\epsilon)}
	\end{pmatrix} (\partial\phi)
	+ \begin{pmatrix}
	r\slashed{\D}_{\Lbar} \chi_{(\text{small})} + r\slashed{\D}_L \chibar_{(\text{small})} 
	\\
	r\slashed{\nabla}^2 \log \mu
	\end{pmatrix} (\slashed{\nabla}\phi)
	+ \bm{\Gamma}^{(1)}_{(-1-\delta)} (\slashed{\D} \mathscr{Z} \phi)
	+ \bm{\Gamma}^{(1)}_{(-1+C_{(1)}\epsilon)} (\overline{\slashed{\D}} \mathscr{Z} \phi)
	\\
	&\phantom{=}
	+ \bm{\Gamma}^{(1)}_{(-1+C_{(1)}\epsilon)} L(r L \phi)
	\end{split}
	\end{equation*}
	As above, we can use proposition \ref{proposition div chihat} to write
	\begin{equation*}
	r\slashed{\Div} \hat{\chi} = \mathscr{Z} \tr_{\slashed{g}}\chi_{(\text{small})} + (\slashed{\D}\mathscr{Z} h)_{(\text{frame})} + \bm{\Gamma}^{(0)}_{(-1-\delta)}
	\end{equation*}
	We can also use propositions \ref{proposition transport chibar} and \ref{proposition Lbar chi} to express, schematically,
	\begin{equation*}
	\begin{split}
	r\slashed{\D}_{\Lbar} \chi_{(\text{small})} + r\slashed{\D}_L \chibar_{(\text{small})}
	&=
	\slashed{\Pi}_\mu^{\phantom{\mu}a}\slashed{\Pi}_\nu^{\phantom{\mu}a} R_{L\mu\Lbar\nu}
	+ r\slashed{\nabla}\zeta
	+ r\slashed{\nabla}^2 \log \mu
	+ \bm{\Gamma}^{(1)}_{(-1+2C_{(1)}\epsilon)}
	\\
	&=
	r(\tilde{\Box}_g h)_{(\text{frame})}
	+ (\slashed{\D} \mathscr{Z} h)_{(\text{frame})}
	+ r\slashed{\nabla}^2 \log \mu
	+ \bm{\Gamma}^{(1)}_{(-1+2C_{(1)}\epsilon)}
	\end{split}
	\end{equation*}
	
	Unlike in the case of commuting with a vector field, there are some additional error terms arising even in the case that $\phi$ is a scalar. These are produced by the term
	\begin{equation*}
	[\slashed{\D}_\mu, \slashed{\D}_\nu] \left( r\slashed{\Pi}_{\slashed{\alpha}}^{\phantom{\slashed{\alpha}}\nu} \D^\mu \phi \right)
	=
	-\frac{1}{2} r\Omega_{L\slashed{\beta} \slashed{\alpha}}^{\phantom{L\slashed{\beta} \slashed{\alpha}} \slashed{\beta}} (\Lbar \phi)
	-\frac{1}{2} r\Omega_{\Lbar\slashed{\beta} \slashed{\alpha}}^{\phantom{L\slashed{\beta} \slashed{\alpha}} \slashed{\beta}} (L \phi)
	+ r\Omega_{\slashed{\gamma}\slashed{\beta} \slashed{\alpha}}^{\phantom{L\slashed{\beta} \slashed{\alpha}} \slashed{\beta}} (\slashed{\nabla}^{\slashed{\gamma}} \phi)
	\end{equation*}
	Note that the terms involving the Riemann curvature $R$ (rather than the curvature $\Omega$ of the vector bundle $B$) vanish by the symmetries of the Riemann tensor. The terms above can be given schematically, using the expression in proposition \ref{proposition expression for Omega}, as
	\begin{equation*}
	[\slashed{\D}_\mu, \slashed{\D}_\nu] \left( r\slashed{\Pi}_{\slashed{\alpha}}^{\phantom{\slashed{\alpha}}\nu} \D^\mu \phi \right)
	= 
	\begin{pmatrix}
		(\overline{\slashed{\D}} \mathscr{Z} h) \\
		\bm{\Gamma}^{(1)}_{(-1-\delta)}
	\end{pmatrix} (\partial\phi)
	+ \bm{\Gamma}^{(1)}_{(-1+2C_{(1)}\epsilon)} (\bar{\partial} \phi)
	\end{equation*}

	On the other hand, if $\phi$ is a higher rank $S_{\tau,r}$-tangent tensor field, then we have
	\begin{equation*}
	\begin{split}
	&
	\slashed{\D}_\mu \left( {^{(r\slashed{\Pi})}\mathscr{J}[\phi]_{\slashed{\alpha}}^{\phantom{\slashed{\alpha}}\mu}} \right)
	+ \frac{1}{2}\left( {^{(r\slashed{\Pi})}\pi_{\slashed{\alpha}\mu}^{\phantom{\slashed{\alpha}\mu}\mu}}\right) \tilde{\slashed{\Box}}_g \phi
	- (r\slashed{\nabla}_{\slashed{\alpha}} \omega)\slashed{\D}_{\Lbar}\phi
	+ \omega \left( \left(\slashed{\D}_{\Lbar} (r\slashed{\Pi}_{\slashed{\alpha}}^{\phantom{\slashed{\alpha}}\mu}) \right)
	- r\slashed{\nabla}_{\slashed{\alpha}} \Lbar^\mu \right) \slashed{\D}_\mu \phi
	\\
	&
	- \frac{1}{2} \omega \left({^{(r\slashed{\Pi})}\pi_{\slashed{\alpha}\mu}^{\phantom{\slashed{\alpha}\mu}\mu}}\right) \slashed{\D}_{\Lbar} \phi
	\\
	&=
	\begin{pmatrix}
		\mathscr{Z} \tr_{\slashed{g}}\chi_{(\text{small})} \\
		r (\slashed{\Div} \hat{\chi}) \\
		(\slashed{\D} \mathscr{Z} h)_{(\text{frame})} \\
		\bm{\Gamma}^{(1)}_{(-1+2C_{(1)}\epsilon)}
	\end{pmatrix} (\slashed{\D}\phi)
	+ \begin{pmatrix}
		r\slashed{\D}_{\Lbar} \chi_{(\text{small})} + r\slashed{\D}_L \chibar_{(\text{small})} 
		\\
		r\slashed{\nabla}^2 \log \mu
	\end{pmatrix} (\slashed{\nabla}\phi)
	+ \bm{\Gamma}^{(1)}_{(-1-\delta)} (\slashed{\D} \mathscr{Z} \phi)
	+ \bm{\Gamma}^{(1)}_{(-1+C_{(1)}\epsilon)} (\overline{\slashed{\D}} \mathscr{Z} \phi)
	\\
	&\phantom{=}
	+ \bm{\Gamma}^{(1)}_{(-1+C_{(1)}\epsilon)} \slashed{\D}_L(r\slashed{\D}_L \phi)
	+ \bm{\Gamma}^{(1)}_{(-1+C_{(1)}\epsilon)} \slashed{\D}\phi
	+ \bm{\Gamma}^{(1)}_{(C_{(1)}\epsilon)} \slashed{\Omega}_{L} \phi
	+ \bm{\Gamma}^{(1)}_{(-\delta)} \slashed{\Omega}_{\Lbar} \phi
	\end{split}
	\end{equation*}
	The only terms which differ in their schematic form from the ones already present if $\phi$ is a scalar field are the final two terms, involving the field $\phi$ itself rather than its derivatives. Again using proposition \ref{proposition expression for Omega} these are given by
	\begin{equation*}
	\bm{\Gamma}^{(1}_{(C_{(1)}\epsilon)} \slashed{\Omega}_L + 	\bm{\Gamma}^{(1)}_{-\delta} \slashed{\Omega}_L
	=
	\bm{\Gamma}^{(1)}_{(-1 + C_{(1)}\epsilon)} (\slashed{\D} \mathscr{Z} h)_{(\text{frame})} + \bm{\Gamma}^{(1)}_{(-2 - \delta + C_{(1)}\epsilon)}
	\end{equation*}
	
	There are also new error terms arising from the ``additional error terms'' in the case that $\phi$ is a higher order tensor field. These arise from the terms
	\begin{equation*}
	\slashed{\D}^\mu\left(r\slashed{\Pi}_{\slashed{\alpha}}^{\phantom{\slashed{\alpha}}\nu} [\slashed{\D}_\mu \, , \slashed{\D}_\nu]\phi \right)
	+ r\omega \Lbar^\mu \Pi_{\slashed{\alpha}}^{\phantom{\slashed{\alpha}}\nu} [\slashed{\D}_\mu \, , \slashed{\D}_\nu] \phi
	\end{equation*}
	Proposition \ref{proposition additional error term commuting with rnabla} gives the bound for these terms, proving the proposition.
		
\end{proof}

\begin{proposition}[The structure of the inhomogeneous terms after commuting once with $r\slashed{\D}_L$]
	\label{proposition inhomogeneous terms commute with rL}
	Suppose that $\phi$ is an $S_{\tau,r}$-tangent tensor field satisfying the equation
	\begin{equation*}
	\tilde{\Box}_g \phi = F
	\end{equation*}
	
	Then, if $\phi$ is a scalar field, $rL\phi$ satisfies the following equation, given schematically except for the first four terms:
	\begin{equation}
	\begin{split}
	\tilde{\Box}_g (rL \phi)
	&= r(LF)
	+ \slashed{\Delta}\phi
	+ r^{-1} L(rL\phi)
	+ r^{-1} (L\phi)
	+ \begin{pmatrix}
		\overline{\slashed{\D}}(\mathscr{Z}h)_{(\text{frame})} \\
		\left(L(rLh)\right)_{(\text{frame})} \\
		r(\tilde{\Box}_g h)_{LL} \\
		\bm{\Gamma}^{(0)}_{(-1)}
	\end{pmatrix} (\partial \phi)
	+ \begin{pmatrix}
		r(\tilde{\Box}_g h)_{(\text{frame})} \\
		(\slashed{\D}\mathscr{Z} h)_{(\text{frame})} \\
		r\slashed{\Delta} \log \mu \\
		\bm{\Gamma}^{(0)}_{(-1+2C_{(0)}\epsilon)}
	\end{pmatrix} (\bar{\partial} \phi)
	\\
	&\phantom{=}
	+ \bm{\Gamma}^{(0)}_{(-1+C_{(0)}\epsilon)} \overline{\slashed{\D}}(\mathscr{Z}\phi)
	+ \bm{\Gamma}^{(0)}_{(-1)} L(rL\phi)
	+ r\bm{\Gamma}^{(0)}_{(-1)} F
	\end{split}
	\end{equation}
	
	On the other hand, if $\phi$ is a higher order tensor field, then $r\slashed{\D}_L \phi$ satisfies the following equation, also given schematically except for the first four terms:
	\begin{equation*}
	\begin{split}
	\tilde{\slashed{\Box}}_g (r\slashed{\D}_L \phi)
	&= r(\slashed{\D}_L F)
	+ \slashed{\Delta}\phi
	+ r^{-1} \slashed{\D}_L(r\slashed{\D}_L \phi)
	+ r^{-1} (\slashed{\D}_L \phi)
	+ \begin{pmatrix}
		\overline{\slashed{\D}}(\mathscr{Z}h)_{(\text{frame})} \\
		\left(L(rLh)\right)_{(\text{frame})} \\
		r(\tilde{\Box}_g h)_{LL} \\
		\bm{\Gamma}^{(0)}_{(-1)}
	\end{pmatrix} (\slashed{\D} \phi)
	\\
	&\phantom{=}
	+ \begin{pmatrix}
		r(\tilde{\Box}_g h)_{(\text{frame})} \\
		(\slashed{\D}\mathscr{Z} h)_{(\text{frame})} \\
		r\slashed{\Delta} \log \mu \\
		\bm{\Gamma}^{(0)}_{(-1+2C_{(0)}\epsilon)}
	\end{pmatrix} (\overline{\slashed{\D}} \phi)
	+ \begin{pmatrix}
	r^{-1} \left(\slashed{\D}_L (r\slashed{\D}_L \mathscr{Z}h)\right)_{(\text{frame})} \\
	\bm{\Gamma}^{(0)}_{(-1+C_{(0)}\epsilon)} \left(L(rLh)\right)_{(\text{frame})} \\
	(\tilde{\slashed{\Box}}_g \mathscr{Z} h)_{(\text{frame})} \\
	\bm{\Gamma}^{(0)}_{(C_{(0)}\epsilon)} (\tilde{\Box}_g h)_{(\text{frame})} \\
	r^{-1} (\slashed{\D}\mathscr{Z}^2 h)_{(\text{frame})} \\
	\bm{\Gamma}^{(0)}_{(-1+C_{(0)}\epsilon)} (\slashed{\D}\mathscr{Z} h)_{(\text{frame})} \\
	\bm{\Gamma}^{(0)}_{(-1+C_{(0)}\epsilon)} \left(\mathscr{Z} \chi_{(\text{small})} \right) \\
	\bm{\Gamma}^{(0)}_{(-2+C_{(0)}\epsilon)} \mathscr{Z}^2 \log \mu \\
	\bm{\Gamma}^{(0)}_{(-2+3C_{(0)}\epsilon)}
	\end{pmatrix} \cdot \phi	
	+ \bm{\Gamma}^{(0)}_{(-1+C_{(0)}\epsilon)} \overline{\slashed{\D}}(\mathscr{Z}\phi)
	\\
	&\phantom{=}
	+ \bm{\Gamma}^{(0)}_{(-1)} \slashed{\D}_L(r\slashed{\D}_L\phi)
	+ r\bm{\Gamma}^{(0)}_{(-1)} F
	\end{split}
	\end{equation*}
	
\end{proposition}

\begin{proof}
	
	Again, we use proposition \ref{proposition commute vector field}. First, suppose that $\phi$ is a scalar field. From proposition \ref{proposition commute rL} we find, schematically,
	\begin{equation*}
	\begin{split}
	&\slashed{\D}_\mu \left( {^{(rL)}\mathscr{J}[\phi]} \right)^\mu - (rL\omega)\slashed{\D}_{\Lbar}\phi  + \frac{1}{2}\left( \tr_{\slashed{g}} {^{(rL)}\pi} \right) \tilde{\slashed{\Box}}_g \phi + \omega\slashed{\D}_{[\Lbar, rL]}\phi - \frac{1}{2}\omega ( \tr {^{(rL)}\pi})\slashed{\D}_{\Lbar}\phi
	\\
	&=
	r(L F)
	+ \slashed{\Delta}\phi
	+ r^{-1} L(rL\phi)
	+ r^{-1} (L\phi)
	+ \begin{pmatrix}
		rL(\tr_{\slashed{g}}\chi_{(\text{small})}) \\
		rL\omega \\
		\bm{\Gamma}^{(0)}_{(-1)}
	\end{pmatrix} (\Lbar \phi)
	+ \begin{pmatrix}
		r\Lbar \tr_{\slashed{g}}\chi_{(\text{small})} \\
		r\slashed{\Div}\zeta \\
		r\slashed{\Delta}\log\mu \\
		rL\omega \\
		r\slashed{\D}_L \zeta \\
		r\slashed{\nabla} \omega \\
		r\slashed{\Div}\hat{\chi} \\
		\bm{\Gamma}^{(1)}_{(-1+2C_{(1)})}
	\end{pmatrix} (\bar{\partial}\phi)
	\\
	&\phantom{=}
	+ \bm{\Gamma}^{(0)}_{(-1+C_{(1)}\epsilon)} \overline{\slashed{\D}}(\mathscr{Z}\phi)
	+ \bm{\Gamma}^{(0)}_{(-1)} L(rL\phi)
	+ r\bm{\Gamma}^{(0)}_{(-1)} \tilde{\Box}_g \phi	
	\end{split}
	\end{equation*}
	Using proposition \ref{proposition transport chi chismall} we can write (again schematically)
	\begin{equation*}
	r\slashed{\D}_L \chi_{(\text{small})}
	= r \slashed{\Pi}_\mu^{\phantom{\mu}a} \slashed{\Pi}_\nu^{\phantom{\nu}b} R_{LaLb}
	+ \bm{\Gamma}^{(0)}_{(-1)}
	\end{equation*}
	and using proposition \ref{proposition explicit expression for Riemann} we find that
	\begin{equation*}
	r \slashed{\Pi}_\mu^{\phantom{\mu}a} \slashed{\Pi}_\nu^{\phantom{\nu}b} R_{LaLb}
	=
	\overline{\slashed{\D}}(\mathscr{Z} h)_{(\text{frame})}
	+ \left( L(rL h) \right)_{(\text{frame})}
	+ \bm{\Gamma}^{(0)}_{(-1-\delta+C_{(0)}\epsilon)}
	\end{equation*}
	
	Next, using propositions \ref{proposition transport mu}, \ref{proposition transport La} and \ref{proposition transport rectangular} we find that
	\begin{equation*}
	\begin{split}
	rL\omega
	&=
	\left(L(r\Lbar h)\right)_{LL}
	+ \left(L(rL h)\right)_{(\text{frame})}
	+ \bm{\Gamma}^{(0)}_{(-1)}
	\\
	&=
	r(\tilde{\Box}_g h)_{LL}
	+ \left(L(rL h)\right)_{(\text{frame})}
	+ (\overline{\slashed{\D}}\mathscr{Z}h)_{(\text{frame})}
	+ \bm{\Gamma}^{(0)}_{(-1)}
	\end{split}
	\end{equation*}
	and also, using proposition \ref{proposition angular rectangular} as well, we have
	\begin{equation*}
	r\slashed{\nabla} \omega
	=
	\left(\slashed{\D}(\mathscr{Z} h)\right)_{(\text{frame})}
	+ \bm{\Gamma}^{(0)}_{(-1)}
	\end{equation*}
	
	Additionally, using proposition \ref{proposition Lbar chi} we have
	\begin{equation*}
	\begin{split}
	r\Lbar \tr_{\slashed{g}}\chibar
	&= 
	r(\slashed{g}^{-1})^{ab} R_{La\Lbar b}
	+ r\slashed{\Div}\zeta
	+ r\slashed{\Delta} \log \mu
	+ \bm{\Gamma}^{(1)}_{(-1+2C_{(1)}\epsilon)}
	\\
	&=
	r(\tilde{\Box}_g h)_{(\text{frame})}
	+ (\slashed{\D}\mathscr{Z}h)_{(\text{frame})}
	+ r\slashed{\Delta} \log \mu
	+ \bm{\Gamma}^{(1)}_{(-1+2C_{(1)}\epsilon)}
	\end{split}
	\end{equation*}
	and using proposition \ref{proposition zeta} we have
	\begin{equation*}
	r\slashed{\Div}\zeta
	=
	(\slashed{\D} \mathscr{Z}h)_{(\text{frame})}
	+ \bm{\Gamma}^{(1)}_{(-1+2C_{(1)}\epsilon)}
	\end{equation*}
	and also
	\begin{equation*}
	r\slashed{\D}_L\zeta
	=
	(\slashed{\D} \mathscr{Z}h)_{(\text{frame})}
	+ \left(L (rLh)\right)_{(\text{frame})}
	+ r(\tilde{\Box}_g h)_{(\text{frame})}
	+ \bm{\Gamma}^{(1)}_{(-1+2C_{(1)}\epsilon)}
	\end{equation*}
	
	Finally, using proposition \ref{proposition div chihat} we have, schematically,
	\begin{equation*}
	r\slashed{\Div}\hat{\chi}
	=
	\mathscr{Z}\tr_{\slashed{g}}\chi_{(\text{small})}
	+ (\slashed{g}^{-1})^{\slashed{\beta}\slashed{\gamma}} R_{\slashed{\alpha}\slashed{\beta} L\slashed{\gamma}}
	+ \bm{\Gamma}^{(0)}_{(-1 + C_{(0)}\epsilon)}
	\end{equation*}
	and, using proposition \ref{proposition explicit expression for Riemann} we have
	\begin{equation*}
	(\slashed{g}^{-1})^{\slashed{\beta}\slashed{\gamma}} R_{\slashed{\alpha}\slashed{\beta} L\slashed{\gamma}}
	=
	(\overline{\slashed{\D}} \mathscr{Z} h)_{(\text{frame})}
	+ \bm{\Gamma}^{(0)}_{(-1-\frac{1}{2}\delta)}
	\end{equation*}
	finishing the calculations needed in the case that $\phi$ is a scalar field.

	If $\phi$ is a higher rank tensor field, then there are additional error terms: in the ``main error terms'' there is the additional error term
	\begin{equation*}
	r(\zeta^\mu + \slashed{\nabla}^\mu \log \mu) L^\nu [\slashed{\D}_\mu \, , \slashed{\D}_\nu] \phi
	= \bm{\Gamma}^{(1)}_{(C_{(1)}\epsilon)} \slashed{\Omega}_L \cdot \phi
	\end{equation*}
	and we can bound $\slashed{\Omega}_L$ using proposition \ref{proposition expression for Omega}, finding that
	\begin{equation*}
	r(\zeta^\mu + \slashed{\nabla}^\mu \log \mu) L^\nu [\slashed{\D}_\mu \, , \slashed{\D}_\nu] \phi
	=
	\begin{pmatrix}
	\bm{\Gamma}^{(1)}_{(-1+ C_{(1)}\epsilon)} (\overline{\slashed{\D}}\mathscr{Z} h)_{(\text{frame})} \\
	\bm{\Gamma}^{(1)}_{(-1-\delta + 2C_{(1)}\epsilon)} (\overline{\slashed{\D}}\mathscr{Z} h)_{(\text{frame})}
	\end{pmatrix}
	\cdot \phi
	\end{equation*}
	
	The remaining ``additional error terms'' were computed in proposition \ref{proposition additional error term commuting with rL}, and are of the form
	\begin{equation*}
	\begin{pmatrix}
		(\slashed{\D}\mathscr{Z}h)_{(\text{frame})} \\
		\bm{\Gamma}^{(1)}_{(-1+C_{(1)}\epsilon)}
	\end{pmatrix} \cdot (\overline{\slashed{\D}}\phi)
	+ \begin{pmatrix}
		r^{-1} \left(\slashed{\D}_L (r\slashed{\D}_L \mathscr{Z}h)\right)_{(\text{frame})} \\
		\bm{\Gamma}^{(1)}_{(-1+C_{(1)}\epsilon)} \left(L(rLh)\right)_{(\text{frame})} \\
		(\tilde{\slashed{\Box}}_g \mathscr{Z} h)_{(\text{frame})} \\
		\bm{\Gamma}^{(1)}_{(C_{(1)}\epsilon)} (\tilde{\Box}_g h)_{(\text{frame})} \\
		r^{-1} (\slashed{\D}\mathscr{Z}^2 h)_{(\text{frame})} \\
		\bm{\Gamma}^{(1)}_{(-1+C_{(1)}\epsilon)} (\slashed{\D}\mathscr{Z} h)_{(\text{frame})} \\
		\bm{\Gamma}^{(1)}_{(-1+C_{(1)}\epsilon)} \left(\mathscr{Z} \chi_{(\text{small})} \right) \\
		\bm{\Gamma}^{(1)}_{(-2+C_{(1)}\epsilon)} \mathscr{Z}^2 \log \mu \\
		\bm{\Gamma}^{(1)}_{(-2+3C_{(1)}\epsilon)}
	\end{pmatrix} \cdot \phi
	\end{equation*}
	
\end{proof}

\begin{proposition}[The structure of the inhomogeneous terms after commuting $n$ times with the operators $\mathscr{Z}$]
	\label{proposition inhomogeneous terms after commuting n times with Z}
	
	Suppose that $\phi$ is a scalar field satisfying the equation
	\begin{equation*}
	\tilde{\Box}_g \phi = F
	\end{equation*}
	
	Then $\mathscr{Z}^n \phi$ satisfies the following schematic equation:
	\begin{equation*}
	\begin{split}
	\tilde{\slashed{\Box}}_g \mathscr{Z}^n \phi 
	&=
	\mathscr{Z}^n F
	+ \bm{\Gamma}^{(0)}_{(-1)} (\slashed{\D} \mathscr{Z}^n \phi)
	+ \bm{\Gamma}^{(1)}_{(-1+ C_{(1)}\epsilon)} (\overline{\slashed{\D}} \mathscr{Z}^n \phi)
	+ \bm{\Gamma}^{(1)}_{(-1+ C_{(1)}\epsilon)} \left(\slashed{\D}_L ( r \slashed{\D}_L \mathscr{Z}^{n-1} \phi \right)
	\\
	& \phantom{=}
	+ (\partial \phi) (\slashed{\D} \mathscr{Z}^n h)_{LL}
	+ \begin{pmatrix} r^{-1}(\partial \phi) \\ (\bar{\partial}\phi) \\ r^{-1} \mathscr{Z} \phi \end{pmatrix} (\slashed{\D} \mathscr{Z}^n h)_{(\text{frame})}
	+ (\partial \phi) (\overline{\slashed{\D}} \mathscr{Z}^n h)_{(\text{frame})}
	\\
	& \phantom{=}
	+ \begin{pmatrix} (\partial \phi) \\ r(\bar{\partial} \phi) \\ \mathscr{Z} \phi \end{pmatrix} (\tilde{\slashed{\Box}}_g \mathscr{Z}^{n-1} h)_{(\text{frame})}
	+ (\partial \phi) \mathscr{Z}^n \tr_{\slashed{g}}\chi_{(\text{small})}
	+ (\bar{\partial} \phi) (r\slashed{\nabla}^2 \mathscr{Z}^{n-1} \log \mu)
	\\
	& \phantom{=}
	+ \begin{pmatrix} r^{-1} (\overline{\slashed{\D}} \mathscr{Z} \phi) \\ \bm{\Gamma}^{(1)}_{(-2+ C_{(1)}\epsilon)} (\mathscr{Z}\phi) \end{pmatrix} (\mathscr{Z}^n \log \mu)
	+ \sum_{\substack{ j+k \leq n+1 \\ j \leq n-1 \\ k \leq n-1}} \bm{\Gamma}^{(j)}_{(-1 + C_{(j)}\epsilon)} (\slashed{\D} \mathscr{Z}^k \phi)
	\\
	& \phantom{=}
	+ \sum_{\substack{ j+k \leq n \\ j \leq n-1 \\ k \leq n-1}} \tilde{\slashed{\Box}}_g (\mathscr{Z}^j h)_{(\text{frame})}  (\slashed{\D} \mathscr{Z}^k \phi)
	+ \sum_{\substack{ j+k \leq n \\ j \leq n-1 \\ k \leq n-1}} r\tilde{\slashed{\Box}}_g (\mathscr{Z}^j h)_{(\text{frame})}  (\overline{\slashed{\D}} \mathscr{Z}^k \phi)
	\\
	& \phantom{=}
	+ \sum_{\substack{ j+k \leq n+1 \\ j \leq n-1 \\ k \leq n-1}} \bm{\Gamma}^{(j)}_{(-1+C_{(j)}\epsilon)} \left(\slashed{\D}_L (\mathscr{Y} \mathscr{Z}^{k-1} \phi) \right)
	+ \sum_{\substack{j+k \leq n+1 \\ j \leq n-1 \\ k \leq n-1}} \bm{\Gamma}^{(j)}_{(-2 + 2C_{(j)}\epsilon)} (\mathscr{Z}^k \phi)
	\end{split}
	\end{equation*}

\end{proposition}

\begin{proof}
	The proof is by induction on $n$. First, note that the proposition clearly holds for $n = 0$. In fact, by propositions \ref{proposition inhomogeneous terms commute with T} and \ref{proposition inhomogeneous terms commute with rnabla}, it also holds for $n = 1$. It is important to remember that $\phi$ is a \emph{scalar field}, so that no terms of the form $(\slashed{\D}\mathscr{Z}^2 h)_{(\text{frame})}$ arise until the second time we commute.
	
	First, we perform a preliminary calculation. Let
	\begin{equation*}
	\begin{split}
	\tilde{F}_{(n)} &:=
	\mathscr{Z}^n F
	+ \bm{\Gamma}^{(0)}_{(-1)} (\slashed{\D} \mathscr{Z}^n \phi)
	+ \bm{\Gamma}^{(1)}_{(-1+ C_{(1)}\epsilon)} (\overline{\slashed{\D}} \mathscr{Z}^n \phi)
	+ \bm{\Gamma}^{(1)}_{(-1+ C_{(1)}\epsilon)} \left(\slashed{\D}_L ( r \slashed{\D}_L \mathscr{Z}^{n-1} \phi \right)
	\\
	& \phantom{=}
	+ (\partial \phi) (\slashed{\D} \mathscr{Z}^n h)_{LL}
	+ \begin{pmatrix} r^{-1}(\partial \phi) \\ (\bar{\partial}\phi) \\ r^{-1} \mathscr{Z} \phi \end{pmatrix} (\slashed{\D} \mathscr{Z}^n h)_{(\text{frame})}
	+ (\partial \phi) (\overline{\slashed{\D}} \mathscr{Z}^n h)_{(\text{frame})}
	\\
	& \phantom{=}
	+ \begin{pmatrix} (\partial \phi) \\ r(\bar{\partial} \phi) \\ \mathscr{Z} \phi \end{pmatrix} (\tilde{\slashed{\Box}}_g \mathscr{Z}^{n-1} h)_{(\text{frame})}
	+ (\partial \phi) \mathscr{Z}^n \tr_{\slashed{g}}\chi_{(\text{small})}
	+ (\bar{\partial} \phi) (r\slashed{\nabla}^2 \mathscr{Z}^{n-1} \log \mu)
	\\
	& \phantom{=}
	+ \begin{pmatrix} r^{-1} (\overline{\slashed{\D}} \mathscr{Z} \phi) \\ \bm{\Gamma}^{(1)}_{(-2+ C_{(1)}\epsilon)} (\mathscr{Z}\phi) \end{pmatrix} (\mathscr{Z}^n \log \mu)
	+ \sum_{\substack{ j+k \leq n+1 \\ j \leq n-1 \\ k \leq n-1}} \bm{\Gamma}^{(j)}_{(-1 + C_{(j)}\epsilon)} (\slashed{\D} \mathscr{Z}^k \phi)
	\\
	& \phantom{=}
	+ \sum_{\substack{ j+k \leq n \\ j \leq n-1 \\ k \leq n-1}} \tilde{\slashed{\Box}}_g (\mathscr{Z}^j h)_{(\text{frame})}  (\slashed{\D} \mathscr{Z}^k \phi)
	+ \sum_{\substack{ j+k \leq n \\ j \leq n-1 \\ k \leq n-1}} r\tilde{\slashed{\Box}}_g (\mathscr{Z}^j h)_{(\text{frame})}  (\overline{\slashed{\D}} \mathscr{Z}^k \phi)
	\\
	& \phantom{=}
	+ \sum_{\substack{ j+k \leq n +1 \\ j \leq n-1 \\ k \leq n-1}} \bm{\Gamma}^{(j)}_{(-1+C_{(j)}\epsilon)} \left(\slashed{\D}_L (\mathscr{Y} \mathscr{Z}^{k-1} \phi) \right)
	+ \sum_{\substack{j+k \leq n+1 \\ j \leq n-1 \\ k \leq n-1}} \bm{\Gamma}^{(j)}_{(-2 + 2C_{(j)}\epsilon)} (\mathscr{Z}^k \phi)
	\end{split}
	\end{equation*}
	Then we can compute $\mathscr{Z} \tilde{F}_{(n)}$. The only term which is a little bit difficult to estimate is
	\begin{equation*}
	\bm{\Gamma}^{(1)}_{(-1+C_{(1)}\epsilon)} \left( \mathscr{Z} \slashed{\D}_L (\mathscr{Y} \mathscr{Z}^{n-1} \phi) \right)
	\end{equation*}
	Using propositions \ref{proposition commuting DT with first order operators} and \ref{proposition commuting rnabla with first order operators} we have, schematically,
	\begin{equation*}
	\begin{split}
	\bm{\Gamma}^{(1)}_{(-1+C_{(1)}\epsilon)} \left( \mathscr{Z} \slashed{\D}_L (r\slashed{\D}_L \mathscr{Z}^{n-1} \phi) \right)
	&=
	\bm{\Gamma}^{(1)}_{(-1+C_{(1)}\epsilon)} \left(  \slashed{\D}_L (r\mathscr{Z} \slashed{\D}_L \mathscr{Z}^{n-1} \phi) \right)
	+ \bm{\Gamma}^{(1)}_{(-1+C_{(1)}\epsilon)} \slashed{\D}_{\Lbar} (\slashed{\D}_L \mathscr{Z}^{n-1} \phi)
	\\
	&\phantom{=}
	+ \bm{\Gamma}^{(1)}_{(-1+2C_{(1)}\epsilon)} \slashed{\D}_L (\slashed{\D}_L \mathscr{Z}^{n-1} \phi)
	+ \bm{\Gamma}^{(1)}_{(-1-\delta+C_{(1)}\epsilon)} r\slashed{\nabla} (\slashed{\D}_L \mathscr{Z}^{n-1} \phi)
	\\
	&\phantom{=}
	+ \bm{\Gamma}^{(1)}_{(C_{(1)}\epsilon)} \Omega_{L\Lbar} \cdot \slashed{\D}_L (\mathscr{Z}^{n-1} \phi)
	+ \bm{\Gamma}^{(1)}_{(C_{(1)}\epsilon)} r\slashed{\Omega}_{L} \cdot \slashed{\D}_L (\mathscr{Z}^{n-1} \phi)
	\\
	\\
	&=
	\bm{\Gamma}^{(1)}_{(-1+C_{(1)}\epsilon)} \left(  \slashed{\D}_L (r \slashed{\D}_L \mathscr{Z}^{n} \phi) \right)
	\\
	&\phantom{=}
	+ \bm{\Gamma}^{(1)}_{(-1+C_{(1)}\epsilon)} \left( r\bm{\Gamma}^{(0)}_{(-1)} (\slashed{\D}_{\Lbar} \slashed{\D}_L \mathscr{Z}^{n-1} \phi) \right)
	\\
	&\phantom{=}
	+ \bm{\Gamma}^{(1)}_{(-1+C_{(1)}\epsilon)} \left( \bm{\Gamma}^{(1)}_{(-1+C_{(1)}\epsilon)} (\slashed{\D}_L \slashed{\D}_L \mathscr{Z}^{n-1} \phi) \right)
	\\
	&\phantom{=}
	+ \bm{\Gamma}^{(1)}_{(-1+C_{(1)}\epsilon)} \left( \bm{\Gamma}^{(1)}_{(-1 - \delta)} (r\slashed{\nabla} \slashed{\D}_L \mathscr{Z}^{n-1} \phi) \right)
	\\
	&\phantom{=}
	+ \bm{\Gamma}^{(1)}_{(-1+C_{(1)}\epsilon)} \slashed{\D}_{\Lbar} (\slashed{\D}_L \mathscr{Z}^{n-1} \phi)
	+ \bm{\Gamma}^{(1)}_{(-2+2C_{(1)}\epsilon)} \slashed{\D}_L (r\slashed{\D}_L \mathscr{Z}^{n-1} \phi)
	\\
	&\phantom{=}
	+ \bm{\Gamma}^{(1)}_{(-1-\delta+C_{(1)}\epsilon)} (\slashed{\D}_L \mathscr{Z}^{n} \phi)
	+ \bm{\Gamma}^{(1)}_{(-1+C_{(1)}\epsilon)} \slashed{\D}_L (\mathscr{Z}^{n-1} \phi)
	\\
	&\phantom{=}
	+ \text{ (lower order terms)}
	\\
	\\
	&=
	\bm{\Gamma}^{(1)}_{(-1+C_{(1)}\epsilon)} \left( \slashed{\D}_L (r\slashed{\D}_L \mathscr{Z}^n \phi) \right)
	+ \bm{\Gamma}^{(1)}_{(-1+C_{(1)}\epsilon - \delta)} (\slashed{\D}_L \mathscr{Z}^n \phi)
	\\
	&\phantom{=}
	+ \text{ (lower order terms)}
	\end{split}
	\end{equation*}
	where the lower order terms are already controlled by $F_{(n)}$.
	
	It is now easy to see that we have	
	\begin{equation*}
	\mathscr{Z} \tilde{F}_{(n)}
	=
	\tilde{F}_{(n+1)}
	\end{equation*}
	where the equality holds ``schematically'', that is, the schematic expression for $\mathscr{Z} \tilde{F}_{(n)}$ is given by $\tilde{F}_{(n+1)}$.

	Now, we suppose that the proposition holds for some value of $n$, say $n = n_0$. Commuting one more time with $\slashed{\D}_T$, and making use of proposition \ref{proposition inhomogeneous terms commute with T}, we find
	\begin{equation}
	\begin{split}
	\tilde{\slashed{\Box}}_g \slashed{\D}_T \mathscr{Z}^{n_0} \phi
	&=
	\slashed{\D}_T \tilde{F}_{(n_0)}
	+ \begin{pmatrix}
		r^{-1} (\slashed{\D} \mathscr{Z} h)_{(\text{frame})} \\
		(\overline{\slashed{\D}} \mathscr{Z} h)_{(\text{frame})} \\
		(\tilde{\Box}_g h)_{(\text{frame})} \\
		\bm{\Gamma}^{(1)}_{(-2+2C_{(1)}\epsilon)}
	\end{pmatrix} (\slashed{\D} \mathscr{Z}^{n_0} \phi)
	+ \begin{pmatrix}
		\mathscr{Z} \tr_{\slashed{g}}\chi_{(\text{small})} \\
		(\slashed{\D}\mathscr{Z}h)_{(\text{frame})}
	\end{pmatrix} (\overline{\slashed{\D}} \mathscr{Z}^{n_0} \phi)
	\\
	&\phantom{=}
	+ \bm{\Gamma}^{(0)}_{(-1)} (\slashed{\D} \slashed{\D}_T \mathscr{Z}^{n_0} \phi)
	+ \bm{\Gamma}^{(1)}_{(-1-\delta)} (\slashed{\D} \mathscr{Z}^{n_0 + 1} \phi)
	+ \begin{pmatrix}
		r^{-1} \left( \tilde{\slashed{\Box}}_g (\mathscr{Z} h) \right)_{(\text{frame})} \\
		r^{-1} \left(\slashed{\D}(\mathscr{Z}^2 h)\right)_{(\text{frame})} \\
		\bm{\Gamma}^{(1)}_{(-1+C_{(1)}\epsilon)} \cdot (\slashed{\D} \mathscr{Z} h)_{(\text{frame})} \\
		\bm{\Gamma}^{(1)}_{(-1+C_{(1)}\epsilon)} \cdot (\tilde{\Box}_g h)_{(\text{frame})} \\
		\bm{\Gamma}^{(1)}_{(-1+C_{(1)}\epsilon)} \cdot \mathscr{Z} \chi_{(\text{small})} \\
		\bm{\Gamma}^{(1)}_{(-2+C_{(1)}\epsilon)} \mathscr{Z}^2 \log \mu \\
		\bm{\Gamma}^{(1)}_{(-3 + 3C_{(1)}\epsilon)}
	\end{pmatrix} (\mathscr{Z}^{n_0}\phi)
	\end{split}
	\end{equation}
	and so we can easily see that
	\begin{equation*}
	\tilde{\slashed{\Box}}_g \slashed{\D}_T \mathscr{Z}^{n_0} \phi = \tilde{F}_{(n_0+1)}
	\end{equation*}
	Indeed, we have already seen that $\mathscr{Z} \tilde{F}_{(n_0)} = \tilde{F}_{(n_0+1)}$. The only other terms which are leading order (in the sense that they involve $n_0 + 1$ commutation operators) are given by
	\begin{equation*}
	\bm{\Gamma}^{(0)}_{(-1)} \slashed{\D} (\mathscr{Z}^{n_0 + 1} \phi)
	+ \bm{\Gamma}^{(0)}_{(-1-\delta)} \slashed{\D} (\mathscr{Z}^{n_0 + 1} \phi)
	\end{equation*}
	which are included in the schematic expression for $\tilde{F}_{(n_0 + 1)}$.
	
	We need to perform a similar calculation after commuting with $r\slashed{\nabla}$. Using proposition \ref{proposition inhomogeneous terms commute with rnabla} we obtain, schematically,
	\begin{equation*}
	\begin{split}
	\tilde{\slashed{\Box}}_g (r\slashed{\nabla} \mathscr{Z}^{n_0}\phi) 
	&= 
	r\slashed{\nabla}\tilde{F}_{(n_0)}
	+ \begin{pmatrix}
		\mathscr{Z} \tr_{\slashed{g}}\chi_{(\text{small})} \\
		(\slashed{\D} \mathscr{Z} h)_{LL} \\
		(\overline{\slashed{\D}} \mathscr{Z} h)_{(\text{frame})} \\
		\bm{\Gamma}^{(1)}_{(-1+2C_{(1)}\epsilon)}
	\end{pmatrix} (\slashed{\D} \mathscr{Z}^{n_0}\phi)
	+ \begin{pmatrix}
		r(\tilde{\Box}_g h)_{(\text{frame})} \\
		r\slashed{\nabla}^2 \log \mu \\
		\bm{\Gamma}^{(1)}_{(-1+2C_{(1)}\epsilon)}
	\end{pmatrix} (\slashed{\nabla}\mathscr{Z}^{n_0}\phi)
	\\
	&\phantom{=}
	+ \bm{\Gamma}^{(1)}_{(-1-\delta)} (\slashed{\D} \mathscr{Z}^{n_0 + 1}\phi) 
	+ \bm{\Gamma}^{(1)}_{(-1+C_{(1)}\epsilon)} (\overline{\slashed{\D}} \mathscr{Z}^{n_0 + 1}\phi)
	+ \bm{\Gamma}^{(1)}_{(-1+C_{(1)}\epsilon)} (\slashed{\D}_L(r\slashed{\D}_L \mathscr{Z}^{n_0}\phi))
	\\
	&\phantom{=}
	+ \begin{pmatrix}
	r^{-1} \slashed{\D}(\mathscr{Z}^2 h)_{(\text{frame})} \\
		\tilde{\slashed{\Box}}_g (\mathscr{Z}h)_{(\text{frame})} \\
		\bm{\Gamma}^{(1)}_{(-1+C_{(1)}\epsilon)} (\tilde{\Box}_g h)_{(\text{frame})} \\
		\bm{\Gamma}^{(1)}_{(-1+C_{(1)}\epsilon)} \mathscr{Z}^2 \log \mu \\
		\bm{\Gamma}^{(1)}_{(-1+C_{(1)}\epsilon)} \mathscr{Z} \chi_{(\text{small})} \\
		\bm{\Gamma}^{(1)}_{(-2+3C_{(1)}\epsilon)} 
	\end{pmatrix} \cdot (\mathscr{Z}^{n_0}\phi)
	\end{split}
	\end{equation*}
	Now, the only new terms which are leading order and which are not already included in the schematic expression for $r\slashed{\nabla} \tilde{F}_{(n_0)} \sim \tilde{F}_{(n_0+1)}$ are
	\begin{equation*}
	\bm{\Gamma}^{(1)}_{(-1-\delta)} (\slashed{\D} \mathscr{Z}^{n_0 + 1}\phi)
	+ \bm{\Gamma}^{(1)}_{(-1+ C_{(1)}\epsilon)} (\overline{\slashed{\D}} \mathscr{Z}^{n_0 + 1}\phi)
	+ \bm{\Gamma}^{(1)}_{(-1+ C_{(1)}\epsilon)} \left(\slashed{\D}_L (r\slashed{\D}_L\mathscr{Z}^{n_0}\phi) \right)
	\end{equation*}
	and these terms also appear in the schematic expression for $\tilde{F}_{(n_0 + 1)}$.
	
	Putting these calculations together, we have shown that, if, schematically,
	\begin{equation*}
	\tilde{\slashed{\Box}}_g \mathscr{Z}^{n_0} \phi = \tilde{F}_{(n_0)}
	\end{equation*}
	then
	\begin{equation*}
	\tilde{\slashed{\Box}}_g \mathscr{Z}^{n_0+1} \phi = \tilde{F}_{(n_0+1)}
	\end{equation*}
	proving the inductive step.	
	
\end{proof}

\begin{proposition}[The structure of the inhomogeneous terms after commuting $n$ times with the operators $\mathscr{Y}$]
\label{proposition inhomogeneous terms after commuting n times with Y}
	
	Suppose that $\phi$ is a \emph{scalar} field satisfying the equation
	\begin{equation*}
	\tilde{\Box}_g \phi = F
	\end{equation*}
	
	Then, if the operator $r\slashed{\D}_L$ appears $k$ times in the expansion of the operator $\mathscr{Y}^n$ (with $k \geq 1$), we have the following schematic equation:
	\begin{equation}
	\begin{split}
	\tilde{\slashed{\Box}}_g \mathscr{Y}^{n} \phi 
	&- k\slashed{\Delta} \mathscr{Y}^{n-1} \phi
	- (2^k-1)r^{-1} \slashed{\D}_L (r\slashed{\D}_L \mathscr{Y}^{n-1} \phi)
	- (2^k-1)r^{-1} \slashed{\D}_L ( \mathscr{Y}^{n-1} \phi)
	\\
	&=
	\mathscr{Y}^{n} F
	+ r^{-1} \overline{\slashed{\D}} (\mathscr{Y}^{\leq n-1} \phi)
	+ \bm{\Gamma}^{(1)}_{(-1+C_{(1)}\epsilon)} \overline{\slashed{\D}}(\mathscr{Y}^n \phi)
	+ \bm{\Gamma}^{(0)}_{(-1)} (\slashed{\D} \mathscr{Y}^n \phi)
	\\
	\\
	&\phantom{=}
	+ \begin{pmatrix} r^{-1} (\partial \phi) \\ (\bar{\partial} \phi) \\ r^{-1} \mathscr{Y}\phi \end{pmatrix} \left( \slashed{\D} \mathscr{Y}^{n} h \right)_{(\text{frame})}
	+ (\partial \phi) (\slashed{\D} \mathscr{Y}^n h)_{LL}
	+ (\partial \phi)\left( \overline{\slashed{\D}} \mathscr{Y}^{n} h \right)_{(\text{frame})}
	\\
	&\phantom{=}
	+  \begin{pmatrix} (\partial \phi) \\ r (\bar{\partial} \phi) \\ \mathscr{Y} \phi \end{pmatrix} \left( \tilde{\slashed{\Box}}_g \mathscr{Y}^{n-1} h \right)_{(\text{frame})}
	+ (\partial \phi) \mathscr{Y}^n \tr_{\slashed{g}}\chi_{(\text{small})}
	+ (\bar{\partial} \phi) (r\slashed{\nabla}^2 \mathscr{Z}^{n-1} \log \mu)
	\\
	&\phantom{=}
	+ \begin{pmatrix} r^{-1} (\overline{\slashed{\D}} \mathscr{Y} \phi) \\ \bm{\Gamma}^{(1)}_{(-2+ C_{(1)}\epsilon)} (\mathscr{Y}\phi) \end{pmatrix} (\mathscr{Y}^n \log \mu)
	+ \sum_{\substack{ j+k \leq n+1 \\ j \leq n-1 \\ k \leq n-1}} \bm{\Gamma}^{(j)}_{(-1 + C_{(j)}\epsilon)} (\slashed{\D} \mathscr{Y}^k \phi)
	\\
	& \phantom{=}
	+ \sum_{\substack{ j+k \leq n \\ j \leq n-1 \\ k \leq n-1}} r\tilde{\slashed{\Box}}_g (\mathscr{Y}^j h)_{(\text{frame})}  (\slashed{\D} \mathscr{Y}^k \phi)
	+ \sum_{\substack{j+k \leq n+1 \\ j \leq n-1 \\ k \leq n-1}} \bm{\Gamma}^{(j)}_{(-2 + 2C_{(j)}\epsilon)} (\mathscr{Y}^k \phi)
	\end{split}
	\end{equation}
	
\end{proposition}

\begin{proof}
	As above, we shall prove this proposition by induction on $n$. Using proposition \ref{proposition inhomogeneous terms commute with rL} it is clear that the proposition is true if $n = 1$.
	
	Now, assume that the proposition is true for all $n \leq n_0$. We aim to show that it is also true for $n = n_0 + 1$. There are two cases to consider: either we are applying the operator $r\slashed{\D}_L$ to the quantity $\mathscr{Z}^{n_0}\phi$, or we are applying some operator $\mathscr{Y}$ to a quantity $\mathscr{Y}^{n_0} \phi$, where the operator $r\slashed{\D}_L$ already appears in the expansion of the operator $\mathscr{Y}^{n_0}$. We consider each of these cases in turn.
	
	First, we compute the quantity
	\begin{equation*}
	\tilde{\slashed{\Box}}_g \left( r\slashed{\D}_L \mathscr{Z}^{n_0} \phi \right)
	- \slashed{\Delta} \left( r\slashed{\D}_L \mathscr{Z}^{n_0-1} \phi \right)
	- r^{-1} \slashed{\D}_L \left((r\slashed{\D}_L \left( r\slashed{\D}_L \mathscr{Z}^{n_0-1} \phi \right) \right)
	- r^{-1} \slashed{\D}_L \left( r\slashed{\D}_L \left( \mathscr{Z}^{n_0-1} \phi \right) \right)
	\end{equation*}
	i.e.\ the relevant quantity in the case that the only time $r\slashed{\D}_L$ appears in the operator $\mathscr{Y}^{n_0 + 1}$ is as the first factor on the left hand side.
	
	We make use of the schematic equation given in proposition \ref{proposition inhomogeneous terms after commuting n times with Y} for $\tilde{\slashed{\Box}}_g \left( \mathscr{Z}^{n_0} \phi \right)$. Commuting one more time with $r\slashed{\D}_L$ and using proposition \ref{proposition inhomogeneous terms commute with rL} we find
	\begin{equation}
	\begin{split}
	&\tilde{\slashed{\Box}}_g (r\slashed{\D}_L \mathscr{Z}^{n_0} \phi)
	- \slashed{\Delta} \mathscr{Z}^{n_0} \phi
	- r^{-1} \slashed{\D}_L \left( r\slashed{\D}_L \mathscr{Z}^{n_0} \phi \right)
	- r^{-1} \slashed{\D}_L \left( \mathscr{Z}^{n_0} \phi \right)
	\\
	& =
	r\slashed{\D}_L \tilde{F}_{(n_0)}
	+ r\bm{\Gamma}^{(0)}_{(-1)} \tilde{F}_{(n_0)}
	+ \begin{pmatrix}
	\overline{\slashed{\D}}(\mathscr{Y}h)_{(\text{frame})} \\
	r(\tilde{\Box}_g h)_{LL} \\
	\bm{\Gamma}^{(0)}_{(-1)}
	\end{pmatrix} (\slashed{\D} \mathscr{Z}^{n_0} \phi)
	\\
	&\phantom{=}
	+ \begin{pmatrix}
	r(\tilde{\Box}_g h)_{(\text{frame})} \\
	(\slashed{\D}\mathscr{Z} h)_{(\text{frame})} \\
	r\slashed{\Delta} \log \mu \\
	\bm{\Gamma}^{(1)}_{(-1+2C_{(1)}\epsilon)}
	\end{pmatrix} (\overline{\slashed{\D}} \mathscr{Z}^{n_0}\phi)
	+ \begin{pmatrix}
	r^{-1} \left(\slashed{\D}_L (r\slashed{\D}_L \mathscr{Z}h)\right)_{(\text{frame})} \\
	\bm{\Gamma}^{(1)}_{(-1+C_{(1)}\epsilon)} \left(L(rLh)\right)_{(\text{frame})} \\
	(\tilde{\slashed{\Box}}_g \mathscr{Z} h)_{(\text{frame})} \\
	\bm{\Gamma}^{(1)}_{(C_{(1)}\epsilon)} (\tilde{\Box}_g h)_{(\text{frame})} \\
	r^{-1} (\slashed{\D}\mathscr{Z}^2 h)_{(\text{frame})} \\
	\bm{\Gamma}^{(1)}_{(-1+C_{(1)}\epsilon)} (\slashed{\D}\mathscr{Z} h)_{(\text{frame})} \\
	\bm{\Gamma}^{(1)}_{(-1+C_{(1)}\epsilon)} \left(\mathscr{Z} \chi_{(\text{small})} \right) \\
	\bm{\Gamma}^{(1)}_{(-2+C_{(1)}\epsilon)} \mathscr{Z}^2 \log \mu \\
	\bm{\Gamma}^{(1)}_{(-2+3C_{(1)}\epsilon)}
	\end{pmatrix} \cdot (\mathscr{Z}^{n_0}\phi)
	+ \bm{\Gamma}^{(1)}_{(-1+C_{(1)}\epsilon)} \overline{\slashed{\D}}(\mathscr{Z}^{n_0 + 1}\phi)
	\\
	&\phantom{=}
	+ \bm{\Gamma}^{(0)}_{(-1)} \slashed{\D}_L(\mathscr{Y} \mathscr{Z}^{n_0}\phi)
	\end{split}
	\end{equation}
	where the schematic expression for $\tilde{F}_{(n_0)}$ is given in proposition \ref{proposition inhomogeneous terms after commuting n times with Y}. Using the schematic notation, we therefore have
	\begin{equation}
	\begin{split}
	&\tilde{\slashed{\Box}}_g (r\slashed{\D}_L \mathscr{Z}^{n_0} \phi)
	- \slashed{\Delta} \mathscr{Z}^{n_0} \phi
	- r^{-1} \slashed{\D}_L \left( r\slashed{\D}_L \mathscr{Z}^{n_0} \phi \right)
	- r^{-1} \slashed{\D}_L \left( \mathscr{Z}^{n_0} \phi \right)
	\\
	& =
	r\slashed{\D}_L \tilde{F}_{(n_0)}
	+ r\bm{\Gamma}^{(0)}_{(-1)} \tilde{F}_{(n_0)}
	+ \begin{pmatrix}
	r(\tilde{\Box}_g h)_{LL} \\
	\bm{\Gamma}^{(1)}_{(-1 - \delta)} \\
	\bm{\Gamma}^{(0)}_{(-1)}
	\end{pmatrix} (\slashed{\D} \mathscr{Z}^{n_0} \phi)
	\\
	&\phantom{=}
	+ \begin{pmatrix}
	r(\tilde{\Box}_g h)_{(\text{frame})} \\
	\bm{\Gamma}^{(2)}_{(-1+C_{(2)}\epsilon)}
	\end{pmatrix} (\overline{\slashed{\D}} \mathscr{Z}^{n_0}\phi)
	+ \begin{pmatrix}
	\bm{\Gamma}^{(2)}_{-2 + 2C_{(2)}\epsilon} \\
	(\tilde{\slashed{\Box}}_g \mathscr{Z} h)_{(\text{frame})} \\
	\bm{\Gamma}^{(1)}_{(C_{(1)}\epsilon)} (\tilde{\Box}_g h)_{(\text{frame})} \\
	\end{pmatrix} \cdot \phi	
	+ \bm{\Gamma}^{(1)}_{(-1+C_{(1)}\epsilon)} \overline{\slashed{\D}}(\mathscr{Z}^{n_0 + 1}\phi)
	\\
	&\phantom{=}
	+ \bm{\Gamma}^{(0)}_{(-1)} \slashed{\D}_L(\mathscr{Y} \mathscr{Z}^{n_0}\phi)
	\end{split}
	\end{equation}
	Note that all of these terms are clearly of the required form, except, perhaps, for $r\slashed{\D}_L \tilde{F}_{(n_0)}$. However, it is fairly easy to see that the schematic form of $r\slashed{\D}_L \tilde{F}_{(n_0)}$ is \emph{also} of the required form: note that acting with $r\slashed{\D}_L$ raises the index $n$ in quantities like $\bm{\Gamma}^{(n)}_{(-1+C_{(n)}\epsilon)}$, while, to commute $r\slashed{\D}_L$ with the derivatives $\slashed{\D}$ we can use proposition \ref{proposition commuting rL with first order operators}.
	
	One term which deserves further consideration is
	\begin{equation*}
	(\bar{\partial} \phi) r\slashed{\D}_L \left( \slashed{\nabla}^2 \mathscr{Z}^{n_0-1} \log \mu \right)
	\end{equation*}
	We can write this term as
	\begin{equation*}
	(\bar{\partial} \phi) r\slashed{\D}_L \left( r^{-2}\mathscr{Z}^{n_0+1} \log \mu \right)
	=
	(\bar{\partial} \phi) r^{-1} \slashed{\D}_L \left( \mathscr{Z}^{n_0+1} \log \mu \right)
	- 2(\bar{\partial} \phi) r^{-2} (\mathscr{Z}^{n_0+1} \log \mu)
	\end{equation*}
	and, using proposition \ref{proposition transport Yn mu} we have, schematically,
	\begin{equation*}
	\begin{split}
	\slashed{\D}_L \left(\mathscr{Z}^{n_0+1} \log \mu \right)
	&=
	\bm{\Gamma}^{(0)}_{(-1)} (\mathscr{Z}^{n_0+1} \log \mu)
	+ \sum_{j+k\leq n_0} \bm{\Gamma}^{(j+1)}_{(-1-\delta)} (\mathscr{Z}^k \log \mu)
	\\
	&\phantom{=}
	+ \bm{\Gamma}^{(0)}_{(-1, \text{large})} \left( \mathscr{Z}^{n_0} \bar{X}_{(\text{frame})} \right)
	+ (\slashed{\D} \mathscr{Z}^{n_0} h)_{(\text{frame})}
	+ \bm{\Gamma}^{(0)}_{(-1 -\delta)} (\mathscr{Z}^{n_0} X_{(\text{frame})})
	\\
	&\phantom{=}
	+ \bm{\Gamma}^{(n_0)}_{(-1+3C_{(n_0)}\epsilon)}
	\end{split}
	\end{equation*}
	
	Note that using this computation means that we do not have to include a term of the form
	\begin{equation*}
	(\bar{\partial} \phi) \left( r^{-2}\mathscr{Y}^{n} \log \mu \right)
	\end{equation*}
	but only a term of the form
	\begin{equation*}
	(\bar{\partial} \phi) \left( r^{-2}\mathscr{Z}^{n} \log \mu \right)
	\end{equation*}
	We could do a similar computation for some of the other terms, but it will not be necessary.

	Next, we deal with the case in which the operator $r\slashed{\D}_L$ has already been applied at least once. Let us suppose that it has been applied $k$ times, where $k \geq 1$. We set
	\begin{equation*}
	\begin{split}
	\tilde{F}_{(rL, n)}
	&:=
	\mathscr{Y}^{\leq n} F
	+ r^{-1} \overline{\slashed{\D}} (\mathscr{Y}^{\leq n-1} \phi)
	+ \bm{\Gamma}^{(1)}_{(-1+C_{(1)}\epsilon)} \overline{\slashed{\D}}(\mathscr{Y}^n \phi)
	+ \bm{\Gamma}^{(0)}_{(-1)} (\slashed{\D} \mathscr{Y}^n \phi)
	\\
	\\
	&\phantom{=}
	+ \begin{pmatrix} r^{-1} (\partial \phi) \\ (\bar{\partial} \phi) \\ r^{-1} \mathscr{Y}\phi \end{pmatrix} \left( \slashed{\D} \mathscr{Y}^{n} h \right)_{(\text{frame})}
	+ (\partial \phi) (\slashed{\D} \mathscr{Y}^n h)_{LL}
	+ (\partial \phi)\left( \overline{\slashed{\D}} \mathscr{Y}^{n} h \right)_{(\text{frame})}
	\\
	&\phantom{=}
	+  \begin{pmatrix} (\partial \phi) \\ r (\bar{\partial} \phi) \\ \mathscr{Y} \phi \end{pmatrix} \left( \tilde{\slashed{\Box}}_g \mathscr{Y}^{n-1} h \right)_{(\text{frame})}
	+ (\partial \phi) \mathscr{Y}^n \tr_{\slashed{g}}\chi_{(\text{small})}
	+ (\bar{\partial} \phi) (r\slashed{\nabla}^2 \mathscr{Z}^{n-1} \log \mu)
	\\
	&\phantom{=}
	+ \begin{pmatrix} r^{-1} (\overline{\slashed{\D}} \mathscr{Y} \phi) \\ \bm{\Gamma}^{(1)}_{(-2+ C_{(1)}\epsilon)} (\mathscr{Y}\phi) \end{pmatrix} (\mathscr{Y}^n \log \mu)
	+ \sum_{\substack{ j+k \leq n+1 \\ j \leq n-1 \\ k \leq n-1}} \bm{\Gamma}^{(j)}_{(-1 + C_{(j)}\epsilon)} (\slashed{\D} \mathscr{Y}^k \phi)
	\\
	& \phantom{=}
	+ \sum_{\substack{ j+k \leq n \\ j \leq n-1 \\ k \leq n-1}} r\tilde{\slashed{\Box}}_g (\mathscr{Z}^j h)_{(\text{frame})}  (\slashed{\D} \mathscr{Y}^k \phi)
	+ \sum_{\substack{j+k \leq n+1 \\ j \leq n-1 \\ k \leq n-1}} \bm{\Gamma}^{(j)}_{(-2 + 2C_{(j)}\epsilon)} (\mathscr{Y}^k \phi)
	\end{split}
	\end{equation*}
	Then, by the inductive hypothesis we have
	\begin{equation*}
	\begin{split}
	\tilde{\slashed{\Box}}_g \mathscr{Y}^{n_0} \phi 
	- k\slashed{\Delta} \mathscr{Y}^{n_0-1} \phi
	- (2^k-1)r^{-1} \slashed{\D}_L (r\slashed{\D}_L \mathscr{Y}^{n_0-1} \phi)
	- (2^k-1)r^{-1} \slashed{\D}_L ( \mathscr{Y}^{n_0-1} \phi)
	=
	\tilde{F}_{(rL, n_0)}
	\end{split}
	\end{equation*}
	since at least one factor of $r\slashed{\D}_L$ appears in the operator $\mathscr{Y}^{n_0}$. Now, we commute once more with $\mathscr{Y}$. This time, $\mathscr{Y}$ might be $\slashed{\D}_T$, $r\slashed{\nabla}$ or $r\slashed{\D}_L$, and we consider each of these in turn.
	
	We have
	\begin{equation*}
	\begin{split}
	&\tilde{\slashed{\Box}}_g \slashed{\D}_T \mathscr{Y}^{n_0} \phi 
	- k\slashed{\Delta} \slashed{\D}_T \mathscr{Y}^{n_0-1} \phi
	- r^{-1} \slashed{\D}_L (r\slashed{\D}_L \slashed{\D}_T \mathscr{Y}^{n_0-1} \phi)
	- r^{-1} \slashed{\D}_L ( \slashed{\D}_T \mathscr{Y}^{n_0-1} \phi)
	\\
	& =
	[\tilde{\slashed{\Box}}_g \, , \slashed{\D}_T] \mathscr{Y}^{n_0}\phi
	- [\slashed{\Delta} \, , \slashed{\D}_T] \mathscr{Y}^{n_0}\phi
	- (2^k-1)[r^{-1} \slashed{\D}_L \left( r \slashed{\D}_L (\cdot) \right) \, , \slashed{\D}_T] \mathscr{Y}^{n_0}\phi
	- (2^k-1)[r^{-1} \slashed{\D}_L \, , \slashed{\D}_T] \mathscr{Y}^{n_0}\phi
	\\
	&\phantom{=}
	+ \slashed{\D}_T \tilde{F}_{(rL, n_0)}
	\end{split}
	\end{equation*}
	
	Since $\slashed{\D}_T(r) = 0$, we can use proposition \ref{proposition commuting DT with first order operators} to find (schematically)
	\begin{equation*}
	[r^{-1} \slashed{\D}_L \, , \slashed{\D}_T] \mathscr{Y}^{n_0}\phi
	=
	\bm{\Gamma}^{(1)}_{(-2 + C_{(1)}\epsilon)} \slashed{\D} \mathscr{Y}^{n_0} \phi
	+ \bm{\Gamma}^{(1)}_{(-3 + C_{(1)}\epsilon)} \mathscr{Y}^{n_0} \phi
	\end{equation*}
	These terms are included in the schematic expression for $\tilde{F}_{(rL, n_0 + 1)}$.
	
	Next, we compute
	\begin{equation*}
	\begin{split}
	[r^{-1} \slashed{\D}_L \left( r\slashed{\D}_L (\cdot) \right) \, \slashed{\D}_T] \mathscr{Y}^{n_0}\phi
	&=
	\slashed{\D}_L \left( [\slashed{\D}_L , \slashed{\D}_T]\phi \right)
	+ [\slashed{\D}_L , \slashed{\D}_T] \slashed{\D}_L \phi
	+ r^{-1} [\slashed{\D}_L , \slashed{\D}_T] \phi
	\\
	&=
	\omega (\slashed{\D}_L \slashed{\D}_{\Lbar} \mathscr{Y}^{n_0}\phi + \slashed{\D}_{\Lbar} \slashed{\D}_L \mathscr{Y}^{n_0}\phi)
	+ \bm{\Gamma}^{(1)}_{(-2 + C_{(1)}\epsilon)} \slashed{\D} \mathscr{Y}^{n_0+1}\phi
	\\
	&=
	\bm{\Gamma}^{(0)}_{(-1)} \tilde{\slashed{\Box}}_g \mathscr{Y}^{n_0} \phi
	+ \bm{\Gamma}^{(1)}_{(-2 + C_{(1)}\epsilon)} \slashed{\D} \mathscr{Y}^{n_0+1}\phi
	\end{split}
	\end{equation*}
	
	Next, using the fact that $\slashed{\D}$ is a metric connection on $\mathcal{B}$ with fibre metric $\slashed{g}$, together with proposition \ref{proposition commuting DT with first order operators} we have
	\begin{equation*}
	\begin{split}
	[\slashed{\Delta} , \slashed{\D}_T] \mathscr{Y}^{n_0} \phi
	&=
	(\slashed{g}^{-1})^{\mu\nu} \left( \slashed{\nabla}_\mu [\slashed{\nabla}_\nu , \slashed{\D}_T] \mathscr{Y}^{n_0} \phi + [\slashed{\nabla}_\mu , \slashed{\D}_T] \slashed{\nabla}_\nu \mathscr{Y}^{n_0} \phi \right)
	\\
	&=
	\bm{\Gamma}^{(1)}_{(-2+2C_{(1)}\epsilon)} \slashed{\D} \mathscr{Y}^{n_0+1} \phi
	+ \bm{\Gamma}^{(2)}_{(-1+C_{(2)}\epsilon)} \slashed{\D} \mathscr{Y}^{n_0} \phi
	+ \bm{\Gamma}^{(2)}_{(-2+C_{(2)}\epsilon)} \mathscr{Y}^{n_0} \phi
	\end{split}
	\end{equation*}
	and all of these terms are in $F_{(rL, n_0 + 1)}$.

	Next, we need to make use of proposition \ref{proposition inhomogeneous terms commute with T} to calculate $[\tilde{\slashed{\Box}}_g , \slashed{\D}_T] \mathscr{Y}^{n_0}\phi$. We find
	\begin{equation*}
	\begin{split}
	[\tilde{\slashed{\Box}}_g , \slashed{\D}_T] \mathscr{Y}^{n_0}\phi
	&=
	\begin{pmatrix}
		\bm{\Gamma}^{(1)}_{(-2 + 2C_{(1)}\epsilon)} \\
		(\tilde{\Box}_g h)_{(\text{frame})}
	\end{pmatrix} (\slashed{\D} \mathscr{Y}^{n_0} \phi)
	+ \bm{\Gamma}^{(0)}_{(-1)} (\slashed{\D} \mathscr{Y}^{n_0 + 1} \phi)
	+ \bm{\Gamma}^{(1)}_{(-1-\delta)} (\slashed{\D} \mathscr{Y}^{n_0 + 1} \phi)
	\\
	&\phantom{=}
	+ \begin{pmatrix}
		\bm{\Gamma}^{(2)}_{(-2 + 2C_{(2)}\epsilon)} \\
		r^{-1} \left( \tilde{\slashed{\Box}}_g (\mathscr{Z}h) \right)_{(\text{frame})} \\
		\bm{\Gamma}^{(1)}_{(-1+C_{(1)}\epsilon)} \left( \tilde{\slashed{\Box}}_g h \right)_{(\text{frame})} \\
	\end{pmatrix} \mathscr{Y}^{n_0} \phi
	\end{split}
	\end{equation*}
	Once again, all of these terms are in $\tilde{F}_{(rL, n_0 + 1)}$
	
	Finally, we need to compute $\slashed{\D}_T \tilde{F}_{(rL, n_0)}$. It is straightforward to see that most of the terms produced in this computation are in $\tilde{F}_{(rL, n_0 + 1)}$ due to the fact that, for example, $\slashed{\D}_T \bm{\Gamma}^{(n)}_{(-1 + C_{(n)})} = \bm{\Gamma}^{(n+1)}_{(-1 + C_{(n+1)})}$. The remaining terms are also easy to compute: when commuting $\slashed{\D}_T$ with $\slashed{\D}$ or $\overline{\slashed{\D}}$ we only produce error terms with better decay, and we certainly do not encounter any terms involving a higher number of derivatives. Finally, we can use proposition \ref{proposition inhomogeneous terms commute with T} again to compute, for example,
	\begin{equation*}
	\begin{split}
	[\tilde{\slashed{\Box}}_g , \slashed{\D}_T] \mathscr{Y}^{n_0-1}h_{(\text{rect})}
	&=
	\begin{pmatrix}
	\bm{\Gamma}^{(1)}_{(-2 + 2C_{(1)}\epsilon)} \\
	(\tilde{\Box}_g h)_{(\text{frame})}
	\end{pmatrix} (\slashed{\D} \mathscr{Y}^{n_0-1} h_{(\text{rect})})
	+ \bm{\Gamma}^{(0)}_{(-1)} (\slashed{\D} \mathscr{Y}^{n_0} h_{(\text{rect})})
	\\
	&\phantom{=}
	+ \bm{\Gamma}^{(1)}_{(-1-\delta)} (\slashed{\D} \mathscr{Y}^{n_0 } h_{(\text{rect})})
	+ \begin{pmatrix}
	\bm{\Gamma}^{(2)}_{(-2 + 2C_{(2)}\epsilon)} \\
	r^{-1} \left( \tilde{\slashed{\Box}}_g (\mathscr{Z}h) \right)_{(\text{frame})} \\
	\bm{\Gamma}^{(1)}_{(-1+C_{(1)}\epsilon)} \left( \tilde{\slashed{\Box}}_g h \right)_{(\text{frame})} \\
	\end{pmatrix} \mathscr{Y}^{n_0-1} h_{(\text{rect})}
	\end{split}
	\end{equation*}
	These terms appear multiplying the quantities $(\partial \phi)$, $r\bar{\partial}\phi$ or $\mathscr{Y} \phi$, and can all be found in $\tilde{F}_{(n_0 + 1)}$.

	We need to repeat the steps above for the cases $\mathscr{Y} = r\slashed{\nabla}$ and $\mathscr{Y} = r\slashed{\D}_L$. First consider $\mathscr{Y} = r\slashed{\nabla}$. Then we have
	\begin{equation*}
	\begin{split}
	&\tilde{\slashed{\Box}}_g \left( r\slashed{\nabla} \mathscr{Y}^{n_0} \phi \right)
	- k\slashed{\Delta} \left( r\slashed{\nabla} \mathscr{Y}^{n_0-1} \phi \right)
	- (2^k-1)r^{-1} \slashed{\D}_L \left(r\slashed{\D}_L \left( r\slashed{\nabla} \mathscr{Y}^{n_0-1} \phi \right)\right)
	- (2^k-1)r^{-1} \slashed{\D}_L \left( r\slashed{\nabla} \mathscr{Y}^{n_0-1} \phi \right)
	\\
	& =
	[\tilde{\slashed{\Box}}_g \, , r\slashed{\nabla}] \mathscr{Y}^{n_0}\phi
	- k[\slashed{\Delta} \, , r\slashed{\nabla}] \mathscr{Y}^{n_0}\phi
	- (2^k-1)[r^{-1} \slashed{\D}_L \left( r \slashed{\D}_L (\cdot) \right) \, , r\slashed{\nabla}] \mathscr{Y}^{n_0}\phi
	- (2^k-1)[r^{-1} \slashed{\D}_L \, , r\slashed{\nabla}] \mathscr{Y}^{n_0}\phi
	\\
	&\phantom{=}
	+ r\slashed{\nabla} \tilde{F}_{(rL, n_0)}
	\end{split}
	\end{equation*}
	
	Now, using proposition \ref{proposition commuting rnabla with first order operators} we have
	\begin{equation*}
	\begin{split}
	[r^{-1} \slashed{\D}_L , r\slashed{\nabla}] \mathscr{Y}^{n_0} \phi
	&=
	-(\chi_{(\text{small})}) \cdot \slashed{\nabla} \mathscr{Y}^{n_0} \phi
	+ \slashed{\Omega}_L \cdot \phi
	\\
	&=
	\bm{\Gamma}^{(0)}_{(-1-\delta)} \overline{\slashed{\D}} \mathscr{Y}^{n_0} \phi
	+ \bm{\Gamma}^{(1)}_{(-2 + C_{(1)}\epsilon)} \phi
	\end{split}
	\end{equation*}
	
	Next, we calculate
	\begin{equation*}
	\begin{split}
	[r^{-1} \slashed{\D}_L \left( r \slashed{\D}_L (\cdot) \right) \, , r\slashed{\nabla}] \mathscr{Y}^{n_0}\phi
	&=
	r^{-1} [\slashed{\D}_L , r\slashed{\nabla}]\mathscr{Y}^{n_0} \phi
	+ \slashed{\D}_L \left( [\slashed{\D}_L , r\slashed{\nabla}] \mathscr{Y}^{n_0} \phi \right)
	+ [\slashed{\D}_L , \slashed{\nabla}] (r\slashed{\D}_L \mathscr{Y}^{n_0} \phi)
	\\
	&=
	\chi_{(\text{small})} \cdot \slashed{\nabla}\mathscr{Y}^{n_0} \phi
	+ \slashed{\Omega}_L \cdot \mathscr{Y}^{n_0} \phi
	+ \slashed{\D}_L \left( r\chi_{(\text{small})} \cdot \slashed{\nabla}\mathscr{Y}^{n_0} \phi
		+ r\slashed{\Omega}_L \cdot \mathscr{Y}^{n_0} \phi \right)
	\\
	&\phantom{=}
	+ \left(r^{-1} + \chi_{(\text{small})}\right) \cdot \left(  
		\slashed{\nabla}(r\slashed{\D}_L\mathscr{Y}^{n_0} \phi)
		+ \slashed{\Omega}_L \cdot r\slashed{\D}_L\mathscr{Y}^{n_0} \phi
	\right)
	\\
	&=
	\begin{pmatrix}
		\bm{\Gamma}^{(1)}_{(-1-\delta)} \\
		r^{-1}	
	\end{pmatrix} \overline{\slashed{\D}} \mathscr{Y}^{n_0 + 1} \phi
	+ \bm{\Gamma}^{(2)}_{(-2+C_{(2)}\epsilon)} \mathscr{Y}^{n_0} \phi
	\end{split}
	\end{equation*}
	
	Next, we can compute, schematically,
	\begin{equation*}
	\begin{split}
	[\slashed{\Delta}, r\slashed{\nabla}] \mathscr{Y}^{n_0} \phi
	&=
	r[\slashed{\Delta}, \slashed{\nabla}] \mathscr{Y}^{n_0} \phi
	\\
	&=
	r\slashed{\nabla} \cdot ([\slashed{\nabla}, \slashed{\nabla}] \mathscr{Y}^{n_0} \phi)
	+ [\slashed{\nabla}, \slashed{\nabla}] (r\slashed{\nabla}\mathscr{Y}^{n_0} \phi)
	\\
	&=
	\slashed{\Omega} \cdot r\slashed{\nabla}\mathscr{Y}^{n_0} \phi
	+ \left( r\slashed{\nabla} \slashed{\Omega} \right) \cdot \mathscr{Y}^{n_0} \phi
	\end{split}
	\end{equation*}
	Note that, in fact, $\slashed{\Omega}$ can be expressed in terms of the Gauss curvature of the spheres $K$, together with the metric $\slashed{g}$. In any case, the expressions given in proposition \ref{proposition expression for Omega} yield (schematically)
	\begin{equation*}
	\slashed{\Omega} = r^{-2} + \bm{\Gamma}^{(0)}_{(-2-2\delta)}
	\end{equation*}
	so that we find
	\begin{equation*}
	[\slashed{\Delta}, r\slashed{\nabla}]
	=
	\begin{pmatrix}
	r^{-1} \\ \bm{\Gamma}^{(0)}_{(-1 -\delta)}
	\end{pmatrix}\slashed{\nabla}\phi
	+ \bm{\Gamma}^{1}_{(-2-\delta)} \phi
	\end{equation*}
	Again, it is easy to see that these terms appear in $\tilde{F}_{(rL, n_0 + 1)}$.
	
	Next, we can use proposition \ref{proposition inhomogeneous terms commute with rnabla} to compute the term $[\tilde{\slashed{\Box}}_g , r\slashed{\nabla}] \mathscr{Y}^{n_0} \phi$. We find
	\begin{equation*}
	\begin{split}
	[\tilde{\slashed{\Box}}_g , r\slashed{\nabla}] \mathscr{Y}^{n_0} \phi
	&=
	\begin{pmatrix}
	\mathscr{Z} \tr_{\slashed{g}}\chi_{(\text{small})} \\
	(\slashed{\D} \mathscr{Z} h)_{LL} \\
	(\overline{\slashed{\D}} \mathscr{Z} h)_{(\text{frame})} \\
	\bm{\Gamma}^{(1)}_{(-1+2C_{(1)}\epsilon)}
	\end{pmatrix} (\slashed{\D} \mathscr{Y}^{n_0}\phi)
	+ \begin{pmatrix}
	r(\tilde{\Box}_g h)_{(\text{frame})} \\
	r\slashed{\nabla}^2 \log \mu \\
	\bm{\Gamma}^{(1)}_{(-1+2C_{(1)}\epsilon)}
	\end{pmatrix} (\slashed{\nabla} \mathscr{Y}^{n_0}\phi)
	+ \bm{\Gamma}^{(1)}_{(-1-\delta)} (\slashed{\D} \mathscr{Y}^{n_0+1 }\phi) 
	\\
	&\phantom{=}
	+ \bm{\Gamma}^{(1)}_{(-1+C_{(1)}\epsilon)} (\overline{\slashed{\D}} \mathscr{Y}^{n_0 + 1}\phi)
	+ \begin{pmatrix}
	r^{-1} \slashed{\D}(\mathscr{Z}^2 h)_{(\text{frame})} \\
	\tilde{\slashed{\Box}}_g (\mathscr{Z}h)_{(\text{frame})} \\
	\bm{\Gamma}^{(1)}_{(-1+C_{(1)}\epsilon)} (\tilde{\Box}_g h)_{(\text{frame})} \\
	\bm{\Gamma}^{(1)}_{(-1+C_{(1)}\epsilon)} \mathscr{Z}^2 \log \mu \\
	\bm{\Gamma}^{(1)}_{(-1+C_{(1)}\epsilon)} \mathscr{Z} \chi_{(\text{small})} \\
	\bm{\Gamma}^{(1)}_{(-2+3C_{(1)}\epsilon)} 
	\end{pmatrix} \cdot \mathscr{Y}^{n_0} \phi
	\\
	\\
	&=
	\bm{\Gamma}^{(1)}_{(-1 + 2C_{(1)}\epsilon)} (\slashed{\D} \mathscr{Y}^{n_0} \phi)
	+ \begin{pmatrix}
		r(\tilde{\Box}_g h)_{(\text{frame})} \\
		\bm{\Gamma}^{(2)}_{(-1+C_{(2)}\epsilon)}
	\end{pmatrix} (\overline{\slashed{\D}}\mathscr{Y}^{n_0} \phi)
	+ \bm{\Gamma}^{(1)}_{(-1-\delta)} (\slashed{\D} \mathscr{Y}^{n_0 + 1} \phi)
	\\
	&\phantom{=}
	+ \bm{\Gamma}^{(1)}_{(-1+C_{(1)}\epsilon)} (\overline{\slashed{\D}} \mathscr{Y}^{n_0 + 1} \phi)
	+ \begin{pmatrix}
		\tilde{\slashed{\Box}}_g (\mathscr{Z}h)_{(\text{frame})} \\
		\bm{\Gamma}^{(1)}_{(-1+C_{(1)}\epsilon)} (\tilde{\Box}_g h)_{(\text{frame})} \\
		\bm{\Gamma}^{(2)}_{(-2 + C_{(2)}\epsilon)}
	\end{pmatrix} \mathscr{Y}^{n_0} \phi
	\end{split}
	\end{equation*}
	again, these terms can be found in $\tilde{F}_{(rL, n_0 + 1)}$
	
	Finally, we need to compute $r\slashed{\nabla} \tilde{F}_{(rL, n_0)}$. However, the details here are almost identical to the computation of $r\slashed{\nabla} \tilde{F}_{(rL, n_0)}$, so we will not repeat them here.

	We must deal with the final case, in which $\mathscr{Y} = r\slashed{\D}_L$. This is slightly different from the other cases, since $r\slashed{\D}_L$ produces ``large'' error terms when commuted through the wave operator. As before, we begin with the computation
	\begin{equation*}
	\begin{split}
	&\tilde{\slashed{\Box}}_g \left( r\slashed{\D}_L \mathscr{Y}^{n_0} \phi \right)
	- k\slashed{\Delta} \left( r\slashed{\D}_L \mathscr{Y}^{n_0-1} \phi \right)
	- r^{-1} (2^k-1)\slashed{\D}_L \left(r\slashed{\D}_L \left( r\slashed{\D}_L \mathscr{Y}^{n_0-1} \phi \right)\right)
	\\
	&
	- r^{-1} (2^k-1)\slashed{\D}_L \left( r\slashed{\D}_L \mathscr{Y}^{n_0-1} \phi \right)
	\\
	& =
	[\tilde{\slashed{\Box}}_g \, , r\slashed{\D}_L] \mathscr{Y}^{n_0}\phi
	- k[\slashed{\Delta} \, , r\slashed{\D}_L] \mathscr{Y}^{n_0}\phi
	- [r^{-1} \slashed{\D}_L \left( r \slashed{\D}_L (\cdot) \right) \, , r\slashed{\D}_L] \mathscr{Y}^{n_0}\phi
	- [r^{-1} \slashed{\D}_L \, , r\slashed{\D}_L] \mathscr{Y}^{n_0}\phi
	\\
	&\phantom{=}
	+ r\slashed{\D}_L \tilde{F}_{(rL, n_0)}
	\end{split}
	\end{equation*}
	
	We also have
	\begin{equation*}
	[r^{-1} \slashed{\D}_L , r\slashed{\D}_L] \mathscr{Y}^{n_0} \phi
	= 2r^{-1} \slashed{\D}_L \mathscr{Y}^{n_0} \phi
	\end{equation*}
	
	Next, we have
	\begin{equation*}
	\begin{split}
	[r^{-1} \slashed{\D}_L \left( r \slashed{\D}_L (\cdot) \right) \, , r\slashed{\D}_L] \mathscr{Y}^{n_0}\phi
	&=
	2r^{-1} \slashed{\D}_L \left( r\slashed{\D}_L \mathscr{Y}^{n_0}\phi \right) 
	\\
	&=
	2r^{-1} \slashed{\D}_L \left( \mathscr{Y}^{n_0 + 1}\phi \right) 
	\end{split}
	\end{equation*}
	
	Using proposition \ref{proposition commuting rL with first order operators} we also compute
	\begin{equation*}
	\begin{split}
	[\slashed{\Delta}, r\slashed{\D}_L] \mathscr{Y}^{n_0} \phi
	&=
	(\slashed{g}^{-1})^{\mu\nu} \left(
		\slashed{\nabla}_\mu [\slashed{\nabla}_\nu, r\slashed{\D}_L] \mathscr{Y}^{n_0} \phi
		+ [\slashed{\nabla}_\mu, r\slashed{\D}_L] \slashed{\nabla}_\nu \mathscr{Y}^{n_0} \phi
	\right)
	\\
	&=
	\begin{pmatrix}
		r^{-1} \\
		\chi_{(\text{small})}
	\end{pmatrix} \slashed{\nabla}\left( r\slashed{\nabla}_\nu \mathscr{Y}^{n_0} \phi \right)
	+ \begin{pmatrix}
		r\slashed{\Div} \chi_{(\text{small})} \\
		r\slashed{\Omega}_L
	\end{pmatrix} \left( \slashed{\nabla} \mathscr{Y}^{n_0} \phi \right)
	+ r(\slashed{\Div} \slashed{\Omega}_L) \mathscr{Y}^{n_0} \phi
	\\
	&= \begin{pmatrix}
		r^{-1} \\
		\bm{\Gamma}^{(0)}_{(-1-\delta)}
	\end{pmatrix} \left( \overline{\slashed{\D}} \mathscr{Y}^{n_0 + 1} \phi \right)
	+ \bm{\Gamma}^{(1)}_{(-1-\delta)} \left( \overline{\slashed{\D}} \mathscr{Y}^{n_0} \phi \right)
	+ \bm{\Gamma}^{(2)}_{(-2-\delta)} \mathscr{Y}^{n_0} \phi
	\end{split}
	\end{equation*}
	
	To compute $[\tilde{\slashed{\Box}}_g, r\slashed{\D}_L] \mathscr{Y}^{n_0} \phi$ we make use of proposition \ref{proposition inhomogeneous terms commute with rL}. We find
	\begin{equation}
	\begin{split}
	[\tilde{\slashed{\Box}}_g, r\slashed{\D}_L] \mathscr{Y}^{n_0} \phi
	&=
	\slashed{\Delta}\mathscr{Y}^{n_0} \phi
	+ r^{-1} \slashed{\D}_L(r\slashed{\D}_L \mathscr{Y}^{n_0} \phi)
	+ r^{-1} (\slashed{\D}_L \mathscr{Y}^{n_0} \phi)
	+ \begin{pmatrix}
		\overline{\slashed{\D}}(\mathscr{Z}h)_{(\text{frame})} \\
		\left(L(rLh)\right)_{(\text{frame})} \\
		r(\tilde{\Box}_g h)_{LL} \\
		\bm{\Gamma}^{(0)}_{(-1)}
	\end{pmatrix} (\slashed{\D} \mathscr{Y}^{n_0} \phi)
	\\
	&\phantom{=}
	+ \begin{pmatrix}
		r(\tilde{\Box}_g h)_{(\text{frame})} \\
		(\slashed{\D}\mathscr{Z} h)_{(\text{frame})} \\
		r\slashed{\Delta} \log \mu \\
		\bm{\Gamma}^{(0)}_{(-1+2C_{(0)}\epsilon)}
	\end{pmatrix} (\overline{\slashed{\D}} \mathscr{Y}^{n_0} \phi)
	+ \begin{pmatrix}
	r^{-1} \left(\slashed{\D}_L (r\slashed{\D}_L \mathscr{Z}h)\right)_{(\text{frame})} \\
	\bm{\Gamma}^{(0)}_{(-1+C_{(0)}\epsilon)} \left(L(rLh)\right)_{(\text{frame})} \\
	(\tilde{\slashed{\Box}}_g \mathscr{Z} h)_{(\text{frame})} \\
	\bm{\Gamma}^{(0)}_{(C_{(0)}\epsilon)} (\tilde{\Box}_g h)_{(\text{frame})} \\
	r^{-1} (\slashed{\D}\mathscr{Z}^2 h)_{(\text{frame})} \\
	\bm{\Gamma}^{(0)}_{(-1+C_{(0)}\epsilon)} (\slashed{\D}\mathscr{Z} h)_{(\text{frame})} \\
	\bm{\Gamma}^{(0)}_{(-1+C_{(0)}\epsilon)} \left(\mathscr{Z} \chi_{(\text{small})} \right) \\
	\bm{\Gamma}^{(0)}_{(-2+C_{(0)}\epsilon)} \mathscr{Z}^2 \log \mu \\
	\bm{\Gamma}^{(0)}_{(-2+3C_{(0)}\epsilon)}
	\end{pmatrix} \cdot \mathscr{Y}^{n_0} \phi
	\\
	&\phantom{=}
	+ \bm{\Gamma}^{(0)}_{(-1+C_{(0)}\epsilon)} \overline{\slashed{\D}}(\mathscr{Z} \mathscr{Y}^{n_0} \phi)
	+ \bm{\Gamma}^{(0)}_{(-1)} \slashed{\D}_L(r\slashed{\D}_L\mathscr{Y}^{n_0} \phi)
	+ r\bm{\Gamma}^{(0)}_{(-1)} \tilde{F}_{(rL, n_0)}
	\\
	\\
	&=
	r^{-1} \slashed{\D}_L(r\slashed{\D}_L \mathscr{Y}^{n_0} \phi)
	+ r^{-1} (\slashed{\D}_L \mathscr{Y}^{n_0} \phi)
	+ r^{-1} (\overline{\slashed{\D}} \mathscr{Y}^{n_0} \phi)
	+ \bm{\Gamma}^{(0)}_{(-1)} \slashed{\D}_L(\mathscr{Y}^{n_0+1} \phi)
	\\
	&\phantom{=}
	+ \bm{\Gamma}^{(0)}_{(-1+C_{(0)}\epsilon)} \overline{\slashed{\D}}(\mathscr{Z} \mathscr{Y}^{n_0} \phi)
	+ r\bm{\Gamma}^{(0)}_{(-1)} \tilde{F}_{(rL, n_0)}
	+ \begin{pmatrix}
		\bm{\Gamma}^{(0)}_{(-1)} \\
		\bm{\Gamma}^{(1)}_{(-1-\delta)} \\
		r(\tilde{\Box}_g h)_{LL}
	\end{pmatrix} (\slashed{\D} \mathscr{Y}^{n_0} \phi)
	\\
	&\phantom{=}
	+ \begin{pmatrix}
		\bm{\Gamma}^{(2)}_{(-1 + C_{(2)}\epsilon)} \\
		r(\tilde{\Box}_g h)_{(\text{frame})}
	\end{pmatrix} (\overline{\slashed{\D}} \mathscr{Y}^{n_0} \phi)
	+ \begin{pmatrix} 
		\bm{\Gamma}^{(2)}_{(-2 + 2C_{(2)}\epsilon)} \\
		\bm{\Gamma}^{(0)}_{(C_{(0)}\epsilon)} (\tilde{\Box}_g h)_{(\text{frame})} \\
		(\tilde{\slashed{\Box}}_g \mathscr{Z} h)_{(\text{frame})} \\
	\end{pmatrix} (\mathscr{Y}^{n_0} \phi) 
	\end{split}
	\end{equation}
	Again, all of these terms appear in the schematic expression for $\tilde{F}_{(rL, n_0+1)}$.

	Finally, we need to compute $r\slashed{\D}_L \tilde{F}_{(rL, n_0+1)}$. Using the inductive hypothesis, we have
	\begin{equation*}
	\begin{split}
	r\slashed{\D}_L \tilde{F}_{(rL, n_0)}
	&=
	r\slashed{\D}_L \mathscr{Y}^{\leq n_0} F
	+ r\slashed{\D}_L \left( r^{-1} \overline{\slashed{\D}} (\mathscr{Y}^{\leq n_0-1} \phi) \right)
	+ r\slashed{\D}_L \left( \bm{\Gamma}^{(1)}_{(-1+C_{(1)}\epsilon)} \overline{\slashed{\D}}(\mathscr{Y}^{n_0} \phi) \right)
	\\
	&\phantom{=}
	+ r\slashed{\D}_L \left( \bm{\Gamma}^{(0)}_{(-1)} (\slashed{\D} \mathscr{Y}^{n_0} \phi) \right)
	+ \begin{pmatrix} r\slashed{\D}_L \left( r^{-1} (\partial \phi) \right) \\ r\slashed{\D}_L \left( (\bar{\partial} \phi) \right) \\ r\slashed{\D}_L \left( r^{-1} \mathscr{Y}\phi \right) \end{pmatrix} \left( \slashed{\D} \mathscr{Y}^{n_0} h \right)_{(\text{frame})}
	\\
	&\phantom{=}
	+ \begin{pmatrix} r^{-1} (\partial \phi) \\ (\bar{\partial} \phi) \\ r^{-1} \mathscr{Y}\phi \end{pmatrix} \left(r\slashed{\D}_L \left( \slashed{\D} \mathscr{Y}^{n_0} h \right) \right)_{(\text{frame})}		
	+ r\slashed{\D}_L \left( (\partial \phi) (\slashed{\D} \mathscr{Y}^{n_0} h)_{LL} \right)
	\\
	&\phantom{=}
	+ r\slashed{\D}_L \left( (\partial \phi)\left( \overline{\slashed{\D}} \mathscr{Y}^{n_0} h \right)_{(\text{frame})} \right)
	+  \begin{pmatrix} r\slashed{\D}_L(\partial \phi) \\ r\slashed{\D}_L \left(r (\bar{\partial} \phi)\right) \\ r\slashed{\D}_L \mathscr{Y} \phi \end{pmatrix} \left( \tilde{\slashed{\Box}}_g \mathscr{Y}^{n_0-1} h \right)_{(\text{frame})}
	\\
	&\phantom{=}
	+  \begin{pmatrix} (\partial \phi) \\ \left(r (\bar{\partial} \phi)\right) \\ \mathscr{Y} \phi \end{pmatrix} \left( r\slashed{\D}_L \tilde{\slashed{\Box}}_g \mathscr{Y}^{n_0-1} h \right)_{(\text{frame})}	
	+ r\slashed{\D}_L \left( (\partial \phi) \mathscr{Y}^{n_0} \tr_{\slashed{g}}\chi_{(\text{small})} \right)
	\\
	&\phantom{=}
	+ r\slashed{\D}_L \left( (\bar{\partial} \phi) (r\slashed{\nabla}^2 \mathscr{Y}^{n_0-1} \log \mu) \right)
	+ \begin{pmatrix} r\slashed{\D}_L \left( r^{-1} (\overline{\slashed{\D}} \mathscr{Y} \phi) \right) \\ r\slashed{\D}_L \left( \bm{\Gamma}^{(1)}_{(-2+ C_{(1)}\epsilon)} (\mathscr{Y}\phi) \right) \end{pmatrix} (\mathscr{Y}^{n_0} \log \mu)
	\\
	&\phantom{=}
	+ \begin{pmatrix} r^{-1} (\overline{\slashed{\D}} \mathscr{Y} \phi) \\ \bm{\Gamma}^{(1)}_{(-2+ C_{(1)}\epsilon)} (\mathscr{Y}\phi) \end{pmatrix} (r\slashed{\D}_L\mathscr{Y}^{n_0} \log \mu)
	+ \sum_{\substack{ j+k \leq n_0+1 \\ j \leq n_0-1 \\ k \leq n_0 -1}} r\slashed{\D}_L \left( \bm{\Gamma}^{(j)}_{(-1 + C_{(j)}\epsilon)} (\slashed{\D} \mathscr{Y}^k \phi) \right)
	\\
	& \phantom{=}
	+ \sum_{\substack{ j+k \leq n_0 \\ j \leq n_0-1 \\ k \leq n-1}} r\slashed{\D}_L \left( r\tilde{\slashed{\Box}}_g (\mathscr{Z}^j h)_{(\text{frame})}  (\slashed{\D} \mathscr{Y}^k \phi) \right)
	+ \sum_{\substack{j+k \leq n_0+1 \\ j \leq n_0-1 \\ k \leq n-1}} r\slashed{\D}_L \left( \bm{\Gamma}^{(j)}_{(-2 + 2C_{(j)}\epsilon)} (\mathscr{Y}^k \phi) \right)
	\end{split}
	\end{equation*}
	
	These terms can all be bounded by terms found in $\tilde{F}_{(rL, n_0 + 1)}$. For example, we have
	\begin{equation*}
		\begin{split}
		r\slashed{\D}_L \left( r^{-1} \overline{\slashed{\D}} (\mathscr{Y}^{\leq n_0-1} \phi) \right)
		&=
		r\slashed{\D}_L \left( r^{-2} \cdot r \overline{\slashed{\D}}(\mathscr{Y}^{\leq n_0-1} \phi) \right)
		\\
		&=
		r^{-1}\slashed{\D}_L \left( \mathscr{Y}^{\leq n_0} \phi \right)
		- 2 r^{-1} \overline{\slashed{\D}}(\mathscr{Y}^{\leq n_0-1} \phi)
		\end{split}
	\end{equation*}
	To take another example, we find that
	\begin{equation*}
		\begin{split}
		r\slashed{\D}_L \left( (\bar{\partial} \phi) (r\slashed{\nabla}^2 \mathscr{Y}^{n_0-1} \log \mu) \right)
		&=
		\begin{pmatrix}
			\overline{\slashed{\D}} \mathscr{Y} \phi \\
			\overline{\slashed{\D}} \phi
		\end{pmatrix} \cdot \bm{\Gamma}^{(n_0)}_{(-1+C_{(n_0)}\epsilon)}
		+ (\bar{\partial} \phi) \cdot \begin{pmatrix}
			r\slashed{\nabla}^2 \mathscr{Y}^{n_0} \log \mu \\
			[r\slashed{\D}_L , r\slashed{\nabla}^2] \mathscr{Y}^{n_0-1} \log \mu
		\end{pmatrix}
		\\
		&=
		\begin{pmatrix}
			\overline{\slashed{\D}} \mathscr{Y} \phi \\
			\overline{\slashed{\D}} \phi
		\end{pmatrix} \cdot \bm{\Gamma}^{(n_0)}_{(-1+C_{(n_0)}\epsilon)}
		+ (\bar{\partial} \phi) \cdot \begin{pmatrix}
			r\slashed{\nabla}^2 \mathscr{Y}^{n_0} \log \mu \\
			r[r\slashed{\D}_L , \slashed{\nabla}] \slashed{\nabla} \mathscr{Y}^{n_0-1} \log \mu \\
			r\slashed{\nabla} [r\slashed{\D}_L , \slashed{\nabla}] \mathscr{Y}^{n_0-1} \log \mu \\
		\end{pmatrix}
		\end{split}
	\end{equation*}
	now, using proposition \ref{proposition commuting rL with first order operators} we have
	\begin{equation*}
	\begin{split}
		r[r\slashed{\D}_L , \slashed{\nabla}] \slashed{\nabla} \mathscr{Y}^{n_0-1} \log \mu
		&=
		r^2 \chi \cdot \slashed{\nabla}^2 \mathscr{Y}^{n_0-1} \log \mu
		+ r^2 \slashed{\Omega}_L \cdot \slashed{\nabla} \mathscr{Y}^{n_0-1} \log \mu
		\\
		&=
		\left( r^{-1} + \bm{\Gamma}^{(0)}_{(-1-\delta)} \right) \mathscr{Y}^{n_0+1} \log \mu
		+ \bm{\Gamma}^{(1)}_{(-1-\delta)} \bm{\Gamma}^{(n_0)}_{(C_{(n_0)}\epsilon)}
	\end{split}
	\end{equation*}
	and similarly
	\begin{equation*}
	\begin{split}
		r\slashed{\nabla}[r\slashed{\D}_L , \slashed{\nabla}] \mathscr{Y}^{n_0-1} \log \mu
		&=
		r^2 \slashed{\nabla} \left( \chi \cdot \slashed{\nabla} \mathscr{Y}^{n_0-1} \log \mu \right)
		+ r^2\slashed{\nabla} \left( \slashed{\Omega}_L \cdot  \mathscr{Y}^{n_0-1} \log \mu \right)
		\\
		&=
		\left( r^{-1} + \bm{\Gamma}^{(0)}_{(-1-\delta)} \right) \mathscr{Y}^{n_0+1} \log \mu
		+ \bm{\Gamma}^{(1)}_{(-1-\delta)} \bm{\Gamma}^{(n_0)}_{(C_{(n_0)}\epsilon)}
		+ \bm{\Gamma}^{(2)}_{(-1-\delta)} \bm{\Gamma}^{(n_0-1)}_{(C_{(n_0-1)}\epsilon)}
	\end{split}
	\end{equation*}
	Again, all of these terms can be found in $\tilde{F}_{(rL, n_0 + 1)}$. The remaining terms can be dealt with similarly.
	
\end{proof}

In order to obtain \emph{improved} decay estimates (see appendix \ref{appendix improved energy decay}) we will also need the form of the inhomogeneous term after commuting with $\slashed{\D}_T \mathscr{Y}^{n-1}$. The important point here is that the operator $\slashed{\D}_T$ appears at least once in every quadratic term.

\begin{proposition}[The structure of the inhomogeneous term after commuting with $(\slashed{\D}_T)^j \mathscr{Y}^{n-j}$]
	\label{proposition inhomogeneous terms after commuting with T and Y}
	
	Suppose that $\phi$ is a \emph{scalar} field satisfying the equation
	\begin{equation*}
	\tilde{\Box}_g \phi = F
	\end{equation*}
	
	Then, if the operator $r\slashed{\D}_L$ appears $k$ times in the expansion of the operator $\mathscr{Y}^{n-j}$, and if $j \geq 1$ then we have the following schematic equation:
	\begin{equation}
	\begin{split}
	\tilde{\slashed{\Box}}_g \slashed{\D}_T^j \mathscr{Y}^{n-j} \phi 
	&- k\slashed{\Delta} \slashed{\D}_T^j \mathscr{Y}^{n-1} \phi
	- (2^k-1)r^{-1} \slashed{\D}_L (r\slashed{\D}_L \mathscr{Y}^{n-1} \phi)
	- (2^k-1)r^{-1} \slashed{\D}_L ( \mathscr{Y}^{n-1} \phi)
	\\
	&=
	\slashed{\D}_T^j \mathscr{Y}^{n-j} F
	+ r^{-1} \overline{\slashed{\D}} (\mathscr{Y}^{\leq n-1} \phi)
	+ \begin{pmatrix}
		(\slashed{\nabla} \log \mu) \\
		\zeta \\
		(\chi_{(\text{small})} + \chibar_{(\text{small})})
	\end{pmatrix} \overline{\slashed{\D}}(\mathscr{Y}^n \phi)
	\\
	&\phantom{=}
	+ \bm{\Gamma}^{(0)}_{(-1)} (\slashed{\D} \slashed{\D}_T^j \mathscr{Y}^{n-j} \phi)
	+ \bm{\Gamma}^{(1)}_{(-1-\delta)} (\overline{\slashed{\D}} \slashed{\D}_T^j \mathscr{Y}^{n-j} \phi)
	+ \begin{pmatrix} r^{-1} (\partial \phi) \\ (\bar{\partial} \phi) \\ r^{-1} \mathscr{Y}\phi \end{pmatrix} \left( \slashed{\D} \mathscr{Y}^{n} h \right)_{(\text{frame})}
	\\
	&\phantom{=}
	+ (\partial \phi) (\slashed{\D} \mathscr{Y}^n h)_{LL}
	+ (\partial \phi)\left( \overline{\slashed{\D}} \mathscr{Y}^{n} h \right)_{(\text{frame})}
	+  \begin{pmatrix} (\partial \phi) \\ r (\bar{\partial} \phi) \\ \mathscr{Y} \phi \end{pmatrix} \left( \tilde{\slashed{\Box}}_g \mathscr{Y}^{n-1} h \right)_{(\text{frame})}
	\\
	&\phantom{=}
	+ (\partial \phi) \mathscr{Y}^n \tr_{\slashed{g}}\chi_{(\text{small})}
	+ (\bar{\partial} \phi) (r\slashed{\nabla}^2 \slashed{\D}_T^j \mathscr{Y}^{n-j-1} \log \mu)
	+ \begin{pmatrix} r^{-1} (\overline{\slashed{\D}} \mathscr{Y} \phi) \\ \bm{\Gamma}^{(1)}_{(-2+ C_{(1)}\epsilon)} (\mathscr{Y}\phi) \end{pmatrix} (\mathscr{Y}^n \log \mu)
	\\
	&\phantom{=}
	+ \sum_{\substack{ j+k \leq n+1 \\ j \leq n-1 \\ k \leq n-1}} \bm{\Gamma}^{(j)}_{(-1 + C_{(j)}\epsilon)} (\slashed{\D} \mathscr{Y}^k \phi)
	+ \sum_{\substack{ j+k \leq n \\ j \leq n-1 \\ k \leq n-1}} r\tilde{\slashed{\Box}}_g (\mathscr{Y}^j h)_{(\text{frame})}  (\slashed{\D} \mathscr{Y}^k \phi)
	\\
	&\phantom{=}
	+ \sum_{\substack{j+k \leq n+1 \\ j \leq n-1 \\ k \leq n-1}} \bm{\Gamma}^{(j)}_{(-2 + 2C_{(j)}\epsilon)} (\mathscr{Y}^k \phi)
	\end{split}
	\end{equation}

\end{proposition}

\begin{proof}
	Note that this expression differs only from that given in proposition \ref{proposition inhomogeneous terms after commuting n times with Y} only by certain top order terms. We will prove this proposition by induction on $j$. First, we show it for $j = 1$. Then, from proposition \ref{proposition inhomogeneous terms after commuting n times with Y} we have
	
	\begin{equation*}
	\begin{split}
		\tilde{\slashed{\Box}}_g \mathscr{Y}^{n-1} \phi 
		&- k\slashed{\Delta} \mathscr{Y}^{n-2} \phi
		- (2^k-1)r^{-1} \slashed{\D}_L (r\slashed{\D}_L \mathscr{Y}^{n-2} \phi)
		- (2^k-1)r^{-1} \slashed{\D}_L ( \mathscr{Y}^{n-2} \phi)
		\\
		&=
		\mathscr{Y}^{n-1} F
		+ r^{-1} \overline{\slashed{\D}} (\mathscr{Y}^{\leq n-2} \phi)
		+ \bm{\Gamma}^{(1)}_{(-1+C_{(1)}\epsilon)} \overline{\slashed{\D}}(\mathscr{Y}^{n-1} \phi)
		+ \bm{\Gamma}^{(0)}_{(-1)} (\slashed{\D} \mathscr{Y}^{n-1} \phi)
		\\
		\\
		&\phantom{=}
		+ \begin{pmatrix} 
			r^{-1} (\partial \phi) \\
			(\bar{\partial} \phi) \\
			r^{-1} \mathscr{Y}\phi
		\end{pmatrix} \left( \slashed{\D} \mathscr{Y}^{n-1} h \right)_{(\text{frame})}
		+ (\partial \phi) (\slashed{\D} \mathscr{Y}^{n-1} h)_{LL}
		+ (\partial \phi)\left( \overline{\slashed{\D}} \mathscr{Y}^{n-1} h \right)_{(\text{frame})}
		\\
		&\phantom{=}
		+  \begin{pmatrix}
			(\partial \phi) \\
			r (\bar{\partial} \phi) \\
			\mathscr{Y} \phi
		\end{pmatrix} \left( \tilde{\slashed{\Box}}_g \mathscr{Y}^{n-2} h \right)_{(\text{frame})}
		+ (\partial \phi) \mathscr{Y}^{n-1} \tr_{\slashed{g}}\chi_{(\text{small})}
		+ (\bar{\partial} \phi) (r\slashed{\nabla}^2 \mathscr{Z}^{n-2} \log \mu)
		\\
		&\phantom{=}
		+ \begin{pmatrix} 
			r^{-1} (\overline{\slashed{\D}} \mathscr{Y} \phi) \\
			\bm{\Gamma}^{(1)}_{(-2+ C_{(1)}\epsilon)} (\mathscr{Y}\phi)
		\end{pmatrix} (\mathscr{Y}^{n-1} \log \mu)
		+ \sum_{\substack{ j+k \leq n \\ j,k \leq n-2 }} \bm{\Gamma}^{(j)}_{(-1 + C_{(j)}\epsilon)} (\slashed{\D} \mathscr{Y}^k \phi)
		\\
		& \phantom{=}
		+ \sum_{\substack{ j+k \leq n-1 \\ j, k \leq n-2 }} r\tilde{\slashed{\Box}}_g (\mathscr{Y}^j h)_{(\text{frame})}  (\slashed{\D} \mathscr{Y}^k \phi)
		+ \sum_{\substack{j+k \leq n \\ j, k \leq n-2 }} \bm{\Gamma}^{(j)}_{(-2 + 2C_{(j)}\epsilon)} (\mathscr{Y}^k \phi)
	\end{split}
	\end{equation*}
	
	Now, commuting one more time with $\slashed{\D}_T$ we have
	\begin{equation*}
	\begin{split}
		\tilde{\slashed{\Box}}_g \slashed{\D}_T \mathscr{Y}^{n-1} \phi 
		&- k\slashed{\Delta} \slashed{\D}_T\mathscr{Y}^{n-2} \phi
		- (2^k-1)r^{-1} \slashed{\D}_L (r\slashed{\D}_L \slashed{\D}_T\mathscr{Y}^{n-2} \phi)
		- (2^k-1)r^{-1} \slashed{\D}_L ( \slashed{\D}_T\mathscr{Y}^{n-2} \phi)
		\\
		&=
		[\tilde{\slashed{\Box}}_g , \slashed{\D}_T] \mathscr{Y}^{n-1} \phi 
		- k[\slashed{\Delta}, \slashed{\D}_T]\mathscr{Y}^{n-2} \phi
		- (2^k-1) [r^{-1}\slashed{\D}_L (r\slashed{\D}_L \cdot) , \slashed{\D}_T] (\mathscr{Y}^{n-2} \phi)
		\\
		&\phantom{=}
		- (2^k-1) [r^{-1} \slashed{\D}_L, \slashed{\D}_T](\mathscr{Y}^{n-2} \phi)
		+ \slashed{\D}_T \mathscr{Y}^{n-1} F
		+ r^{-1} \overline{\slashed{\D}} (\slashed{\D}_T \mathscr{Y}^{\leq n-2} \phi)
		\\
		&\phantom{=}
		+ \bm{\Gamma}^{(1)}_{(-1+C_{(1)}\epsilon)} \overline{\slashed{\D}}(\slashed{\D}_T \mathscr{Y}^{n-1} \phi)
		+ \bm{\Gamma}^{(0)}_{(-1)} (\slashed{\D} \slashed{\D}_T \mathscr{Y}^{n-1} \phi)
		+ \begin{pmatrix} 
			r^{-1} (\partial \phi) \\
			(\bar{\partial} \phi) \\
			r^{-1} \mathscr{Y}\phi
		\end{pmatrix} \left( \slashed{\D} \slashed{\D}_T \mathscr{Y}^{n-1} h \right)_{(\text{frame})}
		\\
		&\phantom{=}
		+ (\partial \phi) (\slashed{\D} \slashed{\D}_T \mathscr{Y}^{n-1} h)_{LL}
		+ (\partial \phi)\left( \overline{\slashed{\D}} \mathscr{Y}^{n-1} h \right)_{(\text{frame})}
		+  \begin{pmatrix}
			(\partial \phi) \\
			r (\bar{\partial} \phi) \\
			\mathscr{Y} \phi
		\end{pmatrix} \left( \tilde{\slashed{\Box}}_g \slashed{\D}_T \mathscr{Y}^{n-2} h \right)_{(\text{frame})}
		\\
		&\phantom{=}
		+ (\partial \phi) \slashed{\D}_T \mathscr{Y}^{n-1} \tr_{\slashed{g}}\chi_{(\text{small})}
		+ (\bar{\partial} \phi) (r\slashed{\nabla}^2 \slashed{\D}_T \mathscr{Z}^{n-2} \log \mu)
		\\
		&\phantom{=}
		+ \begin{pmatrix} 
			r^{-1} (\overline{\slashed{\D}} \mathscr{Y} \phi) \\
			\bm{\Gamma}^{(1)}_{(-2+ C_{(1)}\epsilon)} (\mathscr{Y}\phi)
		\end{pmatrix} (\slashed{\D}_T \mathscr{Y}^{n-1} \log \mu)
		+ \sum_{\substack{ j+k \leq n+1 \\ j,k \leq n-1 }} \bm{\Gamma}^{(j)}_{(-1 + C_{(j)}\epsilon)} (\slashed{\D} \mathscr{Y}^k \phi)
		\\
		& \phantom{=}
		+ \sum_{\substack{ j+k \leq n \\ j, k \leq n-1 }} r\tilde{\slashed{\Box}}_g (\mathscr{Y}^j h)_{(\text{frame})}  (\slashed{\D} \mathscr{Y}^k \phi)
		+ \sum_{\substack{j+k \leq n+1 \\ j, k \leq n-1 }} \bm{\Gamma}^{(j)}_{(-2 + 2C_{(j)}\epsilon)} (\mathscr{Y}^k \phi)
	\end{split}
	\end{equation*}
	where we have used propositions \ref{proposition commuting DT with first order operators} to handle some of the lower order terms.
	
	Now, we use proposition \ref{proposition inhomogeneous terms commute with T} to compute $[\tilde{\slashed{\Box}}_g , \slashed{\D}_T] \mathscr{Y}^{n-1} \phi$. We have also computed the quantities $[\slashed{\Delta}, \slashed{\D}_T]\mathscr{Y}^{n-2} \phi$, $ [r^{-1}\slashed{\D}_L (r\slashed{\D}_L \cdot) , \slashed{\D}_T] (\mathscr{Y}^{n-2} \phi)$ and $(2^k-1) [r^{-1} \slashed{\D}_L, \slashed{\D}_T](\mathscr{Y}^{n-2} \phi)$ in the proof of proposition \ref{proposition inhomogeneous terms after commuting n times with Y}. In fact, these terms only involve lower order terms, with the exception of $[\tilde{\slashed{\Box}}_g , \slashed{\D}_T] \mathscr{Y}^{n-1} \phi$. This leads to
	\begin{equation*}
	\begin{split}
	\tilde{\slashed{\Box}}_g \slashed{\D}_T \mathscr{Y}^{n-1} \phi 
	&- k\slashed{\Delta} \slashed{\D}_T\mathscr{Y}^{n-2} \phi
	- (2^k-1)r^{-1} \slashed{\D}_L (r\slashed{\D}_L \slashed{\D}_T\mathscr{Y}^{n-2} \phi)
	- (2^k-1)r^{-1} \slashed{\D}_L ( \slashed{\D}_T\mathscr{Y}^{n-2} \phi)
	\\
	&=
	\bm{\Gamma}^{(0)}_{(-1)} (\slashed{\D} \slashed{\D}_T \mathscr{Y}^{n-1} \phi)
	+ r^{-1} \begin{pmatrix}
		\slashed{\nabla} \log \mu \\
		\zeta \\
		(\chi_{(\text{small})} + \chibar_{(\text{small})})
	\end{pmatrix} (\slashed{\D} \mathscr{Y}^{n} \phi)
	+ \slashed{\D}_T \mathscr{Y}^{n-1} F
	\\
	&\phantom{=}
	+ r^{-1} \overline{\slashed{\D}} (\slashed{\D}_T \mathscr{Y}^{\leq n-2} \phi)
	+ \bm{\Gamma}^{(1)}_{(-1+C_{(1)}\epsilon)} \overline{\slashed{\D}}(\slashed{\D}_T \mathscr{Y}^{n-1} \phi)
	+ \bm{\Gamma}^{(0)}_{(-1)} (\slashed{\D} \slashed{\D}_T \mathscr{Y}^{n-1} \phi)
	\\
	&\phantom{=}
	+ \begin{pmatrix} 
		r^{-1} (\partial \phi) \\
		(\bar{\partial} \phi) \\
		r^{-1} \mathscr{Y}\phi
	\end{pmatrix} \left( \slashed{\D} \slashed{\D}_T \mathscr{Y}^{n-1} h \right)_{(\text{frame})}
	+ (\partial \phi) (\slashed{\D} \slashed{\D}_T \mathscr{Y}^{n-1} h)_{LL}
	+ (\partial \phi)\left( \overline{\slashed{\D}} \mathscr{Y}^{n-1} h \right)_{(\text{frame})}
	\\
	&\phantom{=}
	+  \begin{pmatrix}
		(\partial \phi) \\
		r (\bar{\partial} \phi) \\
		\mathscr{Y} \phi
	\end{pmatrix} \left( \tilde{\slashed{\Box}}_g \slashed{\D}_T \mathscr{Y}^{n-2} h \right)_{(\text{frame})}
	+ (\partial \phi) \slashed{\D}_T \mathscr{Y}^{n-1} \tr_{\slashed{g}}\chi_{(\text{small})}
	\\
	&\phantom{=}
	+ (\bar{\partial} \phi) (r\slashed{\nabla}^2 \slashed{\D}_T \mathscr{Z}^{n-2} \log \mu)
	+ \begin{pmatrix} 
		r^{-1} (\overline{\slashed{\D}} \mathscr{Y} \phi) \\
		\bm{\Gamma}^{(1)}_{(-2+ C_{(1)}\epsilon)} (\mathscr{Y}\phi)
	\end{pmatrix} (\slashed{\D}_T \mathscr{Y}^{n-1} \log \mu)
	\\
	& \phantom{=}
	+ \sum_{\substack{ j+k \leq n+1 \\ j,k \leq n-1 }} \bm{\Gamma}^{(j)}_{(-1 + C_{(j)}\epsilon)} (\slashed{\D} \mathscr{Y}^k \phi)
	+ \sum_{\substack{ j+k \leq n \\ j, k \leq n-1 }} r\tilde{\slashed{\Box}}_g (\mathscr{Y}^j h)_{(\text{frame})}  (\slashed{\D} \mathscr{Y}^k \phi)
	\\
	& \phantom{=}
	+ \sum_{\substack{j+k \leq n+1 \\ j, k \leq n-1 }} \bm{\Gamma}^{(j)}_{(-2 + 2C_{(j)}\epsilon)} (\mathscr{Y}^k \phi)
	\end{split}
	\end{equation*}
	which is of the correct form.
	
	Now, the proof for larger values of $j$ follows from an almost identical calculation, where we begin with the expression given in the proposition and, using this as the inductive hypothesis, we commute one more time with $\slashed{\D}_T$.

\end{proof}

\section{\texorpdfstring{$L^2$}{L2} bounds for geometric error terms}
\label{section L2 bounds for geometric error terms}

In the previous section we found expressions for the inhomogeneous terms after commuting with $\mathscr{Z}^n$ or $\mathscr{Y}^n$. To make use of these, we have to show that the inhomogeneities produced by commuting in this way satisfy suitable $L^2$ based estimates, so that we can apply the energy estimates in chapter \ref{chapter boundedness and energy decay} and appendix \ref{appendix improved energy decay}, and that is the purpose of this section.

We first note that bounds for quantities which are given directly in terms of the fields $\phi_{(A)}$ or their derivatives follow immediately from the $L^2$ bootstrap bounds of chapter \ref{chapter bootstrap}. Similarly, since the rectangular components of the metric perturbations can be expressed algebraically (and, effectively, linearly) in terms of the fields $\phi_{(A)}$, these terms can also be bounded directly in terms of the $L^2$ bootstrap bounds.

On the other hand, there are many ``auxiliary quantities'' which we have introduced, and which are linked by certain equations to the metric perturbations. Examples of such quantities are the foliation density $\mu$, the rectangular components of the null frame $X_{(\text{frame})}$ and the extrinsic curvatures of the sphere $\chi$ and $\chibar$.

Since many of the bounds we can establish on these auxiliary quantities involve the $L^2$ norms of other auxiliary quantities, we will establish these bounds in the context of another bootstrap argument. Specifically, under the assumption that appropriate bootstrap bounds hold, we will bound the $L^2$ norms of the auxiliary quantities in terms of suitable $L^2$ norms of the metric components $h$, in such a way that, if the (higher order) energies associated with the metric components $h$ obey suitable energy estimates, then these bootstrap bounds on the auxiliary quantities can be improved.

\subsection{\texorpdfstring{$L^2$}{L2} bootstrap assumptions for geometric quantities}

From now on we will assume that the following bootstrap assumptions hold for the $L^2$ norms of various geometric quantities:

For all $n \leq N_2-1$,
\begin{equation}
\label{equation bootstrap L2 geometric}
\begin{split}
	\int_{\mathcal{M}_{\tau}^{\tau_1} \cap \{r \geq r_0\}} \bigg(
		C_{[n]}\epsilon r^{-1-C_{[n]}\epsilon} |\bm{\Gamma}^{(n)}_{(-1+C_{(n)}\epsilon)}|^2
	\bigg) \dVol_g
	&\lesssim
	\epsilon^{2(N_2 + 1 - n)} (1+\tau)^{-1 + C_{(n)}\delta}
	\\
	\\
	\int_{\mathcal{M}_{\tau}^{\tau_1} \cap \{r \geq r_0\}} \bigg(
		c\delta r^{-1-c\delta} |\bm{\Gamma}^{(n)}_{(-1+C_{(n)}\epsilon)}|^2
	\bigg) \dVol_g
	&\lesssim
	\epsilon^{2(N_2 + 1 - n)} (1+\tau)^{-1 + C_{(n)}\delta}
	\\
	\\
	\int_{\mathcal{M}_{\tau}^{\tau_1} \cap \{r \geq r_0\}} \bigg(
		\delta r^{-1+(\frac{1}{2}-c_{[n]})\delta} |\bm{\Gamma}^{(n)}_{(-1-\delta)}|^2
	\bigg) \dVol_g
	&\lesssim
	\epsilon^{2(N_2 + 1 - n)} (1+\tau)^{-1 + C_{(n)}\delta}
\end{split}
\end{equation}
and, at the top order, we have
\begin{equation}
\label{equation bootstrap L2 geometric top order}
\begin{split}
\int_{\mathcal{M}_{\tau}^{\tau_1} \cap \{r \geq r_0\}} \bigg(
	\delta r^{-1-\delta} |\bm{\Gamma}^{(N_2)}_{(-1+C_{(N_2)})}|^2
\bigg) \dVol_g
&\lesssim
\epsilon^2 (1+\tau)^{-1 + C_{(N_2)}\delta}
\end{split}
\end{equation}

\subsection{\texorpdfstring{$L^2$}{L2} bounds for geometric quantities}

Some of the auxiliary quantities mentioned above (the ``easy'' auxiliary quantities) can be expressed directly in terms of derivatives of the rectangular components of the metric, or else in terms of other auxiliary quantities. Many other auxiliary quantities satisfy transport equations in the $L$ direction, and these can be used to provide $L^2$ bounds. Note that these bounds often involve a ``loss of derivatives''. As such, they are not suitable for bounding top order error terms, however, they can be used to provide bounds on lower order error terms. Finally, for these quantities, we must find a way to establish $L^2$ bounds at top-order \emph{without} losing derivatives.

\begin{proposition}[$L^2$ bounds on the rectangular components of the frame fields]
	\label{proposition L2 bounds rectangular}
	
	Suppose that the pointwise bootstrap bounds of chapter \ref{chapter bootstrap} hold, as well as the $L^2$ bootstrap assumptions of equation \eqref{equation bootstrap L2 geometric}.
	
	Then the commuted rectangular components of the frame fields $\mathscr{Y}^n X_{(\text{frame})}$ satisfy the following $L^2$ based bounds: for all $\tau_1 \geq \tau \geq \tau_0$, and for all $n \leq N_2$,
	
	\begin{equation*}
	\begin{split}
	&\int_{\mathcal{M}_\tau^{\tau_1} \cap \{r \geq r_0\}} C_{[n]}\epsilon r^{-3-\frac{1}{2}C_{[n]}\epsilon} |\mathscr{Y}^n X_{(\text{frame, small})}|^2\dVol_{g}
	\\
	&\lesssim
	\left( \frac{\epsilon}{\delta C_{[n]}} + \frac{1}{C_{[n-1]} C_{[n]}} \right) \epsilon^{2(N_2 + 1 - n)}(1+\tau)^{-1+C_{(n)}\delta}
	\\
	&\phantom{\lesssim}
	+ \int_{\mathcal{M}_\tau^{\tau_1} \cap \{r \geq r_0\}} \bigg(
		\frac{1}{C_{[n]}\epsilon} r^{-1+ (\frac{1}{2} - c_{[n]})\delta}|\overline{\slashed{\D}}\mathscr{Y}^n h|^2_{(\text{frame})}
		+ \frac{1}{C_{[n]}\epsilon} r^{-3-\delta}|\mathscr{Y}^n h_{(\text{rect})}|^2
	\bigg) \dVol_{g}
	\end{split}
	\end{equation*}
	also
	\begin{equation*}
	\begin{split}
	&\int_{\mathcal{M}_\tau^{\tau_1} \cap \{r \geq r_0\}} C_{[n]}\epsilon r^{-3-\frac{1}{2}C_{[n]}\epsilon} |\mathscr{Y}^n X_{(\text{frame})}|^2 \dVol_{g2}
	\\
	&\lesssim
	C_{[n]}\epsilon(\tau_1 - \tau)
	+ \left( \frac{\epsilon}{\delta C_{[n]}} + \frac{1}{C_{[n-1]} C_{[n]}} \right) \epsilon^{2(N_2 + 1 - n)}(1+\tau)^{-1+C_{(n)}\delta}
	\\
	&\phantom{\lesssim}
	+ \int_{\mathcal{M}_\tau^{\tau_1} \cap \{r \geq r_0\}} \bigg(
	\frac{1}{C_{[n]}\epsilon} r^{-1+ (\frac{1}{2} - c_{[n]})\delta}|\overline{\slashed{\D}}\mathscr{Y}^n h|^2_{(\text{frame})}
	+ \frac{1}{C_{[n]}\epsilon} r^{-3-\delta}|\mathscr{Y}^n h_{(\text{rect})}|^2
	\bigg) \dVol_{g}
	\end{split}
	\end{equation*}
	and
	\begin{equation*}
	\begin{split}
	&\int_{\mathcal{M}_\tau^{\tau_1}} c\delta r^{-3-c\delta} |\mathscr{Y}^n X_{(\text{frame, small})}|^2 \dVol_{g}
	\\
	&\lesssim
	\left( \frac{\epsilon^2}{c\delta^2} + \frac{\epsilon}{c \delta C_{[n-1]} } \right) \epsilon^{2(N_2 + 1 - n)}(1+\tau)^{-1+C_{(n)}\delta}
	\\
	&\phantom{\lesssim}
	+ \int_{\mathcal{M}_\tau^{\tau_1} \cap \{r \geq r_0\}} \bigg(
	\frac{1}{c\delta} r^{-1+ (\frac{1}{2} - c_{[n]})\delta}|\overline{\slashed{\D}}\mathscr{Y}^n h|^2_{(\text{frame})}
	+ \frac{1}{c\delta} r^{-3-\delta}|\mathscr{Y}^n h_{(\text{rect})}|^2
	\bigg) \dVol_{g}
	\end{split} 
	\end{equation*}
	and finally,
	\begin{equation*}
	\begin{split}
	&\int_{\mathcal{M}_\tau^{\tau_1} \cap \{r \geq r_0\}} \delta r^{-3 + (\frac{1}{2}-c_{[n]})\delta} |\mathscr{Y}^n \bar{X}_{(\text{frame, small})}|^2 \dVol_{g}
	\\
	&\lesssim
	\epsilon^{2(N_2 +2 - n)} (1+\tau)^{-1+C_{(n)}\delta}
	+ \int_{\mathcal{M}_\tau^{\tau_1} \cap \{r \geq r_0\}} \bigg(
	\delta r^{-1 + (\frac{1}{2}-c_{[n]})\delta}|\overline{\slashed{\D}}\mathscr{Y}^n h|^2_{(\text{frame})}
	+ \delta r^{-3-\delta} |\mathscr{Y}^n h_{(\text{rect})}|^2
	\bigg) \dVol_{g}
	\end{split}
	\end{equation*}

\end{proposition}

\begin{proof}
	We begin with the bound for $\bar{X}_{(\text{frame})}$. Proposition \ref{proposition transport Yn rectangular} gives, schematically
	\begin{equation*}
	\slashed{\D}_L \left( r\mathscr{Y}^n \bar{X}_{(\text{frame})} \right)
	=
	\bm{\Gamma}^{(0)}_{(-1)} \mathscr{Y}^n (r\bar{X}_{(\text{frame})})
	+ \bm{\Gamma}^{(n)}_{(-1+C_{(n)}\epsilon)} \bm{\Gamma}^{(0)}_{(-\frac{3}{2}\delta)}
	+ r \bm{\Gamma}^{(0)}_{(0, \text{large})} (\overline{\slashed{\D}} \mathscr{Y}^n h)_{(\text{frame})}
	+ r\bm{\Gamma}^{(n-1)}_{(-1-\delta)}
	\end{equation*}
	so, using the third part of proposition \ref{proposition Hardy} we have
	\begin{equation}
	\label{equation L2 Xbar internal}
	\begin{split}
	& \int_{\Sigma_\tau \cap \{r \geq r_0\}} r^{-3 + (\frac{1}{2}-c_{[n]})\delta} |r\mathscr{Y}^n \bar{X}_{(\text{frame})}|^2 \upd r \wedge \dVol_{\mathbb{S}^2}
	\\
	&\lesssim
	\int_{\Sigma_\tau \cap \{r \geq r_0\}} r^{-1 + (\frac{1}{2}-c_{[n]})\delta} |\slashed{\D}_L \left( r\mathscr{Y}^n \bar{X}_{(\text{frame})}\right)|^2 \upd r \wedge \dVol_{\mathbb{S}^2}
	+ \int_{S_{\tau, r_0}} |\mathscr{Y}^n \bar{X}_{(\text{frame})}|^2 \dVol_{\mathbb{S}^2}
	\\
	\\
	&\lesssim
	\int_{\Sigma_\tau \cap \{r \geq r_0\}} \bigg(
	\epsilon^2 r^{-3 + (\frac{1}{2}-c_{[n]})\delta} |r\mathscr{Y}^n \bar{X}_{(\text{frame})}|^2
	+ \epsilon^2 r^{-1-\delta} |\bm{\Gamma}^{(n)}_{(-1+C_{(n)}\epsilon)}|^2
	+ r^{-1 + (\frac{1}{2}-c_{[n]})\delta}|\overline{\slashed{\D}} \mathscr{Y}^n h|^2_{(\text{frame})}
	\\
	&\phantom{\lesssim \int_{\Sigma_\tau \cap \{r \geq r_0\}} \bigg(}
	+ r^{-1 + (\frac{1}{2}-c_{[n]})\delta}|\bm{\Gamma}^{(n-1)}_{(-1-\delta)}|
	\bigg)\upd r \wedge \dVol_{\mathbb{S}^2}
	+ \int_{S_{\tau, r_0}} |\mathscr{Y}^n \bar{X}_{(\text{frame})}|^2 \dVol_{\mathbb{S}^2}
	\end{split}
	\end{equation}
	Note that the first term on the right hand side can be absorbed on the left.
	
	To bound the final term, note that the rectangular components of the frame fields can be expressed (at $r = r_0$) in terms of the metric components $h$ - see section \ref{section geometric quantities in r leq r0}. Using this, together with the $L^2$ bootstrap we can show
	\begin{equation*}
	\begin{split}
	&\int_{\mathcal{M}_\tau^{\tau_1} \cap \{r \geq r_0\}} r^{-3 + (\frac{1}{2}-c_{[n]})\delta} |r\mathscr{Y}^n \bar{X}_{(\text{frame, small})}|^2 \dVol_{g}
	\\
	&\lesssim
	\delta^{-1} \epsilon^{2(N_2 + 2 - n)} (1+\tau)^{-1+C_{(n)}\delta}
	+ \int_{\mathcal{M}_\tau^{\tau_1} \cap \{r \geq r_0\}} \bigg(
		r^{-1 + (\frac{1}{2}-c_{[n]})\delta}|\overline{\slashed{\D}}\mathscr{Y}^n h|^2_{(\text{frame})}
	\bigg) \dVol_{g}
	\\
	&\phantom{=}
	+ \int_{\tau'=\tau}^{\tau_1} \left( \int_{S_{\tau',r_0}} \left( |\mathscr{Y}^n h_{(\text{rect})}|^2 \right)\dVol_{\mathbb{S}^2} \right) \upd \tau'
	\end{split}
	\end{equation*}
	
	Now, this last term can be estimated using proposition \ref{proposition spherical mean in terms of energy} applied to the field $r^{-\frac{1}{2}-\frac{1}{2}\delta}\mathscr{Y}^n h_{(\text{rect})}$ and with the choice $\alpha = 0$, giving
	\begin{equation*}
	\int_{S_{\tau,r_0}} |\mathscr{Y}^n h_{(\text{rect})}|^2 \dVol_{\mathbb{S}^2}
	\lesssim
	\int_{\Sigma_\tau\cap\{r \geq r_0\}} \left(
		r^{-1 - \delta} |\overline{\slashed{\D}} \mathscr{Y}^n h_{(\text{rect})}|^2
		+ r^{-3-\delta} |\mathscr{Y}^n h_{(\text{rect})}|^2
	\right) r^2 \upd r \wedge \dVol_{\mathbb{S}^2}
	\end{equation*}
	Note that this same calculation can be performed for $X_{(\text{frame})}$ as well as for $\bar{X}_{(\text{frame})}$. Note that, by the pointwise bounds on the frame fields, we can estimate this in terms of the $L^2$ norm of the ``frame'' components rather than the rectangular components if we prefer. Returning to equation \eqref{equation L2 Xbar internal}, we have
	\begin{equation*}
	\begin{split}
	&\int_{\mathcal{M}_\tau^{\tau_1} \cap \{r \geq r_0\}} \delta r^{-3 + (\frac{1}{2}-c_{[n]})\delta} |\mathscr{Y}^n \bar{X}_{(\text{frame, small})}|^2 \dVol_{g}
	\\
	&\lesssim
	\epsilon^{2(N_2 +2 - n)} (1+\tau)^{-1+C_{(n)}\delta}
	+ \int_{\mathcal{M}_\tau^{\tau_1} \cap \{r \geq r_0\}} \bigg(
		\delta r^{-1 + (\frac{1}{2}-c_{[n]})\delta}|\overline{\slashed{\D}}\mathscr{Y}^n h|^2_{(\text{frame})}
		+ \delta r^{-3-\delta} |\mathscr{Y}^n h_{(\text{rect})}|^2
	\bigg) \dVol_{g}
	\end{split}
	\end{equation*}
	
%	
%	The first term on the right hand side can be absorbed by the left hand side. Integrating over $\tau$, using the $L^2$ bootstrap bounds for the other terms, and controlling the boundary term on $S_{\tau,r_0}$ as before, we have
%	\begin{equation}
%	\label{equation L2 Xbar internal 2}
%	\begin{split}
%	& \int_{\mathcal{M}_\tau^{\tau_1} \cap \{r \geq r_0\}} r^{-3 + (\frac{1}{2}-c_{[n]}\delta)} |\mathscr{Y}^n \bar{X}_{(\text{frame})}|^2 \dVol_{g}
%	\\
%	&\lesssim
%	(1+\delta^{-2})\epsilon^{2(N_2 + 3 -n)}
%	\\
%	&\phantom{\lesssim}
%	+ \int_{\mathcal{M}_\tau^{\tau_1} \cap \{r \geq r_0\}} \bigg(
%	r^{-1 + (\frac{1}{2}-c_{[n]}\delta)}|\overline{\slashed{\D}} \mathscr{Y}^n h|^2_{(\text{frame})}
%	+ r^{-1 - \delta}|\slashed{\D} \mathscr{Y}^n h|^2_{(\text{frame})}
%	+ r^{-3 - \delta}|\mathscr{Y}^n h|^2_{(\text{frame})}
%	\bigg)\dVol_{g}
%	\end{split}
%	\end{equation}
%	

	Next, we turn to the bounds for $X_{(\text{frame, small})}$, for which we use proposition \ref{proposition transport Yn rectangular}. Using this, we find that
	\begin{equation*}
	\begin{split}
	\slashed{\D}_L \left(\mathscr{Y}^n X_{(\text{frame, small})} \right)
	&=
	\bm{\Gamma}^{(0)}_{(-1)} \left( \mathscr{Y}^k X_{(\text{frame, small})} \right)
	+ \bm{\Gamma}^{(0)}_{(0, \text{large})} (\overline{\slashed{\D}} \mathscr{Y}^n h)_{(\text{frame})}
	+ \bm{\Gamma}^{(0)}_{(-1, \text{large})} (\mathscr{Y}^n \bar{X}_{(\text{frame})})
	\\
	&\phantom{=}
	+ \bm{\Gamma}^{(n-1)}_{(-1+C_{(n-1)}\epsilon)}
	\end{split}
	\end{equation*}
	
	Now, using the third part of proposition \ref{proposition Hardy} we have
	\begin{equation*}
	\begin{split}
	&\int_{\Sigma_\tau \cap \{r \geq r_0\}} r^{-1-\frac{1}{2}C_{[n]}\epsilon} |\mathscr{Y}^n X_{(\text{frame, small})}|^2 \upd r \wedge \dVol_{\mathbb{S}^2}
	\\
	&\lesssim
	\frac{1}{(C_{[n]})^2 \epsilon^2} \int_{\Sigma_\tau \cap \{r \geq r_0\}} r^{1-\frac{1}{2}C_{[n]}\epsilon} |\slashed{\D}_L \mathscr{Y}^n X_{(\text{frame}, small)}|^2 \upd r \wedge \dVol_{\mathbb{S}^2}
	+ \frac{1}{C_{[n]} \epsilon} \int_{S_{\tau,r_0}} |\mathscr{Y}^n X_{(\text{frame})}|^2 \dVol_{\mathbb{S}^2}
	\end{split}
	\end{equation*}
	and so
	\begin{equation*}
	\begin{split}
	&\int_{\Sigma_\tau \cap \{r \geq r_0\}} r^{-1-\frac{1}{2}C_{[n]}\epsilon} |\mathscr{Y}^n X_{(\text{frame, small})}|^2 \upd r \wedge \dVol_{\mathbb{S}^2}
	\\
	&\lesssim
	\int_{\Sigma_\tau \cap \{r \geq r_0\}} \bigg(
		\frac{1}{(C_{[n]})^2} r^{-1-\frac{1}{2}C_{[n]}\epsilon} |\mathscr{Y}^n X_{(\text{frame, small})}|^2
		+ \frac{1}{(C_{[n]})^2\epsilon^2} r^{1-\frac{1}{2}C_{[n]}\epsilon}|\overline{\slashed{\D}}\mathscr{Y}^n h|^2_{(\text{frame})}
		\\
		&\phantom{\lesssim \int_{\Sigma_\tau \cap \{r \geq r_0\}} \bigg(}
		+ \frac{1}{(C_{[n]})^2\epsilon^2} r^{-1-\frac{1}{2}C_{[n]}\epsilon} |\mathscr{Y}^n \bar{X}_{(\text{frame, small})}|^2
		\\
		&\phantom{\lesssim \int_{\Sigma_\tau \cap \{r \geq r_0\}} \bigg(}
		+ \frac{1}{(C_{[n]})^2\epsilon^2} r^{1-\frac{1}{2}C_{[n]}\epsilon} |\bm{\Gamma}^{(n-1)}_{(-1+C_{(n-1)}\epsilon)}|^2		
	\bigg) \upd r \wedge \dVol_{\mathbb{S}^2}
	\\
	&\phantom{\lesssim}
	+ \frac{1}{C_{[n]} \epsilon} \int_{S_{\tau,r_0}} |\mathscr{Y}^n X_{(\text{frame, small})}|^2 \upd r \wedge \dVol_{\mathbb{S}^2}
	\end{split}
	\end{equation*}
	
	We can absorb the first term on the right hand side by the left hand side, as long as $C_{[n]}$ is sufficiently large. Integrating over $\tau$ up to $\tau_1$, using the pointwise bootstrap assumption $\Omega \sim r$ and the coarea formula, using the bound already established on $\mathscr{Y}^n \bar{X}_{(\text{frame})}$ and finally using the bootstrap assumptions to bound the other terms, we have
	\begin{equation*}
	\begin{split}
	&\int_{\mathcal{M}_\tau^{\tau_1} \cap \{r \geq r_0\}} r^{-1-\frac{1}{2}C_{[n]}\epsilon} |\mathscr{Y}^n X_{(\text{frame, small})}|^2 r^{-2} \dVol_{g}
	\\
	&\lesssim
	\frac{1}{\delta(C_{[n]})^2} \epsilon^{2(N_2 + 1 - n)}(1+\tau)^{-1+C_{(n)}\delta}
	+ \frac{1}{C_{[n-1]} (C_{[n]})^2} \epsilon^{-1} \epsilon^{2(N_2 + 1 - n)}(1+\tau)^{-1+C_{(n-1)}\delta}
	\\
	&\phantom{\lesssim}
	+ \int_{\mathcal{M}_\tau^{\tau_1} \cap \{r \geq r_0\}} \bigg(
		\frac{1}{C_{[n]}^2\epsilon^2} r^{-1+ (\frac{1}{2} - c_{[n]})\delta}|\overline{\slashed{\D}}\mathscr{Y}^n h|^2_{(\text{frame})}
		+ \frac{1}{C_{[n]}^2\epsilon^2} r^{-3-\delta}|\mathscr{Y}^n h_{(\text{rect})}|^2
	\bigg) \dVol_{g}
	\\
	&\phantom{=}
	+ \frac{1}{C_{[n]} \epsilon} \int_{\tau' = \tau}^{\tau_1} \left( \int_{S_{\tau',r_0}} |\mathscr{Y}^n X_{(\text{frame, small})}|^2 r^{-2} \dVol_{g} \right) \upd \tau'
	\end{split}
	\end{equation*}
	Dealing with the final term as before, and using the fact that $C_{(N-1)} \ll C_{(n)}$ we have
	\begin{equation*}
	\begin{split}
	&\int_{\mathcal{M}_\tau^{\tau_1} \cap \{r \geq r_0\}} C_{[n]}\epsilon r^{-1-\frac{1}{2}C_{[n]}\epsilon} |\mathscr{Y}^n X_{(\text{frame, small})}|^2 r^{-2} \dVol_{g}
	\\
	&\lesssim
	\left( \frac{\epsilon}{\delta C_{[n]}} + \frac{1}{C_{[n-1]} C_{[n]}} \right) \epsilon^{2(N_2 + 1 - n)}(1+\tau)^{-1+C_{(n)}\delta}
	\\
	&\phantom{\lesssim}
	+ \int_{\mathcal{M}_\tau^{\tau_1} \cap \{r \geq r_0\}} \bigg(
		\frac{1}{C_{[n]}\epsilon} r^{-1+ (\frac{1}{2} - c_{[n]})\delta}|\overline{\slashed{\D}}\mathscr{Y}^n h|^2_{(\text{frame})}
		+ \frac{1}{C_{[n]}\epsilon} r^{-3-\delta}|\mathscr{Y}^n h_{(\text{rect})}|^2
	\bigg) \dVol_{g}
	\end{split}
	\end{equation*}

	Next we note that the rectangular components of the frame fields differ from the ``small'' versions by either a constant $1$ or by the addition (or subtraction) of a term $\frac{x^i}{r}$. Then, since $\mathscr{Y}(1) = 0$ and $\mathscr{Y}(\frac{x^i}{r}) = X_{(\text{frame})}$, we find
	\begin{equation*}
	\begin{split}
	&\int_{\mathcal{M}_\tau^{\tau_1} \cap \{r \geq r_0\}} C_{[n]}\epsilon r^{-1-\frac{1}{2}C_{[n]}\epsilon} |\mathscr{Y}^n X_{(\text{frame})}|^2 \upd r \wedge \dVol_{\mathbb{S}^2}
	\\
	&\lesssim
	C_{[n]}\epsilon(\tau_1 - \tau)
	+ \left( \frac{\epsilon}{\delta C_{[n]}} + \frac{1}{C_{[n-1]} C_{[n]}} \right) \epsilon^{2(N_2 + 1 - n)}(1+\tau)^{-1+C_{(n)}\delta}
	\\
	&\phantom{\lesssim}
	+ \int_{\mathcal{M}_\tau^{\tau_1} \cap \{r \geq r_0\}} \bigg(
	\frac{1}{C_{[n]}\epsilon} r^{-1+ (\frac{1}{2} - c_{[n]})\delta}|\overline{\slashed{\D}}\mathscr{Y}^n h|^2_{(\text{frame})}
	+ \frac{1}{C_{[n]}\epsilon} r^{-3-\delta}|\mathscr{Y}^n h_{(\text{rect})}|^2
	\bigg) \dVol_{g}
	\end{split}
	\end{equation*}
	where the only difference is due to the fact that the frame fields $X_{(\text{frame})}$ have a ``large'' ($\mathcal{O}(1)$) component at $r = r_0$, which we have integrated over time.

	Now, by following almost identical calculations, for any $c$ we have
	\begin{equation*}
	\begin{split}
	&\int_{\mathcal{M}_\tau^{\tau_1}} c\delta r^{-1-c\delta} |\mathscr{Y}^n X_{(\text{frame, small})}|^2 \upd r \wedge \dVol_{\mathbb{S}^2}
	\\
	&\lesssim
	\left( \frac{\epsilon^2}{c\delta^2} + \frac{\epsilon}{c \delta C_{[n-1]} } \right) \epsilon^{2(N_2 + 1 - n)}(1+\tau)^{-1+C_{(n)}\delta}
	\\
	&\phantom{\lesssim}
	+ \int_{\mathcal{M}_\tau^{\tau_1} \cap \{r \geq r_0\}} \bigg(
		\frac{1}{c\delta} r^{-1+ (\frac{1}{2} - c_{[n]})\delta}|\overline{\slashed{\D}}\mathscr{Y}^n h|^2_{(\text{frame})}
		+ \frac{1}{c\delta} r^{-3-\delta}|\mathscr{Y}^n h_{(\text{rect})}|^2
	\bigg) \dVol_{g}
	\end{split}
	\end{equation*}
	Proving the proposition.

\end{proof}

\begin{remark}[$L^2$ bounds on quantities like $(\mathscr{Y}^n X_{(\text{frame})}) (\partial \phi)$]
	\label{remark L2 bounds involving Yn Xframe}
	When commuting, we will encounter error terms in which many of the commutation operators fall on the rectangular components of the frame fields. As such, we will need to estimate quantities like
	\begin{equation*}
	\int_{\mathcal{M}_\tau^{\tau_1} \cap \{r \geq r_0\}} \epsilon r^{-1-C\epsilon} |\mathscr{Y}^n X_{(\text{frame})}|^2 |\partial \phi|^2 \dVol_g
	\end{equation*}
	The obvious way to estimate this kind of quantity is to first bound $\partial \phi$ in $L^\infty$ and then bound the resulting quantity in $L^2$. However, since we expect $\mathscr{Y}^n X_{(\text{frame})} \sim 1$, this would require pointwise bootstrap assumptions on $\partial \phi$ with \emph{a priori} decay in $\tau$, which we prefer to avoid. Instead, we note that, schematically
	\begin{equation*}
	\begin{split}
	\mathscr{Y}^n X_{(\text{frame})}
	&=
	\mathscr{Y}^n X_{(\text{frame, small})}
	+ \mathscr{Y}^n \left(\frac{x^i}{r}\right)
	\\
	&=
	\mathscr{Y}^n X_{(\text{frame, small})}
	+ \begin{pmatrix} 1 \\ r^{-1} \end{pmatrix} \mathscr{Y}^{n-1} X_{(\text{frame})}
	\end{split}	
	\end{equation*}
	Iterating this argument we find (schematically)
	\begin{equation*}
	\mathscr{Y}^n X_{(\text{frame})}
	=
	\sum_{j + k \leq n} r^{-j} \mathscr{Y}^k X_{(\text{frame, small})}
	+ X_{(\text{frame})}
	\end{equation*}
	Thus we can estimate
	\begin{equation*}
	\begin{split}
	&\int_{\mathcal{M}_\tau^{\tau_1} \cap \{r \geq r_0\}} \epsilon r^{-1-C\epsilon} |\mathscr{Y}^n X_{(\text{frame})}|^2 |\partial \phi|^2 \dVol_g
	\\
	&\lesssim
	\int_{\mathcal{M}_\tau^{\tau_1} \cap \{r \geq r_0\}} \left(
	\sum_{j \leq n} \epsilon r^{-1-C\epsilon} |\mathscr{Y}^j X_{(\text{frame, small})}|^2 |\partial \phi|^2
	+ \epsilon r^{-1-C\epsilon} |\partial \phi|^2 \right)
	\dVol_g
	\end{split}
	\end{equation*}
	Now, we can bound $|\partial \phi|$ in $L^\infty$ to control the first term, and the second term can be bounded using $L^2$ based bounds for $\phi$.
	
\end{remark}

\begin{remark}
	Note that we have bounded $\int_{\mathcal{M}} r^{-3-\frac{1}{2}C_{[n]}\epsilon} |\mathscr{Y}^n X_{(\text{frame, small})}|^2 \dVol_g$. The exponent in the $r$-weight has an extra factor of $1/2$ relative to what might be expected given the bootstrap bounds of \eqref{equation bootstrap L2 geometric}. This is included because it is possible to obtain this improved weight (a fact which can be traced back to the fact that only \emph{good} derivatives appear in the transport equations for the rectangular components of the frame fiels), and this gives us more flexibility in the following estimates. Of course, these estimates also hold without the factor of $1/2$.
\end{remark}

\begin{proposition}[$L^2$ bounds on $\mathscr{Y}^n \omega$ below top-order]
	\label{proposition L2 omega low}
	Suppose that all the bootstrap bounds hold.
	
	Then for all $n \leq N_2-1$ we have
	\begin{equation*}
	\begin{split}
	&\int_{\mathcal{M}_\tau^{\tau_1} \cap \{r \geq r_0\}} C_{[n]}\epsilon r^{-1-C_{[n]}\epsilon} |\mathscr{Y}^n \omega|^2 \dVol_{g}
	\\
	&\lesssim
	C_{[n]}\epsilon^2 \epsilon^{2(N + 1 - n)}(1+\tau)^{-1 + C_{[n]}\delta}
	\\
	&\phantom{\lesssim}
	+ \int_{\mathcal{M}_\tau^{\tau_1} \cap \{r \geq r_0\}} \bigg(
	C_{[n]}\epsilon r^{-1-C_{[n]}\epsilon} |\slashed{\D} \mathscr{Y}^n h|^2_{(\text{frame})}
	+ C_{[n]}\epsilon^2 r^{-1-\delta} |\slashed{\D} \mathscr{Y}^n h_{(\text{rect})}|^2
	\\
	&\phantom{\lesssim + \int_{\mathcal{M}_\tau^{\tau_1} \cap \{r \geq r_0\}} \bigg(}
	+ C_{[n]}\epsilon^2 r^{-3-\delta} |\mathscr{Y}^n h_{(\text{rect})}|^2
	\bigg) \dVol_{g}
	\end{split}
	\end{equation*}
	and, for any $c > 0$,
	\begin{equation*}
	\begin{split}
	&\int_{\mathcal{M}_\tau^{\tau_1} \cap \{r \geq r_0\}} c\delta r^{-1-c\delta} |\mathscr{Y}^n \omega|^2 \dVol_{g}
	\\
	&\lesssim
	\epsilon^2 \left( c\delta^{-1} + 1 \right) \epsilon^{2(N + 1 - n)}(1+\tau)^{-1 + C_{[n]}\delta}
	\\
	&\phantom{\lesssim}
	+ \int_{\mathcal{M}_\tau^{\tau_1} \cap \{r \geq r_0\}} \bigg(
	r^{-1-c\delta} |\slashed{\D} \mathscr{Y}^n h|^2_{(\text{frame})}
	+ \epsilon r^{-1-\delta} |\slashed{\D} \mathscr{Y}^n h_{(\text{rect})}|^2
	+ \epsilon r^{-3-\delta} |\mathscr{Y}^n h_{(\text{rect})}|^2
	\\
	&\phantom{\lesssim \int_{\mathcal{M}_\tau^{\tau_1} \cap \{r \geq r_0\}} \bigg(}
	+ \delta^{-1} r^{-1+2\delta} |\overline{\slashed{\D}} \mathscr{Y}^n h|^2_{(\text{frame})}
	\bigg) \dVol_{g}
	\end{split}
	\end{equation*}
		
\end{proposition}

\begin{proof}
	Recall proposition \ref{proposition Yn omega}, which gives
	\begin{equation*}
	\begin{split}
	&\int_{\mathcal{M}_{\tau}^{\tau_1} \cap \{r \geq r_0\}} r^{-1-C_{[n]}\epsilon} |\mathscr{Y}^n \omega|^2 \dVol_{g}
	\\
	&\lesssim
	\int_{\mathcal{M}_{\tau}^{\tau_1} \cap \{r \geq r_0\}} r^{-1-C_{[n]}\epsilon} \bigg(
	\left( r^{-2} + \epsilon^2 r^{-2 + 2C_{(0)}\epsilon} \right) |\mathscr{Y}^n \bar{X}_{(\text{frame})}|^2
	+ |\bm{\Gamma}^{(0)}_{(-1-\delta)}|^2 |\mathscr{Y}^n X_{(\text{frame})}|^2
	+ |\slashed{\D} \mathscr{Y}^n h|^2_{LL}
	\\
	&\phantom{\lesssim \int_{\Sigma_\tau \cap \{r \geq r_0\}} r w_n \bigg(}		
	+ \epsilon^2 r^{-2\delta} |\bm{\Gamma}^{(n)}_{(-1+C_{(n)}\epsilon)}|^2
	+ |\bm{\Gamma}^{(n-1)}_{(-1+3C_{(n-1)}\epsilon)}|^2
	\bigg) \upd r \wedge \dVol_{\mathbb{S}^2}
	\end{split}
	\end{equation*}
	and so, using proposition \ref{proposition L2 bounds rectangular} together with the remark \ref{remark L2 bounds involving Yn Xframe} we have
	\begin{equation*}
	\begin{split}
	&\int_{\mathcal{M}_\tau^{\tau_1} \cap \{r \geq r_0\}} r^{-1-C_{[n]}\epsilon} |\mathscr{Y}^n \omega|^2 \dVol_{g}
	\\
	&\lesssim
	\int_{\mathcal{M}_\tau^{\tau_1} \cap \{r \geq r_0\}} \bigg(
	r^{-1-C_{[n]}\epsilon} |\slashed{\D} \mathscr{Y}^n h|^2_{(\text{frame})}
	+ \epsilon r^{-1-\delta} |\slashed{\D} \mathscr{Y}^n h_{(\text{rect})}|^2
	+ \epsilon r^{-3-\delta} |\mathscr{Y}^n h_{(\text{rect})}|^2
	\\
	&\phantom{\lesssim \int_{\mathcal{M}_\tau^{\tau_1} \cap \{r \geq r_0\}} \bigg(}
	+ \delta^{-1} r^{-1+2\delta} |\overline{\slashed{\D}} \mathscr{Y}^n h|^2_{(\text{frame})}
	+ \epsilon^2 r^{-1-\delta} |\bm{\Gamma}^{(n)}_{(-1+C_{(n)}\epsilon)}|^2
	\\
	&\phantom{\lesssim \int_{\mathcal{M}_\tau^{\tau_1} \cap \{r \geq r_0\}} \bigg(}
	+ \delta^{-1} r^{-1+(2-\frac{1}{4}c_{[n]})\delta} |\bm{\Gamma}^{(n-1)}_{(-1-\delta)}|^2
	+ r^{-1-\frac{1}{2}C_{[n]}\epsilon} |\bm{\Gamma}^{(n-1)}_{-1+C_{(n-1)}\epsilon}|^2 
	\bigg) \dVol_{g}
	\end{split}
	\end{equation*}
	Finally, using the $L^2$ bootstrap bounds in equation \eqref{equation bootstrap L2 geometric}, we find
	\begin{equation*}
	\begin{split}
	&\int_{\mathcal{M}_\tau^{\tau_1} \cap \{r \geq r_0\}} r^{-1-C_{[n]}\epsilon} |\mathscr{Y}^n \omega|^2 \dVol_{g}
	\\
	&\lesssim
	\left( \epsilon + \frac{\epsilon^2}{\delta} + \frac{\epsilon^2}{\delta}\right) \epsilon^{2(N + 1 - n)}(1+\tau)^{-1 + C_{[n]}\delta}
	\\
	&\phantom{\lesssim}
	+ \int_{\mathcal{M}_\tau^{\tau_1} \cap \{r \geq r_0\}} \bigg(
	r^{-1-C_{[n]}\epsilon} |\slashed{\D} \mathscr{Y}^n h|^2_{(\text{frame})}
	+ \epsilon r^{-1-\delta} |\slashed{\D} \mathscr{Y}^n h_{(\text{rect})}|^2
	+ \epsilon r^{-3-\delta} |\mathscr{Y}^n h_{(\text{rect})}|^2
	\bigg) \dVol_{g}
	\end{split}
	\end{equation*}
	In particular, using $\epsilon \ll \delta$ we have
	\begin{equation*}
	\begin{split}
	&\int_{\mathcal{M}_\tau^{\tau_1} \cap \{r \geq r_0\}} C_{[n]}\epsilon r^{-1-C_{[n]}\epsilon} |\mathscr{Y}^n \omega|^2 \dVol_{g}
	\\
	&\lesssim
	C_{[n]}\epsilon^2 \epsilon^{2(N + 1 - n)}(1+\tau)^{-1 + C_{[n]}\delta}
	\\
	&\phantom{\lesssim}
	+ \int_{\mathcal{M}_\tau^{\tau_1} \cap \{r \geq r_0\}} \bigg(
		C_{[n]}\epsilon r^{-1-C_{[n]}\epsilon} |\slashed{\D} \mathscr{Y}^n h|^2_{(\text{frame})}
		+ C_{[n]}\epsilon^2 r^{-1-\delta} |\slashed{\D} \mathscr{Y}^n h_{(\text{rect})}|^2
		\\
		&\phantom{\lesssim + \int_{\mathcal{M}_\tau^{\tau_1} \cap \{r \geq r_0\}} \bigg(}
		+ C_{[n]}\epsilon^2 r^{-3-\delta} |\mathscr{Y}^n h_{(\text{rect})}|^2
	\bigg) \dVol_{g}
	\end{split}
	\end{equation*}

	Following along very similar lines, we have
	\begin{equation*}
	\begin{split}
	&\int_{\mathcal{M}_\tau^{\tau_1} \cap \{r \geq r_0\}} r^{-1-c\delta} |\mathscr{Y}^n \omega|^2 \dVol_{g}
	\\
	&\lesssim
	\int_{\mathcal{M}_\tau^{\tau_1} \cap \{r \geq r_0\}} \bigg(
	r^{-1-c\delta} |\slashed{\D} \mathscr{Y}^n h|^2_{(\text{frame})}
	+ \epsilon r^{-1-\delta} |\slashed{\D} \mathscr{Y}^n h_{(\text{rect})}|^2
	+ \epsilon r^{-3-\delta} |\mathscr{Y}^n h_{(\text{rect})}|^2
	\\
	&\phantom{\lesssim \int_{\mathcal{M}_\tau^{\tau_1} \cap \{r \geq r_0\}} \bigg(}
	+ \delta^{-1} r^{-1+2\delta} |\overline{\slashed{\D}} \mathscr{Y}^n h|^2_{(\text{frame})}
	+ \epsilon^2 r^{-1-\delta} |\bm{\Gamma}^{(n)}_{(-1+C_{(n)}\epsilon)}|^2
	\\
	&\phantom{\lesssim \int_{\mathcal{M}_\tau^{\tau_1} \cap \{r \geq r_0\}} \bigg(}
	+ \delta^{-1} r^{-1+(2-\frac{1}{4}c_{[n]})\delta} |\bm{\Gamma}^{(n-1)}_{(-1-\delta)}|^2
	+ r^{-1-\frac{1}{2}c\delta} |\bm{\Gamma}^{(n-1)}_{-1+C_{(n-1)}\epsilon}|^2 
	\bigg) \dVol_{g}
	\end{split}
	\end{equation*}
	and again using the pointwise bootstrap bounds we find
	\begin{equation*}
	\begin{split}
	&\int_{\mathcal{M}_\tau^{\tau_1} \cap \{r \geq r_0\}} r^{-1-c\delta} |\mathscr{Y}^n \omega|^2 \dVol_{g}
	\\
	&\lesssim
	\epsilon^2 \left( \delta^{-1} + \delta^{-2} + (c\delta)^{-1} \right) \epsilon^{2(N + 1 - n)}(1+\tau)^{-1 + C_{[n]}\delta}
	\\
	&\phantom{\lesssim}
	+ \int_{\mathcal{M}_\tau^{\tau_1} \cap \{r \geq r_0\}} \bigg(
	r^{-1-c\delta} |\slashed{\D} \mathscr{Y}^n h|^2_{(\text{frame})}
	+ \epsilon r^{-1-\delta} |\slashed{\D} \mathscr{Y}^n h_{(\text{rect})}|^2
	+ \epsilon r^{-3-\delta} |\mathscr{Y}^n h_{(\text{rect})}|^2
	\\
	&\phantom{\lesssim \int_{\mathcal{M}_\tau^{\tau_1} \cap \{r \geq r_0\}} \bigg(}
	+ \delta^{-1} r^{-1+2\delta} |\overline{\slashed{\D}} \mathscr{Y}^n h|^2_{(\text{frame})}
	\bigg) \dVol_{g}
	\end{split}
	\end{equation*}
	so, using the facts that $\delta \ll 1$ proves the proposition.
	
\end{proof}

\begin{proposition}[$L^2$ bounds for $\omega$ at top order]
	\label{proposition L2 omega high}
	Suppose that the bootstrap bounds hold. Then we have
	\begin{equation*}
	\begin{split}
	&\int_{\mathcal{M}_\tau^{\tau_1} \cap \{r \geq r_0\}} \delta r^{-1-\delta} |\mathscr{Y}^{N_2} \omega|^2 \dVol_{g}
	\\
	&\lesssim
	\delta^{-1} \epsilon^4 (1+\tau)^{-1+C_{[N_2]}\delta}
	+ \int_{\mathcal{M}_\tau^{\tau_1} \cap \{r \geq r_0\}} \bigg(
	\delta r^{-1-\delta} |\slashed{\D} \mathscr{Y}^{N_2} h|^2_{(\text{frame})}
	+ \epsilon \delta r^{-1-\delta} |\slashed{\D} \mathscr{Y}^{N_2} h|^2_{(\text{frame})}
	\\
	&\phantom{\lesssim \delta^{-1} \epsilon^4 (1+\tau)^{-1+C_{[N_2]}\delta} + \int_{\mathcal{M}_\tau^{\tau_1} \cap \{r \geq r_0\}} \bigg(}
	+ \epsilon \delta r^{-3-\delta} |\mathscr{Y}^{N_2} h|^2_{(\text{frame})}
	+ r^{-1+2\delta} |\overline{\slashed{\D}} \mathscr{Y}^{N_2} h|^2_{(\text{frame})}
	\bigg) \dVol_{g}
	\end{split}
	\end{equation*}
\end{proposition}

\begin{proof}
	Following identical calculations to those in proposition \ref{proposition L2 omega low} we have
	\begin{equation*}
	\begin{split}
	&\int_{\mathcal{M}_\tau^{\tau_1} \cap \{r \geq r_0\}} r^{-1-\delta} |\mathscr{Y}^{N_2} \omega|^2 \dVol_{g}
	\\
	&\lesssim
	\int_{\mathcal{M}_\tau^{\tau_1} \cap \{r \geq r_0\}} \bigg(
		r^{-1-\delta} |\slashed{\D} \mathscr{Y}^{N_2} h|^2_{(\text{frame})}
		+ \epsilon r^{-1-\delta} |\slashed{\D} \mathscr{Y}^{N_2} h_{(\text{rect})}|^2
		+ \epsilon r^{-3-\delta} |\mathscr{Y}^{N_2} h_{(\text{rect})}|^2
		\\
		&\phantom{\lesssim \int_{\mathcal{M}_\tau^{\tau_1} \cap \{r \geq r_0\}} \bigg(}
		+ \delta^{-1} r^{-1+2\delta} |\overline{\slashed{\D}} \mathscr{Y}^{N_2} h|^2_{(\text{frame})}
		+ \epsilon^2 r^{-1-\delta} |\bm{\Gamma}^{(N_2)}_{(-1+C_{(N_2)}\epsilon)}|^2
		\\
		&\phantom{\lesssim \int_{\mathcal{M}_\tau^{\tau_1} \cap \{r \geq r_0\}} \bigg(}
		+ \delta^{-1} r^{-1+(2-\frac{1}{4}c_{[N_2]})\delta} |\bm{\Gamma}^{(N_2-1)}_{(-1-\delta)}|^2
		+ r^{-1-\frac{1}{2}c\delta} |\bm{\Gamma}^{(N_2-1)}_{-1+C_{(N_2-1)}\epsilon}|^2 
	\bigg) \dVol_{g}
	\end{split}
	\end{equation*}
	so, using the $L^2$ bootstrap bounds from equations \eqref{equation bootstrap L2 geometric} and \eqref{equation bootstrap L2 geometric top order} we find
	\begin{equation*}
	\begin{split}
	&\int_{\mathcal{M}_\tau^{\tau_1} \cap \{r \geq r_0\}} r^{-1-\delta} |\mathscr{Y}^{N_2} \omega|^2 \dVol_{g}
	\\
	&\lesssim
	\epsilon^2 (\delta^{-1} + \delta^{-2}) \epsilon^2 (1+\tau)^{-1+C_{[N_2]}\delta}
	\left( 1 + \delta^{-1} \epsilon^2 + c^{-2} \delta^{-2} \epsilon^2 \right) \epsilon^4 (1+\tau)^{-1+C_{[N_2]}\delta}
	\\
	&\phantom{\lesssim}
	+ \int_{\mathcal{M}_\tau^{\tau_1} \cap \{r \geq r_0\}} \bigg(
		r^{-1-\delta} |\slashed{\D} \mathscr{Y}^{N_2} h|^2_{(\text{frame})}
		+ \epsilon r^{-1-\delta} |\slashed{\D} \mathscr{Y}^{N_2} h_{(\text{rect})}|^2
		+ \epsilon r^{-3-\delta} |\mathscr{Y}^{N_2} h_{(\text{rect})}|^2
		\\
		&\phantom{\lesssim + \int_{\mathcal{M}_\tau^{\tau_1} \cap \{r \geq r_0\}} \bigg(}
		+ \delta^{-1} r^{-1+2\delta} |\overline{\slashed{\D}} \mathscr{Y}^{N_2} h|^2_{(\text{frame})}
	\bigg) \dVol_{g}
	\end{split}
	\end{equation*}
	Finally, using the fact that $\delta \ll 1$ proves the proposition.
\end{proof}

\begin{proposition}[$L^2$ bounds on $\mathscr{Y}^n \zeta$ below top order]
	\label{proposition L2 zeta low}
	Suppose that all of the bootstrap bounds hold.
	
	Then, for $n \leq N_2 -1$ we have
	\begin{equation*}
	\begin{split}
	&\int_{\mathcal{M}_{\tau}^{\tau_1} \cap \{r \geq r_0\}} C_{[n]}\epsilon r^{-1-2C_{[n]}\epsilon} |\mathscr{Y}^n \zeta|^2 \dVol_g
	\\
	&\lesssim
	C_{[n]}\epsilon^2 \epsilon^{2(N_2 + 1 - n)} (1+\tau)^{-1+C_{[n]}\delta}
	\\
	&\phantom{\lesssim}	
	+ \int_{\mathcal{M}_\tau^{\tau_1} \cap \{r \geq r_0\}} \bigg(
		C_{[n]}\epsilon r^{-1-C_{[n]}\epsilon} |\slashed{\D} \mathscr{Y}^n h|^2_{(\text{frame})}
		+ C_{[n]}\epsilon r^{-1-\frac{1}{2}C_{[n]}\epsilon} |\overline{\slashed{\D}} \mathscr{Y}^n h|^2_{(\text{frame})}
		\\
		&\phantom{\lesssim + \int_{\mathcal{M}_\tau^{\tau_1} \cap \{r \geq r_0\}} \bigg(}
		+ C_{[n]}\epsilon^2 r^{-1-\delta} |\slashed{\D} \mathscr{Y}^n h|^2_{(\text{rect})}
		+ C_{[n]}\epsilon (\epsilon + (C_{[n]})^{-2}) r^{-3-\delta} |\mathscr{Y}^n h|^2_{(\text{frame})}
		\\
		&\phantom{\lesssim + \int_{\mathcal{M}_\tau^{\tau_1} \cap \{r \geq r_0\}} \bigg(}
		+ C_{[n]}\epsilon (\delta^{-1} + (C_{[n]})^{-2}) r^{-1+2\delta} |\overline{\slashed{\D}} \mathscr{Y}^n h|^2_{(\text{frame})}
	\bigg) \dVol_{g}
	\end{split}
	\end{equation*}
	and
	\begin{equation*}
	\begin{split}
	&\int_{\mathcal{M}_{\tau}^{\tau_1} \cap \{r \geq r_0\}} c\delta r^{-1-c\delta} |\mathscr{Y}^n \zeta|^2 \dVol_g
	\\
	&\lesssim
	\delta^{-1} \epsilon^2 \epsilon^{2(N_2 + 1 - n)} (1+\tau)^{-1+C_{(n)}\delta}
	\\
	&\phantom{\lesssim}
	+ \int_{\mathcal{M}_\tau^{\tau_1} \cap \{r \geq r_0\}} \bigg(
		c\delta r^{-1-\frac{1}{2}c\delta} |\overline{\slashed{\D}} \mathscr{Y}^n h|^2_{(\text{frame})}
		+ c\delta r^{-1-c\delta} |\slashed{\D} \mathscr{Y}^n h|^2_{(\text{frame})}
		+ c\delta \epsilon r^{-3-\delta} |\mathscr{Y}^n h|^2_{(\text{frame})}
		\\
		&\phantom{\lesssim+  \int_{\mathcal{M}_\tau^{\tau_1} \cap \{r \geq r_0\}} \bigg(}
		+ c r^{-1+2\delta} |\overline{\slashed{\D}} \mathscr{Y}^n h|^2_{(\text{frame})}
	\bigg) \dVol_{g}
	\end{split}
	\end{equation*}
	
\end{proposition}

\begin{proof}
	Recall proposition \ref{proposition Yn zeta}, which gives us
	\begin{equation*}
	\begin{split}
	&\int_{\mathcal{M}_{\tau}^{\tau_1} \cap \{r \geq r_0\}} r^{-1-2C_{[n]}\epsilon} |\mathscr{Y}^n \zeta|^2 \dVol_g
	\\
	&\lesssim
	\int_{\mathcal{M}_{\tau}^{\tau_1} \cap \{r \geq r_0\}} \bigg(
		r^{-1-C_{[n]}\epsilon} |\slashed{\D} \mathscr{Y}^n h|^2_{(\text{frame})}
		+ r^{-3-C_{[n]}\epsilon} |\mathscr{Y}^n \bar{X}_{(\text{frame})}|^2
		\\
		&\phantom{\lesssim 	\int_{\mathcal{M}_{\tau}^{\tau_1} \cap \{r \geq r_0\}} \bigg(}
		+ r^{-1-C_{[n]}\epsilon} |\bm{\Gamma}^{(0)}_{(-1+C_{(0)}\epsilon)}|^2 |\mathscr{Y}^n X_{(\text{frame})}|^2
		+ \epsilon^2 r^{-1-2\delta} |\bm{\Gamma}^{(n)}_{(-1+C_{(n)}\epsilon)}|^2
		\\
		&\phantom{\lesssim 	\int_{\mathcal{M}_{\tau}^{\tau_1} \cap \{r \geq r_0\}} \bigg(}
		+ r^{-1-C_{[n]}\epsilon} |\bm{\Gamma}^{(n-1)}_{(-1+C_{(n-1)}\epsilon)}|^2
	\bigg) \dVol_{g}
	\end{split}
	\end{equation*}
	so, again making use of remark \ref{remark L2 bounds involving Yn Xframe} and, subsequently, proposition \ref{proposition L2 bounds rectangular} we have
	\begin{equation*}
	\begin{split}
	&\int_{\mathcal{M}_{\tau}^{\tau_1} \cap \{r \geq r_0\}} r^{-1-C_{[n]}\epsilon} |\mathscr{Y}^n \zeta|^2 \dVol_g
	\\
	&\lesssim
	\int_{\mathcal{M}_\tau^{\tau_1} \cap \{r \geq r_0\}} \bigg(
		r^{-1-C_{[n]}\epsilon} |\slashed{\D} \mathscr{Y}^n h|^2_{(\text{frame})}
		+ \epsilon r^{-1-\delta} |\slashed{\D} \mathscr{Y}^n h_{(\text{rect})}|^2
		+ \epsilon r^{-3-\delta} |\mathscr{Y}^n h_{(\text{rect})}|^2
		\\
		&\phantom{\lesssim \int_{\mathcal{M}_\tau^{\tau_1} \cap \{r \geq r_0\}} \bigg(}
		+ \delta^{-1} r^{-1+2\delta} |\overline{\slashed{\D}} \mathscr{Y}^n h|^2_{(\text{frame})}
		+ \epsilon^2 r^{-3-\frac{1}{2}C_{[n]}\epsilon} |\mathscr{Y}^{n} X_{(\text{frame, small})}|^2
		\\
		&\phantom{\lesssim \int_{\mathcal{M}_\tau^{\tau_1} \cap \{r \geq r_0\}} \bigg(}
		+ \epsilon^2 r^{-1-2\delta} |\bm{\Gamma}^{(n)}_{(-1+C_{(n)}\epsilon)}|^2
		+ \delta^{-1} r^{-1+(2-\frac{1}{4}c_{[n]})\delta} |\bm{\Gamma}^{(n-1)}_{(-1-\delta)}|^2
		\\
		&\phantom{\lesssim \int_{\mathcal{M}_\tau^{\tau_1} \cap \{r \geq r_0\}} \bigg(}
		+ r^{-1-C_{(n-1)}\epsilon} |\bm{\Gamma}^{(n-1)}_{-1+C_{(n-1)}\epsilon}|^2 
	\bigg) r^2 \upd r \wedge \dVol_{\mathbb{S}^2}
	\\
	\\
	&\lesssim
	\left(\frac{\epsilon^2}{\delta (C_{[n]})^2} + \frac{\epsilon^2}{C_{[n]}C_{[n-1]}\epsilon} \right) \epsilon^{2(N_2 + 1 - n)} (1+\tau)^{-1+C_{[n]}\delta}
	\\
	&\phantom{\lesssim}	
	+ \int_{\mathcal{M}_\tau^{\tau_1} \cap \{r \geq r_0\}} \bigg(
		r^{-1-C_{[n]}\epsilon} |\slashed{\D} \mathscr{Y}^n h|^2_{(\text{frame})}
		+ r^{-1-\frac{1}{2}C_{[n]}\epsilon} |\overline{\slashed{\D}} \mathscr{Y}^n h|^2_{(\text{frame})}
		+ \epsilon r^{-1-\delta} |\slashed{\D} \mathscr{Y}^n h|^2_{(\text{rect})}
		\\
		&\phantom{\lesssim + \int_{\mathcal{M}_\tau^{\tau_1} \cap \{r \geq r_0\}} \bigg(}
		+ (\epsilon + (C_{[n]})^{-2}) r^{-3-\delta} |\mathscr{Y}^n h|^2_{(\text{frame})}
		+ (\delta^{-1} + (C_{[n]})^{-2}) r^{-1+2\delta} |\overline{\slashed{\D}} \mathscr{Y}^n h|^2_{(\text{frame})}
		\\
		&\phantom{\lesssim + \int_{\mathcal{M}_\tau^{\tau_1} \cap \{r \geq r_0\}} \bigg(}
		+ \epsilon^2 r^{-1-2\delta} |\bm{\Gamma}^{(n)}_{(-1+C_{(n)}\epsilon)}|^2
		+ \delta^{-1} r^{-1+(2-\frac{1}{4}c_{[n]})\delta} |\bm{\Gamma}^{(n-1)}_{(-1-\delta)}|^2
		\\
		&\phantom{\lesssim + \int_{\mathcal{M}_\tau^{\tau_1} \cap \{r \geq r_0\}} \bigg(}
		+ r^{-1-C_{(n-1)}\epsilon} |\bm{\Gamma}^{(n-1)}_{-1+C_{(n-1)}\epsilon}|^2 
	\bigg) \dVol_{g}
	\end{split}
	\end{equation*}
	So, using the $L^2$ bootstrap bounds from equation \eqref{equation bootstrap L2 geometric} to bound the terms which are expressed schematically, we have
	\begin{equation*}
	\begin{split}
	&\int_{\mathcal{M}_{\tau}^{\tau_1} \cap \{r \geq r_0\}} r^{-1-2C_{[n]}\epsilon} |\mathscr{Y}^n \zeta|^2 \dVol_g
	\\
	&\lesssim
	\left(\epsilon + \frac{\epsilon^2}{\delta} + \frac{\epsilon^2}{\delta^2} + \frac{\epsilon^2}{\delta (C_{[n]})^2} + \frac{\epsilon}{C_{[n]}C_{[n-1]}} \right) \epsilon^{2(N_2 + 1 - n)} (1+\tau)^{-1+C_{[n]}\delta}
	\\
	&\phantom{\lesssim}	
	+ \int_{\mathcal{M}_\tau^{\tau_1} \cap \{r \geq r_0\}} \bigg(
		r^{-1-C_{[n]}\epsilon} |\slashed{\D} \mathscr{Y}^n h|^2_{(\text{frame})}
		+ r^{-1-\frac{1}{2}C_{[n]}\epsilon} |\overline{\slashed{\D}} \mathscr{Y}^n h|^2_{(\text{frame})}
		+ \epsilon r^{-1-\delta} |\slashed{\D} \mathscr{Y}^n h|^2_{(\text{rect})}
		\\
		&\phantom{\lesssim + \int_{\mathcal{M}_\tau^{\tau_1} \cap \{r \geq r_0\}} \bigg(}
		+ (\epsilon + (C_{[n]})^{-2}) r^{-3-\delta} |\mathscr{Y}^n h|^2_{(\text{frame})}
		+ (\delta^{-1} + (C_{[n]})^{-2}) r^{-1+2\delta} |\overline{\slashed{\D}} \mathscr{Y}^n h|^2_{(\text{frame})}
	\bigg) \dVol_{g}
	\end{split}
	\end{equation*}
	Finally, using the fact that $\epsilon \ll \delta$ and $C_{[n]} \gg 1$ proves the first part of the proposition.

	Following almost identical calculations, we can show that
	\begin{equation*}
	\begin{split}
	&\int_{\mathcal{M}_{\tau}^{\tau_1} \cap \{r \geq r_0\}} r^{-1-c\delta} |\mathscr{Y}^n \zeta|^2 \dVol_g
	\\
	&\lesssim
	\int_{\mathcal{M}_\tau^{\tau_1} \cap \{r \geq r_0\}} \bigg(
		r^{-1-c\delta} |\slashed{\D} \mathscr{Y}^n h|^2_{(\text{frame})}
		+ \epsilon r^{-1-\delta} |\slashed{\D} \mathscr{Y}^n h_{(\text{rect})}|^2
		+ \epsilon r^{-3-\delta} |\mathscr{Y}^n h_{(\text{rect})}|^2
		\\
		&\phantom{\lesssim \int_{\mathcal{M}_\tau^{\tau_1} \cap \{r \geq r_0\}} \bigg(}
		+ \delta^{-1} r^{-1+2\delta} |\overline{\slashed{\D}} \mathscr{Y}^n h|^2_{(\text{frame})}
		+ \epsilon^2 r^{-3-\frac{1}{2}c\delta} |\mathscr{Y}^{n} X_{(\text{frame, small})}|^2
		+ \epsilon^2 r^{-1-2\delta} |\bm{\Gamma}^{(n)}_{(-1+C_{(n)}\epsilon)}|^2
		\\
		&\phantom{\lesssim \int_{\mathcal{M}_\tau^{\tau_1} \cap \{r \geq r_0\}} \bigg(}
		+ \delta^{-1} r^{-1+(2-\frac{1}{4}c_{[n]})\delta} |\bm{\Gamma}^{(n-1)}_{(-1-\delta)}|^2
		+ r^{-1-\frac{1}{2}c\delta} |\bm{\Gamma}^{(n-1)}_{-1+C_{(n-1)}\epsilon}|^2 
	\bigg) \dVol_{g}
	\end{split}
	\end{equation*}
	so, again using proposition \ref{proposition L2 bounds rectangular} and the fact that $\epsilon \ll c\delta$ we find
	\begin{equation*}
	\begin{split}
	&\int_{\mathcal{M}_{\tau}^{\tau_1} \cap \{r \geq r_0\}} r^{-1-c\delta} |\mathscr{Y}^n \zeta|^2 \dVol_g
	\\
	&\lesssim
	\epsilon^3 c^{-2} \delta^{-2} \left( \frac{\epsilon}{\delta} + \frac{1}{C_{[n-1]}}\right) \epsilon^{2(N_2 + 1 - n)} (1+\tau)^{-1+C_{(n)}\delta}
	\\
	&\phantom{\lesssim}
	+ \int_{\mathcal{M}_\tau^{\tau_1} \cap \{r \geq r_0\}} \bigg(
		r^{-1-\frac{1}{2}c\delta} |\overline{\slashed{\D}} \mathscr{Y}^n h|^2_{(\text{frame})}
		+ r^{-1-c\delta} |\slashed{\D} \mathscr{Y}^n h|^2_{(\text{frame})}
		+ \epsilon r^{-3-\delta} |\mathscr{Y}^n h_{(\text{rect})}|^2
		\\
		&\phantom{\lesssim + \int_{\mathcal{M}_\tau^{\tau_1} \cap \{r \geq r_0\}} \bigg(}
		+ \delta^{-1} r^{-1+2\delta} |\overline{\slashed{\D}} \mathscr{Y}^n h|^2_{(\text{frame})}
		+ \epsilon^2 r^{-1-2\delta} |\bm{\Gamma}^{(n)}_{(-1+C_{(n)}\epsilon)}|^2
		\\
		&\phantom{\lesssim + \int_{\mathcal{M}_\tau^{\tau_1} \cap \{r \geq r_0\}} \bigg(}
		+ \delta^{-1} r^{-1+(2-\frac{1}{4}c_{[n]})\delta} |\bm{\Gamma}^{(n-1)}_{(-1-\delta)}|^2
		+ r^{-1-\frac{1}{2}c\delta} |\bm{\Gamma}^{(n-1)}_{-1+C_{(n-1)}\epsilon}|^2 
	\bigg) \dVol_{g}
	\end{split}
	\end{equation*}
	Finally, substituting the $L^2$ bounds from equation \eqref{equation bootstrap L2 geometric} we have
	\begin{equation*}
	\begin{split}
	&\int_{\mathcal{M}_{\tau}^{\tau_1} \cap \{r \geq r_0\}} r^{-1-c\delta} |\mathscr{Y}^n \zeta|^2 \dVol_g
	\\
	&\lesssim
	\left( \epsilon^3 c^{-2} \delta^{-2} \left( \frac{\epsilon}{\delta} + \frac{1}{C_{[n-1]}}\right) + \epsilon^2 \delta^{-1} + \epsilon^2 \delta^{-2} + \epsilon^2 (c\delta)^{-1} \right) \epsilon^{2(N_2 + 1 - n)} (1+\tau)^{-1+C_{(n)}\delta}
	\\
	&\phantom{\lesssim}
		+ \int_{\mathcal{M}_\tau^{\tau_1} \cap \{r \geq r_0\}} \bigg(
		r^{-1-\frac{1}{2}c\delta} |\overline{\slashed{\D}} \mathscr{Y}^n h|^2_{(\text{frame})}
		+ r^{-1-c\delta} |\slashed{\D} \mathscr{Y}^n h|^2_{(\text{frame})}
		+ \epsilon r^{-3-\delta} |\mathscr{Y}^n h|^2_{(\text{frame})}
		\\
		&\phantom{\lesssim+  \int_{\mathcal{M}_\tau^{\tau_1} \cap \{r \geq r_0\}} \bigg(}
		+ \delta^{-1} r^{-1+2\delta} |\overline{\slashed{\D}} \mathscr{Y}^n h|^2_{(\text{frame})}
	\bigg) \dVol_{g}
	\end{split}
	\end{equation*}
	Finally, using the fact that $\epsilon \ll c\delta \ll \delta \ll 1$ proves the proposition.
	
\end{proof}

\begin{proposition}[$L^2$ bounds on $\mathscr{Y}^n \zeta$ at top order]
	\label{proposition L2 zeta high}
	Suppose the bootstrap bounds hold. Then we have
	\begin{equation}
	\begin{split}
	&\int_{\mathcal{M}_{\tau}^{\tau_1} \cap \{r \geq r_0\}} \delta r^{-1-\delta} |\mathscr{Y}^{N_2} \zeta|^2 \dVol_g
	\\
	&\lesssim
	\epsilon^4 (1+\tau)^{-1+C_{(N_2)} \delta}
	+ \int_{\mathcal{M}_{\tau}^{\tau_1} \cap \{r \geq r_0\}} \bigg(
		\epsilon^2 \delta^{-1} r^{-1+\frac{1}{2}\delta} |\overline{\slashed{\D}}\mathscr{Y}^{N_2} h|^2_{(\text{frame})}
		+ \delta r^{-1-\delta} |\slashed{\D}\mathscr{Y}^{N_2} h|^2_{(\text{frame})}
		\\
		&\phantom{ \lesssim \delta^{-2} \epsilon^6 (1+\tau)^{-1+C_{(N_2)} \delta} + \int_{\mathcal{M}_{\tau}^{\tau_1} \cap \{r \geq r_0\}} \bigg(}
		+ \epsilon^2 r^{-3-\delta} |\mathscr{Y}^{N_2} h|^2_{(\text{frame})}
	\bigg) \dVol_{g}
	\end{split}
	\end{equation}
\end{proposition}

\begin{proof}
	From proposition \ref{proposition Yn zeta} we have
	\begin{equation*}
	\begin{split}
	&\int_{\mathcal{M}_{\tau}^{\tau_1} \cap \{r \geq r_0\}} r^{-1-\delta} |\mathscr{Y}^{N_2} \zeta|^2 \dVol_g
	\\
	&\lesssim
	\int_{\mathcal{M}_{\tau}^{\tau_1} \cap \{r \geq r_0\}} \bigg(
		r^{-1-\delta} |\slashed{\D} \mathscr{Y}^{N_2} h|^2_{(\text{frame})}
		+ r^{-1-\frac{1}{2}\delta} |\mathscr{Y}^{N_2} \bar{X}_{(\text{frame})}|^2
		\\
		&\phantom{\lesssim \int_{\Sigma_\tau \cap \{r \geq r_0\}} \bigg(}
		+ r^{-1-\delta} |\bm{\Gamma}^{(0)}_{(-1+C_{(0)}\epsilon)}|^2 |\mathscr{Y}^{N_2} X_{(\text{frame})}|^2
		+ \epsilon^2 r^{-1-2\delta} |\bm{\Gamma}^{(N_2)}_{(-1+C_{(N_2)}\epsilon)}|^2
		\\
		&\phantom{\lesssim \int_{\Sigma_\tau \cap \{r \geq r_0\}} \bigg(}
		+ r^{-1-\frac{1}{2}\delta} |\bm{\Gamma}^{(N_2-1)}_{(-1+C_{(N_2-1)}\epsilon)}|^2
	\bigg) \dVol_{g}
	\end{split}
	\end{equation*}
	Using remark \ref{remark L2 bounds involving Yn Xframe} and then substituting for $X_{(\text{frame, small})}$ from proposition \ref{proposition L2 bounds rectangular} (and also using $\epsilon \ll \delta$) we find
	\begin{equation*}
	\begin{split}
	&\int_{\mathcal{M}_{\tau}^{\tau_1} \cap \{r \geq r_0\}} r^{-1-\delta} |\mathscr{Y}^{N_2} \zeta|^2 \dVol_g
	\\
	&\lesssim
	\epsilon^6 (1+\tau)^{-1+C_{(N_2)} \delta}
	\\
	&\phantom{\lesssim}
	+ \int_{\mathcal{M}_{\tau}^{\tau_1} \cap \{r \geq r_0\}} \bigg(
	\epsilon^2 \delta^{-2} r^{-1+\frac{1}{2}\delta} |\overline{\slashed{\D}}\mathscr{Y}^{N_2} h|^2_{(\text{frame})}
	+ \epsilon^2 \delta^{-1} r^{-3-\delta} |\mathscr{Y}^{N_2} h|^2_{(\text{frame})}
	+ r^{-1-\delta} |\slashed{\D} \mathscr{Y}^{N_2} h|^2_{(\text{frame})}
	\\
	&\phantom{\lesssim + \int_{\Sigma_\tau \cap \{r \geq r_0\}} \bigg(}
	+ \epsilon^2 r^{-1-2\delta} |\bm{\Gamma}^{(N_2)}_{(-1+C_{(N_2)}\epsilon)}|^2
	+ r^{-1-\frac{1}{2}\delta} |\bm{\Gamma}^{(N_2-1)}_{(-1+C_{(N_2-1)}\epsilon)}|^2
	\bigg) \dVol_{g}
	\end{split}
	\end{equation*}
	Now substituting the $L^2$ bootstrap bounds in equation \eqref{equation bootstrap L2 geometric} and using $1 \gg \delta \gg \epsilon$ we have
	\begin{equation}
	\begin{split}
	&\int_{\mathcal{M}_{\tau}^{\tau_1} \cap \{r \geq r_0\}} r^{-1-\delta} |\mathscr{Y}^{N_2} \zeta|^2 \dVol_g
	\\
	&\lesssim
	\delta^{-1} \epsilon^4 (1+\tau)^{-1+C_{(N_2)} \delta}
	+ \int_{\mathcal{M}_{\tau}^{\tau_1} \cap \{r \geq r_0\}} \bigg(
	\epsilon^2 \delta^{-2} r^{-1+\frac{1}{2}\delta} |\overline{\slashed{\D}}\mathscr{Y}^{N_2} h|^2_{(\text{frame})}
	+ r^{-1-\delta} |\slashed{\D}\mathscr{Y}^{N_2} h|^2_{(\text{frame})}
	\\
	&\phantom{ \lesssim \delta^{-2} \epsilon^6 (1+\tau)^{-1+C_{(N_2)} \delta} + \int_{\mathcal{M}_{\tau}^{\tau_1} \cap \{r \geq r_0\}} \bigg(}
	+ \epsilon^2 \delta^{-1} r^{-3-\delta} |\mathscr{Y}^{N_2} h|^2_{(\text{frame})}
	\bigg) \dVol_{g}
	\end{split}
	\end{equation}
	 
\end{proof}

\begin{proposition}[$L^2$ bounds on $\tr_{\slashed{g}}\chi_{(\text{small})}$ below top-order]
	\label{proposition L2 tr chi low}
	Suppose that the bootstrap assumptions hold. Then for all $n \leq N_2 - 1$ we have
	\begin{equation*}
	\begin{split}
	&\int_{\mathcal{M}_\tau^{\tau_1} \cap \{r \geq r_0\}}
	\delta r^{-1+(\frac{1}{2}-c_{[n]})\delta} |\mathscr{Y}^n \tr_{\slashed{g}}\chi_{(\text{small})}|^2
	\dVol_{g}
	\\
	&\lesssim
	\epsilon^2 (1 + c_{[n]}^{-1}) \epsilon^{2(N_2+1-n)} (1+\tau)^{-1 + C_{(n)}\delta}
	\\
	&\phantom{\lesssim}
	+ \int_{\mathcal{M}_\tau^{\tau_1} \cap \{r \geq \frac{1}{2}r_0\}} \bigg(
		\delta r^{-1+(\frac{1}{2}-c_{[n]})\delta} |\overline{\slashed{\D}} \mathscr{Y}^{n+1} h|^2_{(\text{frame})}
		+ \delta r^{-1+\frac{1}{2}\delta} |\overline{\slashed{\D}} \mathscr{Y}^{n} h|^2_{(\text{frame})}
		\\
		&\phantom{\lesssim + \int_{\mathcal{M}_\tau^{\tau_1} \cap \{r \geq \frac{1}{2}r_0\}} \bigg(}
		+ \epsilon^2 \delta  r^{-1-\delta} |\slashed{\D} \mathscr{Y}^{n} h|^2_{(\text{frame})}
		+ \delta r^{-1-\delta} |\slashed{\D} \mathscr{Y}^n h|^2_{(\text{frame})}
		+ \delta r^{-3-\delta} |\mathscr{Y}^n h|^2_{(\text{frame})}
	\bigg) \dVol_{g}
	\end{split}
	\end{equation*}

\end{proposition}

\begin{proof}
	
	We apply the third part of proposition \ref{proposition Hardy}. This time we apply it to the field $\chi_{(r_0)}(r) r^2 \mathscr{Y}^n \mathcal{X}_{(\text{low})}$ and choose $\alpha = -3 + (\frac{1}{2}-c_{[n]})\delta$, obtaining
	\begin{equation*}
	\begin{split}
	&\int_{\Sigma_\tau \cap \{r \geq r_0\}}
	r^{1+(\frac{1}{2}-c_{[n]})\delta} |\mathscr{Y}^n \mathcal{X}_{(\text{low})}|^2
	\upd r \wedge \dVol_{\mathbb{S}^2}
	\\
	&\lesssim
	\int_{\Sigma_\tau \cap \{r \geq \frac{1}{2}r_0\}}
	r^{-1+(\frac{1}{2}-c_{[n]})\delta} |\slashed{\D}_L \left(r^2 \mathscr{Y}^n \mathcal{X}_{(\text{low})}\right)|^2
	\upd r \wedge \dVol_{\mathbb{S}^2}
	+ \int_{\Sigma_\tau \cap \{\frac{1}{2} r_0 \leq r \leq r_0\}}
	|\mathscr{Y}^n \mathcal{X}_{(\text{low})}|^2
	\upd r \wedge \dVol_{\mathbb{S}^2}
	\end{split}
	\end{equation*}
	We deal with each of these terms in turn. First, using proposition \ref{proposition Yn tr chi low} we have
	\begin{equation*}
	\begin{split}
	\slashed{\D}_L \left( r^2 \mathscr{Y}^n \mathcal{X}_{(\text{low})} \right)
	&=
	\bm{\Gamma}^{(0)}_{(-1)} r^2 \mathscr{Y}^n \mathcal{X}_{(\text{low})}
	+ r \bm{\Gamma}^{(0)}_{(0, \text{large})}\left( \overline{\slashed{\D}}\mathscr{Y}^{n+1} h \right)_{(\text{frame})}
	+ r\bm{\Gamma}^{(1)}_{(C_{(1)}\epsilon, \text{large})} \left( \overline{\slashed{\D}}\mathscr{Y}^{n} h \right)_{(\text{frame})}
	\\
	&\phantom{=}
	+ \bm{\Gamma}^{(0)}_{(1-\delta)} (\slashed{\D} \mathscr{Y}^n h)_{(\text{frame})}
	+ \bm{\Gamma}^{(1)}_{(-\delta)} (\mathscr{Y}^n X_{(\text{frame})})
	+ \bm{\Gamma}^{(0)}_{(-\frac{1}{2}+\delta)} \bm{\Gamma}^{(n)}_{(C_{(n)}\epsilon)}
	\\
	&\phantom{=}
	+ r^2 \bm{\Gamma}^{(0)}_{(-1-\delta)}\bm{\Gamma}^{(n)}_{(-1-\delta)}
	+ r \bm{\Gamma}^{(0)}_{(kC_{(0)}\epsilon, \text{large})}\bm{\Gamma}^{(n-1)}_{(-1-\delta)}
	\end{split}
	\end{equation*}
	and so we have the estimate
	\begin{equation*}
	\begin{split}
	&\int_{\Sigma_\tau \cap \{r \geq r_0\}}
	r^{1+(\frac{1}{2}-c_{[n]})\delta} |\mathscr{Y}^n \mathcal{X}_{(\text{low})}|^2
	\upd r \wedge \dVol_{\mathbb{S}^2}
	\\
	&\lesssim
	\int_{\Sigma_\tau \cap \{r \geq \frac{1}{2}r_0\}} \bigg(
		\epsilon^2 r^{-1+(\frac{1}{2}-c_{[n]})\delta} |\mathscr{Y}^n \mathcal{X}_{(\text{low})}|^2
		+ r^{-1+(\frac{1}{2}-c_{[n]})\delta} |\overline{\slashed{\D}} \mathscr{Y}^{n+1} h|^2_{(\text{frame})}
		\\
		&\phantom{\lesssim \int_{\Sigma_\tau \cap \{r \geq \frac{1}{2}r_0\}} \bigg(}
		+ r^{-1+\frac{1}{2}\delta} |\overline{\slashed{\D}} \mathscr{Y}^{n} h|^2_{(\text{frame})}
		+ \epsilon^2 r^{-1-\delta} |\slashed{\D} \mathscr{Y}^{n} h|^2_{(\text{frame})}
		\\
		&\phantom{\lesssim \int_{\Sigma_\tau \cap \{r \geq \frac{1}{2}r_0\}} \bigg(}
		+ \epsilon^2 r^{-3-c_{[n]}\delta} |\mathscr{Y}^n X_{(\text{frame, small})}|^2
		+ \epsilon^2 r^{-2 + 4\delta} |\bm{\Gamma}^{(n)}_{(-1 + C_{(n)}\epsilon)}|^2
		+ \epsilon^2 r^{-1 -\delta} |\bm{\Gamma}^{(n)}_{(-1-\delta)}|^2
		\\
		&\phantom{\lesssim \int_{\Sigma_\tau \cap \{r \geq \frac{1}{2}r_0\}} \bigg(}
		+ r^{-1 + (\frac{1}{2}-c_{[n-1]})\delta} |\bm{\Gamma}^{(n-1)}_{(-1-\delta)}|^2		
	\bigg) r^2 \upd r \wedge \dVol_{\mathbb{S}^2}
	\\
	&\phantom{\lesssim}
	+ \int_{\Sigma_\tau \cap \{\frac{1}{2} r_0 \leq r \leq r_0\}}
	|\mathscr{Y}^n \mathcal{X}_{(\text{low})}|^2
	\upd r \wedge \dVol_{\mathbb{S}^2}
	\end{split}
	\end{equation*}
	so now, absorbing the first term on the right hand side by the left hand side, integrating over $\tau$ and  using the $L^2$ bootstrap bounds from equation \eqref{equation bootstrap L2 geometric} we have
	\begin{equation*}
	\begin{split}
	&\int_{\mathcal{M}_\tau^{\tau_1} \cap \{r \geq r_0\}}
	r^{-1+(\frac{1}{2}-c_{[n]})\delta} |\mathscr{Y}^n \mathcal{X}_{(\text{low})}|^2
	\dVol_{g}
	\\
	&\lesssim
	\epsilon^2 \delta^{-1} (1 + c_{[n]}^{-1}) \epsilon^{2(N_2+1-n)} (1+\tau)^{-1 + C_{(n)}\delta}
	\\
	&\phantom{\lesssim}
	+ \int_{\Sigma_\tau \cap \{r \geq \frac{1}{2}r_0\}} \bigg(
	r^{-1+(\frac{1}{2}-c_{[n]})\delta} |\overline{\slashed{\D}} \mathscr{Y}^{n+1} h|^2_{(\text{frame})}
	+ r^{-1+\frac{1}{2}\delta} |\overline{\slashed{\D}} \mathscr{Y}^{n} h|^2_{(\text{frame})}
	\\
	&\phantom{\lesssim + \int_{\Sigma_\tau \cap \{r \geq \frac{1}{2}r_0\}} \bigg(}
	+ \epsilon^2 r^{-1-\delta} |\slashed{\D} \mathscr{Y}^{n} h|^2_{(\text{frame})}
	\bigg) \dVol_{g}
	\\
	&\phantom{\lesssim}
	+ \int_{\mathcal{M}_\tau^{\tau_1} \cap \{\frac{1}{2} r_0 \leq r \leq r_0\}}
	|\mathscr{Y}^n \mathcal{X}_{(\text{low})}|^2
	\upd r \wedge \dVol_{\mathbb{S}^2}
	\end{split}
	\end{equation*}
	
	Next we need to deal with the final integral, which is an integral only over the region $\frac{1}{2} r_0 \leq r \leq r_0$. We note that, in this region (and in fact in the whole region $r \leq r_0$) we can express $\chi_{(\text{small})}$ (and hence $\mathcal{X}_{(\text{low})}$) directly in terms of $h$ and $\partial h$. Indeed, we note that $\chi$ can expressed as
	\begin{equation*}
	\chi_\mu^{\phantom{\mu}\nu}
	=
	\slashed{\Pi}_a^{\phantom{a}\nu} (\slashed{\nabla}_\mu L^a)
	+ \slashed{\Pi}_\mu^{\phantom{\mu}b} \slashed{\Pi}_c^{\phantom{c}\nu} L^a \Gamma^{c}_{ab}
	\end{equation*}
	and so
	\begin{equation*}
	\tr_{\slashed{g}}\chi_{(\text{small})}
	=
	(\slashed{\nabla}_a L^a_{(\text{small})})
	+ \slashed{\Pi}_c^{\phantom{c}b} L^a \Gamma^{c}_{ab}
	\end{equation*}
	For $r \leq r_0$, the rectangular components of the vector fields $L^a$ and $\Lbar^a$ can be expressed in terms of $h_{ab}$ as mentioned above, and so for $r \leq r_0$,
	\begin{equation*}
	|\tr_{\slashed{g}} \chi_{\text{small})}| \lesssim |h_{(\text{rect})}| + |\partial h_{(\text{rect})}|
	\end{equation*}
	and similarly for higher derivatives. Hence we obtain
	\begin{equation*}
	\begin{split}
	&\int_{\mathcal{M}_\tau^{\tau_1} \cap \{r \geq r_0\}}
	r^{-1+(\frac{1}{2}-c_{[n]})\delta} |\mathscr{Y}^n \mathcal{X}_{(\text{low})}|^2
	\dVol_{g}
	\\
	&\lesssim
	\epsilon^2 \delta^{-1} (1 + c_{[n]}^{-1}) \epsilon^{2(N_2+1-n)} (1+\tau)^{-1 + C_{(n)}\delta}
	\\
	&\phantom{\lesssim}
	+ \int_{\mathcal{M}_\tau^{\tau_1} \cap \{r \geq \frac{1}{2}r_0\}} \bigg(
	r^{-1+(\frac{1}{2}-c_{[n]})\delta} |\overline{\slashed{\D}} \mathscr{Y}^{n+1} h|^2_{(\text{frame})}
	+ r^{-1+\delta} |\overline{\slashed{\D}} \mathscr{Y}^{n} h|^2_{(\text{frame})}
	\\
	&\phantom{\lesssim + \int_{\Sigma_\tau \cap \{r \geq \frac{1}{2}r_0\}} \bigg(}
	+ \epsilon^2 r^{-1-\delta} |\slashed{\D} \mathscr{Y}^{n} h|^2_{(\text{frame})}
	+ r^{-1-\delta} |\slashed{\D} \mathscr{Y}^n h|^2_{(\text{frame})}
	+ r^{-3-\delta} |\mathscr{Y}^n h|^2_{(\text{frame})}
	\bigg) \dVol_{g}
	\end{split}
	\end{equation*}
	
	Finally, we need to relate $\mathcal{X}_{(\text{low})}$ to $\tr_{\slashed{g}}\chi_{(\text{small})}$. Proposition \ref{proposition Yn tr chi low} gives
	\begin{equation*}
	\begin{split}
	|\mathscr{Y}^n \tr_{\slashed{g}}\chi_{(\text{small})}|
	&\lesssim
	|\mathscr{Y}^n \mathcal{X}_{(\text{low})}|
	+ |\overline{\slashed{\D}}\mathscr{Y}^n h|_{(\text{frame})}
	+ |\bm{\Gamma}^{(0)}_{(-1-\delta)}| |\mathscr{Y}^n X_{(\text{frame})}|
	+ \epsilon (1+r)^{-\frac{1}{2} + \delta} |\mathscr{Y}^n \bm{\Gamma}^{(n)}_{(-1+C_{(n)}\epsilon)}|
	\\
	&\phantom{\lesssim}
	+ |\bm{\Gamma}^{(n-1)}_{-1-\delta}|
	\end{split}
	\end{equation*}
	This leads to error terms in the $L^2$ bound for $\tr_{\slashed{g}}\chi_{(\text{small})}$ which are similar to those already encountered. Handling these in the same way as before, we prove the proposition.
	
\end{proof}

\begin{proposition}[$L^2$ bounds on $\tr_{\slashed{g}}\chi_{(\text{small})}$ at top order]
	\label{proposition L2 tr chi high}
	Suppose that the bootstrap bounds hold. Additionally, suppose that $N_2 \leq 2N_1$, and that the inhomogeneous terms $F$ satisfy the pointwise bounds
	\begin{equation*}
	|\mathscr{Y}^n F|_{(\text{frame})} \lesssim \epsilon (1+r)^{-2+2C_{(n)}\epsilon}
	\end{equation*}
	for all $n \leq N_1$.
	
	Then we have
	\begin{equation*}
	\begin{split}
	&\int_{\mathcal{M}_\tau^{\tau_1} \cap \{r \geq r_0\}}
	C_{[N_2]}\epsilon r^{-1-C_{[N_2]}\epsilon} |\mathscr{Y}^n \tr_{\slashed{g}}\chi_{(\text{small})}|^2
	\dVol_{g}
	\\
	&\lesssim
	C_{[N_2]}\epsilon^5 (1+\tau)^{-1 + C_{(N_2)}\delta}
	\\
	&\phantom{\lesssim}
	+ \int_{\mathcal{M}_\tau^{\tau_1} \cap \{r \geq \frac{1}{2}r_0\}}\bigg(
		C_{[N_2]}\epsilon r^{-1-C_{[N_2]}\epsilon} |\slashed{\D} \mathscr{Y}^{N_2} h|^2_{LL}
		+ C_{[N_2]}\epsilon^3 r^{-1-\frac{1}{2}C_{[N_2]}\epsilon} |\overline{\slashed{\D}} \mathscr{Y}^{N_2} h|^2_{(\text{frame})}
		\\
		&\phantom{\lesssim + \int_{\mathcal{M}_\tau^{\tau_1} \cap \{r \geq \frac{1}{2}r_0\}}\bigg(}
		+ C_{[N_2]}\epsilon r^{-1-\delta} |\slashed{\D} \mathscr{Y}^{N_2} h|^2_{(\text{frame})}
		+ C_{[N_2]}\epsilon r^{-3-\delta} |\mathscr{Y}^{N_2} h|^2_{(\text{frame})}
		\\
		&\phantom{\lesssim + \int_{\mathcal{M}_\tau^{\tau_1} \cap \{r \geq \frac{1}{2}r_0\}}\bigg(}
		+ C_{[N_2]}\epsilon r^{1-C_{[N_2]}\epsilon} |\mathscr{Y}^{N_2} F|^2_{LL}
		+ \sum_{j \leq N_2 - 1} C_{[N_2]}\epsilon^3 r^{1-\frac{1}{2}C_{[N_2]}\epsilon} |\mathscr{Y}^j F|_{(\text{frame})}^2
	\bigg) \dVol_{g}
	\end{split}
	\end{equation*}

\end{proposition}

\begin{proof}
	Following calculations almost identical to those at the start of the proof of proposition \ref{proposition L2 tr chi low}, we find
	\begin{equation*}
	\begin{split}
	\int_{\Sigma_\tau \cap \{r \geq r_0\}}
	r^{1-C_{[N_2]}\epsilon} |\mathscr{Y}^{N_2} \mathcal{X}_{(\text{high})}|^2
	\upd r \wedge \dVol_{\mathbb{S}^2}
	&\lesssim
	\int_{\Sigma_\tau \cap \{r \geq \frac{1}{2}r_0\}}
	r^{-1-C_{[N_2]}\epsilon} |\slashed{\D}_L \left(r^2 \mathscr{Y}^{N_2} \mathcal{X}_{(\text{high})}\right)|^2
	\upd r \wedge \dVol_{\mathbb{S}^2}
	\\
	&\phantom{\lesssim}
	+ \int_{\Sigma_\tau \cap \{\frac{1}{2} r_0 \leq r \leq r_0\}}
	|\mathscr{Y}^{N_2} \mathcal{X}_{(\text{high})}|^2
	\upd r \wedge \dVol_{\mathbb{S}^2}
	\end{split}
	\end{equation*}
	where, this time, we have applied proposition \ref{proposition Hardy} to the field $\chi_{(r_0)} r^2 \mathcal{X}_{(\text{high})}$ with the choice $\alpha = 3-C_{[N_2]}$.
	
	Now, using propositions \ref{proposition Yn tr chi high} together with proposition \ref{proposition Yn omega} we obtain
	\begin{equation*}
	\begin{split}
	\slashed{\D}_L \left( r^2 \mathscr{Y}^{N_2} \mathcal{X}_{(\text{high})} \right)
	&= 
	\bm{\Gamma}^{(0)}_{(-1)} \left( r^2 \mathscr{Y}^{N_2} \mathcal{X}_{(\text{high})} \right)
	+ r^2 (\mathscr{Y}^{N_2} F)_{LL}
	+ r^2 (F)_{(\text{frame})} (\mathscr{Y}^{N_2} \bar{X}_{(\text{frame})})
	\\
	&\phantom{=}
	+ r^2 \sum_{\substack{j+k \leq N_2 \\ j,k \leq N_2-1}} (\mathscr{Y}^j X_{(\text{frame})}) (\mathscr{Y}^{k} F)_{(\text{frame})}
	+ r (\slashed{\D} \mathscr{Y}^{N_2} h)_{LL}
	\\
	&\phantom{=}
	+ r^2 \bm{\Gamma}^{(0)}_{(-1-\delta)} (\slashed{\D} \mathscr{Y}^{N_2} h)_{(\text{frame})}
	+ r^2 \bm{\Gamma}^{(0)}_{(-1+C_{(0)}\epsilon)} (\overline{\slashed{\D}} \mathscr{Y}^{N_2} h)_{(\text{frame})}
	+ \bm{\Gamma}^{(0)}_{(-\frac{1}{2} + \delta)} \bm{\Gamma}^{(N_2)}_{(C_{(N_2)}\epsilon)}
	\\
	&\phantom{=}
	+ r^2 \bm{\Gamma}^{(n)}_{(-1-\delta)} \bm{\Gamma}^{(0)}_{(-1-\delta)}
	+ r^2 \sum_{\substack{j+k \leq N_2 \\ j,k \leq N_2-1}} \bm{\Gamma}^{(j)}_{(-1+C_{(j)}\epsilon)} \bm{\Gamma}^{(k)}_{(-1+C_{(k)}\epsilon)}
	\end{split}
	\end{equation*}
	and so, using remark \ref{remark L2 bounds involving Yn Xframe} to bound some of the ``frame'' terms, we find
	\begin{equation*}
	\begin{split}
	&\int_{\Sigma_\tau \cap \{r \geq r_0\}}
	r^{1-C_{[N_2]}\epsilon} |\mathscr{Y}^n \mathcal{X}_{(\text{high})}|^2
	\upd r \wedge \dVol_{\mathbb{S}^2}
	\\
	&\lesssim
	\int_{\Sigma_\tau \cap \{r \geq \frac{1}{2}r_0\}}\bigg(
	\epsilon^2 r^{1-C_{[N_2]}\epsilon} |\mathscr{Y}^n \mathcal{X}_{(\text{high})}|^2
	+ r^{3-C_{[N_2]}\epsilon} |\mathscr{Y}^{N_2} F|^2_{LL}
	+ r^{3-C_{[N_2]}\epsilon} |F|_{(\text{frame})}^2 |\mathscr{Y}^n \bar{X}_{(\text{frame})}|^2
	\\
	&\phantom{\lesssim \int_{\Sigma_\tau \cap \{r \geq \frac{1}{2}r_0\}}\bigg(}
	+ r^{5-C_{[N_2]}\epsilon} \sum_{\substack{j+k \leq N_2 \\ j,k \leq N_2-1}} |\bm{\Gamma}^{(j)}_{(-1+C_{(j)}\epsilon)}|^2 |\mathscr{Y}^k F|_{(\text{frame})}^2
	+ r^{1-C_{[N_2]}\epsilon} |\slashed{\D} \mathscr{Y}^{N_2} h|^2_{LL}
	\\
	&\phantom{\lesssim \int_{\Sigma_\tau \cap \{r \geq \frac{1}{2}r_0\}}\bigg(}
	+ \epsilon^2 r^{1-2\delta} |\slashed{\D} \mathscr{Y}^{N_2} h|^2_{(\text{frame})}
	+ \epsilon^2 r^{1-\frac{1}{2}C_{[N_2]}\epsilon} |\overline{\slashed{\D}} \mathscr{Y}^{N_2} h|^2_{(\text{frame})}
	+ \epsilon^2 r^{2\delta} |\bm{\Gamma}^{(N_2)}_{(-1+C_{(N_2)}\epsilon)}|^2
	\\
	&\phantom{\lesssim \int_{\Sigma_\tau \cap \{r \geq \frac{1}{2}r_0\}}\bigg(}
	+ \epsilon^2 r^{1-2\delta} |\bm{\Gamma}^{(N_2)}_{(-1-\delta)}|^2
	+ \sum_{\substack{j+k \leq n \\ j,k \leq N_2-1}} r^{3-C_{[N_2]}} |\bm{\Gamma}^{(j)}_{(-1+C_{(j)}\epsilon)}|^2 |\bm{\Gamma}^{(k)}_{(-1+C_{(k)}\epsilon)}|^2
	\bigg)\upd r \wedge \dVol_{\mathbb{S}^2}
	\\
	&\phantom{\lesssim}
	+ \int_{\Sigma_\tau \cap \{\frac{1}{2} r_0 \leq r \leq r_0\}}
	|\mathscr{Y}^{N_2} \mathcal{X}_{(\text{high})}|^2
	\upd r \wedge \dVol_{\mathbb{S}^2}
	\end{split}
	\end{equation*}
	
	Next, we use the fact that $\epsilon \ll 1$ to absorb the first term on the right hand side on the left hand side. We also use the assumption that $N_2 \leq 2N_1$, which means that we can bound some of the quadratic terms in $L^\infty$. Finally, we use the pointwise assumptions on $F$, which allow us to obtain the bound
	\begin{equation*}
	\begin{split}
	&\int_{\Sigma_\tau \cap \{r \geq r_0\}}
	r^{1-C_{[N_2]}\epsilon} |\mathscr{Y}^n \mathcal{X}_{(\text{high})}|^2
	\upd r \wedge \dVol_{\mathbb{S}^2}
	\\
	&\lesssim
	\int_{\Sigma_\tau \cap \{r \geq \frac{1}{2}r_0\}}\bigg(
	r^{3-C_{[N_2]}\epsilon} |\mathscr{Y}^{N_2} F|^2_{LL}
	+ \epsilon^2 r^{-1-\frac{1}{2}C_{[N_2]}\epsilon} |\mathscr{Y}^{N_2}	\bar{X}_{(\text{frame})}|^2	
	+ r^{1-C_{[N_2]}\epsilon} |\slashed{\D} \mathscr{Y}^{N_2} h|^2_{LL}
	\\
	&\phantom{\lesssim \int_{\Sigma_\tau \cap \{r \geq \frac{1}{2}r_0\}}\bigg(}	
	+ \epsilon^2 r^{1-2\delta} |\slashed{\D} \mathscr{Y}^{N_2} h|^2_{(\text{frame})}
	+ \epsilon^2 r^{1-\frac{1}{2}C_{[N_2]}\epsilon} |\overline{\slashed{\D}} \mathscr{Y}^{N_2} h|^2_{(\text{frame})}
	\\
	&\phantom{\lesssim \int_{\Sigma_\tau \cap \{r \geq \frac{1}{2}r_0\}}\bigg(}		
	+ \epsilon^2 r^{2\delta} |\bm{\Gamma}^{(N_2)}_{(-1+C_{(N_2)}\epsilon)}|^2
	+ \epsilon^2 r^{1-2\delta} |\bm{\Gamma}^{(N_2)}_{(-1-\delta)}|^2
	\\
	&\phantom{\lesssim \int_{\Sigma_\tau \cap \{r \geq \frac{1}{2}r_0\}}\bigg(}		
	+ \sum_{j \leq N_2 - 1} \epsilon^2 r^{1-\frac{1}{2}C_{[N_2]}\epsilon} |\bm{\Gamma}^{(j)}_{(-1+C_{(j)}\epsilon)}|^2
	+ \sum_{j \leq N_2 - 1} \epsilon^2 r^{3-\frac{1}{2}C_{[N_2]}\epsilon} |\mathscr{Y}^j F|_{(\text{frame})}^2
	\bigg)\upd r \wedge \dVol_{\mathbb{S}^2}
	\\
	&\phantom{\lesssim}
	+ \int_{\Sigma_\tau \cap \{\frac{1}{2} r_0 \leq r \leq r_0\}}
	|\mathscr{Y}^{N_2} \mathcal{X}_{(\text{high})}|^2
	\upd r \wedge \dVol_{\mathbb{S}^2}
	\end{split}
	\end{equation*}
	
	As before, we can bound the final integral (which is an integral only over the region $\frac{1}{2}r_0 \leq r \leq r_0$) by using the fact that $\chi$ (and hence $\mathcal{X}_{(\text{high})}$) can be expressed in terms of $h$ and $\partial h$ in this region. Finally, after integrating over $\tau$, we are in a position to use the $L^2$ bounds. Noting especially that the $L^2$ bound for $\bar{X}_{(\text{frame})}$ established in proposition \ref{proposition L2 bounds rectangular} holds also in the case $n = N_2$, we find
	\begin{equation*}
	\begin{split}
	&\int_{\mathcal{M}_\tau^{\tau_1} \cap \{r \geq r_0\}}
	r^{-1-C_{[N_2]}\epsilon} |\mathscr{Y}^n \mathcal{X}_{(\text{high})}|^2
	\dVol_{g}
	\\
	&\lesssim
	\epsilon^4 (1+\tau)^{-1 + C_{[N_2]}\delta}
	\\
	&\phantom{\lesssim}
	+ \int_{\mathcal{M}_\tau^{\tau_1} \cap \{r \geq \frac{1}{2}r_0\}}\bigg(
	r^{-1-C_{[N_2]}\epsilon} |\slashed{\D} \mathscr{Y}^{N_2} h|^2_{LL}
	+ \epsilon^2 r^{-1-\frac{1}{2}C_{[N_2]}\epsilon} |\overline{\slashed{\D}} \mathscr{Y}^{N_2} h|^2_{(\text{frame})}
	\\
	&\phantom{\lesssim \int_{\mathcal{M}_\tau^{\tau_1} \cap \{r \geq \frac{1}{2}r_0\}}\bigg(}
	+ \sum_{j\leq N_2} r^{-1-\delta} |\slashed{\D} \mathscr{Y}^j h|^2_{(\text{frame})}
	+ \sum_{j\leq N_2} r^{-3-\delta} |\mathscr{Y}^j h|^2_{(\text{frame})}
	\\
	&\phantom{\lesssim \int_{\mathcal{M}_\tau^{\tau_1} \cap \{r \geq \frac{1}{2}r_0\}}\bigg(}
	+ r^{1-C_{[N_2]}\epsilon} |\mathscr{Y}^{N_2} F|^2_{LL}
	+ \sum_{j \leq N_2 - 1} \epsilon^2 r^{1-\frac{1}{2}C_{[N_2]}\epsilon} |\mathscr{Y}^j F|_{(\text{frame})}^2
	\bigg) \dVol_{g}
	\end{split}
	\end{equation*}
	
\end{proof}

\begin{proposition}[$L^2$ bounds on $\mathscr{Y}^n \tr_{\slashed{g}}\chibar_{(\text{small})}$ below top-order]
	\label{proposition L2 tr chibar low}
	Suppose that all the bootstrap bounds hold, and also that $n \leq N_2 - 1$. Then
	\begin{equation*}
	\begin{split}
	&\int_{\mathcal{M}_{\tau}^{\tau_1} \cap \{r \geq r_0\}} C_{[n]}\epsilon r^{-1-C_{[n]}\epsilon} |\mathscr{Y}^n \tr_{\slashed{g}} \chibar_{(\text{small})}|^2 \dVol_g
	\\
	&\lesssim
	\left(\frac{\epsilon^2}{C_{[n-1]}} + \frac{C_{[n]}\epsilon^3}{\delta} \right)\epsilon^{2(N_2 + 1 - n)} (1+\tau)^{-1+C_{(n)}\delta}
	\\
	&\phantom{\lesssim}
	+ \int_{\mathcal{M}_\tau^{\tau_1} \cap \{r \geq r_0\}} \bigg(
		C_{[n]}\epsilon r^{-1-C_{[n]}\epsilon} |\slashed{\D} \mathscr{Y}^n h|^2_{(\text{frame})}
		+ C_{[n]}\epsilon r^{-1+\delta} |\overline{\slashed{\D}} \mathscr{Y}^{n+1} h|^2_{(\text{frame})}
		\\
		&\phantom{\lesssim + \int_{\mathcal{M}_\tau^{\tau_1} \cap \{r \geq r_0\}} \bigg(}
		+ C_{[n]}\epsilon r^{-1+\delta} |\overline{\slashed{\D}} \mathscr{Y}^{n} h|^2_{(\text{frame})}
		+ C_{[n]}\epsilon r^{-3-\delta} |\mathscr{Y}^n h|^2_{(\text{frame})}
	\bigg) \dVol_{g}
	\end{split}
	\end{equation*}
	and also, for any $c > 0$,
	\begin{equation*}
	\begin{split}
	&\int_{\mathcal{M}_{\tau}^{\tau_1} \cap \{r \geq r_0\}} c\delta r^{-1-c\delta} |\mathscr{Y}^n \tr_{\slashed{g}} \chibar_{(\text{small})}|^2 \dVol_g
	\\
	&\lesssim
	c \epsilon^2 \epsilon^{2(N_2 + 1 - n)} (1+\tau)^{-1+C_{(n)}\delta}
	\\
	&\phantom{\lesssim}
	+ \int_{\mathcal{M}_\tau^{\tau_1} \cap \{r \geq r_0\}} \bigg(
		c\delta r^{-1-c\delta} |\slashed{\D} \mathscr{Y}^n h|^2_{(\text{frame})}
		+ c\delta r^{-1+\frac{1}{2}\delta} |\overline{\slashed{\D}} \mathscr{Y}^{n+1} h|^2_{(\text{frame})}
		+ c\delta r^{-1+\delta} |\overline{\slashed{\D}} \mathscr{Y}^{n} h|^2_{(\text{frame})}
		\\
		&\phantom{\lesssim + \int_{\mathcal{M}_\tau^{\tau_1} \cap \{r \geq r_0\}} \bigg(}
		+ c\delta r^{-3-\delta} |\mathscr{Y}^n h|^2_{(\text{frame})}
		+ \epsilon^2 r^{-3-\frac{1}{2}c\delta} |\mathscr{Y}^n h|^2_{(\text{frame})}
	\bigg) \dVol_{g}
	\end{split}
	\end{equation*}
\end{proposition}

\begin{proof}
	Recall proposition \ref{proposition Yn tr chibar}. This gives
	\begin{equation*}
	\begin{split}
	&\int_{\mathcal{M}_{\tau}^{\tau_1} \cap \{r \geq r_0\}} r^{-1-C_{[n]}\epsilon} |\mathscr{Y}^n \tr_{\slashed{g}} \chibar_{(\text{small})}|^2 \dVol_g
	\\
	&\lesssim
	\int_{\mathcal{M}_\tau^{\tau_1} \cap \{r \geq r_0\}} \bigg(
		r^{-1-C_{[n]}\epsilon} |\mathscr{Y}^n \tr_{\slashed{g}} \chi_{(\text{small})}|^2
		+ r^{-1-C_{[n]}\epsilon} |\slashed{\D} \mathscr{Y}^n h|^2_{(\text{frame})}
		\\
		&\phantom{\lesssim \int_{\mathcal{M}_\tau^{\tau_1} \cap \{r \geq r_0\}} \bigg(}
		+ \epsilon^2 r^{-3-\frac{1}{2}C_{[n]}\epsilon} |\mathscr{Y}^n X_{(\text{frame, small})}|^2
		+ \epsilon^2 r^{-1-\frac{1}{2}C_{[n]}\epsilon} |\overline{\slashed{\D}} \mathscr{Y}^n h|^2_{(\text{frame})}
		\\
		&\phantom{\lesssim \int_{\mathcal{M}_\tau^{\tau_1} \cap \{r \geq r_0\}} \bigg(}
		+ \epsilon^2 r^{-1-2\delta} |\bm{\Gamma}^{(n)}_{(-1+C_{(n)}\epsilon)}|^2
		+ r^{-1-C_{[n]}\epsilon} |\bm{\Gamma}^{(n-1)}_{(-1-\delta)}|^2
	\bigg) \dVol_{g}
	\end{split}
	\end{equation*}
	where, as usual, we have followed remark \ref{remark L2 bounds involving Yn Xframe} to handle terms involving $\mathscr{Y}^n X_{(\text{frame})}$.
	
	Now, we can bound the norms of $\mathscr{Y}^n X_{(\text{frame, small})}$ and $\tr_{\slashed{g}}\chi_{(\text{small})}$ using propositions \ref{proposition L2 bounds rectangular} and \ref{proposition L2 tr chi low} respectively, to find
	\begin{equation*}
	\begin{split}
	&\int_{\mathcal{M}_{\tau}^{\tau_1} \cap \{r \geq r_0\}} r^{-1-C_{[n]}\epsilon} |\mathscr{Y}^n \tr_{\slashed{g}} \chibar_{(\text{small})}|^2 \dVol_g
	\\
	&\lesssim
	\left( \frac{\epsilon^2}{\delta (C_{[n]})^2} + \frac{\epsilon}{C_{[n]}C_{[n-1]}} + \frac{\epsilon^2}{\delta} \right)\epsilon^{2(N_2 + 1 - n)} (1+\tau)^{-1+C_{(n)}\delta}
	\\
	&\phantom{\lesssim}
	+ \int_{\mathcal{M}_\tau^{\tau_1} \cap \{r \geq r_0\}} \bigg(
		r^{-1-C_{[n]}\epsilon} |\slashed{\D} \mathscr{Y}^n h|^2_{(\text{frame})}
		+ r^{-1+\frac{1}{2}\delta} |\overline{\slashed{\D}} \mathscr{Y}^{n+1} h|^2_{(\text{frame})}
		+ r^{-1+\frac{1}{2}\delta} |\overline{\slashed{\D}} \mathscr{Y}^{n} h|^2_{(\text{frame})}
		\\
		&\phantom{\lesssim + \int_{\mathcal{M}_\tau^{\tau_1} \cap \{r \geq r_0\}} \bigg(}
		+ r^{-3-\delta} |\mathscr{Y}^n h|^2_{(\text{frame})}
		+ \epsilon^2 r^{-1-2\delta} |\bm{\Gamma}^{(n)}_{(-1+C_{(n)}\epsilon)}|^2
		+ r^{-1-C_{[n]}\epsilon} |\bm{\Gamma}^{(n-1)}_{(-1-\delta)}|^2	
	\bigg) \dVol_{g}
	\end{split}
	\end{equation*}
	
	Next, substituting the bootstrap bounds from equation \eqref{equation bootstrap L2 geometric top order} and using $C_{[n]} \gg 1$ leads to
	\begin{equation*}
	\begin{split}
	&\int_{\mathcal{M}_{\tau}^{\tau_1} \cap \{r \geq r_0\}} r^{-1-C_{[n]}\epsilon} |\mathscr{Y}^n \tr_{\slashed{g}} \chibar_{(\text{small})}|^2 \dVol_g
	\\
	&\lesssim
	\left(\frac{\epsilon}{C_{[n]}C_{[n-1]}} + \frac{\epsilon^2}{\delta} \right)\epsilon^{2(N_2 + 1 - n)} (1+\tau)^{-1+C_{(n)}\delta}
	\\
	&\phantom{\lesssim}
	+ \int_{\mathcal{M}_\tau^{\tau_1} \cap \{r \geq r_0\}} \bigg(
		r^{-1-C_{[n]}\epsilon} |\slashed{\D} \mathscr{Y}^n h|^2_{(\text{frame})}
		+ r^{-1+\delta} |\overline{\slashed{\D}} \mathscr{Y}^{n+1} h|^2_{(\text{frame})}
		+ r^{-1+\delta} |\overline{\slashed{\D}} \mathscr{Y}^{n} h|^2_{(\text{frame})}
		\\
		&\phantom{\lesssim + \int_{\mathcal{M}_\tau^{\tau_1} \cap \{r \geq r_0\}} \bigg(}
		+ r^{-3-\delta} |\mathscr{Y}^n h|^2_{(\text{frame})}
	\bigg) \dVol_{g}
	\end{split}
	\end{equation*}
	Proving the first part of the proposition.

	Next, following almost identical calculations, we have
	\begin{equation*}
	\begin{split}
	&\int_{\mathcal{M}_{\tau}^{\tau_1} \cap \{r \geq r_0\}} r^{-1-c\delta} |\mathscr{Y}^n \tr_{\slashed{g}} \chibar_{(\text{small})}|^2 \dVol_g
	\\
	&\lesssim
	\int_{\mathcal{M}_\tau^{\tau_1} \cap \{r \geq r_0\}} \bigg(
		r^{-1-c\delta} |\mathscr{Y}^n \tr_{\slashed{g}} \chi_{(\text{small})}|^2
		+ r^{-1-c\delta} |\slashed{\D} \mathscr{Y}^n h|^2_{(\text{frame})}
		+ \epsilon^2 r^{-3-\frac{1}{2}c\delta} |\mathscr{Y}^n X_{(\text{frame, small})}|^2
		\\
		&\phantom{\lesssim \int_{\mathcal{M}_\tau^{\tau_1} \cap \{r \geq r_0\}} \bigg(}
		+ \epsilon^2 r^{-1-\frac{1}{2}c\delta} |\overline{\slashed{\D}} \mathscr{Y}^n h|^2_{(\text{frame})}
		+ \epsilon^2 r^{-1-2\delta} |\bm{\Gamma}^{(n)}_{(-1+C_{(n)}\epsilon)}|^2
		+ r^{-1-c\delta} |\bm{\Gamma}^{(n-1)}_{(-1-\delta)}|^2
	\bigg) \dVol_{g}
	\end{split}
	\end{equation*}
	
	Again, we substitute the bounds from proposition \ref{proposition L2 tr chi low} for the term involving $\mathscr{Y}^n \tr_{\slashed{g}} \chi_{(\text{small})}$, and the bounds from proposition \ref{proposition L2 bounds rectangular} for the rectangular terms, and finally we use the bootstrap bounds for the other terms, to find
	\begin{equation*}
	\begin{split}
	&\int_{\mathcal{M}_{\tau}^{\tau_1} \cap \{r \geq r_0\}} r^{-1-c\delta} |\mathscr{Y}^n \tr_{\slashed{g}} \chibar_{(\text{small})}|^2 \dVol_g
	\\
	&\lesssim
	\frac{\epsilon^2}{\delta} \left(1 + \frac{\epsilon^2}{c^2 \delta^2} + \frac{\epsilon}{c^2 \delta C_{[n-1]}} \right) \epsilon^{2(N_2 + 1 - n)} (1+\tau)^{-1+C_{(n)}\delta}
	\\
	&\phantom{\lesssim}
	+ \int_{\mathcal{M}_\tau^{\tau_1} \cap \{r \geq r_0\}} \bigg(
		r^{-1-c\delta} |\slashed{\D} \mathscr{Y}^n h|^2_{(\text{frame})}
		+ r^{-1+\frac{1}{2}\delta} |\overline{\slashed{\D}} \mathscr{Y}^{n+1} h|^2_{(\text{frame})}
		\\
		&\phantom{\lesssim + \int_{\mathcal{M}_\tau^{\tau_1} \cap \{r \geq r_0\}} \bigg(}
		+ (1+c^{-2}\delta^{-2}\epsilon^2) r^{-1+\delta} |\overline{\slashed{\D}} \mathscr{Y}^{n} h|^2_{(\text{frame})}
		+ r^{-3-\delta} |\mathscr{Y}^n h|^2_{(\text{frame})}
		\\
		&\phantom{\lesssim + \int_{\mathcal{M}_\tau^{\tau_1} \cap \{r \geq r_0\}} \bigg(}
		+ \epsilon^2 c^{-1} \delta^{-1} r^{-3-\frac{1}{2}c\delta} |\mathscr{Y}^n h|^2_{(\text{frame})}
	\bigg) \dVol_{g}
	\end{split}
	\end{equation*}
	Finally, using $\epsilon \ll c\delta \ll 1$ proves the second part of the proposition.

\end{proof}

\begin{proposition}[$L^2$ bounds on $\mathscr{Y}^n \tr_{\slashed{g}}\chibar_{(\text{small})}$ at top-order]
	\label{proposition L2 tr chibar high}
	Suppose that all the bootstrap bounds hold, and also that $N_2 \leq 2N_1$. Also suppose that, for $n \leq N_1$ we have
	\begin{equation*}
	|\mathscr{Y}^n F|_{(\text{frame})} \lesssim \epsilon (1+r)^{-2+2C_{(n)}\epsilon}
	\end{equation*}
	Then
	\begin{equation*}
	\begin{split}
	&\int_{\mathcal{M}_{\tau}^{\tau_1} \cap \{r \geq r_0\}} C_{[N_2]} \epsilon r^{-1-C_{[N_2]}\epsilon} |\mathscr{Y}^n \tr_{\slashed{g}} \chibar_{(\text{small})}|^2 \dVol_g
	\\
	&\lesssim
	C_{[N_2]} \delta^{-1} \epsilon^5 (1+\tau)^{-1+C_{(n)}\delta}
	\\
	&\phantom{\lesssim}
	+ \int_{\mathcal{M}_\tau^{\tau_1} \cap \{r \geq r_0\}} \bigg(
		C_{[N_2]}\epsilon r^{-1-C_{[N_2]}\epsilon} |\slashed{\D} \mathscr{Y}^{N_2} h|^2_{LL}
		+ C_{[N_2]}\epsilon r^{-1-\delta} |\slashed{\D} \mathscr{Y}^{N_2} h|^2_{(\text{frame})}
		\\
		&\phantom{\lesssim + \int_{\mathcal{M}_\tau^{\tau_1} \cap \{r \geq r_0\}} \bigg(}
		+ C_{[N_2]}\epsilon \left(\epsilon^2 + (C_{[N_2]})^{-2} \right) r^{-1-\frac{1}{2}C_{[N_2]}\epsilon} |\overline{\slashed{\D}} \mathscr{Y}^{N_2} h|^2_{(\text{frame})}
		+ C_{[N_2]}\epsilon r^{-3-\delta} |\mathscr{Y}^n h|^2_{(\text{frame})}
		\\
		&\phantom{\lesssim + \int_{\mathcal{M}_\tau^{\tau_1} \cap \{r \geq r_0\}} \bigg(}
		+ C_{[N_2]}\epsilon r^{1-C_{[N_2]}\epsilon} |\mathscr{Y}^n F|^2_{LL}
		+ \sum_{j \leq N_2 - 1} C_{[N_2]} \epsilon^3 r^{1-\frac{1}{2}C_{[N_2]}\epsilon} |\mathscr{Y}^j F|_{(\text{frame})}^2
		\\
		&\phantom{\lesssim + \int_{\mathcal{M}_\tau^{\tau_1} \cap \{r \geq r_0\}} \bigg(}
		+ C_{[N_2]}\epsilon^3 r^{-1-2\delta} |\bm{\Gamma}^{(N_2)}_{(-1+C_{(n)}\epsilon)}|^2
		+ C_{[N_2]}\epsilon r^{-1-C_{[n]}\epsilon} |\bm{\Gamma}^{(N_2-1)}_{(-1-\delta)}|^2
	\bigg) \dVol_{g}
	\end{split}
	\end{equation*}
\end{proposition}

\begin{proof}
	As in the proof of proposition \ref{proposition L2 tr chibar low}, we have
	\begin{equation*}
	\begin{split}
	&\int_{\mathcal{M}_{\tau}^{\tau_1} \cap \{r \geq r_0\}} r^{-1-C_{[N_2]}\epsilon} |\mathscr{Y}^n \tr_{\slashed{g}} \chibar_{(\text{small})}|^2 \dVol_g
	\\
	&\lesssim
	\int_{\mathcal{M}_\tau^{\tau_1} \cap \{r \geq r_0\}} \bigg(
		r^{-1-C_{[N_2]}\epsilon} |\mathscr{Y}^n \tr_{\slashed{g}} \chi_{(\text{small})}|^2
		+ r^{-1-C_{[N_2]}\epsilon} |\slashed{\D} \mathscr{Y}^n h|^2_{(\text{frame})}
		\\
		&\phantom{\lesssim \int_{\mathcal{M}_\tau^{\tau_1} \cap \{r \geq r_0\}} \bigg(}
		+ \epsilon^2 r^{-3-\frac{1}{2}C_{[N_2]}\epsilon} |\mathscr{Y}^n X_{(\text{frame, small})}|^2
		+ \epsilon^2 r^{-1-\frac{1}{2}C_{[N_2]}\epsilon} |\overline{\slashed{\D}} \mathscr{Y}^n h|^2_{(\text{frame})}
		\\
		&\phantom{\lesssim \int_{\mathcal{M}_\tau^{\tau_1} \cap \{r \geq r_0\}} \bigg(}
		+ \epsilon^2 r^{-1-2\delta} |\bm{\Gamma}^{(N_2)}_{(-1+C_{(n)}\epsilon)}|^2
		+ r^{-1-C_{[n]}\epsilon} |\bm{\Gamma}^{(N_2-1)}_{(-1-\delta)}|^2
	\bigg) \dVol_{g}
	\end{split}
	\end{equation*}
	
	As before, we can bound $\mathscr{Y}^{N_2} X_{(\text{frame, small})}$ using proposition \ref{proposition L2 bounds rectangular}, but this time we bound the connection coefficient $\mathscr{Y}^n \tr_{\slashed{g}}\chi_{(\text{small})}$ using proposition \ref{proposition L2 tr chi high} instead of proposition \ref{proposition L2 tr chi low}. This gives
	\begin{equation*}
	\begin{split}
	&\int_{\mathcal{M}_{\tau}^{\tau_1} \cap \{r \geq r_0\}} r^{-1-C_{[N_2]}\epsilon} |\mathscr{Y}^n \tr_{\slashed{g}} \chibar_{(\text{small})}|^2 \dVol_g
	\\
	&\lesssim
	\epsilon^4 (1+\tau)^{-1+C_{(n)}\delta}
	\\
	&\phantom{\lesssim}
	+ \int_{\mathcal{M}_\tau^{\tau_1} \cap \{r \geq r_0\}} \bigg(
		r^{-1-C_{[N_2]}\epsilon} |\slashed{\D} \mathscr{Y}^{N_2} h|^2_{LL}
		+ r^{-1-\delta} |\slashed{\D} \mathscr{Y}^{N_2} h|^2_{(\text{frame})}
		\\
		&\phantom{\lesssim + \int_{\mathcal{M}_\tau^{\tau_1} \cap \{r \geq r_0\}} \bigg(}
		+ \left(\epsilon^2 + (C_{[N_2]})^{-2} \right) r^{-1-\frac{1}{2}C_{[N_2]}\epsilon} |\overline{\slashed{\D}} \mathscr{Y}^{N_2} h|^2_{(\text{frame})}
		+ r^{-3-\delta} |\mathscr{Y}^n h|^2_{(\text{frame})}
		\\
		&\phantom{\lesssim + \int_{\mathcal{M}_\tau^{\tau_1} \cap \{r \geq r_0\}} \bigg(}
		+ r^{1-C_{[N_2]}\epsilon} |\mathscr{Y}^n F|^2_{LL}
		+ \sum_{j \leq N_2 - 1} \epsilon^2 r^{1-\frac{1}{2}C_{[N_2]}\epsilon} |\mathscr{Y}^j F|_{(\text{frame})}^2
		\\
		&\phantom{\lesssim + \int_{\mathcal{M}_\tau^{\tau_1} \cap \{r \geq r_0\}} \bigg(}
		+ \epsilon^2 r^{-1-2\delta} |\bm{\Gamma}^{(N_2)}_{(-1+C_{(n)}\epsilon)}|^2
		+ r^{-1-C_{[n]}\epsilon} |\bm{\Gamma}^{(N_2-1)}_{(-1-\delta)}|^2
	\bigg) \dVol_{g}
	\end{split}
	\end{equation*}
	
	Finally, using the bootstrap bounds from equations \eqref{equation bootstrap L2 geometric} and \eqref{equation bootstrap L2 geometric top order} proves the proposition.

\end{proof}

\begin{proposition}[$L^2$ bounds on the foliation density]
	\label{proposition L2 foliation density}
	Assuming that all the bootstrap bounds hold, the foliation density satisfies the following $L^2$ estimate: for all $n \leq N_2$
	\begin{equation*}
	\begin{split}
	&\int_{\mathcal{M}_\tau^{\tau_1} \cap \{r \geq r_0\}} C_{[n]}\epsilon r^{-3-C_{[n]}\epsilon} |\mathscr{Y}^n \log\mu|^2 \dVol_{g}
	\\
	&\lesssim
	\frac{\epsilon}{(C_{[n]}) \delta} \left( 1 + \frac{1}{\delta C_{[n-1]}} \right) \epsilon^{2(N_2 + 1 - n)} (1+\tau)^{-1+C_{(n)}\delta}
	\\
	&\phantom{\lesssim}
	+ \frac{1}{(C_{[n]})^2 \epsilon^2} \int_{\mathcal{M}_\tau^{\tau_1} \cap \{r \geq r_0\}} \bigg(
	C_{[n]}\epsilon r^{-1-C_{[n]}\epsilon} |\slashed{\D}\mathscr{Y}^n h|_{(\text{frame})}^2
	+ C_{[n]}\epsilon r^{-1+\frac{1}{2}\delta} |\overline{\slashed{\D}}\mathscr{Y}^n h|_{(\text{frame})}^2
	\\
	&\phantom{\lesssim + \frac{1}{(C_{[n]})^2 \epsilon^2} \int_{\mathcal{M}_\tau^{\tau_1} \cap \{r \geq r_0\}} \bigg(}
	+ C_{[n]}\epsilon r^{-3-\delta} |\mathscr{Y}^n h|_{(\text{frame})}^2
	\bigg)\dVol_{g}
	\end{split}
	\end{equation*}
	and also
	\begin{equation*}
	\begin{split}
	&\int_{\mathcal{M}_\tau^{\tau_1} \cap \{r \geq r_0\}} c\delta r^{-3-c\delta} |\mathscr{Y}^n \log\mu|^2 \dVol_{g}
	\\
	&\lesssim
	\epsilon c^{-3} \delta^{-3} \left( \frac{\epsilon}{\delta} + \frac{1}{C_{[n-1]}} \right) \epsilon^{2(N_2 + 1 - n)} (1+\tau)^{-1+C_{(n)}\delta}
	\\
	&\phantom{\lesssim}
	+ \frac{1}{c\delta} \int_{\mathcal{M}_\tau^{\tau_1} \cap \{r \geq r_0\}} \bigg(
	c\delta r^{1-c\delta} |\slashed{\D}\mathscr{Y}^n h|_{(\text{frame})}^2
	+ c^{-3} \delta^{-3} r^{-1+\frac{1}{2}\delta} |\overline{\slashed{\D}}\mathscr{Y}^n h|_{(\text{frame})}^2
	\\
	&\phantom{\lesssim + \frac{1}{c\delta} \int_{\mathcal{M}_\tau^{\tau_1} \cap \{r \geq r_0\}} \bigg(}
	+ c^{-3} \delta^{-3} r^{-3-\delta} |\mathscr{Y}^n h|_{(\text{frame})}^2
	\bigg) \dVol_{g}
	\end{split}
	\end{equation*}

\end{proposition}

\begin{proof}
	Using the third part of proposition \ref{proposition Hardy} we have
	\begin{equation*}
	\begin{split}
	\int_{\Sigma_\tau \cap \{r \geq r_0\}} r^{-1-C_{[n]}\epsilon} |\mathscr{Y}^n \log\mu|^2 \upd r \wedge \dVol_{\mathbb{S}^2}
	&\lesssim
	\frac{1}{C_{[n]}^2 \epsilon^2} \int_{\Sigma_\tau \cap \{r \geq r_0\}} r^{1-C_{[n]}\epsilon} |\slashed{\D}_L \left( \mathscr{Y}^n \log\mu \right)|^2 \upd r \wedge \dVol_{\mathbb{S}^2}
	\\
	&\phantom{\lesssim}
	+ \frac{1}{C_{[n]} \epsilon} \int_{S_{\tau,r_0}}|\mathscr{Y}^n \log\mu|^2 \upd r \wedge \dVol_{\mathbb{S}^2}
	\end{split}
	\end{equation*}
	So, using proposition \ref{proposition transport Yn mu} we have
	\begin{equation*}
	\begin{split}
	&\int_{\Sigma_\tau \cap \{r \geq r_0\}} r^{-1-C_{[n]}\epsilon} |\mathscr{Y}^n \log\mu|^2 \upd r \wedge \dVol_{\mathbb{S}^2}
	\\
	&\lesssim
	\frac{1}{(C_{[n]})^2 \epsilon^2} \int_{\Sigma_\tau \cap \{r \geq r_0\}} \bigg(
		\epsilon^2 r^{-1-C_{[n]}\epsilon} |\mathscr{Y}^n \log \mu|^2
		+ r^{-1-C_{[n]}\epsilon} |\mathscr{Y}^n \bar{X}_{(\text{frame})}|^2
		\\
		&\phantom{\lesssim \frac{1}{(C_{[n]})^2 \epsilon^2} \int_{\Sigma_\tau \cap \{r \geq r_0\}} \bigg(}
		+ \epsilon^2 r^{-1-2\delta} |\mathscr{Y}^n X_{(\text{frame, small})}|^2
		+ r^{1-C_{[n]}\epsilon} |\slashed{\D}\mathscr{Y}^n h|_{(\text{frame})}^2
		\\
		&\phantom{\lesssim \frac{1}{(C_{[n]})^2 \epsilon^2} \int_{\Sigma_\tau \cap \{r \geq r_0\}} \bigg(}
		+ r^{1-\frac{1}{2}C_{[n]}\epsilon} |\bm{\Gamma}^{(n-1)}_{(-1+C_{(n-1)}\epsilon)}|^2
		\bigg) \upd r \wedge \dVol_{\mathbb{S}^2}
	\\
	&\phantom{\lesssim}
	+ \frac{1}{C_{[n]} \epsilon} \int_{S_{\tau,r_0}}|\mathscr{Y}^n \log\mu|^2 \dVol_{\mathbb{S}^2}
	\end{split}
	\end{equation*}
	The first term on the right hand side can be absorbed on the left hand side, for sufficiently large $C_{[n]}$. 
	
	Next, we integrate over $\tau$ and use the fact that $\Omega \sim r$. If we then use proposition \ref{proposition L2 bounds rectangular} to bound the second and third terms on the right hand side we find
	\begin{equation*}
	\begin{split}
	&\int_{\mathcal{M}_\tau^{\tau_1} \cap \{r \geq r_0\}} r^{-3-C_{[n]}\epsilon} |\mathscr{Y}^n \log\mu|^2 \dVol_{g}
	\\
	&\lesssim
	\frac{1}{(C_{[n]})^2 \delta} \left( 1 + \frac{\epsilon^2}{\delta^2} + \frac{1}{\delta C_{[n-1]}} \right) \epsilon^{2(N_2 + 1 - n)}(1+\tau)^{-1+C_{(n)}\delta}
	+ \frac{1}{C_{[n]} \epsilon} \int_{S_{\tau,r_0}}|\mathscr{Y}^n \log\mu|^2 \upd r \wedge \dVol_{\mathbb{S}^2}
	\\
	&\phantom{\lesssim}
	+ \frac{1}{(C_{[n]})^2 \epsilon^2} \int_{\mathcal{M}_\tau^{\tau_1} \cap \{r \geq r_0\}} \bigg(
		r^{-1-C_{[n]}\epsilon} |\slashed{\D}\mathscr{Y}^n h|_{(\text{frame})}^2
		+ r^{-1+\frac{1}{2}\delta} |\overline{\slashed{\D}}\mathscr{Y}^n h|_{(\text{frame})}^2
		+ r^{-3-\delta} |\mathscr{Y}^n h|_{(\text{frame})}^2
		\\
		&\phantom{\lesssim + \frac{1}{(C_{[n]})^2 \epsilon^2} \int_{\mathcal{M}_\tau^{\tau_1} \cap \{r \geq r_0\}} \bigg(}
		+ r^{1-\frac{1}{2}C_{[n]}\epsilon} |\bm{\Gamma}^{(n-1)}_{(-1+C_{(n-1)}\epsilon)}|^2
		\bigg)\dVol_{g}
	\end{split}
	\end{equation*}	
	We can use the $L^2$ bootstrap bounds of equation \eqref{equation bootstrap L2 geometric} to handle the last term in the integral over $\Sigma_\tau$.
	
	Note that $\log \mu = 0$ on the surface $r = r_0$. Hence $T \log \mu = r\slashed{\nabla} \log \mu = 0$ here. On the other hand, $rL\log\mu = r\omega$, which can itself be estimated in terms of $h_{(\text{frame})}$ and $(\mathscr{Y} h)_{(\text{frame})}$ (or $(\partial h_{(\text{frame})})$) at $r = r_0$, so we can write
	\begin{equation*}
	\mathscr{Y}^n \log \mu \big|_{r=r_0} = \sum_{j \leq n} \mathscr{Y}^j h_{(\text{rect})}
	\end{equation*}
	and so,  controlling the term on the cylinder $r = r_0$ as in proposition \ref{proposition L2 bounds rectangular} and using $\epsilon \ll \delta$, we have
	\begin{equation*}
	\begin{split}
	&\int_{\mathcal{M}_\tau^{\tau_1} \cap \{r \geq r_0\}} r^{-3-C_{[n]}\epsilon} |\mathscr{Y}^n \log\mu|^2 \dVol_{g}
	\\
	&\lesssim
	\frac{1}{(C_{[n]})^2 \delta} \left( 1 + \frac{1}{\delta C_{[n-1]}} \right) \epsilon^{2(N_2 + 1 - n)} (1+\tau)^{-1+C_{(n)}\delta}
	\\
	&\phantom{\lesssim}
	+ \frac{1}{(C_{[n]})^2 \epsilon^2} \int_{\mathcal{M}_\tau^{\tau_1} \cap \{r \geq r_0\}} \bigg(
	r^{-1-C_{[n]}\epsilon} |\slashed{\D}\mathscr{Y}^n h|_{(\text{frame})}^2
	+ r^{-1+\frac{1}{2}\delta} |\overline{\slashed{\D}}\mathscr{Y}^n h|_{(\text{frame})}^2
	+ r^{-3-\delta} |\mathscr{Y}^n h|_{(\text{frame})}^2
	\\
	&\phantom{\lesssim + \frac{1}{(C_{[n]})^2 \epsilon^2} \int_{\mathcal{M}_\tau^{\tau_1} \cap \{r \geq r_0\}} \bigg(}
	+ r^{1-\frac{1}{2}C_{[n]}\epsilon} |\bm{\Gamma}^{(n-1)}_{(-1+C_{(n-1)}\epsilon)}|^2
	\bigg)\dVol_{g}
	\end{split}
	\end{equation*}
	
	Finally, using the fact that $\epsilon \ll \delta$ and the $L^2$ bootstrap bounds proves the first part of the proposition.
	
	Following very similar calculations, we have
	\begin{equation*}
	\begin{split}
	\int_{\Sigma_\tau \cap \{r \geq r_0\}} r^{-1-c\delta} |\mathscr{Y}^n \log\mu|^2 \upd r \wedge \dVol_{\mathbb{S}^2}
	&\lesssim
	\frac{1}{c^2 \delta^2} \int_{\Sigma_\tau \cap \{r \geq r_0\}} r^{1-c\delta} |\slashed{\D}_L \left( \mathscr{Y}^n \log\mu \right)|^2 \upd r \wedge \dVol_{\mathbb{S}^2}
	\\
	&\phantom{\lesssim}
	+ \frac{1}{c\delta} \int_{S_{\tau,r_0}}|\mathscr{Y}^n \log\mu|^2 \upd r \wedge \dVol_{\mathbb{S}^2}
	\end{split}
	\end{equation*}
	and so
	\begin{equation*}
	\begin{split}
	&\int_{\mathcal{M}_\tau^{\tau_1} \cap \{r \geq r_0\}} r^{-3-c\delta} |\mathscr{Y}^n \log\mu|^2 \dVol_{g}
	\\
	&\lesssim
	\frac{1}{c^2\delta^2} \int_{\mathcal{M}_\tau^{\tau_1} \cap \{r \geq r_0\}} \bigg(
		\epsilon^2 r^{-3-c\delta} |\mathscr{Y}^n \log \mu|^2
		+ r^{-3-c\delta} |\mathscr{Y}^n \bar{X}_{(\text{frame})}|^2
		+ \epsilon^2 r^{-3-2\delta} |\mathscr{Y}^n X_{(\text{frame, small})}|^2
		\\
		&\phantom{\lesssim \frac{1}{C_{[n]}^2 \epsilon^2} \int_{\Sigma_\tau \cap \{r \geq r_0\}} \bigg(}
		+ r^{-1-c\delta} |\slashed{\D}\mathscr{Y}^n h|_{(\text{frame})}^2
		+ r^{-1-\frac{1}{2}c\delta} |\bm{\Gamma}^{(n-1)}_{(-1+C_{(n-1)}\epsilon)}|^2
	\bigg) \dVol_{g}
	\\
	&\phantom{\lesssim}
	+ \frac{1}{c\delta} \int_{\tau}^{\tau_1} \left( \int_{S_{\tau',r_0}}|\mathscr{Y}^n \log\mu|^2 \dVol_{\mathbb{S}^2} \right) \upd \tau'
	\end{split}
	\end{equation*}
	Absorbing the first term by the left hand side, bounding the final term as before, and using the bootstrap bounds and proposition \ref{proposition L2 bounds rectangular} for the other terms, we have
	\begin{equation*}
	\begin{split}
	&\int_{\mathcal{M}_\tau^{\tau_1} \cap \{r \geq r_0\}} r^{-3-c\delta} |\mathscr{Y}^n \log\mu|^2 \dVol_{g}
	\\
	&\lesssim
	\epsilon c^{-4} \delta^{-4} \left( \frac{\epsilon}{\delta} + \frac{1}{C_{[n-1]}} \right) \epsilon^{2(N_2 + 1 - n)} (1+\tau)^{-1+C_{(n)}\delta}
	\\
	&\phantom{\lesssim}
	+ \frac{1}{c^2\delta^2} \int_{\mathcal{M}_\tau^{\tau_1} \cap \{r \geq r_0\}} \bigg(
		r^{1-c\delta} |\slashed{\D}\mathscr{Y}^n h|_{(\text{frame})}^2
		+ c^{-4} \delta^{-4} r^{-1+\frac{1}{2}\delta} |\overline{\slashed{\D}}\mathscr{Y}^n h|_{(\text{frame})}^2
		\\
		&\phantom{\lesssim + \frac{1}{c^2\delta^2} \int_{\mathcal{M}_\tau^{\tau_1} \cap \{r \geq r_0\}} \bigg(}
		+ c^{-4} \delta^{-4} r^{-3-\delta} |\mathscr{Y}^n h|_{(\text{frame})}^2
	\bigg) \dVol_{g}
	\end{split}
	\end{equation*}

\end{proof}

\begin{remark}
	As in  the case of the $L^2$ bounds on the frame fields, the bounds established above for the foliation density hold in the case $n = N_2$, as well as for lower orders.
\end{remark}

We can use proposition \ref{proposition L2 foliation density} to control up to $N_2$ derivatives of the foliation density. Note, however, that we cannot use this to bound $\slashed{\nabla}^2 \mathscr{Y}^{N_2-1} \log \mu$ -- a na\"ive estimate would be to estimate this in terms of $\mathscr{Y}^{N_2+1} \log \mu$, but using proposition \ref{proposition L2 foliation density} to bound this would involve $N_2+1$ derivatives of $h$, which we do not control. Instead, we use the following proposition.

\begin{proposition}[$L^2$ bounds on $\slashed{\nabla}^2 \log \mu$ at top order]
	\label{proposition L2 nabla log mu}
	Suppose that the bootstrap bounds hold. Suppose also that the inhomogeneous terms $F_{(\text{frame})}$ satisfy the pointwise bounds
	\begin{equation*}
	|\mathscr{Y}^n F|_{(\text{frame})} \lesssim \epsilon^2 (1+r)^{-2 + 2C_{(n)}\epsilon}
	\end{equation*}
	for all $n \leq N_1$. Also, suppose that $N_2 \leq 2N_1$. Finally, let $C_{[N_2(j)]} > 0$ be any positive constant.
	
	Then we have
	\begin{equation*}
	\begin{split}
	&\int_{\mathcal{M}_\tau^{\tau_1} \cap \{r \geq r_0\}} C_{[N_2(j)]}\epsilon r^{-1-C_{[N_2(j)]}\epsilon} |\slashed{\nabla}^2 \mathscr{Y}^{N_2-1} \log \mu|^2 \dVol_{g}
	\\
	&\lesssim
	C_{[N_2(j)]}\delta^{-1} \epsilon^5 (1+\tau)^{-1+C_{(N_2)}\delta}
	\\
	&\phantom{\lesssim}
	+ \int_{\mathcal{M}_\tau^{\tau_1} \cap \{r \geq r_0\}} \bigg(
		C_{[N_2(j)]}\epsilon r^{-1-(C_{[N_2(j)]} - 2C_{(0)})\epsilon} |\slashed{\D} \slashed{\D}_T \mathscr{Y}^{N_2-1} h|^2_{LL}
		+ C_{[N_2(j)]}\epsilon r^{-1-\delta} |\slashed{\D} \mathscr{Y}^{N_2} h|^2_{(\text{frame})}
		\\
		&\phantom{\lesssim + \int_{\mathcal{M}_\tau^{\tau_1} \cap \{r \geq r_0\}} \bigg(}
		+ C_{[N_2(j)]}\epsilon^3 r^{-1-\frac{1}{2}C_{[N_2(j)]}\epsilon} |\overline{\slashed{\D}} \mathscr{Y}^{N_2} h|^2_{(\text{frame})}
		+ C_{[N_2(j)]}\epsilon^3 \delta^{-2} r^{-1+\frac{1}{2}\delta} |\overline{\slashed{\D}} \mathscr{Y}^{N_2} h|^2_{(\text{frame})}
		\\
		&\phantom{\lesssim + \int_{\mathcal{M}_\tau^{\tau_1} \cap \{r \geq r_0\}} \bigg(}
		+ C_{[N_2(j)]}\epsilon r^{-3-\delta} |\mathscr{Y}^{N_2} h|^2_{(\text{frame})}
		+ C_{[N_2(j)]}\epsilon r^{1-C_{[N_2(j)]}\epsilon} |\slashed{\D}_T \mathscr{Y}^{N_2-1} F|^2_{LL}
		\\
		&\phantom{\lesssim + \int_{\mathcal{M}_\tau^{\tau_1} \cap \{r \geq r_0\}} \bigg(}
		+ \sum_{j \leq N_2 - 1} C_{[N_2(j)]}\epsilon^3 r^{1-\frac{1}{2}C_{[N_2(j)]}\epsilon} |\mathscr{Y}^j F|_{(\text{frame})}^2
		+ \sum_{k \leq N_2} C_{[N_2(j)]}\epsilon r^{-1} |\mathscr{Y}^k F|_{(\text{frame})}
	\bigg) \dVol_{g}
	\end{split}
	\end{equation*}
	
\end{proposition}

\begin{proof}	
	Recall proposition \ref{proposition commute laplacian mu}, which gives
	\begin{equation}
	\begin{split}
	\slashed{\Delta} \mathscr{Z}^{n-1} \log \mu
	&=
	\bm{\Gamma}^{(0)}_{(-1-\delta)} \mathscr{Y}^n \log \mu
	+ (1+\bm{\Gamma}^{(0)}_{(C_{(0)}\epsilon)}) \slashed{\D}_T \mathscr{Y}^{n-1} \tr_{\slashed{g}} \chi_{(\text{small})}
	+ r^{-1} \mathscr{Y}^n \zeta
	+ r^{-1} (\slashed{\D} \mathscr{Y}^n h)_{(\text{frame})}
	\\
	&\phantom{=}
	+ \sum_{j+k \leq n-1} (\mathscr{Y}^{j} F)_{(\text{rect})} \bm{\Gamma}^{(k)}_{(C_{(k)}\epsilon, \text{large})}
	+ \sum_{j+k \leq n} \bm{\Gamma}^{(j)}_{(-1+C_{(j)}\epsilon)}\bm{\Gamma}^{(k)}_{(-1+C_{(k)}\epsilon)}
	\end{split}
	\end{equation}
	
	Now, proposition \ref{proposition Poisson estimate} gives us
	\begin{equation*}
	\begin{split}
	&\int_{S_{\tau_r}} |\slashed{\nabla}^2 \mathscr{Z}^{n-1} \log \mu|^2 \dVol_{\mathbb{S}^2}
	\\
	&\lesssim
	\int_{S_{\tau_r}} \left( 
	|\slashed{\Delta} \mathscr{Z}^{n-1} \log \mu|^2
	+ K |\slashed{\nabla} \mathscr{Z}^{n-1} \log \mu|^2
	+ |\slashed{\nabla} K| |\slashed{\nabla} \mathscr{Z}^{n-1} \log \mu| | \mathscr{Z}^{n-1} \log \mu|
	\right) \dVol_{\mathbb{S}^2}
	\end{split}
	\end{equation*}
	where we recall that $K$ is the Gauss curvature of the sphere $S_{\tau,r}$ with respect to the metric $\slashed{g}$. Using the pointwise bounds on up to three derivatives of the metric, we find that
	\begin{equation*}
	\begin{split}
	&\int_{S_{\tau,r}} |\slashed{\nabla}^2 \mathscr{Z}^{n-1} \log \mu|^2 \dVol_{\mathbb{S}^2}
	\\
	&\lesssim
	\int_{S_{\tau,r}} \left( 
	|\slashed{\Delta} \mathscr{Z}^{n-1} \log \mu|^2
	+ r^{-4} |\mathscr{Y}^{n} \log \mu|^2
	+ \epsilon r^{-3} |\slashed{\nabla} \mathscr{Z}^{n-1} \log \mu| | \mathscr{Z}^{n-1} \log \mu|
	\right) \dVol_{\mathbb{S}^2}
	\\
	&\lesssim
	\int_{S_{\tau,r}} \left( 
	|\slashed{\Delta} \mathscr{Z}^{n-1} \log \mu|^2
	+ r^{-4} |\mathscr{Y}^{n} \log \mu|^2
	+ \epsilon r^{-4} | \mathscr{Z}^{n-1} \log \mu|
	\right) \dVol_{\mathbb{S}^2}
	\end{split}
	\end{equation*}
	and so, using remark \ref{remark L2 bounds involving Yn Xframe} to control the derivatives of the frame fields, we find
	\begin{equation*}
	\begin{split}
	&\int_{\mathcal{M}_\tau^{\tau_1} \cap \{r \geq r_0\}} r^{-1-\frac{1}{2}C_{[N_2]}\epsilon} |\slashed{\nabla}^2 \mathscr{Y}^{N_2-1} \log \mu|^2 \dVol_{g}
	\\
	&\lesssim
	\int_{\mathcal{M}_\tau^{\tau_1} \cap \{r \geq r_0\}} r^{-1-C_{[N_2(j)]}\epsilon} \bigg(
	\epsilon^2 r^{-2-2\delta} |\mathscr{Y}^{N_2} \log \mu|^2
	+ \left( 1+\epsilon^2 r^{2C_{(0)}\epsilon} \right) |\slashed{\D}_T \mathscr{Y}^{N_2-1} \tr_{\slashed{g}} \chi_{(\text{small})}|^2
	\\
	&\phantom{\lesssim \int_{\mathcal{M}_\tau^{\tau_1} \cap \{r \geq r_0\}} r^{-1-C_{[N_2(j)]}\epsilon} \bigg(}
	+ r^{-2} |\mathscr{Y}^{N_2} \zeta|^2
	+ r^{-2} |\slashed{\D} \mathscr{Y}^{N_2} h|^2_{(\text{frame})}
	\\
	&\phantom{\lesssim \int_{\mathcal{M}_\tau^{\tau_1} \cap \{r \geq r_0\}} r^{-1-\frac{1}{2}C_{[N_2]}\epsilon} \bigg(}
	+ r^2 \sum_{j+k \leq N_2-1} |\mathscr{Y}^{j} F|^2_{(\text{frame})} |\bm{\Gamma}^{(k)}_{(-1 + C_{(k)}\epsilon)}|^2
	+ \sum_{j \leq N_2 - 1} |\mathscr{Y}^j F|^2_{(\text{frame})}
	\\
	&\phantom{\lesssim \int_{\mathcal{M}_\tau^{\tau_1} \cap \{r \geq r_0\}} r^{-1-\frac{1}{2}C_{[N_2]}\epsilon} \bigg(}
	+ \sum_{k+\ell \leq N_2} |\bm{\Gamma}^{(k)}_{(-1+C_{(k)}\epsilon)}|^2 |\bm{\Gamma}^{(\ell)}_{(-1+C_{(\ell)}\epsilon)}|^2
	\bigg) \dVol_{g}
	\end{split}
	\end{equation*}
	Now, using the fact that $N_2 \leq 2N_1$, we can use pointwise bootstrap bounds for $\bm{\Gamma}^{(\leq N_1)}_{(-1 + C_{(\leq N_1)}\epsilon)}$ along with the bounds assumed in the proposition for $F$ to bound the quadratic terms, to find
	\begin{equation*}
	\begin{split}
	&\int_{\mathcal{M}_\tau^{\tau_1} \cap \{r \geq r_0\}} r^{-1-C_{[N_2(j)]}\epsilon} |\slashed{\nabla}^2 \mathscr{Y}^{N_2-1} \log \mu|^2 \dVol_{g}
	\\
	&\lesssim
	\int_{\mathcal{M}_\tau^{\tau_1} \cap \{r \geq r_0\}} \bigg(
	r^{-1-\left(C_{[N_2(j)]} - 2C_{(0)}\right)\epsilon} |\slashed{\D}_T \mathscr{Y}^{N_2-1} \tr_{\slashed{g}} \chi_{(\text{small})}|^2
	+ \epsilon^2 r^{-3} |\Gamma^{(N_2)}_{(-1+C_{(N_2)}\epsilon)}|^2
	\\
	& \phantom{\lesssim \int_{\mathcal{M}_\tau^{\tau_1} \cap \{r \geq r_0\}} \bigg( }
	+ r^{-3} |\slashed{\D} \mathscr{Y}^{N_2} h|^2
	+ \sum_{j \leq N_2} r^{-1} |\mathscr{Y}^j F|_{(\text{frame})}
	+ r^{-3} |\mathscr{Y}^{N_2} \zeta|^2
	\bigg) \dVol_{g}
	\end{split}
	\end{equation*}

	Next, we use proposition \ref{proposition L2 zeta high} to bound the term involving $\mathscr{Y}^{N_2} \zeta$, and the bootstrap bounds from equation \eqref{equation bootstrap L2 geometric} to obtain
	\begin{equation*}
	\begin{split}
	&\int_{\mathcal{M}_\tau^{\tau_1} \cap \{r \geq r_0\}} r^{-1-C_{[N_2(j)]}\epsilon} |\slashed{\nabla}^2 \mathscr{Y}^{N_2-1} \log \mu|^2 \dVol_{g}
	\\
	&\lesssim
	\delta^{-1} \epsilon^4 (1+\tau)^{-1+C_{(N_2)}\delta}
	\\
	&\phantom{\lesssim}
	+ \int_{\mathcal{M}_\tau^{\tau_1} \cap \{r \geq r_0\}} \bigg(
		r^{-1-\left( C_{[N_2(j)]} - 2C_{(0)}\right)\epsilon} |\slashed{\D}_T \mathscr{Y}^{N_2-1} \tr_{\slashed{g}} \chi_{(\text{small})}|^2
		+ \epsilon^2 \delta^{-2} r^{-1+\frac{1}{2}\delta} |\overline{\slashed{\D}} \mathscr{Y}^{N_2} h|^2_{(\text{frame})}
		\\
		& \phantom{\lesssim + \int_{\mathcal{M}_\tau^{\tau_1} \cap \{r \geq r_0\}} \bigg( }
		+ r^{-1-\delta} |\slashed{\D} \mathscr{Y}^{N_2} h|^2_{(\text{frame})}
		+ \epsilon^2 \delta^{-1} r^{-3-\delta} |\mathscr{Y}^{N_2} h|^2_{(\text{frame})}
		+ \sum_{j \leq N_2} r^{-1} |\mathscr{Y}^j F|_{(\text{frame})}
	\bigg) \dVol_{g}
	\end{split}
	\end{equation*}

	Finally, we use proposition \ref{proposition L2 tr chi high} to find
	\begin{equation*}
	\begin{split}
	&\int_{\mathcal{M}_\tau^{\tau_1} \cap \{r \geq r_0\}} r^{-1-C_{[N_2]}\epsilon} |\slashed{\nabla}^2 \mathscr{Y}^{N_2-1} \log \mu|^2 \dVol_{g}
	\\
	&\lesssim
	(1+\delta^{-1}) \epsilon^4 (1+\tau)^{-1+C_{(N_2)}\delta}
	\\
	&\phantom{\lesssim}
	+ \int_{\mathcal{M}_\tau^{\tau_1} \cap \{r \geq r_0\}} \bigg(
		r^{-1-(C_{[N_2(j)]} - 2C_{(0)})\epsilon} |\slashed{\D} \slashed{\D}_T \mathscr{Y}^{N_2-1} h|^2_{LL}
		+ r^{-1-\delta} |\slashed{\D} \mathscr{Y}^{N_2} h|^2_{(\text{frame})}
		\\
		&\phantom{\lesssim + \int_{\mathcal{M}_\tau^{\tau_1} \cap \{r \geq r_0\}} \bigg(}
		+ \epsilon^2 r^{-1-\frac{1}{2}C_{[N_2(j)]}\epsilon} |\overline{\slashed{\D}} \mathscr{Y}^{N_2} h|^2_{(\text{frame})}
		+ \epsilon^2 \delta^{-2} r^{-1+\frac{1}{2}\delta} |\overline{\slashed{\D}} \mathscr{Y}^{N_2} h|^2_{(\text{frame})}
		\\
		&\phantom{\lesssim + \int_{\mathcal{M}_\tau^{\tau_1} \cap \{r \geq r_0\}} \bigg(}
		+ (1+\epsilon^2 \delta^{-1}) r^{-3-\delta} |\mathscr{Y}^{N_2} h|^2_{(\text{frame})}
		+ r^{1-C_{[N_2(j)]}\epsilon} |\slashed{\D}_T \mathscr{Y}^{N_2-1} F|^2_{LL}
		\\
		&\phantom{\lesssim + \int_{\mathcal{M}_\tau^{\tau_1} \cap \{r \geq r_0\}} \bigg(}
		+ \sum_{j \leq N_2 - 1} \epsilon^2 r^{1-\frac{1}{2}C_{[N_2(j)]}\epsilon} |\mathscr{Y}^j F|_{(\text{frame})}^2
		+ \sum_{k \leq N_2} r^{-1} |\mathscr{Y}^k F|_{(\text{frame})}
	\bigg) \dVol_{g}
	\end{split}
	\end{equation*}
	Now, $\delta \ll 1$ proves the proposition.

\end{proof}

\begin{proposition}[$L^2$ bounds on $\hat{\chi}$ below top-order]
	\label{proposition L2 chihat low}
	Suppose that the bootstrap assumptions hold. Then for all $n \leq N_2 - 1$ we have
	\begin{equation*}
	\begin{split}
	&\int_{\mathcal{M}_\tau^{\tau_1} \cap \{r \geq r_0\}}
		\delta r^{-1+(\frac{1}{2}-c_{[n]})\delta} |\mathscr{Y}^n \hat{\chi}|^2
	\dVol_{g}
	\\
	&\lesssim
	\epsilon^2 (1 + c_{[n]}^{-1}) \epsilon^{2(N_2+1-n)} (1+\tau)^{-1 + C_{(n)}\delta}
	\\
	&\phantom{\lesssim}
	+ \int_{\mathcal{M}_\tau^{\tau_1} \cap \{r \geq \frac{1}{2}r_0\}} \bigg(
		\delta r^{-1+(\frac{1}{2}-c_{[n]})\delta} |\overline{\slashed{\D}} \mathscr{Y}^{n+1} h|^2_{(\text{frame})}
		+ \delta r^{-1+\frac{1}{2}\delta} |\overline{\slashed{\D}} \mathscr{Y}^{n} h|^2_{(\text{frame})}
		\\
		&\phantom{\lesssim + \int_{\mathcal{M}_\tau^{\tau_1} \cap \{r \geq \frac{1}{2}r_0\}} \bigg(}
		+ \epsilon^2 \delta r^{-1-\delta} |\slashed{\D} \mathscr{Y}^{n} h|^2_{(\text{frame})}
		+ \delta r^{-1-\delta} |\slashed{\D} \mathscr{Y}^n h|^2_{(\text{frame})}
		+ \delta r^{-3-\delta} |\mathscr{Y}^n h|^2_{(\text{frame})}
	\bigg) \dVol_{g}
	\end{split}
	\end{equation*}
\end{proposition}

\begin{proof}
	Recall proposition \ref{proposition transport Yn chihat}. This gives us
	\begin{equation*}
	\begin{split}
	\slashed{\D}_L \left( r^2 \mathscr{Y}^n \hat{\mathcal{X}} \right)
	&=
	\bm{\Gamma}^{(0)}_{(-1)} r^2 \mathscr{Y}^n \mathcal{X}_{(\text{low})}
	+ r \bm{\Gamma}^{(0)}_{(C_{(0)}\epsilon)}(\overline{\slashed{\D}} \mathscr{Y}^{n+1} h)_{(\text{frame})}
	\\
	&\phantom{=}
	+ r \bm{\Gamma}^{(0)}_{(C_{(0)}, \text{large})} (\overline{\slashed{\D}}\mathscr{Z}h_{(\text{rect})}) (\mathscr{Y}^n X_{(\text{frame})})
	+ r \bm{\Gamma}^{(1)}_{(C_{(1)}\epsilon + C_{(0)}\epsilon, \text{large})} (\overline{\slashed{\D}} \mathscr{Y}^n h_{(\text{rect})})
	\\
	&\phantom{=}
	+ r^2 \bm{\Gamma}^{(0)}_{(-1-\delta + C_{(0)}\epsilon)} (\slashed{\D} \mathscr{Y}^n h)_{(\text{frame})}
	+ \bm{\Gamma}^{(0)}_{(\frac{1}{2}+\delta)} \bm{\Gamma}^{(n)}_{(-1+C_{(n)}\epsilon)}
	+ r^2 \bm{\Gamma}^{(0)}_{(-1-\delta + C_{(0)}\epsilon)} \bm{\Gamma}^{(n)}_{(-1-\delta)}
	\\
	&\phantom{=}
	+ r \bm{\Gamma}^{(n-1)}_{(-1-\delta + C_{(0)}\epsilon)}
	\end{split}
	\end{equation*}
	
	Note that, apart from an additional factor of $r^{C_{(0)}\epsilon}$ in some of the error terms, $\mathscr{Y}^n\hat{\chi}$ satisfies a transport equation with an identical schematic form to that satisfied by $\mathscr{Y}^n \tr_{\slashed{g}}\chi_{(\text{low})}$ (see proposition \ref{proposition L2 tr chi low}). Indeed, these additional factors only multiply terms involving the fields $h_{(\text{frame})}$ and its derivatives. Hence, following an identical set of computations, we find

	\begin{equation*}
	\begin{split}
	&\int_{\mathcal{M}_\tau^{\tau_1} \cap \{r \geq r_0\}}
	r^{-1+(\frac{1}{2}-c_{[n]})\delta} |\mathscr{Y}^n \hat{\mathcal{X}}|^2
	\dVol_{g}
	\\
	&\lesssim
	\epsilon^2 \delta^{-1} (1 + c_{[n]}^{-1}) \epsilon^{2(N_2+1-n)} (1+\tau)^{-1 + C_{(n)}\delta}
	\\
	&\phantom{\lesssim}
	+ \int_{\mathcal{M}_\tau^{\tau_1} \cap \{r \geq \frac{1}{2}r_0\}} \bigg(
		r^{-1+(\frac{1}{2}-c_{[n]})\delta} |\overline{\slashed{\D}} \mathscr{Y}^{n+1} h|^2_{(\text{frame})}
		+ r^{-1+\frac{1}{2}\delta} |\overline{\slashed{\D}} \mathscr{Y}^{n} h|^2_{(\text{frame})}
		\\
		&\phantom{\lesssim + \int_{\mathcal{M}_\tau^{\tau_1} \cap \{r \geq \frac{1}{2}r_0\}} \bigg(}
		+ \epsilon^2 r^{-1-\delta} |\slashed{\D} \mathscr{Y}^{n} h|^2_{(\text{frame})}
		+ r^{-1-\delta} |\slashed{\D} \mathscr{Y}^n h|^2_{(\text{frame})}
		+ r^{-3-\delta} |\mathscr{Y}^n h|^2_{(\text{frame})}
	\bigg) \dVol_{g}
	\end{split}
	\end{equation*}

	Now, we need to relate $\mathscr{Y}^n \hat{\mathcal{X}}$ to $\mathscr{Y}^n \hat{\chi}$. Proposition \ref{proposition transport Yn chihat} gives
	\begin{equation*}
	\begin{split}
	\mathscr{Y}^n \hat{\mathcal{X}}
	&=
	\mathscr{Y}^n \hat{\chi}
	+ (\overline{\slashed{\D}} \mathscr{Y}^n h)_{(\text{frame})}
	+ (\bar{\partial} h_{(\text{rect})}) \bm{\Gamma}^{(0)}_{(0, \text{large})} (\mathscr{Y}^n X_{(\text{frame})})
	+ \bm{\Gamma}^{(n)}_{(-\delta)} \bm{\Gamma}^{(0)}_{(-1+C_{(0)}\epsilon)}
	\\
	&\phantom{=}
	+ \bm{\Gamma}^{(0)}_{(-\delta)} \bm{\Gamma}^{(n)}_{(-1+C_{(n)}\epsilon)}
	+ \bm{\Gamma}^{(n-1)}_{(-1-\delta)}
	\end{split}
	\end{equation*}
	and so, in the region $r \geq r_0$
	\begin{equation*}
	\begin{split}
	|\mathscr{Y}^n \hat{\chi}|
	&\lesssim
	|\mathscr{Y}^n \hat{\mathcal{X}}|
	+ |\overline{\slashed{\D}} \mathscr{Y}^n h|_{(\text{frame})}
	+ |\bar{\partial} h_{(\text{rect})}| |\mathscr{Y}^n X_{(\text{frame})}|
	+ \epsilon r^{-1+C_{(0)}\epsilon} |\bm{\Gamma}^{(n)}_{(-\delta)}|
	\\
	&\phantom{=}
	+ \epsilon r^{-\delta} |\bm{\Gamma}^{(n)}_{(-1+C_{(n)}\epsilon)}|
	+ |\bm{\Gamma}^{(n-1)}_{(-1-\delta)}|
	\end{split}
	\end{equation*}
	Using this, it is easy to see that $\mathscr{Y}^n \hat{\chi}$ satisfies the same bound in $L^2$ as $\hat{\mathcal{X}}$.
	
\end{proof}

Note that the proposition above also ``loses'' a derivative, so it cannot be used to bound $\mathscr{Y}^{N_2} \hat{\chi}$. For this, we need the following proposition.

\begin{proposition}[$L^2$ bounds on $\hat{\chi}$ at top order]
	\label{proposition L2 chihat high}
	Suppose that the bootstrap bounds hold. Suppose also that the inhomogeneous terms $F_{(\text{frame})}$ satisfy the pointwise bounds
	\begin{equation*}
	|\mathscr{Y}^n F|_{(\text{frame})} \lesssim \epsilon^2 (1+r)^{-2 + 2C_{(n)}\epsilon}
	\end{equation*}
	for all $n \leq N_1$. Also, suppose that $N_2 \leq 2N_1$. Finally, suppose that
	\begin{equation*}
	C_{[N_2(T)]} \leq C_{[N_2]} - 2C_{(0)}
	\end{equation*}
	
	Then we have
	\begin{equation*}
	\begin{split}
	&\int_{\mathcal{M}_\tau^{\tau_1}\cap\{r \geq r_0\}} C_{[N_2(j)]}\epsilon r^{-1-C_{[N_2(j)]}\epsilon} \left( |\mathscr{Y}^{N_2} \hat{\chi}|^2 \right) \dVol_{g}
	\\
	&\lesssim
	C_{[N_2(j)]} \delta^{-1} \epsilon^5 (1+\tau)^{-1+C_{(N_2)}\delta}
	\\
	&\phantom{=}
	\int_{\mathcal{M}_\tau^{\tau_1} \cap \{r \geq r_0\}} \bigg(
		C_{[N_2(j)]}\epsilon r^{-1-C_{[N_2(j+1)]}\epsilon} |\slashed{\D} \slashed{\D}_T \mathscr{Y}^{N_2-1} h|^2_{LL}
		+ C_{[N_2(j)]}\epsilon r^{-1-C_{[N_2(j)]}\epsilon} |\slashed{\D} \mathscr{Y}^{N_2} h|^2_{LL}
		\\
		&\phantom{=	\int_{\mathcal{M}_\tau^{\tau_1} \cap \{r \geq r_0\}} \bigg(}
		+ C_{[N_2(j)]}\epsilon r^{-1-\delta} |\slashed{\D} \mathscr{Y}^{N_2} h|^2_{(\text{frame})}
		+ C_{[N_2(j)]} \epsilon^3 r^{-1-\frac{1}{2}C_{[N_2(j)]}\epsilon} |\overline{\slashed{\D}} \mathscr{Y}^{N_2} h|^2_{(\text{frame})}
		\\
		&\phantom{=	\int_{\mathcal{M}_\tau^{\tau_1} \cap \{r \geq r_0\}} \bigg(}
		+ C_{[N_2(j)]}\epsilon r^{-1-C_{[N_2(j)]}\epsilon} |\overline{\slashed{\D}} \mathscr{Y}^{n} h|^2_{(\text{frame})}
		+ C_{[N_2(j)]}\epsilon r^{-3-\delta} |\mathscr{Y}^{N_2} h|^2_{(\text{frame})}
		\\
		&\phantom{=	\int_{\mathcal{M}_\tau^{\tau_1} \cap \{r \geq r_0\}} \bigg(}
		+ C_{[N_2(j)]}\epsilon r^{1-C_{[N_2(j)]}\epsilon} |\slashed{\D}_T \mathscr{Y}^{N_2-1} F|^2_{LL}
		+ \sum_{j \leq N_2 - 1} C_{[N_2(j)]}\epsilon^3 r^{1-\frac{1}{2}C_{[N_2(T)]}\epsilon} |\mathscr{Y}^j F|_{(\text{frame})}^2
		\\
		&\phantom{=	\int_{\mathcal{M}_\tau^{\tau_1} \cap \{r \geq r_0\}} \bigg(}
		+ \sum_{j \leq N_2} C_{[N_2(j)]}\epsilon r^{-1} |\mathscr{Y}^j F|_{(\text{frame})}
	\bigg)\dVol_g
	\end{split}
	\end{equation*}
	
\end{proposition}

\begin{proof}

	The proof of this proposition is perhaps the most complicated of any of the top-order $L^2$ estimates which we need to prove, since we have to consider each of the commutator operators $\slashed{\D}_T$, $r\slashed{\D}_L$ and $r\slashed{\nabla}$ in turn.
	
	First, we examine $\slashed{\D}_T \mathscr{Y}^{n-1} \hat{\chi}$. Using proposition \ref{proposition DT Yn chi hat} we have
	\begin{equation*}
	\begin{split}
	&\int_{\mathcal{M}_\tau^{\tau_1} \cap \{r \geq r_0\}} r^{-1-C_{[N_2(j)]}\epsilon} |\slashed{\D}_T \mathscr{Y}^{n-1} \hat{\chi}|^2 \dVol_g
	\\
	&\lesssim
	\int_{\mathcal{M}_\tau^{\tau_1} \cap \{r \geq r_0\}} r^{-1-C_{[N_2(j)]}\epsilon} \bigg(
		r^{-2} |\slashed{\D} \mathscr{Y}^{N_2} h|^2_{(\text{frame})}
		+ |\slashed{\nabla}^2 \mathscr{Y}^{N_2-1} \log \mu|^2
		+ \epsilon^2 r^{-2+2C_{(1)}\epsilon} |\bm{\Gamma}^{(N_2)}_{(-1+C_{(N_2)}\epsilon)}|^2
		\\
		&\phantom{\lesssim \int_{\mathcal{M}_\tau^{\tau_1} \cap \{r \geq r_0\}} r^{-1-C_{[N_2(j)]}\epsilon} \bigg(}
		+ |\bm{\Gamma}^{(n-1)}_{(-1+C_{(N_2-1)}\epsilon)}|^2
	\bigg)\dVol_g
	\end{split}
	\end{equation*}
	and we can use proposition \ref{proposition L2 nabla log mu} to bound the term involving $\slashed{\nabla}^2 \mathscr{Z}^{n-1} \log \mu$, together with the pointwise bootstrap bounds, to obtain
	\begin{equation}
	\label{equation L2 chihat high internal 1}
	\begin{split}
	&\int_{\mathcal{M}_\tau^{\tau_1} \cap \{r \geq r_0\}} r^{-1-C_{[N_2]}\epsilon} |\slashed{\D}_T \mathscr{Y}^{n-1} \hat{\chi}|^2 \dVol_g
	\\
	&\lesssim
	\delta^{-1}\epsilon^4 (1+\tau)^{-1+C_{(N_2)}\delta}
	\\
	&\phantom{=}
	\int_{\mathcal{M}_\tau^{\tau_1} \cap \{r \geq r_0\}} \bigg(
		r^{-1-C_{[N_2(j+1)]}\epsilon} |\slashed{\D} \slashed{\D}_T \mathscr{Y}^{N_2-1} h|^2_{LL}
		+ r^{-1-\delta} |\slashed{\D} \mathscr{Y}^{N_2} h|^2_{(\text{frame})}
		\\
		&\phantom{=	\int_{\mathcal{M}_\tau^{\tau_1} \cap \{r \geq r_0\}} \bigg(}
		+ \epsilon^2 r^{-1-\frac{1}{2}C_{[N_2(j)]}\epsilon} |\overline{\slashed{\D}} \mathscr{Y}^{N_2} h|^2_{(\text{frame})}
		+ \epsilon^2 \delta^{-2} r^{-1+\frac{1}{2}\delta} |\overline{\slashed{\D}} \mathscr{Y}^{N_2} h|^2_{(\text{frame})}
		\\
		&\phantom{=	\int_{\mathcal{M}_\tau^{\tau_1} \cap \{r \geq r_0\}} \bigg(}
		+ r^{-3-\delta} |\mathscr{Y}^{N_2} h|^2_{(\text{frame})}
		+ r^{1-C_{[N_2(j)]}\epsilon} |\slashed{\D}_T \mathscr{Y}^{N_2-1} F|^2_{LL}
		\\
		&\phantom{=	\int_{\mathcal{M}_\tau^{\tau_1} \cap \{r \geq r_0\}} \bigg(}
		+ \sum_{j \leq N_2 - 1} \epsilon^2 r^{1-\frac{1}{2}C_{[N_2(j)]}\epsilon} |\mathscr{Y}^j F|_{(\text{frame})}^2
		+ \sum_{j \leq N_2} r^{-1} |\mathscr{Y}^j F|_{(\text{frame})}
	\bigg)\dVol_g
	\end{split}
	\end{equation}

	Next, we examine $r\slashed{\D}_L \mathscr{Y}^{n-1} \hat{\chi}$. Using proposition \ref{proposition transport Yn chihat} we have, schematically,
	\begin{equation*}
	\begin{split}
	r\slashed{\D}_L (\mathscr{Y}^{n-1} \hat{\chi})
	=
	(\overline{\slashed{\D}} \mathscr{Y}^{n} h)_{(\text{frame})}
	+ \bm{\Gamma}^{(n-1)}_{(-1-\delta)}
	\end{split}
	\end{equation*}
	and so
	\begin{equation}
	\label{equation L2 chihat high internal 2}
	\begin{split}
	&\int_{\mathcal{M}_\tau^{\tau_1} \cap \{r \geq r_0\}} r^{-1-C_{[N_2]}\epsilon} |r\slashed{\D}_L \mathscr{Y}^{N_2-1} \hat{\chi}|^2 \dVol_g
	\\
	&\lesssim
	\int_{\mathcal{M}_\tau^{\tau_1} \cap \{r \geq r_0\}} r^{-1-C_{[N_2]}\epsilon} \left(
		|\overline{\slashed{\D}} \mathscr{Y}^{N_2} h|^2_{(\text{frame})}
		+ |\bm{\Gamma}^{(N_2-1)}_{(-1-\delta)}|^2
	\right) \dVol_g
	\\
	\\
	&\lesssim
	\epsilon^4 (1+\tau)^{-1+C_{(N_2)}\delta}
	+\int_{\mathcal{M}_\tau^{\tau_1} \cap \{r \geq r_0\}} 
		r^{-1-C_{[N_2]}\epsilon} 
		|\overline{\slashed{\D}} \mathscr{Y}^{n} h|^2_{(\text{frame})}
	\dVol_g
	\end{split}
	\end{equation}

	Finally, we treat the case where the final commutation operator applied is an angular derivative. In this case, we can first use proposition \ref{proposition elliptic estimates Hodge} (together with the pointwise bounds on $\Omega$ and $K$) to obtain the bound
	\begin{equation*}
	\int_{S_{\tau,r}} \left( |\slashed{\nabla} \mathscr{Y}^{N_2-1} \hat{\chi}|^2 \right) \dVol_{\mathbb{S}^2}
	\lesssim
	\int_{S_{\tau,r}} \left(
	|\slashed{\Div} \mathscr{Y}^{N_2-1} \hat{\chi}|^2
	+ r^{-2}|\mathscr{Y}^{N_2-1} \hat{\chi}|^2
	\right) \dVol_{\mathbb{S}^2}
	\end{equation*}
	and then we can use proposition \ref{proposition commuted div chihat} to bound $(\slashed{\Div}\mathscr{Y}^{N_2-1} \hat{\chi})$, which shows
	\begin{equation*}
	\begin{split}
	&\int_{\Sigma_{\tau}\cap\{r \geq r_0\}} r^{1-C_{[N_2]}\epsilon} \left( |r\slashed{\nabla} \mathscr{Y}^{n-1} \hat{\chi}|^2 \right) \upd r \wedge \dVol_{\mathbb{S}^2}
	\\
	&\lesssim
	\int_{\Sigma_{\tau}\cap\{r \geq r_0\}} r^{1-C_{[N_2]}\epsilon} \bigg(
		r^{-2} |\overline{\slashed{\D}} \mathscr{Y}^{N_2} h|^2_{(\text{frame})}
		+ r^{-2} |\mathscr{Y}^{N_2} \tr_{\slashed{g}}\chi_{(\text{small})}|^2
		+ \sum_{j \leq N_2-1} r^{-2} |\bm{\Gamma}^{(j)}_{(-1+C_{(j)}\epsilon)}|^2
		\\
		&\phantom{\lesssim \int_{\Sigma_{\tau}\cap\{r \geq r_0\}} r^{1-C_{[N_2]}\epsilon} \bigg(}
		+ \sum_{j+k \leq N_2-1} |\bm{\Gamma}^{(j)}_{(-1+C_{(j)}\epsilon)}|^2 |\bm{\Gamma}^{(k)}_{(-1-\delta)}|^2
		\bigg) r^2 \upd r \wedge \dVol_{\mathbb{S}^2}
		\\
	&\lesssim
	\int_{\Sigma_{\tau}\cap\{r \geq r_0\}} \bigg(
		r^{-1-C_{[N_2]}\epsilon} |\overline{\slashed{\D}} \mathscr{Y}^{N_2} h|^2_{(\text{frame})}
		+ r^{-1-C_{[N_2]}\epsilon} |\mathscr{Y}^{N_2} \tr_{\slashed{g}} \chi_{(\text{small})}|^2
		+ r^{-1-C_{[N_2]}\epsilon} |\bm{\Gamma}^{(N_2-1)}_{(-1+C_{(n-1)}\epsilon)}|^2
		\\
		&\phantom{\lesssim \int_{\Sigma_{\tau}\cap\{r \geq r_0\}} \bigg(}		
		+ \epsilon^2 r^{-1-\frac{1}{2}C_{[N_2]}} |\bm{\Gamma}^{(N_2-1)}_{(-1-\delta)}|^2
	\bigg) r^2 \upd r \wedge \dVol_{\mathbb{S}^2}
	\end{split}
	\end{equation*}
	
	Integrating over $\tau$ and using $\Omega \sim r$, we can then use proposition \ref{proposition L2 tr chi high} to bound the term involving $\mathscr{Y}^N \tr_{\slashed{g}}\chi_{(\text{small})}$ and the $L^2$ bootstrap bounds from equation \eqref{equation bootstrap L2 geometric} to find
	\begin{equation}
	\label{equation L2 chihat high internal 3}
	\begin{split}
	&\int_{\mathcal{M}_{\tau}^{\tau_1}\cap\{r \geq r_0\}} r^{-1-C_{[N_2]}\epsilon} \left( |r\slashed{\nabla} \mathscr{Y}^{n-1} \hat{\chi}|^2 \right) \dVol_{g}
	\\
	&\lesssim
	\epsilon^4 (1+\tau)^{-1+C_{(N_2)}\delta}
	\\
	&\phantom{\lesssim}
	+ \int_{\mathcal{M}_{\tau}^{\tau_1}\cap\{r \geq r_0\}} \bigg(
		r^{-1-C_{[N_2]}\epsilon} |\slashed{\D} \mathscr{Y}^{N_2} h|^2_{LL}
		+ \epsilon^2 r^{-1-\frac{1}{2}C_{[N_2]}\epsilon} |\overline{\slashed{\D}} \mathscr{Y}^{N_2} h|^2_{(\text{frame})}
		+ r^{-1-C_{[N_2]}\epsilon} |\overline{\slashed{\D}} \mathscr{Y}^{N_2} h|^2_{(\text{frame})}
		\\
		&\phantom{\lesssim + \int_{\mathcal{M}_\tau^{\tau_1} \cap \{r \geq \frac{1}{2}r_0\}}\bigg(}
		+ r^{-1-\delta} |\slashed{\D} \mathscr{Y}^{N_2} h|^2_{(\text{frame})}
		+ r^{-3-\delta} |\mathscr{Y}^{N_2} h|^2_{(\text{frame})}
		\\
		&\phantom{\lesssim + \int_{\mathcal{M}_\tau^{\tau_1} \cap \{r \geq \frac{1}{2}r_0\}}\bigg(}
		+ r^{1-C_{[N_2]}\epsilon} |\mathscr{Y}^{N_2} F|^2_{LL}
		+ \sum_{j \leq N_2 - 1} \epsilon^2 r^{1-\frac{1}{2}C_{[N_2]}\epsilon} |\mathscr{Y}^j F|_{(\text{frame})}^2
	\bigg) \dVol_{g}
	\end{split}
	\end{equation}
	
	Combining equations \eqref{equation L2 chihat high internal 1}, \eqref{equation L2 chihat high internal 2} and \eqref{equation L2 chihat high internal 3} proves the proposition.

\end{proof}

\begin{proposition}[$L^2$ bounds on $\mathscr{Y}^n \hat{\chibar}$ below top-order]
	\label{proposition L2 chibar hat low}
	Suppose that the bootstrap bounds hold.
	
	Then for all $n \leq N_2$ the quantity $\mathscr{Y}^n \hat{\chibar}$ obeys the bound
	\begin{equation*}
	\begin{split}
	&\int_{\mathcal{M}_\tau^{\tau_1} \cap \{r \geq r_0\}} r^{-1-C_{[n]}\epsilon} |\mathscr{Y}^n \hat{\chi}|^2 \dVol_g
	\\
	&\lesssim
	\epsilon^{2(N_2 + 2 - n)} (1+\tau)^{-1+C_{(n)}\delta}
	\\
	&\phantom{\lesssim}
	+ \int_{\mathcal{M}_\tau^{\tau_1} \cap \{r \geq r_0\}} \bigg(
	r^{-1-C_{[n]}\epsilon} |\slashed{\D} \mathscr{Y}^n h|^2_{(\text{frame})}
	+ r^{-1+\delta} |\overline{\slashed{\D}} \mathscr{Y}^{n+1} h|^2
	+ r^{-1+\delta} |\overline{\slashed{\D}} \mathscr{Y}^{n} h|^2
	\\
	&\phantom{\lesssim + \int_{\mathcal{M}_\tau^{\tau_1} \cap \{r \geq r_0\}} \bigg(}
	+ r^{-3-C_{[n]}\epsilon} |\mathscr{Y}^n h|^2_{(\text{frame})}
	\bigg) \dVol_g
	\end{split}
	\end{equation*}
	and also
	\begin{equation*}
	\begin{split}
	&\int_{\mathcal{M}_\tau^{\tau_1} \cap \{r \geq r_0\}} c\delta r^{-1-c\delta} |\mathscr{Y}^n \hat{\chibar}|^2 \dVol_g
	\\
	&\lesssim
	c \epsilon^2 (1+c_{[n]}^{-1}) \epsilon^{2(N_2 + 1 - n)} (1+\tau)^{-1+C_{(n)}\delta}
	\\
	&\phantom{\lesssim}
	+ \int_{\mathcal{M}_\tau^{\tau_1} \cap \{r \geq r_0\}} \bigg(
		c\delta r^{-1-C_{[n]}\epsilon} |\slashed{\D} \mathscr{Y}^n h|^2_{(\text{frame})}
		+ c\delta r^{-1+\frac{1}{2}\delta} |\overline{\slashed{\D}} \mathscr{Y}^{n+1} h|^2
		+ c\delta r^{-1+\frac{1}{2}\delta} |\overline{\slashed{\D}} \mathscr{Y}^{n} h|^2
		\\
		&\phantom{\lesssim + \int_{\mathcal{M}_\tau^{\tau_1} \cap \{r \geq r_0\}} \bigg(}
		+ c\delta r^{-3-C_{[n]}\epsilon} |\mathscr{Y}^n h|^2_{(\text{frame})}
		+ c\delta \epsilon^2 r^{-1-2\delta} |\bm{\Gamma}^{(n)}_{(-1+C_{(n)}\epsilon)}|^2
		+ c\delta r^{-1-C_{[n]}\epsilon} |\bm{\Gamma}^{(n-1)}_{(-1-\delta)}|^2
	\bigg) \dVol_g
	\end{split}
	\end{equation*}
\end{proposition}

\begin{proof}
	Recall proposition \ref{proposition Yn chibar hat}, which provides the schematic expression
	\begin{equation*}
	\begin{split}
	\mathscr{Y}^n \hat{\chibar}
	&=
	\mathscr{Y}^n \hat{\chi}
	+ (\slashed{\D} \mathscr{Y}^n h)_{(\text{frame})}
	+ \bm{\Gamma}^{(0)}_{(-1, \text{large})} (\mathscr{Y}^n h)_{(\text{frame})}
	+ \bm{\Gamma}^{(0)}_{(-1 + C_{(0)}\epsilon)} (\mathscr{Y}^n X_{(\text{frame})})
	+ \bm{\Gamma}^{(0)}_{(-\delta)} \bm{\Gamma}^{(n)}_{(-1+C_{(n)}\epsilon)}
	\\
	&\phantom{=}
	+ \bm{\Gamma}^{(n-1)}_{(-1 + 2C_{(n-1)}\epsilon)}
	\end{split}
	\end{equation*}
	and so we have
	\begin{equation*}
	\begin{split}
	&\int_{\mathcal{M}_\tau^{\tau_1} \cap \{r \geq r_0\}} r^{-1-C_{[n]}\epsilon} |\mathscr{Y}^n \hat{\chibar}|^2 \dVol_g
	\\
	&\lesssim
	\int_{\mathcal{M}_\tau^{\tau_1} \cap \{r \geq r_0\}} \bigg(
		r^{-1-C_{[n]}\epsilon} |\mathscr{Y}^n \hat{\chi}|^2
		+ r^{-1-C_{[n]}\epsilon} |\slashed{\D} \mathscr{Y}^n h|^2_{(\text{frame})}
		+ r^{-3-C_{[n]}\epsilon} |\mathscr{Y}^n h|^2_{(\text{frame})}
		\\
		&\phantom{\lesssim \int_{\mathcal{M}_\tau^{\tau_1} \cap \{r \geq r_0\}} \bigg(}
		+ \epsilon^2 r^{-3-\frac{1}{2}C_{[n]}\epsilon} |\mathscr{Y}^n X_{(\text{frame, small})}|^2
		+ \epsilon^2 r^{-1-2\delta} |\bm{\Gamma}^{(n)}_{(-1+C_{(n)}\epsilon)}|^2
		\\
		&\phantom{\lesssim \int_{\mathcal{M}_\tau^{\tau_1} \cap \{r \geq r_0\}} \bigg(}
		+ r^{-1-C_{[n]}\epsilon} |\bm{\Gamma}^{(n-1)}_{(-1-\delta)}|^2
	\bigg) \dVol_g
	\end{split}
	\end{equation*}
	
	Now, bounding the term involving $\mathscr{Y}^n X_{(\text{frame, small})}$ using proposition \ref{proposition L2 bounds rectangular} and the term involving $\mathscr{Y}^n \hat{\chi}$ using proposition \ref{proposition L2 chihat low}, and using the fact that $\epsilon \ll \delta$, we have
	\begin{equation*}
	\begin{split}
	&\int_{\mathcal{M}_\tau^{\tau_1} \cap \{r \geq r_0\}} r^{-1-C_{[n]}\epsilon} |\mathscr{Y}^n \hat{\chibar}|^2 \dVol_g
	\\
	&\lesssim
	\delta^{-1} \epsilon^2 (1+c_{[n]}^{-1}) \epsilon^{2(N_2 + 1 - n)} (1+\tau)^{-1+C_{(n)}\delta}
	\\
	&\phantom{\lesssim}
	+ \int_{\mathcal{M}_\tau^{\tau_1} \cap \{r \geq r_0\}} \bigg(
		r^{-1-C_{[n]}\epsilon} |\slashed{\D} \mathscr{Y}^n h|^2_{(\text{frame})}
		+ r^{-1+\frac{1}{2}\delta} |\overline{\slashed{\D}} \mathscr{Y}^{n+1} h|^2
		+ r^{-1+\frac{1}{2}\delta} |\overline{\slashed{\D}} \mathscr{Y}^{n} h|^2
		\\
		&\phantom{\lesssim + \int_{\mathcal{M}_\tau^{\tau_1} \cap \{r \geq r_0\}} \bigg(}
		+ r^{-3-C_{[n]}\epsilon} |\mathscr{Y}^n h|^2_{(\text{frame})}
		+ \epsilon^2 r^{-1-2\delta} |\bm{\Gamma}^{(n)}_{(-1+C_{(n)}\epsilon)}|^2
		+ r^{-1-C_{[n]}\epsilon} |\bm{\Gamma}^{(n-1)}_{(-1-\delta)}|^2
	\bigg) \dVol_g
	\end{split}
	\end{equation*}
	Finally, using the $L^2$ bootstrap bounds for the last two terms proves the first part of the proposition.
	
	Following very similar computations, we have
	\begin{equation*}
	\begin{split}
	&\int_{\mathcal{M}_\tau^{\tau_1} \cap \{r \geq r_0\}} r^{-1-c\delta} |\mathscr{Y}^n \hat{\chi}|^2 \dVol_g
	\\
	&\lesssim
	\int_{\mathcal{M}_\tau^{\tau_1} \cap \{r \geq r_0\}} \bigg(
		r^{-1-c\delta} |\mathscr{Y}^n \hat{\chi}|^2
		+ r^{-1-c\delta} |\slashed{\D} \mathscr{Y}^n h|^2_{(\text{frame})}
		+ r^{-3-c\delta} |\mathscr{Y}^n h|^2_{(\text{frame})}
		\\
		&\phantom{\lesssim \int_{\mathcal{M}_\tau^{\tau_1} \cap \{r \geq r_0\}} \bigg(}
		+ \epsilon^2 r^{-3-\frac{1}{2}c\delta} |\mathscr{Y}^n X_{(\text{frame, small})}|^2
		+ \epsilon^2 r^{-1-2\delta} |\bm{\Gamma}^{(n)}_{(-1+C_{(n)}\epsilon)}|^2
		+ r^{-1-c\delta} |\bm{\Gamma}^{(n-1)}_{(-1-\delta)}|^2
	\bigg) \dVol_g
	\end{split}
	\end{equation*}
	again, substituting for $\mathscr{Y}^n X_{(\text{frame, small})}$ using proposition \ref{proposition L2 bounds rectangular} and the term involving $\mathscr{Y}^n \hat{\chi}$ using proposition \ref{proposition L2 chihat low}, and using the fact that $\epsilon \ll c\delta \ll 1$ we have
	\begin{equation*}
	\begin{split}
	&\int_{\mathcal{M}_\tau^{\tau_1} \cap \{r \geq r_0\}} r^{-1-c\delta} |\mathscr{Y}^n \hat{\chibar}|^2 \dVol_g
	\\
	&\lesssim
	\delta^{-1} \epsilon^2 (1+c_{[n]}^{-1}) \epsilon^{2(N_2 + 1 - n)} (1+\tau)^{-1+C_{(n)}\delta}
	\\
	&\phantom{\lesssim}
	+ \int_{\mathcal{M}_\tau^{\tau_1} \cap \{r \geq r_0\}} \bigg(
		r^{-1-C_{[n]}\epsilon} |\slashed{\D} \mathscr{Y}^n h|^2_{(\text{frame})}
		+ r^{-1+\frac{1}{2}\delta} |\overline{\slashed{\D}} \mathscr{Y}^{n+1} h|^2
		+ r^{-1+\frac{1}{2}\delta} |\overline{\slashed{\D}} \mathscr{Y}^{n} h|^2
		\\
		&\phantom{\lesssim + \int_{\mathcal{M}_\tau^{\tau_1} \cap \{r \geq r_0\}} \bigg(}
		+ r^{-3-C_{[n]}\epsilon} |\mathscr{Y}^n h|^2_{(\text{frame})}
		+ \epsilon^2 r^{-1-2\delta} |\bm{\Gamma}^{(n)}_{(-1+C_{(n)}\epsilon)}|^2
		+ r^{-1-C_{[n]}\epsilon} |\bm{\Gamma}^{(n-1)}_{(-1-\delta)}|^2
	\bigg) \dVol_g
	\end{split}
	\end{equation*}
	Again, bounding the final two terms using the $L^2$ bootstrap bounds in equation \eqref{equation bootstrap L2 geometric} proves the second part of the proposition.
	
\end{proof}

\begin{proposition}[$L^2$ bounds on $\mathscr{Y}^n \hat{\chibar}$ at top-order]
	\label{proposition L2 chibar hat high}
	Suppose that the bootstrap bounds hold. Suppose also that the inhomogeneous terms $F_{(\text{frame})}$ satisfy the pointwise bounds
	\begin{equation*}
	|\mathscr{Y}^n F|_{(\text{frame})} \lesssim \epsilon^2 (1+r)^{-2 + 2C_{(n)}\epsilon}
	\end{equation*}
	for all $n \leq N_1$. Also, suppose that $N_2 \leq 2N_1$.
	
	Then we have
	\begin{equation*}
	\begin{split}
	&\int_{\mathcal{M}_\tau^{\tau_1} \cap \{r \geq r_0\}} C_{[N_2]}\epsilon r^{-1-C_{[N_2]}\epsilon} |\mathscr{Y}^{N_2} \hat{\chibar}|^2 \dVol_g
	\\
	&\lesssim
	C_{[N_2]}\delta^{-1} \epsilon^5 (1+\tau)^{-1+C_{(N_2)}\delta}
	\\
	&\phantom{\lesssim}
	+ \int_{\mathcal{M}_\tau^{\tau_1} \cap \{r \geq r_0\}} \bigg(
		C_{[N_2]}\epsilon r^{-1-(C_{[N_2]} - 2C_{(0)})\epsilon} |\slashed{\D} \slashed{\D}_T \mathscr{Y}^{N_2-1} h|^2_{LL}
		+ C_{[N_2]}\epsilon r^{-1-C_{[N_2]}\epsilon} |\slashed{\D} \mathscr{Y}^{N_2} h|^2_{(\text{frame})}
		\\
		&\phantom{\lesssim + \int_{\mathcal{M}_\tau^{\tau_1} \cap \{r \geq r_0\}} \bigg(}
		+ C_{[N_2]}\epsilon r^{-1-\delta} |\slashed{\D} \mathscr{Y}^{N_2} h|^2_{(\text{frame})}
		+ C_{[N_2]}\epsilon^3 r^{-1-\frac{1}{2}C_{[N_2(T)]}\epsilon} |\overline{\slashed{\D}} \mathscr{Y}^{N_2} h|^2_{(\text{frame})}
		\\
		&\phantom{\lesssim + \int_{\mathcal{M}_\tau^{\tau_1} \cap \{r \geq r_0\}} \bigg(}
		+ C_{[N_2]}\epsilon r^{-1-\frac{1}{2}C_{[N_2]}\epsilon} |\overline{\slashed{\D}} \mathscr{Y}^{n} h|^2_{(\text{frame})}
		+ C_{[N_2]}\epsilon r^{-3-C_{[N_2]}\epsilon} |\mathscr{Y}^{N_2} h|^2_{(\text{frame})}
		\\
		&\phantom{\lesssim + \int_{\mathcal{M}_\tau^{\tau_1} \cap \{r \geq r_0\}} \bigg(}
		+ C_{[N_2]}\epsilon r^{1-C_{[N_2]}\epsilon} |\slashed{\D}_T \mathscr{Y}^{N_2-1} F|^2_{LL}
		+ \sum_{j \leq N_2 - 1} C_{[N_2]}\epsilon^3 r^{1-\frac{1}{2}C_{[N_2]}\epsilon} |\mathscr{Y}^j F|_{(\text{frame})}^2
		\\
		&\phantom{\lesssim + \int_{\mathcal{M}_\tau^{\tau_1} \cap \{r \geq r_0\}} \bigg(}
		+ \sum_{j \leq N_2} C_{[N_2]}\epsilon r^{-1} |\mathscr{Y}^j F|_{(\text{frame})}
	\bigg) \dVol_g
	\end{split}
	\end{equation*}

\end{proposition}

\begin{proof}
	As in proposition \ref{proposition L2 chibar hat low} we have
	\begin{equation*}
	\begin{split}
	&\int_{\mathcal{M}_\tau^{\tau_1} \cap \{r \geq r_0\}} r^{-1-C_{[N_2(j)]}\epsilon} |\mathscr{Y}^{N_2} \hat{\chi}|^2 \dVol_g
	\\
	&\lesssim
	\int_{\mathcal{M}_\tau^{\tau_1} \cap \{r \geq r_0\}} \bigg(
		r^{-1-C_{[N_2(j)]}\epsilon} |\mathscr{Y}^{N_2} \hat{\chi}|^2
		+ r^{-1-C_{[N_2(j)]}\epsilon} |\slashed{\D} \mathscr{Y}^{N_2} h|^2_{(\text{frame})}
		+ r^{-3-C_{[N_2(j)]}\epsilon} |\mathscr{Y}^{N_2} h|^2_{(\text{frame})}
		\\
		&\phantom{\lesssim \int_{\mathcal{M}_\tau^{\tau_1} \cap \{r \geq r_0\}} \bigg(}
		+ \epsilon^2 r^{-3-\frac{1}{2}C_{[N_2(j)]}\epsilon} |\mathscr{Y}^{N_2} X_{(\text{frame, small})}|^2
		+ \epsilon^2 r^{-1-2\delta} |\bm{\Gamma}^{(N_2)}_{(-1+C_{(N_2)}\epsilon)}|^2
		\\
		&\phantom{\lesssim \int_{\mathcal{M}_\tau^{\tau_1} \cap \{r \geq r_0\}} \bigg(}
		+ r^{-1-C_{[N_2(j)]}\epsilon} |\bm{\Gamma}^{(N_2-1)}_{(-1-\delta)}|^2
	\bigg) \dVol_g
	\end{split}
	\end{equation*}
	This time, we substitute for $\mathscr{Y}^{N_2}\hat{\chi}$ from proposition \ref{proposition L2 chihat high}. We also substitute for $\mathscr{Y}^{N_2} X_{(\text{frame, small})}$ from proposition \ref{proposition L2 bounds rectangular} and for the final two terms using the bootstrap bounds in equations \eqref{equation bootstrap L2 geometric} and \eqref{equation bootstrap L2 geometric top order} to obtain
	\begin{equation*}
	\begin{split}
	&\int_{\mathcal{M}_\tau^{\tau_1} \cap \{r \geq r_0\}} r^{-1-C_{[N_2]}\epsilon} |\mathscr{Y}^{N_2} \hat{\chibar}|^2 \dVol_g
	\\
	&\lesssim
	\delta^{-1}\epsilon^4 (1+\tau)^{-1+C_{(N_2)}\delta}
	\\
	&\phantom{\lesssim}
	+ \int_{\mathcal{M}_\tau^{\tau_1} \cap \{r \geq r_0\}} \bigg(
		r^{-1-(C_{[N_2]} - 2C_{(0)})\epsilon} |\slashed{\D} \slashed{\D}_T \mathscr{Y}^{N_2-1} h|^2_{LL}
		+ r^{-1-C_{[N_2]}\epsilon} |\slashed{\D} \mathscr{Y}^{N_2} h|^2_{(\text{frame})}
		\\
		&\phantom{\lesssim + \int_{\mathcal{M}_\tau^{\tau_1} \cap \{r \geq r_0\}} \bigg(}
		+ r^{-1-\delta} |\slashed{\D} \mathscr{Y}^{N_2} h|^2_{(\text{frame})}
		+ \epsilon^2 r^{-1-\frac{1}{2}C_{[N_2(T)]}\epsilon} |\overline{\slashed{\D}} \mathscr{Y}^{N_2} h|^2_{(\text{frame})}
		\\
		&\phantom{\lesssim + \int_{\mathcal{M}_\tau^{\tau_1} \cap \{r \geq r_0\}} \bigg(}
		+ r^{-1-\frac{1}{2}C_{[N_2]}\epsilon} |\overline{\slashed{\D}} \mathscr{Y}^{n} h|^2_{(\text{frame})}
		+ r^{-3-C_{[N_2]}\epsilon} |\mathscr{Y}^{N_2} h|^2_{(\text{frame})}
		+ r^{1-C_{[N_2]}\epsilon} |\slashed{\D}_T \mathscr{Y}^{N_2-1} F|^2_{LL}
		\\
		&\phantom{\lesssim + \int_{\mathcal{M}_\tau^{\tau_1} \cap \{r \geq r_0\}} \bigg(}
		+ \sum_{j \leq N_2 - 1} \epsilon^2 r^{1-\frac{1}{2}C_{[N_2]}\epsilon} |\mathscr{Y}^j F|_{(\text{frame})}^2
		+ \sum_{j \leq N_2} r^{-1} |\mathscr{Y}^j F|_{(\text{frame})}
	\bigg) \dVol_g
	\end{split}
	\end{equation*}

\end{proof}

\section{Putting together the \texorpdfstring{$L^2$}{L2} estimates}

The bounds in the previous subsection allow us to relate $L^2$ bounds for the various geometric quantities to $L^2$ bounds for the components of $h$. As we saw in section \ref{section expressions for the inhomogeneous terms after commuting}, these are the kinds of terms which appear when commuting systems wave equations with the operators $\mathscr{Y}^n$. 

In this section, we will return to the system of wave equations \eqref{equation system of wave equations} and finally produce bounds on the inhomogenous terms after commuting. In other words, we will start with a known system of wave equations (satisfying the weak null condition) and produce bounds (in $L^2$) for the inhomogeneity after commuting $n$ times with the commutation operators $\mathscr{Y}$, where $n$ ranges from $0$ to $N_2$. Note that we have already produced such a bound in the case $n = 0$ in section \ref{section bounds on the inhomogeneous terms before commuting}.

\begin{proposition}[$L^2$ bounds for the inhomogeneous terms after commuting]
	\label{proposition L2 bounds for F after commuting}
	Let $\phi_{(A)}$ be a set of scalar fields satisfying the equations
	\begin{equation*}
	\begin{split}
	\tilde{\Box}_g\phi_{(A)} &= F_{(A,0)}
	\\
	F_{(A,0)} &= F_{(A,0)}^{(0)} + \left(F_{(A,0)}^{(BC)}\right)^{\mu\nu}(\partial_\mu \phi_B)(\partial_\nu \phi_C) + \mathcal{O}\left(\phi (\partial\phi)^2\right)
	\end{split}
	\end{equation*}
	where we further decompose
	\begin{equation*}
	\begin{split}
	F_{(A,0)}^{(0)}
	&=
	F_{(A,0,1)}^{(0)}
	+ F_{(A,0,2)}^{(0)}
	+ F_{(A,0,3)}^{(0)}
	\\
	&=
	F_{(A,0,4)}^{(0)}
	+ F_{(A,0,5)}^{(0)}
	+ F_{(A,0,6)}^{(0)}
	\end{split}	
	\end{equation*}
	and we define (schematically)
	\begin{equation*}
	F^{(0)}_{(A,n)} := \mathscr{Y}^n F^{(0)}_{(A,0)}
	\end{equation*}
	
	We require that the tensor fields $F_{(A,0)}^{(BC)}$ have \emph{constant rectangular components}\footnote{Note that this condition can be weakened: see the footnotes to proposition \ref{proposition L2 bounds for F before commuting}}. Also, we suppose that they satisfy the structural equations
	\begin{equation*}
	\begin{split}
	\left(F_{(A,0)}^{(BC)}\right)^{\mu\nu} 
	&= 	\left(F_{(A,0)}^{(CB)}\right)^{\mu\nu}
	\\
	\left(F_{(A,0)}^{(BC)}\right)_{LL}
	&= 0 \quad \text{if}\quad \phi_{(A)} \in \Phi_{[0]}
	\\
	\left(F_{(A,0)}^{(BC)}\right)_{LL}
	&= 0 \quad \text{if}\quad \phi_{(A)} \in \Phi_{[n]} \text{\, and either } \begin{cases}
	\phi_{(B)} \in \Phi_{[n+1]} \\
	\phi_{(B)} \in \Phi_{[n]} \text{\, and \,} \phi_{(C)} \in \Phi_{[m]} \text{ , } m \geq 1
	\end{cases}
	\end{split}
	\end{equation*}
	
	Suppose moreover that the terms $F^{(0)}_{(A,n)}$ satisfy the following conditions: if $\phi_{(A)} \in \Phi_{[0]}$, then 
	\begin{equation*}
	\begin{split}
	&\int_{\mathcal{M}_{\tau_0}^\tau} \epsilon^{-1} \bigg(
		(1+r)^{1-C_{[0,0]}\epsilon}|F^{(0)}_{(A,0)}|^2
		+ (1+r)^{\frac{1}{2}\delta}(1+\tau)^{1+\delta} |F^{(0)}_{(A,0,1)}|^2
		+ (1+r)^{1-3\delta}(1+\tau)^{2\beta} |F^{(0)}_{(A,0,2)}|^2
		\\
		&\phantom{\int_{\mathcal{M}_{\tau_0}^\tau} \epsilon^{-1} \bigg(}
		+ (1+r)^{1+\frac{1}{2}\delta} |F^{(0)}_{(A,0,3)}|^2
	\bigg)\dVol_g
	\lesssim
	\frac{1}{C_{[0,0]}} \epsilon^{2(N_2 + 2)}(1+\tau)^{-1}
	\\
	\\
	&\int_{\mathcal{M}_{\tau_0}^\tau \cap \{r \geq r_0\}} \epsilon^{-1}\bigg(
		r^{1-C_{[0,0]}\epsilon} (1+\tau)^{1+\delta} |F^{(0)}_{(A,0,4)}|^2
		+ r^{2-C_{[0,0]}\epsilon-2\delta} (1+\tau)^{2\beta} |F^{(0)}_{(A,0,5)}|^2
		\\
		&\phantom{\int_{\mathcal{M}_{\tau_0}^\tau \cap \{r \geq r_0\}} \epsilon^{-1}\bigg(}
		+ r^{2-C_{[0.0]}\epsilon} |F^{(0)}_{(A,0,6)}|^2
	\bigg)\dVol_g
	\lesssim
	\frac{1}{C_{[0,0]}} \epsilon^{2(N_2 + 2)}
	\end{split}
	\end{equation*}
	and in general, if $\phi_{(A)} \in \Phi_{[m]}$, then for $0 \leq n \leq N_2$
	\begin{equation*}
	\begin{split}
	&\int_{\mathcal{M}_{\tau_0}^\tau} \epsilon^{-1} \bigg(
		(1+r)^{1-C_{[n,m]}\epsilon}|F^{(0)}_{(A,n)}|^2
		+ (1+r)^{\frac{1}{2}\delta}(1+\tau)^{1+\delta} |F^{(0)}_{(A,n,1)}|^2
		+ (1+r)^{1-3\delta}(1+\tau)^{2\beta} |F^{(0)}_{(A,n,2)}|^2
		\\
		&\phantom{\int_{\mathcal{M}_{\tau_0}^\tau} \epsilon^{-1} \bigg(}
		+ (1+r)^{1+\frac{1}{2}\delta} |F^{(0)}_{(A,n,3)}|^2
	\bigg)\dVol_g
	\lesssim
	\frac{1}{C_{[n,m]}} \epsilon^{2(N_2 - n + 2)}(1+\tau)^{-1+C_{(n,m)}\delta}
	\\
	\\
	&\int_{\mathcal{M}_{\tau_0}^\tau \cap \{r \geq r_0\}} \epsilon^{-1}\bigg(
		r^{1-C_{[n,m]}\epsilon} (1+\tau)^{1+\delta} |F^{(0)}_{(A,n,4)}|^2
		+ r^{2-C_{[n,m]}\epsilon-2\delta} (1+\tau)^{2\beta} |F^{(0)}_{(A,n,5)}|^2
		\\
		&\phantom{\int_{\mathcal{M}_{\tau_0}^\tau \cap \{r \geq r_0\}} \epsilon^{-1}\bigg(}
		+ r^{2-C_{[n,m]}\epsilon} |F^{(0)}_{(A,n,6)}|^2
	\bigg)\dVol_g
	\lesssim
	\frac{1}{C_{[n,m]}} \epsilon^{2(N_2 - n + 2)}(1+\tau)^{C_{(n,m)}\delta}
	\end{split}
	\end{equation*}

	Furthermore, suppose that both the pointwise bounds \emph{and} the $L^2$ bounds of chapter \ref{chapter bootstrap} hold.
	
	Define $F_{(A,n)}$ as follows: if $\mathscr{Y}^n$ contains no factors of the operator $r\slashed{\D}_L$, then we define
	\begin{equation*}
	\tilde{\slashed{\Box}}_g \mathscr{Y}^n \phi_{(A)} = F_{(A,n)}
	\end{equation*}
	otherwise, if $\mathscr{Y}^n$ contains $k$ factors of the operator $r\slashed{\D}_L$, $k \geq 1$, then we define
	\begin{equation*}
	\tilde{\slashed{\Box}}_g \mathscr{Y}^n \phi_{(A)} 
	- k\slashed{\Delta} \mathscr{Y}^{n-1} \phi_{(A)}
	- (2^k - 1)r^{-1} \slashed{\D}_L (r\slashed{\D}_L \mathscr{Y}^{n-1}\phi)
	- (2^k - 1)r^{-1} \slashed{\D}_L (\mathscr{Y}^{n-1} \phi)
	= F_{(A,n)}
	\end{equation*}

	Then, for all sufficiently small $\epsilon$, we can decompose $F_{(A,n)}$ as
	\begin{equation*}
	\begin{split}
	F_{(A,n)}
	&=
	F_{(A,n,1)} + F_{(A,n,2)} + F_{(A,n,3)}
	\\
	&=
	F_{(A,n,4)} + F_{(A,n,5)} + F_{(A,n,6)}
	\end{split}	
	\end{equation*}
	where, if $\phi_{(A)} \in \Phi_{[m]}$, then we have
	\begin{equation*}
	\begin{split}
	&\int_{\mathcal{M}_{\tau_0}^\tau} \epsilon^{-1} \bigg(
		(1+r)^{1-C_{[n,m]}\epsilon}|F_{(A,n)}|^2
		+ (1+r)^{\frac{1}{2}\delta}(1+\tau)^{1+\delta} |F_{(A,n,1)}|^2
		+ (1+r)^{1-3\delta}(1+\tau)^{2\beta} |F_{(A,n,2)}|^2
		\\
		&\phantom{\int_{\mathcal{M}_{\tau_0}^\tau} \epsilon^{-1} \bigg(}
		+ (1+r)^{1+\frac{1}{2}\delta} |F_{(A,n,3)}|^2
	\bigg)\dVol_g
	\lesssim
	\left( \frac{1}{C_{[n,m]}} + \frac{\epsilon^2}{\delta^6} \right) \epsilon^{2(N_2 + 2 - n)} (1+\tau)^{-1+C_{(n)}\delta}
	\\
	\\
	&\int_{\mathcal{M}_{\tau_0}^\tau \cap \{r \geq r_0\}} \epsilon^{-1}\bigg(
		r^{1-C_{[n,m]}\epsilon} (1+\tau)^{1+\delta} |F_{(A,n,4)}|^2
		+ r^{2-C_{[n,m]}\epsilon-2\delta} (1+\tau)^{2\beta} |F_{(A,n,5)}|^2
		\\
		&\phantom{\int_{\mathcal{M}_{\tau_0}^\tau \cap \{r \geq r_0\}} \epsilon^{-1}\bigg(}
		+ r^{2-C_{[n,m]}\epsilon} |F_{(A,n,6)}|^2
	\bigg)\dVol_g
	\lesssim
	\left( \frac{1}{C_{[n,m]}} + \frac{\epsilon^2}{\delta^6} \right) \epsilon^{2(N_2 + 2 - n)} (1+\tau)^{C_{(n)}\delta}
	\end{split}
	\end{equation*}
	
	Finally, we define $F_{(A,N_2(j))}$ as follows: if $\mathscr{Y}^{N_2 - j}$ contains $k$ factors of the operator $r\slashed{\D}_L$, $k \geq 1$, then we define
	\begin{equation*}
	\begin{split}
	&\tilde{\slashed{\Box}}_g (\slashed{\D}_T)^j \mathscr{Y}^{N_2 - j} \phi_{(A)} 
	- k\slashed{\Delta} (\slashed{\D}_T)^j \mathscr{Y}^{N_2 - 1 - j} \phi_{(A)}
	- (2^k - 1)r^{-1} \slashed{\D}_L (r\slashed{\D}_L (\slashed{\D}_T)^j \mathscr{Y}^{N_2 - 1 - j}\phi)
	\\
	&
	- (2^k - 1)r^{-1} \slashed{\D}_L (\slashed{\D}_T)^j \mathscr{Y}^{N_2 - 1 - j} \phi)
	= F_{(A,N_2(j))}
	\end{split}
	\end{equation*}
	then we have
	\begin{equation*}
	\begin{split}
	&\int_{\mathcal{M}_{\tau_0}^\tau} \epsilon^{-1} \bigg(
		(1+r)^{1-C_{[N_2(j),m]}\epsilon}|F_{(A,N_2(j))}|^2
		+ (1+r)^{\frac{1}{2}\delta}(1+\tau)^{1+\delta} |F_{(A,N_2(j),1)}|^2
		\\
		&\phantom{\int_{\mathcal{M}_{\tau_0}^\tau} \epsilon^{-1} \bigg(}
		+ (1+r)^{1-3\delta}(1+\tau)^{2\beta} |F_{(A,N_2(j),2)}|^2
		+ (1+r)^{1+\frac{1}{2}\delta} |F_{(A,N_2(j),3)}|^2
	\bigg)\dVol_g
	\\
	&\lesssim
	\left( \frac{1}{C_{[N_2(j+1),m]}} + \frac{\epsilon^2}{\delta^6} \right) \epsilon^4 (1+\tau)^{-1+C_{(N_2(j))}\delta}
	\\
	\\
	&\int_{\mathcal{M}_{\tau_0}^\tau \cap \{r \geq r_0\}} \epsilon^{-1}\bigg(
		r^{1-C_{[N_2(j),m]}\epsilon} (1+\tau)^{1+\delta} |F_{(A,[N_2(j),4)}|^2
		+ r^{2-C_{[N_2(j),m]}\epsilon-2\delta} (1+\tau)^{2\beta} |F_{(A,[N_2(j),5)}|^2
		\\
		&\phantom{\int_{\mathcal{M}_{\tau_0}^\tau \cap \{r \geq r_0\}} \epsilon^{-1}\bigg(}
		+ r^{2-C_{[N_2(j),m]}\epsilon} |F_{(A,[N_2(j),6)}|^2
	\bigg)\dVol_g
	\lesssim
	\left( \frac{1}{C_{[N_2(j+1),m]}} + \frac{\epsilon^2}{\delta^6} \right) \epsilon^4 (1+\tau)^{C_{(N_2(j))}\delta}
	\end{split}
	\end{equation*}
	and
	\begin{equation*}
	\begin{split}
	&\int_{\mathcal{M}_{\tau_0}^\tau} \epsilon^{-1} \bigg(
		(1+r)^{1-C_{[N_2(j),m]}\epsilon}|F_{(A,N_2(N_2))}|^2
		+ (1+r)^{\frac{1}{2}\delta}(1+\tau)^{1+\delta} |F_{(A,N_2(N_2),1)}|^2
		\\
		&\phantom{\int_{\mathcal{M}_{\tau_0}^\tau} \epsilon^{-1} \bigg(}
		+ (1+r)^{1-3\delta}(1+\tau)^{2\beta} |F_{(A,N_2(N_2),2)}|^2
		+ (1+r)^{1+\frac{1}{2}\delta} |F_{(A,N_2(N_2),3)}|^2
	\bigg)\dVol_g
	\\
	&\lesssim
	\left( \frac{1}{C_{[N_2(N_2),m]}} + \frac{\epsilon^2}{\delta^6} \right) \epsilon^4 (1+\tau)^{-1+C_{(N_2(N_2))}\delta}
	\\
	\\
	&\int_{\mathcal{M}_{\tau_0}^\tau \cap \{r \geq r_0\}} \epsilon^{-1}\bigg(
		r^{1-C_{[N_2(N_2),m]}\epsilon} (1+\tau)^{1+\delta} |F_{(A,[N_2(N_2),4)}|^2
		\\
		&\phantom{\int_{\mathcal{M}_{\tau_0}^\tau \cap \{r \geq r_0\}} \epsilon^{-1}\bigg(}
		+ r^{2-C_{[N_2(N_2),m]}\epsilon-2\delta} (1+\tau)^{2\beta} |F_{(A,[N_2(N_2),5)}|^2
		+ r^{2-C_{[N_2(N_2),m]}\epsilon} |F_{(A,[N_2(N_2),6)}|^2
	\bigg)\dVol_g
	\\
	&\lesssim
	\left( \frac{1}{C_{[N_2(N_2),m]}} + \frac{\epsilon^2}{\delta^6} \right) \epsilon^4 (1+\tau)^{C_{(N_2(N_2))}\delta}
	\end{split}
	\end{equation*}
	
\end{proposition}

\begin{proof}
	Consider the field $\phi_{(A)} \in \Phi_{[m]}$, which satisfies the equation $\tilde{\Box}_g \phi_{(A)} = F_{(A,0)}$. We claim that $\mathscr{Y}^n F_{(A,0)}$ satisfies the following bound, if $n \leq 2N_1$:
	
	\begin{equation*}
	\begin{split}
	|\mathscr{Y}^n F_{(A,0)}|
	&\lesssim
	\epsilon (1+r)^{-1} |\slashed{\D} \mathscr{Y}^n \phi_{[m]}|
	+ \epsilon (1+r)^{-1+C_{(0,m)}\epsilon} |\slashed{\D} \mathscr{Y}^n \phi_{[0]}|
	+ \epsilon (1+r)^{-1+C_{(0,m-1)}\epsilon} |\slashed{\D} \mathscr{Y}^n \phi_{[m-1]}|
	\\
	&\phantom{\lesssim}
	+ \epsilon (1+r)^{-1-\delta} |\slashed{\D} \mathscr{Y}^n \phi|
	+ \epsilon (1+r)^{-1+C_{[n,N]}} |\overline{\slashed{\D}} \mathscr{Y}^n \phi|
	+ \epsilon (1+r)^{-2+\delta} |\mathscr{Y}^n \phi|
	\\
	&\phantom{\lesssim}
	+ \sum_{j \leq n/2} \epsilon (1+r)^{-1+C_{(j,N)}} |\bm{\Gamma}^{(j-n)}_{-1+C_{(j-1)}\epsilon}|
	+ |F^{(0)}_{(A,n)}|
	\end{split}
	\end{equation*}
	where we recall that $\phi_{[m]}$ stands for any field in the set $\Phi_{[m]}$, and we also recall that we write $\phi$ to stand for \emph{any} of the fields $\phi \in \Phi_{[N]}$, where $\Phi_{[N]}$ represents the ``highest'' level of the hierarchy.
	
	First, for $n = 0$, by the conditions on $F$ we have, schematically,
	\begin{equation}
	\begin{split}
	\label{equation L2 put together internal}
	F_{(A,0)} 
	&=
	F^{(0)}_{(A,0)}
	+ (F_{(A)}^{(BC)}) (\partial \phi_{[m]})(\partial \phi_{[0]})
	+ (F_{(A)}^{(BC)}) (\partial \phi_{[m-1]})(\partial \phi_{[m-1]})
	+ (F_{(A)}^{(BC)}) (\partial \phi)(\bar{\partial} \phi)
	\\
	&\phantom{=}
	+ (\partial \phi)^2 (\phi)
	+ \text{cubic terms}
	\end{split}
	\end{equation}
	where the cubic terms are easier to bound than those which we have displayed explicitly.
	
	Therefore, using the pointwise bounds on the fields together with the fact that the frame components $F^{(BC)}_{(A)}$ are uniformly bounded, we find that
	\begin{equation*}
	\begin{split}
	|F_{(A,0)}| 
	&\lesssim
	|F^{(0)}_{(A,0)}|
	+ \epsilon r^{-1} |\partial \phi_{[m]}|
	+ \epsilon r^{-1 + C_{(0,m)}\epsilon} |\partial \phi_{[0]}|
	+ \epsilon r^{-1 + C_{(0,m-1)}\epsilon} |\partial \phi_{[m-1]}|
	+ \epsilon r^{-1-\delta} |\partial \phi|
	\\
	&\phantom{\lesssim}
	+ \epsilon r^{-1 + C_{(0,N)}\epsilon} |\bar{\partial} \phi|
	+ \epsilon^2 r^{-2 + C_{(0,N)}\epsilon} |\phi|
	\end{split}
	\end{equation*}
	
	Now, we return to equation \eqref{equation L2 put together internal} and apply the operator $\mathscr{Y}^n$. Note that, since the \emph{rectangular} components of $F_{(A)}^{(BC)}$ are constants, we have
	\begin{equation*}
	\mathscr{Y}^{n} \left((F_{(A)}^{(BC)})_{(\text{frame})} \right)
	\lesssim
	\sum_{j+k \leq n} |\mathscr{Y}^j X_{(\text{frame})}| |\mathscr{Y}^k X_{(\text{frame})}|
	\end{equation*}
	Hence we obtain, schematically,
	\begin{equation*}
	\begin{split}
	\mathscr{Y}^n F_{(A, 0)}
	&=
	F^{(0)}_{(A, n)}
	+ (F_{(A)}^{(BC)}) (\partial \phi_{[0]}) (\slashed{\D} \mathscr{Y}^n \phi_{[m]})
	+ (F_{(A)}^{(BC)}) (\partial \phi_{[m]}) (\slashed{\D} \mathscr{Y}^n \phi_{[0]})
	\\
	&\phantom{=}
	+ (F_{(A)}^{(BC)}) (\partial \phi_{[m-1]}) (\slashed{\D} \mathscr{Y}^n \phi_{[m-1]})
	+ (F_{(A)}^{(BC)}) (\partial \phi) (\overline{\slashed{\D}} \mathscr{Y}^n \phi)
	\\
	&\phantom{=}
	+ (F_{(A)}^{(BC)}) (\bar{\partial} \phi) (\slashed{\D} \mathscr{Y}^n \phi)
	+ (\partial \phi) (\phi) (\slashed{\D} \mathscr{Y}^n \phi)
	+ (\partial \phi)^2 (\mathscr{Y}^n \phi)
	\\
	&\phantom{=}
	+ (\mathscr{Y}^n X_{(\text{frame, small})}) (\partial \phi)^2
	+ \sum_{\substack{j+k \leq n \\ j,k \leq n-1}} \left(\bm{\Gamma}^{(j)}_{(-1+C_{(j)}\epsilon)} \right) \left(\bm{\Gamma}^{(k)}_{(-1+C_{(k)}\epsilon)} \right)
	\end{split}
	\end{equation*}

	Before we go any further, we should note that \emph{it is not the case that} $\mathscr{Y}^n F_{(A,0)} = F_{(A,n)}$. Indeed, this would only be the case if the operators $\mathscr{Y}$ commuted with the (modified) wave operator $\tilde{\slashed{\Box}}$, which they do not. Instead, recall proposition \ref{proposition inhomogeneous terms after commuting n times with Y}. Combining this with the equation above, we have
	\begin{equation}
	\label{equation L2 put together internal 2}
	\begin{split}
	F_{(A, n)}
	&=
	F^{(0)}_{(A, n)}
	+ (F_{(A)}^{(BC)}) (\partial \phi_{[0]}) (\slashed{\D} \mathscr{Y}^n \phi_{[m]})
	+ (F_{(A)}^{(BC)}) (\partial \phi_{[m]}) (\slashed{\D} \mathscr{Y}^n \phi_{[0]})
	\\
	&\phantom{=}
	+ (F_{(A)}^{(BC)}) (\partial \phi_{[m-1]}) (\slashed{\D} \mathscr{Y}^n \phi_{[m-1]})
	+ (F_{(A)}^{(BC)}) (\partial \phi) (\overline{\slashed{\D}} \mathscr{Y}^n \phi)
	\\
	&\phantom{=}
	+ (F_{(A)}^{(BC)}) (\bar{\partial} \phi) (\slashed{\D} \mathscr{Y}^n \phi)
	+ (\partial \phi) (\phi) (\slashed{\D} \mathscr{Y}^n \phi)
	+ (\partial \phi)^2 (\mathscr{Y}^n \phi)
	\\
	&\phantom{=}
	+ (\mathscr{Y}^n X_{(\text{frame, small})}) (\partial \phi)^2
	+ \sum_{\substack{j+k \leq n \\ j,k \leq n-1}} \left(\bm{\Gamma}^{(j)}_{(-1+C_{(j)}\epsilon)} \right) \left(\bm{\Gamma}^{(k)}_{(-1+C_{(k)}\epsilon)} \right)
	\\
	&\phantom{=}
	+ r^{-1} \overline{\slashed{\D}} (\mathscr{Y}^{\leq n-1} \phi_{(A)})
	+ \bm{\Gamma}^{(1)}_{(-1+C_{(1)}\epsilon)} \overline{\slashed{\D}}(\mathscr{Y}^n \phi_{(A)})
	+ \bm{\Gamma}^{(0)}_{(-1)} (\slashed{\D} \mathscr{Y}^n \phi_{(A)})
	\\
	\\
	&\phantom{=}
	+ \begin{pmatrix} r^{-1} (\partial \phi_{(A)}) \\ (\bar{\partial} \phi_{(A)}) \\ r^{-1} \mathscr{Y}\phi_{(A)} \end{pmatrix} \left( \slashed{\D} \mathscr{Y}^{n} h \right)_{(\text{frame})}
	+ (\partial \phi_{(A)}) (\slashed{\D} \mathscr{Y}^n h)_{LL}
	+ (\partial \phi_{(A)}) \left( \overline{\slashed{\D}} \mathscr{Y}^{n} h \right)_{(\text{frame})}
	\\
	&\phantom{=}
	+  \begin{pmatrix} (\partial \phi_{(A)}) \\ r (\bar{\partial} \phi_{(A)}) \\ \mathscr{Y} \phi_{(A)} \end{pmatrix} \left( \tilde{\slashed{\Box}}_g \mathscr{Y}^{n-1} h \right)_{(\text{frame})}
	+ (\partial \phi_{(A)}) \mathscr{Y}^n \tr_{\slashed{g}}\chi_{(\text{small})}
	+ (\bar{\partial} \phi_{(A)}) (r\slashed{\nabla}^2 \mathscr{Z}^{n-1} \log \mu)
	\\
	&\phantom{=}
	+ \begin{pmatrix} r^{-1} (\overline{\slashed{\D}} \mathscr{Y} \phi_{(A)}) \\ \bm{\Gamma}^{(1)}_{(-2+ C_{(1)}\epsilon)} (\mathscr{Y}\phi) \end{pmatrix} (\mathscr{Y}^n \log \mu)
	+ \sum_{\substack{ j+k \leq n+1 \\ j \leq n-1 \\ k \leq n-1}} \bm{\Gamma}^{(j)}_{(-1 + C_{(j)}\epsilon)} (\slashed{\D} \mathscr{Y}^k \phi_{(A)})
	\\
	& \phantom{=}
	+ \sum_{\substack{ j+k \leq n \\ j \leq n-1 \\ k \leq n-1}} r\tilde{\slashed{\Box}}_g (\mathscr{Z}^j h)_{(\text{frame})}  (\slashed{\D} \mathscr{Y}^k \phi_{(A)})
	+ \sum_{\substack{j+k \leq n+1 \\ j \leq n-1 \\ k \leq n-1}} \bm{\Gamma}^{(j)}_{(-2 + 2C_{(j)}\epsilon)} (\mathscr{Y}^k \phi_{(A)})
	\end{split}
	\end{equation}

	The first bound we will consider is the one which does not allow for any decomposition of $F_{(A,n)}$: that is, we will produce a bound on
	\begin{equation*}
	\int_{\mathcal{M}_{\tau_0}^\tau \cap \{ r \geq r_0\}} \epsilon^{-1} (1+r)^{1-C_{[n,m]}\epsilon} |F_{(A,n)}|^2 \dVol_g
	\end{equation*}
	For this purpose, we use pointwise bounds with the maximum decay in $r$, for example, we use the pointwise bounds
	\begin{equation*}
	|\bm{\Gamma}^{(j)}_{(-1+C_{(j)}\epsilon)}| \lesssim \epsilon (1+r)^{(-1+C_{(j)}\epsilon)}
	\end{equation*}
	which holds for $j \leq N_1$.	Since $2N_1 \leq N_2$, we can use the pointwise bounds on at least one of each of the factors in every quadratic term in equation \eqref{equation L2 put together internal 2} (in particular, we use the improved pointwise bounds on $\partial \mathscr{Y}^{\leq 1} \phi_{(A)}$ and $\bar{\partial} \mathscr{Y}^{\leq 1} \phi_{(A)}$), to obtain the bound
	\begin{equation*}
	\begin{split}
	|F_{(A,n)}|
	&\lesssim
	|F^{(0)}_{(A,n)}|
	+ \epsilon (1+r)^{-1} |\slashed{\D} \mathscr{Y}^n \phi_{[m]}|
	+ \epsilon (1+r)^{-1+C_{(0,m)}\epsilon} |\slashed{\D} \mathscr{Y}^n \phi_{[0]}|
	\\
	&\phantom{\lesssim}
	+ \epsilon (1+r)^{-1+C_{(0,m-1)}\epsilon} |\slashed{\D} \mathscr{Y}^n \phi_{[m-1]}|
	+ \epsilon (1+r)^{-1+C_{(0,N)}\epsilon} |\overline{\slashed{\D}} \mathscr{Y}^n \phi|
	+ \epsilon (1+r)^{-1-\delta} |\slashed{\D} \mathscr{Y}^n \phi|
	\\
	&\phantom{\lesssim}
	+ \epsilon^2 (1+r)^{-2+2C_{(0,N)}\epsilon} |\mathscr{Y}^n \phi|
	+ \epsilon^2 (1+r)^{-2+2C_{(0,N)}\epsilon} |\mathscr{Y}^n X_{(\text{frame, small})}|
	\\
	&\phantom{\lesssim}
	+ \epsilon (1+r)^{-1+C_{(N_1)}\epsilon} |\bm{\Gamma}^{(n-1)}_{(-1+C_{(n-1)}\epsilon)}|
	+ r^{-1} |\overline{\slashed{\D}} \mathscr{Y}^{\leq n-1} \phi_{(A)}|
	\\
	&\phantom{\lesssim}
	+ \epsilon^3 (1+r)^{-1+C_{(1)}\epsilon} |\overline{\slashed{\D}} \mathscr{Y}^{n} \phi_{(A)}|
	+ \epsilon^3 (1+r)^{-1} |\slashed{\D} \mathscr{Y}^{n} \phi_{(A)}|
	+ \epsilon^3 (1+r)^{-1-\delta} |\slashed{\D} \mathscr{Y}^n h|_{(\text{frame})}
	\\
	&\phantom{\lesssim}
	+ \epsilon^3 (1+r)^{-1 + C_{(0,m)}\epsilon} |\slashed{\D} \mathscr{Y}^n h|_{LL}
	+ \epsilon^3 (1+r)^{-1 + C_{(0,m)}\epsilon} |\overline{\slashed{\D}} \mathscr{Y}^n h|_{(\text{frame})}
	\\
	&\phantom{\lesssim}
	+ \epsilon^3 (1+r)^{-\frac{1}{2} + \delta} |\tilde{\slashed{\Box}}_g \mathscr{Y}^{n-1} h|_{(\text{frame})}
	+ \epsilon^5 (1+r)^{-1 + C_{(0,m)}\epsilon} |\mathscr{Y}^n \tr_{\slashed{g}}\chi_{(\text{small})}|
	\\
	&\phantom{\lesssim}
	+ |\mathscr{Y} \phi_{(A)}| |\slashed{\nabla}^2 \mathscr{Z}^{n-1} \log \mu|
	+ \epsilon^5 (1+r)^{-2 - \delta} |\mathscr{Y}^n \log \mu|
	+ \sum_{\substack{j+k \leq n \\ j,k\leq n-1}} r |\tilde{\slashed{\Box}}_g (\mathscr{Z}^j h)_{(\text{frame})}| |\slashed{\D} \mathscr{Y}^k \phi_{(A)}|
	\\
	&\phantom{\lesssim}
	+ \epsilon^3 (1+r)^{-2+2C_{(N)}\epsilon} |\mathscr{Y}^{\leq n-1} \phi_{(A)}|	
	\end{split}
	\end{equation*}
	
	We need to further estimate a couple of these terms. First, recall that we write
	\begin{equation*}
	\left(\tilde{\slashed{\Box}}_g (\mathscr{Z}^j h)\right)_{(\text{frame})}
	= (F^{(j)})_{(\text{frame})}
	\end{equation*}
	and, for $n \leq N_1$ we have
	\begin{equation*}
	|F^{(j)}|_{(\text{frame})} \lesssim \epsilon^2 (1+r)^{-2+2C_{(j)}\epsilon}
	\end{equation*}
	
	Similarly, we have, schematically,
	\begin{equation*}
	\begin{split}
	&(\tilde{\slashed{\Box}}_g \mathscr{Y}^{n-1} h)_{(\text{frame})}
	- (\slashed{\Delta} \mathscr{Y}^{n-2} h)_{(\text{frame})}
	- (2^k -1) r^{-1} \slashed{\D}_L (\mathscr{Y}^{n-1} h)_{(\text{frame})}
	- (2^k -1) r^{-1} \slashed{\D}_L (\mathscr{Y}^{n-2} h)_{(\text{frame})}
	\\
	&= (F^{(n-1)})_{(\text{frame})}
	\end{split}
	\end{equation*}
	where $k$ is the number of times the operator $r\slashed{\D}_L$ appears in the expansion of $\mathscr{Y}^{n-1}$ (c.f.\ proposition \ref{proposition inhomogeneous terms after commuting n times with Y}). Hence, we have
	\begin{equation*}
	\begin{split}
	&|\tilde{\slashed{\Box}}_g \mathscr{Y}^{n-1} h|_{(\text{frame})}
	\lesssim 
	r^{-1} |\overline{\slashed{\D}} \mathscr{Y}^{n-1} h|_{(\text{frame})}
	+ r^{-1} |\overline{\slashed{\D}} \mathscr{Y}^{n-2} h|_{(\text{frame})}
	+ |F^{(n-1)}|_{(\text{frame})}
	\end{split}
	\end{equation*}

	Next, we bound the term $|\mathscr{Y} \phi_{(A)}| |\slashed{\nabla}^2 \mathscr{Z}^{n-1} \log \mu|$. We do this in two different ways, depending on whether $n = N_2$ or $n \leq N_2 - 1$. If $n \leq N_2 - 1$ then we can use the bound
	\begin{equation*}
	|\mathscr{Y} \phi_{(A)}| \lesssim \epsilon^5 (1+r)^{-\frac{1}{2} + \delta} (1+\tau)^{-\beta}
	\end{equation*}
	to bound this term as
	\begin{equation*}
	|\mathscr{Y} \phi_{(A)}| |\slashed{\nabla}^2 \mathscr{Z}^{n-1} \log \mu|
	\lesssim
	\epsilon^5 (1+r)^{-\frac{5}{2} + \delta} (1+\tau)^{-\beta} |\mathscr{Y}^{n+1} \log \mu|^2
	\end{equation*}
	On the other hand, in the case $n = N_2$ we bound this term using the bound
	\begin{equation*}
	|\mathscr{Y} \phi_{(A)}| \lesssim \epsilon^5 (1+r)^{-1 + C_{[N_1]}\epsilon} (1+\tau)^{C_{(N_1)}\delta}
	\end{equation*}	
	to bound it as
	\begin{equation*}
	|\mathscr{Y} \phi_{(A)}| |\slashed{\nabla}^2 \mathscr{Y}^{N_2-1} \log \mu|
	\lesssim
	\epsilon^5 (1+r)^{-1+C_{[N_1]}\epsilon} (1+\tau)^{C_{(N_1)}\delta} |\slashed{\nabla}^2 \mathscr{Y}^{N_2-1} \log \mu|
	\end{equation*}
	
	Putting these calculations together, we have, for $n \leq N_2 - 1$,
	\begin{equation*}
	\begin{split}
		&\int_{\mathcal{M}_\tau^{\tau_1}} \epsilon^{-1} (1+r)^{1-C_{[n,m]}\epsilon} |F_{(A,n)}|^2 \dVol_g
		\\
		&\lesssim
		\int_{\mathcal{M}_\tau^{\tau_1}} \bigg(
			\epsilon^{-1} (1+r)^{1-C_{[n,m]}\epsilon} |F^{(0)}_{(A,n)}|^2
			+ \epsilon (1+r)^{-1-C_{[n,m]}\epsilon} |\slashed{\D} \mathscr{Y}^n \phi_{[m]}|^2
			\\
			&\phantom{\lesssim \int_{\mathcal{M}_\tau^{\tau_1}} \bigg(}
			+ \epsilon (1+r)^{-1-(C_{[n,m]} - 2C_{(0,m)})\epsilon} |\slashed{\D} \mathscr{Y}^n \phi_{[0]}|^2
			+ \epsilon (1+r)^{-1-(C_{[n,m]} - 2C_{(0,m-1)})\epsilon} |\slashed{\D} \mathscr{Y}^n \phi_{[m-1]}|^2
			\\
			&\phantom{\lesssim \int_{\mathcal{M}_\tau^{\tau_1}} \bigg(}
			+ \epsilon (1+r)^{-1-(C_{[n,m]} - 2C_{(0,N)})\epsilon} |\overline{\slashed{\D}} \mathscr{Y}^n \phi|^2
			+ \epsilon (1+r)^{-1-2\delta} |\slashed{\D} \mathscr{Y}^n \phi|^2
			\\
			&\phantom{\lesssim \int_{\mathcal{M}_\tau^{\tau_1}} \bigg(}
			+ \epsilon^3 (1+r)^{-3-(C_{[n,m]} - 4C_{(0,N)})\epsilon} |\mathscr{Y}^n \phi|^2
			+ \epsilon^5 (1+r)^{-3-(C_{[n,m]} - 4C_{(0,N)})\epsilon} |\mathscr{Y}^n X_{(\text{frame, small})}|^2
			\\
			&\phantom{\lesssim \int_{\mathcal{M}_\tau^{\tau_1}} \bigg(}
			+ \epsilon (1+r)^{-1-(C_{[n,m]} - 2C_{(N_1)})\epsilon} |\bm{\Gamma}^{(n-1)}_{(-1+C_{(n-1)}\epsilon)}|^2
			+ \epsilon^{-1} r^{-1-C_{[n,m]}\epsilon} |\overline{\slashed{\D}} \mathscr{Y}^{\leq n-1} \phi_{(A)}|^2
			\\
			&\phantom{\lesssim \int_{\mathcal{M}_\tau^{\tau_1}} \bigg(}
			+ \epsilon (1+r)^{-1-(C_{[n,m]} - 2C_{(1)})\epsilon} |\overline{\slashed{\D}} \mathscr{Y}^{n} \phi_{(A)}|^2
			+ \epsilon^3 (1+r)^{-1-2\delta} |\slashed{\D} \mathscr{Y}^n h|^2_{(\text{frame})}
			\\
			&\phantom{\lesssim \int_{\mathcal{M}_\tau^{\tau_1}} \bigg(}
			+ \epsilon^3 (1+r)^{-1 -(C_{[n,m]} - 2C_{(0,m)})\epsilon} |\slashed{\D} \mathscr{Y}^n h|^2_{LL}
			+ \epsilon^3 (1+r)^{-1 -(C_{[n,m]} - 2C_{(0,m)})\epsilon} |\overline{\slashed{\D}} \mathscr{Y}^n h|^2_{(\text{frame})}
			\\
			&\phantom{\lesssim \int_{\mathcal{M}_\tau^{\tau_1}} \bigg(}
			+ \epsilon^3 (1+r)^{2\delta} |F^{(n-1)}|^2_{(\text{frame})}
			+ \epsilon^5 (1+r)^{-1 -(C_{[n,m]} - 2C_{(0,m)})\epsilon} |\mathscr{Y}^n \tr_{\slashed{g}}\chi_{(\text{small})}|^2
			\\
			&\phantom{\lesssim \int_{\mathcal{M}_\tau^{\tau_1}} \bigg(}
			+ \epsilon^5 (1+r)^{-4+2\delta} (1+\tau)^{-2\beta} |\mathscr{Y}^{n+1} \log \mu|^2
			+ \epsilon^5 (1+r)^{-3-2\delta} |\mathscr{Y}^n \log \mu|^2
			\\
			&\phantom{\lesssim \int_{\mathcal{M}_\tau^{\tau_1}} \bigg(}
			+ \epsilon^3 (1+r)^{1 - (C_{[n,m]} - 2C_{(N_1,m)})\epsilon} |F^{(\leq n-1)}|^2_{(\text{frame})}
			\\
			&\phantom{\lesssim \int_{\mathcal{M}_\tau^{\tau_1}} \bigg(}
			+ \epsilon (1+r)^{-3 - (C_{[n,m]} - 4C_{(N)})\epsilon} |\mathscr{Y}^{\leq n-1} \phi_{(A)}|^2
		\bigg) \dVol_g
	\end{split}
	\end{equation*}
	
	Using the conditions on the constants $C_{(n,m)}$, $C_{[n,m]}$, $C_{(n)}$ and $C_{[n]}$ (and re-ordering the terms) we find
	\begin{equation*}
	\begin{split}
	&\int_{\mathcal{M}_\tau^{\tau_1}} \epsilon^{-1} (1+r)^{1-C_{[n,m]}\epsilon} |F_{(A,n)}|^2 \dVol_g
	\\
	&\lesssim
	\int_{\mathcal{M}_\tau^{\tau_1}} \bigg( 
		\epsilon^{-1} (1+r)^{1-C_{[n,m]}\epsilon} |F^{(0)}_{(A,n)}|^2
		+ \epsilon (1+r)^{2\delta} |F^{(n-1)}|^2_{(\text{frame})}
		+ \epsilon (1+r)^{1 - C_{[n-1,m]}\epsilon} |F^{(\leq n-1)}|^2_{(\text{frame})}
		\\
		&\phantom{\lesssim \int_{\mathcal{M}_\tau^{\tau_1}} \bigg(}
		+ \epsilon (1+r)^{-1-C_{[n,m]}\epsilon} |\slashed{\D} \mathscr{Y}^n \phi_{[m]}|^2
		+ \epsilon (1+r)^{-1-C_{[n,0]}\epsilon} |\slashed{\D} \mathscr{Y}^n \phi_{[0]}|^2
		\\
		&\phantom{\lesssim \int_{\mathcal{M}_\tau^{\tau_1}} \bigg(}
		+ \epsilon (1+r)^{-1-C_{[n,m-1]}\epsilon} |\slashed{\D} \mathscr{Y}^n \phi_{[m-1]}|^2
		+ \epsilon (1+r)^{-1-2\delta} |\slashed{\D} \mathscr{Y}^n \phi|^2
		+ \epsilon (1+r)^{-1-C_{[n-1,0]}\epsilon} |\overline{\slashed{\D}} \mathscr{Y}^n \phi|^2
		\\
		&\phantom{\lesssim \int_{\mathcal{M}_\tau^{\tau_1}} \bigg(}
		+ \epsilon^{-1} r^{-1-C_{[n,m]}\epsilon} |\overline{\slashed{\D}} \mathscr{Y}^{\leq n-1} \phi_{(A)}|^2
		+ \epsilon^3 (1+r)^{-3-C_{[n-1,0]}\epsilon} |\mathscr{Y}^n \phi|^2
		\\
		&\phantom{\lesssim \int_{\mathcal{M}_\tau^{\tau_1}} \bigg(}
		+ \epsilon (1+r)^{-3 - C_{[n-1,N]}\epsilon} |\mathscr{Y}^{\leq n-1} \phi_{(A)}|^2
		+ \epsilon^5 (1+r)^{-3-\frac{1}{2}C_{[n,m]}\epsilon} |\mathscr{Y}^n X_{(\text{frame, small})}|^2
		\\
		&\phantom{\lesssim \int_{\mathcal{M}_\tau^{\tau_1}} \bigg(}
		+ \epsilon^5 (1+r)^{-1 -(C_{[n,m]} - 2C_{(0,m)})\epsilon} |\mathscr{Y}^n \tr_{\slashed{g}}\chi_{(\text{small})}|^2
		+ \epsilon^5 (1+r)^{-4+2\delta} (1+\tau)^{-2\beta} |\mathscr{Y}^{n+1} \log \mu|^2
		\\
		&\phantom{\lesssim \int_{\mathcal{M}_\tau^{\tau_1}} \bigg(}
		+ \epsilon^5 (1+r)^{-3-2\delta} |\mathscr{Y}^n \log \mu|^2
		+ \epsilon (1+r)^{-1-C_{[n-1,N]}\epsilon} |\bm{\Gamma}^{(n-1)}_{(-1+C_{(n-1)}\epsilon)}|^2
		\\
		&\phantom{\lesssim \int_{\mathcal{M}_\tau^{\tau_1}} \bigg(}
		+ \epsilon^3 (1+r)^{-1-2\delta} |\slashed{\D} \mathscr{Y}^n h|^2_{(\text{frame})}
		+ \epsilon^3 (1+r)^{-1 -(C_{[n,m]} - 2C_{(0,m)})\epsilon} |\slashed{\D} \mathscr{Y}^n h|^2_{LL}
		\\
		&\phantom{\lesssim \int_{\mathcal{M}_\tau^{\tau_1}} \bigg(}
		+ \epsilon^3 (1+r)^{-1 -C_{[n-1,0]}\epsilon} |\overline{\slashed{\D}} \mathscr{Y}^n h|^2_{(\text{frame})}
	\bigg) \dVol_g 
	\end{split}
	\end{equation*}

	Now, we can use proposition \ref{proposition L2 bounds rectangular} to bound the term involving $\mathscr{Y}^n X_{(\text{frame, small})}$, proposition \ref{proposition L2 tr chi low} to bound the term involving $\mathscr{Y}^n \tr_{\slashed{g}}\chi_{(\text{small})}$, proposition \ref{proposition L2 foliation density} to bound the terms involving $\mathscr{Y}^n \log \mu$ and $\mathscr{Y}^{n+1} \log \mu$, and finally the bootstrap bounds in equation \eqref{equation bootstrap L2 geometric} to bound the lower order terms. Using $\epsilon \ll C_{[n]}\delta$ for all $n$, this leads to the bound
	\begin{equation*}
	\begin{split}
	&\int_{\mathcal{M}_\tau^{\tau_1}} \epsilon^{-1} (1+r)^{1-C_{[n,m]}\epsilon} |F_{(A,n)}|^2 \dVol_g
	\\
	&\lesssim
	\epsilon \epsilon^{2(N_2 + 2 - n)} \left( (1+\tau)^{-1+C_{(n+1)}\delta - 2\beta} + (1+\tau)^{-1+C_{(n)}\delta} \right)
	\\
	&+ \int_{\mathcal{M}_\tau^{\tau_1}} \bigg( 
		\epsilon^{-1} (1+r)^{1-C_{[n,m]}\epsilon} |F^{(0)}_{(A,n)}|^2
		+ \epsilon (1+r)^{2\delta} |F^{(n-1)}|^2_{(\text{frame})}
		+ \epsilon (1+r)^{1 - C_{[n-1,m]}\epsilon} |F^{(\leq n-1)}|^2_{(\text{frame})}
		\\
		&\phantom{\lesssim \int_{\mathcal{M}_\tau^{\tau_1}} \bigg(}
		+ \epsilon (1+r)^{-1-C_{[n,m]}\epsilon} |\slashed{\D} \mathscr{Y}^n \phi_{[m]}|^2
		+ \epsilon (1+r)^{-1-C_{[n,0]}\epsilon} |\slashed{\D} \mathscr{Y}^n \phi_{[0]}|^2
		\\
		&\phantom{\lesssim \int_{\mathcal{M}_\tau^{\tau_1}} \bigg(}
		+ \epsilon (1+r)^{-1-C_{[n,m-1]}\epsilon} |\slashed{\D} \mathscr{Y}^n \phi_{[m-1]}|^2
		+ \epsilon (1+r)^{-1-2\delta} |\slashed{\D} \mathscr{Y}^n \phi|^2
		+ \epsilon (1+r)^{-1-C_{[n-1,0]}\epsilon} |\overline{\slashed{\D}} \mathscr{Y}^n \phi|^2
		\\
		&\phantom{\lesssim \int_{\mathcal{M}_\tau^{\tau_1}} \bigg(}
		+ \epsilon^{-1} r^{-1-C_{[n,m]}\epsilon} |\overline{\slashed{\D}} \mathscr{Y}^{\leq n-1} \phi_{(A)}|^2
		+ \epsilon^3 (1+r)^{-3-C_{[n-1,0]}\epsilon} |\mathscr{Y}^n \phi|^2
		\\
		&\phantom{\lesssim \int_{\mathcal{M}_\tau^{\tau_1}} \bigg(}
		+ \epsilon (1+r)^{-3 - C_{[n-1,N]}\epsilon} |\mathscr{Y}^{\leq n-1} \phi_{(A)}|^2
		+ \epsilon^5 (1+r)^{-1-\delta} (1+\tau)^{-2\beta} |\slashed{\D} \mathscr{Y}^{n+1} h|^2_{(\text{frame})}
		\\
		&\phantom{\lesssim \int_{\mathcal{M}_\tau^{\tau_1}} \bigg(}
		+ \epsilon^3 (1+r)^{-1-\delta} |\slashed{\D} \mathscr{Y}^{n} h|^2_{(\text{frame})}
		+ \epsilon^3 (1+r)^{-1 -(C_{[n,m]} - 2C_{(0,m)})\epsilon} |\slashed{\D} \mathscr{Y}^n h|^2_{LL}
		\\
		&\phantom{\lesssim \int_{\mathcal{M}_\tau^{\tau_1}} \bigg(}
		+ \epsilon^5 (1+r)^{-1+\frac{1}{2}\delta} (1+\tau)^{-2\beta} |\overline{\slashed{\D}} \mathscr{Y}^{n+1} h|^2_{(\text{frame})}
		+ \left( \epsilon^3 + \frac{\epsilon^5}{\delta^5} \right) (1+r)^{-1+\frac{1}{2}\delta} |\overline{\slashed{\D}} \mathscr{Y}^{n} h|^2_{(\text{frame})}
		\\
		&\phantom{\lesssim \int_{\mathcal{M}_\tau^{\tau_1}} \bigg(}
		+ \frac{\epsilon^5}{\delta^5} (1+r)^{-3-\delta} |\mathscr{Y}^{n+1} h|^2_{(\text{frame})}
		+ \epsilon^3 (1+r)^{-3-\delta} |\mathscr{Y}^{n} h|^2_{(\text{frame})}
	\bigg) \dVol_g 
	\end{split}
	\end{equation*}
	
	Finally, we want to apply the $L^2$ bootstrap bounds for derivatives of $h$ and $\phi$. We can consider two possible cases: either $\phi_{(A)} \in \Phi_{[0]}$ or $\phi_{(A)} \in \Phi_{[m]}$ for $m \geq 1$. In the first case, we have $C_{(0,0)} = 0$, so
	\begin{equation*}
	\begin{split}
	\int_{\mathcal{M}_\tau^{\tau_1}} \bigg( 
		\epsilon^3 (1+r)^{-1 -(C_{[n,m]} - 2C_{(0,m)})\epsilon} |\slashed{\D} \mathscr{Y}^n h|^2_{LL}
	\bigg) \dVol_g
	&\lesssim
	\int_{\mathcal{M}_\tau^{\tau_1}} \bigg( 
		\epsilon^3 (1+r)^{-1 -C_{[n,0]}\epsilon} |\slashed{\D} \mathscr{Y}^n h|^2_{LL}
	\bigg) \dVol_g
	\\
	&\lesssim
	\frac{1}{C_{[n,0]}} \epsilon^{2(N_2 + 2 - n)} (1+\tau)^{-1+C_{(n)}\epsilon}
	\end{split}	
	\end{equation*}
	where in the second line we have used the bootstrap assumptions. On the other hand, in the second case we have
	\begin{equation*}
	\begin{split}
	\int_{\mathcal{M}_\tau^{\tau_1}} \bigg( 
		\epsilon^3 (1+r)^{-1 -(C_{[n,m]} - 2C_{(0,m)})\epsilon} |\slashed{\D} \mathscr{Y}^n h|^2_{LL}
	\bigg) \dVol_g
	&\lesssim
	\int_{\mathcal{M}_\tau^{\tau_1}} \bigg( 
		\epsilon^3 (1+r)^{-1 -C_{[n,m-1]}\epsilon} |\slashed{\D} \mathscr{Y}^n h|^2_{LL}
	\bigg) \dVol_g
	\\
	&\lesssim
	\frac{1}{C_{[n,0]}} \epsilon^{2(N_2 + 2 - n)} (1+\tau)^{-1+C_{(n)}\epsilon}
	\end{split}	
	\end{equation*}
	so in both cases we are led to the same bound. Note that these are the only terms where we need to consider the change-of-frame from the rectangular frame to the null frame, each of which defines different sections of the vector bundle of symmetric matrices over $\mathcal{M}$. The reason for this is that, for the other terms, we can simply use
	\begin{equation*}
	|\partial h|_{(\text{frame})} \sim |\partial h|_{(\text{rect})} \sim |\partial \phi_{[N]}|
	\end{equation*}
	and similarly for higher order quantities. If we estimate these terms in this way, then we lose information about the hierarchy in the fields: we must estimate all terms as if they are the ``worst'' fields $\phi_{[N]}$. Despite this fact, the resulting estimates are sufficient to control the other error terms.

	Now, we are finally in a position to apply the $L^2$ bootstrap bounds for $h$ and $\phi$ (see section \ref{section L2 bootstrap bounds}). Combining this with the fact that $\epsilon \ll \delta \ll  1$ and that $\delta \ll \beta$ we have
	\begin{equation*}
	\int_{\mathcal{M}_\tau^{\tau_1}} \epsilon^{-1} (1+r)^{1-C_{[n,m]}\epsilon} |F_{(A,n)}|^2 \dVol_g
	\lesssim
	\left( \frac{1}{C_{[n,m]}} + \frac{\epsilon^2}{\delta^6} \right) \epsilon^{2(N_2 + 2 - n)} (1+\tau)^{-1+C_{(n)}\delta}
	\end{equation*}

	We need to provide a similar bound in the case $n = N_2$. In this case, we need to distinguish between the cases where the final commutation operator applied is $\slashed{\D}_T$, and the case where it is an arbitrary commutation operator. In fact, we will establish another hierarchy of estimates where, as we ascend the hierarchy there are fewer and fewer factors of $\slashed{\D}_T$ in the expansion of $\mathscr{Y}^n$, and the decay in $r$ is also progressively worse.

	First, we consider the case where the final commutation operator applied is arbitrary, i.e.\ it could be $r\slashed{\nabla}$, $r\slashed{\D}_L$ or $\slashed{\D}_T$. Following almost identical calculations to those above, we find
	\begin{equation}
	\label{equation L2 put together internal 3}
	\begin{split}
	&\int_{\mathcal{M}_\tau^{\tau_1}} \epsilon^{-1} (1+r)^{1-C_{[N_2,m]}\epsilon} |F_{(A,N_2)}|^2 \dVol_g
	\\
	&\lesssim
	\int_{\mathcal{M}_\tau^{\tau_1}} \bigg( 
		\epsilon^{-1} (1+r)^{1-C_{[N_2,m]}\epsilon} |F^{(0)}_{(A,N_2)}|^2
		+ \epsilon (1+r)^{2\delta} |F^{(N_2-1)}|^2_{(\text{frame})}
		\\
		&\phantom{\lesssim \int_{\mathcal{M}_\tau^{\tau_1}} \bigg(}
		+ \epsilon (1+r)^{1 - C_{[N_2-1,m]}\epsilon} |F^{(\leq N_2-1)}|^2_{(\text{frame})}
		+ \epsilon (1+r)^{-1-C_{[N_2,m]}\epsilon} |\slashed{\D} \mathscr{Y}^{N_2} \phi_{[m]}|^2
		\\
		&\phantom{\lesssim \int_{\mathcal{M}_\tau^{\tau_1}} \bigg(}
		+ \epsilon (1+r)^{-1-C_{[N_2,0]}\epsilon} |\slashed{\D} \mathscr{Y}^{N_2} \phi_{[0]}|^2
		+ \epsilon (1+r)^{-1-C_{[N_2,m-1]}\epsilon} |\slashed{\D} \mathscr{Y}^{N_2} \phi_{[m-1]}|^2
		\\
		&\phantom{\lesssim \int_{\mathcal{M}_\tau^{\tau_1}} \bigg(}
		+ \epsilon (1+r)^{-1-2\delta} |\slashed{\D} \mathscr{Y}^{N_2} \phi|^2
		+ \epsilon (1+r)^{-1-C_{[N_2-1,0]}\epsilon} |\overline{\slashed{\D}} \mathscr{Y}^{N_2} \phi|^2
		\\
		&\phantom{\lesssim \int_{\mathcal{M}_\tau^{\tau_1}} \bigg(}
		+ \epsilon^{-1} r^{-1-C_{[N_2,m]}\epsilon} |\overline{\slashed{\D}} \mathscr{Y}^{\leq N_2-1} \phi_{(A)}|^2
		+ \epsilon^3 (1+r)^{-3-C_{[N_2-1,0]}\epsilon} |\mathscr{Y}^{N_2} \phi|^2
		\\
		&\phantom{\lesssim \int_{\mathcal{M}_\tau^{\tau_1}} \bigg(}
		+ \epsilon (1+r)^{-3 - C_{[N_2-1,N]}\epsilon} |\mathscr{Y}^{\leq N_2-1} \phi_{(A)}|^2
		+ \epsilon^5 (1+r)^{-3-\frac{1}{2}C_{[N_2,m]}\epsilon} |\mathscr{Y}^{N_2} X_{(\text{frame, small})}|^2
		\\
		&\phantom{\lesssim \int_{\mathcal{M}_\tau^{\tau_1}} \bigg(}
		+ \epsilon^5 (1+r)^{-1 -(C_{[N_2,m]} - 2C_{(0,m)})\epsilon} |\mathscr{Y}^{N_2} \tr_{\slashed{g}}\chi_{(\text{small})}|^2
		\\
		&\phantom{\lesssim \int_{\mathcal{M}_\tau^{\tau_1}} \bigg(}
		+ |\mathscr{Y} \phi| |\slashed{\nabla}^2 \mathscr{Z}^{N_2-1} \log \mu|^2
		+ \epsilon^5 (1+r)^{-3-2\delta} |\mathscr{Y}^{N_2} \log \mu|^2
		\\
		&\phantom{\lesssim \int_{\mathcal{M}_\tau^{\tau_1}} \bigg(}
		+ \epsilon (1+r)^{-1-C_{[N_2-1,0]}\epsilon} |\bm{\Gamma}^{(N_2-1)}_{(-1+C_{(N_2-1)}\epsilon)}|^2
		+ \epsilon^3 (1+r)^{-1-2\delta} |\slashed{\D} \mathscr{Y}^{N_2} h|^2_{(\text{frame})}
		\\
		&\phantom{\lesssim \int_{\mathcal{M}_\tau^{\tau_1}} \bigg(}
		+ \epsilon^3 (1+r)^{-1 -(C_{[N_2,m]} - 2C_{(0,m)})\epsilon} |\slashed{\D} \mathscr{Y}^{N_2} h|^2_{LL}
		+ \epsilon^3 (1+r)^{-1 -C_{[N_2-1,0]}\epsilon} |\overline{\slashed{\D}} \mathscr{Y}^{N_2} h|^2_{(\text{frame})}
	\bigg) \dVol_g 
	\end{split}
	\end{equation}

	Using proposition \ref{proposition commute laplacian mu} and the pointwise bounds we have
	\begin{equation*}
	\begin{split}
	&\int_{\mathcal{M}_{\tau}^{\tau_1}} \bigg(
	\epsilon^{-1} (1+r)^{1-C_{[N_2(j), m]}\epsilon} r^2 |\bar{\partial} \phi|^2 |\slashed{\nabla}^2  \mathscr{Z}^{N_2-1} \log \mu|^2
	\bigg) \dVol_g
	\\
	&\lesssim
	\int_{\tau_0}^{\tau_1} \Bigg( \int_{r=0}^\infty \epsilon^{-1} (1+r)^{1-C_{[N_2(j), m]}\epsilon} r^4 \sup_{S_{\tau,r}} \left( |\bar{\partial} \phi|^2 \right) \int_{\mathbb{S}^2} \bigg(
		(1+r)^{2C_{(0)}\epsilon} |\slashed{\D}_T \mathscr{Y}^{N_2-1} \tr_{\slashed{g}}\chi_{(\text{small})}|^2
		\\
		&\phantom{\lesssim \int_{\tau_0}^{\tau_1} \Bigg( \int_{r=0}^\infty \epsilon^{-1} r^4 \sup_{S_{\tau,r}} \left( |\bar{\partial} \phi|^2 \right) \int_{\mathbb{S}^2} \bigg(}
		+ \epsilon^2 r^{-2} (1+r)^{-2+2C_{(1)}\epsilon} |\mathscr{Y}^{N_2} \log \mu|^2
		+ r^{-2} |\mathscr{Y}^n \zeta|^2
		\\
		&\phantom{\lesssim \int_{\tau_0}^{\tau_1} \Bigg( \int_{r=0}^\infty \epsilon^{-1} r^4 \sup_{S_{\tau,r}} \left( |\bar{\partial} \phi|^2 \right) \int_{\mathbb{S}^2} \bigg(}
		+ r^{-2} |\slashed{\D} \mathscr{Y}^n h|^2_{(\text{frame})}
		+ \epsilon^2 (1+r)^{-1 + C_{(N_1)}\epsilon} |\bm{\Gamma}^{(N_2 - 1)}_{(-1+C_{(N_2 - 1)}\epsilon)}|^2
		\\
		&\phantom{\lesssim \int_{\tau_0}^{\tau_1} \Bigg( \int_{r=0}^\infty \epsilon^{-1} r^4 \sup_{S_{\tau,r}} \left( |\bar{\partial} \phi|^2 \right) \int_{\mathbb{S}^2} \bigg(}
		+ \epsilon^2 (1+r)^{-1 + C_{(N_1)}\epsilon} |F^{(\leq N_2 - 1)}|_{(\text{frame})}^2		
	\bigg) \upd r \wedge \dVol_{\mathbb{S}^2} \Bigg) \upd \tau
	\\
	\end{split}
	\end{equation*}
	For the first and last two terms on the right hand side, we use the bound $r|\bar{\partial} \phi| \lesssim |\mathscr{Y} \phi| \lesssim \epsilon^3 (1+r)^{-1+C_{[N_1]}\epsilon}(1+\tau)^{C_{(N_1)}\delta}$, while for the other terms we use the bound $r|\bar{\partial} \phi| \lesssim \epsilon^3 (1+r)^{-\delta} (1+\tau)^{-\beta}$. This leads to the bound
	\begin{equation*}
	\begin{split}
	&\int_{\mathcal{M}_{\tau}^{\tau_1}} \bigg(
	\epsilon^{-1} r^2 |\bar{\partial} \phi|^2 |\slashed{\nabla}^2 \slashed{\D}_T^j \mathscr{Z}^{N_2-j-1} \log \mu|^2
	\bigg) \dVol_g
	\\
	&\lesssim
	\int_{\mathcal{M}_{\tau}^{\tau_1}} \bigg(
	\epsilon^3 (1+r)^{-1-(C_{[N_2]} - 2C_{[N_1]} - 2C_{(0)})\epsilon} (1+\tau)^{C_{(N_1)}\delta} |\slashed{\D}_T^{j+1} \mathscr{Y}^{N_2-j-2} \tr_{\slashed{g}}\chi_{(\text{small})}|^2
	\\
	&\phantom{\lesssim \int_{\mathcal{M}_{\tau}^{\tau_1}} \bigg(}
	+ \epsilon^3 (1+r)^{-3-2\delta} (1+\tau)^{-2\beta} |\mathscr{Y}^{N_2} \log \mu|^2
	+ \epsilon^3 (1+r)^{-3-2\delta} (1+\tau)^{-2\beta} |\mathscr{Y}^n \zeta|^2
	\\
	&\phantom{\lesssim \int_{\mathcal{M}_{\tau}^{\tau_1}} \bigg(}
	+ \epsilon^3 (1+r)^{-3-2\delta} (1+\tau)^{-2\beta} |\slashed{\D} \mathscr{Y}^n h|^2_{(\text{frame})}
	+ \epsilon^3 (1+r)^{-1-\frac{1}{2}C_{[N_2]}\epsilon} (1+\tau)^{C_{(N_1)}\delta} |\bm{\Gamma}^{(N_2 - 1)}_{(-1+C_{(N_2 - 1)}\epsilon)}|^2
	\\
	&\phantom{\lesssim \int_{\mathcal{M}_{\tau}^{\tau_1}} \bigg(}
	+ \epsilon^3 (1+r)^{-1-\frac{1}{2}C_{[N_2]}\epsilon} (1+\tau)^{C_{(N_1)}\delta} |F^{(\leq N_2 - 1)}|_{(\text{frame})}^2
	\bigg) \dVol_g
	\end{split}
	\end{equation*}

	Now, we return to equation \eqref{equation L2 put together internal 3}. This time, we use proposition \ref{proposition L2 bounds rectangular} to bound the term involving $\mathscr{Y}^{N_2} X_{(\text{frame, small})}$, proposition \ref{proposition L2 tr chi high} to bound the term involving $\mathscr{Y}^{N_2} \tr_{\slashed{g}}\chi_{(\text{small})}$, proposition \ref{proposition L2 foliation density} to bound the term involving $\mathscr{Y}^{N_2} \log \mu$, the bound derived above to bound the term involving $\slashed{\nabla}^2 \mathscr{Z}^{N_2-1} \log \mu$ and the bootstrap bounds to bound the lower order terms. In this way we obtain
	\begin{equation*}
	\begin{split}
	&\int_{\mathcal{M}_\tau^{\tau_1}} \epsilon^{-1} (1+r)^{1-C_{[N_2,m]}\epsilon} |F_{(A,N_2)}|^2 \dVol_g
	\\
	&\lesssim
	\epsilon \epsilon^4(1+\tau)^{-1+C_{(N_2)}\delta}
	\\
	&\phantom{\lesssim}
	+ \int_{\mathcal{M}_\tau^{\tau_1}} \bigg( 
		\epsilon^{-1} (1+r)^{1-C_{[N_2,m]}\epsilon} |F^{(0)}_{(A,N_2)}|^2
		+ \epsilon (1+r)^{2\delta} |F^{(N_2-1)}|^2_{(\text{frame})}
		\\
		&\phantom{\lesssim + \int_{\mathcal{M}_\tau^{\tau_1}} \bigg(}
		+ \epsilon (1+r)^{1 - C_{[N_2-1,m]}\epsilon} |F^{(\leq N_2-1)}|^2_{(\text{frame})}
		+ \epsilon (1+r)^{-1-C_{[N_2,m]}\epsilon} |\slashed{\D} \mathscr{Y}^{N_2} \phi_{[m]}|^2
		\\
		&\phantom{\lesssim + \int_{\mathcal{M}_\tau^{\tau_1}} \bigg(}
		+ \epsilon (1+r)^{-1-C_{[N_2,0]}\epsilon} |\slashed{\D} \mathscr{Y}^{N_2} \phi_{[0]}|^2
		+ \epsilon (1+r)^{-1-C_{[N_2,m-1]}\epsilon} |\slashed{\D} \mathscr{Y}^{N_2} \phi_{[m-1]}|^2
		\\
		&\phantom{\lesssim + \int_{\mathcal{M}_\tau^{\tau_1}} \bigg(}
		+ \epsilon (1+r)^{-1-2\delta} |\slashed{\D} \mathscr{Y}^{N_2} \phi|^2
		+ \epsilon (1+r)^{-1-C_{[N_2-1,0]}\epsilon} |\overline{\slashed{\D}} \mathscr{Y}^{N_2} \phi|^2
		\\
		&\phantom{\lesssim + \int_{\mathcal{M}_\tau^{\tau_1}} \bigg(}
		+ \epsilon^{-1} r^{-1-C_{[N_2,m]}\epsilon} |\overline{\slashed{\D}} \mathscr{Y}^{\leq N_2-1} \phi_{(A)}|^2
		+ \epsilon^3 (1+r)^{-3-C_{[N_2-1,0]}\epsilon} |\mathscr{Y}^{N_2} \phi|^2
		\\
		&\phantom{\lesssim + \int_{\mathcal{M}_\tau^{\tau_1}} \bigg(}
		+ \epsilon (1+r)^{-3 - C_{[N_2-1,N]}\epsilon} |\mathscr{Y}^{\leq N_2-1} \phi_{(A)}|^2
		+ \epsilon^3 (1+r)^{-1 -(C_{[N_2,m]} - 2C_{(0,m)})\epsilon} |\slashed{\D} \mathscr{Y}^{N_2} h|^2_{LL}
		\\
		&\phantom{\lesssim + \int_{\mathcal{M}_\tau^{\tau_1}} \bigg(}
		+ \epsilon^3 (1+r)^{-1 -(C_{[N_2,m]} - 2C_{[N_1]})\epsilon} (1+\tau)^{C_{(N_1)}\delta}|\slashed{\D}_T \mathscr{Y}^{N_2-1} \tr_{\slashed{g}} \chi_{(\text{small})}|^2
		\\
		&\phantom{\lesssim + \int_{\mathcal{M}_\tau^{\tau_1}} \bigg(}
		+ \epsilon^5 \delta^{-2} (1+r)^{-1-\delta} |\slashed{\D} \mathscr{Y}^{N_2} h|^2_{(\text{frame})}
		+ \epsilon^5 \delta^{-4} (1+r)^{-1+\frac{1}{2}\delta} |\overline{\slashed{\D}} \mathscr{Y}^{N_2} h|^2_{(\text{frame})}
		\\
		&\phantom{\lesssim + \int_{\mathcal{M}_\tau^{\tau_1}} \bigg(}
		+ \epsilon^5 \delta^{-3} (1+r)^{-3-\delta}|\mathscr{Y}^{N_2}h|^2_{(\text{frame})}
		+ \epsilon^3 (1+r)^{1 -(C_{[N_2,m]} - 2C_{(0,m)})\epsilon} |\mathscr{Y}^{N_2} F|^2_{LL}
		\\
		&\phantom{\lesssim + \int_{\mathcal{M}_\tau^{\tau_1}} \bigg(}
		+ \epsilon^3 (1+r)^{1 -(C_{[N_2,m]} - 2C_{[N_1]})\epsilon}(1+\tau)^{C_{(N_1)}\delta} |\slashed{\D}_T \mathscr{Y}^{N_2-1} F|^2_{LL}
		\\
		&\phantom{\lesssim + \int_{\mathcal{M}_\tau^{\tau_1}} \bigg(}
		+ \epsilon \sum_{k\leq N_2 -1} (1+r)^{1 -\frac{1}{2}C_{[N_2,m]}\epsilon} |\mathscr{Y}^{k} F|^2_{LL}
	\bigg) \dVol_g 
	\end{split}
	\end{equation*}

	Now, we can use proposition \ref{proposition L2 tr chi high} again to bound the term involving $\slashed{\D}_T \mathscr{Y}^{N_2-1} \tr_{\slashed{g}}\chi_{(\text{small})}$. Noting the formula for $[\slashed{\D}_L , \slashed{\D}_T]$ (see proposition \ref{proposition commuting DT with first order operators}) we find that the top order terms either include a $\slashed{\D}_T$ operator or a quantity that can be estimated in $L^\infty$ and that gives additional decay in $\tau$. Hence we obtain
	\begin{equation*}
	\begin{split}
	&\int_{\mathcal{M}_\tau^{\tau_1}} \epsilon^{-1} (1+r)^{1-C_{[N_2,m]}\epsilon} |F_{(A,N_2)}|^2 \dVol_g
	\\
	&\lesssim
	\epsilon \epsilon^4(1+\tau)^{-1+C_{(N_2)}\delta}
	\\
	&\phantom{\lesssim}
	+ \int_{\mathcal{M}_\tau^{\tau_1}} \bigg( 
	\epsilon^{-1} (1+r)^{1-C_{[N_2,m]}\epsilon} |F^{(0)}_{(A,N_2)}|^2
	+ \epsilon (1+r)^{2\delta} |F^{(N_2-1)}|^2_{(\text{frame})}
	\\
	&\phantom{\lesssim + \int_{\mathcal{M}_\tau^{\tau_1}} \bigg(}
	+ \epsilon (1+r)^{1 - C_{[N_2-1,m]}\epsilon} |F^{(\leq N_2-1)}|^2_{(\text{frame})}
	+ \epsilon (1+r)^{-1-C_{[N_2,m]}\epsilon} |\slashed{\D} \mathscr{Y}^{N_2} \phi_{[m]}|^2
	\\
	&\phantom{\lesssim + \int_{\mathcal{M}_\tau^{\tau_1}} \bigg(}
	+ \epsilon (1+r)^{-1-C_{[N_2,0]}\epsilon} |\slashed{\D} \mathscr{Y}^{N_2} \phi_{[0]}|^2
	+ \epsilon (1+r)^{-1-C_{[N_2,m-1]}\epsilon} |\slashed{\D} \mathscr{Y}^{N_2} \phi_{[m-1]}|^2
	\\
	&\phantom{\lesssim + \int_{\mathcal{M}_\tau^{\tau_1}} \bigg(}
	+ \epsilon (1+r)^{-1-2\delta} |\slashed{\D} \mathscr{Y}^{N_2} \phi|^2
	+ \epsilon (1+r)^{-1-C_{[N_2-1,0]}\epsilon} |\overline{\slashed{\D}} \mathscr{Y}^{N_2} \phi|^2
	\\
	&\phantom{\lesssim + \int_{\mathcal{M}_\tau^{\tau_1}} \bigg(}
	+ \epsilon^{-1} r^{-1-C_{[N_2,m]}\epsilon} |\overline{\slashed{\D}} \mathscr{Y}^{\leq N_2-1} \phi_{(A)}|^2
	+ \epsilon^3 (1+r)^{-3-C_{[N_2-1,0]}\epsilon} |\mathscr{Y}^{N_2} \phi|^2
	\\
	&\phantom{\lesssim + \int_{\mathcal{M}_\tau^{\tau_1}} \bigg(}
	+ \epsilon (1+r)^{-3 - C_{[N_2-1,N]}\epsilon} |\mathscr{Y}^{\leq N_2-1} \phi_{(A)}|^2
	+ \epsilon^3 (1+r)^{-1 -(C_{[N_2,m]} - 2C_{(0,m)})\epsilon} |\slashed{\D} \mathscr{Y}^{N_2} h|^2_{LL}
	\\
	&\phantom{\lesssim + \int_{\mathcal{M}_\tau^{\tau_1}} \bigg(}
	+ \epsilon^3 (1+r)^{-1 -(C_{[N_2,m]} - 2C_{[N_1]})\epsilon} (1+\tau)^{C_{(N_1)}\delta}|\slashed{\D} \slashed{\D}_T \mathscr{Y}^{N_2-1} h|^2_{LL}
	\\
	&\phantom{\lesssim + \int_{\mathcal{M}_\tau^{\tau_1}} \bigg(}
	+ \epsilon^5 \delta^{-2} (1+r)^{-1-\delta} |\slashed{\D} \mathscr{Y}^{N_2} h|^2_{(\text{frame})}
	+ \epsilon^5 \delta^{-4} (1+r)^{-1+\frac{1}{2}\delta} |\overline{\slashed{\D}} \mathscr{Y}^{N_2} h|^2_{(\text{frame})}
	\\
	&\phantom{\lesssim + \int_{\mathcal{M}_\tau^{\tau_1}} \bigg(}
	+ \epsilon^5 \delta^{-3} (1+r)^{-3-\delta}|\mathscr{Y}^{N_2}h|^2_{(\text{frame})}
	+ \epsilon^3 (1+r)^{1 -(C_{[N_2,m]} - 2C_{(0,m)})\epsilon} |\mathscr{Y}^{N_2} F|^2_{LL}
	\\
	&\phantom{\lesssim + \int_{\mathcal{M}_\tau^{\tau_1}} \bigg(}
	+ \epsilon^3 (1+r)^{1 -(C_{[N_2,m]} - 2C_{[N_1]})\epsilon} (1+\tau)^{C_{(N_1)}\delta} |\slashed{\D}_T \mathscr{Y}^{N_2-1} F|^2_{LL}
	\\
	&\phantom{\lesssim + \int_{\mathcal{M}_\tau^{\tau_1}} \bigg(}
	+ \epsilon \sum_{k\leq N_2 -1} (1+r)^{1 -\frac{1}{2}C_{[N_2,m]}\epsilon} |\mathscr{Y}^{k} F|^2_{LL}
	\bigg) \dVol_g 
	\end{split}
	\end{equation*}

	As before, we note that, either $m = 0$, in which case
	\begin{equation*}
	\begin{split}
	&\int_{\mathcal{M}_{\tau}^{\tau_1}} \epsilon^3 (1+r)^{-1-(C_{[N_2,m]} - 2C_{(0,m)}\epsilon)} (1+\tau)^ |\slashed{\D} \mathscr{Y}^{N_2} h|^2_{LL} \dVol_g
	\\
	&=
	\int_{\mathcal{M}_{\tau}^{\tau_1}} \epsilon^3 (1+r)^{-1-C_{[N_2,m]}\epsilon} |\slashed{\D} \mathscr{Y}^{N_2} h|^2_{LL} \dVol_g
	\\
	&\lesssim
	\frac{1}{C_{[N_2,m]}} \epsilon^4 (1+\tau)^{-1+C_{(N_2)}\delta}
	\end{split}
	\end{equation*}
	while if $m \geq 1$ then we have
	\begin{equation*}
	\begin{split}
	\int_{\mathcal{M}_{\tau}^{\tau_1}} \epsilon^3 (1+r)^{-1-(C_{[N_2,m]} - 2C_{(0,m)}\epsilon)} |\slashed{\D} \mathscr{Y}^{N_2} h|^2_{LL} \dVol_g
	&\lesssim
	\int_{\mathcal{M}_{\tau}^{\tau_1}} \epsilon^3 (1+r)^{-1-C_{[N_2,m-1]}\epsilon} |\slashed{\D} \mathscr{Y}^{N_2} \phi_{[0]}|^2 \dVol_g
	\\
	&\lesssim
	\frac{1}{C_{[N_2,m]}} \epsilon^4 (1+\tau)^{-1+C_{(N_2)}\delta}
	\end{split}
	\end{equation*}
	
	On the other hand, we choose the constants so that
	\begin{equation*}
	\begin{split}
	C_{[N_2(j-1),m]} \gg C_{[N_2(j),m]}
	\\
	C_{(N_2(j-1))} \gg C_{(N_2(j))}
	\end{split}	
	\end{equation*}
	and $C_{[N_2(0),m]} = C_{[N_2,m]}$, $C_{(N_2(0))} = C_{(N_2)}$. The idea is that the index $j$ will label the number of times that $\slashed{\D}_T$ appears on the left hand side of the expansion of $\mathscr{Y}^n$. Hence we have
	\begin{equation*}
	\begin{split}
	&\int_{\mathcal{M}_\tau^{\tau_1}} \epsilon^3 (1+r)^{-1-(C_{[N_2,m]} - 2C_{[N_1]})\epsilon} (1+\tau)^{C_{(N_1)}\delta} |\slashed{\D} \slashed{\D}_T \mathscr{Y}^{N_2 - 1} h|^2_{LL} \dVol_g
	\\
	&\lesssim
	\int_{\mathcal{M}_\tau^{\tau_1}} \epsilon^3 (1+r)^{-1-C_{[N_2(1),m]}\epsilon} (1+\tau)^{C_{(N_1)}\delta} |\slashed{\D} \slashed{\D}_T \mathscr{Y}^{N_2 - 1} h|^2_{LL} \dVol_g
	\end{split}
	\end{equation*}
	
	So, using the bootstrap assumptions (see section \ref{section L2 bootstrap bounds}) we have that
	\begin{equation*}
	\begin{split}
	&\int_{\mathcal{M}_\tau^{\tau_1}} \epsilon^3 (1+r)^{-1-(C_{[N_2,m]} - 2C_{[N_1]})\epsilon} (1+\tau)^{C_{(N_1)}\delta} |\slashed{\D} \slashed{\D}_T \mathscr{Y}^{N_2 - 1} h|^2_{LL} \dVol_g
	\\
	&\lesssim
	\frac{1}{C_{[N_2(1),m]}}\epsilon^4 (1+\tau)^{-1+C_{(N_2(1))}\delta + C_{(N_1)}\delta}
	\\
	&\lesssim
	\frac{1}{C_{[N_2(1),m]}}\epsilon^4 (1+\tau)^{-1+C_{(N_2(0))}\delta}
	\end{split}
	\end{equation*}
	where we have used proposition \ref{proposition integral inequality} in the first line, and the fact that $N_{(2)}(0) \gg N_2(1)$ in the second line. Hence, putting these last few calculations together, we have
	\begin{equation*}
	\int_{\mathcal{M}_\tau^{\tau_1}} \epsilon^{-1} (1+r)^{1-C_{[N_2(0),m]}\epsilon} |F_{(A,N_2(0))}|^2 \dVol_g
	\lesssim
	\frac{1}{C_{[N_2(1),m]}}\epsilon^4 (1+\tau)^{-1+C_{(N_2(0))}\delta}
	\end{equation*}

	We also need to provide bounds for $F_{(A, N_2(j))}$ with $j \geq 1$. This time we must be more careful with the decay in $\tau$, since we must also show the improved decay in $\tau$. We claim that, schematically, we have
	\begin{equation}
	\begin{split}
		F_{(A, N_2(j))}
		&=
		F^{(0)}_{(A, N_2(j))}
		+ (F_{(A)}^{(BC)}) (\partial \phi_{[0]}) (\slashed{\D} \mathscr{Y}^{N_2} \phi_{[m]})
		+ (F_{(A)}^{(BC)}) (\partial \phi_{[m]}) (\slashed{\D} \mathscr{Y}^{N_2} \phi_{[0]})
		\\
		&\phantom{=}
		+ (F_{(A)}^{(BC)}) (\partial \phi_{[m-1]}) (\slashed{\D} \mathscr{Y}^{N_2} \phi_{[m-1]})
		+ (F_{(A)}^{(BC)}) (\partial \phi) (\overline{\slashed{\D}} \mathscr{Y}^{N_2} \phi)
		\\
		&\phantom{=}
		+ (F_{(A)}^{(BC)}) (\bar{\partial} \phi) (\slashed{\D} \mathscr{Y}^{N_2} \phi)
		+ (\partial \phi) (\phi) (\slashed{\D} \mathscr{Y}^{N_2} \phi)
		+ (\partial \phi)^2 (\mathscr{Y}^{N_2} \phi)
		\\
		&\phantom{=}
		+ (\mathscr{Y}^{N_2} X_{(\text{frame, small})}) (\partial \phi)^2
		+ \sum_{\substack{j+k \leq N_2 \\ j,k \leq N_2 - 1}} \left(\bm{\Gamma}^{(j)}_{(-1+C_{(j)}\epsilon)} \right) \left(\bm{\Gamma}^{(k)}_{(-1+C_{(k)}\epsilon)} \right)
		+ r^{-1} \overline{\slashed{\D}} (\mathscr{Y}^{\leq N_2 - 1} \phi)
		\\
		&\phantom{=}
		+ \begin{pmatrix}
			(\slashed{\nabla} \log \mu) \\
			\zeta \\
			(\chi_{(\text{small})} + \chibar_{(\text{small})})
		\end{pmatrix} \overline{\slashed{\D}}(\mathscr{Y}^{N_2} \phi_{(A)})
		+ \bm{\Gamma}^{(0)}_{(-1)} (\slashed{\D} \slashed{\D}_T^j \mathscr{Y}^{N_2-j} \phi_{(A)})
		\\
		&\phantom{=}
		+ \bm{\Gamma}^{(1)}_{(-1-\delta)} (\overline{\slashed{\D}} \slashed{\D}_T^j \mathscr{Y}^{N_2-j} \phi_{(A)})
		+ \begin{pmatrix}
			r^{-1} (\partial \phi_{(A)}) \\
			(\bar{\partial} \phi_{(A)}) \\
			r^{-1} \mathscr{Y}\phi_{(A)}
		\end{pmatrix} \left( \slashed{\D} \mathscr{Y}^{N_2} h \right)_{(\text{frame})}
		+ (\partial \phi_{(A)}) (\slashed{\D} \mathscr{Y}^{N_2} h)_{LL}
		\\
		&\phantom{=}
		+ (\partial \phi_{(A)})\left( \overline{\slashed{\D}} \mathscr{Y}^{N_2} h \right)_{(\text{frame})}
		+  \begin{pmatrix}
			(\partial \phi_{(A)}) \\
			r (\bar{\partial} \phi_{(A)}) \\
			\mathscr{Y} \phi_{(A)}
		\end{pmatrix} \left( \tilde{\slashed{\Box}}_g \mathscr{Y}^{N_2-1} h \right)_{(\text{frame})}
		+ (\partial \phi_{(A)}) \mathscr{Y}^{N_2} \tr_{\slashed{g}}\chi_{(\text{small})}
		\\
		&\phantom{=}
		+ (\bar{\partial} \phi_{(A)}) (r\slashed{\nabla}^2 \slashed{\D}_T^j \mathscr{Z}^{N_2-j-1} \log \mu)
		+ \begin{pmatrix}
			r^{-1} (\overline{\slashed{\D}} \mathscr{Y} \phi_{(A)}) \\
			\bm{\Gamma}^{(1)}_{(-2+ C_{(1)}\epsilon)} (\mathscr{Y}\phi_{(A)})
		\end{pmatrix} (\mathscr{Y}^{N_2} \log \mu)
		\\
		&\phantom{=}
		+ \sum_{\substack{ j+k \leq N_2+1 \\ j,k \leq N_2-1}} \bm{\Gamma}^{(j)}_{(-1 + C_{(j)}\epsilon)} (\slashed{\D} \mathscr{Y}^k \phi_{(A)})
		+ \sum_{\substack{ j+k \leq N_2 \\ j,k \leq N_2-1}} r\tilde{\slashed{\Box}}_g (\mathscr{Y}^j h)_{(\text{frame})}  (\slashed{\D} \mathscr{Y}^k \phi_{(A)})
		\\
		&\phantom{=}
		+ \sum_{\substack{j+k \leq N_2+1 \\ j,k \leq N_2-1 }} \bm{\Gamma}^{(j)}_{(-2 + 2C_{(j)}\epsilon)} (\mathscr{Y}^k \phi_{(A)})
	\end{split}
	\end{equation}
	In fact, proposition \ref{proposition inhomogeneous terms after commuting with T and Y} provides us with most of the terms, while the first few terms follow from the form of $F_{(A, 0)}$.
	
	Now, using the pointwise bounds we have
	\begin{equation}
	\begin{split}
	|F_{(A, N_2(j))}|
	&\lesssim
	|F^{(0)}_{(A, N_2(j))}|
	+ \epsilon (1+r)^{-1} (1+\tau)^{-\beta} |\slashed{\D} \mathscr{Y}^{N_2} \phi_{[m]}|
	+ \epsilon (1+r)^{-1+C_{(0,m)}\epsilon} (1+\tau)^{-\beta} |\slashed{\D} \mathscr{Y}^{N_2} \phi_{[0]}|
	\\
	&\phantom{\lesssim}
	+ \epsilon (1+r)^{-1+C_{(0,m-1)}\epsilon} (1+\tau)^{-\beta} |\slashed{\D} \mathscr{Y}^{N_2} \phi_{[m-1]}|
	+ \epsilon (1+r)^{-1+C_{(0,N)}\epsilon} (1+\tau)^{-\beta} |\overline{\slashed{\D}} \mathscr{Y}^{N_2} \phi|
	\\
	&\phantom{\lesssim}
	+ \epsilon (1+r)^{-1-\delta} (1+\tau)^{-\beta} |\slashed{\D} \mathscr{Y}^{N_2} \phi|
	+ \epsilon^2 (1+r)^{-2+2C_{(N)}\epsilon} (1+\tau)^{-2\beta} |\mathscr{Y}^{N_2} \phi|
	\\
	&\phantom{\lesssim}
	+ \epsilon (1+r)^{-1+C_{(0,N_1)}\epsilon} \left|\bm{\Gamma}^{(N_2 - 1)}_{(-1+C_{(N_2 - 1)}\epsilon)} \right|
	+ r^{-1} |\overline{\slashed{\D}} (\mathscr{Y}^{\leq N_2 - 1} \phi)|
	\\
	&\phantom{\lesssim}
	+ \epsilon (1+r)^{-1 + C_{(1)}\epsilon)} (1+\tau)^{-C^* \delta} |\overline{\slashed{\D}} (\mathscr{Y}^{N_2} \phi)|
	+ \epsilon (1+r)^{-1}(1+\tau)^{-\beta} |\slashed{\D} \slashed{\D}_T^j \mathscr{Y}^{N_2-j} \phi|
	\\
	&\phantom{\lesssim}
	+ \epsilon^3 (1+r)^{-1-\delta} (1+\tau)^{-\beta} |\slashed{\D} \mathscr{Y}^{N_2} h |_{(\text{frame})}
	+ \epsilon^3 (1+r)^{-1+C_{(0,m)}\epsilon}(1+\tau)^{-\beta} |\slashed{\D} \mathscr{Y}^{N_2} h|_{LL}
	\\
	&\phantom{\lesssim}
	+ \epsilon^3 (1+r)^{-1+C_{(0,m)}\epsilon}(1+\tau)^{-\beta} |\overline{\slashed{\D}} \mathscr{Y}^{N_2} h|_{(\text{frame})}
	+ \epsilon (1+r)^{-\delta} (1+\tau)^{-\beta} |\tilde{\slashed{\Box}}_g \mathscr{Y}^{N_2-1} h |_{(\text{frame})}
	\\
	&\phantom{\lesssim}
	+ \epsilon^5 (1+r)^{-1+C_{(0,m)}\epsilon}(1+\tau)^{-\beta} |\mathscr{Y}^{N_2} \tr_{\slashed{g}}\chi_{(\text{small})}|
	+ r|\bar{\partial} \phi| |\slashed{\nabla}^2 \slashed{\D}_T^j \mathscr{Z}^{N_2-j-1} \log \mu|
	\\
	&\phantom{\lesssim}
	+ \epsilon^5 (1+r)^{-2-\delta}(1+\tau)^{-\beta} |\mathscr{Y}^{N_1} \log \mu|
	+ \sum_{\substack{ j+k \leq N_2 \\ j,k \leq N_2-1}} r|\tilde{\slashed{\Box}}_g (\mathscr{Y}^j h)|_{(\text{frame})} |\slashed{\D} \mathscr{Y}^k \phi|
	\\
	&\phantom{\lesssim}
	+ \epsilon (1+r)^{-2+2C_{(N_1)}\epsilon} (\mathscr{Y}^{\leq N_2 - 1} \phi) 
	\end{split}
	\end{equation}
	
	In particular, this gives us 
	\begin{equation*}
	\begin{split}
	&\int_{\mathcal{M}_{\tau}^{\tau_1}} \epsilon^{-1} (1+r)^{1-C_{[N_2(j), m]}\epsilon} 	|F_{(A, N_2(j))}| \dVol_g
	\\
	&\lesssim
	\int_{\mathcal{M}_{\tau}^{\tau_1}} \bigg(
		\epsilon^{-1} (1+r)^{-1-C_{[N_2(j), m]}\epsilon} |F^{(0)}_{(A, N_2(j))}|^2
		+ \epsilon (1+r)^{-1-C_{[N_2(j), m]}\epsilon} (1+\tau)^{-2\beta} |\slashed{\D} \mathscr{Y}^{N_2} \phi_{[m]}|^2
		\\
		&\phantom{\lesssim \int_{\mathcal{M}_{\tau}^{\tau_1}} \bigg(}
		+ \epsilon (1+r)^{-1-(C_{[N_2(j), m]}-2C_{(0,m)})\epsilon} (1+\tau)^{-2\beta} |\slashed{\D} \mathscr{Y}^{N_2} \phi_{[0]}|^2
		\\
		&\phantom{\lesssim \int_{\mathcal{M}_{\tau}^{\tau_1}} \bigg(}
		+ \epsilon (1+r)^{-1-(C_{[N_2(j), m]}-C_{(m-1)})\epsilon} (1+\tau)^{-2\beta} |\slashed{\D} \mathscr{Y}^{N_2} \phi_{[m-1]}|^2
		\\
		&\phantom{\lesssim \int_{\mathcal{M}_{\tau}^{\tau_1}} \bigg(}
		+ \epsilon (1+r)^{-1-\frac{1}{2}C_{[N_2(j), m]}\epsilon} (1+\tau)^{-2\beta} |\overline{\slashed{\D}} \mathscr{Y}^{N_2} \phi|^2
		+ \epsilon (1+r)^{-2-2\delta} (1+\tau)^{-2\beta} |\slashed{\D} \mathscr{Y}^{N_2} \phi|^2
		\\
		&\phantom{\lesssim \int_{\mathcal{M}_{\tau}^{\tau_1}} \bigg(}
		+ \epsilon^3 (1+r)^{-3-\frac{1}{2}C_{[N_2(j), m]}\epsilon} (1+\tau)^{-4\beta} |\mathscr{Y}^{N_2} \phi|		
		+ \epsilon (1+r)^{-1-\frac{1}{2}C_{[N_2(j), m]}\epsilon} \left|\bm{\Gamma}^{(N_2 - 1)}_{(-1+C_{(N_2 - 1)}\epsilon)} \right|^2
		\\
		&\phantom{\lesssim \int_{\mathcal{M}_{\tau}^{\tau_1}} \bigg(}
		+ \epsilon^{-1} (1+r)^{-1-C_{[N_2(j), m]}\epsilon} |\overline{\slashed{\D}} (\mathscr{Y}^{\leq N_2 - 1} \phi)|^2
		+ \epsilon (1+r)^{-1 -\frac{1}{2}C_{[N_2(j), m]}\epsilon} (1+\tau)^{-C^* \delta} |\overline{\slashed{\D}} (\mathscr{Y}^{N_2} \phi)|^2
		\\
		&\phantom{\lesssim \int_{\mathcal{M}_{\tau}^{\tau_1}} \bigg(}
		+ \epsilon (1+r)^{-1 -C_{[N_2(j), m]}\epsilon} (1+\tau)^{-2\beta} |\slashed{\D} \slashed{\D}_T^j \mathscr{Y}^{N_2-j} \phi|^2
		+ \epsilon^3 (1+r)^{-1-2\delta} (1+\tau)^{-2\beta} |\slashed{\D} \mathscr{Y}^{N_2} h |_{(\text{frame})}^2
		\\
		&\phantom{\lesssim \int_{\mathcal{M}_{\tau}^{\tau_1}} \bigg(}
		+ \epsilon^3 (1+r)^{-1-(C_{[N_2(j), m]}-2C_{(0,m)})\epsilon}(1+\tau)^{-2\beta} |\slashed{\D} \mathscr{Y}^{N_2} h|_{LL}
		\\
		&\phantom{\lesssim \int_{\mathcal{M}_{\tau}^{\tau_1}} \bigg(}
		+ \epsilon^5 (1+r)^{-1-\frac{1}{2}C_{[N_2(j), m]}\epsilon}(1+\tau)^{-2\beta} |\overline{\slashed{\D}} \mathscr{Y}^{N_2} h|_{(\text{frame})}^2
		\\
		&\phantom{\lesssim \int_{\mathcal{M}_{\tau}^{\tau_1}} \bigg(}
		+ \epsilon^3 (1+r)^{1-2\delta} (1+\tau)^{-2\beta} |\tilde{\slashed{\Box}}_g \mathscr{Y}^{N_2-1} h |_{(\text{frame})}^2
		\\
		&\phantom{\lesssim \int_{\mathcal{M}_{\tau}^{\tau_1}} \bigg(}
		+ \epsilon^5 (1+r)^{-1-(C_{[N_2(j), m]}-2C_{(0,m)})\epsilon}(1+\tau)^{-2\beta} |\mathscr{Y}^{N_2} \tr_{\slashed{g}}\chi_{(\text{small})}|^2
		\\
		&\phantom{\lesssim \int_{\mathcal{M}_{\tau}^{\tau_1}} \bigg(}
		+ \epsilon^{-1} r^2|\bar{\partial} \phi|^2 |\slashed{\nabla}^2 \slashed{\D}_T^j \mathscr{Z}^{N_2-j-1} \log \mu|^2
		+ \epsilon^5 (1+r)^{-3-2\delta}(1+\tau)^{-2\beta} |\mathscr{Y}^{N_1} \log \mu|^2
		\\
		&\phantom{\lesssim \int_{\mathcal{M}_{\tau}^{\tau_1}} \bigg(}
		+ \epsilon^4 (1+r)^{1-\frac{1}{2}C_{[N_2(j), m]}\epsilon} (1+\tau)^{-2\beta} |\tilde{\slashed{\Box}}_g (\mathscr{Y}^{\leq N_2 - 1} h)|_{(\text{frame})}^2
		\\
		&\phantom{\lesssim \int_{\mathcal{M}_{\tau}^{\tau_1}} \bigg(}
		+ \epsilon (1+r)^{-3-\frac{1}{2}C_{[N_2(j), m]}\epsilon} |\mathscr{Y}^{\leq N_2 - 1} \phi|^2
	\bigg) \dVol_g
	\end{split}
	\end{equation*}
	
	As above, we can bound some of the terms involving \emph{critical} decay in $r$ by considering separately the cases $m = 0$ (in which case $C_{(0,0)} = 0$) and $m \geq 0$. This time, however, many of the terms have additional decay in $\tau$. We can use proposition \ref{proposition L2 tr chi high} to bound the term involving $\mathscr{Y}^{N_2} \tr_{\slashed{g}}\chi_{(\text{small})}$, obtaining
	\begin{equation*}
	\begin{split}
		&\int_{\mathcal{M}_{\tau}^{\tau_1}} \bigg(
			\epsilon^3 (1+r)^{-1-C_{[N_2(j), 0]}\epsilon}(1+\tau)^{-2\beta} |\mathscr{Y}^{N_2} \tr_{\slashed{g}}\chi_{(\text{small})}|^2
		\bigg) \dVol_g
		\\
		&\lesssim
		\epsilon^7 (1+\tau)^{-1 -2\beta + C_{(N_2)}\delta}
		\\
		&\phantom{\lesssim}
		+ (1+\tau)^{-2\beta}\int_{\mathcal{M}_\tau^{\tau_1} \cap \{r \geq \frac{1}{2}r_0\}}\bigg(
			\epsilon^3 r^{-1-C_{[N_2]}\epsilon} |\slashed{\D} \mathscr{Y}^{N_2} h|^2_{LL}
			+ \epsilon^3 r^{-1-\frac{1}{2}C_{[N_2]}\epsilon} |\overline{\slashed{\D}} \mathscr{Y}^{N_2} h|^2_{(\text{frame})}
			\\
			&\phantom{\lesssim + (1+\tau)^{-2\beta}\int_{\mathcal{M}_\tau^{\tau_1} \cap \{r \geq \frac{1}{2}r_0\}}\bigg(}
			+ \epsilon^3 r^{-1-\delta} |\slashed{\D} \mathscr{Y}^{N_2} h|^2_{(\text{frame})}
			+ \epsilon^3 r^{-3-\delta} |\mathscr{Y}^{N_2} h|^2_{(\text{frame})}
			\\
			&\phantom{\lesssim + (1+\tau)^{-2\beta}\int_{\mathcal{M}_\tau^{\tau_1} \cap \{r \geq \frac{1}{2}r_0\}}\bigg(}
			+ \epsilon^3 r^{1-C_{[N_2]}\epsilon} |\mathscr{Y}^{N_2} F|^2_{LL}
			+ \sum_{j \leq N_2 - 1} \epsilon^3 r^{1-\frac{1}{2}C_{[N_2]}\epsilon} |\mathscr{Y}^j F|_{(\text{frame})}^2
		\bigg) \dVol_{g}
	\end{split}
	\end{equation*}
	
	Next, we can use the coarea formula (proposition \ref{proposition coarea}) to obtain
	\begin{equation*}
	\begin{split}
		&\int_{\mathcal{M}_{\tau}^{\tau_1}} \bigg(
			\epsilon^{-1} (1+r)^{1-C_{[N_2(j), m]}\epsilon} r^2 |\bar{\partial} \phi|^2 |\slashed{\nabla}^2 \slashed{\D}_T^j \mathscr{Z}^{N_2-j-1} \log \mu|^2
		\bigg) \dVol_g
		\\
		\lesssim
		&\int_{\tau_0}^{\tau_1} \Bigg( \int_{\Sigma_\tau} \bigg(
			\epsilon^{-1} (1+r)^{1-C_{[N_2(j), m]}\epsilon} r^2 |\bar{\partial} \phi|^2 |\slashed{\nabla}^2 \slashed{\D}_T^j \mathscr{Y}^{N_2-j-1} \log \mu|^2
		\bigg) r^2 \upd r \wedge \dVol_{\mathbb{S}^2} \Bigg) \upd \tau
		\\
		\lesssim
		&\int_{\tau_0}^{\tau_1} \Bigg( \int_{r=0}^\infty \epsilon^{-1} (1+r)^{1-C_{[N_2(j), m]}\epsilon} r^4 \sup_{S_{\tau,r}} \left( |\bar{\partial} \phi|^2 \right) \int_{\mathbb{S}^2} \bigg(
			 |\slashed{\nabla}^2 \slashed{\D}_T^j \mathscr{Y}^{N_2-j-1} \log \mu|^2
		\bigg) \upd r \wedge \dVol_{\mathbb{S}^2} \Bigg) \upd \tau
		\\
	\end{split}
	\end{equation*}
	Now, using proposition \ref{proposition Poisson estimate} we obtain
	\begin{equation*}
	\begin{split}
	&\int_{\mathcal{M}_{\tau}^{\tau_1}} \bigg(
		\epsilon^{-1} (1+r)^{1-C_{[N_2(j), m]}\epsilon} r^2 |\bar{\partial} \phi|^2 |\slashed{\nabla}^2 \slashed{\D}_T^j \mathscr{Z}^{N_2-j-1} \log \mu|^2
	\bigg) \dVol_g
	\\
	&\lesssim
	\int_{\tau_0}^{\tau_1} \Bigg( \int_{r=0}^\infty \epsilon^{-1} (1+r)^{1-C_{[N_2(j), m]}\epsilon} r^4 \sup_{S_{\tau,r}} \left( |\bar{\partial} \phi|^2 \right) \int_{\mathbb{S}^2} \bigg(
		|\slashed{\Delta} \slashed{\D}_T^j \mathscr{Y}^{N_2-j-1} \log \mu|^2
		\\
		&\phantom{\lesssim \int_{\tau_0}^{\tau_1} \Bigg( \int_{r=0}^\infty \epsilon^{-1} }
		+ r^{-4}|\slashed{\nabla} \slashed{\D}_T^j \mathscr{Y}^{N_2-j-1} \log \mu|^2
		+ \epsilon r^{-4}|\slashed{\D}_T^j \mathscr{Y}^{N_2-j-1} \log \mu|^2
	\bigg) \upd r \wedge \dVol_{\mathbb{S}^2} \Bigg) \upd \tau
	\\
	\end{split}
	\end{equation*}
	(c.f.\ proposition \ref{proposition L2 nabla log mu}). Next, and similarly to before, we have
	\begin{equation*}
	\begin{split}
	&\int_{\mathcal{M}_{\tau}^{\tau_1}} \bigg(
	\epsilon^{-1} r^2 |\bar{\partial} \phi|^2 |\slashed{\nabla}^2 \slashed{\D}_T^j \mathscr{Z}^{N_2-j-1} \log \mu|^2
	\bigg) \dVol_g
	\\
	&\lesssim
	\int_{\mathcal{M}_{\tau}^{\tau_1}} \bigg(
	\epsilon^3 (1+r)^{-1-(C_{[N_2]} - 2C_{[N_1]} - 2C_{(0)})\epsilon} (1+\tau)^{C_{(N_1)}\delta} |\slashed{\D}_T^{j+1} \mathscr{Y}^{N_2-j-1} \tr_{\slashed{g}}\chi_{(\text{small})}|^2
	\\
	&\phantom{\lesssim \int_{\mathcal{M}_{\tau}^{\tau_1}} \bigg(}
	+ \epsilon^3 (1+r)^{-3-2\delta} (1+\tau)^{-2\beta} |\mathscr{Y}^{N_2} \log \mu|^2
	+ \epsilon^3 (1+r)^{-3-2\delta} (1+\tau)^{-2\beta} |\mathscr{Y}^n \zeta|^2
	\\
	&\phantom{\lesssim \int_{\mathcal{M}_{\tau}^{\tau_1}} \bigg(}
	+ \epsilon^3 (1+r)^{-3-2\delta} (1+\tau)^{-2\beta} |\slashed{\D} \mathscr{Y}^n h|^2_{(\text{frame})}
	+ \epsilon^3 (1+r)^{-1-\frac{1}{2}C_{[N_2]}\epsilon} (1+\tau)^{C_{(N_1)}\delta} |\bm{\Gamma}^{(N_2 - 1)}_{(-1+C_{(N_2 - 1)}\epsilon)}|^2
	\\
	&\phantom{\lesssim \int_{\mathcal{M}_{\tau}^{\tau_1}} \bigg(}
	+ \epsilon^3 (1+r)^{-1-\frac{1}{2}C_{[N_2]}\epsilon} (1+\tau)^{C_{(N_1)}\delta} |F^{(\leq N_2 - 1)}|_{(\text{frame})}^2
	\bigg) \dVol_g
	\end{split}
	\end{equation*}

	Putting all of these calculations together, we find that
	\begin{equation*}
	\begin{split}
	&\int_{\mathcal{M}_{\tau}^{\tau_1}} \epsilon^{-1} (1+r)^{1-C_{[N_2(j), m]}\epsilon} |F_{(A, N_2(j))}| \dVol_g
	\\
	&\lesssim
	\epsilon^7 (1+\tau)^{-1-2\beta + C_{(N_2)}\delta}
	\\
	&+ \int_{\mathcal{M}_{\tau}^{\tau_1}} \bigg(
		\epsilon^{-1} (1+r)^{-1-C_{[N_2(j), m]}\epsilon} |F^{(0)}_{(A, N_2(j))}|^2
		+ \epsilon (1+r)^{-1-C_{[N_2(j), m]}\epsilon} (1+\tau)^{-2\beta} |\slashed{\D} \mathscr{Y}^{N_2} \phi_{[m]}|^2
		\\
		&\phantom{\lesssim \int_{\mathcal{M}_{\tau}^{\tau_1}} \bigg(}
		+ \epsilon (1+r)^{-1-(C_{[N_2(j), m]}-2C_{(0,m)})\epsilon} (1+\tau)^{-2\beta} |\slashed{\D} \mathscr{Y}^{N_2} \phi_{[0]}|^2
		\\
		&\phantom{\lesssim \int_{\mathcal{M}_{\tau}^{\tau_1}} \bigg(}
		+ \epsilon (1+r)^{-1-(C_{[N_2(j), m]}-C_{(m-1)})\epsilon} (1+\tau)^{-2\beta} |\slashed{\D} \mathscr{Y}^{N_2} \phi_{[m-1]}|^2
		\\
		&\phantom{\lesssim \int_{\mathcal{M}_{\tau}^{\tau_1}} \bigg(}
		+ \epsilon (1+r)^{-1-\frac{1}{2}C_{[N_2(j), m]}\epsilon} (1+\tau)^{-2\beta} |\overline{\slashed{\D}} \mathscr{Y}^{N_2} \phi|^2
		+ \epsilon (1+r)^{-2-2\delta} (1+\tau)^{-2\beta} |\slashed{\D} \mathscr{Y}^{N_2} \phi|^2
		\\
		&\phantom{\lesssim \int_{\mathcal{M}_{\tau}^{\tau_1}} \bigg(}
		+ \epsilon^3 (1+r)^{-3-\frac{1}{2}C_{[N_2(j), m]}\epsilon} (1+\tau)^{-4\beta} |\mathscr{Y}^{N_2} \phi|		
		+ \epsilon (1+r)^{-1-\frac{1}{2}C_{[N_2(j), m]}\epsilon} \left|\bm{\Gamma}^{(N_2 - 1)}_{(-1+C_{(N_2 - 1)}\epsilon)} \right|^2
		\\
		&\phantom{\lesssim \int_{\mathcal{M}_{\tau}^{\tau_1}} \bigg(}
		+ \epsilon^{-1} (1+r)^{-1-C_{[N_2(j), m]}\epsilon} |\overline{\slashed{\D}} (\mathscr{Y}^{\leq N_2 - 1} \phi)|^2
		+ \epsilon (1+r)^{-1 -\frac{1}{2}C_{[N_2(j), m]}\epsilon} (1+\tau)^{-C^* \delta} |\overline{\slashed{\D}} (\mathscr{Y}^{N_2} \phi)|^2
		\\
		&\phantom{\lesssim \int_{\mathcal{M}_{\tau}^{\tau_1}} \bigg(}
		+ \epsilon (1+r)^{-1 -C_{[N_2(j), m]}\epsilon} (1+\tau)^{-2\beta} |\slashed{\D} \slashed{\D}_T^j \mathscr{Y}^{N_2-j} \phi|^2
		+ \epsilon^3 (1+r)^{-1-\delta} (1+\tau)^{-2\beta} |\slashed{\D} \mathscr{Y}^{N_2} h |_{(\text{frame})}^2
		\\
		&\phantom{\lesssim \int_{\mathcal{M}_{\tau}^{\tau_1}} \bigg(}
		+ \epsilon^3 (1+r)^{-1-(C_{[N_2(j), m]}-2C_{(0,m)})\epsilon}(1+\tau)^{-2\beta} |\slashed{\D} \mathscr{Y}^{N_2} h|_{LL}
		\\
		&\phantom{\lesssim \int_{\mathcal{M}_{\tau}^{\tau_1}} \bigg(}
		+ \epsilon^3 (1+r)^{-1-\frac{1}{2}C_{[N_2(j), m]}\epsilon}(1+\tau)^{-2\beta} |\overline{\slashed{\D}} \mathscr{Y}^{N_2} h|_{(\text{frame})}^2
		\\
		&\phantom{\lesssim \int_{\mathcal{M}_{\tau}^{\tau_1}} \bigg(}
		+ \epsilon^3 (1+r)^{1-2\delta} (1+\tau)^{-2\beta} |\tilde{\slashed{\Box}}_g \mathscr{Y}^{N_2-1} h |_{(\text{frame})}^2
		+ \epsilon^3 (1+\tau)^{-2\beta} (1+r)^{-3-\delta} |\mathscr{Y}^{N_2} h|^2_{(\text{frame})}
		\\
		&\phantom{\lesssim \int_{\mathcal{M}_{\tau}^{\tau_1}} \bigg(}
		+ \epsilon^3 (1+\tau)^{-2\beta} (1+r)^{1-C_{[N_2]}\epsilon} |\mathscr{Y}^{N_2} F|^2_{LL}
		+ \sum_{j \leq N_2 - 1} \epsilon^5 (1+\tau)^{-2\beta} (1+r)^{1-\frac{1}{2}C_{[N_2]}\epsilon} |\mathscr{Y}^j F|_{(\text{frame})}^2
		\\
		&\phantom{\lesssim \int_{\mathcal{M}_{\tau}^{\tau_1}} \bigg(}
		+ \epsilon^3 (1+r)^{-1-(C_{[N_2]} - 2C_{[N_1]})\epsilon} (1+\tau)^{C_{(N_1)}\delta}|\slashed{\D}_T^{j+1} \mathscr{Y}^{N_2-j-2} \tr_{\slashed{g}}\chi_{(\text{small})}|^2
		\\
		&\phantom{\lesssim \int_{\mathcal{M}_{\tau}^{\tau_1}} \bigg(}
		+ \epsilon^3 (1+r)^{-3-2\delta} (1+\tau)^{-2\beta} |\mathscr{Y}^{N_2} \log \mu|^2
		+ \epsilon^3 (1+r)^{-3-2\delta} (1+\tau)^{-2\beta} |\mathscr{Y}^n \zeta|^2
		\\
		&\phantom{\lesssim \int_{\mathcal{M}_{\tau}^{\tau_1}} \bigg(}
		+ \epsilon^3 (1+r)^{-1-\frac{1}{2}C_{[N_2]}\epsilon} (1+\tau)^{C_{(N_1)}\delta}|F^{(\leq N_2 - 1)}|_{(\text{frame})}^2
		+ \epsilon^3 (1+r)^{-3-2\delta} (1+\tau)^{-2\beta} |\mathscr{Y}^{N_2} \log \mu|^2
		\\
		&\phantom{\lesssim \int_{\mathcal{M}_{\tau}^{\tau_1}} \bigg(}
		+ \epsilon^3 (1+r)^{-3-\frac{1}{2}C_{[N_2(j), m]}\epsilon} |\mathscr{Y}^{\leq N_2 - 1} \phi|^2
	\bigg) \dVol_g 
	\end{split}
	\end{equation*}

	The terms involving $\mathscr{Y}^{N_2} \zeta$ and $\mathscr{Y}^{N_2} \log \mu$ can be controlled as before, using propositions \ref{proposition L2 zeta high} and \ref{proposition L2 foliation density} respectively. Importantly, every term either has additional decay in $\tau$ (at least $(1+\tau)^{-C^* \delta}$) \emph{or} has additional factors of $\slashed{\D}_T$ \emph{or} is lower order. In the latter two cases there might actually be \emph{growth} in $\tau$, but only at the rate $(1+\tau)^{C_{(N_1)}\delta}$.
	
	Using the $L^2$ bootstrap bounds we therefore have
	\begin{equation*}
	\begin{split}
	&\int_{\mathcal{M}_{\tau}^{\tau_1}} \epsilon^{-1} (1+r)^{1-C_{[N_2(j), m]}\epsilon} |F_{(A, N_2(j))}| \dVol_g
	\lesssim
	\frac{1}{C_{[N_2(j+1),m]}}\epsilon^4 (1+\tau)^{-1 + C_{(N_2(j+1))}\delta}
	\end{split}
	\end{equation*}

	Finally, we need to deal with the case in which $j = N_2$, i.e.\ $\mathscr{Y}^{N_2} = (\slashed{\D}_T)^{N_2}$. The point is that this case generates \emph{better} error terms than the general case.
	
	We claim that, if $\phi_{(A)}$ is a scalar field satisfying $\tilde{\Box}_g \phi_{(A)} = F_{(A)}$, then $(\slashed{\D}_T)^n \phi_{(A)} = T^n \phi_{(A)}$ satisfies
	\begin{equation*}
	\begin{split}
	\tilde{\Box}_g (T^n \phi_{(A)})
	&=
	T^n F_{(A)}
	+ \bm{\Gamma}^{(0)}_{(-1)} (\slashed{\D} T^n \phi_{(A)})
	+ \bm{\Gamma}^{(1)}_{(-1-\delta)} (\slashed{\D} \mathscr{Z}^n h)_{(\text{frame})}
	+ \bm{\Gamma}^{(0)}_{(-1+C_{(0)}\epsilon)} (\overline{\slashed{\D}} \mathscr{Z}^n h)_{(\text{frame})}
	\\
	&\phantom{=}
	+ \sum_{j+k \leq n-1} \bm{\Gamma}^{(j)}_{(-1+C_{(j)}\epsilon)} (F^{(k)})_{(\text{frame})}
	+ \bm{\Gamma}^{(0)}_{(-1-\delta)} \mathscr{Z}^{n} \tr_{\slashed{g}}\chi_{(\text{small})}
	+ \sum_{j+k \leq n} \bm{\Gamma}^{(j)}_{(-1-\delta)} \bm{\Gamma}^{(k+1)}_{(-1+C_{(k+1)}\epsilon)}
	\end{split}
	\end{equation*}
	We can prove this by induction on $n$. First, we note that it is evidently true if $n = 0$. Now, suppose that this is true for some value $n = n_1$. Then we have
	\begin{equation*}
	\begin{split}
	\tilde{\Box}_g (T^{n_1 + 1} \phi_{(A)})
	&=
	[\tilde{\Box}_g, T] T^{n_1} \phi_{(A)}
	+ T^{n_1 + 1} F_{(A)}
	+ \bm{\Gamma}^{(0)}_{(-1)} (\slashed{\D} T^{n_1 + 1} \phi_{(A)})
	+ \bm{\Gamma}^{(1)}_{(-1+C_{(1)}\epsilon)} (\slashed{\D} T^{n_1} \phi_{(A)})
	\\
	&\phantom{=}
	+ \bm{\Gamma}^{(0)}_{(-1)} ([\slashed{\D}, T] T^{n_1} \phi_{(A)})
	+ \bm{\Gamma}^{(0)}_{(-1+C_{(0)}\epsilon)} (\overline{\slashed{\D}} \mathscr{Z}^{n_1 + 1} h)_{(\text{frame})}
	+ \bm{\Gamma}^{(1)}_{(-1+C_{(1)}\epsilon)} (\overline{\slashed{\D}} \mathscr{Z}^{n_1} h)_{(\text{frame})}
	\\
	&\phantom{=}
	+ \bm{\Gamma}^{(0)}_{(-1+C_{(0)}\epsilon)} ([\overline{\slashed{\D}}, \slashed{\D}_T] \mathscr{Z}^{n_1} h)_{(\text{frame})}
	+ \sum_{j+k \leq n_1-1} \bm{\Gamma}^{(j+1)}_{(-1+C_{(j+1)}\epsilon)} (F^{(k)})_{(\text{frame})}
	\\
	&\phantom{=}
	+ \sum_{j+k \leq n_1-1} \bm{\Gamma}^{(j)}_{(-1+C_{(j)}\epsilon)} (F^{(k+1)})_{(\text{frame})}
	+ \sum_{j+k \leq n_1-1} \bm{\Gamma}^{(j)}_{(-1+C_{(j)}\epsilon)} (F^{(k)})_{(\text{frame})} (TX_{(\text{frame})})
	\\
	&\phantom{=}
	+ \bm{\Gamma}^{(1)}_{(-1-\delta)} \mathscr{Z}^{n_1} \tr_{\slashed{g}}\chi_{(\text{small})}
	+ \bm{\Gamma}^{(0)}_{(-1-\delta)} \mathscr{Z}^{n_1+1} \tr_{\slashed{g}}\chi_{(\text{small})}
	+ \sum_{j,k \leq n-1} \bm{\Gamma}^{(j+1)}_{(-1-\delta)} \bm{\Gamma}^{(k+1)}_{(-1+C_{(k+1)}\epsilon)}
	\\
	&\phantom{=}
	+ \sum_{j,k \leq n-1} \bm{\Gamma}^{(j)}_{(-1-\delta)} \bm{\Gamma}^{(k+2)}_{(-1+C_{(k+2)}\epsilon)}
	\end{split}
	\end{equation*}
	Simplifying these expressions proves the inductive step. In particular, we use proposition \ref{proposition inhomogeneous terms commute with T} to compute $[\tilde{\Box}_g, T] T^{n_1} \phi_{(A)}$, and note that the term given there as $\bm{\Gamma}^{(1)}_{(-1-\delta)} (\slashed{\D} \mathscr{Z} \phi)$ is actually $\bm{\Gamma}^{(0)}_{(-1-\delta)} (\slashed{\D} \mathscr{Z} \phi) + r^{-1}\bm{\Gamma}^{(1)}_{(-1+C_{(1)}\epsilon)} (\slashed{\D} \mathscr{Z} \phi)$.
	
	Notice that this time there is no term involving $\slashed{\nabla}^2 \mathscr{Y}^{n-1} \phi$. In particular, this means that we have sufficient decay in $r$ for all coefficients. Now, proceeding exactly as before, we find that
	\begin{equation*}
	\begin{split}
	&\int_{\mathcal{M}_{\tau}^{\tau_1}} \epsilon^{-1} (1+r)^{1-C_{[N_2(N_2), m]}\epsilon} |F_{(A, N_2(j))}| \dVol_g
	\\
	&\lesssim
	\frac{1}{C_{[N_2(N_2)]}} \epsilon^4 \left( (1+\tau)^{-1 + (C_{(N_2(N_2))}-8)\delta} + (1+\tau)^{-1 + (C_{(N_2-1)})\delta} \right)
	\end{split}
	\end{equation*}

	This gives us all of the required bounds on terms of the form $\int_{\mathcal{M}_\tau^{\tau_1}} \epsilon^{-1} (1+r)^{1-C_{[N_2(N_2),m]}\epsilon} |F|^2 \dVol_g$, where we do not split the inhomogeneous term into different pieces. 
	
	When splitting the inhomogeneous terms into different parts, we can place all of the error terms into $F_{(A, n, 2)}$ and $F_{(A,n,4)}$ except for the terms with the superscript $(0)$. In other words, we set
	\begin{equation*}
	\begin{split}
	F_{(A, n, 2)} &:= F_{(A, n)} - F_{(A, n, 2)}^{(0)} - F_{(A, n, 3)}^{(0)} \\
	F_{(A, n, 4)} &:= F_{(A, n)} - F_{(A, n, 4)}^{(0)} - F_{(A, n, 5)}^{(0)}
	\end{split}
	\end{equation*}
	then, following the calculations above it should be clear that, for all $n \leq N_2 - 1$ we have
	\begin{equation*}
	\begin{split}
	&\int_{\mathcal{M}_{\tau}^{\tau_1}} \epsilon^{-1} (1+r)^{1-3\delta}(1+\tau)^{6\delta} |F_{(A, n, 2)}|^2 \dVol_g
	\\
	&\lesssim
	\frac{1}{C_{[N_2(n+1), m]}} \epsilon^{2(N_2 + 2 - n)} \Big(
	(1+\tau)^{-1 + (C_{(n)}-2)\delta}
	+ (1+\tau)^{-1 + (C_{(n-1)} + 6)\delta} \Big)
	\end{split}
	\end{equation*}
	In fact, we have more than enough decay in $r$ this time. Similarly, we have
	\begin{equation*}
	\begin{split}
	&\int_{\mathcal{M}_{\tau}^{\tau_1}} \epsilon^{-1} (1+r)^{1-3\delta}(1+\tau)^{6\delta} |F_{(A, N_2(j), 2)}|^2 \dVol_g
	\\
	&\lesssim
	\frac{1}{C_{[N_2(j+1), m]}} \epsilon^4 \Big(
		(1+\tau)^{-1 + (C_{(N_2(j))}-2)\delta}
		+ (1+\tau)^{-1 + (C_{(N_2(j+1))} + 6)\delta}
		+ (1+\tau)^{-1 + (C_{(N_2-1)} + 6)\delta}
	\Big)
	\end{split}
	\end{equation*}
	and
	\begin{equation*}
	\begin{split}
	&\int_{\mathcal{M}_{\tau}^{\tau_1}} \epsilon^{-1} (1+r)^{1-3\delta}(1+\tau)^{6\delta} |F_{(A, N_2(N_2), 2)}|^2 \dVol_g
	\\
	&\lesssim
	\frac{1}{C_{[N_2(N_2), m]}} \epsilon^4 \Big(
		(1+\tau)^{-1 + (C_{(N_2(N_2))}-2)\delta}
		+ (1+\tau)^{-1 + (C_{(N_2-1)} + 6)\delta}
	\Big)
	\end{split}
	\end{equation*}

	Next, for $n \leq N_2$ we have
	\begin{equation*}
	\begin{split}
	&\int_{\mathcal{M}_{\tau}^{\tau_1}} \epsilon^{-1} (1+r)^{1-C_{[n]}\epsilon}(1+\tau)^{1+\delta} |F_{(A, n, 2)}|^2 \dVol_g
	\\
	&\lesssim
	\frac{1}{C_{[N_2(n+1), m]}}\epsilon^{2(N_2 + 2 - n)} \Big(
	(1+\tau)^{(C_{(n)} - 7)\delta}
	+ (1+\tau)^{(C_{(n-1)} + 1)\delta} \Big)
	\end{split}
	\end{equation*}
	and similarly
	\begin{equation*}
	\begin{split}
	&\int_{\mathcal{M}_{\tau}^{\tau_1}} \epsilon^{-1} (1+r)^{1-C_{[N_2(j)]}\epsilon}(1+\tau)^{1+\delta} |F_{(A, N_2(j), 2)}|^2 \dVol_g
	\\
	&\lesssim
	\frac{1}{C_{[N_2(j+1), m]}}\epsilon^4 \Big(
	(1+\tau_1)^{(C_{(N_2(j))} -7)\delta}
	+ (1+\tau_1)^{(C_{(N_2(j+1))} + 1)\delta}
	+ (1+\tau_1)^{(C_{(N_2-1)} + 1)\delta} \Big)
	\end{split}
	\end{equation*}
	and finally
	\begin{equation*}
	\begin{split}
	&\int_{\mathcal{M}_{\tau}^{\tau_1}} \epsilon^{-1} (1+r)^{1-C_{[N_2(j)]}\epsilon}(1+\tau)^{1+\delta} |F_{(A, N_2(N_2), 2)}|^2 \dVol_g
	\\
	&\lesssim
	\frac{1}{C_{[N_2(N_2), m]}} \epsilon^4 \Big(
	(1+\tau_1)^{(C_{(N_2(N_2))} -7)\delta}
	+ (1+\tau_1)^{(C_{(N_2-1)} + 1)\delta} \Big)
	\end{split}
	\end{equation*}

\end{proof}

\section{Dealing with a point-dependent change of basis}

We also need to consider the inhomogeneous term which appears on the right hand side of an equation for some quantity which is obtained after a point-dependent change of basis, for example, the quantity $h_{LL}$. In general, suppose that the quantities $\phi_{(a)}$ satisfy
\begin{equation*}
\tilde{\slashed{\Box}}_g \phi_{(a)} = F_{(a)}
\end{equation*}
then we define $F_{(A)}$ by
\begin{equation*}
M_{(A)}^{(a)}\tilde{\slashed{\Box}}_g \phi_{(a)} = F_{(A)}
\end{equation*}
where $M_{(A)}^{(a)}$ is the change-of-basis matrix (see section \ref{section energy estimates involving a point-dependent change of basis}).

Note also from section \ref{section energy estimates involving a point-dependent change of basis} that there are additional error terms in this case, if the change of basis matrix $M_{(A)}^{(a)}$ is a function (and not just a constant) on the manifold $\mathcal{M}$. We do not include these error terms in our definition of $F_{(A)}$, because when we commute we do not apply the commutation operators to the change-of-basis matrix: in other words, we set
\begin{equation*}
F_{(A, \mathscr{Z}\phi)} = M_{(A)}^{(a)}\tilde{\slashed{\Box}}_g \mathscr{Z} \phi_{(a)} 
\end{equation*}

From these considerations, it is easy to see that we have the following proposition, which is almost identical to proposition \ref{proposition L2 bounds for F after commuting}.

\begin{proposition}[$L^2$ bounds for the inhomogeneous terms after commuting and a change of basis]
	\label{proposition L2 bounds for F after commuting and change of basis}
	Let $\phi_{(a)}$ be a set of scalar fields satisfying the equations
	\begin{equation*}
	\begin{split}
	\tilde{\Box}_g\phi_{(a)} &= F_{(a,0)}
	\\
	F_{(a,0)} &= F_{(a,0)}^{(0)} + \left(F_{(a,0)}^{(bc)}\right)^{\mu\nu}(\partial_\mu \phi_B)(\partial_\nu \phi_C) + \mathcal{O}\left(\phi (\partial\phi)^2\right)
	\end{split}
	\end{equation*}
	where we further decompose
	\begin{equation*}
	\begin{split}
	F_{(a,0)}^{(0)}
	&=
	F_{(a,0,1)}^{(0)}
	+ F_{(a,0,2)}^{(0)}
	+ F_{(a,0,3)}^{(0)}
	\\
	&=
	F_{(a,0,4)}^{(0)}
	+ F_{(a,0,5)}^{(0)}
	+ F_{(a,0,6)}^{(0)}
	\end{split}	
	\end{equation*}
	and we define (schematically)
	\begin{equation*}
	F^{(0)}_{(a,n)} := \mathscr{Y}^n F^{(0)}_{(a,0)}
	\end{equation*}
	
	Now, let $M_{(A)}^{(a)}$ be a (possibly point-dependent) change-of-basis matrix, and let us define
	\begin{equation*}
	\begin{split}
	F_{(A,0)}^{(0)} &:= M_{(A)}^{(a)} F_{(a,0)}^{(0)}
	\\
	F_{(A,0,n)}^{(0)} &:= M_{(A)}^{(a)} F_{(a,0,n)}^{(0)} \quad \quad \text{for} \quad n = 1, \ldots 6
	\\
	F_{(A,0)}^{(BC)} &:= M_{(A)}^{(a)} (M^{-1})_{(b)}^{(B)} (M^{-1})_{(c)}^{(C)} F_{(a,0)}^{(bc)}
	\end{split}	
	\end{equation*}

	We require that the tensor fields $F_{(A,0)}^{(BC)}$ have \emph{constant rectangular components} (if not, then it should be the case that we can include any non-constant parts in the lower order terms)\footnote{Note that this condition can be weakened: see the footnotes to proposition \ref{proposition L2 bounds for F before commuting}}. Also, we suppose that they satisfy the structural equations
	\begin{equation*}
	\begin{split}
	\left(F_{(A,0)}^{(BC)}\right)^{\mu\nu} 
	&= 	\left(F_{(A,0)}^{(CB)}\right)^{\mu\nu}
	\\
	\left(F_{(A,0)}^{(BC)}\right)_{LL}
	&= 0 \quad \text{if}\quad \phi_{(A)} \in \Phi_{[0]}
	\\
	\left(F_{(A,0)}^{(BC)}\right)_{LL}
	&= 0 \quad \text{if}\quad \phi_{(A)} \in \Phi_{[n]} \text{\, and either } \begin{cases}
	\phi_{(B)} \in \Phi_{[n+1]} \\
	\phi_{(B)} \in \Phi_{[n]} \text{\, and \,} \phi_{(C)} \in \Phi_{[m]} \text{ , } m \geq 1
	\end{cases}
	\end{split}
	\end{equation*}
	
	Suppose moreover that the terms $F^{(0)}_{(A,n)}$ satisfy the following conditions: if $\phi_{(A)} \in \Phi_{[0]}$, then 
	\begin{equation*}
	\begin{split}
	&\int_{\mathcal{M}_{\tau_0}^\tau} \epsilon^{-1} \bigg(
	(1+r)^{1-C_{[0,0]}\epsilon}|F^{(0)}_{(A,0)}|^2
	+ (1+r)^{\frac{1}{2}\delta}(1+\tau)^{1+\delta} |F^{(0)}_{(A,0,1)}|^2
	+ (1+r)^{1-3\delta}(1+\tau)^{2\beta} |F^{(0)}_{(A,0,2)}|^2
	\\
	&\phantom{\int_{\mathcal{M}_{\tau_0}^\tau} \epsilon^{-1} \bigg(}
	+ (1+r)^{1+\frac{1}{2}\delta} |F^{(0)}_{(A,0,3)}|^2
	\bigg)\dVol_g
	\lesssim
	\frac{1}{C_{[0,0]}} \delta^{-1} \epsilon^{2(N_2 + 2)}(1+\tau)^{-1}
	\\
	\\
	&\int_{\mathcal{M}_{\tau_0}^\tau \cap \{r \geq r_0\}} \epsilon^{-1}\bigg(
	r^{1-C_{[0,0]}\epsilon} (1+\tau)^{1+\delta} |F^{(0)}_{(A,0,4)}|^2
	+ r^{2-C_{[0,0]}\epsilon-2\delta} (1+\tau)^{2\beta} |F^{(0)}_{(A,0,5)}|^2
	\\
	&\phantom{\int_{\mathcal{M}_{\tau_0}^\tau \cap \{r \geq r_0\}} \epsilon^{-1}\bigg(}
	+ r^{2-C_{[0.0]}\epsilon} |F^{(0)}_{(A,0,6)}|^2
	\bigg)\dVol_g
	\lesssim
	\frac{1}{C_{[0,0]}} \delta^{-1} \epsilon^{2(N_2 + 2)}
	\end{split}
	\end{equation*}
	and in general, if $\phi_{(A)} \in \Phi_{[m]}$, then
	\begin{equation*}
	\begin{split}
	&\int_{\mathcal{M}_{\tau_0}^\tau} \epsilon^{-1} \bigg(
	(1+r)^{1-C_{[n,m]}\epsilon}|F^{(0)}_{(A,n)}|^2
	+ (1+r)^{\frac{1}{2}\delta}(1+\tau)^{1+\delta} |F^{(0)}_{(A,n,1)}|^2
	+ (1+r)^{1-3\delta}(1+\tau)^{2\beta} |F^{(0)}_{(A,n,2)}|^2
	\\
	&\phantom{\int_{\mathcal{M}_{\tau_0}^\tau} \epsilon^{-1} \bigg(}
	+ (1+r)^{1+\frac{1}{2}\delta} |F^{(0)}_{(A,n,3)}|^2
	\bigg)\dVol_g
	\lesssim
	\frac{1}{C_{[n,m]}} \delta^{-1} \epsilon^{2(N_2 + 2)}(1+\tau)^{-1+C_{(n,m)}\delta}
	\\
	\\
	&\int_{\mathcal{M}_{\tau_0}^\tau \cap \{r \geq r_0\}} \epsilon^{-1}\bigg(
	r^{1-C_{[n,m]}\epsilon} (1+\tau)^{1+\delta} |F^{(0)}_{(A,n,4)}|^2
	+ r^{2-C_{[n,m]}\epsilon-2\delta} (1+\tau)^{2\beta} |F^{(0)}_{(A,n,5)}|^2
	\\
	&\phantom{\int_{\mathcal{M}_{\tau_0}^\tau \cap \{r \geq r_0\}} \epsilon^{-1}\bigg(}
	+ r^{2-C_{[n,m]}\epsilon} |F^{(0)}_{(A,n,6)}|^2
	\bigg)\dVol_g
	\lesssim
	\frac{1}{C_{[0,m]}} \delta^{-1} \epsilon^{2(N_2 + 2)}(1+\tau)^{C_{(n,m)}\delta}
	\end{split}
	\end{equation*}

	Furthermore, suppose that both the pointwise bounds \emph{and} the $L^2$ bounds of chapter \ref{chapter bootstrap} hold.
	
	Define $F_{(A,n)}$ as follows: if $\mathscr{Y}^n$ contains no factors of the operator $r\slashed{\D}_L$, then we define
	\begin{equation*}
	\tilde{\slashed{\Box}}_g \mathscr{Y}^n \phi_{(A)} = F_{(A,n)}
	\end{equation*}
	otherwise, if $\mathscr{Y}^n$ contains $k$ factors of the operator $r\slashed{\D}_L$, $k \geq 1$, then we define
	\begin{equation*}
	\tilde{\slashed{\Box}}_g \mathscr{Y}^n \phi_{(A)} 
	- k\slashed{\Delta} \mathscr{Y}^{n-1} \phi_{(A)}
	- (2^k - 1)r^{-1} \slashed{\D}_L (r\slashed{\D}_L \mathscr{Y}^{n-1}\phi)
	- (2^k - 1)r^{-1} \slashed{\D}_L (\mathscr{Y}^{n-1} \phi)
	= F_{(A,n)}
	\end{equation*}

	Then, for all sufficiently small $\epsilon$, we can decompose $F_{(A,n)}$ as
	\begin{equation*}
	\begin{split}
	F_{(A,n)}
	&=
	F_{(A,n,1)} + F_{(A,n,2)} + F_{(A,n,3)}
	\\
	&=
	F_{(A,n,4)} + F_{(A,n,5)} + F_{(A,n,6)}
	\end{split}	
	\end{equation*}
	where, if $\phi_{(A)} \in \Phi_{[m]}$, then we have
	\begin{equation*}
	\begin{split}
	&\int_{\mathcal{M}_{\tau_0}^\tau} \epsilon^{-1} \bigg(
	(1+r)^{1-C_{[n,m]}\epsilon}|F_{(A,n)}|^2
	+ (1+r)^{\frac{1}{2}\delta}(1+\tau)^{1+\delta} |F_{(A,n,1)}|^2
	+ (1+r)^{1-3\delta}(1+\tau)^{2\beta} |F_{(A,n,2)}|^2
	\\
	&\phantom{\int_{\mathcal{M}_{\tau_0}^\tau} \epsilon^{-1} \bigg(}
	+ (1+r)^{1+\frac{1}{2}\delta} |F_{(A,n,3)}|^2
	\bigg)\dVol_g
	\lesssim
	\left( \frac{1}{C_{[n,m]}} + \frac{\epsilon^2}{\delta^6} \right) \delta^{-1} \epsilon^{2(N_2 + 2 - n)} (1+\tau)^{-1+C_{(n)}\delta}
	\\
	\\
	&\int_{\mathcal{M}_{\tau_0}^\tau \cap \{r \geq r_0\}} \epsilon^{-1}\bigg(
	r^{1-C_{[n,m]}\epsilon} (1+\tau)^{1+\delta} |F_{(A,n,4)}|^2
	+ r^{2-C_{[n,m]}\epsilon-2\delta} (1+\tau)^{2\beta} |F_{(A,n,5)}|^2
	\\
	&\phantom{\int_{\mathcal{M}_{\tau_0}^\tau \cap \{r \geq r_0\}} \epsilon^{-1}\bigg(}
	+ r^{2-C_{[n,m]}\epsilon} |F_{(A,n,6)}|^2
	\bigg)\dVol_g
	\lesssim
	\left( \frac{1}{C_{[n,m]}} + \frac{\epsilon^2}{\delta^6} \right) \delta^{-1} \epsilon^{2(N_2 + 2 - n)} (1+\tau)^{C_{(n)}\delta}
	\end{split}
	\end{equation*}
	
	Finally, we define $F_{(A,N_2(j))}$ as follows: if $\mathscr{Y}^{N_2 - j}$ contains $k$ factors of the operator $r\slashed{\D}_L$, $k \geq 1$, then we define
	\begin{equation*}
	\begin{split}
	&\tilde{\slashed{\Box}}_g (\slashed{\D}_T)^j \mathscr{Y}^{N_2 - j} \phi_{(A)} 
	- k\slashed{\Delta} (\slashed{\D}_T)^j \mathscr{Y}^{N_2 - 1 - j} \phi_{(A)}
	- (2^k - 1)r^{-1} \slashed{\D}_L (r\slashed{\D}_L (\slashed{\D}_T)^j \mathscr{Y}^{N_2 - 1 - j}\phi)
	\\
	&
	- (2^k - 1)r^{-1} \slashed{\D}_L (\slashed{\D}_T)^j \mathscr{Y}^{N_2 - 1 - j} \phi)
	= F_{(A,N_2(j))}
	\end{split}
	\end{equation*}
	then we have
	\begin{equation*}
	\begin{split}
	&\int_{\mathcal{M}_{\tau_0}^\tau} \epsilon^{-1} \bigg(
	(1+r)^{1-C_{[N_2(j),m]}\epsilon}|F_{(A,N_2(j))}|^2
	+ (1+r)^{\frac{1}{2}\delta}(1+\tau)^{1+\delta} |F_{(A,N_2(j),1)}|^2
	\\
	&\phantom{\int_{\mathcal{M}_{\tau_0}^\tau} \epsilon^{-1} \bigg(}
	+ (1+r)^{1-3\delta}(1+\tau)^{2\beta} |F_{(A,N_2(j),2)}|^2
	+ (1+r)^{1+\frac{1}{2}\delta} |F_{(A,N_2(j),3)}|^2
	\bigg)\dVol_g
	\\
	&\lesssim
	\left( \frac{1}{C_{[N_2(j+1),m]}} + \frac{\epsilon^2}{\delta^6} \right) \delta^{-1} \epsilon^4 (1+\tau)^{-1+C_{(N_2(j))}\delta}
	\\
	\\
	&\int_{\mathcal{M}_{\tau_0}^\tau \cap \{r \geq r_0\}} \epsilon^{-1}\bigg(
	r^{1-C_{[N_2(j),m]}\epsilon} (1+\tau)^{1+\delta} |F_{(A,[N_2(j),4)}|^2
	+ r^{2-C_{[N_2(j),m]}\epsilon-2\delta} (1+\tau)^{2\beta} |F_{(A,[N_2(j),5)}|^2
	\\
	&\phantom{\int_{\mathcal{M}_{\tau_0}^\tau \cap \{r \geq r_0\}} \epsilon^{-1}\bigg(}
	+ r^{2-C_{[N_2(j),m]}\epsilon} |F_{(A,[N_2(j),6)}|^2
	\bigg)\dVol_g
	\lesssim
	\left( \frac{1}{C_{[N_2(j+1),m]}} + \frac{\epsilon^2}{\delta^6} \right) \delta^{-1} \epsilon^4 (1+\tau)^{C_{(N_2(j))}\delta}
	\end{split}
	\end{equation*}
	and
	\begin{equation*}
	\begin{split}
	&\int_{\mathcal{M}_{\tau_0}^\tau} \epsilon^{-1} \bigg(
	(1+r)^{1-C_{[N_2(j),m]}\epsilon}|F_{(A,N_2(N_2))}|^2
	+ (1+r)^{\frac{1}{2}\delta}(1+\tau)^{1+\delta} |F_{(A,N_2(N_2),1)}|^2
	\\
	&\phantom{\int_{\mathcal{M}_{\tau_0}^\tau} \epsilon^{-1} \bigg(}
	+ (1+r)^{1-3\delta}(1+\tau)^{2\beta} |F_{(A,N_2(N_2),2)}|^2
	+ (1+r)^{1+\frac{1}{2}\delta} |F_{(A,N_2(N_2),3)}|^2
	\bigg)\dVol_g
	\\
	&\lesssim
	\left( \frac{1}{C_{[N_2(N_2),m]}} + \frac{\epsilon^2}{\delta^6} \right) \delta^{-1} \epsilon^4 (1+\tau)^{-1+C_{(N_2(N_2))}\delta}
	\\
	\\
	&\int_{\mathcal{M}_{\tau_0}^\tau \cap \{r \geq r_0\}} \epsilon^{-1}\bigg(
	r^{1-C_{[N_2(N_2),m]}\epsilon} (1+\tau)^{1+\delta} |F_{(A,[N_2(N_2),4)}|^2
	\\
	&\phantom{\int_{\mathcal{M}_{\tau_0}^\tau \cap \{r \geq r_0\}} \epsilon^{-1}\bigg(}
	+ r^{2-C_{[N_2(N_2),m]}\epsilon-2\delta} (1+\tau)^{2\beta} |F_{(A,[N_2(N_2),5)}|^2
	+ r^{2-C_{[N_2(N_2),m]}\epsilon} |F_{(A,[N_2(N_2),6)}|^2
	\bigg)\dVol_g
	\\
	&\lesssim
	\left( \frac{1}{C_{[N_2(N_2),m]}} + \frac{\epsilon^2}{\delta^6} \right) \delta^{-1} \epsilon^4 (1+\tau)^{C_{(N_2(N_2))}\delta}
	\end{split}
	\end{equation*}
	
\end{proposition}

\chapter{Proving the theorem}
\label{chapter proving the theorem}

Now we are finally in a position to prove a detailed version of theorem \ref{theorem main theorem rough statement}, which is the main result of this work.

\begin{theorem}[Global existence for small initial data for systems of wave equations obeying the hierarchical weak null condition]
	\label{theorem main theorem}
	
	Let $\mathcal{M}$ be a smooth, four-dimensional Lorentzian manifold, equipped with a set of smooth coordinate functions (which we call the \emph{rectangular coordinates}) $(x^0, x^1, x^2, x^3)$ with respect to which the metric $g$ has components
	\begin{equation*}
	g(\partial_{x^a}, \partial_{x^b})
	= g_{ab}
	= m_{ab} + h_{ab}
	\end{equation*}
	where $m_{ab}$ are the rectangular components of the Minkowksi metric: that is,
	\begin{equation*}
	\begin{split}
	m_{00} &= -1 \\
	m_{11} = m_{22} = m_{33} &= 1 \\
	m_{ab} &= 0 \quad \quad \text{if} \quad a \neq b
	\end{split}	
	\end{equation*}

	Let $\mathcal{V}$ be a finite-dimensional, trivial vector bundle over $\mathcal{M}$, so that $\mathcal{V} = \mathcal{M} \times V$ for some finite-dimensional vector space $V$. Let $(v^{(1)}, v^{(2)}, \ldots)$ be a basis of smooth sections for $\mathcal{V}$, so that, for every point $x \in \mathcal{M}$, the vectors $(v^{(1)}(x), v^{(2)}(x), \ldots)$ span the vector space $V$.
	
	Let $\bm{\phi} = \phi_{(a)} v^{(a)}$ be a section of $\mathcal{V}$, such that the components of $\phi_{(a)}$ satisfy the system of wave equations
	\begin{equation*}
	\tilde{\Box}_g \phi_{(a)} = F_{(a)}(x, \bm{\phi}, \partial\bm{\phi})
	\end{equation*}
	where $\tilde{\Box}_g$ is the operator\footnote{See chapter \ref{chapter preliminaries}, and in particular definition \ref{definition L and Lbar} for the definition of the vector field $\Lbar$, and see chapter \ref{chapter null frame connection coefficients}, definition \ref{definition omega} for the definition of the scalar field $\omega$.}
	\begin{equation*}
	\tilde{\Box}_g \phi := \Box_g \phi + \omega \Lbar \phi
	\end{equation*}
	Here, $\Box_g = (g^{-1})^{\mu\nu} \D_\mu \D_\nu $ is the ``geometric'' wave operator, where $\D$ is the covariant derivative with respect to the metric $g$.

	Let $M$ be a section of the vector bundle with fibres that are invertible linear maps\footnote{Note, in particular, that $M$ could be the identity map, in which case we say that we do not need to change the basis of sections for $\mathcal{V}$. Alternatively, $M$ could be a \emph{constant} linear map (i.e.\ independent of the point on the manifold $\mathcal{M}$), in which case we say that we do not need to perform a \emph{point-dependent} change of basis.} from $V$ to $V$. Then we can realise $M$ as a section of the vector bundle $\mathcal{M} \times V \times V^*$,  where $V^*$ is the dual to $V$. We write this section as $M_{(A)}^{(a)}(x, \bm{\phi})$. Suppose that the components of $M$ with respect to the basis $(v^{(1)}, v^{(2)}, \ldots)$ are given by
	\begin{equation*}
	M_{(A)}^{\phantom{(A)}(a)}(x, \bm{\phi})
	=
	M_{(A)}^{\phantom{(A)}(a)}(x, \bm{\phi})
	\end{equation*}
	where, \emph{if the pointwise bounds in section \ref{section pointwise bootstrap}} hold, then these components satisfy the pointwise bounds\footnote{see chapter \ref{chapter preliminaries}, section \ref{section null frame} for the definitions of the norms $|\partial \phi|$ and $|\bar{\partial} \phi|$.}
	\begin{equation}
	\begin{split}
	|M_{(A)}^{\phantom{(A)}(a)}| &\lesssim 1 \\
	|\bar{\partial} M_{(A)}^{\phantom{(A)}(a)}| &\lesssim (1+r)^{-1 - \delta} \\
	|\partial M_{(A)}^{\phantom{(A)}(a)}| &\lesssim (1+r)^{-1} + (1+r)^{-1 + \delta}(1+\tau)^{-\beta}
	\end{split}
	\end{equation}
	for some $\frac{1}{2} > \beta > \delta > 0$

	Define
	\begin{equation*}
	F_{(A)}(x, \bm{\phi}, \partial\bm{\phi}) := M_{(A)}^{\phantom{(A)}(a)}(x, \bm{\phi}) F_{(a)}(x, \bm{\phi}, \partial\bm{\phi})
	\end{equation*}
	
	Suppose that the fields $F_{(A)}$ can be decomposed as
	\begin{equation*}
	F_{(A)} = F_{(A)}^{(0)}(x) + \left(F_{(A)}^{(BC)}(x)\right)^{\mu\nu}(\partial_\mu \phi)_B (\partial_\nu \phi)_C + \mathcal{O}\left(\phi (\partial\phi)^2\right) + \mathcal{O}\left( (\partial\phi)^3\right)
	\end{equation*}
	where the tensor fields $F_{(A)}^{(BC)}$ have \emph{constant rectangular components}\footnote{Note that this condition can be weakened: see the footnotes to proposition \ref{proposition L2 bounds for F before commuting}}.
	
	Next, we suppose that the vector space $V$ admits a decomposition as a direct sum
	\begin{equation*}
	V = V_{[1]} \oplus V_{[2]} \oplus \ldots \oplus V_{[N]}
	\end{equation*}
	such that each of the basis sections $v^{(A)}(x) = (M^{-1})_{(a)}^{\phantom{(a)}(A)} v^{(a)}(x)$ is in precisely \emph{one} of the subspaces $V_{[n]}$ (in particular, $v^{(A)}(x)$ is in the same subspace regardless of the point $x$ on the manifold). We write $\phi_{(A)} \in \Phi_{[n]}$ if $\phi_{(A)}(x) v^{(A)}(x) \in V_{[n]}$ for all $x \in \mathcal{M}$.
	
	For any tensor field $F^{\mu\nu}$ we define $F_{LL} := F^{\mu\nu} L_\mu L_{\nu}$ (again, see definition \ref{definition L and Lbar} for the definition of the vector field $L$). Then we require that the tensor fields $(F_{(A)}^{(BC)})^{\mu\nu}$ satisfy the structural equations
	\begin{equation*}
	\begin{split}
	\left(F_{(A)}^{(BC)}\right)^{\mu\nu} 
	&= 	\left(F_{(A)}^{(CB)}\right)^{\mu\nu}
	\\
	\left(F_{(A)}^{(BC)}\right)_{LL}
	&= 0 \quad \text{if}\quad \phi_{(A)} \in \Phi_{[0]}
	\\
	\left(F_{(A)}^{(BC)}\right)_{LL}
	&= 0 \quad \text{if}\quad \phi_{(A)} \in \Phi_{[n]} \text{\, and either } \begin{cases}
	\phi_{(B)} \in \Phi_{[n+1]} \\
	\phi_{(B)} \in \Phi_{[n]} \text{\, and \,} \phi_{(C)} \in \Phi_{[m]} \text{ , } m \geq 1
	\end{cases}
	\end{split}
	\end{equation*}
	which we will refer to as the \emph{hierarchical weak null condition}.

	We require the following bounds for the inhomogeneous terms which are independent of the field $\bm{\phi}$ and its derivatives. Define $F^{(0)}_{(A, n)} = \mathscr{Y}^n F_{(A)}^{(0)}$, where each operator $\mathscr{Y}$ is any operator from the set $\{ \slashed{\D}_T , r\slashed{\nabla}, r\slashed{\D}_L \}$ (see chapter \ref{chapter preliminaries} for this notation). Then we require the following bounds for these inhomogeneous terms $F_{(A)}^{(0)}$:
	
	Pointwise bounds: for all $n \leq N_1$, if $\phi_{(A)} \in \Phi_{[m]}$ then we have
	\begin{equation*}
	|F^{(0)}_{(A,n)}| \leq \epsilon^{(N_2 - 2 - n_0)} (1+r)^{-2 + 2C_{(n_0, m)}\epsilon} (1+\tau)^{-\beta}
	\end{equation*}
	for some sufficiently large constants $C_{(n,m)}$, which satisfy
	\begin{equation*}
	C_{(n, m)} \gg C_{(n, m-1)}
	\end{equation*}
	and for all $m_1, m_2$
	\begin{equation*}
	C_{(n, m_1)} \gg C_{(n-1, m_2)}
	\end{equation*}
	
	We also require the following special pointwise bound: if $\phi_{(A)} \in \Phi_{[0]}$ then
	\begin{equation*}
	|F^{(0)}_{(A,0)}| = |F^{(0)}_{(A)}| \leq \epsilon^{(N_2 - 2)} (1+r)^{-2} (1+\tau)^{-\beta}
	\end{equation*}

	We also require some $L^2$ bounds on these homogeneous terms. In order to state these, we first decompose the homogeneous terms as follows:
	\begin{equation*}
	\begin{split}
	F_{(A,n)}^{(0)}
	&=
	F_{(A,n,1)}^{(0)}
	+ F_{(A,n,2)}^{(0)}
	+ F_{(A,n,3)}^{(0)}
	\\
	&=
	F_{(A,n,4)}^{(0)}
	+ F_{(A,n,5)}^{(0)}
	+ F_{(A,n,6)}^{(0)}
	\end{split}	
	\end{equation*}
	then we require the $L^2$ bounds	
	\begin{equation*}
	\begin{split}
	&\int_{\mathcal{M}_{\tau_0}^\tau} \epsilon^{-1} \bigg(
	(1+r)^{1-C_{[n,m]}\epsilon}|F^{(0)}_{(A,n)}|^2
	+ (1+r)^{\frac{1}{2}\delta}(1+\tau)^{1+\delta} |F^{(0)}_{(A,n,1)}|^2
	+ (1+r)^{1-3\delta}(1+\tau)^{2\beta} |F^{(0)}_{(A,n,2)}|^2
	\\
	&\phantom{\int_{\mathcal{M}_{\tau_0}^\tau} \epsilon^{-1} \bigg(}
	+ (1+r)^{1+\frac{1}{2}\delta} |F^{(0)}_{(A,n,3)}|^2
	\bigg)\dVol_g
	\lesssim
	\frac{1}{C_{[n,m]}} \epsilon^{2(N_2 - n + 2)}(1+\tau)^{-1+C_{(n,m)}\delta}
	\\
	\\
	&\int_{\mathcal{M}_{\tau_0}^\tau \cap \{r \geq r_0\}} \epsilon^{-1}\bigg(
	r^{1-C_{[n,m]}\epsilon} (1+\tau)^{1+\delta} |F^{(0)}_{(A,n,4)}|^2
	+ r^{2-C_{[n,m]}\epsilon-2\delta} (1+\tau)^{2\beta} |F^{(0)}_{(A,n,5)}|^2
	\\
	&\phantom{\int_{\mathcal{M}_{\tau_0}^\tau \cap \{r \geq r_0\}} \epsilon^{-1}\bigg(}
	+ r^{2-C_{[n,m]}\epsilon} |F^{(0)}_{(A,n,6)}|^2
	\bigg)\dVol_g
	\lesssim
	\frac{1}{C_{[n,m]}} \epsilon^{2(N_2 - n + 2)}(1+\tau)^{C_{(n,m)}\delta}
	\end{split}
	\end{equation*}
	where $\delta$ is some sufficiently small constant and the $C_{[n,m]}$ are some sufficiently large constants, satisfying
	\begin{equation*}
	C_{[n, m]} \gg C_{[n, m-1]}
	\end{equation*}
	and, for all $m_1, m_2$
	\begin{equation*}
	C_{[n, m_1]} \gg C_{[n-1, m_2]}
	\end{equation*}
	and finally, for all $n_1$, $n_2$, $m_1$, $m_2$
	\begin{equation*}
	C_{[n_1, m_1]} \gg C_{(n_2, m_2)}
	\end{equation*}

	Next, suppose that the rectangular components of the metric components can be expressed as
	\begin{equation*}
	h_{ab}(x, \bm{\phi}) = h^{(0)}_{ab}(x) + h^{(1)}_{ab}(x, \bm{\phi})
	\end{equation*}
	such that the following bounds hold: for all $n \leq N_1$
	\begin{equation*}
	\begin{split}
	|\mathscr{Y}^n h^{(1)}_{ab}| 
	&\lesssim
	\sum_{m \leq n} \sum_{(a)} |\mathscr{Y}^m \phi_{(a)}| + \mathcal{O}(|\mathscr{Y}^{\leq n} \bm{\phi}|^2)
	\\
	|\slashed{\D} \mathscr{Y}^n h^{(1)}_{ab}| 
	&\lesssim
	\sum_{m \leq n} \sum_{(a)} |\slashed{\D} \mathscr{Y}^m \phi_{(a)}| + \mathcal{O}\left( \sum_{j+k \leq n}|\slashed{\D}\mathscr{Y}^{j} \bm{\phi}||\mathscr{Y}^{k} \bm{\phi}| \right)
	\\
	|\overline{\slashed{\D}} \mathscr{Y}^n h^{(1)}_{ab}| 
	&\lesssim
	\sum_{m \leq n} \sum_{(a)} |\overline{\slashed{\D}} \mathscr{Y}^m \phi_{(a)}| + \mathcal{O}\left( \sum_{j+k \leq n}|\slashed{\D}\mathscr{Y}^{j} \bm{\phi}||\mathscr{Y}^{k} \bm{\phi}| \right)
	\end{split}
	\end{equation*}
	(see chapter \ref{chapter preliminaries} for the notation used here) and also such that we have the following bound
	\begin{equation*}
	|\partial h^{(1)}|_{LL} = |\partial h^{(1)}_{ab}|L^a L^b \lesssim \sum_{(A) \, | \, \phi_{(A)} \in \Phi_{[0]}} |\partial \phi|_{(A)} + \mathcal{O}(|\bm{\phi}||\partial \bm{\phi}|)
	\end{equation*}
	Additionally, the lower order terms in the metric perturbations are required to satisfy the following pointwise bounds, for all $n \leq N_1$
	\begin{equation*}
	\begin{split}
	|\mathscr{Y}^n h^{(0)}_{ab}| &\leq \frac{1}{2}\epsilon (1+r)^{-\frac{1}{2} + \delta}
	\\
	|\slashed{\D}\mathscr{Y}^n h^{(0)}_{ab}| &\leq \begin{cases}
	\frac{1}{2}\epsilon \left( (1+r)^{-1} + (1+r)^{-1+\delta} (1+\tau)^{-\beta} \right) \\
	\frac{1}{2}\epsilon (1+r)^{-1+C_{(n)}\epsilon}
	\end{cases}
	\\
	|\slashed{\D}\mathscr{Y}^n h^{(0)}_{ab}| &\leq \frac{1}{2}\epsilon (1+r)^{-1-\delta}
	\end{split}
	\end{equation*}
	where $C_{(n)}$ is some sufficiently large constant satisfying
	\begin{equation*}
	C_{(n)} \geq C_{(n,m)} \quad \text{for all } m
	\end{equation*}
	We also require the following bounds, giving additional control over lower order terms:
	\begin{equation*}
	\begin{split}
	|\partial h^{(0)}_{ab}| &\leq \frac{1}{2}\epsilon (1+r)^{-1 + C_{(0)}}(1+\tau)^{-C^* \delta}
	\\
	|\slashed{\D} \mathscr{Y} h^{(0)}_{ab}| &\leq \frac{1}{2}\epsilon (1+r)^{-1 + C_{(1)}}(1+\tau)^{-C^* \delta} 
	\\
	|\partial h^{(0)}_{ab}|L^a L^b &\leq \frac{1}{2} \epsilon (1+r)^{-1}
	\end{split}
	\end{equation*}

	Finally, we suppose that the initial data for the fields $\phi_{(a)}$ is posed on the hypersurface $\Sigma_{\tau_0}$, which consists of two parts: the hypersurface\footnote{Here, as elsewhere, $r$ is defined relative to the \emph{rectangular coordinates} by $r = \sqrt{(x^1)^2 + (x^2)^2 + (x^3)^2}$.} $\{ x^0 = t_0 = \text{constant}\} \cap \{r \leq r_0\}$ together with an outgoing \emph{characteristic hypersurface} emanating from the sphere $r = r_1$, $t = t_0$ and normal to this sphere.
	
	The initial data is required to satisfy the following bounds: for $\phi_{(A)} \in \Phi_{[m]}$ and for all $n \leq N_2$, we have
	\begin{equation*}
	\int_{\Sigma_{\tau_0}} \bigg(
	(1+r)^{-C_{[n,m]}\epsilon} |\overline{\slashed{\D}} \mathscr{Y}^n \phi|^2_{(A)} 
	+ (1+r)^{1-C_{[n,m]}\epsilon} |\slashed{\D}_L \mathscr{Y}^n \phi|^2_{(A)} 
	\bigg)\dVol_g \leq \epsilon^{2(N_2 + 3 - n)}
	\end{equation*}
	as well as the pointwise bounds
	\begin{equation*}
	\int_{\bar{S}_{t,r}} |\mathscr{Y}^n \phi|_{(A)} \dVol_{\mathbb{S}^2} \lesssim \epsilon^{2(N_2 + 3 - n)} (t + 1 - t_0)^{-1 + \frac{1}{2}C_{[n,m]}\epsilon}
	\end{equation*}
	Again, see chapter \ref{chapter preliminaries} for the definitions of the volume forms and the spheres $\bar{S}_{t,r}$.

	Then, if $N_2 \geq 8$ and $N_2 - 4 \geq N_1 \geq 4$, for all sufficiently small $\epsilon$ the system of wave equations $\tilde{\Box}_g \phi_{(a)} = F_{(a)}$ has a \emph{unique, global solution}, i.e.\ a unique solution in the region to the future of $\Sigma_{\tau_0}$. Furthermore, this solution will obey the pointwise bounds and $L^2$ bounds given in chapter \ref{chapter bootstrap}, as well as the $L^2$ bounds\footnote{In fact, the system will obey all of these bounds with an additional factor of, say $1/2$ on the right hand side.} in section \ref{section L2 bounds for geometric error terms}.

\end{theorem}

\begin{proof}
	This is the main result of this work, and as such we will essentially be appealing to all of the results that we have proved so far.
	
	We proceed under the bootstrap assumptions of chapter \ref{chapter bootstrap}. We also make the bootstrap assumptions on the $L^2$ norms of various geometric quantities, as given in section \ref{section L2 bounds for geometric error terms}. These latter bootstrap bounds are to be regarded as, in a sense, auxiliary: we will very quickly be able to improve these bounds by using the other bootstrap bounds.
	
	Let $\tau_{(\text{max})}$ be the largest time such that these bootstrap bounds hold. We will first establish that that all of the bootstrap bounds hold for some (potentially very short) time, by using the conditions on the initial data.
	
	First, we note that the required $L^2$ bounds on the fields $\bm{\phi}$ hold, at least up to some time $\tau_1 > \tau_0$, by the semi-global existence result (see appendix \ref{appendix semi-global existence}). Together with the relationship between the metric fields $h_{ab}$ and the fields $\bm{\phi}$ that is assumed in the proposition, this implies that the remaining $L^2$ bootstrap bounds from section \ref{section L2 bounds for geometric error terms} hold up to time $\tau_1 > \tau_0$. Then, using the bounds established in section \ref{section L2 bounds for geometric error terms} we find that the $L^2$ bootstrap bounds for the various geometric terms (see section \eqref{equation bootstrap L2 geometric} and \eqref{equation bootstrap L2 geometric top order}) also hold for some short time.
	
	Next, we note that the semi-global existence result of appendix \ref{appendix semi-global existence} also provides us with pointwise bounds on the fields $\bm{\phi}$. Together with the conditions on the initial data, this can be used to show that the pointwise bootstrap bounds on the fields $\bm{\phi}$ hold, at least up to some time $\tau_1 > \tau_0$. Again, noting the relationships between the fields $\bm{\phi}$ and the metric components $h_{ab}$, we conclude that the required bounds on the metric components also hold for this short time. Finally, the calculations in chapter \ref{chapter pointwise bounds} show that the pointwise bounds we require on the other geometric quantities also hold for this short time.
	
	In summary, we have established that the time $\tau_{(\text{max})} > \tau_0$. In other words, there is some small time during which all of the bootstrap bounds hold.

	Next, we aim to show that, in fact, at all times $\tau$ such that $\tau_0 \leq \tau \leq \tau_{(\text{max})}$, the bootstrap bounds can be improved. That is, the bootstrap bounds which we have assumed actually hold with an additional factor of, say, $1/2$ on the right hand side. Since we are restricting to the region $\tau_0 \leq \tau \leq \tau_{(\text{max})}$, we can attempt this task while assuming that the bootstrap bounds hold (without, of course, assuming the improvement of a factor of $1/2$!).
	
	\vspace{4mm}
	
	\textbf{Pointwise bounds on the inhomogeneous terms}
	\vspace{3mm}
	
	Many of the bounds we need to improve involve bounds on the inhomogeneous terms in the various wave equations we are considering. The structure of these terms, together with some $L^2$ bounds relating them to other quantities, were established in chapter \ref{chapter bounds for the inhomogeneous term}. Using these expressions, together with the bootstrap bounds, we can establish various pointwise and $L^2$ bounds for the inhomogeneous terms. 
	
	Specifically, using proposition \ref{proposition inhomogeneous terms after commuting n times with Y} together with the pointwise bootstrap bounds, we easily conclude that
	\begin{equation*}
	|\mathscr{Y}^n F|_{(A)}
	\lesssim
	\epsilon^2 (1+r)^{-2+2C_{(n,m)}\epsilon}
	\end{equation*}
	if $\phi_{(A)} \in \Phi_{[m]}$ and $n \leq N_1 - 1$.

	\vspace{4mm}

	\textbf{Improving the} $\bm{L^2}$ \textbf{bounds on geometric quantities}
	
	\vspace{3mm}
	
	The first (and easiest) of the bootstrap bounds which we can improve are the $L^2$ bootstrap bounds for the geometric quantities, which are stated in section \ref{section L2 bounds for geometric error terms}.
	
	If we combine propositions \ref{proposition L2 bounds rectangular}, \ref{proposition L2 omega low}, \ref{proposition L2 zeta low}, \ref{proposition L2 tr chi low}, \ref{proposition L2 tr chibar low}, \ref{proposition L2 foliation density}, \ref{proposition L2 chihat low} and \ref{proposition L2 chibar hat low} then we obtain the following bound: for all $n \leq N_2 - 1$, and for small enough $\epsilon$ and $\delta$,
	\begin{equation*}
	\begin{split}
	&\int_{\mathcal{M}_\tau^{\tau_1}} \bigg(
	C_{[n]}\epsilon (1+r)^{-1-C_{[n]}\epsilon} |\bm{\Gamma}^{(n)}_{(-1+C_{(n)}\epsilon)}|^2
	+ c\delta (1+r)^{-1-c\delta} |\bm{\Gamma}^{(n)}_{(-1+C_{(n)}\epsilon)}|^2
	\\
	&\phantom{\int_{\mathcal{M}_\tau^{\tau_1}} \bigg(}
	+ \delta (1+r)^{-1+(\frac{1}{2} - c_{[n]})\delta} |\bm{\Gamma}^{(n)}_{(-1-\delta)}|^2
	\bigg)\dVol_g
	\\
	&\lesssim
	\left( \frac{\epsilon}{\delta C_{[n]}} + \frac{1}{C_{[n-1]}C_{[n]}} + \frac{\epsilon^2}{c^3 \delta^2} + \frac{\epsilon}{c^3 \delta^3 C_{[n-1]}} \right) \epsilon^{2(N_2 + 1 - n)}(1+\tau)^{-1+C_{(n)}\delta}
	\\
	&\phantom{\lesssim}
	+ \frac{1}{C_{[n]}\epsilon} \int_{\mathcal{M}_\tau^{\tau_1}} \bigg(
	(1+r)^{-1+\frac{1}{2}\delta} |\overline{\slashed{\D}} \mathscr{Y}^n h|^2_{(\text{frame})}
	+ (1+r)^{-3-\delta} |\mathscr{Y}^n h_{(\text{rect})}|^2
	\\
	&\phantom{\lesssim + \frac{1}{C_{[n]}\epsilon} \int_{\mathcal{M}_\tau^{\tau_1}} \bigg(}
	+ C_{[n]} \epsilon^3 \delta (1+r)^{-1+\frac{1}{2}\delta} |\overline{\slashed{\D}} \mathscr{Y}^{n+1} h^2_{(\text{frame})}
	\bigg) \dVol_g
	\end{split}
	\end{equation*}
	so, using the bootstrap bounds for the metric fields we have
	\begin{equation*}
	\begin{split}
	&\int_{\mathcal{M}_\tau^{\tau_1}} \bigg(
	C_{[n]}\epsilon (1+r)^{-1-C_{[n]}\epsilon} |\bm{\Gamma}^{(n)}_{(-1+C_{(n)}\epsilon)}|^2
	+ c\delta (1+r)^{-1-c\delta} |\bm{\Gamma}^{(n)}_{(-1+C_{(n)}\epsilon)}|^2
	\bigg)\dVol_g
	\\
	&\lesssim
	\left( \frac{\epsilon}{\delta C_{[n]}} + \frac{1}{C_{[n-1]}C_{[n]}} + \frac{\epsilon^2}{c^3 \delta^2} + \frac{\epsilon}{c^3 \delta^3 C_{[n-1]}} \right) \epsilon^{2(N_2 + 1 - n)}(1+\tau)^{-1+C_{(n)}\delta}
	\end{split}
	\end{equation*}
	In particular, if the $C_{[n]}$'s are sufficiently large, and if $\epsilon$ is sufficiently small compared to $\delta$ and $c$, then we can improve the bootstrap bound to find
	\begin{equation*}
	\begin{split}
	&\int_{\mathcal{M}_\tau^{\tau_1}} \bigg(
	C_{[n]}\epsilon (1+r)^{-1-C_{[n]}\epsilon} |\bm{\Gamma}^{(n)}_{(-1+C_{(n)}\epsilon)}|^2
	+ c\delta (1+r)^{-1-c\delta} |\bm{\Gamma}^{(n)}_{(-1+C_{(n)}\epsilon)}|^2
	\bigg)\dVol_g
	\\
	&\leq \frac{1}{2} \epsilon^{2(N_2 + 1 - n)}(1+\tau)^{-1+C_{(n)}\delta}
	\end{split}
	\end{equation*}

	Next, we we combine propositions \ref{proposition L2 bounds rectangular}, \ref{proposition L2 omega high}, \ref{proposition L2 zeta high}, \ref{proposition L2 tr chi high}, \ref{proposition L2 tr chibar high}, \ref{proposition L2 foliation density}, \ref{proposition L2 chihat high} and \ref{proposition L2 chibar hat high}, and use the pointwise bounds on the inhomogeneous terms which were obtained above. Then we obtain the following bound:
	\begin{equation*}
	\begin{split}
	&\int_{\mathcal{M}_\tau^{\tau_1}} \bigg(
	\delta (1+r)^{-1-\delta} |\bm{\Gamma}^{(n)}_{(-1+C_{(n)}\epsilon)}|^2
	\bigg)\dVol_g
	\\
	&\lesssim
	\left( \frac{\epsilon^2}{\delta} + \frac{\epsilon}{\delta C_{[N_2]}} + \frac{1}{C_{[N_2-1]}C_{[N_2]}} + \frac{\epsilon^2}{\delta^4} + \frac{\epsilon}{\delta^3 C_{[N-2 - 1]}} \right) \epsilon^2 (1+\tau)^{-1+C_{(N_2)}\delta}
	\\
	& + \frac{1}{\delta^4} \int_{\mathcal{M}_\tau^{\tau_1}} \bigg(
	(1+r)^{-1-\delta} |\slashed{\D} \mathscr{Y}^{N_2} h|^2_{(\text{frame})}
	+ (1+r)^{-3-\delta} |\mathscr{Y}^{N_2} h|^2_{(\text{frame})}
	+ (1+r)^{-1+ \frac{1}{2}\delta} |\overline{\slashed{\D}} \mathscr{Y}^{N_2} h|^2_{(\text{frame})}
	\\
	&\phantom{ + \frac{1}{\delta^4} \int_{\mathcal{M}_\tau^{\tau_1}} \bigg(}
	+ \delta^4 (1+r)^{1-C_{[N_2]}\epsilon} |\mathscr{Y}^{N_2} F|^2_{LL}
	+ \delta^4 \sum_{n = 0}^{N_2 - 1 } (1+r)^{1-C_{[n]}\epsilon} |\mathscr{Y}^{n} F|^2_{LL}
	\bigg)\dVol_g
	\end{split}
	\end{equation*}
	
	Now, using the bounds obtained in proposition \ref{proposition L2 bounds for F after commuting and change of basis} we can bound these last terms as
	\begin{equation*}
	\int_{\mathcal{M}_\tau^{\tau_1}} \bigg(
	+ (1+r)^{1-C_{[N_2]}\epsilon} |\mathscr{Y}^{N_2} F|^2_{LL}
	+ \sum_{n = 0}^{N_2 - 1 } (1+r)^{1-C_{[n]}\epsilon} |\mathscr{Y}^{n} F|^2_{LL}
	\bigg)\dVol_g
	\lesssim
	\delta^{-1} \epsilon^4 (1+\tau)^{-1+C_{(N_2)}\delta}
	\end{equation*}
	Using the $L^2$ bootstrap bounds it is again easy to see that, for $\epsilon$ sufficiently small compared with $\delta$ and for all sufficiently large constants $C_{[n]}$, we can improve the bounds to give us
	\begin{equation*}
	\int_{\mathcal{M}_\tau^{\tau_1}} \bigg(
	\delta (1+r)^{-1-\delta} |\bm{\Gamma}^{(n)}_{(-1+C_{(n)}\epsilon)}|^2
	+ c\delta (1+r)^{-1-c\delta} |\bm{\Gamma}^{(n)}_{(-1+C_{(n)}\epsilon)}|^2
	\bigg)\dVol_g
	\\
	\leq \frac{1}{2}\epsilon^2 (1+\tau)^{-1+C_{(N_2)}\delta}
	\end{equation*}
	
	Thus all of the bounds in equations \eqref{equation bootstrap L2 geometric} and \eqref{equation bootstrap L2 geometric} have been improved.

	\vspace{4mm}
	\textbf{First improvements of the pointwise bounds on the fields and the metric components}
	\vspace{3mm}
	
	Some quantities can be controlled in $L^\infty$ directly using the energy estimates together with propositions \ref{proposition spherical mean in terms of energy} or \ref{proposition higher weighted spherical integral} and the Sobolev embedding on the sphere, or the elliptic estimates (see section \ref{section pointwise bounds in r < r0}). In other words, we mean the fields $\phi_{(a)}$, their derivatives, and other quantities which are related \emph{algebraically} to them.
	
	Using proposition \ref{proposition spherical mean in terms of energy} together with the Sobolev inequalities on the sphere (proposition \ref{proposition Sobolev}) and the $L^2$ bootstrap bounds, we obtain, in the region $r \geq r_0$,
	\begin{equation*}
	\sup_{S_{\tau, r}} |\mathscr{Y}^n \phi_{(a)}|
	\lesssim
	\epsilon^{N_2 - n} r^{-\frac{1}{2} + \frac{1}{2} C_{[n+2, m]}\epsilon} (1+\tau)^{-\frac{1}{2} + \frac{1}{2}C_{(n+2)}\delta}
	\end{equation*}
	if $\phi_{(a)} \in \Phi_{[m]}$ and for all $n \leq N_2 - 2$. Similarly, if $\phi_{(A)} \in \Phi_{[m]}$ then, using the pointwise bounds on the matrices $M_{(A)}^{\phantom{(a)}(a)}$ we also obtain
	\begin{equation*}
	\sup_{S_{\tau, r}} |\mathscr{Y}^n \phi|_{(A)}
	\lesssim
	\epsilon^{N_2 - n} r^{-\frac{1}{2} + \frac{1}{2} C_{[n+2, m]}\epsilon} (1+\tau)^{-\frac{1}{2} + \frac{1}{2}C_{(n+2)}\delta}
	\end{equation*}

	Alternatively, we can use proposition \ref{proposition higher weighted spherical integral} together with the Sobolev inequalities to obtain
	\begin{equation*}
	\begin{split}
	\sup_{S_{\tau, r}} |\mathscr{Y}^n \phi_{(a)}|
	&\lesssim
	r^{-1} \left( \int_{S_{\tau, r_0}} \sum_{j \leq n+2} |\mathscr{Y}^j \phi|^2 \dVol_{\mathbb{S}^2} \right)^{\frac{1}{2}}
	\\
	&\phantom{\lesssim}
	+ \frac{1}{C_{[n]}\epsilon} r^{-1 + \frac{1}{2}C_{[n]}\epsilon} \left( \int_{\Sigma_\tau \cap \{r \geq r_0\}} \sum_{j \leq n+2} (r')^{1-C_{[j]}\epsilon} |\slashed{\D}_L (r\mathscr{Y}^j \phi)|^2 \upd r \wedge \dVol_{\mathbb{S}^2} \right)^{\frac{1}{2}}
	\end{split}
	\end{equation*}
	so then, using the $L^2$ bootstrap bounds for the second term together with the bounds we have already established for the first term, we have
	\begin{equation*}
	\sup_{S_{\tau, r}} |\mathscr{Y}^n \phi_{(a)}|
	\lesssim
	\epsilon^{N_2 - n} r^{-1} (1+\tau)^{-\frac{1}{2} + \frac{1}{2}C_{(n+2)}\delta}
	+ \frac{1}{C_{[n]}} \epsilon^{N_2 - n - 1} r^{-1 + \frac{1}{2}C_{[n]}\epsilon} (1+\tau)^{\frac{1}{2}C_{(n)}\delta}
	\end{equation*}
	The first term is strictly smaller than the second, for sufficiently small $\epsilon$. Hence
	\begin{equation*}
	\sup_{S_{\tau, r}} |\mathscr{Y}^n \phi_{(a)}|
	\lesssim
	\frac{1}{C_{[n]}} \epsilon^{N_2 - n - 1} r^{-1 + \frac{1}{2}C_{[n]}\epsilon} (1+\tau)^{\frac{1}{2}C_{(n)}\delta}
	\end{equation*}

	%	Interpolating between the two bounds we have obtained for the quantity $|\mathscr{Y}^n \phi_{(a)}|$ we can obtain the bound
	%	\begin{equation*}
	%	\sup_{S_{\tau, r}} |\mathscr{Y}^n \phi_{(a)}|
	%	\lesssim
	%	\frac{1}{(C_{[n]})^2} \epsilon^{N_2 - n - 2} r^{-1 + C_{(n+2)}\delta} 
	%	\end{equation*}
	%	where we have used the fact that $\epsilon \ll \delta$. We note again that these bound holds only in the region $r \geq r_0$ and for $n \leq N_2 - 2$. Moreover, they also hold for fields $\phi_{(A)} = M_{(a)}^{\phantom{(a)}(A)} \phi_{(a)}$, by using the bounds on $M_{(a)}^{\phantom{(a)}(A)}$.

	Next, by using the fact that the ``good'' derivatives can be expressed in terms of the commutation operators, we have
	\begin{equation*}
	\sup_{S_{\tau, r}} |\overline{\slashed{\D}} \mathscr{Y}^n \phi_{(a)}|
	\lesssim
	\epsilon^{N_2 - n - 1} r^{-\frac{1}{2} + \frac{1}{2} C_{[n+3, m]}\epsilon} (1+\tau)^{-\frac{1}{2} + \frac{1}{2}C_{(n+3)}\delta}
	\end{equation*}
	in the region $r \geq r_0$, and for $\phi_{(a)} \in \Phi_{[m]}$ and $n \leq N_2 - 3$. Similarly, if $\phi_{(A)} \in \Phi_{[m]}$ then under the same conditions
	\begin{equation*}
	\sup_{S_{\tau, r}} |\overline{\slashed{\D}} \mathscr{Y}^n \phi|_{(A)}
	\lesssim
	\epsilon^{N_2 - n - 1} r^{-\frac{1}{2} + \frac{1}{2} C_{[n+3, m]}\epsilon} (1+\tau)^{-\frac{1}{2} + \frac{1}{2}C_{(n+3)}\delta}
	\end{equation*}

	Next, we can use proposition \ref{proposition pointwise bound lbar} to obtain a bound on the $\Lbar$ derivatives. Before we can use this proposition, we need to obtain a bound on the quantity
	\begin{equation*}
	\int_{\Sigma_{\tau}\cap\{r \geq r'\}}
	\bigg(
	r^{1 - 2\delta} \sum_{0 \leq m \leq 2} |(r\slashed{\nabla})^m F_{(A, n)}|^2
	\bigg) \dVol_{\Sigma_\tau}
	\end{equation*}
	We can bound this by the following calculation: for any $\tau_1 \geq \tau$ we have
	\begin{equation*}
	\begin{split}
	\int_{\Sigma_{\tau}\cap\{r \geq r'\}} \bigg(
	r^{1 - 2\delta} \sum_{0 \leq m \leq 2} |(r\slashed{\nabla})^m F_{(A, n)}|^2
	\bigg) \dVol_{\Sigma_\tau}
	&\lesssim
	\int_{\Sigma_{\tau_1}\cap\{r \geq r'\}} \bigg(
	r^{1 - 2\delta} \sum_{0 \leq m \leq 2} |(r\slashed{\nabla})^m F_{(A, n)}|^2
	\bigg) \dVol_{\Sigma_\tau}
	\\
	&\phantom{\lesssim}
	+ \int_{\mathcal{M}_\tau^{\tau_1}} \bigg(
	r^{1 - 2\delta} \partial_\tau \big|_{r, \vartheta^A} \sum_{0 \leq m \leq 2} |(r\slashed{\nabla})^m F_{(A, n)}|^2
	\bigg) \dVol_{g}
	\end{split}
	\end{equation*}
	and, using equation \eqref{equation d dtau and ddr on sigma} and the pointwise bootstrap bounds we have, for any function $\phi$,
	\begin{equation*}
	|\partial_\tau \big|_{r, \vartheta^A} \phi|
	\lesssim
	\mu\left( |T \phi| + |r\slashed{\nabla} \phi| \right)
	\lesssim (1+r)^{C_{(0)}\epsilon} |\mathscr{Y}\phi|
	\end{equation*}
	Putting this together, we have
	\begin{equation*}
	\begin{split}
	\int_{\Sigma_{\tau}\cap\{r \geq r'\}} \bigg(
	r^{1 - 2\delta} \sum_{0 \leq m \leq 2} |(r\slashed{\nabla})^m F_{(A, n)}|^2
	\bigg) \dVol_{\Sigma_\tau}
	&\lesssim
	\int_{\Sigma_{\tau_1}\cap\{r \geq r'\}} \bigg(
	r^{1 - 2\delta} \sum_{0 \leq m \leq 2} |\mathscr{Y}^m F_{(A, n)}|^2
	\bigg) \dVol_{\Sigma_\tau}
	\\
	&\phantom{\lesssim}
	+ \int_{\mathcal{M}_\tau^{\tau_1}} \bigg(
	r^{1 - \frac{3}{2}\delta} \sum_{0 \leq m_1 \leq 3} |\mathscr{Y}^m F_{(A, n)}|^2
	\bigg) \dVol_{g}
	\end{split}
	\end{equation*}
	Now, the $L^2$ bootstrap bounds together with the bound on the inhomogeneous terms in proposition \ref{proposition L2 bounds for F after commuting} imply that, for all $n \leq N_2$, we easily find that
	\begin{equation*}
	\begin{split}
	\int_{\mathcal{M}_{\tau}^{\tau_1}\cap\{r \geq r'\}} \bigg(
	r^{1 - \delta} |F_{(A, n)}|^2
	\bigg) \dVol_{\Sigma_\tau}
	&\lesssim
	\delta^{-1} \epsilon^{2(N_2 + 2 - n)} (1+\tau)^{-1+C_{(n)}\delta}
	\end{split}
	\end{equation*}
	
	Now, we note that $|\mathscr{Y}^2 F_{(A, n)}|$ and $|\mathscr{Y}^3 F_{(A, n)}|$ are not precisely the same as $|F_{(A, n+3)}|$ and $|F_{(A, n+2)}|$, because they do not include terms from commuting the $\mathscr{Y}^n$ through the wave operator. However, as should be clear from the proof of proposition \ref{proposition inhomogeneous terms after commuting n times with Y}, these terms do have the same \emph{schematic} form as $|F_{(A, n+3)}|$ and $|F_{(A, n+2)}|$ respectively. Hence, following exactly as in the proof of proposition \ref{proposition L2 bounds for F after commuting}, we have, for any $\tau_2 \geq \tau$,
	\begin{equation*}
	\int_{\mathcal{M}_\tau^{\tau_2}\cap\{r \geq r'\}} \bigg(
	r^{1 - 2\delta} \sum_{0 \leq m \leq 2} |\mathscr{Y}^m F_{(A, n)}|^2
	\bigg) \dVol_{\Sigma_\tau}
	\lesssim
	\delta^{-1} \epsilon^{2(N_2 - n)}
	\end{equation*}
	and so we can find a sequence of times $\tau_n$ such that, at these times,
	\begin{equation*}
	\int_{\Sigma_{\tau_n}} \bigg(
	r^{1 - 2\delta} \sum_{0 \leq m \leq 2} |\mathscr{Y}^m F_{(A, n)}|^2
	\bigg) \dVol_{\Sigma_\tau}
	\lesssim
	\delta^{-1} \epsilon^{2(N_2 - n)} (1+\tau_n)^{-1}
	\end{equation*}
	Now, we choose $\tau_1$ to be one of these times, such that $\tau_1$ is comparable to $\tau$. Then we also have
	\begin{equation*}
	\int_{\mathcal{M}_\tau^{\tau_2}\cap\{r \geq r'\}} \bigg(
	r^{1 - 2\delta} \sum_{0 \leq m \leq 3} |\mathscr{Y}^m F_{(A, n)}|^2
	\bigg) \dVol_{\Sigma_\tau}
	\lesssim
	\delta^{-1} \epsilon^{2(N_2 -1 - n)} (1+\tau)^{-1+C_{(n+3)}\delta}
	\end{equation*}
	So, putting these calculations together, we find that
	\begin{equation*}
	\begin{split}
	\int_{\Sigma_{\tau}\cap\{r \geq r'\}} \bigg(
	r^{1 - 2\delta} \sum_{0 \leq m \leq 2} |(r\slashed{\nabla})^m F_{(A, n)}|^2
	\bigg) \dVol_{\Sigma_\tau}
	&\lesssim
	\delta^{-1} \epsilon^{2(N_2 -1 - n)} (1+\tau)^{-1+C_{(n+3)}\delta}
	\end{split}
	\end{equation*}
	
	Now we can use proposition \ref{proposition pointwise bound lbar}. Together with the $L^2$ bounds, this gives us the pointwise bound
	\begin{equation*}
	|\slashed{\D}_{\Lbar} \mathscr{Y}^n \phi_{(a)}|
	\lesssim
	\delta^{-\frac{1}{2}} \epsilon^{(N_2 -1 - n)} (1+\tau)^{-\frac{1}{2} + \frac{1}{2}C_{(n+3)}\delta}
	\end{equation*}
	for all $n \leq N_2 - 3$, from which it follows (using the pointwise bounds we are assuming on $M_{(a)}^{\phantom{(a)}(A)}$) that we also have
	\begin{equation*}
	|\slashed{\D}_{\Lbar} \mathscr{Y}^n \phi|_{(A)}
	\lesssim
	\delta^{-\frac{1}{2}} \epsilon^{(N_2 -1 - n)} r^{-1+\delta} (1+\tau)^{-\frac{1}{2} + \frac{1}{2}C_{(n+3)}\delta}
	\end{equation*}

	Putting together the calculations so far in this part of the proof, and using the fact that the rectangular components of the metric are related to the fields, we can improve the pointwise bounds on the fields, metric perturbations and their derivatives as follows: if $r \geq r_0$ then we have the improved bounds
	\begin{equation*}
	\begin{split}
	|\phi_{(a)}| &\lesssim \epsilon^{N_2} r^{-\frac{1}{2} + \frac{1}{2}C_{[2, m]}\epsilon} (1+\tau)^{-\frac{1}{2} + \frac{1}{2}C_{(2)}\delta}
	\\
	|\phi|_{(A)} &\lesssim \epsilon^{N_2} r^{-\frac{1}{2} + \frac{1}{2}C_{[2, m]}\epsilon} (1+\tau)^{-\frac{1}{2} + \frac{1}{2}C_{(2)}\delta}
	\\
	|\mathscr{Y}\phi_{(a)}| &\lesssim \epsilon^{N_2 - 1} r^{-\frac{1}{2} + \frac{1}{2}C_{[3, m]}\epsilon} (1+\tau)^{-\frac{1}{2} + \frac{1}{2}C_{(3)}\delta}
	\\
	|\mathscr{Y}\phi|_{(A)} &\lesssim \epsilon^{N_2 - 1} r^{-\frac{1}{2} + \frac{1}{2}C_{[3, m]}\epsilon} (1+\tau)^{-\frac{1}{2} + \frac{1}{2}C_{(3)}\delta}
	\\
	|h_{(\text{rect})}| &\lesssim \epsilon^{N_2} r^{-\frac{1}{2} + \frac{1}{2}C_{[2]}\epsilon} (1+\tau)^{-\frac{1}{2} + \frac{1}{2}C_{(2)}\delta}
	\\
	|h|_{(\text{frame})} &\lesssim \epsilon^{N_2} r^{-\frac{1}{2} + \frac{1}{2}C_{[2]}\epsilon} (1+\tau)^{-\frac{1}{2} + \frac{1}{2}C_{(2)}\delta}
	\\
	|\mathscr{Y}h_{(rect)}| &\lesssim \epsilon^{N_2 - 1} r^{-\frac{1}{2} + \frac{1}{2}C_{[3]}\epsilon} (1+\tau)^{-\frac{1}{2} + \frac{1}{2}C_{(3)}\delta}
	\\
	|\mathscr{Y}h|_{(\text{frame})} &\lesssim \epsilon^{N_2 - 1} r^{-\frac{1}{2} + \frac{1}{2}C_{[3]}\epsilon} (1+\tau)^{-\frac{1}{2} + \frac{1}{2}C_{(3)}\delta}
	\end{split}
	\end{equation*}
	additionally, for $r \geq r_0$ we have the following bounds, giving improved decay in $r$:
	\begin{equation*}
	\begin{split}
	|\phi_{(a)}| &\lesssim \epsilon^{N_2 - 1} r^{-1 + \frac{1}{2}C_{[2]}\epsilon} (1+\tau)^{\frac{1}{2}C_{(2)}\delta}
	\\
	|\phi|_{(A)} &\lesssim \epsilon^{N_2 - 1} r^{-1 + \frac{1}{2}C_{[2]}\epsilon} (1+\tau)^{\frac{1}{2}C_{(2)}\delta}
	\\
	|\mathscr{Y}\phi_{(a)}| &\lesssim \epsilon^{N_2 - 2} r^{-1 + \frac{1}{2}C_{[3]}\epsilon} (1+\tau)^{\frac{1}{2}C_{(3)}\delta}
	\\
	|\mathscr{Y}\phi|_{(A)} &\lesssim \epsilon^{N_2 - 2} r^{-1 + \frac{1}{2}C_{[3]}\epsilon} (1+\tau)^{\frac{1}{2}C_{(3)}\delta}
	\\
	|h_{(\text{rect})}| &\lesssim \epsilon^{N_2 - 1} r^{-1 + \frac{1}{2}C_{[2]}\epsilon} (1+\tau)^{\frac{1}{2}C_{(2)}\delta}
	\\
	|h|_{(\text{frame})} &\lesssim \epsilon^{N_2 - 1} r^{-1 + \frac{1}{2}C_{[2]}\epsilon} (1+\tau)^{\frac{1}{2}C_{(2)}\delta}
	\\
	|\mathscr{Y}h_{(\text{rect})}| &\lesssim \epsilon^{N_2 - 2} r^{-1 + \frac{1}{2}C_{[3]}\epsilon} (1+\tau)^{\frac{1}{2}C_{(3)}\delta}
	\\
	|\mathscr{Y}h|_{(\text{frame})} &\lesssim \epsilon^{N_2 - 2} r^{-1 + \frac{1}{2}C_{[3]}\epsilon} (1+\tau)^{\frac{1}{2}C_{(3)}\delta}
	\\
	\end{split}
	\end{equation*}
	Also, for all $n \leq N_2 - 2$ we have the bounds
	\begin{equation*}
	\begin{split}
	|\mathscr{Y}^n \phi_{(a)}| &\lesssim \epsilon^{N_2 - n} r^{-1 + \frac{1}{2}C_{[n+2]}\epsilon} (1+\tau)^{\frac{1}{2}C_{(n+2)}\delta}
	\\
	|\mathscr{Y}^n \phi|_{(A)} &\lesssim \epsilon^{N_2 - n} r^{-1 + \frac{1}{2}C_{[n+2]}\epsilon} (1+\tau)^{\frac{1}{2}C_{(n+2)}\delta}
	\\		
	|\mathscr{Y}^n h_{(\text{rect})}| &\lesssim \epsilon^{N_2 - n} r^{-1 + \frac{1}{2}C_{[n+2]}\epsilon} (1+\tau)^{\frac{1}{2}C_{(n+2)}\delta}
	\\
	|\mathscr{Y}^n h|_{(\text{frame})} &\lesssim \epsilon^{N_2 - n} r^{-1 + \frac{1}{2}C_{[n+2]}\epsilon} (1+\tau)^{\frac{1}{2}C_{(n+2)}\delta}
	\end{split}
	\end{equation*}
	and for all $n \leq N_2 - 3$ we have the improved bounds
	\begin{equation*}
	\begin{split}
	|\overline{\slashed{\D}}\mathscr{Y}^n \phi_{(a)}| &\lesssim \epsilon^{N_2 -1 - n} r^{-\frac{3}{2} + \frac{1}{2}C_{[n+3]}\epsilon} (1+\tau)^{-\frac{1}{2} + \frac{1}{2}C_{(n+3)}\delta}
	\\
	|\overline{\slashed{\D}}\mathscr{Y}^n \phi|_{(A)} &\lesssim \epsilon^{N_2 -1 - n} r^{-\frac{3}{2} + \frac{1}{2}C_{[n+3]}\epsilon} (1+\tau)^{-\frac{1}{2} + \frac{1}{2}C_{(n+3)}\delta}
	\\		
	|\overline{\slashed{\D}}\mathscr{Y}^n h_{(\text{rect})}| &\lesssim \epsilon^{N_2 -1 - n} r^{-\frac{3}{2} + \frac{1}{2}C_{[n+3]}\epsilon} (1+\tau)^{-\frac{1}{2} + \frac{1}{2}C_{(n+3)}\delta}
	\\
	|\overline{\slashed{\D}}\mathscr{Y}^n h|_{(\text{frame})} &\lesssim \epsilon^{N_2 -1 - n} r^{-\frac{3}{2} + \frac{1}{2}C_{[n+3]}\epsilon} (1+\tau)^{-\frac{1}{2} + \frac{1}{2}C_{(n+3)}\delta}
	\end{split}
	\end{equation*}
	and finally, for all $n \leq N_2 - 3$ we have the improved bounds
	\begin{equation*}
	\begin{split}
	|\slashed{\D}\mathscr{Y}^n \phi_{(a)}| &\lesssim \delta^{-\frac{1}{2}} \epsilon^{N_2 -1 - n} r^{-1 + \delta} (1+\tau)^{-\frac{1}{2} + \frac{1}{2}C_{(n+3)}\delta}
	\\
	|\slashed{\D}\mathscr{Y}^n \phi|_{(A)} &\lesssim \delta^{-\frac{1}{2}} \epsilon^{N_2 -1 - n} r^{-1 + \delta} (1+\tau)^{-\frac{1}{2} + \frac{1}{2}C_{(n+3)}\delta}
	\\
	|\slashed{\D}\mathscr{Y}^n h_{(\text{rect})}| &\lesssim \delta^{-\frac{1}{2}} \epsilon^{N_2 -1 - n} r^{-1 + \delta} (1+\tau)^{-\frac{1}{2} + \frac{1}{2}C_{(n+3)}\delta}
	\\
	|\slashed{\D}\mathscr{Y}^n h|_{(\text{frame})} &\lesssim \delta^{-\frac{1}{2}} \epsilon^{N_2 -1 - n} r^{-1 + \delta} (1+\tau)^{-\frac{1}{2} + \frac{1}{2}C_{(n+3)}\delta}
	\end{split}
	\end{equation*}

	Note that all of these improved bounds improve the constant relative to the bootstrap constants in section \ref{section pointwise bootstrap} by at least a factor\footnote{There is also a (potentially large) numerical constant in the improved bounds we have derived. However, for sufficiently small $\epsilon$ we still obtain an improvement.} of $\delta^{-\frac{1}{2}} \epsilon$, as long as we take $N_2 \geq 7$. Note also that we must assume that $N_1 \leq N_2 - 3$ for all of these bounds to hold.

	The bounds derived above improve almost all of the pointwise bootstrap bounds for the fields $\phi$, their derivatives and the metric fields $h$ that were made in section \ref{section pointwise bootstrap}, with the important exception of the \emph{improved} pointwise bootstraps. By this, we mean those bootstrap assumptions giving sharp decay in $r$ (often at the expense of decay in $\tau$). Recovering these bounds is a more delicate matter, and in fact it will involve simultaneously improving the pointwise bootstrap bounds on \emph{geometric} quantities, so we postpone it for now. Instead, we note that we are already in a position to improve the pointwise bounds of the fields at the lowest order, so we turn to this next.

	\vspace{4mm}
	\textbf{Improving the pointwise bounds on lowest-order quantities}
	\vspace{3mm}

	We can use these bounds to improve the \emph{sharp} pointwise bounds on the bad derivatives $\Lbar \mathscr{Y}^n \phi$ at lowest order. First, we note that, if $\phi_{(a)}$ (or $\phi_{(A)}$) is in $\Phi_{[0]}$, then the improved pointwise bounds above give us the bounds
	\begin{equation*}
	\begin{split}
	|F_{(a,0)}| &\lesssim \delta^{-\frac{1}{2}} \epsilon^{2(N_2 - 1)} r^{-\frac{5}{2} + \frac{1}{2}C_{[3]}\epsilon + \delta} (1+\tau)^{-1 + C_{(3)}\delta}
	\\
	|F|_{(A,0)} &\lesssim \delta^{-\frac{1}{2}} \epsilon^{2(N_2 - 1)} r^{-\frac{5}{2} + \frac{1}{2}C_{[3]}\epsilon + \delta} (1+\tau)^{-1 + C_{(3)}\delta}
	\end{split}
	\end{equation*}
	Next, we cam appeal to proposition \ref{proposition improved pointwise bound Lbar}, which gives us that, in the region $r \geq r_0$, we actually have the improved bound
	\begin{equation*}
	\begin{split}
	|\Lbar \phi_{(a)}| &\lesssim \delta^{-\frac{1}{2}} \epsilon^{(N_2 - 1)} r^{-1} (1+\tau)^{-\frac{1}{2} + \frac{1}{2}C_{(3)}\delta}
	\\
	|\Lbar \phi|_{(A)} &\lesssim \delta^{-\frac{1}{2}} \epsilon^{(N_2 - 1)} r^{-1} (1+\tau)^{-\frac{1}{2} + \frac{1}{2}C_{(3)}\delta}	
	\end{split}
	\end{equation*}
	Note that, to prove this second bound, we can follow the construction of proposition \ref{proposition improved pointwise bound Lbar}, and combine this with the bound
	\begin{equation*}
	|L M_{(A)}^{\phantom{(A)}(a)}| \lesssim r^{-1-\delta}
	\end{equation*}

	On the other hand, if $\phi_{(a)}$ (or $\phi_{(A)}$) is in $\phi_{[m]}$, for $m \geq 1$, then we can obtain a sharp bound by an induction argument. Suppose that, for all $m \leq m_1$, we have already obtained the improved pointwise bounds
	\begin{equation*}
	\begin{rcases*}
	|\Lbar \phi_{(a)}| \\
	|\Lbar \phi|_{(A)} 
	\end{rcases*} \lesssim \delta^{-\frac{1}{2}} \epsilon^{(N_2 - 1)} r^{-1 + \frac{1}{2}C_{(m)}\epsilon}
	\quad \text{for } \phi_{(a)} \, , \, \phi_{(A)} \in \Phi_{[m]} 
	\end{equation*}
	Then, using the structure of the inhomogeneous term $F_{(0,m_1 + 1)}$ we have
	\begin{equation*}
	|F_{(0, m_1 + 1)}|
	\lesssim
	\epsilon r^{-1} |\Lbar \phi_{[m_1 + 1]}|
	+ \delta^{-1} \epsilon^{2(N_2 - 1)} r^{-2 + C_{(m_1)}\epsilon}
	\end{equation*}
	Now, using proposition\footnote{Technically, the field $\phi_{[m_1 + 1]}$ might not be the same as the field $\phi_{(a)}$ that we are estimating, but it is at the same level of the semilinear hierarchy, and the proof of proposition \ref{proposition improved pointwise bound Lbar} can easily be adapted to this case.} \ref{proposition improved pointwise bound Lbar}, we obtain the bounds
	\begin{equation*}
	\begin{rcases*}
	|\Lbar \phi_{(a)}| \\
	|\Lbar \phi|_{(A)} 
	\end{rcases*} \lesssim \delta^{-\frac{1}{2}} \epsilon^{(N_2 - 1)} r^{-1 + \frac{1}{2}C_{(m_1 + 1)}\epsilon}
	\quad \text{for } \phi_{(a)} \, , \, \phi_{(A)} \in \Phi_{[m_1 + 1]} 
	\end{equation*}
	as long as $C_{(m_1 + 1)}$ is sufficiently large compared with $C_{(m_1)}$. This proves the inductive step.

	Now, before we can improve the rest of the sharp pointwise bounds on the bad derivatives (that is, the pointwise bounds giving improved decay in $r$ \emph{after} commuting) we need to first obtain some improved bounds on the geometric error terms. This is because, in order to utilise proposition \ref{proposition improved pointwise bound Lbar}, we need some improved bounds on the inhomogeneous terms $F_{(a,n)}$. Apart from the lowest-order terms, these error terms will involve geometric quantities (such as $\tr_{\slashed{g}}\chi$ etc.) and so in order to obtain an improvement on the $\Lbar$ derivatives we first need some improvement on these geometric quantities.

	\vspace{4mm}
	\textbf{First improvements of the pointwise bounds on geometric quantities}
	\vspace{3mm}
	
	Using the improved bounds that we have already obtained (summarized directly above) we can begin to improve the pointwise bounds on the various geometric quantities. We will prove the improved bounds on geometric quantities using an induction on the number of commutation operators. We make the inductive hypothesis:
	\begin{equation}
	\label{equation induction pointwise geometric bounds}
	\begin{split}
	\bm{\Gamma}^{(n)}_{(-1+C_{(n)}\delta)} &\lesssim \epsilon^{N_2 - 4 - n} (1+r)^{-1 + C_{(n)}\epsilon} (1+\tau)^{-\beta}
	\\
	\bm{\Gamma}^{(n)}_{(-1-\delta)} &\lesssim \epsilon^{N_2 - 4 - n} (1+r)^{-1 -(2-c_{(n)})\delta} (1+\tau)^{-\beta}
	\end{split}
	\end{equation}
	where, as was previously the case, $\beta$ is some constant in the range $(0, \frac{1}{2})$, satisfying $\frac{1}{2} - \beta \gg C^* \delta$.
	
	This means that, for each of the various geometric quantities, we will need to show both the inductive step and the base case. The latter means showing the estimate before applying any commutation operators, the former will be split up into two cases: we will first show that, assuming the inductive hypothesis holds for all $n \leq n_0 - 1$, we have bounds of the form
	\begin{equation*}
	\begin{split}
	\bm{\Gamma}^{(n_0)}_{(-1+C_{(n_0)}\epsilon)} &\lesssim \epsilon^{N_2 - 4 - n_0} (1+r)^{-1 + \delta} (1+\tau)^{-\beta}
	\\
	\bm{\Gamma}^{(n_0)}_{(-1-\delta)} &\lesssim \delta^{-2}\epsilon^{N_2 - 3 - n_0} (1+r)^{-1 -(2-c_{(n_0)})\delta} (1+\tau)^{-\beta}
	\end{split}
	\end{equation*}
	Now, these bounds will be sufficient to recover the \emph{sharp} pointwise decay of the fields $\phi$ and the metric components $h$. Finally, after this is done, we can then return to the pointwise decay for the geometric quantities, and finish proving the inductive step.

	\vspace{2mm}
	\textit{Improving the pointwise bounds on the rectangular components of the frame fields}
	\vspace{2mm}
	
	We can now apply proposition \ref{proposition pointwise bounds rectangular}. With the pointwise bounds we have now obtained on the metric and its first derivatives, we can choose $\sqrt{\mathcal{E}} = \epsilon^{N_2 - 1}(1+\tau)^{-\frac{1}{2} + \frac{1}{2}C_{(2)}\delta}$. This leads to the bounds
	\begin{equation*}
	|X_{(\text{frame})}|
	\leq
	1 + C \delta^{-1} \epsilon^{N_2 - 1} (1+\tau)^{-\frac{1}{2} + \frac{1}{2}C_{(3)}\delta}
	\end{equation*}
	For some numerical constant $C$. Now, for $\epsilon$ sufficiently small, this allows us to obtain the bound
	\begin{equation*}
	|X_{(\text{frame})}|
	\leq
	1 + C_{(0)}\delta^{-1} \epsilon^{N_2 - 1} (1+\tau)^{-\frac{1}{2} + \frac{1}{2}C_{(3)}\delta}
	\end{equation*}
	From this proposition, we also obtain the bound
	\begin{equation*}
	|\bar{X}_{(\text{frame})}|
	\lesssim
	\delta^{-1} \epsilon^{N_2 - 1} (1+\tau)^{-\frac{1}{2} + \frac{1}{2}C_{(3)}\delta}
	\end{equation*}
	which, for small enough $\epsilon$, allows us to improve the bound on $\bar{X}_{(\text{frame})}$ to
	\begin{equation*}
	|\bar{X}_{(\text{frame})}|
	\leq
	\frac{1}{2}\epsilon (1+r)^{-2\beta} (1+\tau)^{-\beta}
	\end{equation*}
	
	Similarly, proposition \ref{proposition pointwise bound rectangular small} allows us to improve the bounds $X_{(\text{frame, small})}$ to
	\begin{equation*}
	|X_{(\text{frame})}|
	\leq
	C_{(0)}\delta^{-1} \epsilon^{N_2 - 1} (1+\tau)^{-\frac{1}{2} + \frac{1}{2}C_{(3)}\delta}
	\end{equation*}
	
	This finishes the proof of the ``base case'' in the inductive argument that establishes equations \eqref{equation induction pointwise geometric bounds} in the case of the rectangular components of the null frame fields.

	Next, we use the first part of proposition \ref{proposition pointwise bounds Yn Xframe}. Now, using the bounds on the metric perturbations $h$ and their derivatives that we have already obtained, together with the inductive hypothesis for the lower order geometric error terms, we can apply the first part of proposition \ref{proposition pointwise bounds Yn Xframe} with the choice $\sqrt{\mathcal{E}} = \epsilon^{N_2 - 3 - n_0}(1+\tau)^{-\beta}$, we find that, for all $n_0 \leq N_2 - 3$ we have
	\begin{equation*}
	\begin{split}
	|\mathscr{Y}^{n_0} X_{(\text{frame})}|
	&\lesssim
	(1+r)^{C_{(n_0-1)}\epsilon}
	\\
	|\mathscr{Y}^{n_0} \bar{X}_{(\text{frame})}|
	&\lesssim
	\delta^{-1} \epsilon^{N_2 - 3 - n_0} (1+r)^{-2\delta} (1+\tau)^{-\beta}
	\end{split}
	\end{equation*}
	
	Similarly, using proposition \ref{proposition pointwise bound rectangular small} with we obtain the bound
	\begin{equation*}
	|\mathscr{Y}^{n_0} X_{(\text{frame, small})}|
	\lesssim
	\epsilon^{N_2 - 4 - n_0}(1+r)^{C_{(n_0-1)}\epsilon} (1+\tau)^{-\beta}
	\end{equation*}
	
	Note that these bounds are actually already strong enough to prove the inductive step for the quantities $X_{(\text{frame, small})}$ and $\bar{X}_{(\text{frame, small})}$ (see equation \eqref{equation induction pointwise geometric bounds}).

	\vspace{2mm}
	\textit{First improvements of the pointwise bounds on the foliation density}
	\vspace{2mm}
	
	Now we will obtain our first improvements on the foliation density. As above, we will be able to prove the base case for the induction, but this time we will not be able to obtain the inductive step straight away. Instead, we will be able to obtain an intermediate bound.

	First, we can use proposition \ref{proposition pointwise bound mu} with the choice $\sqrt{\mathcal{E}} = \epsilon^{N_2 - 1}(1+\tau)^{-\frac{1}{2} + \frac{1}{2}C_{(3)}\delta}$. This gives use the bounds
	\begin{equation*}
	\begin{split}
	|\mu| &\leq 2(1+r)^{C \epsilon^{N_2 - 1}(1+\tau)^{-\frac{1}{2} + \frac{1}{2}C_{(3)}\delta}}
	\\
	|\mu^{-1}| &\leq 2(1+r)^{C \epsilon^{N_2 - 1}(1+\tau)^{-\frac{1}{2} + \frac{1}{2}C_{(3)}\delta}}
	\end{split}
	\end{equation*}

	Next, the bounds we have recovered so far, along with the inductive hypothesis, allow us to use the first part of proposition \ref{proposition pointwise bound Yn mu} with the choice $\sqrt{\mathcal{E}} = \epsilon^{N_2 - 3 - n_0}(1+\tau)^{-\frac{1}{2} + \frac{1}{2}C_{(n_0+3)}\delta}$. This yields the bound
	\begin{equation*}
	|\mathscr{Y}^{n_0} \log \mu| \lesssim \delta^{-2} \epsilon^{N_2 - 3 - n_0} (1+r)^{-1+\delta} (1+\tau)^{-\beta}
	\end{equation*}

	\vspace{2mm}
	\textit{First improvements of the pointwise bounds on the connection coefficient $\omega$}
	\vspace{2mm}
	
	Next, we will now attempt to improve the pointwise behaviour of the connection coefficient $\omega$. Using proposition \ref{proposition pointwise bound omega} with the choice $\sqrt{\mathcal{E}} = \epsilon^{N_2 - 1}(1+\tau)^{-\frac{1}{2} + \frac{1}{2}C_{(3)}\delta}$ together with the bounds we have obtained on $h$ and $\partial h$, we obtain the bound
	\begin{equation*}
	|\omega| \lesssim \epsilon^{N_2 - 1}(1+r)^{-1} (1+\tau)^{-\frac{1}{2} + \frac{1}{2}C_{(3)}\delta}
	\end{equation*}
	proving the base case for $\omega$.
	
	Now, we can use proposition \ref{proposition pointwise bound Yn omega} to bound higher derivatives of $\omega$. Specifically, using the pointwise bounds on the metric components and their derivatives, together with the inductive hypothesis, we can apply the first part of proposition \ref{proposition pointwise bound Yn omega} with the choice $\sqrt{\mathcal{E}} = \epsilon^{N_2 - 3 - n_0}(1+\tau)^{-\beta}$, which yields
	\begin{equation*}
	|\mathscr{Y}^n \omega| \lesssim \epsilon^{N_2 - 3 - n_0}(1+r)^{-1+\delta} (1+\tau)^{-\beta}
	\end{equation*}

	\vspace{2mm}
	\textit{First improvements of the pointwise bounds on the connection coefficient $\zeta$}
	\vspace{2mm}
	
	We now examine the connection coefficient $\zeta$. Using proposition \ref{proposition pointwise bound zeta} together with the bounds we have obtained on $h$ and $\partial h$, we find that we can choose $\sqrt{\mathcal{E}} = \epsilon^{N_2 - 1}(1+\tau)^{-\frac{1}{2} + \frac{1}{2}C_{(3)}\delta}$. This leads to the bound
	\begin{equation*}
	|\zeta| \lesssim \epsilon^{N_2 - 1}(1+r)^{-1+C_{(0)}\epsilon} (1+\tau)^{-\frac{1}{2} + \frac{1}{2}C_{(3)}\delta}
	\end{equation*}
	proving the base case for $\zeta$.

	Next, we use proposition \ref{proposition pointwise bound Yn zeta} to improve the bounds on $\mathscr{Y}^n \zeta$. Substituting the pointwise bounds on the metric components and their derivatives, together with the inductive hypothesis, we can apply the first part of proposition \ref{proposition pointwise bound Yn zeta} with the choice $\sqrt{\mathcal{E}} = \epsilon^{N_2 - 3 - n_0}(1+\tau)^{-\beta}$ to obtain
	\begin{equation*}
	|\mathscr{Y}^n \zeta| \lesssim \epsilon^{N_2 - 3 - n_0}(1+r)^{-1+\delta} (1+\tau)^{-\beta}
	\end{equation*}

	\vspace{2mm}
	\textit{Improving the pointwise bounds on the connection coefficient $\tr_{\slashed{g}}\chi_{(\text{small})}$}
	\vspace{2mm}
	
	Using proposition \ref{proposition pointwise bound trchi small} together with the bounds that we already possess for the metric perturbations and their derivatives, we find that
	\begin{equation*}
	|\tr_{\slashed{g}}\chi_{(\text{small})}|
	\lesssim
	\delta^{-1} \epsilon^{N_2 - 1}(1+r)^{-1-2\delta} (1+\tau)^{-\frac{1}{2} + \frac{1}{2}C_{(3)}\delta}
	\end{equation*}
	
	Next, we apply proposition \ref{proposition pointwise bound Yn tr chi}. Using the bounds on the metric components together with the inductive hypothesis, this leads to the bound
	\begin{equation*}
	|\mathscr{Y}^n \tr_{\slashed{g}} \chi_{(\text{small})}| \lesssim \delta^{-1} \epsilon^{N_2 - 3 - n_0}(1+r)^{-1-2\delta + C_{(n-1)}\epsilon} (1+\tau)^{-\beta}
	\end{equation*}

	\vspace{2mm}
	\textit{Improving the pointwise bounds on the connection coefficient $\hat{\chi}$}
	\vspace{2mm}
	
	Using proposition \ref{proposition pointwise bound chihat} together with the bounds that we already possess for the metric perturbations and their derivatives, we find that
	\begin{equation*}
	|\tr_{\slashed{g}}\chi_{(\text{small})}|
	\lesssim
	\delta^{-1} \epsilon^{N_2 - 1}(1+r)^{-1-2\delta} (1+\tau)^{-\frac{1}{2} + \frac{1}{2}C_{(3)}\delta}
	\end{equation*}
	
	Next, we apply proposition \ref{proposition pointwise bound Yn chihat}. Using the bounds on the metric components together with the inductive hypothesis, this leads to the bound
	\begin{equation*}
	|\mathscr{Y}^n \tr_{\slashed{g}} \chi_{(\text{small})}| \lesssim \delta^{-1} \epsilon^{N_2 - 3 - n_0}(1+r)^{-1-2\delta + C_{(n-1)}\epsilon} (1+\tau)^{-\beta}
	\end{equation*}

	\vspace{2mm}
	\textit{First improvements of the pointwise bounds on the connection coefficient $\tr_{\slashed{g}} \chibar_{(\text{small})}$}
	\vspace{2mm}
	
	Using proposition \ref{proposition pointwise bound tr chibar} together with the bounds that we already possess for the metric perturbations and their derivatives, we find that
	\begin{equation*}
	|\tr_{\slashed{g}}\chibar_{(\text{small})}|
	\lesssim
	\delta^{-1} \epsilon^{N_2 - 1}(1+r)^{-1+\delta} (1+\tau)^{-\frac{1}{2} + \frac{1}{2}C_{(3)}\delta}
	\end{equation*}
	
	Next, we apply proposition \ref{proposition pointwise bound Yn tr chibar}. Using the bounds on the metric components together with the inductive hypothesis, this leads to the bound
	\begin{equation*}
	|\mathscr{Y}^n \tr_{\slashed{g}} \chibar_{(\text{small})}| \lesssim \delta^{-1} \epsilon^{N_2 - 3 - n_0}(1+r)^{-1+\delta} (1+\tau)^{-\beta}
	\end{equation*}

	\vspace{2mm}
	\textit{First improvements of the pointwise bounds on the connection coefficient $\hat{\chibar}$}
	\vspace{2mm}
	
	Using proposition \ref{proposition pointwise bound chibar hat} together with the bounds that we already possess for the metric perturbations and their derivatives, we find that
	\begin{equation*}
	|\tr_{\slashed{g}}\chibar_{(\text{small})}|
	\lesssim
	\delta^{-1} \epsilon^{N_2 - 1}(1+r)^{-1+\delta} (1+\tau)^{-\frac{1}{2} + \frac{1}{2}C_{(3)}\delta}
	\end{equation*}
	
	Next, we apply proposition \ref{proposition pointwise bound Yn chibar hat}. Using the bounds on the metric components together with the inductive hypothesis, this leads to the bound
	\begin{equation*}
	|\mathscr{Y}^n \tr_{\slashed{g}} \chibar_{(\text{small})}| \lesssim \delta^{-1} \epsilon^{N_2 - 3 - n_0}(1+r)^{-1+\delta} (1+\tau)^{-\beta}
	\end{equation*}

	\vspace{4mm}
	\textbf{The sharp decay estimates on the field and metric components}
	\vspace{3mm}

	Now that we have obtained preliminary bounds for all of the geometric quantities, we are in a position to improve the pointwise bounds on the fields $\phi$ (and thereby on the metric perturbations $h$), in order to recover the sharp pointwise bounds on these quantities. We will then be able to return to the geometric quantities that need further improvement, and recover this additional improvement, which will finally complete the proof of the inductive step.

	We first make an additional inductive hypothesis: we suppose that, if \emph{either} $n \leq n_0 - 1$ \emph{or} if $m \leq m_0 - 1$, then we have already obtained the bounds
	\begin{equation}
	\label{equation induction pointwise bad derivs}
	|\slashed{\D} \mathscr{Y}^{n} \phi_{(A)}|
	\lesssim
	\epsilon^{(N_2 - 3 - n)} (1+r)^{-1 + C_{(n, m)}\epsilon} (1+\tau)^{-\beta}
	\end{equation}
	for all $\phi_{(A)} \in \Phi_{[m]}$. Note that we have already proved this in the case $n = 0$ (i.e.\ before commuting).

	We first show how to prove the ``base case'' for this induction. That is, we will prove this bound in the case that $m = 0$, so that $\phi \in \Phi_{[0]}$, assuming that we have already obtained the bound for smaller values of $n$. In this case, if we substitute the bounds that we have already obtained into the expression for the inhomogeneous term given in proposition \ref{proposition inhomogeneous terms after commuting n times with Y} then we easily obtain the bound
	\begin{equation*}
	\begin{split}
	|F_{(A, n_0)}|
	&\lesssim
	\epsilon(1+r)^{-1} |\slashed{\D} \mathscr{Y}^n \phi_{(A)}|
	+ \epsilon (1+r)^{-1} |\slashed{\D} \mathscr{Y}^n h|_{LL}
	\\
	&\phantom{\lesssim}
	+ \epsilon^{(N_2 - 2 - n_0)} (1+r)^{-2 + 2C_{(n_0-1)}\epsilon} (1+\tau)^{-\beta}
	\end{split}
	\end{equation*}
	as long as $n_0 \leq N_2 - 4$. Indeed, we could actually obtain a better bound, with slightly faster decay in $\tau$ and $r$. For example, we now have the bound
	\begin{equation*}
	\begin{split}
	|\bar{\partial} \phi_{(A)}| |r\slashed{\nabla}^2 \mathscr{Y}^{n_0 - 1} \log \mu|
	&\lesssim
	\epsilon^{N_2 - 1} (1+r)^{-\frac{3}{2} + C_{(3)}} (1+\tau)^{-1 + \frac{1}{2}C_{(3)}\delta}
	r^{-1}|\mathscr{Y}^{n_0 + 1} \log \mu|
	\\
	&\lesssim
	\delta^{-2} \epsilon^{2N_2 - 4 - n_0} (1+r)^{-\frac{5}{2} + \delta + C_{(3)}} (1+\tau)^{-1 + \frac{1}{2}C_{(3)}\delta -\beta}
	\end{split}
	\end{equation*}
	which has additional decay in $r$ and $\tau$ compared to what we require, and, if $N_2 \geq 2$ it also has additional factors of $\epsilon$ (which more than compensate for the factor of $\delta^{-2}$). Similar bounds can be places on the other terms, assuming that $N_2 \geq 3$. The worst decay in $r$ comes from terms of the form
	\begin{equation*}
	\begin{split}
	|\partial \mathscr{Y} \phi| |\slashed{\D} \mathscr{Y}^{n_0 - 1} \phi|
	&\lesssim
	\epsilon^{2N_2 - 6 - n_0} (1+r)^{-2 + C_{(1)} + C_{(n_0-1)}} (1+\tau)^{-\frac{1}{2} + \frac{1}{2}C_{(4)}\delta - \beta}
	\end{split}
	\end{equation*}
	which, again, has more than enough decay in $r$ and $\tau$, and additional factors of $\epsilon$ if $N_2 > 3$.

	Now, we can follow the calculations in proposition \ref{proposition improved pointwise bound Lbar}, but applied to all the fields $\phi_{(A)} \in \Phi_{[0]}$ and simultaneously to the field $h_{LL}$. Specifically, if we set
	\begin{equation*}
	\tilde{\phi}_{(0)} = \left(  (\mathscr{Y}^{n_0} \phi_{(A)} | \phi_{(A)} \in \Phi_{[0]}) \, ,  \, (\mathscr{Y}^{n_0}h)_{LL} \right)
	\end{equation*}
	then, treating this as a vector of fields, we have the vector differential inequality\footnote{In order to obtain this inequality for the field $(\mathscr{Y}^n h)_{LL}$, we have to commute the rectangular components $L^a$ through the differential operator $L$, but this can be achieved without adding any error terms that cannot be easily controlled using the bounds that we already have on the derivatives of the rectangular components of the frame fields.}
	\begin{equation*}
	|\slashed{\D}_L (r\tilde{\phi}_{(0)})|
	\lesssim
	\epsilon (1+r)^{-1} |r\tilde{\phi}_{(0)}| + \epsilon^{(N_2 - 2 - n_0)} (1+r)^{-1 + 2C_{(n_0-1)}\epsilon} (1+\tau)^{-\beta}
	\end{equation*}
	where the norm of these vector quantities can be taken to be the sum of the absolute values of their components\footnote{Note that, if the operator $\mathscr{Y}^n$ includes any factors of $r\slashed{\D}_L$, then we do not actually have $\tilde{\slashed{\Box}}_g \mathscr{Y}^n \phi_{(A)} = F_{(A,n)}$, but there are additional terms of the form $\overline{\slashed{\D}}\mathscr{Y}^{n+1} \phi$ on the right hand side. However, the ``unimproved bounds'' are sufficient to show that this quantity behaves like $\epsilon^{N_2 - 2 - n_0} r^{-1-\delta}\tau^{-\beta}$. }. We can integrate this from $r = r_0$ along the integral curves of $L$. Note that, at $r = r_0$, we already have the bounds $|\slashed{\D} \mathscr{Y}^n \phi_{(A)}| \lesssim \epsilon^{N_2 - 1 - n} (1+\tau)^{-\beta}$, and $|\slashed{\D} \mathscr{Y}^n h|_{LL} \lesssim \epsilon^{N_2 - 1 - n} (1+\tau)^{-\beta}$ from the ``unimproved'' bounds. Hence, using the Gronwall inequality, we arrive at the bounds
	\begin{equation*}
	|\tilde{\phi}_{(0)}| \lesssim \epsilon^{(N_2 - 3 - n_0)} (1+r)^{-1 + (2C_{(n_0-1)} + C)\epsilon} (1+\tau)^{-\beta}
	\end{equation*}
	for some numerical constant $C$. Since $C_{(n_0, 0)} \gg C_{(n_0 - 1)}$ this proves that the inductive hypothesis \eqref{equation induction pointwise bad derivs} holds in the ``base case'' $m = 0$, assuming that it  holds for all $n \leq n_0 - 1$ and (with this value of $n$) for any value of $m$.

	Next, we will show that, if the inductive hypothesis \eqref{equation induction pointwise bad derivs} holds for all $n \leq n_0 - 1$ \emph{and} it holds for all $m \leq m_0$, then it also holds for $m = m_0 + 1$. Let $\phi_{(A)} \in \Phi_{[m_0 + 1]}$. Then, again by substituting the bounds that we have already obtained (and those given by the inductive hypothesis) into the expression for the inhomogeneous term given in proposition \ref{proposition inhomogeneous terms after commuting n times with Y}, we find the pointwise bound
	\begin{equation*}
	\begin{split}
	|F_{(A, n_0)}|
	&\lesssim
	\epsilon(1+r)^{-1} |\slashed{\D} \mathscr{Y}^n \phi_{(A)}|
	\\
	&\phantom{\lesssim}
	+ \epsilon^{(N_2 - 2 - n_0)} (1+r)^{-2 + 2C_{(n_0, m)}\epsilon} (1+\tau)^{-\beta}
	\end{split}
	\end{equation*}
	where, again, we must take $n_0 \leq N_2 - 4$. This time, the term with the worst decay in $r$ is of the form
	\begin{equation*}
	\begin{split}
	|\partial \phi_{[m]}| |\slashed{\D} \mathscr{Y}^{n_0} \phi_{[m]}|
	&\lesssim
	\epsilon^{2N_2 - 4 - n_0} (1+r)^{-2 + C_{(0, m)} + C_{(n, m)}} (1+\tau)^{-1 + \frac{1}{2}C_{(4)}\delta - \beta}
	\end{split}
	\end{equation*}
	Again, we have more than enough decay in $r$ and $\tau$, and we also have additional factors of $\epsilon$ if $N_2 \geq 2$ (in general, we need $N_2 \geq 3$ as before).
	
	Now, applying proposition \ref{proposition improved pointwise bound Lbar} to the field $\phi_{(A)}$ and using the bound for $|F_{(A, n_0)}|$ obtained above, we find that
	\begin{equation*}
	\begin{split}
	|\slashed{\D}_{\Lbar} \mathscr{Y}^{n_0} \phi_{(A)} | 
	&\lesssim \epsilon^{(N_2 - 3 - n_0)} (1+r)^{-1 + (2C_{(n_0, m)} + C)\epsilon} (1+\tau)^{-\beta}
	\\
	&\lesssim \epsilon^{(N_2 - 3 - n_0)} (1+r)^{-1 + C_{(n_0, m+1)}\epsilon} (1+\tau)^{-\beta}
	\end{split}
	\end{equation*}
	where $C$ is some numerical constant, and the second line follows from the fact that $C_{(n_0, m+1)} \gg C_{(n_0, m)}$. This finishes the proof of the inductive step: we can now conclude that the bounds in equation \eqref{equation induction pointwise bad derivs} holds for all values of $m$ and for all $n_0 \leq N_2 - 4$.

	Now, using the algebraic relationship between the metric and the fields $\phi_{(A)}$, we find that we can improve the pointwise bounds on the bad derivatives of the metric. Specifically, for all $n \leq N_2 - 4$ we obtain the bounds
	\begin{equation*}
	\begin{split}
	|\slashed{\D} \mathscr{Y}^{n} h|_{LL}
	&\lesssim
	\epsilon^{(N_2 - 3 - n)} (1+r)^{-1 + C_{(n, 0)}\epsilon} (1+\tau)^{-\beta}
	\\
	|\slashed{\D} \mathscr{Y}^{n} h_{(\text{rect})}|
	&\lesssim
	\epsilon^{(N_2 - 3 - n)} (1+r)^{-1 + C_{(n)}\epsilon} (1+\tau)^{-\beta}
	\\
	|\slashed{\D} \mathscr{Y}^{n} h|_{(\text{frame})}
	&\lesssim
	\epsilon^{(N_2 - 3 - n)} (1+r)^{-1 + C_{(n)}\epsilon} (1+\tau)^{-\beta}
	\end{split}
	\end{equation*}
	where the first inequality follows from the fact that $h_{LL} \in \Phi_{[0]}$, and the second and third follow from the fact that $C_{(n)} = C_{(n, N_1)}$, the fact that $h_{(\text{rect})}$ can be expressed in terms of the fields $\phi_{(A)}$, and the pointwise bounds on the rectangular components of the frame fields.

	\vspace{4mm}
	\textbf{The sharp decay estimates on the geometric quantities}
	\vspace{3mm}
	
	Now that we have obtained sharp decay estimates on the metric components, we can use these bounds to obtain sharp decay estimates on various other geometric quantities. Specifically, we can upgrade the pointwise decay that we have obtained above, to obtain
	\begin{equation*}
	\bm{\Gamma}^{(n)}_{(-1+C_{(n)}\delta)} \lesssim \epsilon^{N_2 - 4 - n} (1+r)^{-1 + C_{(n)}\epsilon} (1+\tau)^{-\beta}
	\end{equation*}
	matching the inductive hypothesis \eqref{equation induction pointwise geometric bounds}. We do this using the estimates in section \ref{section pointwise bounds for geometric quantities}. Since most of the details are identical to those given above, we will not provide them here. For example, using proposition \ref{proposition pointwise bound Yn mu}, together with the bound that we have now obtained on the metric, we immediately obtain
	\begin{equation*}
	|\mathscr{Y}^n \log \mu|
	\lesssim
	C_{(n, 0)}^{-1} \epsilon^{N_2 - 4 - n} (1+r)^{-1 + C_{(n, 0)}\epsilon} (1+\tau)^{-\beta}
	\end{equation*}
	The other bounds follow from the estimates in section \ref{section pointwise bounds for geometric quantities} in just as straightforward a manner.

	Additionally we have the improved bounds on the scalar field $\Omega$ and the Gauss curvature $K$, which follow from propositions \ref{proposition pointwise bound on Omega} and \ref{proposition pointwise bound on Gauss curvature}:
	\begin{equation*}
	\begin{split}
	|\Omega| + |\Omega^{-1}| &\lesssim \exp{\delta^{-2} \epsilon^{N_2-4}} \\
	|r^2 K| &\lesssim 1 + \delta^{-2} \epsilon^{N_2 - 6}
	\end{split}
	\end{equation*}

	This finally allows us to prove the inductive step corresponding to the inductive hypothesis \eqref{equation induction pointwise geometric bounds}. To recap our argument: we showed, using the unimproved bounds, the sharp pointwise bounds on the fields $\phi$, the metric components $h$, and the geometric quantities. We then made the inductive hypothesis \eqref{equation induction pointwise geometric bounds}. This was sufficient to show some intermediate bounds on the geometric fields, which correspond to the unimproved bounds on the fields and the metric components, but not to prove the inductive step straight away. However, these intermediate bounds were then used to improve the bounds on the fields and the metric components, and to obtain sharp decay estimates for these quantities. Finally, these improved bounds were used to improve the bounds on the geometric quantities, finishing the proof of the inductive step.

	\vspace{4mm}
	\textbf{Improving the decay estimates in the region $\frac{1}{2}r_0 \leq r \leq r_0$}
	\vspace{3mm}
	
	In the region $\frac{1}{2}r_0 \leq r \leq r_0$, the geometric components can be related algebraically to the rectangular components of the metric $h$ and their derivatives (see section \ref{section geometric quantities in r leq r0}). Furthermore, $r$ is bounded, so we do not need to worry about decay in the $r$ direction.
	
	Using proposition \ref{proposition pointwise bounds 1/2 r0 leq r leq r0}, together with the $L^2$ bootstrap bounds, and bounding the inhomogeneous terms in a similar way to above, we find that in the region $\frac{1}{2}r_0 \leq r \leq r_0$, for all $n \leq N_2 - 3$ we have the bounds
	\begin{equation*}
	\begin{split}
	|\mathscr{Y}^n \phi_{(a)}| &\lesssim  \epsilon^{N_2 - n} (1+\tau)^{\frac{1}{2} C_{(n+2)}\delta}
	\\
	|\mathscr{Y}^n \phi|_{(A)} &\lesssim  \epsilon^{N_2 - n} (1+\tau)^{\frac{1}{2} C_{(n+2)}\delta}
	\\
	|\mathscr{Y}^n h_{(\text{rect})}| &\lesssim  \epsilon^{N_2 - n} (1+\tau)^{\frac{1}{2} C_{(n+2)}\delta}
	\\
	|\mathscr{Y}^n h|_{(\text{frame})} &\lesssim  \epsilon^{N_2 - n} (1+\tau)^{\frac{1}{2} C_{(n+2)}\delta}
	\\
	|\slashed{\D}\mathscr{Y}^n \phi_{(a)}| &\lesssim  \epsilon^{N_2 - 1 - n} (1+\tau)^{\frac{1}{2} C_{(n+2)}\delta}
	\\
	|\slashed{\D}\mathscr{Y}^n \phi|_{(A)} &\lesssim  \epsilon^{N_2 - 1 - n} (1+\tau)^{\frac{1}{2} C_{(n+2)}\delta}
	\\
	|\slashed{\D}\mathscr{Y}^n h_{(\text{rect})}| &\lesssim  \epsilon^{N_2 - 1 - n} (1+\tau)^{\frac{1}{2} C_{(n+2)}\delta}
	\\
	|\slashed{\D}\mathscr{Y}^n h|_{(\text{frame})} &\lesssim  \epsilon^{N_2 - 1 - n} (1+\tau)^{\frac{1}{2} C_{(n+2)}\delta}
	\end{split}
	\end{equation*}

	\vspace{4mm}
	\textbf{Improving the decay estimates in the region $r \leq \frac{1}{2}r_0$}
	\vspace{3mm}
	
	In the region $r \leq \frac{1}{2}r_0$ we can no longer rely on the angular derivatives $r\slashed{\nabla}$ together with Sobolev embedding on the sphere of radius $r$ to obtain pointwise bounds, since $r$ is not bounded from below in this region. Instead, we shall use the elliptic estimates given in section \ref{section pointwise bounds in r < r0}. Note that $r$ is also bounded in this region, so we do not need to obtain decay in $r$.
	
	First, using the expressions for the inhomogeneous term $F$ (before commuting) given in the theorem, for any given field $\phi_{(a)}$ we have
	\begin{equation*}
	\tilde{\Box}_g \phi_{(a)} = F_{(a)}
	\end{equation*}
	where $F_{(A)}$ satisfies the bound
	\begin{equation*}
	|F_{(A)}| \lesssim |F^{(0)}_{(a)}| + |\partial \phi|^2
	\end{equation*}
	
	Using both the pointwise bootstrap assumptions (section \ref{section pointwise bootstrap}) and the $L^2$ bounds (section \ref{section L2 bootstrap bounds}), we find that
	\begin{equation*}
	\begin{split}
	\int_{\Sigma_\tau \cap \{r \leq r_0\}} F^2 \dVol_{\Sigma_\tau}
	&\lesssim \epsilon^2 \int_{\Sigma_\tau \cap \{r \leq r_0\}} |\partial \phi|^2 \dVol_{\Sigma_\tau}
	\\
	&\lesssim \epsilon^{2(N_2 + 3)} (1+\tau)^{-1 + C_{(0)}\delta}
	\end{split}
	\end{equation*}
	The $L^2$ bootstrap bounds also give us the bound
	\begin{equation*}
	\begin{split}
	\int_{\Sigma_\tau \cap \{r \leq r_0\}} \left( (\partial T \phi)^2 + (T\phi)^2 + |\phi|^2 \right)\dVol_{\Sigma_\tau}
	\lesssim \epsilon^{2(N_2 + 1)} (1+\tau)^{-1 + C_{(1)}\delta}
	\end{split}
	\end{equation*}
	where we have used a Hardy inequality to bound the last term. In total we have the bound
	\begin{equation*}
	\begin{split}
	\int_{\Sigma_\tau \cap \{r \leq r_0\}} \left( (\partial T \phi)^2 + (T\phi)^2 + |\phi|^2 + |F|^2 \right)\dVol_{\Sigma_\tau}
	\lesssim \epsilon^{2(N_2 + 1)} (1+\tau)^{-1 + C_{(1)}\delta}
	\end{split}
	\end{equation*}
	
	Now, proposition \ref{proposition H2 and C0 bounds} immediately gives us the bounds
	\begin{equation*}
	\begin{split}
	||\phi_{(a)}||_{H^2{[\Sigma_\tau \cap \{r \leq \frac{3}{4}r_0\}]}} &\lesssim \epsilon^{N_2 + 1} (1+\tau)^{-\frac{1}{2} + \frac{1}{2}C_{(1)}\delta}
	\\
	||\phi_{(a)}||_{C^{0, \frac{1}{2}}{[\Sigma_\tau \cap \{r \leq \frac{3}{4}r_0\}]}} &\lesssim \epsilon^{N_2 + 1} (1+\tau)^{-\frac{1}{2} + \frac{1}{2}C_{(1)}\delta}
	\end{split}
	\end{equation*}
	The second bound gives us a pointwise bound on the field $\phi_{(a)}$. Consequently, we also obtain the corresponding bounds for the metric components:
	\begin{equation*}
	\begin{split}
	||h_{(\text{rect})}||_{H^2{[\Sigma_\tau \cap \{r \leq \frac{3}{4}r_0\}]}} &\lesssim \epsilon^{N_2 + 1} (1+\tau)^{-\frac{1}{2} + \frac{1}{2}C_{(1)}\delta}
	\\
	||h_{(\text{rect})}||_{C^{0, \frac{1}{2}}{[\Sigma_\tau \cap \{r \leq \frac{3}{4}r_0\}]}} &\lesssim \epsilon^{N_2 + 1} (1+\tau)^{-\frac{1}{2} + \frac{1}{2}C_{(1)}\delta}
	\end{split}
	\end{equation*}

	Repeating these calculations with $\phi$ replaced by $T\phi$, we obtain\footnote{Note that in the region $r \leq r_0$, the rectangular components of $T$ can be expressed directly in terms of the rectangular components of $h$. Thus, the error terms found when commuting with $T$ are all easily controlled, and can be written directly in terms of lower order quantities.} the bound
	\begin{equation*}
	||T\phi_{(a)}||_{H^2{[\Sigma_\tau \cap \{r \leq \frac{3}{4}r_0\}]}} \lesssim \epsilon^{N_2} (1+\tau)^{-\frac{1}{2} + \frac{1}{2}C_{(2)}\delta}
	\end{equation*}
	
	We also find
	\begin{equation*}
	\begin{split}
	\int_{\Sigma_\tau \cap \{r \leq r_0\}} \left( (\partial T^2 \phi_{(a)})^2 + (T^2\phi_{(a)})^2 \right)\dVol_{\Sigma_\tau}
	\lesssim \epsilon^{2N_2} (1+\tau)^{-1 + C_{(2)}\delta}
	\end{split}
	\end{equation*}
	using the $L^2$ bounds from section \ref{section L2 bootstrap bounds}. In particular, this gives us that
	\begin{equation*}
	|| T^2 \phi_{(a)} ||_{H^1{[\Sigma_\tau \cap \{r \leq \frac{3}{4}r_0\}]}}
	\lesssim \epsilon^{N_2} (1+\tau)^{-\frac{1}{2} + \frac{1}{2}C_{(2)}\delta}
	\end{equation*}
	
	Additionally, we have
	\begin{equation*}
	\begin{split}
	|| F ||_{H^1{[\Sigma_\tau \cap \{r \leq \frac{3}{4}r_0\}]}}
	&\lesssim
	|| F^{(0)} ||_{H^1{[\Sigma_\tau \cap \{r \leq \frac{3}{4}r_0\}]}}
	+ || (\partial \phi)(\partial \phi) ||_{H^1{[\Sigma_\tau \cap \{r \leq \frac{3}{4}r_0\}]}}
	\\
	&\lesssim
	|| F^{(0)} ||_{H^1{[\Sigma_\tau \cap \{r \leq \frac{3}{4}r_0\}]}}
	+ || \partial \phi ||_{H^1{[\Sigma_\tau \cap \{r \leq \frac{3}{4}r_0\}]}} || \partial \phi ||^2_{L^\infty{[\Sigma_\tau \cap \{r \leq \frac{3}{4}r_0\}]}}
	\\
	&\lesssim
	|| F^{(0)} ||_{H^1{[\Sigma_\tau \cap \{r \leq \frac{3}{4}r_0\}]}}
	+ || \phi ||_{H^2{[\Sigma_\tau \cap \{r \leq \frac{3}{4}r_0\}]}} || \partial \phi ||^2_{L^\infty{[\Sigma_\tau \cap \{r \leq \frac{3}{4}r_0\}]}}
	\end{split}
	\end{equation*}
	where in the second line we have used a product estimate. Using the bounds we already have (together with the bounds on $F^{(0)}$) we have
	\begin{equation*}
	|| F ||_{H^1{[\Sigma_\tau \cap \{r \leq \frac{3}{4}r_0\}]}}
	\lesssim
	\epsilon^{N_2} (1+\tau)^{-\frac{1}{2} + \frac{1}{2}C_{(2)}\delta}
	\end{equation*}

	In summary, we have
	\begin{equation*}
	\begin{split}
	||\phi||^2_{H^1{[\Sigma_\tau \cap \{r \leq \frac{3}{4}r_0\}]}}
	&\lesssim
	\epsilon^{2(N_2 + 2)} (1+\tau)^{-\frac{1}{2} + \frac{1}{2}C_{(0)}\delta}
	\\
	||T\phi||^2_{H^2{[\Sigma_\tau \cap \{r \leq \frac{3}{4}r_0\}]}}
	&\lesssim
	\epsilon^{2N_2} (1+\tau)^{-\frac{1}{2} + \frac{1}{2}C_{(2)}\delta}
	\\
	||T^2\phi||^2_{H^1{[\Sigma_\tau \cap \{r \leq \frac{3}{4}r_0\}]}}
	&\lesssim
	\epsilon^{2N_2} (1+\tau)^{-\frac{1}{2} + \frac{1}{2}C_{(2)}\delta}
	\\
	||F||^2_{H^1{[\Sigma_\tau \cap \{r \leq \frac{3}{4}r_0\}]}}
	&\lesssim
	\epsilon^{N_2} (1+\tau)^{-\frac{1}{2} + \frac{1}{2}C_{(2)}\delta}
	\\
	||h||^2_{H^2{[\Sigma_\tau \cap \{r \leq \frac{3}{4}r_0\}]}}
	&\lesssim
	\epsilon^{2(N_2 + 1)} (1+\tau)^{-\frac{1}{2} + \frac{1}{2}C_{(1)}\delta}
	\\
	||Th||^2_{H^1{[\Sigma_\tau \cap \{r \leq \frac{3}{4}r_0\}]}}
	&\lesssim
	\epsilon^{2(N_2 + 1)} (1+\tau)^{-\frac{1}{2} + \frac{1}{2}C_{(1)}\delta}
	\end{split}
	\end{equation*}
	
	This, together with the bootstrap bound on $h$ and $\partial h$, allows to to use proposition \ref{proposition H3 and C1 bounds} to obtain
	\begin{equation*}
	\begin{split}
	||\phi||_{H^3{[\Sigma_\tau \cap \{r \leq \frac{2}{3}r_0\}]}}
	&\lesssim
	\epsilon^{N_2} (1+\tau)^{-\frac{1}{2} + \frac{1}{2}C_{(2)}\delta}
	\\
	||\phi||_{C^{1, \frac{1}{2}}{[\Sigma_\tau \cap \{r \leq \frac{2}{3}r_0\}]}}
	&\lesssim
	\epsilon^{N_2} (1+\tau)^{-\frac{1}{2} + \frac{1}{2}C_{(2)}\delta}
	\end{split}
	\end{equation*}
	Again, using the relationship between the metric components and the fields gives us the bounds
	\begin{equation*}
	\begin{split}
	||h_{(\text{rect})}||_{H^3{[\Sigma_\tau \cap \{r \leq \frac{2}{3}r_0\}]}}
	&\lesssim
	\epsilon^{N_2} (1+\tau)^{-\frac{1}{2} + \frac{1}{2}C_{(2)}\delta}
	\\
	||h_{(\text{rect})}||_{C^{1, \frac{1}{2}}{[\Sigma_\tau \cap \{r \leq \frac{2}{3}r_0\}]}}
	&\lesssim
	\epsilon^{N_2} (1+\tau)^{-\frac{1}{2} + \frac{1}{2}C_{(2)}\delta}
	\end{split}
	\end{equation*}

	Now, we are in a position to apply the $C^{k, \frac{1}{2}}$ estimates from proposition \ref{proposition Ck estimates}. First we note that, by repeating the above calculations but with $\phi$ replaced by $T^n \phi$, we can conclude the following: if $n \leq N_2 - 2$, then
	\begin{equation*}
	\begin{split}
	||T^{n} \phi||_{C^{1, \frac{1}{2}}[\Sigma_\tau \cap \{r \leq \frac{2}{3}r_0\}]}
	\lesssim
	\epsilon^{N_2 - n} (1+\tau)^{-\frac{1}{2} + \frac{1}{2}C_{(n + 2)}\delta}
	\end{split}
	\end{equation*}
	and if $n \leq N_2 - 1$ then
	\begin{equation*}
	\begin{split}
	||T^{n} \phi||_{C^{0, \frac{1}{2}}[\Sigma_\tau \cap \{r \leq \frac{2}{3}r_0\}]}
	\lesssim
	\epsilon^{N_2 - n} (1+\tau)^{-\frac{1}{2} + \frac{1}{2}C_{(n + 2)}\delta}
	\end{split}
	\end{equation*}

	Now, suppose that, for all $k \leq k_0$ we have already obtained the following bounds: if $j \leq N_2 - k$ and $j \leq N_2 - 2$, then
	\begin{equation*}
	||T^{j} \phi||_{C^{k, \frac{1}{2}}[\Sigma_\tau \cap \{r \leq (\frac{2}{3} - k_0 \delta) r_0\}]}
	\lesssim
	\epsilon^2 \delta^{-2k} (1+\tau)^{-\frac{1}{2} + \frac{1}{2}C_{(N_2)}\delta}
	\end{equation*}
	By the results above we can take $k_0 = 1$.

	Then, we also have
	\begin{equation*}
	||T^{j} h_{(\text{rect})}||_{C^{k, \frac{1}{2}}[\Sigma_\tau \cap \{r \leq (\frac{2}{3} - k_0 \delta) r_0\}]}
	\lesssim
	\epsilon^2 \delta^{-2k} (1+\tau)^{-\frac{1}{2} + \frac{1}{2}C_{(N_2)}\delta}
	\end{equation*}
	since the metric components are linear in the fields $\phi_{(a)}$. Furthermore, we can compute
	\begin{equation*}
	\begin{split}
	||F||_{C^{k_0-2, \frac{1}{2}}[\Sigma_\tau \cap \{r \leq (\frac{2}{3} - k_0 \delta) r_0\}]}
	&\lesssim
	||F^{(0)}||_{C^{k_0-2, \frac{1}{2}}[\Sigma_\tau \cap \{r \leq (\frac{2}{3} - k_0 \delta) r_0\}]}
	+ ||(\partial \phi)^2||_{C^{k_0-2, \frac{1}{2}}[\Sigma_\tau \cap \{r \leq (\frac{2}{3} - k_0 \delta) r_0\}]}
	\\
	\\
	&\lesssim
	||F^{(0)}||_{C^{k_0-2, \frac{1}{2}}[\Sigma_\tau \cap \{r \leq (\frac{2}{3} - k_0 \delta) r_0\}]}
	\\
	&\phantom{\lesssim}
	+ ||\partial \phi||_{L^{\infty}[\Sigma_\tau \cap \{r \leq (\frac{2}{3} - k_0 \delta) r_0\}]} ||\phi||_{C^{k_0-1, \frac{1}{2}}[\Sigma_\tau \cap \{r \leq (\frac{2}{3} - k_0 \delta) r_0\}]}
	\\
	\\
	&
	\lesssim
	\epsilon^2 \delta^{-2k} (1+\tau)^{-\frac{1}{2} + \frac{1}{2}C_{(N_2)}\delta}
	\end{split}
	\end{equation*}
	
	Hence we can apply proposition \ref{proposition Ck estimates}, which immediately gives us the bound
	\begin{equation*}
	||T^{j} \phi||_{C^{k_0 + 1, \frac{1}{2}}[\Sigma_\tau \cap \{r \leq (\frac{2}{3} - k_0 \delta) r_0\}]}
	\lesssim
	\epsilon^2 \delta^{-2(k_0 + 1)} (1+\tau)^{-\frac{1}{2} + \frac{1}{2}C_{(N_2)}\delta}
	\end{equation*}
	for all $j \leq N_2 - k - 1$. In other words, if we know these bounds hold for some $k_0$ then they also hold for $k_0 + 1$, providing that $k_0 + 1 \leq N_2 - 2$. Since we know that they hold for $k_0 = 1$, they in fact hold for all $k_0 \leq N_2 - 2$.
	
	In summary, for sufficiently small $\epsilon$, we have obtained the improved bounds
	\begin{equation*}
	\begin{split}
	||T^j \phi_{(a)}||_{C^{k, \frac{1}{2}}[\Sigma_\tau \cap \{r \leq \frac{1}{2}r_0\}]}
	&\leq
	\frac{1}{2} \epsilon (1+\tau)^{-\frac{1}{2} + \frac{1}{2}C_{(N_2)}\delta}
	\\
	||T^j h_{(\text{rect})}||_{C^{k, \frac{1}{2}}[\Sigma_\tau \cap \{r \leq \frac{1}{2}r_0\}]}
	&\leq
	\frac{1}{2} \epsilon (1+\tau)^{-\frac{1}{2} + \frac{1}{2}C_{(N_2)}\delta}
	\end{split}
	\end{equation*}
	for all $j + k \leq N_2$, and for $k \leq N_2 - 2$.

	\vspace{4mm}
	\textbf{Summary of the improved pointwise bounds}
	\vspace{3mm}
	
	In summary, we see that, as long as $N_2 \geq 8$ and $\epsilon$ is sufficiently small, we can improve all of the pointwise bootstrap assumptions in section \ref{section pointwise bootstrap} by at least a factor of $1/2$.

	\textbf{Improving the} $\bm{L^2}$ \textbf{bounds on the fields and metric components}
	
	Finally, we come to the improvements of the $L^2$ bootstrap bounds made in section \ref{section L2 bootstrap bounds}. First we note that, using proposition \ref{proposition L2 bounds for F after commuting}, for sufficiently small $\epsilon$ we already have $L^2$ bounds on the inhomogeneous terms. We also suppose that the initial data for $\phi$ satisfies
	\begin{equation*}
	\begin{split}
	\mathcal{E}^{(\tilde{w}T)}[\mathscr{Y}^n\phi_{(a)}](\tau_0)
	&\lesssim
	\frac{1}{C_{[n,m]}} \epsilon^{2(N_2 + 2 - n)}
	\\
	\mathcal{E}^{(\tilde{w}T)}[\slashed{\D}_T^j \mathscr{Y}^{N_2-j} \phi_{(a)}](\tau_0)
	&\lesssim
	\frac{1}{C_{[N_2(j),m]}} \epsilon^4
	\\
	\mathcal{E}^{(L, 1-C_{[n,m]}\epsilon)}[\mathscr{Y}^n\phi_{(a)}](\tau_0)
	&\lesssim
	\frac{1}{C_{[n,m]}} \epsilon^{2(N_2 + 2 - n)}
	\\
	\mathcal{E}^{(L, 1-C_{[N_2(j),m]}\epsilon)}[\slashed{\D}_T^j \mathscr{Y}^{N_2-j} \phi_{(a)}](\tau_0)
	&\lesssim
	\frac{1}{C_{[N_2(j),m]}} \epsilon^4
	\\
	\int_{\bar{S}_t,r} |\mathscr{Y}^n \phi_{(a)}|^2 \dVol_{\mathbb{S}^2}
	&\lesssim
	\frac{1}{C_{[n,m]}} \epsilon^{2(N_2 + 2 - n)} (t-\tau_0)^{-1 + \frac{1}{2}C_{[n,m]}\epsilon}
	\\
	\int_{\bar{S}_t,r} |\slashed{\D}_T^j \mathscr{Y}^{N_2 - j} \phi_{(a)}|^2 \dVol_{\mathbb{S}^2}
	&\lesssim
	\frac{1}{C_{[N_2(j),m]}} \epsilon^{2(N_2 + 2 - n)} (t-\tau_0)^{-1 + \frac{1}{2}C_{[N_2(j),m]}\epsilon}
	\\
	\\
	\mathcal{E}^{(\tilde{w}T)}[\mathscr{Y}^n\phi]_{(A)}(\tau_0)
	&\lesssim
	\frac{1}{C_{[n,m]}} \epsilon^{2(N_2 + 2 - n)}
	\\
	\mathcal{E}^{(\tilde{w}T)}[\slashed{\D}_T^j \mathscr{Y}^{N_2-j} \phi]_{(A)}(\tau_0)
	&\lesssim
	\frac{1}{C_{[N_2(j),m]}} \epsilon^4
	\\
	\mathcal{E}^{(L, 1-C_{[n,m]}\epsilon)}[\mathscr{Y}^n\phi]_{(A)}(\tau_0)
	&\lesssim
	\frac{1}{C_{[n,m]}} \epsilon^{2(N_2 + 2 - n)}
	\\
	\mathcal{E}^{(L, 1-C_{[N_2(j),m]}\epsilon)}[\slashed{\D}_T^j \mathscr{Y}^{N_2-j} \phi]_{(A)}(\tau_0)
	&\lesssim
	\frac{1}{C_{[N_2(j),m]}} \epsilon^4
	\\
	\int_{\bar{S}_t,r} |\mathscr{Y}^n \phi|_{(A)}^2 \dVol_{\mathbb{S}^2}
	&\lesssim
	\frac{1}{C_{[n,m]}} \epsilon^{2(N_2 + 2 - n)} (t-\tau_0)^{-1 + \frac{1}{2}C_{[n,m]}\epsilon}
	\\
	\int_{\bar{S}_t,r} |\slashed{\D}_T^j \mathscr{Y}^{N_2 - j} \phi|_{(A)}^2 \dVol_{\mathbb{S}^2}
	&\lesssim
	\frac{1}{C_{[N_2(j),m]}} \epsilon^{2(N_2 + 2 - n)} (t-\tau_0)^{-1 + \frac{1}{2}C_{[N_2(j),m]}\epsilon}
	\end{split}
	\end{equation*}
	if $\phi_{(a)}$ (or $\phi_{(A)}$) is in $\Phi_{[m]}$.
	
	Substituting these bound (for the inhomogeneous terms and the initial data) into lemma \ref{lemma p weighted}, and using corollaries \ref{corollary small p decay} and \ref{corollay energy decay}, we obtain identical bounds to those given as the $L^2$ bootstrap bounds in section \ref{section L2 bootstrap bounds}, but with an additional factor of $1/C_{[n,m]}$ or $1/C_{[N_2(j),m]}$ on the right hand side\footnote{Note that there is also an implicit numerical constant in these bounds, which might be very large. Nevertheless, choosing the constants $C_{[n,m]}$ and $C_{[N_2(j),m]}$ sufficiently large, we can overwhelm this numerical constant.}. Hence, if we choose these constants sufficiently large, then we can improve the $L^2$ bootstrap bounds.
	
	Now all the bootstrap bounds in chapter \ref{chapter bootstrap} are improved, finishing the proof of the theorem.
	
\end{proof}

\appendix

\chapter{Improved energy decay}
\label{appendix improved energy decay}

If we commute once with the operator $r\slashed{\D}_L$, then we can improve the decay in $r$ for $L$ derivatives, which was used in the main body of the paper. Additionally, we can improve the decay in $\tau$ for the energy of $\slashed{\D}_T \phi$. The improved decay in $r$ was used extensively in the main body of the paper, but in contrast we did not need the additional decay in $\tau$.

Note that, with the help of this improved energy decay, we will be able to obtain both improved decay for the field $\slashed{\D}_T \phi$, as well as improved decay in the spacetime region $r \leq r_{(\text{max})}$ for some fixed $r_{(\text{max})}$.

The idea of obtaining improved decay for $T$ derivatives in the context of the $p$-weighted estimates first appears in \cite{Schlue2010}, where the commutation vector field is $L$ (rather than $rL$), but this still allows higher values of $p$ to be taken in the $p$-weightes estimates, which translates into improved decay in $\tau$. This approach was expanded in \cite{Angelopoulos2018a} and \cite{Angelopoulos2018b}, where the authors commute with both $rL$ and $r^2L$ and obtain greatly improved decay in $\tau$. Note, however, that these papers deal with the \emph{linear} wave equation, albeit on black hole spacetimes. Thus, even without commuting, $p$ can be taken to be $2$, which would actually be sufficient for our purposes (although we would have to lose a little bit of $\tau$ decay each time we commute with $r\slashed{\nabla}$ in this case).

\section{Improved \texorpdfstring{$p$}{p}-weighted estimates for \texorpdfstring{$r\slashed{\D}_L \phi$}{rDL phi}}

\begin{proposition}[The basic $p$-weighted energy estimates after commuting with $r\slashed{\D}_L \phi$]
	\label{proposition basic p weighted rL}
	Suppose $\phi$ is an $S_{\tau,r}$-tangent tensor field such that, in the region $r \geq r_0$, $r\slashed{\D}_L \phi$ satisfies the equation
	\begin{equation}
	\tilde{\slashed{\Box}}_g (r\slashed{\D}_L \phi) 
	=
	\slashed{\Delta} \phi
	+ (2^k - 1)r^{-1} \slashed{\D}_L \left( r\slashed{\D}_L \phi \right)
	+ (2^k - 1)r^{-1} \slashed{\D}_L \phi
	+ F_{(rL, 1)}
	+ F_{(rL, 2)}
	+ F_{(rL, 3)}
	\end{equation}
	where $k \in \mathbb{N}$ satisfies $k \geq 1$.
	
	Suppose additionally that $\phi$ satisfies the equation
	\begin{equation*}
	\tilde{\slashed{\Box}}_g \phi = F = F_1 + F_2 + F_3
	\end{equation*}
	
	Assume that all the bootstrap assumptions from chapter \ref{chapter bootstrap} hold. Let $0 < p \leq \delta$, and choose $t$ sufficiently large relative to $R$ and $\tau$ so that 
	\begin{equation*}
	\text{supp}(\chi_{r_0, R}) \cap {_{\tau_0}^{\tau}\bar{\Sigma}}_t = \emptyset
	\end{equation*}
	
	Then, for all positive values of $p$, if $t$ is sufficiently large compared to $R$ and $(\tau - \tau_0)$, then for all sufficiently small $\epsilon$ (depending on $\delta$) we have
	
	\begin{equation}
	\begin{split}
	&\mathcal{E}^{(L,p)}[r\slashed{\D}_L\phi](\tau, R) 
	+ \int_{^t\mathcal{M}_{\tau_0}^{\tau}} \chi_{(2r_0, R)}\Big( pr^{p-1}|\slashed{\D}_L (r\slashed{\D}_L\phi)|^2 + (2-p)r^{p-1}|\slashed{\nabla}(r\slashed{\D}_L\phi)|^2 
	\\
	&\phantom{\mathcal{E}^{(L,p)}[r\slashed{\D}_L\phi](\tau, R) 
		+ \int_{^t\mathcal{M}_{\tau_0}^{\tau}} \chi_{(2r_0, R)}\Big(}
	+ p(1-p)r^{p-3}|r\slashed{\D}_L\phi|^2 + (2^k-1)r^{p-3} \left| \slashed{\D}_L (r^2\slashed{\D}_L \phi) \right|^2 \Big) \dVol_g \\
	&\lesssim
	\mathcal{E}^{(L,p)}[r\slashed{\D}_L\phi](\tau_0, R)
	\\
	&\phantom{\lesssim}
	+ \int_{^t\mathcal{M}_{\tau_0}^{\tau}} \bigg(
	\epsilon \chi_{(2r_0, R)} r^{p-2} (1+\tau)^{-1-\delta} |\slashed{\D}_L (r\slashed{\D}_L \phi)|^2
	+ \epsilon \chi_{(2r_0, R)} r^{p-1}|\overline{\slashed{\D}}(r\slashed{\D}_L\phi)|^2
	\\
	&\phantom{\lesssim + \int_{^t\mathcal{M}_{\tau_0}^{\tau}} \bigg(}
	+ \epsilon \chi_{(2r_0, R)} r^{p-1}|\overline{\slashed{\D}}(\mathscr{Z}\phi)|^2
	+ \epsilon \chi_{(2r_0, R)} r^{p-2+2C_{(0)}\epsilon}(1+\tau)^{1-\delta}|\overline{\slashed{\D}}(\mathscr{Z}\phi)|^2
	\\
	&\phantom{\lesssim + \int_{^t\mathcal{M}_{\tau_0}^{\tau}} \bigg(}
	+ \epsilon \chi_{(2r_0, R)} r^{p-1+2C_{(0)}\epsilon}|\overline{\slashed{\D}}\phi|^2
	+ \epsilon \chi_{(2r_0, R)} r^{p}(1+\tau)^{-1-\delta} |\slashed{\D}_L \phi|^2
	\\
	&\phantom{\lesssim + \int_{^t\mathcal{M}_{\tau_0}^{\tau}} \bigg(}
	+ \epsilon \chi_{(2r_0, R)} r^{p-2}(1+\tau)^{1-\delta}|\slashed{\D} \phi|^2
	+ \epsilon \chi_{(2r_0, R)} r^{p-3+4C_{(0)}\epsilon}|\phi|^2
	+ \epsilon \chi_{(2r_0, R)} r^{p} (1+\tau)^{1-\delta}|F_1|^2
	\\
	&\phantom{\lesssim + \int_{^t\mathcal{M}_{\tau_0}^{\tau}} \bigg(}
	+ \epsilon \chi_{(2r_0, R)} r^{p+1-2\delta} |F_2|^2
	+ \epsilon \chi_{(2r_0, R)} r^{p+1} |F_3|^2
	+ \epsilon^{-1} \chi_{(2r_0, R)} r^{p} (1+\tau)^{1+\delta}|F_{(rL,1)}|^2
	\\
	&\phantom{\lesssim + \int_{^t\mathcal{M}_{\tau_0}^{\tau}} \bigg(}
	+ \epsilon^{-1} \chi_{(2r_0, R)} r^{p+1-2\delta} (1+\tau)^{6\delta}|F_{(rL,2)}|^2
	+ \epsilon^{-1} \chi_{(2r_0, R)} r^{p+1}|F_{(rL,3)}|^2
	\bigg) \dVol_g
	\\
	&\phantom{\lesssim}
	+ \int_{^t\mathcal{M}_{\tau_0}^{\tau} \cap \{r_0 \leq r \leq 2r_0\}} \left( |\slashed{\D} (\mathscr{Z}\phi)|^2 + |\slashed{\D}\phi|^2 + |\phi|^2 \right) \dVol_g \\
	&\phantom{\lesssim}
	+ \int_{^t\Sigma_{\tau}} \epsilon (1+r)^{p-\delta} |\slashed{\D}_L (r\slashed{\D}_L\phi)|^2 r^2 \upd r \wedge \dVol_{\mathbb{S}^2}
	+ \int_{\bar{S}_{t, \tau}} \epsilon(1+r)^{p + 1 - \delta } |r\slashed{\D}_L\phi|^2 \dVol_{\mathbb{S}^2} \\
	&\phantom{\lesssim}
	+ \int_{^t\Sigma_{\tau_0}} \epsilon \frac{1}{(1-p+\delta)^2}(1+r)^{p-\delta} |\slashed{\D}_L (r\slashed{\D}_L\phi)|^2 r^2 \upd r \wedge \dVol_{\mathbb{S}^2}
	\\
	&\phantom{\lesssim}
	+ \int_{\bar{S}_{t, \tau_0}} \epsilon \frac{1}{(1-p+\delta)}(1+r)^{p + 1 - \delta } |r\slashed{\D}_L\phi|^2 \dVol_{\mathbb{S}^2} \\
	&\phantom{\lesssim} 
	+ \int_{\tau_0}^\tau \left( \int_{^t\Sigma_{\tau'}\cap\left\{ \frac{1}{2}R \leq r \leq R \right\}} p^{-2} r^{p-1}|\overline{\slashed{\D}}(r\slashed{\D}_L\phi)|^2 r^2 \upd r \wedge \dVol_{\mathbb{S}^2} 
	+ \int_{S_{\tau,R}} p^{-1} r^p |r\slashed{\D}_L \phi|^2 \dVol_{\mathbb{S}^2} \right) \upd \tau
	\end{split}
	\end{equation}
	Where, if $\phi$ is a scalar field then the final term in the spacetime integral, involving the term 
	\begin{equation*}
	\epsilon^{-1}\chi_{(r_0)} r^{p-3}|\phi|^2
	\end{equation*}
	is not present.

\end{proposition}

\begin{proof}
	We repeat the calculations of proposition \ref{proposition basic p weighted estimate} applied to the field $r\slashed{\D}_L \phi$. The additional terms are 
	\begin{equation*}
	\begin{split}
	&\int_{^t\mathcal{M}_{\tau_0}^{\tau}} \chi_{(2r_0, R)} \Big( 
	\slashed{\Delta} \phi
	+ r^{-1} (2^k-1)\slashed{\D}_L \left( r\slashed{\D}_L \phi \right)
	+ r^{-1} (2^k-1)\slashed{\D}_L \phi
	\Big) \left( r^p \slashed{\D}_L \left( r\slashed{\D}_L \phi \right) + r^p \slashed{\D}_L \phi \right) \dVol_g
	\\
	&=
	\int_{^t\mathcal{M}_{\tau_0}^{\tau}} \chi_{(2r_0, R)} \Big( 
	(2^k-1)r^{p-3} \left| \slashed{\D}_L (r^2\slashed{\D}_L \phi) \right|^2
	+ r^{p-2} \slashed{\nabla}^\mu (r\slashed{\nabla}_\mu \phi) \left( \slashed{\D}_L (r^2\slashed{\D}_L \phi) \right) \Big) \dVol_g
	\end{split}
	\end{equation*}
	
	The first term is positive and so gives an additional good term in the energy estimate, while the second term can be controlled assuming that we \emph{already} have control over $r\slashed{\nabla} \phi$. 
	
	When applied to the field $r\slashed{\D}_L \phi$, the error terms in the $p$ weighted energy can be estimated as
	\begin{equation*}
	\begin{split}
	|\textit{Err}_{(L,p,\textrm{bulk})}|
	&\lesssim
	r^p f'_L |\overline{\slashed{\D}}(r\slashed{\D}_L \phi)|^2
	+ (r^{p-1}f_L'' + r^{p-2}f_L') |r\slashed{\D}_L\phi|^2
	\\
	&\phantom{\lesssim}
	+ \epsilon f_L r^{p-1}(1+\tau)^{-\beta}|\slashed{\D}_L (r\slashed{\D}_L\phi)| |\slashed{\D}_{\Lbar} (r\slashed{\D}_L\phi)|
	+ \epsilon f_L r^{p-1-\delta} |\slashed{\D}_L (r\slashed{\D}_L\phi)| |\slashed{\D}_{\Lbar} (r\slashed{\D}_L\phi)|
	\\
	&\phantom{\lesssim}
	+ \epsilon f_L r^{p-1+\delta} (1+\tau)^{-\beta} |\slashed{\D}_L (r\slashed{\D}_L\phi)| |\slashed{\nabla} (r\slashed{\D}_L\phi)|
	+ \epsilon f_L r^{p-1} |\slashed{\nabla} (r\slashed{\D}_L\phi)|^2
	\\
	&\phantom{\lesssim}
	+ \epsilon f_L (p-1) r^{p-3+\delta}(1+\tau)^{-\beta} |r\slashed{\D}_L\phi|^2
	+ \epsilon f_L (p-1)r^{p-3}|(r\slashed{\D}_L\phi)|^2
	\end{split}
	\end{equation*}
	
	Now, we note that
	\begin{equation*}
	\begin{split}
	\slashed{\D}_{\Lbar} (r\slashed{\D}_L \phi)
	&=
	- r\tilde{\slashed{\Box}}_g \phi
	+ r\slashed{\Delta} \phi
	- \left( 1 + \frac{1}{2}r\tr_{\slashed{g}}\chi_{(\text{small})} - \frac{1}{2}r\omega \right) \slashed{\D}_{\Lbar} \phi
	+ \left(- \frac{1}{2}r\tr_{\slashed{g}}\chibar_{(\text{small})} + r\omega\right) \slashed{\D}_L \phi
	\\
	&\phantom{=}
	+ r\left(\zeta^\alpha + 2\slashed{\nabla}^{\alpha} \log \mu \right) \slashed{\nabla}_\alpha \phi
	- \frac{1}{2}rL^\mu \Lbar^\nu [\slashed{\D}_\mu \, , \slashed{\D}_\nu] \phi
	\end{split}
	\end{equation*}
	So, using the bootstrap bounds, we have
	\begin{equation*}
	|\slashed{\D}_{\Lbar} (r\slashed{\D}_L \phi)|^2
	\lesssim
	r^2|\tilde{\slashed{\Box}}_g \phi|^2
	+ | \overline{\slashed{\D}} (r\slashed{\nabla} \phi) |^2
	+ |\slashed{\D} \phi|^2 
	+ r^{2C_{(0)}\epsilon} |\overline{\slashed{\D}} \phi|^2
	+ \epsilon^2 r^{-2 + 4C_{(0)}\epsilon} |\phi|^2
	\end{equation*}
	%	Also, we have
	%	\begin{equation*}
	%	\slashed{\nabla}_{\alpha}(r\slashed{\D}_L \phi)
	%	=
	%	\slashed{\D}_L (r\slashed{\nabla}\phi)
	%	+ r(\chi_{(\text{small})})_\alpha^{\phantom{\alpha}\beta} \slashed{\nabla}_\beta \phi
	%	+ rL^\mu \slashed{\Pi}_\alpha^{\phantom{\alpha}\nu}[\slashed{\D}_\mu , \slashed{\D}_\nu]\phi
	%	\end{equation*}
	%	and so
	%	\begin{equation*}
	%	|\slashed{\nabla}_{\alpha}(r\slashed{\D}_L \phi)|
	%	\lesssim
	%	|\slashed{\D}_L (r\slashed{\nabla}\phi)|
	%	+ \epsilon r^{-\delta} |\slashed{\nabla}_\beta \phi|
	%	+ \epsilon r^{-3+2C_{(0)}\epsilon} |\phi|
	%	\end{equation*}
	
	Hence, we can estimate
	\begin{equation*}
	\begin{split}
	&\epsilon r^{p-1}(1+\tau)^{-\beta} |\slashed{\D}_L (r\slashed{\D}_L \phi)| |\slashed{\D}_{\Lbar} (r\slashed{\D}_L \phi)|
	\\
	&\lesssim
	\epsilon r^{p}(1+\tau)^{-1-\delta} |\slashed{\D}_L (r\slashed{\D}_L \phi)|^2
	+ \epsilon r^{p-1+2\delta}(1+\tau)^{-\beta} |\slashed{\D}_L (r\slashed{\D}_L \phi)|^2
	+ \epsilon r^{p-1} |\slashed{\D}_L (r\slashed{\D}_L \phi)|^2
	\\
	&\phantom{\lesssim}
	+ \epsilon r^{p}(1+\tau)^{1+\delta-2\beta}|F_1|^2
	+ \epsilon r^{p+1-2\delta}(1+\tau)^{-\beta}|F_2|^2
	+ \epsilon r^{p+1}(1+\tau)^{-2\beta}|F_3|^2	
	\\
	&\phantom{\lesssim}
	+ \epsilon r^{p-2}(1+\tau)^{1+\delta-2\beta}|\slashed{\D} \phi|^2
	+ \epsilon r^{p-1} |\overline{\slashed{\D}} \mathscr{Z} \phi|^2
	+ \epsilon r^{p-1+2C_{(0)}\epsilon} |\overline{\slashed{\D}}\phi|^2
	+ \epsilon^3 r^{p-3+4C_{(0)}\epsilon} |\phi|^2
	\\
	\\
	&\lesssim
	\epsilon r^{p}(1+\tau)^{-1-\delta} |\slashed{\D}_L (r\slashed{\D}_L \phi)|^2
	+ \epsilon r^{p-1} |\slashed{\D}_L (r\slashed{\D}_L \phi)|^2
	\\
	&\phantom{\lesssim}
	+ \epsilon r^{p}(1+\tau)^{1+\delta-2\beta}|F_1|^2
	+ \epsilon r^{p+1-2\delta}(1+\tau)^{-\beta}|F_2|^2
	+ \epsilon r^{p+1}(1+\tau)^{-2\beta}|F_3|^2	
	\\
	&\phantom{\lesssim}
	+ \epsilon r^{p-2}(1+\tau)^{1+\delta-2\beta}|\slashed{\D} \phi|^2
	+ \epsilon r^{p-1} |\overline{\slashed{\D}} \mathscr{Z} \phi|^2
	+ \epsilon r^{p-1+2C_{(0)}\epsilon} |\overline{\slashed{\D}}\phi|^2
	+ \epsilon^3 r^{p-3+4C_{(0)}\epsilon} |\phi|^2
	\end{split}
	\end{equation*}
	and
	\begin{equation*}
	\begin{split}
	&\epsilon r^{p-1-\delta}|\slashed{\D}_L (r\slashed{\D}_L \phi)| |\slashed{\D}_{\Lbar} (r\slashed{\D}_L \phi)|
	\\
	&\lesssim
	\epsilon r^{p-2\delta}(1+\tau)^{-1+\delta}|\slashed{\D}_L (r\slashed{\D}_L \phi)|^2
	+ \epsilon r^{p-1}|\slashed{\D}_L (r\slashed{\D}_L \phi)|^2
	+ \epsilon r^p (1+\tau)^{1-\delta} |F_1|^2
	\\
	&\phantom{\lesssim}
	+ \epsilon r^{p+1-2\delta} \left( |F_2|^2 + |F_3|^2 \right)	
	+ \epsilon r^{p-1-2\delta} |\overline{\slashed{\D}} \mathscr{Z} \phi|^2
	+ \epsilon r^{p-2} (1+\tau)^{1-\delta} |\slashed{\D} \phi|^2
	\\
	&\phantom{\lesssim}
	+ \epsilon r^{p-2+ 2C_{(0)}\epsilon} (1+\tau)^{1-\delta} |\overline{\slashed{\D}} \phi|^2
	+ \epsilon^3 r^{p-4+ 4C_{(0)}\epsilon} (1+\tau)^{1-\delta} |\phi|^2
	\\
	\\
	&\lesssim
	\epsilon r^{p}(1+\tau)^{-1-\delta}|\slashed{\D}_L (r\slashed{\D}_L \phi)|^2
	+ \epsilon r^{p-1}|\slashed{\D}_L (r\slashed{\D}_L \phi)|^2
	+ \epsilon r^p (1+\tau)^{1-\delta} |F_1|^2
	\\
	&\phantom{\lesssim}
	+ \epsilon r^{p+1-2\delta} \left( |F_2|^2 + |F_3|^2 \right)	
	+ \epsilon r^{p-1-2\delta} |\overline{\slashed{\D}} \mathscr{Z} \phi|^2
	+ \epsilon r^{p-2} (1+\tau)^{1-\delta} |\slashed{\D} \phi|^2
	\\
	&\phantom{\lesssim}
	+ \epsilon r^{p-2+ 2C_{(0)}\epsilon} (1+\tau)^{1-\delta} |\overline{\slashed{\D}} \phi|^2
	+ \epsilon^3 r^{p-3-2\delta+ 4C_{(0)}\epsilon} |\phi|^2
	\end{split}
	\end{equation*}
	Finally, we note that
	\begin{equation*}
	\begin{split}
	r^{p-1+\delta}(1+\tau)^{-\beta} |\slashed{\D}_L (r\slashed{\D}_L \phi)| |\slashed{\nabla} (r\slashed{\D}_L \phi)|
	&\lesssim
	r^{p-1+2\delta}(1+\tau)^{-2\beta} |\slashed{\D}_L (r\slashed{\D}_L \phi)|^2
	+ r^{p-1}|\slashed{\nabla} (r\slashed{\D}_L \phi)|^2
	\\
	&\lesssim
	r^{p}(1+\tau)^{-1-\delta} |\slashed{\D}_L (r\slashed{\D}_L \phi)|^2
	+ r^{p-1}|\overline{\slashed{\D}} (r\slashed{\D}_L \phi)|^2
	\end{split}
	\end{equation*}
	
	Putting the above calculations together, we conclude that
	\begin{equation*}
	\begin{split}
	|\textit{Err}_{(L,p,\textrm{bulk})}|
	&\lesssim
	r^p f'_L |\overline{\slashed{\D}}(r\slashed{\D}_L \phi)|^2
	+ (r^{p+1}f_L'' + r^{p}f_L') |\slashed{\D}_L\phi|^2
	\\
	&\phantom{\lesssim}
	+ \epsilon r^p (1+\tau)^{-1-\delta} |\slashed{\D}_L (r\slashed{\D}_L \phi)|^2
	+ \epsilon f_L r^{p-1}|\overline{\slashed{\D}}(r\slashed{\D}_L\phi)|^2
	+ \epsilon f_L r^p (1+\tau)^{1-\delta} |F_1|^2
	\\
	&\phantom{\lesssim}
	+ \epsilon f_L r^{p+1-2\delta} |F_2|^2
	+ \epsilon f_L r^{p+1} |F_3|^2
	+ \epsilon f_L r^{p-1}|\overline{\slashed{\D}}(\mathscr{Z}\phi)|^2
	+ \epsilon f_L r^{p-1+2C_{(0)}\epsilon}|\overline{\slashed{\D}}\phi|^2
	\\
	&\phantom{\lesssim}
	+ \epsilon f_L r^{p-2+2C_{(0)}\epsilon}(1+\tau)^{1-\delta}|\overline{\slashed{\D}}(\mathscr{Z}\phi)|^2
	+ \epsilon f_L r^{p-2}(1+\tau)^{1-\delta}|\slashed{\D} \phi|^2
	+ \epsilon f_L r^{p-3+4C_{(0)}\epsilon}|\phi|^2
	\\
	&\phantom{\lesssim}
	+ \epsilon f_L r^{p}(1+\tau)^{-1-\delta} |\slashed{\D}_L \phi|^2
	\end{split}
	\end{equation*}
	
	Next, we can estimate
	\begin{equation*}
	\begin{split}
	&\int_{^t\mathcal{M}_{\tau_0}^{\tau}} \chi_{(2r_0, R)}
	r^{p-2} \slashed{\nabla}^\mu (r\slashed{\nabla}_\mu \phi) \left( \slashed{\D}_L (r^2\slashed{\D}_L \phi) \right)\dVol_g
	\\
	&\lesssim
	\int_{^t\mathcal{M}_{\tau_0}^{\tau}} \chi_{(2r_0, R)}
	\Big( \delta^{-1} r^{p-1} |\slashed{\nabla}^\mu (r\slashed{\nabla}_\mu \phi)|^2 +  \delta r^{p-3} | \slashed{\D}_L (r^2\slashed{\D}_L \phi) |^2 \Big) \dVol_g
	\end{split}
	\end{equation*}
	
	Putting these calculations together, and following steps similar to those followed in proposition \ref{proposition basic p weighted estimate}, we can prove the proposition.

\end{proof}

\begin{lemma}[$p$-weighted energy estimates for $r\slashed{\D}_L \phi$]
	\label{lemma p weighted for rL}
	Suppose that the conditions of lemma \ref{lemma p weighted} hold. In particular, the pointwise bootstrap bounds of chapter \ref{chapter bootstrap} are assumed to hold, as are the relevant $L^2$ based bounds for the inhomogeneity $F$. Suppose additionally that the same bounds hold with $\phi$ replaced by $\mathscr{Z}\phi$.
	
	Suppose also that $r\slashed{\D}_L \phi$ satisfies the equation
	\begin{equation*}
	\tilde{\slashed{\Box}}_g (r\slashed{\D}_L\phi) = F_{(rL)} = F_{(rL,1)} + F_{(rL,2)}
	\end{equation*}
	where $F_{(rL)}$, $F_{(rL,1)}$ and $F_{(rL,2)}$ satisfy the bounds
	\begin{equation*}
	\begin{split}
	&\int_{\mathcal{M}^{\tau_1}_{\tau}} \left(
	\epsilon^{-1} \chi_{(2r_0)} r^{1-\frac{1}{2}\delta} (1+\tau)^{2\beta}|F_{(rL)}|^2
	\right) \dVol_g 
	\lesssim \mathcal{E}_0 (1+\tau)^{-1+C_{(\phi)}\delta}
	\\
	&\int_{\mathcal{M}^\tau_{\tau_0}} \left(
	\epsilon^{-1} \chi_{(2r_0)} r^{1-C_{(\phi)}\epsilon} (1+\tau)^{1+\delta}|F_{(rL,1)}|^2
	+ \epsilon^{-1} \chi_{(2r_0)} r^{2-2\delta+C_{(\phi)}\epsilon} (1+\tau)^{2\beta}|F_{(rL,2)}|^2
	\right) \dVol_g 
	\lesssim \mathcal{E}_0
	\end{split}
	\end{equation*}
	
	Then, for $-1 + \frac{1}{10}\delta \leq p \leq \frac{3}{2}\delta$ we have
	\begin{equation}
	\begin{split}
	&\mathcal{E}^{(L,p)}[r\slashed{\D}_L \phi](\tau_1)
	\\
	&
	+ \int_{\mathcal{M}_{\tau}^{\tau_1}} \chi_{(2r_0)} \left( 
	(1+p)r^{p-1}|\slashed{\D}_L (r\slashed{\D}_L \phi)|^2
	+ (2-p)r^{p-1}|\slashed{\D}\phi|^2
	+ p(1-p)r^{p-1}|\slashed{\D}_L\phi|^2
	\right) \dVol_g
	\\ \\
	&\lesssim
	\mathcal{E}^{(L,p)}[r\slashed{\D}_L \phi](\tau)
	+ \delta^{-1} \mathcal{E}_0 (1+\tau)^{-1+ C_{(\phi)}\delta}
	\end{split}
	\end{equation}
	
	On the other hand, for $-1 + \frac{1}{10}\delta \leq p \leq 1 - C_{(\phi)}\epsilon$, we have
	\begin{equation}
	\begin{split}
	&\mathcal{E}^{(L,p)}[r\slashed{\D}_L \phi](\tau)
	+ \int_{\mathcal{M}_{\tau_0}^{\tau}} \chi_{(2r_0)} \bigg( 
	(1+p-\frac{1}{10}\delta)r^{p-1}|\slashed{\D}_L (r\slashed{\D}_L \phi)|^2
	+ (2-p)r^{p-1}|\slashed{\D}\phi|^2
	\\
	& \phantom{\mathcal{E}^{(L,p)}[r\slashed{\D}_L \phi](\tau) + \int_{\mathcal{M}_{\tau_0}^{\tau}} \chi_{(2r_0)} \bigg(}
	+ p(1-p)r^{p-1}|\slashed{\D}_L\phi|^2
	\bigg) \dVol_g
	\\
	&\lesssim
	\mathcal{E}^{(L,p)}[r\slashed{\D}_L \phi](\tau_0)
	+ \delta^{-1} \mathcal{E}_0
	\end{split}
	\end{equation}	
	
\end{lemma}

\begin{proof}
	We begin with proposition \ref{proposition basic p weighted rL}. We first note that we can take $p$ to be \emph{negative} and still obtain a positive bulk term. In fact, we have
	\begin{equation*}
	pr^{p-1} |\slashed{\D}_L (r\slashed{\D}_L\phi)|^2
	+ r^{p-3} |\slashed{\D}_L (r^2\slashed{\D}_L\phi)|^2
	\gtrsim
	(p+1 - c\delta)r^{p-1} |\slashed{\D}_L (r\slashed{\D}_L\phi)|^2
	- c^{-1} \delta^{-1} r^{p-1} |\slashed{\D}_L \phi|^2
	\end{equation*}
	for any small constant $c$. Note that the spacetime integral of the second term can already be controlled, for $p \leq 1-C_{(\phi)}\epsilon$. Hence, if we choose $p \geq -1 + c\delta$ (for any small $c$) then we still obtain a positive bulk term for the spacetime integral of $r^{p-1}|\slashed{\D}_L (r^2\slashed{\D}_L\phi)|^2$

	Now we move on to the error terms. We first use the $p$-weighted estimate to control the term
	\begin{equation*}
	\int_{^t\mathcal{M}^\tau_{\tau_1}} \epsilon \chi_{(2r_0, R)} r^{p-2-\delta} (1+\tau)^{-1-\beta} |\slashed{\D}_L (r\phi)|^2
	\end{equation*}
	Using lemma \ref{lemma p weighted} we have the following: for $p \leq \frac{3}{2}\delta$,
	\begin{equation*}
	\mathcal{E}^{(L, p-\delta)}[\phi](\tau) \lesssim \delta^{-1} \mathcal{E}_0 (1+\tau)^{-1+C_{(\phi)}\delta}
	\end{equation*}
	and so, if $p \leq \frac{3}{2}\delta$, 
	\begin{equation*}
	\int_{^t\mathcal{M}^\tau_{\tau_1}} \epsilon \chi_{(2r_0, R)} r^{p-2-\delta} (1+\tau)^{-1-\beta} |\slashed{\D}_L (r\phi)|^2
	\lesssim 
	\epsilon \delta^{-1} \mathcal{E}_0 (1+\tau_1)^{-1-\beta+C_{(\phi)}\delta}
	\end{equation*}
	
	On the other hand, if $p \leq 1+\delta-C_{(\phi)}\epsilon$, then we have
	\begin{equation*}
	\mathcal{E}^{(L, p-\delta)}[\phi](\tau) \lesssim \delta^{-1} \mathcal{E}_0
	\end{equation*}
	and so, in this case,
	\begin{equation*}
	\int_{^t\mathcal{M}^\tau_{\tau_1}} \epsilon \chi_{(2r_0, R)} r^{p-2-\delta} (1+\tau)^{-1-\beta} |\slashed{\D}_L (r\phi)|^2
	\lesssim 
	\epsilon \delta^{-1} \beta^{-1} \mathcal{E}_0 (1+\tau_1)^{-\beta}
	\end{equation*}
	
	Next, we bound the terms involving the ``good'' derivatives. For $p \leq \frac{1}{2}\delta$ we have
	\begin{equation*}
	\int_{^t\mathcal{M}^\tau_{\tau_1}} \delta^{-1}\chi_{(2r_0,R)} r^{p-1} \left( |\overline{\slashed{\D}}(r\slashed{\nabla}\phi)|^2 + |\overline{\slashed{\D}}\phi|^2 \right) \dVol_g
	\lesssim
	p^{-1} \delta^{-1} \mathcal{E}_0 (1+\tau)^{-1+C_{(\phi)}\delta}
	\end{equation*}
	where, again, we have used lemma \ref{lemma p weighted}.
	
	On the other hand, if $p \leq 1-C_{(\phi)}\epsilon$, we have
	\begin{equation*}
	\int_{^t\mathcal{M}^\tau_{\tau_1}} \delta^{-1}\chi_{(2r_0,R)} r^{p-1} |\overline{\slashed{\D}}(r\slashed{\nabla}\phi)|^2 \dVol_g
	\lesssim
	p^{-1} \delta^{-1} \mathcal{E}_0
	\end{equation*}
	
	Next, we bound the terms involving the ``bad'' derivatives, and the lower order terms. For $p \leq 2 - \delta - C_{(\phi)}\epsilon$ we have
	\begin{equation*}
	\int_{^t\mathcal{M}^\tau_{\tau_1}} \epsilon\chi_{(2r_0,R)} \left( r^{p-3+2C_{(0)}\epsilon} |\slashed{\D}(r\slashed{\nabla}\phi)|^2 + r^{p-5+2C_{(0)}\epsilon} |\slashed{\D}(r\slashed{\nabla}\phi)|^2 \right)\dVol_g
	\lesssim
	\epsilon \delta^{-3} \mathcal{E}_0 (1+\tau)^{-2+2\delta}
	\end{equation*}
	
	Similarly, we can bound the terms in the region $r \leq 2r_0$ by
	\begin{equation*}
	\int_{^t\mathcal{M}_{\tau_0}^{\tau} \cap \{r_0 \leq r \leq 2r_0\}} \left( |\slashed{\D} (\mathscr{Z}\phi)|^2 + |\slashed{\D}\phi|^2 + |\phi|^2 \right) \dVol_g
	\lesssim
	\delta^{-1} \mathcal{E}_0 (1+\tau)^{-2+2\delta}
	\end{equation*}
	
	Next, we need to bound the error terms on the spheres $\bar{S}$. First, we consider $p \leq -1 + \frac{1}{4}\delta - 2C_{(0)}\epsilon$
	
	Following calculations similar to those in proposition \ref{proposition spherical mean in terms of energy}, we have
	\begin{equation}
	\label{equation internal bound on spherical error terms}
	\begin{split}
	&\int_{\bar{S}_{t,\tau}} r^{p+3-\delta} |\slashed{\D}_L \phi|^2 \dVol_{\mathbb{S}^2}
	\\
	&\lesssim
	(r_{(\text{min})}(t))^{-\frac{1}{4}\delta} \Bigg( \int_{r = r_{(\text{min})}(t)}^{r_{(\text{max})}} \left( \int_{S_{r,\tau}} r^{-1+\frac{1}{4}\delta} \left(r^{p+3-\delta}|\slashed{\D}_L \slashed{\D}_L \phi|^2 + r^{p+1-\delta}|\slashed{\D}_L \phi|^2 \right) r^2 \dVol_{\mathbb{S}^2} \right) \upd r
	\\
	&\phantom{\lesssim (1+r_{(\text{min})}(t))^{-\frac{1}{4}\delta} \Bigg(}
	+ \int_{S_{r_{(\text{max})},\tau}} r^{p+3-\frac{1}{2}\delta}|\slashed{\D}_L \phi|^2 \dVol_{\mathbb{S}^2}
	\Bigg)
	\end{split}
	\end{equation}
	where 
	\begin{equation*}
	r_{(\text{min})}(t) = \min_{x \in \bar{S}_{t,\tau}}r(x) \geq r_0
	\end{equation*}
	and $r_{(\text{max})}$ satisfies $r_{(\text{max})} > r_{(\text{min})}(t)$ but is otherwise arbitrary.
	
	For $p \leq \frac{3}{4}\delta - C(\phi)\epsilon$ we have
	\begin{equation*}
	\int_{r = r_0}^\infty \left( \int_{S_{r,\tau}} r^{p+3-\frac{3}{4}\delta}|\slashed{\D}_L \phi|^2 \dVol_{\mathbb{S}^2} \right) \upd r \leq \delta^{-1}\mathcal{E}_0
	\end{equation*}
	Hence we can pick $r_{(\text{max})}$ so that
	\begin{equation*}
	\int_{S_{r_{(\text{max})},\tau}} r^{p+3-\frac{1}{2}\delta}|\slashed{\D}_L \phi|^2 \dVol_{\mathbb{S}^2}
	\lesssim
	\delta^{-1}\mathcal{E}_0 (r_{(\text{max})})^{-1}
	\end{equation*}
	Moreover, there is a sequence of such $r_{(\text{max})}$ such that $r_{(\text{max})} \rightarrow \infty$.
	
	Next, we note that, for all $p \leq 1 + \frac{3}{4}\delta - C_{(\phi)}$ we have
	\begin{equation*}
	\int_{r = r_{(\text{min})}(t)}^\infty \left( \int_{S_{\tau,r}} \left( r^{p-\frac{3}{4}\delta} |\slashed{\D}_L \phi|^2 \right) r^2\dVol_{\mathbb{S}^2} \right) \upd r \lesssim \mathcal{E}_0
	\end{equation*}
	
	Finally, we have
	\begin{equation*}
	\begin{split}
	& \int_{r=r_{(\text{min})}}^{\infty} \left( \int_{S_{\tau,r}} r^{p+2-\frac{3}{4}\delta} |\slashed{\D}_L \slashed{\D}_L \phi|^2 r^2 \dVol_{\mathbb{S}^2} \right) \upd r
	\\
	& \lesssim
	\int_{r=r_{(\text{min})}}^{\infty} \left( \int_{S_{\tau,r}} r^{p+2-\frac{3}{4}\delta} |\slashed{\D}_L \slashed{\D}_T \phi|^2 r^2 \dVol_{\mathbb{S}^2} \right) \upd r
	+ \int_{r=r_{(\text{min})}}^{\infty} \left( \int_{S_{\tau,r}} r^{p+2-\frac{3}{4}\delta} |\slashed{\D}_L \slashed{\D}_{\Lbar} \phi|^2 r^2 \dVol_{\mathbb{S}^2} \right) \upd r
	\\
	& \lesssim
	\int_{r=r_{(\text{min})}}^{\infty} \left( \int_{S_{\tau,r}} r^{p+2-\frac{3}{4}\delta} |\slashed{\D}_L \slashed{\D}_T \phi|^2 r^2 \dVol_{\mathbb{S}^2} \right) \upd r
	\\
	&\phantom{\lesssim}
	+ \int_{r=r_{(\text{min})}}^{\infty} \left( \int_{S_{\tau,r}} r^{p+2-\frac{3}{4}\delta} \left( r^{-2}|\overline{\slashed{\D}} (r\slashed{\nabla}\phi)|^2
	+ r^{-2+2C_{(0)}\epsilon}|\slashed{\D}\phi|^2 \right)
	r^2 \dVol_{\mathbb{S}^2} \right) \upd r
	\end{split}
	\end{equation*}
	
	For $p < 1 + \frac{3}{4}\delta - C(\phi)\epsilon$ we have
	\begin{equation*}
	\int_{r=r_{(\text{min})}}^{\infty} \left( \int_{S_{\tau,r}} r^{p+2-\frac{3}{4}\delta} |\slashed{\D}_L \slashed{\D}_T \phi|^2 r^2 \dVol_{\mathbb{S}^2} \right) \upd r
	\lesssim \mathcal{E}_0
	\end{equation*}
	
	For $p \leq \frac{3}{4}\delta - C_{(\phi)}\epsilon$ we have
	\begin{equation*}
	\int_{r=r_{(\text{min})}}^{\infty} \left( \int_{S_{\tau,r}} r^{p+2-\frac{3}{4}\delta} \left( r^{-2}|\overline{\slashed{\D}} (r\slashed{\nabla}\phi)|^2 \right)
	r^2 \dVol_{\mathbb{S}^2} \right) \upd r \lesssim \mathcal{E}_0
	\end{equation*}
	
	Finally, for $p \leq -1 + \frac{1}{4}\delta - 2C_{(0)}\epsilon$ we have
	\begin{equation*}
	\int_{r=r_{(\text{min})}}^{\infty} \left( \int_{S_{\tau,r}} r^{p-\frac{3}{4}\delta + 2C_{(0)}\epsilon} |\slashed{\D}\phi|^2 r^2 \dVol_{\mathbb{S}^2} \right) \upd r
	\lesssim
	\int_{\Sigma_\tau} (1+r)^{-1-\frac{1}{2}\delta} |\slashed{\D}\phi|^2 r^2 \upd r \wedge \dVol_{\mathbb{S}^2}
	\end{equation*}
	and we also have
	\begin{equation*}
	\begin{split}
	\int_{\Sigma_\tau} (1+r)^{-1-\frac{1}{2}\delta} &|\slashed{\D}\phi|^2 r^2 \upd r \wedge \dVol_{\mathbb{S}^2}
	\\
	&\lesssim
	\int_{\Sigma_{\tau_{(\text{max})}}} (1+r)^{-1-\frac{1}{2}\delta} |\slashed{\D}\phi|^2 r^2 \upd r \wedge \dVol_{\mathbb{S}^2}
	\\
	&\phantom{\lesssim}
	+ \int_{\mathcal{M}_\tau^{\tau_{(\text{max})}}} (1+r)^{-1-\frac{1}{2}\delta} \left( \mu|\slashed{\D}\phi|^2 + \mu^{-1}\left| \frac{\partial}{\partial \tau}\Big|_{(r,\vartheta)} \slashed{\D}\phi\right|^2 \right) r^2 \upd \tau \wedge \upd r \wedge \dVol_{\mathbb{S}^2}
	\end{split}
	\end{equation*}
	Now, recall the expression for $\frac{\partial}{\partial \tau}\Big|_{(r,\vartheta)}$ given in \eqref{equation d dtau and ddr on sigma}. Using this, we find
	\begin{equation*}
	\begin{split}
	\left| \frac{\partial}{\partial \tau}\Big|_{(r,\vartheta)} \slashed{\D}\phi\right|^2
	&\lesssim
	\mu^2 |\slashed{\D} \slashed{\D}\phi|^2
	\\
	&\lesssim
	\mu^2 \left( |\slashed{\D} \mathscr{Z}\phi|^2 + |\slashed{\D} \phi|^2 + |F|^2 \right) 
	\end{split}
	\end{equation*}
	and so we find
	\begin{equation*}
	\begin{split}
	\int_{\Sigma_\tau} (1+r)^{-1-\frac{1}{2}\delta} &|\slashed{\D}\phi|^2 r^2 \upd r \wedge \dVol_{\mathbb{S}^2}
	\\
	&\lesssim
	\int_{\Sigma_{\tau_{(\text{max})}}} (1+r)^{-1-\frac{1}{2}\delta} |\slashed{\D}\phi|^2 r^2 \upd r \wedge \dVol_{\mathbb{S}^2}
	\\
	&\phantom{\lesssim}
	+ \int_{\mathcal{M}_\tau^{\tau_{(\text{max})}}} (1+r)^{-1-\frac{1}{2}\delta} \left( |\slashed{\D}\mathscr{Z}\phi|^2 + |\slashed{\D}\phi|^2 + |F|^2 \right) \dVol_g
	\end{split}
	\end{equation*}
	By the integrated local energy inequality (see corollary \ref{corollary ILED decay}), we can choose some $\tau_{(\text{max})}$ so that the first term is arbitrarily small. Moreover, the second term is bounded by $\mathcal{E}_0$ (in fact it decays in $\tau$).
	
	Putting together the calculations above, we see that, for $p \leq -1 + \frac{1}{4}\delta - 2C_{(0)}\epsilon$ we can bound
	\begin{equation*}
	\begin{split}
	\int_{r=r_{(\text{min})}}^{\infty} \left( \int_{S_{\tau,r}} r^{p+2-\frac{3}{4}\delta} |\slashed{\D}_L \slashed{\D}_L \phi|^2 r^2 \dVol_{\mathbb{S}^2} \right) \upd r
	\lesssim
	\mathcal{E}_0
	\end{split}
	\end{equation*}
	
	Returning to equation \eqref{equation internal bound on spherical error terms} we have, for $p \leq -1 + \frac{1}{4}\delta - 2C_{(0)}\epsilon$,
	\begin{equation*}
	\int_{\bar{S}_{t,\tau}} r^{p+3-\delta} |\slashed{\D}_L \phi|^2 \dVol_{\mathbb{S}^2}
	\lesssim
	\delta^{-1}\mathcal{E}_0 \left(r_{(min)}(t)\right)^{-\frac{1}{4}\delta}
	\end{equation*}
	In particular, as $t \rightarrow \infty$ this term tends to zero.
	
	We also need to bound this term for higher values of $p$. In this case, we write
	\begin{equation*}
	\begin{split}
	&\int_{\bar{S}_{t,\tau}} r^{p+3-\delta} |\slashed{\D}_L \phi|^2 \dVol_{\mathbb{S}^2}
	\\
	&\lesssim
	(r_{(\text{min})}(t))^{-\frac{1}{4}\delta} \Bigg( \int_{r = r_{(\text{min})}(t)}^{r_{(\text{max})}} \bigg( \int_{S_{r,\tau}} r^{-1+\frac{1}{4}\delta} \bigg(r^{p+1-\delta}|\slashed{\D}_L (r\slashed{\D}_L \phi)|^2
	\\
	&\phantom{\lesssim (r_{(\text{min})}(t))^{-\frac{1}{4}\delta} \Bigg( \int_{r = r_{(\text{min})}(t)}^{r_{(\text{max})}} \bigg( }
	+ \left(p+1-\frac{3}{4}\delta\right) r^{p+1-\delta}|\slashed{\D}_L \phi|^2 \bigg) r^2 \dVol_{\mathbb{S}^2} \bigg) \upd r
	\\
	&\phantom{\lesssim (1+r_{(\text{min})}(t))^{-\frac{1}{4}\delta} \Bigg(}
	+ \int_{S_{r_{(\text{max})},\tau}} r^{p+3-\frac{1}{2}\delta}|\slashed{\D}_L \phi|^2 \dVol_{\mathbb{S}^2}
	\Bigg)
	\end{split}
	\end{equation*}
	So, if we suppose that we \emph{already} know that
	\begin{equation*}
	\int_{\Sigma_\tau \cap \{r \geq r_0\}} r^{p - \frac{3}{4}\delta} |\slashed{\D}_L (r\slashed{\D}_L \phi)|^2 \dVol_{\Sigma_\tau} \lesssim \delta^{-1}\mathcal{E}_0
	\end{equation*}
	then we can conclude that
	\begin{equation*}
	\int_{\bar{S}_{t,\tau}} r^{p+3-\delta}|\slashed{\D}_L\phi|^2 \dVol_{\mathbb{S}^2}
	\lesssim
	\delta^{-1}\mathcal{E}_0 \left(r_{(\text{min})}(t)\right)^{-\frac{1}{4}\delta}
	\end{equation*}
	
	Next, we turn to the error terms in the region $\frac{1}{2}R \leq r \leq R$. We have
	\begin{equation*}
	\begin{split}
	&\int_{\tau_0}^\tau \left( \int_{^t\Sigma_{\tau'} \cap \{\frac{1}{2}R \leq r \leq R \}} r^{p-1}|\overline{\slashed{\D}}(r\slashed{\D}_L \phi)|^2 \dVol_{\Sigma_\tau} \right) \upd \tau
	\\
	&\lesssim
	\int_{\tau_0}^\tau \bigg( \int_{^t\Sigma_{\tau'} \cap \{\frac{1}{2}R \leq r \leq R \}} r^{p-1} \Big( r^2 |\slashed{\D}_L \mathscr{Z}\phi|^2 + r^2 |F|^2 + |\slashed{\D}\phi|^2 + r^{C_{(0)}\epsilon}  |\overline{\slashed{\D}}\phi|^2
	\\
	&\phantom{ \lesssim \int_{\tau_0}^\tau \bigg( \int_{^t\Sigma_{\tau'} \cap \{\frac{1}{2}R \leq r \leq R \}} r^{p-1} \Big(}
	+ r^{-2+2C_{(0)}\epsilon} |\phi|^2 \Big) \dVol_{\Sigma_\tau} \bigg) \upd \tau
	\end{split}
	\end{equation*}
	And so, for $p \leq -\delta$ we can bound this term by
	\begin{equation*}
	\int_{\tau_0}^\tau \left( \int_{^t\Sigma_{\tau'} \cap \{\frac{1}{2}R \leq r \leq R \}} r^{p-1}|\overline{\slashed{\D}}(r\slashed{\D}_L \phi)|^2 \dVol_{\Sigma_\tau} \right) \upd \tau
	\lesssim
	(\tau-\tau_0) \delta^{-1} \mathcal{E}_0 R^{-1}
	\end{equation*}
	In particular, as we take the limit $R \rightarrow \infty$ this term tends to zero.
	
	On the other hand, for larger values of $p$ we can estimate
	\begin{equation*}
	\begin{split}
	&\int_{\tau_0}^\tau \left( \int_{^t\Sigma_{\tau'} \cap \{\frac{1}{2}R \leq r \leq R \}} r^{p-1}|\overline{\slashed{\D}}(r\slashed{\D}_L \phi)|^2 \dVol_{\Sigma_\tau} \right) \upd \tau
	\\
	&\lesssim
	\int_{\tau_0}^\tau \left( \int_{^t\Sigma_{\tau'} \cap \{\frac{1}{2}R \leq r \leq R \}}
	\left( r^{p-1}|\slashed{\D}_L (r\slashed{\D}_L \phi)|^2 
	+ r^{p-1}|r\slashed{\nabla}\slashed{\D}_L \phi)|^2 
	\right)
	\dVol_{\Sigma_\tau} \right) \upd \tau
	\\
	&\lesssim
	\int_{\tau_0}^\tau \bigg( \int_{^t\Sigma_{\tau'} \cap \{\frac{1}{2}R \leq r \leq R \}}
	\Big( r^{p-1}|\slashed{\D}_L (r\slashed{\D}_L \phi)|^2 
	+ r^{p-1}|\slashed{\D}_L (r\slashed{\nabla}\phi)|^2
	+ \epsilon^2 r^{p-1-2\delta}|\slashed{\nabla}\phi|^2
	\\
	&\phantom{\lesssim \int_{\tau_0}^\tau \bigg( \int_{^t\Sigma_{\tau'} \cap \{\frac{1}{2}R \leq r \leq R \}} \Big(}
	+ \epsilon^2 r^{p-3-2\delta}|\phi|^2
	\Big)
	\dVol_{\Sigma_\tau} \bigg) \upd \tau
	\end{split}
	\end{equation*}
	where we have used proposition \ref{proposition commute rnabla} to commute $r\slashed{\nabla}$ and $\slashed{\D}_L$, as well as proposition \ref{proposition expression for Omega} to handle the lower order terms. Hence, \emph{if we already have}
	\begin{equation*}
	\int_{^t\Sigma_{\tau'} \cap \{r \geq r_0\}}
	\Big( r^{p-1}|\slashed{\D}_L (r\slashed{\D}_L \phi)|^2 \Big)\dVol_{\Sigma_\tau}
	\lesssim \delta^{-1}\mathcal{E}_0
	\end{equation*}
	then we can bound
	\begin{equation*}
	\int_{\tau_0}^\tau \left( \int_{^t\Sigma_{\tau'} \cap \{\frac{1}{2}R \leq r \leq R \}} r^{p-1}|\overline{\slashed{\D}}(r\slashed{\D}_L \phi)|^2 \dVol_{\Sigma_\tau} \right) \upd \tau
	\lesssim
	\delta^{-1}\mathcal{E}_0 (\tau - \tau_0) R^{-1}
	\end{equation*}
	Despite the fact that this grows in $\tau$, it also vanishes in the limit $R \rightarrow \infty$ as required.
	
	Finally, we encounter an error term involving the integral over the cylinder $r = R$. We can bound this term in a similar way to the previous bounds of error terms on the spheres $\bar{S}_{\tau,t}$. Specifically, for $p \leq -\frac{3}{4}\delta - 2C_{(0)}\epsilon$ we have
	\begin{equation*}
	\int_{S_{\tau,R}} r^{p+2}|\slashed{\D}_L \phi|^2 \dVol_{\mathbb{S}^2}
	\lesssim
	\delta^{-1}\mathcal{E}_0 R^{-\frac{1}{4}\delta}
	\end{equation*}
	
	On the other hand, \emph{if we already know that}
	\begin{equation*}
	\int_{\Sigma_\tau \cap \{r \geq r_0\}} r^{p -1 + \frac{1}{4}\delta} |\slashed{\D}_L (r\slashed{\D}_L \phi)|^2 \dVol_{\Sigma_\tau} \lesssim \delta^{-1}\mathcal{E}_0
	\end{equation*}
	then we can conclude that
	\begin{equation*}
	\int_{\bar{S}_{t,\tau}} r^{p+2}|\slashed{\D}_L\phi|^2 \dVol_{\mathbb{S}^2}
	\lesssim
	\delta^{-1}\mathcal{E}_0 R^{-\frac{1}{4}\delta}
	\end{equation*}
	Note that these terms also have to be integrated over $\tau$, but, importantly, they all tend to zero as $R \rightarrow \infty$.

\end{proof}

\begin{corollary}[Decay of the $p$-weighted energy with $p = \delta$]
	Suppose that the conditions of lemma \ref{lemma p weighted for rL} hold.
	
	Then we have
	\begin{equation}
	\mathcal{E}^{(L, \delta)}[r\slashed{\D}_L \phi](\tau) \lesssim \delta^{-1} \mathcal{E}_{(0)} (1+\tau_n)^{-1+C_{(\phi)}\delta}
	\end{equation}
\end{corollary}

\begin{proof}
	Applying lemma \ref{lemma p weighted for rL} with the choice $p = 1-C_{(\phi)}\epsilon$ we obtain, in particular, for all $\tau \geq \tau_0$,
	\begin{equation*}
	\int_{\mathcal{M}_{\tau_0}^{\tau}} \chi_{(2r_0)} r^{-C_{(\phi)}\epsilon} \left|\slashed{\D}_L (r\slashed{\D}_L \phi) \right|^2 \dVol_g \lesssim (1+\delta^{-1}) \mathcal{E}_0
	\end{equation*}
	so we can pick a diadic sequence of times $\tau_n \rightarrow \infty$ such that
	\begin{equation*}
	\int_{\Sigma_{\tau_n}} \chi_{(2r_0)} r^{-C_{(\phi)}\epsilon} |\slashed{\D}_L (r\slashed{\D}_L \phi)|^2 \dVol_{\Sigma_\tau}
	\lesssim
	(1+\delta^{-1})\mathcal{E}_0 (1+\tau_n)^{-1}
	\end{equation*}
	Additionally, at these times (and in fact at all times $\tau \geq \tau_0$) we have
	\begin{equation*}
	\mathcal{E}^{(L, 1-C_{(\phi)}\epsilon)}[r\slashed{\D}_L \phi](\tau) \lesssim (1+\delta^{-1})\mathcal{E}_0
	\end{equation*}
	Interpolating between these two inequalities using H\"older's inequality, and using the fact that $\delta^{-1} \geq 1$ we have
	\begin{equation*}
	\mathcal{E}^{(L, \delta)}[r\slashed{\D}_L \phi](\tau_n) \lesssim \delta^{-1}\mathcal{E}_0 (1+\tau_n)^{-1+\delta+C_{(\phi)}\epsilon}
	\end{equation*}
	
	Next, we appeal to lemma \ref{lemma p weighted for rL} again, this time with the choice $p = \delta$, for some time $\tau$ satisfying $\tau_n \leq \tau \leq \tau_{n+1}$. We find
	\begin{equation*}
	\begin{split}
	\mathcal{E}^{(L, \delta)}[r\slashed{\D}_L \phi](\tau)
	&\lesssim
	\mathcal{E}^{(L, \delta)}[r\slashed{\D}_L \phi](\tau_n)
	+ \delta^{-1} \mathcal{E}_{(0)} (1+\tau_n)^{-1+C_{(\phi)}\delta}
	\\
	&\lesssim
	\delta^{-1} \mathcal{E}_{(0)} (1+\tau_n)^{-1+C_{(\phi)}\delta}
	\end{split}
	\end{equation*}
	proving the corollary.
	
\end{proof}

\section{Improved energy estimates for \texorpdfstring{$\slashed{\D}_T \phi$}{DT phi}}

The lemma above can be used to obtain improved bounds for quantities involving $L$ derivatives. In particular, it will play an important role in proving improved pointwise decay in $r$ for $L$ derivatives, which in turn is necessary to control certain error terms arising when commuting with $rL$. On the other hand, these estimates also lead to improved decay in $\tau$ for $T$ derivatives, established in the following lemma.

\begin{lemma}[Improved energy decay for $\slashed{\D}_T \phi$]
	\label{lemma improved decay DT phi}
	Suppose that the same conditions as those assumed in lemma \ref{lemma p weighted for rL} hold. Suppose, in addition, that $\slashed{\D}_T \phi$ satisfies
	\begin{equation*}
	\tilde{\slashed{\Box}}_g (\slashed{\D}_T \phi) 
	=
	\slashed{\Delta} \slashed{\D}_T \phi
	+ (2^k - 1) r^{-1} \slashed{\D}_L (r\slashed{\D}_L \slashed{\D}_T  \phi)
	+ (2^k - 1) r^{-1} \slashed{\D}_L (\slashed{\D}_T \phi)
	+ F_{(T)}
	\end{equation*}
	where 
	\begin{equation*}
	F_{(T)} = F_{(T,1)} + F_{(T,2)}
	\end{equation*}	
	and the $S_{\tau,r}$-tangent tensor fields $F_{(1,T)}$ and $F_{(2,T)}$ satisfy	
	\begin{equation*}
	\begin{split}
	\int_{\mathcal{M}_\tau^{\tau_1}} \left( \epsilon^{-1} \chi_{(r_0)} r^{1-C_{(\phi)}\epsilon}(1+\tau)^{1+\delta}|F_{(T)}|^2 \right) \dVol_g 
	&\lesssim \mathcal{E}_0 (1+\tau)^{1-K} \\ \\
	\int_{\mathcal{M}^{\tau_1}_{\tau}} \epsilon^{-1}\left( (1+r)^{1-C_{(\phi)}\epsilon}|F_{(T,1)}|^2 + (1+r)^{1-\delta}(1+\tau)^{6\delta}|F_{(T,2)}|^2 \right) \dVol_g
	&\lesssim \mathcal{E}_0(1+\tau))^{-K}
	\end{split}
	\end{equation*}
	for some $C_{(\phi)} > 0$ sufficiently large, and for some constant $K > 0$.
	
	Finally, suppose that the $\mathscr{Z}^{(i)}\phi$ (for $i = 0, 1$)satisfies
	\begin{equation*}
	\tilde{\slashed{\Box}}_g \mathscr{Z}^{(i)}\phi
	=
	\slashed{\Delta} \mathscr{Z}^{(i)}\phi
	+ (2^k - 1) r^{-1} \slashed{\D}_L (r\slashed{\D}_L \mathscr{Z}^{(i)}\phi)
	+ (2^k - 1) r^{-1} \slashed{\D}_L (\mathscr{Z}^{(i)}\phi)
	+ F_{(\mathscr{Z}^{(i)})}
	\end{equation*}
	where $F_{(\mathscr{Z}^{(i)})}$ satisfies 
	\begin{equation*}
	\begin{split}
	\int_{\mathcal{M}_\tau^{\tau_1}} \left( \epsilon^{-1} \chi_{(r_0)} r^{1-C_{(\phi)}\epsilon}(1+\tau)^{1+\delta}|F_{(\mathscr{Z}^{(i)})}|^2 \right) \dVol_g 
	&\lesssim \mathcal{E}_0 (1+\tau)^{C_{(\mathscr{Z}^{(i)}\phi)}\delta} \\ \\
	\int_{\mathcal{M}^{\tau_1}_{\tau}} \epsilon^{-1}\left( (1+r)^{1-C_{(\phi)}\epsilon}|F_{(T,1)}|^2 + (1+r)^{1-\delta}(1+\tau)^{2\beta}|F_{(T,2)}|^2 \right) \dVol_g
	&\lesssim \mathcal{E}_0(1+\tau))^{-1+C_{(\mathscr{Z}^{(i)}\phi)}\delta}
	\end{split}
	\end{equation*}
	
	Then we can improve the decay of the energy of $\slashed{\D}_T \phi$ to find
	\begin{equation*}
	\mathcal{E}^{(wT)}[\slashed{\D}_T \phi](\tau) + \mathcal{E}^{(L, \frac{1}{2}\delta)}[\slashed{\D}_T \phi](\tau) \lesssim \delta^{-9}\mathcal{E}_0 (1+\tau)^{-1 -K^*}
	\end{equation*}
	where $K^*$ is defined by
	\begin{equation*}
	K^* := \min\{ K, 2-C_{(\slashed{\D}_T \phi)}\delta\}
	\end{equation*}

\end{lemma}

\begin{proof}
	Using lemma \ref{lemma p weighted for rL} with the choice $p = 1-C_{(r\slashed{\D}_L \phi)}\epsilon$, we find, in particular, that for all $\tau \geq \tau_0$
	\begin{equation*}
	\int_{\mathcal{M}_{\tau_0}^\tau} \chi_{(2r_0)} r^{-C_{(r\slashed{\D}_L \phi)}\epsilon} |\slashed{\D}_L (r\slashed{\D}_L \phi)|^2 \dVol_g \lesssim \delta^{-1} \mathcal{E}_0 (1+\tau)^{C_{(r\slashed{\D}_L\phi)}\delta}
	\end{equation*}
	So, we can pick a sequence of times $\tau_n \rightarrow \infty$ such that
	\begin{equation*}
	\int_{\Sigma_{\tau_n}} \chi_{(2r_0)} \left( r^{-C_{(r\slashed{\D}_L \phi)}\epsilon} |\slashed{\D}_L (r\slashed{\D}_L \phi)|^2 \right) r^2 \upd r \wedge \dVol_{\mathbb{S}^2} \lesssim \delta^{-1} \mathcal{E}_0 (1+\tau_n)^{-1+C_{(r\slashed{\D}_L\phi)}\delta}
	\end{equation*}
	
	Now, we can also apply lemma \ref{lemma p weighted for rL} with the choice $p = -C_{(r\slashed{\D}_L \phi)}\epsilon$. Choosing the initial time for this estimate to be one of the $\tau_n$ and repeating the calculation above, we obtain a new sequence of times $\tau'_n \rightarrow \infty$ such that
	\begin{equation*}
	\int_{\Sigma_{\tau'_n}} \chi_{(2r_0)} \left( r^{-1-C_{(r\slashed{\D}_L \phi)}\epsilon} |\slashed{\D}_L (r\slashed{\D}_L \phi)|^2 \right) r^2 \upd r \wedge \dVol_{\mathbb{S}^2} \lesssim \delta^{-1}\mathcal{E}_0 (1+\tau'_n)^{-2 + \frac{3}{4}\delta}
	\end{equation*}
	
	Next, we use $L = 2T - \Lbar$ to write
	\begin{equation*}
	\begin{split}
	r^{-1-C_{(r\slashed{\D}_L \phi)}\epsilon} |\slashed{\D}_L (r\slashed{\D}_L \phi)|^2 
	&\gtrsim
	r^{-1-C_{(r\slashed{\D}_L \phi)}\epsilon}|\slashed{\D}_L (r\slashed{\D}_T \phi)|^2
	- r^{-1-C_{(r\slashed{\D}_L \phi)}\epsilon}|\slashed{\D}_L \left(r \slashed{\D}_{\Lbar} \phi\right)|^2
	\\
	&\gtrsim
	r^{-1-C_{(r\slashed{\D}_L \phi)}\epsilon}|\slashed{\D}_L (r\slashed{\D}_T \phi)|^2
	- r^{-1-C_{(r\slashed{\D}_L \phi)}\epsilon}|\slashed{\nabla}\mathscr{Z} \phi|^2
	\\
	&\phantom{\gtrsim}
	- r^{-1-C_{(r\slashed{\D}_L \phi)}\epsilon + 2C_{(0)}\epsilon} |\overline{\slashed{\D}} \phi|^2
	- \epsilon^2 r^{-1 - 2\delta + 2C_{(0)}\epsilon} |\slashed{\D} \phi|^2
	- \epsilon^2 r^{-3 - 2\delta + 2C_{(0)}\epsilon} |\phi|^2
	\\
	&\phantom{\gtrsim}
	- r^{1-C_{(r\slashed{\D}_L \phi)}\epsilon} |\tilde{\slashed{\Box}}_g \phi|^2
	\end{split}
	\end{equation*}
	and so we have a sequence of times such that
	\begin{equation*}
	\begin{split}
	&\int_{\Sigma_{\tau'_n}} \chi_{(2r_0)} \left( 
	r^{1-C_{(r\slashed{\D}_L \phi)}\epsilon } |\slashed{\D}_L \slashed{\D}_T \phi|^2 
	\right) r^2 \upd r \wedge \dVol_{\mathbb{S}^2}
	\\
	&\lesssim
	\delta^{-1}\mathcal{E}_0 (1+\tau'_n)^{-2+\frac{3}{4}\delta}
	\\
	&\phantom{\lesssim}
	+ \int_{\Sigma_{\tau'_n}} \chi_{(2r_0)} \bigg(
	r^{-1-C_{(r\slashed{\D}_L \phi)}\epsilon} |\overline{\slashed{\D}} \mathscr{Z}\phi|^2
	r^{-1-C_{(r\slashed{\D}_L \phi)}\epsilon} |\overline{\slashed{\D}} \phi|^2
	+ r^{-1-C_{(r\slashed{\D}_L \phi)}\epsilon + 2C_{(0)}\epsilon }|\overline{\slashed{\D}} \phi|^2
	\\
	&\phantom{\lesssim + \int_{\Sigma_{\tau'_n}} \chi_{(2r_0)} \bigg(}
	+ \epsilon^2 r^{-1-\delta}|\slashed{\D}\phi|^2 + r^{-3-\delta }|\phi|^2
	+ r^{1-C_{(r\slashed{\D}_L \phi)}\epsilon}|F|^2
	\bigg) r^2 \upd r \wedge \dVol_{\mathbb{S}^2}
	\end{split}
	\end{equation*}
	
	Now, if we use corollaries \ref{corollary ILED decay}  and \ref{corollary small p decay} applied to the fields $\phi$ and $\mathscr{Z}\phi$, as well as the assumptions on $F$ and $F_{(\mathscr{Z})}$, we find that
	\begin{equation*}
	\begin{split} 
	& \int_{\mathcal{M}_{\tau}^{\tau_1}} \bigg( 
	(1+r)^{-1-\frac{1}{2}\delta}|\slashed{\D}\phi|^2
	+ (1+r)^{-3-\frac{1}{2}\delta}|\phi|^2
	+ (1+r)^{-1+\frac{1}{2}\delta}|\overline{\slashed{\D}}\phi|^2
	+ (1+r)^{1-\frac{1}{2}C_{(T)}\epsilon} |F|^2
	\\
	&\phantom{ \int_{\mathcal{M}_{\tau}^{\tau_1}} \bigg( }
	+ (1+r)^{-1-\frac{1}{2}\delta}|\slashed{\D}\mathscr{Z}\phi|^2
	+ (1+r)^{-3-\frac{1}{2}\delta}|\mathscr{Z}\phi|^2
	+ (1+r)^{-1+\frac{1}{2}\delta}|\overline{\slashed{\D}}\mathscr{Z}\phi|^2
	\\
	&\phantom{ \int_{\mathcal{M}_{\tau}^{\tau_1}} \bigg( }
	+ (1+r)^{1-\frac{1}{2}C_{(T)}\epsilon} |F_{(\mathscr{Z})}|^2
	\bigg) \dVol_g
	\\
	& \lesssim
	\delta^{-1} \mathcal{E}_0 (1+\tau)^{-1 + C_{(\mathscr{Z}\phi)}\delta}	
	\end{split}
	\end{equation*}
	where we have assumed that $C_{(\mathscr{Z}\phi)} \geq C_{(\phi)}$. Hence there is a sequence of times $\tau''_n$ such that
	\begin{equation*}
	\begin{split}
	& \int_{\Sigma_{\tau''_n}} \bigg( 
	(1+r)^{-1-\frac{1}{2}\delta}|\slashed{\D}\phi|^2
	+ (1+r)^{-3-\frac{1}{2}\delta}|\phi|^2
	+ (1+r)^{-1+\frac{1}{2}\delta}|\overline{\slashed{\D}}\phi|^2
	+ (1+r)^{1-\frac{1}{2}C_{(T)}\epsilon} |F|^2
	\\
	&\phantom{\int_{\Sigma_{\tau''_n}} \bigg( }
	+ (1+r)^{-1-\frac{1}{2}\delta}|\slashed{\D}\mathscr{Z}\phi|^2
	+ (1+r)^{-3-\frac{1}{2}\delta}|\mathscr{Z}\phi|^2
	+ (1+r)^{-1+\frac{1}{2}\delta}|\overline{\slashed{\D}}\mathscr{Z}\phi|^2
	\\
	&\phantom{\int_{\Sigma_{\tau''_n}} \bigg( }
	+ (1+r)^{1-\frac{1}{2}C_{(T)}\epsilon} |F_{(\mathscr{Z})}|^2
	\bigg) r^2 \upd r \wedge \dVol_{\mathbb{S}^2}
	\\
	&\lesssim
	\delta^{-1} \mathcal{E}_0 (1+\tau'')^{-2 + C_{(\mathscr{Z}\phi)}\delta}	
	\end{split}
	\end{equation*}
	Moreover (for example, by adding together the two spacetime integrals that led to these inequalities) we can arrange things so that the times $\tau''_n$ to coincide with the times $\tau'_n$. Hence, we have
	\begin{equation*}
	\int_{\Sigma_{\tau'_n}} \chi_{(2r_0)} \left( r^{-1-C_{(T)}\epsilon} |\slashed{\D}_L (r \slashed{\D}_T \phi)|^2 \right) r^2 \upd r \wedge \dVol_{\mathbb{S}^2}
	\lesssim
	\delta^{-3}\mathcal{E}_0 (1+\tau'_n)^{-2+C_{(\mathscr{Z}\phi)}\delta}
	\end{equation*}
	
	Next, we want to apply the $p$-weighted energy estimate to the field $\slashed{\D}_T \phi$, with the choice $p = 1-C_{(T)}\epsilon$. From lemma \ref{lemma p weighted} applied to the field $\slashed{\D}_T \phi$, we obtain
	\begin{equation*}
	\mathcal{E}^{(wT)}[\slashed{\D}_T \phi](\tau) + \mathcal{E}^{(L, \frac{1}{2}\delta)}[\slashed{\D}_T \phi](\tau) \lesssim \delta^{-3}\mathcal{E}_0 (1+\tau)^{-1 + C_{(\slashed{\D}_T \phi)}\delta}
	\end{equation*}
	where
	\begin{equation*}
	w = (1+r)^{-C_{(\slashed{\D}_T \phi)}\epsilon}
	\end{equation*}
	Note that his holds for all $\tau \geq \tau_0$, so \emph{a fortiori} it holds for $\tau = \tau'_n$. In particular, at the times $\tau'_n$ we have
	\begin{equation*}
	\mathcal{E}^{(L, 1-C_{(T)}\epsilon)}[\slashed{\D}_T \phi](\tau'_n)
	+ \mathcal{E}^{(wT)}[\slashed{\D}_T \phi](\tau'_n) 
	+ \mathcal{E}^{(L, \frac{1}{2}\delta)}[\slashed{\D}_T \phi](\tau'_n)
	\lesssim
	\delta^{-3} \mathcal{E}_0 (1+\tau'_n)^{-1+C_{(\slashed{\D}_T \phi)} \delta} 
	\end{equation*}

	Following exactly the same steps as in lemma \ref{lemma p weighted}, but now using the fact that
	\begin{equation*}
	\begin{split}
	\int_{\mathcal{M}_\tau^{\tau_1}} \left( \epsilon^{-1} \chi_{(r_0)} r^{1-C_{(T)}\epsilon}(1+\tau)^{1+\delta}|F_{(T)}|^2 \right) \dVol_g 
	&\lesssim \mathcal{E}_0 (1+\tau)^{-K} \\ \\
	\int_{\mathcal{M}^{\tau_1}_{\tau}} \epsilon^{-1}\left( (1+r)^{1-C_{(T)}\epsilon}|F_{(T,1)}|^2 + (1+r)^{1-\delta}(1+\tau)^{2\beta}|F_{(T,2)}|^2 \right) \dVol_g
	&\lesssim \mathcal{E}_0(1+\tau))^{-1-K}
	\end{split}
	\end{equation*}
	we can improve the decay of the energy of $\slashed{\D}_T \phi$ to find, in analogy with corollaries \ref{corollay energy decay} and \ref{corollary small p decay}
	\begin{equation*}
	\mathcal{E}^{(wT)}[\slashed{\D}_T \phi](\tau) + \mathcal{E}^{(L, \frac{1}{2}\delta)}[\slashed{\D}_T \phi](\tau) \lesssim \delta^{-6}\mathcal{E}_0 (1+\tau)^{-1-K^*}
	\end{equation*}
	where
	\begin{equation*}
	K^* := \min\{ K, 1-C_{(\slashed{\D}_T \phi)}\delta\}
	\end{equation*}
	
	Note also that these calculations allows us to drop the reliance on a subsequence, and to show that \emph{for all} $\tau \geq \tau_0$ we have
	\begin{equation*}
	\mathcal{E}^{(L, 1-C_{(T)}\epsilon)}[\slashed{\D}_T \phi](\tau) \lesssim \delta^{-7}\mathcal{E}_0 (1+\tau)^{-1-K^*}
	\end{equation*}
	
	If $K > 1-C_{(T)}\epsilon$, so that $K^* = 1-C_{\slashed{\D}_T \phi}$, then we can repeat this argument, with the new and improved boundary estimates. In fact, this means that we can take
	\begin{equation*}
	K^* := \min\{ K, 2-C_{(\slashed{\D}_T \phi)}\delta\}
	\end{equation*}
	so that if $K$ is sufficiently large, we have energy decay like $\tau^{-3 + \delta}$.
	
\end{proof}

Note that, if the inhomogeneous terms after commuting with $T$ decay sufficiently fast so that $K > 1$, then the square root of the energy is integrable in $\tau$. If the same statement is true for higher derivatives, then we could obtain \emph{pointwise} estimates on the $T$ derivatives of the fields which are integrable in $\tau$.

Note also that the proposition above can be used to provide improved decay in $\tau$ in a bounded region, but since this requires different bootstrap assumptions (specifically, additional decay in $\tau$ for certain quantities) we will not perform this estimate here.

\chapter{Semi-global existence and uniqueness}
\label{appendix semi-global existence}

In this appendix we outline the proof of semi-global existence and uniqueness for the kind of wave equations we have considered in this work. By ``semi-global'' existence, we mean \emph{local in} $\tau$ but global in $r$. Due to the nature of the foliation by leaves of constant $\tau$, this ``local'' existence should already be seen as ``semi-global'' - it is global in $r$ but local in $\tau$, and certain kinds of wave equations (e.g.\ those which do not satisfy the weak null condition) do not admit a result of this kind. Note also that we have to work with finite degenerate energy \emph{and} finite $p$-weighted energy, rather than the more standard finite initial energy.

The first step of the proof is to establish existence and uniqueness for a suitable \emph{linear} equation. The linear equation that we need to solve is
\begin{equation*}
\tilde{\slashed{\Box}}_g \phi = F
\end{equation*}
where we now consider $g$ and $F$ as a \emph{fixed} metric and inhomogeneous term respectively. It is important that we use the \emph{reduced} wave operator $\tilde{\slashed{\Box}}_g$ rather than the standard wave operator $\slashed{\Box}_g$ here, since we will later iterate these solutions, and it is only solutions to the linear equations \emph{involving the reduced wave operator} which will have the right asymptotics.

\begin{figure}[htb]
	\centering
	\includegraphics[width = 0.9\linewidth, keepaspectratio]{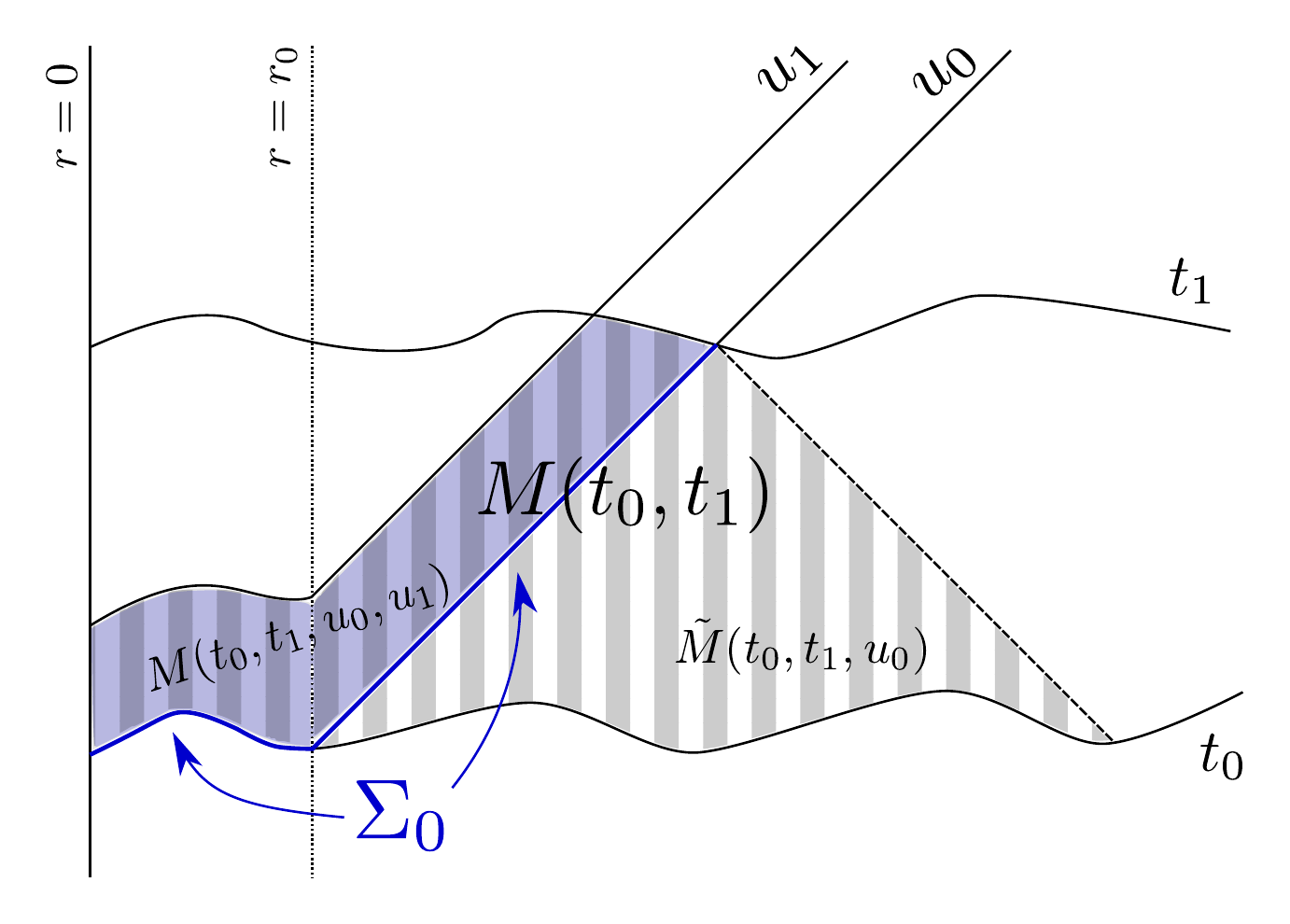}
	\caption{
	A figure showing the various spacetime regions for use in the linear/semi-global existence problems. Initial data is posed on the hypersurface $\Sigma_0$, and we seek a solution in the (shaded blue) region $M(t_0, t_1, u_0, u_1)$. We produce this solution by first defining the functions $\phi_2$ and $\rho$ in the entire striped region, which includes a region to the past of $\Sigma_0$. We then solve a wave equation for a quantitiy $\chi$, with initial data posed on the \emph{spacelike} surface $t = t_0$, again finding a solution in the striped region. Finally, the solution we are seeking is constructed from both $\chi$ and $\phi_2$.
	}
	\label{figure spacetime regions local}
\end{figure}

To obtain semi-global existence and uniqueness for the \emph{linear} equations, we can (very slightly) adapt the approach of \cite{Rendall1990}. Consider the following situation: we are given ``initial'' data on a spacelike hypersurface surface $t = t_0 \cap \{u_0 \leq u \leq u_1\}$, together with data on the null surface $u = u_0 \cap \{t_0 \leq t \leq t_1\}$. Here, $u_0$, $u_1$, $t_0$ and $t_1$ are constants. Then we will produce a solution to the wave equation in the region $M(t_0, t_1, u_0, u_1) := \{p \in \mathcal{M} | t_0 \leq t(p) \leq t_1 \ \& \  u_0 \leq u(p) \leq u_1\}$. See figure \ref{figure spacetime regions local} for a diagram of these hypersurfaces and the spacetime region.

Consider the integral curves of $\Lbar$ through the surface $u = u_0 \cap \{t_0 \leq t \leq t_1\}$. Following these backwards (i.e.\ in the direction of decreasing $t$, that is, increasing $r$), we find that they intersect the hypersurface $t = t_0$ at some points. Again, see figure \ref{figure spacetime regions local}. The idea, following Rendall \cite{Rendall1990}, is to find a function $\phi_1$ such that
\begin{itemize}
	\item $\phi_1 = \phi$ on the hypersurface $u = u_0 \cap \{t_0 \leq t \leq t_1\}$
	\item transverse derivatives of $\phi_1$ match the transverse derivatives of $\phi$ on the hypersurface $u = u_0 \cap \{t_0 \leq t \leq t_1\}$, if $\phi$ satisfies the wave equation near this hypersurface
\end{itemize}

Let us explain this process in more detail. Define three functions\footnote{Here we are assuming that we are given \emph{smooth} initial data. If this is not the case, and the initial data is only in some Sobolev space, then we can simply approximate the data by smooth data and take the limit.}
\begin{equation*}
\begin{split}
	f_1 &\in C^\infty(\{t = t_0 \cap \{u_0 \leq u \leq u_1\}\} \rightarrow \mathbb{R})
	\\
	f_2 &\in C^\infty(\{t = t_0 \cap \{u_0 \leq u \leq u_1\}\} \rightarrow \mathbb{R})
	\\
	f_3 &\in C^\infty(\{u = u_0 \cap \{t_0 \leq t \leq t_1\}\} \rightarrow \mathbb{R})
\end{split}
\end{equation*}
which are such that
\begin{equation*}
f_1\big|_{\{t = t_0, u = u_0\}}
= f_3\big|_{\{t = t_0, u = u_0\}}
\end{equation*}
We want to find a function $\phi$ such that
\begin{equation*}
\begin{cases}
	\phi\big|_{\{t = t_0 \cap \{u_0 \leq u \leq u_1\}\}} &= f_1
	\\
	(\partial_t \phi)\big|_{\{t = t_0 \cap \{u_0 \leq u \leq u_1\}\}} &= f_2
	\\
	\phi\big|_{\{u = u_0 \cap \{t_0 \leq t \leq t_1\}\}} &= f_3
\end{cases}
\end{equation*}
Note that we do not specify the derivatives of $\phi$ on the null part of the initial data surface.

Using the expression for the wave operator given in proposition \ref{proposition scalar wave operator}, if $\phi$ satisfies the wave equation $\tilde{\Box}_g \phi = F$ then we have
\begin{equation}
L\left( r\Lbar\phi\right)
	+ \frac{1}{2} (\tr_{\slashed{g}}\chi_{(\text{small})}) r\Lbar\phi 
	=
	r\slashed{\Delta}\phi 
	- \frac{1}{2}(\tr_{\slashed{g}}\chibar )r L\phi 
	- \zeta^\alpha r \slashed{\nabla}_\alpha \phi
	- rF
\end{equation}
We can use this equation to find the initial value that $r\Lbar\phi$ should take along the null part of the initial data. Specifically, we can define the quantity $r\Lbar\phi_1$ \emph{on the hypersurface $u = u_0$} by the linear ODE
\begin{equation*}
L\left( r\Lbar\phi_1\right)
	+ \frac{1}{2} (\tr_{\slashed{g}}\chi_{(\text{small})}) (r\Lbar\phi_1 )
	=
	r\slashed{\Delta}f_3 
	- \frac{1}{2}(\tr_{\slashed{g}}\chibar )r Lf_3 
	- \zeta^\alpha r \slashed{\nabla}_\alpha f_3
	- rF
\end{equation*}
Note that we only take derivatives of $f_3$ that are \emph{tangent} to the initial data surface, so these quantities are well defined. The initial values for these ODEs is defined as follows:
\begin{equation*}
r\Lbar \phi_1 \big|_{t = t_0, u = u_0} = 
r\Lbar^0 f_2 \big|_{t = t_0, u = u_0}
+ \lim_{u_\epsilon \searrow u_0} r\Lbar^i (\partial_i f_1) \big|_{t = t_0, u = u_\epsilon}
\end{equation*}
Note that the spatial derivatives of $f_1$ have to be defined as ``one-sided'' derivatives since $f_1$ is only defined for $u \geq u_0$.

Importantly, if we have the bounds
\begin{equation*}
\begin{split}
	|rF| &\lesssim \epsilon (1+r)^{-1+C\epsilon}
	\\
	|rL f_3| &\lesssim \epsilon (1+r)^{-\delta}
	\\
	|r\slashed{\Delta} f_3| &\lesssim \epsilon (1+r)^{-1-\delta}
\end{split}
\end{equation*}
and if the usual pointwise bootstrap bounds hold on the geometric quantities, then the solution to the linear ODE above satisfies the bound
\begin{equation*}
|r\Lbar \phi_1| \lesssim \left( r|f_2| + r|\underline{\partial} f_1| \right)\big|_{\{t = t_0\} \cap \{u = u_0\}} (1+r)^{C\epsilon}
\end{equation*}
where we recall that $\underline{\partial}f_1$ refers to the \emph{spatial} derivatives of $f_1$. Note that $r$ is bounded on the ``sphere'' $\{t = t_0\} \cap \{u = u_0\}$, so this gives us that $\Lbar \phi_1 \sim r^{-1+C\epsilon}$.

Higher transverse derivatives of $\phi_1$ can be determined by commuting with the vector field $T$ and then using the fact that $T = \frac{1}{2}(L + \Lbar)$. For example, suppose that we set $\phi_1\big|_{\{u = u_0\}} = f_3$ and we have also found $\Lbar \phi_1 \big|_{\{u = u_0\}}$ by following the process above. Then, to find $r\Lbar T \phi_1$ we have to solve the ODE
\begin{equation*}
	L\left( r\Lbar T\phi_1\right)
	+ \frac{1}{2} (\tr_{\slashed{g}}\chi_{(\text{small})}) (r\Lbar T\phi_1 )
	=
	r\slashed{\Delta} T\phi_1 
	- \frac{1}{2}(\tr_{\slashed{g}}\chibar )r L T\phi_1
	- \zeta^\alpha r \slashed{\nabla}_\alpha T\phi_1
	+ r[T, \tilde{\Box}_g] \phi_1
	- rT(F)
\end{equation*}
Note that all of the quantities on the right hand side involve at most one transverse derivative of $\phi_1$, so these are all quantities which we already know. To be even more explicit, we have
\begin{equation*}
\begin{split}
	L\left( r\Lbar T\phi_1\right)
	+ \frac{1}{2} (\tr_{\slashed{g}}\chi_{(\text{small})}) (r\Lbar T\phi_1 )
	&=
	\frac{1}{2} r\slashed{\Delta} L f_3
	- \frac{1}{4}(\tr_{\slashed{g}}\chibar )r L L f_3
	- \frac{1}{2}\zeta^\alpha r \slashed{\nabla}_\alpha L f_3
	+ \frac{1}{2} r\slashed{\Delta} (\Lbar \phi_1)
	\\
	&\phantom{=}
	- \frac{1}{4}(\tr_{\slashed{g}}\chibar )r L \Lbar \phi_1
	- \frac{1}{2}\zeta^\alpha r \slashed{\nabla}_\alpha \Lbar \phi_1
	+ r[T, \tilde{\Box}_g] \phi_1
	- rT(F)
\end{split}
\end{equation*}
This expression only involves derivatives of $f_3$ or $\Lbar \phi_1$ which are \emph{tangential} to the hypersurface $u = u_0$. Since we know these quantities everywhere on this hypersurface, we can compute these quantities. Note that $[T, \tilde{\Box}_g] \phi_1$ is some second order operator (see proposition \ref{proposition commute T}) which we can express as
\begin{equation*}
[T, \tilde{\Box}_g] \phi_1
=
- \omega \slashed{\D}_{\Lbar} \slashed{\D}_T \phi \, + \, (\text{terms with at most one transversal derivative of $\phi$})
\end{equation*}
so, schematically, we have
\begin{equation*}
\begin{split}
L\left( r\Lbar T\phi_1\right)
+ \left( \omega + \frac{1}{2}\tr_{\slashed{g}}\chi_{(\text{small})} \right) (r\Lbar T\phi_1 )
&=
- rT(F)
\, + \, (\text{good terms})
\end{split}
\end{equation*}
To find the initial conditions for this, we can use the wave equation \ref{proposition scalar wave operator} to find the required value of $\phi_1$ at the ``corner'' $\{u = u_0\} \cap \{t = t_0\}$. Note that, using the wave equation, \emph{all} the second derivatives on the spacelike part of the initial data surface can be determined from the functions $f_1$ and $f_2$. Indeed, the only terms which are not immediately given are the second \emph{transverse} derivatives, but these can found using the equation as
\begin{equation*}
(g^{-1})^{00} \partial_t^2 \phi_1
=
-2(g^{-1})^{0i} \partial_i f_2
- (g^{-1})^{ij} \partial_i \partial_j f_1
+ (g^{-1})^{ab} \Gamma^0_{ab} f_2
+ (g^{-1})^{ab} \Gamma^i_{ab} \partial_i f_1
- \omega \Lbar^0 f_2
- \omega \Lbar^i \partial_i f_1
+ F
\end{equation*}

After we have found $\Lbar T \phi_1$, we can find $\Lbar \Lbar \phi_1$ by using
\begin{equation*}
\Lbar\Lbar\phi_1 = 2\Lbar T \phi_1 - \Lbar L \phi_1
\end{equation*}
and we already know both the terms on the right hand side. Importantly, if $|rT(F)| \lesssim \epsilon (1+r)^{C_1 \epsilon}$ for large enough $C_1$, and we have the usual bootstrap bounds for the geometric quantities, then the propagation equation for $\Lbar T \phi_1$ (with suitable initial conditions) gives the bound
\begin{equation*}
|\Lbar \Lbar \phi_1| \lesssim \epsilon (1+r)^{-1+2C_1 \epsilon}
\end{equation*}

Since we are only presenting a sketch of this proof, we will not proceed any further, but it is clear that, in this way, \emph{all} transverse derivatives of $\phi_1$ on the initial data surface can be found. Moreover, with the bounds we can assume on $F$ and its derivatives, the transverse derivatives of $\phi_1$ decay at rates $\sim r^{-1+C\epsilon}$. We can also repeat this entire process if $\phi$ is not a scalar field, but rather a section of the vector bundle $\mathcal{B}$, etc.

Now, we have seen that $\phi_1$ and all of its transverse derivatives can be defined on the initial data surface. We then find a function $\phi_2$ which matches $\phi_1$ and all of its derivatives, but is defined on the entire spacetime region $\mathcal{M}(t_0,t_1)$ (see figure \ref{figure spacetime regions local}) - that is, over the intersection of the causal past of $\{t = t_1 \}\cap\{u \geq u_0\}$ with the causal future of $\{t = t_0\}$. Note that the only requirement on this function is that it matches $\phi_1$ and all of its derivatives on the surface $\{u = u_0\} \cap \{t_0 \leq t \leq t_1\}$. Such a function can be constructed using the Whitney embedding theorem \cite{Rendall1990}.

We now define a function $\rho$ by
\begin{equation*}
\rho := \begin{cases}
	\tilde{\Box}_g \phi_2 \quad &\text{ if} \quad u \leq u_0
	\\
	F \quad &\text{ if} \quad u > u_0
\end{cases}
\end{equation*}
By construction, this function is smooth. Now, we attempt to solve the equation
\begin{equation}
\label{equation appendix box chi}
\tilde{\Box}_g \chi = -\tilde{\Box}_g \phi_2 + \rho
\end{equation}
with vanishing initial data for $\chi$ on the spacelike hypersurface $t = t_0$. In the region $\tilde{M}(t_0,t_1, u_0)$ (see figure \ref{figure spacetime regions local}), by a domain-of-dependence argument this is just the solution to a linear wave equation with vanishing initial data, so by standard results $\chi = 0$ in this region. Hence, the function $(\chi + \phi_2)$ matches the original initial data on the surface $\Sigma_0$ (see figure \ref{figure spacetime regions local}). Furthermore, in the causal future of this hypersurface (i.e.\ in the region ${M}(t_0,t_1) \setminus \tilde{M}(t_0,t_1, u_0)$) we have
\begin{equation*}
\tilde{\Box}_g (\chi + \phi_2) = F
\end{equation*}
so $(\chi + \phi_2)$ is the solution we seek. We have already seen how to construct the function $\phi_2$; we construct $\chi$ by solving the wave equation \eqref{equation appendix box chi} which has initial data posed on a \emph{spacelike} hypersurface. Moreover, since this equation is linear, and the initial data has \emph{locally} bounded energy (in fact, it has bounded ``degenerate energy''), the solution can be constructed in as large a region as we like.

Now that we have seen that a solution to the \emph{linear} wave equation exists (and is unique), we can set up an iteration scheme to find the local solution to the \emph{nonlinear} problem. We seek a solution in the region $\mathcal{M}_0^{\epsilon}$, which we regard as the space
\begin{equation*}
\mathcal{M}_0^{\epsilon} = \Big( [0, \epsilon] \times [0, r_0] \times \mathbb{S}^2 \Big) \cup \Big( [0, \epsilon] \times (r_0, \infty) \times \mathbb{S}^2 \Big)
\end{equation*}
with points on the boundaries $r = r_0$ identified in the obvious way, so that the region $r \leq r_0$ is provided with coordinates $t \in [0, \epsilon]$, $r \in [0, r_0]$, while the region $r > r_0$ is provided with coordinates $\tau \in [0, \epsilon]$, $r \in (r_0, \infty)$, and both regions are provided with a map from the region in question to the unit sphere. Note that our iteration scheme will ensure that the surfaces of constant $u$ are null for \emph{every} iterate, and not just in the limit. Since we are only sketching the proof, we will focus on the region $r \geq r_0$, since this is where all the difficulties lie.

Now, we \emph{define} the null frame fields $L$, $\Lbar$, $\slashed{\nabla}$ in terms of these coordinates by
\begin{equation*}
\begin{split}
	L &:= \partial_r
	\\
	\Lbar_{(n)} &:= 2(\mu_{(n)})^{-1} \partial_\tau - \partial_r + (b_{(n)})
	\\
	\left(\slashed{\nabla}_{(n)}\right)_{\alpha} &:= \left(\slashed{\Pi}_{(n)}\right)_\alpha^{\phantom{\alpha}\beta} \slashed{\D}_\beta
\end{split}
\end{equation*}
where the projection operator $\slashed{\Pi}_{(n)}$ is defined as
\begin{equation*}
\left(\slashed{\Pi}_{(n)}\right)_\alpha^{\phantom{\alpha}\beta}
=
\delta_\alpha^\beta + \frac{1}{2}L_\alpha (\Lbar_{(n)})^\beta + \frac{1}{2}(\Lbar_{(n)})_\alpha L^\beta 
\end{equation*}
Note that we do not need to include an index $(n)$ on the vector field $L$, since it is always the same irrespective of the iteration (i.e.\ it is independent of $n$). Similarly, we \emph{define} the wave operator by
\begin{equation*}
\tilde{\Box}_{g_{(n)}} \phi
:=
-L\Lbar_{(n)} \phi
+ \slashed{\Delta}_{(n)} \phi
- \frac{1}{2} \tr_{\slashed{g}} \chi_{(n)} \Lbar_{(n)} \phi
- \frac{1}{2} \tr_{\slashed{g}} \chibar_{(n)} L \phi
- (\zeta_{(n)})^\alpha (\slashed{\nabla}_{(n)})\phi
\end{equation*}
and in general we can define various geometric quantities with an index $(n)$ by including an index $(n)$ on all of the relevant geometric quantities\footnote{Note that, for clarity of notation, we do not insert an index $n$ into the subscript $\slashed{g}$ in terms like $\tr_{\slashed{g}} \chi_{(n)}$. Nevertheless, the trace is taken here with respect to the metric $\slashed{g}_{(n)}$!}. 

We define the quantities with index $(0)$ as those that correspond to the vanishing of all of the fields $\phi_{(a)}$. In many cases, this means that the metric $g$ is the Minkowski metric $m$, and the various quantities take their Minkowski values. For example, in this case we would have $\mu_{(0)} = 1$ and $\tr_{\slashed{g}}\chi_{(0)} = \frac{2}{r}$.

We then iterate by solving linear system of wave equations
\begin{equation*}
\tilde{\Box}_{g_{(n)}} \phi_{(a, n+1)} = F_{(a, n)}
\end{equation*}
where $F_{(a,n)}$ is the inhomogenous term corresponding to the field $\phi_{(a)}$, with the fields $\phi_{(a, n)}$ entered instead of $\phi_{(a)}$. Since both the wave operator and the inhomogenous terms depend only on the fields $\phi_{(a, n-1)}$ and not on the field $\phi_{(a, n)}$, this is a set of \emph{linear} equations. We can therefore construct a local solution (and in fact, a global solution) using the method outlined above.

Note that the (linear) solutions constructed in this way will satisfy all of the bounds (both in $L^2$ and $L^\infty$) used in the main part of the proof\footnote{This is why it is important to iterate using the \emph{reduced} wave operator $\tilde{\Box}_{g_{(n)}}$ rather than the standard wave operator $\Box_{g_{(n)}}$}. Indeed, up until the point that we begin using the bootstrap bounds, all of the proof may be taken as establishing energy bounds and pointwise bounds on \emph{linear} solutions on manifolds with suitable metrics. This also means that the metric in the next iteration also obeys all of the required bounds, since the metric components and other geometric quantities (such as the null frame connection components) are related to the solutions of the wave equation in the right kinds of ways, as shown in the main body of the paper.

We need to show that this sequence of solutions converges, at least for sufficiently small $\epsilon$. The difference between the $(n+1)$-th iterate and the $n$-th iterate satisfies the wave equation
\begin{equation*}
\tilde{\Box}_{g_{(n)}} \left( \phi_{(a, n+1)} - \phi_{(a, n)} \right)
= 
\left( \tilde{\Box}_{g_{(n-1)}} - \tilde{\Box}_{g_{(n)}} \right) \phi_{(a,n)}
+ F_{(a,n)} - F_{(a,n-1)}
\end{equation*}
with vanishing initial data. We can expand this, first by expanding the wave operators, to obtain
\begin{equation*}
\begin{split}
\tilde{\Box}_{g_{(n)}} \left( \phi_{(a, n+1)} - \phi_{(a, n)} \right)
&= 
-L\left( \Lbar_{(n-1)} \phi_{(a,n)} - \Lbar_{(n)} \phi_{(a,n)} \right)
+ \left( \slashed{\Delta}_{(n-1)} - \slashed{\Delta}_{(n)} \right) \phi_{(a,n)}
\\
&\phantom{=}
- \frac{1}{2} \tr_{\slashed{g}} \chi_{(n-1)} \Lbar_{(n-1)} \phi_{(a,n)}
+ \frac{1}{2} \tr_{\slashed{g}} \chi_{(n)} \Lbar_{(n)} \phi_{(a,n)}
\\
&\phantom{=}
- \frac{1}{2} \left( \tr_{\slashed{g}} \chibar_{(n-1)} - \tr_{\slashed{g}} \chibar_{(n)} \right) L\phi_{(a,n)}
\\
&\phantom{=}
- (\zeta_{(n-1)})^\alpha (\slashed{\nabla}_{(n-1)})\phi_{(a,n)}
+ (\zeta_{(n)})^\alpha (\slashed{\nabla}_{(n)})\phi_{(a,n)}
+ F_{(a,n)} - F_{(a,n-1)}
\end{split}
\end{equation*}
and then by writing everthing in terms of the null frame associated\footnote{Note that, for a \emph{scalar} field $\phi$, we have $(\slashed{\nabla}_{(n)}) \phi = (\slashed{\nabla}_{(n-1)}) \phi$, since both of these quantities are equal to the restriction to the tangent space of the spheres of the covector $\upd \phi$, and this notion is independent of the metric (and the metric on the spheres). Of course, this is not the case for higher order quantities.} with the metric $g_{(n-1)}$, we obtain
\begin{equation*}
\begin{split}
\tilde{\Box}_{g_{(n)}} \left( \phi_{(a, n+1)} - \phi_{(a, n)} \right)
&=
-\frac{1}{2}(h_{(n-1)} - h_{(n)})_{\Lbar_{(n)} L} \left( L\Lbar_{(n-1)} \phi_{(a,n)} \right)
\\
&\phantom{=}
- \frac{1}{2}(g_{(n-1)})_{ab} (\Lbar_{(n-1)})^a (\Lbar_{(n)} - \Lbar_{(n-1)})^b (LL\phi_{(a,n)})
\\
&\phantom{=}
- (\slashed{\Pi}_{(n-1)})_a^{\phantom{a}\alpha} (\Lbar_{(n)} - \Lbar_{(n-1)})^a L\left( (\slashed{\nabla}_{(n-1)} \phi_{(a,n)}) \right)
\\
&\phantom{=}
+ \left( (\slashed{g}_{(n-1)}^{-1}) - (\slashed{g}_{(n)}^{-1}) \right)^{\alpha\beta} \left( (\slashed{\nabla}_{(n-1)})_\alpha (\slashed{\nabla}_{(n-1)})_\beta \phi_{(a,n)} \right)
\\
&\phantom{=}
- \frac{1}{2} \bigg(
	L \left( (g_{(n-1)})_{ab} \Lbar_{(n-1)}^a (\Lbar_{(n)} - \Lbar_{(n-1)})^b \right)
	\\
	&\phantom{=	- \frac{1}{2} \bigg(}
	+ \tr_{\slashed{g}} \chibar_{(n-1)} - \tr_{\slashed{g}} \chibar_{(n)}
\bigg) (L \phi_{(a,n-1)})
\\
&\phantom{=}
- \frac{1}{2} \bigg( 
	L \left( \left( h_{(n-1)} - h_{(n)} \right)_{\Lbar_{(n)} L} \right)
	+ \tr_{\slashed{g}} \chi_{(n-1)} - \tr_{\slashed{g}} \chi_{(n)}
	\\
	&\phantom{= - \frac{1}{2} \bigg( }
	+ \frac{1}{2} \left( (h_{(n-1)} - h_{(n)})_{\Lbar_{(n)} L} \right) \tr_{\slashed{g}} \chi_{(n)}
\bigg) (\Lbar_{(n-1)} \phi_{(a,n)})
\\
&\phantom{=}
- \left( \zeta_{(n-1)}^{\alpha} - \zeta_n^\alpha
	- (\slashed{g}_{(n)}^{-1})^{\beta\gamma} \slashed{\omega}^\alpha_{\phantom{\alpha}\beta\gamma}
\right) \left( (\slashed{\nabla}_{(n-1)})_{\alpha} \phi_{(a,n)} \right)
\\
&\phantom{=}
+ F_{(a,n)} - F_{(a,n-1)}
\end{split}
\end{equation*}
where here, we note that $\slashed{g}_{(n)}$ and $\slashed{g}_{(n-1)}$ define two \emph{different} metrics on the spheres. The difference between the connections associated with these two metrics can be used to define a tensor field: we define $\slashed{\omega}$ by
\begin{equation*}
(\slashed{\nabla}_{(n-1)})_\alpha X^\beta - (\slashed{\nabla}_{(n)})_\alpha X^\beta = \slashed{\omega}^\beta_{\alpha\gamma} X^\gamma
\end{equation*}
for all vector fields $X^\beta$ on the sphere.

We can continue to expand this expression in more detail. In particular, using the wave equation satisfied by $\phi_{(a,n)}$, we can replace the term $L\Lbar_{(n-1)}\phi_{(a,n)}$. Schematically, we obtain an expression of the form
\begin{equation*}
\begin{split}
&\tilde{\Box}_{g_{(n)}} \left( \phi_{(a, n+1)} - \phi_{(a, n)} \right)
\\
&=
\left( 
	r^{-1} \left(\overline{\slashed{\D}} \mathscr{Y}_{(n-1)} \phi_{(n)} \right)
	+ r^{-1} \left(\slashed{\D} \phi_{(n)} \right)
\right) \left( \left( \phi_{(n)} - \phi_{(n-1)} \right) + \left( (X_{(frame)})_{(n)} - (X_{(\text{frame})})_{(n-1)} \right) \right)
\\
&\phantom{=}
+ \left(
	\bar{\partial} (\phi_{(n)} - \phi_{(n-1)}) 
	+ \left( (\bm{\Gamma}^{(0)}_{(-1-\delta)})_{(n)} - (\bm{\Gamma}^{(0)}_{(-1-\delta)})_{(n-1)} \right)
	\right) \slashed{\D} \phi_{(n)}
\\
&\phantom{=}
+ \left(
	\bar{\partial} \left( (X_{(\text{frame})})_{(n)} - (X_{(\text{frame})})_{(n-1)} \right)
	+ \left( (\bm{\Gamma}^{(0)}_{(-1+C_{(0)}\epsilon)})_{(n)} - (\bm{\Gamma}^{(0)}_{(-1+C_{(0)}\epsilon)})_{(n-1)} \right)
\right) \overline{\slashed{\D}} \phi_{(n)}
\\
&\phantom{=}
+ F_{(a,n)} - F_{(a, n-1)}
\end{split}
\end{equation*}
where now all of the covariant derivatives are taken with respect to the metric $g_{(n-1)}$, and where $\phi_{(n)}$ stands for any of the fields $\phi_{(n, a)}$. Note that similar inequalities can be obtained after commuting, but, since we are only sketching the proof here, we will not enter into that issue.

Ignoring, for a moment, the final terms involving the inhomogeneous terms $F_{(a,n)}$ and $F_{(a,n-1)}$, note each of the extra ``error terms'' on the right hand sides has good properties: each of them involves at least one ``good derivative'' or ``good'' connection coefficient, so, morally, these error terms have the ``classical null condition''. Note also that every error term in the equation for the difference $(\phi_{(n+1)} - \phi_{(n)})$ involves the difference between the previous iterates $(\phi_{(n)} - \phi_{(n-1)})$.

Now, suppose that $\phi_{(a)} \in \Phi_{[m]}$. Then, schematically, and neglecting easier terms, we can write
\begin{equation*}
\begin{split}
F_{(a,n)} - F_{(a,n-1)}
&\sim
(\partial \phi_{[0]})_{(n)} \left( (\partial \phi_{[m]})_{(n)} - (\partial \phi_{[m]})_{(n-1)} \right)
+ (\partial \phi_{[m]})_{(n-1)} \left( (\partial \phi_{[0]})_{(n)} - (\partial \phi_{[0]})_{(n-1)} \right)
\\
&\phantom{\sim}
+ \left( (\partial \phi_{[m-1]})_{(n)} + (\partial \phi_{[m-1]})_{(n-1)} \right) \left( (\partial \phi_{[m-1]})_{(n)} - (\partial \phi_{[m-1]})_{(n-1)} \right)
\end{split}
\end{equation*}
Note that the coefficient of the term $\left( (\partial \phi_{[m]})_{(n)} - (\partial \phi_{[m]})_{(n-1)} \right)$ decays like $r^{-1}$. The coefficients of the terms of the form $\left(\partial \phi_{[m-1]})_{(n)} - (\partial \phi_{[m-1]})_{(n-1)} \right)$ do not decay at this sharp rate, however, these involve fields which are lower in the hierarchy.

Now, let us define the norm
\begin{equation*}
\begin{split}
||\phi||_{(C\epsilon)} &:=
\sup_{\tau \in [0, \epsilon]} \int_{\Sigma_\tau} \left( 
	(1+r)^{-C\epsilon} |\overline{\partial} \phi|^2
	+ \chi_{r_0} r^{-1-C\epsilon} |L (r\phi)|^2
\right) r^2 \upd r \dVol_{\mathbb{S}^2}
\\
&\phantom{:=}
+ \int_{\mathcal{M}_0^{\epsilon}} \left( 
	C\epsilon (1+r)^{-1-C\epsilon} |\partial \phi|^2
	+ \chi_{r_0} r^{-C\epsilon} |\overline{\partial} \phi|^2
\right) \dVol_g
\end{split}
\end{equation*}
and further let us define
\begin{equation*}
||\phi||_{([m], C\epsilon)} := \sup_{\phi_{(a)} \in \Phi_{[m]}} ||\phi_{(a)}||_{(C\epsilon)}
\end{equation*}
and finally, we define
\begin{equation*}
||\phi|| := \sup_{[m]} ||\phi||_{([m], C_{[m]}\epsilon)}
\end{equation*}

If we use the energy estimates from the main body of the paper, and pay attention to the error terms above (together with the fact that the differences between the iterates have vanishing initial data) we can obtain the bound
\begin{equation*}
||(\phi)_{(n+1)} - (\phi)_{(n)}||
\lesssim
\frac{1}{C_{[m]}} ||(\phi)_{(n)} - (\phi)_{(n-1)}||
\end{equation*}
If we take $C_{[m]}$ sufficiently large, then we see that our sequence of iterates converges as required.

Finally, note that, \emph{after} a solution is found on the manifold $\Big( [0, \epsilon] \times [0, r_0] \times \mathbb{S}^2 \Big) \cup \Big( [0, \epsilon] \times (r_0, \infty) \times \mathbb{S}^2 \Big)$, we can map this solution onto some subset of $\mathbb{R}^4$ as follows: first, the subset $\Big( [0, \epsilon] \times [0, r_0] \times \mathbb{S}^2 \Big)$ is mapped onto $\{x \in \mathbb{R}^4 \; | \; 0 \leq t(x) \leq \epsilon \; , \; r(x) \leq r_0\}$ in the obvious way. That is, the point $(t_1, r_1, \sigma)$ (where $\sigma \in \mathbb{S}^2$) is mapped to the point in $\mathbb{R}^4$ with $t$ coordinate $t_1$, $r$ coordinate $r_1$ and with spherical coordinate corresponding to the point $\sigma$. Note that, by pulling back the coordinate functions $t, x^1, x^2, x^3$ by this map, we obtain ``rectangular coordinates'' on the manifold $\Big( [0, \epsilon] \times [0, r_0] \times \mathbb{S}^2 \Big)$.

In the region $r \geq r_0$ things are not so simple. However, we can find functions $x^a$ on the manifold $\Big( [0, \epsilon] \times (r_0, \infty) \times \mathbb{S}^2 \Big)$ as follows: on the ``inner boundary'', where $r = r_0$, we begin with coordinates corresponding to the rectangular coordinates constructed above (recall that points on the boundaries of the two manifolds, i.e.\ at $r = r_0$, are identified). Then, we extend these functions by solving the ODE along an outgoing integral curve of $L$:
\begin{equation*}
L(x^a) = L^a
\end{equation*}
where the ``rectangular components'' $L^a$ are themselves constructed by solving the transport equations given in proposition \ref{proposition transport La}. The coordinates $x^a$ then give the required map onto a subset of $\mathbb{R}^4$. Note that, if we do this during the iteration process, then we obtain a different map for \emph{each} iterate. For example, if we use the solution $\phi_{(n)}$ to construct this map, then the push-forward of the vector field $L$ will be null \emph{with respect to the metric} $g_{(n)}$. Similarly, the push-forward of the function $u$ will satisfy the eikonal equation with respect to the metric $g_{(n)}$.

\chapter{An explicit example of shock formation at infinity}
\label{appendix explicit shock formation}

In the main body of this work we have mentioned, several times, that our estimates are consistent with a phenomenon that we have called ``shock formation at infinity''. However, it is not clear that this phenomenon actually \emph{does} take place for any set of wave equations - for example, we can actually rule out this kind of behaviour for the Einstein equations. However, in this appendix we will present a simple equation, and an explicit solution, which exibits shock formation at infinity.

In fact, we do not need to construct the solution: we can use our control over the initial data to choose data so that the shock forms \emph{instantaneously}. To be more explicit: we can choose initial data so that the inverse foliation density $\mu$ (which is not prescribed directly in the initial data, but which can be computed from the data) tends to zero as $r \rightarrow 0$. The data that we choose is smooth, compactly supported and spherically symmetric. In fact, it will vanish in the region $r \geq r_0$

Consider the scalar quasilinear wave equation
\begin{equation*}
\tilde{\Box}_g \phi = 0
\end{equation*}
where the metric is given by
\begin{equation*}
g := \frac{1}{\left( 1 + \chi_{(\frac{1}{2}r_0)}\phi+\frac{1}{4}\chi_{(\frac{1}{2}r_0)}^2\phi^2 \right)} \left( -\upd t^2 + \chi_{(\frac{1}{2}r_0)}\phi \, \upd t \upd r + (1+\chi_{(\frac{1}{2}r_0)} \phi) \upd r^2 \right) + \slashed{g}_{\mathbb{S}^2}
\end{equation*}
where the cut-off function $\chi_{(\frac{1}{2}r_0)}$ is defined in equation \eqref{equation cut off functions}, and $\slashed{g}_{\mathbb{S}^2}$ is the standard round metric on the sphere of radius $r$. Note that this metric is identical to the Minkowksi metric in the region $r \leq \frac{1}{4}r_0$.

With this metric\footnote{The chosen metric might appear to be ``artificial'', but in fact we found this metric by first choosing $L$ and $\Lbar$ to be particularly simple, which then specifies $g^{-1}$. We then invert $g^{-1}$ to find $g$.}, if we choose axisymmetric initial data then the null frame vectors must lie in the span of the coordinate vector fields $\partial_t$ and $\partial_r$. Hence, in the axisymmetric case the null frame vector fields are found to be
\begin{equation*}
\begin{split}
L &= (1+ \chi_{(\frac{1}{2}r_0)} \phi)\partial_t + \partial_r
\\
\Lbar &= \partial_t - \partial_r
\end{split}
\end{equation*}
Note that the metric satisfies the radial normalisation condition (as can also be seen directly by computing $g^{-1}$), and that the rectangular components of $L$ are given by $L^i = \frac{x^i}{r}$.

The metric component $h_{LL}$ is found to be
\begin{equation*}
h_{LL} = g_{LL} - m_{LL} = \phi(2 + \phi)
\end{equation*}
Furthermore, using proposition \ref{proposition initial data for mu}, we can compute the inverse foliation density as $r \searrow r_0$. We find that
\begin{equation*}
\mu = \frac{1}{1+\phi}
\end{equation*}

We choose the initial data for $\phi$ as follows: in the region $r \leq r_0$, on the surface $t = t_0$ we fix
\begin{equation*}
\begin{split}
\phi &= \epsilon \chi_{(\frac{1}{2}r_0)}(r) (r - r_0)
\\
\partial_t \phi &= 0
\end{split}
\end{equation*}
for some $\epsilon > 0$. Note that, with these prescriptions, in the region $r \geq \frac{1}{2}r_0$ we have
\begin{equation*}
\begin{split}
L\phi &= \epsilon \\
\Lbar \phi &= -\epsilon
\end{split}
\end{equation*}

For $r \geq r_0$, on the surface $u = t_0 - r_0$, we choose the initial data
\begin{equation*}
\phi = \epsilon(1-\chi_{(3r_0)}(r))(r-r_0)
\end{equation*}
Note that $\phi$ is continuous at $r = r_0$. In fact, $\phi$ is smooth: this is easy to see, since it is the restriction of a smooth function to the initial data surface. This smooth function can be chosen to be
\begin{equation*}
	\tilde{\phi} := \epsilon(1-\chi_{(3r_0)}(r)) \chi_{(\frac{1}{2}r_0)}(r) (r - r_0) 
\end{equation*}

We can compute the derivatives of $\phi$ in the region $r \geq r_0$, we have
\begin{equation*}
	L\phi = \epsilon(1-\chi_{(3r_0)}(r)) - \epsilon \chi'_{(3r_0)}(r)(r-r_0)
\end{equation*}
Note that, for $r \geq 3r_0$, both $\phi$ and $L\phi$ vanish.

Now, we can compute the transverse derivative $\Lbar \phi$ on the initial data, in the region $r \geq r_0$. Using the wave equation given in proposition \ref{proposition scalar wave operator} and the fact that the data is spherically symmetric, together with the form of the metric, we find that $r\Lbar \phi$ satisfies the transport equation
\begin{equation*}
L\left( r\Lbar \phi \right) = -\frac{1}{2}\tr_{\slashed{g}} \chi_{(\text{small})} \left( r\Lbar \phi \right) - \frac{1}{2}r \tr_{\slashed{g}} \chibar (L\phi)
\end{equation*}
In fact, it turns out that, with our prescribed metric, $\tr_{\slashed{g}} \chi_{(\text{small})}$ and $\tr_{\slashed{g}}\chibar_{(\text{small})}$ vanish identically. This is because the metric on the spheres is the standard round metric on the sphere of radius $r$, and the induced metric on the spheres is independent of $t$. Hence $\Lbar \phi$ satisfies
\begin{equation*}
L\left( r\Lbar \phi \right) = L\phi
\end{equation*}
In other words, $\left( r\Lbar \phi - \phi \right)$ is \emph{conserved} along the initial data surface. At $r = r_0$, we have $\Lbar \phi = -\epsilon$, we actually have
\begin{equation*}
\Lbar \phi = -\frac{\epsilon r_0}{r} + \epsilon(1-\chi_{(3r_0)}(r))\left(1-\frac{r_0}{r} \right) \quad \text{for } r \geq r_0
\end{equation*}
In particular, we see that, for $r \geq 3r_0$, we have $\Lbar \phi = -\frac{\epsilon r_0}{r}$.

We can also compute the initial value of $\mu$ (see proposition \ref{proposition initial data for mu}) as $r \searrow r_0$. Recall that, for $r \geq r_0$, $u$ solves the eikonal equation, so $g^{-1}(\upd u, \upd u) = 0$. At $r = r_0$, $g$ is the Minkowski metric, so we find that
\begin{equation*}
	\lim_{r \searrow r_0} \upd u \big|_{r} = \upd t - \upd r
\end{equation*}
Hence, using the fact that $g$ approaches the Minkowski metric as $r \rightarrow r_0$, we also find that
\begin{equation*}
	g^{-1}(\upd u, \upd r) \rightarrow -1 \text{ \, as } r \searrow r_0
\end{equation*}
and so $\mu \rightarrow 1$ as $r$ tends to $r_0$ from above.

Now, a slightly detailed computation shows that we find that, in the region $r \geq r_0$, we have
\begin{equation*}
(\Lbar h)_{LL} = 2 \Lbar \phi
	+ 2\phi (\Lbar \phi)
	- 2\phi (\Lbar \phi) \left( \frac{ 2 + \phi + \frac{1}{4}\phi^2 + \frac{1}{4}\phi^3 }{ 1 + \phi + \frac{1}{4}\phi^2 } \right)
\end{equation*}
Importantly, since $\phi = 0$ for $r \geq 3r_0$, we have
\begin{equation*}
(\Lbar h)_{LL} = -\frac{2\epsilon r_0}{r} \text{ \, for } r \geq 3r_0
\end{equation*}

If we recall the transport equation for the inverse foliation density $\mu$ given in proposition \ref{proposition transport mu}, and we substitute for $L^i$ (recalling that $L^i = \frac{x^i}{r}$) then we find that $\mu$ satisfies the transport equation\footnote{Note that $(\Lbar h)_{LL} = \Lbar(h_{LL})$, since in this particular case $\Lbar L^a = 0$.}
\begin{equation*}
L\log \mu = -\frac{1}{2} (Lh)_{L\Lbar} + \frac{1}{4} (Lh)_{LL} + \frac{1}{4} (\Lbar h)_{LL}
\end{equation*}

The first two terms vanish for $r \leq 3r_0$. Hence we can integrate the equation for $\mu$, finding
\begin{equation*}
\log \mu = C(r) \epsilon + \log \left( \left( \frac{r}{r_0} \right)^{-\frac{1}{2}\epsilon r_0} \right)
\end{equation*}
where $C(r)$ is given by
\begin{equation*}
	C(R) := \frac{1}{\epsilon} \int_{r = r_0}^{R} \left(
		- \frac{1}{2} (Lh)_{L\Lbar}
		+ \frac{1}{4} (Lh)_{LL}
		+ \frac{1}{2} \phi (\Lbar \phi)
		- \frac{1}{2} \phi (\Lbar \phi)  \left( \frac{ 2 + \phi + \frac{1}{4}\phi^2 + \frac{1}{4}\phi^3 }{ 1 + \phi + \frac{1}{4}\phi^2 } \right)
		\right) \upd r
\end{equation*}
Note that all of these terms are compactly supported in $r \leq 3r_0$, and moreover they are at least $\mathbb{O}(\epsilon)$. Hence $C(r)$ is constant for $r \geq 3r_0$. In particular, $C(3r_0)$ is some bounded numerical constant. Moreover, if we consider $C(3r_0)$ as a function of $\epsilon$, then for all $0 < \epsilon < 1$, $C(3r_0)$ is uniformly bounded by some other numerical constant $\tilde{C}$.

In particular, for \emph{all} $\epsilon$, we have
\begin{equation*}
\mu \leq e^{ \tilde{C} \epsilon} \left( \frac{r_0}{r} \right)^{\frac{1}{2}\epsilon r_0}
\end{equation*}

In particular, for \emph{any} positive value of $\epsilon$, we have $\mu \rightarrow 0$ as $r \rightarrow \infty$, signalling shock formation at infinity. Note also that, if we choose $\epsilon < 0$ then we form an ``anti-shock'', where $\mu \rightarrow \infty$ as $r \rightarrow \infty$.

In summary, we have exhibited a family of smooth, compactly supported initial data for a wave equation that satisfies the weak null condition, depending on a small parameter $\epsilon$. The data is such that, as $\epsilon \rightarrow 0$, the initial data becomes trivial. In particular, for all sufficiently small values of $\epsilon$, our main theorem \ref{theorem main theorem} guarantees a global solution. However, for \emph{every} value of $\epsilon$, the corresponding solution exhibits \emph{immediate shock formation at infinity}.

If we want to find an example of a set of wave equations where the derivatives of the \emph{fields} actually exhibit different asymptotic behaviour from the linear case, then we can study the exact same scenario, but introduce a second field $\phi_2$ satisfying
\begin{equation*}
\tilde{\Box}_g \phi_2 = (T \phi)^2
\end{equation*}
and where $\phi_2$ is given \emph{trivial} initial data over the entire initial data surface (including the spacelike portion $r \leq r_0$). Then, using the wave equation, we find that on the initial data surface, in the region $r \geq r_0$, $\Lbar \phi_2$ satisfies the transport equation
\begin{equation*}
L\left( r\Lbar \phi_2 \right) = r(T \phi)^2 = \frac{1}{4} r (\Lbar \phi)^2
\end{equation*}
substituting for $\Lbar \phi$ from above, we have that, for $r \geq r_0$, 
\begin{equation*}
L\left( r\Lbar \phi_2 \right) = \frac{\epsilon^2 r_0^2}{4r} 
\end{equation*}
and so $\Lbar \phi_2$ is given by
\begin{equation*}
\Lbar \phi_2 = \frac{\epsilon^2 r_0^2}{4 r} \log \left( \frac{r}{r_0} \right)
\end{equation*}
so that this does \emph{not} have the $\frac{1}{r}$ decay of the linear equation. Similarly, if we consider the field $\phi_3$ satisfying
\begin{equation*}
\tilde{\Box}_g \phi_2 = (T \phi)(T\phi_3)
\end{equation*}
and again give $\phi_3$ trivial initial data, then we find that $\Lbar\phi_3$ is given, on the initial data surface, in the region $r \geq r_0$, by
\begin{equation*}
\Lbar \phi_3 = \frac{1}{r} \left( \frac{r_0}{r} \right)^{\frac{1}{4}\epsilon}
\end{equation*}
Both the fields $\phi_2$ and $\phi_3$ can be said to exhibit immediate ``blowup at infinity''.

\sloppy
\printbibliography

\end{document}